\documentclass[12pt,a4paper]{amsart}
\usepackage{amssymb}
\usepackage[all,cmtip]{xy}
\usepackage{hyperref}

\calclayout

\newcommand{\+}{\protect\nobreakdash-}
\renewcommand{\:}{\colon}

\newcommand{\rarrow}{\longrightarrow}
\newcommand{\larrow}{\longleftarrow}
\newcommand{\ot}{\otimes}
\newcommand{\ocn}{\odot}
\newcommand{\st}{\star}
\newcommand{\bec}{\natural}

\DeclareFontFamily{U}{mathb}{\hyphenchar\font45}
\DeclareFontShape{U}{mathb}{m}{n}{
      <5> <6> <7> <8> <9> <10> gen * mathb
      <10.95> mathb10 <12> <14.4> <17.28> <20.74> <24.88> mathb12
      }{}
\DeclareSymbolFont{mathb}{U}{mathb}{m}{n}
\DeclareFontSubstitution{U}{mathb}{m}{n}
\DeclareMathSymbol{\blackdiamond}{0}{mathb}{"0C}

\newcommand{\bu}{{\text{\smaller\smaller$\scriptstyle\bullet$}}}
\newcommand{\cu}{{\text{\smaller$\scriptstyle\blackdiamond$}}}
\newcommand{\subcu}{{\text{\smaller$\scriptscriptstyle\blackdiamond$}}}

\newcommand{\lrarrow}{\mskip.5\thinmuskip\relbar\joinrel\relbar\joinrel
 \rightarrow\mskip.5\thinmuskip\relax}
\newcommand{\llarrow}{\mskip.5\thinmuskip\leftarrow\joinrel\relbar
 \joinrel\relbar\mskip.5\thinmuskip\relax}

\DeclareMathOperator{\Spec}{Spec}
\DeclareMathOperator{\Hom}{Hom}
\DeclareMathOperator{\Ext}{Ext}
\DeclareMathOperator{\Tor}{Tor}
\DeclareMathOperator{\cone}{cone}
\DeclareMathOperator{\coker}{coker}
\DeclareMathOperator{\im}{im}
\DeclareMathOperator{\Sym}{Sym}

\DeclareMathOperator{\Cohom}{\mathfrak{Cohom}}
\DeclareMathOperator{\fHom}{\mathfrak{Hom}}
\DeclareMathOperator{\cHom}{\mathcal{H}\mskip-.3\thinmuskip\textit{om}}

\newcommand{\Modl}{{\operatorname{\mathsf{--Mod}}}}
\newcommand{\Modr}{{\operatorname{\mathsf{Mod--}}}}
\newcommand{\biMod}{{\operatorname{\mathsf{--Mod--}}}}
\newcommand{\QMod}{{\operatorname{\mathsf{--QMod}}}}
\newcommand{\Ctrh}{{\operatorname{\mathsf{--Ctrh}}}}
\newcommand{\Lcth}{{\operatorname{\mathsf{--Lcth}}}}
\newcommand{\Qcoh}{{\operatorname{\mathsf{--Qcoh}}}}
\newcommand{\Qcohr}{{\operatorname{\mathsf{Qcoh--}}}}
\newcommand{\biQcoh}{{\operatorname{\mathsf{--Qcoh--}}}}
\newcommand{\QQcoh}{{\operatorname{\mathsf{--QQcoh}}}}
\newcommand{\Sh}{{\operatorname{\mathsf{--Sh}}}}
\newcommand{\Cosh}{{\operatorname{\mathsf{--Cosh}}}}

\newcommand{\bModl}{{\operatorname{\mathbf{--Mod}}}}
\newcommand{\bModr}{{\operatorname{\mathbf{Mod--}}}}
\newcommand{\bCtrh}{{\operatorname{\mathbf{--Ctrh}}}}
\newcommand{\bLcth}{{\operatorname{\mathbf{--Lcth}}}}
\newcommand{\bQcoh}{{\operatorname{\mathbf{--Qcoh}}}}
\newcommand{\bQcohr}{{\operatorname{\mathbf{Qcoh--}}}}

\newcommand{\Ab}{\mathsf{Ab}}
\newcommand{\Ac}{\mathsf{Ac}}
\newcommand{\Hot}{\mathsf{Hot}}
\newcommand{\Com}{\mathsf{Com}}
\newcommand{\bCom}{\mathbf{Com}}

\newcommand{\rop}{\mathrm{op}}
\newcommand{\sop}{\mathsf{op}}
\newcommand{\id}{\mathrm{id}}
\newcommand{\gr}{\mathrm{gr}}
\newcommand{\sk}{\mathrm{sk}}
\newcommand{\sy}{\mathrm{sy}}
\newcommand{\cry}{\mathsf{cr}}
\newcommand{\str}{\mathsf{st}}
\newcommand{\qu}{\mathsf{qu}}

\newcommand{\bb}{\mathsf{b}}
\newcommand{\abs}{\mathsf{abs}}
\newcommand{\co}{\mathsf{co}}
\newcommand{\ctr}{\mathsf{ctr}}
\newcommand{\bco}{\mathsf{bco}}
\newcommand{\bctr}{\mathsf{bctr}}
\newcommand{\si}{\mathsf{si}}
\newcommand{\bsi}{\mathsf{bsi}}

\newcommand{\inj}{\mathsf{inj}}
\newcommand{\proj}{\mathsf{proj}}
\newcommand{\binj}{\mathbf{inj}}
\newcommand{\bproj}{\mathbf{proj}}
\renewcommand{\cot}{\mathsf{cot}}
\newcommand{\cta}{\mathsf{cta}}
\newcommand{\lct}{\mathsf{lct}}
\newcommand{\lin}{\mathsf{lin}}
\newcommand{\lvfl}{\mathsf{lvfl}}
\newcommand{\lfl}{\mathsf{lfl}}
\newcommand{\rfl}{\mathsf{rfl}}
\newcommand{\al}{\mathsf{al}}
\newcommand{\qc}{\mathsf{qc}}
\newcommand{\ct}{\mathsf{ct}}
\newcommand{\thk}{\mathsf{th}}
\newcommand{\bth}{\mathbf{th}}
\newcommand{\ff}{\mathsf{ff}}
\newcommand{\vflp}{\mathsf{vflp}}
\newcommand{\flp}{\mathsf{flp}}
\newcommand{\fl}{\mathsf{fl}}
\newcommand{\alfp}{\mathsf{alfp}}

\newcommand{\dfl}{{\operatorname{\mathsf{-fl}}}}
\newcommand{\dlct}{{\operatorname{\mathsf{-lct}}}}
\newcommand{\dlin}{{\operatorname{\mathsf{-lin}}}}
\newcommand{\dcta}{{\operatorname{\mathsf{-cta}}}}
\newcommand{\dcot}{{\operatorname{\mathsf{-cot}}}}
\newcommand{\dinj}{{\operatorname{\mathsf{-inj}}}}
\newcommand{\dvflp}{{\operatorname{\mathsf{-vflp}}}}
\newcommand{\dflp}{{\operatorname{\mathsf{-flp}}}}
\newcommand{\drflp}{{\operatorname{\mathsf{-rflp}}}}

\newcommand{\red}{{\operatorname{\mathsf{-red}}}}
\newcommand{\dqc}{{\operatorname{\mathsf{-qc}}}}

\newcommand{\sA}{\mathsf A}
\newcommand{\sC}{\mathsf C}
\newcommand{\sD}{\mathsf D}
\newcommand{\sE}{\mathsf E}
\newcommand{\sF}{\mathsf F}
\newcommand{\sH}{\mathsf H}
\newcommand{\sK}{\mathsf K}
\newcommand{\sL}{\mathsf L}
\newcommand{\sZ}{\mathsf Z}

\newcommand{\cA}{\mathcal A}
\newcommand{\cB}{\mathcal B}
\newcommand{\cC}{\mathcal C}
\newcommand{\cD}{\mathcal D}
\newcommand{\cE}{\mathcal E}
\newcommand{\cO}{\mathcal O}
\newcommand{\cT}{\mathcal T}
\newcommand{\cU}{\mathcal U}
\newcommand{\cV}{\mathcal V}

\newcommand{\C}{\mathcal C}
\newcommand{\E}{\mathcal E}
\newcommand{\F}{\mathcal F}
\newcommand{\G}{\mathcal G}
\newcommand{\cH}{\mathcal H}
\newcommand{\J}{\mathcal J}
\newcommand{\K}{\mathcal K}
\newcommand{\cL}{\mathcal L}
\newcommand{\M}{\mathcal M}
\newcommand{\N}{\mathcal N}
\newcommand{\cR}{\mathcal R}

\newcommand{\gA}{\mathfrak A}
\newcommand{\gB}{\mathfrak B}
\newcommand{\gC}{\mathfrak C}
\newcommand{\gF}{\mathfrak F}
\newcommand{\gI}{\mathfrak I}
\newcommand{\gJ}{\mathfrak J}
\newcommand{\gK}{\mathfrak K}
\newcommand{\gM}{\mathfrak M}
\newcommand{\gN}{\mathfrak N}
\renewcommand{\P}{\mathfrak P}
\newcommand{\Q}{\mathfrak Q}
\newcommand{\R}{\mathfrak R}

\newcommand{\g}{\mathfrak g}
\newcommand{\h}{\mathfrak h}
\newcommand{\p}{\mathfrak p}

\newcommand{\vect}{\mathfrak{vect}}

\newcommand{\bA}{\mathbf A}
\newcommand{\bB}{\mathbf B}
\newcommand{\bC}{\mathbf C}
\newcommand{\bE}{\mathbf E}
\newcommand{\bF}{\mathbf F}
\newcommand{\bW}{\mathbf W}
\newcommand{\bT}{\mathbf T}

\newcommand{\boZ}{\mathbb Z}
\newcommand{\boQ}{\mathbb Q}
\newcommand{\boR}{\mathbb R}
\newcommand{\boL}{\mathbb L}

\newcommand{\Section}[1]{\bigskip\section{#1}\medskip}
\theoremstyle{plain}
\newtheorem{thm}{Theorem}[section]
\newtheorem{prop}[thm]{Proposition}
\newtheorem{lem}[thm]{Lemma}
\newtheorem{cor}[thm]{Corollary}
\theoremstyle{definition}
\newtheorem{rem}[thm]{Remark}
\newtheorem{rems}[thm]{Remarks}
\newtheorem{quest}[thm]{Question}
\newtheorem{ex}[thm]{Example}
\newtheorem{exs}[thm]{Examples}

\begin{document}

\title{$\cD$-$\Omega$ duality on the contra side}

\author{Leonid Positselski}

\address{Institute of Mathematics, Czech Academy of Sciences \\
\v Zitn\'a~25, 115~67 Prague~1 \\ Czech Republic} 

\email{positselski@math.cas.cz}

\begin{abstract}
 Given a smooth morphism of schemes $X\rarrow T$, denote by
$\cD_{X/T}^\cry$ the sheaf of rings of fiberwise crystalline
differential operators on $X$ relative to $T$ and by $\Omega^\bu_{X/T}$
the de~Rham sheaf of DG\+algebras of relative differential forms on
$X$ over~$T$.
 Assume that the scheme $X$ is quasi-compact and semi-separated.
 We construct a commutative square diagram of triangulated equivalences
between four triangulated categories: the derived category of
quasi-coherent sheaves of $\cD_{X/T}^\cry$\+modules, the reduced
coderived category of quasi-coherent DG\+modules over
$\Omega_{X/T}^\bu$, the derived category of contraherent cosheaves of
$\cD_{X/T}^\cry$\+modules, and the reduced contraderived category of
contraherent DG\+modules over~$\Omega_{X/T}^\bu$.
 The equivalence involving the contraderived category was previously
known for affine varieties only; we use contraherent cosheaves in
order to obtain a nonaffine generalization of the ``contra side''
of the story.
 The exposition is written in the generality of finite locally free
twisted Lie algebroids $(\g,\widetilde\g)$ over quasi-compact
semi-separated schemes $X$, the quasi-coherent twisted universal
enveloping quasi-algebras of $(\g,\widetilde\g)$, and
the Chevalley--Eilenberg quasi-coherent CDG\+quasi-algebras
of~$(\g,\widetilde\g)$.
 The equivalence between the derived categories of quasi-coherent
and contraherent $\cA$\+modules, called the ``na\"\i ve co-contra 
correspondence'', is proved quite generally for any quasi-coherent
quasi-algebra $\cA$ over~$X$.
\end{abstract}

\maketitle

\tableofcontents

\section*{Introduction}
\medskip

 The concept of \emph{Koszul duality} goes back to
Chevalley--Eilenberg~\cite{ChE}, Koszul~\cite{Kos},
Eilenberg--Mac~Lane~\cite{EML}, Adams~\cite{Ad1,Ad2},
Quillen~\cite{Quil}, and Priddy~\cite{Prid}.
 \emph{Derived Koszul duality} means triangulated equivalences of
module categories over Koszul dual algebras; this topic goes back to
Bernstein--Gelfand--Gelfand~\cite{BGG} and Beilinson~\cite{Beil}.
 For recent surveys penned by the present author,
see~\cite[Prologue]{Prel} and~\cite{Pksurv}.

 The duality between rings of differential operators and DG\+algebras
of differential forms is one of the primary examples of Koszul duality
in algebraic and differential geometry.
 The $\cD$\+$\Omega$ duality is an instance of \emph{nonhomogeneous
Koszul duality}, in that the rings of differential operators are not
graded, but only filtered.

 The relations defining the rings of differential operators involve
functions and vector fields.
 These relations are \emph{nonhomogeneous quadratic}, in that they have
principal quadratic parts (with respect to the vector fields) and
additional linear and scalar terms.
 The differential on the de~Rham DG\+algebra corresponds to
the lower-degree summands in the relations defining the ring of
differential operators.

 Attempts to work out triangulated equivalences of the derived
$\cD$\+$\Omega$ duality were made by Kapranov~\cite{Kap} and
Beilinson--Drinfeld~\cite[Section~7.2]{BD2}.
 The current state of the art can be found in the present author's
memoir~\cite[Appendix~B]{Pkoszul} and book~\cite{Prel}.

 As established by Hinich~\cite{Hin2}, Lef\`evre-Hasegawa~\cite{Lef}, 
Keller~\cite{Kel}, and the present author~\cite{PV},
\cite[Section~6]{Pkoszul}, Koszul duality is naturally viewed as
connecting \emph{algebras} with \emph{coalgebras}.
 Following an observation going back to Eilenberg and Moore~\cite{EM}
and developed in~\cite{Psemi,Prev}, there are \emph{two} kinds of
natural module categories over a coalgebra: in addition to the more
familiar \emph{comodules}, there are also \emph{contramodules}.
 Accordingly, the Koszul duality as worked out in~\cite{Pkoszul}
happens on the comodule and contramodule side, connected with each
other by the \emph{comodule-contramodule correspondence}.
 We refer to~\cite{Prev} and~\cite[Sections~8\+-9]{Pksurv} for
introductory survey expositions.

 The main innovation of the modern theory of derived nonhomogeneous
Koszul duality is the \emph{derived categories of the second kind}.
 The basic idea and the terminology of the \emph{first} vs.\
\emph{second kind} in this context goes back to Husemoller, Moore,
and Stasheff~\cite{HMS}, but the actual definitions of weak equivalences
of the second kind are due to Hinich~\cite{Hin2},
Lef\`evre-Hasegawa~\cite{Lef}, Keller~\cite{Kel}, and the present
author~\cite{Psemi,Prev}.
 The terminology is a bit misleading, as ``derived categories of
the second kind'' is a generic name for several different constructions,
the most important of them known as the \emph{absolute derived},
\emph{coderived}, and \emph{contraderived categories}.
 We refer to~\cite[Remark~9.2]{PS4} and~\cite[Section~7]{Pkoszul} for
historical and philosophical discussions.

 Following the philosophy elaborated upon in~\cite{Psemi,Pkoszul}, one
is generally supposed to consider the conventional \emph{derived}
categories of modules, the \emph{coderived} categories of comodules,
and the \emph{contraderived} categories of contramodules (though there
are many exceptions to this rule).
 So the ``co'' side of derived Koszul duality features the coderived
category of comodules, while on the ``contra'' side one considers
the contraderived category of contramodules.

 In addition to the distinction between the \emph{homogeneous} Koszul
duality (as in~\cite{BGG,Beil,BGS}) and the \emph{nonhomogeneous} one
(as in~\cite{Hin2,Lef,Kel,PV,Pkoszul,Pksurv}), there is also
a relevant distinction between \emph{absolute} and \emph{relative}
Koszul duality.
 Here ``absolute'' means ``over a field'' (or in the most complicated
cases, over a semisimple ring, as in~\cite{BGS}, or over a semisimple
coalgebra, as in~\cite{HL}).
 ``Relative'' means over a more or less arbitrary (certainly
nonsemisimple) ring~\cite[Section~0.4]{Psemi}, \cite{Prel} or
coalgebra~\cite[Chapter~11]{Psemi}.
 The $\cD$\+$\Omega$ duality is a species of relative nonhomogeneous
Koszul duality, in that it happens over the base ring of functions.

 The derived nonhomogeneous Koszul duality theory developed in the main
body of the memoir~\cite{Pkoszul} and surveyed in~\cite{Pksurv} is
absolute over a field.
 The relative case over a base ring is more complicated in many
respects, and first of all, in that it needs to be properly understood
what constitutes the coalgebra side.
 The answer is that the ring of differential operators $\cD$ stands
on the algebra side, while the de~Rham DG\+ring $\Omega^\bu$ plays
the role of the coalgebra.
 Accordingly, one needs to consider the coderived category of
DG\+modules over $\Omega^\bu$ on the ``comodule'' side, and
the contraderived category of the ``contramodule'' side.

 Generally speaking, in the context of relative nonhomogeneous Koszul
duality over a base ring, one has to consider comodules and
contramodules over graded rings~\cite[Chapters~5\+7]{Prel}.
 This distinction does \emph{not} manifest itself in the case of
the $\cD$\+$\Omega$ duality, as the de~Rham DG\+algebra $\Omega^\bu$
is bounded in its cohomological grading.
 But the distinction (and the connections) between the \emph{coderived}
and the \emph{contraderived} category of DG\+modules over $\Omega^\bu$
is important.

 The relative nonhomogeneous Koszul duality theory spelled out in
the book~\cite{Prel} is affine; there are \emph{no} nonaffine
schemes in~\cite{Prel}.
 The earlier exposition on $\cD$\+$\Omega$ duality
in~\cite[Appendix~B]{Pkoszul} covers smooth nonaffine algebraic
varieties, but only on the ``co'' side.
 All the assertions about the contraderived categories
in~\cite[Appendix~B]{Pkoszul} apply to affine varieties only.
 The reason is that the definition of the contraderived category
does \emph{not} make sense for quasi-coherent sheaves.

 In fact, there are two such definitions.
 The approach of the present author, as per~\cite{Psemi,Pkoszul},
uses infinite products, and makes sense when the infinite product
functors are exact.
 But infinite products of quasi-coherent sheaves are not exact.
 The alternative approach, named after Becker with the reference to
his paper~\cite{Bec} and going back to J\o rgensen~\cite{Jor},
Krause~\cite{Kra}, and Neeman~\cite{Neem2}, uses projective objects
or DG\+modules whose underlying graded modules are projective.
 But nonzero projective quasi-coherent sheaves usually do not exist.
 
 \emph{Contraherent cosheaves} were invented in~\cite{Pcosh} in order
to rectify the situation.
 For any scheme $X$, the exact category of contraherent cosheaves
$X\Ctrh$ has exact functors of infinite products.
 Under mild assumptions (e.~g., for quasi-compact semi-separated
schemes~$X$), this exact category also has enough projective objects.
 So both the Positselski and Becker constructions of contraderived
categories are applicable to contraherent cosheaves.

 In the context of $\cD$\+$\Omega$ duality, the constructions of
the derived categories of the second kind are applied to DG\+modules
over the de~Rham DG\+algebra $\Omega^\bu$, while for modules over
the sheaf of rings of differential operators $\cD$ one is mostly
interested in the conventional derived categories (or in their
modifications called the ``semiderived categories'').
 In the first part of this paper, we work with what we call
\emph{quasi-coherent quasi-algebras} $\cA$ over quasi-compact
semi-separated schemes~$X$.
 These are natural generalizations of the sheaves of rings of
differential operators~$\cD$.

 In this setting, we construct a triangulated equivalence
\begin{equation} \label{introd-naive-co-contra}
 \sD(\cA\Qcoh)\simeq\sD(\cA\Ctrh)
\end{equation}
between the conventional derived category of the abelian category
$\cA\Qcoh$ of quasi-coherent sheaves of $\cA$\+modules and
the conventional derived category of the exact category $\cA\Ctrh$ of
contraherent cosheaves of $\cA$\+modules on~$X$.
 The equivalence holds for the bounded and one-sided bounded as well as
for the unbounded conventional derived categories.
 This is an instance of what we call the ``na\"\i ve co-contra
correspondence'' phenomenon and a far-reaching generalization
of the similar result for quasi-coherent and contraherent
$\cO_X$\+modules~\cite[Theorem~4.8.1]{Pcosh}.
 See Theorems~\ref{quasi-algebra-bounded-derived-naive-co-contra}
and~\ref{quasi-algebra-unbounded-derived-naive-co-contra}
(and formula~\eqref{bounded-derived-naive-co-contra-diagram-II})
below in the main body of this paper.

 The word ``na\"\i ve'' means precisely that this version of
the co-contra correspondence connects the conventional derived
categories rather than the coderived and contraderived categories.
 We refer to the introduction to the paper~\cite{Pmgm} for a discussion
of the philosophy of co-contra correspondence, and to~\cite[Section~0.11
of the Introduction]{Pcosh} for further comments on the na\"\i ve
co-contra correspondence and an explanation of how it arises.

 In the second part of the paper, we consider (what we call)
quasi-coherent twisted Lie algebroids $(\g,\widetilde\g)$ over
quasi-compact semi-separated schemes~$X$.
 The Chevalley--Eilenberg quasi-coherent curved DG\+quasi-algebras
$\cB^\cu=\cC^\cu_X(\g,\widetilde\g)$ of twisted Lie algebroids
$(\g,\widetilde\g)$ with a finite locally free underlying quasi-coherent
sheaf~$\g$ are natural generalizations of the sheaves of de~Rham
DG\+rings of differential forms~$\Omega^\bu$.

 We construct a triangulated equivalence
\begin{equation} \label{introd-genuine-co-contra}
 \sD^\co_{X\red}(\cB^\cu\bQcoh)\simeq\sD^\ctr_{X\red}(\cB^\cu\bCtrh)
\end{equation}
between the reduced coderived category of quasi-coherent sheaves of
curved DG\+mod\-ules and the reduced contraderived category of
contraherent cosheaves of curved DG\+modules over~$\cB^\cu$ on~$X$.
 This is an instance of ``non-na\"\i ve'' or genuine co-contra
correspondence, in that it connects some versions of the coderived
and the contraderived category.
 See Corollary~\ref{cdg-module-co-contra-correspondence-cor} below.

 Notice that a comparable equivalence between the coderived category
of quasi-coherent sheaves and the contraderived category of
contraherent cosheaves~\cite[Theorem~6.7.1 or~6.8.2]{Pcosh} assumes
the scheme $X$ to be Noetherian and depends on a dualizing complex
on~$X$.
 Our result in this paper is different in that it applies to arbitrary
quasi-compact semi-separated schemes $X$ and does not need a dualizing
complex.
 Instead, we consider the \emph{reduced} coderived and contraderived
categories, which are a kind of mixtures of the conventional derived
category along the structure sheaf $\cO_X$ of the scheme $X$ and
the coderived or contraderived category in the direction of $\Omega^\bu$
or $\cC^\cu_X(\g,\widetilde\g)$ relative to~$\cO_X$.

 For a semi-separated regular Noetherian scheme $X$ of finite Krull
dimension, the reduced co/contraderived categories are no different
from the usual (unreduced) ones.
 However, the reduced co/contraderived categories play a crucial role
in the $\cD$\+$\Omega$ duality and derived co-contra correspondence
theory for DG\+modules over $\Omega^\bu$ in the case of a scheme $X$
that is \emph{not} necessarily (Noetherian and) regular of finite
Krull dimension.
 For example, if one is interested in relative/fiberwise differential
operators and differential forms for a smooth morphism of schemes
$X\rarrow T$, then the singular (or even non-Noetherian) schemes $X$
and $T$ are a natural setting.

 Let us now describe the content of the paper section by section.
 The paper consists of three parts.
 Sections~\ref{module-prelim-secn}\+-\ref{naive-co-contra-secn}
form the first part, leading to the construction of the derived
na\"\i ve co-contra correspondence over an arbitrary quasi-coherent
quasi-algebra~$\cA$.
 Sections~\ref{contraherent-cdg-modules-secn}\+-%
\ref{cdg-co-contra-secn} constitute the second part, where we work
out the co-contra correspondence over the Chevalley--Eilenberg
quasi-coherent CDG\+quasi-algebra $\cC^\bu_X(\g,\widetilde\g)$ of
a finite locally free twisted Lie algebroid~$(\g,\widetilde\g)$.
 Section~\ref{becker-semicontraderived-secn} stands apart a bit in
that it returns to the setting of an arbitrary quasi-coherent
quasi-algebra (rather than a quasi-coherent CDG\+quasi-algebra or
a twisted locally free Lie algebroid) typical to
Sections~\ref{module-prelim-secn}\+-\ref{naive-co-contra-secn}.
 The $\cD$\+$\Omega$ duality and quadrality theorems are presented
in the third part of the paper, which consists of
Sections~\ref{D-Omega-duality-secn}\+-\ref{well-behaved-schemes-secn}.
 Once again, the final Section~\ref{well-behaved-schemes-secn} stands
a bit apart in that (with the exception of some technical results)
it makes more restrictive assumptions on the scheme~$X$.

 The module-theoretic Section~\ref{module-prelim-secn}, besides known
results on cotorsion modules~\cite{En2}, contraadjusted and very flat
modules~\cite{Pcosh,Pcta,PSl1} etc., includes an in-depth discussion
of \emph{quasi-modules}~\cite{BB,Ptd} and modules over
\emph{quasi-algebras}~\cite[Section~2.3]{Pedg} over commutative rings,
as well as what we call \emph{flaprojective modules} (a relative
mixture of the flat and projective ones).
 Quasi-modules, known classically as ``differential
bimodules''~\cite[Section~1.1]{BB}, are basically bimodules over
a commutative rings $R$ in which the left and right actions of $R$
almost, but not quite, agree with each other; technically speaking,
these are modules over $R\otimes_\boZ R$ supported set-theoretically
in the diagonal of the Cartesian square $\Spec R\times_{\Spec\boZ}
\Spec R$ \,\cite[Corollary~2.3]{Ptd}.
 Quasi-modules localize well over schemes, so one can speak of
\emph{quasi-coherent quasi-modules}.
 Quasi-algebras are basically quasi-modules with an associative ring
structure.

 Section~\ref{contraherent-cosheaves-of-O-modules-secn} is a reminder
of the background material on contraherent cosheaves from
the long manuscript~\cite{Pcosh}.
 Section~\ref{contraherent-cosheaves-of-A-modules-secn} works out
an extension of this background from the case of contraherent cosheaves
of modules over the structure sheaf $\cO_X$ to that of contraherent
cosheaves of $\cA$\+modules for a quasi-coherent quasi-algebra $\cA$
over~$X$.

 The important technical results concerning the \emph{antilocality}
property of contraadjusted and cotorsion quasi-coherent sheaves
and the concept of antilocal contraherent cosheaves, going back
to~\cite[Sections~4.1 and~4.3]{Pcosh} and elaborated upon in
the paper~\cite{Pal} and the appendix~\cite[Appendix~B]{Pcosh},
are extended to the realm of $\cA$\+modules in
Section~\ref{antilocal-classes-secn}.
 That is where the concept of flaprojective modules, introduced
in Section~\ref{prelim-flaprojective-subsecn} and discussed in
the paper~\cite{PBas}, plays a role.
 The derived na\"\i ve co-contra correspondence between quasi-coherent
and locally contraherent $\cA$\+modules, for a quasi-coherent
quasi-algebra $\cA$ over a quasi-compact semi-separated scheme $X$,
is established in Section~\ref{naive-co-contra-secn}.

 In the second part of the paper, in
Section~\ref{contraherent-cdg-modules-secn}, we introduce
quasi-coherent curved differential graded (CDG) quasi-algebras over
schemes and construct the DG\+categories of quasi-coherent and
locally contraherent CDG\+modules over such CDG\+quasi-algebras.
 In Section~\ref{contraderived-and-semiderived-secn}, we define and
discuss the coderived categories of quasi-coherent CDG\+modules,
the contraderived categories of locally contraherent CDG\+modules,
and the semiderived categories of (complexes of) quasi-coherent and
locally contraherent modules over a quasi-algebra.

 The advantage of the constructions of the Becker co/contraderived
categories is that they tend to provably ``give the correct answer''
or ``produce the correct category'' in a greater generality than
the respective Positselski constructions.
 However, the Becker coderived categories, and particularly the Becker
contraderived categories, are often more difficult to work with
than the Positselski ones.
 In particular, establishing the locality of the Becker
contraacyclicity property for complexes of contraherent cosheaves
is a nontrivial task~\cite[Section~5.4]{Pcosh}.
 In Section~\ref{becker-semicontraderived-secn}, we prove the locality
of the Becker contraacyclicity for complexes of contraherent
modules over a quasi-coherent quasi-algebra~$\cA$.
 This result is relevant in the context of the definition and
basic properties of the Becker semicontraderived categories of
contraherent $\cA$\+modules, which is the main aim of
Section~\ref{becker-semicontraderived-secn}.

 Section~\ref{twisted-lie-algebroids-secn} offers a discussion of
twisted Lie algebroids $(\g,\widetilde\g)$ over schemes $X$, their
twisted universal enveloping quasi-algebras $\cA_X(\g,\widetilde\g)$
(including twisted differential operators), and the Chevalley--Eilenberg
quasi-coherent CDG\+quasi-algebras $\cB^\cu=\cC^\cu_X(\g,\widetilde\g)$.
 These form a natural setting for the $\cD$\+$\Omega$ duality theory
over schemes.
 Here we assume that the underlying quasi-coherent sheaf~$\g$ of
the quasi-coherent twisted Lie algebroid $(\g,\widetilde\g)$ is
a finite locally free sheaf on~$X$.
 In the case of the (nontwisted) Lie algebroid of vector fields
$\g=\vect_{X/T}$ for a smooth morphism of schemes $X\rarrow T$,
the quasi-coherent universal enveloping quasi-algebra $\cA=\cA_X(\g)=
\cD_{X/T}^\cry$ is the quasi-algebra of fiberwise crystalline
differential operators on $X$ relative to $T$, while
the Chevalley--Eilenberg quasi-coherent CDG\+quasi-algebra
$\cB^\cu=\cC^\bu_X(\g)=\Omega^\bu_{X/T}$ is the de~Rham sheaf of
DG\+algebras of relative differential forms on $X$ over~$T$.

 In the longish Section~\ref{thick-graded-modules-secn},
we study a special (adjusted) class of graded modules over the exterior
algebra of a finite locally free sheaf over a scheme, called
the \emph{thick graded modules}.
 We discuss thick quasi-coherent graded modules, thick locally
contraherent graded modules, and thick CDG\+modules.
 These are important as a technical instrument for the $\Omega$ side of
$\cD$\+$\Omega$ duality.
 In Section~\ref{cdg-reduced-contraderived-secn}, which is even longer,
we define and discuss the \emph{reduced coderived categories} of
quasi-coherent CDG\+modules and \emph{reduced contraderived categories}
of locally contraherent CDG\+modules over the Chevalley--Eilenberg
quasi-coherent CDG\+quasi-algebra.

 The derived co-contra correspondence between quasi-coherent sheaves
and contraherent cosheaves of CDG\+modules over
the Chevalley--Eilenberg quasi-coherent CDG\+quasi-algebra
$\cB^\cu=\cC^\cu_X(\g,\widetilde\g)$ is established in
Section~\ref{cdg-co-contra-secn}.
 As mentioned above, this triangulated equivalence involves
the reduced coderived category of quasi-coherent CDG\+modules
$\sD^\co_{X\red}(\cB^\cu\bQcoh)$ and the reduced contraderived category
of contraherent CDG\+modules $\sD^\ctr_{X\red}(\cB^\cu\bCtrh)$.

 In Section~\ref{D-Omega-duality-secn}, we establish the triangulated
equivalences of $\cD$\+$\Omega$ duality on the co and contra sides.
 Let $\cA=\cA_X(\g,\widetilde\g)$ be the quasi-coherent twisted
universal enveloping quasi-algebra of a twisted Lie algebroid
$(\g,\widetilde\g)$ with a finite locally free underlying
quasi-coherent sheaf~$\g$, and let
$\cB^\cu=\cC^\cu_X(\g,\widetilde\g)$ be the Chevalley--Eilenberg
quasi-coherent CDG\+quasi-algebra.
 On the co side, we construct the triangulated equivalences
describing the (conventional and semico)derived categories
of quasi-coherent \emph{right} $\cA$\+modules,
\begin{align}
 \sD^\si(\Qcohr\cA) &\simeq \sD^\co(\bQcohr\cB^\cu),
 \label{introd-semiderived-koszul-duality-right-co-side} \\
 \sD(\Qcohr\cA) &\simeq \sD^\co_{X\red}(\bQcohr\cB^\cu).
 \label{introd-reduced-koszul-duality-right-co-side}
\end{align}
 On the contra side, we construct the triangulated equivalences
describing the (conventional and semicontra)derived categories
of contraherent \emph{left} $\cA$\+modules,
\begin{align}
 \sD^\si(\cA\Ctrh) &\simeq \sD^\ctr(\cB^\cu\bCtrh),
 \label{introd-semiderived-koszul-duality-contra-side} \\
 \sD(\cA\Ctrh) &\simeq \sD^\ctr_{X\red}(\cB^\cu\bCtrh).
 \label{introd-reduced-koszul-duality-contra-side}
\end{align}
 The triangulated
equivalences~\eqref{introd-semiderived-koszul-duality-right-co-side}
and~\eqref{introd-reduced-koszul-duality-right-co-side} are the results
of Theorems~\ref{semiderived-koszul-duality-right-co-side}
and~\ref{reduced-koszul-duality-right-co-side}.
 The triangulated
equivalences~\eqref{introd-semiderived-koszul-duality-contra-side}
and~\eqref{introd-reduced-koszul-duality-contra-side} are the results
of Theorems~\ref{semiderived-koszul-duality-contra-side}
and~\ref{reduced-koszul-duality-contra-side}.

 The triangulated
equivalences~(\ref{introd-semiderived-koszul-duality-right-co-side}\+-%
\ref{introd-reduced-koszul-duality-right-co-side}) are more general
than the co side $\cD$\+$\Omega$ duality result
of~\cite[Theorem~B.2(a)]{Pkoszul}
in that the scheme $X$ is not assumed to be regular
in~(\ref{introd-semiderived-koszul-duality-right-co-side}\+-%
\ref{introd-reduced-koszul-duality-right-co-side}).
 The triangulated
equivalences~(\ref{introd-semiderived-koszul-duality-contra-side}\+-%
\ref{introd-reduced-koszul-duality-contra-side}) are more general
than the contra side $D$\+$\Omega$ duality result
of~\cite[Theorem~B.2(b)]{Pkoszul} in that the scheme $X$ is not
assumed to be either affine or regular
in~(\ref{introd-semiderived-koszul-duality-contra-side}\+-%
\ref{introd-reduced-koszul-duality-contra-side}).

 In the next-to-final Section~\ref{quadrality-diagram-secn}, we
construct a commutative pentagonal diagram of triangulated equivalences
\begin{equation} \label{introd-pentagonal-diagram}
\begin{gathered}
 \xymatrix{
  \sD(\Qcohr\cA^\circ)
  \ar@{=}[rd]^-{\eqref{introd-reduced-koszul-duality-right-co-side}}
  \ar@{=}[d] \\
  \sD(\cA\Qcoh) \ar@{=}[r] \ar@{=}[d]_{\eqref{introd-naive-co-contra}} &
  \sD^\co_{X\red}(\cB^\cu\bQcoh)
  \ar@{=}[d]^{\eqref{introd-genuine-co-contra}} \\
  \sD(\cA\Ctrh)
  \ar@{=}[r]_-{\eqref{introd-reduced-koszul-duality-contra-side}}
  & \sD^\ctr_{X\red}(\cB^\cu\bCtrh)
 }
\end{gathered}
\end{equation}
 Here $\cA^\circ=\cA_X(\g,\widetilde\g^\circ)$ is the quasi-coherent
twisted enveloping quasi-algebra of the opposite quasi-coherent twisted
Lie algebroid $(\g,\widetilde\g^\circ)$ to the quasi-coherent twisted
Lie algebroid~$(\g,\widetilde\g)$.
 One has $\cB^\circ{}^\cu=\cC^\cu_X(\g,\widetilde\g^\circ)\simeq
\cC^\cu_X(\g,\widetilde\g)^\rop=\cB^\rop{}^\cu$, so there is a natural
equivalence (in fact, isomorphism) of DG\+categories
$\bQcohr\cB^\circ{}^\cu\simeq\cB^\cu\bQcoh$.
 The upper diagonal line in~\eqref{introd-pentagonal-diagram} is
the equivalence~\eqref{introd-reduced-koszul-duality-right-co-side}
for $(\g,\widetilde\g^\circ)$, while the lower horizontal line
in~\eqref{introd-pentagonal-diagram} is
the equivalence~\eqref{introd-reduced-koszul-duality-contra-side}
for~$(\g,\widetilde\g)$.
 The upper leftmost vertical triangulated equivalence
$\sD(\Qcohr\cA^\circ)\simeq\sD(\cA\Qcoh)$ is induced by the classical
``conversion'' equivalence of abelian categories $\Qcohr\cA^\circ
\simeq\cA\Qcoh$.

 In the case of a nontwisted quasi-coherent Lie algebroid~$\g$
(e.~g., the usual nontwisted crystalline differential operators
$\cA=\cD_{X/T}^\cry$), one has $(\g,\widetilde\g^\circ)\simeq
(\g,\widetilde\g)$, hence $\cB^\circ{}^\cu=\cB^\cu$ and
$\cA^\circ=\cA$.
 The lower commutative square in~\eqref{introd-pentagonal-diagram}
is called the ``$\cD$\+$\Omega$ quadrality''.
 This is the main result of this paper.

 The upper commutative triangle in~\eqref{introd-pentagonal-diagram}
is the result of
Corollary~\ref{left-right-co-side-reduced-koszul-duality}.
 The lower commutative square in~\eqref{introd-pentagonal-diagram}
is the result of Theorem~\ref{main-quadrality-theorem} (see also
Theorem~\ref{main-lct-quadrality-theorem} for locally cotorsion
versions).
 The whole commutative diagram of triangulated
equivalences~\eqref{introd-pentagonal-diagram} is stated as
Corollary~\ref{main-hexagonality-corollary} (see also
Corollary~\ref{main-lct-hexagonality-corollary}).

 In the final Section~\ref{well-behaved-schemes-secn}, we discuss
the simplifications and improvements of the results and arguments
from the previous sections that can be achieved by making more
restrictive assumptions on the scheme~$X$.
 For the most part of the section, this means assuming that
the scheme is semi-separated and Noetherian of finite Krull dimension
(though some technical results are still proved for quasi-compact
semi-separated schemes quite generally, as in the rest of the paper).
 At the end of the section, we specialize even further to \emph{regular}
Noetherian schemes of finite Krull dimension, and explain that one can
drop the adjective ``$X$\+reduced'' in the formulations of the paper's
results for such schemes.

 Before we finish this introduction, let us mention that the most
technical aspects of the exposition in this paper are only worked out
for the Positselski co/contraderived categories.
 We believe that the similar results should hold for the Becker
co/contraderived categories, but, as we already mentioned, these are
harder to work with (mainly because proving that a functor preserves
the Becker co/contraacyclicity is often much more difficult than for
the Positselski co/contraacyclicity).
 Nevertheless, we prove the semiderived Koszul duality theorems on
the co and contra sides
(Theorems~\ref{semiderived-koszul-duality-right-co-side},
\ref{semiderived-koszul-duality-contra-side},
and~\ref{semiderived-koszul-duality-becker-contra-side})
for the Becker co/contraderived categories alongside with
the Positselski ones.

 At the end of the day, in the main results of this paper as per
the exposition in Section~\ref{quadrality-diagram-secn}, there is
no difference between the Positselski and Becker versions of
the exotic derived categories, as one can see from the formulations
of the main results.
 See Remark~\ref{becker=positselski-at-the-end-remark} for a discussion.

\subsection*{Acknowledgement}
 I~am grateful to Boris Tsygan for stimulating interest.
 I~also wish to thank Jan \v St\!'ov\'\i\v cek, Michal Hrbek, and
Sergio Pavon for helpful conversations.
 The author is supported by the GA\v CR project 23-05148S
and the Institute of Mathematics, Czech Academy of Sciences
(research plan RVO:~67985840).

\Section{Module-Theoretic Preliminaries}
\label{module-prelim-secn}

 The first half of this section is mostly an extraction
from~\cite[Chapter~1]{Pcosh} and~\cite[Section~4]{Pphil}.
 The second half of the section, starting from
Section~\ref{prelim-quasi-modules-subsecn}, expands and improves
upon the first half of~\cite[Section~2.3]{Pedg}.
 One of the most important new results is
Lemma~\ref{very-flat-quasi-module-lemma}.

\subsection{Homological formula of derived Hom-tensor adjunction}
 The following ``homological formula'' is useful for proving many
lemmas in Section~\ref{module-prelim-secn}.

\begin{lem} \label{Ext-homological-formula}
 Let $R$ and $S$ be two associative rings.
 Let $E$ be a left $R$\+module, $F$ be an $R$\+$S$\+bimodule
and $G$ be a left $S$\+module.
 Let $n\ge0$ be an integer such that\/ $\Ext^i_R(F,E)=0=\Tor^S_i(F,G)$
for all\/ $0<i\le n$.
 Then there is a natural isomorphism of abelian groups
$$
 \Ext^n_S(G,\Hom_R(F,E))\simeq\Ext^n_R(F\ot_S G,\>E).
$$
\end{lem}

\begin{proof}
 This is a variation on the theme of~\cite[formula~(4) in
Section~XVI.4]{CE} suggested in~\cite[formula~(1.1)
in Section~1.2]{Pcosh}.
 The formula is easily provable by a spectral sequence argument
similar to the one in~\cite{CE}.
 Alternatively, here is a more elementary proof.
 Pick a projective resolution
$$
 \dotsb\lrarrow Q_2\lrarrow Q_1\lrarrow Q_0\lrarrow G\lrarrow0
$$
of the left $S$\+module $G$ and an injective coresolution
$$
 0\lrarrow E\lrarrow J^0\lrarrow J^1\lrarrow J^2\lrarrow\dotsb
$$
of the left $R$\+module~$E$.

 By assumption, the complex $F\ot_SQ_\bu$ is exact at all
the homological degrees $0<i\le n$; so its terms sitting in
the homological degrees from~$0$ to $n+1$ form an initial
fragment of a resolution of the left $R$\+module $F\ot_SG$.
 Extending this initial fragment to an actual resolution $C_\bu$,
we obtain an exact complex of left $R$\+modules $C_\bu\rarrow F\ot_SG
\rarrow0$ together with a commutative diagram of morphisms of
left $R$\+modules
$$
 \xymatrixcolsep{1.2em}\xymatrix{
  & & F\ot_RQ_{n+1} \ar[r] \ar@{=}[d] & F\ot_RQ_n \ar[r] \ar@{=}[d]
  & \dotsb \ar[r] & F\ot_RQ_0 \ar[r]
  \ar@{=}[d] & F\ot_RG \ar[r] \ar@{=}[d] & 0 \\
  \dotsb \ar[r] & C_{n+2} \ar[r] & C_{n+1} \ar[r] & C_n \ar[r]
  & \dotsb \ar[r] & C_0 \ar[r] & F\ot_RG \ar[r] & 0
 }
$$

 Dual-analogously, the complex $\Hom_R(F,J^\bu)$ is exact at all
the cohomological degrees $0<i\le n$; so its terms sitting in
the cohomological degrees from~$0$ to $n+1$ form an initial
fragment of a coresolution of the left $S$\+module $\Hom_R(F,E)$.
  Extending this initial fragment to an actual coresolution $K^\bu$,
we obtain an exact complex of left $S$\+modules $0\rarrow
\Hom_R(F,E)\rarrow K^\bu$ together with a commutative diagram
of morphisms of left $S$\+modules
$$
 \xymatrixcolsep{1.2em}\xymatrix{
 0 \ar[r] & \Hom_R(F,E) \ar[r] \ar@{=}[d] & \Hom_R(F,J^0)
 \ar[r] \ar@{=}[d] & \dotsb \ar[r] & \Hom_R(F,J^{n+1}) \ar@{=}[d] \\
 0 \ar[r] & \Hom_R(F,E) \ar[r] & K^0 \ar[r] & \dotsb \ar[r] & K^{n+1}
 \ar[r] & K^{n+2} \ar[r] & \dotsb
 }
$$

 Now we have natural isomorphisms of abelian groups
\begin{multline*}
 \Ext^n_S(G,\Hom_R(F,E))=H^n\Hom_S(Q_\bu,\Hom_R(F,E))\\
 \simeq H^n\Hom_S(Q_\bu,K^\bu)=H^n\Hom_S(Q_\bu,\Hom_R(F,J^\bu)) \\
 \simeq H^n\Hom_R(F\ot_R Q_\bu,\>J^\bu)=H^n\Hom_R(C_\bu,J^\bu) \\
 \simeq H^n\Hom_R(F\ot_R G,\>J^\bu)=\Ext^n_R(F\ot_RG,\>E).
\end{multline*}
 Here the total complexes of the bicomplexes of $\Hom$ are presumed in
the notation.
 The first and last equalities hold by the definition of $\Ext^n$, while
the two middle equalities follow from the commutative diagrams above.
 The inner isomorphism is the usual Hom-tensor adjunction, while
the two outer isomorphisms hold due to the quasi-isomorphisms
$\Hom_R(F,E)\rarrow K^\bu$ and $C_\bu\rarrow F\ot_SG$.
\end{proof}

\subsection{Flat and cotorsion modules}
\label{prelim-cotorsion-subsecn}
 Let $R$ be an associative ring.
 We will denote by $R\Modl$ the category of left $R$\+modules and by
$\Modr R$ the category of right $R$\+modules.

 A left $R$\+module $C$ is said to be \emph{cotorsion} (\emph{in
the sense of Enochs}~\cite{En2}) if $\Ext^1_R(F,C)=0$ for all flat
left $R$\+modules~$F$.
 Clearly, the class of cotorsion $R$\+modules is closed under
extensions, direct summands, and direct products.

\begin{lem} \label{flat-cotorsion-pair-hereditary}
\textup{(a)} One has\/ $\Ext^n_R(F,C)=0$ for all flat $R$\+modules $F$,
all cotorsion $R$\+modules $C$, and all $n\ge1$. \par
\textup{(b)} The class of cotorsion $R$\+modules is closed under
cokernels of monomorphisms in $R\Modl$.
\end{lem}

\begin{proof}
 Part~(a) follows from the fact that, in any projective resolution
of the $R$\+module $F$, the $R$\+module of cycles are flat.
 In other words, the point is that the class of flat $R$\+modules
is closed under kernels of epimorphisms.
 Part~(b) follows easily from part~(a) for $n=2$.
\end{proof}

\begin{lem} \label{flat-cotorsion-is-a-cotorsion-pair}
 A left $R$\+module $F$ is flat if and only if\/ $\Ext_R^1(F,C)=0$
for all cotorsion left $R$\+modules~$C$.
\end{lem}

\begin{proof}
 The assertion can be deduced from
Theorem~\ref{flat-cotorsion-pair-complete}(a)
below; see, e.~g., \cite[Lemma~B.1.2]{Pcosh}.
 However, in fact, Lemma~\ref{flat-cotorsion-is-a-cotorsion-pair} is
a much easier (and historically much earlier) result than the much more
substantial Theorem~\ref{flat-cotorsion-pair-complete}.
 See~\cite[Lemma~3.4.1]{Xu}.
 The point is that all left $R$\+modules of the form
$C=\Hom_\boZ(N,\boQ/\boZ)$, where $N$ ranges over the right
$R$\+modules, are cotorsion (in fact, they are pure-injective,
which is a stronger property).
 Then one computes that $\Ext^1_R(F,\Hom_\boZ(N,\boQ/\boZ))\simeq
\Hom_\boZ(\Tor_1^R(N,F),\boQ/\boZ)$ for all left $R$\+modules~$F$.
\end{proof}

 The following theorem says that there are enough cotorsion modules,
in a strong sense.

\begin{thm} \label{flat-cotorsion-pair-complete}
 For any $R$\+module $M$, there exist \par
\textup{(a)} a short exact sequence of $R$\+modules\/
$0\rarrow C'\rarrow F\rarrow M\rarrow0$ with a flat $R$\+module $F$
and a cotorsion $R$\+module~$C'$; \par
\textup{(b)} a short exact sequence of $R$\+modules\/
$0\rarrow M\rarrow C\rarrow F'\rarrow0$ with a cotorsion $R$\+module $C$
and a flat $R$\+module~$F'$.
\end{thm}

\begin{proof}
 This one of the formulations of what became known as
the \emph{Flat Cover Conjecture}, posed by Enochs in~\cite{En}
and proved by Bican--El~Bashir--Enochs in~\cite{BBE}.
 In fact, there are two proofs in the paper~\cite{BBE}; the one
more closely related to the formulation above is based on
the preceding work of Eklof and Trlifaj~\cite{ET}.

 The short exact sequence in part~(a) is called a \emph{special
flat precover sequence}.
 The short exact sequence in part~(b) is called a \emph{special
cotorsion preenvelope sequence}.
 A further discussion can be found in~\cite[Theorem~8.1(a)]{GT},
\cite[Theorem~1.3.1]{Pcosh}, and~\cite[Theorem~4.4]{Pphil}.
 For explicit examples of special flat precovers and special cotorsion
preenvelopes of abelian groups, see~\cite[Section~12]{Pcta}.
\end{proof}

 The results of Lemmas~\ref{flat-cotorsion-pair-hereditary}\+-%
\ref{flat-cotorsion-is-a-cotorsion-pair} and
Theorem~\ref{flat-cotorsion-pair-complete} can be summarized by saying
that the pair of classes (flat modules, cotorsion modules) is
a \emph{hereditary complete cotorsion pair} in $R\Modl$.
 The notion of a cotorsion pair (or a \emph{cotorsion theory}) goes
back to Salce~\cite{Sal}.
 See~\cite[Chapters~5\+-6]{GT}, \cite[Section~4]{Pphil},
or~\cite[Appendix~B]{Pcosh} for a discussion.

\begin{lem} \label{hom-injective-cotorsion}
 Let $R$ and $S$ be two associative rings, $F$ be an $R$\+$S$\+bimodule,
and $E$ be a left $R$\+module.
 In this context: \par
\textup{(a)} if $F$ is a flat left $R$\+module and $E$ is a cotorsion
left $R$\+module, then\/ $\Hom_R(F,E)$ is a cotorsion left
$S$\+module; \par
\textup{(b)} if $F$ is a flat right $S$\+module and $E$ is an injective
left $R$\+module, then\/ $\Hom_R(F,E)$ is an injective left $S$\+module;
\par
\textup{(c)} if $E$ is an injective left $R$\+module, then\/
$\Hom_R(F,E)$ is a cotorsion left $S$\+module.
\end{lem}

\begin{proof}
 All the three assertions can be obtained from the homological
formula of Lemma~\ref{Ext-homological-formula}; this is
the approach in~\cite[Lemma~1.3.3]{Pcosh}.
 In all the three cases, the following more elementary argument also
works.
 Let us explain part~(a).
 It suffices to show that the functor
$G\longmapsto\Hom_S(G,\Hom_R(F,E))$ takes short exact sequences of
flat left $S$\+modules $G$ to short exact sequences of abelian groups
(since all projective modules are flat and the syzygy modules of
flat modules are flat).
 Indeed, we have $\Hom_S(G,\Hom_R(F,E))\simeq\Hom_R(F\ot_SG,\>E)$,
the functor $F\ot_S{-}$ is exact on the category of flat left
$S$\+modules $G$, the tensor product $F\ot_SG$ of an $R$\+flat
$R$\+$S$\+bimodule $F$ and a flat left $S$\+module $G$ is a flat
left $R$\+module, and the functor $\Hom_R({-},E)$ into a cotorsion
left $R$\+module $E$ is exact on the category flat left $R$\+modules.
\end{proof}

\begin{lem} \label{restriction-coextension-injective-cotorsion}
 Let $R\rarrow S$ be a homomorphism of associative rings.
 In this context: \par
\textup{(a)} if $D$ is a cotorsion $S$\+module, then $D$ is also
cotorsion as an $R$\+module; \par
\textup{(b)} if $K$ is an injective left $S$\+module and $S$ is
flat as a right $R$\+module, then $K$ is also injective as a left
$R$\+module; \par
\textup{(c)} if $C$ is a cotorsion left $R$\+module and $S$ is
flat as a left $R$\+module, then\/ $\Hom_R(S,C)$ is a cotorsion
left $S$\+module; \par
\textup{(d)} if $J$ is a an injective $R$\+module, then\/
$\Hom_R(S,J)$ is an injective $S$\+module.
\end{lem}

\begin{proof}
 This is~\cite[Lemmas~1.3.4 and~1.3.5]{Pcosh}.
 Both an argument based on Lemma~\ref{Ext-homological-formula} and
an argument similar to the one spelled out in the proof of
Lemma~\ref{hom-injective-cotorsion} above are applicable in all
the cases.
\end{proof}

\begin{prop} \label{relative-bass-theorem-cotorsion-prop}
 Let $R\rarrow A$ be a homomorphism of associative rings such that
$A$ is finitely generated and projective as a right $R$\+module.
 Then a left $A$\+module $D$ is cotorsion if and only if $D$ is
cotorsion as a left $R$\+module.
\end{prop}

\begin{proof}
 The easy ``only if'' assertion is provided by
Lemma~\ref{restriction-coextension-injective-cotorsion}(a) (and
does not depend on the assumptions about $A$ as right $R$\+module).
 The nontrivial ``if'' assertion is~\cite[Theorem~2.8]{PBas}.
\end{proof}

\subsection{Very flat and contraadjusted modules}
\label{prelim-very-flat-subsecn}
 Let $R$ be a commutative ring.
 Given an element $s\in R$, we denote by $R[s^{-1}]$ the ring $R$
with the element~$s$ inverted, i.~e., the localization of $R$ at
the multiplicative subset $\{1,s,s^2,s^3,\dotsc\}$.
 For an $R$\+module $M$, we put $M[s^{-1}]=R[s^{-1}]\ot_RM$.

 An $R$\+module $C$ is said to be
\emph{contraadjusted}~\cite[Section~1.1]{Pcosh},
\cite[Section~2]{Pcta} if $\Ext^1_R(R[s^{-1}],C)=0$ for all
elements $s\in R$.
 (The terminology ``contraadjusted'' means ``adjusted to
contraherent cosheaves''.)
 An $R$\+module $F$ is said to be
\emph{very flat}~\cite[Section~1.1]{Pcosh}, \cite{PSl1} if
$\Ext^1_R(F,C)=0$ for all contraadjusted $R$\+modules~$C$.

 It is clear that the class of contraadjusted $R$\+modules is closed
under extensions, direct summands, and infinite products; while
the class of very flat $R$\+modules is closed under extensions,
direct summands, and infinite direct sums.
 The following lemma lists slightly more nontrivial properties.

\begin{lem} \label{very-flat-cotorsion-pair-hereditary}
\textup{(a)} Any quotient $R$\+module of a contraadjusted
$R$\+module is contraadjusted. \par
\textup{(b)} The projective dimension of any very flat $R$\+module
does not exceed\/~$1$. \par
\textup{(c)} The class of very flat $R$\+modules is closed under
kernels of epimorphisms in $R\Modl$.
\end{lem}

\begin{proof}
 This is~\cite[beginning of Section~1.1]{Pcosh}.
 Part~(a) follows from the fact that the projective dimension of
the $R$\+module $R[s^{-1}]$ never exceeds~$1$ \,\cite[proof of
Lemma~2.1]{Pcta}.
 Both parts~(b) and~(c) follow from part~(a).
\end{proof}

 The following theorem says that there are enough very flat and
contraadjusted modules.

\begin{thm} \label{very-flat-cotorsion-pair-complete}
 For any $R$\+module $M$, there exist \par
\textup{(a)} a short exact sequence of $R$\+modules\/
$0\rarrow C'\rarrow F\rarrow M\rarrow0$ with a very flat $R$\+module $F$
and a contraadjusted $R$\+module~$C'$; \par
\textup{(b)} a short exact sequence of $R$\+modules\/
$0\rarrow M\rarrow C\rarrow F'\rarrow0$ with a contraadjusted
$R$\+module $C$ and a very flat $R$\+module~$F'$.
\end{thm}

\begin{proof}
 This is~\cite[Theorem~1.1.1]{Pcosh} or~\cite[Theorem~4.5]{Pphil},
or a special case of~\cite[Theorems~2 and~10]{ET}.
 In the terminology of cotorsion pairs~\cite{Sal,ET,GT}, the pair
of classes (very flat modules, contraadjusted modules) is, by
the definition, the cotorsion pair \emph{generated by} the set of
$R$\+modules $R[s^{-1}]$, \,$s\in R$, in $R\Modl$.
 The theorem says that this cotorsion pair is \emph{complete}
(while Lemma~\ref{very-flat-cotorsion-pair-hereditary} implies
that it is \emph{hereditary}).
 For any associative ring $R$, any cotorsion pair generated by
a set of modules in complete in $R\Modl$ \,\cite{ET},
\cite[Theorem~6.11]{GT}.
\end{proof}

 The proof of Theorem~\ref{very-flat-cotorsion-pair-complete}
provides more information than the theorem says.
 Specifically, the class of very flat $R$\+modules can be described
as the class of all direct summands of \emph{transfinitely iterated
extensions} (\emph{in the sense of the inductive limit}) of
the $R$\+modules $R[s^{-1}]$, \,$s\in R$.
 See~\cite[Corollary~1.1.4]{Pcosh}, \cite[Theorem~4.5(b)]{Pphil},
or~\cite[Corollary~6.14]{GT}.

\begin{exs} \label{very-flat-ring-examples}
 (1)~Let $U$ be an affine scheme and $V$ be an affine open
subscheme in~$U$.
 Then the ring of functions $\cO(V)$ is a very flat module over
the ring of functions $\cO(V)$.
 This is~\cite[Lemma~1.2.4]{Pcosh}.

 (2)~Much more generally, let $R\rarrow S$ be a morphism of
commutative rings such that $S$ is a finitely presented commutative
$R$\+algebra and a flat $R$\+module.
 Then $S$ is a very flat $R$\+module.
 This is the assertion of the \emph{Very Flat Conjecture}, which was
proved in~\cite[Main Theorem~1.1]{PSl1}.

 See also the discussion in~\cite[Remark~5.1]{Pphil}.
\end{exs}

\begin{lem} \label{vfl-cta-tensor-hom-lemma}
\textup{(a)} For any two very flat $R$\+modules $F$ and $G$,
the tensor product $F\ot_RG$ is a very flat $R$\+module. \par
\textup{(b)} For any very flat $R$\+module $F$ and contraadjusted
$R$\+module $C$, the $R$\+module\/ $\Hom_R(F,C)$ is contraadjusted.
\end{lem}

\begin{proof}
 This is~\cite[Lemma~1.2.1]{Pcosh}.
\end{proof}

\begin{lem} \label{restriction-co-extension-vfl-cta}
 Let $R\rarrow S$ be a homomorphism of commutative rings.
 In this context: \par
\textup{(a)} if $D$ is a contraadjusted $S$\+module, then $D$ is
also contraadjusted as an $R$\+module; \par
\textup{(b)} for any very flat $R$\+module $F$, the $S$\+module
$S\ot_RF$ is very flat; \par
\textup{(c)} if the $R$\+module $S[s^{-1}]$ is very flat for
every element $s\in S$ and $G$ is a very flat $S$\+module, then $G$
is also very flat as an $R$\+module; \par
\textup{(d)} if the $R$\+module $S[s^{-1}]$ is very flat for
every element $s\in S$ and $C$ is a contraadjusted $R$\+module,
then\/ $\Hom_R(S,C)$ is a contraadjusted $S$\+module.
\end{lem}

\begin{proof}
 This is~\cite[Lemma~1.2.2(a\+-b) and Lemma~1.2.3(a\+-b)]{Pcosh}.
\end{proof}

 Notice that the $R$\+module $S$ itself being very flat is \emph{not}
enough for the validity of
Lemma~\ref{restriction-co-extension-vfl-cta}(c\+-d), and the stronger
assumption that $S[s^{-1}]$ is a very flat $R$\+module for all $s\in S$
is needed.
 See~\cite[Lemma~9.3(a) and Example~9.7]{PSl1}
or~\cite[Lemma~6.1]{Pphil}.
 When the latter assumption is satisfied, one says that \emph{$S$~is
a very flat commutative $R$\+algebra} or \emph{the commutative ring
homomorphism $R\rarrow S$ is very flat}.

\subsection{Quasi-modules} \label{prelim-quasi-modules-subsecn}
 The terminology \emph{quasi-modules} in this section agrees
with the one in~\cite[Section~2.1]{Ptd} and does \emph{not} agree
with the terminology in~\cite[Section~2.3]{Pedg}.
 The more narrow class of quasi-modules considered
in~\cite[Section~2.3]{Pedg} is called the class of \emph{strong
quasi-modules} in~\cite[Section~2.2]{Ptd}, and we follow the terminology
of~\cite{Ptd} in this paper.
 A more common terminology for what we,
following~\cite[Section~2.2]{Ptd}, call a strong quasi-module (and what
is called simply a ``quasi-module'' in~\cite[Section~2.3]{Pedg}) is
a \emph{differential bimodule}~\cite[Section~1.1]{BB}.
 The point is that, out of the three definitions of a quasi-module (with
various qualifiers) discussed in~\cite[Section~2]{Ptd}, the least
restrictive one still words well for our purposes.

 Let $R$ be a commutative ring and $B$ be an $R$\+$R$\+bimodule.
 Given an element $r\in R$, let us define the natural increasing
filtration $F^{(r)}$ on $B$ by the rules
\begin{itemize}
\item $F^{(r)}_nB=0$ for all integers $n<0$;
\item for each integer $n\ge0$, \ $F^{(r)}_nB$ is the set of all
elements $b\in B$ such that $rb-br\in F^{(r)}_{n-1}B$.
\end{itemize}
 One can easily see that $F^{(r)}B$ is an $R$\+$R$\+subbimodule in $B$
and $F^{(r)}_{n-1}B\subset F^{(r)}_n$ for all $n\in\boZ$.
 We say that $B$ is a \emph{quasi-module}~\cite[Section~2.1]{Ptd}
over $R$ if, for every element $r\in R$, the filtration $F^{(r)}$
on $B$ is exhaustive, i.~e., $B=\bigcup_{n\ge0}F_n^{(r)}B$.

 For a geometric interpretation of quasi-modules as bimodules supported
in the diagonal $\Spec R\subset\Spec R\times_{\Spec\boZ}\Spec R$,
see~\cite[Corollary~2.3]{Ptd}.

 Clearly, the full subcategory of quasi-modules $R\QMod\subset
R\biMod R$ is closed under subobjects, quotients, and infinite direct
sums in the abelian category of $R$\+$R$\+bimodules $R\biMod R$.
 According to the discussion in~\cite[Sections~1.1 and~2.1]{Ptd}
(see specifically~\cite[Lemma~1.1 and Corollary~2.2]{Ptd}),
the full subcategory $R\QMod$ is also closed under extensions in
$R\biMod R$.
 In other words, one can say that $R\QMod$ is a \emph{hereditary
torsion class} (in the terminology of~\cite[Sections~VI.2\+-3]{Ste}),
or equivalently, a \emph{localizing subcategory} in $R\biMod R$.
 Hence $R\QMod$ is an abelian category (in fact, a Grothendieck
category) with an exact inclusion functor $R\QMod\rarrow R\biMod R$
preserving infinite direct sums.

 The natural ordinal-indexed increasing filtration $F$ on
an $R$\+$R$\+bimodule $B$ is defined by the rules
\begin{itemize}
\item $F_0B\subset B$ is the $R$\+$R$\+subbimodule consisting
of all elements $b\in B$ such that $rb-br=0$ in $B$ for all $r\in R$;
\item for every ordinal~$\alpha>0$, the $R$\+$R$\+subbimodule
$F_\alpha B\subset B$ consists of all elements $b\in B$ such that
$rb-br\in \bigcup_{\beta<\alpha}F_\beta B$ for all $r\in R$.
\end{itemize}
 Obviously, one has $F_\beta B\subset F_\alpha B$ for all ordinals
$\beta\le\alpha$.

 We say that an $R$\+$R$\+bimodule $B$ is a \emph{quite quasi-module
over~$R$} \,\cite[Section~2.3]{Ptd} if $B=\bigcup_\beta F_\beta B$,
where the union is taken over all the ordinals~$\beta$.
 Equivalently, $B$ is a quite quasi-module over $R$ if and only if
there exists an ordinal~$\alpha$ such that $B=\bigcup_{\beta<\alpha}
F_\beta B$.
 All quite quasi-modules are quasi-modules over $R$, but the converse
is \emph{not} true in general~\cite[Examples~1.10, 2.8, and~3.7]{Ptd}.

 We say that an $R$\+$R$\+bimodule $B$ is a \emph{strong quasi-module
over~$R$} \,\cite[Section~2.2]{Ptd} if $B=\bigcup_{n<\omega} F_nB$
(where the union is taken over the nonnegative integers~$n$).
So, by the definition, any strong quasi-module is a quite quasi-module.
 The converse is \emph{not} true in general, for the basic reason that
the class of strong quasi-modules \emph{need not} be closed under
extensions in $R\biMod R$ \,\cite[Examples~1.6, 2.5, and~3.3]{Ptd}.

 In particular, $F^{(r)}_0B\subset B$ is the subbimodule of all
elements in $B$ on which the left and right actions of
the element $r\in R$ agree.
 Therefore, any $R$\+module $M$, viewed as an $R$\+$R$\+bimodule in
which the left and right actions of $R$ agree, is a quasi-module
over~$R$ (with $F^{(r)}_0M=M$ for all $r\in R$).
 Moreover, one has $F_0M=M$; hence $M$ is even a strong quasi-module
over~$R$.
 So we have $R\Modl\subset R\QMod\subset R\biMod R$.

 Given two $R$\+$R$\+bimodules $A$ and~$B$ (in particular, two
quasi-modules over~$R$) we will use the notation $A\ot_RB$ in
the usual sense of the tensor product of the right $R$\+module
structure on $A$ and the left $R$\+module structure on~$B$.
 So the right $R$\+module structure on $A$ and the left and
$R$\+module structure on $B$ get eaten up in the construction of
the tensor product, which becomes an $R$\+$R$\+bimodule with
the left action of $R$ coming from the left action of $R$ in $A$
and the right action of $R$ coming from the right action of $R$ in~$B$.

 Similarly, given an $R$\+$R$\+bimodule $B$ and an $R$\+module $M$,
the notation $B\ot_RM$ stands for the tensor product of the right
$R$\+module structure on $B$ with the $R$\+module structure of
$M$, while $M\ot_RB$ is the tensor product of the $R$\+module
structure on $M$ with the left $R$\+module structure on~$B$.
 Viewing $M$ as an $R$\+$R$\+bimodule in which the left and right
actions of $R$ agree, we obtain $R$\+$R$\+bimodule structures on
$B\ot_RM$ and $M\ot_RB$, as per the previous paragraph.

 One can check that the full subcategory $R\QMod$ is closed under
tensor products over~$R$ in $R\biMod R$ \,\cite[Proposition~2.4]{Ptd}.
 Furthermore, the unit object of the monoidal category $R\biMod R$
(with respect to the tensor product functor~$\ot_R$) is
the diagonal bimodule $R$ in which the left and right actions of $R$
agree, so $R\in R\Modl\subset R\QMod$ according to the discussion
above.
 Thus $R\QMod$ is a monoidal full subcategory in $R\biMod R$.

\begin{lem} \label{quasi-module-restriction-of-scalars}
 Let $f\:R\rarrow S$ be a homomorphism of commutative rings and $B$ be
a quasi-module over~$S$.
 Then, viewed as an $R$\+$R$\+bimodule, $B$ is also a quasi-module
over~$R$.
\end{lem}

\begin{proof}
 Let $r\in R$ be an element and $s=f(r)\in S$ be the image of~$r$
in~$S$.
 Denoting by $F^{(r)}$ the related increasing filtration on $B$ as
an $R$\+$R$\+bimodule and by $G^{(s)}$ the related increasing filtration
on $B$ as an $S$\+$S$\+bimodule, one notices that the two filtrations
coincide, $G^{(s)}_nB = F^{(r)}_nB$ for all $n\ge0$.
 Thus, if the filtration $G^{(s)}$ is exhaustive for all $s\in f(R)
\subset S$, then so is the filtration $F^{(r)}$ for all $r\in R$.
\end{proof}

\begin{lem} \label{quasi-module-element-localization-lemma}
 Let $R$ be a commutative ring and $r\in R$ be an element.
 Let $B$ be a quasi-module over~$R$.
 Then \par
\textup{(a)} the natural maps $R[r^{-1}]\ot_RB\rarrow
R[r^{-1}]\ot_RB\ot_RR[r^{-1}]\larrow B\ot_RR[r^{-1}]$ are
isomorphisms; \par
\textup{(b)} the $R[r^{-1}]$\+$R[r^{-1}]$\+bimodule $R[r^{-1}]\ot_RB$
is a quasi-module over $R[r^{-1}]$.
\end{lem}

\begin{proof}
 The argument is similar to~\cite[proof of Lemma~2.2]{Pedg}.
 Part~(a): it suffices to show that the $R$\+$R[r^{-1}]$\+bimodule
$B\ot_RR[r^{-1}]$, viewed as an $R$\+module with respect to the left
action of $R$, belongs to the essential image of the fully
faithful functor of restriction of scalars $R[r^{-1}]\Modl\rarrow
R\Modl$.
 Indeed, consider the natural filtration $F^{(r)}$ on the
$R$\+$R$\+bimodule $B$ as defined above in this section,
and put $G_n^{(r)}(B\ot_RR[r^{-1}])=(F_n^{(r)}B)\ot_RR[r^{-1}]
\subset B\ot_RR[r^{-1}]$.
 Then $G$ is an exhaustive increasing filtration of
the $R$\+$R[r^{-1}]$\+bimodule $B\ot_RR[r^{-1}]$ by subbimodules
such that the left and right actions of the element $r\in R$
in the successive quotient bimodules
$$
 G_n^{(r)}(B\ot_RR[r^{-1}])/G_{n-1}^{(r)}(B\ot_RR[r^{-1}])
 \simeq (F_n^{(r)}B/F_{n-1}^{(r)}B)\ot_RR[r^{-1}]
$$
agree.
 Hence, viewed as $R$\+modules with respect to the left action of $R$,
the successive quotient bimodules belong to the essential image of
the functor of restriction of scalars $R[r^{-1}]\Modl\rarrow
R\Modl$.
 It remains to point out that this essential image is closed under
extensions and inductive limits in $R\Modl$.

 Part~(b): in view of the filtration argument in part~(a) and the fact
that the class of all quasi-modules over $R[r^{-1}]$ is closed under
extensions and inductive limits in $R[r^{-1}]\biMod R[r^{-1}]$, it
suffices to consider the case when the left and right actions of
the element~$r$ in the bimodule~$B$ coincide.
 In this case, one easily checks that, for any element $s\in R$,
the rule $G_n^{(s)}(R[r^{-1}]\ot_RB)=R[r^{-1}]\ot_R F_n^{(s)}B$ defines
a filtration of the $R[r^{-1}]$\+$R[r^{-1}]$\+bimodule $R[r^{-1}]\ot_RB$
by $R[r^{-1}]$\+$R[r^{-1}]$\+subbimodules such that the left and right
actions of the element~$s$ in the successive quotient bimodules
coincide.
 This proves that the natural increasing filtration $F^{(s)}$ of
the $R[r^{-1}]$\+$R[r^{-1}]$\+bimodule $R[r^{-1}]\ot_RB$ is exhaustive.
 The same assertion in the general case of an element $s\in R[r^{-1}]$
also follows easily.
\end{proof}

 We refer to~\cite[Chapter~XI]{Ste}, \cite[Section~4]{GL},
and~\cite[Section~5]{Ptd} for a discussion of \emph{flat epimorphisms}
of commutative rings.
 Following these references, if $\sigma\:R\rarrow S$ is a flat
epimorphism of commutative rings, then the natural maps
$S\rightrightarrows S\ot_RS\rarrow S$ are isomorphisms, the functor of
restriction of scalars $\sigma_*\:S\Modl\rarrow R\Modl$ is fully
faithful, and its essential image $\sigma_*(S\Modl)$ is closed under
extensions, infinite direct sums, infinite direct products, kernels,
and cokernels in $R\Modl$.

 In particular, a homomorphism of commutative rings $\sigma\:R\rarrow S$
induces an open immersion of affine schemes $\Spec S\rarrow\Spec R$
if and only if $\sigma$~is a \emph{flat epimorphism of finite
presentation}, i.~e., a flat epimorphism of commutative rings
making $S$ a finitely presented
$R$\+algebra~\cite[Th\'eor\`eme~IV.17.9.1]{EGA4}.

 The following lemma is a generalization of
Lemma~\ref{quasi-module-element-localization-lemma}.

\begin{lem} \label{quasi-module-localization-lemma}
 Let $R\rarrow S$ be a flat epimorphism of commutative rings and
$B$ be a quasi-module over~$R$.
 Then \par
\textup{(a)} the natural maps $S\ot_RB\rarrow S\ot_RB\ot_RS\larrow
B\ot_RS$ are isomorphisms; \par
\textup{(b)} the $S$\+$S$\+bimodule $S\ot_RB$ is a quasi-module
over~$S$.
\end{lem}

\begin{proof}
 This is~\cite[Proposition~5.2]{Ptd}.
 Let us spell out a separate proof in the special case when the morphism
of affine schemes $\Spec S\rarrow\Spec R$ is an open immersion.

 Let $r_1$,~\dots, $r_N\in R$ be a finite family of elements such
that $\Spec S=\bigcup_{\alpha=1}^N\Spec R[r_\alpha^{-1}]\subset\Spec R$.
 Then, for any $S$\+module $M$, the \v Cech coresolution
\begin{multline} \label{R-module-M-cech-coresolution-over-Spec-S}
 0\lrarrow S\ot_RM\lrarrow\bigoplus\nolimits_{\alpha=1}^N
 R[r_\alpha^{-1}]\ot_RM \\
 \lrarrow\bigoplus\nolimits_{1\le\alpha<\beta\le N}
 R[r_\alpha^{-1},r_\beta^{-1}]\ot_RM \\
 \lrarrow\dotsb\lrarrow R[r_1^{-1},\dotsc,r_N^{-1}]\ot_RM
 \lrarrow0
\end{multline}
is an exact sequence.

 Returning to our quasi-module $B$ over~$R$, consider three
complexes: the \v Cech complex of the left $S$\+module $S\ot_RB$
\begin{multline} \label{bimodule-B-left-cech-coresolution}
 0\lrarrow S\ot_RB\lrarrow\bigoplus\nolimits_{\alpha=1}^N
 R[r_\alpha^{-1}]\ot_RB\lrarrow
 \bigoplus\nolimits_{1\le\alpha<\beta\le N}
 R[r_\alpha^{-1},r_\beta^{-1}]\ot_RB \\
 \lrarrow\dotsb\lrarrow R[r_1^{-1},\dotsc,r_N^{-1}]\ot_RB
 \lrarrow0,
\end{multline}
the \v Cech complex of the right $S$\+module $B\ot_RS$
\begin{multline} \label{bimodule-B-right-cech-coresolution}
 0\lrarrow B\ot_RS\lrarrow\bigoplus\nolimits_{\alpha=1}^N
 B\ot_RR[r_\alpha^{-1}]\lrarrow
 \bigoplus\nolimits_{1\le\alpha<\beta\le N}
 B\ot_RR[r_\alpha^{-1},r_\beta^{-1}] \\
 \lrarrow\dotsb\lrarrow B\ot_RR[r_1^{-1},\dotsc,r_N^{-1}]
 \lrarrow0,
\end{multline}
and the \v Cech complex of the $S$\+$S$\+bimodule $S\ot_RB\ot_RS$
\begin{multline} \label{bimodule-B-two-sided-cech-coresolution}
 0\lrarrow S\ot_RB\ot_RS\lrarrow\bigoplus\nolimits_{\alpha=1}^N
 R[r_\alpha^{-1}]\ot_RB\ot_RR[r_\alpha^{-1}] \\
 \lrarrow\bigoplus\nolimits_{1\le\alpha<\beta\le N}
 R[r_\alpha^{-1},r_\beta^{-1}]\ot_RB\ot_RR[r_\alpha^{-1},r_\beta^{-1}]
 \\ \lrarrow\dotsb\lrarrow R[r_1^{-1},\dotsc,r_N^{-1}]\ot_R
 B\ot_RR[r_1^{-1},\dotsc,r_N^{-1}]\lrarrow0.
\end{multline}

 The complexes~\eqref{bimodule-B-left-cech-coresolution}
and~\eqref{bimodule-B-right-cech-coresolution} are exact as
particular cases of~\eqref{R-module-M-cech-coresolution-over-Spec-S}.
 The complex~\eqref{bimodule-B-two-sided-cech-coresolution} can be
obtained by applying the functor ${-}\ot_RS$ to
the complex~\eqref{bimodule-B-left-cech-coresolution}, as one can
see from Lemma~\ref{quasi-module-element-localization-lemma}(a).
 So the complex~\eqref{bimodule-B-two-sided-cech-coresolution} is
exact as well.

 Now there are natural morphisms of complexes of $R$\+$R$\+bimodules
from both the complexes~\eqref{bimodule-B-left-cech-coresolution}
and~\eqref{bimodule-B-right-cech-coresolution} into
the complex~\eqref{bimodule-B-two-sided-cech-coresolution}.
 Applying Lemma~\ref{quasi-module-element-localization-lemma}(a)
again, we see that both the morphisms of complexes are isomorphisms
on all the terms except perhaps the leftmost ones.
 It follows that both the morphisms of compexes are also isomorphisms
on the leftmost terms, as desired.

 Part~(b): according to~\eqref{bimodule-B-left-cech-coresolution},
\,$S\ot_RB$ is an $S$\+$S$\+subbimodule of the $S$\+$S$\+bimodule
$\bigoplus_{\alpha=1}^N R[r_\alpha^{-1}]\ot_RB$.
 Each direct summand $R[r_\alpha^{-1}]\ot_RB$ is a quasi-module
over $R[r_\alpha^{-1}]$ by
Lemma~\ref{quasi-module-element-localization-lemma}(b), hence
(by Lemma~\ref{quasi-module-restriction-of-scalars}) also
a quasi-module over~$S$.
 It follows that $S\ot_RB$ is a quasi-module over~$S$.
\end{proof}

\begin{cor} \label{for-quasi-module-quasi-coherence-aux-cor}
 Let $R\rarrow S$ be a flat epimorphism of commutative rings.
 Let $A$ be an $R$\+$R$\+bimodule, $B$ be an $S$\+$S$\+bimodule,
and $A\rarrow B$ be a homomorphism of $R$\+$R$\+bimodules.
 Assume that the $R$\+$R$\+bimodule $A$ is a quasi-module over~$R$.
 Then the induced morphism of $S$\+$R$\+bimodules $S\ot_RA\rarrow B$
is an isomorphism if and only if the induced morphism of
$R$\+$S$\+bimodules $A\ot_RS\rarrow B$ is an isomorphism.
\end{cor}

\begin{proof}
 This is our version of~\cite[Corollary~2.3]{Pedg}.
 The assertion follows from
Lemma~\ref{quasi-module-localization-lemma}(a).
\end{proof}

\begin{lem} \label{quasi-module-locality}
 Let $R\rarrow S_\alpha$ be a finite collection of homomorphisms of
commutative rings such that the collection of induced maps of
the spectra\/ $\Spec S_\alpha\rarrow\Spec R$ is an affine open covering
of an affine scheme.
 Let $B$ be an $R$\+$R$\+bimodule.
 Then the following two conditions are equivalent:
\begin{enumerate}
\item $B$ is a quasi-module over~$R$;
\item for every index~$\alpha$, the right $R$\+module structure on
the tensor product $S_\alpha\ot_R B$ comes from a (necessarily unique)
$S_\alpha$\+module structure, and the $S_\alpha$\+$S_\alpha$\+bimodule
$S_\alpha\ot_RB$ is a quasi-module over~$S_\alpha$.
\end{enumerate}
\end{lem}

\begin{proof}
 If $f\:R\rarrow S$ is a commutative ring homomorphism such that
the induced map of the spectra $\Spec S\rarrow\Spec R$ is an open
immersion, then $f$~is a ring epimorphism.
 So the functor of restriction of scalars $S\Modl\rarrow R\Modl$ is
fully faithful, and the $S$\+module structure of any $S$\+module
is uniquely determined by its underlying $R$\+module structure.
 This explains the uniqueness assertion in~(2).

 Now the implication (1)\,$\Longrightarrow$\,(2) holds by
Lemma~\ref{quasi-module-localization-lemma}.
 To prove that (2)\,$\Longrightarrow$\,(1), notice that the natural
map $B\rarrow\bigoplus_\alpha S_\alpha\ot_R B$ is an $R$\+module
monomorphism.
 By Lemma~\ref{quasi-module-restriction-of-scalars}, the condition
that $S_\alpha\ot_R B$ is a quasi-module over $S_\alpha$ implies that
$S_\alpha\ot_R B$ is a quasi-module over~$R$.
 It remains to point out that the (finite) direct sums of quasi-modules
are quasi-modules and any subbimodule of a quasi-module is
a quasi-module.
\end{proof}

 Given an $R$\+$R$\+bimodule $B$ and an $R$\+module $E$, the notation
$\Hom_R(B,E)$ in part~(b) of the following lemma stands for the abelian
group of all maps $B\rarrow E$ that are $R$\+linear with respect to
the left action of $R$ in $B$ and the action of $R$ in~$E$.
 The abelian group $\Hom_R(B,E)$ is endowed with the $R$\+module
structure induced by the right action of $R$ in~$B$.
 The notation $\Ext^*_R(B,E)$ is understood similarly.

\begin{lem} \label{very-flat-quasi-module-lemma}
 Let $R$ be a commutative ring, let $B$ be a quasi-module over $R$,
and let $G$ and $E$ be $R$\+modules.  In this context: \par
\textup{(a)} if the left $R$\+module $B$ is very flat and
the $R$\+module $G$ is very flat, then the left $R$\+module $B\ot_RG$
is very flat; \par
\textup{(b)} if the left $R$\+module $B$ is very flat and
the $R$\+module $E$ is contraadjusted, then the $R$\+module\/
$\Hom_R(B,E)$ is contraadjusted.
\end{lem}

\begin{proof}
 This is a generalization of
Lemma~\ref{vfl-cta-tensor-hom-lemma}.
 Notice that the condition that $B$ is a quasi-module cannot be dropped:
the assertions of the lemma are \emph{not} true for an arbitrary
$R$\+$R$\+bimodule~$B$.
 See~\cite[Example~6.2]{Pphil} and
Example~\ref{quasi-module-assumption-necessary-counterex} below.

 To prove part~(b), we need to show that
$\Ext^1_R(R[s^{-1}],\Hom_R(B,E))=0$ for all $s\in R$.
 One approach is to compute that
\begin{multline*}
 \Ext^1_R(R[s^{-1}],\Hom_R(B,E))\simeq
 \Ext^1_R(B\ot_R R[s^{-1}],\>E) \\
 \simeq\Ext^1_R(R[s^{-1}]\ot_R B,\>E)
 \simeq\Ext^1_R(B,\Hom_R(R[s^{-1}],E))=0
\end{multline*}
using Lemma~\ref{Ext-homological-formula} twice,
Lemma~\ref{quasi-module-element-localization-lemma}(a),
and Lemma~\ref{vfl-cta-tensor-hom-lemma}(b).

 Alternatively, here is a more elementary argument in
the spirit of the proof of Lemma~\ref{hom-injective-cotorsion}.
 Let $0\rarrow Q\rarrow P\rarrow R[s^{-1}]\rarrow0$ be a two-term
projective resolution of the $R$\+module $R[s^{-1}]$, as
in~\cite[proof of Lemma~2.1]{Pcta}.
 We have to show that the map $\Hom_R(P,\Hom_R(B,E))\rarrow
\Hom_R(Q,\Hom_R(B,E))$ is surjective.
 The latter maps is isomorphic to the map $\Hom_R(B\ot_RP,\>E)
\rarrow\Hom_R(B\ot_RQ,\>E)$.
 Now $0\rarrow B\ot_RQ\rarrow B\ot_RP\rarrow B\ot_RR[s^{-1}]\rarrow0$
is a two-term resolution of the $R$\+module $B\ot_RR[s^{-1}]$ by
very flat $R$\+modules $B\ot_RP$ and $B\ot_RQ$ (in the left
$R$\+module structures).
 It remains to point out that $B\ot_RR[s^{-1}]\simeq R[s^{-1}]\ot_RB$
is a very flat $R$\+module (in the left $R$\+module structure) by
Lemmas~\ref{quasi-module-element-localization-lemma}(a)
and~\ref{vfl-cta-tensor-hom-lemma}(a).

 To deduce part~(a), we need to show that $\Ext^1_R(B\ot_RG,\>E)=0$
for any contraadjusted $R$\+module~$E$.
 Once again,
$$
 \Ext^1_R(B\ot_RG,\>E)\simeq\Ext^1_R(G,\Hom_R(B,E))=0
$$
by Lemma~\ref{Ext-homological-formula} and part~(b) of the present
lemma.
 Alternatively, one can use
Lemmas~\ref{quasi-module-element-localization-lemma}(a)
and~\ref{vfl-cta-tensor-hom-lemma}(a) together with the description of
the $R$\+module $G$ as a direct summand of a transfinitely iterated
extension of the $R$\+modules $R[s^{-1}]$, \,$s\in R$ (as mentioned in
Section~\ref{prelim-very-flat-subsecn}).
\end{proof}

\begin{ex} \label{quasi-module-assumption-necessary-counterex}
 Let us present an explicit counterexample showing that
Lemma~\ref{very-flat-quasi-module-lemma} it \emph{not} true for
an arbitrary $R$\+$R$\+bimodule~$B$.
 As one of the proofs of Lemma~\ref{very-flat-quasi-module-lemma}(a)
deduces it from Lemma~\ref{very-flat-quasi-module-lemma}(b),
any refutation of part~(a) of the lemma with the quasi-module
assumption on $B$ dropped would also refute part~(b) without
the quasi-module assumption on~$B$.

 According to~\cite[Example~6.2]{Pphil} (based
on~\cite[Example~9.7]{PSl1}), there exists a pair of commutative
rings $R$ and $S$, an $R$\+very flat $R$\+$S$\+bimodule $F$, and
a very flat $S$\+module $G$ such that the tensor product $F\ot_SG$
is not very flat as a left $R$\+module.
 Let us consider the Cartesian product ring $T=R\times S$, and view
$F$ as a $T$\+$T$\+bimodule with the left and right $T$\+module
structures induced from the $R$\+module and $S$\+module structures
via the projection morphisms of rings $T\rarrow R$ and $T\rarrow S$.
 Similarly, let us consider $G$ as a $T$\+module.
 Notice that an $R$\+module (or an $S$\+module) is very flat if and
only if it is a very flat as a $T$\+module.
 Therefore, in the situation at hand $F$ is very flat as a left
$T$\+module and $G$ is very flat as a $T$\+module, but the $T$\+module
$F\ot_TG=F\ot_SG$ is \emph{not} very flat.
\end{ex}

\subsection{Quasi-algebras} \label{prelim-quasi-algebras-subsecn}
 Let $R$ be a commutative ring.
 A \emph{quasi-algebra} $A$ over $R$ is an associative, unital ring
endowed with a ring homomorphism $R\rarrow A$ such that, viewed as
an $R$\+$R$\+bimodule with respect to the induced $R$\+$R$\+bimodule
structure, $A$ becomes a quasi-module over~$R$.
 Alternatively, a quasi-algebra over $R$ is the same thing as a monoid
object in the monoidal category $R\QMod$ defined in
Section~\ref{prelim-quasi-modules-subsecn}.
 So a quasi-algebra $A$ over $R$ is a quasi-module endowed with
$R$\+$R$\+bimodule morphisms of \emph{multiplication} $A\ot_RA\rarrow A$
and \emph{unit} $R\rarrow A$ satisfying the usual associativity
and unitality axioms.

\begin{exs} \label{diffoperators-quasi-algebras-quasi-modules-exs}
 Let $K\rarrow R$ be a homomorphism of commutative rings and $U$ be
an $R$\+module.
 Consider the $K$\+algebra $E=\Hom_K(U,U)$ of all $K$\+linear maps
$U\rarrow U$.
 Then the action of $R$ in itself provides a $K$\+algebra homomorphism
$R\rarrow E$.
 In particular, $E$ is an $R$\+$R$\+bimodule.
 
 Let $\cD_{R/K}(U,U)=\bigcap_{r\in R}\bigcup_{n\ge0}F^{(r)}_nE\subset E$
be the maximal $R$\+$R$\+subbimodule of $E$ that is a quasi-module
over~$R$.
 By~\cite[Corollary~3.9 and the subsequent paragraph]{Ptd},
\,$\cD_{R/K}(U,U)$ is a $K$\+subalgebra in~$E$ (obviously, containing
the image of $R$ in~$E$, which is even contained in the $K$\+subalgebra
$\Hom_R(U,U)=\bigcap_{r\in R}F^{(r)}_0E$ of all $R$\+linear maps
$U\rarrow U$).
 Thus $\cD_{R/K}(U,U)$ is a quasi-algebra over~$R$.
 We call the quasi-algebra $\cD_{R/K}(U,U)$ the \emph{quasi-algebra
of $K$\+linear $R$\+differential operators $U\rarrow U$}.
 These are the $R$\+differential operators as defined
in~\cite[Section~3.3]{Ptd}, i.~e., in the most general sense among
the three notions of a differential operator discussed
in~\cite[Section~3]{Ptd}.

 Furthermore, let $\cD_{R/K}^\qu(U,U)=\bigcup_\beta F_\beta E\subset E$
be the maximal $R$\+$R$\+subbimodule of $E$ that is a quite quasi-module
over~$R$.
 Here $F$ is the natural ordinal-indexed increasing filtration on
$R$\+$R$\+bimodules defined in
Section~\ref{prelim-quasi-modules-subsecn}, and the direct union
is taken over all ordinals~$\beta$.
 Finally, let $\cD_{R/K}^\str(U,U)=\bigcup_{n<\omega} F_n E\subset E$
be the maximal $R$\+$R$\+subbimodule of $E$ that is a strong quasi-module
over~$R$ (here the direct union is taken over the nonnegative
integers~$n$).
 The elements of $\cD_{R/K}^\str(U,U)$ and $\cD_{R/K}^\qu(U,U)$ are
called the $K$\+linear \emph{strongly $R$\+differential operators} and
\emph{quite $R$\+differential operators} $U\rarrow U$, respectively.

 According to~\cite[the paragraphs after Corollaries~3.2 and~3.6]{Ptd},
all the three $R$\+$R$\+subbimodules $\cD_{R/K}^\str(U,U)\subset
\cD_{R/K}^\qu(U,U)\subset\cD_{R/K}(U,U)$ are $K$\+subalgebras in $E$
(containing the image of $R$ in~$E$); so these are quasi-algebras
over~$R$.
 By construction, the ring $\cD_{R/K}^\str(U,U)$ is endowed with
a natural increasing filtration~$F$ by the orders of the differential
operators, indexed by the nonnegative integers; while the similar
filtration $F$ on the ring $\cD_{R/K}^\qu(U,U)$ is indexed by
the ordinals~\cite[Sections~3.1\+-3.2]{Ptd}.

 More generally, let $U$ and $V$ be two $R$\+modules.
 Then the $K$\+module $\Hom_K(U,V)$ has a natural structure of
$R$\+$R$\+bimodule (where the left and right actions of $K$ agree,
but the left and ring actions of $R$ are different).
 The left action of $R$ on $\Hom_K(U,V)$ is induced by the left
action of $R$ on $V$, while the right action of $R$ on
$\Hom_K(U,V)$ is induced by the right action of $R$ on~$U$.

 Then, as above, $\cD_{R/K}(U,V)=\bigcap_{r\in R}\bigcup_{n\ge0}
F^{(r)}_nE\subset E$ is the maximal $R$\+$R$\+sub\-bi\-mod\-ule of $E$
that is a quasi-module over~$R$.
 We call $\cD_{R/K}(U,V)$ the \emph{quasi-module of $K$\+linear
$R$\+differential operators $U\rarrow V$} \,\cite[Section~3.3]{Ptd}.
 The natural $R$\+$R$\+subbimodules $\cD^\str_{R/K}(U,V)\subset
\cD^\qu_{R/K}(U,V)\subset\cD_{R/K}(U,V)$ of \emph{$K$\+linear
strongly $R$\+differential operators} and \emph{$K$\+linear quite
$R$\+differential operators} are constructed similarly to the above
as the maximal $R$\+$R$\+subbimodules in $E$ that are a strong
quasi-module and a quite quasi-module over $R$, respectively.
 See~\cite[Sections~3.1\+-3.2]{Ptd}.

 When $R$ is a finitely generated $K$\+algebra, all the three classes
of differential operators described above coincide, i.~e., one has
$\cD^\str_{R/K}(U,V)=\cD^\qu_{R/K}(U,V)=\cD_{R/K}(U,V)$.
 This follows from~\cite[Lemmas~1.3 and~2.1]{Ptd}; see the discussion
in~\cite[Sections~1.3, 2.3, and~3.3]{Ptd}.
\end{exs}

\begin{lem} \label{quasi-algebra-co-extension-of-scalars}
 Let $R\rarrow S$ be a flat epimorpism of commutative rings and
$A$ be a quasi-algebra over~$R$.
 In this context: \par
\textup{(a)} The tensor product $S\ot_RA\simeq S\ot_RA\ot_RS\simeq
A\ot_RS$ (as per Lemma~\ref{quasi-module-localization-lemma})
has a unique structure of associative ring for which the natural map
$S\rarrow S\ot_RA$ is a ring homomorphism inducing the given
$S$\+$S$\+bimodule structure on $S\ot_RA$ and $A\rarrow S\ot_RA$
is a ring homomorphism.
 Endowed with this ring structure, $S\ot_RA$ becomes a quasi-algebra
over~$S$. \par
\textup{(b)} For any left $A$\+module $M$, the left $S$\+module of
tensor product $S\ot_RM$ admits a unique extension of its $S$\+module
structure to a structure of left module over $S\ot_RA$ such that
$M\rarrow S\ot_RM$ is an $A$\+module morphism. \par
\textup{(c)} For any left $A$\+module $P$, the left $S$\+module of
homomorphisms\/ $\Hom_R(S,P)$ admits a unique extension of its
$S$\+module structure to a structure of left module over $S\ot_RA$
such that $\Hom_R(S,P)\rarrow P$ is an $A$\+module morphism.
\end{lem}

\begin{proof}
 Part~(a): the uniqueness holds because an arbitrary $S$\+$S$\+bimodule
map $(S\ot_R\nobreak A)\ot_\boZ(A\ot_R\nobreak S)\rarrow S\ot_RA\ot_RS$
is determined by its composition with the map $A\ot_\boZ A\rarrow
(S\ot_R\nobreak A)\ot_\boZ(A\ot_R\nobreak S)$.
 To prove existence, one observes that the functor $B\longmapsto S\ot_RB
\simeq S\ot_RB\ot_RS\simeq B\ot_RS$ is a monoidal functor $R\QMod\rarrow
S\QMod$ (as one can easily check); so it takes monoid objects to
monoid objects.

 Part~(b): the uniqueness holds because an arbitrary $S$\+module map
$(S\ot_RA)\ot_S(S\ot_RM)\simeq S\ot_R(A\ot_RM)\rarrow S\ot_RM$ is
determined by its composition with the map $A\ot_RM\rarrow
S\ot_R(A\ot_RM)$.
 To prove existence, one can say that, conversely, the composition
$(S\ot_RA)\ot_S(S\ot_RM)\simeq S\ot_R(A\ot_RM)\rarrow S\ot_RM$, with
the map $S\ot_R(A\ot_RM)\rarrow S\ot_RM$ induced by the left action map
$A\ot_RM\rarrow M$, defines the desired module structure on $S\ot_RM$.
 Alternatively, to establish the existence it suffices to point out
the natural isomorphism $(S\ot_RA)\ot_AM\simeq S\ot_RM$, which allows
to define the desired $(S\ot_R\nobreak A)$\+module $S\ot_RM$ as
the extension of scalars of the $A$\+module $M$ with respect to the ring
homomorphism $A\rarrow S\ot_RA$.

 Part~(c): the uniqueness holds because an arbitrary $S$\+module map
$\Hom_R(S,P)\rarrow\Hom_S(A\ot_RS,\>\Hom_R(S,P))\simeq
\Hom_R(A\ot_RS,\>P)\simeq\Hom_R(S,\Hom_R(A,P))$ is determined by its
composition with the map $\Hom_R(S,\Hom_R(A,P))\rarrow\Hom_R(A,P)$.
 To prove existence, one can say that, conversely, the composition
$\Hom_R(S,P)\rarrow\Hom_R(S,\Hom_R(A,P))\simeq\Hom_R(A\ot_RS,\>P)
\simeq\Hom_S(A\ot_RS,\>\Hom_R(S,P))$, with the map $\Hom_R(S,P)
\rarrow\Hom_R(S,\allowbreak\Hom_R(A,P))$ induced by the left action
map $P\rarrow\Hom_R(A,P)$, defines the desired module structure
on $\Hom_R(S,P)$.
 Alternatively, to establish the existence it suffices to point out
the natural isomorphism $\Hom_A(A\ot_RS,\>P)\simeq\Hom_R(S,P)$, which
allows to define the desired $(S\ot_R\nobreak A)$\+module $\Hom_R(S,P)$
as the coextension of scalars of the $A$\+module $P$ with respect to
the ring homomorphism $A\rarrow A\ot_RS$. \hfuzz=1.6pt
\end{proof}

\begin{lem} \label{quasi-algebra-change-of-scalars-tensor-hom}
 Let $R\rarrow S$ be a flat epimorphism of commutative rings and
$A$ be a quasi-algebra over~$R$.
 Then \par
\textup{(a)} for any left $A$\+module $M$ and any right module $N$
over $S\ot_RA$, the map of tensor products $N\ot_AM\rarrow
N\ot_{S\ot_RA}(S\ot_RM)$ induced by the ring and module maps
$A\rarrow S\ot_RA$ and $M\rarrow S\ot_RM$ is an isomorphism; \par
\textup{(b)} for any left $A$\+module $M$ and any left module $Q$
over $S\ot_RA$, the map of Hom modules $\Hom_{S\ot_RA}(S\ot_RM,\>Q)
\rarrow\Hom_A(M,Q)$ induced by the ring and module maps $A\rarrow
S\ot_RA$ and $M\rarrow S\ot_RM$ is an isomorphism; \par
\textup{(c)} for any left module $N$ over $S\ot_RA$ and any left
$A$\+module $P$, the map of Hom modules $\Hom_{S\ot_RA}(N,\Hom_R(S,P))
\rarrow\Hom_A(N,P)$ induced by the ring and module maps $A\rarrow
S\ot_RA$ and $\Hom_R(S,P)\rarrow P$ is an isomorphism.
\end{lem}

\begin{proof}
 All the assertions follow from the alternative proofs of
Lemma~\ref{quasi-algebra-co-extension-of-scalars}(b\+-c).
\end{proof}

\begin{cor} \label{quasi-algebra-change-of-scalars-adjustedness}
 Let $R\rarrow S$ be a flat epimorphism of commutative rings and
$A$ be a quasi-algebra over~$R$.
 Then \par
\textup{(a)} the ring $S\ot_RA\simeq A\ot_RS$ is flat as a left and
right module over~$A$; \par
\textup{(b)} for any flat left $A$\+module $F$, the left module
$S\ot_RF$ is flat over $S\ot_RA$; \par
\textup{(c)} for any projective left $A$\+module $F$, the left module
$S\ot_RF$ is projective over $S\ot_RA$; \par
\textup{(d)} for any cotorsion left $A$\+module $C$, the left module\/
$\Hom_R(S,C)$ is cotorsion over $S\ot_RA$; \par
\textup{(e)} for any injective left $A$\+module $J$, the left module\/
$\Hom_R(S,J)$ is injective over $S\ot_RA$.
\end{cor}

\begin{proof}
 Part~(a) holds since $S$ is a flat $R$\+module.
 Parts~(b\+-c) and~(e) follow either from the respective parts of
Lemma~\ref{quasi-algebra-change-of-scalars-tensor-hom}, or directly
from the alternative proofs of
Lemma~\ref{quasi-algebra-co-extension-of-scalars}(b\+-c).
 In the latter approach, in part~(e),
Lemma~\ref{restriction-coextension-injective-cotorsion}(d) is
relevant.
 In part~(d), one has to use part~(a), which implies that any flat
module over $S\ot_RA$ is also flat as an $A$\+module.
 One can also refer to
Lemma~\ref{restriction-coextension-injective-cotorsion}(c).
\end{proof}

 Given a ring homomorphism $R\rarrow A$ with a commutative ring $R$,
an $A$\+module is said to be \emph{$R$\+contraadjusted} if it is
contraadjusted as an $R$\+module.
 \emph{$R$\+cotorsion $A$\+modules} and \emph{$R$\+injective
$A$\+modules} (for an arbitrary homomorphism of associative rings
$R\rarrow A$) are defined similarly.

\begin{lem} \label{quasi-algebra-adjustedness-co-locality}
 Let $R\rarrow S_\alpha$, \ $1\le\alpha\le N$, be a finite collection
of homomorphisms of commutative rings such that the collection of
induced maps of the spectra\/ $\Spec S_\alpha\rarrow\Spec R$ is
an affine open covering of an affine scheme.
 Let $A$ be a quasi-algebra over~$R$.  Then \par
\textup{(a)} a left $A$\+module $F$ is flat if and only if the left
module $S_\alpha\ot_RF$ is flat over $S_\alpha\ot_RA$ for all
indices~$\alpha$; \par
\textup{(b)} an $R$\+contraadjusted left $A$\+module $C$ is cotorsion
over $A$ if and only if the left module\/ $\Hom_R(S_\alpha,C)$ is
cotorsion over $S_\alpha\ot_RA$ for all indices~$\alpha$; \par
\textup{(c)} an $R$\+contraadjusted left $A$\+module $J$ is injective
over $A$ if and only if the left module\/ $\Hom_R(S_\alpha,J)$ is
injective over $S_\alpha\ot_RA$ for all indices~$\alpha$.
\end{lem}

\begin{proof}
 Part~(a): the ``only if'' assertion holds by
Corollary~\ref{quasi-algebra-change-of-scalars-adjustedness}(b).
 To prove the ``if'', consider the \v Cech coresolution
\begin{multline} \label{module-over-quasi-algebra-cech-coresolution}
 0\lrarrow F\lrarrow\bigoplus\nolimits_{\alpha=1}^N S_\alpha\ot_RF
 \lrarrow\bigoplus\nolimits_{1\le\alpha<\beta\le N}
 S_\alpha\ot_RS_\beta\ot_RF \\
 \lrarrow\dotsb\lrarrow S_1\ot_R S_2\ot_R\dotsb\ot_R S_N\ot_RF
 \lrarrow0
\end{multline}
(cf.~\cite[formula~(1.4) in the proof of Lemma~1.2.6(a)]{Pcosh}).
 The complex~\eqref{module-over-quasi-algebra-cech-coresolution} is
an exact sequence of left $A$\+modules (see
Lemma~\ref{quasi-algebra-co-extension-of-scalars}(b)).
 For any subsequence of indices $1\le\alpha_1<\dotsb<\alpha_k\le N$,
\,$k\ge1$, the tensor product $S_{\alpha_1}\ot_R\dotsb\ot_R S_{\alpha_k}
\ot_RF$ is a flat left module over the ring $S_{\alpha_1}\ot_R\dotsb
\ot_R S_{\alpha_k}\ot_RA$ by
Corollary~\ref{quasi-algebra-change-of-scalars-adjustedness}(b).
 Since $S_{\alpha_1}\ot_R\dotsb\ot_R S_{\alpha_k}\ot_RA$ is a flat left
$A$\+module (see
Corollary~\ref{quasi-algebra-change-of-scalars-adjustedness}(a)), it
follows that $S_{\alpha_1}\ot_R\dotsb\ot_R S_{\alpha_k}\ot_RF$ is a flat
left $A$\+module.
 Now the assertion that $F$ is a flat left $A$\+module is provable
by moving along the exact
sequence~\eqref{module-over-quasi-algebra-cech-coresolution} from
its right end to the left one and using the fact that the class of flat
$A$\+modules is closed under the kernels of epimorphisms.

 Part~(b): the ``only if'' assertion holds by
Corollary~\ref{quasi-algebra-change-of-scalars-adjustedness}(d).
 To prove the ``if'', consider the \v Cech resolution
\begin{multline} \label{module-over-quasi-algebra-cech-resolution}
 0\lrarrow\Hom_R(S_1\ot_R\dotsb\ot_RS_N,\>C)\lrarrow\dotsb \\
 \lrarrow\bigoplus\nolimits_{\alpha<\beta}\Hom_R(S_\alpha\ot_R
 S_\beta,\>C)\lrarrow\bigoplus\nolimits_\alpha\Hom_R(S_\alpha,C)
 \lrarrow C\lrarrow0
\end{multline}
as per~\cite[formula~(1.3) in Lemma~1.2.6(b)]{Pcosh}.
 The complex~\eqref{module-over-quasi-algebra-cech-resolution} is
an exact sequence of left $A$\+modules (see
Lemma~\ref{quasi-algebra-co-extension-of-scalars}(c)).
 For any subsequence of indices $1\le\alpha_1<\dotsb<\alpha_k\le N$,
\,$k\ge1$, the Hom module $\Hom_R(S_{\alpha_1}\ot_R\dotsb\ot_R 
S_{\alpha_k},\>C)$ is a cotorsion left module over the ring
$S_{\alpha_1}\ot_R\dotsb \ot_R S_{\alpha_k}\ot_RA$ by
Corollary~\ref{quasi-algebra-change-of-scalars-adjustedness}(d).
 By Lemma~\ref{restriction-coextension-injective-cotorsion}(a),
it follows that $\Hom_R(S_{\alpha_1}\ot_R\dotsb\ot_R S_{\alpha_k},\>C)$
is a cotorsion left $A$\+module.
 Now the assertion that $C$ is a cotorsion left $A$\+module is provable
by moving along the exact
sequence~\eqref{module-over-quasi-algebra-cech-resolution} from
its left end to the right one and using the fact that the class of
cotorsion $A$\+modules is closed under cokernels of monomorphisms
(Lemma~\ref{flat-cotorsion-pair-hereditary}(b)).

 The proof of part~(c) is similar to that of part~(b) and also based
on the \v Cech
resolution~\eqref{module-over-quasi-algebra-cech-resolution}.
 One has to use
Corollary~\ref{quasi-algebra-change-of-scalars-adjustedness}(a,e)
and Lemma~\ref{restriction-coextension-injective-cotorsion}(b).
\end{proof}

\begin{quest} \label{quasi-algebra-relative-flatness-localization-quest}
 Let $R$ be a commutative ring and $A$ be a quasi-algebra over~$R$.
 Let $R\rarrow S$ be a homomorphism of commutative rings such that
the induced morphism of affine schemes $\Spec S\rarrow\Spec R$ is
an open immersion (e.~g., $S=R[r^{-1}]$, where $r\in R$ is an element).

 Let $M$ be a left $A$\+module.
 Then, by Lemma~\ref{quasi-algebra-co-extension-of-scalars}(a\+-b),
\,$S\ot_RA$ is a quasi-algebra over $S$ and $S\ot_RM$ is a left module
over $S\ot_RA$.
 By Corollary~\ref{quasi-algebra-change-of-scalars-adjustedness}(a\+-b),
if $F$ is a flat left $A$\+module, then $S\ot_RF$ is a flat left
module over $S\ot_RA$, and consequently also a flat left $A$\+module.

 Now let $E$ be an associative ring and $E\rarrow A$ be a homomorphism
of associative rings.
 Let $M$ be a left $A$\+module that is flat as a module over~$E$.
 Does it follow that the left $A$\+module $S\ot_RM$ is flat as a module
over~$E$\,?
\end{quest}

 Nevertheless, the following lemma is obvious.

\begin{lem} \label{quasi-algebra-change-of-scalars-ring-flatness}
 Let $R\rarrow S$ be a flat epimorphism of commutative rings and
$A$ be a quasi-algebra over~$R$.
 Let $A$ be a quasi-algebra over~$R$.
 Let $E$ be an associative ring and $E\rarrow A$ be a ring homomorphism.
 Assume that the ring $A$ is flat as a left $E$\+module.
 Then the ring $S\ot_RA\simeq A\ot_RS$ is also flat as a left
$E$\+module.
\end{lem}

\begin{proof}
 Follows immediately from
Corollary~\ref{quasi-algebra-change-of-scalars-adjustedness}(a).
\end{proof}

\begin{lem} \label{quasi-algebra-ring-flatness-locality}
 Let $R\rarrow S_\alpha$, \ $1\le\alpha\le N$, be a finite collection
of homomorphisms of commutative rings such that the collection of
induced maps of the spectra\/ $\Spec S_\alpha\rarrow\Spec R$ is an
affine open covering of an affine scheme.
 Let $A$ be a quasi-algebra over $R$, and let $E\rarrow A$ be a ring
homomorphism.
 Then the left $E$\+module $A$ is flat if and only if the left
$E$\+module $S_\alpha\ot_RA\simeq A\ot_RS_\alpha$ is flat for
every~$\alpha$.
\end{lem}

\begin{proof}
 The ``only if'' assertion holds by
Lemma~\ref{quasi-algebra-change-of-scalars-ring-flatness}.
 To prove the ``if'', consider the \v Cech
coresolution~\eqref{module-over-quasi-algebra-cech-coresolution}
for $F=A$,
\begin{multline} \label{quasi-algebra-itself-cech-coresolution}
 0\lrarrow A\lrarrow\bigoplus\nolimits_{\alpha=1}^N S_\alpha\ot_RA
 \lrarrow\bigoplus\nolimits_{1\le\alpha<\beta\le N}
 S_\alpha\ot_RS_\beta\ot_RA \\
 \lrarrow\dotsb\lrarrow S_1\ot_R S_2\ot_R\dotsb\ot_R S_N\ot_RA
 \lrarrow0.
\end{multline}
 The complex~\eqref{quasi-algebra-itself-cech-coresolution} is
an exact sequence of $A$\+$A$\+bimodules, hence also an exact
sequence of left $E$\+modules.
 For any subsequence of indices $1\le\alpha_1<\dotsb<\alpha_k\le N$,
\,$k\ge1$, the tensor product $S_{\alpha_1}\ot_R\dotsb\ot_R S_{\alpha_k}
\ot_RA$ is a flat left $E$\+module by
Lemma~\ref{quasi-algebra-change-of-scalars-ring-flatness}.
 The argument finishes similarly to the proof of
Lemma~\ref{quasi-algebra-adjustedness-co-locality}(a).
\end{proof}

\subsection{Flaprojective and very flaprojective modules}
\label{prelim-flaprojective-subsecn}
 Let $R\rarrow A$ be a homomorphism of associative rings.
 The definitions of $R$\+cotorsion $A$\+modules and (if $R$ is
commutative) $R$\+contraadjusted $A$\+modules were given in
the previous Section~\ref{prelim-quasi-algebras-subsecn}.

 We will say that a left $A$\+module $F$ is \emph{flaprojective relative
to~$R$} (or \emph{$A/R$\+fla\-pro\-jec\-tive}) if $\Ext^1_A(F,C)=0$ for
all $R$\+cotorsion left $A$\+modules~$C$.
 If the ring $R$ is commutative, we will say that a left $A$\+module $F$
is \emph{very flaprojective relative to~$R$} (or \emph{$A/R$\+very
flaprojective}) if $\Ext^1_A(F,C)=0$ for all $R$\+contraadjusted left
$A$\+modules~$C$.
 Clearly, the classes of $A/R$\+flaprojective and $A/R$\+very 
flaprojective $A$\+modules are closed under extensions, direct summands,
and infinite direct sums.

\begin{exs} \label{relative-flaprojectivity-degenerated-cases}
 (1)~If $A=R$ and $R\rarrow A$ is the identity map, then the class
of $A/R$\+flaprojective $A$\+modules modules coincides with the class
of flat $R$\+modules.
 The class of $A/R$\+very flaprojective $A$\+modules coincides with
the class of very flat $R$\+modules in this case.

 (2)~If $R=k$ is a field, then the both the classes of
$A/R$\+flaprojective $A$\+modules and $A/R$\+very flaprojective
$A$\+modules coincide with the class of projective $A$\+modules.
\end{exs}

\begin{rem} \label{very-flaprojective-terminology-explained}
 Explaining the terminology ``flaprojective'', one can say informally,
based on Examples~\ref{relative-flaprojectivity-degenerated-cases},
that the $A/R$\+flaprojective $A$\+modules are ``flat along~$R$'' and
``projective in the direction of $A$ relative to~$R$''.
 Similarly, the $A/R$\+very flaprojective $A$\+modules are
``very flat along~$R$'' and ``projective in the direction of $A$
relative to~$R$''.

 The informal explanation above should be taken with caution.
 For example, all projective $A$\+modules are obviously
$A/R$\+flaprojective by the definition (and even $A/R$\+very 
flaprojective, if $R$ is commutative).
 Still, the projective $A$\+module $A$ need not be flat as
an $R$\+module in general.

 On the other hand, all cotorsion $A$\+modules are also cotorsion
as $R$\+modules by
Lemma~\ref{restriction-coextension-injective-cotorsion}(a) (hence
also contraadjusted as $R$\+modules if $R$ is commutative).
 Therefore, Lemma~\ref{flat-cotorsion-is-a-cotorsion-pair} implies
that all $A/R$\+flaprojective $A$\+modules are flat as $A$\+modules.
 In particular, if $R$ is commutative, then all
$A/R$\+very flaprojective $A$\+modules are flat as $A$\+modules.
 If the left $R$\+module $A$ is flat, it follows that all
$A/R$\+(very) flaprojective $A$\+modules are flat as $R$\+modules.

 See also the discussion in~\cite[Remarks~3.3]{PBas}.
\end{rem}

\begin{lem} \label{flaprojective-cotorsion-pair-hereditary}
\textup{(a)} One has\/ $\Ext^n_A(F,C)=0$ for all $A/R$\+flaprojective
$A$\+modules $F$, all $R$\+cotorsion $A$\+modules $C$, and all $n\ge1$.
\par
\textup{(b)} The class of $A/R$\+flaprojective $A$\+modules is closed
under kernels of epimorphisms in $A\Modl$.
\end{lem}

\begin{proof}
 Part~(a): the point is that all injective $A$\+modules are cotorsion
by definition, hence also $R$\+cotorsion by
Lemma~\ref{restriction-coextension-injective-cotorsion}(a).
 Furthermore, the class of $R$\+cotorsion $A$\+modules is closed
under cokernels of monomorphisms in $A\Modl$ by
Lemma~\ref{flat-cotorsion-pair-hereditary}(b).
 With these observations in mind, the assertion of part~(a) follows
by applying~\cite[Lemma~7.1]{PS6}.
 Part~(b) follows immediately from part~(a).
\end{proof}

\begin{lem} \label{very-flaprojective-cotorsion-pair-hereditary}
 Assume that the ring $R$ is commutative. \par
\textup{(a)} One has\/ $\Ext^n_A(F,C)=0$ for all $A/R$\+very
flaprojective $A$\+modules $F$, all $R$\+con\-tra\-ad\-justed
$A$\+modules $C$, and all $n\ge1$. \par
\textup{(b)} The class of $A/R$\+very flaprojective $A$\+modules is
closed under kernels of epimorphisms in $A\Modl$.
\end{lem}

\begin{proof}
 Similar to the proof of
Lemma~\ref{flaprojective-cotorsion-pair-hereditary}.
 Lemma~\ref{very-flat-cotorsion-pair-hereditary}(a) is relevant.
\end{proof}

 In fact, all $A/R$\+very flaprojective $A$\+modules have projective
dimensions at most~$1$.
 So one has $\Ext_A^n(F,C)=0$ for all $A/R$\+very flaprojective
$A$\+modules $F$, \emph{all} $A$\+modules $C$, and all $n\ge2$.
 See~\cite[Example~6.11(2)]{Pfltp}.

\begin{lem} \label{very-flaprojective-is-a-cotorsion-pair}
\textup{(a)} A left $A$\+module $C$ is $R$\+cotorsion if and only if\/
$\Ext_A^1(F,C)=0$ for all $A/R$\+flaprojective left $A$\+modules~$F$.
\par
\textup{(b)} Assuming that the ring $R$ is commutative, a left
$A$\+module $C$ is $R$\+contraadjusted if and only if\/
$\Ext_A^1(F,C)=0$ for all $A/R$\+very flaprojective $A$\+modules~$F$.
\end{lem}

\begin{proof}
 The ``only if'' assertions in~(a) and~(b) hold by the definition.
 To prove the ``if'' implications, one argues as follows.
 For any flat left $R$\+module $F$ and any left $A$\+module $C$,
the natural isomorphism $\Ext^n_A(A\ot_RF,\>C)\simeq\Ext^n_R(F,C)$
holds for all $n\ge0$ by Lemma~\ref{Ext-homological-formula}
(see also~\cite[Lemma~4.1(a)]{PSl1} or~\cite[Lemma~1.7(e)]{Pal}).
 Therefore, for any flat left $R$\+module $F$, the left $A$\+module
$A\ot_RF$ is $A/R$\+flaprojective; and for any very flat $R$\+module
$F$, the left $A$\+module $A\ot_RF$ is $A/R$\+very flaprojective.
 Now the Ext isomorphism above implies that a left $A$\+module $C$ is
$R$\+cotorsion if and only if $\Ext^1_A(A\ot_RF,\>C)=0$ for all
flat left $R$\+modules $F$, and $C$ is $R$\+contraadjusted if and only
if $\Ext^1_A(A\ot_RF,\>C)=0$ for all very flat $R$\+modules~$F$.
\end{proof}

\begin{rem} \label{very-flaprojective-cotorsion-pair-generated-by}
 In the terminology of~\cite{Sal,ET,GT}, the proof of
Lemma~\ref{very-flaprojective-is-a-cotorsion-pair} tells us that
the pair of classes ($A/R$\+flaprojective $A$\+modules, $R$\+cotorsion
$A$\+modules) is the cotorsion pair in $A\Modl$ generated by the class
of all modules of the form $A\ot_RF$, where $F$ ranges over flat
left $R$\+modules~\cite[proof of Theorem~3.5]{PBas}.
 Similarly, the pair of classes ($A/R$\+very flaprojective $A$\+modules,
$R$\+contraadjusted $A$\+modules) is the cotorsion pair in $A\Modl$
generated by the class of all modules of the form $A\ot_RF$, where
$F$ ranges over very flat $R$\+modules.
 As both the classes of flat and very flat $R$\+modules are
deconstructible (in the case of flat modules, this is~\cite[Lemma~1
and Proposition~2]{BBE}), both the cotorsion pairs are actually
generated by suitable subsets of the mentioned two classes of modules.
 The next Theorems~\ref{flaprojective-cotorsion-pair-complete}
and~\ref{very-flaprojective-cotorsion-pair-complete} follow
immediately from these observations.
 It also follows that the classes of $A/R$\+(very) flaprojective
$A$\+modules can be described as the classes of all direct summands
of transfinitely iterated extensions (in the sense of the inductive
limit) of the $A$\+modules $A\ot_RF$, where the $R$\+module $F$ is
(very) flat~\cite[Theorems~2 and~10]{ET}, \cite[Corollary~6.14]{GT}.
 See also~\cite[Proposition~3.3]{Pctrl} for a discussion much
more specific than in~\cite{ET,GT}, still considerably more general
than in the present section.
\end{rem}

\begin{thm} \label{flaprojective-cotorsion-pair-complete}
 For any $A$\+module $M$, there exist \par
\textup{(a)} a short exact sequence of $A$\+modules\/
$0\rarrow C'\rarrow F\rarrow M\rarrow0$ with an $A/R$\+flaprojective
$A$\+module $F$ and an $R$\+cotorsion $A$\+module~$C'$; \par
\textup{(b)} a short exact sequence of $A$\+modules\/
$0\rarrow M\rarrow C\rarrow F'\rarrow0$ with an $R$\+cotorsion
$A$\+module $C$ and an $A/R$\+flaprojective $A$\+module~$F'$.
\end{thm}

\begin{proof}
 This is~\cite[Theorem~3.5(a)]{PBas}.
\end{proof}

\begin{thm} \label{very-flaprojective-cotorsion-pair-complete}
 Assume that the ring $R$ is commutative.
 Then, for any $A$\+module $M$, there exist \par
\textup{(a)} a short exact sequence of $A$\+modules\/
$0\rarrow C'\rarrow F\rarrow M\rarrow0$ with an $A/R$\+very
flaprojective $A$\+module $F$ and an $R$\+contraadjusted
$A$\+module~$C'$; \par
\textup{(b)} a short exact sequence of $A$\+modules\/
$0\rarrow M\rarrow C\rarrow F'\rarrow0$ with an $R$\+contraadjusted
$A$\+module $C$ and an $A/R$\+very flaprojective $A$\+module~$F'$.
\end{thm}

\begin{proof}
 In view of the discussion in
Remark~\ref{very-flaprojective-cotorsion-pair-generated-by},
this is a special cases of~\cite[Theorems~2 and~10]{ET}
or~\cite[Theorem~6.11]{GT}, similarly to
Theorem~\ref{flaprojective-cotorsion-pair-complete}.
 This is also a special case of~\cite[Proposition~3.3]{Pctrl}.
\end{proof}

\begin{cor} \label{relative-bass-theorem-flaprojective-cor}
 Let $R\rarrow A$ be a homomorphism of associative rings such that
$A$ is a finitely generated projective right $R$\+module.
 Then the classes of flat left $A$\+modules and $A/R$\+flaprojective
left $A$\+modules coincide.
\end{cor}

\begin{proof}
 This is a corollary of
Proposition~\ref{relative-bass-theorem-cotorsion-prop}.
 The assertion follows immediately from the definition of
an $A/R$\+flaprojective $A$\+module together with
Proposition~\ref{relative-bass-theorem-cotorsion-prop} and
Lemma~\ref{flat-cotorsion-is-a-cotorsion-pair}.
 See~\cite[Corollary~3.7]{PBas} for a discussion.
\end{proof}

\begin{lem} \label{flaprojective-restriction-extension-of-scalars}
 Let $R\rarrow S$ be a flat epimorphism of commutative rings.
 Let $A$ be a quasi-algebra over $R$ and $B=S\ot_RA=A\ot_RS$ be
the related quasi-algebra over $S$, as per
Lemma~\ref{quasi-algebra-co-extension-of-scalars}(a).
 In this context: \par
\textup{(a)} any $B/S$\+flaprojective $B$\+module is also
$A/R$\+flaprojective as an $A$\+module; \par
\textup{(b)} if $F$ is an $A/R$\+flaprojective left $A$\+module,
then $S\ot_RF$ is a $B/S$\+flaprojective left $B$\+module.
\end{lem}

\begin{proof}
 Part~(a): given an $R$\+cotorsion left $A$\+module $C$, the related
left $B$\+module $\Hom_A(B,C)\simeq\Hom_R(S,C)$ is $S$\+cotorsion by
Lemma~\ref{restriction-coextension-injective-cotorsion}(c).
 Furthermore, one has $\Ext^n_A(B,C)=\Ext^n_A(A\ot_RS,\>C)\simeq
\Ext^n_R(S,C)=0$ for all $n\ge1$ by Lemma~\ref{Ext-homological-formula}
(see also~\cite[Lemma~4.1(a)]{PSl1}).
 Consequently, the natural isomorphism $\Ext^n_A(G,C)\simeq
\Ext^n_B(G,\Hom_A(B,C))$ holds for all left $A$\+modules $G$ and all
$n\ge0$ by Lemma~\ref{Ext-homological-formula} (see
also~\cite[Lemma~1.7(e)]{Pal}).
 Now it is clear that $G$ is $A/R$\+flaprojective whenever it is
$B/S$\+flaprojective.

 Part~(b): recall the isomorphism $S\ot_RF\simeq B\ot_AF$ from
the alternative proof of
Lemma~\ref{quasi-algebra-co-extension-of-scalars}(b).
 Now for any left $B$\+module $D$, one has $\Ext^n_B(B\ot_AF,\>D)\simeq
\Ext_A^n(F,D)$ for all $n\ge1$ by Lemma~\ref{Ext-homological-formula}
(see also~\cite[Lemma~4.1(a)]{PSl1}).
 It remains to point out that any $S$\+cotorsion left $B$\+module $D$
is also $R$\+cotorsion as an $A$\+module by
Lemma~\ref{restriction-coextension-injective-cotorsion}(a).
\end{proof}

\begin{lem} \label{very-flaprojective-restriction-extension-of-scalars}
 Let $R\rarrow S$ be a flat epimorphism of commutative rings.
 Let $A$ be a quasi-algebra over $R$, and $B=S\ot_RA=A\ot_RS$ be
the related quasi-algebra over $S$ as per
Lemma~\ref{quasi-algebra-co-extension-of-scalars}(a).
 In this context: \par
\textup{(a)} if the morphism of affine schemes\/ $\Spec S\rarrow\Spec R$
is an open immersion, then any $B/S$\+very flaprojective $B$\+module
is also $A/R$\+very flaprojective as an $A$\+module; \par
\textup{(b)} if $F$ is an $A/R$\+very flaprojective left $A$\+module,
then $S\ot_RF$ is a very $B/S$\+fla\-pro\-jective left $B$\+module.
\end{lem} 

\begin{proof}
 Similar to the proof of
Lemma~\ref{flaprojective-restriction-extension-of-scalars}.
 One needs to use Lemma~\ref{restriction-co-extension-vfl-cta}(d)
with Example~\ref{very-flat-ring-examples}(1) in part~(a) and
Lemma~\ref{restriction-co-extension-vfl-cta}(a) in part~(b).
\end{proof}

\begin{lem} \label{quasi-algebra-very-flaprojectivity-locality}
 Let $R\rarrow S_\alpha$ be a finite collection of homomorphisms of
commutative rings such that the collection of induced maps of
the spectra\/ $\Spec S_\alpha\rarrow\Spec R$ is an affine open covering
of an affine scheme.
 Let $A$ be a quasi-algebra over~$R$.  Then \par
\textup{(a)} a left $A$\+module $F$ is $A/R$\+flaprojective if and only
if the left module $S_\alpha\ot_RF$ over the ring $B_\alpha=
S_\alpha\ot_RA$ is $B_\alpha/S_\alpha$\+flaprojective for all
indices~$\alpha$; \par
\textup{(b)} a left $A$\+module $F$ is $A/R$\+very flaprojective if and
only if the left module $S_\alpha\ot_RF$ over the ring $B_\alpha=
S_\alpha\ot_RA$ is $B_\alpha/S_\alpha$\+very flaprojective for all
indices~$\alpha$.
\end{lem}

\begin{proof}
 The argument from the proof of
Lemma~\ref{quasi-algebra-adjustedness-co-locality}(a), based on
the \v Cech
coresolution~\eqref{module-over-quasi-algebra-cech-coresolution},
is applicable.
 In part~(a), one needs to use
Lemmas~\ref{flaprojective-restriction-extension-of-scalars}(a\+-b)
and~\ref{flaprojective-cotorsion-pair-hereditary}(b).
 For part~(b), use
Lemmas~\ref{very-flaprojective-restriction-extension-of-scalars}(a\+-b)
and~\ref{very-flaprojective-cotorsion-pair-hereditary}(b).
\end{proof}

\subsection{Robustly flaprojective bimodules}
\label{prelim-robustly-flaprojective-subsecn}
 This section is a bit sketchy.
 It sheds some light on the exposition below in this paper, but
using its results can be avoided (as we will point out).

 Let $R\rarrow A$ be a homomorphism of associative rings and $K$ be
an associative ring.
 We will say that an $A$\+$K$\+bimodule $F$ is \emph{robustly 
flaprojective relative to~$R$} (or \emph{$(A/R,K)$\+robustly
flaprojective}) if the left $A$\+module $F\ot_KH$ is
$A/R$\+flaprojective for every flat left $K$\+module~$H$.
 The class of $(A/R,K)$\+robustly flaprojective $A$\+$K$\+bimodules
is closed under transfinitely iterated extensions and kernels of
epimorphisms (since so is the class of $A/R$\+flaprojective
left $A$\+modules).

\begin{ex}
 For any $R$\+$K$\+bimodule $U$ that is flat as a left $R$\+module,
the $A$\+$K$\+bi\-mod\-ule $A\ot_RU$ is $(A/R,K)$\+robustly
flaprojective.
 Indeed, one has $(A\ot_R\nobreak U)\ot_KH\simeq A\ot_R(U\ot_KH)$,
and $U\ot_KH$ is a flat left $R$\+module (cf.\
Remark~\ref{very-flaprojective-cotorsion-pair-generated-by}).
\end{ex}

\begin{lem} \label{robustly-flaprojective-restrict-exten-of-scalars}
  Let $R\rarrow S$ be a flat epimorphism of commutative rings and
$A$ be a quasi-algebra over~$R$.
 Let $A$ and $K$ be quasi-algebras over $R$, and $B=S\ot_RA$ and
$L=S\ot_RK$ be the related quasi-algebras over~$S$.
 In this context: \par
\textup{(a)} Let $G$ be a $(B/S,L)$\+robustly flaprojective
$B$\+$L$\+bimodule.
 Then $G$ is also $(A/R,K)$\+robustly flaprojective as
an $A$\+$K$\+bimodule. \par
\textup{(b)} Let $F$ be an $A$\+$K$\+bimodule whose underlying
$R$\+$R$\+bimodule is a quasi-module over~$R$.
 Assume that the $A$\+$K$\+bimodule $F$ is $(A/R,K)$\+robustly
flaprojective.
 Then the $B$\+$L$\+bimodule $S\ot_RF\simeq S\ot_RF\ot_RS\simeq
F\ot_RS$ is $(B/S,L)$\+robustly flaprojective.
\end{lem}

\begin{proof}
 Part~(a): given a flat left $K$\+module $H$, the left $L$\+module
$S\ot_RH\simeq L\ot_KH$ is flat by
Corollary~\ref{quasi-algebra-change-of-scalars-adjustedness}(b).
 Now the left $B$\+module $G\ot_KH\simeq G\ot_L(L\ot_KH)$ is
$B/S$\+flaprojective by assumption.
 By Lemma~\ref{flaprojective-restriction-extension-of-scalars}(a),
it follows that $G\ot_KH$ is also $A/R$\+flaprojective as
a left $A$\+module.

 Part~(b): any flat left $L$\+module $H$ is also flat as a left
$K$\+module by
Corollary~\ref{quasi-algebra-change-of-scalars-adjustedness}(a).
 Now we have $(S\ot_RF\ot_RS)\ot_LH\simeq(S\ot_RF\ot_KL)\ot_LH
\simeq S\ot_RF\ot_KH$.
 The left $A$\+module $F\ot_KH$ is $A/R$\+flaprojective by assumption.
 By Lemma~\ref{flaprojective-restriction-extension-of-scalars}(b),
it follows that $S\ot_RF\ot_KH$ is $B/S$\+flaprojective as
a left $B$\+module.
\end{proof}

\begin{lem} \label{quasi-algebras-robust-flaprojectivity-locality}
 Let $R\rarrow S_\alpha$ be a finite collection of homomorphisms of
commutative rings such that the collection of induced maps of
the spectra\/ $\Spec S_\alpha\rarrow\Spec R$ is an affine open covering
of an affine scheme.
 Let $A$ and $K$ be quasi-algebras over $R$, and let $F$ be
an $A$\+$K$\+bimodule whose underlying $R$\+$R$\+bimodule is
a quasi-module over~$R$.
 Let $B_\alpha=S_\alpha\ot_RA$ and $L_\alpha=S_\alpha\ot_RK$ be
the related rings.
 Then the $A$\+$K$\+bimodule $F$ is $(A/R,K)$\+robustly flaprojective
if and only if the $B_\alpha$\+$L_\alpha$\+bimodule $S_\alpha\ot_RF$
is $(B_\alpha/S_\alpha,L_\alpha)$\+robustly flaprojective for
all indices~$\alpha$. 
\end{lem}

\begin{proof}
 This is similar to
Lemmas~\ref{quasi-algebra-adjustedness-co-locality}(a)
and~\ref{quasi-algebra-very-flaprojectivity-locality}, and based on
Lemma~\ref{robustly-flaprojective-restrict-exten-of-scalars}(a\+b).
 The fact that the class of $(A/R,K)$\+robustly flaprojective modules
is closed under kernels of epimorphisms has to be used.
\end{proof}

 Let us define two more cotorsion pairs in the category of modules
over a commutative ring~$R$.
 Let us say that an $R$\+module $C$ is \emph{flepi cotorsion} if for
every flat epimorphism of commutative rings $R\rarrow S$ one has
$\Ext^1_R(S,C)=0$.
 We will say that an $R$\+module $F$ is \emph{flepi flat} if for
every flepi cotorsion $R$\+module $C$ one has $\Ext^1_R(F,C)=0$.
 Furthermore, let us say that an $R$\+module $C$ is \emph{strongly
flepi cotorsion} if for every flat epimorphism of commutative rings
$R\rarrow S$ one has $\Ext^n_R(S,C)=0$ for all $n\ge1$.
 We will say that an $R$\+module $F$ is \emph{weakly flepi flat} if
for every strongly flepi cotorsion $R$\+module $C$ one has
$\Ext^1_R(F,C)=0$ (or equivalently, for every strongly flepi cotorsion
$R$\+module $C$ one has $\Ext^n_R(F,C)=0$ for all $n\ge1$).
 Obviously, all cotorsion $R$\+modules are strongly flepi cotorsion,
hence all weakly flepi flat $R$\+modules are flat.

 For any fixed ring $R$, there exists, up to isomorphism, at most a set
of flat ring epimorphisms $R\rarrow S$ (as one can see
from~\cite[Theorem~XI.2.1.]{Ste}.
 So the theorem of Eklof and Trlifaj~\cite[Theorems~2 and~10]{ET},
\cite[Theorem~6.11 and Corollary~6.14]{GT} is applicable, and it tells
us that the pair of classes (flepi flat modules, flepi cotorsion
modules) is a complete cotorsion pair in $R\Modl$, while the pair of
classes (weakly flepi flat modules, strongly flepi cotorsion modules)
is a hereditary complete cotorsion pair in $R\Modl$.
 The flepi flat $R$\+modules are precisely all the direct summands
of transfinitely iterated extensions of the $R$\+modules $S$, where
$R\rarrow S$ are flat epimorphisms of commutative rings.
 The weakly flepi flat $R$\+modules are precisely all the direct
summands of transfinitely iterated extensions of the $R$\+modules $S$
as above and their \emph{syzygy modules} (i.~e., the $R$\+modules
of cocycles in projective resolutions of the $R$\+modules~$S$).

 Part~(b) of the next lemma explains why did not define ``robustly
very flaprojective bimodules''.
 Abusing the terminology, one can say that \emph{all very flaprojective
bimodules are robustly very flaprojective} (in the setting we are
interested in).

\begin{lem} \label{flepi-very-flat-modules-unproblematic-lemma}
 Let $R$ be a commutative ring, $A$ be a quasi-algebra over $R$, and
$F$ be an $A$\+$R$\+bimodule whose underlying $R$\+$R$\+bimodule is
a quasi-module over~$R$.
 In this context: \par
\textup{(a)} if the left $A$\+module $F$ is $A/R$\+flaprojective and
$H$ is a weakly flepi flat $R$\+module, then the left $A$\+module
$F\ot_RH$ is $A/R$\+flaprojective; \par
\textup{(b)} if the left $A$\+module $F$ is $A/R$\+very flaprojective
and $H$ is a very flat $R$\+module, then the left $A$\+module
$F\ot_RH$ is $A/R$\+very flaprojective.
\end{lem}

\begin{proof}
 In part~(a), one observes that the class of $A/R$\+flaprojective
left $A$\+modules is closed under transfinitely iterated extensions
(by Remark~\ref{very-flaprojective-cotorsion-pair-generated-by}),
kernels of epimorphisms in $A\Modl$ (by
Lemma~\ref{flaprojective-cotorsion-pair-hereditary}(b)),
and direct summands.
 Therefore, it suffices to consider the case of $H=S$, where
$R\rarrow S$ is a flat epimorphism of commutative rings.
 In this case, the left $A$\+module $F\ot_RS\simeq S\ot_RF$ is
$A/R$\+flaprojective by
Lemma~\ref{flaprojective-restriction-extension-of-scalars}(a\+-b).

 Part~(b) is similar but simpler.
 One observes that the class of $A/R$\+very flaprojective
left $A$\+modules is closed under transfinitely iterated extensions
(by Remark~\ref{very-flaprojective-cotorsion-pair-generated-by})
and direct summands.
 Therefore, it suffices to consider the case of $H=R[r^{-1}]$, where
$r\in R$ is an element.
 In this case, the left $A$\+module $F\ot_RR[r^{-1}]\simeq
R[r^{-1}]\ot_RF$ is $A/R$\+very flaprojective by
Lemma~\ref{very-flaprojective-restriction-extension-of-scalars}(a\+-b).
\end{proof}

\begin{rems}
 We \emph{do not know} whether the assertion of
Lemma~\ref{flepi-very-flat-modules-unproblematic-lemma}(a) holds for
all flat $R$\+modules~$H$ (in other words, whether every
left $A/R$\+flaprojective $(A,R)$\+bimodule whose underlying
$R$\+$R$\+bimodule is a quasi-module over $R$ needs to be
$(A/R,R)$\+robustly flaprojective).
 We also do not know how ubiquitous the weakly flepi flat $R$\+modules
are in general.

 The following two results are known, however:
\begin{itemize}
\item Over any Noetherian commutative ring $R$, all countably generated
flat modules are flepi flat (and even quite flat, which is a stronger
property) \cite[Theorem~2.4]{HPS}.
\item Over any Noetherian commutative ring $R$ with (at most)
countable spectrum, all flat modules are flepi flat (and even quite
flat) \cite[Theorem~1.17]{PSl2}, \cite[Theorems~2.4 and~3.8]{HPS}.
\end{itemize}

 So, for a quasi-algebra $A$ over a Noetherian commutative ring $R$
with countable spectrum, an $(A,R)$\+bi\-mod\-ule whose underlying
$R$\+$R$\+bimodule is a quasi-module over $R$ is $(A/R,R)$\+ro\-bustly
flaprojective whenever it is $A/R$\+flaprojective as a left $A$\+module.
\end{rems}

\begin{lem} \label{robustly-flaprojective-Hom-tensor-lemma}
 Let $R\rarrow A$ and $S\rarrow B$ be two homomorphisms of associative
rings, $F$ be an $A$\+$B$\+bimodule, $G$ be a left $B$\+module, and
$E$ be a left $A$\+module.
 In this context: \par
\textup{(a)} if the $A$\+$S$\+bimodule $F$ is $(A/R,S)$\+robustly
flaprojective and the left $B$\+module $G$ is $B/S$\+flaprojective,
then the left $A$\+module $F\ot_BG$ is a $A/R$\+flaprojective; \par
\textup{(b)} if the $A$\+$S$\+bimodule $F$ is $(A/R,S)$\+robustly
flaprojective and the left $R$\+module $E$ is cotorsion, then
the left $S$\+module\/ $\Hom_A(F,E)$ is cotorsion.
\end{lem}

\begin{proof}
 Part~(b): let $H$ be a flat left $S$\+module.
 Then $\Ext^1_A(F,E)=0$ (since the $A$\+module $F$ is
$A/R$\+flaprojective) and $\Tor_1^S(F,H)=0$ (sinse the $S$\+module $H$
is flat).
 By Lemma~\ref{Ext-homological-formula}, it follows that
$\Ext^1_S(H,\Hom_A(F,E))\simeq\Ext_A^1(F\ot_SH,\>E)$.
 It remains to point out that the left $A$\+module $F\ot_SH$ is
$A/R$\+flaprojective, so $\Ext_A^1(F\ot_SH,\>E)=0$.

 Part~(a): we need to show that $\Ext^1_A(F\ot_BG,\>E)=0$.
 Indeed, we have $\Ext^1_A(F,E)=0$ (as mentioned above) and
$\Tor_1^B(F,G)=0$ (since all $B/S$\+flaprojective left $B$\+modules
are flat as left $B$\+modules by
Remark~\ref{very-flaprojective-terminology-explained}).
 By Lemma~\ref{Ext-homological-formula}, it follows that
$\Ext_A^1(F\ot_BG,\>E)\simeq\Ext^1_B(G,\Hom_A(F,E))$.
 It remains to recall that the left $B$\+module $\Hom_A(F,E)$
is cotorsion as a left $S$\+module by part~(b).

 Alternatively, in part~(a) one can say that the left $B$\+module $G$
is a direct summand of a $B$\+module filtered by the $B$\+modules
$B\ot_SH$, where $H$ ranges over flat left $S$\+modules, by
Remark~\ref{very-flaprojective-cotorsion-pair-generated-by}.
 So it suffices to consider the case of $G=B\ot_SH$, and then
$F\ot_BG\simeq F\ot_SH$ is an $A/R$\+flaprojective left $A$\+module
by assumption.
\end{proof}

\begin{cor} \label{robustly-flaprojective-tensor-product-cor}
 Let $R\rarrow A$ and $S\rarrow B$ be two homomorphisms of associative
rings, $K$ be an associative ring, $F$ be a $A$\+$B$\+bimodule whose
underlying $A$\+$S$\+bimodule is $(A/R,S)$\+robustly flaprojective,
and $G$ be a $(B/S,K)$\+robustly flaprojective $B$\+$K$\+bimodule.
 Then the $A$\+$K$\+bimodule $F\ot_BG$ is $(A/R,K)$\+robustly
flaprojective.
\end{cor}

\begin{proof}
 Given a flat left $K$\+module $H$, we need to show that the left
$A$\+module $F\ot_BG\ot_KH$ is $A/R$\+flaprojective.
 Indeed, the left $B$\+module $G'=G\ot_KH$ is $B/S$\+flaprojective
by assumption, and it remains to apply
Lemma~\ref{robustly-flaprojective-Hom-tensor-lemma}(a) to
the $A$\+$B$\+bimodule $F$ and the left $B$\+module~$G'$.
\end{proof}

\begin{lem} \label{very-flaprojective-Hom-tensor-lemma}
 Let $R$ be a commutative ring, and let $A$ and $B$ be two
quasi-algebras over~$R$.
 Let $F$ be an $A$\+$B$\+bimodule whose underlying $R$\+$R$\+bimodule
is a quasi-module over~$R$.
 Let $G$ be a left $B$\+module and $E$ be a left $A$\+module.
 In this context: \par
\textup{(a)} if the left $A$\+module $F$ is $A/R$\+very flaprojective
and the left $B$\+module $G$ is $B/R$\+very flaprojective, then
the left $A$\+module $F\ot_BG$ is $A/R$\+very flaprojective; \par
\textup{(b)} if the left $A$\+module $F$ is $A/R$\+very flaprojective
and the $R$\+module $E$ is contraadjusted, then the $R$\+module\/
$\Hom_A(F,E)$ is contraadjusted.
\end{lem}

\begin{proof}
 Part~(b): let $H$ be a very flat left $R$\+module.
 Then $\Ext_A^1(F,E)=0$ and $\Tor^R_1(F,H)=0$.
 By Lemma~\ref{Ext-homological-formula}, it follows that
$\Ext^1_R(H,\Hom_A(F,E))\simeq\Ext^1_A(F\ot_RH,\>E)$.
 It remains to point out that the left $A$\+module $F\ot_RH$ is
$A/R$\+very flaprojective by
Lemma~\ref{flepi-very-flat-modules-unproblematic-lemma}(b), so
$\Ext^1_A(F\ot_RH,\>E)=0$.

 Part~(a): we need to show that $\Ext_A^1(F\ot_BG,\>E)=0$.
 Indeed, we have $\Ext_A^1(F,E)=0$ and $\Tor^B_1(F,G)=0$ (since
all $B/S$\+very flaprojective left $B$\+modules are flat as left
$B$\+modules by Remark~\ref{very-flaprojective-terminology-explained}).
 By Lemma~\ref{Ext-homological-formula}, it follows that
$\Ext_A^1(F\ot_BG,\>E)\simeq\Ext_B^1(G,\Hom_A(F,E))$.
 It remains to recall that the left $B$\+module $\Hom_A(F,E)$ is
contraadjusted as an $R$\+module by part~(b).

 Alternatively, in part~(a) one can say that the left $B$\+module $G$
is a direct summand of a $B$\+module filtered by the $B$\+modules
$B\ot_RH$, where $H$ ranges over very flat $R$\+modules, by
Remark~\ref{very-flaprojective-cotorsion-pair-generated-by}.
 So it suffices to consider the case of $G=B\ot_RH$, and then
$F\ot_BG\simeq F\ot_RH$ is an $A/R$\+very flaprojective left
$A$\+module by
Lemma~\ref{flepi-very-flat-modules-unproblematic-lemma}(b).
\end{proof}

\Section{Contraherent Cosheaves of $\cO$-Modules}
\label{contraherent-cosheaves-of-O-modules-secn}

 There are no new results is this section, which is mostly
an extraction from~\cite[Chapters~1\+-3]{Pcosh}.
 The preprint~\cite{Pphil} can be used as a supplementary reference.

\subsection{Cosheaves of $\cO$-modules}
\label{cosheaves-of-modules-subsecn}
 In this section we mostly follow~\cite[Section~2.1]{Pcosh}.

 Let $X$ be a topological space.
 Recall that a \emph{presheaf of abelian groups} $\M$ on $X$ is
a contravariant functor from the category of open subsets in~$X$
(with identity inclusions as morphisms) to the category of
abelian groups~$\Ab$.
 So, to every open subset $U\subset X$, a presheaf $\M$ assigns
an abelian group $\M(U)$ of \emph{sections of $\M$ over~$U$};
and to every pair of open subsets $V\subset U\subset X$, it assigns
the \emph{restriction map} $\M(U)\rarrow\M(V)$, which must be
an abelian group homomorphism.

 A presheaf of abelian groups $\M$ on $X$ is called a \emph{sheaf} if
it satisfies the following \emph{sheaf axiom}.
 For every open covering $U=\bigcup_\alpha V_\alpha$ of an open
subset $U\subset X$, the short sequence of abelian groups
\begin{equation} \label{sheaf-axiom}
 0\lrarrow\M(U)\lrarrow\prod\nolimits_\alpha\M(V_\alpha)\lrarrow
 \prod\nolimits_{\alpha,\beta}\M(V_\alpha\cap V_\beta)
\end{equation}
must be left exact.

 Dual-analogously, a \emph{copresheaf of abelian groups} $\P$ on $X$ is
a covariant functor from the category of open subsets in $X$ to~$\Ab$.
 To every open subset $U\subset X$, a copresheaf $\P$ assigns
an abelian group $\P[U]$ of \emph{cosections of\/ $\P$ over~$U$};
and to every pair of open subsets $V\subset U\subset X$, it assigns
the \emph{corestriction map} $\P[V]\rarrow\P[U]$, which must be
an abelian group homomorphism.

 A copresheaf of abelian groups $\P$ on $X$ is called a \emph{cosheaf}
if it satisfies the following \emph{cosheaf axiom}.
 For every open covering $U=\bigcup_\alpha U_\alpha$ of an open subset
$U\subset X$, the short sequence of abelian groups
\begin{equation} \label{cosheaf-axiom}
 0\llarrow\P[U]\llarrow\coprod\nolimits_\alpha\P[V_\alpha]\llarrow
 \coprod\nolimits_{\alpha,\beta}\P[V_\alpha\cap V_\beta]
\end{equation}
must be right exact.

 Now let $(X,\cO)$ be a ringed space.
 So $\cO$ is a sheaf of (associative and unital, not necessarily
commutative) rings on~$X$.
 In the rest of this Section~\ref{cosheaves-of-modules-subsecn}, all
\emph{modules} will be presumed to be left modules.

 A \emph{presheaf of $\cO$\+modules} $\M$ on $X$ is a presheaf of
abelian groups together with an additional structure of
an $\cO(U)$\+module on every group of sections $\M(U)$ such that
the restriction maps $\M(U)\rarrow\M(V)$ are $\cO(U)$\+module morphisms.
 Here the $\cO(U)$\+module structure on $\M(V)$ is defined using
the restriction map $\cO(U)\rarrow\cO(V)$ in the sheaf of rings~$\cO$.

 A presheaf of $\cO$\+modules on $X$ is called a \emph{sheaf of
$\cO$\+modules} if its underlying presheaf of abelian groups is
a sheaf of abelian groups.
 So the sequence~\eqref{sheaf-axiom} is a left exact sequence of
$\cO(U)$\+modules in this case.
 We will denote the category of sheaves of $\cO$\+modules on $X$
(with the obvious morphisms) by $(X,\cO)\Sh$.

 Dual-analogously, a \emph{copresheaf of $\cO$\+modules} $\P$ on $X$
is a copresheaf of abelian groups together with an additional structure
of an $\cO(U)$\+module on the group of cosections $\P[U]$ over every
open subset $U\subset X$ such that, for every pair of open subsets
$V\subset U\subset X$, the corestriction map $\P[V]\rarrow\P[U]$ is
an $\cO(U)$\+module morphism.
 Here the $\cO(U)$\+module structure on $\P[V]$ is defined in terms of
the $\cO(V)$\+module structure on $\P[V]$ and the ring homomorphism
$\cO(U)\rarrow\cO(V)$ of restriction of sections in the sheaf of
rings~$\cO$.

 A copresheaf of $\cO$\+modules on $X$ is called a \emph{cosheaf of
$\cO$\+modules} if its underlying copresheaf of abelian groups is
a cosheaf of abelian groups.
 So the sequence~\eqref{cosheaf-axiom} is a right exact sequence of
$\cO(U)$\+modules in this case.
 We will denote the category of cosheaves of $\cO$\+modules on $X$
(with the obvious morphisms) by $(X,\cO)\Cosh$.

 Now let $\bB$ be a base of open subsets in~$X$.
 The aim of the subsequent paragraphs is to generalize the discussion
above from the case of the base of \emph{all} open subsets in $X$ to
the case of arbitrary topology base~$\bB$.

 We will view $\bB$ as a category with identity inclusion maps as
morphisms.
 Then a \emph{presheaf of abelian groups} $\N$ on $\bB$ is simply defined
as a functor $\bB^\sop\rarrow\Ab$.
 Once again, to every open subset $U\in\bB$, a presheaf $\N$ assigns
an abelian group $\N(U)$ of \emph{sections of $\N$ over~$U$}; and to
every pair of open subsets $V\subset U$, \ $U$, $V\in\bB$, it assigns
the \emph{restriction map} $\N(U)\rarrow\N(V)$.

 A presheaf of abelian groups $\N$ on $\bB$ is called a \emph{sheaf} if
it satisfies the following version of the sheaf axiom.
 For every open covering $U=\bigcup_\alpha V_\alpha$ of an open subset
$U\in\bB$ by open subsets $V_\alpha\in\bB$, and every family of open
coverings $V_\alpha\cap V_\beta=\bigcup_\gamma W_{\alpha,\beta,\gamma}$
of all the intersections $V_\alpha\cap V_\beta$ by open subsets
$W_{\alpha,\beta,\gamma}\in\bB$, the short sequence of abelian groups
\begin{equation} \label{sheaf-axiom-topology-base}
 0\lrarrow\N(U)\lrarrow\prod\nolimits_\alpha\N(V_\alpha)\lrarrow
 \prod\nolimits_{\alpha,\beta,\gamma}\N(W_{\alpha,\beta,\gamma})
\end{equation}
must be left exact.

 Dual-analogously, a \emph{copresheaf of abelian groups} $\Q$ on $\bB$
is a functor $\bB\rarrow\Ab$.
 To every open subset $U\in\bB$, a copresheaf $\Q$ assigns an abelian
group $\Q[U]$ of \emph{cosections of\/ $\Q$ over~$U$}; and to every
pair of open subsets $V\subset U$, \ $U$, $V\in\bB$, it assigns
the \emph{corestriction map} $\Q[V]\rarrow\Q[U]$.

 A copresheaf of abelian groups $\Q$ on $\bB$ is called a \emph{cosheaf}
if it satisfies the following version of the cosheaf axiom.
 For every open covering $U=\bigcup_\alpha V_\alpha$ of an open subset
$U\in\bB$ by open subsets $V_\alpha\in\bB$, and every family of open
coverings $V_\alpha\cap V_\beta=\bigcup_\gamma W_{\alpha,\beta,\gamma}$
of all the intersections $V_\alpha\cap V_\beta$ by open subsets
$W_{\alpha,\beta,\gamma}\in\bB$, the short sequence of abelian groups
\begin{equation} \label{cosheaf-axiom-topology-base}
 0\llarrow\Q[U]\llarrow\coprod\nolimits_\alpha\Q[V_\alpha]\llarrow
 \coprod\nolimits_{\alpha,\beta}\Q[W_{\alpha,\beta,\gamma}]
\end{equation}
must be right exact.

 As above, let $(X,\cO)$ be a ringed space.
 A \emph{presheaf of $\cO$\+modules} $\N$ on $\bB$ is a presheaf of
abelian groups on $\bB$ such that the abelian group $\N(U)$ is endowed
with an additional structure of an $\cO(U)$\+module for every $U\in\bB$
and the restriction map $\N(U)\rarrow\N(V)$ is an $\cO(U)$\+module
morphism for all $V\subset U$, \ $U$, $V\in\bB$.
 A presheaf of $\cO$\+modules on $\bB$ is called a \emph{sheaf of
$\cO$\+modules} on $\bB$ if its underlying presheaf of abelian groups
is a sheaf of abelian groups on~$\bB$.
 We will denote the category of sheaves of $\cO$\+modules on $\bB$
(with the obvious morphisms) by $(\bB,\cO)\Sh$.

 Dual-analogously, a \emph{copresheaf of $\cO$\+modules} $\Q$ on $\bB$
is a copresheaf of abelian groups on $\bB$ such that the abelian group
$\Q[U]$ is endowed with an $\cO(U)$\+module structure for every
$U\in\bB$ and the corestriction map $\Q[V]\rarrow\Q[U]$ is
an $\cO(U)$\+module morphism for all $V\subset U$, \ $U$, $V\in\bB$.
 Here, as above, the $\cO(U)$\+module structure on $\Q[V]$ is defined
in terms of the $\cO(V)$\+module structure on $\Q[V]$ and
the restriction map $\cO(U)\rarrow\cO(V)$.
 A copresheaf of $\cO$\+modules on $\bB$ is called a \emph{cosheaf of
$\cO$\+modules} on $\bB$ if its underlying copresheaf of abelian groups
is a cosheaf of abelian groups on~$\bB$.
 We will denote the category of cosheaves of $\cO$\+modules on $\bB$
(with the obvious morphisms) by $(\bB,\cO)\Cosh$.

\begin{thm} \label{extension-of-co-sheaves-from-topology-base}
 Let $(X,\cO)$ be a ringed space and\/ $\bB$ be a topology base in~$X$.
 In this setting: \par
\textup{(a)} Any sheaf of $\cO$\+modules on\/ $\bB$ can be extended to
a sheaf of $\cO$\+modules on $X$ (i.~e., on all open subsets of~$X$)
in a unique way.
 The functor of restriction from all open subsets of $X$ to the open
subsets belonging to\/ $\bB$ is an equivalence of categories
$$
 (X,\cO)\Sh\simeq(\bB,\cO)\Sh.
$$ \par
\textup{(b)} Any cosheaf of $\cO$\+modules on\/ $\bB$ can be extended to
a cosheaf of $\cO$\+modules on $X$ (i.~e., on all open subsets of~$X$)
in a unique way.
 The functor of restriction from all open subsets of $X$ to the open
subsets belonging to\/ $\bB$ is an equivalence of categories
$$
 (X,\cO)\Cosh\simeq(\bB,\cO)\Cosh.
$$
\end{thm}

\begin{proof}
 Part~(a) is essentially~\cite[Section~0.3.2]{EGA1}; see
also~\cite[Lemma Tag~009U]{SP}.
 Part~(b) is~\cite[Theorem~2.1.2]{Pcosh}.
\end{proof}

 Given an open subset $Y\subset X$, denote by $\cO|_Y$ the restriction
of the sheaf of rings $\cO$ on $X$ to the open subset~$Y$.
 Then, for any sheaf of $\cO$\+modules $\M$ on $X$, the sheaf of
$\cO|_Y$\+modules $\M|_Y$ on $Y$ is defined by the rule $\M|_Y(U)=\M(U)$
for all open subsets $U\subset Y$.
 For any cosheaf of $\cO$\+modules $\P$ on $X$, the cosheaf of
$\cO|_Y$\+modules $\P|_Y$ on $Y$ is defined by the rule $\P|_Y[U]=\P[U]$
for all open subsets $U\subset Y$.

\subsection{Locally contraherent cosheaves}
\label{locally-contraherent-cosheaves-subsecn}
 This section is largely based on~\cite[Sections~2.2 and~3.1]{Pcosh}.

 Let $X$ be a scheme with the structure sheaf $\cO_X$ and a fixed open
covering~$\bW$.
 We will say that an open subscheme $U\subset X$ is \emph{subordinate
to\/~$\bW$} if there exists an open subscheme $W\in\bW$ such that
$U\subset W$.
 An open covering $X=\bigcup_\alpha U_\alpha$ is said to be
\emph{subordinate to\/~$\bW$} if every open subscheme $U_\alpha$
belonging to this covering is subordinate to~$\bW$.

 Denote by $\bB$ the base of open subsets in $X$ consisting of all
the affine open subschemes $U\subset X$ subordinate to~$\bW$.

 Let $\M$ be a presheaf of $\cO_X$\+modules on~$\bB$.
 One says that the presheaf $\M$ is \emph{quasi-coherent} if
the following \emph{quasi-coherence axiom} holds:
\begin{itemize}
\item For every pair of affine open subschemes $V\subset U\subset X$
subordinate to $\bW$, the $\cO_X(V)$\+module map
$$
 \cO_X(V)\ot_{\cO_X(U)}\M(U)\lrarrow\M(V)
$$
corresponding by adjunction to the restriction map $\M(U)\rarrow\M(V)$
is an isomorphism.
\end{itemize}

\begin{lem} \label{quasi-coherence-implies-sheaf}
 The quasi-coherence axiom implies the sheaf axiom: Any quasi-coherent
presheaf on\/ $\bB$ is a sheaf of $\cO_X$\+modules on\/~$\bB$.
\end{lem}

\begin{proof}
 It suffices to consider the case of an affine scheme $X=U$ with
the open covering $\bW=\{U\}$.
 Then the point is that the rule $V\longmapsto
\cO_X(V)\ot_{\cO_X(U)}\M(U)$ defines the restriction to $\bB$ of
the quasi-coherent sheaf on $\Spec U$ corresponding to
the $\cO_X(U)$\+module~$\M(U)$.
 Alternatively, one can observe that, in the context of the definition
of a sheaf of $\cO_X$\+modules on $\bB$, it suffices to check the sheaf
axiom for one's favorite choice of the coverings $V_\alpha\cap V_\beta=
\bigcup_\gamma W_{\alpha,\beta,\gamma}$ of the intersections
$V_\alpha\cap V_\beta$, rather than for all such coverings (this is
a part of the formulation of~\cite[Proposition~2.1.3]{Pcosh}).
 In the situation at hand, we have $V_\alpha\cap V_\beta\in\bB$, since
the intersections of affine open subschemes in an affine scheme $U$
are affine.
 So let $V_\alpha\cap V_\beta$ be its own covering.
 Furthermore, as affine schemes are quasi-compact, it suffices to
consider \emph{finite} affine open coverings of $U$
only~\cite[Remark~2.1.4]{Pcosh}.
 This reduces the question to an exercise in commutative algebra, which
we leave to the reader (cf.~\cite[Lemma~1.2.6]{Pcosh}).
\end{proof}

 We denote the full subcategory of quasi-coherent (pre)sheaves on $\bB$
by $\bB\Qcoh\subset(\bB,\cO_X)\Sh$.
 The full subcategory of quasi-coherent sheaves in $(X,\cO_X)\Sh$ is
denoted by $X\Qcoh\subset(X,\cO_X)\Sh$.

\begin{cor} \label{enochs-estrada}
 Let $X$ be a scheme with an open covering\/ $\bW$ and\/ $\bB$ be
the topology base of $X$ formed by all the affine open subschemes
subordinate to\/~$\bW$.
 Then the category $X\Qcoh$ of quasi-coherent sheaves on $X$ is
equivalent to the category\/ $\bB\Qcoh$ of quasi-coherent presheaves
on\/~$\bB$.
\end{cor}

\begin{proof}
 This result is essentially due to Enochs and
Estrada~\cite[Section~2]{EE}.
 By Lemma~\ref{quasi-coherence-implies-sheaf}, any quasi-coherent
presheaf $\M$ on $\bB$ is a sheaf on~$\bB$; so
Theorem~\ref{extension-of-co-sheaves-from-topology-base}(a) tells us
that $\M$ can be uniquely extended to a sheaf of $\cO_X$\+modules
on~$X$.
 It remains to explain that the class of sheaves of $\cO_X$\+modules
one obtains in this way coincides with the class of all quasi-coherent
sheaves on~$X$.
 In other words, a sheaf of $\cO_X$\+modules on $X$ is quasi-coherent
if and only if it satisfies the quasi-coherence axiom above for all
pairs of affine open subschemes $V\subset U\subset X$ subordinate to
any fixed open covering $\bW$ of the scheme~$X$.
 This is clear from the standard local definition of a quasi-coherent
sheaf on a ringed space~\cite[Section Tag~01BD]{SP} together with
the classical description of quasi-coherent sheaves on affine schemes
in terms of modules over the ring of global
functions~\cite[Section Tag~01I6]{SP}.
 As a bonus, these arguments prove that the quasi-coherence condition is
local in the sense that the category $\bB\Qcoh$ of quasi-coherent
presheaves on $\bB$ does not depend on the choice of an open covering
$\bW$ of a scheme~$X$.
\end{proof}

 Let $\P$ be a copresheaf of $\cO_X$\+modules on~$\bB$.
 One says that the copresheaf $\P$ is \emph{contraherent} (on~$\bB$) if
the following two axioms hold for every pair of affine open subschemes
$V\subset U\subset X$ subordinate to $\bW$:
\begin{enumerate}
\renewcommand{\theenumi}{\roman{enumi}}
\item The $\cO_X(V)$\+module map
$$
 \P[V]\lrarrow\Hom_{\cO_X(U)}(\cO_X(V),\P[U])
$$
corresponding by adjunction to the corestriction map $\P[V]\rarrow\P[U]$
is an isomorphism.
\item One has
$$
 \Ext^1_{\cO_X(U)}(\cO_X(V),\P[U])=0.
$$
\end{enumerate}
 The condition~(i) is called the \emph{contraherence axiom}.
 The condition~(ii) is called the \emph{contraadjustedness axiom}.

 Concerning condition~(ii), notice first of all that the projective
dimension of the $\cO_X(U)$\+module $\cO_X(V)$ never exceeds~$1$
\,\cite[Remark~5.1]{Pphil}.
 Moreover, the $\cO_X(U)$\+mod\-ule $\cO_X(V)$ is \emph{very flat}
(in the sense of Section~\ref{prelim-very-flat-subsecn} above);
see Example~\ref{very-flat-ring-examples}(1).
 On the other hand, all principal affine open subschemes in $U$ occur
as the schemes $V$ in the context of condition~(ii).
 So, setting $R=\cO_X(U)$, for any element $s\in R$ one has
$R[s^{-1}]=\cO_X(V)$ for a certain affine open subscheme $V\subset U$.

 It follows that, as $U$ is fixed and $V$ varies, condition~(ii) is
equivalent to the $\cO_X(U)$\+module $\P[U]$ being
\emph{contraadjusted} in the sense of
Section~\ref{prelim-very-flat-subsecn}.
 Hence the terminology ``contraadjustedness axiom''.

\begin{lem} \label{contraherence+contraadjustedness-imply-cosheaf}
 The contraherence and contraadjustedness axioms imply the cosheaf
axiom: Any contraherent copresheaf on\/ $\bB$ is a cosheaf of
$\cO_X$\+modules on\/~$\bB$.
\end{lem}

\begin{proof}
 One observes that, in the context of the definition of a cosheaf of
$\cO_X$\+modules on $\bB$, it suffices to check the cosheaf axiom
for one's favorite choice of the coverings $V_\alpha\cap V_\beta=
\bigcup_\gamma W_{\alpha,\beta,\gamma}$ of the intersections
$V_\alpha\cap V_\beta$, rather than for all such coverings.
 This is a part of the formulation of~\cite[Theorem~2.1.2]{Pcosh}.
 As mentioned in the proof of Lemma~\ref{quasi-coherence-implies-sheaf},
in the situation at hand we have $V_\alpha\cap V_\beta\in\bB$; so let
$V_\alpha\cap V_\beta$ be its own covering.
 Furthermore, it suffices to consider \emph{finite} affine open
coverings $U=\bigcup_{\alpha=1}^N V_\alpha$
only~\cite[Remark~2.1.4]{Pcosh}.
 This reduces the question to a commutative algebra
lemma~\cite[Lemma~1.2.6(b)]{Pcosh}.
\end{proof}

 A cosheaf of $\cO_X$\+modules $\P$ on $X$ is said to be
\emph{$\bW$\+locally contraherent} if the restriction of $\P$ to $\bB$
is a contraherent co(pre)sheaf on~$\bB$.
 In other words, a cosheaf $\P$ is called a $\bW$\+locally contraherent
cosheaf on $X$ if the contraherence axiom~(i) and the contraadjustedness
axiom~(ii) hold for $\P$ for all pairs of affine open subschemes
$V\subset U\subset X$ subordinate to~$\bW$.
 We denote the full subcategory of $\bW$\+locally contraherent cosheaves
by $X\Lcth_\bW\subset(X,\cO_X)\Cosh$.
 The full subcategory of contraherent co(pre)sheaves on $\bB$ is
denoted by $\bB\Ctrh\subset(\bB,\cO_X)\Cosh$.

\begin{cor} \label{loc-contraherent-cosheaves-recovered-from-affines}
 Let $X$ be a scheme with an open covering\/ $\bW$ and\/ $\bB$ be
the topology base of $X$ formed by all the affine open subschemes
subordinate to\/~$\bW$.
 Then the category $X\Lcth_\bW$ of\/ $\bW$\+locally contraherent
cosheaves on $X$ is equivalent to the category\/ $\bB\Ctrh$ of
contraherent copresheaves on\/~$\bB$.
\end{cor}

\begin{proof}
 This is~\cite[Theorem~2.2.1]{Pcosh} (in the case of the trivial open
covering $\bW=\{X\}$) or~\cite[Theorem~3.1.1]{Pcosh} (in the general
case of an arbitrary open covering~$\bW$).
 The assertion follows from
Lemma~\ref{contraherence+contraadjustedness-imply-cosheaf}
and Theorem~\ref{extension-of-co-sheaves-from-topology-base}(b).
\end{proof}

 A cosheaf of $\cO_X$\+modules $\P$ on $X$ is said to be
(\emph{globally}) \emph{contraherent} if it is locally contraherent
for the trivial covering $\{X\}$ of the scheme~$X$.
 In other words, a cosheaf $\P$ is called contraherent if
the contraherence axiom~(i) and the contraadjustedness axiom~(ii)
hold for $\P$ for all pairs of affine open subschemes $V\subset U
\subset X$.
 We denote the full subcategory of contraherent cosheaves on $X$ by
$X\Ctrh\subset(X,\cO_X)\Cosh$.

 A cosheaf of $\cO_X$\+modules $\P$ on $X$ is said to be \emph{locally
contraherent} if there \emph{exists} an open covering $\bW$ such
that $\P$ is $\bW$\+locally contraherent.
 The full subcategory of locally contraherent cosheaves is denoted
by $X\Lcth=\bigcup_\bW X\Lcth_\bW\subset(X,\cO_X)\Cosh$.
 By the definition, we have inclusions of full subcategories
$$
 X\Ctrh\subset X\Lcth_\bW\subset X\Lcth\subset(X,\cO_X)\Cosh
$$
for any open covering $\bW$ of a scheme~$X$.

 Unlike the quasi-coherence property of a sheaf of $\cO_X$\+modules,
the contraherence property of a cosheaf of $\cO_X$\+modules is
\emph{not} local: a $\bW$\+locally contraherent cosheaf on $X$
\emph{need not} be (globally) contraherent.
 In fact, the contraadjustedness axiom~(ii) is local, but
the contraherence axiom~(i) is not.
 A counterexample can be found in~\cite[Example~3.2.1]{Pcosh}.

 The $\bW$\+locally contraherent cosheaves in the sense of
the definition above are sometimes called \emph{locally contraadjusted},
to distinguish them from the following two narrower classes of
$\bW$\+locally contraherent cosheaves.
 A $\bW$\+locally contraherent cosheaf $\P$ on $X$ is said to be
\emph{locally cotorsion} if the $\cO_X(U)$\+module $\P[U]$ is
cotorsion (in the sense of Section~\ref{prelim-cotorsion-subsecn})
for every affine open subscheme $U\subset X$ subordinate to~$\bW$.
 A $\bW$\+locally contraherent cosheaf $\gJ$ on $X$ is said to be
\emph{locally injective} if the $\cO_X(U)$\+module $\gJ[U]$ is
injective for every affine open subscheme $U\subset X$ subordinate
to~$\bW$.
 (Similarly one can speak about \emph{locally cotorsion} and
\emph{locally injective} contraherent copresheaves of
$\cO_X$\+modules on~$\bB$.)

 The local cotorsion and local injectivity properties of locally
contraherent cosheaves \emph{are} local (just as the terminology
suggests).
 This means that, given an affine open covering $X=\bigcup_\alpha
U_\alpha$ of the scheme $X$ subordinate to $\bW$, a $\bW$\+locally
contraherent cosheaf $\P$ on $X$ is locally cotorsion if and only if
the $\cO_X(U_\alpha)$\+module $\P[U_\alpha]$ is cotorsion for
every~$\alpha$ (see~\cite[Lemma~1.3.6(a) and Section~3.1]{Pcosh}).
 Similarly, a $\bW$\+locally contraherent cosheaf $\gJ$ on $X$ is
locally injective if and only if the $\cO_X(U_\alpha)$\+module
$\gJ[U_\alpha])$ is injective for every~$\alpha$
(see~\cite[Lemma~1.3.6(b) and Section~3.1]{Pcosh}).

 We denote the full subcategories of locally cotorsion and locally
injective $\bW$\+locally contraherent cosheaves by
$$
 X\Lcth_\bW^\lin\subset X\Lcth_\bW^\lct\subset X\Lcth_\bW.
$$
 Specializing to the case of the trivial open covering $\bW$, we
obtain the full subcategories of locally cotorsion and locally injective
(globally) contraherent cosheaves
$$
 X\Ctrh^\lin\subset X\Ctrh^\lct\subset X\Ctrh.
$$
 Finally, we put $X\Lcth^\lin=\bigcup_\bW X\Lcth_\bW^\lin$
and $X\Lcth^\lct=\bigcup_\bW X\Lcth_\bW^\lct$; so
$$
 X\Lcth^\lin\subset X\Lcth^\lct\subset X\Lcth.
$$
 It follows from the previous paragraph that one has
$$
 X\Lcth_\bW^\lct= X\Lcth^\lct\cap X\Lcth_\bW\subset X\Lcth
$$
and 
$$
 X\Lcth_\bW^\lin= X\Lcth^\lin\cap X\Lcth_\bW\subset X\Lcth.
$$

\subsection{Exact categories of locally contraherent cosheaves}
\label{exact-categories-of-contrah-subsecn}
 In this section we mostly follow~\cite[Sections~3.1\+-3.2]{Pcosh}.

 Let $X$ be a scheme with an open covering~$\bW$.
 A short sequence of $\bW$\+locally contraherent cosheaves
$0\rarrow\P\rarrow\Q\rarrow\R\rarrow0$ on $X$ is said to be
\emph{exact} if, for every affine open subscheme $U\subset X$
subordinate to $\bW$, the short sequence of $\cO_X(U)$\+modules
$0\rarrow\P[U]\rarrow\Q[U]\rarrow\R[U]\rarrow0$ is exact.
 This rule defines an exact category structure (in the sense of
Quillen~\cite{Bueh}) on the additive category $X\Lcth_\bW$ of
$\bW$\+locally contraherent cosheaves on~$X$
\,\cite[Section~3.1]{Pcosh}.

 Specializing to the case of the trivial open covering $\{X\}$ of
the scheme $X$, we obtain the natural exact category structure on
the additive category $X\Ctrh$ of (globally) contraherent cosheaves
on~$X$.
 Passing to the inductive limit over refinements of coverings $\bW$,
we obtain a natural exact category structure on the additive category
$X\Lcth$ of locally contraherent cosheaves on~$X$.
 The latter sentence means that a short sequence in $X\Lcth$ is
called \emph{exact} if it is exact in $X\Lcth_\bW$ for \emph{some}
open covering~$\bW$.
 The full subcategories $X\Ctrh$ and $X\Lcth_\bW$ are closed under
extensions in the exact category $X\Lcth$
\,\cite[Corollary~3.2.4]{Pcosh}.

 The property of a short sequence of $\bW$\+locally contraherent
cosheaves on $X$ to be exact is local.
 This means that, given an affine open covering $X=\bigcup_\alpha
U_\alpha$ of the scheme $X$ subordinate to $\bW$, a short sequence
of $\bW$\+locally contraherent cosheaves $0\rarrow\P\rarrow\Q\rarrow\R
\rarrow0$ on $X$ is exact if and only if the short sequence of
$\cO_X(U_\alpha)$\+modules $0\rarrow\P[U_\alpha]\rarrow\Q[U_\alpha]
\rarrow\R[U_\alpha]\rarrow0$ is exact for every~$\alpha$
\,\cite[Lemma~1.4.1(a) or~3.1.2(a)]{Pcosh}.
 It follows that a short sequence in $X\Lcth_\bW$ is exact if and only
if it is exact in $X\Lcth$ \,\cite[Section~3.2]{Pcosh}.

 A morphism in $X\Ctrh$ or in $X\Lcth_\bW$ is an admissible epimorphism
if and only if it is an admissible epimorphism in $X\Lcth$
\,\cite[Lemma~1.4.1(b)]{Pcosh}.
 In other words, the full subcategories $X\Ctrh$ and $X\Lcth_\bW$ are 
closed under kernels of admissible epimorphisms in $X\Lcth$
\,\cite[Section~3.2]{Pcosh}.
 For admissible monomorphisms, the similar assertions are \emph{not}
true~\cite[Example~3.2.1]{Pcosh}.

 A short sequence of locally cotorsion (resp., locally injective)
$\bW$\+locally contraherent cosheaves on $X$ is said to be
\emph{exact} if it is exact in $X\Lcth_\bW$.
 This rule defines exact category structures on the additive categories
$X\Lcth_\bW^\lct$ and $X\Lcth_\bW^\lin$ of locally cotorsion and
locally injective $\bW$\+locally contraherent cosheaves on~$X$.
 The full subcategories $X\Lcth_\bW^\lct$ and $X\Lcth_\bW^\lin$ are
closed under extensions and cokernels of admissible monomorphisms
in $X\Lcth_\bW$ (by Lemma~\ref{flat-cotorsion-pair-hereditary}(b)),
but not under the kernels of admissible epimorphisms.

 Specializing to the case of the trivial open covering $\{X\}$ of
the scheme $X$, we obtain natural exact category structures on
the additive categories $X\Ctrh^\lct$ and $X\Ctrh^\lin$ of locally
cotorsion and locally injective (globally) contraherent cosheaves
on~$X$.
 Passing to the inductive limits over refinements of coverings $\bW$,
we obtain natural exact category structures on the additive categories
$X\Lcth^\lct$ and $X\Lcth^\lin$ of locally cotorsion and locally
injective locally contraherent cosheaves on~$X$.

 The exact categories $X\Ctrh$ and $X\Lcth_\bW$ have exact functors
of infinite products.
 The functor of cosections ${-}[U]\:X\Lcth_\bW\rarrow\cO_X(U)\Modl$
preserves infinite products for every quasi-compact quasi-separated open
subscheme $U\subset X$ \,\cite[Remark~2.1.4 and Section~3.1]{Pcosh}.
 The full subcategories $X\Lcth_\bW^\lct$ and $X\Lcth_\bW^\lin$ are
closed under infinite products in $X\Lcth_\bW$ (in particular,
the full subcategories $X\Ctrh^\lct$ and $X\Ctrh^\lin$ are closed
under infinite products in $X\Ctrh$).

 For any affine scheme $U=\Spec R$, the functor of global cosections
${-}[U]$ provides an equivalence between the exact category $U\Ctrh$ of
contraherent cosheaves on $U$ and the exact category $R\Modl^\cta$ of
contraadjusted $R$\+modules,
$$
 U\Ctrh\simeq R\Modl^\cta.
$$
 Here the exact category structure on the category of contraadjusted
$R$\+modules $R\Modl^\cta$ is inherited from the abelian exact
structure on the ambient abelian category $R\Modl$.
 Lemma~\ref{restriction-co-extension-vfl-cta}(d) together with
Example~\ref{very-flat-ring-examples}(1) is relevant to this category
equivalence.
 See~\cite[Section~5.6]{Pphil} for a discussion.

 Similarly, the same functor of global cosections ${-}[U]$ provides
an equivalence between the exact category $U\Ctrh^\lct$ of locally
cotorsion contraherent cosheaves on $U$ and the exact category
$R\Modl^\cot$ of cotorsion $R$\+modules,
$$
 U\Ctrh^\lct\simeq R\Modl^\cot.
$$
 Here the exact category structure on the category of cotorsion
$R$\+modules $R\Modl^\cot$ is inherited from the abelian exact structure
on the ambient abelian category $R\Modl$.
 Lemma~\ref{restriction-coextension-injective-cotorsion}(c)
is relevant to this category equivalence.

 Finally, on an affine scheme $U=\Spec R$, the exact structure on
the category of locally injective contraherent cosheaves $U\Ctrh^\lin$
is split (i.~e., all short exact sequences are split).
 The very same functor of global cosections ${-}[U]$ provides
an equivalence between the additive category $U\Ctrh^\lin$ and
the additive category $R\Modl^\inj$ of injective $R$\+modules,
$$
 U\Ctrh^\lin\simeq R\Modl^\inj.
$$
 Lemma~\ref{restriction-coextension-injective-cotorsion}(d)
is relevant to this category equivalence.

\subsection{Direct images} \label{direct-images-of-O-co-sheaves-subsecn}
 In this section and the next one we mostly follow~\cite[Sections~2.3
and~3.3]{Pcosh}.

 Let $f\:(Y,\cO_Y)\rarrow(X,\cO_X)$ be a morphism of ringed spaces.
 Given a copresheaf of $\cO_Y$\+modules $\Q$ on $Y$, the \emph{direct
image} $f_!\Q$ of the copresheaf $\Q$ is a copresheaf of
$\cO_X$\+modules on $X$ defined by the rule
$$
 (f_!\Q)[U]=\Q[f^{-1}(U)]
$$
for every open subset $U\subset X$.
 Here $f^{-1}(U)\subset Y$ is the preimage of $U$ in~$Y$ (which is
an open subset in~$Y$).
 As a part of the datum of a morphism of ringed spaces~$f$, we have
a ring homomorphism $\cO_X(U)\rarrow\cO_Y(f^{-1}(U))$, which allows
to define an $\cO_X(U)$\+module structure on
the $\cO_Y(f^{-1}(U))$\+module $\Q[f^{-1}(U)]$.
 The corestriction maps in the copresheaf $f_!\Q$ are defined as
the respective corestriction maps in the copresheaf~$\Q$.

 For any cosheaf of $\cO_Y$\+modules $\Q$ on $Y$, the copresheaf
of $\cO_X$\+modules $f_!\Q$ on $X$ is a cosheaf (since the preimage
of any open covering of $U$ is an open covering of $f^{-1}(U)$).
 We have constructed an additive functor of direct image
$f_!\:(Y,\cO_Y)\Cosh\rarrow(X,\cO_X)\Cosh$.

 Let $j\:(Y,\cO_Y)\rarrow(X,\cO_X)$ be an identity open inclusion of
ringed spaces.
 So $Y$ is an open subset in $X$ and the sheaf of rings $\cO_Y=\cO_X|_Y$
is the restriction of the sheaf of rings $\cO_X$ to $Y\subset X$.
 Then the functor $j_!\:(Y,\cO_Y)\Cosh\rarrow(X,\cO_X)\Cosh$ is
left adjoint to the restriction functor $\P\longmapsto\P|_Y\:
(X,\cO_X)\Cosh\rarrow(Y,\cO_Y)\Cosh$ defined at the end of
Section~\ref{cosheaves-of-modules-subsecn}.
 In other words, for any cosheaves $\P\in(X,\cO_X)\Cosh$ and
$\Q\in(Y,\cO_Y)\Cosh$ there is a natural adjunction isomorphism
of abelian groups~\cite[formula~(2.7) in Section~2.3]{Pcosh}
\begin{equation} \label{open-immers-direct-inverse-cosheaf-adjunction}
 \Hom^{\cO_X}(j_!\Q,\P)\simeq\Hom^{\cO_Y}(\Q,\P|_Y)
\end{equation}
 Here we denote by $\Hom^{\cO_X}$ and $\Hom^{\cO_Y}$ the groups of
morphisms in the categories of cosheaves of modules $(X,\cO_X)\Cosh$
and $(Y,\cO_Y)\Cosh$.

 Now let $f\:Y\rarrow X$ be a morphism of schemes.
 The morphism~$f$ is called \emph{affine} if, for every affine open
subscheme $U\subset X$, the preimage $f^{-1}(U)\subset Y$ is an affine
open subscheme in~$Y$.
 It suffices to check this condition for affine open subschemes $U$
belonging to any chosen affine open covering $X=\bigcup_\alpha U_\alpha$
of the scheme~$X$ \,\cite[D\'efinition~II.1.2.1 and
Corollaire~II.1.3.2]{EGA2}, \cite[Lemma Tag~01S8]{SP}.

 For any affine morphism of schemes $f\:Y\rarrow X$, the direct image
functor $f_!\:(Y,\cO_Y)\Cosh\rarrow(X,\cO_X)\Cosh$ takes contraherent
cosheaves on $Y$ to contraherent cosheaves on~$X$.
 In particular, given a contraherent cosheaf $\Q$ on $Y$,
the $\cO_X(U)$\+module $(f_!\Q)[U]=\Q[f^{-1}(V)]$ is contraadjusted
for all affine open subschemes $U\subset X$ by
Lemma~\ref{restriction-co-extension-vfl-cta}(a).
 So we have the direct image functor between the categories of
contraherent cosheaves~\cite[Section~2.3]{Pcosh}
$$
 f_!\:Y\Ctrh\lrarrow X\Ctrh.
$$
 This assertion does \emph{not} hold without the assumption that
the morphism~$f$ is affine (\emph{not} even for open immersions~$f$
of smooth algebraic varieties over a field~\cite[Remark~2.3.1
and Example~2.3.2]{Pcosh}).

 Furthermore, for any affine morphism of schemes $f\:Y\rarrow X$,
the direct image functor $f_!\:(Y,\cO_Y)\Cosh\rarrow(X,\cO_X)\Cosh$
takes locally cotorsion contraherent cosheaves on $Y$ to locally
cotorsion contraherent cosheaves on~$X$.
 In particular, given a locally cotorsion contraherent cosheaf $\Q$
on $Y$, the $\cO_X(U)$\+module $(f_!\Q)[U]=\Q[f^{-1}(V)]$ is cotorsion
for all affine open subschemes $U\subset X$ by
Lemma~\ref{restriction-coextension-injective-cotorsion}(a).
 So we have the direct image functor between the categories of
locally cotorsion contraherent cosheaves~\cite[Section~2.3]{Pcosh}
$$
 f_!\:Y\Ctrh^\lct\lrarrow X\Ctrh^\lct.
$$

 See the beginning of the next
Section~\ref{inverse-images-of-O-co-sheaves-subsecn}
for a brief discussion of \emph{flat} morphisms of schemes.
 For any flat affine morphism of schemes $f\:Y\rarrow X$, the direct
image functor $f_!\:(Y,\cO_Y)\Cosh\rarrow(X,\cO_X)\Cosh$ takes locally
injective contraherent cosheaves on $Y$ to locally injective
contraherent cosheaves on~$X$.
 In particular, given a locally injective contraherent cosheaf $\gI$
on $Y$, the $\cO_X(U)$\+module $(f_!\gI)[U]=\gI[f^{-1}(V)]$ is injective
for all affine open subschemes $U\subset X$ by
Lemma~\ref{restriction-coextension-injective-cotorsion}(b).
 So we have the direct image functor between the categories of
locally injective contraherent cosheaves~\cite[Section~2.3]{Pcosh}
$$
 f_!\:Y\Ctrh^\lin\lrarrow X\Ctrh^\lin.
$$

 More generally, let $\bW$ be an open covering of the scheme $X$ and
$\bT$ be an open covering of the scheme~$Y$.
 A morphism of schemes $f\:Y\rarrow X$ is said to be
\emph{($\bW,\bT)$\+affine} if, for every affine open subscheme
$U\subset X$ subordinate to $\bW$, the preimage $f^{-1}(U)\subset Y$
an affine open subscheme subordinate to~$\bT$.
 So all $(\bW,\bT)$\+affine morphisms of schemes are affine.

 For any $(\bW,\bT)$\+affine morphism of schemes $f\:Y\rarrow X$,
the direct image functor $f_!\:(Y,\cO_Y)\Cosh\rarrow(X,\cO_X)\Cosh$
takes $\bT$\+locally contraherent cosheaves on $Y$ to $\bW$\+locally
contraherent cosheaves on~$X$.
 The resulting functor between the exact categories of locally
contraherent cosheaves
$$
 f_!\:Y\Lcth_\bT\lrarrow X\Lcth_\bW
$$
is exact and preserves infinite products~\cite[Section~3.3]{Pcosh}.

 Similarly, for any $(\bW,\bT)$\+affine morphism of schemes
$f\:Y\rarrow X$, the direct image functor $f_!\:(Y,\cO_Y)\Cosh\rarrow
(X,\cO_X)\Cosh$ takes locally cotorsion $\bT$\+locally contraherent
cosheaves on $Y$ to locally cotorsion $\bW$\+locally contraherent
cosheaves on~$X$.
 The resulting functor between the exact categories of locally cotorsion
locally contraherent cosheaves
$$
 f_!\:Y\Lcth_\bT^\lct\lrarrow X\Lcth_\bW^\lct
$$
is exact and preserves infinite products~\cite[Section~3.3]{Pcosh}.

 Finally, for any flat $(\bW,\bT)$\+affine morphism of schemes
$f\:Y\rarrow X$, the direct image functor $f_!\:(Y,\cO_Y)\Cosh\rarrow
(X,\cO_X)\Cosh$ takes locally injective $\bT$\+locally contraherent
cosheaves on $Y$ to locally injective $\bW$\+locally contraherent
cosheaves on~$X$.
 The resulting functor between the exact categories of locally injective
locally contraherent cosheaves
$$
 f_!\:Y\Lcth_\bT^\lin\lrarrow X\Lcth_\bW^\lin
$$
is exact and preserves infinite products~\cite[Section~3.3]{Pcosh}.

\subsection{Inverse images}
\label{inverse-images-of-O-co-sheaves-subsecn}
 We refer to~\cite[Section~0.10 of the Introduction]{Pcosh}
and~\cite[Section~5.6]{Pphil} for a introductory discussion of partially
defined functors between exact categories, which sheds some light
on the material of this section.

 A morphism of schemes $f\:Y\rarrow X$ is said to be \emph{flat} if,
for every pair of affine open subschemes $U\subset X$ and $V\subset Y$
such that $f(V)\subset U$, the commmutative ring homomorphism
$\cO_X(U)\rarrow\cO_Y(V)$ makes $\cO_Y(V)$ a flat $\cO_X(U)$\+module.
 It suffices to check this condition for affine open subschemes $U$
belonging to any chosen affine open covering $X=\bigcup_\alpha U_\alpha$
of the scheme $X$ and affine open subschemes $V$ belonging to any
chosen affine open covering $f^{-1}(U)=\bigcup_\beta V_\beta$ of
the open subscheme $f^{-1}(U)\subset Y$ \,\cite[Lemma Tag~01U5]{SP}.

 A morphism of schemes $f\:Y\rarrow X$ is said to be \emph{very flat}
if, for every pair of affine open subschemes $U\subset X$ and
$V\subset Y$ such that $f(V)\subset U$, the commmutative ring
homomorphism $\cO_X(U)\rarrow\cO_Y(V)$ makes $\cO_Y(V)$ a very
flat $\cO_X(U)$\+module (in the sense of
Section~\ref{prelim-very-flat-subsecn}).
 It suffices to check this condition for affine open subschemes $U$
belonging to any chosen affine open covering $X=\bigcup_\alpha U_\alpha$
of the scheme~$X$ \,\cite[Corollary~1.2.5(b), Lemma~1.2.6(a),
and Section~1.10]{Pcosh}.

 Equivalently, a morphism of schemes $f\:Y\rarrow X$ is very flat if and
only if, for every pair of affine open subschemes $U\subset X$ and
$V\subset Y$ such that $f(V)\subset U$, the homomorphism of commutative
rings $\cO_X(U)\rarrow\cO_Y(V)$ is very flat, or in other words,
$\cO_Y(V)$ is a very flat commutative algebra over $\cO_X(U)$
(in the sense of the definition at the end of
Section~\ref{prelim-very-flat-subsecn}).
 It suffices to check this condition for affine open subschemes $U$
belonging to any chosen affine open covering $X=\bigcup_\alpha U_\alpha$
of the scheme $X$ and affine open subschemes $V$ belonging to any
chosen affine open covering $f^{-1}(U)=\bigcup_\beta V_\beta$ of
the open subscheme $f^{-1}(U)\subset Y$ \,\cite[Lemma~1.8.6 and
Section~1.10]{Pcosh}.

 A morphism of schemes $f\:Y\rarrow X$ is said to be \emph{coaffine} if
for every affine open subscheme $V\subset Y$ there exists an affine
open subscheme $U\subset X$ such that $f(V)\subset U$.
 Interesting examples of \emph{noncoaffine} morphisms of algebraic
varieties (in particular, with an affine variety~$Y$) are provided by
the \emph{Jouanolou trick} and its generalization by
Thomason~\cite[Proposition~4.3 or~4.4]{Wei}.
 Any morphism into an affine scheme $X$ is coaffine.
 Any open immersion of schemes is coaffine.

 Now let $\bW$ be an open covering of the scheme $X$ and $\bT$ be
an open covering of the scheme~$Y$.
 A morphism of schemes $f\:Y\rarrow X$ is said to be
\emph{($\bW,\bT)$\+coaffine} if, for every affine open subscheme
$V\subset Y$ subordinate to $\bT$, there exists an affine open
subscheme $U\subset X$ subordinate to $\bW$ such that $f(V)\subset U$.

 Given an open covering $\bW$ of the scheme $X$ and a morphism of
schemes $f\:Y\rarrow X$, one can consider the open covering $\bT$ of
the scheme $Y$ formed by the preimages of all affine open subschemes
$U\subset X$ subordinate to~$\bW$.
 Then the morphism~$f$ is ($\bW,\bT)$\+coaffine.
 Moreover, the morphism~$f$ is ($\bW,\bT)$\+affine whenever
it is affine.

 Let $X$ be a scheme with an open covering $\bW$ and $Y$ be a scheme
with an open covering~$\bT$.
 Let $f\:Y\rarrow X$ be a $(\bW,\bT)$\+coaffine morphism of schemes.

 Assume that the morphism~$f$ is very flat, and let $\P$ be
a $\bW$\+locally contraherent cosheaf on~$X$.
 Then the $\bT$\+locally contraherent cosheaf $f^!\P$ on $Y$ is
defined by the rule
$$
 (f^!\P)[V]=\Hom_{\cO_X(U)}(\cO_Y(V),\P[U])
$$
for any affine open subscheme $V\subset Y$ subordinate to $\bT$ and
any chosen affine open subscheme $U\subset X$ subordinate to $\bW$
such that $f(V)\subset U$.
 The $\cO_Y(V)$\+module $\Hom_{\cO_X(U)}(\cO_Y(V),\P[U])$ is
contraadjusted by Lemma~\ref{restriction-co-extension-vfl-cta}(d).

 We refer to~\cite[Sections~2.3 and~3.3]{Pcosh} for the details of
this construction, including an explanation of why the $\bT$\+locally
contraherent cosheaf $f^!\P$ on $Y$ is well-defined.
 Proving this becomes a little bit more involved when the scheme $X$ is
\emph{not} assumed to be semi-separated; see~\cite[Section~3.3]{Pcosh}.
 The resulting inverse image functor between the exact categories of
locally contraherent cosheaves
$$
 f^!\:X\Lcth_\bW\lrarrow Y\Lcth_\bT
$$
is exact and preserves infinite products.

 Assume that the morphism~$f$ is flat, and let $\P$ be
a locally cotorsion $\bW$\+locally contraherent cosheaf on~$X$.
 Then the locally cotorsion $\bT$\+locally contraherent cosheaf $f^!\P$
on $Y$ is defined by the same rule
$$
 (f^!\P)[V]=\Hom_{\cO_X(U)}(\cO_Y(V),\P[U])
$$
for any affine open subscheme $V\subset Y$ subordinate to $\bT$ and
any chosen affine open subscheme $U\subset X$ subordinate to $\bW$
such that $f(V)\subset U$.
 The $\cO_Y(V)$\+module $\Hom_{\cO_X(U)}(\cO_Y(V),\P[U])$ is
cotorsion by Lemma~\ref{restriction-coextension-injective-cotorsion}(c).
 The resulting inverse image functor between the exact categories of
locally cotorsion locally contraherent cosheaves
$$
 f^!\:X\Lcth_\bW^\lct\lrarrow Y\Lcth_\bT^\lct
$$
is exact and preserves infinite products~\cite[Sections~2.3
and~3.3]{Pcosh}.

 Finally, let $f\:Y\rarrow X$ be an arbitrary $(\bW,\bT)$\+coaffine
morphism of schemes, and let $\gJ$ be a locally injective
$\bW$\+locally contraherent cosheaf on~$X$.
 Then the locally injective $\bT$\+locally contraherent cosheaf $f^!\gJ$
on $Y$ is defined by the same rule
$$
 (f^!\gJ)[V]=\Hom_{\cO_X(U)}(\cO_Y(V),\gJ[U])
$$
for any affine open subscheme $V\subset Y$ subordinate to $\bT$ and
any chosen affine open subscheme $U\subset X$ subordinate to $\bW$
such that $f(V)\subset U$.
 The $\cO_Y(V)$\+module $\Hom_{\cO_X(U)}(\cO_Y(V),\gJ[U])$ is
injective by Lemma~\ref{restriction-coextension-injective-cotorsion}(d).
 The resulting inverse image functor between the exact categories of
locally injective locally contraherent cosheaves
$$
 f^!\:X\Lcth_\bW^\lct\lrarrow Y\Lcth_\bT^\lct
$$
is exact and preserves infinite products~\cite[Sections~2.3
and~3.3]{Pcosh}.

 Passing to the inductive limit of categories over the refinements
of the coverings, we obtain an exact functor of inverse image of
locally contraherent cosheaves
$$
 f^!\:X\Lcth\lrarrow Y\Lcth
$$
for any very flat morphism of schemes~$f$, an exact functor of
inverse image of locally cotorsion locally contraherent cosheaves
$$
 f^!\:X\Lcth^\lct\lrarrow Y\Lcth^\lct
$$
for any flat morphism of schemes~$f$, and an exact functor of
inverse image of locally injective contraherent cosheaves
$$
 f^!\:X\Lcth^\lin\lrarrow Y\Lcth^\lin
$$
for an arbitrary morphism of schemes $f\:Y\rarrow X$.

 Let $f\:Y\rarrow X$ be a morphism of schemes, and let $\P$ be a locally
contraherent cosheaf on $X$ such that the locally contraherent cosheaf
of inverse image $f^!\P$ on $Y$ is well-defined as per one of
the constructions above.
 Let $\Q$ be a cosheaf of $\cO_Y$\+modules on~$Y$.
 Then there is a natural (adjunction) isomorphism of abelian
groups~\cite[formulas~(3.4) and~(3.5) in Section~3.3]{Pcosh}
\begin{equation} \label{direct-inverse-cosheaf-general-adjunction}
 \Hom^{\cO_X}(f_!\Q,\P)\simeq\Hom^{\cO_Y}(\Q,f^!\P).
\end{equation}

 Given a scheme $X$ with an open covering $\bW$ and an open subscheme
$Y\subset X$, denote by $\bW|_Y$ the open covering of $Y$ consisting
of all the open subschemes $Y\cap W$ with $W\in\bW$.
 Then the open immersion $j\:Y\rarrow X$ is a very flat
$(\bW,\bW|_Y)$\+coaffine morphism.
 In this case, the functor of inverse image $j^!\:X\Lcth_\bW\rarrow
Y\Lcth_{\bW|_Y}$ defined above in this section agrees with the functor
of restriction to an open subscheme $\P\longmapsto\P|_Y$ from
the end of Section~\ref{cosheaves-of-modules-subsecn} for
$\bW$\+locally contraherent cosheaves $\P$ on~$X$.

 We refer to Section~\ref{inverse-images-of-A-co-sheaves-subsecn} below
for a further discussion with some details spelled out.

\Section{Contraherent Cosheaves of $\cA$-Modules}
\label{contraherent-cosheaves-of-A-modules-secn}

 This section offers several points of view on locally contraherent
cosheaves of modules over quasi-coherent quasi-algebras $\cA$ over
schemes $X$, and spells out the basics of the theory, including
various exact categories of locally contraherent $\cA$\+modules,
direct and inverse image functors, etc.
 This is new material, not yet presented elsewhere, even though parts
of the expositions are parallel to the theory of locally contraherent
cosheaves (of $\cO_X$\+modules) elaborated upon in~\cite{Pcosh} and
discussed in Section~\ref{contraherent-cosheaves-of-O-modules-secn}
above.

\subsection{Quasi-coherent quasi-algebras}
\label{quasi-coherent-quasi-algebras-subsecn}
 This section is largely based on~\cite[second half of
Section~2.3]{Pedg}, with a generalization worked out
in~\cite[Section~2.1]{Ptd} and
Sections~\ref{prelim-quasi-modules-subsecn}\+-%
\ref{prelim-quasi-algebras-subsecn} above.

 Let $X$ be a topological space endowed with a sheaf of commutative
rings~$\cO$.
 We are interested in sheaves of $\cO$\+$\cO$\+bimodules on~$X$.
 For the definition of a \emph{quasi-module} over a commutative ring,
see Section~\ref{prelim-quasi-modules-subsecn}.

\begin{lem} \label{sheaf-of-quasi-modules-locality}
 Let $\cB$ be a sheaf of $\cO$\+$\cO$\+bimodules on~$X$, and let
$U\subset X$ be an open subset with a finite open covering
$U=\bigcup_{\alpha=1}^N V_\alpha$.
 Suppose that, for every index~$\alpha$,
the $\cO(V_\alpha)$\+$\cO(V_\alpha)$\+bimodule $\cB(V_\alpha)$
is a quasi-module over $\cO(V_\alpha)$.
 Then the $\cO(U)$\+$\cO(U)$\+bimodule $\cB(U)$ is a quasi-module
over $\cO(U)$.
\end{lem}

\begin{proof}
 This is similar to the proof of
Lemma~\ref{quasi-module-locality}\,(2)\,$\Rightarrow$\,(1).
 By the sheaf axiom, the natural map $\cB(U)\rarrow
\bigoplus_{\alpha=1}^N\cB(V_\alpha)$ is injective.
 It is also a homomorphism of $\cO(U)$\+$\cO(U)$\+bimodules.
 By Lemma~\ref{quasi-module-restriction-of-scalars}, the assumption
that the $\cO(V_\alpha)$\+$\cO(V_\alpha)$\+bimodule $\cB(V_\alpha)$
is a quasi-module over $\cO(V_\alpha)$ implies that, viewed as
a $\cO(U)$\+$\cO(U)$\+bimodule, $\cB(V_\alpha)$ is a quasi-module
over $\cO(U)$.
 Finally, the direct sums of quasi-modules are quasi-modules and
any subbimodule of a quasi-module is a quasi-module, as mentioned in
the beginning of Section~\ref{prelim-quasi-modules-subsecn}. 
\end{proof}

 Assume that quasi-compact open subsets in the topological space $X$
form a base of the topology.
 Then we will say that a sheaf of $\cO$\+$\cO$\+bimodules $\cB$ on $X$
is a \emph{sheaf of quasi-modules over~$\cO$} if, for every
quasi-compact open subset $U\subset X$, the $\cO(U)$\+$\cO(U)$\+bimodule
$\cB(U)$ is a quasi-module over $\cO(U)$.
 By Lemma~\ref{sheaf-of-quasi-modules-locality}, it suffices to
check this condition for open subsets $U$ belonging to any chosen
base of neighborhoods of zero $\bB$ in $X$ consisting of quasi-compact
open subsets.

 A \emph{sheaf of quasi-algebras} $\cA$ over the sheaf of rings~$\cO$
on $X$ is defined as a sheaf of associative rings on $X$ endowed with
a morphism of sheaves of rings $\cO\rarrow\cA$ such that, viewed with
the induced structure of a sheaf of $\cO$\+$\cO$\+bimodules, the sheaf
$\cA$ is a sheaf of quasi-modules over~$\cO$ on~$X$.

 Now let $X$ be a scheme.
 As a particular case of the discussion above, we obtain the definitions
of a \emph{sheaf of quasi-modules over~$\cO_X$} and a \emph{sheaf of
quasi-algebras over~$\cO_X$} on~$X$.

\begin{lem} \label{quasi-coherence-of-sheaves-of-quasi-modules}
 Let $X$ be a scheme and $\cB$ be a sheaf of quasi-modules over~$\cO_X$.
 Then the sheaf $\cB$, viewed as a sheaf of $\cO_X$\+modules with
respect to the left $\cO_X$\+module structure, is quasi-coherent if
and only if it is quasi-coherent as a sheaf of $\cO_X$\+modules with
respect to the right $\cO_X$\+module structure.
\end{lem}

\begin{proof}
 This assertion is based on
Corollary~\ref{for-quasi-module-quasi-coherence-aux-cor}.
 The argument from~\cite[proof of Proposition~2.4]{Pedg} applies.
\end{proof}

 A \emph{quasi-coherent sheaf of quasi-modules on~$X$} (or just
a \emph{quasi-coherent quasi-module over~$X$} for brevity) is defined
as a sheaf of quasi-modules over $\cO_X$ satisfying the equivalent
conditions of Lemma~\ref{quasi-coherence-of-sheaves-of-quasi-modules}.
 \emph{Morphisms of quasi-coherent quasi-modules} are just the usual
morphisms of sheaves of $\cO_X$\+$\cO_X$\+bimodules.
 We denote the category of quasi-coherent quasi-modules on $X$ by
$X\QQcoh$.
 One can easily check that $X\QQcoh$ is a Grothendieck abelian category.
 The two forgetful functors $X\QQcoh\rightrightarrows X\Qcoh$ (assigning
to a quasi-coherent quasi-module its underlying quasi-coherent sheaves
of with the left and right $\cO_X$\+module structures) are exact and
preserve infinite direct sums.

 Any quasi-coherent sheaf on $X$ can be viewed as a quasi-coherent
quasi-module over $X$ in which the left and right actions of $\cO_X$
agree.
 The resulting inclusion functor $X\Qcoh\rarrow X\QQcoh$ is exact and
fully faithful, and preserves infinite direct sums.

 For an affine scheme $U=\Spec R$, the category of quasi-coherent
quasi-modules over $U$ is naturally equivalent to the category
of quasi-modules over $R$, that is $U\QQcoh\simeq R\QMod$.
 Lemma~\ref{quasi-module-localization-lemma}(b) is relevant here.

\begin{lem} \label{direct-image-quasi-coherent-quasi-module}
 Let $f\:Y\rarrow X$ be a quasi-compact quasi-separated morphism of
schemes and $\cB$ be a quasi-coherent quasi-module on~$Y$.
 Then the direct image $f_*\cB$ is a quasi-coherent quasi-module on~$X$.
\end{lem}

\begin{proof}
 Recall that, by the definition, for any morphism of schemes (or ringed
spaces), the functor of direct image of sheaves of $\cO$\+modules
$f_*\:(Y,\cO_Y)\Sh\rarrow(X,\cO_X)\Sh$ agrees with the functor of
direct image of sheaves of abelian groups under the map of topological
spaces $f\:Y\rarrow X$ (see the beginning of
Section~\ref{direct-images-of-O-co-sheaves-subsecn}).
 Accordingly, for any sheaf of $\cO_Y$\+$\cO_Y$\+bimodules $\cB$ on $Y$,
the direct image $f_*\cB$ is naturally a sheaf of
$\cO_X$\+$\cO_X$\+bimodules on~$X$.
 Furthermore, for any quasi-compact quasi-separated morphism of
schemes~$f$, the functor of direct image~$f_*$ preserves quasi-coherence
of sheaves of $\cO$\+modules~\cite[Proposition~I.6.7.1]{EGA1},
\cite[Lemma Tag~01LC]{SP}.

 It remains to show that, for every affine open subscheme $U\subset X$,
the $\cO_X(U)$\+$\cO_X(U)$\+bimodule $(f_*\cB)(U)=\cB(f^{-1}(U))$ is
a quasi-module over $\cO_X(U)$.
 This is another application of the argument from the proofs of
Lemmas~\ref{quasi-module-locality}\,(2)\,$\Rightarrow$\,(1)
and~\ref{sheaf-of-quasi-modules-locality}.
 Since the morphism~$f$ is quasi-compact by assumption, the open
subset $f^{-1}(U)\subset Y$ is quasi-compact.
 So $f^{-1}(U)$ is a finite union of affine open subschemes $V_\alpha$
in~$Y$; that is $f^{-1}(U)=\bigcup_{\alpha=1}^N V_\alpha$.
 Now the $\cO_X(U)$\+$\cO_X(U)$\+bimodule $\cB(V_\alpha)$ is
a quasi-module over $\cO_X(U)$ by
Lemma~\ref{quasi-module-restriction-of-scalars};
the $\cO_X(U)$\+$\cO_X(U)$\+bimodule $\cB(f^{-1}(U))$ is a subbimodule
of $\bigoplus_{\alpha=1}^N\cB(V_\alpha)$, and any subbimodule
of a quasi-module is a quasi-module.  \hbadness=1200
\end{proof}

 Let $j\:Y\rarrow X$ be an open immersion of schemes and $\cA$ be
a quasi-coherent quasi-module on~$X$.
 Viewing $\cA$ as a quasi-coherent sheaf on $X$ in the left or right
$\cO_X$\+module structure, one can consider the inverse image~$j^*\cA$.
 Both the inverse images so obtained are, of course, the same things
as the restriction $\cA|_Y$ of the sheaf of $\cO_X$\+$\cO_X$\+bimodules
$\cA$ to the open subscheme $Y\subset X$.
 So we write $j^*\cA=\cA|_Y$.
 The sheaf of $\cO_Y$\+$\cO_Y$\+bimodules $\cA|_Y$ is a quasi-coherent
quasi-module on~$Y$.

 Recall from Section~\ref{prelim-quasi-modules-subsecn} that, given
a commutative ring $R$ and two $R$\+$R$\+bimodules $A$ and $B$,
the notation $A\ot_RB$ stands for the tensor product of the right
$R$\+module structure on $A$ and the left $R$\+module structure on $B$,
as usual.
 Given an $R$\+module $M$, the notation $B\ot_RM$ and $M\ot_RB$ is
understood similarly.

\begin{lem} \label{quasi-coherent-quasi-modules-tensor-product}
 Let $X$ be a scheme. \par
\textup{(a)} Given two quasi-coherent quasi-modules $\cA$ and $\cB$
over $X$, their tensor product $\cA\ot_{\cO_X}\cB$, defined by the rule
$$
 (\cA\ot_{\cO_X}\cB)(U)=\cA(U)\ot_{\cO_X(U)}\cB(U),
$$
for all affine open subschemes $U\subset X$, is again a quasi-coherent
quasi-module over~$X$. \par
\textup{(b)} Given a quasi-coherent quasi-module $\cB$ and
a quasi-coherent sheaf $\M$ on $X$, the two tensor products
$\cB\ot_{\cO_X}\M$ and $\M\ot_{\cO_X}\cB$, defined by
the rules
$$
 (\cB\ot_{\cO_X}\M)(U)=\cB(U)\ot_{\cO_X(U)}\M(U)
 \text{ and }
 (\M\ot_{\cO_X}\cB)(U)=\M(U)\ot_{\cO_X(U)}\cB(U),
$$
for all affine open subschemes $U\subset X$, are quasi-coherent sheaves
on~$X$.
\end{lem}

\begin{proof}
 This is our version of~\cite[Lemma~2.5]{Pedg}.
 In part~(a), the point is that $\cA(U)\ot_{\cO_X(U)}\cB(U)$ is
a quasi-module over $\cO_X(U)$ by~\cite[Proposition~2.4]{Ptd}, and for
every pair of affine open subschemes $V\subset U\subset X$ one has
\begin{multline*}
 \cA(V)\ot_{\cO_X(V)}\cB(V)\simeq
 (\cO_X(V)\ot_{\cO_X(U)}\cA(U))\ot_{\cO_X(V)}
 (\cO_X(V)\ot_{\cO_X(U)}\cB(U)) \\ \simeq
 \cO_X(V)\ot_{\cO_X(U)}(\cA(U)\ot_{\cO_X(U)}\cB(U)).
\end{multline*}
 Part~(b) can be viewed as a particular case of part~(a).
\end{proof}

 The abelian category of quasi-coherent quasi-modules $X\QQcoh$ is
an (associative, noncommutative) monoidal category with respect to
the tensor product operation~$\ot_{\cO_X}$.
 The structure sheaf $\cO_X$, viewed as an object of $X\QQcoh$ via
the fully faithful inclusion functor $X\Qcoh\rarrow X\QQcoh$, is
the unit object of this monoidal category.
 Viewing $X\Qcoh$ as a (commutative) monoidal category with respect
to the usual tensor product functor~$\ot_X$ makes the inclusion
functor $X\Qcoh\rarrow X\QQcoh$ monoidal.
 Furthermore, the category $X\Qcoh$ is both a left and a right module
category over the monoidal category $X\QQcoh$ with respect to
the tensor product functors~$\ot_{\cO_X}$ defined above.

 An (associative, unital) \emph{quasi-coherent quasi-algebra} $\cA$
over $X$ is defined as a monoid object in the monoidal category
$X\QQcoh$.
 So a quasi-coherent quasi-algebra $\cA$ is a quasi-coherent
quasi-module endowed with quasi-coherent quasi-module morphisms of
\emph{multiplication} $\cA\ot_{\cO_X}\cA\rarrow\cA$ and \emph{unit}
$\cO_X\rarrow\cA$ satisfying the usual associativity and unitality
axioms.
 Alternatively, a quasi-coherent quasi-algebra $\cA$ over $X$ is
the same thing as a sheaf of quasi-algebras over $\cO_X$ whose
underlying sheaf of quasi-modules over $\cO_X$ is quasi-coherent.

 For an affine scheme $U=\Spec R$, the category of quasi-coherent
quasi-algebras over $U$ is naturally equivalent to the category
of quasi-algebras over~$R$ (in the sense of
Section~\ref{prelim-quasi-algebras-subsecn}).
 Lemma~\ref{quasi-algebra-co-extension-of-scalars}(a) is relevant here.

\begin{rems}
 (1)~Let $K\rarrow R$ be a homomorphism of commutative rings, and
let $B$ be an $R$\+$R$\+bimodule over~$K$ (i.~e., the left and right
actions of $K$ in $B$ agree).
 Consider the ring $T=R\ot_KR$.
 Then \cite[Corollary~2.3]{Ptd} says that $B$ is a quasi-module over $R$
if and only if, viewed as a $T$\+module, $B$ is supported
set-theoretically in the image of the diagonal closed immersion of
schemes $\Spec R\rarrow\Spec R\times_{\Spec K}\Spec R=\Spec T$.

 (2)~Consider a morphism of schemes $f\:X\rarrow S$.
 Let us say that a sheaf of $\cO_X$\+$\cO_X$\+bimodules $\cB$ on $X$
is a \emph{module over~$S$} if, for every pair of
affine open subschemes $U\subset S$ and $V\subset X$ such that
$f(V)\subset U$, the left and right actions of $\cO_X(V)$ in
$\cB(V)$ restrict to one and the same action of $\cO_S(U)$ in $\cB(V)$,
i.~e., the left and right actions of $\cO_S(U)$ in $\cB(V)$ agree.

 (3)~Assume that the morphism~$f$ is separated, i.~e., the diagonal
morphism $X\rarrow X\times_SX$ is a closed immersion.
 Globalizing the assertion of~(1) from affine to arbitrary schemes, one
can say that the abelian category of quasi-coherent quasi-modules over
$X$ that are modules over $S$ is equivalent to the abelian category of
quasi-coherent sheaves on $W=X\times_SX$ supported set-theoretically
in the image $Z\subset W$ of the diagonal closed immersion
$X\rarrow X\times_SX$.
\end{rems}

\begin{exs} \label{fiberwise-differential-operators-examples}
 The following examples of sheaves of quasi-modules, sheaves of
quasi-algebras, quasi-coherent quasi-modules, and quasi-coherent
quasi-algebras are relevant for the purposes of this paper.

\smallskip
 (1)~Let $\tau\:X\rarrow T$ be a morphism of schemes and $\cU$, $\cV$
be quasi-coherent sheaves on~$X$.
 Let $Y\subset X$ and $W\subset T$ be affine open subschemes such that
$\tau(Y)\subset W$.
 Then the rule $\cD_{X/T}(\cU,\cV)(Y)=
\cD_{\cO_X(Y)/\cO_T(W)}(\cU(Y),\cV(Y))$ (where the notation in
the right-hand side comes from
Examples~\ref{diffoperators-quasi-algebras-quasi-modules-exs})
defines a sheaf of quasi-modules $\cD_{X/T}(\cU,\cV)$ over~$X$
\,\cite[Section~0.5 of the Introduction]{Ptd}.
 The restriction maps from affine open subschemes $Y\subset X$ to
smaller affine open subschemes $Y'\subset Y$ are constructed
using~\cite[Theorem~6.2]{Ptd}, and the result
of~\cite[Proposition~6.4]{Ptd} implies that the sheaf axiom holds.
 The sheaf $\cD_{X/T}(\cU,\cV)$ is called the \emph{sheaf of}
(\emph{fiberwise}) \emph{$X/T$\+differential operators}
acting from $\cU$ to~$\cV$.

\smallskip
 (2)~In the same setting as in~(1), the rules
$\cD^\qu_{X/T}(\cU,\cV)(Y)=
\cD^\qu_{\cO_X(Y)/\cO_T(W)}(\cU(Y),\allowbreak\cV(Y))$ and
$\cD^\str_{X/T}(\cU,\cV)(Y)=
\cD^\str_{\cO_X(Y)/\cO_T(W)}(\cU(Y),\cV(Y))$ define sheaves of
quasi-modules $\cD_{X/T}^\qu(\cU,\cV)$ and $\cD_{X/T}^\str(\cU,\cV)$
over~$X$.
 The restriction maps from affine open subschemes $Y\subset X$ to
smaller affine open subschemes are constructed
using~\cite[Proposition~6.6]{Ptd}, and the result
of~\cite[Proposition~6.8]{Ptd} implies the sheaf axiom.
 The sheaf $\cD_{X/T}^\qu(\cU,\cV)$ is called the \emph{sheaf of
quite $X/T$\+differential operators} acting from $\cU$ to $\cV$,
while the sheaf $\cD_{X/T}^\str(\cU,\cV)$ is called the \emph{sheaf of
strongly $X/T$\+differential operators}.
 Obviously, one has $\cD_{X/T}^\str(\cU,\cV)\subset
\cD_{X/T}^\qu(\cU,\cV)\subset\cD_{X/T}(\cU,\cV)$.

\smallskip
 (3)~In particular, when $\cU=\cV$, all the three sheaves
$\cD_{X/T}^\str(\cU,\cU)\subset\cD_{X/T}^\qu(\cU,\cU)\subset
\cD_{X/T}(\cU,\cU)$ have natural structures of sheaves of quasi-algebras
over~$X$ (with respect to the composition of differential operators).
 See the discussion in
Examples~\ref{diffoperators-quasi-algebras-quasi-modules-exs}.

\smallskip
 (4)~A morphism of schemes $\tau\:X\rarrow T$ is said to be
\emph{locally of finite type} if for every pair of affine open
subschemes $V\subset X$ and $U\subset T$ such that $\tau(V)\subset U$,
the commutative $\cO_T(U)$\+algebra $\cO_X(V)$ is finitely generated.
 It suffices to check this condition for affine open subschemes $U$
belonging to any given affine open covering $T=\bigcup_\alpha U_\alpha$
of the scheme $T$ and affine open subschemes $V$ belonging to any
given affine open covering $\tau^{-1}(U)=\bigcup_\beta V_\beta$ of
the open subscheme $\tau^{-1}(U)\subset X$
\,\cite[D\'efinition~I.6.2.1 and Proposition~I.6.2.5]{EGA1},
\cite[Definition Tag~01T1 and Lemma Tag~01T2]{SP}.

 For a morphism of schemes $\tau\:X\rarrow T$ locally of finite type
and any quasi-coherent sheaves $\cU$ and $\cV$ on $X$, our three
sheaves of fiberwise differential operators coincide, that is,
$\cD_{X/T}^\str(\cU,\cV)=\cD_{X/T}^\qu(\cU,\cV)=\cD_{X/T}(\cU,\cV)$.
 This follows from~\cite[Lemmas~1.3 and~2.1]{Ptd}; see the last
paragraph of
Examples~\ref{diffoperators-quasi-algebras-quasi-modules-exs}.

\smallskip
 The following examples will be discussed in more detail in
Section~\ref{twisted-lie-algebroids-secn}.

 (4)~A morphism of schemes $\tau\:X\rarrow T$ is said to be
\emph{locally of finite presentation} if for every pair of affine open
subschemes $V\subset X$ and $U\subset T$ such that $\tau(V)\subset U$,
the commutative $\cO_T(U)$\+algebra $\cO_X(V)$ is finitely presented
(as an algebra, not as a module!).
 It suffices to check this condition for affine open subschemes $U$
belonging to any chosen affine open covering $T=\bigcup_\alpha U_\alpha$
of the scheme $T$ and affine open subschemes $V$ belonging to any
chosen affine open covering $\tau^{-1}(U)=\bigcup_\beta V_\beta$ of
the open subscheme $\tau^{-1}(U)\subset X$
\,\cite[D\'efinition~I.6.2.1 and Proposition~I.6.2.9]{EGA1},
\cite[Definition Tag~01TP and Lemma Tag~01TQ]{SP}.

 The definition of a \emph{locally finitely presented} quasi-coherent
sheaf on a scheme $X$ can be found below in
Section~\ref{loc-free-sheaves-subsecn}.
 In the terminology of~\cite[0.5.2.5]{EGA1}
and~\cite[Definition Tag~01BN]{SP}, such quasi-coherent
sheaves are said to be ``of finite presentation''.

 For a morphism of schemes $\tau\:X\rarrow T$ locally of finite 
presentation, a locally finitely presented quasi-coherent sheaf $\cU$
on $X$, and any quasi-coherent sheaf $\cV$ on $X$, the sheaf of
quasi-modules $\cD_{X/T}(\cU,\cV)$ on $X$ is
quasi-coherent~\cite[Propositions~IV.16.8.6 and~IV.16.8.8]{EGA4}.
 Accordingly, the sheaf of quasi-algebras $\cD_{X/T}(\cU,\cU)$ is
a quasi-coherent quasi-algebra over $X$ in this case.

\smallskip
 (5)~For any quasi-coherent twisted Lie algebroid $(\g,\widetilde\g)$
over a scheme~$X$ (in the sense of
Section~\ref{twisted-lie-algebroids-subsecn}), the twisted universal
enveloping quasi-algebra $\cA_X(\g,\widetilde\g)$ is a quasi-coherent
quasi-algebra over~$X$.
 See Section~\ref{enveloping-algebra-subsecn} below.

\smallskip
 (6)~In particular, for any weakly smooth morphism of schemes
$\tau\:X\rarrow T$ (in the sense of
Section~\ref{crystalline-diffoperators-subsecn}), the sheaf of rings of
crystalline differential operators $\cD^\cry_{X/T}=\cA_X(\vect_{X/T})$
is a quasi-coherent quasi-algebra over~$X$.
 Here $\vect_{X/T}$ is the Lie algebroid of vector fields on $X$
over~$T$ (as defined in Section~\ref{vector-fields-subsecn}).
 See Section~\ref{crystalline-diffoperators-subsecn} for a discussion
of this example.

\smallskip
 (7)~As another particular case of~(5) and a generalization of~(6),
for any weakly smooth morphism of schemes $\tau\:X\rarrow T$ and
any quasi-coherent twisted Lie algebroid $(\g,\widetilde\g)=
(\vect_{X/T},\widetilde\vect_{X/T})$ whose underlying quasi-coherent
Lie algebroid~$\g$ is the Lie algebroid of vector fields
$\g=\vect_{X/T}$, the sheaf of rings of twisted differential operators
$\cA_X(\g,\widetilde\g)$ is a quasi-coherent quasi-algebra over~$X$.
 See Section~\ref{twisted-diffoperators-subsecn} for a discussion
of this class of examples.
\end{exs}

\subsection{$\Cohom$ from quasi-coherent quasi-modules to
contraherent cosheaves} \label{cohom-from-quasi-modules-subsecn}
 The material of this section is a (partial) generalization
of~\cite[Sections~2.4 and~3.6]{Pcosh}.
 We refer to~\cite[Section~0.10 of the Introduction]{Pcosh}
and~\cite[Sections~5.5\+-5.6]{Pphil} for a relevant discussion of
partially defined functors between exact categories.

 Let $X$ be a scheme with an open covering~$\bW$.
 Denote by $\bB$ the base of open subsets in $X$ consisting of all
the affine open subschemes $U\subset X$ subordinate to~$\bW$.

 Let $\cB$ be a quasi-coherent quasi-module and $\P$ be
a $\bW$\+locally contraherent cosheaf on~$X$.
 The copresheaf of $\cO_X$\+modules $\Cohom_X(\cB,\P)$ on $\bB$
is defined by the rule
$$
 \Cohom_X(\cB,\P)[U]=\Hom_{\cO_X(U)}(\cB(U),\P[U])
$$
for all affine open subschemes $U\subset X$ subordinate to~$\bW$.

 Here, as usual, the $\Hom$ of $\cO_X(U)$\+modules is taken with
respect to the left $\cO_X(U)$\+module structure on $\cB(U)$ (and
the only given $\cO_X(U)$\+module structure on $\P[U]$).
 The abelian group $\Hom_{\cO_X(U)}(\cB(U),\P[U])$ is endowed with
the $\cO_X(U)$\+module structure induced by the right
$\cO_X(U)$\+module structure on $\cB(U)$ (as per the convention
in Section~\ref{prelim-quasi-modules-subsecn}).
{\hbadness=1075\par}

 The $\cO_X(U)$\+module $\cB(U)$ depends contravariantly on
$U\in\bB$, while the $\cO_X(U)$\+module $\P[U]$ depends covariantly
on $U\in\bB$.
 Hence the $\cO_X(U)$\+module $\Hom_{\cO_X(U)}(\cB(U),\P[U])$ depends
covariantly on $U\in\bB$, so $\Cohom_X(\cB,\P)$ is indeed
a copresheaf of $\cO_X$\+modules on~$\bB$.
{\hbadness=1350\par}

\begin{lem} \label{cohom-loc-contraherent-lemma}
\textup{(a)} If a quasi-coherent quasi-module $\cB$ on $X$ is very
flat as a quasi-coherent sheaf on $X$ with respect to the left
$\cO_X$\+module structure, and\/ $\P$ is a\/ $\bW$\+locally
contraherent cosheaf on $X$, then\/ $\Cohom_X(\cB,\P)$ is
a contraherent copresheaf on\/~$\bB$. \par
\textup{(b)} If a quasi-coherent quasi-module $\cB$ on $X$ is
flat as a quasi-coherent sheaf on $X$ with respect to the left
$\cO_X$\+module structure, and\/ $\P$ is a locally cotorsion\/
$\bW$\+locally contraherent cosheaf on $X$, then\/ $\Cohom_X(\cB,\P)$ is
a locally cotorsion contraherent copresheaf on\/~$\bB$. \par
\textup{(c)} If a quasi-coherent quasi-module $\cB$ on $X$ is
flat as a quasi-coherent sheaf on $X$ with respect to the right
$\cO_X$\+module structure, and\/ $\gJ$ is a locally injective\/
$\bW$\+locally contraherent cosheaf on $X$, then\/ $\Cohom_X(\cB,\gJ)$
is a locally injective contraherent copresheaf on\/~$\bB$. \par
\textup{(d)} For any quasi-coherent quasi-module $\cB$ on $X$, if\/
$\gJ$ is a locally injective\/ $\bW$\+locally contraherent cosheaf
on $X$, then\/ $\Cohom_X(\cB,\gJ)$ is a locally cotorsion contraherent
copresheaf on\/~$\bB$.
\end{lem}

\begin{proof}
 The copresheaf of $\cO_X$\+modules $\Cohom_X(\cB,\P)$ on $\bB$
satisfies the contraherence axiom~(i) from
Section~\ref{locally-contraherent-cosheaves-subsecn} for all
quasi-coherent quasi-modules $\cB$ and all $\bW$\+locally contraherent
cosheaves $\P$ on~$X$.
 Indeed, we have
\begin{multline*}
 \Hom_{\cO_X(U)}(\cO_X(V),\Cohom_X(\cB,\P)[U]) \\
 = \Hom_{\cO_X(U)}(\cO_X(V),\Hom_{\cO_X(U)}(\cB(U),\P[U])) \\
 \simeq\Hom_{\cO_X(U)}(\cB(U)\ot_{\cO_X(U)}\cO_X(V),\>\P[U]) \\
 \simeq\Hom_{\cO_X(U)}(\cO_X(V)\ot_{\cO_X(U)}\cB(U)\ot_{\cO_X(U)}
 \cO_X(V),\>\P[U]) \\
 \simeq\Hom_{\cO_X(U)}(\cB(U)\ot_{\cO_X(U)}\cO_X(V),\>
 \Hom_{\cO_X(U)}(\cO_X(V),\P[U])) \\
 \simeq\Hom_{\cO_X(U)}(\cB(V),\P[V])
 \simeq\Hom_{\cO_X(V)}(\cB(V),\P[V])
\end{multline*}
for every pair of affine open subschemes $V\subset U$, \ $U$, $V\in\bB$,
by Lemma~\ref{quasi-module-localization-lemma}(a).

 It remains to point out that the contraadjustedness axiom~(ii) holds
in part~(a) by Lemma~\ref{very-flat-quasi-module-lemma}(b).
 In the assumptions of part~(b), the $\cO_X(U)$\+module
$\Cohom_X(\cB,\P)[U]$ is cotorsion by
Lemma~\ref{hom-injective-cotorsion}(a).
 In the assumptions of part~(c), the $\cO_X(U)$\+module
$\Cohom_X(\cB,\P)[U]$ is injective by
Lemma~\ref{hom-injective-cotorsion}(b).
 In the assumptions of part~(d), the $\cO_X(U)$\+module
$\Cohom_X(\cB,\P)[U]$ is cotorsion by
Lemma~\ref{hom-injective-cotorsion}(c).
\end{proof}

 Using
Corollary~\ref{loc-contraherent-cosheaves-recovered-from-affines},
in each of the cases~(a\+-d) of Lemma~\ref{cohom-loc-contraherent-lemma}
we extend the contraherent copresheaf $\Cohom_X(\cB,\P)$ or
$\Cohom_X(\cB,\gJ)$ on $\bB$ to a $\bW$\+locally contraherent cosheaf
on (all open subsets of) the whole scheme~$X$.
 The resulting $\bW$\+locally contraherent cosheaf on $X$ is also
denoted by $\Cohom_X(\cB,\P)$ or $\Cohom_X(\cB,\gJ)$, respectively.
 To summarize:
\begin{itemize}
\item $\Cohom_X(\cB,\P)$ is $\bW$\+locally contraherent if $\cB$ is
very flat on the left and $\P$ is $\bW$\+locally contraherent;\
\item $\Cohom_X(\cB,\P)$ is locally cotorsion $\bW$\+locally
contraherent if $\cB$ is flat on the left and $\P$ is locally cotorsion
$\bW$\+locally contraherent;
\item $\Cohom_X(\cB,\gJ)$ is locally injective $\bW$\+locally
contraherent if $\cB$ is flat on the right and $\gJ$ is locally
injective $\bW$\+locally contraherent;
\item $\Cohom_X(\cB,\gJ)$ is locally cotorsion $\bW$\+locally
contraherent if $\gJ$ is locally injective $\bW$\+locally contraherent.
\end{itemize}

 The following version of the $\Cohom$ projection
formula~\cite[formula~(3.21) in Section~3.8]{Pcosh} for quasi-coherent
quasi-modules will be needed in
Section~\ref{koszul-resolutions-of-trivial-ctrh-cdg-mods}.
 Let $X$ be a scheme with an open covering $\bW$, let $Y\subset X$
be an open subscheme whose open immersion morphism $j\:Y\rarrow X$
is affine, and let $\bW|_Y$ be the restriction of the open covering
$\bW$ to the open subscheme~$Y$.
 Let $\cB$ be a quasi-coherent quasi-module on $X$ and $\Q$ be
a $\bW|_Y$\+locally contraherent cosheaf on~$Y$.

 Assume that $\cB$ is very flat as a quasi-coherent sheaf on $X$ with
respect to its left $\cO_X$\+module structure.
 Then there is a natural isomorphism of $\bW$\+locally contraherent
cosheaves
\begin{equation} \label{qcoh-quasi-mod-cohom-projection-formula}
 \Cohom_X(\cB,j_!\Q)\simeq j_!\Cohom_Y(j^*\cB,\Q)
\end{equation}
on the scheme~$X$.
 The same computation as in~\cite[Section~3.8]{Pcosh} is applicable.
 The similar isomorphism of locally cotorsion $\bW$\+locally
contraherent cosheaves on $X$ holds for any quasi-coherent quasi-module
$\cB$ on $X$ that is flat as a quasi-coherent sheaf with respect to its
left $\cO_X$\+module structure and any locally cotorsion
$\bW|_Y$\+locally contraherent cosheaf $\Q$ on~$Y$.

\subsection{Locally contraherent cosheaves of $\cA$-modules}
\label{cosheaves-of-A-modules-subsecn}
 Let $X$ be a scheme and $\cA$ be a sheaf of rings on $X$ endowed with
a morphism of sheaves of rings $\cO_X\rarrow\cA$.
 Then a \emph{quasi-coherent sheaf of left $\cA$\+modules} on~$X$ (or
a \emph{quasi-coherent left $\cA$\+module} for brevity) can be simply
defined as a sheaf of left $\cA$\+modules whose underlying sheaf of
$\cO_X$\+modules is quasi-coherent.
 \emph{Quasi-coherent right $\cA$\+modules} are defined similarly.
 The category of quasi-coherent left $\cA$\+modules $\cA\Qcoh$ and
the category of quasi-coherent right $\cA$\+modules $\Qcohr\cA$ are
Grothendieck abelian categories.
 The forgetful functor $\cA\Qcoh\rarrow X\Qcoh$ is exact and preserves
infinite direct sums.

 Let $\bW$ be an open covering of~$X$.
 Then a \emph{$\bW$\+locally contraherent cosheaf of} (\emph{left})
\emph{$\cA$\+modules} on~$X$ (or a \emph{$\bW$\+locally contraherent
$\cA$\+module} for brevity) can be defined as a cosheaf of
$\cA$\+modules on~$X$ (in the sense of
Section~\ref{cosheaves-of-modules-subsecn}) whose underlying
cosheaf of $\cO_X$\+modules is $\bW$\+locally contraherent.
 A $\bW$\+locally contraherent $\cA$\+module is said to be
\emph{$X$\+locally cotorsion} (respectively, \emph{$X$\+locally
injective}) if its underlying $\bW$\+locally contraherent cosheaf
of $\cO_X$\+modules is locally cotorsion (resp., locally injective).

 We will denote the category of $\bW$\+locally contraherent left
$\cA$\+modules by $\cA\Lcth_\bW$, the category of $X$\+locally
cotorsion $\bW$\+locally contraherent left $\cA$\+modules by
$\cA\Lcth_\bW^{X\dlct}$, and the category of $X$\+locally
cotorsion $\bW$\+locally contraherent left $\cA$\+modules by
$\cA\Lcth_\bW^{X\dlin}$.

 We also put $\cA\Ctrh=\cA\Lcth_{\{X\}}$, \ $\cA\Ctrh^{X\dlct}=
\cA\Lcth_{\{X\}}^{X\dlct}$, and $\cA\Ctrh^{X\dlin}=
\cA\Lcth_{\{X\}}^{X\dlin}$.
 The objects of $\cA\Ctrh$ are called \emph{contraherent cosheaves of
$\cA$\+modules} (or \emph{contraherent $\cA$\+modules} for brevity),
the objects of $\cA\Ctrh^{X\dlct}$ are called \emph{$X$\+locally
cotorsion contraherent $\cA$\+modules}, and similarly for the objects
of $\cA\Ctrh^{X\dlin}$.
 Finally, we put $\cA\Lcth=\bigcup_\bW\cA\Lcth_\bW$, \
$\cA\Lcth^{X\dlct}=\bigcup_\bW\cA\Lcth_\bW^{X\dlct}$, and
$\cA\Lcth^{X\dlin}=\bigcup_\bW\cA\Lcth_\bW^{X\dlin}$.
 The objects of $\cA\Lcth$ are called \emph{locally contraherent
cosheaves of $\cA$\+modules}, and similarly for the objects of
$\cA\Lcth^{X\dlct}$ and $\cA\Lcth^{X\dlin}$.

 A short sequence of $\bW$\+locally contraherent cosheaves of
$\cA$\+modules is said be \emph{exact} in $\cA\Lcth_\bW$ if its
underlying short sequence of $\bW$\+locally contraherent cosheaves
of $\cO_X$\+modules is exact as a short sequence in $X\Lcth_\bW$.
 The similar rules define short exact sequence in all the other
categories of (locally) contraherent cosheaves of $\cA$\+modules
mentioned in the previous two paragraphs.
 Endowed with these classes of short exact sequences, all of these
categories become exact categories.

 The exact categories $\cA\Lcth_\bW$, \ $\cA\Lcth_\bW^{X\dlct}$, and
$\cA\Lcth_\bW^{X\dlin}$ have exact functors of infinite products.
 The forgetful functors from these exact categories to the respective
exact categories of locally contraherent cosheaves (of $\cO_X$\+modules)
on $X$ are exact and preserve infinite products.

 The forgetful functors $\cA\Lcth_\bW\rarrow X\Lcth_\bW$, \ 
$\cA\Lcth_\bW^{X\dlct}\rarrow X\Lcth_\bW^\lct$, and
$\cA\Lcth_\bW^{X\dlin}\rarrow X\Lcth_\bW^\lin$ also (preserve and)
reflect admissible monomorphisms and admissible epimorphisms.
 So a morphism in $\cA\Lcth_\bW$ is an admissible monomorphism
(respectively, an admissible epimorphisms) if and only if it is
an admissible monomorphism (resp., admissible epimorphism) in
$X\Lcth_\bW$; and the same applies to the two other forgetful
functors above.

 Now let $\cA$ be a quasi-coherent quasi-algebra over~$X$.
 As explained in Section~\ref{quasi-coherent-quasi-algebras-subsecn},
the category of quasi-coherent quasi-modules $X\QQcoh$ is a monoidal
category, the category of quasi-coherent sheaves $X\Qcoh$ is a (left
and right) module category over $X\QQcoh$, and $\cA$ is a monoid
object in $X\QQcoh$.
 In this context, quasi-coherent left $\cA$\+modules can be equivalently
defined as left module objects in the left module category $X\Qcoh$
over the monoid object $\cA$ in the monoidal category $X\QQcoh$.
 So a quasi-coherent left $\cA$\+module $\M$ is the same thing as
a quasi-coherent sheaf on $X$ endowed with a \emph{left action}
morphism of quasi-coherent sheaves $\cA\ot_{\cO_X}\M\rarrow\M$
satisfying the usual associativity and unitality axioms together
with the multiplication and unit morphisms of~$\cA$.

 Given a quasi-coherent sheaf $\M$ on $X$, the tensor product
$\cA\ot_{\cO_X}\M$ has a natural structure of a quasi-coherent
left $\cA$\+module.
 The functor $\cA\ot_{\cO_X}{-}\,\:X\Qcoh\rarrow\cA\Qcoh$ is left
adjoint to the forgetful functor $\cA\Qcoh\rarrow X\Qcoh$.

 For a quasi-coherent quasi-algebra $\cA$ over an affine scheme
$U=\Spec R$, the abelian category of quasi-coherent left $\cA$\+modules
is naturally equivalent to the abelian category of left modules over
the ring $A=\cA(U)$, that is $\cA\Qcoh\simeq A\Modl$.
 Lemma~\ref{quasi-algebra-co-extension-of-scalars}(b) is relevant here.

 Denote by $X\QQcoh_\lvfl\subset X\QQcoh$ the full subcategory of
all quasi-coherent quasi-modules that are very flat as quasi-coherent
sheaves on $X$ with respect to the left $\cO_X$\+module structure.
 Similarly, denote by $X\QQcoh_\lfl\subset X\QQcoh$ the full
subcategory of all quasi-coherent quasi-modules that are flat as
quasi-coherent sheaves on $X$ with respect to the left $\cO_X$\+module
structure, and by $X\QQcoh_\rfl\subset X\QQcoh$ the full subcategory
of all quasi-coherent quasi-modules that are flat as quasi-coherent
sheaves on $X$ with respect to the right $\cO_X$\+module structure.
 Clearly, the structure sheaf $\cO_X\in X\Qcoh\subset X\QQcoh$ belongs
to all the three full subcategories $X\QQcoh_\lvfl$, $X\Qcoh_\lfl$,
and $X\Qcoh_\rfl$.

 It is clear from the definitions that the $X\QQcoh_\lfl$ and
$X\QQcoh_\rfl$ are monoidal full subcategories in $X\QQcoh$.
 By Lemma~\ref{very-flat-quasi-module-lemma}(a), \,$X\QQcoh_\lvfl$
is also a monoidal full subcategory in $X\QQcoh$.

 The functor
$$
 \Cohom_X({-},{-})\:(X\QQcoh_\lvfl)^\sop\times X\Lcth_\bW\lrarrow
 X\Lcth_\bW
$$
from Lemma~\ref{cohom-loc-contraherent-lemma}(a) and the subsequent
discussion in Section~\ref{cohom-from-quasi-modules-subsecn} defines
a structure of right module category over $X\QQcoh_\lvfl$ on
the opposite category to the category $X\Lcth_\bW$ of $\bW$\+locally
contraherent cosheaves on~$X$.
 Similarly, the functor
$$
 \Cohom_X({-},{-})\:(X\QQcoh_\lfl)^\sop\times X\Lcth_\bW^\lct\lrarrow
 X\Lcth_\bW^\lct
$$
from Lemma~\ref{cohom-loc-contraherent-lemma}(b) and the subsequent
discussion in Section~\ref{cohom-from-quasi-modules-subsecn} defines
a structure of right module category over $X\QQcoh_\lfl$ on the opposite
category to the category $X\Lcth_\bW^\lct$ of locally cotorsion
$\bW$\+locally contraherent cosheaves on~$X$.
 Finally, the functor
$$
 \Cohom_X({-},{-})\:(X\QQcoh_\rfl)^\sop\times X\Lcth_\bW^\lin\lrarrow
 X\Lcth_\bW^\lin
$$
from Lemma~\ref{cohom-loc-contraherent-lemma}(c) and the subsequent
discussion in Section~\ref{cohom-from-quasi-modules-subsecn} defines
a structure of right module category over $X\QQcoh_\rfl$ on the opposite
category to the category $X\Lcth_\bW^\lin$ of locally injective
$\bW$\+locally contraherent cosheaves on~$X$.

 These assertions are based on the observations that
\begin{itemize}
\item for all $\cA$, $\cB\in X\QQcoh_\lvfl$ and $\P\in X\Lcth_\bW$,
a natural associativity isomorphism $\Cohom_X(\cB,\Cohom_X(\cA,\P))
\simeq\Cohom_X(\cA\ot_{\cO_X}\cB,\>\P)$ holds in $X\Lcth_\bW$;
\item for all $\cA$, $\cB\in X\QQcoh_\lfl$ and $\P\in X\Lcth_\bW^\lct$,
a natural associativity isomorphism $\Cohom_X(\cB,\Cohom_X(\cA,\P))
\simeq\Cohom_X(\cA\ot_{\cO_X}\cB,\>\P)$ holds in $X\Lcth_\bW^\lct$;
\item for all $\cA$, $\cB\in X\QQcoh_\rfl$ and $\gJ\in X\Lcth_\bW^\lin$,
a natural associativity isomorphism $\Cohom_X(\cB,\Cohom_X(\cA,\gJ))
\simeq\Cohom_X(\cA\ot_{\cO_X}\cB,\>\gJ)$ holds in $X\Lcth_\bW^\lin$.
\end{itemize}
 The other relevant observation is that of a natural isomorphism
$\Cohom_X(\cO_X,\P)\simeq\P$ for any $\bW$\+locally contraherent
cosheaf $\P$ on~$X$.

 Let $\cA$ be a quasi-coherent quasi-algebra over $X$ that is very
flat with respect to its left $\cO_X$\+module structure.
 Then $\bW$\+locally contraherent left $\cA$\+modules on $X$ can be
equivalently defined as right module objects in the right module
category $X\Lcth_\bW$ over the monoid object $\cA$ in the monoidal
category $X\QQcoh_\lvfl$.
 So a $\bW$\+locally contraherent left $\cA$\+module $\P$ is the same
thing as a $\bW$\+locally contraherent cosheaf on $X$ endowed with
a \emph{left action} morphism of $\bW$\+locally contraherent cosheaves
$\P\rarrow\Cohom_X(\cA,\P)$ satisfying natural associativity and
unitality axioms together with the multiplication and unit morphisms
of~$\cA$.

 Specifically, the two compositions
$$
 \P\rarrow\Cohom_X(\cA,\P)\rightrightarrows
 \Cohom_X(\cA\ot_{\cO_X}\cA,\>\P)\simeq\Cohom_X(\cA,\Cohom_X(\cA,\P))
$$
of the left action morphism $\P\rarrow\Cohom_X(\cA,\P)$ with
the morphism $\Cohom_X(\cA,\P)\allowbreak\rarrow
\Cohom_X(\cA\ot_{\cO_X}\cA,\>\P)$ induced by the multiplication morphism
$\cA\ot_{\cO_X}\cA\rarrow\cA$ and with the morphism $\Cohom_X(\cA,\P)
\rarrow\Cohom_X(\cA,\Cohom_X(\cA,\P))$ induced by the left action
morphism $\P\rarrow\Cohom_X(\cA,\P)$ must be equal to each other
(associativity).
 The composition
$$
 \P\rarrow\Cohom_X(\cA,\P)\rarrow\P
$$
of the left action morphism $\P\rarrow\Cohom_X(\cA,\P)$ with
the morphism $\Cohom_X(\cA,\P)\allowbreak\rarrow\P$ induced by the unit
morphism $\cO_X\rarrow\cA$ must be equal to the identity endomorphism
of~$\P$ (unitality).

 Similarly, let $\cA$ be a quasi-coherent quasi-algebra over $X$ that
is flat with respect to its left $\cO_X$\+module structure.
 Then $X$\+locally cotorsion $\bW$\+locally contraherent left
$\cA$\+modules on $X$ can be equivalently defined as right module
objects in the right module category $X\Lcth_\bW^\lct$ over the monoid
object $\cA$ in the monoidal category $X\QQcoh_\lfl$.
 So an $X$\+locally cotorsion $\bW$\+locally contraherent left
$\cA$\+module $\P$ is the same thing as a locally cotorsion
$\bW$\+locally contraherent cosheaf on $X$ endowed with a \emph{left
action} morphism of locally cotorsion $\bW$\+locally contraherent
cosheaves $\P\rarrow\Cohom_X(\cA,\P)$ satisfying the associativity
and unitality axioms above.

 Finally, let $\cA$ be a quasi-coherent quasi-algebra over $X$ that
is flat with respect to its right $\cO_X$\+module structure.
 Then $X$\+locally cotorsion $\bW$\+locally contraherent left
$\cA$\+modules on $X$ can be equivalently defined as right module
objects in the right module category $X\Lcth_\bW^\lin$ over the monoid
object $\cA$ in the monoidal category $X\QQcoh_\rfl$.
 So an $X$\+locally injective $\bW$\+locally contraherent left
$\cA$\+module $\gJ$ is the same thing as a locally injective
$\bW$\+locally contraherent cosheaf on $X$ endowed with a \emph{left
action} morphism of locally injective $\bW$\+locally contraherent
cosheaves $\gJ\rarrow\Cohom_X(\cA,\gJ)$ satisfying the same
associativity and unitality axioms.

 Let $\cA$ be a quasi-coherent quasi-algebra over~$X$.
 For any $\bW$\+locally contraherent cosheaf $\Q$ on $X$ such that
a $\bW$\+locally contraherent cosheaf $\Cohom_X(\cA,\Q)$ is
well-defined as per one of the conditions in
Section~\ref{cohom-from-quasi-modules-subsecn}, the $\bW$\+locally
contraherent cosheaf $\Cohom_X(\cA,\Q)$ on $X$ has a natural left
$\cA$\+module structure.
 Under these assumptions, for any $\bW$\+locally contraherent
$\cA$\+module $\P$ on $X$ there is a natural adjunction isomorphism of
abelian groups
\begin{equation} \label{cohom-A-adjunction}
 \Hom^X(\P,\Q)\simeq\Hom^\cA(\P,\Cohom_X(\cA,\Q)).
\end{equation}
 Here we denote by $\Hom^X$ the groups of morphisms in the category
of locally contraherent cosheaves on $X$, and by $\Hom^\cA$
the groups of morphisms in the category of locally contraherent
$\cA$\+modules on~$X$.

 Let $U=\Spec R$ be an affine scheme and $\cA$ be a quasi-coherent
quasi-algebra over~$U$ with the ring of global sections $A=\cA(U)$;
so $A$ is a quasi-algebra over~$R$.
 Then the exact category of contraherent $\cA$\+modules $\cA\Ctrh=
\cA\Lcth_{\{U\}}$ is naturally equivalent to the exact category
$A\Modl^{R\dcta}$ of $R$\+contraadjusted left $A$\+modules (see
Section~\ref{prelim-quasi-algebras-subsecn} for the terminology).
 Here the exact category structure on the full subcategory
$A\Modl^{R\dcta}\subset A\Modl$ is inherited from the abelian
exact structure of the ambient abelian category $A\Modl$.

 Similarly, the exact category of $U$\+locally cotorsion contraherent
$\cA$\+modules $\cA\Ctrh^{U\dlct}=\cA\Lcth_{\{U\}}^{U\dlct}$ is
naturally equivalent to the exact category $A\Modl^{R\dcot}$ of
$R$\+cotorsion left $A$\+modules.
 Once again, the exact category structure on the full subcategory
$A\Modl^{R\dcot}\subset A\Modl$ is inherited from the abelian
exact structure of the ambient abelian category $A\Modl$.

 Finally, the exact category of $U$\+locally injective contraherent
$\cA$\+modules $\cA\Ctrh^{U\dlin}=\cA\Lcth_{\{U\}}^{U\dlin}$ is
naturally equivalent to the exact category $A\Modl^{R\dinj}$ of
$R$\+injective left $A$\+modules.
 As above, the exact category structure on the full subcategory
$A\Modl^{R\dinj}\subset A\Modl$ is inherited from the abelian
exact structure of the ambient abelian category $A\Modl$.
 Lemma~\ref{quasi-algebra-co-extension-of-scalars}(c) is relevant to
these assertions.
{\hbadness=2200\par}

\subsection{$\cA$-locally cotorsion and $\cA$-locally injective
cosheaves of $\cA$-modules}
\label{A-loc-cotors-loc-inj-cosheaves-subsecn}
 Let $X$ be a scheme and $\cA$ be a quasi-coherent quasi-algebra
over~$X$.
 A quasi-coherent left $\cA$\+module $\F$ on $X$ is said to be
\emph{flat} (or \emph{$\cA$\+flat}) if, for every affine open subscheme
$U\subset X$, the left $\cA(U)$\+module $\F(U)$ is flat.
 By Lemma~\ref{quasi-algebra-adjustedness-co-locality}(a), it suffices
to check this condition for affine open subschemes $U$ belonging to
any chosen affine open covering $X=\bigcup_\alpha U_\alpha$ of
the scheme~$X$.
 If $\cA$ is flat as a quasi-coherent sheaf on $X$ with respect to
its left $\cO_X$\+module structure, then any $\cA$\+flat left
$\cA$\+module on $X$ is also flat as a quasi-coherent sheaf of
$\cO_X$\+modules.

 The full subcategory $\cA\Qcoh^{\cA\dfl}$ of $\cA$\+flat
quasi-coherent left $\cA$\+modules is closed under kernels of
epimorphisms, extensions, and infinite direct sums in the abelian
category of quasi-coherent left $\cA$\+modules $\cA\Qcoh$.
 So the full subcategory $\cA\Qcoh^{\cA\dfl}\subset\cA\Qcoh$ inherits
an exact category structure from the abelian exact structure of
$\cA\Qcoh$.

 Let $\bW$ be an open covering of~$X$.
 A $\bW$\+locally contraherent $\cA$\+module $\P$ on $X$ is said to be
\emph{$\cA$\+locally cotorsion} if, for every affine open subscheme
$U\subset X$ subordinate to $\bW$, the left $\cA(U)$\+module $\P[U]$ is
cotorsion.
By Lemma~\ref{quasi-algebra-adjustedness-co-locality}(b), it suffices
to check this condition for affine open subschemes $U$ belonging to
any chosen affine open covering $X=\bigcup_\alpha U_\alpha$ of
the scheme $X$ subordinate to~$\bW$.
 By Lemma~\ref{restriction-coextension-injective-cotorsion}(a), any
$\cA$\+locally cotorsion $\bW$\+locally contraherent $\cA$\+module
on $X$ is $X$\+locally cotorsion.

 The full subcategory of $\cA\Lcth_\bW^{\cA\dlct}$ of $\cA$\+locally
cotorsion $\bW$\+locally contraherent $\cA$\+modules is closed under
cokernels of admissible monomorphisms, extensions, and infinite products
in the exact category of $\bW$\+locally contraherent $\cA$\+modules
$\cA\Lcth_\bW$ (or in the exact category of $X$\+locally cotorsion
$\bW$\+locally contraherent $\cA$\+modules $\cA\Lcth_\bW^{X\dlct}$).
 See, in particular, Lemma~\ref{flat-cotorsion-pair-hereditary}(b).
 So the full subcategory $\cA\Lcth_\bW^{\cA\dlct}\subset
\cA\Lcth_\bW$ inherits an exact category structure from $\cA\Lcth_\bW$.

 A $\bW$\+locally contraherent $\cA$\+module $\gJ$ on $X$ is said to be
\emph{$\cA$\+locally injective} if, for every affine open subscheme
$U\subset X$ subordinate to $\bW$, the left $\cA(U)$\+module $\gJ[U]$ is
injective.
 By Lemma~\ref{quasi-algebra-adjustedness-co-locality}(c), it suffices
to check this condition for affine open subschemes $U$ belonging to
any chosen affine open covering $X=\bigcup_\alpha U_\alpha$ of
the scheme $X$ subordinate to~$\bW$.
 If $\cA$ is flat as a quasi-coherent sheaf on $X$ with respect to
its right $\cO_X$\+module structure then,
by Lemma~\ref{restriction-coextension-injective-cotorsion}(b), any
$\cA$\+locally injective $\bW$\+locally contraherent $\cA$\+module
on $X$ is $X$\+locally injective.

 The full subcategory of $\cA\Lcth_\bW^{\cA\dlin}$ of $\cA$\+locally
injective $\bW$\+locally contraherent $\cA$\+modules is closed under
cokernels of admissible monomorphisms, extensions, and infinite products
in the exact category of $\bW$\+locally contraherent $\cA$\+modules
$\cA\Lcth_\bW$.
 So the full subcategory $\cA\Lcth_\bW^{\cA\dlin}\subset
\cA\Lcth_\bW$ inherits an exact category structure from $\cA\Lcth_\bW$.

 As usual, we put $\cA\Ctrh^{\cA\dlct}=\cA\Lcth_{\{X\}}^{\cA\dlct}$
and $\cA\Ctrh^{\cA\dlin}=\cA\Lcth_{\{X\}}^{\cA\dlin}$, as well as
$\cA\Lcth^{\cA\dlct}=\bigcup_\bW\cA\Lcth_\bW^{\cA\dlct}$ and
$\cA\Lcth^{\cA\dlin}=\bigcup_\bW\cA\Lcth_\bW^{\cA\dlin}$.

 Let $U=\Spec R$ and $\cA$ be a quasi-coherent quasi-algebra over~$U$
with the ring of global sections $A=\cA(U)$.
 Then the exact category of $\cA$\+locally cotorsion contraherent
$\cA$\+modules $\cA\Ctrh^{\cA\dlct}=\cA\Lcth_{\{U\}}$ is naturally
equivalent to the exact category $A\Modl^\cot$ of cotorsion left
$A$\+modules.
 Here the exact category structure on the full subcategory
$A\Modl^\cot\subset A\Modl$ is inherited from the abelian
exact structure of the ambient abelian category $A\Modl$.

 The exact category of $\cA$\+locally injective contraherent
$\cA$\+modules $\cA\Ctrh^{\cA\dlin}=\cA\Lcth_{\{U\}}^{\cA\dlin}$ over
an affine scheme $U$ has a split exact structure (all short exact
sequences are split).
 The additive category $\cA\Ctrh^{\cA\dlin}$ is naturally equivalent
to the additive category $A\Modl^\inj$ of injective left $A$\+modules.

\subsection{Direct images} \label{direct-images-of-A-co-sheaves-subsecn}
 Let us introduce one extra piece of terminology:
 A \emph{quasi-ringed scheme} $(X,\cA)$ is a pair consisting of a scheme
$X$ and a quasi-coherent quasi-algebra $\cA$ over~$X$.
 Given a quasi-ringed scheme $(Y,\cB)$ and a morphism of schemes
$f\:Y\rarrow X$, the direct image $f_*\cB$ is, generally speaking,
just a sheaf of rings on $X$ endowed with a morphism of sheaves
of rings $f_*\cO_Y\rarrow f_*\cB$.
 A \emph{morphism of quasi-ringed schemes} $f\:(Y,\cB)\rarrow(X,\cA)$
is a pair consisting of a morphism of schemes $f\:Y\rarrow X$ and
a morphism of sheaves of rings $\cA\rarrow f_*\cB$ on $X$ making
the square diagram of morphisms of sheaves of rings on~$X$
\begin{equation} \label{morphism-of-ringed-schemes-diagram}
\begin{gathered}
 \xymatrix{
  f_*\cO_Y \ar[r] & f_*\cB \\
  \cO_X \ar[r] \ar[u] & \cA \ar[u]
 }
\end{gathered}
\end{equation}
commutative.

 When the morphism of schemes $f\:Y\rarrow X$ is quasi-compact and
quasi-separated, the direct image $f_*\cB$ is a quasi-coherent
quasi-module on $X$ by
Lemma~\ref{direct-image-quasi-coherent-quasi-module}.
 On the other hand, $f_*\cB$ is a sheaf of rings on $X$ endowed with
a morphism of sheaves of rings $\cO_X\rarrow f_*\cB$ (constructed
as the composition $\cO_X\rarrow f_*\cO_Y\rarrow f_*\cB$).
 One can easily check that the two resulting structures of a sheaf
of $\cO_X$\+$\cO_X$\+bimodules on $f_*\cB$ agree; so $f_*\cB$ is
a quasi-coherent quasi-algebra over~$X$.
 Thus, in the context of a morphism of quasi-ringed schemes
$f\:(Y,\cB)\rarrow(X,\cA)$, the morphism of sheaves of rings
$\cA\rarrow f_*\cB$ is a morphism of quasi-coherent quasi-algebras
over $X$ whenever $f\:X\rarrow Y$ is quasi-compact and quasi-separated.

 More generally, let $f\:Y\rarrow X$ be a morphism of schemes, $\cA$ be
a sheaf of rings on $X$ endowed with a morphism of sheaves of rings
$\cO_X\rarrow\cA$, and $\cB$ be a sheaf of rings on $Y$ endowed with
a morphism of sheaves of rings $\cO_Y\rarrow\cB$.
 Assume that we are given a morphism of sheaves of rings
$\cA\rarrow f_*\cB$ on $X$ making the square
diagram~\eqref{morphism-of-ringed-schemes-diagram} commutative.

 Assume that the morphism $f\:Y\rarrow X$ is quasi-compact and
quasi-separated.
 Let $\N$ be a quasi-coherent left $\cB$\+module on~$Y$.
 Then the direct image $f_*\N$ has a natural structure of a sheaf
of $f_*\cB$\+modules on $X$ whose underlying sheaf of $\cO_X$\+modules
is quasi-coherent.
 Using the morphism of sheaves of rings $\cA\rarrow f_*\cB$, one defines
the induced structure of a quasi-coherent left $\cA$\+module on~$f_*\N$.
 So we have the direct image functor
\begin{equation} \label{qcoh-A-modules-direct-image}
 f_*\:\cB\Qcoh\lrarrow\cA\Qcoh.
\end{equation}

 Let $\bW$ be an open covering of $X$ and $\bT$ be an open covering
of $Y$ such that the morphism of schemes $f\:Y\rarrow X$ is
$(\bW,\bT)$\+affine (in the sense of
Section~\ref{direct-images-of-O-co-sheaves-subsecn}).
 Let $\Q$ be a $\bT$\+locally contraherent $\cA$\+module on~$Y$.
 Then the direct image $f_!\Q$ is a $\bW$\+locally contraherent
cosheaf on $X$ with a natural structure of a locally contraherent
cosheaf of $f_*\cB$\+modules.
 Using the morphism of sheaves of rings $\cA\rarrow f_*\cB$, one defines
the induced structure of a $\bW$\+locally contraherent $\cA$\+module
on~$f_!\Q$.
 So we obtain the direct image functor
$$
 f_!\:\cB\Lcth_\bT\lrarrow\cA\Lcth_\bW,
$$
which restricts to a direct image functor between the categories
of $Y$\+locally cotorsion $\bT$\+locally contraherent $\cB$\+modules
and $X$\+locally cotorsion $\bW$\+locally contraherent $\cA$\+modules
$$
 f_!\:\cB\Lcth_\bT^{Y\dlct}\lrarrow\cA\Lcth_\bW^{X\dlct}.
$$

 Assuming additionally that $f\:Y\rarrow X$ is a flat morphism of
schemes, the two direct image functors~$f_!$ above restrict to a direct
image functor between the categories of $Y$\+locally injective
$\bT$\+locally contraherent $\cB$\+modules and $X$\+locally injective
$\bW$\+locally contraherent $\cA$\+modules,
$$
 f_!\:\cB\Lcth_\bT^{Y\dlin}\lrarrow\cA\Lcth_\bW^{X\dlin}.
$$

 In particular, let $j\:Y\rarrow X$ be an open immersion of schemes,
and let $\cA$ be a sheaf of rings on $X$ endowed with a morphism of
sheaves of rings $\cO_X\rarrow\cA$.
 Put $\cB=\cA|_Y$; then there are natural morphisms of sheaves of
rings $\cO_Y\rarrow\cB$ on $Y$ and $\cA\rarrow j_*\cB$ on $X$
making the square diagram~\eqref{morphism-of-ringed-schemes-diagram}
commutative (for $f=j$).
 Given an open covering $\bW$ of $X$, consider the open covering
$\bW|_Y$ of $Y$ as at the end of
Section~\ref{inverse-images-of-O-co-sheaves-subsecn};
then $j\:Y\rarrow X$ is a $(\bW,\bW|_Y)$\+coaffine morphism of schemes.
 If the open immersion morphism~$j$ is affine, then it is also
a $(\bW,\bW|_Y)$\+affine morphism of schemes.

 In this context, for any quasi-coherent left $\cA$\+module $\M$ on $X$,
the restriction $j^*\M=\M|_Y$ of $M$ to $Y$ is naturally
a quasi-coherent left $\cB$\+module on~$Y$; and for any $\bW$\+locally
contraherent $\cA$\+module $\P$ on $X$, the restriction
$j^!\P=\P|_Y$ of $\P$ to $Y$ is naturally a $\bW|_Y$\+locally
contraherent $\cB$\+module on~$Y$.
 Assuming that the open immersion morphism~$j$ is affine,
as a special case of the adjunction
isomorphism~\eqref{open-immers-direct-inverse-cosheaf-adjunction},
we obtain an adjunction isomorphism
\begin{equation} \label{affine-open-immers-contrah-module-adjunction}
 \Hom^\cA(j_!\Q,\P)\simeq\Hom^\cB(\Q,\P|_Y)
\end{equation}
for any $\bW|_Y$\+locally contraherent $\cB$\+module $\Q$ on $Y$ and
any $\bW$\+locally contraherent $\cA$\+module $\P$ on~$X$.
 Here we denote by $\Hom^\cA$ and $\Hom^\cB$ the groups of morphisms in
the categories of locally contraherent cosheaves of $\cA$\+modules and
$\cB$\+modules on $X$ and~$Y$.

 In this paper, we are mostly interested in the direct images
of quasi-coherent and locally contraherent $\cB$\+modules in
the context of the previous two paragraphs (i.~e., for an affine
open immersion morphism $j\:Y\rarrow X$ and the sheaf of rings
$\cB=\cA|_Y$).
 Moreover, we will almost always assume that $\cA$ is a quasi-coherent
quasi-algebra over~$X$ (hence $\cB$ is a quasi-coherent quasi-algebra
over~$Y$).

 Now let $f\:(Y,\cB)\rarrow(X,\cA)$ be a morphism of quasi-ringed
schemes (so $\cB$ is a quasi-coherent quasi-algebra over $X$ and
$\cA$ is a quasi-coherent quasi-algebra over~$Y$).
 Let us say that $f$~is \emph{right flat} (as a morphism of
quasi-ringed schemes) if, for every pair of affine open subschemes
$U\subset X$ and $V\subset Y$ such that $f(V)\subset U$, the ring
homomorphism $\cA(U)\rarrow\cB(V)$ makes $\cB(V)$ a flat right
$\cA(U)$\+module.
 It suffices to check this condition for affine open subschemes $U$
belonging to any chosen affine open covering $X=\bigcup_\alpha U_\alpha$
of the scheme $X$ and affine open subschemes $V$ belonging to any
chosen affine open covering $f^{-1}(U)=\bigcup_\beta V_\beta$ of
the open subscheme $f^{-1}(U)\subset Y$ (by
Lemmas~\ref{quasi-algebra-adjustedness-co-locality}(a)
and~\ref{quasi-algebra-ring-flatness-locality}).

 Let $f\:(Y,\cB)\rarrow(X,\cA)$ be a morphism of quasi-ringed schemes,
and let $\bW$ and $\bT$ be open coverings of $X$ and $Y$ such that
the morphism of schemes $f\:Y\rarrow X$ is $(\bW,\bT)$\+affine.
 Then the direct image functor $f_!\:\cB\Lcth_\bT^{Y\dlct}\rarrow
\cA\Lcth_\bW^{X\dlct}$ takes $\cB$\+locally cotorsion $\bT$\+locally
contraherent cosheaves on $Y$ to $\cA$\+locally cotorsion
$\bW$\+locally contraherent cosheaves on~$X$ (by
Lemma~\ref{restriction-coextension-injective-cotorsion}(a)).
 So we have the direct image functor
$$
 f_!\:\cB\Lcth_\bT^{\cB\dlct}\lrarrow\cA\Lcth_\bW^{\cA\dlct}.
$$

 Assume that the morphism of quasi-ringed schemes
$f\:(Y,\cB)\rarrow(X,\cA)$ is right flat.
 Then the direct image functor $f_!\:\cB\Lcth_\bT\rarrow\cA\Lcth_\bW$
takes $\cB$\+locally injective $\bT$\+locally contraherent cosheaves
on $Y$ to $\cA$\+locally injective $\bW$\+locally contraherent cosheaves
on~$X$ (by Lemma~\ref{restriction-coextension-injective-cotorsion}(b)).
 So we obtain the direct image functor
$$
 f_!\:\cB\Lcth_\bT^{\cB\dlin}\lrarrow\cA\Lcth_\bW^{\cA\dlin}.
$$

\subsection{Inverse images}
\label{inverse-images-of-A-co-sheaves-subsecn}
 The definition of a right flat morphism of quasi-ringed schemes
$f\:(Y,\cB)\rarrow(X,\cA)$ was given in the previous
Section~\ref{direct-images-of-A-co-sheaves-subsecn}.
 The definition of a \emph{left flat} morphism of quasi-ringed schemes
is similar.

 Given a morphism of quasi-ringed schemes $f\:(Y,\cB)\rarrow(X,\cA)$
consider the topology base $\bB$ of $Y$ consisting of all the affine
open subschemes $V\subset Y$ for which there exists an affine open
subscheme $U\subset X$ such that $f(V)\subset U$.
 For any quasi-coherent $\cA$\+module $\M$ on $X$, one defines
the presheaf of $\cB$\+modules $f^*\M$ on $\bB$ by the rule
$$
 (f^*\M)(V)=\cB(V)\ot_{\cA(U)}\M(U)
$$
for any affine open subscheme $V\subset Y$, \ $V\in\bB$ and any chosen
affine open subscheme $U\subset X$ such that $f(V)\subset U$.

 The tricky part of this definition is to check that
the $\cB(V)$\+module $\cB(V)\ot_{\cA(U)}\M(U)$ does not depend on
the choice of an affine open subscheme $U\subset X$.
 For the beginning, assume that the scheme $X$ is \emph{semi-separated},
i.~e., finite intersections of affine open subschemes in $X$
are affine.
 Then, for any two affine open subschemes $U'\subset X$ and
$U''\subset X$ such that $f(V)\subset U'$ and $f(V)\subset U''$,
one has $f(V)\subset U'\cap U''$, where $U'\cap U''$ is a smaller
affine open subscheme in~$X$.
 It remains to compute that, given a pair of affine open subschemes
$U\subset U'\subset X$ such that $f(V)\subset U$, one has
$$
 \cB(V)\ot_{\cA(U)}\M(U)\simeq\cB(V)\ot_{\cA(U)}
 (\cA(U)\ot_{\cA(U')}\M(U'))\simeq\cB(V)\ot_{\cA(U')}\M(U')
$$
because
\begin{multline*}
 \M(U)\simeq\cO_X(U)\ot_{\cO_X(U')}\M(U') \\ \simeq
 (\cO_X(U)\ot_{\cO_X(U')}\cA(U'))\ot_{\cA(U')}\M(U')\simeq
 \cA(U)\ot_{\cA(U')}\M(U).
\end{multline*}

 To check the quasi-coherence axiom from
Section~\ref{locally-contraherent-cosheaves-subsecn} for
the presheaf of $\cB$\+modules $f^*\M$ on $\bB$, one computes that,
for any affine open subschemes $V\subset V'\subset Y$ and any affine
open subscheme $U\subset X$ such that $f(V')\subset U$,
\begin{multline*}
 \cO_Y(V)\ot_{\cO_Y(V')}(f^*\M)(V')=
 \cO_Y(V)\ot_{\cO_Y(V')}(\cB(V')\ot_{\cA(U)}\M(U)) \\
 \simeq(\cO_Y(V)\ot_{\cO_Y(V')}\cB(V'))\ot_{\cA(U)}\M(U)
 \simeq\cB(V)\ot_{\cA(U)}\M(U)=(f^*\M)(V).
\end{multline*}
 Now, by Theorem~\ref{extension-of-co-sheaves-from-topology-base}(a)
and Lemma~\ref{quasi-coherence-implies-sheaf}, the quasi-coherent
presheaf of $\cB$\+modules $f^*\M$ on $\bB$ extends uniquely to
a quasi-coherent sheaf of $\cB$\+modules on $Y$, which we also denote
by~$f^*\M$.

 It remains to explain what to do if the scheme $X$ is not
semi-separated.
 In this case, one just covers $X$ by semi-separated (or affine)
open subschemes $X_\alpha$, puts $Y_\alpha=f^{-1}(X_\alpha)$,
and restricts the quasi-coherent quasi-algebras $\cA$ and $\cB$ to
$X_\alpha$ and $Y_\alpha$, respectively.
 Furthermore, one also restricts the quasi-coherent $\cA$\+module $\M$
to the open subschemes $X_\alpha\subset X$, obtaining quasi-coherent
$\cA|_{X_\alpha}$\+modules $\M|_{X_\alpha}$.
 Then one takes the inverse images of the quasi-coherent
$\cA|_{X_\alpha}$\+modules $\M|_{X_\alpha}$ along the (obvious)
morphisms of quasi-ringed schemes $f_\alpha\:(Y_\alpha,\cB|_{Y_\alpha})
\rarrow(X_\alpha,\cA|_{X_\alpha})$.
 Finally, one observes that the resulting quasi-coherent
$\cB|_{Y_\alpha}$\+modules $f_\alpha^*(\M|_{X_\alpha})$ on $Y_\alpha$
agree on the intersections $Y_\alpha\cap Y_\beta$ (e.~g., since the open
subschemes $X_\alpha\cap X_\beta\subset X$ are also semi-separated); so
(e.~g., by Theorem~\ref{extension-of-co-sheaves-from-topology-base}(a)) 
they can be glued together to obtain the desired quasi-coherent
$\cB$\+module $f^*\M$ on~$Y$.
 (Cf.\ the discussion of the non-semi-separatedness issue
in~\cite[Section~3.3]{Pcosh}.)

 We have constructed the inverse image functor
\begin{equation} \label{qcoh-A-modules-inverse-image}
 f^*\:\cA\Qcoh\lrarrow\cB\Qcoh.
\end{equation}
 The functor~$f^*$ is right exact and preserves infinite direct sums.
 When the morphism of quasi-ringed schemes $f\:(Y,\cB)\rarrow(X,\cA)$
is right flat, the functor~$f^*$ is exact.
 When the morphism of schemes $f\:Y\rarrow X$ is quasi-compact and
quasi-separated, the functor~$f^*$ is left adjoint to the direct
image functor~$f_*$~\eqref{qcoh-A-modules-direct-image}.

 Having spelled out the quasi-coherent case in terms instructive for
our purposes, let us pass to the locally contraherent context, which is
only slightly more complicated.
 Let $f\:(Y,\cB)\rarrow(X,\cA)$ be a left flat morphism of quasi-ringed
schemes.
 Let $\bW$ and $\bT$ be open coverings of $X$ and $Y$ such that
the morphism of schemes $f\:Y\rarrow X$ is $(\bW,\bT)$\+coaffine.
 Consider the topology base $\bB$ of $Y$ consisting of all the affine
open subschemes $V\subset Y$ subordinate to~$\bT$.

 Given an $\cA$\+locally cotorsion $\bW$\+locally contraherent
$\cA$\+module $\P$ on $X$, we define the copresheaf of $\cB$\+modules
$f^!\P$ on $\bB$ by the rule
$$
 (f^!\P)[V]=\Hom_{\cA(U)}(\cB(V),\P[U])
$$
for any affine open subscheme $V\subset Y$ subordinate to $\bT$ and
any chosen affine open subscheme $U\subset X$ subordinate to $\bW$
such that $f(V)\subset U$.
 Notice that the $\cB(V)$\+module $\Hom_{\cA(U)}(\cB(V),\P[U])$ is
cotorsion by Lemma~\ref{restriction-coextension-injective-cotorsion}(c).

 Once again, we assume that the scheme $X$ is semi-separated for
the beginning, and check that the $\cB(V)$\+module
$\Hom_{\cA(U)}(\cB(V),\P[U])$ does not depend on the choice of
an affine open subscheme $U\subset X$.
 Indeed, given a pair of affine open subschemes $U\subset U'\subset X$
such that $f(V)\subset U$, one has
\begin{multline*}
 \Hom_{\cA(U)}(\cB(V),\P[U])\simeq
 \Hom_{\cA(U)}(\cB(V),\Hom_{\cA(U')}(\cA(U),P[U'])) \\
 \simeq\Hom_{\cA(U')}(\cB(V),\P[U'])
\end{multline*}
because
\begin{multline*}
 \P[U]\simeq\Hom_{\cO_X(U')}(\cO_X(U),\P[U']) \\
 \simeq\Hom_{\cA(U')}(\cA(U')\ot_{\cO_X(U')}\cO_X(U),\>\P[U'])
 \simeq\Hom_{\cA(U')}(\cA(U),\P[U']).
\end{multline*}

 Let us check that the contraherence axiom~(i) from
Section~\ref{locally-contraherent-cosheaves-subsecn} for
the copresheaf of $\cB$\+modules $f^!\P$ on~$\bB$.
 Indeed, for any affine open subschemes $V\subset V'\subset Y$ and any
affine open subscheme $U\subset X$ such that $f(V')\subset U$, we have
\begin{multline*}
 \Hom_{\cO_Y(V')}(\cO_Y(V),(f^!\P)[V'])=
 \Hom_{\cO_Y(V')}(\cO_Y(V),\Hom_{\cA(U)}(\cB(V'),\P[U])) \\
 \simeq\Hom_{\cA(U)}(\cB(V')\ot_{\cO_Y(V')}\cO_Y(V),\>\P[U])
 \simeq\Hom_{\cA(U)}(\cB(V),\P[U])=(f^!\P)[V].
\end{multline*}

 Furthermore, for any $V\in\bB$, the module $(f^!\P)[V]$ is cotorsion
over $\cB(V)$ (as mentioned above), hence also cotorsion over $\cO_Y(V)$
(by Lemma~\ref{restriction-coextension-injective-cotorsion}(a)); so
the contraadjustedness axiom~(ii) holds (for the underlying copresheaf
of $\cO_Y$\+modules of the copresheaf of $\cB$\+modules $f^!\P$)
as well.
 Now, by Theorem~\ref{extension-of-co-sheaves-from-topology-base}(b)
and Lemma~\ref{contraherence+contraadjustedness-imply-cosheaf},
the contraherent copresheaf of $\cB$\+modules $f^!\P$ on $\bB$ extends
uniquely to a $\cB$\+locally cotorsion $\bT$\+locally contraherent
cosheaf of $\cB$\+modules on $Y$, which we also denote by~$f^!\P$.

 The case of a non-semi-separated scheme $X$ is dealt with similarly
to the way it is treated in the discussion above in the context of
quasi-coherent sheaves; see also~\cite[Section~3.3]{Pcosh}.

 We have constructed the inverse image functor
\begin{equation} \label{ctrh-lct-A-modules-inverse-image}
 f^!\:\cA\Lcth_\bW^{\cA\dlct}\lrarrow\cB\Lcth_\bT^{\cB\dlct}.
\end{equation}
 The functor~$f^!$ is exact (as a functor of exact categories) and
preserves infinite products.

 Now let $f\:(Y,\cB)\rarrow(X,\cA)$ be an arbitrary morphism of
quasi-ringed schemes.
 Then, given an $\cA$\+locally injective $\bW$\+locally contraherent
$\cA$\+module $\gJ$ on $X$, we define the copresheaf of $\cB$\+modules
$f^!\gJ$ on $\bB$ by the same rule
$$
 (f^!\gJ)[V]=\Hom_{\cA(U)}(\cB(V),\gJ[U])
$$
for $V\in\bB$ and $U$ as above.
 The $\cB(V)$\+module $\Hom_{\cA(U)}(\cB(V),\gJ[U])$ is injective by
by Lemma~\ref{restriction-coextension-injective-cotorsion}(d),
hence cotorsion by the definition; hence it is also cotorsion
as an $\cO_X(V)$\+module by
Lemma~\ref{restriction-coextension-injective-cotorsion}(a).

 So the contraadjustedness axiom~(ii) from
Section~\ref{locally-contraherent-cosheaves-subsecn} holds for
(the underlying copresheaf of $\cO_Y$\+modules of) the copresheaf
of $\cB$\+modules $f^!\gJ$ on~$\bB$.
 The rest of the arguments, proving that the $\cB(V)$\+module
$(f^!\gJ)[V]$ is well-defined and the contraherence axiom~(i)
is satisfied for $f^!\gJ$, is the same as above.
 Hence the contraherent copresheaf of $\cB$\+modules $f^!\gJ$ on $\bB$
extends uniquely to a $\cB$\+locally injective $\bT$\+locally
contraherent cosheaf of $\cB$\+modules on $Y$, which we also denote
by~$f^!\gJ$.

 We have constructed the inverse image functor
\begin{equation} \label{ctrh-lin-A-modules-inverse-image}
 f^!\:\cA\Lcth_\bW^{\cA\dlin}\lrarrow\cB\Lcth_\bT^{\cB\dlin}.
\end{equation}
 The functor~$f^!$ is exact and preserves infinite products.

 Let $f\:(Y,\cB)\rarrow(X,\cA)$ be a morphism of quasi-ringed schemes,
and let $\P$ be locally contraherent $\cA$\+module on $X$ such that
the locally contraherent $\cB$\+module $f^!\P$ on $Y$ is well-defined
as per one of the two constructions above.
 Let $\Q$ be a cosheaf of $\cB$\+modules on~$Y$.
 Then there is a natural (adjunction) isomorphism of abelian groups
\begin{equation} \label{direct-inverse-A-module-general-adjunction}
 \Hom^\cA(f_!\Q,\P)\simeq\Hom^\cB(\Q,f^!\P).
\end{equation}
 Here we denote by $\Hom^\cA$ and $\Hom^\cB$ the groups of morphisms
in the categories of cosheaves of $\cA$\+modules on $X$ and cosheaves
of $\cB$\+modules on~$Y$.

 Indeed, let $\bW$ be an open covering of $X$ such that the cosheaf
$\P$ is $\bW$\+locally contraherent, and let $\bT$ be an open covering
of $Y$ such that the morphism~$f$ is $(\bW,\bT$)\+coaffine (so
the cosheaf $f^!\P$ is $\bT$\+locally contraherent).
 Then both the abelian groups
in~\eqref{direct-inverse-A-module-general-adjunction} are naturally
isomorphic to the group of all compatible systems of $\cA(U)$\+module
maps $\Q[V]\rarrow\P[U]$ defined for affine open subschemes $U\subset X$
subordinate to $\bW$ and $V\subset Y$ subordinate to $\bT$ such that
$f(V)\subset U$.
 Here the $\cA(U)$\+module structure on $\Q[V]$ is obtained from
the $\cB(V)$\+module structure by restriction of scalars via
the ring homomorphism $\cA(U)\rarrow\cB(V)$.

 Given a quasi-ringed scheme $(X,\cA)$ and an open subscheme
$Y\subset X$ with the open immersion morphism $j\:Y\rarrow X$,
put $\cB=j^*\cA=\cA|_Y$.
 Then $(Y,\cB)$ is a quasi-ringed scheme and $j\:(Y,\cB)\rarrow(X,\cA)$
is naturally a left and right flat morphism of quasi-ringed schemes.
 Given an open covering $\bW$ of $X$, consider the open covering
$\bW|_Y$ of $Y$ as at the end of
Section~\ref{inverse-images-of-O-co-sheaves-subsecn} and in
Section~\ref{direct-images-of-A-co-sheaves-subsecn}.
 In this context, the inverse image
functor~$j^*$~\eqref{qcoh-A-modules-inverse-image} defined in
this section agrees with the restriction functor $\M\longmapsto\M|_Y$
discussed in Section~\ref{direct-images-of-A-co-sheaves-subsecn}.
 The inverse image functors~$j^!$
\eqref{ctrh-lct-A-modules-inverse-image}
and~\eqref{ctrh-lin-A-modules-inverse-image} defined in this section
agree with the restriction functor $\P\longmapsto\P|_Y$ from
Section~\ref{direct-images-of-A-co-sheaves-subsecn} wherever
the former are defined.

 In this paper, we are mostly interested in the inverse images of
quasi-coherent and locally contraherent $\cA$\+modules in the context
of the previous paragraph (i.~e., for an open immersion morphism
$j\:Y\rarrow X$ and the quasi-coherent quasi-algebra
$\cB=j^*\cA=\cA|_Y$).
 We will use the convenient notation $j^*$ and~$j^!$ for the functors
of restriction of quasi-coherent and locally contraherent
$\cA$\+modules to the open subscheme $Y\subset X$ (even though
the more restrictive assumptions on a locally contraherent
$\cA$\+module $\P$ required for the general definition of the inverse
image $f^!\P$ above in this section may be not satisfied).
 So we just put $j^*\M=\M|_Y$ and $j^!\P=\P|_Y$, by abuse of notation.

 Let us empasize, however, that, generally speaking, for a morphism
of quasi-ringed schemes $f\:(Y,\cB)\rarrow(X,\cA)$, the inverse
image functor~\eqref{qcoh-A-modules-inverse-image} from this section
does \emph{not} agree with the functor of inverse image of
quasi-coherent sheaves of $\cO$\+modules $f^*\:X\Qcoh\rarrow Y\Qcoh$.
 Similarly, the inverse image
functors~\eqref{ctrh-lct-A-modules-inverse-image}
and~\eqref{ctrh-lin-A-modules-inverse-image} from this section
do \emph{not} agree with the functors of inverse image of
locally contraherent cosheaves of $\cO$\+modules from
Section~\ref{inverse-images-of-O-co-sheaves-subsecn}.
 On ther other hand, the direct image functors $f_*$ and~$f_!$
from Section~\ref{direct-images-of-A-co-sheaves-subsecn}, by
the definition, \emph{do} agree with the direct image functors
from Section~\ref{direct-images-of-O-co-sheaves-subsecn}.

\subsection{\v Cech (co)resolutions} \label{cech-subsecn}
 Let $X$ be a quasi-compact semi-separated scheme and
$X=\bigcup_{\alpha=1}^N U_\alpha$ be a finite affine open covering
of~$X$.
 For every subsequence of indices $1\le\alpha_1<\dotsb<\alpha_k\le N$,
denote by $j_{\alpha_1,\dotsc,\alpha_k}\:\bigcap_{s=1}^kU_{\alpha_s}
\rarrow X$ the open immersion morphism.
 Let $\cA$ be a sheaf of rings on $X$ endowed with a morphism of
sheaves of rings $\cO_X\rarrow\cA$.

 Let $\M$ be a quasi-coherent $\cA$\+module on~$X$.
 Then the \v Cech coresolution
\begin{multline} \label{qcoh-sheaf-of-rings-cech-coresolution}
 0\lrarrow\M\lrarrow\bigoplus\nolimits_{\alpha=1}^N
 j_\alpha{}_*j_\alpha^*\M\lrarrow
 \bigoplus\nolimits_{1\le\alpha<\beta\le N}
 j_{\alpha,\beta}{}_*j_{\alpha,\beta}^*\M \\
 \lrarrow\dotsb\lrarrow j_{1,\dotsc,N}{}_*j_{1,\dotsc,N}^*\M
 \lrarrow0
\end{multline}
is a finite exact sequence of quasi-coherent $\cA$\+modules on~$X$.
 Here, as suggested in the previous
Section~\ref{inverse-images-of-A-co-sheaves-subsecn}, we just use
the symbol~$j^*$ as a simplified notation for the functor
$\M\longmapsto\M|_Y\:\cA\Qcoh\rarrow\cA|_Y\Qcoh$ of restriction of
quasi-coherent sheaves of $\cA$\+modules to an open subscheme
$Y\subset X$ with the open immersion morphism $j\:Y\rarrow X$.

 Now let $\bW$ be an open covering of~$X$, and let
$X=\bigcup_{\alpha=1}^N U_\alpha$ be a finite affine open covering
subordinate to~$\bW$.
 Let $\P$ be a $\bW$\+locally contraherent $\cA$\+module on~$X$.
 Then the \v Cech resolution~\cite[formula~(3.6) in Section~3.3]{Pcosh}
\begin{multline} \label{lcth-sheaf-of-rings-cech-resolution}
 0\lrarrow j_{1,\dotsc,N}{}_!j_{1,\dotsc,N}^!\P \lrarrow\dotsb \\
 \lrarrow\bigoplus\nolimits_{1\le\alpha<\beta\le N}
 j_{\alpha,\beta}{}_!j_{\alpha,\beta}^!\P\lrarrow
 \bigoplus\nolimits_{\alpha=1}^N j_\alpha{}_!j_\alpha^!\P
 \lrarrow\P\lrarrow0
\end{multline}
is a finite exact sequence of $\bW$\+locally contraherent
$\cA$\+modules on~$X$.
 In fact, \eqref{lcth-sheaf-of-rings-cech-resolution}~is a finite
resolution of $\P$ by (globally) contraherent $\cA$\+modules on~$X$.
 Once again, we use the symbol~$j^!$ as simplified notation for
the functor $\P\longmapsto\P|_Y\:\cA\Lcth_\bW\rarrow
\cA|_Y\Lcth_{\bW|_Y}$ of restriction of locally contraherent cosheaves
of $\cA$\+modules to an open subscheme $Y\subset X$ with the open
immersion morphism $j\:Y\rarrow X$.

 So~\eqref{lcth-sheaf-of-rings-cech-resolution} is a finite exact
sequence in the exact category $\cA\Lcth_\bW$.
 Whenever the $\bW$\+locally contraherent $\cA$\+module $\P$ belongs
to one of the exact subcategories $\cA\Lcth_\bW^{X\dlct}$,
$\cA\Lcth_\bW^{X\dlin}$, $\cA\Lcth_\bW^{\cA\dlct}$, or
$\cA\Lcth_\bW^{\cA\dlin}\subset\cA\Lcth_\bW$, the
\v Cech resolution~\eqref{lcth-sheaf-of-rings-cech-resolution}
is a finite exact sequence in the respective exact category.

\subsection{Projection formulas}
 Let $X$ be a scheme and $Y\subset X$ be an open subscheme such that
the open immersion morphism $j\:Y\rarrow X$ is affine.
 Then, for any quasi-coherent quasi-module $\cA$ on $X$ and any
quasi-coherent sheaf $\N$ on $Y$ there is a natural isomorphism
\begin{equation} \label{quasi-module-tensor-projection-formula}
 \cA\ot_{\cO_X}j_*\N\simeq j_*(j^*\cA\ot_{\cO_Y}\N)
\end{equation}
of quasi-coherent sheaves on~$X$.
 If $\cA$ is a quasi-coherent quasi-algebra on $X$,
then~\eqref{quasi-module-tensor-projection-formula} is an isomorphism
of quasi-coherent left $\cA$\+modules on~$X$.

 Let $\bW$ be an open covering of $X$ and $\bW|_Y$ be the restriction
of $\bW$ to~$Y$, as in
Sections~\ref{inverse-images-of-O-co-sheaves-subsecn}
and~\ref{direct-images-of-A-co-sheaves-subsecn}.
 Let $\cA$ be a quasi-coherent quasi-module on $X$ that is very flat
as a quasi-coherent sheaf in the left $\cO_X$\+module structure,
and let $\Q$ be a $\bW|_Y$\+locally contraherent cosheaf on~$Y$.
 Then there is a natural isomorphism
\begin{equation} \label{flat-quasi-module-cohom-projection-formula}
 \Cohom_X(\cA,j_!\Q)\simeq j_!\Cohom_Y(j^*\cA,\Q)
\end{equation}
of $\bW$\+locally contraherent cosheaves on~$X$
(cf.~\cite[formula~(3.21) in Section~3.8]{Pcosh}).

 Indeed, for any affine open subscheme $U\subset X$ subordinate to
$\bW$, one has
\begin{multline*}
 \Cohom_X(\cA,f_!\Q)[U]=\Hom_{\cO_X(U)}(\cA(U),\Q[Y\cap U])\\
 \simeq\Hom_{\cO_X(Y\cap U)}(\cO_X(Y\cap U)\ot_{\cO_X(U)}\cA(U),
 \>\Q[Y\cap U])\\
 \simeq\Hom_{\cO_X(Y\cap U)}(\cA(Y\cap U),\Q[Y\cap U])=
 (j_!\Cohom_Y(j^*\cA,\Q))[U].
\end{multline*}

 Similarly, a natural
isomorphism~\eqref{flat-quasi-module-cohom-projection-formula} of
locally cotorsion $\bW$\+locally contraherent cosheaves on~$X$ holds
for any quasi-coherent quasi-module $\cA$ on $X$ that is flat as
a quasi-coherent sheaf in the left $\cO_X$\+module structure and any
locally cotorsion $\bW|_Y$\+contraherent cosheaf $\Q$ on~$Y$.
 Finally, for any quasi-coherent quasi-module $\cA$ on $X$ and any
locally injective $\bW|_Y$\+locally contraherent cosheaf $\gI$ on $Y$,
there is a natural isomorphism
\begin{equation} \label{quasi-module-cohom-into-lin-projection-formula}
 \Cohom_X(\cA,j_!\gI)\simeq j_!\Cohom_Y(j^*\cA,\gI)
\end{equation}
of locally cotorsion $\bW$\+locally contraherent cosheaves on~$X$
(cf.~\cite[formula~(3.22) in Section~3.8]{Pcosh}).

 If $\cA$ is a quasi-coherent quasi-algebra on $X$,
then~\eqref{flat-quasi-module-cohom-projection-formula}
and~\eqref{quasi-module-cohom-into-lin-projection-formula} are
isomorphisms of $\bW$\+locally contraherent left $\cA$\+modules on~$X$.

\Section{Some Antilocal Classes} \label{antilocal-classes-secn}

 The concept of \emph{antilocality} was defined in
the paper~\cite{Pal} in the context of commutative algebra, i.~e.,
affine schemes.
 For quasi-compact semi-separated schemes, it was worked out
in~\cite[Sections~B.3\+-B.6]{Pcosh}.
 The construction on which this concept is based goes back
to~\cite[proof of Lemma~A.1]{EP}.

 The results of
Sections~\ref{antilocality-of-X-contraadjusted-subsecn}\+-%
\ref{antilocality-of-A-cotorsion-subsecn}, concerning quasi-coherent
sheaves, provide three different $\cA$\+module versions/generalizations
of most results of~\cite[Section~4.1]{Pcosh}.
 The results of
Sections~\ref{antilocal-contrah-A-modules-subsecn}\+-%
\ref{antilocal-contrah-A-lct-A-modules-subsecn}, concerning locally
contraherent cosheaves, provide three different $\cA$\+module
versions/generalizations of most results of~\cite[Section~4.3]{Pcosh}.

\subsection{Antilocality of $X$-contraadjusted quasi-coherent
$\cA$-modules} \label{antilocality-of-X-contraadjusted-subsecn}
 Let $X$ be a scheme.
 We denote by $\Ext_X^*({-},{-})$ the Ext groups computed in
the abelian category of quasi-coherent sheaves $X\Qcoh$.

 A quasi-coherent sheaf $\F$ on $X$ is said to be
\emph{very flat}~\cite[Section~1.10]{Pcosh} if the $\cO_X(U)$\+module
$\F(U)$ is very flat for all affine open subschemes $U\subset X$.
 It suffices to check this condition for affine open subschemes $U$
belonging to any chosen affine open covering of the scheme~$X$
\,\cite[Lemma~1.2.6(a) and Section~1.10]{Pcosh}.

 The class of very flat quasi-coherent sheaves on $X$ is closed under
extensions, kernels of epimorphisms, direct summands, and infinite
direct sums (see Section~\ref{prelim-very-flat-subsecn} and
Lemma~\ref{very-flat-cotorsion-pair-hereditary}(c)).
 The class of very flat quasi-coherent sheaves on schemes is preserved
by the inverse images with respect to all morphisms of schemes
(by Lemma~\ref{restriction-co-extension-vfl-cta}(b)) and the direct
images with respect to very flat affine morphisms of schemes
(in the sense of Section~\ref{inverse-images-of-O-co-sheaves-subsecn};
by Lemma~\ref{restriction-co-extension-vfl-cta}(c)).

 A quasi-coherent sheaf $\C$ on $X$ is said to be
\emph{contraadjusted}~\cite[Section~2.5]{Pcosh} if $\Ext_X^1(\F,\C)=0$
for all very flat quasi-coherent sheaves $\F$ on~$X$.
 The class of contraadjusted quasi-coherent sheaves on $X$ is closed
under extensions, direct summands, and infinite
products~\cite[Corollary 8.3]{CoFu}, \cite[Corollary A.2]{CoSt}.
 On a quasi-compact semi-separated scheme $X$, it is also closed
under cokernels of monomorphisms~\cite[Corollary~4.1.2(c)]{Pcosh}.

 The class of contraadjusted quasi-coherent sheaves on schemes is
preserved by the direct images with respect to quasi-compact
quasi-separated morphisms of schemes~\cite[Section~2.5]{Pcosh},
\cite[Lemma~1.7(c)]{Pal}.
 For more substantial results about very flat and contraadjusted
quasi-coherent sheaves on quasi-compact semi-separated schemes,
see~\cite[Corollary~4.1.4]{Pcosh}.

 Let $\cA$ be a quasi-coherent quasi-algebra over~$X$
(in the sense of the definition in
Section~\ref{quasi-coherent-quasi-algebras-subsecn}).
 We will say that a quasi-coherent $\cA$\+module $\C$ on~$X$
(in the sense of Section~\ref{cosheaves-of-A-modules-subsecn})
is \emph{$X$\+contraadjusted} if the underlying quasi-coherent sheaf
(of $\cO_X$\+modules) of $\C$ is contraadjusted.

 Furthermore, we will say that a quasi-coherent $\cA$\+module $\F$
on $X$ is \emph{very flaprojective} if the $\cA(U)$\+module $\F(U)$
is $\cA(U)/\cO_X(U)$\+very flaprojective (in the sense of
Section~\ref{prelim-flaprojective-subsecn}) for all affine open
subschemes $U\subset X$.
 By Lemma~\ref{quasi-algebra-very-flaprojectivity-locality}(b),
it suffices to check this condition for affine open subschemes $U$
belonging to any chosen affine open covering of the scheme~$X$.
 The class of very flaprojective quasi-coherent $\cA$\+modules is
closed under extensions, kernels of epimorphisms, direct summands,
and infinite direct sums in $\cA\Qcoh$
(see Section~\ref{prelim-flaprojective-subsecn}
and Lemma~\ref{very-flaprojective-cotorsion-pair-hereditary}(b)).

 Now let $Y\subset X$ be an open subscheme with the open immersion
morphism $j\:Y\rarrow X$.
 Then, by the definition, the inverse image functor $\F\longmapsto
\F|_Y=j^*\F\:\cA\Qcoh\rarrow\cA|_Y\Qcoh$ takes very flaprojective
quasi-coherent $\cA$\+modules on $X$ to very flaprojective
quasi-coherent $\cA|_Y$\+modules on~$Y$.
 Assuming that the open immersion morphism $j\:Y\rarrow X$ is affine,
the direct image functor $j_*\:\cA|_Y\Qcoh\rarrow\cA\Qcoh$ takes
very flaprojective quasi-coherent $\cA|_Y$\+modules on $Y$ to very
flaprojective quasi-coherent $\cA$\+modules on $X$ (by
Lemma~\ref{very-flaprojective-restriction-extension-of-scalars}(a)).

 Given a quasi-coherent quasi-algebra $\cA$ over a scheme $X$, we
denote by $\Ext_\cA^*({-},{-})$ the Ext groups computed in the abelian
category $\cA\Qcoh$ of quasi-coherent left $\cA$\+modules on~$X$.

\begin{lem} \label{direct-image-of-vflp-from-affine-Ext-orthogonality}
 Let $X$ be a semi-separated scheme, $U\subset X$ be an affine
open subscheme with the open immersion morphism $j\:U\rarrow X$,
and $\cA$ be a quasi-coherent quasi-algebra over~$X$.
 Let $\G$ be a very flaprojective quasi-coherent left $\cA|_U$\+module
on~$U$ and $\C$ be an $X$\+contraadjusted quasi-coherent left
$\cA$\+module on~$X$.
 Then one has \par
\textup{(a)} $\Ext^1_\cA(j_*\G,\C)=0$; \par
\textup{(b)} assuming that $X$ is also quasi-compact,
$\Ext^n_\cA(j_*\G,\C)=0$ for all $n\ge1$.
\end{lem}

\begin{proof}
 Part~(a): the equivalence of categories $\cA|_U\Qcoh\simeq\cA(U)\Modl$
identifies very flaprojective quasi-coherent left $\cA|_U$\+modules
on $U$ with $\cA(U)/\cO(U)$\+very flaprojective left $\cA(U)$\+modules.
 By Remark~\ref{very-flaprojective-cotorsion-pair-generated-by},
the $\cA(U)/\cO(U)$\+very flaprojective $\cA(U)$\+module $\G(U)$
is a direct summand of a transfinitely iterated extension of
$\cA(U)$\+modules of the form $\cA(U)\ot_{\cO(U)}F$, where $F$ ranges
over very flat $\cO(U)$\+modules.
 The open immersion morphism $j\:U\rarrow X$ is affine, so the direct
image functor $j_*\:\cA|_Y\Qcoh\rarrow\cA\Qcoh$ preserves extensions
and inductive limits, hence also transfinitely iterated extensions.
 In view of the Eklof lemma~\cite[Lemma~1]{ET}, \cite[Lemma~7.5]{PS6},
it suffices to consider the case when $\G(U)=\cA(U)\ot_{\cO(U)}F$.

 Denoting by $\F$ the very flat quasi-coherent sheaf on $U$
corresponding to the $\cO(U)$\+module $F$, we have
$\G\simeq\cA|_U\ot_{\cO_U}\F$.
 By the projection
formula~\eqref{quasi-module-tensor-projection-formula}, we have
$j_*(\G)=j_*(j^*\cA\ot_{\cO_U}\F)\simeq\cA\ot_{\cO_X}j_*\F$.
 Notice that $j_*\F$ is a very flat quasi-coherent sheaf on~$X$
(as $j$~is a very flat affine morphism).

 It remains to point out that the functor $\cA\ot_{\cO_X}{-}\,\:X\Qcoh
\rarrow\cA\Qcoh$, which is well-defined by
Lemma~\ref{quasi-coherent-quasi-modules-tensor-product}(b),
is left adjoint to the forgetful functor $\cA\Qcoh\rarrow X\Qcoh$
(as mentioned in Section~\ref{cosheaves-of-A-modules-subsecn}).
 By~\cite[Lemma~1.7(e)]{Pal}, we have an $\Ext^1$ adjunction
isomorphism $\Ext^1_\cA(\cA\ot_{\cO_X}\cH,\>\C)\simeq
\Ext^1_X(\cH,\C)$ for any flat quasi-coherent sheaf $\cH$ and any
quasi-coherent left $\cA$\+module $\C$ on~$X$.
 Taking $\cH=j_*\F$, we arrive to the desired $\Ext^1_\cA$ vanishing
conclusion.

 Part~(b): the category $\cA\Qcoh$ is a Grothendieck abelian category,
so it has enough injective objects.
 The $\Ext^1$ adjunction isomorphism mentioned in the previous paragraph
implies that all injective quasi-coherent $\cA$\+modules are
$X$\+contraadjusted (see
Section~\ref{antilocality-of-A-cotorsion-subsecn} below for
a stronger assertion).
 The class of all $X$\+contraadjusted quasi-coherent $\cA$\+modules is
closed under cokernels of monomorphisms in $\cA\Qcoh$, since the class
of all contraadjusted quasi-coherent sheaves is closed under cokernels
of monomorphisms in $X\Qcoh$ \,\cite[Corollary~4.1.2(c)]{Pcosh}.
 Therefore, part~(b) follows from part~(a) by~\cite[Lemma~7.1]{PS6}.
\end{proof}

\begin{cor} \label{vflp-Ext-orthogonal-to-X-cta}
 Let $X$ be a quasi-compact semi-separated scheme and $\cA$ be
a quasi-coherent quasi-algebra over~$X$.
 Then, for any very flaprojective quasi-coherent left $\cA$\+module $\F$
and any $X$\+contraadjusted quasi-coherent left $\cA$\+module $\C$
on $X$, one has\/ $\Ext^n_\cA(\F,\C)=0$ for all $n\ge1$.
\end{cor}

\begin{proof}
 The argument is based on the \v Cech coresolution and
Lemma~\ref{direct-image-of-vflp-from-affine-Ext-orthogonality}(b).
 Let $X=\bigcup_{\alpha=1}^N U_\alpha$ be a finite affine open covering
of the scheme~$X$.
 Then the \v Cech
coresolution~\eqref{qcoh-sheaf-of-rings-cech-coresolution}
\begin{multline*}
 0\lrarrow\F\lrarrow\bigoplus\nolimits_{\alpha=1}^N
 j_\alpha{}_*j_\alpha^*\F\lrarrow
 \bigoplus\nolimits_{1\le\alpha<\beta\le N}
 j_{\alpha,\beta}{}_*j_{\alpha,\beta}^*\F \\
 \lrarrow\dotsb\lrarrow j_{1,\dotsc,N}{}_*j_{1,\dotsc,N}^*\F
 \lrarrow0
\end{multline*}
is a finite coresolution of $\F$ by the direct images of
very flaprojective quasi-coherent $\cA|_U$\+modules from
affine open subschemes $U=U_{\alpha_1}\cap\dotsb\cap U_{\alpha_k}
\subset X$ in $\cA\Qcoh$.
 By Lemma~\ref{direct-image-of-vflp-from-affine-Ext-orthogonality}(b),
we have $\Ext^n_\cA(j_*j^*\F,\C)=0$ for all $n\ge1$, where
$j\:U\rarrow X$ is the open immersion morphism.
 The assertion of the corollary follows.
\end{proof}

\begin{lem} \label{qcomp-qsep-very-flaprojective-precover}
 Let $X$ be a quasi-compact semi-separated scheme with a finite
affine open covering $X=\bigcup_\alpha U_\alpha$.
 Let $\cA$ be a quasi-coherent quasi-algebra over~$X$ and $\M$ be
a quasi-coherent $\cA$\+module.
 Then there exists a short exact sequence of quasi-coherent
$\cA$\+modules\/ $0\rarrow\C'\rarrow\F\rarrow\M\rarrow0$ on $X$
such that $\F$ is a very flaprojective quasi-coherent $\cA$\+module
and $\C'$ is a finitely iterated extension of the direct images of
$U_\alpha$\+contraadjusted quasi-coherent $\cA|_{U_\alpha}$\+modules
in the abelian category $\cA\Qcoh$.
\end{lem}

\begin{proof}
 This is a version of~\cite[proof of Lemma~A.1]{EP},
\cite[Proposition~4.3]{Pal}, \cite[Lemma~4.1.1 or
Proposition~B.4.4]{Pcosh}, etc.
 The same construction, proceeding (roughly speaking) by induction on
the number of affine open subschemes in the finite affine open covering
$X=\bigcup_\alpha U_\alpha$ of the scheme $X$, is applicable in
the situation at hand as well.
 The local nature of the notion of a very flaprojective quasi-coherent
$\cA$\+module plays a key role.
 It is also important that the class of very flaprojective
quasi-coherent modules is preserved by the direct images with respect
to affine open immersions of schemes, as well as extensions and
kernels of epimorphisms.
 Of course, one also needs to use
Theorem~\ref{very-flaprojective-cotorsion-pair-complete}(a).
\end{proof}

\begin{lem} \label{X-contraadjusted-modules-cogenerating-class}
 Let $X$ be a quasi-compact semi-separated scheme with a finite
affine open covering $X=\bigcup_\alpha U_\alpha$.
 Denote by $j_\alpha\:U_\alpha\rarrow X$ the open immersion morphisms.
 Let $\cA$ be a quasi-coherent quasi-algebra over~$X$.
 Then any quasi-coherent $\cA$\+module $\M$ on $X$ is a submodule of
a finite direct sum\/ $\bigoplus_\alpha j_\alpha{}_*\C_\alpha$ for
some $U_\alpha$\+contraadjusted quasi-coherent
$\cA|_{U_\alpha}$\+modules $\C_\alpha$ on~$U_\alpha$.
\end{lem}

\begin{proof}
 For every index~$\alpha$, consider the quasi-coherent
$\cA|_{U_\alpha}$\+module $j_\alpha^*\M$ on~$U_\alpha$, and pick
a $U_\alpha$\+contraadjusted quasi-coherent $\cA|_{U_\alpha}$\+module
$\C_\alpha$ together with an injective morphism of quasi-coherent
$\cA|_{U_\alpha}$\+modules $j_\alpha^*\M\rarrow\C_\alpha$
on~$U_\alpha$.
 For example, it suffices to take $\C_\alpha$ to be an injective
object of the module category $\cA|_{U_\alpha}\Qcoh\simeq
\cA(U_\alpha)\Modl$, and recall that all injective
$\cA(U_\alpha)$\+modules are contraadjusted over $\cO(U_\alpha)$ (cf.\
the proof of Lemma~\ref{flaprojective-cotorsion-pair-hereditary}(a)).
 By adjunction, we obtain a morphism of quasi-coherent $\cA$\+modules
$\M\rarrow\bigoplus_\alpha j_\alpha{}_*\C_\alpha$, which is injective
because it is injective in restriction to every~$U_\alpha$.
\end{proof}

\begin{thm} \label{qcomp-qsep-very-flaproj-complete-cotorsion-pair-thm}
 Let $X$ be a quasi-compact semi-separated scheme with a finite
affine open covering $X=\bigcup_\alpha U_\alpha$.
 Let $\cA$ be a quasi-coherent quasi-algebra over~$X$ and $\M$ be
a quasi-coherent $\cA$\+module.
 In this context: \par
\textup{(a)} there exists a short exact sequence of quasi-coherent
$\cA$\+modules\/ $0\rarrow\C'\rarrow\F\rarrow\M\rarrow0$ on $X$
such that $\F$ is a very flaprojective quasi-coherent $\cA$\+module
and $\C'$ is an $X$\+contraadjusted quasi-coherent $\cA$\+module; \par
\textup{(b)} there exists a short exact sequence of quasi-coherent
$\cA$\+modules\/ $0\rarrow\M\rarrow\C\rarrow\F'\rarrow0$ on $X$
such that $\C$ is an $X$\+contraadjusted quasi-coherent $\cA$\+module
and $\F'$ is a very flaprojective quasi-coherent $\cA$\+module; \par
\textup{(c)} a quasi-coherent $\cA$\+module on $X$ is
$X$\+contraadjusted if and only if it is a direct summand of a finitely
iterated extension of the direct images of $U_\alpha$\+contraadjusted
quasi-coherent $\cA|_{U_\alpha}$\+modules in the abelian category
$\cA\Qcoh$.
\end{thm}

\begin{proof}
 Part~(a) follows immediately from
Lemma~\ref{qcomp-qsep-very-flaprojective-precover} (since the class
of contraadjusted quasi-coherent sheaves is preserved by extensions and
direct images with respect to quasi-compact quasi-separated morphisms).
 Part~(b) follows from part~(a) and
Lemma~\ref{X-contraadjusted-modules-cogenerating-class} by virtue
of a classical argument known as the \emph{Salce lemma}~\cite{Sal},
\cite[proof of Theorem~10]{ET}, \cite[Lemma~B.1.1(b) or proof of
Lemma~4.1.3]{Pcosh}.

 To prove part~(c), one needs to observe that the proof of part~(b)
actually provides a more precise result: in part~(b), one can choose
the short exact sequence so that $\C$ is a a finitely iterated extension
of the direct images of $U_\alpha$\+contraadjusted quasi-coherent
$\cA|_{U_\alpha}$\+modules.
 Now let $\cD$ be an $X$\+contraadjusted quasi-coherent $\cA$\+module
on~$X$.
 Choose a short exact sequence $0\rarrow\cD\rarrow\C\rarrow\F'\rarrow0$,
where $\C$ is a finitely iterated extension of the direct images from
$U_\alpha$ as above, while $\F'$ is a very flaprojective
quasi-coherent $\cA$\+module.
 Then $\Ext^1_\cA(\F',\cD)=0$ by
Corollary~\ref{vflp-Ext-orthogonal-to-X-cta}, hence $\cD$ is a direct
summand of $\C$, as desired.
\end{proof}

\begin{cor} \label{qcomp-qsep-vflp-X-cta-is-a-cotorsion-pair-cor}
 Let $X$ be a quasi-compact semi-separated scheme, and let $\cA$ be
a quasi-coherent quasi-algebra over~$X$.
 In this context: \par
\textup{(a)} a quasi-coherent left $\cA$\+module $\F$ is very
flaprojective if and only if\/ $\Ext^1_\cA(\F,\C)\allowbreak=0$ for all
$X$\+contraadjusted quasi-coherent left $\cA$\+modules~$\C$; \par
\textup{(b)} a quasi-coherent left $\cA$\+module $\C$ is
$X$\+contraadjusted if and only if\/ $\Ext^1_\cA(\F,\C)\allowbreak=0$
for all very flaprojective quasi-coherent left $\cA$\+modules~$\F$.
\end{cor}

\begin{proof}
 The ``only if'' assertions are provided by
Corollary~\ref{vflp-Ext-orthogonal-to-X-cta}.
 The ``if'' implications follow from the same corollary and
Theorem~\ref{qcomp-qsep-very-flaproj-complete-cotorsion-pair-thm}(a\+-b)
by the argument similar to the one in the proof of
Theorem~\ref{qcomp-qsep-very-flaproj-complete-cotorsion-pair-thm}(c);
see~\cite[Lemma~B.1.2]{Pcosh}.
 The point is that the classes of very flaprojective and
$X$\+contraadjusted quasi-coherent $\cA$\+modules are closed under
direct summands in $\cA\Qcoh$.
\end{proof}

\begin{rem} \label{X-contraadjusted-A-modules-concluding-remark}
 The assertions of
Theorem~\ref{qcomp-qsep-very-flaproj-complete-cotorsion-pair-thm}(a\+-b)
and Corollaries~\ref{vflp-Ext-orthogonal-to-X-cta}
and~\ref{qcomp-qsep-vflp-X-cta-is-a-cotorsion-pair-cor}
can be summarized by saying that the pair of classes (very flaprojective
quasi-coherent $\cA$\+modules, $X$\+contraadjusted quasi-coherent
$\cA$\+modules) is a hereditary complete cotorsion pair in $\cA\Qcoh$.
 On the other hand,
Theorem~\ref{qcomp-qsep-very-flaproj-complete-cotorsion-pair-thm}(c) is
an extra piece of information which one obtains from the proof of
Theorem~\ref{qcomp-qsep-very-flaproj-complete-cotorsion-pair-thm}(a\+-b)
presented above.

 The assertion of
Theorem~\ref{qcomp-qsep-very-flaproj-complete-cotorsion-pair-thm}(c)
can be expressed by saying that the class of $X$\+contraadjusted
quasi-coherent $\cA$\+modules on a quasi-compact semi-separated
scheme $X$ is \emph{antilocal} in the sense of~\cite[Section~4]{Pal}.
 This is the main result of this
Section~\ref{antilocality-of-X-contraadjusted-subsecn} for our purposes
in this paper.
 The whole machinery of very flaprojective modules was developed in
Section~\ref{prelim-flaprojective-subsecn} with the sole aim to prove
this result (which does not mention very flaprojective modules).
\end{rem}

\subsection{Antilocality of $X$-cotorsion quasi-coherent $\cA$-modules}
\label{antilocality-of-X-cotorsion-subsecn}
 This section is a slight variation upon the previous one.
 All the proofs are similar, and therefore omitted or replaced by
references.

 Let $X$ be a scheme.
 A quasi-coherent sheaf $\F$ on $X$ is said to be \emph{flat} if
the $\cO_X(U)$\+module $\F(U)$ is flat for all affine open subschemes
$U\subset X$.
 It suffices to check this condition for affine open subschemes $U$
belonging to any chosen affine open covering of the scheme~$X$
(see~\cite[Lemma~1.8.3]{Pcosh} for a more general relative version
of this assertion).

 The class of flat quasi-coherent sheaves on $X$ is closed under
extensions, kernels of epimorphisms, direct summands, and infinite
direct sums.
 The class of flat quasi-coherent sheaves on schemes is preserved
by the inverse images with respect to all morphisms of schemes
and the direct images with respect to flat affine morphisms of schemes.

 A quasi-coherent sheaf $\C$ on $X$ is said to be
\emph{cotorsion}~\cite{EE} if $\Ext_X^1(\F,\C)=0$ for all flat
quasi-coherent sheaves $\F$ on~$X$.
 The class of cotorsion quasi-coherent sheaves on $X$ is closed
under extensions, direct summands, and infinite
products~\cite[Corollary 8.3]{CoFu}, \cite[Corollary A.2]{CoSt}.
 On a quasi-compact semi-separated scheme $X$, it is also closed
under cokernels of monomorphisms~\cite[Corollary~4.1.9(c)]{Pcosh}.

 The class of cotorsion quasi-coherent sheaves on schemes is preserved
by the direct images with respect to quasi-compact quasi-separated
morphisms of schemes~\cite[Section~2.5]{Pcosh},
\cite[Lemma~1.7(c)]{Pal}.
 For more substantial results about flat and cotorsion quasi-coherent
sheaves on quasi-compact semi-separated schemes,
see~\cite[Corollary~4.2]{EE}, \cite[Section~2.4]{M-n},
\cite[Lemma~A.1]{EP}, or~\cite[Corollary~4.1.11]{Pcosh}.

 Let $\cA$ be a quasi-coherent quasi-algebra over~$X$.
 We will say that a quasi-coherent $\cA$\+module $\C$ on $X$ is
\emph{$X$\+cotorsion} if the underlying quasi-coherent sheaf
(of $\cO_X$\+modules) of $\C$ is cotorsion.

 Furthermore, we will say that a quasi-coherent $\cA$\+module $\F$
on $X$ is \emph{flaprojective} if the $\cA(U)$\+module $\F(U)$
is $\cA(U)/\cO_X(U)$\+flaprojective (in the sense of
Section~\ref{prelim-flaprojective-subsecn}) for all affine open
subschemes $U\subset X$.
 By Lemma~\ref{quasi-algebra-very-flaprojectivity-locality}(a),
it suffices to check this condition for affine open subschemes $U$
belonging to any chosen affine open covering of the scheme~$X$.
 The class of flaprojective quasi-coherent $\cA$\+modules is
closed under extensions, kernels of epimorphisms, direct summands,
and infinite direct sums in $\cA\Qcoh$
(see Section~\ref{prelim-flaprojective-subsecn}
and Lemma~\ref{flaprojective-cotorsion-pair-hereditary}(b)).

 Now let $Y\subset X$ be an open subscheme with the open immersion
morphism $j\:Y\rarrow X$.
 Then, by the definition, the inverse image functor $\F\longmapsto
\F|_Y=j^*\F\:\cA\Qcoh\rarrow\cA|_Y\Qcoh$ takes flaprojective
quasi-coherent $\cA$\+modules on $X$ to flaprojective
quasi-coherent $\cA|_Y$\+modules on~$Y$.
 Assuming that the open immersion morphism $j\:Y\rarrow X$ is affine,
the direct image functor $j_*\:\cA|_Y\Qcoh\rarrow\cA\Qcoh$ takes
flaprojective quasi-coherent $\cA|_Y$\+modules on $Y$ to flaprojective
quasi-coherent $\cA$\+modules on $X$ (by
Lemma~\ref{flaprojective-restriction-extension-of-scalars}(a)).

\begin{lem} \label{direct-image-of-flp-from-affine-Ext-orthogonality}
 Let $X$ be a semi-separated scheme, $U\subset X$ be an affine
open subscheme with the open immersion morphism $j\:Y\rarrow X$,
and $\cA$ be a quasi-coherent quasi-algebra over~$X$.
 Let $\G$ be a flaprojective quasi-coherent left $\cA|_U$\+module on~$U$
and $\C$ be an $X$\+cotorsion quasi-coherent left $\cA$\+module on~$X$.
 Then one has \par
\textup{(a)} $\Ext^1_\cA(j_*\G,\C)=0$; \par
\textup{(b)} assuming that $X$ is also quasi-compact,
$\Ext^n_\cA(j_*\G,\C)=0$ for all $n\ge1$.
\end{lem}

\begin{proof}
 Similar to
Lemma~\ref{direct-image-of-vflp-from-affine-Ext-orthogonality}.
 In part~(b), one observes that all injective quasi-coherent
$\cA$\+modules are actually $X$\+cotorsion (on any scheme~$X$;
see Section~\ref{antilocality-of-A-cotorsion-subsecn} below for
a more general assertion).
 Then one also has to use~\cite[Lemma~4.1.9(c)]{Pcosh}.
\end{proof}

\begin{cor} \label{flp-Ext-orthogonal-to-X-cot}
 Let $X$ be a quasi-compact semi-separated scheme and $\cA$ be
a quasi-coherent quasi-algebra over~$X$.
 Then, for any flaprojective quasi-coherent left $\cA$\+module $\F$
and any $X$\+cotorsion quasi-coherent left $\cA$\+module $\C$ on $X$,
one has\/ $\Ext^n_\cA(\F,\C)=0$ for all $n\ge1$.
\end{cor}

\begin{proof}
 Similar to Corollary~\ref{vflp-Ext-orthogonal-to-X-cta}, and based
on Lemma~\ref{direct-image-of-flp-from-affine-Ext-orthogonality}(b).
\end{proof}

\begin{lem} \label{qcomp-qsep-flaprojective-precover}
 Let $X$ be a quasi-compact semi-separated scheme with a finite
affine open covering $X=\bigcup_\alpha U_\alpha$.
 Let $\cA$ be a quasi-coherent quasi-algebra over~$X$ and $\M$ be
a quasi-coherent $\cA$\+module.
 Then there exists a short exact sequence of quasi-coherent
$\cA$\+modules\/ $0\rarrow\C'\rarrow\F\rarrow\M\rarrow0$ on $X$
such that $\F$ is a flaprojective quasi-coherent $\cA$\+module
and $\C'$ is a finitely iterated extension of the direct images of
$U_\alpha$\+cotorsion quasi-coherent $\cA|_{U_\alpha}$\+modules
in the abelian category $\cA\Qcoh$.
\end{lem}

\begin{proof}
 Similar to Lemma~\ref{qcomp-qsep-very-flaprojective-precover}
and using Theorem~\ref{flaprojective-cotorsion-pair-complete}(a).
\end{proof}

\begin{lem} \label{X-cotorsion-modules-cogenerating-class}
 Let $X$ be a quasi-compact semi-separated scheme with a finite
affine open covering $X=\bigcup_\alpha U_\alpha$.
 Denote by $j_\alpha\:U_\alpha\rarrow X$ the open immersion morphisms.
 Let $\cA$ be a quasi-coherent quasi-algebra over~$X$.
 Then any quasi-coherent $\cA$\+module $\M$ on $X$ is a submodule of
a finite direct sum\/ $\bigoplus_\alpha j_\alpha{}_*\C_\alpha$ for
some $U_\alpha$\+cotorsion quasi-coherent $\cA|_{U_\alpha}$\+modules
$\C_\alpha$ on~$U_\alpha$.
\end{lem}

\begin{proof}
 This is a slightly stronger version of
Lemma~\ref{X-contraadjusted-modules-cogenerating-class},
provable in the same way.
\end{proof}

\begin{thm} \label{qcomp-qsep-flaproj-complete-cotorsion-pair-thm}
 Let $X$ be a quasi-compact semi-separated scheme with a finite
affine open covering $X=\bigcup_\alpha U_\alpha$.
 Let $\cA$ be a quasi-coherent quasi-algebra over~$X$ and $\M$ be
a quasi-coherent $\cA$\+module.
 In this context: \par
\textup{(a)} there exists a short exact sequence of quasi-coherent
$\cA$\+modules\/ $0\rarrow\C'\rarrow\F\rarrow\M\rarrow0$ on $X$
such that $\F$ is a flaprojective quasi-coherent $\cA$\+module
and $\C'$ is an $X$\+cotorsion quasi-coherent $\cA$\+module; \par
\textup{(b)} there exists a short exact sequence of quasi-coherent
$\cA$\+modules\/ $0\rarrow\M\rarrow\C\rarrow\F'\rarrow0$ on $X$
such that $\C$ is an $X$\+cotorsion quasi-coherent $\cA$\+module
and $\F'$ is a flaprojective quasi-coherent $\cA$\+module; \par
\textup{(c)} a quasi-coherent $\cA$\+module on $X$ is $X$\+cotorsion
if and only if it is a direct summand of a finitely iterated extension
of the direct images of $U_\alpha$\+cotorsion quasi-coherent
$\cA|_{U_\alpha}$\+modules in the abelian category $\cA\Qcoh$.
\end{thm}

\begin{proof}
 Similar to
Theorem~\ref{qcomp-qsep-very-flaproj-complete-cotorsion-pair-thm},
and based on Corollary~\ref{flp-Ext-orthogonal-to-X-cot} and
Lemmas~\ref{qcomp-qsep-flaprojective-precover}\+-%
\ref{X-cotorsion-modules-cogenerating-class}.
\end{proof}

\begin{cor} \label{qcomp-qsep-flp-X-cot-is-a-cotorsion-pair-cor}
 Let $X$ be a quasi-compact semi-separated scheme, and let $\cA$ be
a quasi-coherent quasi-algebra over~$X$.
 In this context: \par
\textup{(a)} a quasi-coherent left $\cA$\+module $\F$ is flaprojective
if and only if\/ $\Ext^1_\cA(\F,\C)=0$ for all $X$\+cotorsion
quasi-coherent left $\cA$\+modules~$\C$; \par
\textup{(b)} a quasi-coherent left $\cA$\+module $\C$ is
$X$\+cotorsion if and only if\/ $\Ext^1_\cA(\F,\C)=0$ for all
flaprojective quasi-coherent left $\cA$\+modules~$\F$.
\end{cor}

\begin{proof}
 Similar to
Corollary~\ref{qcomp-qsep-vflp-X-cta-is-a-cotorsion-pair-cor}.
\end{proof}

\begin{rem} \label{X-cotorsion-A-modules-concluding-remark}
 Similarly to Remark~\ref{X-contraadjusted-A-modules-concluding-remark},
the assertions of
Theorem~\ref{qcomp-qsep-flaproj-complete-cotorsion-pair-thm}(a\+-b)
and Corollaries~\ref{flp-Ext-orthogonal-to-X-cot}
and~\ref{qcomp-qsep-flp-X-cot-is-a-cotorsion-pair-cor} can be summarized
by saying that the pair of classes (flaprojective quasi-coherent
$\cA$\+modules, $X$\+cotorsion quasi-coherent $\cA$\+modules) is
a hereditary complete cotorsion pair in $\cA\Qcoh$.
 On the other hand,
Theorem~\ref{qcomp-qsep-flaproj-complete-cotorsion-pair-thm}(c) is
an extra piece of information which one obtains from our proof of
Theorem~\ref{qcomp-qsep-flaproj-complete-cotorsion-pair-thm}(a\+-b).

 The assertion of
Theorem~\ref{qcomp-qsep-flaproj-complete-cotorsion-pair-thm}(c)
can be expressed by saying that the class of $X$\+cotorsion
quasi-coherent $\cA$\+modules on a quasi-compact semi-separated
scheme $X$ is \emph{antilocal} in the sense of~\cite[Section~4]{Pal}.
 This is the main result of this
Section~\ref{antilocality-of-X-cotorsion-subsecn} for our purposes
in this paper.
 The whole machinery of flaprojective modules was developed in
Section~\ref{prelim-flaprojective-subsecn} with the sole aim to prove
this result (which does not mention flaprojective modules).
\end{rem}

\subsection{Antilocality of $\cA$-cotorsion quasi-coherent
$\cA$-modules} \label{antilocality-of-A-cotorsion-subsecn}
 Let $X$ be a scheme and $\cA$ be a quasi-coherent quasi-algebra
over~$X$.
 Recall that, according to
Section~\ref{A-loc-cotors-loc-inj-cosheaves-subsecn},
a quasi-coherent left $\cA$\+module is said to be \emph{flat}
(or \emph{$\cA$\+flat}) if the $\cA(U)$\+module $\F(U)$ is flat
for all affine open subschemes $U\subset X$.
 It suffices to check this condition for affine open subschemes $U$
belonging to any chosen affine open covering of the scheme~$X$
(see Lemma~\ref{quasi-algebra-adjustedness-co-locality}(a)).

 The class of flat quasi-coherent $\cA$\+modules is closed under
extensions, kernels of epimorphisms, direct summands, and infinite
direct sums in $\cA\Qcoh$.
 Given an open subscheme $Y\subset X$ with the open immersion morphism
$j\:Y\rarrow X$, the inverse image functor $j^*\:\cA\Qcoh\rarrow
\cA|_Y\Qcoh$ takes flat quasi-coherent $\cA$\+modules to flat
quasi-coherent $\cA|_Y$\+modules.
 If the open immersion morphism~$j$ is affine, then the direct image
functor $j_*\:\cA|_Y\Qcoh\rarrow\cA\Qcoh$ takes flat quasi-coherent
$\cA|_Y$\+modules to flat quasi-coherent $\cA$\+modules (as one can see
from Corollary~\ref{quasi-algebra-change-of-scalars-adjustedness}(a)).

 A quasi-coherent left $\cA$\+module $\C$ is said to be
\emph{$\cA$\+cotorsion} (or just \emph{cotorsion}) if
$\Ext^1_\cA(\F,\C)=0$ for all flat quasi-coherent left $\cA$\+modules
$\F$ on~$X$.
 The class of $\cA$\+cotorsion quasi-coherent $\cA$\+modules is closed
under extensions, direct summands, and direct products in $\cA\Qcoh$
\,\cite[Corollary 8.3]{CoFu}, \cite[Corollary A.2]{CoSt}.
 For a quasi-coherent quasi-algebra $\cB$ over an affine scheme $U$,
a quasi-coherent $\cB$\+module $\C$ on $U$ is $\cB$\+cotorsion if
and only if the $\cB(U)$\+module $\C(U)$ is cotorsion.

 Given an open subscheme $Y\subset X$ with a quasi-compact open
immersion morphism $j\:Y\rarrow X$, the direct image functor
$j_*\:\cA|_Y\Qcoh\rarrow\cA\Qcoh$ takes $\cA|_Y$\+cotorsion
quasi-coherent $\cA|_Y$\+modules to $\cA$\+cotorsion quasi-coherent
$\cA$\+modules.
 Indeed, the direct image functor~$j_*$ is right adjoint to the inverse
image functor $j^*\:\cA\Qcoh\rarrow\cA|_Y\Qcoh$, which is exact.
 Therefore, for any quasi-coherent left $\cA$\+module $\M$ on $X$ and
any quasi-coherent left $\cA|_Y$\+module $\N$ on $Y$ there is
a natural injective map of abelian groups $\Ext^1_\cA(\M,j_*\N)\rarrow
\Ext^1_{\cA|_Y}(j^*\M,\N)$ \,\cite[Lemma~1.7(b)]{Pal}.
 So $\Ext^1_{\cA|_Y}(j^*\M,\N)=0$ implies $\Ext^1_\cA(\M,j_*\N)=0$.

 All $\cA$\+cotorsion quasi-coherent $\cA$\+modules are $X$\+cotorsion.
 Indeed, the functor $\cA\ot_{\cO_X}{-}\,\:X\Qcoh\rarrow\cA\Qcoh$
takes flat quasi-coherent sheaves on $X$ to flat quasi-coherent
$\cA$\+modules.
 The functor $\cA\ot_{\cO_X}{-}$ is also left adjoint to the exact
forgetful functor $\cA\Qcoh\rarrow X\Qcoh$.
 Therefore, for any flat quasi-coherent sheaf $\cH$ and any
quasi-coherent left $\cA$\+module $\C$ on $X$ there is a natural
isomorphism of abelian groups $\Ext^1_\cA(\cA\ot_{\cO_X}\cH,\>\C)
\simeq\Ext^1_X(\cH,\C)$ \,\cite[Lemma~1.7(e)]{Pal}, as mentioned in
the proof of
Lemma~\ref{direct-image-of-vflp-from-affine-Ext-orthogonality}.
 So $\Ext^1_\cA(\cA\ot_{\cO_X}\cH,\>\C)=0$ implies
$\Ext^1_X(\cH,\C)=0$.

\begin{lem} \label{qcomp-qsep-quasi-algebra-flat-precover}
 Let $X$ be a quasi-compact semi-separated scheme with a finite
affine open covering $X=\bigcup_\alpha U_\alpha$.
 Let $\cA$ be a quasi-coherent quasi-algebra over~$X$ and $\M$ be
a quasi-coherent $\cA$\+module.
 Then there exists a short exact sequence of quasi-coherent
$\cA$\+modules\/ $0\rarrow\C'\rarrow\F\rarrow\M\rarrow0$ on $X$
such that $\F$ is a flat quasi-coherent $\cA$\+module and $\C'$ is
a finitely iterated extension of the direct images of
$\cA|_{U_\alpha}$\+cotorsion quasi-coherent $\cA|_{U_\alpha}$\+modules
in the abelian category $\cA\Qcoh$.
\end{lem}

\begin{proof}
 Similar to the proofs of
Lemmas~\ref{qcomp-qsep-very-flaprojective-precover}
and~\ref{qcomp-qsep-flaprojective-precover}.
 One needs to use various properties of flat quasi-coherent
$\cA$\+modules listed above, such as the local nature of the definition,
closedness under extensions and kernels of epimorphisms, and 
preservation by the direct images with respect to affine open
immersions.
 One also has to use Theorem~\ref{flat-cotorsion-pair-complete}(a)
for the rings $R=\cA(U_\alpha)$.
\end{proof}

\begin{cor} \label{A-cotorsion-hereditariness-properties}
 Let $X$ be quasi-compact semi-separated scheme and $\cA$ be
a quasi-coherent quasi-algebra over~$X$.
 In this context: \par
\textup{(a)} A quasi-coherent left $\cA$\+module $\C$ on $X$ is
$\cA$\+cotorsion if and only if the functor\/ $\Hom_\cA({-},\C)$ takes
takes short exact sequences of flat quasi-coherent left $\cA$\+modules
to short exact sequences of abelian groups. \par
\textup{(b)} A quasi-coherent left $\cA$\+module $\C$ on $X$ is
$\cA$\+cotorsion if and only if\/ $\Ext_\cA^n(\F,\C)=0$ for all
all flat quasi-coherent left $\cA$\+modules $\F$ on~$X$. \par
\textup{(c)} The class of $\cA$\+cotorsion quasi-coherent $\cA$\+modules
is closed under cokernels of monomorphisms in $\cA\Qcoh$.
\end{cor}

\begin{proof}
 This is similar to~\cite[Corollaries~4.1.2 and~4.1.9]{Pcosh}
and based on the fact that every quasi-coherent $\cA$\+module is
a quotient of a flat one (which is a weak version of
Lemma~\ref{qcomp-qsep-quasi-algebra-flat-precover}) together with
the fact that the class of flat quasi-coherent $\cA$\+modules is
closed under kernels of epimorphisms.
\end{proof}

\begin{lem} \label{A-cotorsion-modules-cogenerating-class}
 Let $X$ be a quasi-compact semi-separated scheme with a finite
affine open covering $X=\bigcup_\alpha U_\alpha$.
 Denote by $j_\alpha\:U_\alpha\rarrow X$ the open immersion morphisms.
 Let $\cA$ be a quasi-coherent quasi-algebra over~$X$.
 Then any quasi-coherent $\cA$\+module $\M$ on $X$ is a submodule of
a finite direct sum\/ $\bigoplus_\alpha j_\alpha{}_*\C_\alpha$ for
some $\cA|_{U_\alpha}$\+cotorsion quasi-coherent
$\cA|_{U_\alpha}$\+modules $\C_\alpha$ on~$U_\alpha$.
\end{lem}

\begin{proof}
 Similar to
Lemmas~\ref{X-contraadjusted-modules-cogenerating-class}
and~\ref{X-cotorsion-modules-cogenerating-class}.
 Notice that all injective $\cA(U_\alpha)$\+modules are cotorsion
by the definition.
\end{proof}

\begin{thm} \label{qcomp-qsep-A-flat-complete-cotorsion-pair-thm}
 Let $X$ be a quasi-compact semi-separated scheme with a finite
affine open covering $X=\bigcup_\alpha U_\alpha$.
 Let $\cA$ be a quasi-coherent quasi-algebra over~$X$ and $\M$ be
a quasi-coherent $\cA$\+module.
 In this context: \par
\textup{(a)} there exists a short exact sequence of quasi-coherent
$\cA$\+modules\/ $0\rarrow\C'\rarrow\F\rarrow\M\rarrow0$ on $X$
such that $\F$ is a flat quasi-coherent $\cA$\+module and $\C'$ is
an $\cA$\+cotorsion quasi-coherent $\cA$\+module; \par
\textup{(b)} there exists a short exact sequence of quasi-coherent
$\cA$\+modules\/ $0\rarrow\M\rarrow\C\rarrow\F'\rarrow0$ on $X$
such that $\C$ is an $\cA$\+cotorsion quasi-coherent $\cA$\+module
and $\F'$ is a flat quasi-coherent $\cA$\+module; \par
\textup{(c)} a quasi-coherent $\cA$\+module on $X$ is $\cA$\+cotorsion
if and only if it is a direct summand of a finitely iterated extension
of the direct images of $\cA|_{U_\alpha}$\+cotorsion quasi-coherent
$\cA|_{U_\alpha}$\+modules in the abelian category $\cA\Qcoh$.
\end{thm}

\begin{proof}
 Similar to
Theorems~\ref{qcomp-qsep-very-flaproj-complete-cotorsion-pair-thm}
and~\ref{qcomp-qsep-flaproj-complete-cotorsion-pair-thm}.
 Part~(a) follows immediately from
Lemma~\ref{qcomp-qsep-quasi-algebra-flat-precover}, since the class
of $\cA$\+cotorsion quasi-coherent $\cA$\+modules is preserved by
extensions and direct images from affine open subschemes.
 Part~(b) follows from part~(a) and
Lemma~\ref{A-cotorsion-modules-cogenerating-class} by virtue of
the Salce lemma.
 Concerning part~(c), one has $\Ext^1_\cA(\F',\cD)=0$ for all
flat quasi-coherent $\cA$\+modules $\F'$ and all $\cA$\+cotorsion
quasi-coherent $\cA$\+modules $\cD$ on $X$ by the definition.
\end{proof}

\begin{cor} \label{qcomp-qsep-A-flat-is-a-cotorsion-pair-cor}
 Let $X$ be a quasi-compact semi-separated scheme, and let $\cA$ be
a quasi-coherent quasi-algebra over~$X$.
 Then a quasi-coherent left $\cA$\+module $\F$ is flat if and only if\/
$\Ext^1_\cA(\F,\C)=0$ for all $\cA$\+cotorsion quasi-coherent left
$\cA$\+modules~$\C$.
\end{cor}

\begin{proof}
 Similar to
Corollaries~\ref{qcomp-qsep-vflp-X-cta-is-a-cotorsion-pair-cor}(a)
and~\ref{qcomp-qsep-flp-X-cot-is-a-cotorsion-pair-cor}(a).
\end{proof}

\begin{rem} \label{A-cotorsion-A-modules-concluding-remark}
 Similarly to Remarks~\ref{X-contraadjusted-A-modules-concluding-remark}
and~\ref{X-cotorsion-A-modules-concluding-remark} above, the assertions
of Theorem~\ref{qcomp-qsep-A-flat-complete-cotorsion-pair-thm}(a\+-b)
and Corollaries~\ref{A-cotorsion-hereditariness-properties}
and~\ref{qcomp-qsep-A-flat-is-a-cotorsion-pair-cor} can be summarized by
saying that the pair of classes (flat quasi-coherent $\cA$\+modules,
$\cA$\+cotorsion quasi-coherent $\cA$\+modules) is a hereditary complete
cotorsion pair in $\cA\Qcoh$.
 On the other hand,
Theorem~\ref{qcomp-qsep-A-flat-complete-cotorsion-pair-thm}(c) is
an extra piece of information which one obtains from our proof of
Theorem~\ref{qcomp-qsep-A-flat-complete-cotorsion-pair-thm}(a\+-b).

 The assertion of
Theorem~\ref{qcomp-qsep-A-flat-complete-cotorsion-pair-thm}(c)
can be expressed by saying that the class of $\cA$\+cotorsion
quasi-coherent $\cA$\+modules on a quasi-compact semi-separated
scheme $X$ is \emph{antilocal} in the sense of~\cite[Section~4]{Pal}.
 This is the main result of this
Section~\ref{antilocality-of-A-cotorsion-subsecn} for our purposes
in this paper.
\end{rem}

\subsection{Antilocal locally contraherent $\cA$\+modules}
\label{A-antilocal=X-antilocal-subsecn}
 Let $X$ be a scheme with an open covering~$\bW$.
 A $\bW$\+locally contraherent cosheaf $\P$ on~$X$ is said to be
\emph{antilocal}~\cite[Section~4.3]{Pcosh} if the functor
$\Hom^X(\P,{-})$ takes short exact sequences of locally injective
$\bW$\+locally contraherent cosheaves on $X$ to short exact
sequences of abelian groups.
 In particular, any contraherent cosheaf on an affine scheme $U$
with the open covering $\{U\}$ is antilocal.

 We refer to~\cite[Section~4.3]{Pcosh} for substantial results about
antilocal $\bW$\+locally contraherent cosheaves on quasi-compact
semi-separated schemes~$X$.
 In particular, the class of all such cosheaves does not depend on
the open covering $\bW$, and all of them are (globally) contraherent
on~$X$ \,\cite[Corollaries~4.3.5(c) and~4.3.7]{Pcosh}.
 The class of all antilocal contraherent cosheaves on $X$ is closed
under extensions, kernels of admissible epimorphisms, direct summands,
and direct products in $X\Ctrh$ \,\cite[Corollaries~4.3.3(b)
and~4.3.10]{Pcosh}.
 The class of antilocal contraherent cosheaves is preserved by
the direct images with respect to all morphisms of quasi-compact
semi-separated schemes~\cite[Corollary~4.6.3(a)]{Pcosh}.

 Let $\cA$ be a quasi-coherent quasi-algebra over~$X$.
 We will say that a $\bW$\+locally contraherent $\cA$\+module $\P$
on~$X$ (in the sense of the definition in
Section~\ref{cosheaves-of-A-modules-subsecn})
is \emph{$X$\+antilocal} if the underlying $\bW$\+locally contraherent
cosheaf (of $\cO_X$\+modules) of $\P$ is antilocal.

 We will also say that a $\bW$\+locally contraherent $\cA$\+module $\P$
is \emph{$\cA$\+antilocal} if the functor $\Hom^\cA(\P,{-})$ takes short
exact sequences of $\cA$\+locally injective $\bW$\+locally contraherent
$\cA$\+modules on~$X$ (in the sense of
Section~\ref{A-loc-cotors-loc-inj-cosheaves-subsecn}) to short exact
sequences of abelian groups.
 In particular, any contraherent $\cB$\+module on an affine scheme $U$
with the open covering $\{U\}$ and a quasi-coherent quasi-algebra $\cB$
over $U$ is antilocal, because all short exact sequences of
$\cB$\+locally injective contraherent $\cB$\+modules on $U$ are split
(the category $\cB\Ctrh^{\cB\dlin}$ being equivalent to the category
$\cB(U)\Modl^\inj$ of injective modules over the ring $\cB(U)$;
see Section~\ref{A-loc-cotors-loc-inj-cosheaves-subsecn}).

\begin{lem} \label{A-antilocal-are-X-antilocal}
 For any scheme $X$ with an open covering\/ $\bW$ and any quasi-coherent
quasi-algebra $\cA$ over $X$, all $\cA$\+antilocal\/ $\bW$\+locally
contraherent $\cA$\+modules on $X$ are $X$\+antilocal.
\end{lem}

\begin{proof}
 The argument uses the functor $\Cohom_X$ from quasi-coherent
quasi-modules to $\bW$\+locally contraherent cosheaves on $X$, which
was defined in Section~\ref{cohom-from-quasi-modules-subsecn}.
 For any locally injective $\bW$\+locally contraherent cosheaf $\gJ$
on $X$, we have an $\cA$\+locally injective $\bW$\+locally contraherent
$\cA$\+module $\Cohom_X(\cA,\gJ)$ on~$X$ (by
Lemmas~\ref{cohom-loc-contraherent-lemma}(d)
and~\ref{restriction-coextension-injective-cotorsion}(d),
and the discussion in Section~\ref{cosheaves-of-A-modules-subsecn}).
 For any $\bW$\+locally contraherent cosheaf $\P$ on $X$,
the adjunction isomorphism of abelian groups $\Hom^X(\P,\gJ)\simeq
\Hom^\cA(\P,\Cohom_X(\cA,\gJ))$ holds by
formula~\eqref{cohom-A-adjunction}.
 
 Now if $0\rarrow\gJ'\rarrow\gJ\rarrow\gJ''\rarrow0$ is a short exact
sequence in the exact category $X\Lcth_\bW^\lin$, then $0\rarrow
\Cohom_X(\cA,\gJ')\rarrow\Cohom_X(\cA,\gJ)\rarrow\Cohom_X(\cA,\gJ'')
\rarrow0$ is a short exact sequence in the exact category
$\cA\Lcth_\bW^{\cA\dlin}$.
 The assertion of the lemma follows from these observations.
\end{proof}

 Let $X$ be a scheme with an open covering $\bW$ and $\cA$ be
a quasi-coherent quasi-algebra over~$X$.
 Let us denote temporarily by $\Ext^{X,*}({-},{-})$ and
$\Ext^{\cA,*}({-},{-})$ the Ext groups computed in the exact
categories $X\Lcth_\bW$ and $\cA\Lcth_\bW$, respectively
(see~\cite[Section~4.3]{Pcosh} and
Remark~\ref{Ext-unambiguous-remark} below for a discussion).

 The following lemma and corollary are dual-analogous to
Lemma~\ref{direct-image-of-vflp-from-affine-Ext-orthogonality} and
Corollary~\ref{vflp-Ext-orthogonal-to-X-cta} (as well as 
to Lemma~\ref{direct-image-of-flp-from-affine-Ext-orthogonality} and
Corollary~\ref{flp-Ext-orthogonal-to-X-cot}).

\begin{lem} \label{direct-image-of-A-lin-from-affine-Ext-orthogonality}
 Let $X$ be a semi-separated scheme with an open covering\/ $\bW$,
let $U\subset X$ be an affine open subscheme subordinate to\/ $\bW$
with the open immersion morphism $j\:U\rarrow X$, and let $\cA$ be
a quasi-coherent quasi-algebra over~$X$.
 Let $\P$ be a\/ $\bW$\+locally contraherent $\cA$\+module on $X$ such
that its underlying\/ $\bW$\+locally contraherent cosheaf (of
$\cO_X$\+modules) has the property that\/ $\Ext^{X,1}(\P,\gK)=0$
for all locally injective\/ $\bW$\+locally contraherent cosheaves\/
$\gK$ on~$X$.
 Let\/ $\gI$ be an $\cA|_U$\+locally injective contraherent
$\cA|_U$\+module on~$U$.
 Then one has \par
\textup{(a)} $\Ext^{\cA,1}(\P,j_!\gI)=0$; \par
\textup{(b)} assuming that $X$ is also quasi-compact,
$\Ext^{\cA,n}(\P,j_!\gI)=0$ for all $n\ge1$.
\end{lem}

\begin{proof}
 Part~(a): the equivalence of categories $\cA|_U\Ctrh\simeq
\cA(U)\Modl^{\cO(U)\dcta}$ identifies $\cA|_U$\+locally injective
contraherent $\cA$\+modules on $U$ with injective
left $\cA(U)$\+modules.
 For any associative ring homomorphism $R\rarrow A$, all injective
left $A$\+modules are direct summands of $A$\+modules of the form
$\Hom_R(A,J)$, where $J$ ranges over injective left $R$\+modules.
 In particular, all injective left $\cA(U)$\+module are direct summands
of $\cA(U)$\+modules of the form $\Hom_{\cO(U)}(\cA(U),J)$, where
$J$ ranges over injective $\cO_X(U)$\+modules.
 Therefore, it suffices to consider the case when $\gI[U]=
\Hom_{\cO(U)}(\cA(U),J)$.

 Denoting by $\gJ$ the locally injective contraherent cosheaf on $U$
corresponding to the $\cO(U)$\+module $J$, we have
$\gI=\Cohom_U(\cA|_U,\gJ)$.
 By the projection
formula~\eqref{quasi-module-cohom-into-lin-projection-formula},
we have $j_!\gI=j_!\Cohom_U(j^*\cA,\gJ)\simeq\Cohom_X(\cA,j_!\gJ)$.
 Notice that $\gK=j_!\gJ$ is a locally injective contraherent cosheaf
on~$X$ (as $j$~is a flat affine morphism; see
Section~\ref{direct-images-of-O-co-sheaves-subsecn}).

 Now we argue similarly to~\cite[proof of Lemma~1.7(c)]{Pal}.
 Suppose given a short exact sequence $0\rarrow\Cohom_X(\cA,\gK)
\rarrow\Q\rarrow\P\rarrow0$ in the exact category $\cA\Lcth_\bW$.
 Consider its underlying short exact sequence in $X\Lcth_\bW$.
 The unit morphism $\cO_X\rarrow\cA$ of the quasi-coherent
quasi-algebra $\cA$ induces a natural (adjunction) morphism
$\Cohom_X(\cA,\gK)\rarrow\gK$ in $X\Lcth_\bW$.
 Consider the pushout diagram
$$
 \xymatrix{
  0 \ar[r] & \Cohom_X(\cA,\gK) \ar[r] \ar[d]
  & \Q \ar[rr]\ar[d] \ar@{-->}[ld]
  && \P \ar[rr] \ar@{=}[d] && 0 \\
  0 \ar[r] & \gK \ar[r] & \R \ar[rr] && \P \ar[rr] && 0
 }
$$
in the exact category $X\Lcth$.
 By assumption, we have $\Ext^{X,1}(\P,\gK)=0$, so the short
exact sequence $0\rarrow\gK\rarrow\R\rarrow\P\rarrow0$ in
$X\Lcth_\bW$ splits.
 Let $\R\rarrow\gK$ be the splitting morphism.
 The composition $\Q\rarrow\R\rarrow\gK$ provides the dotted diagonal
arrow on the diagram, making both the triangles on the left-hand side
of the diagram commutative.
 By the adjunction formula~\eqref{cohom-A-adjunction}, the morphism
$\Q\rarrow\gK$ in $X\Lcth_\bW$ corresponds to a morphism
$\Q\rarrow\Cohom_X(\cA,\gK)$ in $\cA\Lcth_\bW$.
 The latter morphism splits the original short exact sequence
$0\rarrow\Cohom_X(\cA,\gK)\rarrow\Q\rarrow\P\rarrow0$ in the exact
category $\cA\Lcth_\bW$, as desired.

 Part~(b): for a quasi-compact semi-separated scheme $X$, the condition
on the $\bW$\+locally contraherent $\cA$\+module $\P$ in the formulation
of the lemma just means that $\P$ is
$X$\+antilocal~\cite[Corollary~4.3.3(a)]{Pcosh}.

 Choose a finite affine open covering $X=\bigcup_\alpha U_\alpha$ of
the scheme $X$ subordinate to $\bW$, and denote by $j_\alpha\:U_\alpha
\rarrow X$ the open immersion morphisms.
 For any $\bW$\+locally contraherent $\cA$\+module $\Q$ on $X$,
the natural morphism $\bigoplus_\alpha j_\alpha{}_!j_\alpha^!\Q
\rarrow\Q$ (given by
the adjunction~\eqref{affine-open-immers-contrah-module-adjunction})
is an admissible epimorphism in $\cA\Lcth_\bW$, as one can see from
the \v Cech resolution~\eqref{lcth-sheaf-of-rings-cech-resolution}.
 The $\bW$\+locally contraherent (in fact, contraherent) $\cA$\+module
$\bigoplus_\alpha j_\alpha{}_!j_\alpha^!\Q$ on $X$ is $X$\+antilocal
according to the discussion in~\cite[Section~4.3]{Pcosh}.

 The class of all $X$\+antilocal ($\bW$\+locally) contraherent
$\cA$\+modules is closed under kernels of admissible epimorphisms
in $\cA\Lcth_\bW$, since the class of all antilocal contraherent
cosheaves on $X$ is closed under kernels of admissible epimorphisms
in $X\Lcth_\bW$ by~\cite[Corollary~4.3.3(b)]{Pcosh} (see
also~\cite[Corollaries~4.3.5(c) and~4.3.7]{Pcosh}).
 For all the reasons mentioned above, part~(b) follows from part~(a)
by~\cite[Lemma~7.1]{PS6}.
\end{proof}

\begin{cor} \label{X-antilocal-equivalent-to-A-antilocal}
 Let $X$ be a quasi-compact semi-separated scheme with an open
covering\/ $\bW$ and $\cA$ be a quasi-coherent quasi-algebra over~$X$.
 Let\/ $\P$ be a $\bW$\+locally contraherent $\cA$\+module on~$X$.
 Then the following conditions are equivalent:
\begin{enumerate}
\item $\P$ is $X$\+antilocal;
\item $\P$ is $\cA$\+antilocal;
\item $\Ext^{\cA,1}(\P,\gJ)=0$ for all $\cA$\+locally injective
$\bW$\+locally contraherent $\cA$\+modules\/ $\gJ$ on~$X$;
\item $\Ext^{\cA,n}(\P,\gJ)=0$ for all $\cA$\+locally injective
$\bW$\+locally contraherent $\cA$\+modules\/ $\gJ$ on $X$ and
all $n\ge1$.
\end{enumerate}
\end{cor}

\begin{proof}
 The implications (4)\,$\Longrightarrow$\,(3)\,$\Longrightarrow$\,(2)
are obvious (and hold for any scheme~$X$).
 The implication (2)\,$\Longrightarrow$\,(1) holds by
Lemma~\ref{A-antilocal-are-X-antilocal} (for any scheme~$X$).

 It remains to prove (1)\,$\Longrightarrow$\,(4).
 The argument is based on the \v Cech resolution and
Lemma~\ref{direct-image-of-A-lin-from-affine-Ext-orthogonality}(b).
 Once again, for a quasi-compact semi-separated scheme $X$,
condition~(1) implies the condition on $\P$ in the formulation of
Lemma~\ref{direct-image-of-A-lin-from-affine-Ext-orthogonality},
by~\cite[Corollary~4.3.3(a)]{Pcosh}.
 Let $X=\bigcup_{\alpha=1}^N U_\alpha$ be a finite affine open
covering of $X$ subordinate to~$\bW$.
 Then the \v Cech
resolution~\eqref{lcth-sheaf-of-rings-cech-resolution}
\begin{multline*}
 0\lrarrow j_{1,\dotsc,N}{}_!j_{1,\dotsc,N}^!\gJ \lrarrow\dotsb \\
 \lrarrow\bigoplus\nolimits_{1\le\alpha<\beta\le N}
 j_{\alpha,\beta}{}_!j_{\alpha,\beta}^!\gJ\lrarrow
 \bigoplus\nolimits_{\alpha=1}^N j_\alpha{}_!j_\alpha^!\gJ
 \lrarrow\gJ\lrarrow0
\end{multline*}
is a finite resolution of $\gJ$ by the direct images of
$\cA|_U$\+locally injective contraherent $\cA|_U$\+modules from
affine open subschemes $U=U_{\alpha_1}\cap\dotsb\cap U_{\alpha_k}
\subset X$ in $\cA\Lcth_\bW$.
 By Lemma~\ref{direct-image-of-A-lin-from-affine-Ext-orthogonality}(b),
we have $\Ext^{\cA,n}(\P,j_!j^!\gJ)=0$ for all $n\ge1$, where
$j\:U\rarrow X$ is the open immersion morphism.
 The desired Ext vanishing~(4) follows.
\end{proof}

\subsection{Antilocality of antilocal contraherent $\cA$-modules}
\label{antilocal-contrah-A-modules-subsecn}
 Let $X$ be a quasi-compact semi-separated scheme with an open covering
$\bW$ and $\cA$ be a quasi-coherent quasi-algebra over~$X$.

 By Corollary~\ref{X-antilocal-equivalent-to-A-antilocal}, the classes
of $\cA$\+antilocal and $X$\+antilocal $\bW$\+locally contraherent
$\cA$\+modules on $X$ coincide.
 Moreover, by~\cite[Corollaries~4.3.5(c) and~4.3.7]{Pcosh}, all
$X$\+antilocal $\bW$\+locally contraherent $\cA$\+modules on $X$
are (globally) contraherent, and the class of such $\cA$\+modules
does not depend on the open covering~$\bW$.

 Accordingly, from now on, in the context of quasi-compact
semi-separated schemes, we will speak simply of \emph{antilocal
contraherent $\cA$\+modules} instead of ``$\cA$\+antilocal
$\bW$\+locally contraherent $\cA$\+modules'' or ``$X$\+antilocal
$\bW$\+locally contraherent $\cA$\+modules''.

\begin{lem} \label{qcomp-qsep-A-loc-inj-preenvelope}
 Let $X$ be a quasi-compact semi-separated scheme with an open
covering\/ $\bW$ and $X=\bigcup_\alpha U_\alpha$ be a finite affine
open covering of $X$ subordinate to\/~$\bW$.
 Let $\cA$ be a quasi-coherent quasi-algebra over $X$ and\/ $\gM$ be
a\/ $\bW$\+locally contraherent $\cA$\+module on~$X$.
 Then there exists a short exact sequence of\/ $\bW$\+locally
contraherent $\cA$\+modules\/ $0\rarrow\gM\rarrow\gJ\rarrow\P\rarrow0$
on $X$ such that\/ $\gJ$ is an $\cA$\+locally injective\/
$\bW$\+locally contraherent $\cA$\+module and\/ $\P$ is a finitely
iterated extension of the direct images of contraherent
$\cA|_{U_\alpha}$\+modules from $U_\alpha$ in the exact category
$\cA\Lcth_\bW$ (or $\cA\Ctrh$).
\end{lem}

\begin{proof}
 This is a version of~\cite[Proposition~5.2]{Pal}
and~\cite[Lemma~4.3.2 or Proposition~B.6.4]{Pcosh}.
 The same construction, which is the locally contraherent cosheaf
(dual-analogous) version of the construction of
Lemma~\ref{qcomp-qsep-very-flaprojective-precover}, is applicable
in the situation at hand.
 The local nature of the $\cA$\+local injectivity property of
locally contraherent $\cA$\+modules (see
Section~\ref{A-loc-cotors-loc-inj-cosheaves-subsecn}) plays a key role.
 It is also important that the class of $\cA$\+locally injective
$\bW$\+locally contraherent $\cA$\+modules is preserved by
the direct images with respect to $(\bW,\bT)$\+affine open immersions of
schemes (see Sections~\ref{direct-images-of-A-co-sheaves-subsecn}\+-%
\ref{inverse-images-of-A-co-sheaves-subsecn}), as well as by extensions
and cokernels of admissible monomorphisms in $\cA\Lcth_\bW$.
 Of course, one also has to use the fact that every
$\cA(U_\alpha)$\+module is a submodule of an injective
$\cA(U_\alpha)$\+module.
\end{proof}

\begin{lem} \label{antilocal-modules-generating-class}
 Let $X$ be a quasi-compact semi-separated scheme with an open
covering\/ $\bW$ and $X=\bigcup_\alpha U_\alpha$ be a finite affine
open covering of $X$ subordinate to\/~$\bW$.
 Denote by $j_\alpha\:U_\alpha\rarrow X$ the open immersion morphisms.
 Let $\cA$ be a quasi-coherent quasi-algebra over $X$ and\/ $\gM$ be
a\/ $\bW$\+locally contraherent $\cA$\+module on~$X$.
 Then there exists an admissible epimorphism onto\/ $\gM$ in
$\cA\Lcth_\bW$ from a finite direct sum\/ $\bigoplus_\alpha
j_\alpha{}_!\gM_\alpha$ for some contraherent
$\cA|_{U_\alpha}$\+modules\/ $\gM_\alpha$ on~$U_\alpha$.
\end{lem}

\begin{proof}
 This was explained in the proof of
Lemma~\ref{direct-image-of-A-lin-from-affine-Ext-orthogonality}(b).
\end{proof}

\begin{rem} \label{Ext-unambiguous-remark}
 Now we can present the discussion of the notation
$\Ext^{\cA,*}({-},{-})$ that was promised in
Section~\ref{A-antilocal=X-antilocal-subsecn}
(cf.~\cite[Section~4.3]{Pcosh}).
 Let $X$ be a quasi-compact semi-separated scheme with an open
covering $\bW$ and $\cA$ be a quasi-coherent quasi-algebra over~$X$.
 First of all, the Ext groups in the category $\cA\Lcth_\bW$ do not
depend on the covering $\bW$ and agree with the Ext groups in
the whole category $\cA\Lcth$ of locally contraherent $\cA$\+modules
on~$X$.

 Indeed, the full exact subcategory $\cA\Lcth_\bW$ is closed under
extensions and kernels of admissible epimorphisms in $\cA\Lcth$
(by the definition of the exact categories of locally contraherent
$\cA$\+modules, this assertion holds just because it holds for the exact
categories of locally contraherent cosheaves of $\cO_X$\+modules;
see Sections~\ref{exact-categories-of-contrah-subsecn}
and~\ref{cosheaves-of-A-modules-subsecn}).
 Furthermore, every object of $\cA\Lcth$ is the target of an admissible
epimorphism from an object of $\cA\Ctrh\subset\cA\Lcth_\bW$ (and even
from an antilocal contraherent $\cA$\+module on~$X$)
by Lemma~\ref{antilocal-modules-generating-class}.
 So the dual version of~\cite[Theorem~12.1(b)]{Kel0} or the dual
version of~\cite[Proposition~13.2.2(i)]{KS} is applicable.
 See also~\cite[Proposition~A.2.1 or~A.3.1(a)]{Pcosh}
(or Proposition~\ref{infinite-resolutions} below) for further details.

 For the same reasons (up to inverting the arrows), the Ext groups
computed in the exact categories of $X$\+locally cotorsion,
$\cA$\+locally cotorsion, and $\cA$\+locally injective $\bW$\+locally
contraherent $\cA$\+modules $\cA\Lcth_\bW^{X\dlct}$, \
$\cA\Lcth_\bW^{\cA\dlct}$, and $\cA\Lcth_\bW^{\cA\dlin}$ (and also in
$\cA\Lcth^{X\dlct}$, \ $\cA\Lcth^{\cA\dlct}$, and
$\cA\Lcth^{\cA\dlin}$).
 Indeed, the full exact subcategories $\cA\Lcth_\bW^{X\dlct}$, \
$\cA\Lcth_\bW^{\cA\dlct}$, and $\cA\Lcth_\bW^{\cA\dlin}$ are closed
under extensions and cokernels of admissible monomorphisms in
$\cA\Lcth_\bW$, and we have just shown in
Lemma~\ref{qcomp-qsep-A-loc-inj-preenvelope} that any $\bW$\+locally
contraherent $\cA$\+module has an admissible monomorphism into
an $\cA$\+locally injective $\bW$\+locally contraherent $\cA$\+module
in $\cA\Lcth_\bW$.
 Furthermore, one has $\cA\Lcth_\bW^{\cA\dlin}\subset
\cA\Lcth_\bW^{\cA\dlct}\subset\cA\Lcth_\bW^{X\dlct}$ (see
Section~\ref{A-loc-cotors-loc-inj-cosheaves-subsecn}).

 When the quasi-coherent quasi-algebra $\cA$ is flat as a quasi-coherent
sheaf on $X$ with respect to its right $\cO_X$\+module structure,
the assertions of the previous paragraph also apply to the exact
categories of $X$\+locally injective locally contraherent $\cA$\+modules
$\cA\Lcth_\bW^{X\dlin}$ and $\cA\Lcth^{X\dlin}$, because one has
$\cA\Lcth_\bW^{\cA\dlin}\subset\cA\Lcth_\bW^{X\dlin}$ in this case.

 The conclusion of this discussion is that the notation
$\Ext^{\cA,*}({-},{-})$ for the Ext groups in the exact categories
of locally contraherent $\cA$\+modules is largely unambiguous in
the case of a quasi-compact semi-separated scheme~$X$.
\end{rem}

\begin{thm} \label{qcomp-qsep-antilocal-complete-cotorsion-pair-thm}
 Let $X$ be a quasi-compact semi-separated scheme with an open
covering\/ $\bW$ and $X=\bigcup_\alpha U_\alpha$ be a finite affine
open covering of $X$ subordinate to\/~$\bW$.
 Let $\cA$ be a quasi-coherent quasi-algebra over $X$ and\/ $\gM$ be
a\/ $\bW$\+locally contraherent $\cA$\+module on~$X$.
 In this context: \par
\textup{(a)} there exists a short exact sequence of\/ $\bW$\+locally
contraherent $\cA$\+modules\/ $0\rarrow\gJ'\rarrow\P\rarrow\gM\rarrow0$
on $X$ such that\/ $\P$ is an antilocal contraherent $\cA$\+module and\/
$\gJ'$ is an $\cA$\+locally injective\/ $\bW$\+locally contraherent
$\cA$\+module; \par
\textup{(b)} there exists a short exact sequence of\/ $\bW$\+locally
contraherent $\cA$\+modules\/ $0\rarrow\gM\rarrow\gJ\rarrow\P'\rarrow0$
on $X$ such that\/ $\gJ$ is an $\cA$\+locally injective\/ $\bW$\+locally
contraherent $\cA$\+module and\/ $\P'$ is an antilocal contraherent
$\cA$\+module; \par
\textup{(c)} a (\/$\bW$\+locally) contraherent $\cA$\+module on $X$ is
antilocal if and only if it is a direct summand of a finitely
iterated extension of the direct images of contraherent
$\cA|_{U_\alpha}$\+modules from $U_\alpha$ in the exact category
$\cA\Lcth_\bW$ or $\cA\Ctrh$.
\end{thm}

\begin{proof}
 Part~(b) follows immediately from
Lemma~\ref{qcomp-qsep-A-loc-inj-preenvelope} (because the class of
antilocal contraherent $\cA$\+modules is closed under extensions in
$\cA\Lcth_\bW$ by~\cite[Corollary~4.3.3(b)]{Pcosh} and preserved by
the direct images by~\cite[Section~4.3 or Corollary~4.6.3(a)]{Pcosh}).
 Part~(a) follows from part~(b) and
Lemma~\ref{antilocal-modules-generating-class} by virtue of the Salce
lemma~\cite{Sal}, \cite[Lemma~B.1.1(a)]{Pcosh} (cf.~\cite[proof
of Lemma~4.3.4]{Pcosh}).

 To prove part~(c), one needs to observe that the proof of part~(a)
actually provides a more precise result: in part~(a), one can choose
the short exact sequence so that $\P$ is a finitely iterated extension
of the direct images of contraherent $\cA|_{U_\alpha}$\+modules
from~$U_\alpha$.
 Now let $\Q$ be an antilocal contraherent $\cA$\+module on~$X$.
 Choose a short exact sequence $0\rarrow\gJ'\rarrow\P\rarrow\Q\rarrow0$,
where $\P$ is a finitely iterated extension of the direct images from
$U_\alpha$ as above, while $\gJ'$ is an $\cA$\+locally injective
$\bW$\+locally contraherent $\cA$\+module.
 Then $\Ext^{\cA,1}(\Q,\gJ')=0$ by
Corollary~\ref{X-antilocal-equivalent-to-A-antilocal}, hence $\Q$ is
a direct summand of $\P$, as desired.
\end{proof}

\begin{cor} \label{qcomp-qsep-antiloc-A-lin-is-a-cotorsion-pair-cor}
 Let $X$ be a quasi-compact semi-separated scheme with an open
covering\/ $\bW$ and $X=\bigcup_\alpha U_\alpha$ be a finite affine
open covering of $X$ subordinate to\/~$\bW$.
 Let $\cA$ be a quasi-coherent quasi-algebra over~$X$.
 Then a\/ $\bW$\+locally contraherent $\cA$\+module\/ $\gJ$ is
$\cA$\+locally injective if and only if\/ $\Ext^{\cA,1}(\P,\gJ)=0$
for all antilocal contraherent $\cA$\+modules\/~$\P$.
\end{cor}

\begin{proof}
 The ``only if'' assertion is provided by
Corollary~\ref{X-antilocal-equivalent-to-A-antilocal}.
 The ``if'' implication follows from the same corollary and
Theorem~\ref{qcomp-qsep-antilocal-complete-cotorsion-pair-thm}(b)
by the argument dual to the proof of
Theorem~\ref{qcomp-qsep-antilocal-complete-cotorsion-pair-thm}(c);
see~\cite[Lemma~B.1.2]{Pcosh}.
 The point is that the class of $\cA$\+locally injective
$\bW$\+locally contraherent $\cA$\+modules is closed under direct
summands in $\cA\Lcth_\bW$.
\end{proof}

\begin{rem} \label{antilocal-A-modules-concluding-remark}
 The assertions of
Theorem~\ref{qcomp-qsep-antilocal-complete-cotorsion-pair-thm}(a\+-b)
and Corollaries~\ref{X-antilocal-equivalent-to-A-antilocal}
and~\ref{qcomp-qsep-antiloc-A-lin-is-a-cotorsion-pair-cor} can be
summarized by saying that the pair of classes (antilocal contraherent
$\cA$\+modules, $\cA$\+locally injective $\bW$\+locally contraherent
$\cA$\+modules) is a hereditary complete cotorsion pair in the exact
category $\cA\Lcth_\bW$.
 On the other hand,
Theorem~\ref{qcomp-qsep-antilocal-complete-cotorsion-pair-thm}(c)
is an extra piece of information which one obtains from our proof of
Theorem~\ref{qcomp-qsep-antilocal-complete-cotorsion-pair-thm}(a\+-b).

 The assertion of
Theorem~\ref{qcomp-qsep-antilocal-complete-cotorsion-pair-thm}(c)
can be expressed by saying that the class of all antilocal contraherent
$\cA$\+modules on a quasi-compact semi-separated scheme $X$ is
(just as the terminology suggests) \emph{antilocal} in the sense
of~\cite[Section~4]{Pal}.
 Together with the claim that every $\bW$\+locally contraherent
$\cA$\+module on $X$ is an admissible subobject of an $\cA$\+locally
injective one (a weak version of
Lemma~\ref{qcomp-qsep-A-loc-inj-preenvelope} or
Theorem~\ref{qcomp-qsep-antilocal-complete-cotorsion-pair-thm}(b)),
this is the main result of this
Section~\ref{antilocal-contrah-A-modules-subsecn} for our purposes
in this paper.
\end{rem}

\subsection{Antilocality of antilocal $X$-locally cotorsion
contraherent $\cA$-modules}
\label{antilocal-contrah-X-lct-A-modules-subsecn}
 Let $X$ be a quasi-compact semi-separated scheme and $\cA$ be
a quasi-coherent quasi-algebra over~$X$.
 By an ``antilocal $X$\+locally cotorsion contraherent $\cA$\+module''
we mean a ($\bW$\+locally) contraherent $\cA$\+module that is
simultaneously antilocal and $X$\+locally cotorsion.

\begin{lem} \label{qcomp-qsep-X-lct-A-loc-inj-preenvelope}
 Let $X$ be a quasi-compact semi-separated scheme with an open
covering\/ $\bW$ and $X=\bigcup_\alpha U_\alpha$ be a finite affine
open covering of $X$ subordinate to\/~$\bW$.
 Let $\cA$ be a quasi-coherent quasi-algebra over $X$ and\/ $\gM$ be
an $X$\+locally cotorsion\/ $\bW$\+locally contraherent $\cA$\+module
on~$X$.
 Then there exists a short exact sequence of $X$\+locally cotorsion\/
$\bW$\+locally contraherent $\cA$\+modules\/ $0\rarrow\gM\rarrow\gJ
\rarrow\P\rarrow0$ on $X$ such that\/ $\gJ$ is an $\cA$\+locally
injective\/ $\bW$\+locally contraherent $\cA$\+module and\/ $\P$ is
a finitely iterated extension of the direct images of
$U_\alpha$\+locally cotorsion contraherent $\cA|_{U_\alpha}$\+modules
in the exact category $\cA\Lcth_\bW^{X\dlct}$ (or $\cA\Ctrh^{X\dlct}$).
\end{lem}

\begin{proof}
 Similar to the proof of Lemma~\ref{qcomp-qsep-A-loc-inj-preenvelope}.
 Notice that the class of $X$\+locally cotorsion $\bW$\+locally
contraherent $\cA$\+modules is preserved by the direct images with
respect to $(\bW,\bT)$\+affine open immersions of schemes
(see Section~\ref{direct-images-of-O-co-sheaves-subsecn} and
the end of Section~\ref{inverse-images-of-A-co-sheaves-subsecn}).
\end{proof}

\begin{lem} \label{antilocal-X-lct-modules-generating-class}
 Let $X$ be a quasi-compact semi-separated scheme with an open
covering\/ $\bW$ and $X=\bigcup_\alpha U_\alpha$ be a finite affine
open covering of $X$ subordinate to\/~$\bW$.
 Denote by $j_\alpha\:U_\alpha\rarrow X$ the open immersion morphisms.
 Let $\cA$ be a quasi-coherent quasi-algebra over $X$ and\/ $\gM$ be
an $X$\+locally cotorsion\/ $\bW$\+locally contraherent $\cA$\+module
on~$X$.
 Then there exists an admissible epimorphism onto\/ $\gM$ in
$\cA\Lcth_\bW^{X\dlct}$ from a finite direct sum\/ $\bigoplus_\alpha
j_\alpha{}_!\gM_\alpha$ for some $U_\alpha$\+locally cotorsion
contraherent $\cA|_{U_\alpha}$\+modules\/~$\gM_\alpha$.
\end{lem}

\begin{proof}
 This is similar to Lemma~\ref{antilocal-modules-generating-class}.
\end{proof}

\begin{thm} \label{qcomp-qsep-X-lct-aloc-complete-cotorsion-pair-thm}
 Let $X$ be a quasi-compact semi-separated scheme with an open
covering\/ $\bW$ and $X=\bigcup_\alpha U_\alpha$ be a finite affine
open covering of $X$ subordinate to\/~$\bW$.
 Let $\cA$ be a quasi-coherent quasi-algebra over $X$ and\/ $\gM$ be
an $X$\+locally cotorsion\/ $\bW$\+locally contraherent $\cA$\+module
on~$X$.
 In this context: \par
\textup{(a)} there exists a short exact sequence of $X$\+locally
cotorsion\/ $\bW$\+locally contraherent $\cA$\+modules\/ $0\rarrow
\gJ'\rarrow\P\rarrow\gM\rarrow0$ on $X$ such that\/ $\P$ is an antilocal
$X$\+locally cotorsion contraherent $\cA$\+module and\/ $\gJ'$ is
an $\cA$\+locally injective\/ $\bW$\+locally contraherent $\cA$\+module;
\par
\textup{(b)} there exists a short exact sequence of $X$\+locally
cotorsion\/ $\bW$\+locally contraherent $\cA$\+modules\/ $0\rarrow\gM
\rarrow\gJ\rarrow\P'\rarrow0$ on $X$ such that\/ $\gJ$ is
an $\cA$\+locally injective\/ $\bW$\+locally contraherent $\cA$\+module
and\/ $\P'$ is an antilocal $X$\+locally cotorsion contraherent
$\cA$\+module; \par
\textup{(c)} an $X$\+locally cotorsion (\/$\bW$\+locally) contraherent
$\cA$\+module on $X$ is antilocal if and only if it is a direct summand
of a finitely iterated extension of the direct images of
$U_\alpha$\+locally cotorsion contraherent $\cA|_{U_\alpha}$\+modules
in the exact category $\cA\Lcth_\bW^{X\dlct}$ or $\cA\Ctrh^{X\dlct}$.
\end{thm}

\begin{proof}
 Similar to
Theorem~\ref{qcomp-qsep-antilocal-complete-cotorsion-pair-thm}, and
based on Corollary~\ref{X-antilocal-equivalent-to-A-antilocal}
and Lemmas~\ref{qcomp-qsep-X-lct-A-loc-inj-preenvelope}\+-%
\ref{antilocal-X-lct-modules-generating-class}.
\end{proof}

\begin{cor} \label{qcomp-qsep-X-lct-aloc-A-lin-is-a-cotorsion-pair-cor}
 Let $X$ be a quasi-compact semi-separated scheme with an open
covering\/ $\bW$ and $X=\bigcup_\alpha U_\alpha$ be a finite affine
open covering of $X$ subordinate to\/~$\bW$.
 Let $\cA$ be a quasi-coherent quasi-algebra over~$X$.
 Then an $X$\+locally cotorsion\/ $\bW$\+locally contraherent
$\cA$\+module\/ $\gJ$ is $\cA$\+locally injective if and only if\/
$\Ext^{\cA,1}(\P,\gJ)=0$ for all antilocal $X$\+locally cotorsion
contraherent $\cA$\+modules\/~$\P$.
\end{cor}

\begin{proof}
 Similar to
Corollary~\ref{qcomp-qsep-antiloc-A-lin-is-a-cotorsion-pair-cor}.
\end{proof}

\begin{rem} \label{antilocal-X-lct-A-modules-concluding-remark}
 The assertions of
Theorem~\ref{qcomp-qsep-X-lct-aloc-complete-cotorsion-pair-thm}(a\+-b)
together with
Corollaries~\ref{qcomp-qsep-X-lct-aloc-A-lin-is-a-cotorsion-pair-cor}
and~\ref{X-antilocal-equivalent-to-A-antilocal}
(with Remark~\ref{Ext-unambiguous-remark}) can be summarized by saying
that the pair of classes (antilocal $X$\+locally cotorsion contraherent
$\cA$\+modules, $\cA$\+locally injective $\bW$\+locally contraherent
$\cA$\+modules) is a hereditary complete cotorsion pair in the exact
category $\cA\Lcth_\bW^{X\dlct}$.
 This assertion does not need a separate proof: it can be deduced from
the similar result for the exact category $\cA\Lcth_\bW$ as per
Remark~\ref{antilocal-A-modules-concluding-remark} by applying
the lemma about restricting (hereditary) complete cotorsion pairs to
exact subcategories~\cite[Lemmas~1.5(b) and~1.6]{Pal},
\cite[Lemmas~B.1.3 and~B.1.4(b)]{Pcosh}.

 On the other hand,
Theorem~\ref{qcomp-qsep-X-lct-aloc-complete-cotorsion-pair-thm}(c)
is an extra piece of information which one obtains from the proof of
Theorem~\ref{qcomp-qsep-X-lct-aloc-complete-cotorsion-pair-thm}(a\+-b)
sketched above.
 Theorem~\ref{qcomp-qsep-X-lct-aloc-complete-cotorsion-pair-thm}(c)
does \emph{not} follow from
Theorem~\ref{qcomp-qsep-antilocal-complete-cotorsion-pair-thm}(c) (or
any other results of Sections~\ref{A-antilocal=X-antilocal-subsecn}\+-%
\ref{antilocal-contrah-A-modules-subsecn}) in any simple way;
instead, its proof requires going through all the rounds of
the argument again in the new context, as described in the present
Section~\ref{antilocal-contrah-X-lct-A-modules-subsecn}.

 The assertion of
Theorem~\ref{qcomp-qsep-X-lct-aloc-complete-cotorsion-pair-thm}(c)
can be expressed by saying that the class of all antilocal
$X$\+locally cotorsion contraherent $\cA$\+modules on a quasi-compact
semi-separated scheme $X$ is \emph{antilocal} in the sense
of~\cite[Section~4]{Pal}.
 This is the main result of this
Section~\ref{antilocal-contrah-X-lct-A-modules-subsecn} for our purposes
in this paper.
\end{rem}

\subsection{Antilocality of antilocal $\cA$-locally cotorsion
contraherent $\cA$-modules}
\label{antilocal-contrah-A-lct-A-modules-subsecn}
 Let $X$ be a quasi-compact semi-separated scheme with an open
covering $\bW$ and $\cA$ be a quasi-coherent quasi-algebra over~$X$.
 By an ``antilocal $\cA$\+locally cotorsion contraherent $\cA$\+module''
we mean a ($\bW$\+locally) contraherent $\cA$\+module that is
simultaneously antilocal and $\cA$\+locally cotorsion.

\begin{rem} \label{antilocal-contrah-A-lct-A-modules-remark}
 Everything that was said in the previous
Section~\ref{antilocal-contrah-X-lct-A-modules-subsecn} remains valid
with the $X$\+locally cotorsion condition on $\bW$\+locally
contraherent $\cA$\+modules replaced by the $\cA$\+locally cotorsion
condition.

 It suffices to replace the words ``$X$\+locally cotorsion'' by
the words ``$\cA$\+locally cotorsion'', the words
``$U_\alpha$\+locally cotorsion'' by the words
``$\cA|_{U_\alpha}$\+locally cotorsion'', and the upper index
``$X\dlct$'' with the upper index ``$\cA\dlct$'' in the notation for
the exact categories throughout
the whole Section~\ref{antilocal-contrah-X-lct-A-modules-subsecn}
in order to obtain a correct exposition of the theory of antilocality
of antilocal $\cA$\+locally cotorsion contraherent $\cA$\+modules
and the related cotorsion pair.
 One can also point out that the class of $\cA$\+locally cotorsion
$\bW$\+locally contraherent $\cA$\+modules is preserved by
the direct images with respect to $(\bW,\bT)$\+affine open immersions
of schemes, as explained in
Section~\ref{direct-images-of-A-co-sheaves-subsecn}.
 To avoid unnecessary repetition, we mostly restrict ourselves to this
short remark in the present section.
\end{rem}

 The following corollary, which is an $\cA$\+locally cotorsion version
of~\cite[Corollary~4.6.3(a\+-b)]{Pcosh}, illustrates the power of
the main result obtained in this way (viz., the $\cA$\+locally
cotorsion version of
Theorem~\ref{qcomp-qsep-X-lct-aloc-complete-cotorsion-pair-thm}(c)).

 Let us introduce some notation first.
 The full exact subcategory of antilocal contraherent cosheaves on
a scheme~$X$ (defined in~\cite[Section~4.3]{Pcosh}, see
Section~\ref{A-antilocal=X-antilocal-subsecn}) will be denoted by
$X\Ctrh_\al\subset X\Ctrh$.

 We will denote the full subcategory of antilocal contraherent
$\cA$\+modules (defined in
Sections~\ref{A-antilocal=X-antilocal-subsecn}\+-%
\ref{antilocal-contrah-A-modules-subsecn}) by
$\cA\Ctrh_\al\subset\cA\Lcth_\bW$.
 The full subcategory of antilocal $X$\+locally cotorsion contraherent
$\cA$\+modules (defined in
Section~\ref{antilocal-contrah-X-lct-A-modules-subsecn}) will be denoted
by $\cA\Ctrh_\al^{X\dlct}\subset\cA\Lcth_\bW^{X\dlct}$.
 The full subcategory of antilocal $\cA$\+locally cotorsion contraherent
$\cA$\+modules (defined above in this section) will be denoted
by $\cA\Ctrh_\al^{\cA\dlct}\subset\cA\Lcth_\bW^{\cA\dlct}$.
 The full subcategories $\cA\Ctrh_\al^{\cA\dlct}\subset
\cA\Ctrh_\al^{X\dlct}\subset\cA\Ctrh_\al$ inherit exact category
structures from $\cA\Lcth_\bW^{\cA\dlct}\subset\cA\Lcth_\bW^{X\dlct}
\subset\cA\Lcth_\bW$. {\emergencystretch=1em\par}

\begin{cor} \label{contrah-al-A-lct-direct-image}
 Let $X$ be a quasi-compact semi-separated scheme and $Y$ be
a quasi-compact open subscheme in $X$ with the open immersion morphism
$f\:Y\rarrow X$.
 Let $\cA$ be a quasi-coherent quasi-algebra over~$X$.
 Then the direct image functor $f_!\:(Y,\cA|_Y)\Cosh\rarrow
(X,\cA)\Cosh$ takes antilocal $\cA|_Y$\+locally cotorsion contraherent
cosheaves on $Y$ to antilocal $\cA$\+locally cotorsion contraherent
cosheaves on~$X$.
 The resulting functor of exact categories
$f_!\:\cA|_Y\Ctrh^{\cA|_Y\dlct}_\al\rarrow\cA\Ctrh^{\cA\dlct}_\al$
is exact.
\end{cor}

\begin{proof}
 Notice first of all that the direct image functor
$f_!\:(Y,\cO_Y)\Cosh\rarrow(X,\cO_X)\Cosh$ takes antilocal contraherent
cosheaves on $Y$ to antilocal contraherent cosheaves on $X$, and
the resulting functor $f_!\:Y\Ctrh_\al\rarrow X\Ctrh_\al$ is exact
(for any morphism of quasi-compact semi-separated schemes
$f\:Y\rarrow X$) by~\cite[Corollary~4.6.3(a)]{Pcosh}, which was
already mentioned in Section~\ref{A-antilocal=X-antilocal-subsecn}.
 On the other hand, for an \emph{affine} open immersion morphism
$f\:Y\rarrow X$, the functor $f_!\:(Y,\cA|_Y)\Cosh\rarrow
(X,\cA)\Cosh$ takes all (and not only antilocal) $\cA|_Y$\+locally
cotorsion contraherent cosheaves on $Y$ to $\cA$\+locally cotorsion
contrarent cosheaves on $X$
by Lemma~\ref{restriction-coextension-injective-cotorsion}(a),
as mentioned in Section~\ref{direct-images-of-A-co-sheaves-subsecn}.
 This suffices to prove the corollary in the case when the open
immersion morphism~$f$ is affine.

 In the full generality of the assumptions of the corollary, it follows
from~\cite[Corollary~4.6.3(a)]{Pcosh} that the functor
$f_!\:(Y,\cA|_Y)\Cosh\rarrow(X,\cA)\Cosh$ takes antilocal contraherent
$\cA|_Y$\+modules on $Y$ to antilocal contraherent $\cA$\+modules on
$X$, and the resulting functor $f_!\:\cA|_Y\Ctrh_\al\rarrow\cA\Ctrh_\al$
is exact.
 By Theorem~\ref{qcomp-qsep-X-lct-aloc-complete-cotorsion-pair-thm}(c)
and Remark~\ref{antilocal-contrah-A-lct-A-modules-remark}, all
the antilocal $\cA|_Y$\+locally cotorsion contraherent $\cA|_Y$\+modules
on $Y$ are direct summands of finitely iterated extensions of the direct
images of $\cA|_V$\+locally cotorsion contraherent $\cA|_V$\+modules
from affine open subschemes $V\subset Y$.
 One can consider such finitely iterated extensions as happening in
the exact category of antilocal $\cA|_Y$\+modules $\cA|_Y\Ctrh_\al$.
 The functor $f_!\:\cA|_Y\Ctrh_\al\rarrow\cA\Ctrh_\al$ preserves
finitely iterated extensions, since it is exact.
 The assertions of the corollary follow.
\end{proof}

 For a discussion of a yet another cotorsion pair in $\cA\Lcth_\bW$
somewhat similar to those described in
Sections~\ref{antilocal-contrah-A-modules-subsecn}\+-%
\ref{antilocal-contrah-A-lct-A-modules-subsecn}, see
Section~\ref{antilocally-flaprojective-subsecn} below.

\Section{Na\"\i ve Co-Contra Correspondence for $\cA$-Modules}
\label{naive-co-contra-secn}

 The aim of this section is to construct a triangulated equivalence
$\sD(\cA\Qcoh)\simeq\sD(\cA\Lcth_\bW)$ between the derived category
of the abelian category of quasi-coherent sheaves of $\cA$\+modules and
the derived category of the exact category of $\bW$\+locally
contraherent cosheaves of $\cA$\+modules for a quasi-compact
semi-separated scheme $X$ with an open covering $\bW$ and
a quasi-coherent quasi-algebra $\cA$ over~$X$.

\subsection{Resolution and coresolution dimensions}
\label{co-resolution-dimensions-subsecn}
 Let $\sE$ be an exact category (in the sense of Quillen~\cite{Bueh}).
 A full subcategory $\sF\subset\sE$ is said to be \emph{generating}
if every object of $\sE$ is an admissible quotient object of
an object from~$\sF$.
 A generating full subcategory $\sF\subset\sE$ is said to be
\emph{resolving} if $\sF$ is closed under extensions and kernels of
admissible epimorphisms in~$\sE$.

 Dually, a full subcategory $\sC\subset\sE$ is said to be
\emph{cogenerating} if every object of $\sE$ is an admissible
suboboject of and object from~$\sF$.
 A cogenerating full subcategory $\sC\subset\sE$ is said to be
\emph{coresolving} if $\sC$ is closed under extensions and cokernels
of admissible subobjects in~$\sE$.

 Any (co)resolving full subcategory in $\sE$ is, by definition,
closed under extensions; so it inherits an exact category structure.
 We will consider the (co)resolving full subcategories in $\sE$ as
exact categories with the inherited exact category structure.

 Given an additive category $\sA$, we denote by $\sK^\st(\sA)$
the category of complexes in $\sA$ with morphisms up to the cochain
homotopy.
 We refer to~\cite[Section~10.4]{Bueh} for the definition of
the derived category $\sD^\st(\sE)$ of an exact category~$\sE$.
 Here $\st$~is one of the conventional derived category symbols
$\st=\bb$, $+$, $-$, or~$\varnothing$.

\begin{prop} \label{infinite-resolutions}
 Let\/ $\sF$ be a resolving subcategory in an exact category\/ $\sE$.
 Then the inclusion of the exact categories\/ $\sF\rarrow\sE$ induces
an equivalence of the derived categories
$$
 \sD^-(\sF)\simeq\sD^-(\sE).
$$
\end{prop}

\begin{proof}
 This is the dual version of~\cite[Theorem~12.1]{Kel0},
or the dual version of~\cite[Proposition~13.2.2(i)]{KS},
or~\cite[Proposition~A.3.1(a)]{Pcosh}.
\end{proof}

 The classical homological algebra concepts of projective, injective,
and flat dimensions can be extended to exact categories with
(co)resolving subcategories as follows.
 Assume that the exact category $\sE$ is weakly
idempotent-complete~\cite[Section~7]{Bueh}; then so is any
(co)resolving subcategory in~$\sE$.

 Let $\sF\subset\sE$ be a resolving subcategory and $d\ge-1$ be
an integer.
 An object $E\in\sE$ is said to have \emph{resolution dimension~$\le d$}
with respect to $\sF$ if there exists an exact sequence
$0\rarrow F_d\rarrow F_{d-1}\rarrow\dotsb\rarrow F_0\rarrow E\rarrow0$
in~$\sE$.

 Dually, let $\sC\subset\sE$ be a coresolving subcategory.
 An object $E\in\sE$ is said to have \emph{coresolution
dimension~$\le d$} with respect to $\sC$ if there exists an exact
sequence $0\rarrow E\rarrow C^0\rarrow C^1\rarrow\dotsb\rarrow C^d
\rarrow0$ in~$\sE$.

 Just as the projective, injective, and flat dimensions,
the (co)resolution dimension does not depend on the choice of
a (co)resolution.
 A precise formulation and proof this result can be found,
on various generality levels, in~\cite[Lemma~2.1]{Zhu},
\cite[Proposition~2.3(1)]{Sto}, or~\cite[Corollary~A.5.2]{Pcosh}.

\begin{prop} \label{finite-resolutions}
 Let\/ $\sF$ be a resolving subcategory in a weakly idempotent-complete
exact category\/ $\sE$ such that the resolution dimensions of all
the objects of\/ $\sE$ with respect to\/ $\sF$ are finite and bounded
by a fixed constant~$d$.
 Then, for any conventional derived category symbol\/ $\st=\bb$, $+$,
$-$, or\/~$\varnothing$, the inclusion of exact categories\/
$\sF\rarrow\sE$ induces an equivalence of the derived categories
$$
 \sD^\st(\sF)\simeq\sD^\st(\sE).
$$
\end{prop}

\begin{proof}
 The case $\st=-$ is covered by the more general
Proposition~\ref{infinite-resolutions}.
 The case $\st=\bb$ is the dual version
of~\cite[Proposition~13.2.2(ii)]{KS}.
 The case $\st=\varnothing$ is~\cite[Proposition~5.14]{Sto}.
 The result of~\cite[Proposition~A.5.8]{Pcosh} covers all
the cases.
\end{proof}

\subsection{Cotorsion periodicity}
 \emph{Periodicity theorems} in homological algebra are linked to
properties of the objects of cocycles in acyclic complexes.
 These results and their corollaries describe the phenomena when certain
(co)resolving subcategories exhibit behavior typical of such
subcategories with finite (co)resolution dimensions, even though
the (co)resolution dimensions are actually infinite.
 We refer to the papers~\cite{BHP,PS6} and
the preprint~\cite[Section~7]{Pcosh} for a general discussion of
the periodicity theorems from the perspective relevant to
the present paper.

 The following \emph{cotorsion periodicity theorem} for modules over
associative rings is due to Bazzoni, Cort\'es-Izurdiaga, and
Estrada~\cite[Theorem~1.2(2), Proposition~4.8(2), or
Theorem~5.1(2)]{BCE}.
 The definition of a cotorsion $R$\+module can be found in
Section~\ref{prelim-cotorsion-subsecn}.

\begin{thm} \label{module-cotorsion-periodicity}
 Let $R$ be an associative ring and $C^\bu$ be an acyclic complex
of $R$\+modules.
 Assume that all the terms $C^n$, \,$n\in\boZ$, of the complex $C^\bu$
are cotorsion $R$\+modules.
 Then all the $R$\+modules of cocycles of the acyclic complex $C^\bu$
are also cotorsion $R$\+modules. \qed
\end{thm}

 The following theorem is a quasi-coherent sheaf version of
the preceding one.
 The definition of a cotorsion quasi-coherent sheaf on a scheme can be
found in Section~\ref{antilocality-of-X-cotorsion-subsecn}.

\begin{thm} \label{qcoh-cotorsion-periodicity}
 Let $X$ be a quasi-compact semi-separated scheme and $\C^\bu$ be
an acyclic complex of quasi-coherent sheaves on~$X$.
 Assume that all the terms $\C^n$, \,$n\in\boZ$, of the complex $\C^\bu$
are cotorsion quasi-coherent sheaves.
 Then all the quasi-coherent sheaves of cocycles of the acyclic complex
$\C^\bu$ are also cotorsion.
\end{thm}

\begin{proof}
 This is~\cite[Corollary~10.4]{PS6}.
 A version of this approach was developed by Estrada with collaborators;
see~\cite[Remark~10.3]{PS6} for a discussion.
 Another proof can be found in~\cite[Section~8]{Pphil}.
\end{proof}

 The next theorem is a common generalization of
Theorems~\ref{module-cotorsion-periodicity}
and~\ref{qcoh-cotorsion-periodicity}.
 The definition of an $\cA$\+cotorsion quasi-coherent $\cA$\+module
can be found in Section~\ref{antilocality-of-A-cotorsion-subsecn}.

\begin{thm} \label{qcoh-quasi-algebra-cotorsion-periodicity}
 Let $X$ be a quasi-compact semi-separated scheme and $\cA$ be
a quasi-coherent quasi-algebra over~$X$.
 Let $\C^\bu$ be an acyclic complex of quasi-coherent $\cA$\+modules
on~$X$.
 Assume that all the terms $\C^n$, \,$n\in\boZ$, of the complex $\C^\bu$
are $\cA$\+cotorsion quasi-coherent $\cA$\+modules.
 Then all the quasi-coherent $\cA$\+modules of cocycles of the acyclic
complex $\C^\bu$ are also $\cA$\+cotorsion.
\end{thm}

\begin{proof}
 The argument from~\cite[proof of Corollary~10.4]{PS6}, based as it is
on numerous preceding results of~\cite{PS6} and the category-theoretic
references therein, works in the situation at hand without any changes.
 The result of
Lemma~\ref{quasi-algebra-adjustedness-co-locality}(a), establishing
the local nature of the class of flat quasi-coherent $\cA$\+modules,
is helpful.
 One proves that all flat quasi-coherent $\cA$\+modules are filtered
inductive limits of locally countably presented flat quasi-coherent
$\cA$\+modules, shows that the latter ones have finite projective
dimensions in $\cA\Qcoh$, etc.
 We omit the straightforward but tedious details.
 The contraherent cosheaf-based argument from~\cite[Section~8]{Pphil}
can be also adopted to the situation at hand, using the next
Corollary~\ref{loc-contrah-quasi-algebra-cotorsion-periodicity},
the constructions of
Sections~\ref{fHom-over-qcoh-quasi-algebra-subsecn}
and~\ref{contratensor-over-qcoh-quasi-algebra-subsecn},
and Lemma~\ref{quasi-algebra-underived-naive-co-contra} below.
\end{proof}

 The locally contraherent version of the preceding results follows
from the module-theoretic version.
 The definitions of $X$\+locally cotorsion and $\cA$\+locally
cotorsion $\bW$\+locally contraherent $\cA$\+modules can be found
in Sections~\ref{cosheaves-of-A-modules-subsecn}
and~\ref{A-loc-cotors-loc-inj-cosheaves-subsecn}.

\begin{cor} \label{loc-contrah-quasi-algebra-cotorsion-periodicity}
 Let $X$ be a scheme with an open covering\/ $\bW$, and let $\cA$ be
a quasi-coherent quasi-algebra over~$X$.
 Let\/ $\gC^\bu$ be an acyclic complex in the exact category\/
$\cA\Lcth_\bW$.
 In this context: \par
\textup{(a)} If all the terms\/ $\gC^n$, \,$n\in\boZ$, of the complex\/
$\gC^\bu$ are $X$\+locally cotorsion\/ $\bW$\+locally contraherent
$\cA$\+modules, then the complex\/ $\gC^\bu$ is acyclic in the exact
category of $X$\+locally cotorsion\/ $\bW$\+locally contraherent
$\cA$\+modules $\cA\Lcth_\bW^{X\dlct}$. \par
\textup{(b)} If all the terms\/ $\gC^n$, \,$n\in\boZ$, of the complex\/
$\gC^\bu$ are $\cA$\+locally cotorsion\/ $\bW$\+locally contraherent
$\cA$\+modules, then the complex\/ $\gC^\bu$ is acyclic in the exact
category of $\cA$\+locally cotorsion\/ $\bW$\+locally contraherent
$\cA$\+modules $\cA\Lcth_\bW^{\cA\dlct}$.
\end{cor}

\begin{proof}
 Both the assertions follow immediately from
Theorem~\ref{module-cotorsion-periodicity} (applied to the rings
$R=\cO_X(U)$ in the case~(a) and to the rings $R=\cA(U)$ in
the case~(b), where $U$ ranges over the affine open subschemes
$U\subset X$ subordinate to~$\bW$).
\end{proof}

\subsection{Background derived equivalences for quasi-coherent
$\cA$-modules} \label{background-derived-quasi-coherent}
 In this section we prove $\cA$\+module generalizations of some
results about the conventional derived categories of quasi-coherent
sheaves from~\cite[Section~4.7]{Pcosh}.

 Let $X$ be a scheme and $\cA$ be a quasi-coherent quasi-algebra
over~$X$.
 We will denote the full subcategory of $X$\+contraadjusted
quasi-coherent $\cA$\+modules (defined in
Section~\ref{antilocality-of-X-contraadjusted-subsecn})
by $\cA\Qcoh^{X\dcta}\subset\cA\Qcoh$.
 The full subcategory of $X$\+cotorsion quasi-coherent $\cA$\+modules
(defined in Section~\ref{antilocality-of-X-cotorsion-subsecn}) will be
denoted by $\cA\Qcoh^{X\dcot}\subset\cA\Qcoh$.
 The full subcategory of $\cA$\+cotorsion quasi-coherent $\cA$\+modules
(defined in Section~\ref{antilocality-of-A-cotorsion-subsecn}) will be
denoted by $\cA\Qcoh^{\cA\dcot}\subset\cA\Qcoh$.

\begin{lem} \label{X-A-cta-cot-coresolving-coresolution-dimension}
 Let $X$ be a quasi-compact semi-separated scheme with an affine
open covering $X=\bigcup_{\alpha=1}^N U_\alpha$ and $\cA$ be
a quasi-coherent quasi-algebra over~$X$.
 Then \par
\textup{(a)} the full subcategories of $X$\+contraadjusted,
$X$\+cotorsion, and $\cA$\+cotorsion quasi-coherent sheaves
$\cA\Qcoh^{X\dcta}$, $\cA\Qcoh^{X\dcot}$, and $\cA\Qcoh^{\cA\dcot}$
are coresolving in the abelian category $\cA\Qcoh$; \par
\textup{(b)} all quasi-coherent $\cA$\+modules have finite coresolution
dimensions not exceeding $N$ with respect to the coresolving
subcategory $\cA\Qcoh^{X\dcta}\subset\cA\Qcoh$.
\end{lem}

\begin{proof}
 Part~(a): recall that one has $\cA\Qcoh^{\cA\dcot}\subset
\cA\Qcoh^{X\dcot}\subset\cA\Qcoh^{X\dcta}$, as explained in
Section~\ref{antilocality-of-A-cotorsion-subsecn}.
 Clearly, all injective quasi-coherent $\cA$\+modules are
$\cA$\+cotorsion.
 So all the three full subcategories in question are cogenerating
in $\cA\Qcoh$.
 All the three full subcategories are also closed under extensions
(by the definition) and under cokernels of monomorphisms
(by~\cite[Corollary~4.1.2(c)]{Pcosh}, \cite[Corollary~4.1.9(c)]{Pcosh},
and Corollary~\ref{A-cotorsion-hereditariness-properties}(c) above,
respectively).

 Part~(b): the assertion follows from the facts that the coresolution
dimension does not depend on the choice of a coresolution
(see Section~\ref{co-resolution-dimensions-subsecn}) and that all
quasi-coherent sheaves on $X$ have finite coresolution dimensions
not exceeding $N$ with respect to the coresolving subcategory
of contraadjusted quasi-coherent sheaves $X\Qcoh^\cta\subset
X\Qcoh$ \,\cite[Lemma~4.7.1(b)]{Pcosh}.
\end{proof}

 For a discussion of the coresolution dimension with respect to
the coresolving subcategory $\cA\Qcoh^{X\dcot}\subset\cA\Qcoh$ under
more restrictive assumptions on the scheme $X$, see
Lemma~\ref{A-qcoh-X-cot-coresolution-dimension} below.

\begin{cor} \label{qcoh-X-cta-derived-equivalence}
 Let $X$ be a quasi-compact semi-separated scheme and $\cA$ be
a quasi-coherent quasi-algebra over~$X$.
 Then, for any conventional derived category symbol\/ $\st=\bb$, $+$,
$-$, or\/~$\varnothing$, the inclusion of exact/abelian categories
$\cA\Qcoh^{X\dcta}\rarrow\cA\Qcoh$ induces an equivalence of
the derived categories
$$
 \sD^\st(\cA\Qcoh^{X\dcta})\simeq\sD^\st(\cA\Qcoh).
$$
\end{cor}

\begin{proof}
 Follows from
Lemma~\ref{X-A-cta-cot-coresolving-coresolution-dimension}(a\+-b)
in view of the dual version of Proposition~\ref{finite-resolutions}.
\end{proof}

 For a similar equivalence involving also the derived category
$\sD^\st(\cA\Qcoh^{X\dcot})$ for any conventional derived category
symbol $\st=\bb$, $+$, $-$, or~$\varnothing$ (holding under more
restrictive assumptions on the scheme~$X$), see
Corollary~\ref{qcoh-X-cot-derived-equivalence} below.

 The following theorem is a generalization
of~\cite[Corollary~10.7]{PS6}.

\begin{thm} \label{qcoh-X-cot-A-cot-derived-equivalence}
 Let $X$ be a quasi-compact semi-separated scheme and $\cA$ be
a quasi-coherent quasi-algebra over~$X$.
 Then, for any derived category symbol\/ $\st=+$ or\/~$\varnothing$,
the inclusions of exact/abelian categories
$\cA\Qcoh^{\cA\dcot}\rarrow X\Qcoh^{X\dcot}\rarrow\cA\Qcoh$ induce
equivalences of the derived categories
$$
 \sD^\st(\cA\Qcoh^{\cA\dcot})\simeq
 \sD^\st(\cA\Qcoh^{X\dcot})\simeq\sD^\st(\cA\Qcoh).
$$
\end{thm}

\begin{proof}
 It suffices to show that the triangulated functors
$\sD^\st(\cA\Qcoh^{\cA\dcot})\rarrow\sD^\st(\cA\Qcoh)$ and
$\sD^\st(\cA\Qcoh^{X\dcot})\rarrow\sD^\st(\cA\Qcoh)$ are
triangulated equivalences.
 The case $\st=+$ is easy and follows from
Lemma~\ref{X-A-cta-cot-coresolving-coresolution-dimension}(a)
by virtue of the dual version of Proposition~\ref{infinite-resolutions}.
 The case of the doubly unbounded derived categories
($\st=\varnothing$) is interesting and based on the cotorsion
periodicity (Theorems~\ref{qcoh-cotorsion-periodicity}
and~\ref{qcoh-quasi-algebra-cotorsion-periodicity}).

 It is convenient to refer to~\cite[Proposition~10.6]{PS6}.
 In view of this proposition, in order to prove that the functor
$\sD(\cA\Qcoh^{X\dcot})\rarrow\sD(\cA\Qcoh)$ is a triangulated
equivalence, it suffices to check that a complex in $\cA\Qcoh^{X\dcot}$
is acyclic in $\cA\Qcoh^{X\dcot}$ whenever it is acyclic in $\cA\Qcoh$.
 This assertion follows from the fact that a complex in the exact
category of cotorsion quasi-coherent sheaves $X\Qcoh^\cot$ is
acyclic in $X\Qcoh^\cot$ whenever it is acyclic in $X\Qcoh$.
 This is Theorem~\ref{qcoh-cotorsion-periodicity}.

 By the same proposition from~\cite{PS6}, in order to prove that
$\sD(\cA\Qcoh^{\cA\dcot})\rarrow\sD(\cA\Qcoh)$ is a triangulated
equivalence, it suffices to check that a complex in $\cA\Qcoh^{\cA\dcot}$
is acyclic in $\cA\Qcoh^{\cA\dcot}$ whenever it is acyclic
in $\cA\Qcoh$.
 This is Theorem~\ref{qcoh-quasi-algebra-cotorsion-periodicity}.
\end{proof}

\subsection{Background derived equivalences for locally
contraherent $\cA$-modules} \label{background-derived-contraherent}
 In this section we prove $\cA$\+module generalizations of some
results about the conventional derived categories of locally
contraherent cosheaves from~\cite[Section~4.7]{Pcosh}.

 Let $X$ be a quasi-compact semi-separated scheme with an open
covering $\bW$ and $\cA$ be a quasi-coherent quasi-algebra over~$X$.
 We refer to Section~\ref{antilocal-contrah-A-lct-A-modules-subsecn}
for the notation $\cA\Ctrh^{\cA\dlct}_\al\subset\cA\Ctrh^{X\dlct}_\al
\subset\cA\Ctrh_\al$ for the exact categories of antilocal
contraherent $\cA$\+modules.
 The similar notation $X\Ctrh_\al$ for the exact category of
antilocal contraherent cosheaves on $X$ was also mentioned there.

\begin{lem} \label{antilocal-resolving-resolution-dimension}
 Let $X$ be a quasi-compact semi-separated scheme with an open
covering\/ $\bW$ and a finite affine open covering
$X=\bigcup_{\alpha=1}^N U_\alpha$ subordinate to\/~$\bW$.
 Let $\cA$ be a quasi-coherent quasi-algebra over~$X$.
 In this context: \par
\textup{(a)} The full subcategories $\cA\Ctrh_\al\subset\cA\Ctrh
\subset\cA\Lcth_\bW$ of contraherent $\cA$\+modules and antilocal
contraherent $\cA$\+modules are resolving in the exact category
of\/ $\bW$\+locally contraherent $\cA$\+modules $\cA\Lcth_\bW$.
 All the objects of $\cA\Lcth_\bW$ have finite resolution dimensions
not exceeding $N-1$ with respect to the resolving subcategories
$\cA\Ctrh_\al$ and $\cA\Ctrh$. \par
\textup{(b)} The full subcategories $\cA\Ctrh_\al^{X\dlct}\subset
\cA\Ctrh^{X\dlct}\subset\cA\Lcth_\bW^{X\dlct}$ of $X$\+locally cotorsion
contraherent $\cA$\+modules and antilocal $X$\+locally cotorsion
contraherent $\cA$\+modules are resolving in the exact category
of $X$\+locally cotorsion\/ $\bW$\+locally contraherent $\cA$\+modules
$\cA\Lcth_\bW^{X\dlct}$.
 All the objects of $\cA\Lcth_\bW^{X\dlct}$ have finite resolution
dimensions not exceeding $N-1$ with respect to the resolving
subcategories $\cA\Ctrh_\al^{X\dlct}$ and $\cA\Ctrh^{X\dlct}$. \par
\textup{(c)} The full subcategories $\cA\Ctrh_\al^{\cA\dlct}\subset
\cA\Ctrh^{\cA\dlct}\subset\cA\Lcth_\bW^{\cA\dlct}$ of $\cA$\+locally
cotorsion contraherent $\cA$\+modules and antilocal $\cA$\+locally
cotorsion contraherent $\cA$\+modules are resolving in the exact
category of $\cA$\+locally cotorsion\/ $\bW$\+locally contraherent
$\cA$\+modules $\cA\Lcth_\bW^{\cA\dlct}$.
 All the objects of $\cA\Lcth_\bW^{\cA\dlct}$ have finite resolution
dimensions not exceeding $N-1$ with respect to the resolving
subcategories $\cA\Ctrh_\al^{\cA\dlct}$ and $\cA\Ctrh^{\cA\dlct}$.
\end{lem}

\begin{proof}
 The full subcategory $\cA\Ctrh$ is closed under extensions and
kernels of admissible epimorphisms in the exact category $\cA\Lcth_\bW$
in view of the discussion in
Section~\ref{exact-categories-of-contrah-subsecn}.
 It follows that the full subcategory $\cA\Ctrh^{X\dlct}$ is closed
under extensions and kernels of admissible epimorphisms in
the exact category $\cA\Lcth_\bW^{X\dlct}$, and the full subcategory
$\cA\Ctrh^{\cA\dlct}$ is closed under extensions and kernels of
admissible epimorphisms in the exact category $\cA\Lcth_\bW^{\cA\dlct}$.

 The full subcategory $\cA\Ctrh_\al$ is closed under extensions and
kernels of admissible epimorphisms in the exact category $\cA\Lcth_\bW$,
since the full subcategory of antilocal contraherent cosheaves
$X\Ctrh_\al$ is closed under extensions and kernels of admissible
epimorphisms in the exact category $X\Lcth_\bW$.
 See~\cite[Corollary~4.3.3(b)]{Pcosh} and the discussion in
Section~\ref{A-antilocal=X-antilocal-subsecn}.
 The full subcategory $\cA\Ctrh_\al^{X\dlct}$ is closed under extensions
and kernels of admissible epimorphisms in the exact category
$\cA\Lcth_\bW^{X\dlct}$, and the full subcategory
$\cA\Ctrh_\al^{\cA\dlct}$ is closed under extensions and kernels of
admissible epimorphisms in the exact category $\cA\Lcth_\bW^{\cA\dlct}$,
for the same reasons.

 Every object of $\cA\Lcth_\bW$ is an admissible quotient of an object
from $\cA\Ctrh_\al$, every object of $\cA\Lcth_\bW^{X\dlct}$ is
an admissible quotient of an object from $\cA\Ctrh_\al^{X\dlct}$, and
every object of $\cA\Lcth_\bW^{\cA\dlct}$ is an admissible quotient of
an object from $\cA\Ctrh_\al^{\cA\dlct}$ by
Lemmas~\ref{antilocal-modules-generating-class}
and~\ref{antilocal-X-lct-modules-generating-class},
and Remark~\ref{antilocal-contrah-A-lct-A-modules-remark}.
 Finally, the resolution dimensions do not exceed $N-1$ because
the \v Cech resolution~\eqref{lcth-sheaf-of-rings-cech-resolution}
is a resolution by antilocal contraherent $\cA$\+modules.
 (Cf.~\cite[Lemmas~4.7.1(c) and~4.7.2(a\+-b)]{Pcosh}.)
\end{proof}

\begin{cor} \label{A-lcth-ctrh-al-derived-equivalences}
 Let $X$ be a quasi-compact semi-separated scheme with an open
covering\/ $\bW$, and let $\cA$ be a quasi-coherent quasi-algebra
over~$X$.
 Let\/ $\st=\bb$, $+$, $-$, or\/~$\varnothing$ be a conventional
derived category symbol.
 Then \par
\textup{(a)} the inclusions of exact categories\/ $\cA\Ctrh_\al
\rarrow\cA\Ctrh\rarrow\cA\Lcth_\bW$ induce equivalences of
the derived categories
$$
 \sD^\st(\cA\Ctrh_\al)\simeq\sD^\st(\cA\Ctrh)
 \simeq\sD^\st(\cA\Lcth_\bW);
$$ \par
\textup{(b)} the inclusions of exact categories\/ $\cA\Ctrh_\al^{X\dlct}
\rarrow\cA\Ctrh^{X\dlct}\rarrow\cA\Lcth_\bW^{X\dlct}$ induce
equivalences of the derived categories
$$
 \sD^\st(\cA\Ctrh_\al^{X\dlct})\simeq\sD^\st(\cA\Ctrh^{X\dlct})
 \simeq\sD^\st(\cA\Lcth_\bW^{X\dlct});
$$ \par
\textup{(c)} the inclusions of exact categories\/
$\cA\Ctrh_\al^{\cA\dlct}\rarrow\cA\Ctrh^{\cA\dlct}\rarrow
\cA\Lcth_\bW^{\cA\dlct}$ induce equivalences of the derived categories
$$
 \sD^\st(\cA\Ctrh_\al^{\cA\dlct})\simeq\sD^\st(\cA\Ctrh^{\cA\dlct})
 \simeq\sD^\st(\cA\Lcth_\bW^{\cA\dlct}).
$$
\end{cor}

\begin{proof}
 All the assertions follow from the respective assertions of
Lemma~\ref{antilocal-resolving-resolution-dimension} in view of
Proposition~\ref{finite-resolutions}.
\end{proof}

\begin{lem} \label{X-lct-A-lct-coresolving}
 Let $X$ be a quasi-compact semi-separated scheme with an open
covering\/ $\bW$, and let $\cA$ be a quasi-coherent quasi-algebra
over~$X$.
 Then the full subcategories of $X$\+locally cotorsion and
$\cA$\+locally cotorsion\/ $\bW$\+locally contraherent $\cA$\+modules
$\cA\Lcth_\bW^{X\dlct}$ and $\cA\Lcth_\bW^{\cA\dlct}$ are coresolving
in the exact category of\/ $\bW$\+locally contraherent $\cA$\+modules
$\cA\Lcth_\bW$.
\end{lem}

\begin{proof}
 The full subcategories in question are cogenerating in
$\cA\Lcth_\bW$ by Lemma~\ref{qcomp-qsep-A-loc-inj-preenvelope}
(recall that $\cA\Lcth_\bW^{\cA\dlin}\subset\cA\Lcth_\bW^{\cA\dlct}
\subset\cA\Lcth_\bW^{X\dlct}$ by
Lemma~\ref{restriction-coextension-injective-cotorsion}(a), as pointed
out in Section~\ref{A-loc-cotors-loc-inj-cosheaves-subsecn}).
 The two full subcategories are also closed under extensions and
cokernels of admissible monomorphisms in $\cA\Lcth_\bW$ by
Lemma~\ref{flat-cotorsion-pair-hereditary}(b), as mentioned in
Sections~\ref{exact-categories-of-contrah-subsecn}
and~\ref{A-loc-cotors-loc-inj-cosheaves-subsecn}.
\end{proof}

\begin{lem} \label{monic-quis-into-complex-of-injective-modules}
 Let $R$ be an associative ring and $C^\bu$ be a complex of left
$R$\+modules.
 Then there exists a complex of injective left $R$\+modules $J^\bu$
together with a termwise injective morphism of complexes of $R$\+modules
$C^\bu\rarrow J^\bu$ such that the quotient complex $J^\bu/C^\bu$
is acyclic.
\end{lem}

\begin{proof}
 In fact, the requirements in the lemma are weak and can be strengthened
in various ways: in particular, one can require the complex $J^\bu$ to
be a \emph{homotopy injective} complex of injective modules, or
alternatively, one can require the complex $J^\bu/C^\bu$ to be
\emph{coacyclic in the sense of Becker} (but not both simultaneously).
 One can start by noticing an easy construction of a termwise injective
morphism of any complex $C^\bu$ into a contractible complex of injective
modules~$K^\bu$.
 Then one can find a morphism from $C^\bu$ into a complex of injective
modules $I^\bu$ such that the cone of the morphism $C^\bu\rarrow I^\bu$
is acyclic, as per~\cite[proof of Proposition~10.6]{PS6}.
 Finally, put $J^\bu=I^\bu\oplus K^\bu$, and consider the diagonal
morphism $C^\bu\rarrow I^\bu\oplus K^\bu$.
 The same arguments work in any Grothendieck category in place of
$R\Modl$.
\end{proof}

\begin{lem} \label{contrah-quasi-alg-quis-into-A-lin}
 Let $X$ be a quasi-compact semi-separated scheme with an open
covering\/ $\bW$, and let $\cA$ be a quasi-coherent quasi-algebra
over~$X$.
 Then, for any complex of\/ $\bW$\+locally contraherent $\cA$\+modules\/
$\gM^\bu$ on $X$, there exists a complex of $\cA$\+locally injective\/
$\bW$\+locally contraherent $\cA$\+modules\/ $\gJ^\bu$ on $X$ together
with a morphism of complexes\/ $\gM^\bu\rarrow\gJ^\bu$ that is
a termwise admissible monomorphism in $\cA\Lcth_\bW$ such that
the cokernel\/ $\gJ^\bu/\gM^\bu$ is acyclic in $\cA\Lcth_\bW$.
\end{lem}

\begin{proof}
 This is a (weak) locally contraherent $\cA$\+module version
of~\cite[Lemmas~4.7.8 and~4.9.4]{Pcosh}.
 The idea is to plug in whole complexes of cosheaves
(rather than single cosheaves) into the construction of 
Lemmas~\ref{qcomp-qsep-A-loc-inj-preenvelope}
and~\ref{qcomp-qsep-X-lct-A-loc-inj-preenvelope}.
 One has to use locality of the notion of an $\cA$\+locally
injective locally contraherent $\cA$\+module, together with the facts
that the class of such modules is closed under extensions and
cokernels of monomorphisms in $\cA\Lcth_\bW$ and preserved by
the direct images with respect to $(\bW,\bT)$\+affine open
immersions of schemes.
 One also needs to use
Lemma~\ref{monic-quis-into-complex-of-injective-modules}.
 The construction produces a (termwise admissible) short exact sequence
$0\rarrow\gM^\bu\rarrow\gJ^\bu\rarrow\P^\bu\rarrow0$ of complexes
in the exact category $\cA\Lcth_\bW$, where
$\gJ^\bu\in\cA\Lcth_\bW^{\cA\dlin}$ and $\P^\bu$ is a finitely iterated
extension of the direct images of acyclic complexes in the exact
categories $\cA|_{U_\alpha}\Ctrh$ of contraherent
$\cA|_{U_\alpha}$\+modules on~$U_\alpha$.
 Here $X=\bigcup_\alpha U_\alpha$ is some finite affine open covering
of $X$ subordinate to~$\bW$.
\end{proof}

 The following theorem is a generalization
of~\cite[Theorem~4.7.9]{Pcosh};
see also~\cite[Corollary~6.4.4(b)]{Pcosh}.

\begin{thm} \label{A-lcth-X-lct-A-lct-derived-equivalence}
 Let $X$ be a quasi-compact semi-separated scheme with an open
covering\/ $\bW$, and let $\cA$ be a quasi-coherent quasi-algebra
over~$X$.
 Then, for any derived category symbol\/ $\st=+$ or\/~$\varnothing$,
the inclusions of exact categories $\cA\Lcth_\bW^{\cA\dlct}\rarrow
\cA\Lcth_\bW^{X\dlct}\rarrow\cA\Lcth_\bW$ induce equivalences
of the derived categories
$$
 \sD^\st(\cA\Lcth_\bW^{\cA\dlct})\simeq\sD^\st(\cA\Lcth_\bW^{X\dlct})
 \simeq\sD^\st(\cA\Lcth_\bW).
$$
\end{thm}

\begin{proof}
 The easy case $\st=+$ follows from Lemma~\ref{X-lct-A-lct-coresolving}
by virtue of the dual version of Proposition~\ref{infinite-resolutions}.
 The interesting case $\st=\varnothing$ follows from
Corollary~\ref{loc-contrah-quasi-algebra-cotorsion-periodicity}
and Lemma~\ref{contrah-quasi-alg-quis-into-A-lin} by virtue of
a well-known lemma about localizations of categories by multiplicative
classes of morphisms~\cite[Proposition~10.2.7(i)]{KS},
\cite[Lemma~1.6(b)]{Pkoszul}, or~\cite[Lemma~A.3.3(b)]{Pcosh}.
 The point is that Lemma~\ref{contrah-quasi-alg-quis-into-A-lin}
provides a quasi-isomorphism in $\cA\Lcth_\bW$ from any given complex
in $\cA\Lcth_\bW$ to a complex in $\cA\Lcth_\bW^{\cA\dlin}$, and one
has $\cA\Lcth_\bW^{\cA\dlin}\subset\cA\Lcth_\bW^{\cA\dlct}\subset
\cA\Lcth_\bW^{X\dlct}$.
\end{proof}

 For a similar equivalence involving the derived categories
$\sD^\st(\cA\Lcth_\bW^{X\dlct})$ and $\sD^\st(\cA\Lcth_\bW)$
for any conventional derived category symbol $\st=\bb$, $+$, $-$,
or~$\varnothing$ (valid under more restrictive assumptions
on the scheme~$X$), see
Corollary~\ref{A-lcth-al-X-lct-X-lcta-derived-equivalences}(a) below.

\subsection{Contraherent $\fHom$ over a quasi-coherent quasi-algebra}
\label{fHom-over-qcoh-quasi-algebra-subsecn}
 This section is a partial extension of~\cite[Section~2.5]{Pcosh}
from $\cO_X$\+modules to $\cA$\+modules.
 Let $X$ be a scheme, and let $\cA$ and $\cB$ be two quasi-coherent
quasi-algebras over~$X$.

 A \emph{quasi-coherent $\cA$\+$\cB$\+bimodule} $\cE$ on $X$ is
a quasi-coherent quasi-module on $X$ endowed with a structure of
a sheaf of $\cA$\+$\cB$\+bimodules satisfying the following condition.
 The $\cO_X$\+module structure underlying the left $\cA$\+module
structure on $\cE$ must coincide with the left $\cO_X$\+module structure
on the quasi-coherent quasi-module $\cE$, and the $\cO_X$\+module
structure underlying the right $\cB$\+module structure on $\cE$ must
coincide with the right $\cO_X$\+module structure on the quasi-coherent
quasi-module $\cE$ on~$X$.
 Equivalently, a structure of quasi-coherent $\cA$\+$\cB$\+bimodule
on a quasi-coherent quasi-module $\cE$ on $X$ is given by
an \emph{action map} $\cA\ot_{\cO_X}\cE\ot_{\cO_X}\cB\rarrow\cE$,
which must be a morphism of quasi-coherent quasi-modules on~$X$
(cf.\ Lemma~\ref{quasi-coherent-quasi-modules-tensor-product}(a))
satisfying the usual associativity and unitality axioms.
 The category of quasi-coherent $\cA$\+$\cB$\+bimodules on $X$ is
a Grothendieck abelian category.

 In the case of two quas-coherent quasi-algebras $\cA$ and $\cB$ over
an affine scheme $U$, the global sections functor provides
an equivalence between the category of quasi-coherent
$\cA$\+$\cB$\+bimodules $\cE$ on $X$ and the category of
$\cA(U)$\+$\cB(U)$\+bimodules $E=\cE(U)$ whose underlying
$\cO(U)$\+$\cO(U)$\+bimodule is a quasi-module over~$\cO(U)$.

 Let $\bB$ denote the topology base of $X$ consisting of all
the affine open subschemes $U\subset X$.
 Let $\F$ be a quasi-coherent left $\cA$\+module on $X$ and $R$ be
an associative ring.
 Assume that $R$ acts on $\F$ on the right by left $\cA$\+module
endomorphisms.
 Given a left $R$\+module $M$, define a presheaf of $\cA$\+modules
$\F\ot_RM$ on $\bB$ by the rule
$$
 (\F\ot_RM)(U)=\F(U)\ot_RM
$$
for all $U\in\bB$.
 One can easily check that (the underlying presheaf of $\cO_X$\+modules
of) $\F\ot_RM$ satisfies the quasi-coherence axiom from
Section~\ref{locally-contraherent-cosheaves-subsecn}.
 By Lemma~\ref{quasi-coherence-implies-sheaf}
and Theorem~\ref{extension-of-co-sheaves-from-topology-base}(a),
the presheaf $\F\ot_RM$ extends uniquely to a quasi-coherent sheaf
of $\cA$\+modules on $X$, which we will denote also by $\F\ot_RM$.

 Let $\C$ be a quasi-coherent left $\cA$\+module.
 Then the Hom group $\Hom_\cA(\F,\C)$ has a natural left $R$\+module
structure, and the adjunction isomorphism
$$
 \Hom_R(M,\Hom_\cA(\F,\C))\simeq\Hom_\cA(\F\ot_RM,\>\C)
$$
holds for any left $R$\+module~$M$.

\begin{lem} \label{Ext-over-quasi-algebra-homological-formula}
 Let $\F$ be a quasi-coherent left $\cA$\+module with a right action of
a ring $R$ by left $\cA$\+module endomorphisms, $\C$ be a quasi-coherent
left $\cA$\+module, and $M$ be a left $R$\+module.
 Let $n\ge0$ be an integer such that\/ $\Ext^i_\cA(\F,\C)=0$ for all\/
$0<i\le n$.
 Assume that\/ $\Tor_i^R(\F(U),M)=0$ for all affine open subschemes
$U\subset X$ and all\/ $0<i\le n$.
 Then there is a natural isomorphism of abelian groups
$$
 \Ext^n_R(M,\Hom_X(\F,\C))\simeq\Ext_\cA^n(\F\ot_RM,\>\C).
$$
\end{lem}

\begin{proof}
 This is a quasi-coherent sheaf version of
Lemma~\ref{Ext-homological-formula} and a slight generalization
of~\cite[Lemma~2.5.1]{Pcosh}, provable in the same way.
\end{proof}

 Assume that the scheme $X$ is quasi-separated.
 Let $\cE$ be a quasi-coherent $\cA$\+$\cB$\+bimodule and $\C$ be
a quasi-coherent left $\cA$\+module on~$X$.
 The copresheaf of $\cB$\+modules $\fHom_\cA(\cE,\C)$ on $\bB$ is
defined by the rule
$$
 \fHom_\cA(\cE,\C)[U]=\Hom_\cA(j_*j^*\cE,\C)
$$
for all affine open subschemes $U\subset X$ with the open immersion
morphism $j\:U\rarrow X$.
 Notice that $j^*\cE$ is a quasi-coherent left $\cA|_U$\+module
with a right action of the ring $R=\cB(U)$, hence $j_*j^*\cE$ is
a quasi-coherent left $\cA$\+module with a right action of
the same ring.
 This right action of the ring $\cB(U)$ on the quasi-coherent
left $\cA$\+module $j_*j^*\C$ induces a left action of $\cB(U)$ on
the group $\Hom_\cA(j_*j^*\cE,\C)$, making $\fHom_\cA(\cE,\C)[U]$
a left $\cB(U)$\+module.

 Furthermore, let $V\subset U\subset X$ be a pair of affine open
subschemes with open immersion morphisms $j\:U\rarrow X$, \
$h\:V\rarrow U$, and $k=jh\:V\rarrow X$.
 Then there is a natural adjunction morphism
$$
 j_*j^*\cE\lrarrow j_*h_*h^*j^*\cE\simeq k_*k^*\cE.
$$
 This morphism of quasi-coherent $\cA$\+modules $j_*j^*\cE\rarrow
k_*k^*\cE$ induces a map of abelian groups
$\Hom_\cA(k_*k^*\cE,\C)\rarrow\Hom_\cA(j_*j^*\cE,\C)$, providing
the corestriction map $\fHom_\cA(\cE,\C)[V]\rarrow\fHom_\cA(\cE,\C)[U]$
in the copresheaf $\fHom_\cA(\cE,\C)$ on~$\bB$.

 The contraherence axiom~(i) from
Section~\ref{locally-contraherent-cosheaves-subsecn} always holds
for (the underlying copresheaf of $\cO_X$\+modules of) the copresheaf
of $\cB$\+modules $\fHom_\cA(\cE,\C)$ on~$\bB$.
 Indeed, one computes
\begin{multline*}
 \Hom_{\cO_X(U)}(\cO_X(V),\fHom_\cA(\cE,\C)[U])=
 \Hom_{\cO_X(U)}(\cO_X(V),\Hom_\cA(j_*j^*\cE,\C)) \\
 \simeq\Hom_\cA((j_*j^*\cE)\ot_{\cO_X(U)}\cO_X(V),\>\C)
 \simeq\Hom_\cA(k_*k^*\cE,\C)=\fHom_\cA(\cE,\C)[V].
\end{multline*}
 Here we are using the natural isomorphism of quasi-coherent
$\cA$\+modules
\begin{equation} \label{direct-inverse-image-tensor-isomorphism}
 (j_*j^*\cE)\ot_{\cO_X(U)}\cO_X(V)\simeq
 j_*(j^*\cE\ot_{\cO_X(U)}\cO_X(V))\simeq
 j_*h_*h^*j^*\cE\simeq k_*k^*\cE
\end{equation} 
on~$X$.
 The isomorphism $(j_*\G)\ot_RH\simeq j_*(\G\ot_RH)$ holds
for any quasi-coherent sheaf $\G$ on $U$ with a right action
of a ring $R$ and any flat left $R$\+module $H$, and only requires
the morphism~$j$ to be quasi-compact and quasi-separated
(see~\cite[Section~2.5]{Pcosh}); when the morphism~$j$ is affine,
this isomorphism holds for an arbitrary left $R$\+module~$H$.
 On the other hand, the isomorphism $j^*\cE\ot_{\cO_X(U)}\cO_X(V)
\simeq h_*h^*j^*\cE$ holds due to quasi-coherence of the sheaf
$\cE$ with respect to its right $\cO_X$\+module structure.

 The lemma below provides sufficient conditions for
the contraadjustedness axiom~(ii) from
Section~\ref{locally-contraherent-cosheaves-subsecn} and its stronger
versions to be satisfied for the copresheaf $\fHom_\cA(\cE,\C)$.
 See~\cite[Section~0.10 of the Introduction]{Pcosh}
and~\cite[Sections~5.6 and~8.7]{Pphil} for a relevant discussion of
partially defined functors between exact categories.

 We will say that a quasi-coherent $\cA$\+$\cB$\+bimodule is
\emph{$\cA$\+very flaprojective} if it is very flaprojective (in
the sense of Section~\ref{antilocality-of-X-contraadjusted-subsecn})
as a quasi-coherent left $\cA$\+module.
 The terminology \emph{$\cA$\+flaprojective quasi-coherent
$\cA$\+$\cB$\+bimodule} and \emph{$\cA$\+flat quasi-coherent
$\cA$\+$\cB$\+bimodule} has a similar meaning (with the references to
the definitions in Sections~\ref{antilocality-of-X-cotorsion-subsecn}
and~\ref{antilocality-of-A-cotorsion-subsecn}).
 In particular, any quasi-coherent quasi-algebra $\cA$ over $X$ is
naturally an $\cA$\+very flaprojective (hence also $\cA$\+flaprojective
and $\cA$\+flat) quasi-coherent $\cA$\+$\cA$\+bimodule
(cf.\ Remark~\ref{very-flaprojective-terminology-explained}).

 Now let $\cE$ be a quasi-coherent $\cA$\+$\cO_X$\+bimodule.
 We will say that $\cE$ is \emph{robustly flaprojective} if, for every
affine open subscheme $U\subset X$, the $A$\+$R$\+bimodule $\cE(U)$
is $(A/R,R)$\+robustly flaprojective in the sense of the definition
in Section~\ref{prelim-robustly-flaprojective-subsecn}.
 Here we use the notation $A=\cA(U)$ and $R=\cO_X(U)$.
 In view of Lemma~\ref{quasi-algebras-robust-flaprojectivity-locality}
(for $K=R$), it suffices to check this condition for affine open
subschemes $U\subset X$ belonging to any chosen affine open covering
of the scheme~$X$.
 For example, for any quasi-coherent quasi-module $\F$ on $X$ that
is flat as a quasi-coherent sheaf with respect to its left
$\cO_X$\+module structure, the quasi-coherent $\cA$\+$\cO_X$\+bimodule
$\cA\ot_{\cO_X}\F$ is robustly flaprojective.
 We will say that a quasi-coherent $\cA$\+$\cB$\+bimodule $\cE$ is
\emph{$\cA$\+robustly flaprojective} if $\cE$ is robustly flaprojective
as a quasi-coherent $\cA$\+$\cO_X$\+bimodule.

 When no confusion is likely to arise between the left and right module
structures, we will also speak of \emph{$\cB$\+flat} quasi-coherent
$\cA$\+$\cB$\+bimodules, meaning that the underlying quasi-coherent
right $\cB$\+module is flat.
 Finally, we will say that a quasi-coherent left $\cA$\+module $\J$
is \emph{injective} if it is injective as an object of the Grothendieck
abelian category $\cA\Qcoh$.
 The definitions of \emph{$X$\+contraadjusted}, \emph{$X$\+cotorsion},
and \emph{$\cA$\+cotorsion} quasi-coherent $\cA$\+modules can be also
found in Sections~\ref{antilocality-of-X-contraadjusted-subsecn}\+-%
\ref{antilocality-of-A-cotorsion-subsecn}.

 A simpler affine version of the following lemma is stated as
Lemma~\ref{module-Hom-contraadjusted-and-more-lemma} below.

\begin{lem} \label{fHom-contraadjusted-and-more-lemma}
 Let $\cA$ and $\cB$ be quasi-coherent quasi-algebras over
a scheme $X$, let $\cE$ be a quasi-coherent $\cA$\+$\cB$\+bimodule,
and let $\C$ and $\J$ be a quasi-coherent left $\cA$\+modules.
 Let $U\subset X$ be an affine open subscheme.
 In this context: \par
\textup{(a)} if the scheme $X$ is quasi-compact and semi-separated,
the quasi-coherent $\cA$\+$\cB$\+bimodule $\cE$ is $\cA$\+very
flaprojective, and the quasi-coherent left $\cA$\+module $\C$ is
$X$\+contraadjusted, then the $\cO_X(U)$\+module\/
$\fHom_\cA(\cE,\C)[U]$ is contraadjusted; \par
\textup{(b)} if the scheme $X$ is quasi-compact and semi-separated,
the quasi-coherent $\cA$\+$\cB$\+bimodule $\cE$ is $\cA$\+robustly
flaprojective, and the quasi-coherent left $\cA$\+module $\C$ is
$X$\+cotorsion, then the $\cO_X(U)$\+module\/ $\fHom_\cA(\cE,\C)[U]$
is cotorsion; \par
\textup{(c)} if the scheme $X$ is semi-separated, the quasi-coherent
$\cA$\+$\cB$\+bimodule $\cE$ is $\cA$\+flat, and the quasi-coherent
left $\cA$\+module $\C$ is $\cA$\+cotorsion, then the $\cB(U)$\+module\/
$\fHom_\cA(\cE,\C)[U]$ is cotorsion; \par
\textup{(d)} if the scheme $X$ is semi-separated, the quasi-coherent
$\cA$\+$\cB$\+bimodule $\cE$ is $\cB$\+flat and the quasi-coherent
left $\cA$\+module $\J$ is injective, then the $\cB(U)$\+module\/
$\fHom_\cA(\cE,\J)[U]$ is injective; \par
\textup{(e)} if the scheme $X$ is quasi-separated and the quasi-coherent
left $\cA$\+module $\J$ is injective, then the $\cB(U)$\+module\/
$\fHom_\cA(\cE,\J)[U]$ is cotorsion.
\end{lem}

\begin{proof}
 Notice first of all that flatness of quasi-coherent $\cA$\+modules
and $\cB$\+modules is preserved by the restrictions to open subschemes
and the direct images from affine immersions of open subschemes.
 The same applies to flaprojectivity and very flaprojectivity
(see Sections~\ref{antilocality-of-X-contraadjusted-subsecn}\+-%
\ref{antilocality-of-A-cotorsion-subsecn}).
 To make the open immersion morphism $j\:U\rarrow X$ affine for
an affine open subscheme $U\subset X$, the semi-separatedness assumption
is imposed in parts~(a\+-d).

 In part~(a), we let $G$ be a very flat $\cO_X(U)$\+module and compute
\begin{multline*}
 \Ext^1_{\cO_X(U)}(G,\Hom_\cA(j_*j^*\cE,\C)) \\ \simeq
 \Ext^1_\cA((j_*j^*\cE)\ot_{\cO_X(U)}G,\>\C) \simeq
 \Ext^1_\cA(j_*(j^*\cE\ot_{\cO_X(U)}G),\>\C)=0.
\end{multline*}
 Here the first isomorphism holds by
Lemma~\ref{Ext-over-quasi-algebra-homological-formula} (for $n=1$
and $R=\cO_X(U)$), the second isomorphism is true because $j$~is
an affine morphism (or because $G$ is a flat $\cO_X(U)$\+module),
and the final vanishing holds since $\cE(U)\ot_{\cO_X(U)}G$ is
an $\cA(U)/\cO_X(U)$\+very flaprojective $\cA(U)$\+module
(by Lemma~\ref{flepi-very-flat-modules-unproblematic-lemma}(b)).
 It is helpful to notice that the class of very flaprojective
quasi-coherent $\cA$\+modules is preserved by the direct images from
affine open immersion morphisms (by
Lemma~\ref{very-flaprojective-restriction-extension-of-scalars}(a), as
mentioned in Section~\ref{antilocality-of-X-contraadjusted-subsecn}),
and one has $\Ext^1_\cA(\F,\C)=0$ for any very flaprojective
quasi-coherent $\cA$\+module $\F$ and any $X$\+contraadjusted
quasi-coherent $\cA$\+module $\C$ on $X$ by
Corollary~\ref{vflp-Ext-orthogonal-to-X-cta}.

 Alternatively, instead of referring to
Lemma~\ref{flepi-very-flat-modules-unproblematic-lemma}(b), one can say
that it suffices to consider the $\cO_X(U)$\+modules $G=\cO_X(V)$
for affine open subschemes $V\subset U$.
 Then $\cE(U)\ot_{\cO_X(U)}\cO_X(V)\simeq\cE(V)\simeq \cO_X(V)
\ot_{\cO_X(U)}\cE(U)$ is an $\cA(U)/\cO_X(U)$\+very flaprojective
$\cA(U)$\+module by
Lemma~\ref{flaprojective-restriction-extension-of-scalars}(a).

 The same computation works in part~(b) for a flat
$\cO_X(U)$\+module~$G$.
 In this context, instead of referring to
Lemma~\ref{flepi-very-flat-modules-unproblematic-lemma}, one says
that the $\cA(U)$\+module $\cE(U)\ot_{\cO_X(U)}G$ is
$\cA(U)/\cO_X(U)$\+flaprojective since the
$\cA(U)$\+$\cO_X(U)$\+bimodule $\cE(U)$ is
$(\cA(U)/\cO_X(U),\cO_X(U))$\+robustly flaprojective by assumption.
 The class of flaprojective quasi-coherent $\cA$\+modules is preserved
by the direct images from affine open immersion morphisms
(by Lemma~\ref{flaprojective-restriction-extension-of-scalars}(a), as
mentioned in Section~\ref{antilocality-of-X-cotorsion-subsecn}), and
one has $\Ext^1_\cA(\F,\C)=0$ for any flaprojective quasi-coherent
$\cA$\+module $\F$ and any $X$\+cotorsion quasi-coherent $\cA$\+module
$\C$ on $X$ by Corollary~\ref{flp-Ext-orthogonal-to-X-cot}.

 In part~(c), we let $G$ be a flat left $\cB(U)$\+module and compute
\begin{multline*}
 \Ext^1_{\cB(U)}(G,\Hom_\cA(j_*j^*\cE,\C)) \\ \simeq
 \Ext^1_\cA((j_*j^*\cE)\ot_{\cB(U)}G,\>\C) \simeq
 \Ext^1_\cA(j_*(j^*\cE\ot_{\cB(U)}G),\>\C)=0.
\end{multline*}
 Here the first isomorphism holds by
Lemma~\ref{Ext-over-quasi-algebra-homological-formula} (for $n=1$
and $R=\cB(U)$), the second isomorphism is true because $j$~is
an affine morphism (or because $G$ is a flat $\cB(U)$\+module),
and the final vanishing holds since $\cE(U)\ot_{\cB(U)}G$ is a flat
$\cA(U)$\+module.
 Notice that the class of flat quasi-coherent $\cA(U)$\+modules is
preserved by the direct images from affine open immersion morphisms
(as mentioned in Section~\ref{antilocality-of-A-cotorsion-subsecn}).

 In part~(d), we let $M$ be an arbitrary left $\cB(U)$\+module
and compute
\begin{multline*}
 \Ext^1_{\cB(U)}(M,\Hom_\cA(j_*j^*\cE,\J)) \\ \simeq
 \Ext^1_\cA((j_*j^*\cE)\ot_{\cB(U)}M,\>\J) \simeq
 \Ext^1_\cA(j_*(j^*\cE\ot_{\cB(U)}M),\>\J)=0.
\end{multline*}
 Here the first isomorphism holds by
Lemma~\ref{Ext-over-quasi-algebra-homological-formula} (for $n=1$
and $R=\cB(U)$) and the second isomorphism is true because $j$~is
an affine morphism.
 The condition of $\cB$\+flatness of the quasi-coherent
$\cA$\+$\cB$\+bimodule $\cE$ is needed in order to satisfy
the Tor vanishing assumption in
Lemma~\ref{Ext-over-quasi-algebra-homological-formula}.

 In part~(e), we let $G$ be a flat left $\cB(U)$\+module and compute
\begin{multline*}
 \Ext^1_{\cB(U)}(G,\Hom_\cA(j_*j^*\cE,\J)) \\ \simeq
 \Ext^1_\cA((j_*j^*\cE)\ot_{\cB(U)}G,\>\J) \simeq
 \Ext^1_\cA(j_*(j^*\cE\ot_{\cB(U)}G),\>\J)=0.
\end{multline*}
 Here the first isomorphism holds by
Lemma~\ref{Ext-over-quasi-algebra-homological-formula} (for $n=1$
and $R=\cB(U)$) and the second isomorphism is true because $G$ is
a flat $\cB(U)$\+module.
 The condition of flatness of the $\cB(U)$\+module $G$ is also needed
in order to satisfy the Tor vanishing assumption in
Lemma~\ref{Ext-over-quasi-algebra-homological-formula}.
\end{proof}

 Under any one of the assumptions~(a\+-e) of
Lemma~\ref{fHom-contraadjusted-and-more-lemma}, we have proved, in
particular, that the copresheaf $\fHom_\cA(\cE,\C)$ or
$\fHom_\cA(\cE,\J)$ on $\bB$, viewed as a copresheaf of
$\cO_X$\+modules, satisfies the contraadjustedness
axiom~(ii) from Section~\ref{locally-contraherent-cosheaves-subsecn}.
 Recall that any cotorsion $\cB(U)$\+module is also cotorsion
(hence contraadjusted) as a $\cO_X(U)$\+module by
Lemma~\ref{restriction-coextension-injective-cotorsion}(a).
 The contraherence axiom~(i) was established above.
 By Lemma~\ref{contraherence+contraadjustedness-imply-cosheaf}
and Theorem~\ref{extension-of-co-sheaves-from-topology-base}(b),
the copresheaf $\fHom_\cA(\cE,\C)$ or $\fHom_\cA(\cE,\J)$ on $\bB$
extends uniquely to a contraherent cosheaf of $\cB$\+modules on $X$,
which we will denote also by $\fHom_\cA(\cE,\C)$ or
$\fHom_\cA(\cE,\J)$.

 To summarize:
\begin{itemize}
\item $\fHom_\cA(\cE,\C)$ is a contraherent $\cB$\+module whenever $\cE$
is $\cA$\+very flaprojective, $\C$ is $X$\+contraadjusted, and $X$ is
quasi-compact and semi-separated;
\item $\fHom_\cA(\cE,\C)$ is an $X$\+locally cotorsion contraherent
$\cB$\+module whenever $\cE$ is $\cA$\+robustly flaprojective, $\C$
is $X$\+cotorsion, and $X$ is quasi-compact and semi-separated;
\item $\fHom_\cA(\cE,\C)$ is a $\cB$\+locally cotorsion contraherent
$\cB$\+module whenever $\cE$ is $\cA$\+flat, $\C$ is $\cA$\+cotorsion,
and $X$ is semi-separated;
\item $\fHom_\cA(\cE,\J)$ is a $\cB$\+locally injective contraherent
$\cB$\+module whenever $\cE$ is $\cB$\+flat, $\J$ is injective,
and $X$ is semi-separated;
\item $\fHom_\cA(\cE,\J)$ is a $\cB$\+locally cotorsion contraherent
$\cB$\+module whenever $\J$ is injective and $X$ is quasi-separated.
\end{itemize}
 Here the definitions of a contraherent $\cB$\+module and
an $X$\+locally cotorsion contraherent $\cB$\+module were given in
Section~\ref{cosheaves-of-A-modules-subsecn}, while the definitions
of a $\cB$\+locally cotorsion contraherent $\cB$\+module and
a $\cB$\+locally injective contraherent $\cB$\+module can be found
in Section~\ref{A-loc-cotors-loc-inj-cosheaves-subsecn}.

 Denote by $\cA\biQcoh\cB$ the Grothendieck abelian category of
quasi-coherent $\cA$\+$\cB$\+bi\-mod\-ules on~$X$.
 Denote further by $\cA\biQcoh\cB^{\cA\dvflp}$ the full subcategory
of $\cA$\+very flaprojective quasi-coherent $\cA$\+$\cB$\+bimodules,
by $\cA\biQcoh\cB^{\cA\drflp}$ the full subcategory of $\cA$\+robustly
flaprojective quasi-coherent $\cA$\+$\cB$\+bimodules,
by $\cA\biQcoh\cB^{\cA\dfl}$ the full subcategory of $\cA$\+flat
quasi-coherent $\cA$\+$\cB$\+bimodules, and by
$\cA\biQcoh\cB^{\cB\dfl}$ the full subcategory of $\cB$\+flat
quasi-coherent $\cA$\+$\cB$\+bimodules in $\cA\biQcoh\cB$.
 All the four full subcategories $\cA\biQcoh\cB^{\cA\dvflp}$,
\ $\cA\biQcoh\cB^{\cA\drflp}$, \ $\cA\biQcoh\cB^{\cA\dfl}$, and
$\cA\biQcoh\cB^{\cB\dfl}\subset\cA\biQcoh\cB$ are closed under
extensions and kernels of epimorphisms in $\cA\biQcoh\cB$, so they
inherit exact category structures from the abelian exact structure
of $\cA\biQcoh\cB$.

 For any abelian category $\sA$, we denote by $\sA^\inj\subset\sA$
the full subcategory of injective objects in~$\sA$.
 The additive category $\sA^\inj$ is endowed with the split exact
category structure.

 It is clear from the constructions and arguments above that, in
the context of each of the parts~(a\+-d) of
Lemma~\ref{fHom-contraadjusted-and-more-lemma}, the related functor
$\fHom_\cA$ is exact as a functor of two arguments between
the respective exact categories of quasi-coherent (bi)modules and
contraherent modules over quasi-algebras over~$X$.
 So, under the respective assumptions on the scheme $X$, we have
exact functors of two arguments, all of them denoted by
$\fHom_\cA({-},{-})$:
\begin{align}
 (\cA\biQcoh\cB^{\cA\dvflp})^\sop\times\cA\Qcoh^{X\dcta}
 &\lrarrow\cB\Ctrh;
 \label{fHom-A-very-flaproj-X-cta-exact} \\
 (\cA\biQcoh\cB^{\cA\drflp})^\sop\times\cA\Qcoh^{X\dcot}
 &\lrarrow\cB\Ctrh^{X\dlct};
 \label{fHom-A-rob-flaproj-X-cot-X-lct-exact} \\
 (\cA\biQcoh\cB^{\cA\dfl})^\sop\times\cA\Qcoh^{\cA\dcot}
 &\lrarrow\cB\Ctrh^{\cB\dlct};
 \label{fHom-A-flat-A-cot-B-lct-exact} \\
 (\cA\biQcoh\cB^{\cB\dfl})^\sop\times\cA\Qcoh^\inj
 &\lrarrow\cB\Ctrh^{\cB\dlin}. \label{fHom-B-flat-inj-B-dlin-exact}
\end{align}

\subsection{Cosections of the $\fHom$ cosheaf}
\label{cosections-of-fHom-subsecn}
 The aim of this section is prove the following version
of~\cite[Lemma~2.5.2]{Pcosh}.
 For affine open subschemes $Y\subset X$, the assertions of the lemma
hold by the definition of $\fHom$; the point of the lemma is
that there is often a similar isomorphism for \emph{nonaffine} open
subschemes $Y\subset X$.

\begin{lem} \label{cosections-of-fHom-lemma}
 Let $X$ be a scheme and $Y\subset X$ be an open subscheme in $X$ with
the open immersion morphism $j\:Y\rarrow X$.
 Let $\cA$ and $\cB$ be quasi-coherent quasi-algebras over $X$,
let $\cE$ be a quasi-coherent $\cA$\+$\cB$\+bimodule, and let $\C$
and $\J$ be quasi-coherent left $\cA$\+modules on~$X$.
 In this context: \par
\textup{(a)} if the scheme $X$ is quasi-compact and semi-separated,
the open immersion morphism~$j$ is affine, the quasi-coherent
$\cA$\+$\cB$\+bimodule $\cE$ is $\cA$\+very flaprojective, and
the quasi-coherent left $\cA$\+module $\C$ is $X$\+contraadjusted,
then there is a natural isomorphism of left $\cB(Y)$\+modules\/
$\fHom_\cA(\cE,\C)[Y]\simeq\Hom_\cA(j_*j^*\cE,\C)$; \par
\textup{(b)} if the scheme $X$ is quasi-compact and semi-separated,
the open immersion morphism~$j$ is affine, the quasi-coherent
$\cA$\+$\cB$\+bimodule $\cE$ is $\cA$\+robustly flaprojective, and
the quasi-coherent left $\cA$\+module is $X$\+cotorsion, then
there is a natural isomorphism of left $\cB(Y)$\+modules\/
$\fHom_\cA(\cE,\C)[Y]\simeq\Hom_\cA(j_*j^*\cE,\C)$; \par
\textup{(c)} if the scheme $X$ is semi-separated, the scheme $Y$
is quasi-compact, the open immersion morphism~$j$ is affine,
the quasi-coherent $\cA$\+$\cB$\+bimodule $\cE$ is $\cA$\+flat, and
the quasi-coherent left $\cA$\+module is $\cA$\+cotorsion, then
there is a natural isomorphism of left $\cB(Y)$\+modules\/
$\fHom_\cA(\cE,\C)[Y]\simeq\Hom_\cA(j_*j^*\cE,\C)$; \par
\textup{(d)} if the scheme $X$ is quasi-separated, the scheme $Y$
is quasi-compact, and the quasi-coherent left $\cA$\+module $\J$ is
injective, then there is a natural isomorphism of left
$\cB(Y)$\+modules\/ $\fHom_\cA(\cE,\J)[Y]\simeq\Hom_\cA(j_*j^*\cE,\J)$.
\end{lem}

\begin{proof}
 First of all, the $\Hom$ group $\Hom_\cA(j_*j^*\cE,\C)$ or
$\Hom_\cA(j_*j^*\cE,\J)$ is a left $\cB(Y)$\+module because
the ring $\cB(Y)$ acts on the right on the quasi-coherent left
$\cA|_Y$\+module $j^*\cE$ on $Y$, hence also on the quasi-coherent
left $\cA$\+module $j_*j^*\cE$ on~$X$.
 The group of cosections $\fHom_\cA(\cE,\C)[Y]$ or
$\fHom_\cA(\cE,\J)[Y]$ is a left $\cB(Y)$\+module since
$\fHom_\cA(\cE,\C)$ or $\fHom_\cA(\cE,\J)$ is a contraherent
$\cB$\+module, hence, in particular, a cosheaf of $\cB$\+modules
on~$X$ (according to the previous
Section~\ref{fHom-over-qcoh-quasi-algebra-subsecn}).

 The argument is similar to the proof of~\cite[Lemma~2.5.2]{Pcosh}.
 In all cases~(a\+-d), the scheme $X$ is, at least, quasi-separated,
and the open subscheme $Y$ is, at least, quasi-compact.
 Let $Y=\bigcup_{\alpha=1}^N V_\alpha$ be a finite affine open
covering of~$Y$.
 Denote by $k_{\alpha_1,\dotsc,\alpha_i}$ the open immersions
$V_{\alpha_1}\cap\dotsb\cap V_{\alpha_i}\rarrow Y$.
 For any quasi-coherent $\cA$\+module $\G$ on $Y$, the \v Cech
coresolution~\eqref{qcoh-sheaf-of-rings-cech-coresolution}
\begin{multline} \label{qcoh-sheaf-G-on-scheme-Y-chech-coresolution}
 0\lrarrow\G\lrarrow\bigoplus\nolimits_{\alpha=1}^N
 k_\alpha{}_*k_\alpha^*\G\lrarrow
 \bigoplus\nolimits_{1\le\alpha<\beta\le N}
 k_{\alpha,\beta}{}_*k_{\alpha,\beta}^*\G \\
 \lrarrow\dotsb\lrarrow k_{1,\dotsc,N}{}_*k_{1,\dotsc,N}^*\G
 \lrarrow0
\end{multline}
is a finite exact sequence of quasi-coherent $\cA|_Y$\+modules on~$Y$.
 Put $\G=j^*\cE$.

 When the open immersion morphism $j\:Y\rarrow X$ is affine
(cases~(a\+-c)), the functor $j_*\:\cA|_Y\Qcoh\rarrow\cA\Qcoh$ is exact.
 Denote by $h_{\alpha_1,\dotsc,\alpha_i}=jk_{\alpha_1,\dotsc,\alpha_i}$
the open immersions $V_{\alpha_1}\cap\dotsb\cap V_{\alpha_i}\rarrow X$.
 Applying the functor~$j_*$ to the exact
sequence~\eqref{qcoh-sheaf-G-on-scheme-Y-chech-coresolution}, we obtain
a finite exact sequence of quasi-coherent $\cA$\+modules with a right
action of the ring~$\cB(Y)$
\begin{multline} \label{direct-image-of-cech-coresolution}
 0\lrarrow j_*j^*\cE\lrarrow\bigoplus\nolimits_{\alpha=1}^N
 h_\alpha{}_*h_\alpha^*\cE\lrarrow
 \bigoplus\nolimits_{1\le\alpha<\beta\le N}
 h_{\alpha,\beta}{}_*h_{\alpha,\beta}^*\cE \\
 \lrarrow\dotsb\lrarrow h_{1,\dotsc,N}{}_*h_{1,\dotsc,N}^*\cE
 \lrarrow0.
\end{multline}

 In the case~(a), the terms of
the sequence~\eqref{direct-image-of-cech-coresolution} are
very flaprojective quasi-coherent $\cA$\+modules on $X$ by
Lemma~\ref{very-flaprojective-restriction-extension-of-scalars}(a),
as per the discussion in
Section~\ref{antilocality-of-X-contraadjusted-subsecn}.
 Using Lemma~\ref{very-flaprojective-cotorsion-pair-hereditary}(b)
and arguing by induction moving from the rightmost end of the exact
sequence ot the leftmost end, one shows that the $\cA$\+modules of
cocycles in~\eqref{direct-image-of-cech-coresolution} are also
very flaprojective.
 In view of Corollary~\ref{vflp-Ext-orthogonal-to-X-cta}, the functor
$\Hom_\cA({-},\C)$ transforms~\eqref{direct-image-of-cech-coresolution}
into an exact sequence of left $\cB(Y)$\+modules.

 In the case~(b), the terms of
the sequence~\eqref{direct-image-of-cech-coresolution} are
flaprojective quasi-coherent $\cA$\+modules on $X$ by
Lemma~\ref{flaprojective-restriction-extension-of-scalars}(a),
as per the discussion in
Section~\ref{antilocality-of-X-cotorsion-subsecn}.
 Using Lemma~\ref{flaprojective-cotorsion-pair-hereditary}(b)
and Corollary~\ref{flp-Ext-orthogonal-to-X-cot}, one similarly shows
that the functor $\Hom_\cA({-},\C)$
takes~\eqref{direct-image-of-cech-coresolution}
to an exact sequence of left $\cB(Y)$\+modules.

 In the case~(c), one arrives to the very same conclusion using
Corollary~\ref{quasi-algebra-change-of-scalars-adjustedness}(a),
the discussion in Section~\ref{antilocality-of-A-cotorsion-subsecn},
and the definition of an $\cA$\+cotorsion quasi-coherent $\cA$\+module.
 Indeed, the terms of
the sequence~\eqref{direct-image-of-cech-coresolution} are flat
quasi-coherent $\cA$\+modules on $X$ in this case.

 In all the cases~(a\+-c), we have obtained a finite exact sequence of
left $\cB(Y)$\+modules ending in the terms
\begin{multline} \label{cosections-of-fHom-computation-ends}
 \bigoplus\nolimits_{1\le\alpha<\beta\le N}
 \fHom_\cA(\cE,\C)[V_\alpha\cap V_\beta]\lrarrow
 \bigoplus\nolimits_{\alpha=1}^N\fHom_\cA(\cE,\C)[V_\alpha] \\
 \lrarrow\Hom_\cA(j_*j^*\cE,\C)\lrarrow0.
\end{multline}
 Thus we have computed the $\cB(Y)$\+module $\Hom_\cA(j_*j^*\cE,\C)$ in
a way that agrees with the description of the module of cosections of
the cosheaf $\fHom_\cA(\cE,\C)$ over the open subset $Y\subset X$
provided by the cosheaf axom~\eqref{cosheaf-axiom}.

 The proof of part~(d) is very similar to~\cite[proof of
Lemma~2.5.2(c)]{Pcosh}.
 We refrain from reiterating the details.
\end{proof}

\subsection{Contratensor product over a quasi-coherent quasi-algebra}
\label{contratensor-over-qcoh-quasi-algebra-subsecn}
 This section is a partial generalization of~\cite[Section~2.6]{Pcosh}.

 Let $X$ be a quasi-separated scheme, and let $\cA$ and $\cB$ be two
quasi-coherent quasi-algebras over~$X$.
 Let $\cE$ be a quasi-coherent $\cA$\+$\cB$\+bimodule on~$X$ (as defined
in Section~\ref{fHom-over-qcoh-quasi-algebra-subsecn}, and let $\P$
be a cosheaf of $\cB$\+modules on~$X$.
 The quasi-coherent left $\cA$\+module $\cE\ocn_\cB\P$ on $X$, called
the \emph{contratensor product} of $\cE$ and $\P$ over $\cA$,
is constructed as the following (nonfiltered) inductive limit in
the abelian category $\cA\Qcoh$.

 The diagram is indexed by the poset of all affine open subschemes
$U\subset X$.
 To every such affine open subscheme $U$ with the open immersion
morphism $j\:U\rarrow X$, the diagram assigns the quasi-coherent
$\cA$\+module
$$
 (j_*j^*\cE)\ot_{\cB(U)}\P[U].
$$
 Notice that the ring $\cB(U)$ acts naturally on the right on
the quasi-coherent $\cA|_U$\+module $j^*\cE$ on $U$, hence also
on the quasi-coherent $\cA$\+module $j_*j^*\cE$ on~$X$.
 So the tensor product, as defined in
Section~\ref{fHom-over-qcoh-quasi-algebra-subsecn}, makes sense.

 To a pair of affine open subschemes $V\subset U\subset X$ with open
immersion morphisms $j\:U\rarrow X$ and $k\:V\rarrow X$, the diagram
assigns the morphism of quasi-coherent $\cA$\+modules
\begin{equation} \label{contratensor-diagram-transition-morphism}
 (k_*k^*\cE)\ot_{\cB(V)}\P[V]\lrarrow(j_*j^*\cE)\ot_{\cB(U)}\P[U]
\end{equation}
constructed as follows.
 By~\eqref{direct-inverse-image-tensor-isomorphism}, we have
$k_*k^*\cE\simeq(j_*j^*\cE)\ot_{\cO_X(U)}\cO_X(V)$, hence
\begin{multline} \label{for-contratensor-computation}
 (k_*k^*\cE)\ot_{\cB(V)}\P[V]\simeq((j_*j^*\cE)\ot_{\cO_X(U)}\cO_X(V))
 \ot_{\cB(V)}\P[V] \\
 \simeq((j_*j^*\cE)\ot_{\cB(U)}(\cB(U)\ot_{\cO_X(U)}\cO_X(V))
 \ot_{\cB(V)}\P[V] \\
 \simeq((j_*j^*\cE)\ot_{\cB(U)}\cB(V))\ot_{\cB(V)}\P[V]
 \simeq(j_*j^*\cE)\ot_{\cB(U)}\P[V].
\end{multline}
 The desired morphism~\eqref{contratensor-diagram-transition-morphism}
is now induced by the corestriction morphism of $\cB(U)$\+modules
$\P[V]\rarrow\P[U]$.

 So we put
\begin{equation} \label{contratensor-product-defined-formula}
 \cE\ocn_\cB\P=\varinjlim\nolimits_U((j_*j^*\cE)\ot_{\cB(U)}\P[U]).
\end{equation}
 Now let $\C$ be a quasi-coherent $\cA$\+module on $X$ such that
the contraherent $\cB$\+module $\fHom_\cA(\cE,\C)$ is well-defined
as per one of the sufficient conditions in
Section~\ref{fHom-over-qcoh-quasi-algebra-subsecn}
(see Lemma~\ref{fHom-contraadjusted-and-more-lemma}).
 Then there is a natural (adjunction) isomorphism of abelian groups
\begin{equation} \label{fHom-contratensor-adjunction}
 \Hom^\cB(\P,\fHom_\cA(\cE,\C))\simeq
 \Hom_\cA(\cE\ocn_\cB\P,\>\C).
\end{equation}
 Indeed, both the abelian groups in~\eqref{fHom-contratensor-adjunction}
are naturally isomorphic to the group of all compatible collections
of morphisms of quasi-coherent $\cA$\+modules
$$
 (j_*j^*\cE)\ot_{\cB(U)}\P[U]\lrarrow\C
$$
on $X$, or equivalently, all the compatible collections of morphisms
of $\cB(U)$\+modules
$$
 \P[U]\rarrow\Hom_\cA(j_*j^*\cE,\C)
$$
defined for all the affine open subschemes $U\subset X$ with open
immersion morphisms $j\:U\rarrow X$.
 Notice that affine open subschemes form a topology base $\bB$ on $X$,
so a morphism of cosheaves of $\cB$\+modules on $X$ is uniquely
determined by an arbitrary morphism of the related copresheaves of
$\cB$\+modules on~$\bB$ (by
Theorem~\ref{extension-of-co-sheaves-from-topology-base}(b)).

 Now let $\bW$ be an open covering of $X$, and let $\bB$ be the topology
base consisting of all affine open subschemes subordinate to~$\bW$.
 Then, instead of computing the contratensor product $\E\ocn_\cB\P$
as the inductive limit over all affine open subschemes $U\subset X$
as above, one can restrict the construction to affine open subschemes
$U$ subordinate to $\bW$, and put
\begin{equation} \label{W-contratensor-product-definition-formula}
 \E\ocn_\cB\P=
 \varinjlim\nolimits_{U\in\bB}((j_*j^*\cE)\ot_{\cB(U)}\P[U]).
\end{equation}
 According to the argument above,
the adjunction~\eqref{fHom-contratensor-adjunction} still holds for
such more general definition of the contratensor product depending on
an open covering~$\bW$.
 In particular, \eqref{fHom-contratensor-adjunction}~holds for injective
quasi-coherent sheaves $\cC=\J$ on $X$, for which the contraherent
$\cB$\+module $\fHom_\cA(\E,\J)$ is well-defined by
Lemma~\ref{fHom-contraadjusted-and-more-lemma}(e).

 As there are enough injective quasi-coherent sheaves on any scheme $X$,
we can conclude that the contratensor product $\E\ocn_\cB\P$ as
defined by formula~\eqref{W-contratensor-product-definition-formula}
actually does not depend on an open covering $\bW$ of the scheme $X$,
and the two
definitions~\eqref{contratensor-product-defined-formula}
and~\eqref{W-contratensor-product-definition-formula} agree.
 We refer to~\cite[Section~2.6]{Pcosh} for a further discussion
with a futher generalization.

 For a quasi-coherent quasi-algebra $\cA$ over a quasi-compact
semi-separated scheme $X$, the notion of an antilocal contraherent
$\cA$\+module on $X$ was defined in
Sections~\ref{A-antilocal=X-antilocal-subsecn}\+-%
\ref{antilocal-contrah-A-modules-subsecn}.
 The notation $\cA\Ctrh_\al$ for the full exact subcategory of antilocal
contraherent $\cA$\+modules in the exact category $\cA\Ctrh\subset
\cA\Lcth_\bW$ was introduced in
Section~\ref{antilocal-contrah-A-lct-A-modules-subsecn}.

\begin{lem} \label{contratensor-exact-on-antilocal}
 Let $X$ be a quasi-compact semi-separated scheme, and let $\cA$ and
$\cB$ be quasi-coherent quasi-algebras over~$X$.
 Then the functor of contratensor product\/~$\ocn_\cB$ is exact
on the Cartesian product of the full subcategories of $\cB$\+flat
quasi-coherent $\cA$\+$\cB$\+bimodules $\cA\biQcoh\cB^{\cB\dfl}
\subset\cA\biQcoh\cB$ and antilocal contraherent $\cB$\+modules
$\cB\Ctrh_\al\subset(X,\cB)\Cosh$.
 So the functor
\begin{equation} \label{contratensor-B-flat-antilocal-exact}
 \ocn_\cB\:\cA\biQcoh\cB^{\cB\dfl}\times
 \cB\Ctrh_\al\lrarrow\cA\Qcoh
\end{equation}
is exact as a functor of two arguments acting between exact/abelian
categories.
\end{lem}

\begin{proof}
 The Grothendieck abelian category $\cA\Qcoh$ has enough injective
objects.
 Hence it suffices to show that, for any fixed injective quasi-coherent
$\cA$\+module $\J\in\cA\Qcoh^\inj$, the functor $(\cE,\P)\longmapsto
\Hom_\cA(\cE\ocn_\cB\P,\>\J)$ is exact as a functor of two arguments
$(\cA\biQcoh\cB^{\cB\dfl})^\sop\times(\cB\Ctrh_\al)^\sop\rarrow\Ab$.
 According to the adjunction
formula~\eqref{fHom-contratensor-adjunction} and taking into account
Lemma~\ref{fHom-contraadjusted-and-more-lemma}(d), we have
$$
 \Hom_\cA(\cE\ocn_\cB\P,\>\J)\simeq\Hom^\cB(\P,\fHom_\cA(\cE,\J)),
$$
where $\fHom_\cA(\cE,\J)$ is a $\cB$\+locally injective contraherent
$\cB$\+module.
 It remains to point out that the functor $\cE\longmapsto
\fHom_\cA(\cE,\J)\:(\cA\biQcoh\cB^{\cB\dfl})^\sop\rarrow
\cB\Ctrh^{\cB\dlin}$ is exact as mentioned
in~\eqref{fHom-B-flat-inj-B-dlin-exact}, while the functor
$\Hom^\cB\:(\cB\Ctrh_\al)^\sop\times\cB\Lcth_\bW^{\cB\dlin}\rarrow\Ab$
is exact (for any open covering $\bW$ of the scheme~$X$) by
Corollary~\ref{X-antilocal-equivalent-to-A-antilocal}(3).
\end{proof}

\subsection{Projection formulas}
\label{fHom-contratensor-projection-formulas}
 This section is a very partial generalization of the results
of~\cite[Section~3.8]{Pcosh} concerning the contraherent $\fHom$
and contratensor product functors.
 We restrict ourselves to \emph{open immersions} $f\:Y\rarrow X$.

 Let $X$ be a scheme and $Y\subset X$ be an open subscheme such that
the open immersion morphism $f\:Y\rarrow X$ is quasi-compact.
 Let $\cA$ and $\cB$ be quasi-coherent quasi-algebras over $X$, let
$\cE$ be a quasi-coherent $\cA$\+$\cB$\+bimodule on $X$, and let $\cD$
be a quasi-coherent left $f^*\cA$\+module on~$Y$.

 Assume that one of the five sufficient conditions of
Section~\ref{fHom-over-qcoh-quasi-algebra-subsecn}
(see Lemma~\ref{fHom-contraadjusted-and-more-lemma})
for the existence of a contraherent cosheaf $\fHom$ is satisfied
simultaneously for $\cE$ on $X$, for $\cD$ on $Y$, and for
the scheme $X$ itself.
 In this case, the same conditions are also satisfied for
the quasi-coherent $f^*\cA$\+$f^*\cB$\+module $f^*\cE$ on $Y$,
for the quasi-coherent $\cA$\+module $f_*\cD$ on $X$, and for
the scheme~$Y$.
 In particular:
\begin{itemize}
\item If the quasi-coherent $f^*\cA$\+module $\cD$ on $Y$ is
$Y$\+contraadjusted, then the quasi-coherent $\cA$\+module $f_*\cD$
on $X$ is $X$\+contraadjusted, as mentioned
in Section~\ref{antilocality-of-X-contraadjusted-subsecn}
and explained in~\cite[Section~2.5]{Pcosh}.
\item If the quasi-coherent $f^*\cA$\+module $\cD$ on $Y$ is
$Y$\+cotorsion, then the quasi-coherent $\cA$\+module $f_*\cD$
on $X$ is $X$\+cotorsion, as mentioned
in Section~\ref{antilocality-of-X-cotorsion-subsecn}
and explained in~\cite[Section~2.5]{Pcosh}.
\item If the quasi-coherent $f^*\cA$\+module $\cD$ on $Y$ is
$f^*\cA$\+cotorsion, then the quasi-coherent $\cA$\+module $f_*\cD$
on $X$ is $\cA$\+cotorsion, as explained
in Section~\ref{antilocality-of-A-cotorsion-subsecn}.
\item If the quasi-coherent $f^*\cA$\+module $\cD$ on $Y$ is
injective, then the quasi-coherent $\cA$\+module $f_*\cD$ on $X$
is injective, since the functor $f_*\:(f^*\cA)\Qcoh\rarrow\cA\Qcoh$
is right adjoint to an exact functor $f^*\:\cA\Qcoh\rarrow
(f^*\cA)\Qcoh$.
\end{itemize}

 Under the assumptions above, there is a natural isomorphism of
cosheaves of $\cB$\+modules on~$X$
\begin{equation} \label{fHom-projection-formula}
 \fHom_\cA(\cE,f_*\cD)\simeq f_!\fHom_{f^*\cA}(f^*\cE,\cD),
\end{equation}
where the cosheaf of $\cB$\+modules $\fHom_\cA(\cE,f_*\cD)$ on $X$
is contraherent by Section~\ref{fHom-over-qcoh-quasi-algebra-subsecn} 
(hence it follows from~\eqref{fHom-projection-formula} that the cosheaf
of $\cB$\+modules $f_!\fHom_{f^*\cA}(f^*\cE,\cD)$ on $X$
is contraherent, too).
 The formula~\eqref{fHom-projection-formula} is our version
of~\cite[formulas~(3.23) and~(3.24)]{Pcosh}.

 Indeed, let $U\subset X$ be an affine open subscheme.
 Put $V=Y\cap U$, and notice that the scheme $V$ is quasi-compact
(since the morphism~$f$ is quasi-compact by assumption).
 Furthermore, the open immersion morphism $V\rarrow Y$ is affine
whenever the scheme $X$ is semi-separated (since the open immersion
morphism $U\rarrow X$ is affine in this case).
 Denote the open immersion morphisms by $j\:U\rarrow X$, \
$j'\:V\rarrow Y$, and $f'\:V\rarrow U$.
 Then we have natural isomorphisms of $\cB(U)$\+modules
\begin{multline*}
 \fHom_\cA(\cE,f_*\cD)[U]=\Hom_\cA(j_*j^*\cE,f_*\cD)\simeq
 \Hom_{f^*\cA}(f^*j_*j^*\cE,\cD) \\
 \simeq\Hom_{f^*\cA}(j'_*f'{}^*j^*\cE,\cD)\simeq
 \Hom_{f^*\cA}(j'_*j'{}^*f^*\cE,\cD) \\
 \simeq\fHom_{f^*\cA}(f^*\cE,\cD)[V]
 =(f_!\fHom_{f^*\cA}(f^*\cE,\cD))[U],
\end{multline*}
as desired.
 Here the base change isomorphism $f^*j_*\M\simeq j'_*f'{}^*\M$ holds
for all quasi-coherent sheaves $\M$ on $U$, since $j$~is a quasi-compact
and quasi-separated morphism of schemes, while $f$~is a flat morphism
of schemes.
 The isomorphism $\Hom_{f^*\cA}(j'_*j'{}^*f^*\cE,\cD)\simeq
\fHom_{f^*\cA}(f^*\cE,\cD)[V]$ holds by
Lemma~\ref{cosections-of-fHom-lemma}.

 We refer to Sections~\ref{A-antilocal=X-antilocal-subsecn}\+-%
\ref{antilocal-contrah-A-modules-subsecn} for the definition of
an \emph{antilocal} contraherent $\cA$\+module on a quasi-compact
semi-separated scheme~$X$.

\begin{cor} \label{fHom-antilocal-corollary}
 Let $X$ be a quasi-compact semi-separated scheme, and let $\cA$ and
$\cB$ be quasi-coherent quasi-algebras over~$X$.
 Let $\cE$ be a quasi-coherent $\cA$\+$\cB$\+bimodule and $\C$ be
a quasi-coherent left $\cA$\+module on~$X$.
 Assume that one of the five sufficient conditions of
Section~\ref{fHom-over-qcoh-quasi-algebra-subsecn}
(see Lemma~\ref{fHom-contraadjusted-and-more-lemma})
for the existence of a contraherent cosheaf\/ $\fHom_\cA(\cE,\C)$
is satisfied for $\cE$ and~$\C$.
 Then the contraherent $\cB$\+module\/ $\fHom_\cA(\cE,\C)$ on $X$
is antilocal.
\end{cor}

\begin{proof}
 The argument is based on formula~\eqref{fHom-projection-formula} for
the immersions of affine open subschemes $f=j\:U\rarrow X$ and
the results of Section~\ref{antilocal-classes-secn}.

 Let $X=\bigcup_{\alpha=1}^N U_\alpha$ be a finite affine open
covering of the scheme~$X$.
 Under the assumptions of
Lemma~\ref{fHom-contraadjusted-and-more-lemma}(a), it is important
that every $X$\+contraadjusted quasi-coherent $\cA$\+module $\C$
on $X$ is a direct summand of a finitely iterated extension of
the direct images of $U_\alpha$\+contraadjusted quasi-coherent
$\cA|_{U\alpha}$\+modules from~$U_\alpha$ (see
Theorem~\ref{qcomp-qsep-very-flaproj-complete-cotorsion-pair-thm}(c)).
 Under the assumptions of
Lemma~\ref{fHom-contraadjusted-and-more-lemma}(b), it is important
that every $X$\+cotorsion quasi-coherent $\cA$\+module $\C$
on $X$ is a direct summand of a finitely iterated extension of
the direct images of $U_\alpha$\+cotorsion quasi-coherent
$\cA|_{U\alpha}$\+modules from~$U_\alpha$ (see
Theorem~\ref{qcomp-qsep-flaproj-complete-cotorsion-pair-thm}(c)).
 Under the assumptions of
Lemma~\ref{fHom-contraadjusted-and-more-lemma}(c), it is important
that every $\cA$\+cotorsion quasi-coherent $\cA$\+module $\C$
on $X$ is a direct summand of a finitely iterated extension of
the direct images of $\cA|_{U_\alpha}$\+cotorsion quasi-coherent
$\cA|_{U\alpha}$\+modules from~$U_\alpha$ (see
Theorem~\ref{qcomp-qsep-A-flat-complete-cotorsion-pair-thm}(c)).
 Under the assumptions of
Lemma~\ref{fHom-contraadjusted-and-more-lemma}(d) or~(e), one needs
to use the fact that every injective quasi-coherent $\cA$\+module $\J$
on $X$ is a direct summand of a finite direct sum of the direct images
of injective quasi-coherent $\cA|_{U\alpha}$\+modules from~$U_\alpha$
(see Lemma~\ref{injective-qcoh-modules-strongly-antilocal} below).

 Then one needs to use the facts that
the functors~(\ref{fHom-A-very-flaproj-X-cta-exact}\+-%
\ref{fHom-A-flat-A-cot-B-lct-exact}) are exact; so they preserve
finitely iterated extensions in the second argument.
 Formula~\eqref{fHom-projection-formula} plays the key role; and
then it remains to recall that the direct summands of finitely
iterated extensions of the direct images of contraherent
$\cA|_{U\alpha}$\+modules from~$U_\alpha$ are antilocal contraherent
$\cA$\+modules on~$X$ (this is the easy implication in
Theorem~\ref{qcomp-qsep-antilocal-complete-cotorsion-pair-thm}(c)).
\end{proof}

 Now assume that the scheme $X$ is quasi-separated and the open
immersion morphism $f\:Y\rarrow X$ is affine.
 Let $\Q$ be a cosheaf of $f^*\cB$\+modules on~$Y$.
 Then there is a natural isomorphism of quasi-coherent $\cA$\+modules
on~$X$
\begin{equation} \label{contratensor-projection-formula}
 \cE\ocn_\cB f_!\Q\simeq f_*(f^*\cE\ocn_{f^*\cB}\Q).
\end{equation}
 The formula~\eqref{contratensor-projection-formula} is our version
of~\cite[formula~(3.25)]{Pcosh}.

 Indeed, in the notation above, and denoting also the open immersion
$V\rarrow X$ by $k=jf'=fj'$, we have natural isomorphisms of
quasi-coherent $\cA$\+modules on~$X$
\begin{multline*}
 (j_*j^*\cE)\ot_{\cB(U)}((f_!\Q)[U])=(j_*j^*\cE)\ot_{\cB(U)}\Q[V]
 \simeq(k_*k^*\cE)\ot_{\cB(V)}\Q[V] \\
 \simeq(f_*j'_*j'{}^*f^*\cE)\ot_{\cB(V)}\Q[V]\simeq
 f_*((j'_*j'{}^*f^*\cE)\ot_{\cB(V)}\Q[V])
\end{multline*}
where the second isomorphism is~\eqref{for-contratensor-computation},
and the last isomorphism holds since we are assuming that the morphism
of schemes~$f$ is affine.
 It remains to point out that the direct image functor
$f_*\:(f^*\cA)\Qcoh\rarrow\cA\Qcoh$ for an affine morphism~$f$
preserves all inductive limits.

 The following lemma, which is our version of~\cite[construction in 
Section~3.8 and Lemma~8.4.1]{Pcosh}, extends the natural
isomorphism~\eqref{contratensor-projection-formula} to certain
\emph{not necessarily affine} open immersion morphisms~$f$.
 However, part~(b) requires more restrictive assumptions on
the cosheaf~$\Q$.

 The notation $\cA\Ctrh_\al$ for the exact category of antilocal
contraherent $\cA$\+modules on a scheme $X$ was introduced in
Section~\ref{antilocal-contrah-A-lct-A-modules-subsecn}, and
the notation $X\Ctrh_\al$ for the exact category of antilocal
contraherent cosheaves on $X$ was also mentioned there.

\begin{lem} \label{contratensor-projection-comparison-and-antilocal}
 Let $X$ a scheme and $Y\subset X$ be an open subscheme with the open
immersion morphism $f\:Y\rarrow X$.
 Let $\cA$ and $\cB$ be quasi-coherent quasi-algebras over $X$, let
$\cE$ be a quasi-coherent $\cA$\+$\cB$\+bimodule on $X$, and let
$\Q$ be a cosheaf of $f^*\cB$\+modules on~$Y$.
 In this context: \par
\textup{(a)} If the scheme $X$ is quasi-separated and the open immersion
morphism~$f$ is quasi-compact, then there is a natural morphism
of quasi-coherent $\cA$\+modules on~$X$
\begin{equation} \label{contratensor-projection-comparison-morphism}
 \cE\ocn_\cB f_!\Q\lrarrow f_*(f^*\cE\ocn_{f^*\cB}\Q).
\end{equation} \par
\textup{(b)} If the schemes $X$ and $Y$ are quasi-compact and
semi-separated, the quasi-coherent $\cA$\+$\cB$\+bimodule $\cE$ on $X$
is $\cB$\+flat, and $\Q$ is an antilocal contraherent $f^*\cB$\+module
on $Y$, then
the morphism~\eqref{contratensor-projection-comparison-morphism}
is an isomorphism.
\end{lem}

\begin{proof}
 Part~(a): the construction is similar to the one
in~\cite[Section~3.8]{Pcosh}, but a bit simpler, essentially because
the morphism of quasi-ringed schemes $f\:(Y,f^*\cB)\rarrow
(X,\cB)$ is left flat in our case.
 Let $U\subset X$ and $V\subset Y$ be affine open subschemes such
that $V\subset U$.
 As above, we denote the open immersion morphisms by $j\:U\rarrow X$,
\ $j'\:V\rarrow Y$, and $f'\:V\rarrow U$.
 The desired morphism is produced as the composition of natural
morphism and isomorphisms
\begin{multline*}
 \cE\ocn_\cB f_!\Q=\varinjlim\nolimits_{U\subset X}
 j_*j^*\cE\ocn_{\cB(U)}\Q[Y\cap U] \simeq
 \varinjlim\nolimits_{U\subset X}j_*j^*\cE\ocn_{\cB(U)}
 (\varinjlim\nolimits_{V\subset Y\cap U}\Q[V]) \\
 \simeq\varinjlim\nolimits_{U\subset X}
 \varinjlim\nolimits_{V\subset Y\cap U}j_*j^*\cE\ocn_{\cB(U)}\Q[V]
 \simeq
 \varinjlim\nolimits_{(U,V):V\subset Y\cap U}j_*j^*\cE\ocn_{\cB(U)}\Q[V]
 \\ \simeq \varinjlim\nolimits_{(U,V):V\subset Y\cap U}
 f_*j'_*j'{}^*f^*\cE\ocn_{\cB(V)}\Q[V] \\
 \lrarrow f_*(\varinjlim\nolimits_{(U,V):V\subset Y\cap U}
 j'_*j'{}^*f^*\cE\ocn_{\cB(V)}\Q[V])\\ 
 \simeq f_*(\varinjlim\nolimits_{V\subset Y}
 j'_*j'{}^*f^*\cE\ocn_{\cB(V)}\Q[V]) = f_*(f^*\cE\ocn_{f^*\cB}\Q).
\end{multline*}
 Here the isomorphism in the third line holds by
formula~\eqref{for-contratensor-computation}.

 Part~(b): the argument is similar
to~\cite[proof of Lemma~8.4.1]{Pcosh}.
 The idea is to show that both the left-hand side and the right-hand
side of~\eqref{contratensor-projection-comparison-morphism} are exact
as functors $f^*\cB\Ctrh_\al\rarrow\cA\Qcoh$ (or equivalently,
$f^*\cB\Ctrh_\al\rarrow X\Qcoh$) of the argument
$\Q\in f^*\cB\Ctrh_\al$ for any fixed $\cB$\+flat quasi-coherent
$\cA$\+$\cB$\+bimodule $\cE$ on~$X$.

 The direct image functor $f_!\:(Y,\cO_Y)\Cosh\rarrow(X,\cO_X)\Cosh$
takes antilocal contraherent cosheaves on $Y$ to antilocal contraherent
cosheaves on $X$, and is exact as a functor $f_!\:Y\Ctrh_\al\rarrow
X\Ctrh_\al$, by~\cite[Corollary~4.6.3(a)]{Pcosh}.
 It follows that the direct image functor $f_!\:(Y,f^*\cB)\Cosh\rarrow
(X,\cB)\Cosh$ takes antilocal contraherent $f^*\cB$\+modules on $Y$ to
antilocal contraherent $\cB$\+modules on $X$, and is exact as
a functor $f_!\:f^*\cB\Ctrh_\al\rarrow\cB\Ctrh_\al$.
 On the other hand, the functor $\cE\ocn_\cB{-}\,\:\cB\Ctrh_\al\rarrow
\cA\Qcoh$ is exact by Lemma~\ref{contratensor-exact-on-antilocal}.
 This proves that the left-hand side
of~\eqref{contratensor-projection-comparison-morphism} is exact
as a functor of $\Q\in f^*\cB\Ctrh_\al$.

 The proof of exactness of the right-hand side
of~\eqref{contratensor-projection-comparison-morphism} is based on
the notion of a \emph{dilute} quasi-coherent sheaf on a quasi-compact
semi-separated scheme, which was suggested in~\cite[Section~4.1]{M-n}
and elaborated upon in~\cite[Section~4.2]{Pcosh}.
 The functor $f^*\cE\ocn_{f^*\cB}{-}\,\:f^*\cB\Ctrh_\al\rarrow
f^*\cA\Qcoh$ is exact by Lemma~\ref{contratensor-exact-on-antilocal}.
 As the functor~$f_*$ of direct image of quasi-coherent sheaves with
respect to a morphism~$f$ of quasi-compact semi-separated schemes is
exact on the full exact subcategory of dilute quasi-coherent
sheaves~\cite[Corollary~4.2.7]{Pcosh}, it suffices to show that
the quasi-coherent left $f^*\cA$\+module $f^*\cE\ocn_{f^*\cB}\Q$ is
dilute as a quasi-coherent sheaf on $Y$ for all $\Q\in f^*\cB\Ctrh_\al$.
 The class of dilute quasi-coherent sheaves is closed under
extensions in $Y\Qcoh$ \,\cite[Section~4.2]{Pcosh}, so in view of
Theorem~\ref{qcomp-qsep-antilocal-complete-cotorsion-pair-thm}(c) it
suffices to consider the case of a contraherent $f^*\cB$\+module $\Q$
of the form $\Q=j'_!\R$, where $j'\:V\rarrow Y$ in the open immersion
of an affine open subscheme $V\subset Y$ and $\R$ is a contraherent
$j'{}^*f^*\cB$\+module on~$V$.
 By~\eqref{contratensor-projection-formula}, we have
$f^*\cE\ocn_{f^*\cB}j'_!\R\simeq
j'_*(j'{}^*f^*\cE\ocn_{j'{}^*f^*\cB}\R)$, and it remains to mention
that the direct images of quasi-coherent sheaves from affine open
subschemes are dilute on quasi-compact semi-separated
schemes~\cite[Lemma~4.2.1 or Corollary~4.2.7]{Pcosh}.
 Notice that all quasi-coherent sheaves on affine schemes
are dilute~\cite[Section~4.2]{Pcosh}.

 Now that both sides
of~\eqref{contratensor-projection-comparison-morphism} are known to be
exact as functors of $\Q\in f^*\cB\Ctrh_\al$,
Theorem~\ref{qcomp-qsep-antilocal-complete-cotorsion-pair-thm}(c)
implies that, in order to show
that~\eqref{contratensor-projection-comparison-morphism} is
an isomorphism, it suffices to consider the case of $\Q=j'_!\R$
as above.
 In this case, the assertion follows from
the formula~\eqref{contratensor-projection-formula} applied to two
affine open immersion morphisms $j'\:V\rarrow Y$ and $fj'\:V\rarrow X$.
 One needs to use commutativity of the triangular diagram
$$
 \cE\ocn_\cB f_!j'_!\R\lrarrow f_*(f^*\cE\ocn_{f^*\cB}j'_!\R)
 \lrarrow f_*j'_*(j'{}^*f^*\cE\ocn_{j'{}^*f^*\cB}\R)
$$
of natural morphisms of quasi-coherent $\cA$\+modules on~$X$.
 The rightmost morphism and the composition are isomorphisms
by~\eqref{contratensor-projection-formula}, hence the leftmost morphism
is an isomorphism, too, as desired.
\end{proof}

\subsection{Compatibility with forgetting the action of~$\cA$}
 Let $X$ be a scheme.
 Then the structure sheaf $\cO_X$ of $X$ is an obvious example
of a quasi-coherent quasi-algebra over~$X$.
 The structure sheaf $\cE=\cO_X$ can be also viewed as a quasi-coherent
$\cO_X$\+$\cO_X$\+bimodule.

 One can specialize the discussion
in Sections~\ref{fHom-over-qcoh-quasi-algebra-subsecn}\+-%
\ref{contratensor-over-qcoh-quasi-algebra-subsecn} to the case of
$\cA=\cO_X=\cB$ and a quasi-coherent sheaf $\cE$ on $X$ (viewed as
a quasi-coherent quasi-module on $X$ via the natural inclusion functor
$X\Qcoh\rarrow X\QQcoh$, as per the discussion in
Section~\ref{quasi-coherent-quasi-algebras-subsecn}).
 This leads to the definitions of the partially defined functor
$\fHom_{\cO_X}({-},{-})$ and the functor ${-}\ocn_{\cO_X}{-}$, acting
between the categories of quasi-coherent sheaves and contraherent
cosheaves (or more generally, cosheaves of $\cO_X$\+modules) on~$X$.
 These functors agree with the functors $\fHom_X$ and $\ocn_X$
defined in~\cite[Sections~2.5\+-2.6]{Pcosh}.
 So we will write $\fHom_{\cO_X}=\fHom_X$ and $\ocn_{\cO_X}=\ocn_X$.

 Now let $\cA$ be an arbitrary quasi-coherent quasi-algebra over~$X$.
 Then the quasi-coherent $\cA$\+$\cA$\+bimodule $\cA$ is
$\cA$\+very flaprojective, $\cA$\+robustly flaprojective, and
$\cA$\+flat on both sides, as mentioned in
Section~\ref{fHom-over-qcoh-quasi-algebra-subsecn}.

 The notation $\cA\Qcoh^{\cA\dcot}\subset\cA\Qcoh^{X\dcot}
\subset\cA\Qcoh^{X\dcta}$ for exact subcategories in the abelian
category $\cA\Qcoh$ can be found in
Section~\ref{background-derived-quasi-coherent}.
 The similar notation $X\Qcoh^\cot\subset X\Qcoh^\cta\subset X\Qcoh$
was mentioned in the proofs of
Lemma~\ref{X-A-cta-cot-coresolving-coresolution-dimension}(b)
and Theorem~\ref{qcoh-X-cot-A-cot-derived-equivalence};
see~\cite[Section~4.1]{Pcosh} for a detailed discussion.

\begin{lem} \label{fHom-forgetting-A-module-structure-lemma}
 Let $X$ be a scheme and $\cA$ be a quasi-coherent quasi-algebra
over~$X$.  In this context: \par
\textup{(a)} Assume that the scheme $X$ is quasi-compact and
semi-separated.
 Then there is a commutative diagram of exact functors between
exact categories
\begin{equation} \label{X-cta-fHom-forgetful-diagram}
\begin{gathered}
 \xymatrix{
  \cA\Qcoh^{X\dcta} \ar[rr]^-{\fHom_\cA(\cA,{-})} \ar[d]
  && \cA\Ctrh \ar[d] \\
  X\Qcoh^\cta \ar[rr]_-{\fHom_X(\cO_X,{-})} && X\Ctrh
 }
\end{gathered}
\end{equation} \par
\textup{(b)} Assume that the scheme $X$ is quasi-compact and
semi-separated.
 Then there is a commutative diagram of exact functors between
exact categories
\begin{equation} \label{X-cot-fHom-forgetful-diagram}
\begin{gathered}
 \xymatrix{
  \cA\Qcoh^{X\dcot} \ar[rr]^-{\fHom_\cA(\cA,{-})} \ar[d]
  && \cA\Ctrh^{X\dlct} \ar[d] \\
  X\Qcoh^\cot \ar[rr]_-{\fHom_X(\cO_X,{-})} && X\Ctrh^\lct
 }
\end{gathered}
\end{equation} \par
\textup{(c)} Assume that the scheme $X$ is semi-separated.
 Then there is a commutative diagram of exact functors between
exact categories
\begin{equation} \label{A-cot-fHom-forgetful-diagram}
\begin{gathered}
 \xymatrix{
  \cA\Qcoh^{\cA\dcot} \ar[rr]^-{\fHom_\cA(\cA,{-})} \ar[d]
  && \cA\Ctrh^{\cA\dlct} \ar[d] \\
  X\Qcoh^\cot \ar[rr]_-{\fHom_X(\cO_X,{-})} && X\Ctrh^\lct
 }
\end{gathered}
\end{equation} \par
\textup{(d)} Assume that the scheme $X$ is semi-separated.
 Then there is a commutative diagram of exact functors between
additive/exact categories
\begin{equation} \label{A-inj-A-lin-X-cot-forgetful-diagram}
\begin{gathered}
 \xymatrix{
  \cA\Qcoh^\inj \ar[rr]^-{\fHom_\cA(\cA,{-})} \ar[d]
  && \cA\Ctrh^{\cA\dlin} \ar[d] \\
  X\Qcoh^\cot \ar[rr]_-{\fHom_X(\cO_X,{-})} && X\Ctrh^\lct
 }
\end{gathered}
\end{equation} \par
\textup{(e)} Assume that the scheme $X$ is semi-separated and
the quasi-coherent quasi-algebra $\cA$ is flat as a quasi-coherent
sheaf on $X$ with respect to the right $\cO_X$\+module structure.
 Then there is a commutative diagram of exact functors between
additive/exact categories
\begin{equation} \label{A-inj-A-lin-X-inj-forgetful-diagram}
\begin{gathered}
 \xymatrix{
  \cA\Qcoh^\inj \ar[rr]^-{\fHom_\cA(\cA,{-})} \ar[d]
  && \cA\Ctrh^{\cA\dlin} \ar[d] \\
  X\Qcoh^\inj \ar[rr]_-{\fHom_X(\cO_X,{-})} && X\Ctrh^\lin
 }
\end{gathered}
\end{equation} \par
\textup{(f)} Assume that the scheme $X$ is quasi-separated and
the quasi-coherent quasi-algebra $\cA$ is flat as a quasi-coherent
sheaf on $X$ with respect to the right $\cO_X$\+module structure.
 Then there is a commutative diagram of exact functors between
additive/exact categories
\begin{equation} \label{A-inj-A-lct-X-inj-forgetful-diagram}
\begin{gathered}
 \xymatrix{
  \cA\Qcoh^\inj \ar[rr]^-{\fHom_\cA(\cA,{-})} \ar[d]
  && \cA\Ctrh^{\cA\dlct} \ar[d] \\
  X\Qcoh^\inj \ar[rr]_-{\fHom_X(\cO_X,{-})} && X\Ctrh^\lct
 }
\end{gathered}
\end{equation}
 On all the diagrams~\textup{(\ref{X-cta-fHom-forgetful-diagram}\+-%
\ref{A-inj-A-lct-X-inj-forgetful-diagram})},
the vertical arrows denote the forgetful functors.
\end{lem}

\begin{proof}
 The leftmost vertical (forgetful) functors on all the diagrams are
the respective restrictions of the exact forgetful functor
$\cA\Qcoh\rarrow X\Qcoh$.
 The leftmost vertical functors in~(a) and~(b) are well-defined by
the definitions of $\cA\Qcoh^{X\dcta}$ and $\cA\Qcoh^{X\dcot}$.
 The leftmost vertical functors in~(c) and~(d) are well-defined
according to the argument in
Section~\ref{antilocality-of-A-cotorsion-subsecn}.
 The leftmost vertical functors in~(e) and~(f) are well-defined since,
under the assumptions of~(e\+-f), the forgetful functor
$\cA\Qcoh\rarrow X\Qcoh$ is right adjoint to an exact functor
$\cA\ot_{\cO_X}{-}\,\:X\Qcoh\rarrow\cA\Qcoh$; so this forgetful
functor takes injective objects to injective objects.

 The rightmost vertical (forgetful) functor in part~(a) is provided
by the definition of the exact category $\cA\Ctrh$.
 All the other rightmost vertical functors are the respective
restrictions of the one in part~(a).
 The rightmost vertical functor in part~(b) is well-defined by
the definition of $\cA\Ctrh^{X\dlct}$.
 The rightmost vertical functors in parts~(c\+-d) are well-defined by
Lemma~\ref{restriction-coextension-injective-cotorsion}(a), as
mentioned in Section~\ref{A-loc-cotors-loc-inj-cosheaves-subsecn}.
 The rightmost vertical functors in parts~(e\+-f) are well-defined by
Lemma~\ref{restriction-coextension-injective-cotorsion}(b), as also
mentioned in Section~\ref{A-loc-cotors-loc-inj-cosheaves-subsecn}.

 All the lower horizontal functors are provided by the construction
and discussion of~\cite[Section~2.5]{Pcosh}.
 All the upper horizontal functors are provided by the construction
of Section~\ref{fHom-over-qcoh-quasi-algebra-subsecn} and the respective
parts of Lemma~\ref{fHom-contraadjusted-and-more-lemma}.

 Let us explain why the diagrams are commutative.
 In all the cases, this follows from the constructions of the functors
$\fHom_\cA$ and $\fHom_X$ in view of the natural isomorphism of
quasi-coherent left $\cA$\+modules
\begin{equation} \label{quasi-algebra-itself-projection-formula}
 j_*j^*\cA\simeq j_*(j^*\cA\ot_{\cO_U}\cO_U)\simeq
 \cA\ot_{\cO_X}j_*j^*\cO_X
\end{equation}
for any affine open subscheme $U\subset X$ with the open immersion
morphism $j\:U\rarrow X$ and the fact that the functor
$\cA\ot_{\cO_X}{-}\,\:X\Qcoh\rarrow\cA\Qcoh$ is left adjoint to
the forgetful functor $\cA\Qcoh\rarrow X\Qcoh$.
 In all the cases~(a\+-e),
\,\eqref{quasi-algebra-itself-projection-formula}~is a special case of
the projection formula~\eqref{quasi-module-tensor-projection-formula}.
 In the case~(f), where the open immersion morphism~$j$ need not be
affine, the isomorphism~\eqref{quasi-algebra-itself-projection-formula}
holds due to the assumption that $\cA$ is flat as a quasi-coherent
sheaf in the right $\cO_X$\+module structure.

 Finally, let us point out an alternative argument proving that
the functor $\fHom_\cA(\cA,{-})$ takes $X$\+cotorsion quasi-coherent
$\cA$\+modules to $X$\+locally cotorsion contraherent $\cA$\+modules
in part~(b) without a reference to
Lemma~\ref{fHom-contraadjusted-and-more-lemma}(b) and the complicated
notion of a robustly flaprojective quasi-coherent bimodule.
 Instead of the reference to
Lemma~\ref{fHom-contraadjusted-and-more-lemma}(b), in the proof of
part~(b) one could just refer to the commutative
diagram~\eqref{X-cta-fHom-forgetful-diagram} in part~(a)
and the fact that the functor $\fHom_X(\cO_X,{-})$ takes cotorsion
quasi-coherent sheaves to locally cotorsion contraherent
cosheaves~\cite[Section~2.5]{Pcosh}.
 This is based on Lemma~\ref{fHom-contraadjusted-and-more-lemma}(a),
but avoids using Lemma~\ref{fHom-contraadjusted-and-more-lemma}(b).
\end{proof}

\begin{lem} \label{contratensor-forgetting-A-module-structure-lemma}
 Let $X$ be a scheme and $\cA$ be a quasi-coherent quasi-algebra
over~$X$.
 Assume that either \par
\textup{(a)} the scheme $X$ is semi-separated, or \par
\textup{(b)} the scheme $X$ is quasi-separated and $\cA$ is flat
as a quasi-coherent sheaf on $X$ with respect to the left
$\cO_X$\+module structure. \par
 Then there is a commutative diagram of additive functors
\begin{equation} \label{quasi-algebra-contratensor-forgetful-diagram}
\begin{gathered}
 \xymatrix{
  (X,\cA)\Cosh \ar[rr]^-{\cA\ocn_\cA{-}} \ar[d] && \cA\Qcoh \ar[d] \\
  (X,\cO_X)\Cosh \ar[rr]_-{\cO_X\ocn_X{-}} && X\Qcoh
 }
\end{gathered}
\end{equation}
where the vertical arrows denote the forgetful functors.
\end{lem}

\begin{proof}
 The constructions of the vertical (forgetful) functors are obvious.
 The lower horizontal functor was constructed
in~\cite[Section~2.6]{Pcosh}.
 The upper horizontal functor is provided by the construction
of Section~\ref{contratensor-over-qcoh-quasi-algebra-subsecn}.

 Let us explain why the diagram is commutative.
 In both cases~(a) and~(b), this follows from the constructions of
the functors $\ocn_\cA$ and~$\ocn_X$ in view of the natural isomorphism
\begin{equation} \label{quasi-algebra-itself-tensor-direct-image}
 j_*j^*\cA\simeq j_*(j^*\cO_X\ot_{\cO_X(U)}\cA(U))
 \simeq(j_*j^*\cO_X)\ot_{\cO_X(U)}\cA(U)
\end{equation}
of quasi-coherent sheaves on $X$ with a right action of
the ring~$\cA(U)$.
 Here $U\subset X$ is an affine open subscheme with the open
immersion morphism $j\:U\rarrow X$.
 In case~(a), the isomorphism holds because $j$~is an affine morphism;
in case~(b), because $\cA(U)$ is a flat left $\cO_X(U)$\+module.

 One also needs to use the obvious isomorphism of tensor products
$$
 (\F\ot_RS)\ot_SP\simeq\F\ot_RP
$$
for any quasi-coherent sheaf $\F$ on a scheme $X$ with a right action
of a ring $R$ on $\F$, any associative ring homomorphism $R\rarrow S$,
and any $S$\+module~$P$.
 Finally, one needs to use the fact that the forgetful functor
$\cA\Qcoh\rarrow X\Qcoh$ preserves all inductive limits.
\end{proof}

\subsection{Na\"\i ve co-contra correspondence}
\label{naive-co-contra-subsecn}
 We refer to the introduction to the paper~\cite{Pmgm} for
a general discussion of the co-contra correspondence phenomenon.
 This section is an $\cA$\+module generalization
of~\cite[Section~4.8]{Pcosh}.
 We start with underived versions of the na\"\i ve co-contra
correspondence for $\cA$\+modules before proceeding to deduce
the derived versions.

 Let $X$ be a quasi-compact semi-separated scheme.
 The notation $\cA\Ctrh_\al^{\cA\dlct}\subset\cA\Ctrh_\al^{X\dlct}
\subset\cA\Ctrh_\al$ for exact categories of antilocal contraherent
$\cA$\+modules was introduced in
Section~\ref{antilocal-contrah-A-lct-A-modules-subsecn},
and the similar notation $X\Ctrh_\al$ for the exact category of
antilocal contraherent cosheaves was also mentioned there.
 We refer to
Sections~\ref{locally-contraherent-cosheaves-subsecn}\+-%
\ref{exact-categories-of-contrah-subsecn} for the notation
$X\Ctrh^\lct\subset X\Lcth_\bW^\lct$, and put $X\Ctrh^\lct_\al
=X\Ctrh_\al\cap X\Lcth_\bW^\lct$ for any open covering $\bW$
of the scheme~$X$.
 See~\cite[Section~4.3]{Pcosh} for a detailed discussion.

\begin{lem} \label{quasi-algebra-underived-naive-co-contra}
 Let $X$ be a quasi-compact semi-separated scheme and $\cA$ be
a quasi-coherent quasi-algebra over~$X$.
 In this setting: \par
\textup{(a)} The functors $\cA\ocn_\cA{-}$ and\/ $\fHom_\cA(\cA,{-})$
restrict to mutually inverse equivalences between the exact
categories $\cA\Qcoh^{X\dcta}$ and $\cA\Ctrh_\al$, forming a commutative
square diagram with the similar mutually inverse equivalences
$\cO_X\ocn_X{-}$ and\/ $\fHom_X(\cO_X,{-})$ between the exact
categories $X\Qcoh^\cta$ and $X\Ctrh_\al$,
\begin{equation} \label{X-cta-underived-naive-co-contra}
\begin{gathered}
 \qquad\xymatrix{
  \text{\llap{$\fHom_\cA(\cA,{-})\:$}} \cA\Qcoh^{X\dcta}
  \ar@{=}[rr] \ar[d]
  && \cA\Ctrh_\al \text{\rlap{$\,\,:\!\cA\ocn_\cA{-}$}} \ar[d] \\
  \text{\llap{$\fHom_X(\cO_X,{-})\:$}} X\Qcoh^\cta \ar@{=}[rr] 
  && X\Ctrh_\al \text{\rlap{$\,\,:\!\cO_X\ocn_X{-}$}}
 }
\end{gathered}
\end{equation} \par
\textup{(b)} The functors $\cA\ocn_\cA{-}$ and\/ $\fHom_\cA(\cA,{-})$
restrict to mutually inverse equivalences between the exact
categories $\cA\Qcoh^{X\dcot}$ and $\cA\Ctrh_\al^{X\dlct}$, forming
a commutative square diagram with the similar mutually inverse
equivalences $\cO_X\ocn_X{-}$ and\/ $\fHom_X(\cO_X,{-})$ between
the exact categories $X\Qcoh^\cot$ and $X\Ctrh_\al^\lct$,
\begin{equation} \label{X-cot-X-lct-underived-naive-co-contra}
\begin{gathered}
 \qquad\xymatrix{
  \text{\llap{$\fHom_\cA(\cA,{-})\:$}} \cA\Qcoh^{X\dcot}
  \ar@{=}[rr] \ar[d]
  && \cA\Ctrh_\al^{X\dlct} \text{\rlap{$\,\,:\!\cA\ocn_\cA{-}$}}
  \ar[d] \\
  \text{\llap{$\fHom_X(\cO_X,{-})\:$}} X\Qcoh^\cot \ar@{=}[rr] 
  && X\Ctrh_\al^\lct \text{\rlap{$\,\,:\!\cO_X\ocn_X{-}$}}
 }
\end{gathered}
\end{equation} \par
\textup{(c)} The functors $\cA\ocn_\cA{-}$ and\/ $\fHom_\cA(\cA,{-})$
restrict to mutually inverse equivalences between the exact
categories $\cA\Qcoh^{\cA\dcot}$ and $\cA\Ctrh_\al^{\cA\dlct}$, forming
a commutative square diagram with the similar mutually inverse
equivalences $\cO_X\ocn_X{-}$ and\/ $\fHom_X(\cO_X,{-})$ between
the exact categories $X\Qcoh^\cot$ and $X\Ctrh_\al^\lct$,
\begin{equation} \label{A-cot-A-lct-underived-naive-co-contra}
\begin{gathered}
 \qquad\xymatrix{
  \text{\llap{$\fHom_\cA(\cA,{-})\:$}} \cA\Qcoh^{\cA\dcot}
  \ar@{=}[rr] \ar[d]
  && \cA\Ctrh_\al^{\cA\dlct} \text{\rlap{$\,\,:\!\cA\ocn_\cA{-}$}}
  \ar[d] \\
  \text{\llap{$\fHom_X(\cO_X,{-})\:$}} X\Qcoh^\cot \ar@{=}[rr] 
  && X\Ctrh_\al^\lct \text{\rlap{$\,\,:\!\cO_X\ocn_X{-}$}}
 }
\end{gathered}
\end{equation}
\end{lem}

\begin{proof}
 First of all, the lower horizontal equivalences of exact categories
are established in~\cite[Lemma~4.8.2]{Pcosh} for part~(a)
and~\cite[Lemma~4.8.4(a)]{Pcosh} for parts~(b\+-c).
 Therefore, commutativity of the square
diagram~\eqref{X-cta-fHom-forgetful-diagram} in
Lemma~\ref{fHom-forgetting-A-module-structure-lemma}(a) implies
that the functor $\fHom_\cA(\cA,{-})\:\cA\Qcoh^{X\dcta}\rarrow\cA\Ctrh$
takes values in the full exact subcategory
$\cA\Ctrh_\al\subset\cA\Ctrh$ (alternatively, one can just refer
to Corollary~\ref{fHom-antilocal-corollary}).
 Similarly, commutativity of the square
diagram~\eqref{quasi-algebra-contratensor-forgetful-diagram} in
Lemma~\ref{contratensor-forgetting-A-module-structure-lemma} implies
that the functor $\cA\ocn_\cA{-}\,\:\cA\Ctrh_\al\rarrow\cA\Qcoh$
takes values in the full exact subcategory
$\cA\Qcoh^{X\dcta}\subset\cA\Qcoh$, while the functor
$\cA\ocn_\cA{-}\,\:\cA\Ctrh_\al^{X\dlct}\rarrow\cA\Qcoh$
takes values in the full exact subcategory
$\cA\Qcoh^{X\dcot}\subset\cA\Qcoh$.
 This establishes the existence of all the functors on the diagrams
in parts~(a) and~(b).
 The upper horizontal functors on these diagrams are exact by
formulas~(\ref{fHom-A-very-flaproj-X-cta-exact}\+-%
\ref{fHom-A-rob-flaproj-X-cot-X-lct-exact})
and Lemma~\ref{contratensor-exact-on-antilocal}.

 The upper horizontal functors on the diagrams in~(a) and~(b) are
adjoint to each other by formula~\eqref{fHom-contratensor-adjunction}.
 Moreover, this adjunction agrees with the similar adjunction between
the two lower horizontal functors on the respective
diagrams~\cite[formula~(2.19) in Section~2.6]{Pcosh}, in
the sense that the vertical (forgetful) functors take the adunction
morphisms to the adjunction morphisms.
 As the vertical functors are conservative (take nonisomorphisms to
nonisomorphisms), the fact that the lower horizontal adjunction
is a category equivalence implies that the upper horizontal adjunction
is a ctategory equivalence, too, in both parts~(a) and~(b).

 Alternatively, the upper horizontal equivalences in parts~(a) and~(b)
can be proved directly and independently of the results
of~\cite[Section~4.8]{Pcosh} by arguments similar to the proof of
part~(c) below.
 Then, for part~(a), one needs to use
Theorems~\ref{qcomp-qsep-very-flaproj-complete-cotorsion-pair-thm}(c)
and~\ref{qcomp-qsep-antilocal-complete-cotorsion-pair-thm}(c).
 For part~(b), one needs to use
Theorems~\ref{qcomp-qsep-flaproj-complete-cotorsion-pair-thm}(c)
and~\ref{qcomp-qsep-X-lct-aloc-complete-cotorsion-pair-thm}(c).

 Part~(c): the functors
$\fHom_\cA(\cA,{-})\:\cA\Qcoh^{\cA\dcot}\rarrow\cA\Ctrh^{\cA\dlct}$ and
$\cA\ocn_\cA\nobreak{-}\,\:\allowbreak\cA\Ctrh_\al^{\cA\dlct}\rarrow
\cA\Qcoh$ are exact by formula~\eqref{fHom-A-flat-A-cot-B-lct-exact}
and Lemma~\ref{contratensor-exact-on-antilocal}.
 In order to establish their properties mentioned in the lemma, one
starts with considering the case of an affine scheme $U$ with
a quasi-coherent quasi-algebra $\cB$ over~$U$.

 In this case, one has $\cB\Qcoh\simeq\cB(U)\Modl$ as
$\cB\Ctrh^{\cB\dlct}\simeq\cB(U)\Modl^\cot$, as explained in
Sections~\ref{cosheaves-of-A-modules-subsecn}\+-%
\ref{A-loc-cotors-loc-inj-cosheaves-subsecn}.
 The equivalence of categories $\cB\Qcoh\simeq\cB(U)\Modl$ identifies
the full subcategory $\cB\Qcoh^{\cB\dcot}\subset\cB\Qcoh$ with
the full subcategory $\cB(U)\Modl^\cot\subset\cB(U)\Modl$.
 Furthermore, all contraherent cosheaves on an affine scheme $U$
are antilocal.
 With these identifications in mind, the functor
$\fHom_\cB(\cB,{-})\:\cB\Qcoh^{\cB\dcot}\rarrow\cB\Ctrh^{\cB\dlct}$ is
just the identity functor $\cB(U)\Modl^\cot\allowbreak
\rarrow\cB(U)\Modl^\cot$.
 The functor $\cB\ocn_\cB{-}\,\:\cB\Ctrh_\al^{\cB\dlct}\rarrow\cB\Qcoh$
is just the identity inclusion $\cB(U)\Modl^\cot\rarrow\cB(U)\Modl$.

 Returning to our scheme $X$, the result of
Theorem~\ref{qcomp-qsep-A-flat-complete-cotorsion-pair-thm}(c) describes
the objects of $\cA\Qcoh^{\cA\dcot}$ as the direct summands of finitely
iterated extensions of the direct images of objects of
$\cB_\alpha\Qcoh^{\cB_\alpha\dcot}$ from affine open subschemes
$U_\alpha\subset X$ with the quasi-coherent quasi-algebras
$\cB_\alpha=\cA|_{U_\alpha}$.
 In the same notation, the result of
Theorem~\ref{qcomp-qsep-X-lct-aloc-complete-cotorsion-pair-thm}(c)
and Remark~\ref{antilocal-contrah-A-lct-A-modules-remark} describes
the objects of $\cA\Ctrh^{\cA\dlct}_\al$ as the direct summands of
finitely iterated exensions of the direct images of objects of
$\cB_\alpha\Ctrh^{\cB_\alpha\dlct}$.

 It remains to point out that the functors $\fHom_\cA(\cA,{-})$
and $\cA\ocn_\cA{-}$ form commutative square diagrams with the direct
images from affine open subschemes, as per
the formulas~\eqref{fHom-projection-formula}
and~\eqref{contratensor-projection-formula}.
 It follows from these observations that the functor
$\fHom_\cA(\cA,{-})\:\cA\Qcoh^{\cA\dcot}\rarrow\cA\Ctrh^{\cA\dlct}$
takes values in $\cA\Ctrh^{\cA\dlct}_\al$, the functor
$\cA\ocn_\cA{-}\,\:\cA\Ctrh_\al^{\cA\dlct}\rarrow\cA\Qcoh$ takes values
in $\cA\Qcoh^{\cA\dcot}$, and the resulting pair of adjoint functors
(shown by the upper horizontal double line on
the diagram~\eqref{A-cot-A-lct-underived-naive-co-contra} in part~(c))
are mutually inverse equivalences.
 (Cf.~\cite[proof of Lemma~4.8.2]{Pcosh}.)
\end{proof}

\begin{lem} \label{underived-naive-co-contra-direct-image-lemma}
 Let $X$ be a quasi-compact semi-separated scheme and $Y\subset X$ be
a quasi-compact open subscheme with the open immersion morphism
$f\:Y\rarrow X$.
 Let $\cA$ be a quasi-coherent quasi-algebra over~$X$.
 In this setting: \par
\textup{(a)} The equivalences of exact categories from
Lemma~\ref{quasi-algebra-underived-naive-co-contra}(a) form
a commutative square diagram with the direct image functors $f_*$
and~$f_!$,
\begin{equation} \label{X-cta-under-naive-co-contra-direct-image}
\begin{gathered}
 \qquad\quad\xymatrix{
  \text{\llap{$\fHom_{\cA|_Y}(\cA|_Y,{-})\:$}}\cA|_Y\Qcoh^{Y\dcta}
  \ar@{=}[rr] \ar[d]_{f_*}
  && \cA|_Y\Ctrh_\al \text{\rlap{$\,\,:\!\cA|_Y\ocn_{\cA|_Y}{-}$}} 
  \ar[d]^{f_!} \\
  \text{\llap{$\fHom_\cA(\cA,{-})\:$}} \cA\Qcoh^{X\dcta} \ar@{=}[rr]
  && \cA\Ctrh_\al \text{\rlap{$\,\,:\!\cA\ocn_\cA{-}$}}
 }
\end{gathered}
\end{equation} \par
\textup{(b)} The equivalences of exact categories from
Lemma~\ref{quasi-algebra-underived-naive-co-contra}(b) form
a commutative square diagram with the direct image functors $f_*$
and~$f_!$,
\begin{equation} \label{X-cot-under-naive-co-contra-direct-image}
\begin{gathered}
 \qquad\quad\xymatrix{
  \text{\llap{$\fHom_{\cA|_Y}(\cA|_Y,{-})\:$}}\cA|_Y\Qcoh^{Y\dcot}
  \ar@{=}[rr] \ar[d]_{f_*}
  && \cA|_Y\Ctrh_\al^{Y\dlct}
  \text{\rlap{$\,\,:\!\cA|_Y\ocn_{\cA|_Y}{-}$}} \ar[d]^{f_!} \\
  \text{\llap{$\fHom_\cA(\cA,{-})\:$}}
  \cA\Qcoh^{X\dcot} \ar@{=}[rr] 
  && \cA\Ctrh_\al^{X\dlct}
  \text{\rlap{$\,\,:\!\cA\ocn_\cA{-}$}}
 }
\end{gathered}
\end{equation} \par
\textup{(c)} The equivalences of exact categories from
Lemma~\ref{quasi-algebra-underived-naive-co-contra}(c) form
a commutative square diagram with the direct image functors $f_*$
and~$f_!$,
\begin{equation} \label{A-cot-under-naive-co-contra-direct-image}
\begin{gathered}
 \qquad\quad\xymatrix{
  \text{\llap{$\fHom_{\cA|_Y}(\cA|_Y,{-})\:$}}\cA|_Y\Qcoh^{\cA|_Y\dcot}
  \ar@{=}[rr] \ar[d]_{f_*}
  && \cA|_Y\Ctrh_\al^{\cA|_Y\dlct}
  \text{\rlap{$\,\,:\!\cA|_Y\ocn_{\cA|_Y}{-}$}} \ar[d]^{f_!} \\
  \text{\llap{$\fHom_\cA(\cA,{-})\:$}} \cA\Qcoh^{\cA\dcot} \ar@{=}[rr] 
  && \cA\Ctrh_\al^{\cA\dlct}
  \text{\rlap{$\,\,:\!\cA\ocn_\cA{-}$}}
 }
\end{gathered}
\end{equation} \par
 The direct image functors $f_*$ and~$f_!$ on all
the diagrams~\textup{(\ref{X-cta-under-naive-co-contra-direct-image}\+-%
\ref{A-cot-under-naive-co-contra-direct-image})} are exact.
\end{lem}

\begin{proof}
 Part~(a): the functor $f_*\:\cA|_Y\Qcoh\rarrow\cA\Qcoh$ takes
$Y$\+contraadjusted quasi-coherent $\cA|_Y$\+modules on $Y$ to
$X$\+contraadjusted quasi-coherent $\cA$\+modules on $X$
by~\cite[Section~2.5]{Pcosh}, \cite[Lemma~1.7(c)]{Pal}, as mentioned
in Section~\ref{antilocality-of-X-contraadjusted-subsecn}.
 The functor $f_*\:Y\Qcoh^\cta\rarrow X\Qcoh^\cta$ is exact
by~\cite[Corollary~4.1.14(a)]{Pcosh}; hence the functor
$f_*\:\cA|_Y\Qcoh^{Y\dcta}\rarrow\cA\Qcoh^{X\dcta}$ is exact, too.

 The functor $f_!\:(Y,\cO_Y)\Cosh\rarrow(X,\cO_X)\Cosh$ takes
antilocal contraherent cosheaves on $Y$ to antilocal contraherent
cosheaves, and the functor $f_!\:Y\Ctrh_\al\rarrow X\Ctrh_\al$ is
exact by~\cite[Corollary~4.6.3(a)]{Pcosh}.
 Therefore, the functor $f_!\:(Y,\cA|_Y)\Cosh\allowbreak\rarrow
(X,\cA)\Cosh$ takes antilocal contraherent $\cA|_Y$\+modules to
antilocal contraherent $\cA$\+modules, and the functor
$f_!\:\cA|_Y\Ctrh_\al\rarrow\cA\Ctrh_\al$ is exact (as it was
already mentioned in the proof of
Corollary~\ref{contrah-al-A-lct-direct-image}).

 This proves the existence and exactness of all the functors on
the diagram~\eqref{X-cta-under-naive-co-contra-direct-image}.
 The diagram is commutative by formula~\eqref{fHom-projection-formula},
or alternatively, by
Lemma~\ref{contratensor-projection-comparison-and-antilocal}(b).

 The proof of part~(b) is similar.
 One can refer to Section~\ref{antilocality-of-X-cotorsion-subsecn}
and~\cite[Corollaries~4.1.14(b) and~4.6.3(b)]{Pcosh}.

 Part~(c): the functor $f_*\:\cA|_Y\Qcoh\rarrow\cA\Qcoh$ takes
$\cA|_Y$\+cotorsion quasi-coherent $\cA|_Y$\+modules on $Y$ to
$\cA$\+cotorsion quasi-coherent $\cA$\+modules on $X$
by~\cite[Lemma~1.7(b)]{Pal}, as explained in
Section~\ref{antilocality-of-A-cotorsion-subsecn}.
 It is clear from part~(a) or~(b) that the functor
$f_*\:\cA|_Y\Qcoh^{\cA|_Y\dcot}\rarrow\cA\Qcoh^{\cA\dcot}$ is exact
as a restriction of an exact functor to full exact subcategories of
the source and target categories.

 The functor $f_!\:(Y,\cA|_Y)\Cosh\allowbreak\rarrow(X,\cA)\Cosh$ takes
antilocal $\cA|_Y$\+locally cotorsion contraherent $\cA|_Y$\+modules to
antilocal $\cA|_Y$\+locally cotorsion contraherent $\cA$\+modules, and
the functor $f_!\:\cA|_Y\Ctrh^{\cA|_Y\dlct}_\al\rarrow
\cA\Ctrh^{\cA\dlct}_\al$ is exact by
Corollary~\ref{contrah-al-A-lct-direct-image}.

 This proves the existence and exactness of all the functors on
the diagram~\eqref{A-cot-under-naive-co-contra-direct-image}.
 The diagram is commutative by formula~\eqref{fHom-projection-formula},
or alternatively, by
Lemma~\ref{contratensor-projection-comparison-and-antilocal}(b).
\end{proof}

 The following two theorems form the main result of this
Section~\ref{naive-co-contra-secn}.

\begin{thm} \label{quasi-algebra-bounded-derived-naive-co-contra}
 Let $X$ be a quasi-compact semi-separated scheme with an open
covering\/ $\bW$ and $\cA$ be a quasi-coherent quasi-algebra over~$X$.
 Then, for any conventional derived category symbol\/ $\st=\bb$, $+$,
$-$, or\/~$\varnothing$, there are natural triangulated equivalences
between the derived categories of quasi-coherent and\/ $\bW$\+locally
contraherent $\cA$\+modules on~$X$,
\begin{equation} \label{bounded-derived-naive-co-contra-diagram-I}
 \xymatrix{
  \sD^\st(\cA\Qcoh)\ar@{-}@<2pt>[r]
  & \sD^\st(\cA\Qcoh^{X\dcta}) \ar@<2pt>[l] \ar@{=}[r]
  & \sD^\st(\cA\Ctrh_\al) \ar@<-2pt>[r]
  & \sD^\st(\cA\Lcth_\bW). \ar@{-}@<-2pt>[l]
 }
\end{equation}
\end{thm}

\begin{proof}
 The leftmost equivalence, induced by the exact inclusion of
exact/abelian categories $\cA\Qcoh^{X\dcta}\rarrow\cA\Qcoh$, is
the result of Corollary~\ref{qcoh-X-cta-derived-equivalence}.
 The rightmost equivalence, induced by the exact inclusion of
exact categories $\cA\Ctrh_\al\rarrow\cA\Lcth_\bW$, is the result of
Corollary~\ref{A-lcth-ctrh-al-derived-equivalences}(a).
 The middle equivalence is induced by the equivalence of exact
categories $\cA\Qcoh^{X\dcta}\simeq\cA\Ctrh_\al$ provided by
Lemma~\ref{quasi-algebra-underived-naive-co-contra}(a).
\end{proof}

 The notation for the resulting triangulated equivalence of
Theorem~\ref{quasi-algebra-bounded-derived-naive-co-contra} is
\begin{equation} \label{bounded-derived-naive-co-contra-diagram-II}
 \xymatrix{
  \boR\fHom_\cA(\cA,{-})\:\sD^\st(\cA\Qcoh)\ar@{=}[r]
  & \sD^\st(\cA\Lcth_\bW)\,:\!\cA\ocn_\cA^\boL{-},
 }
\end{equation}
where $\boR\fHom_\cA(\cA,{-})$ means the right derived functor of
the contraherent\/ $\fHom$ functor $\fHom_\cA(\cA,{-})$, while
$\cA\ocn_\cA^\boL{-}$ stands for the left derived functor of
the contratensor product functor $\cA\ocn_\cA{-}$.

 The second theorem is a generalization
of~\cite[Corollary~4.8.7]{Pcosh}.

\begin{thm} \label{quasi-algebra-unbounded-derived-naive-co-contra}
 Let $X$ be a quasi-compact semi-separated scheme with an open
covering\/ $\bW$ and $\cA$ be a quasi-coherent quasi-algebra over~$X$.
 Then, for any derived category symbol\/ $\st=+$ or\/~$\varnothing$,
there is a commutative diagram of natural triangulated equivalences
between the derived categories of various abelian/exact categories of
quasi-coherent and\/ $\bW$\+locally contraherent $\cA$\+modules on~$X$,
\begin{equation} \label{unbounded-derived-naive-co-contra-diagram}
\begin{gathered}
 \xymatrix{
  \sD^\st(\cA\Qcoh)\ar@{-}@<2pt>[r]
  & \sD^\st(\cA\Qcoh^{X\dcta}) \ar@<2pt>[l] \ar@{=}[r] \ar@{-}@<2pt>[d]
  & \sD^\st(\cA\Ctrh_\al) \ar@<-2pt>[r] \ar@{-}@<-2pt>[d]
  & \sD^\st(\cA\Lcth_\bW) \ar@{-}@<-2pt>[l] \ar@{-}@<-2pt>[d] \\
  & \sD^\st(\cA\Qcoh^{X\dcot}) \ar@{=}[r]
  \ar@<2pt>[u] \ar@{-}@<2pt>[d]
  & \sD^\st(\cA\Ctrh_\al^{X\dlct}) \ar@<-2pt>[r]
  \ar@<-2pt>[u] \ar@{-}@<-2pt>[d]
  & \sD^\st(\cA\Lcth_\bW^{X\dlct}) \ar@{-}@<-2pt>[l]
  \ar@<-2pt>[u] \ar@{-}@<-2pt>[d] \\
  & \sD^\st(\cA\Qcoh^{\cA\dcot}) \ar@{=}[r] \ar@<2pt>[u]
  & \sD^\st(\cA\Ctrh_\al^{\cA\dlct}) \ar@<-2pt>[r] \ar@<-2pt>[u]
  & \sD^\st(\cA\Lcth_\bW^{\cA\dlct}) \ar@{-}@<-2pt>[l] \ar@<-2pt>[u]
 }
\end{gathered}
\end{equation}
\end{thm}

\begin{proof}
 The three triangulated equivalences shown by double horizontal lines
without arrows (in the middle column of equivalences) are induced by
the equivalences of exact categories provided by
Lemma~\ref{quasi-algebra-underived-naive-co-contra}(a\+-c).
 All the other triangulated functors (shown by arrows) are induced by
the respective exact inclusions of exact/abelian categories.

 The upper horizontal line of triangulated equivalences is a particular
case of Theorem~\ref{quasi-algebra-bounded-derived-naive-co-contra}.
 The triangulated equivalences in the left-hand part of the diagram
are the result of Corollary~\ref{qcoh-X-cta-derived-equivalence}
and Theorem~\ref{qcoh-X-cot-A-cot-derived-equivalence}.
 The horizontal triangulated equivalences in the right-hand part of
the diagram are the result of
Corollary~\ref{A-lcth-ctrh-al-derived-equivalences}(a\+-c).
 The vertical triangulated equivalences in the rightmost column of
the diagram are the result of
Theorem~\ref{A-lcth-X-lct-A-lct-derived-equivalence}.
 The vertical triangulated equivalences in the next-to-rightmost
column of the diagram follow.
\end{proof}

 For a partial extension of
Theorem~\ref{quasi-algebra-unbounded-derived-naive-co-contra}
to the derived category symbols $\st=\bb$ and~$-$ involving
the $X$\+cotorsion quasi-coherent $\cA$\+modules and
the $X$\+locally cotorsion locally contraherent $\cA$\+modules
(in addition to the $X$\+contraadjusted and $X$\+locally contraadjusted
ones, as in
Theorem~\ref{quasi-algebra-bounded-derived-naive-co-contra}),
which holds under more restrictive assumptions on the scheme~$X$, see
Theorem~\ref{quasi-algebra-cot-bounded-derived-naive-co-contra} below.

\begin{rem} \label{naive-co-contra-applies-to-diffoperators-remark}
 In particular, the results of
Theorems~\ref{quasi-algebra-bounded-derived-naive-co-contra}\+-%
\ref{quasi-algebra-unbounded-derived-naive-co-contra} are applicable
to the sheaf of rings of (fiberwise) differential operators
$\cA=\cD_{X/T}(\cU,\cU)$ for any morphism of schemes $X\rarrow T$
locally of finite presentation (with a quasi-compact semi-separated
scheme~$X$) and any locally finitely presented quasi-coherent
sheaf $\cU$ on $X$, as per
Example~\ref{fiberwise-differential-operators-examples}(4).
 The same applies to the twisted universal enveloping quasi-algebra
$\cA=\cA_X(\g,\widetilde\g)$ of any quasi-coherent twisted Lie algebroid
$(\g,\widetilde\g)$ over $X$ and, in particular, to the sheaf of rings
of crystalline differential operators $\cA=\cD^\cry_{X/T}$ for a weakly
smooth morphism of schemes $X\rarrow T$, as per
Examples~\ref{fiberwise-differential-operators-examples}(5\+-7).
\end{rem}

\Section{Contraherent CDG-Modules}
\label{contraherent-cdg-modules-secn}

 The concepts of quasi-coherent CDG\+quasi-algebras over schemes
and quasi-coherent CDG\+modules over them were defined
in~\cite[Section~2.4]{Pedg} (in the context of the definition of
a quasi-module from~\cite[Section~2.3]{Pedg}, which is more narrow than
in our Section~\ref{prelim-quasi-modules-subsecn}).
 In this section, we define locally contraherent CDG\+modules over
quasi-coherent CDG\+quasi-algebras.

\subsection{Quasi-coherent and contraherent graded modules}
\label{graded-modules-subsecn}
 It is advisable to think of graded objects in a category as collections
of objects indexed by the grading set or the grading group, such as
the group of integers~$\boZ$, \emph{rather than} as direct sum
decompositions indexed by~$\boZ$.
 Even though the direct sum notation is convenient and we will use it,
speaking, e.~g., of ``graded rings $B^*=\bigoplus_{i\in\boZ} B^i$\,'',
it is nevertheless more instructive to think of a graded ring $B^*$ as
the collection of abelian groups $B^i$ together with multiplication maps
$B^i\times B^j\rarrow B^{i+j}$, \ $i$, $j\in\boZ$, and a unit element
$1\in B^0$, than to consider $B^*$ as a single abelian group/ring
decomposed into a direct sum of the grading components~$B^i$.

 One way to explain this point-of-view preference is to mention that
the infinite direct sum is a nontrivial functor which need not commute
with other functors.
 An even more important explanation, particularly in the contexts where
one is interested in the ``contra'' side of one's story (involving,
e.~g., contramodules or contraherent cosheaves), is that one often needs
to consider \emph{direct products} of the grading components instead
of the direct sums.

 So, the proper point of view is that the category of graded abelian
groups $\Ab^\boZ$ is just the Cartesian product of $\boZ$ copies of
the category of abelian groups~$\Ab$.
 There is not one but \emph{two} forgetful functors $\Ab^\boZ\rarrow\Ab$
assigning an ungraded abelian group to a graded one: the direct sum
functor $\Sigma\:\Ab^\boZ\rarrow\Ab$ and the direct product functor
$\Pi\:\Ab^\boZ\rarrow\Ab$.
 So we put~\cite[Section~2.1]{Prel}
$$
 \Sigma((A^i)_{i\in\boZ})=\bigoplus\nolimits_{i\in\boZ}A^i
 \quad\text{and}\quad
 \Pi((A^i)_{i\in\boZ})=\prod\nolimits_{i\in\boZ}A^i.
$$

 For any $\boZ$\+graded ring $B^*$, the ungraded abelian group
$\Sigma B^*$ has a natural ungraded ring structure.
 For any $\boZ$\+graded $B^*$\+module $M^*$, \emph{both} the ungraded
abelian groups $\Sigma M^*$ and $\Pi M^*$ have natural
$\Sigma B^*$\+module structures.
 The tensor product $N^*\ot_{B^*}M^*$ of a graded right $B^*$\+module
$N^*$ and a graded left $B^*$\+module $M^*$ carries a natural grading,
and there is the graded abelian group $\Hom_{B^*}^*(M^*,P^*)$ for any
two graded left $B^*$\+modules $M^*$ and~$P^*$.
 The tensor product and $\Hom$ functors interact with the functors of
forgetting the grading according to the rules
$$
 \Sigma(N^*\ot_{B^*}M^*)\simeq\Sigma N^*\ot_{\Sigma B^*}\Sigma M^*
 \ \ \text{and}\ \
 \Pi\Hom^*_{B^*}(M^*,P^*)\simeq\Hom_{\Sigma B^*}(\Sigma M^*,\Pi P^*).
$$

 The abelian categories of graded left and right modules over a graded
ring $B^*$ will be denoted by $B^*\Modl$ and $\Modr B^*$.
 So we have two exact, faithful functors of forgetting the grading
$\Sigma\:B^*\Modl\rarrow\Sigma B^*\Modl$ and $\Pi\:B^*\Modl\rarrow
\Sigma B^*\Modl$.

 A graded left $B^*$\+module $C^*$ is said to be \emph{cotorsion} if
$\Ext^1_{B^*}(F^*,C^*)=0$ for all flat graded left $B^*$\+modules~$F^*$.
 Here $\Ext^n_{B^*}({-},{-})$, \,$n\ge0$, denote the $\Ext$ groups
computed in the category of graded left $B^*$\+modules $B^*\Modl$.
 The notion of a cotorsion graded $B^*$\+module has all the properties
similar to those of its ungraded version, as spelled out in
Section~\ref{prelim-cotorsion-subsecn}.
 A graded left $B^*$\+module $J^*$ is said to be \emph{injective} if
it is injective as an object of $B^*\Modl$.

 Given a commutative ring $R$, a graded $R$\+$R$\+bimodule $B^*$ is said
to be a \emph{graded quasi-module} over $R$ if its grading component
bimodules $B^i$, \,$i\in\boZ$, are quasi-modules over~$R$ (in the sense
of the definition in Section~\ref{prelim-quasi-modules-subsecn}).
 A graded ring $A^*$ endowed with a morphism of graded rings
$R\rarrow A^*$ (i.~e., a ring homomorphism $R\rarrow A^0$) is said to be
a \emph{graded quasi-algebra} over $R$ if the underlying graded
$R$\+$R$\+bimodule of $A^*$ is a graded quasi-module over~$R$
(cf.\ the discussion in Section~\ref{prelim-quasi-algebras-subsecn}).
 Here the underlying graded $R$\+$R$\+bimodule structure on $A^*$ is
induced by the graded ring structure on $A^*$ and the graded ring
homomorphism $R\rarrow A^*$.

 Essentially all the basic results of commutative and noncommutative
ring and module theory remain true for graded modules over graded rings.
 One just needs to restrict all the arguments and constructions to
homogeneous elements, etc.
 This allows us to extend the theory developed in
Sections~\ref{contraherent-cosheaves-of-A-modules-secn}\+-%
\ref{naive-co-contra-secn} to quasi-coherent and (locally) contraherent
graded modules over quasi-coherent graded quasi-algebras over schemes
without spelling out all the details anew.
 One only needs to be careful with the functors of forgetting
the grading.

 From our perspective, we prefer to think of (co)sheaves of graded
modules rather than graded (co)sheaves of modules.
 So the graded generalization of the theory developed above in
this paper starts with generalizing
Section~\ref{cosheaves-of-modules-subsecn} to the case of a topological
space $X$ endowed with a sheaf of graded rings~$\cO^*$.
 The definitions of sheaves and cosheaves of $\cO$\+modules
in Section~\ref{cosheaves-of-modules-subsecn} are extended
straightforwardly to the case of sheaves and cosheaves of graded modules
over a sheaf of graded rings~$\cO^*$.
 We denote the category of sheaves of $\cO^*$\+modules on $X$ by
$(X,\cO^*)\Sh$ and the category of cosheaves of $\cO^*$\+modules on $X$
by $(X,\cO^*)\Cosh$.

 Then, given a scheme $X$, a sheaf of graded $\cO_X$\+$\cO_X$\+bimodules
$\cB^*$ is said to be a \emph{quasi-coherent graded quasi-module} over
$X$ if its grading component sheaves of $\cO_X$\+$\cO_X$\+bimodules
$\cB^i$, \,$i\in\boZ$, are quasi-coherent quasi-modules over~$X$
(in the sense of the definition in
Section~\ref{quasi-coherent-quasi-algebras-subsecn}).
 A sheaf of graded rings $\cA^*$ on $X$ endowed with a morphism of
sheaves of graded rings $\cO_X\rarrow\cA^*$ is said to be
a \emph{quasi-coherent graded quasi-algebra} over $X$ if the underlying
sheaf of graded $\cO_X$\+$\cO_X$\+bimodules of $\cA^*$ is
a quasi-coherent graded quasi-module.
 Here $\cO_X$ is considered as a sheaf of graded rings concentrated
in the grading~$0$.

 Let $X$ be a scheme and $\cA^*$ be a sheaf of graded rings on $X$
endowed with a morphism of sheaves of graded rings $\cO_X\rarrow\cA^*$.
 A sheaf of graded $\cA^*$\+modules $\M^*$ on $X$ is said to be
a \emph{quasi-coherent graded $\cA^*$\+module} if the grading component
sheaves $\M^i$ of $\M^*$ are quasi-coherent as sheaves of
$\cO_X$\+modules.
 Given an open covering $\bW$ of $X$, a cosheaf of graded
$\cA^*$\+modules $\P^*$ on $X$ is said to be a \emph{$\bW$\+locally
contraherent $\cA^*$\+module} if the grading component cosheaves $\P^i$
of $\P^*$ are $\bW$\+locally contraherent as cosheaves of
$\cO_X$\+modules.

 The functor of forgetting the grading assigns to a quasi-coherent
graded algebra $\cA^*$ its underlying ungraded quasi-coherent
algebra $\Sigma\cA^*=\bigoplus_{i\in\boZ}\cA^i$.
 The similar functor of forgetting the grading assigns to
a quasi-coherent graded $\cA^*$\+module $\M^*$ its underlying
ungraded quasi-coherent $\Sigma\cA^*$\+module $\Sigma\M^*=
\bigoplus_{i\in\boZ}\M^i$.
 To a $\bW$\+locally contraherent graded $\cA^*$\+module $\P^*$,
the functor of forgetting the grading assigns its underlying
ungraded $\bW$\+locally contraherent $\Sigma\cA^*$\+module
$\Pi\P^*=\prod_{i\in\boZ}\P^i$.

 We denote the category of quasi-coherent graded left $\cA^*$\+modules
on $X$ by $\cA^*\Qcoh$, the category of quasi-coherent graded right
$\cA^*$\+modules by $\Qcohr\cA^*$, and the category of $\bW$\+locally 
contraherent graded (left) $\cA^*$\+modules on $X$ by $\cA^*\Lcth_\bW$.
 So the functors of forgetting the grading are
$\Sigma\:\cA^*\Qcoh\rarrow\Sigma\cA^*\Qcoh$
and $\Pi\:\cA^*\Lcth_\bW\rarrow\Sigma\cA^*\Lcth_\bW$.
 The category $\cA^*\Qcoh$ is a Grothendieck abelian category, while
$\cA^*\Lcth_\bW$ is a naturally an exact category with exact functors
of infinite products.
 The exact category structure on $\cA^*\Lcth_\bW$ is defined similarly
to the ungraded version of this definition given in
Section~\ref{cosheaves-of-A-modules-subsecn}.

 The exact category of contraherent graded $\cA^*$\+modules $\cA^*\Ctrh$
on $X$ is $\cA^*\Ctrh=\cA^*\Lcth_{\{X\}}$.
 The exact category of locally contraherent graded $\cA^*$\+modules
$\cA^*\Lcth$ on $X$ is defined as $\cA^*\Lcth=\bigcup_\bW
\cA^*\Lcth_\bW$.
 So, for a cosheaf of $\cA^*$\+modules $\P^*$ to be locally
contraherent, there must exist an open covering $\bW$ of $X$ such that
the cosheaf of $\cO_X$\+modules $\P^i$ is $\bW$\+locally contraherent
for all $i\in\boZ$ (one open covering $\bW$ has to be suitable for all
grading components~$\P^i$).

 A $\bW$\+locally contraherent graded $\cA^*$\+module $\P^*$ on $X$
is called \emph{$X$\+locally cotorsion} if all its grading component
$\bW$\+locally contraherent cosheaves $\P^i$, \,$i\in\boZ$, are locally
cotorsion.
 The \emph{$X$\+locally injective} $\bW$\+locally contraherent
graded $\cA^*$\+modules on $X$ are defined similarly.
 The exact categories of $X$\+locally cotorsion (locally) contraherent
graded $\cA^*$\+modules are defined similarly to
Section~\ref{cosheaves-of-A-modules-subsecn} and denoted by
$\cA^*\Lcth_\bW^{X\dlct}$, \ $\cA^*\Ctrh^{X\dlct}=
\cA^*\Lcth_{\{X\}}^{X\dlct}$, and $\cA^*\Lcth^{X\dlct}=
\bigcup_\bW\cA^*\Lcth_\bW^{X\dlct}$.

 A quasi-coherent graded $\cA^*$\+module $\M^*$ is called \emph{flat}
(or \emph{$\cA^*$\+flat}) if the graded module $\M^*(U)$ over
the graded ring $\cA^*(U)$ is flat for every affine open subscheme
$U\subset X$.
 This definition has all the properties similar to those of its
ungraded version, as per
Section~\ref{A-loc-cotors-loc-inj-cosheaves-subsecn}.
 The full exact subcategory of flat quasi-coherent graded
$\cA^*$\+modules is denoted by $\cA^*\Qcoh^{\cA^*\dfl}\subset
\cA^*\Qcoh$.

 A $\bW$\+locally contraherent graded $\cA^*$\+module $\P^*$ is
called \emph{$\cA^*$\+locally cotorsion} if, for every affine open
subscheme $U\subset X$ subordinate to $\bW$, the graded left
$\cA^*(U)$\+module $\P^*[U]$ is cotorsion.
 Once again, this definition has all the properties similar to those
of its ungraded version, as per
Section~\ref{A-loc-cotors-loc-inj-cosheaves-subsecn}.
 The exact categories of $\cA^*$\+locally cotorsion (locally)
contraherent graded $\cA^*$\+modules are defined similarly to
Section~\ref{A-loc-cotors-loc-inj-cosheaves-subsecn} and denoted by
$\cA^*\Lcth_\bW^{\cA^*\dlct}$, \ $\cA^*\Ctrh^{\cA^*\dlct}=
\cA^*\Lcth_{\{X\}}^{\cA^*\dlct}$, and $\cA^*\Lcth^{\cA^*\dlct}=
\bigcup_\bW\cA^*\Lcth_\bW^{\cA^*\dlct}$.

 A $\bW$\+locally contraherent graded $\cA^*$\+module $\gJ^*$ is
called \emph{$\cA^*$\+locally injective} if, for every affine open
subscheme $U\subset X$ subordinate to $\bW$, the graded left
$\cA^*(U)$\+module $\gJ^*[U]$ is injective.
 As above, the properties of this definition are very similar to those
of its ungraded version discussed in
Section~\ref{A-loc-cotors-loc-inj-cosheaves-subsecn}.
 The exact categories of $\cA^*$\+locally injective (locally)
contraherent graded $\cA^*$\+modules are denoted by
$\cA^*\Lcth_\bW^{\cA^*\dlin}$, \ $\cA^*\Ctrh^{\cA^*\dlin}=
\cA^*\Lcth_{\{X\}}^{\cA^*\dlin}$, and $\cA^*\Lcth^{\cA^*\dlin}=
\bigcup_\bW\cA^*\Lcth_\bW^{\cA^*\dlin}$.

 The direct image and inverse image functors for quasi-coherent
graded $\cA^*$\+modules and locally contraherent graded
$\cA^*$\+modules are constructed similarly to the discussion in
Sections~\ref{direct-images-of-A-co-sheaves-subsecn}\+-%
\ref{inverse-images-of-A-co-sheaves-subsecn}, and denoted in
the same way.
 The contraherent $\fHom$ and contratensor product functors for
quasi-coherent graded (bi)modules and (contraherent) cosheaves of
graded modules are constructed similarly to
Sections~\ref{fHom-over-qcoh-quasi-algebra-subsecn}
and~\ref{contratensor-over-qcoh-quasi-algebra-subsecn}, and denoted
by $(\cE^*,\C^*)\longmapsto\fHom_{\cA^*}^*(\cE^*,\C^*)$ and
$(\cE^*,\P^*)\longmapsto\cE^*\ocn_{\cB^*}\P^*$.

 The results of Sections~\ref{antilocal-classes-secn}
and~\ref{naive-co-contra-secn}, including
Theorems~\ref{quasi-algebra-bounded-derived-naive-co-contra}\+-%
\ref{quasi-algebra-unbounded-derived-naive-co-contra}, remain valid
in the graded setting.
 However, the graded versions of the category equivalences of
Section~\ref{naive-co-contra-secn} connecting the quasi-coherent
and locally contraherent graded $\cA^*$\+modules do \emph{not} form
commutative diagrams with ungraded versions of the same equivalences
and the functors of forgetting the grading.

 It suffices to consider the case of an affine scheme $X=U$ in order
to understand the situation.
 In this case, the na\"\i ve co-contra correspondence reduces
essentially to identity functors or identity inclusion functors
assigning to a graded $\cA^*(U)$\+module the same graded
$\cA^*(U)$\+module, or at worst its coresolution, etc.
 But the functors of forgetting the grading, which are $\Sigma$ on
the quasi-coherent side and $\Pi$ on the contraherent side, still
do \emph{not} agree with each other.

\subsection{Curved DG-rings and curved DG-modules}
\label{cdg-rings-cdg-modules-subsecn}
 The material of this section goes back to~\cite{Pcurv}.
 The most complete exposition can be found in~\cite[Sections~3.2,
4.2, and~6.1]{Prel}; and there are other expositions
in~\cite[Section~3.1]{Pkoszul}, \cite[Section~1.1]{EP}, 
\cite[Sections~2.2, 2.4, and~3.1]{Pedg}, \cite[Section~6.2]{Pksurv},
\cite[Section~1]{PS7}, etc.

 A \emph{curved DG\+ring} (\emph{CDG\+ring}) $B^\cu=(B^*,d,h)$ is
a set of data consisting of
\begin{itemize}
\item a $\boZ$\+graded ring $B^*=\bigoplus_{i\in\boZ}B^i$;
\item an odd derivation~$d$ of degree~$1$ on $B^*$, that is
$d\:B^i\rarrow B^{i+1}$ for all $i\in\boZ$ and
$d(bc)=d(b)c+(-1)^{|b|}bd(c)$ for all homogeneous elements
$b$ and $c\in B^*$ of degrees $|b|$ and~$|c|\in\boZ$;
\item an element $h\in B^2$.
\end{itemize}
 The following two axioms must be satisfied:
\begin{enumerate}
\renewcommand{\theenumi}{\roman{enumi}}
\item $d(d(b))=hb-bh$ for all $b\in B^*$;
\item $d(h)=0$.
\end{enumerate}

 So, the square of the differential~$d$ on $B^*$ is \emph{not} equal
to zero; rather, it is equal to the commutator with the element~$h$.
 The element $h\in B^2$ is called the \emph{curvature element} of
a CDG\+ring~$B^\cu$.
 A \emph{DG\+ring} $A^\bu=(A^*,d)$ can be defined as a CDG\+ring
$(A^*,d,h)$ with a vanishing curvature element, $h=0$.

 Let $A^\cu=(A^*,d_A,h_A)$ and $B^\cu=(B^*,d_B,h_B)$ be two CDG\+rings.
 A \emph{morphism of CDG\+rings} $f=(f,a)\:B^\cu\rarrow A^\cu$ is
a set of data consisting of
\begin{itemize}
\item a homomorphism of $\boZ$\+graded rings $f\:B^*\rarrow A^*$;
\item an element $a\in A^1$.
\end{itemize}
 The following two axioms must be satisfied:
\begin{enumerate}
\renewcommand{\theenumi}{\roman{enumi}}
\setcounter{enumi}{2}
\item $f(d_B(b))=d_A(f(b))+af(b)-(-1)^{|b|}f(b)a$ for all elements
$b\in B^*$ of degree $|b|\in\boZ$;
\item $f(h_B)=h_A+d_A(a)+a^2$.
\end{enumerate}

 So, the morphism of graded rings~$f$ does \emph{not} form a commutative
diagram with the differentials $d_B$ and~$d_A$; rather, the diagram is
commutative up to the graded commutator with the element~$a$.
 The element $a\in A^1$ is called the \emph{change-of-connection
element} of a morphism of CDG\+rings $f=(f,a)$.
 A \emph{morphism of DG\+rings} $f\:A^\bu\rarrow B^\bu$ can be defined
as a morphism of CDG\+rings $(f,a)\:(A^*,d_A,0)\rarrow(B^*,d_B,0)$
with a vanishing change-of-connection element, $a=0$.

 The \emph{composition} of two morphisms of CDG\+rings
$(g,b)\:C^\cu\rarrow B^\cu$ and $(f,a)\:\allowbreak B^\cu\rarrow A^\cu$
is defined by the obvious formula $(f,a)\circ (g,b)=
(f\circ g,\>a+f(b))\:\allowbreak C^\cu\rarrow A^\cu$.
 A morphism of CDG\+rings $(f,a)$ is said to be \emph{strict} if $a=0$.
 Morphisms of CDG\+rings of the form $(\id_{B^*},a)\:(B^*,d',h')\rarrow
(B^*,d,h)$ are called \emph{change-of-connection morphisms}.
 By the definition, one has $d'(b)=d(b)+ab-(-1)^{|b|}ba$ for
all $b\in B^{|b|}$ and $h'=h+d(a)+a^2$ in this case.
 Clearly, all change-of-connection morphisms of CDG\+rings are
isomorphisms (i.~e., invertible).
 Every morphism of CDG\+rings $(f,a)\:B^\cu\rarrow A^\cu$ factorizes
uniquely into the composition of a strict morphism followed by
a change-of-connection isomorphism, $(f,a)=(\id_{A^*},a)\circ(f,0)$.

 A \emph{left CDG\+module} $M^\cu=(M^*,d_M)$ over a CDG\+ring
$B^\cu=(B^*,d,h)$ is a set of data consisting of
\begin{itemize}
\item a $\boZ$\+graded left $B^*$\+module $M^*=\bigoplus_{i\in\boZ}M^i$;
\item an odd derivation~$d_M$ of degree~$1$ on $M^*$ compatible with
the odd derivation~$d$ on $B^*$, that is $d_M\:M^i\rarrow M^{i+1}$
for all $i\in\boZ$ and $d_M(bx)=d(b)x+(-1)^{|b|}bd_M(x)$ for all
homogeneous elements $b\in B$ and $x\in M$ of degrees $|b|$
and $|x|\in\boZ$.
\end{itemize}
 The following axiom must be satisfied:
\begin{enumerate}
\renewcommand{\theenumi}{\roman{enumi}}
\setcounter{enumi}{4}
\item $d_M(d_M(x))=hx$ for all $x\in M^*$.
\end{enumerate}

 For a morphism of CDG\+rings $(f,a)\:B^\cu\rarrow A^\cu$ and
a left CDG\+module $M^\cu=(M^*,d_M)$ over $A^\cu$, the induced
left CDG\+module structure over $B^\cu$ on $M^\cu$ is given by
the pair $(M^*,d'_M)$ with the differential~$d'_M$ on $M^*$ defined
by the formula
\begin{enumerate}
\renewcommand{\theenumi}{\roman{enumi}}
\setcounter{enumi}{5}
\item $d'_M(x)=d_M(x)+ax$ for all $x\in M^*$.
\end{enumerate}

 A \emph{right CDG\+module} $N^\cu=(N^*,d_N)$ over a CDG\+ring
$B^\cu=(B^*,d,h)$ is a set of data consisting of
\begin{itemize}
\item a $\boZ$\+graded right $B^*$\+module
$N^*=\bigoplus_{i\in\boZ}N^i$;
\item an odd derivation~$d_N$ of degree~$1$ on $N^*$ compatible with
the odd derivation~$d$ on $B^*$, that is $d_N\:N^i\rarrow N^{i+1}$
for all $i\in\boZ$ and $d_N(yb)=d_N(y)b+(-1)^{|y|}yd(b)$ for all
homogeneous elements $b\in B$ and $y\in N$ of degrees $|b|$ and
$|y|\in\boZ$.
\end{itemize}
 The following axiom must be satisfied:
\begin{enumerate}
\renewcommand{\theenumi}{\roman{enumi}}
\setcounter{enumi}{6}
\item $d_N(d_N(y))=-yh$ for all $y\in N^*$.
\end{enumerate}

 For a morphism of CDG\+rings $(f,a)\:B^\cu\rarrow A^\cu$ and
a right CDG\+module $N^\cu=(N^*,d_N)$ over $A^\cu$, the induced
right CDG\+module structure over $B^\cu$ on $N^\cu$ is given by
the pair $(N^*,d'_N)$ with the differential~$d'_N$ on $N^*$ defined
by the formula
\begin{enumerate}
\renewcommand{\theenumi}{\roman{enumi}}
\setcounter{enumi}{7}
\item $d'_N(y)=d_M(y)-(-1)^{|y|}ya$ for all $y\in N^{|y|}$.
\end{enumerate}

 Notice that the CDG\+ring $B^\cu=(B^*,d,h)$ does \emph{not} have any
natural structure of a left or right CDG\+module over itself (when
$h\ne0$), because of the mismatch between the formulas~(i), (v),
and~(vii).
 However, any CDG\+ring $B^\cu$ is naturally a \emph{CDG\+bimodule}
over itself, in the sense of the next definition.

 Let $A^\cu=(B^*,d_B,h_B)$ and $C^\cu=(C^*,d_C,h_C)$ be two CDG\+rings.
 A \emph{CDG\+bimodule} $K^\cu$ over $B^\cu$ and $C^\cu$ is a set of
data consisting of
\begin{itemize}
\item a graded $B^*$\+$C^*$\+bimodule $K^*=\bigoplus_{i\in\boZ}K^i$;
\item an odd derivation $d_K$ of degree~$1$ on $K^*$ compatible with
the odd derivations $d_B$ on $B^*$ and $d_C$ on $C^*$, that is
$d_K\:K^i\rarrow K^{i+1}$ for all $i\in\boZ$ and both the identities
$d_K(bz)=d_B(b)z+(-1)^{|b|}bd_K(z)$ and $d_K(zc)=d_K(z)c+
(-1)^{|z|}zd_C(c)$ are satisfied for all $b\in B^{|b|}$, \
$z\in K^{|z|}$, and $c\in C^{|c|}$.
\end{itemize}
 The following axiom must be satisfied:
\begin{enumerate}
\renewcommand{\theenumi}{\roman{enumi}}
\setcounter{enumi}{8}
\item $d_K(d_K(z))=h_Bz-zh_C$ for all $z\in K^*$.
\end{enumerate}

 Let $(f,a_B)\:B^\cu{}'\rarrow B^\cu$ and
$(g,a_C)\:C^\cu{}'\rarrow C^\cu$ be two morphisms of CDG\+rings.
 For a CDG\+bimodule $K^\cu$ over $B^\cu$ and $C^\cu$, the induced
CDG\+bimodule structure over $B^\cu{}'$ and $C^\cu{}'$ on $K^\cu$
is given by the pair $(K^*,d'_K)$ with the differential~$d'_K$ on
$K^*$ defined by the formula
\begin{enumerate}
\renewcommand{\theenumi}{\roman{enumi}}
\setcounter{enumi}{9}
\item $d'_K(z)=d_K(z)+a_Bz-(-1)^{|z|}za_C$ for all $z\in K^{|z|}$.
\end{enumerate}

 Given a graded ring $B^*$ and two graded left $B^*$\+modules $L^*$ and
$M^*$, we denote by $\Hom_{B^*}^*(L^*,M^*)$ the graded abelian group
of homogeneous $B^*$\+module maps $L^*\rarrow M^*$ of various degrees
$n\in\boZ$ with the following sign rule.
 The elements of $\Hom_{B^*}^n(L^*,M^*)$ are homogeneous additive maps
$f\:L^*\rarrow M^*$ of degree~$n$ such that, for every $b\in B^{|b|}$
and $l\in L^{|l|}$ one has $f(bl)=(-1)^{n|b|}bf(l)$ in~$M^*$.

 For any left CDG\+modules $L^\cu=(L^*,d_L)$ and $M^\cu=(M^*,d_M)$ over
a CDG\+ring $B^\cu$, the graded abelian group $\Hom_{B^*}^*(L^*,M^*)$
is endowed with a natural differential~$d$ given by the usual formula
$d(f)(l)=d_M(f(l))-(-1)^{|f|}f(d_L(l))$ for all $l\in L^{|l|}$ and
$f\in\Hom_{B^*}^{|f|}(L^*,M^*)$.
 The differential~$d$ on $\Hom_{B^*}^*(L^*,M^*)$ squares to zero:
one has $d(d(f))=0$ for all $f\in\Hom_{B^*}^*(L^*,M^*)$, because
the two curvature-related terms arising from $d_M^2$ and~$d_L^2$
cancel each other in the expression for~$d^2$.
 So the graded abelian group $\Hom_{B^*}^*(L^*,M^*)$ with
the differential~$d$ is a complex, which we denote by
$\Hom_{B^*}^\bu(L^\cu,M^\cu)$.

 Given a graded ring $B^*$ and two graded right $B^*$\+modules $R^*$ and
$N^*$, we denote by $\Hom_{B^*{}^\rop}^*(R^*,N^*)$ the graded abelian
group of homogeneous $B^*$\+module maps $R^*\rarrow N^*$ of various
degrees $n\in\boZ$.
 In the case of right modules, there is no sign rule.
 The elements of $\Hom_{B^*{}^\rop}^n(R^*,N^*)$ are homogeneous additive
maps $f\:R^*\rarrow N^*$ of degree~$n$ such that, for every
$b\in B^{|b|}$ and $r\in R^{|r|}$ one has $f(rb)=f(r)b$ in~$N^*$.

 For any right CDG\+modules $R^\cu=(R^*,d_R)$ and $N^\cu=(N^*,d_N)$ over
a CDG\+ring $B^\cu$, the graded abelian group
$\Hom_{B^*{}^\rop}^*(R^*,N^*)$ is endowed with a natural
differential~$d$ given by the same formula $d(f)(r)=d_N(f(r))-
(-1)^{|f|}f(d_R(r))$ for all $r\in R^{|r|}$ and
$f\in\Hom_{B^*{}^\rop}^{|f|}(R^*,N^*)$.
 Similarly to the case of left modules, the differential~$d$ on
$\Hom_{B^*{}^\rop}^*(R^*,N^*)$ squares to zero:
one has $d(d(f))=0$ for all $f\in\Hom_{B^*{}^\rop}^*(R^*,N^*)$.
 So the graded abelian group $\Hom_{B^*{}^\rop}^*(R^*,N^*)$ with
the differential~$d$ is a complex, which we denote by
$\Hom_{B^*{}^\rop}^\bu(R^\cu,N^\cu)$.

 Given a graded ring $B^*$, a graded right $B^*$\+module $N^*$, and
a graded left $B^*$\+module $M^*$, the tensor product $N^*\ot_{B^*}M^*$
is naturally a graded abelian group, as mentioned in
Section~\ref{graded-modules-subsecn}.
 The definition of such tensor product is the usual one, and there is
no sign rule involved.

 For any left CDG\+module $N^\cu=(N^*,d_N)$ and any right CDG\+module
$M^\cu=(M^*,d_M)$ over a CDG\+ring $B^\cu$, the graded abelian group
$N^*\ot_{B^*}M^*$ is endowed with a natural differential~$d$ given
by the usual formula $d(y\ot x)=d_N(y)\ot x+(-1)^{|y|}y\ot d_M(x)$
for all $y\in N^{|y|}$ and $x\in M^{|x|}$.
 The differential~$d$ on $N^*\ot_{B^*}M^*$ squares to zero: one has
$d(d(y\ot x))=0$ for all $y\ot x\in N^*\ot_{B^*}M^*$.
 Once again, this happens because the two curvature-related terms
arizing from $d_N^2$ and $d_M^2$ cancel each other in the expression
for~$d^2$.
 So the graded abelian group $N^*\ot_{B^*}M^*$ is a complex, which we
denote by $N^\cu\ot_{B^*}M^\cu$.

 More generally, let $A^\cu$, $B^\cu$, and $C^\cu$ be three
CDG\+rings.
 Let $N^\cu$ be a CDG\+bimodule over $A^\cu$ and $B^\cu$, and let
$M^\cu$ be a CDG\+bimodule over $B^\cu$ and $C^\cu$.
 Then the tensor product of graded bimodules $N^*\ot_{B^*}M^*$,
endowed by the differential given by the usual formula as above,
is a CDG\+bimodule over $A^\cu$ and~$C^\cu$.
 We denote this CDG\+bimodule by $N^\cu\ot_{B^*}M^\cu$.

 Similarly, let $L^\cu$ be a CDG\+bimodule over $B^\cu$ and $A^\cu$,
and let $M^\cu$ be a CDG\+bimodule over $B^\cu$ and $C^\cu$.
 Then the graded abelian group $\Hom_{B^*}^*(L^*,M^*)$, endowed with
the differential given by the usual formula as above, is
a CDG\+bimodule over $A^\cu$ and~$C^\cu$.
 We denote this CDG\+bimodule by $\Hom_{B^*}^\cu(L^\cu,M^\cu)$.
 We refer to~\cite[Section~6.1]{Prel} for the details, including
the sign rule.

 Finally, let $R^\cu$ be a CDG\+bimodule over $A^\cu$ and $B^\cu$,
and let $N^\cu$ be a CDG\+bimodule over $C^\cu$ and $B^\cu$.
 Then the graded abelian group $\Hom_{B^*{}^\rop}^*(R^*,N^*)$,
endowed with the differential given by the usual formula as above,
a CDG\+bimodule over $C^\cu$ and~$A^\cu$.
 We denote this CDG\+bimodule by $\Hom_{B^*{}^\rop}^\cu(R^\cu,N^\cu)$.
 The details can be found in~\cite[Section~6.1]{Prel}
(see also~\cite[Section~1]{PS7}).

\subsection{Quasi-coherent CDG-quasi-algebras}
\label{qcoh-cdg-quasi-algebras-subsecn}
 This section provides numerous supplementary details for a much more
sketchy exposition in~\cite[Section~2.4]{Pedg}.

 The definition of a graded quasi-algebra over a commutative ring $R$
was given in Section~\ref{graded-modules-subsecn}, based on
the discussion in Sections~\ref{prelim-quasi-modules-subsecn}\+-%
\ref{prelim-quasi-algebras-subsecn}.
 The definitions of a CDG\+ring and a morphism of CDG\+rings were
given in Section~\ref{cdg-rings-cdg-modules-subsecn}.

 A \emph{CDG\+quasi-algebra} $A^\cu$ over a commutative ring $R$ is
a CDG\+ring $(A^*,d,h)$ endowed with a morphism of graded rings
$R\rarrow A^*$ making the graded ring $A^*$ a graded quasi-algebra
over~$R$.
 Let us emphasize that the differential $d\:A^*\rarrow A^*$ is
\emph{not} presumed to be $R$\+linear.
 In fact, \emph{no} condition of compatibility between the CDG\+ring
structure and the graded quasi-algebra structure on $A^*$ is imposed.

 A \emph{morphism} of CDG\+quasi-algebras $B^\cu\rarrow A^\cu$ over
$R$ is a morphism of CDG\+rings $(f,a)\:B^\cu\rarrow A^\cu$ such that
the triangular diagram of graded ring homomorphisms $R\rarrow B^*
\rarrow A^*$ is commutative.
 The composition of morphisms of CDG\+quasi-algebras is defined in
the obvious way, based on the construction of the composition of
morphisms of CDG\+rings in Section~\ref{cdg-rings-cdg-modules-subsecn}.

 Similarly, a \emph{DG\+quasi-algebra} $A^\bu$ over $R$ is
a DG\+ring $(A^*,d)$ endowed with a morphism of graded rings
$R\rarrow A^*$ making the graded ring $A^*$ a graded quasi-algebra
over~$R$.
 A \emph{morphism} of DG\+quasi-algebras $B^\bu\rarrow A^\bu$ over
$R$ is a morphism of DG\+rings $f\:B^\bu\rarrow A^\bu$ such that
the triangular diagram of graded ring homomorphisms $R\rarrow B^*
\rarrow A^*$ is commutative.

 The aim of the present section is to work out the localization and
globalization of CDG\+quasi-algebras over Zariski affine open
coverings of schemes.

\begin{lem} \label{localizing-derivation-by-element}
 Let $R$ be a commutative ring and $A^*=\bigoplus_{i\in\boZ}A^i$ be
a graded quasi-algebra over~$R$.
 Let $r\in R$ be an element and $d_A\:A^*\rarrow A^*$ be an odd
derivation of degree~$1$.
 Consider the graded quasi-algebra $A^*[r^{-1}]=R[r^{-1}]\ot_RA^*$
over the commutative ring $R[r^{-1}]$, as per the graded version of
Lemma~\ref{quasi-algebra-co-extension-of-scalars}(a).
 Then there exists a unique extension of~$d_A$ to an odd derivation
$d_{A[r^{-1}]}\:A^*[r^{-1}]\rarrow A^*[r^{-1}]$ of degree~$1$ forming
a commutative square diagram with the derivation~$d_A$ and the natural
(localization) homomorphism of graded rings $A^*\rarrow A^*[r^{-1}]$.
\end{lem}

\begin{proof}
 It suffices to say that the graded ring $A^*[r^{-1}]$ is freely
generated by the graded ring $A^*$ and a formal inverse element~$r^{-1}$
to the element $r\in A^0$.
 Then it follows that the assignment $d_{A[r^{-1}]}(a)=d_A(a)$
for $a\in A^*$ and $d_{A[r^{-1}]}(r^{-1})=-r^{-1}d_A(r)r^{-1}$ extends
uniquely to an odd derivation $d_{A[r^{-1}]}\:A^*[r^{-1}]
\rarrow A^*[r^{-1}]$.
 It suffices to check that the only relations $rr^{-1}=1=r^{-1}r$
are preserved, that is $d_{A[r^{-1}]}(r)r^{-1}+
rd_{A[r^{-1}]}(r^{-1})=0=d_A(1)$ and
$r^{-1}d_{A[r^{-1}]}(r)+d_{A[r^{-1}]}(r^{-1})r=0=d_A(1)$
(see~\cite[Lemma~3.18]{Prel}).

 More explicitly, let $n\ge0$ be an integer.
 The formula
\begin{equation} \label{localization-by-element-derivation-I}
 d_{A[r^{-1}]}(r^{-n}a)=
 r^{-n}d_A(a)-\sum\nolimits_{i=1}^n r^{-i}d_A(r)r^{n-i+1}a
 \quad \text{for all $a\in A^*$}
\end{equation}
where $r\in R$ is viewed as an element of $A^0$, defines
the derivation~$d_{A[r^{-1}]}$.
 Equivalently,
\begin{equation} \label{localization-by-element-derivation-II}
 d_{A[r^{-1}]}(br^{-n})=
 d_A(b)r^{-n}-(-1)^{|b|}\sum\nolimits_{i=1}^n br^{-i}d_A(r)r^{n-i+1}
 \quad \text{for all $b\in A^{|b|}$}.
\end{equation}

 To prove that the two formulas agree with each other and define
an odd derivation $d_{A[r^{-1}]}\:A^*[r^{-1}]\rarrow A^*[r^{-1}]$
compatible with~$d_A$, denote temporarily by~$d'$ the output of
the formula~\eqref{localization-by-element-derivation-I} and by~$d''$
the output of the formula~\eqref{localization-by-element-derivation-II}.
 Clearly, $d'(r^{-0}a)=d_A(a)$ and $d''(br^{-0})=d_A(b)$.
 To show that the maps $d'$ and $d''\:A^*[r^{-1}]\rarrow A^*[r^{-1}]$
are well-defined, one needs to check that $d'(r^{-n-1}(ra))=d'(r^{-n}a)$
and $d''((br)r^{-n-1})=d''(br^{-n})$ in $A^*[r^{-1}]$ for all
$a\in A^{|a|}$ and $b\in A^{|b|}$.
 
 Then one observes that $d_A(r^n)r^{-n}a+r^nd'(r^{-n}a)=d_A(a)$
and $d''(br^{-n})r^n+(-1)^{|b|}br^{-n}d_A(r^n)=d_A(b)$.
 Furthermore, the identities $d'(r^{-n}a_1a_2)=d'(r^{-n}a_1)a_2+
(-1)^{|a_1|}r^{-n}a_1d_A(a_2)$ and $d''(b_1b_2r^{-n})=d_A(b_1)b_2r^{-n}+
(-1)^{|b_1|}b_1d''(b_2r^{-n})$ for all $a_1$, $a_2$, $b_1$,
$b_2\in A^*$ are also easy to check.
 Now one can compute that the equation $r^{-n}a=br^{-n}$ in
$A^*[r^{-1}]$ implies
\begin{multline*}
 r^nd'(r^{-n}a)r^n=d_A(r^nr^{-n}ar^n)-d_A(r^n)r^{-n}ar^n-
 (-1)^{|a|}r^nr^{-n}ad_A(r^n) \\
 =d_A(r^nbr^{-n}r^n)-d_A(r^n)br^{-n}r^n-(-)^{|b|}r^nbr^{-n}d_A(r^n)
 =r^nd''(br^{-n})r^n,
\end{multline*}
hence $d'(r^{-n}a)=d''(br^{-n})$ and $d'=d''$.
 So the notation $d_{A[r^{-1}]}=d'=d''\:A^*[r^{-1}]\allowbreak\rarrow 
A^*[r^{-1}]$ is unambiguous.
 It remains to notice that $d'(r^{-n}x)=d'(r^{-n})x+r^{-n}d'(x)$
and $d''(yr^{-n})=d''(y)r^{-n}+(-1)^{|y|}yd''(r^{-n})$ for all
$x$, $y\in A^*[r^{-1}]$, which allows one to deduce the desired equation
$d_{A[r^{-1}]}(r^{-n}abr^{-n})=d_{A[r^{-1}]}(r^{-n}a)br^{-n}
+(-1)^{|a|}r^{-n}ad_{A[r^{-1}]}(br^{-n})$ for all $a$, $b\in A^*$.

 The uniqueness of $d_{A[r^{-1}]}$ holds because the graded ring
$A^*[r^{-1}]$ is generated by $A^*$ and~$r^{-1}$, and one has
$d(r^{-1})=-r^{-1}d(r)r^{-1}$ for any odd derivation~$d$ on a graded
algebra $B^*$ with an invertible element $r\in B^0$.
\end{proof}

\begin{lem} \label{localizing-derivation-by-affine-open-immersion}
 Let $R$ be a commutative ring and $A^*=\bigoplus_{i\in\boZ}A^i$ be
a graded quasi-algebra over~$R$.
 Let $d_A\:A^*\rarrow A^*$ be an odd derivation of degree~$1$.
 Let $f\:R\rarrow S$ be a homomorphism of commutative rings such that
the related morphism of affine schemes\/ $\Spec S\rarrow\Spec R$ is
an open immersion.
 Consider the graded quasi-algebra $B^*=S\ot_RA^*\simeq A^*\ot_RS$
over the commutative ring $S$, as per the graded version of
Lemma~\ref{quasi-algebra-co-extension-of-scalars}(a).
 Then there exists a unique extension of~$d_A$ to an odd derivation
$d_B\:B^*\rarrow B^*$ of degree\/~$1$ forming a commutative square
diagram with the derivation~$d_A$ and the natural homomorphism of graded
rings $A^*\rarrow B^*$.
\end{lem}

\begin{proof}
 Let $r_1$,~\dots, $r_N\in R$ be a finite family of elements such
that $\Spec S=\bigcup_{\alpha=1}^N\Spec R[r_\alpha^{-1}]\subset\Spec R$.
 Then we have homomorphisms of commutative ring $R\rarrow S\rarrow
R[r_\alpha^{-1}]$ for all $1\le\alpha\le N$, and natural isomorphisms
of commutative rings $R[r_\alpha^{-1}]\simeq S[f(r_\alpha)^{-1}]$.

 To prove uniqueness of~$d_B$, consider the graded quasi-algebra
$C_\alpha^*=S[f(r_\alpha)^{-1}]\ot_SB^*\simeq R[r_\alpha^{-1}]\ot_RA^*$.
 Suppose given two odd derivations $d'_B$ and $d''_B\:B^*\rarrow B^*$
satisfying the conditions in the lemma.
 By Lemma~\ref{localizing-derivation-by-element} applied to the graded
quasi-algebra $B^*$ over $S$ and the element $f(r_\alpha)\in S$, there
exists a unique odd derivation $d'_\alpha\:C_\alpha^*\rarrow C_\alpha^*$
forming a commutative square diagram with~$d'_B$, and a unique odd
derivation $d''_\alpha\:C_\alpha^*\rarrow C_\alpha^*$ forming
a commutative square diagram with~$d''_B$.
 As both the derivations $d'_B$ and~$d''_B$ form commutative square
diagrams with~$d_A$, it follows that both the derivations $d'_\alpha$
and~$d''_\alpha$ form commutative square diagrams with~$d_A$.
 Applying Lemma~\ref{localizing-derivation-by-element} to the graded 
quasi-algebra $A^*$ over $R$ and the element $r_\alpha\in R$, we
conclude that $d'_\alpha=d''_\alpha$.
 Since this holds for all $1\le\alpha\le N$ and the map of graded rings
$B^*\rarrow\bigoplus_{\alpha=1}^N C_\alpha^*$ is injective, it follows
that $d'_B=d''_B$, as desired.

 To prove existence, consider a subsequence of indices
$1\le\alpha_1<\dotsb<\alpha_k\le N$, where $k\ge1$.
 Put $C_{\alpha_1,\dotsc,\alpha_k}^*=
R[r_{\alpha_1}^{-1},\dotsc,r_{\alpha_k}^{-1}]\ot_RA^*$.
 By Lemma~\ref{localizing-derivation-by-element}, there exists
a unique odd derivation $d_{\alpha_1,\dotsc,\alpha_k}\:
C_{\alpha_1,\dotsc,\alpha_k}^*\rarrow C_{\alpha_1,\dotsc,\alpha_k}^*$
on the graded ring $C_{\alpha_1,\dotsc,\alpha_k}^*$ forming
a commutative square diagram with~$d_A$.
 Consider the \v Cech
coresolution~\eqref{bimodule-B-left-cech-coresolution} of
the graded quasi-module $B^*=S\ot_RA^*$ over~$S$,
\begin{multline} \label{graded-bimodule-B-left-cech-coresolution}
 0\lrarrow B^*\lrarrow\bigoplus\nolimits_{\alpha=1}^N
 C_\alpha^*\lrarrow\bigoplus\nolimits_{1\le\alpha<\beta\le N}
 C_{\alpha,\beta}^* \\ \lrarrow\dotsb\lrarrow
 C_{1,\dotsc,N}^*\lrarrow0
\end{multline}
(see the proof of Lemma~\ref{quasi-module-localization-lemma}(a)).
 Endowing every term of the \v Cech
coresolution~\eqref{graded-bimodule-B-left-cech-coresolution}, except
the leftmost one, with its own ``vertical'' differential equal to
the direct sum of the respective derivations
$d_{\alpha_1,\dotsc,\alpha_k}$, we notice that it follows from
Lemma~\ref{localizing-derivation-by-element} that the vertical
and ``horizontal'' (\v Cech) differentials
on~\eqref{graded-bimodule-B-left-cech-coresolution} commute.
 Since the coresolution~\eqref{graded-bimodule-B-left-cech-coresolution}
is exact, we can pass to the kernel and conclude that there exists
a unique differential~$d_B$ on $B^*$ making a commutative square diagram
with the vertical differential on $\bigoplus_{\alpha=1}^N C_\alpha^*$
and the leftmost nontrivial horizontal differential
in~\eqref{graded-bimodule-B-left-cech-coresolution}.
 Since the latter differential $B^*\rarrow\bigoplus_{\alpha=1}^N
C_\alpha^*$ is an injective homomorphism of graded rings, while
the vertical differential on $\bigoplus_{\alpha=1}^N C_\alpha^*$ is
an odd derivation, it follows that $d_B$~is an odd derivation on $B^*$,
as desired.

 An alternative approach is to notice that that every grading component
$d_{A,n}\:A^n\allowbreak\rarrow A^{n+1}$ of the odd derivation~$d_A$ is
a differential operator of order not exceeding~$1$ in the sense
of~\cite[Section~IV.16.8]{EGA4} or~\cite[Section Tag~09CH]{SP}, or
a strongly $R$\+differential operator of order~$\le1$ in
the terminology of~\cite[Section~3.1]{Ptd}.
 By~\cite[Proposition~6.6]{Ptd}, there exists a unique strongly
$S$\+differential operator $d_{B,n}\:B^n\rarrow B^{n+1}$ of
order~$\le1$ forming a commutative square diagram with~$d_{A,n}$ and
the natural maps $A^n\rarrow B^n$ and $A^{n+1}\rarrow B^{n+1}$.
 It remains to use the uniqueness assertion
of~\cite[Proposition~6.6]{Ptd} in order to prove that the collection
of differential operators~$d_{B,n}$ constructed in this way is
an odd derivation $d_B\:B^*\rarrow B^*$.
 The advantage of this argument is that it works more generally for
any flat epimorphism of commutative rings $f\:R\rarrow S$ (not
necessarily of finite presentation).
\end{proof}

\begin{lem} \label{strict-gluing-derivation-on-affine}
 Let $R$ be a commutative ring and $A^*=\bigoplus_{i\in\boZ}A^i$ be
a graded quasi-algebra over~$R$.
 Let $R\rarrow S_\alpha$, \ $1\le\alpha\le N$, be a finite collection
of homomorphisms of commutative rings such that the collection of
induced maps of the spectra\/ $\Spec S_\alpha\rarrow\Spec R$ is
an affine open covering of an affine scheme.
 Consider the graded quasi-algebras $B_\alpha^*=S_\alpha\ot_RA^*$
over the commutative rings $S_\alpha$ and the graded quasi-algebras
$C_{\alpha\beta}^*=S_\alpha\ot_RS_\beta\ot_RA^*$ over the commutative
rings $S_\alpha\ot_RS_\beta$, as per the graded version of
Lemma~\ref{quasi-algebra-co-extension-of-scalars}(a).
 Suppose given, for every index~$\alpha$, an odd derivation $d_\alpha\:
B_\alpha^*\rarrow B_\alpha^*$ of degree\/~$1$ on the graded
ring~$B_\alpha^*$.
 Assume that, for every pair of indices $\alpha$ and~$\beta$,
the two odd derivations on the graded ring $C_{\alpha\beta}^*$
arising from the odd derivations $d_\alpha^*$ and~$d_\beta^*$
on the graded rings $B_\alpha^*$ and $B_\beta^*$, as per
the construction of
Lemma~\ref{localizing-derivation-by-affine-open-immersion}, agree.
 Then there exists a unique odd derivation $d\:A^*\rarrow A^*$ of
degree\/~$1$ on the graded ring $A^*$ such that, for every
index~$\alpha$, the odd derivation~$d_\alpha$ on the graded ring
$B_\alpha^*$ arises from the odd derivation~$d$ on the graded ring
$A^*$ via the construction of
Lemma~\ref{localizing-derivation-by-affine-open-immersion}.
\end{lem}

\begin{proof}
 The argument, based on the \v Cech coresolution
$$
 0\lrarrow A^*\lrarrow\bigoplus\nolimits_{\alpha=1}^N
 B_\alpha^*\lrarrow\bigoplus\nolimits_{1\le\alpha<\beta\le N}
 C_{\alpha\beta}^*\lrarrow\dotsb,
$$
is similar to the proof of the existence assertion in
Lemma~\ref{localizing-derivation-by-affine-open-immersion}.
\end{proof}

\begin{lem} \label{localizing-cdg-ring-by-affine-open-immersion}
 Let $R$ be a commutative ring and $A^*=\bigoplus_{i\in\boZ}A^i$ be
a graded quasi-algebra over~$R$.
 Let $A^\cu=(A^*,d_A,h_A)$ be a CDG\+ring structure on~$A^*$.
 Let $f\:R\rarrow S$ be a homomorphism of commutative rings such that
the related morphism of affine schemes\/ $\Spec S\rarrow\Spec R$ is
an open immersion.
 Consider the graded quasi-algebra $B^*=S\ot_RA^*\simeq A^*\ot_RS$
over the commutative ring $S$, as per the graded version of
Lemma~\ref{quasi-algebra-co-extension-of-scalars}(a).
 Then there exists a unique CDG\+ring structure $B^\cu=(B^*,d_B,h_B)$
on the graded ring $B^*$ such that the natural homomorphism of
graded rings $g=S\ot_Rf\:A^*\rarrow B^*$ is a strict morphism of
CDG\+rings $(g,0)\:A^\cu\rarrow B^\cu$.
\end{lem}

\begin{proof}
 For any strict morphism of CDG\+rings $(g,0)\:(A^*,d_A,h_A)\rarrow
(B^*,d_B,h_B)$, one has $h_B=g(h_A)$.
 Thus the uniqueness assertion of the lemma follows from
the uniqueness assertion of
Lemma~\ref{localizing-derivation-by-affine-open-immersion}.
 To prove the existence, it remains to show that the odd derivation
$d_B\:B^*\rarrow B^*$ provided by
Lemma~\ref{localizing-derivation-by-affine-open-immersion} together
with the curvature element $h_B=g(h_A)$ form a CDG\+ring structure
on~$B^*$.
 There are two equations~(i\+-ii) from
Section~\ref{cdg-rings-cdg-modules-subsecn} that we need to check.
 The equation~(ii), \,$d_B(h_B)=0$, holds because
$d_B(h_B)=d_B(g(h_A))=g(d_A(h_A))=g(0)=0$.
 To check the equation~(i) for the square~$d_B^2$ of
the derivation~$d_B$, notice that both~$d_B^2$ and the commutator
map $[h_B,{-}]$ are even derivations of degree~$2$ on~$B^*$.
 Both of them form a commutative square diagram with one and the same
even derivation $d_A^2=[h_A,{-}]$ on~$A^*$.
 Now one needs to go through the rounds of the proofs of
Lemmas~\ref{localizing-derivation-by-element}
and~\ref{localizing-derivation-by-affine-open-immersion} and check that
the uniqueness assertions of the two lemmas apply to even derivations
of degree~$2$ just as well as to odd derivations of degree~$1$.
\end{proof}

\begin{lem} \label{strict-gluing-cdg-ring-structure-on-affine}
 Let $R$ be a commutative ring and $A^*=\bigoplus_{i\in\boZ}A^i$ be
a graded quasi-algebra over~$R$.
 Let $R\rarrow S_\alpha$, \ $1\le\alpha\le N$, be a finite collection
of homomorphisms of commutative rings such that the collection of
induced maps of the spectra\/ $\Spec S_\alpha\rarrow\Spec R$ is
an affine open covering of an affine scheme.
 Consider the graded quasi-algebras $B_\alpha^*=S_\alpha\ot_RA^*$
over the commutative rings $S_\alpha$ and the graded quasi-algebras
$C_{\alpha\beta}^*=S_\alpha\ot_RS_\beta\ot_RA^*$ over the commutative
rings $S_\alpha\ot_RS_\beta$, as per the graded version of
Lemma~\ref{quasi-algebra-co-extension-of-scalars}(a).
 Suppose given, for every index~$\alpha$, a CDG\+ring structure
$B^\cu_\alpha=(B^*_\alpha,d_\alpha,h_\alpha)$ on the graded
ring~$B_\alpha^*$.
 Assume that, for every pair of indices $\alpha$ and~$\beta$,
the two CDG\+ring structures on the graded ring $C_{\alpha\beta}^*$
arising from the CDG\+ring structures $B^\cu_\alpha$ and $B^\cu_\beta$
on the graded rings $B_\alpha^*$ and $B_\beta^*$, as per
the construction of
Lemma~\ref{localizing-cdg-ring-by-affine-open-immersion}, are
equal to each other (in the strict sense of the equality of the two
differentials and the two change-connection elements).
 Then there exists a unique CDG\+ring structure $A^\cu=(A^*,d,h)$ on
the graded ring $A^*$ such that, for every index~$\alpha$,
the CDG\+ring structure $B^\cu_\alpha$ on the graded ring $B_\alpha^*$
arises from the CDG\+ring structure $A^\cu$ on the graded ring $A^*$
via the construction of
Lemma~\ref{localizing-cdg-ring-by-affine-open-immersion}.
\end{lem}

\begin{proof}
 The existence of an odd derivation $d\:A^*\rarrow A^*$ of degree~$1$
forming commutative square diagrams with the odd derivations
$d_\alpha\:B_\alpha^*\rarrow B_\alpha^*$ is a particular case of
the result of Lemma~\ref{strict-gluing-derivation-on-affine}.
 Using the \v Cech coresolution from the proof of the latter lemma,
one can see that there exists a unique element $h\in A^2$ whose
images in $B_\alpha^2$ are equal to~$h_\alpha$ for all indices~$\alpha$.
 The equations~(i\+-ii) from the definition of a CDG\+ring in
Section~\ref{cdg-rings-cdg-modules-subsecn} hold for $(A^*,d,h)$,
since they hold for $(B_\alpha^*,d_\alpha,h_\alpha)$ and the graded
ring homomorphism $A^*\rarrow\bigoplus_{\alpha=1}^N B_\alpha^*$
is injective.
\end{proof}

 The following definition goes back to~\cite[Section~B.1]{Pkoszul},
\cite[Section~1.2]{EP}, and~\cite[Section~2.4]{Pedg}.
 Let $X$ be a scheme with an open covering~$\bW$.
 A \emph{quasi-coherent CDG\+quasi-algebra} \,$\cB^\cu$ over
the scheme $X$ with the open covering $\bW$ is a set of data
consisting of
\begin{itemize}
\item a quasi-coherent graded quasi-algebra
$\cB^*=\bigoplus_{i\in\boZ}\cB^i$ over~$X$ (see
Sections~\ref{quasi-coherent-quasi-algebras-subsecn}
and~\ref{graded-modules-subsecn} for the definition);
\item an odd derivation $d_U=d_{\cB,U}\:\cB^*(U)\rarrow\cB^*(U)$ of
degree~$1$ on the graded ring $\cB^*(U)$ and a curvature element
$h_U=h_{\cB,U}\in\cB^2(U)$ given for each affine open subscheme
$U\subset X$ subordinate to~$\bW$;
\item a change-of-connection element $a_{VU}=a_{\cB,VU}\in\cB^1(V)$
given for each pair of affine open subschemes $V\subset U\subset X$
subordinate to~$\bW$.
\end{itemize}
 The following three axioms must be satisfied:
\begin{enumerate}
\renewcommand{\theenumi}{\roman{enumi}}
\setcounter{enumi}{10}
\item for each affine open subscheme $U\subset X$ subordinate to $\bW$,
the triple $\cB^\cu(U)=(\cB^*(U),d_{\cB,U},h_{\cB,U})$ is a CDG\+ring
(see Section~\ref{cdg-rings-cdg-modules-subsecn} for the definition);
\item for each pair of affine open subschemes $V\subset U\subset X$
subordinate to  $\bW$, the pair $(\rho_{\cB,VU},a_{\cB,VU})$ is
a morphism of CDG\+rings $(\rho_{\cB,VU},a_{\cB,VU})\:\allowbreak
(\cB^*(U),d_{\cB,U},h_{\cB,U})\rarrow(\cB^*(V),d_{\cB,V},h_{\cB,V})$,
where $\rho_{UV}=\rho_{\cB,UV}\:\cB^*(U)\rarrow\cB^*(V)$ is
the restriction morphism in the sheaf of graded rings $\cB^*$ on~$X$;
\item for each triple of affine open subschemes $T\subset V\subset U
\subset X$ subordinate to~$\bW$, the morphisms of CDG\+rings
$(\rho_{VU},a_{VU})\:(\cB^*(U),d_U,h_U)\rarrow(\cB^*(V),d_V,h_V)$,
\ $(\rho_{TV},a_{TV})\:(\cB^*(V),d_V,h_V)\rarrow(\cB^*(T),d_T,h_T)$,
and $(\rho_{TU},a_{TU})\:(\cB^*(U),d_U,h_U)\rarrow(\cB^*(T),d_T,h_T)$
form a commutative triangular diagram in the category of CDG\+rings;
in other words, the equation $a_{\cB,TU}=a_{\cB,TV}+
\rho_{\cB,TV}(a_{\cB,VU})$ holds in $\cB^1(T)$.  \hbadness=1525
\end{enumerate}

 A \emph{morphism of quasi-coherent CDG\+quasi-algebras} $f^\cu\:\cB^\cu
\rarrow\cA^\cu$ over the scheme $X$ with the open covering $\bW$ is
a set of data consisting of
\begin{itemize}
\item a morphism of quasi-coherent graded quasi-algebras $f\:\cB^*
\rarrow\cA^*$ over~$X$;
\item a change-of-connection element $a_{f,U}\in\cA^1(U)$ given
for each affine open subscheme $U\subset X$ subordinate to~$\bW$.
\end{itemize}
 The following two axioms must be satisfied:
\begin{enumerate}
\renewcommand{\theenumi}{\roman{enumi}}
\setcounter{enumi}{13}
\item for each affine open subscheme $U\subset X$ subordinate to $\bW$,
the pair $(f(U),a_{f,U})$ is a morphism of CDG\+rings
$(f(U),a_{f,U})\:(\cB^*(U),d_{\cB,U},h_{\cB,U})\rarrow
(\cA^*(U),d_{\cA,U},h_{\cA,U})$, where $f(U)\:\cB^*(U)\rarrow\cA^*(U)$
is the homomorphism of graded rings assigned to the affine open
subscheme $U\subset X$ by the morphism of quasi-coherent graded
quasi-algebras $f\:\cB^*\rarrow\cA^*$ over~$X$;
\item for each pair of affine open subschemes $V\subset U\subset X$
subordinate to~$\bW$, the square diagram of morphisms of CDG\+rings
$$
 \xymatrix{
  (\cB^*(U),d_{\cB,U},h_{\cB,U}) \ar[rr]^{(f(U),a_{f,U})}
  \ar[d]_{(\rho_{\cB,VU},a_{\cB,VU})}
  && (\cA^*(U),d_{\cA,U},h_{\cA,U})
  \ar[d]^{(\rho_{\cA,VU},a_{\cA,VU})} \\
  (\cB^*(V),d_{\cB,V},h_{\cB,V}) \ar[rr]^{(f(V),a_{f,V})}
  && (\cA^*(V),d_{\cA,V},h_{\cA,V})
 }
$$
is commutative; in other words, the equation
$a_{\cA,VU}+\rho_{\cA,VU}(a_{f,U})=f(V)(a_{\cB,VU})+a_{f,V}$
holds in $\cA^1(V)$.
\end{enumerate}

 The \emph{composition} of morphisms of quasi-coherent
CDG\+quasi-algebras is defined in the obvious way, based on
the definition of the composition of morphisms of CDG\+rings
in Section~\ref{cdg-rings-cdg-modules-subsecn}.
 We leave it to the reader to fill the straightforward details.

 A morphism of quasi-coherent CDG\+quasi-algebras $f^\cu\:\cB^\cu
\rarrow\cA^\cu$ over a scheme $X$ with an open covering $\bW$ is said
to be \emph{strict} if all the change-of-connection
elements~$a_{f,U}$ vanish.
 We will use this definition in
Section~\ref{opposite-cdg-rings-and-twisted-lie-subsecn} below.

\begin{cor} \label{qcoh-cdg-quasi-algebras-not-depend-on-covering}
 The category of quasi-coherent CDG\+quasi-algebras on a scheme $X$
does not depend on an open covering\/~$\bW$.
 Specifically, the functor of restriction to affine open subschemes
subordinate to\/ $\bW$ is an equivalence between the category of
quasi-coherent CDG\+quasi-algebras on the scheme $X$ with the open
covering\/ $\{X\}$ and the category of quasi-coherent
CDG\+quasi-algebras on the scheme $X$ with the open covering\/~$\bW$.
\end{cor}

\begin{proof}
 Let $\cB^\cu$ be a quasi-coherent CDG\+quasi-algebra on the scheme
$X$ with the open covering~$\bW$.
 We need to explain how to extend $\cB^\cu$ to a quasi-coherent
CDG\+quasi-algebra on the scheme $X$ with the open covering~$\{X\}$.

 Let $U\subset X$ be an affine open subscheme, and let
$U=\bigcup_{\alpha=1}^N V_\alpha$ be a finite covering of $U$ by
affine open subschemes $V_\alpha$ subordinate to~$\bW$.
 For every pair of indices $\alpha$ and~$\beta$, put
$V_{\alpha\beta}=V_\alpha\cap V_\beta\subset U$ and
$c_{\alpha\beta}=a_{V_{\alpha\beta}V_\alpha}-
a_{V_{\alpha\beta}V_\beta}\in\cB^1(V_{\alpha\beta})$.
 It follows from axiom~(xiii) that the collection of elements
$c=(c_{\alpha\beta})_{\alpha,\beta=1}^N$ is a \v Cech $1$\+cocycle
with coefficients in the sheaf of abelian groups $\cB^1|_U$ on
the scheme $U$ with respect to the open covering
$U=\bigcup_\alpha V_\alpha$.
 Since the sheaf $\cB^1|_U$, viewed as a sheaf of $\cO_U$\+modules
with respect either to the left or to the right $\cO_U$\+module
structure, is quasi-coherent, the \v Cech $1$\+cocycle~$c$ is
a coboundary (see, e.~g., \cite[Section Tag~01X8]{SP}).
 This means that there exists a collection of elements
$b_\alpha\in\cB^1(V_\alpha)$, \,$1\le\alpha\le N$, such that
$c_{\alpha\beta}=\rho_{V_{\alpha\beta}V_\alpha}(b_\alpha)-
\rho_{V_{\alpha\beta}V_\beta}(b_\beta)$ for all indices~$\alpha$
and~$\beta$.

 Now, for every index~$\alpha$, consider the CDG\+ring structure
$(\cB^*(V_\alpha),d'_\alpha,h'_\alpha)$ on the graded ring
$\cB^*(V_\alpha)$ such that $(\id,-b_\alpha)$ is a change-of-connection
morphism of CDG\+rings $(\id,-b_\alpha)\:
(\cB^*(V_\alpha),d'_\alpha,h'_\alpha)
\rarrow(\cB^*(V_\alpha),d_{\cB,V_\alpha},h_{\cB,V_\alpha})$.
 Then the collection of CDG\+ring structures
$(\cB^*(V_\alpha),d'_\alpha,h'_\alpha)$ satisfies the assumptions
of Lemma~\ref{strict-gluing-cdg-ring-structure-on-affine}, so it
arises from a unique CDG\+ring structure $(\cB^*(U),d',h')$ on
the graded ring $\cB^*(U)$.
 At this point, we put $(\cB^*(U),d_{\cB,U},h_{\cB,U})=
(\cB^*(U),d',h')$ and $a_{V_\alpha U}=-b_\alpha$.
 The remaining verifications are easier, and left to the reader.
\end{proof}

\begin{lem} \label{qcoh-cdg-quasi-algebras-on-affine-schemes}
 Let $U=\Spec R$ be an affine scheme.
 Then the rule assigning the CDG\+quasi-algebra
$\cB^\cu(U)=(\cB^*(U),d_{\cB,U},h_{\cB,U})$ to a quasi-coherent
CDG\+quasi-algebra $\cB^\cu$ over $U$ defines a natural equivalence
between the category of quasi-coherent CDG\+quasi-algebras over $U$
and the category of CDG\+quasi-algebras over~$R$.
\end{lem}

\begin{proof}
 The functor in one direction is described in the formulation of
the lemma.
 The inverse functor assigns to a CDG\+quasi-algebra
$A^\cu=(A^*,d_A,h_A)$ over $R$ the collection of CDG\+quasi-algebras
$B^\cu=B^\cu_V=(B^*,d_B,h_B)$ over the rings of functions
$S=\cO_U(V)$ on affine open subschemes $V\subset U$ constructed in
Lemma~\ref{localizing-cdg-ring-by-affine-open-immersion}.
 So we put $\cA^\cu(V)=(\cA^*(V),d_{\cA,V},h_{\cA,V})=
B^\cu_V=(B^*,d_B,h_B)$.
 The distinctive property of the resulting quasi-coherent
CDG\+quasi-algebra $\cA^\cu$ on $U$ is that all the change-of-connection
elements $a_{\cA,TV}\in\cA^1(T)$ vanish, i.~e., $a_{\cA,TV}=0$ for all
affine open subschemes $T\subset V\subset U$.
 The point is that this property is, of course, not preserved by
isomorphisms of quasi-coherent CDG\+quasi-algebras.

 Given an arbitrary quasi-coherent CDG\+quasi-algebra $\cB^\cu$ on $U$,
put $A^\cu=\cB^\cu(U)$ and consider the related quasi-coherent
CDG\+quasi-algebra $\cA^\cu$ on $U$, as in the previous paragraph.
 Then, first of all, we have a natural isomorphism of graded
quasi-algebras $f_V\:\cA^*(V)=\cO_U(V)\ot_{\cO(U)}A^*\simeq\cB^*(V)$,
since $\cB^*$ is a quasi-coherent graded quasi-algebra over~$U$.
 Now the collection of (essentially, change-of-connection) isomorphisms
of CDG\+rings $(f_V,a_{\cB,UV})\:\cA^\cu(V)=(B^*,d_B,h_B)\rarrow
(\cB^*(V),d_{\cB,V},h_{\cB,V})=\cB^\cu(V)$ defines a natural
isomorphism of quasi-coherent CDG\+quasi-algebras
$\cA^\cu\rarrow\cB^\cu$ on~$U$.  \hbadness=1100
\end{proof}

 A \emph{quasi-coherent DG\+quasi-algebra} $\cA^\bu$ over a scheme
$X$ can be defined as quasi-coherent CDG\+quasi-algebra with
vanishing curvature elements~$h_U$ and vanishing change-of-connection
elements~$a_{VU}$.
 A \emph{morphism of quasi-coherent DG\+quasi-algebras} is a morphism
of quasi-coherent CDG\+quasi-algebras such that the domain and
the codomain are quasi-coherent DG\+quasi-algebras, and all
the change-of-connection elements~$a_{f,U}$ vanish as well.

 The category of quasi-coherent DG\+quasi-algebras on a scheme $X$
does not depend on an open covering~$\bW$.
 The category of quasi-coherent DG\+quasi-algebras over an affine
scheme $U$ is equivalent to the category of DG\+quasi-algebras over
the commutative ring~$\cO(U)$.
 We leave the proofs of these (simpler) versions of
Corollary~\ref{qcoh-cdg-quasi-algebras-not-depend-on-covering} and
Lemma~\ref{qcoh-cdg-quasi-algebras-on-affine-schemes} to the reader.

\subsection{Quasi-coherent and contraherent CDG-modules}
\label{qcoh-lcth-cdg-modules-subsecn}
 In this section and the next one we use the opportunity to provide
supplementary details to the expositions in~\cite[Section~B.1]{Pkoszul},
\cite[Section~1.2]{EP}, and~\cite[Section~2.4]{Pedg} concerning
quasi-coherent CDG\+modules.
 Then we also introduce locally contraherent CDG\+modules over
quasi-coherent CDG\+quasi-algebras over schemes.

 The following well-known definition was already mentioned in
Section~\ref{cdg-rings-cdg-modules-subsecn}.
 Let $A^*=\bigoplus_{i\in\boZ}A^i$ be a graded ring and
$M^*=\bigoplus_{i\in\boZ}M^i$ be a graded left $A^*$\+module.
 Let $d_A\:A^*\rarrow A^*$ be an odd derivation of degree~$1$.
 A homogeneous additive map $d_M\:M^*\rarrow M^*$ of degree~$1$ is
said to be an \emph{odd derivation on $M$ compatible with
the odd derivation~$d_A$ on~$A$} if one has $d_M(ax)=d_A(a)x+
(-1)^{|a|}ad_M(x)$ for all homogeneous elements $a\in A^{|a|}$
and $x\in M^{|x|}$.

\begin{lem} \label{co-induced-graded-module-derivation}
 Let $A^*$ and $B^*$ be graded rings endowed with odd derivations
$d_A\:A^*\rarrow A^*$ and $d_B\:B^*\rarrow B^*$ of degree\/~$1$.
 Let $A^*\rarrow B^*$ be a homomorphism of graded rings forming
a commutative square diagram with the derivations $d_A$ and~$d_B$.
 In this setting: \par
\textup{(a)} Let $M^*$ be a graded left $A^*$\+module endowed with
an odd derivation $d_M\:M^*\rarrow M^*$ of degree~\/~$1$ compatible
with the odd derivation~$d_A$ on~$A^*$.
 Consider the graded left $B^*$\+module $N^*=B^*\ot_{A^*}M^*$.
 Then there exists a unique odd derivation $d_N\:N^*\rarrow N^*$
compatible with the derivation~$d_B$ on $B^*$ such that the natural
map $\iota\:M^*\rarrow N^*$ forms a commutative square diagram with
the derivations $d_M$ and~$d_N$. \par
\textup{(b)} Let $P^*$ be a graded left $A^*$\+module endowed with
an odd derivation $d_P\:P^*\rarrow P^*$ of degree~\/~$1$ compatible
with the odd derivation~$d_A$ on~$A^*$.
 Consider the graded left $B^*$\+module $Q^*=\Hom_{A^*}^*(B^*,P^*)$.
 Then there exists a unique odd derivation $d_Q\:Q^*\rarrow Q^*$
compatible with the derivation~$d_B$ on $B^*$ such that the natural
map $\pi\:Q^*\rarrow P^*$ forms a commutative square diagram with
the derivations $d_Q$ and~$d_P$.
\end{lem}

\begin{proof}
 Part~(a): to prove the existence, it suffices to put
$d_N(b\ot x)=d_B(b)\ot x+(-1)^{|b|}b\ot d_M(x)$ for all
$b\ot x\in B^{|b|}\ot M^{|x|}$.
 The uniqueness follows from the fact that the graded $B^*$\+module
$N^*$ is generated by the image of the homogenenous map
$\iota\:M^*\rarrow N^*$.

 Part~(b): to prove the existence, it suffices to put
$d_Q(f)(b)=d_P(f(b))-(-1)^{|f|}f(d_B(b))$ for all $b\in B^{|b|}$
and $f\in\Hom_{A^*}^{|f|}(B^*,P^*)$.
 The uniqueness follows from the fact that the graded $B^*$\+module
$Q^*$ is \emph{cogenerated} by the homogeneous map
$\pi\:Q^*\rarrow P^*$.
 The latter statement means that, for every nonzero homogeneous
element $g\in Q^{|g|}$ there exists a homogeneous element $b\in B^{|b|}$
such that $\pi(bg)\ne0$ in~$P^*$ (as one can easily check).
 Now, for every odd derivation~$d_Q$ on $Q^*$ satisfying
the conditions of part~(b) and every element $f\in Q^{|f|}$, we have
$bd_Q(f)=(-1)^{|b|}d_Q(bf)-(-1)^{|b|}d_B(b)f$, hence the element
$\pi(bd_Q(f))=(-1)^{|b|}\pi(d_Q(bf))-(-1)^{|b|}\pi(d_B(b)f)=
(-1)^{|b|}d_P(b\pi(f))-(-1)^{|b|}d_B(b)\pi(f)\in P^*$ is uniquely
determined by $d_B$ and~$d_P$.
\end{proof}

\begin{lem} \label{co-induced-cdg-module}
 Let $A^\cu=(A^*,d_A,h_A)$ and $B^\cu=(B^*,d_B,h_B)$ be two
CDG\+rings, and let $(f,0)\:A^\cu\rarrow B^*$ be a strict morphism
of CDG\+rings.
 In this setting: \par
\textup{(a)} Let $M^\cu=(M^*,d_M)$ be a left CDG\+module over~$A^\cu$.
 Consider the graded left $B^*$\+module $N^*=B^*\ot_{A^*}M^*$.
 Then there exists a unique structure of CDG\+module $N^\cu=(N^*,d_N)$
over $B^\cu$ on the graded module $N^*$ over $B^*$ such that
the natural map $\iota\:M^*\rarrow N^*$ forms a commutative square
diagram with the derivations $d_M$ and~$d_N$. \par
\textup{(b)} Let $P^\cu=(P^*,d_P)$ be a left CDG\+module over~$A^\cu$.
 Consider the graded left $B^*$\+module $Q^*=\Hom_{A^*}^*(B^*,P^*)$.
 Then there exists a unique structure of CDG\+module $Q^\cu=(Q^*,d_Q)$
over $B^\cu$ on the graded module $Q^*$ over $B^*$ such that
the natural map $\pi\:Q^*\rarrow P^*$ forms a commutative square
diagram with the derivations $d_Q$ and~$d_P$.
\end{lem}

\begin{proof}
 Part~(a): the uniqueness follows immediately from
Lemma~\ref{co-induced-graded-module-derivation}(a).
 To prove the existence, notice that $(B^*,d_B)$ is a CDG\+bimodule
over the CDG\+rings $B^\cu$ and $A^\cu$, while $(M^*,d_M)$ is
a CDG\+bimodule over the CDG\+rings $A^\cu$ and $(\boZ,0,0)$.
 The construction of the tensor product of CDG\+bimodules from
Section~\ref{cdg-rings-cdg-modules-subsecn} provides the desired
CDG\+module structure $N^\cu=B^\cu\ot_{A^*}M^\cu$ on the tensor product
$N^*=B^*\ot_{A^*}\nobreak M^*$.

 Part~(b): the uniqueness follows immediately from
Lemma~\ref{co-induced-graded-module-derivation}(b).
 To prove the existence, notice that $(B^*,d_B)$ is a CDG\+bimodule
over the CDG\+rings $A^\cu$ and~$B^\cu$.
 The construction of the $\Hom$ CDG\+module from
Section~\ref{cdg-rings-cdg-modules-subsecn} provides the desired
CDG\+module structure $Q^\cu=\Hom_{A^*}^\cu(B^\cu,M^\cu)$ on
the graded $B^*$\+module $Q^*=\Hom_{A^*}^*(B^*,P^*)$.
\end{proof}

\begin{cor} \label{cdg-quasi-algebra-co-extension-of-scalars}
 Let $R$ be a commutative ring and $A^\cu=(A^*,d_A,h_A)$ be
a CDG\+quasi-algebra over~$R$.
 Let $f\:R\rarrow S$ be a homomorphism of commutative rings such that
the related morphism of affine schemes\/ $\Spec S\rarrow\Spec R$ is
an open immersion.
 Let $B^\cu=(B^*,d_B,h_B)$ be the induced CDG\+quasi-algebra
structure on the graded quasi-algebra $B^*=S\ot_RA^*$ over $S$,
as in Lemma~\ref{localizing-cdg-ring-by-affine-open-immersion}.
 In this setting: \par
\textup{(a)} Let $M^\cu=(M^*,d_M)$ be a left CDG\+module over~$A^\cu$.
 Consider the graded left $B^*$\+module $N^*=S\ot_RM^*$, as per
the graded version of
Lemma~\ref{quasi-algebra-co-extension-of-scalars}(b).
 Then there exists a unique structure of CDG\+module $N^\cu=(N^*,d_N)$
over $B^\cu$ on the graded module $N^*$ over $B^*$ such that
the natural map $\iota\:M^*\rarrow N^*$ forms a commutative square
diagram with the derivations $d_M$ and~$d_N$. \par
\textup{(b)} Let $P^\cu=(P^*,d_M)$ be a left CDG\+module over~$A^\cu$.
 Consider the graded left $B^*$\+module $Q^*=\Hom^*_R(S,P^*)$, as per
the graded version of
Lemma~\ref{quasi-algebra-co-extension-of-scalars}(c).
 Then there exists a unique structure of CDG\+module $Q^\cu=(Q^*,d_Q)$
over $B^\cu$ on the graded module $Q^*$ over $B^*$ such that
the natural map $\pi\:Q^*\rarrow P^*$ forms a commutative square
diagram with the derivations $d_Q$ and~$d_P$.
\end{cor}

\begin{proof}
 In view of the alternative proofs of
Lemma~\ref{quasi-algebra-co-extension-of-scalars}(b\+-c),
this is a special case of Lemma~\ref{co-induced-cdg-module}.
\end{proof}

\begin{lem} \label{strict-gluing-module-derivation-on-affine}
 Let $R$ be a commutative ring and $A^\cu=(A^*,d_A,h_A)$ be
a CDG\+quasi-algebra over~$R$.
 Let $R\rarrow S_\alpha$, \ $1\le\alpha\le k$, be a finite collection
of homomorphisms of commutative rings such that the collection of
induced maps of the spectra\/ $\Spec S_\alpha\rarrow\Spec R$ is
an affine open covering of an affine scheme.
 Consider the induced CDG\+quasi-algebra structures
$B_\alpha^\cu=(B_\alpha^*,d_\alpha,h_\alpha)$ and $C_{\alpha\beta}^\cu
=(C_{\alpha\beta}^*,d_{\alpha\beta},h_{\alpha\beta})$ on the graded
quasi-algebras $B_\alpha^*=S_\alpha\ot_RA^*$ and
$S_{\alpha\beta}^*=S_\alpha\ot_RS_\beta\ot_RA^*$ over the rings
$S_\alpha$ and $S_\alpha\ot_RS_\beta$, as per 
Lemma~\ref{localizing-cdg-ring-by-affine-open-immersion}.
 In this setting: \par
\textup{(a)} Let $M^*$ be a graded left $A^*$\+module.
 Suppose given, for every index~$\alpha$, a CDG\+module structure
$(N^\cu_\alpha,d_{N,\alpha})$ over the CDG\+ring $B_\alpha^\cu$ on
the graded left $B_\alpha^*$\+module $N_\alpha^*=S_\alpha\ot_RM^*$.
 Assume that, for every pair of indices $\alpha$ and~$\beta$,
the two CDG\+module structures over the CDG\+ring $C_{\alpha\beta}^\cu$
on the graded $C_{\alpha\beta}^*$\+module $(S_\alpha\ot_RS_\beta)
\ot_{S_\alpha}N_\alpha^*\simeq(S_\alpha\ot_RS_\beta)\ot_{S_\beta}
N_\beta^*$ arising from the CDG\+module structures $N^\cu_\alpha$
and $N^\cu_\beta$ on the graded modules $N^*_\alpha$ and $N^*_\beta$,
as per the construction of
Corollary~\ref{cdg-quasi-algebra-co-extension-of-scalars}(a),
coincide with each other.
 Then there exists a unique CDG\+module structure $M^\cu=(M^*,d_M)$
over the CDG\+ring $A^\cu$ on the graded $A^*$\+module $M^*$
such that, for every index~$\alpha$, the CDG\+module structure
$N_\alpha^\cu$ over $B_\alpha^\cu$ arises from the CDG\+module
structure $M^\cu$ over $A^\cu$ via the construction of
Corollary~\ref{cdg-quasi-algebra-co-extension-of-scalars}(a). \par
\textup{(b)} Let $P^*$ be a graded left $A^*$\+module such that,
for every degree $i\in\boZ$, the $R$\+module $P^i$ is contraadjusted.
 Suppose given, for every index~$\alpha$, a CDG\+module structure
$(Q^\cu_\alpha,d_{Q,\alpha})$ over the CDG\+ring $B_\alpha^\cu$ on
the graded left $B_\alpha^*$\+module
$Q_\alpha^*=\Hom^*_R(S_\alpha,P^*)$.
 Assume that, for every pair of indices $\alpha$ and~$\beta$,
the two CDG\+module structures over the CDG\+ring $C_{\alpha\beta}^\cu$
on the graded $C_{\alpha\beta}^*$\+module\/ $\Hom_{S_\alpha}^*
(S_\alpha\ot_RS_\beta,\>Q_\alpha^*)\simeq\Hom_{S_\beta}^*
(S_\alpha\ot_RS_\beta,\>Q_\beta^*)$ arising from the CDG\+module
structures $Q^\cu_\alpha$ and $Q^\cu_\beta$ on the graded modules
$Q^*_\alpha$ and $Q^*_\beta$, as per the construction of
Corollary~\ref{cdg-quasi-algebra-co-extension-of-scalars}(b),
coincide with each other.
 Then there exists a unique CDG\+module structure $P^\cu=(P^*,d_M)$
over the CDG\+ring $A^\cu$ on the graded $A^*$\+module $P^*$
such that, for every index~$\alpha$, the CDG\+module structure
$Q_\alpha^\cu$ over $B_\alpha^\cu$ arises from the CDG\+module
structure $P^\cu$ over $A^\cu$ via the construction of
Corollary~\ref{cdg-quasi-algebra-co-extension-of-scalars}(b).
\end{lem}

\begin{proof}
 The proof is similar to that of
Lemma~\ref{strict-gluing-derivation-on-affine}.
 In part~(a), the argument is based on exactness of the \v Cech
coresolution~\eqref{module-over-quasi-algebra-cech-coresolution}
\begin{multline} \label{module-cech-coresolution-frament}
 0\lrarrow M^*\lrarrow\bigoplus\nolimits_{\alpha=1}^k
 S_\alpha\ot_RM^* \\ \lrarrow\bigoplus\nolimits_{1\le\alpha<\beta\le k}
 S_\alpha\ot_RS_\beta\ot_RM^*\lrarrow\dotsb
\end{multline}
 One endows the middle term 
of~\eqref{module-cech-coresolution-frament} with a ``vertical''
differential constructed as the direct sum of
the derivations~$d_{N,\alpha}$.
 Then one endows all the terms
of~\eqref{module-cech-coresolution-frament} going further to the right
with the vertical differentials constructed as the direct sums of
the induced derivations on the tensor products, produced by
the construction of
Corollary~\ref{cdg-quasi-algebra-co-extension-of-scalars}(a).
 The vertical differentials form commutative square diagrams with
the horizontal (\v Cech) differentials, and it remains to construct
the desired differential $d_M\:M^*\rarrow M^*$ by passing to the kernel.
 The pair $(M^*,d_M)$ is a CDG\+module over $A^\cu$ because
the direct sum $\bigoplus_{\alpha=1}^k N_\alpha^\cu$ is a CDG\+module
over $A^\cu$ and the map $M^*\rarrow\bigoplus_{\alpha=1}^k N_\alpha^*
=\bigoplus_{\alpha=1}^k S_\alpha\ot_RM^*$ is injective.

 The proof of part~(b) is dual-analogous and based on exactness of
the \v Cech resolution~\eqref{module-over-quasi-algebra-cech-resolution}
\begin{multline} \label{module-cech-resolution-fragment}
 \dotsb\lrarrow\bigoplus\nolimits_{1\le\alpha<\beta\le k}
 \Hom_R^*(S_\alpha\ot_RS_\beta,\>P^*) \\
 \lrarrow\bigoplus\nolimits_{\alpha=1}^k
 \Hom_R^*(S_\alpha,P^*)\lrarrow P^*\lrarrow0
\end{multline}
(see~\cite[formula~(1.3) in Lemma~1.2.6(b)]{Pcosh}).
 The exactness of~\eqref{module-cech-resolution-fragment} depends on
the assumption that the grading components of $P^*$ are contraadjusted
$R$\+modules; that is where this assumption is used.
\end{proof}

 Let $(f,a)\:B^\cu\rarrow A^\cu$ be a morphism of CDG\+rings.
 Let $M^\cu=(M^*,d_M)$ be a left CDG\+module over $A^\bu$, and let
$N^\cu=(N^*,d_N)$ be a left CDG\+module over $B^\cu$.
 We will say that a homogeneous additive map $g\:N^*\rarrow M^*$
of degree~$0$ is a \emph{map of CDG\+modules compatible with
the morphism of CDG\+rings $(f,a)$} if $g$~is a morphism of
graded $B^*$\+modules and the map~$g$ forms a commutative
square diagram with the differential~$d_N$ on~$N^*$ and
the differential~$d'_M$ on $M^*$ given by the formula~(vi) from
Section~\ref{cdg-rings-cdg-modules-subsecn}.
 Similarly, a homogeneous additive map $g\:M^*\rarrow N^*$ of degree~$0$
is a \emph{map of CDG\+modules compatible with the morphism of
CDG\+rings $(f,a)$} if $g$~is a morphism of graded $B^*$\+modules
and the map~$g$ forms a commutative square diagram with
the differential~$d_N$ on~$N^*$ and the differential~$d'_M$ on $M^*$
given by the same formula~(vi).
 In the context of right CDG\+modules, the notion of a map of
CDG\+modules compatible with a morphism of CDG\+rings is defined
similarly using formula~(viii) from
Section~\ref{cdg-rings-cdg-modules-subsecn}.

 Let $X$ be a scheme with an open covering $\bW$ and $\cA^\cu$ be
a quasi-coherent CDG\+quasi-algebra over~$X$ (as defined in
Section~\ref{qcoh-cdg-quasi-algebras-subsecn}).
 A \emph{quasi-coherent left CDG\+module} $\M^\cu$ over $\cA^\cu$
is a set of data consisting of
\begin{itemize}
\item a quasi-coherent graded left module $\M^*=\bigoplus_{i\in\boZ}
\M^i$ over the quasi-coherent graded quasi-algebra $\cA^*$ over~$X$
(see Sections~\ref{cosheaves-of-A-modules-subsecn}
and~\ref{graded-modules-subsecn} for the definition);
\item an odd derivation $d_{\M,U}\:\M^*(U)\rarrow\M^*(U)$ of
degree~$1$ on the graded $\cA^*(U)$\+module $\M^*(U)$ compatible
with the odd derivation $d_{\cA,U}\:\cA^*(U)\rarrow\cA^*(U)$
on the graded ring $\cA^*(U)$, given for each affine open
subscheme $U\subset X$ subordinate to~$\bW$.
\end{itemize}
 The following two axioms must be satisfied:
\begin{enumerate}
\renewcommand{\theenumi}{\roman{enumi}}
\setcounter{enumi}{15}
\item for each affine open subscheme $U\subset X$ subordinate to $\bW$,
the pair $\M^\cu(U)=(\M^*(U),d_{\M,U})$ is a left CDG\+module over
the CDG\+ring $\cA^\cu(U)=(\cA^*(U),\allowbreak d_{\cA,U},h_{\cA,U})$;
\item for each pair of affine open subschemes $V\subset U\subset X$
subordinate to $\bW$, the map of restriction of sections
$\M^*(U)\rarrow\M^*(V)$ in the quasi-coherent sheaf $\M^*$ is a map
of CDG\+modules compatible with the morphism of CDG\+rings
$(\rho_{\cA,VU},a_{\cA,VU})\:(\cA^*(U),d_{\cA,U},h_{\cA,U})\rarrow
(\cA^*(V),d_{\cA,V},h_{\cA,V})$.
\end{enumerate}
 The definition of a \emph{quasi-coherent right CDG\+module} $\N^\cu$
over $\cA^\cu$ is similar.

 Let $\cB^\cu\rarrow\cA^\cu$ be a morphism of quasi-coherent
CDG\+quasi-algebras over~$X$ (as defined in
Section~\ref{qcoh-cdg-quasi-algebras-subsecn}), and let
$\cB^*\rarrow\cA^*$ be its underlying morphism of quasi-coherent
graded quasi-algebras.
  Given a quasi-coherent left CDG\+module $\M^\cu$ over $\cA^\cu$,
the morphism of quasi-coherent graded quasi-algebras $\cB^*\rarrow\cA^*$
allows one to endow the quasi-coherent graded left $\cA^*$\+module
$\M^*$ with the structure of a quasi-coherent graded left
$\cB^*$\+module.
 In this context, the quasi-coherent graded left $\cB^*$\+module $\M^*$
acquires a natural structure of quasi-coherent left CDG\+module over
$\cB^*$, given by differentials $d'_{\M,U}\:\M^*(U)\rarrow\M^*(U)$
obtained from the differentials~$d_{\M,U}$ using the same formula~(vi)
from Section~\ref{cdg-rings-cdg-modules-subsecn}.
 We leave details to the reader.

\begin{cor} \label{qcoh-cdg-module-not-dependent-on-covering}
 The notion of a quasi-coherent CDG\+module over a quasi-coherent
CDG\+quasi-algebra $\cA^\cu$ over a scheme $X$ does not depend on
an open covering\/~$\bW$.
 Specifically, given a quasi-coherent graded left $\cA^*$\+module
$\M^*$, any structure of a quasi-coherent CDG\+module over $\cA^\cu$
on $\M^*$ defined for the open covering\/ $\bW$ can be uniquely extended
to a structure of quasi-coherent CDG\+module over $\cA^\cu$ on $\M^*$
defined for the open covering\/ $\{X\}$ of the scheme~$X$.
\end{cor}

\begin{proof}
 It suffices to consider the case of an affine scheme $X=U$ with
some open covering~$\bW$.
 In this case, following (the proof of)
Lemma~\ref{qcoh-cdg-quasi-algebras-on-affine-schemes}, one can
assume that all the restriction morphisms of CDG\+rings in
the quasi-coherent CDG\+quasi-algebra $\cA^\cu$ on $U$ are strict,
i.~e., the change-of-connection elements $a_{\cA,TV}\in\cA^1(T)$
vanish for all affine open subschemes $T\subset U\subset X$.
 In this case, it remains to pick a finite affine open covering
$U=\bigcup_{\alpha=1}^k V_\alpha$ of the affine scheme $U$ subordinate
to~$\bW$.
 Then the desired assertion is provided by
Lemma~\ref{strict-gluing-module-derivation-on-affine}(a).
\end{proof}

 Let $X$ be a scheme with an open covering $\bW$ and $\cA^\cu$ be
a quasi-coherent CDG\+quasi-algebra over~$X$.
 A \emph{$\bW$\+locally contraherent CDG\+module} $\P^\cu$ over
$\cA^\cu$ is a set of data consisting of
\begin{itemize}
\item a $\bW$\+locally contraherent graded left module
$\prod_{i\in\boZ}\P^i$ over the quasi-coherent graded quasi-algebra
$\cA^*$ over~$X$ (see Sections~\ref{cosheaves-of-A-modules-subsecn}
and~\ref{graded-modules-subsecn} for the definition);
\item an odd derivation $d_{\P,U}\:\P^*[U]\rarrow\P^*[U]$ of
degree~$1$ on the graded $\cA^*(U)$\+module $\P^*[U]$ compatible with
the odd derivation $d_{\cA,U}\:\cA^*(U)\rarrow\cA^*(U)$ on the graded
ring $\cA^*(U)$, defined for each affine open subscheme $U\subset X$
subordinate to~$\bW$.
\end{itemize}
 The following two axioms must be satisfied:
\begin{enumerate}
\renewcommand{\theenumi}{\roman{enumi}}
\setcounter{enumi}{17}
\item for each affine open subscheme $U\subset X$ subordinate to $\bW$,
the pair $\P^\cu(U)=(\P^*(U),d_{\P,U})$ is a left CDG\+module over
the CDG\+ring $\cA^\cu(U)=(\cA^*(U),\allowbreak d_{\cA,U},h_{\cA,U})$;
\item for each pair of affine open subschemes $V\subset U\subset X$
subordinate to $\bW$, the map of corestriction of cosections
$\P^*[V]\rarrow\P^*[U]$ in the locally contraherent cosheaf $\P^*$ is
a map of CDG\+modules compatible with the morphism of CDG\+rings
$(\rho_{\cA,VU},a_{\cA,VU})\:(\cA^*(U),d_{\cA,U},h_{\cA,U})\rarrow
(\cA^*(V),d_{\cA,V},h_{\cA,V})$.
\end{enumerate}

 Let $\cB^\cu\rarrow\cA^\cu$ be a morphism of quasi-coherent
CDG\+quasi-algebras over~$X$ and let $\cB^*\rarrow\cA^*$ be its
underlying morphism of quasi-coherent graded quasi-algebras.
  Given a $\bW$\+locally contraherent CDG\+module $\P^\cu$ over
$\cA^\cu$, the morphism of quasi-coherent graded quasi-algebras
$\cB^*\rarrow\cA^*$ allows one to endow the locally contraherent
graded $\cA^*$\+module $\P^*$ with the structure of a locally
contraherent graded $\cB^*$\+module.
 In this context, the $\bW$\+locally contraherent graded left
$\cB^*$\+module $\M^*$ acquires a natural structure of
a $\bW$\+locally contraherent CDG\+module over $\cB^*$, given by
differentials $d'_{\P,U}\:\P^*[U]\rarrow\P^*[U]$ obtained from
the differentials~$d_{\P,U}$ using formula~(vi) from
Section~\ref{cdg-rings-cdg-modules-subsecn}.
 We leave details to the reader.

\begin{cor} \label{lcth-cdg-module-not-dependent-on-covering}
 The notion of a locally contraherent CDG\+module structure over
a quasi-coherent CDG\+quasi-algebra $\cA^\cu$ over a scheme $X$ on
a given locally contraherent graded $\cA^*$\+module\/ $\P^*$ does not
depend on an open covering\/~$\bW$.
 Specifically, let\/ $\bW'$ be an open covering of $X$ such that
the open covering\/ $\bW$ is subordinate to\/~$\bW$'.
 Given a\/ $\bW'$\+locally contraherent graded $\cA^*$\+module\/ $\P^*$,
any structure of a\/ $\bW$\+locally contraherent CDG\+module over
$\cA^\cu$ on\/ $\P^*$ can be uniquely extended to a structure of\/
$\bW'$\+locally contraherent CDG\+module over $\cA^\cu$ on\/~$\P^*$.
\end{cor}

\begin{proof}
 This is dual-analogous to
Corollary~\ref{qcoh-cdg-module-not-dependent-on-covering}.
 It suffices to consider the case of an affine scheme $X=U$ with
the open covering $\bW'=\{U\}$ and some open covering~$\bW$.
 In this case, following (the proof of)
Lemma~\ref{qcoh-cdg-quasi-algebras-on-affine-schemes}, one can
assume that all the restriction morphisms of CDG\+rings in
the quasi-coherent CDG\+quasi-algebra $\cA^\cu$ on $U$ are strict,
i.~e., the change-of-connection elements $a_{\cA,TV}\in\cA^1(T)$
vanish for all affine open subschemes $T\subset U\subset X$.
 Then the desired assertion is provided by
Lemma~\ref{strict-gluing-module-derivation-on-affine}(b).
\end{proof}

\subsection{DG-categories of quasi-coherent and locally contraherent
CDG-mod\protect\-ules}  \label{dg-categories-of-cdg-modules-subsecn}
 We refer to~\cite[Section~1.2]{Pkoszul}, \cite[Section~1]{Pedg},
or~\cite[Section~1]{PS5} for a background discussion of DG\+categories
suitable for our context.

 Let us emphasize that, following the approach of~\cite{Pedg,PS5},
we take the \emph{strict point of view} on DG\+categories and
DG\+functors.
 This means that the complexes of morphisms in DG\+categories are
considered up to isomorphism and \emph{not} up to quasi-isomorphism.
 So we work with DG\+categories up to equivalences of DG\+categories
and \emph{not} up to quasi-equivalences.
 See~\cite[Introduction]{Pedg} or~\cite[end of Section~0.4]{PS5}.

 We will use the notation similar to that in~\cite{Pedg,PS5}.
 In particular, for any DG\+category $\bE$, we denote by $\sZ^0(\bE)$
the preadditive category of closed morphisms of degree~$0$ in $\bE$,
and by $\sH^0(\bE)$ the preadditive category of degree~$0$ cohomology
of~$\bE$ (i.~e., closed morphisms of degree~$0$ up to homotopy).
 For any DG\+category $\bE$ with a zero object, shifts, and cones,
the homotopy category $\sH^0(\bE)$ is naturally a triangulated category.
 See~\cite[Section~1.3]{Pedg} or~\cite[Sections~1.3\+-1.5]{PS5}.

 Let $B^\cu=(B^*,d,h)$ be a CDG\+ring.
 For any two left CDG\+modules $L^\cu=(L^*,d_L)$ and $M^\cu=(M^*,d_M)$
over $B^\cu$, the complex of abelian groups
$\Hom_{B^*}^\bu(L^\cu,M^\cu)$ was constructed in
Section~\ref{cdg-rings-cdg-modules-subsecn}.
 Given three left CDG\+modules $K^\cu$, $L^\cu$, and $M^\cu$ over
$B^\cu$, the obvious map of composition of homogeneous maps of
graded $B^*$\+modules
$$
 \Hom_{B^*}^*(L^*,M^*)\ot_\boZ\Hom_{B^*}^*(K^*,L^*)\lrarrow
 \Hom_{B^*}^*(K^*,M^*)
$$
is a morphism of complexes of abelian groups
$$
 \Hom_{B^*}^\bu(L^\cu,M^\cu)\ot_\boZ\Hom_{B^*}^\bu(K^\cu,L^\cu)\lrarrow
 \Hom_{B^*}^\bu(K^\cu,M^\cu).
$$
 This construction defines the \emph{DG\+category of left CDG\+modules
over~$B^\cu$}, which we will denote by $B^\cu\bModl$.

 Similarly, for any two right CDG\+modules $R^\cu=(R^*,d_R)$ and
$N^\cu=(N^*,d_N)$, the construction from
Section~\ref{cdg-rings-cdg-modules-subsecn} produces a complex of
abelian groups $\Hom_{B^*{}^\rop}^\bu(R^\cu,N^\cu)$.
 Given three right CDG\+modules $R^\cu$, $N^\cu$, and $K^\cu$ over
$B^\cu$, the composition of homogeneous maps of graded $B^*$\+modules 
provides a composition morphism of complexes of abelian groups
$$
 \Hom_{B^*{}^\rop}^\bu(N^\cu,K^\cu)\ot_\boZ
 \Hom_{B^*{}^\rop}^\bu(R^\cu,N^\cu) 
 \lrarrow\Hom_{B^*{}^\rop}^\bu(R^\cu,K^\cu).
$$
 This defines the \emph{DG\+category of right CDG\+modules
over~$B^\cu$}, which we will denote by $\bModr B^\cu$.
 For any DG\+ring $B^\cu$, the DG\+categories $B^\cu\bModl$ and
$\bModr B^\cu$ have zero objects, shifts, and cones; in fact, they
even have twists, infinite direct sums, and infinite products.

 Let $X$ be a scheme and $\cB^*$ be a quasi-coherent graded
quasi-algebra over~$X$.
 Let $\cL^*$ and $\M^*$ be two quasi-coherent graded left
$\cB^*$\+modules.
 We denote by $\Hom_{\cB^*}^*(\cL^*,\M^*)$ the graded abelian group
of homogeneous maps of quasi-coherent graded $\cB^*$\+modules $\cL^*
\rarrow\M^*$ of various degrees $n\in\boZ$.
 The sign rule for compatibility with the action of $\cB^*$ is
as written down in Section~\ref{cdg-rings-cdg-modules-subsecn}.
 The notation $\Hom_{\cB^*{}^\rop}^*(\cR^*,\N^*)$ for two quasi-coherent
graded right $\cB^*$\+modules $\cR^*$ and $\N^*$ has similar meaning.

 Similarly, for any locally contraherent graded $\cB^*$\+modules
$\P$ and $\Q$, we denote by $\Hom^{\cB^*,*}(\P^*,\Q^*)$ the graded
abelian group of homogeneous maps of locally contraherent graded
$\cB^*$\+modules $\P^*\rarrow\Q^*$ of various degrees $n\in\boZ$.
 The sign rule for compatibility with the action of $\cB^*$ is
as in Section~\ref{cdg-rings-cdg-modules-subsecn}.

 Let $X$ be a scheme with an open covering $\bW$ and $\cB^\cu$ be
a quasi-coherent CDG\+quasi-algebra over~$X$.
 Let $\cL^\cu$ and $\M^\cu$ be two quasi-coherent left CDG\+modules
over~$\cB^\cu$ (see Section~\ref{qcoh-lcth-cdg-modules-subsecn}
for the definitions).
 Then the differential~$d$ on the graded abelian group
$\Hom_{\cB^*}^*(\cL^*,\M^*)$ is constructed as follows.

 Suppose given a homogeneous map of quasi-coherent graded
$\cB^*$\+modules $f\in\Hom_{\cB^*}^{|f|}(\cL^*,\M^*)$.
 For every affine open subscheme $U\subset X$ subordinate to $\bW$,
we have two left CDG\+modules $(\cL^*(U),d_{\cL,U})$ and
$(\M^*(U),d_{\M,U})$ over the CDG\+ring
$(\cB^*(U),d_{\cB,U},h_{\cB,U})$.
 We also have a homogeneous map of graded left $\cB^*(U)$\+modules
$f(U)\:\cL^*(U)\rarrow\M^*(U)$ of degree~$|f|$.

 Based on these data, the construction of the differential on
the complex of morphisms between two left CDG\+modules, as spelled
out in Section~\ref{cdg-rings-cdg-modules-subsecn}, provides
a homogeneous map of graded left $\cB^*(U)$\+modules
$d(f(U))\:\cL^*(U)\rarrow\M^*(U)$ of degree~$|f|+1$.
 The homogeneous map of quasi-coherent graded $\cB^*$\+modules
$d(f)\in\Hom_{\cB^*}^{|f|+1}(\cL^*,\M^*)$ is defined by the rule
$d(f)(U)=d(f(U))\:\cL^*(U)\rarrow\M^*(U)$ for all affine open
subschemes $U\subset X$ subordinate to~$\bW$.

 It is straightforward to check, using formula~(vi) from
Section~\ref{cdg-rings-cdg-modules-subsecn} and axiom~(xvii)
from Section~\ref{qcoh-lcth-cdg-modules-subsecn},
that the maps $d(f(U))$ and $d(f(V))$ form commutative square diagrams
with the restriction maps $\cL^*(U)\rarrow\cL^*(V)$ and
$\M^*(U)\rarrow\M^*(V)$ for all affine open subschemes
$V\subset U\subset X$ subordinate to~$\bW$.
 So the homogeneous map of quasi-coherent graded $\cB^*$\+modules
$d(f)\:\cL^*\rarrow\M^*$ is well-defined.
 One can see from axiom~(xvi) that the resulting differential~$d$ on
the graded abelian group $\Hom_{\cB^*}^*(\cL^*,\M^*)$ makes it
a complex, which we denote by $\Hom_{\cB^*}^\bu(\cL^\cu,\M^\cu)$.
 Obviously, this construction does not depend on the choice of
an open covering $\bW$ of the scheme~$X$.

 Given three quasi-coherent left CDG\+modules $\K^\cu$, $\cL^\cu$,
and $\M^\cu$ over $\cB^\cu$, the operation of composition of
homogeneous maps of quasi-coherent graded $\cB^*$\+modules provides
the composition morphism of complexes of abelian groups
$$
 \Hom_{\cB^*}^\bu(\cL^\cu,\M^\cu)\ot_\boZ
 \Hom_{\cB^*}^\bu(\K^\cu,\cL^\cu)\lrarrow
 \Hom_{\cB^*}^\bu(\K^\cu,\M^\cu).
$$
 This defines the \emph{DG\+category of quasi-coherent left
CDG\+modules over~$\cB^\cu$}, which we will denote by $\cB^\cu\bQcoh$.
 The \emph{DG\+category of quasi-coherent right CDG\+modules
over~$\cB^\cu$} is defined similarly and denoted by $\bQcohr\cB^\cu$.

 For any quasi-coherent CDG\+quasi-algebra $\cB^\cu$ over
a scheme $X$, the DG\+categories $\cB^\cu\bQcoh$ and $\bQcohr\cB^\cu$
have zero objects, shifts, and cones; in fact, they even have twists
and infinite direct sums.

\begin{lem} \label{qcoh-cdg-modules-on-affine-schemes}
 Let $U=\Spec R$ be an affine scheme, $\cB^\cu$ be a quasi-coherent
CDG\+quasi-algebra over $U$, and $\cB^\cu(U)=
(\cB^*(U),d_{\cB,U},h_{\cB,U})$ be the corresponding CDG\+quasi-algebra
over $R$, as per Lemma~\ref{qcoh-cdg-quasi-algebras-on-affine-schemes}.
 Then the rule assigning the CDG\+module $\M^\cu(U)=(\M^*(U),d_{\M,U})$
over $\cB^\cu(U)$ to a quasi-coherent CDG\+module $\M^\cu$ over
$\cB^\cu$ is a natural equivalence between the DG\+category of
quasi-coherent CDG\+modules over $\cB^\cu$ and the DG\+category of
CDG\+modules over~$\cB^\cu(U)$,
\begin{equation} \label{qcoh-cdg-modules-on-affine-scheme-formula}
 \cB^\cu\bQcoh\simeq\cB^\cu(U)\bModl.
\end{equation}
\end{lem}

\begin{proof}
 This is essentially the same argument as the proof of
Corollary~\ref{qcoh-cdg-module-not-dependent-on-covering}.
 According to (the proof of) 
Lemma~\ref{qcoh-cdg-quasi-algebras-on-affine-schemes}, one can assume
that all the restriction morphisms of CDG\+rings in
the quasi-coherent CDG\+quasi-algebra $\cB^\cu$ on $U$ are strict.
 Then the assertion follows from the equivalence $\cB^*\Qcoh\simeq
\cB^*(U)\Modl$ on the level of the categories of graded modules together
with Corollary~\ref{cdg-quasi-algebra-co-extension-of-scalars}(a).
\end{proof}

 The locally contraherent version of the constructions above is similar.
 Let $X$ be a scheme with an open covering $\bW$ and $\cB^\cu$ be
a quasi-coherent CDG\+quasi-algebra over~$X$.
 Let $\P^\cu$ and $\Q^\cu$ be two $\bW$\+locally contraherent
CDG\+modules over~$\cB^\cu$ (see
Section~\ref{qcoh-lcth-cdg-modules-subsecn} for the definitions).
 Then the differential~$d$ on the graded abelian group
$\Hom^{\cB^*,*}(\P^*,\Q^*)$ is constructed as follows.

 Suppose given a homogeneous map of locally contraherent graded
$\cB^*$\+modules $f\in\Hom^{\cB^*,|f|}(\P^*,\Q^*)$.
 For every affine open subscheme $U\subset X$ subordinate to $\bW$,
we have two left CDG\+modules $(\P^*[U],d_{\P,U})$ and
$(\Q^*[U],d_{\Q,U})$ over the CDG\+ring
$(\cB^*(U),d_{\cB,U},h_{\cB,U})$.
 We also have a homogeneous map of graded left $\cB^*(U)$\+modules
$f[U]\:\P^*[U]\rarrow\Q^*[U]$ of degree~$|f|$.

 Based on these inputs, the construction of the differential on
the complex of morphisms between two left CDG\+modules, as spelled
out in Section~\ref{cdg-rings-cdg-modules-subsecn}, provides
a homogeneous map of graded left $\cB^*(U)$\+modules
$d(f[U])\:\P^*[U]\rarrow\Q^*[U]$ of degree~$|f|+1$.
 The homogeneous map of locally contraherent graded $\cB^*$\+modules
$d(f)\in\Hom^{\cB^*,|f|+1}(\P^*,\Q^*)$ is defined by the rule
$d(f)[U]=d(f[U])\:\P^*[U]\rarrow\Q^*[U]$ for all affine open
subschemes $U\subset X$ subordinate to~$\bW$.

 It is straightforward to check, using formula~(vi) from
Section~\ref{cdg-rings-cdg-modules-subsecn} and axiom~(xix)
from Section~\ref{qcoh-lcth-cdg-modules-subsecn},
that the maps $d(f[V])$ and $d(f[U])$ form commutative square diagrams
with the corestriction maps $\P^*[V]\rarrow\P^*[U]$ and
$\Q^*[V]\rarrow\Q^*[U]$ for all affine open subschemes
$V\subset U\subset X$ subordinate to~$\bW$.
 So the homogeneous map of locally contraherent graded $\cB^*$\+modules
$d(f)\:\P^*\rarrow\Q^*$ is well-defined.
 One can see from axiom~(xviii) that the resulting differential~$d$ on
the graded abelian group $\Hom^{\cB^*,*}(\P^*,\Q^*)$ makes it a complex,
which we denote by $\Hom^{\cB^*,\bu}(\P^\cu,\Q^\cu)$.
 Obviously, this construction does not depend on the choice of
an open covering $\bW$ of the scheme $X$ such that both $\P^\cu$ and
$\Q^\cu$ are $\bW$\+locally contraherent CDG\+modules (i.~e., both
$\P^i$ and $\Q^i$ are $\bW$\+locally contraherent cosheaves on $X$
for all degrees $i\in\boZ$).

 Given three $\bW$\+locally contraherent left CDG\+modules
$\gK^\cu$, $\P^\cu$, and $\Q^\cu$ over $\cB^\cu$, the operation
of composition of homogeneous maps of locally contraherent graded
$\cB^*$\+modules provides the composition morphism of complexes
of abelian groups
$$
 \Hom^{\cB^*,\bu}(\P^\cu,\Q^\cu)\ot_\boZ
 \Hom^{\cB^*,\bu}(\gK^\cu,\P^\cu)\lrarrow
 \Hom^{\cB^*,\bu}(\gK^\cu,\Q^\cu).
$$
 This defines the \emph{DG\+category of\/ $\bW$\+locally contraherent
CDG\+modules over~$\cB^\cu$}, which we will denote by
$\cB^\cu\bLcth_\bW$.
 As usual, we put $\cB^\cu\bCtrh=\cB^\cu\bLcth_{\{X\}}$, and
call the objects of $\cB^\cu\bCtrh$ the \emph{contraherent CDG\+modules
over~$\cB^\cu$}.
 We also put $\cB^\cu\bLcth=\bigcup_\bW\cB^\cu\bLcth_\bW$.

 A $\bW$\+locally contraherent CDG\+module $\P^\cu$ over $\cB^\cu$ is
said to be \emph{$X$\+locally cotorsion}, \emph{$X$\+locally injective},
\emph{$\cB^*$\+locally cotorsion}, or \emph{$\cB^*$\+locally injective}
if its underlying graded $\cB^*$\+module $\P^*$ has the respective
property (defined in Section~\ref{graded-modules-subsecn}).
 We denote the full DG\+subcategory of $X$\+locally cotorsion
$\bW$\+locally contraherent CDG\+modules by $\cB^\cu\bLcth_\bW^{X\dlct}
\subset\cB^\cu\bLcth_\bW$ and the full DG\+subcategory of
$\cB^*$\+locally cotorsion $\bW$\+locally contraherent CDG\+modules by
$\cB^\cu\bLcth_\bW^{\cB^*\dlct}\subset\cB^\cu\bLcth_\bW$.
 The notation $\cB^\cu\bCtrh^{X\dlct}$ and $\cB^\cu\bCtrh^{\cB^*\dlct}$
has similar meaning.

 For any quasi-coherent CDG\+quasi-algebra $\cB^\cu$ over a scheme $X$
and any open covering $\bW$ of $X$, the DG\+category $\cB^\cu\bLcth_\bW$
has a zero object, shifts, and cones; in fact, it even has twists
and infinite direct products.
 The full DG\+subcategories $\cB^\cu\bLcth_\bW^{\cB^*\dlct}\subset
\cB^\cu\bLcth_\bW^{X\dlct}\subset\cB^\cu\bLcth_\bW$ are closed under
shifts, twists, and direct products (hence also under cones).

 Let $R$ be a commutative ring and $B^\cu$ be a CDG\+ring endowed with
a ring homomorphism $R\rarrow B^0$.
 Let us say that a CDG\+module $P^\cu$ over $B^\cu$ is
\emph{$R$\+contraadjusted} (respectively, \emph{$R$\+cotorsion}) if
the $R$\+module $P^i$ is contraadjusted (resp., cotorsion) for every
degree $i\in\boZ$.
 Let us also say that a CDG\+module $P^\cu$ over $B^\cu$ is
\emph{$B^*$\+cotorsion} if the graded $B^*$\+module $P^*$ is cotorsion.
 Let us denote the full DG\+subcategories of $B^*$\+cotorsion,
$R$\+cotorsion, and $R$\+contraadjusted left CDG\+modules over $B^\cu$
by $B^\cu\bModl^{B^*\dcot}\subset B^\cu\bModl^{R\dcot}\subset
B^\cu\bModl^{R\dcta}\subset B^\cu\bModl$.
 Clearly, all the three full DG\+subcategories are closed under shifts,
twists, and direct products (hence also under cones) in $B^\cu\bModl$.

\begin{lem} \label{ctrh-cdg-modules-on-affine-schemes}
 Let $U=\Spec R$ be an affine scheme, $\cB^\cu$ be a quasi-coherent
CDG\+quasi-algebra over $U$, and $B^\cu=\cB^\cu(U)=
(\cB^*(U),d_{\cB,U},h_{\cB,U})$ be the corresponding CDG\+quasi-algebra
over $R$, as per Lemma~\ref{qcoh-cdg-quasi-algebras-on-affine-schemes}.
 Then the rule assigning the CDG\+module\/
$\P^\cu[U]=(\P^*[U],d_{\P,U})$ over $\cB^\cu(U)$ to a contraherent
CDG\+module\/ $\P^\cu$ over $\cB^\cu$ is a natural equivalence between
the DG\+category of contraherent CDG\+modules over $\cB^\cu$ and
the DG\+category of $R$\+contraadjusted CDG\+modules over
the CDG\+ring $B^\cu=\cB^\cu(U)$,
\begin{equation} \label{cta-ctrh-cdg-modules-on-affine-scheme}
 \cB^\cu\bCtrh\simeq B^\cu\bModl^{R\dcta}.
\end{equation}
 The equivalence of
DG\+categories~\eqref{cta-ctrh-cdg-modules-on-affine-scheme} restricts
to equivalences of DG\+categories
\begin{align}
 \cB^\cu\bCtrh^{U\dlct} &\simeq B^\cu\bModl^{R\dcot}, \\
 \cB^\cu\bCtrh^{\cB^*\dlct} &\simeq B^\cu\bModl^{B^*\dcot}.
\end{align}
\end{lem}

\begin{proof}
 This is dual-analogous to
Lemma~\ref{qcoh-cdg-modules-on-affine-schemes}, and essentially
the same argument as the proof of
Corollary~\ref{lcth-cdg-module-not-dependent-on-covering}.
 According to (the proof of) 
Lemma~\ref{qcoh-cdg-quasi-algebras-on-affine-schemes}, one can assume
that all the restriction morphisms of CDG\+rings in
the quasi-coherent CDG\+quasi-algebra $\cB^\cu$ on $U$ are strict.
 Then the assertions follow from the equivalences $\cB^*\Ctrh\simeq
B^*\Modl^{R\dcta}$, \ $\cB^*\Ctrh^{U\dlct}\simeq B^*\Modl^{R\dcot}$,
\ $\cB^*\Ctrh^{\cB^*\dlct}\simeq B^*\Modl^{B^*\dcot}$ on the level
of the categories of graded modules together with
Corollary~\ref{cdg-quasi-algebra-co-extension-of-scalars}(b).
\end{proof}

\Section{Contraderived and Semiderived Categories}
\label{contraderived-and-semiderived-secn}

 In this section we define the contraderived category of $\bW$\+locally
contraherent CDG\+modules over a quasi-coherent CDG\+quasi-algebra
$\cA^\cu$ over a scheme $X$ and prove that this triangulated category
does not depend on the open covering~$\bW$.
 We also prove the similar result for the semiderived category of
$\bW$\+locally contraherent modules over a quasi-coherent quasi-algebra
$\cA$ over~$X$.

\subsection{Freely generated CDG-modules}
 Let $B^\cu=(B^*,d,h)$ be a CDG\+ring.
 The notation $B^*\Modl$ for the abelian category of graded left
$B^*$\+modules appeared in Section~\ref{graded-modules-subsecn},
while the notation $B^\cu\bModl$ for the DG\+category of left
CDG\+modules over $B^\cu$ was introduced in
Section~\ref{dg-categories-of-cdg-modules-subsecn}.

 Furthermore, following~\cite{Pedg,PS5} and the beginning of 
Section~\ref{dg-categories-of-cdg-modules-subsecn}, the notation
$\sZ^0(B^\cu\bModl)$ stands for the abelian category of left
CDG\+modules over $B^\cu$ and closed morphisms of degree~$0$
between them.
 See~\cite[Section~2.2]{Pedg} or~\cite[Section~1.9]{PS5}.

 The forgetful functor $\#\:\sZ^0(B^\cu\bModl)\rarrow B^*\Modl$ assigns
to any CDG\+module $M^\cu=(M^*,d_M)$ over $B^\cu$ its underlying
graded $B^*$\+module $M^*=(M^\cu)^\#$.

\begin{lem} \label{G-plus-G-minus-functors-for-modules}
 The forgetful functor\/ $\#\:\sZ^0(B^\cu\bModl)\rarrow B^*\Modl$ has
a left adjoint functor $G^+\:B^*\Modl\rarrow\sZ^0(B^\cu\bModl)$
and a right adjoint functor $G^-\:B^*\Modl\rarrow\sZ^0(B^\cu\bModl)$.
 For any graded $B^*$\+module $M^*$, there is a natural closed
isomorphism $G^-(M^*)\simeq G^+(M^*)[1]$ of CDG\+modules over $B^\cu$,
where\/ $[1]$~denotes the cohomological degree shift.
 For any graded $B^*$\+module $M^*$, there are natural short exact
sequences of graded $B^*$\+modules \hbadness=1350
\begin{gather}
 0\lrarrow M^*\lrarrow G^+(M^*)^\#\lrarrow M^*[-1]\lrarrow0
 \label{G-plus-sequence} \\
 0\lrarrow M^*[1]\lrarrow G^-(M^*)^\#\lrarrow M^*\lrarrow0,
 \label{G-minus-sequence}
\end{gather}
where $M^*\rarrow G^+(M^*)^\#$ and $G^-(M^*)^\#\rarrow M^*$ are
the adjunction morphisms.
\end{lem}

\begin{proof}
 Explicit constructions of the CDG\+modules $G^+(M^*)$ and $G^-(M^*)$
are spelled out in~\cite[proof of Theorem~3.6]{Pkoszul}.
 The CDG\+module $G^+(M^*)$ is said to be \emph{freely generated}
by the graded module $M^*$, while the CDG\+module $G^-(M^*)$ is
\emph{cofreely cogenerated} by~$M^*$.
 The grading component $G^+(M^*)^i$, \,$i\in\boZ$ of the CDG\+module
$G^+(M^*)$ is the set of all formal expressions $x'+dx''$ with
$x'\in M^i$ and $x''\in M^{i-1}$.
 The action of $B^*$ and the differential on $G^+(M^*)$ are given by
explicit formulas that can be found in~\cite{Pkoszul}.
 The adjunction morphism $M^*\rarrow G^+(M^*)^\#$ is given by
the obvious rule $x\longmapsto x+d(0)$.
 For more fancy constructions of $G^+$ and $G^-$ and the proofs of
the other assertions of the lemma, see~\cite[Proposition~3.1]{Pedg}.
 See also~\cite[Proposition~1.3.2]{Bec} or~\cite[Section~2.6]{PS5}.
\end{proof}

\begin{lem} \label{G-plus-reflects-exactness}
 Let $K^*\rarrow L^*\rarrow N^*$ be a composable pair of morphisms of
graded $B^*$\+modules.
 Then the short sequence of graded $B^*$\+modules\/ $0\rarrow K^*
\rarrow L^*\rarrow N^*\rarrow0$ is exact if and only if the short
sequence of graded $B^*$\+modules\/ $0\rarrow G^+(K^*)^\#\rarrow
G^+(L^*)^\#\rarrow G^+(N^*)^\#\rarrow0$ is exact.
\end{lem}

\begin{proof}
 The explicit construction of the CDG\+module $G^+(M^*)$ mentioned
in the proof of Lemma~\ref{G-plus-G-minus-functors-for-modules}
implies a natural isomorphsm of graded abelian groups (but \emph{not}
$B^*$\+modules) $G^+(M^*)^\#\simeq M^*\oplus M^*[-1]$ for any graded
$B^*$\+module~$M^*$.
 It remains to point out that a short sequence of graded $B^*$\+modules
is exact if and only if it is exact as a short sequence of graded
abelian groups, and that the shifts and finite direct sums preserve
and reflect exactness.
\end{proof}

 Besides the obvious functoriality as adjoint functors,
the constructions $M^*\longmapsto G^+(M^*)$ and $M^*\longmapsto
G^-(M^*)$ are also functorial with respect to morphisms of CDG\+rings,
in the following sense.
 Let $(f,a)\:B^\cu=(B^*,d_B,h_B)\rarrow(A^*,d_A,h_A)=A^\cu$ be
a morphism of CDG\+rings.

 Let $M^*$ be a graded left $A^*$\+module and $N^*$ be a graded left
$B^*$\+module.
 Suppose given a morphism of graded $B^*$\+modules $g\:N^*\rarrow M^*$.
 Consider the composition $N^*\rarrow M^*\rarrow G^+(M^*)^\#$, where
$M^*\rarrow G^+(M^*)^\#$ is the adjunction morphism.
 As any CDG\+module over $A^\cu$, the CDG\+module $G^+(M^*)$ acquires
a structure of CDG\+module over $B^\cu$ using formula~(vi) from
Section~\ref{cdg-rings-cdg-modules-subsecn}.
 By adjunction, the morphism of graded $B^*$\+modules $N^*\rarrow
G^+(M^*)^\#$ corresponds to a morphism $G^+(N^*)\rarrow G^+(M^*)$ of
CDG\+modules over~$B^\cu$.
 We have constructed a map of CDG\+modules $G^+(g)\:G^+(N^*)\rarrow
G^+(M^*)$ compatible with the morphism of CDG\+rings
$(f,a)\:B^\cu\rarrow A^\cu$.
 Explicitly, the map of CDG\+modules $G^+(g)$ is given by the formula
$G^+(g)(y'+dy'')=g(y')+ag(y'')+d(g(y''))$ for all $y'\in N^i$ and
$y''\in N^{i-1}$, \,$i\in\boZ$ (the notation from the proof of
Lemma~\ref{G-plus-G-minus-functors-for-modules}).
 Using the natural isomorphism $G^-\simeq G^+[1]$, we also obtain
a map of CDG\+modules $G^-(g)\:G^-(N^*)\rarrow G^-(M^*)$ compatible
with the morphism of CDG\+rings $B^\cu\rarrow A^\cu$.

 Now suppose given a morphism of graded $B^*$\+modules in the opposite
direction, $g\:M^*\rarrow N^*$.
 Consider the composition $G^-(M^*)^\#\rarrow M^*\rarrow N^*$, where
$G^-(M^*)\rarrow M^*$ is the adjunction morphism.
 The CDG\+module $G^-(M^*)$ over $A^\cu$ acquires a structure of
CDG\+module over $B^\cu$ using formula~(vi) from
Section~\ref{cdg-rings-cdg-modules-subsecn}.
 By adjunction, the morphism of graded $B^*$\+modules $G^-(M^*)^\#
\rarrow N^*$ corresponds to a morphism $G^-(M^*)\rarrow G^-(N^*)$
of CDG\+modules over~$B^\cu$.
 We have constructed a map of CDG\+modules $G^-(g)\:G^-(M^*)\rarrow
G^-(N^*)$ compatible with the morphism of CDG\+rings
$(f,a)\:B^\cu\rarrow A^\cu$.
 Using the natural isomorphism $G^+\simeq G^-[-1]$, we also obtain
a map of CDG\+modules $G^+(g)\:G^+(M^*)\rarrow G^+(N^*)$ compatible
with the morphism of CDG\+rings $B^\cu\rarrow A^\cu$.
 Explicitly, the map of CDG\+modules $G^+(g)$ is given by the formula
$G^+(g)(x'+dx'')=g(x')-g(ax'')+d(g(x''))$ for all $x'\in M^i$ and
$x''\in M^{i-1}$, \,$i\in\boZ$ (the notation from the proof of
Lemma~\ref{G-plus-G-minus-functors-for-modules}).

 Let $X$ be a scheme with an open covering $\bW$, and let $\cB^\cu$ be
a quasi-coherent CDG\+quasi-algebra over~$X$.
 The notation $\cB^*\Qcoh$ for the abelian category of quasi-coherent
graded left $\cB^*$\+modules and $\cB^*\Lcth_\bW$ for the exact
category of $\bW$\+locally contraherent graded $\cB^*$\+modules was
introduced in Section~\ref{graded-modules-subsecn}.

 Following the notation of
Section~\ref{dg-categories-of-cdg-modules-subsecn},
\,$\sZ^0(\cB^\cu\bQcoh)$ denotes the category of quasi-coherent left
CDG\+modules over $\cB^\cu$ and closed morphisms of degree~$0$
between them.
 One can easily see that the category $\sZ^0(\cB^\cu\bQcoh)$ is
abelian (in fact, a Grothendieck category).
 Similarly, the notation $\sZ^0(\cB^\cu\bLcth_\bW)$ stands for
the additive category of $\bW$\+locally contraherent left CDG\+modules
over $\cB^\cu$ and closed morphisms of degree~$0$ between them.
 The category $\sZ^0(\cB^\cu\bLcth_\bW)$ is endowed with the exact
category structure in which a short sequence is exact if and only if
its underlying short sequence of $\bW$\+locally contraherent graded
$\cB^*$\+modules is exact (see
Corollaries~\ref{abelian-dg-category-of-qcoh-cdg-modules}
and~\ref{exact-dg-categories-of-lcth-cdg-modules} in the next
Section~\ref{abelian-and-exact-dg-categs-of-cdg-modules-subsecn}
for the discussion).

 The forgetful functor $\#\:\sZ^0(\cB^\cu\bQcoh)\rarrow\cB^*\Qcoh$
assigns to any quasi-coherent CDG\+module $\M^\cu$ over $\cB^\cu$ its
underlying quasi-coherent graded $\cB^*$\+module $\M^*=(\M^\cu)^\#$.
 The forgetful functor $\#\:\sZ^0(\cB^\cu\bLcth_\bW)\rarrow
\cB^*\Lcth_\bW$ assigns to any $\bW$\+locally contraherent
CDG\+module $\P^\cu$ over $\cB^\cu$ its underlying $\bW$\+locally
contraherent graded $\cB^*$\+module $\P^*=(\P^\cu)^\#$.

\begin{lem} \label{G-plus-G-minus-functors-for-qcoh-lcth}
\textup{(a)} The forgetful functor\/ $\#\:\sZ^0(\cB^\cu\bQcoh)
\rarrow\cB^*\Qcoh$ has a left adjoint functor $G^+\:\cB^*\Qcoh\rarrow
\sZ^0(\cB^\cu\bQcoh)$ and a right adjoint functor $G^-\:\cB^*\Qcoh
\rarrow\sZ^0(\cB^\cu\bQcoh)$.
 For any quasi-coherent graded $\cB^*$\+module $\M^*$, there is
a natural closed isomorphism $G^-(\M^*)\simeq G^+(\M^*)[1]$ of
quasi-coherent CDG\+modules over~$\cB^\cu$.
 For any quasi-coherent graded $\cB^*$\+module $\M^*$, there are natural
short exact sequences of quasi-coherent graded $\cB^*$\+modules
\begin{gather}
 0\lrarrow\M^*\lrarrow G^+(\M^*)^\#\lrarrow\M^*[-1]\lrarrow0
 \label{G-plus-sequence-qcoh} \\
 0\lrarrow\M^*[1]\lrarrow G^-(\M^*)^\#\lrarrow\M^*\lrarrow0,
 \label{G-minus-sequence-qcoh}
\end{gather}
where $\M^*\rarrow G^+(\M^*)^\#$ and $G^-(\M^*)^\#\rarrow\M^*$ are
the adjunction morphisms. \par
\textup{(b)} The forgetful functor\/ $\#\:\sZ^0(\cB^\cu\bLcth_\bW)
\rarrow\cB^*\Lcth_\bW$ has a left adjoint functor $G^+\:\cB^*\Lcth_\bW
\rarrow\sZ^0(\cB^\cu\bLcth_\bW)$ and a right adjoint functor
$G^-\:\cB^*\Lcth_\bW\allowbreak\rarrow\sZ^0(\cB^\cu\bLcth_\bW)$.
 For any\/ $\bW$\+locally contraherent graded $\cB^*$\+module\/ $\P^*$,
there is a natural closed isomorphism $G^-(\P^*)\simeq G^+(\P^*)[1]$
of\/ $\bW$\+locally contraherent CDG\+modules over~$\cB^\cu$.
 For any\/ $\bW$\+locally contraherent graded $\cB^*$\+module\/ $\P^*$,
there are natural admissible short exact sequences
\begin{gather}
 0\lrarrow\P^*\lrarrow G^+(\P^*)^\#\lrarrow\P^*[-1]\lrarrow0
 \label{G-plus-sequence-lcth} \\
 0\lrarrow\P^*[1]\lrarrow G^-(\P^*)^\#\lrarrow\P^*\lrarrow0,
 \label{G-minus-sequence-lcth}
\end{gather}
in the exact category $\cB^*\Lcth_\bW$, where\/
$\P^*\rarrow G^+(\P^*)^\#$ and $G^-(\P^*)^\#\rarrow\P^*$ are
the adjunction morphisms.
\end{lem}

\begin{proof}
 In part~(a), the quasi-coherent CDG\+modules $G^\pm(\M^*)$ are
constructed by the rule $G^\pm(\M^*)(U)=G^\pm((\M^*(U))$ for all
affine open subschemes $U\subset X$, where $G^\pm(\M^*(U))$ is
the CDG\+module over the CDG\+ring $\cB^\cu(U)=
(\cB^*(U),d_{\cB,U},h_{\cB,U})$ (co)freely (co)generated by
the graded $\cB^*$\+module $\M^*(U)$.
 For any pair of affine open subschemes $V\subset U\subset X$, 
the restriction map $G^\pm(\M^*(U))\rarrow G^\pm(\M^*(V))$ in
the (pre)sheaf $G^\pm(\M^*)$ is obtained by applying the construction
$G^\pm$ to the restriction map of graded modules $\M^*(U)\rarrow\M^*(V)$
over the underlying graded rings of the morphism of CDG\+rings
$(\rho_{\cB,VU},a_{\cB,VU})\:\cB^\cu(U)\rarrow\cB^\cu(V)$ (see
the discussion above in this section).
 In view of the natural short exact sequences of graded
modules~(\ref{G-plus-sequence}\+-\ref{G-minus-sequence}), the existence
of natural short exact sequences of presheaves of graded
$\cB^*$\+modules~(\ref{G-plus-sequence-qcoh}\+-%
\ref{G-minus-sequence-qcoh}) on the topology base $\bB$ of all affine
open subschemes of $X$ is straightforward to check.
 As the quasi-coherence condition is preserved by extensions in
the abelian category of presheaves of $\cO_X$\+modules on $\bB$,
the short exact sequences~(\ref{G-plus-sequence-qcoh}\+-%
\ref{G-minus-sequence-qcoh}) of presheaves of graded $\cB^*$\+modules
on $\bB$ obtained in this way imply quasi-coherence (of
the underlying presheaf of graded $\cO_X$\+modules) of the presheaf
of graded $\cB^*$\+modules $G^\pm(\M^*)^\#$ on~$\bB$.
 Then one can refer to Lemma~\ref{quasi-coherence-implies-sheaf}
and Theorem~\ref{extension-of-co-sheaves-from-topology-base}(a).
 (Cf.~\cite[Example~3.18]{Pedg}.)

 The argument in part~(b) is dual-analogous.
 The $\bW$\+locally contraherent CDG\+modules $G^\pm(\P^*)$ are
constructed by the rule $G^\pm(\P^*)[U]=G^\pm((\P^*[U])$ for all
affine open subschemes $U\subset X$ subordinate to $\bW$, where
$G^\pm(\P^*[U])$ is the CDG\+module over the CDG\+ring $\cB^\cu(U)=
(\cB^*(U),d_{\cB,U},h_{\cB,U})$ (co)freely (co)generated by
the graded $\cB^*$\+module $\P^*[U]$.
 For any pair of affine open subschemes $V\subset U\subset X$
subordinate to $\bW$, the corestriction map $G^\pm(\P^*[V])\rarrow
G^\pm(\P^*[U])$ in the co(pre)sheaf $G^\pm(\P^*)$ is obtained by
applying the construction $G^\pm$ to the corestriction map of graded
modules $\P^*[V]\rarrow\P^*[U]$ over the underlying graded rings of
the morphism of CDG\+rings $(\rho_{\cB,VU},a_{\cB,VU})\:\cB^\cu(U)
\rarrow\cB^\cu(V)$ (see the discussion above in this section).
 Using the natural short exact sequences of graded
modules~(\ref{G-plus-sequence}\+-\ref{G-minus-sequence}), the existence
of natural short exact sequences of copresheaves of graded
$\cB^*$\+modules~(\ref{G-plus-sequence-lcth}\+-%
\ref{G-minus-sequence-lcth}) on the topology base $\bB$ of all affine
open subschemes of $X$ subordinate to $\bW$ is straightforward to check.
 As the contraadjustedness and contraherence conditions~(i\+-ii) from
Section~\ref{locally-contraherent-cosheaves-subsecn} are preserved
by extensions in the abelian category of copresheaves of
$\cO_X$\+modules on $\bB$ (see the discussion
in~\cite[Section~3.1]{Pcosh}), the short exact
sequences~(\ref{G-plus-sequence-lcth}\+-%
\ref{G-minus-sequence-lcth}) of copresheaves of graded $\cB^*$\+modules
on $\bB$ obtained in this way imply contraherence (of
the underlying copresheaf of graded $\cO_X$\+modules) of the copresheaf
of graded $\cB^*$\+modules $G^\pm(\P^*)^\#$ on~$\bB$.
 Then one can refer to
Lemma~\ref{contraherence+contraadjustedness-imply-cosheaf}
and Theorem~\ref{extension-of-co-sheaves-from-topology-base}(b).
\end{proof}

\subsection{Abelian and exact DG-categories of CDG-modules}
\label{abelian-and-exact-dg-categs-of-cdg-modules-subsecn}
 An \emph{abelian DG\+cat\-e\-gory} $\bA$ can be simply defined as
a DG\+category with a zero object, shifts, and cones such that
the additive category $\sZ^0(\bA)$ of closed morphisms of degree~$0$
in $\bA$ is abelian~\cite[Section~4.6]{Pedg},
\cite[Proposition~3.8]{PS5}.
 An \emph{exact DG\+category} $\bE$ can be simply defined as
a DG\+category with a zero object, shifts, and cones for which
the additive category $\sZ^0(\bE)$ is endowed with an exact category
structure (in the sense of Quillen~\cite{Bueh}) preserved by
the cohomological degree shift functors~$[1]$ and~$[-1]$ such that
the functor of cone of identity endomorphism $E\longmapsto
\cone(\id_E)\:\sZ^0(\bE)\rarrow\sZ^0(\bE)$ preserves and reflects
admissible short exact sequences~\cite[Proposition~4.7,
Theorem~4.17, and Remark~4.18]{Pedg}.

 Stated in this way, these definitions do not mention
the \emph{underlying abelian/exact category of graded objects},
which plays an important role.
 One starts with the construction of an \emph{almost involution}
assigning to a DG\+category $\bA$ a new DG\+category~$\bA^\bec$
\,\cite[Section~3.2]{Pedg}, \cite[Section~2.1]{PS5}.
 Then the (pre)additive category $\sZ^0(\bA^\bec)$ is viewed as
the underlying category of graded objects for~$\bA$.
 This construction is important for applications of the abelian/exact
DG\+category theory that we are interested in.
 The additive category $\sZ^0(\bA^\bec)$ is abelian for an abelian
DG\+category $\bA$, and carries an exact category structure in
the case of an exact DG\+category~$\bA$.

 For any DG\+category $\bA$ with shifts and cones, there is a natural
faithful additive functor $\Phi_\bA\:\sZ^0(\bA)\rarrow\sZ^0(\bA^\bec)$,
viewed as the functor assigning to a DG\+category object its underlying
graded object.
 The functor $\Phi_\bA$ has a left and a right adjoint functors,
denoted by $\Psi^+_\bA$ and $\Psi^-_\bA\:\sZ^0(\bA^\bec)\rarrow
\sZ^0(\bA)$.
 The functors $\Psi^+_\bA$ and $\Psi^-_\bA$ only differ from each other
by a shift, $\Psi^-_\bA=\Psi^+_\bA[1]$ \,\cite[Lemma~3.4]{Pedg}, \cite[Section~2.2]{PS5}.

 For an exact DG\+category $\bE$, the functors $\Phi_\bE$ and
$\Psi^\pm_\bE$ preserve and refect admissible short exact
sequences~\cite[Section~4.3]{Pedg}.
 A DG\+functor $G\:\bF\rarrow\bE$ between exact DG\+categories $\bF$
and $\bE$ is said to be \emph{exact} if the functor $\sZ^0(G)\:
\sZ^0(\bF)\rarrow\sZ^0(\bE)$ is exact, or equivalently, the functor
$\sZ^0(G^\bec)\:\sZ^0(\bF^\bec)\rarrow\sZ^0(\bE^\bec)$ is exact (as
a functor between exact categories) \cite[Sections~3.4 and~4.4]{Pedg}.

 In particular, for any CDG\+ring $B^\cu$, the DG\+category of
CDG\+modules $B^\cu\bModl$ is an abelian
DG\+category~\cite[Example~4.41]{Pedg}, \cite[Example~3.14]{PS5}.
 According to~\cite[Example~3.17]{Pedg} or~\cite[Section~2.6]{PS5},
the abelian category $\sZ^0((B^\cu\bModl)^\bec)$ is naturally
equivalent to the abelian category of graded $B^*$\+modules.
 The equivalence is provided by the functor denoted by
\begin{equation} \label{CDG-modules-over-CDG-ring-Upsilon-equivalence}
 \Upsilon_{B^\subcu}\:B^*\Modl\simeq\sZ^0((B^\cu\bModl)^\bec).
\end{equation}
 The construction of the functor $\Upsilon_{B^\subcu}$ is based on
the construction of the functor $G^+\:B^*\Modl\rarrow
\sZ^0(B^\cu\bModl)$ from
Lemma~\ref{G-plus-G-minus-functors-for-modules} above.
 The equivalence of categories $\Upsilon_{B^\subcu}$ transforms
the forgetful functor $\#\:\sZ^0(B^\cu\bModl)\rarrow B^*\Modl$
into the functor $\Phi_\bA$ and the functors $G^\pm\:B^*\Modl
\rarrow\sZ^0(B^\cu\bModl)$ into the functor $\Psi^\pm_\bA$ for
the DG\+category $\bA=B^\cu\bModl$ \,\cite[diagrams~(8\+-9)]{Pedg},
\cite[diagrams~(4) and~(6)]{PS5}.

\begin{cor} \label{abelian-dg-category-of-qcoh-cdg-modules}
 Let $X$ be a scheme and $\cB^\cu$ be a quasi-coherent
CDG\+quasi-algebra over~$X$.
 Then the DG\+category $\cB^\cu\bQcoh$ of quasi-coherent left
CDG\+modules over $\cB^\cu$ is an abelian DG\+category.
 The abelian category\/ $\sZ^0((\cB^\cu\bQcoh)^\bec)$ is naturally
equivalent to the abelian category of quasi-coherent
graded left $\cB^*$\+modules $\cB^*\Qcoh$,
\begin{equation} \label{Upsilon-qcoh-equivalence}
 \Upsilon_{\cB^\subcu}^\qc\:
 \cB^*\Qcoh\simeq\sZ^0((\cB^\cu\bQcoh)^\bec).
\end{equation}
 The equivalence of categories\/ $\Upsilon_{\cB^\subcu}^\qc$
transforms the forgetful functor\/ $\#\:\sZ^0(\cB^\cu\bQcoh)\allowbreak
\rarrow\cB^*\Qcoh$ into the functor\/ $\Phi_\bA$, and the functors
$G^\pm\:\cB^*\Qcoh\rarrow \sZ^0(\cB^\cu\bQcoh)$ from
Lemma~\ref{G-plus-G-minus-functors-for-qcoh-lcth}(a) into the functors\/
$\Psi^\pm_\bA$ for the DG\+category\/ $\bA=\cB^\cu\bQcoh$.
 In other words, the following triangular diagrams of faithful
additive functors are commutative:
\begin{equation} \label{qcoh-cdg-modules-Phi-diagram}
\begin{gathered}
 \xymatrix{
  & \sZ^0(\cB^\cu\bQcoh) \ar[ld]_-{\#}
  \ar[rd]^-{\Phi_{\cB^\subcu\bQcoh}} \\
  \cB^*\Qcoh \ar@<-0.4ex>[rr]_-{\Upsilon_{\cB^\subcu}^\qc}
  && \sZ^0((\cB^\cu\bQcoh)^\bec) \ar@<-0.4ex>@{-}[ll]
 }
\end{gathered}
\end{equation}
and
\begin{equation} \label{qcoh-cdg-modules-Psi-diagram}
\begin{gathered}
 \xymatrix{
  \cB^*\Qcoh \ar@<0.4ex>[rr]^-{\Upsilon_{\cB^\subcu}^\qc}
  \ar[rd]_-{G^\pm}
  && \sZ^0((\cB^\cu\bQcoh)^\bec) \ar@<0.4ex>@{-}[ll]
  \ar[ld]^-{\Psi_{\cB^\subcu\bQcoh}^\pm} \\ & \sZ^0(\cB^\cu\bQcoh)
 }
\end{gathered}
\end{equation}
\end{cor}

\begin{proof}
 The construction of the functor $\Upsilon_{\cB^\subcu}^\qc$ is
similar to the one in~\cite[Example~3.18]{Pedg} and based on
the construction of the functor $G^+$ from
Lemma~\ref{G-plus-G-minus-functors-for-qcoh-lcth}(a).
 The reader can start with looking up the similar constructions
for CDG\+modules over CDG\+rings in~\cite[Example~3.17]{Pedg}
or~\cite[Section~2.6]{PS5}, or even with the construction of
the almost involution of DG\+categories $\bA\longmapsto\bA^\bec$
in~\cite[Section~3.2]{Pedg} or~\cite[Section~2.1]{PS5}.
 It is explained in~\cite[Example~3.18]{Pedg} how to construct
the inverse functor to~$\Upsilon_{\cB^\subcu}^\qc$.
 The category $\cB^*\Qcoh$ is obviously abelian, and the DG\+category
$\cB^\cu\bQcoh$ is clearly idempotent-complete with twists and
(even infinite) direct sums; so it is an abelian DG\+category
by~\cite[Corollary~4.38]{Pedg} or~\cite[Proposition~3.8]{PS5}.
 The remaining arguments are straightforward verifications.
\end{proof}

 The proof of the contraherent version of
Corollary~\ref{abelian-dg-category-of-qcoh-cdg-modules}
requires an additional lemma.

\begin{lem} \label{G-plus-periodicity-lemma}
 Let $R^*$ be a graded ring and $B^\cu=(B^*,d,h)$ be a CDG\+ring
endowed with a homomorphism of graded rings $R^*\rarrow B^*$.
 Let $P^*$ be a graded $B^*$\+module.
 Consider the CDG\+module $G^+(P^*)$ over $B^\cu$, as constructed in
Lemma~\ref{G-plus-G-minus-functors-for-modules}.
 In this setting: \par
\textup{(a)} if $R^*=R^0=R$ is an ungraded commutative ring and
the $R$\+module $G^+(P^*)^i$ is contraadjusted for all $i\in\boZ$,
then the $R$\+module $P^i$ is contraadjusted for all $i\in\boZ$; \par
\textup{(b)} if $R^*=R^0=R$ is an ungraded associative ring and
the $R$\+module $G^+(P^*)^i$ is cotorsion for all $i\in\boZ$,
then the $R$\+module $P^i$ is cotorsion for all $i\in\boZ$; \par
\textup{(c)} if the graded $R^*$\+module $G^+(P^*)^\#$ is cotorsion,
then the graded $R^*$\+module $P^*$ is cotorsion.
\end{lem}

\begin{proof}
 Part~(a) holds because, according to the short exact
sequence~\eqref{G-plus-sequence}, the $R$\+module $P^i$ is a quotient
$R$\+module of the $R$\+module $G^+(P)^{i+1}$.
 Recall that the class of contraadjusted $R$\+modules is closed
under quotients by Lemma~\ref{very-flat-cotorsion-pair-hereditary}(a).

 To prove part~(c), we notice that splicing up shifted copies of
the short exact sequence~\eqref{G-plus-sequence} produces an unbounded
acyclic complex of graded $B^*$\+modules
$$
 \dotsb\lrarrow G^+(P^*)[1]^\#\lrarrow G^+(P^*)^\#\lrarrow
 G^+(P^*)[-1]^\#\lrarrow\dotsb,
$$
whose graded $B^*$\+modules of cocycles are $P^*[i]$, \,$i\in\boZ$.
 Now it remains to refer to the graded version of
Theorem~\ref{module-cotorsion-periodicity}.

 Part~(b) is a particular case of part~(c), as a graded module over
an ungraded ring is cotorsion if and only all its grading components
are cotorsion.
\end{proof}

\begin{cor} \label{exact-dg-categories-of-lcth-cdg-modules}
 Let $X$ be a scheme with an open covering\/ $\bW$ and $\cB^\cu$ be
a quasi-coherent CDG\+quasi-algebra over~$X$.
 Then there is a commutative diagram of additive category equivalences
and fully faithful identity inclusion functors
\begin{equation} \label{Upsilon-lcth-diagram}
\begin{gathered}
 \xymatrix{
  \cB^*\Lcth_\bW^{\cB^*\dlct}
  \ar@<0.4ex>[rr]^-{\Upsilon_{\cB^\subcu}^\ct}
  \ar@{>->}[d]
  && \sZ^0((\cB^\cu\bLcth_\bW^{\cB^*\dlct})^\bec) \ar@<0.4ex>@{-}[ll]
  \ar@{>->}[d] \\
  \cB^*\Lcth_\bW^{X\dlct} \ar@<0.4ex>[rr]^-{\Upsilon_{\cB^\subcu}^\ct}
  \ar@{>->}[d]
  && \sZ^0((\cB^\cu\bLcth_\bW^{X\dlct})^\bec) \ar@<0.4ex>@{-}[ll]
  \ar@{>->}[d] \\
  \cB^*\Lcth_\bW \ar@<0.4ex>[rr]^-{\Upsilon_{\cB^\subcu}^\ct}
  && \sZ^0((\cB^\cu\bLcth_\bW)^\bec) \ar@<0.4ex>@{-}[ll]
 }
\end{gathered}
\end{equation}
with horizontal double lines showing category equivalences and
vertical arrows with tails showing fully faithful inclusions.
 The DG\+categories $\cB^\cu\bLcth_\bW$, \ $\cB^\cu\bLcth_\bW^{X\dlct}$,
and $\cB^\cu\bLcth_\bW^{\cB^*\dlct}$ have natural structures of exact
DG\+categories such that the related exact category structures on
the additive categories in the rightmost column of the diagram agree
with the exact category structures on the additive categories in
the leftmost column of the diagram described in
Sections~\ref{cosheaves-of-A-modules-subsecn}\+-%
\ref{A-loc-cotors-loc-inj-cosheaves-subsecn}
and~\ref{graded-modules-subsecn}.

 The equivalence of categories\/ $\Upsilon_{\cB^\subcu}^\ct$
transforms the forgetful functor\/ $\#\:\sZ^0(\cB^\cu\bLcth_\bW)
\allowbreak\rarrow\cB^*\Lcth_\bW$ into the functor\/ $\Phi_\bE$, and
the functors $G^\pm\:\cB^*\Lcth_\bW\rarrow \sZ^0(\cB^\cu\bLcth_\bW)$
from Lemma~\ref{G-plus-G-minus-functors-for-qcoh-lcth}(b) into
the functors\/ $\Psi^\pm_\bE$ for the DG\+category\/
$\bE=\cB^\cu\bLcth_\bW$.
 In other words, the following triangular diagrams of faithful
additive functors are commutative: \hfuzz=10pt \hbadness=1800
\begin{equation} \label{lcth-cdg-modules-Phi-diagram}
\begin{gathered}
 \xymatrix{
  & \sZ^0(\cB^\cu\bLcth_\bW) \ar[ld]_-{\#}
  \ar[rd]^-{\Phi_{\cB^\subcu\bLcth_\bW}} \\
  \cB^*\Lcth_\bW \ar@<-0.4ex>[rr]_-{\Upsilon_{\cB^\subcu}^\ct}
  && \sZ^0((\cB^\cu\bLcth_\bW)^\bec) \ar@<-0.4ex>@{-}[ll]
 }
\end{gathered}
\end{equation}
and
\begin{equation} \label{lcth-cdg-modules-Psi-diagram}
\begin{gathered}
 \xymatrix{
  \cB^*\Lcth_\bW \ar@<0.4ex>[rr]^-{\Upsilon_{\cB^\subcu}^\ct}
  \ar[rd]_-{G^\pm}
  && \sZ^0((\cB^\cu\bLcth_\bW)^\bec) \ar@<0.4ex>@{-}[ll]
  \ar[ld]^-{\Psi_{\cB^\subcu\bLcth_\bW}^\pm} \\ 
  & \sZ^0(\cB^\cu\bLcth_\bW)
 }
\end{gathered}
\end{equation}
 The similar assertions hold for the DG\+categories
$\cB^\cu\bLcth_\bW^{X\dlct}$ and $\cB^\cu\bLcth_\bW^{\cB^*\dlct}$.
\end{cor}

\begin{proof}
 The construction of the functors $\Upsilon_{\cB^\subcu}^\ct$ is
similar to that of the functor $\Upsilon_{\cB^\subcu}^\qc$ in
Corollary~\ref{abelian-dg-category-of-qcoh-cdg-modules}, and
based on the construction of the functor $G^+$ from
Lemma~\ref{G-plus-G-minus-functors-for-qcoh-lcth}(b).
 One needs to use Lemma~\ref{G-plus-periodicity-lemma} in order to
show that the functors $\Upsilon_{\cB^\subcu}^\ct$ are essentially
surjective.
 In the exact category structures on the additive categories
$\sZ^0(\cB^\cu\bLcth_\bW)$, \ $\sZ^0(\cB^\cu\bLcth_\bW^{X\dlct})$,
and $\sZ^0(\cB^\cu\bLcth_\bW^{\cB^*\dlct})$, a short sequence is
exact if and only if the forgetful functor~$\#$ takes it to a short
exact sequence in the respective exact category $\cB^*\Lcth_\bW$,
\ $\cB^*\Lcth_\bW^{X\dlct}$, or $\cB^*\Lcth_\bW^{\cB^*\dlct}$.
 One easily checks that this is a DG\+compatible exact structure
in the sense of~\cite[Section~4.2]{Pedg} inducing the original
exact structure on the underlying category of graded objects,
as per Sections~\ref{cosheaves-of-A-modules-subsecn}\+-%
\ref{A-loc-cotors-loc-inj-cosheaves-subsecn}
and~\ref{graded-modules-subsecn}
(Lemma~\ref{G-plus-reflects-exactness} is helpful here).
 So~\cite[Theorem~4.17]{Pedg} is applicable.
 Checking that the diagrams are commutative is straightforward.
\end{proof}

\subsection{Exact DG-categories of contraadjusted and cotorsion
quasi-coherent CDG-modules} \label{exact-dg-of-cta-cot-qcoh-cdg}
 Let $X$ be a scheme.
 A graded quasi-coherent sheaf $\C^*=\bigoplus_{i\in\boZ}\C^i$ on $X$
is said to be \emph{contraadjusted} if all its grading component
quasi-coherent sheaves $\C^i$, \,$i\in\boZ$, are contraadjusted
(in the sense of
Section~\ref{antilocality-of-X-contraadjusted-subsecn}).
 The graded quasi-coherent sheaf $\C^*$ is said to be \emph{cotorsion}
if all the quasi-coherent sheaves $\C^i$ are cotorsion
(in the sense of
Section~\ref{antilocality-of-X-cotorsion-subsecn}).

 Let $\cA^*$ be a quasi-coherent graded quasi-algebra over~$X$
(in the sense of Section~\ref{graded-modules-subsecn}).
 A quasi-coherent graded $\cA^*$\+module $\C^*$ on $X$ is said to be
\emph{$X$\+contraadjusted} (respectively, \emph{$X$\+cotorsion}) if
the underlying graded quasi-coherent sheaf of $\C^*$ is contraadjusted
(resp., cotorsion).

 Similarly to the definition in
Section~\ref{A-loc-cotors-loc-inj-cosheaves-subsecn},
a quasi-coherent graded $\cA^*$\+module $\F^*$ on $X$ is said to be
\emph{flat} if, for every affine open subscheme $U\subset X$,
the graded $\cA^*(U)$\+module $\F^*(U)$ is flat.
 Similarly to the definition in
Section~\ref{antilocality-of-A-cotorsion-subsecn}, a quasi-coherent
graded left $\cA^*$\+module $\C^*$ on $X$ is said to be
\emph{$\cA^*$\+cotorsion} if $\Ext^1_{\cA^*}(\F^*,\C^*)=0$ for all
flat quasi-coherent graded left $\cA^*$\+modules $\F^*$ on~$X$.
 Here the notation $\Ext^*_{\cA^*}({-},{-})$ stands for the Ext groups
in the abelian category $\cA^*\Qcoh$ of quasi-coherent graded left
$\cA^*$\+modules on~$X$.

 The full subcategories of $X$\+contraadjusted, $X$\+cotorsion,
and $\cA^*$\+cotorsion quasi-coherent graded $\cA^*$\+modules
in the abelian category of quasi-coherent graded left $\cA^*$\+modules
$\cA^*\Qcoh$ on $X$ are denoted by $\cA^*\Qcoh^{X\dcta}$, \
$\cA^*\Qcoh^{X\dcot}$, and $\cA^*\Qcoh^{\cA^*\dcot}$, respectively.
 By the graded version of the argument in
Section~\ref{antilocality-of-A-cotorsion-subsecn}, all
$\cA^*$\+cotorsion quasi-coherent graded $\cA^*$\+modules are
$X$\+cotorsion, that is $\cA^*\Qcoh^{\cA^*\dcot}\subset
\cA^*\Qcoh^{X\dcot}$.
 By the definition, we have $\cA^*\Qcoh^{X\dcot}\subset
\cA^*\Qcoh^{X\dcta}\subset\cA^*\Qcoh$.

 Similarly to the discussion in
Sections~\ref{antilocality-of-X-contraadjusted-subsecn}\+-%
\ref{antilocality-of-A-cotorsion-subsecn}, the full subcategories
$\cA^*\Qcoh^{\cA^*\dcot}\allowbreak\subset\cA^*\Qcoh^{X\dcot}\subset
\cA^*\Qcoh^{X\dcta}\subset\cA^*\Qcoh$ are closed under extensions,
direct summands, and infinite products in $\cA^*\Qcoh$.
 For a quasi-compact semi-separated scheme $X$, these full subcategories
are also closed under cokernels of monomorphisms in $\cA^*\Qcoh$
(cf.\ Corollary~\ref{A-cotorsion-hereditariness-properties}(c)).
 So all the three full subcategories inherit exact category structures
from the abelian exact structure of $\cA^*\Qcoh$.

 Let $\cA^\cu$ be a CDG\+quasi-algebra over~$X$.
 A quasi-coherent CDG\+module $\C^\cu$ over~$\cA^\cu$ (as defined in
Section~\ref{qcoh-lcth-cdg-modules-subsecn}) is said to be
\emph{$X$\+contraadjusted}, \emph{$X$\+cotorsion}, or
\emph{$\cA^*$\+cotorsion} if its underlying graded $\cA^*$\+module
$\C^*$ has the respective property.
 We denote the full DG\+subcategory of $X$\+contraadjusted
quasi-coherent CDG\+modules by $\cA^\cu\bQcoh^{X\dcta}\subset
\cA^\cu\bQcoh$, the full DG\+subcategory of $X$\+cotorsion
quasi-coherent CDG\+modules by $\cA^\cu\bQcoh^{X\dcot}\subset
\cA^\cu\bQcoh$, and the full DG\+subcategory of $\cA^*$\+cotorsion
quasi-coherent CDG\+modules by $\cA^\cu\bQcoh^{\cA^*\dcot}\subset
\cA^\cu\bQcoh$.

 For any quasi-coherent CDG\+quasi-algebra $\cA^\cu$ over a scheme $X$,
the full DG\+sub\-cat\-e\-gories $\cA^\cu\bQcoh^{\cA^*\dcot}\subset
\cA^\cu\bQcoh^{X\dcot}\subset\cA^\cu\bQcoh^{X\dcta}\subset
\cA^\cu\bQcoh$ contain the zero object and are closed under finite
direct sums, shifts, and twists (hence also under cones) in
the DG\+category $\cA^\cu\bQcoh$.

\begin{lem} \label{cta-cot-qcoh-cdg-modules-on-affine-schemes}
 Let $U=\Spec R$ be an affine scheme, $\cB^\cu$ be a quasi-coherent
CDG\+quasi-algebra over $U$, and $B^\cu=\cB^\cu(U)=
(\cB^*(U),d_{\cB,U},h_{\cB,U})$ be the corresponding CDG\+quasi-algebra
over $R$, as per Lemma~\ref{qcoh-cdg-quasi-algebras-on-affine-schemes}.
 Then he equivalence of
DG\+categories~\eqref{qcoh-cdg-modules-on-affine-scheme-formula}
from Lemma~\ref{qcoh-cdg-modules-on-affine-schemes} restricts
to equivalences of DG\+categories
\begin{align}
 \cB^\cu\bQcoh^{U\dcta} &\simeq B^\cu\bModl^{R\dcta}, \\
 \cB^\cu\bQcoh^{U\dcot} &\simeq B^\cu\bModl^{R\dcot}, \\
 \cB^\cu\bQcoh^{\cB^*\dcot} &\simeq B^\cu\bModl^{B^*\dcot}
\end{align}
in the notation from the end of
Section~\ref{dg-categories-of-cdg-modules-subsecn}.
\end{lem}

\begin{proof}
 The assertion follows from the equivalences
$\cB^*\Qcoh^{U\dcta}\simeq B^*\Modl^{R\dcta}$, \
$\cB^*\Qcoh^{U\dcot}\simeq B^*\Modl^{R\dcot}$, \
$\cB^*\Qcoh^{\cB^*\dcot}\simeq B^*\Modl^{B^*\dcot}$
on the level of the categories of graded modules.
\end{proof}

 The functors $G^+$ and $G^-\:\cA^*\Qcoh\rarrow\sZ^0(\cA^\cu\bQcoh)$
were constructed in
Lemma~\ref{G-plus-G-minus-functors-for-qcoh-lcth}(a).
 The following lemma is a quasi-coherent version of
Lemma~\ref{G-plus-periodicity-lemma}.

\begin{lem} \label{qcoh-G-plus-periodicity-lemma}
 Let $X$ be a quasi-compact semi-separated scheme and $\cA^\cu$ be
a quasi-coherent CDG\+quasi-algebra over~$X$.
 Let $\C^*$ be a quasi-coherent graded $\cA^*$\+module.
 In this setting: \par
\textup{(a)} if the quasi-coherent sheaf $G^+(\C^*)^i$ on $X$ is
contraadjusted for all $i\in\boZ$, then the quasi-coherent sheaf $\C^i$
is contraadjusted for all $i\in\boZ$; \par
\textup{(b)} if the quasi-coherent sheaf $G^+(\C^*)^i$ on $X$ is
cotorsion for all $i\in\boZ$, then the quasi-coherent sheaf $\C^i$ is
cotorsion for all $i\in\boZ$; \par
\textup{(c)} if the quasi-coherent graded $\cA^*$\+module
$G^+(\C^*)^\#$ is $\cA^*$\+cotorsion, then the quasi-coherent graded
$\cA^*$\+module $\C^*$ is $\cA^*$\+cotorsion.
\end{lem}

\begin{proof}
 Splicing up shifted copies of the short exact
sequence~\eqref{G-plus-sequence-qcoh} produces an unbounded
acyclic complex of quasi-coherent graded $\cA^*$\+modules
$$
 \dotsb\lrarrow G^+(\C^*)[1]^\#\lrarrow G^+(\C^*)^\#\lrarrow
 G^+(\C^*)[-1]^\#\lrarrow\dotsb,
$$
whose quasi-coherent graded $\cA^*$\+modules of cocycles are
$\C^*[i]$, \,$i\in\boZ$.
 Now part~(b) follows by Theorem~\ref{qcoh-cotorsion-periodicity},
and part~(c) by the graded version of
Theorem~\ref{qcoh-quasi-algebra-cotorsion-periodicity}.
 To prove part~(a), one needs to apply the ``contraadjusted
periodicity theorem for quasi-coherent sheaves'', which is provided
by the dual version of~\cite[Corollary~A.5.4]{Pcosh}, based on the fact
that the contraadjusted coresolution dimension of any quasi-coherent
sheaf on $X$ does not exceed the number of affine open subschemes
in any given affine open covering of~$X$
(see~\cite[Lemma~4.7.1(b)]{Pcosh}, cf.\ the proof of
Lemma~\ref{X-A-cta-cot-coresolving-coresolution-dimension}(b) above).
\end{proof}

\begin{cor} \label{exact-dg-categories-of-qcoh-cta-cot-cdg-modules}
 Let $X$ be a quasi-compact semi-separated scheme and $\cA^\cu$ be
a quasi-coherent CDG\+quasi-algebra over~$X$.
 Then there is a commutative diagram of additive category equivalences
and fully faithful identity inclusion functors
\begin{equation} \label{Upsilon-qcoh-diagram}
\begin{gathered}
 \xymatrix{
  \cA^*\Qcoh^{\cA^*\dcot}
  \ar@<0.4ex>[rr]^-{\Upsilon_{\cA^\subcu}^\qc}
  \ar@{>->}[d]
  && \sZ^0((\cA^\cu\bQcoh^{\cA^*\dcot})^\bec) \ar@<0.4ex>@{-}[ll]
  \ar@{>->}[d] \\
  \cA^*\Qcoh^{X\dcot} \ar@<0.4ex>[rr]^-{\Upsilon_{\cA^\subcu}^\qc}
  \ar@{>->}[d]
  && \sZ^0((\cA^\cu\bQcoh^{X\dcot})^\bec) \ar@<0.4ex>@{-}[ll]
  \ar@{>->}[d] \\
  \cA^*\Qcoh^{X\dcta} \ar@<0.4ex>[rr]^-{\Upsilon_{\cA^\subcu}^\qc}
  \ar@{>->}[d]
  && \sZ^0((\cA^\cu\bQcoh^{X\dcta})^\bec) \ar@<0.4ex>@{-}[ll]
  \ar@{>->}[d] \\
  \cA^*\Qcoh \ar@<0.4ex>[rr]^-{\Upsilon_{\cA^\subcu}^\qc}
  && \sZ^0((\cA^\cu\bQcoh)^\bec) \ar@<0.4ex>@{-}[ll]
 }
\end{gathered}
\end{equation}
with horizontal double lines showing category equivalences and
vertical arrows with tails showing fully faithful inclusions.
 The category equivalence in the lower horizontal line is
the equivalence~\eqref{Upsilon-qcoh-equivalence}
from Corollary~\ref{abelian-dg-category-of-qcoh-cdg-modules}.
 The DG\+categories $\cA^\cu\bQcoh^{X\dcta}$, \
$\cA^\cu\bQcoh^{X\dcot}$, and $\cA^\cu\bQcoh^{\cA^*\dcot}$ have natural
structures of exact DG\+categories such that the related exact category
structures on the additive categories in the rightmost column of
the diagram agree with the exact category structures on the additive
categories in the leftmost column of the diagram described above in
this section.
\end{cor}

\begin{proof}
 The assertions about category equivalences follow immediately
from Corollary~\ref{abelian-dg-category-of-qcoh-cdg-modules} and
Lemma~\ref{qcoh-G-plus-periodicity-lemma}.
 The assertions about the exact DG\+category structures are provable
similarly to the respective assertions of
Corollary~\ref{exact-dg-categories-of-lcth-cdg-modules}.
\end{proof}

\subsection{Exact DG-categories of antilocal contraherent CDG-modules}
\label{exact-dg-categories-of-antilocal-ctrh-subsecn}
 Let $X$ be a scheme with an open covering~$\bW$.
 A graded $\bW$\+locally contraherent cosheaf $\P*=\prod_{i\in\boZ}\P^i$
on $X$ is said to be \emph{antilocal} if all its grading component
contraherent cosheaves $\P^i$, \,$i\in\boZ$, are antilocal (in the sense
of Section~\ref{A-antilocal=X-antilocal-subsecn}).

 Let $\cA^*$ be a quasi-coherent graded quasi-algebra over~$X$
(in the sense of Section~\ref{graded-modules-subsecn}).
 A $\bW$\+locally contraherent graded $\cA^*$\+module $\P^*$ on $X$ is
said to be \emph{$X$\+antilocal} if the underlying graded
$\bW$\+locally contraherent cosheaf (of $\cO_X$\+modules) of $\P^*$ is
antilocal.

 A $\bW$\+locally contraherent graded $\cA^*$\+module $\P^*$ on $X$ is
said to be \emph{$\cA^*$\+antilocal} if the functor
$\Hom^{\cA^*}(\P^*,{-})$ of morphisms in the category of locally
contraherent graded $\cA^*$\+modules on $X$ takes short exact sequences
of $\cA^*$\+locally injective $\bW$\+locally contraherent graded
$\cA^*$\+modules on~$X$ (in the sense of
Sections~\ref{A-loc-cotors-loc-inj-cosheaves-subsecn}
and~\ref{graded-modules-subsecn}) to short exact sequences of
abelian groups.
 This is the graded version of the respective definition from
Section~\ref{A-antilocal=X-antilocal-subsecn}.

 Assume that the scheme $X$ is quasi-compact and semi-separated.
 Then the graded version of
Corollary~\ref{X-antilocal-equivalent-to-A-antilocal} claims that
the classes of $X$\+antilocal and $\cA^*$\+antilocal $\bW$\+locally
contraherent graded $\cA^*$\+modules on $X$ coincde.
 So we will call such $\bW$\+locally contraherent graded
$\cA^*$\+modules simply \emph{antilocal}.
 In view of~\cite[Corollaries~4.3.5(c) and~4.3.7]{Pcosh}, it follows
that all antilocal $\bW$\+locally contraherent graded $\cA^*$\+modules
on $X$ are actually (globally) contraherent, and the class of such
contraherent graded $\cA^*$\+modules does not depend on an open
covering $\bW$ of the scheme~$X$.

 An $X$\+locally cotorsion or $\cA^*$\+locally cotorsion (locally)
contraherent graded $\cA^*$\+module is said to be antilocal if it is
antilocal as an ($X$\+locally contraadjusted) locally contraherent
graded $\cA^*$\+module.
 We will denote the full subcategory of antilocal contraherent
graded $\cA^*$\+modules by $\cA^*\Ctrh_\al\subset\cA^*\Lcth_\bW$.
 The full subcategory of antilocal $X$\+locally cotorsion contraherent
graded $\cA^*$\+modules will be denoted by $\cA^*\Ctrh^{X\dlct}_\al
\subset\cA^*\Lcth_\bW^{X\dlct}$, and the full subcategory of antilocal
$\cA^*$\+locally cotorsion contraherent graded $\cA^*$\+modules will be
denoted by $\cA^*\Ctrh^{\cA^*\dlct}_\al\subset
\cA^*\Lcth_\bW^{\cA^*\dlct}$.
 This is the graded version of the notation introduced in
Section~\ref{antilocal-contrah-A-lct-A-modules-subsecn}.

 By~\cite[Corollaries~4.3.3(b) and~4.3.10]{Pcosh}, the full
subcategory $\cA^*\Ctrh_\al$ is closed under extensions, kernels
of admissible epimorphisms, direct summands, and infinite products
in $X\Lcth_\bW$.
 Accordingly, the full subcategories $\cA^*\Ctrh^{\cA^*\dlct}_\al
\subset\cA^*\Ctrh^{X\dlct}_\al\subset\cA^*\Ctrh_\al$ inherit exact
category structures from $\cA^*\Lcth_\bW^{\cA^*\dlct}\subset
\cA^*\Lcth_\bW^{X\dlct}\subset\cA^*\Lcth_\bW$.

 Let $\cA^\cu$ be a CDG\+quasi-algebra over~$X$.
 A ($\bW$\+locally) contraherent CDG\+module $\P^\cu$ over $\cA^\cu$
(as defined in Sections~\ref{qcoh-lcth-cdg-modules-subsecn}\+-%
\ref{dg-categories-of-cdg-modules-subsecn}) is said to be
\emph{antilocal} if its underlying locally contraherent graded
$\cA^*$\+module $\P^*$ is antilocal.
 We denote the full DG\+subcategory of antilocal contraherent
CDG\+modules by $\cA^\cu\bCtrh_\al\subset\cA^\cu\bLcth_\bW$, the full
DG\+subcategory of antilocal $X$\+locally cotorsion contraherent
CDG\+modules by $\cA^\cu\bCtrh^{X\dlct}_\al\subset
\cA^\cu\bLcth_\bW^{X\dlct}$, and the full subcategory of antilocal
$\cA^*$\+locally cotorsion contraherent CDG\+modules by
$\cA^\cu\bCtrh^{\cA^*\dlct}_\al\subset\cA^\cu\bLcth_\bW^{\cA^*\dlct}$.

 For any quasi-coherent CDG\+quasi-algebra $\cA^\cu$ over
a quasi-compact semi-separated scheme $X$, the full DG\+subcategories
$\cA^\cu\bCtrh^{\cA^*\dlct}_\al\subset\cA^\cu\bCtrh^{X\dlct}_\al
\subset\cA^\cu\bCtrh_\al$ contain the zero object and are closed
under infinite products, shifts, and twists (hence also under cones)
in the DG\+categories $\cA^\cu\bLcth_\bW^{\cA^*\dlct}\subset
\cA^\cu\bLcth_\bW^{X\dlct}\subset\cA^\cu\bLcth_\bW$.
 
 The functors $G^+$ and $G^-\:\cA^*\Lcth_\bW\rarrow
\sZ^0(\cA^\cu\bLcth_\bW)$ were constructed in
Lemma~\ref{G-plus-G-minus-functors-for-qcoh-lcth}(b).
 The following lemma is a version of
Lemmas~\ref{G-plus-periodicity-lemma}
and~\ref{qcoh-G-plus-periodicity-lemma} for the antilocality
property.

\begin{lem} \label{antilocal-G-plus-periodicity-lemma}
 Let $X$ be a quasi-compact semi-separated scheme and $\cA^\cu$ be
a quasi-coherent CDG\+quasi-algebra over~$X$.
 Let\/ $\P^*$ be a\/ $\bW$\+locally contraherent graded $\cA^*$\+module.
 Assume that the\/ $\bW$\+locally contraherent graded
$\cA^*$\+module $G^+(\P^*)^\#$ is antilocal.
 Then the\/ $\bW$\+locally contraherent graded $\cA^*$\+module\/ $\P^*$
is antilocal, too.
\end{lem}

\begin{proof}
 The argument is dual-analogous to the proof of
Lemma~\ref{qcoh-G-plus-periodicity-lemma}(a).
 Splicing up shifted copies of the short exact
sequence~\eqref{G-plus-sequence-lcth} produces an unbounded
acyclic complex of $\bW$\+locally contraherent graded $\cA^*$\+modules
$$
 \dotsb\lrarrow G^+(\P^*)[1]^\#\lrarrow G^+(\P^*)^\#\lrarrow
 G^+(\P^*)[-1]^\#\lrarrow\dotsb,
$$
whose $\bW$\+locally contraherent graded $\cA^*$\+modules of cocycles
are $\P^*[i]$, \,$i\in\boZ$.
 Now we apply the ``antilocal periodicity theorem for $\bW$\+locally
contraherent cosheaves'', which is provided
by~\cite[Corollary~A.5.4]{Pcosh}, based on the fact that
the antilocal resolution dimension of any $\bW$\+locally contraherent
cosheaf on $X$ does not exceed $N-1$, where $N$ is the number of affine
open subschemes in an affine open covering of~$X$ subordinate to~$\bW$
(see~\cite[Lemma~4.7.1(c)]{Pcosh}, cf.\ the proof of
Lemma~\ref{antilocal-resolving-resolution-dimension} above).
\end{proof}

\begin{cor} \label{exact-dg-categories-of-antilocal-ctrh-cdg-modules}
 Let $X$ be a quasi-compact semi-separated scheme with an open
covering\/ $\bW$ and $\cA^\cu$ be a quasi-coherent CDG\+quasi-algebra
over~$X$.
 Then there are commutative diagrams of additive category equivalences
and fully faithful identity inclusion functors
\begin{equation} \label{Upsilon-A-lct-antiloc-diagram}
\begin{gathered}
 \xymatrix{
  \cA^*\Ctrh^{\cA^*\dlct}_\al
  \ar@<0.4ex>[rr]^-{\Upsilon_{\cA^\subcu}^\ct}
  \ar@{>->}[d]
  && \sZ^0((\cA^\cu\bCtrh^{\cA^*\dlct}_\al)^\bec) \ar@<0.4ex>@{-}[ll]
  \ar@{>->}[d] \\
  \cA^*\Lcth_\bW^{\cA^*\dlct}
  \ar@<0.4ex>[rr]^-{\Upsilon_{\cA^\subcu}^\ct}
  && \sZ^0((\cA^\cu\bLcth_\bW^{\cA^*\dlct})^\bec) \ar@<0.4ex>@{-}[ll]
 }
\end{gathered}
\end{equation}
\begin{equation} \label{Upsilon-X-lct-antiloc-diagram}
\begin{gathered}
 \xymatrix{
  \cA^*\Ctrh^{X\dlct}_\al
  \ar@<0.4ex>[rr]^-{\Upsilon_{\cA^\subcu}^\ct}
  \ar@{>->}[d]
  && \sZ^0((\cA^\cu\bCtrh^{X\dlct}_\al)^\bec) \ar@<0.4ex>@{-}[ll]
  \ar@{>->}[d] \\
  \cA^*\Lcth_\bW^{X\dlct}
  \ar@<0.4ex>[rr]^-{\Upsilon_{\cA^\subcu}^\ct}
  && \sZ^0((\cA^\cu\bLcth_\bW^{X\dlct})^\bec) \ar@<0.4ex>@{-}[ll]
 }
\end{gathered}
\end{equation}
\begin{equation} \label{Upsilon-X-lcta-antiloc-diagram}
\begin{gathered}
 \xymatrix{
  \cA^*\Ctrh_\al \ar@<0.4ex>[rr]^-{\Upsilon_{\cA^\subcu}^\ct}
  \ar@{>->}[d]
  && \sZ^0((\cA^\cu\bCtrh_\al)^\bec) \ar@<0.4ex>@{-}[ll]
  \ar@{>->}[d] \\
  \cA^*\Lcth_\bW \ar@<0.4ex>[rr]^-{\Upsilon_{\cA^\subcu}^\ct}
  && \sZ^0((\cA^\cu\bLcth_\bW)^\bec) \ar@<0.4ex>@{-}[ll]
 }
\end{gathered}
\end{equation}
with horizontal double lines showing category equivalences and
vertical arrows with tails showing fully faithful inclusions.
 The category equivalences in the lower horizontal lines of all
the three diagrams are the ones from
diagram~\eqref{Upsilon-lcth-diagram}
in Corollary~\ref{exact-dg-categories-of-lcth-cdg-modules}.
 The DG\+categories $\cA^\cu\bCtrh^{\cA^*\dlct}_\al$, \
$\cA^\cu\bCtrh^{X\dlct}_\al$, and $\cA^\cu\bCtrh_\al$ have natural
structures of exact DG\+categories such that the related exact category
structures on the additive categories in the rightmost column of
the diagram agree with the exact category structures on the additive
categories in the leftmost column of the diagram described above in
this section.
\end{cor}

\begin{proof}
 The assertions about category equivalences follow immediately
from Corollary~\ref{exact-dg-categories-of-lcth-cdg-modules} and
Lemma~\ref{antilocal-G-plus-periodicity-lemma}.
 The assertions about the exact DG\+category structures are provable
similarly to the respective assertions of
Corollary~\ref{exact-dg-categories-of-lcth-cdg-modules}.
\end{proof}

\subsection{Derived categories of the second kind}
\label{derived-second-kind-subsecn}
 In this section we mostly recall the background material going back
to~\cite[Sections~2.1 and~4.1]{Psemi}, \cite[Sections~3.3
and~4.2]{Pkoszul}, \cite{Jor}, \cite[Section~2]{Kra}, \cite{Neem2},
\cite[Section~1.3]{Bec}, \cite[Section~5.1]{Pedg},
\cite[Sections~6\+-8]{PS5}, and~\cite[Appendix~A and
Section~B.7]{Pcosh}.
 See~\cite[Remark~9.2]{PS4} and~\cite[Section~7]{Pksurv} for
a historical and philosophical discussion.

 Given a DG\+category $\bA$, one can consider complexes of objects
of $\bA$ with closed morphisms of degree~$0$ in the role of
the differentials.
 Equivalently, these are just complexes in the  preadditive
category~$\sZ^0(\bA)$ of closed morphisms of degree~$0$ in~$\bA$.
 One can define the notions of the \emph{direct sum totalization}
and the \emph{direct product totalization} of a complex in~$\bA$.
 Both the totalizations, if they exist, are objects of~$\bA$.
 For a finite complex in $\bA$, the direct sum and direct product
totalizations agree, and both exist whenever $\bA$ has shifts and
cones~\cite[Section~1.2]{Pkoszul}, \cite[Section~1.3]{Pedg},
\cite[Section~1.6]{PS5}.

 Let $\bE$ be an exact DG\+category.
 Recall the notation $\sH^0(\bE)$ for the triangulated homotopy category
of~$\bE$ (see Section~\ref{dg-categories-of-cdg-modules-subsecn}).
 We are interested in (admissible) short exact sequences in the exact
category~$\sZ^0(\bE)$.
 Viewed as finite (three-term) complexes in $\sZ^0(\bE)$, such short
exact sequences have totalizations, which are objects of~$\bE$.
 An object $X\in\bE$ is said to be \emph{absolutely acyclic} if it
belongs to the minimal thick subcategory of $\sH^0(\bE)$ containing
the totalizations of short exact sequences in $\sZ^0(\bE)$.
 Equivalently, an object $X\in\bE$ is absolutely acyclic if and only
if it belongs to the minimal subcategory of the exact category
$\sZ^0(\bE)$ containing all contractible objects (i.~e., the objects
that vanish in $\sH^0(\bE)$) and closed under extensions and direct
summands~\cite[Proposition~8.12]{PS5}.
 The thick subcategory of absolutely acyclic objects is denoted by
$\Ac^\abs(\bE)\subset\sH^0(\bE)$.
 The triangulated Verdier quotient category
$$
 \sD^\abs(\bE)=\sH^0(\bE)/\Ac^\abs(\bE)
$$
is called the \emph{absolute derived category} of an exact
DG\+category~$\bE$.

 We refer to~\cite[Section~1.2]{Pkoszul}, \cite[Section~1.3]{Pedg},
or~\cite[Section~1.4]{PS5} for the definitions of infinite direct
sums (coproducts) and infinite products in DG\+categories.
 Any infinite direct sum (respectively, infinite product) in
a DG\+category $\bA$ is also an infinite direct sum (resp., product)
in the preadditive categories $\sZ^0(\bA)$ and $\sH^0(\bA)$;
so infinite direct sums (resp., products) exist in $\sZ^0(\bA)$
and $\sH^0(\bA)$ whenever they exist in~$\bA$.
 One says that infinite direct sums (resp., products) \emph{are
exact} in an exact DG\+category $\bE$ if they are exact in the exact
category $\sZ^0(\bE)$, or equivalently, in the exact category of
underlying graded objects $\sZ^0(\bE^\bec)$ \,\cite[Section~5.1]{Pedg}.

 Given an exact DG\+category $\bE$ where infinite direct sums exist
and are exact, one defines the thick subcategory of \emph{coacyclic
objects} (\emph{in the sense of Positselski}) $\Ac^\co(\bE)\subset
\sH^0(\bE)$ as the minimal triangulated subcategory of $\sH^0(\bE)$
containing the totalizations of short exact sequences in $\sZ^0(\bE)$
and closed under infinite direct sums.
 The triangulated Verdier quotient category
$$
 \sD^\co(\bE)=\sH^0(\bE)/\Ac^\co(\bE)
$$
is called the \emph{coderived category} (\emph{in the sense of
Positselski}) of the exact DG\+cat\-e\-gory~$\bE$.

 Dually, when $\bE$ is an exact DG\+category with exact infinite
products, the thick subcategory of \emph{contraacyclic objects}
(\emph{in the sense of Positselski}) $\Ac^\ctr(\bE)\subset\sH^0(\bE)$
is defined as the minimal triangulated subcategory of $\sH^0(\bE)$
containing the totalizations of short exact sequences in $\sZ^0(\bE)$
and closed under infinite products.
 The triangulated Verdier quotient category
$$
 \sD^\ctr(\bE)=\sH^0(\bE)/\Ac^\ctr(\bE)
$$
is called the \emph{contraderived category} (\emph{in the sense of
Positselski}) of the exact DG\+category~$\bE$.

 Let $\bE$ be an exact DG\+category.
 See Section~\ref{abelian-and-exact-dg-categs-of-cdg-modules-subsecn}
for the notation $\bE^\bec$, \ $\Phi_\bE$, and $\Psi^\pm_\bE$.
 An object $P\in\bE$ is said to be \emph{graded-projective} if
the object $\Phi_\bE(P)$ is projective in the exact category
$\sZ^0(\bE)$.
 Dually, an object $J\in\bE$ is said to be \emph{graded-injective} if
the object $\Phi_\bE(J)$ is injective in the exact category
$\sZ^0(\bE)$.
 The full DG\+subcategories of graded-projective and graded-injective
objects are denoted by $\bE_\bproj$ and $\bE^\binj\subset\bE$.

 An object $X\in\bE$ is said to be \emph{coacyclic in the sense of
Becker} if the complex of morphisms $\Hom^\bu_\bE(X,J)$ is acyclic
for all graded-injective objects $J$ in~$\bE$.
 The full subcategory of Becker-coacyclic objects $\Ac^\bco(\bE)
\subset\sH^0(\bE)$ contains all the totalizations of short exact
sequences in $\sZ^0(\bE)$ and is closed under infinite direct sums
in $\sH^0(\bE)$ (so one has $\Ac^\co(\bE)\subset\Ac^\bco(\bE)$ whenever
$\bE$ has exact infinite direct sums~\cite[Theorem~5.5(a)]{Pedg}).
  The triangulated Verdier quotient category
$$
 \sD^\bco(\bE)=\sH^0(\bE)/\Ac^\bco(\bE)
$$
is called the \emph{coderived category in the sense of Becker} of
the exact DG\+category~$\bE$.

 Dually, an object $Y\in\bE$ is said to be \emph{contraacyclic in
the sense of Becker} if the complex of morphisms $\Hom^\bu_\bE(P,Y)$
is acyclic for all graded-projective objects $P$ in~$\bE$.
 The full subcategory of Becker-contraacyclic objects $\Ac^\bctr(\bE)
\subset\sH^0(\bE)$ contains all the totalizations of short exact
sequences in $\sZ^0(\bE)$ and is closed under infinite products
in $\sH^0(\bE)$ (so one has $\Ac^\ctr(\bE)\subset\Ac^\bctr(\bE)$
whenever $\bE$ has exact infinite products~\cite[Theorem~5.5(b)]{Pedg}).
  The triangulated Verdier quotient category
$$
 \sD^\bctr(\bE)=\sH^0(\bE)/\Ac^\bctr(\bE)
$$
is called the \emph{contraderived category in the sense of Becker} of
the exact DG\+category~$\bE$.

\begin{thm} \label{under-star-condition-Positselski=Becker}
 Let\/ $\bE$ be an exact DG\+category with twists and infinite direct
sums.
 Assume that there are enough injective objects in the exact
category\/ $\sZ^0(\bE^\bec)$, and countable direct sums of injective
objects have finite injective dimensions in\/ $\sZ^0(\bE^\bec)$.
 Then the classes of Positselski-coacyclic and Becker-coacyclic
objects in\/ $\bE$ coincide.
 The composition of the identity inclusion of triangulated categories\/
$\sH^0(\bE^\binj)\rarrow\sH^0(\bE)$ and the triangulated Verdier
quotient functor\/ $\sH^0(\bE)\rarrow\sD^\co(\bE)$ is a triangulated
equivalence
$$
 \sH^0(\bE^\binj)\simeq\sD^\co(\bE)=\sD^\bco(\bE).
$$
\end{thm}

\begin{proof}
 This is~\cite[Theorem~5.10(a)]{Pedg}.
\end{proof}

 An abelian DG\+category $\bA$ is said to be \emph{Grothendieck} if
the abelian category $\sZ^0(\bA)$ is Grothendieck, or equivalently,
$\bA$ is a DG\+category with infinite direct sums and the abelian
category $\sZ^0(\bA^\bec)$ is
Grothendieck~\cite[Proposition~9.4]{Pedg},
\cite[Proposition~7.5]{PS5}.

\begin{thm} \label{coderived-of-Grothendieck-DG-category}
 Let\/ $\bA$ be a Grothendieck abelian DG\+category.
 Then the composition of the identity inclusion of triangulated
categories\/ $\sH^0(\bA^\binj)\rarrow\sH^0(\bA)$ and the triangulated
Verdier quotient functor\/ $\sH^0(\bA)\rarrow\sD^\bco(\bA)$ is
a triangulated equivalence
$$
 \sH^0(\bA^\binj)\simeq\sD^\bco(\bA).
$$
\end{thm}

\begin{proof}
 This is~\cite[Corollary~7.10]{PS5}.
 Notice that the theorem actually tells us more than it says.
 As one can see either from the proof of the theorem or just from
the definitions involved, the actual assertion is that, for
every object $A\in\bA$ there exists a graded-injective object
$J\in\bA$ together with a closed morphism $A\rarrow J$ of degree~$0$
whose cone is Becker-coacyclic in~$\bA$.
\end{proof}

 Given an exact category $\sE$ and a full additive subcategory
$\sF\subset\sE$, one says that the full subcategory $\sF$
\emph{inherits an exact structure} from $\sE$ is the class of all
short sequences in $\sF$ that are (admissible) exact in $\sE$ is
an exact category structure on~$\sF$.
 A characterization of full additive subcategories inheriting
an exact category structure can be found in~\cite[Theorem~2.6]{DS}
or~\cite[Lemma~4.21]{Pedg}.

 Let $\bE$ be an exact DG\+category and $\bF\subset\bE$ be a full
DG\+subcategory containing the zero object and closed under shifts
and cones.
 One says that the full DG\+subcategory $\bF$ \emph{inherits
an exact DG\+category structure} from $\bE$ if the full additive
subcategory $\sZ^0(\bF)\subset\sZ^0(\bE)$ inherits an exact category
structure from $\sZ^0(\bE)$ \emph{and} the full additive subcategory
$\sZ^0(\bF^\bec)\subset\sZ^0(\bE^\bec)$ inherits an exact category
structure from $\sZ^0(\bE^\bec)$.
  In this case, the DG\+category $\bF$ with the inherited exact
category structures on $\sZ^0(\bF)$ and $\sZ^0(\bF^\bec)$ becomes
an exact DG\+category itself~\cite[Section~4.5]{Pedg}.

 The former condition in the previous paragraph actually implies
the latter one~\cite[Proposition~4.24]{Pedg}.
 The latter condition implies the former one whenever the full
DG\+subcategory $\bF$ is closed under twists and direct summands
in~$\bE$ \,\cite[Corollary~4.25]{Pedg}.
 Here a full DG\+subcategory $\bF$ is said to be closed under direct
summands in $\bE$ if the full subcategory $\sZ^0(\bF)$ is closed
under direct summands in the (pre)additive category $\sZ^0(\bE)$
\,\cite[Section~1.3]{Pedg}.

 Let $\bE$ be an exact DG\+category and $\bF\subset\bE$ be a full
DG\+subcategory inheriting an exact DG\+category structure.
 We will say that the exact DG\+subcategory $\bF\subset\bE$ is
\emph{strict} if $F\in\bE$ and $\Phi_\bE(F)\in\sZ^0(\bF^\bec)$ imply
$F\in\bF$ (cf.\ the definition of a ``strict exact DG\+subpair''
in~\cite[Section~6.1]{Pedg}) and in
Section~\ref{exact-dg-pairs-subsecn} below).

 The definitions of a \emph{resolving subcategory} and
the \emph{resolution dimension} were given in
Section~\ref{co-resolution-dimensions-subsecn}.
 The following proposition is a counterpart of
Proposition~\ref{finite-resolutions} for derived categories of
the second kind in the context of exact DG\+categories.

\begin{prop} \label{second-kind-finite-resolutions}
 Let\/ $\bE$ be an exact DG\+category and\/ $\bF\subset\bE$ be
a strict exact DG\+subcategory.
 Assume that the full subcategory\/ $\sZ^0(\bF^\bec)$ is resolving in
the exact category\/ $\sZ^0(\bE^\bec)$, and the resolution dimensions
of all the objects of\/ $\sZ^0(\bE^\bec)$ with respect to\/
$\sZ^0(\bF^\bec)$ are bounded by a fixed constant~$d$.
 In this setting: \par
\textup{(a)} The inclusion of exact DG\+categories\/
$\bF\rarrow\bE$ induces an equivalence of the absolute derived
categories
$$
 \sD^\abs(\bF)\simeq\sD^\abs(\bE).
$$ \par
\textup{(b)} If infinite direct sums exist and are exact in\/ $\bE$,
and the full DG\+subcategory\/ $\bF$ is closed under infinite direct
sums in\/ $\bE$, then the inclusion of exact DG\+categories\/
$\bF\rarrow\bE$ induces an equivalence of the Positselski coderived
categories
$$
 \sD^\co(\bF)\simeq\sD^\co(\bE).
$$ \par
\textup{(c)} If infinite products exist and are exact in\/ $\bE$,
and the full DG\+subcategory\/ $\bF$ is closed under infinite products
in\/ $\bE$, then the inclusion of exact DG\+categories\/ $\bF\rarrow\bE$
induces an equivalence of the Positselski contraderived categories
$$
 \sD^\ctr(\bF)\simeq\sD^\ctr(\bE).
$$ \par
\textup{(d)} If the full DG\+subcategory\/ $\bF$ is closed under
direct summands in\/ $\bE$, then the inclusion of exact DG\+categories\/
$\bF\rarrow\bE$ induces an equivalence of the Becker contraderived
categories
$$
 \sD^\bctr(\bF)\simeq\sD^\bctr(\bE).
$$
\end{prop}

\begin{proof}
 Parts~(a\+-c) are particular cases of the respective assertions
of~\cite[Theorem~6.6(a)]{Pedg}, in view of~\cite[Example~6.1(0)]{Pedg}.
 To prove part~(d), notice that under its assumption the full
subcategory $\sZ^0(\bF^\bec)$ is closed under direct summands in
$\sZ^0(\bE^\bec)$ \,\cite[Section~4.5]{Pedg}.
 Since $\sZ^0(\bF^\bec)$ is assumed to be a resoving subcategory
in $\sZ^0(\bE^\bec)$, it follows that the classes of projective
objects in $\sZ^0(\bE^\bec)$ and $\sZ^0(\bF^\bec)$ coincide; in
particular, all the projective objects of $\sZ^0(\bE^\bec)$ belong
to $\sZ^0(\bF^\bec)$ (see~\cite[first paragaraph of the proof of
Proposition~B.7.9]{Pcosh}).
 Since the exact DG\+subcategory $\bF\subset\bE$ is assumed to be
strict, it follows that the classes of graded-projective objects in
$\bE$ and $\bF$ coincide as well.
 Thus an object of $\bF$ is Becker-contraacyclic in $\bF$ if and only if
it is Becker-contraacyclic in~$\bE$.
 It remains to point out that, by~\cite[beginning of the proof
of Theorem~6.6]{Pedg}, for every object $E\in\bE$ there exists
an object $F\in\bF$ together with a closed morphism $F\rarrow E$ of
degree~$0$ whose cone is absolutely acyclic in~$\bE$.
 Furthermore, all absolutely acyclic objects in $\bE$ are
Becker-contraacyclic, as it is clear from the discussion above
in this section.
 So one can refer to the well-known category-theoretic
lemma~\cite[Proposition~10.2.7(ii)]{KS}, \cite[Lemma~1.6(a)]{Pkoszul},
\cite[Lemma~6.7]{Pedg}, or~\cite[Lemma~A.3.3(a)]{Pcosh}.
\end{proof}

 For a version of Proposition~\ref{second-kind-finite-resolutions}(b)
for the Becker coderived categories (or rather, the dual assertion
to it), see Proposition~\ref{becker-contraderived-finite-coresolutions}
below.

\begin{prop} \label{second-kind-infinite-resolutions}
 Let\/ $\bE$ be an exact DG\+category and\/ $\bF\subset\bE$ be
a DG\+subcategory inheriting an exact DG\+category structure.
 Assume that the full subcategory\/ $\sZ^0(\bF^\bec)$ is resolving
in the exact category\/ $\sZ^0(\bE^\bec)$.
 Assume further that infinite products exist and are exact in\/ $\bE$,
and the full DG\+subcategory\/ $\bF$ is closed under infinite products
in\/~$\bE$.
 Assume also that all twists exist in the DG\+category\/ $\bE$, and
the full DG\+subcategory\/ $\bF\subset\bE$ is closed under twists.
 In this setting: \par
\textup{(a)} The inclusion of exact DG\+categories\/ $\bF\rarrow\bE$
induces an equivalence of the Positselski contraderived categories
$$
 \sD^\ctr(\bF)\simeq\sD^\ctr(\bE).
$$ \par
\textup{(b)} If the full DG\+subcategory\/ $\bF$ is strict and closed
under direct summands in\/ $\bE$, then the inclusion of exact
DG\+categories\/ $\bF\rarrow\bE$ induces an equivalence of the Becker
contraderived categories
$$
 \sD^\bctr(\bF)\simeq\sD^\bctr(\bE).
$$
\end{prop}

\begin{proof}
 Part~(a) is a particular case of~\cite[Theorem~7.11(a)]{Pedg},
in view of~\cite[Example~6.1(0)]{Pedg}.
 The proof of part~(b) is similar to the proof of
Proposition~\ref{second-kind-finite-resolutions}(d) and based on
the fact that all Positselski-contraacyclic objects in $\bE$ are
Becker-contraacyclic together with the construction of the resolution
from~\cite[proof of Theorem~7.11]{Pedg}.
\end{proof}

 Let $\bE$ and $\bF$ be two DG\+categories.
 A DG\+functor $G\:\bF\rarrow\bE$ is said to be \emph{right adjoint}
to a DG\+functor $H\:\bE\rarrow\bF$ if, for every pair of objects
$E\in\bE$ and $F\in\bF$, there is a given isomorphism of complexes of
abelian groups $\Hom^\bu_\bE(E,G(F))\simeq\Hom^\bu_\bF(H(E),F)$
whose underlying isomorphism of graded abelian groups
$\Hom^*_\bE(E,G(F))\simeq\Hom^*_\bF(H(E),F)$ is functorial with
respect to all (and not only closed) homogeneous morphisms of all
cohomological degrees in $\bE$ and~$\bF$.

 Any adjoint pair of DG\+functors $(H,G)$ between DG\+categories
$\bE$ and $\bF$ induces an adjoint pair of additive functors
$\sZ^0(H)$ and $\sZ^0(G)$ between the preadditive categories
$\sZ^0(\bE)$ and $\sZ^0(\bF)$, as well as an adjoint pair of
additive functors $\sH^0(H)$ and $\sH^0(G)$ between the preadditive
categories $\sH^0(\bE)$ and~$\sH^0(\bF)$.
 Any adjoint pair of DG\+functors $(H,G)$ induces also an adjoint
pair of DG\+functors $H^\bec$ and $G^\bec$ between the DG\+categories
$\bE^\bec$ and~$\bF^\bec$.

\begin{lem} \label{DG-functors-preserve-Becker-co-contra-acyclicity}
 Let\/ $\bE$ and\/ $\bF$ be exact DG\+categories, and let
$G\:\bF\rarrow\bE$ be an exact DG\+functor (in the sense of
the definitions in Section~%
\textup{\ref{abelian-and-exact-dg-categs-of-cdg-modules-subsecn}}
and~\cite[Section~4.4]{Pedg}).
 In this context: \par
\textup{(a)} if a right adjoint DG\+functor to $G$ exists, then $G$
takes Becker-coacyclic objects of\/ $\bF$ to Becker-coacyclic objects
of\/~$\bE$; \par
\textup{(b)} if a left adjoint DG\+functor to $G$ exists, then $G$
takes Becker-contraacyclic objects of\/ $\bF$ to Becker-contraacyclic
objects of\/~$\bE$.
\end{lem}

\begin{proof}
 This is a DG\+category generalization of~\cite[Lemma~B.7.5]{Pcosh}.
 Let us prove part~(b).
 Let $H\:\bE\rarrow\bF$ be the left adjoint DG\+functor to~$G$.
 Then the additive functor $\sZ^0(H^\bec)\:\sZ^0(\bE^\bec)\rarrow
\sZ^0(\bF^\bec)$ is left adjoint to the exact functor $\sZ^0(G^\bec)\:
\sZ^0(\bF^\bec)\allowbreak\rarrow\sZ^0(\bE^\bec)$.
 Hence the functor $\sZ^0(H^\bec)$ takes projective objects in
the exact category $\sZ^0(\bE^\bec)$ to projective objects in
the exact category $\sZ^0(\bF^\bec)$.
 In view of commutativity of a suitable diagram
from~\cite[Section~3.4]{Pedg}, it follows the additive functor
$\sZ^0(H)\:\sZ^0(\bE)\allowbreak\rarrow\sZ^0(\bF)$ takes
graded-projective objects to graded-projective objects.

 Now let $Y\in\bF$ be a Becker-contraacyclic object and $P\in\bE$
be a graded-projective object.
 We need to show that the complex $\Hom_\bE^\bu(P,G(Y))$ is acyclic.
 Indeed, we have $\Hom_\bE^\bu(P,G(Y))\simeq\Hom_\bF^\bu(H(P),Y)$,
and the object $H(P)\in\bF$ is graded-projective according to
the previous paragraph.
\end{proof}

\begin{lem} \label{Grothendick-abelian-DG-functor-preserves-Becker-co}
 Let\/ $\bA$ and\/ $\bB$ be Grothendieck abelian DG\+categories, and
let $F\:\bA\rarrow\bB$ be an exact DG\+functor preserving infinite
direct sums.
 Then the DG\+functor $F$ takes Becker-coacyclic objects of\/ $\bA$ to
Becker-coacyclic objects of\/~$\bB$.
\end{lem}

\begin{proof}
 This is a DG\+category generalization of~\cite[Lemma~A.5]{Psemten}.
 The functor of Grothendieck abelian categories $\sZ^0(F)\:
\sZ^0(\bA)\rarrow\sZ^0(\bB)$ is exact and preserves infinite direct
sums; hence it preserves all inductive limits.
 According to~\cite[Corollary~7.17]{PS5}, the class of
Becker-coacyclic objects in $\sZ^0(\bA)$ is the closure of the class
of all contractible objects under extensions and direct limits;
and the same for~$\bB$.
 It remains to point out that any DG\+functor takes contractible objects
to contractible objects.
\end{proof}

\subsection{Coderived category of quasi-coherent CDG-modules}
\label{coderived-qcoh-cdg-modules-subsecn}
 We start with a brief discussion of the absolute derived categories
of quasi-coherent CDG\+modules before turning to the Becker coderived
categories.

\begin{cor} \label{qcoh-cta-cdg-absolute-derived-equivalence}
 Let $X$ be a quasi-compact semi-separated scheme and $\cB^\cu$ be
a quasi-coherent CDG\+quasi-algebra over~$X$.
 Then the inclusion of exact/abelian DG\+categories
$\cB^\cu\bQcoh^{X\dcta}\rarrow\cB^\cu\bQcoh$ induces an equivalence
of the absolute derived categories
$$
 \sD^\abs(\cB^\cu\bQcoh^{X\dcta})\simeq\sD^\abs(\cB^\cu\bQcoh).
$$
\end{cor}

\begin{proof}
 It is clear from the discussion in
Section~\ref{exact-dg-of-cta-cot-qcoh-cdg} (see
Corollary~\ref{exact-dg-categories-of-qcoh-cta-cot-cdg-modules})
that $\cB^\cu\bQcoh^{X\dcta}$ (as well as $\cB^\cu\bQcoh^{X\dcot}$
and $\cB^\cu\bQcoh^{\cB^*\dcot}$) is a strict exact DG\+subcategory
in $\cB^\cu\bQcoh$ in the sense of
Section~\ref{derived-second-kind-subsecn}.
 In view of the graded version of
Lemma~\ref{X-A-cta-cot-coresolving-coresolution-dimension}(b),
the dual version of Proposition~\ref{second-kind-finite-resolutions}(a)
is applicable.
\end{proof}

 For a similar equivalence involving also the absolute derived
category of $X$\+co\-tor\-sion quasi-coherent CDG\+modules
$\sD^\abs(\cB^\cu\bQcoh^{X\dcot})$, which holds under more restrictive
assumption on the scheme $X$, see
Corollary~\ref{qcoh-cot-cdg-absolute-derived-equivalence} below.

 Notice that the following corollary is \emph{not} applicable to
the Positselski coderived categories, as the DG\+categories
$\cB^\cu\bQcoh^{\cB^*\dcot}$, \ $\cB^\cu\bQcoh^{X\dcot}$, and
$\cB^\cu\bQcoh^{X\dcta}$ do \emph{not} have infinite direct sums.
 The point is that the classes of contraadjusted and cotorsion
modules/sheaves are not preserved by infinite direct sums in
the respective Grothendieck categories of modules/sheaves.

\begin{cor} \label{qcoh-cta-cot-cdg-becker-coderived-equivalence}
 Let $X$ be a quasi-compact semi-separated scheme and $\cB^\cu$ be
a quasi-coherent CDG\+quasi-algebra over~$X$.
 Then the inclusions of DG\+category/exact DG\+categories/abelian
DG\+category $\cB^\cu\bQcoh^\binj\rarrow\cB^\cu\bQcoh^{\cB^*\dcot}
\rarrow\cB^\cu\bQcoh^{X\dcot}\allowbreak\rarrow\cB^\cu\bQcoh^{X\dcta}
\rarrow\cB^\cu\bQcoh$ induce triangulated equivalences of the homotopy
category and Becker coderived categories \hfuzz=11pt
\begin{align*}
 \sH^0(\cB^\cu\bQcoh^\binj)\simeq
 \sD^\bco(\cB^\cu\bQcoh^{\cB^*\dcot}) &\simeq
 \sD^\bco(\cB^\cu\bQcoh^{X\dcot}) \\  &\simeq
 \sD^\bco(\cB^\cu\bQcoh^{X\dcta})\simeq\sD^\bco(\cB^\cu\bQcoh).
\end{align*}
\end{cor}

\begin{proof}
 The point is that graded-injective quasi-coherent CDG\+modules over
$\cB^\cu$ are $\cB^*$\+cotorsion (by the definition), hence also
$X$\+cotorsion and $X$\+contraadjusted.
 Furthermore, the classes of graded-injective objects in all the four
exact/abelian DG\+categories $\cB^\cu\bQcoh^{\cB^*\dcot}$, \
$\cB^\cu\bQcoh^{X\dcot}$, \ $\cB^\cu\bQcoh^{X\dcta}$, and
$\cB^\cu\bQcoh$ coincide.
 Indeed, $\cB^*\Qcoh^{\cB^*\dcot}\subset\cB^*\Qcoh^{X\dcot}\subset
\cB^*\Qcoh^{X\dcta}$ are coresolving subcategories closed under direct
summands in the abelian category of quasi-coherent graded modules
$\cB^*\Qcoh$ (cf.\
Lemma~\ref{X-A-cta-cot-coresolving-coresolution-dimension}(a)), so
the classes of injective objects in all the four exact/abelian
categories $\cB^*\Qcoh^{\cB^*\dcot}\subset\cB^*\Qcoh^{X\dcot}\subset
\cB^*\Qcoh^{X\dcta}\subset\cB^*\Qcoh$ coincide (and there are enough
injective objects in each of the four categories).
 Therefore, an object from the respective exact DG\+category is
Becker-coacyclic in $\cB^\cu\bQcoh^{\cB^*\dcot}$, \
$\cB^\cu\bQcoh^{X\dcot}$, \ $\cB^\cu\bQcoh^{X\dcta}$ if and only if
it is Becker-coacyclic in $\cB^\cu\bQcoh$.
 It remains to refer to
Theorem~\ref{coderived-of-Grothendieck-DG-category} and the discussion
in its proof, together with the category-theoretic
lemma~\cite[Proposition~10.2.7(i)]{KS},
\cite[Lemma~1.6(b)]{Pkoszul}, or~\cite[Lemma~A.3.3(b)]{Pcosh}.
\hbadness=1500
\end{proof}

 Let $X$ be a Noetherian scheme.
 A quasi-coherent graded quasi-algebra $\cB^*$ over $X$ is said to be
\emph{left Noetherian} if the graded ring $\cB^*(U)$ is left
Noetherian for every affine open subscheme $U\subset X$.
 It suffices to check this condition for affine open subschemes
belonging to any chosen affine open covering of the scheme~$X$.
 For a left Noetherian quasi-coherent graded quasi-algebra $\cB^*$
over $X$, the Grothendieck category of quasi-coherent graded left
$\cB^*$\+modules $\cB^*\Qcoh$ is locally Noetherian.

\begin{cor} \label{noetherian-quasi-algebra-Positselski=Becker}
 Let $X$ be a Noetherian scheme and $\cB^\cu$ be a quasi-coherent
quasi-algebra over $X$ whose underlying graded quasi-algebra $\cB^*$
is left Noetherian.
 Then the classes of Positselski-coacyclic and Becker-coacyclic
objects in the abelian DG\+category $\cB^\cu\bQcoh$ coincide.
 The inclusion of DG\+categories $\cB^\cu\bQcoh^\binj\allowbreak
\rarrow\cB^\cu\bQcoh$ induces a triangulated equivalence
$$
 \sH^0(\cB^\cu\bQcoh^\binj)\simeq\sD^\co(\cB^\cu\bQcoh)
 =\sD^\bco(\cB^\cu\bQcoh).
$$
\end{cor}

\begin{proof}
 This is a special case of
Theorem~\ref{under-star-condition-Positselski=Becker} and
a slight generalization of~\cite[Lemma~1.7(b)]{EP}.
 Notice that, in a locally Noetherian Grothendieck category, infinite
direct sums of injective objects are injective; so the assumption of
Theorem~\ref{under-star-condition-Positselski=Becker} is satisfied.
\end{proof}

\subsection{Contraderived categories of contraherent CDG-modules}
\label{contraderived-lcth-cdg-modules-subsecn}
 In this section we restrict ourselves to proving that the absolute 
derived and contraderived categories of CDG\+contramodules do not
depend on the open covering $\bW$ of a quasi-compact semi-separated
scheme $X$, and remain unchanged when the antilocality condition is
imposed.

\begin{cor} \label{contraderived-indep-of-covering-or-antilocal}
 Let $X$ be a quasi-compact semi-separated scheme with an open
covering\/ $\bW$ and $\cB^\cu$ be a quasi-coherent CDG\+quasi-algebra
over~$X$.
 Let\/ $\st=\abs$, $\ctr$, or\/~$\bctr$ be an absolute derived or
contraderived category symbol.
 In this context: \par
\textup{(a)} The inclusions of exact DG\+categories
$\cB^\cu\bCtrh_\al\rarrow\cB^\cu\bCtrh\rarrow\cB^\cu\bLcth_\bW$
induce equivalences of triangulated categories \hbadness=1650
$$
 \sD^\st(\cB^\cu\bCtrh_\al)\simeq\sD^\st(\cB^\cu\bCtrh)
 \simeq\sD^\st(\cB^\cu\bLcth_\bW).
$$ \par
\textup{(b)} The inclusions of exact DG\+categories
$\cB^\cu\bCtrh^{X\dlct}_\al\rarrow\cB^\cu\bCtrh^{X\dlct}\rarrow
\cB^\cu\bLcth_\bW^{X\dlct}$ induce equivalences of triangulated
categories
$$
 \sD^\st(\cB^\cu\bCtrh^{X\dlct}_\al)\simeq
 \sD^\st(\cB^\cu\bCtrh^{X\dlct})\simeq
 \sD^\st(\cB^\cu\bLcth_\bW^{X\dlct}).
$$ \par
\textup{(c)} The inclusions of exact DG\+categories
$\cB^\cu\bCtrh^{\cB^*\dlct}_\al\rarrow\cB^\cu\bCtrh^{\cB^*\dlct}\rarrow
\cB^\cu\bLcth_\bW^{\cB^*\dlct}$ induce equivalences of triangulated
categories
$$
 \sD^\st(\cB^\cu\bCtrh^{\cB^*\dlct}_\al)\simeq
 \sD^\st(\cB^\cu\bCtrh^{\cB^*\dlct})\simeq
 \sD^\st(\cB^\cu\bLcth_\bW^{\cB^*\dlct}).
$$
\end{cor}

\begin{proof}
 It is clear from the discussion in
Section~\ref{exact-dg-categories-of-antilocal-ctrh-subsecn} (see
Corollary~\ref{exact-dg-categories-of-antilocal-ctrh-cdg-modules})
that all the full DG\+subcategories in question are strict exact
DG\+subcategories in the sense of
Section~\ref{derived-second-kind-subsecn}.
 Moreover, all the full DG\+subcategories in the situation at hand
are closed under direct summands in each other.
 In view of the graded version of
Lemma~\ref{antilocal-resolving-resolution-dimension},
Proposition~\ref{second-kind-finite-resolutions}(a,c,d) is applicable.
\end{proof}

 For a triangulated equivalence between the absolute
derived/contraderived categories appearing in parts~(a) and~(b)
of Corollary~\ref{contraderived-indep-of-covering-or-antilocal}
(proved under more restrictive assumptions on the scheme~$X$),
see Theorem~\ref{lcth-lcta-lct-cdg-abs-contraderived-equiv-thm} below.

\subsection{Semiderived categories of quasi-coherent $\cA$-modules}
\label{semiderived-quasi-coherent-subsecn}
 The notion of a semiderived category goes back to
the book~\cite{Psemi}, where (as well as in the subsequent
book~\cite{Psemten}) it played a key role.
 See also the paper~\cite[Sections~5\+-6]{Pfp} and the book
manuscript~\cite[Chapter~8]{Pcosh}.
 In the context of relative derived Koszul duality, the semiderived
categories appeared in the book~\cite[Section~6.3 and Theorem~6.14,
Section~7.3 and Theorem~7.11]{Prel}.

 Let $X$ be a scheme and $\cA$ be a quasi-coherent quasi-algebra
over~$X$.
 The semiderived category of (complexes of) quasi-coherent
$\cA$\+modules is an intermediate triangulated Verdier quotient
category between the derived and the coderived category of
quasi-coherent $\cA$\+modules, constructed using the notion of
a coacyclic complex of quasi-coherent sheaves (of $\cO_X$\+modules)
on~$X$.
 Since there are two notions of a coacyclic complex of quasi-coherent
sheaves, the Positselski and the Becker one, there are also two
notions of a semiderived category.

 Let us introduce some notation.
 A quasi-coherent quasi-algebra $\cA$ can be viewed as a special case
of a quasi-coherent CDG\+quasi-algebra $\cA^\cu=\cA$; then
the quasi-coherent CDG\+modules over $\cA^\cu$ are just complexes of
quasi-coherent $\cA$\+modules.
 On the other hand, for any additive category $\sA$, we have
the DG\+category $\bCom(\sA)$ of (unbounded) complexes in~$\sA$.
 For an exact category $\sE$, the DG\+category $\bCom(\sE)$ is naturally
an exact DG\+category (see~\cite[Example~4.40]{Pedg}).
 In the context of a quasi-coherent quasi-algebra $\cA$ over $X$,
we have $\cA^\cu\bQcoh=\bCom(\cA\Qcoh)$, where $\cA\Qcoh$ is
the abelian category of quasi-coherent $\cA$\+modules.
 Similarly, we have $\cA^\cu\bQcoh^{X\dcta}=\bCom(\cA\Qcoh^{X\dcta})$,
\ $\cA^\cu\bQcoh^{X\dcot}=\bCom(\cA\Qcoh^{X\dcot})$, and
$\cA^\cu\bQcoh^{\cA^*\dcot}=\bCom(\cA\Qcoh^{\cA\dcot})$
(where $\cA^*$ denotes the quasi-coherent quasi-algebra $\cA$ viewed
as a quasi-coherent graded quasi-algebra concentrated in degree~$0$).

 Given an additive category $\sA$, we denote by $\Com(\sA)=
\sZ^0(\bCom(\sA))$ the additive category of complexes in $\sA$ and
closed morphisms of degree zero between them.
 We also denote by $\Hot(\sA)=\sH^0(\bCom(\sA))$ the triangulated
homotopy category of complexes in $\sA$, constructed as the additive
quotient category of $\Com(\sA)$ by the ideal of morphisms
cochain homotopic to zero.
 For an exact category $\sE$, as particular cases of the definition
in Section~\ref{derived-second-kind-subsecn}, we obtain the thick
subcategories of absolutely acyclic and coacyclic complexes
$$
 \Ac^\abs(\sE)\subset\Ac^\co(\sE)\subset\Ac^\bco(\sE)\subset
 \Hot(\sE),
$$
where the full subcategory of Positselski-coacyclic complexes
$\Ac^\co(\sE)$ is only well-defined when infinite direct sums exist
and are exact in~$\sE$.
 The related exotic derived categories are the absolute derived
category $\sD^\abs(\sE)=\Hot(\sE)/\Ac^\abs(\sE)$, the Positselski
coderived category $\sD^\co(\sE)=\Hot(\sE)/\Ac^\co(\sE)$, and
the Becker coderived category $\sD^\bco(\sE)=\Hot(\sE)/\Ac^\bco(\sE)$.

 The following three lemmas play an important role.

\begin{lem} \label{forgetful-functor-preserves-coacyclicity}
  Let $X$ be a scheme and $\cA$ be a quasi-coherent quasi-algebra
over~$X$.
 In this context: \par
\textup{(a)} the forgetful functor $\cA\Qcoh\rarrow X\Qcoh$ takes
Positselski-coacyclic complexes in $\cA\Qcoh$ to Positselski-coacyclic
complexes in $X\Qcoh$; \par
\textup{(b)} the forgetful functor $\cA\Qcoh\rarrow X\Qcoh$ takes
Becker-coacyclic complexes in $\cA\Qcoh$ to Becker-coacyclic
complexes in $X\Qcoh$; \par
\textup{(c)} if the scheme $X$ is quasi-compact and semi-separated,
then the forgetful functor $\cA\Qcoh^{X\dcta}\rarrow X\Qcoh^\cta$ takes
Becker-coacyclic complexes in $\cA\Qcoh^{X\dcta}$ to Becker-coacyclic
complexes in $X\Qcoh^\cta$; \par
\textup{(d)} if the scheme $X$ is quasi-compact and semi-separated,
then the forgetful functor $\cA\Qcoh^{X\dcot}\rarrow X\Qcoh^\cot$ takes
Becker-coacyclic complexes in $\cA\Qcoh^{X\dcot}$ to Becker-coacyclic
complexes in $X\Qcoh^\cot$; \par
\textup{(e)} if the scheme $X$ is quasi-compact and semi-separated,
then the forgetful functor $\cA\Qcoh^{\cA\dcot}\rarrow X\Qcoh^\cot$
takes Becker-coacyclic complexes in $\cA\Qcoh^{\cA\dcot}$ to
Becker-coacyclic complexes in $X\Qcoh^\cot$.
\end{lem}

\begin{proof}
 Part~(a) holds because the forgetful functor $\cA\Qcoh\rarrow X\Qcoh$
is exact and preserves infinite direct sums.
 In view of these observations, part~(b) is a special case
of~\cite[Lemma~A.5]{Psemten} (which, in turn, is a special case
of~\cite[Lemma~B.7.5(a)]{Pcosh}); see also
Lemmas~\ref{DG-functors-preserve-Becker-co-contra-acyclicity}(a)
and~\ref{Grothendick-abelian-DG-functor-preserves-Becker-co} above.
 Parts~(c\+-e) follow from part~(b) because a complex in any one of
the categories $\cA\Qcoh^{X\dcta}$, \ $\cA\Qcoh^{X\dcot}$, or
$\cA\Qcoh^{\cA\dcot}$ is Becker-coacyclic if and only if it is
Becker-coacyclic in $\cA\Qcoh$, and similarly for $X\Qcoh$ (see
Corollary~\ref{qcoh-cta-cot-cdg-becker-coderived-equivalence}).
\end{proof}

\begin{lem} \label{becker-coacyclic-complexes-are-acyclic}
 Let $X$ be a scheme and $\cA$ be a quasi-coherent quasi-algebra
over~$X$.
 In this context: \par
\textup{(a)} any Positselski-coacyclic complex in the abelian category
$\cA\Qcoh$ is acyclic in $\cA\Qcoh$; \par
\textup{(b)} any Becker-coacyclic complex in the abelian category
$\cA\Qcoh$ is acyclic in $\cA\Qcoh$; \par
\textup{(c)} if the scheme $X$ is quasi-compact and semi-separated,
then any Becker-coacyclic complex in the exact category
$\cA\Qcoh^{X\dcta}$ is acyclic in $\cA\Qcoh^{X\dcta}$; \par
\textup{(d)} if the scheme $X$ is quasi-compact and semi-separated,
then any Becker-coacyclic complex in the exact category
$\cA\Qcoh^{X\dcot}$ is acyclic in $\cA\Qcoh^{X\dcot}$; \par
\textup{(e)} if the scheme $X$ is quasi-compact and semi-separated,
then any Becker-coacyclic complex in the exact category
$\cA\Qcoh^{\cA\dcot}$ is acyclic in $\cA\Qcoh^{\cA\dcot}$.
\end{lem}

\begin{proof}
 Part~(a) obviously holds for any exact category with exact functors
of infinite direct sums.
 Part~(b), which is a stronger version of part~(a)
(as all Positselski-coacyclic complexes are Becker-coacyclic),
holds for any abelian category with enough injective objects;
see~\cite[Lemma~A.2]{Psemten} and~\cite[Lemma~B.7.3(a) and
Remark~B.7.4]{Pcosh}.
 Parts~(c\+-e) are deduced from part~(b) in the following way.
 A complex in either of the exact categories in question is
Becker-coacyclic if and only if it is Becker-coacyclic in $\cA\Qcoh$
(see Corollary~\ref{qcoh-cta-cot-cdg-becker-coderived-equivalence}).
 Now the point is that any complex in $\cA\Qcoh^{X\dcta}$,
$\cA\Qcoh^{X\dcot}$, or $\cA\Qcoh^{\cA\dcot}$ that is acyclic in
$\cA\Qcoh$ is also acyclic in the respective exact subcategory
$\cA\Qcoh^{X\dcta}$, $\cA\Qcoh^{X\dcot}$, or $\cA\Qcoh^{\cA\dcot}$.
 See Theorems~\ref{qcoh-cotorsion-periodicity}\+-%
\ref{qcoh-quasi-algebra-cotorsion-periodicity} and the proof
of Lemma~\ref{qcoh-G-plus-periodicity-lemma}.
\end{proof}

\begin{lem} \label{qcoh-forgetful-functor-reflects-acyclicity}
 Let $X$ be a scheme and $\cA$ be a quasi-coherent quasi-algebra
over~$X$.
 In this context: \par
\textup{(a)} a complex in $\cA\Qcoh$ is acyclic if and only if it is
acyclic as a complex in $X\Qcoh$; \par
\textup{(b)} a complex in $\cA\Qcoh^{X\dcta}$ is acyclic if and only if
it is acyclic as a complex in $X\Qcoh^\cta$; \par
\textup{(c)} a complex in $\cA\Qcoh^{X\dcot}$ is acyclic if and only if
it is acyclic as a complex in $X\Qcoh^\cot$; \par
\textup{(d)} if the scheme $X$ is quasi-compact and semi-separated,
then a complex in $\cA\Qcoh^{\cA\dcot}$ is acyclic if and only
if it is acyclic as a complex in $X\Qcoh^\cot$.
\end{lem}

\begin{proof}
 Part~(a) holds because the forgetful functor of abelian categories
$\cA\Qcoh\rarrow X\Qcoh$ is exact and faithful.
 Parts~(b\+-c) follow by the definition of the exact categories
$\cA\Qcoh^{X\dcta}$ and $\cA\Qcoh^{X\dcot}$.
 To prove part~(d), one needs to use
Theorem~\ref{qcoh-quasi-algebra-cotorsion-periodicity}; cf.\ the proof
of the previous Lemma~\ref{becker-coacyclic-complexes-are-acyclic}.
\end{proof}

 Now let us spell out the definitions.
 A complex in the abelian category $\cA\Qcoh$ is said to be
\emph{Positselski-semicoacyclic} (relative to $X\Qcoh$) if,
\emph{viewed as a complex in $X\Qcoh$}, it is Positselski-coacyclic
in the abelian category $X\Qcoh$.
 The full subcategory of Positselski-semicoacyclic complexes is denoted
by $\Ac^\si(\cA\Qcoh)\subset\Hot(X\Qcoh)$.
 By Lemma~\ref{forgetful-functor-preserves-coacyclicity}(a),
all Positselski-coacyclic complexes in $\cA\Qcoh$ are
Positselski-semicoacyclic.
  By Lemma~\ref{becker-coacyclic-complexes-are-acyclic}(a)
(for the quasi-coherent quasi-algebra $\cO_X$ over~$X$), all
Positselski-coacyclic complexes in $X\Qcoh$ are acyclic.
 Therefore, it follows from
Lemma~\ref{qcoh-forgetful-functor-reflects-acyclicity}(a) that all
Positselski-semicoacyclic complexes in $\cA\Qcoh$ are acyclic.
 So we have
$$
 \Ac^\co(\cA\Qcoh)\subset\Ac^\si(\cA\Qcoh)\subset\Ac(\cA\Qcoh),
$$
where, for any exact category $\sE$, the notation $\Ac(\sE)\subset
\Hot(\sE)$ stands for the full subcategory of acyclic complexes.
 The triangulated Verdier quotient category
$$
 \sD^\si(\cA\Qcoh)=\Hot(\cA\Qcoh)/\Ac^\si(\cA\Qcoh)
$$
is called the \emph{Positselski semi}(\emph{co})\emph{derived category}
of quasi-coherent left $\cA$\+modules on~$X$.

 Similarly, a complex in the abelian category $\cA\Qcoh$ is said to be
\emph{Becker-semicoacyclic} (relative to $X\Qcoh$) if, \emph{viewed as
a complex in $X\Qcoh$}, it is Becker-coacyclic in the abelian
category $X\Qcoh$.
 The full subcategory of Becker-semicoacyclic complexes is denoted by
$\Ac^\bsi(\cA\Qcoh)\subset\Hot(X\Qcoh)$.
 By Lemma~\ref{forgetful-functor-preserves-coacyclicity}(b),
all Becker-coacyclic complexes in $\cA\Qcoh$ are Becker-semicoacyclic.
 By Lemma~\ref{becker-coacyclic-complexes-are-acyclic}(b)
(for the quasi-coherent quasi-algebra $\cO_X$ over~$X$), all
Becker-coacyclic complexes in $X\Qcoh$ are acyclic.
 Therefore, it follows from
Lemma~\ref{qcoh-forgetful-functor-reflects-acyclicity}(a) that all
Becker-semicoacyclic complexes in $\cA\Qcoh$ are acyclic.
 So we have {\hbadness=1650
$$
 \Ac^\bco(\cA\Qcoh)\subset\Ac^\bsi(\cA\Qcoh)\subset\Ac(\cA\Qcoh).
$$
 The} triangulated Verdier quotient category
$$
 \sD^\bsi(\cA\Qcoh)=\Hot(\cA\Qcoh)/\Ac^\bsi(\cA\Qcoh)
$$
is called the \emph{Becker semi}(\emph{co})\emph{derived category}
of quasi-coherent left $\cA$\+modules on~$X$.

 Assume that $X$ is a quasi-compact semi-separated scheme.
 A complex in the exact category $\cA\Qcoh^{X\dcta}$ is said to be
\emph{Becker-semicoacyclic} (relative to $X\Qcoh^\cta$) if,
\emph{viewed as a complex in $X\Qcoh^\cta$}, it is Becker-coacyclic
in the exact category $X\Qcoh^\cta$.
 The full subcategory of Becker-semicoacyclic complexes is denoted by
$\Ac^\bsi(\cA\Qcoh^{X\dcta})\subset\Hot(X\Qcoh^{X\dcta})$.
 By Lemma~\ref{forgetful-functor-preserves-coacyclicity}(c),
all Becker-coacyclic complexes in $\cA\Qcoh^{X\dcta}$ are
Becker-semicoacyclic.
 By Lemma~\ref{becker-coacyclic-complexes-are-acyclic}(c) (for
the quasi-coherent quasi-algebra $\cO_X$ over~$X$), all Becker-coacyclic
complexes in $X\Qcoh^\cta$ are acyclic in $X\Qcoh^\cta$.
 Therefore, it follows from
Lemma~\ref{qcoh-forgetful-functor-reflects-acyclicity}(b) that all
Becker-semicoacyclic complexes in $\cA\Qcoh^{X\dcta}$ are acyclic
in $\cA\Qcoh^{X\dcta}$.
 So we have
$$
 \Ac^\bco(\cA\Qcoh^{X\dcta})\subset\Ac^\bsi(\cA\Qcoh^{X\dcta})
 \subset\Ac(\cA\Qcoh^{X\dcta}).
$$
 The triangulated Verdier quotient category
$$
 \sD^\bsi(\cA\Qcoh^{X\dcta})=
 \Hot(\cA\Qcoh^{X\dcta})/\Ac^\bsi(\cA\Qcoh^{X\dcta})
$$
is called the \emph{Becker semi}(\emph{co})\emph{derived category}
of $X$\+contraadjusted quasi-coherent left $\cA$\+modules on~$X$.

 Similarly, a complex in the exact category $\cA\Qcoh^{X\dcot}$ is said
to be \emph{Becker-semicoacyclic} (relative to $X\Qcoh^\cot$) if,
\emph{viewed as a complex in $X\Qcoh^\cot$}, it is Becker-coacyclic
in the exact category $X\Qcoh^\cot$.
 The full subcategory of Becker-semicoacyclic complexes is denoted by
$\Ac^\bsi(\cA\Qcoh^{X\dcot})\subset\Hot(X\Qcoh^{X\dcot})$.
 By Lemma~\ref{forgetful-functor-preserves-coacyclicity}(d),
all Becker-coacyclic complexes in $\cA\Qcoh^{X\dcot}$ are
Becker-semicoacyclic.
 By Lemma~\ref{becker-coacyclic-complexes-are-acyclic}(d) (for
the quasi-coherent quasi-algebra $\cO_X$ over~$X$), all Becker-coacyclic
complexes in $X\Qcoh^\cot$ are acyclic in $X\Qcoh^\cot$.
 Therefore, it follows from
Lemma~\ref{qcoh-forgetful-functor-reflects-acyclicity}(c) that all
Becker-semicoacyclic complexes in $\cA\Qcoh^{X\dcot}$ are acyclic
in $\cA\Qcoh^{X\dcot}$.
 So we have
$$
 \Ac^\bco(\cA\Qcoh^{X\dcot})\subset\Ac^\bsi(\cA\Qcoh^{X\dcot})
 \subset\Ac(\cA\Qcoh^{X\dcot}).
$$
 The triangulated Verdier quotient category
$$
 \sD^\bsi(\cA\Qcoh^{X\dcot})=
 \Hot(\cA\Qcoh^{X\dcot})/\Ac^\bsi(\cA\Qcoh^{X\dcot})
$$
is called the \emph{Becker semi}(\emph{co})\emph{derived category}
of $X$\+cotorsion quasi-coherent left $\cA$\+modules on~$X$.

 Finally, a complex in the exact category $\cA\Qcoh^{\cA\dcot}$ is said
to be \emph{Becker-semicoacyclic} (relative to $X\Qcoh^\cot$) if,
\emph{viewed as a complex in $X\Qcoh^\cot$}, it is Becker-coacyclic
in the exact category $X\Qcoh^\cot$.
 The full subcategory of Becker-semicoacyclic complexes is denoted by
$\Ac^\bsi(\cA\Qcoh^{\cA\dcot})\subset\Hot(X\Qcoh^{\cA\dcot})$.
 By Lemma~\ref{forgetful-functor-preserves-coacyclicity}(e),
all Becker-coacyclic complexes in $\cA\Qcoh^{\cA\dcot}$ are
Becker-semicoacyclic.
 By Lemma~\ref{becker-coacyclic-complexes-are-acyclic}(d) (for
the quasi-coherent quasi-algebra $\cO_X$ over~$X$), all Becker-coacyclic
complexes in $X\Qcoh^\cot$ are acyclic in $X\Qcoh^\cot$.
 Therefore, it follows from
Lemma~\ref{qcoh-forgetful-functor-reflects-acyclicity}(d) that all
Becker-semicoacyclic complexes in $\cA\Qcoh^{\cA\dcot}$ are acyclic
in $\cA\Qcoh^{\cA\dcot}$.
 So we have
$$
 \Ac^\bco(\cA\Qcoh^{\cA\dcot})\subset\Ac^\bsi(\cA\Qcoh^{\cA\dcot})
 \subset\Ac(\cA\Qcoh^{\cA\dcot}).
$$
 The triangulated Verdier quotient category
$$
 \sD^\bsi(\cA\Qcoh^{\cA\dcot})=
 \Hot(\cA\Qcoh^{\cA\dcot})/\Ac^\bsi(\cA\Qcoh^{\cA\dcot})
$$
is called the \emph{Becker semi}(\emph{co})\emph{derived category}
of $\cA$\+cotorsion quasi-coherent left $\cA$\+modules on~$X$.

 The definitions of and the notation for the semiderived categories
of quasi-coherent right $\cA$\+modules are similar to those for
quasi-coherent left $\cA$\+modules.
 One just writes $\Qcohr\cA$ instead of $\cA\Qcoh$.

\begin{cor}
 Let $X$ be a quasi-compact semi-separated scheme and $\cA$ be
a quasi-coherent quasi-algebra over~$X$.
 Then the inclusions of exact/abelian categories
$\cA\Qcoh^{\cA\dcot}\rarrow\cA\Qcoh^{X\dcot}\rarrow\cA\Qcoh^{X\dcta}
\rarrow\cA\Qcoh$ and $X\Qcoh^\cot\rarrow X\Qcoh^\cta\rarrow X\Qcoh$
induce triangulated equivalences of the Becker semiderived categories
\hbadness=1100
$$
 \sD^\bsi(\cA\Qcoh^{\cA\dcot})\simeq\sD^\bsi(\cA\Qcoh^{X\dcot})
 \simeq\sD^\bsi(\cA\Qcoh^{X\dcta})\simeq\sD^\bsi(\cA\Qcoh).
$$
\end{cor}

\begin{proof}
 Following the proof of
Lemma~\ref{forgetful-functor-preserves-coacyclicity}(c\+-e) (for
the quasi-coherent quasi-algebra $\cO_X$ over~$X$), a complex in
any one of the categories $\cA\Qcoh^{X\dcta}$, \ $\cA\Qcoh^{X\dcot}$,
or $\cA\Qcoh^{\cA\dcot}$ is semicoacyclic as a complex in the respective
exact category if and only if it is semicoacyclic as a complex in
the abelian category $\cA\Qcoh$.
 Therefore, it remains to show that, for any complex $\M^\bu$ in
$\cA\Qcoh$, there exists a complex $\J^\bu$ in $\cA\Qcoh^{\cA\dcot}$
together with a closed morphism of complexes $\M^\bu\rarrow\J^\bu$
with a cone semiacyclic in $\cA\Qcoh$.
 The latter assertion follows from (the proof of)
Theorem~\ref{coderived-of-Grothendieck-DG-category} or
Corollary~\ref{qcoh-cta-cot-cdg-becker-coderived-equivalence} together
with Lemma~\ref{forgetful-functor-preserves-coacyclicity}(b).
\end{proof}

\begin{cor} \label{noetherian-Positselski=Becker-semicoderived}
 Let $X$ be a locally Noetherian scheme and $\cA$ be a quasi-coherent
quasi-algebra over~$X$.
 Then the classes of Positselski-semicoacyclic and Becker-semicoacyclic
complexes of quasi-coherent $\cA$\+modules on $X$ coincide.
 So one has\/ $\Ac^\si(\cA\Qcoh)=\Ac^\bsi(\cA\Qcoh)$ and\/
$\sD^\si(\cA\Qcoh)=\sD^\bsi(\cA\Qcoh)$.
\end{cor}

\begin{proof}
 For a Noetherian scheme $X$, the assertion follows immediately from
Corollary~\ref{noetherian-quasi-algebra-Positselski=Becker} (for
the quasi-coherent CDG\+quasi-algebra $\cB^\cu=\cO_X$ over~$X$).
 In the general case of a locally Noetherian scheme $X$, the point is
that even though the Grothendieck category $X\Qcoh$ need not be locally
Noetherian~\cite[Example in Section~II.7]{HartRD}, infinite direct
sums of injective objects are still injective in $X\Qcoh$
(because injectivity in $X\Qcoh$ is a local
property~\cite[Lemma~II.7.16 and Theorem~II.7.18]{HartRD}),
so Theorem~\ref{under-star-condition-Positselski=Becker} tells us
that the classes of Positselski-coacyclic and Becker-coacyclic
complexes in $X\Qcoh$ coincide.
 (See also~\cite[Theorem~6.4.10(a)]{Pcosh}.)
\end{proof}

\subsection{Semiderived categories of contraherent $\cA$-modules}
\label{semiderived-contraherent-subsecn}
 Let $X$ be a scheme with an open covering $\bW$ and $\cA$ be
a quasi-coherent quasi-algebra over~$X$.
 The semiderived category of (complexes of) $\bW$\+locally contraherent
$\cA$\+modules is an intermediate triangulated Verdier quotient
category between the derived and the contraderived category of
$\bW$\+locally contraherent $\cA$\+modules, constructed using
the notion of a contraacyclic complex of $\bW$\+locally contraherent
cosheaves (of $\cO_X$\+modules) on~$X$.
 Since there are two notions of a contraacyclic complex of
$\bW$\+locally contraherent cosheaves, the Positselski and the Becker
one, the natural expectation is that there should be are also two
notions of a semi(contra)derived category.

 However, we have seen in
Section~\ref{semiderived-quasi-coherent-subsecn} that the assertions
of Lemma~\ref{forgetful-functor-preserves-coacyclicity}(b\+-e) play
an important role in the theory (or indeed in the definitions)
of semiderived categories of quasi-coherent $\cA$\+modules.
 The similar assertions for $\bW$\+locally contraherent $\cA$\+modules
are harder to prove, because the general result
of~\cite[Lemma~B.7.5(b)]{Pcosh} or
Lemma~\ref{DG-functors-preserve-Becker-co-contra-acyclicity}(b) above
is \emph{not} applicable (as the required left adjoint functors to
the forgetful functors do not exist).
 For this reason, we restrict ourselves to
Positselski semi(contra)derived categories in this
Section~\ref{semiderived-contraherent-subsecn}, and postpone
the discussion of the Becker semicontraderived categories to
Section~\ref{becker-semicontraderived-defined-subsecn}.

 As in Section~\ref{semiderived-quasi-coherent-subsecn}, we start with
introducing some notation.
 A quasi-coherent quasi-algebra $\cA$ can be viewed as a special case
of a quasi-coherent CDG\+quasi-algebra $\cA^\cu=\cA$; then
the $\bW$\+locally contraherent CDG\+modules over $\cA^\cu$ are just
complexes of $\bW$\+locally contraherent $\cA$\+modules.
 In this context, we have $\cA^\cu\bLcth_\bW=\bCom(\cA\Lcth_\bW)$,
and similarly $\cA^\cu\bLcth_\bW^{X\dlct}=
\bCom(\cA\Lcth_\bW^{X\dlct})$ and $\cA^\cu\bLcth_\bW^{\cA\dlct}=
\bCom(\cA\Lcth_\bW^{\cA\dlct})$.

 For an exact category $\sE$, as particular cases of the definition
in Section~\ref{derived-second-kind-subsecn}, we obtain the thick
subcategories of absolutely acyclic and contraacyclic complexes
$$
 \Ac^\abs(\sE)\subset\Ac^\ctr(\sE)\subset\Ac^\bctr(\sE)\subset
 \Hot(\sE),
$$
where the full subcategory of Positselski-contraacyclic complexes
$\Ac^\ctr(\sE)$ is only well-defined when infinite products exist
and are exact in~$\sE$.
 The related exotic derived categories are the absolute derived
category $\sD^\abs(\sE)=\Hot(\sE)/\Ac^\abs(\sE)$, the Positselski
contraderived category $\sD^\ctr(\sE)=\Hot(\sE)/\Ac^\ctr(\sE)$, and
the Becker contraderived category $\sD^\bctr(\sE)=
\Hot(\sE)/\Ac^\bctr(\sE)$.

\begin{lem} \label{forg-functor-preserves-positselski-contraacyclicity}
 Let $X$ be a scheme with an open covering\/ $\bW$ and $\cA$ be
a quasi-coherent quasi-algebra over~$X$.
 In this context: \par
\textup{(a)} the forgetful functor $\cA\Lcth_\bW\rarrow X\Lcth_\bW$
takes Positselski-contraacyclic complexes in $\cA\Lcth_\bW$ to
Positselski-contraacyclic complexes in $X\Lcth_\bW$; \par
\textup{(b)} the forgetful functor $\cA\Lcth_\bW^{X\dlct}\rarrow
X\Lcth_\bW^\lct$ takes Positselski-contraacyclic complexes in
$\cA\Lcth_\bW^{X\dlct}$ to Positselski-contraacyclic complexes
in $X\Lcth_\bW^\lct$; \par
\textup{(c)} the forgetful functor $\cA\Lcth_\bW^{\cA\dlct}\rarrow
X\Lcth_\bW^\lct$ takes Positselski-contraacyclic complexes in
$\cA\Lcth_\bW^{\cA\dlct}$ to Positselski-contraacyclic complexes
in $X\Lcth_\bW^\lct$.
\end{lem}

\begin{proof}
 All the assertions hold because the forgetful functors in question
are exact and preserve infinite products.
\end{proof}

\begin{lem} \label{positselski-contraacyclic-complexes-are-acyclic}
 Let $X$ be a scheme with an open covering\/ $\bW$ and $\cA$ be
a quasi-coherent quasi-algebra over~$X$.
 In this context: \par
\textup{(a)} any Positselski-contraacyclic complex in the exact
category $\cA\Lcth_\bW$ is acyclic in $\cA\Lcth_\bW$; \par
\textup{(b)} any Positselski-contraacyclic complex in the exact
category $\cA\Lcth_\bW^{X\dlct}$ is acyclic in $\cA\Lcth_\bW^{X\dlct}$;
\par
\textup{(c)} any Positselski-contraacyclic complex in the exact
category $\cA\Lcth_\bW^{\cA\dlct}$ is acyclic in
$\cA\Lcth_\bW^{\cA\dlct}$.
\end{lem}

\begin{proof}
 All the assertions hold because the infinite product functors are
exact in the respective exact categories.
\end{proof}

 The following lemma should be compared with~\cite[Lemma~3.1.2]{Pcosh}.

\begin{lem} \label{lcth-forgetful-functor-reflects-acyclicity}
 Let $X$ be a scheme with an open covering\/ $\bW$ and $\cA$ be
a quasi-coherent quasi-algebra over~$X$.
 In this context: \par
\textup{(a)} a complex in $\cA\Lcth_\bW$ is acyclic if and only if it
is acyclic in $X\Lcth_\bW$; \par
\textup{(b)} a complex in $\cA\Lcth_\bW^{X\dlct}$ is acyclic if and only
if it is acyclic in $X\Lcth_\bW^\lct$; \par
\textup{(c)} a complex in $\cA\Lcth_\bW^{\cA\dlct}$ is acyclic if and
only if it is acyclic in $X\Lcth_\bW^\lct$.
\end{lem}

\begin{proof}
 Parts~(a\+-b) follow easily from the definitions.
 To prove part~(c), one needs to use
Corollary~\ref{loc-contrah-quasi-algebra-cotorsion-periodicity}(b).
\end{proof}

 Now we can spell out the definitions.
 A complex in the exact category $\cA\Lcth_\bW$ is said to be
\emph{Positselski-semicontraacyclic} (relative to $X\Lcth_\bW$) if,
\emph{viewed as a complex in $X\Lcth_\bW$}, it is
Positselski-contraacyclic in the exact category $X\Lcth_\bW$.
 The full subcategory of Positselski-semicontraacyclic complexes is
denoted by $\Ac^\si(\cA\Lcth_\bW)\subset\Hot(\cA\Lcth_\bW)$.
 By Lemma~\ref{forg-functor-preserves-positselski-contraacyclicity}(a),
all Positselski-contraacyclic complexes in $\cA\Lcth_\bW$ are
Positselski-semicontraacyclic.
 By Lemma~\ref{positselski-contraacyclic-complexes-are-acyclic}(a)
(for the quasi-coherent quasi-algebra $\cO_X$ over~$X$), all
Positselski-contraacyclic complexes in $X\Lcth_\bW$ are acyclic
in $X\Lcth_\bW$.
 Therefore, it follows from
Lemma~\ref{lcth-forgetful-functor-reflects-acyclicity}(a) that all
Positselski-semicontraacyclic complexes in $\cA\Lcth_\bW$ are acyclic
in $\cA\Lcth_\bW$.
 So we have
$$
 \Ac^\ctr(\cA\Lcth_\bW)\subset\Ac^\si(\cA\Lcth_\bW)\subset
 \Ac(\cA\Lcth_\bW).
$$
 The triangulated Verdier quotient category
$$
 \sD^\si(\cA\Lcth_\bW)=\Hot(\cA\Lcth_\bW)/\Ac^\si(\cA\Lcth_\bW)
$$
is called the Positselski \emph{semi}(\emph{contra})\emph{derived
category} of $\bW$\+locally contraherent $\cA$\+modules on~$X$.

 A complex in the exact category $\cA\Lcth_\bW^{X\dlct}$ is said to be
\emph{Positselski-semi\-con\-tra\-acyclic} (relative to
$X\Lcth_\bW^\lct$) if, \emph{viewed as a complex in $X\Lcth_\bW^\lct$},
it is Positselski-contraacyclic in the exact category $X\Lcth_\bW^\lct$.
 The full subcategory of Positselski-semicontraacyclic complexes is
denoted by $\Ac^\si(\cA\Lcth_\bW^{X\dlct})\subset
\Hot(\cA\Lcth_\bW^{X\dlct})$.
 By Lemma~\ref{forg-functor-preserves-positselski-contraacyclicity}(b),
all Positselski-contraacyclic complexes in $\cA\Lcth_\bW^{X\dlct}$ are
Positselski-semicontraacyclic.
 By Lemma~\ref{positselski-contraacyclic-complexes-are-acyclic}(b)
(for the quasi-coherent quasi-algebra $\cO_X$ over~$X$), all
Positselski-contraacyclic complexes in $X\Lcth_\bW^\lct$ are acyclic
in $X\Lcth_\bW^\lct$.
 Therefore, it follows from
Lemma~\ref{lcth-forgetful-functor-reflects-acyclicity}(b) that all
Positselski-semicontraacyclic complexes in $\cA\Lcth_\bW^{X\dlct}$ are
acyclic in $\cA\Lcth_\bW^{X\dlct}$.
 So we have
$$
 \Ac^\ctr(\cA\Lcth_\bW^{X\dlct})\subset\Ac^\si(\cA\Lcth_\bW^{X\dlct})
 \subset\Ac(\cA\Lcth_\bW^{X\dlct}).
$$
 The triangulated Verdier quotient category
$$
 \sD^\si(\cA\Lcth_\bW^{X\dlct})=
 \Hot(\cA\Lcth_\bW^{X\dlct})/\Ac^\si(\cA\Lcth_\bW^{X\dlct})
$$
is called the Positselski \emph{semi}(\emph{contra})\emph{derived
category} of $X$\+locally cotorsion $\bW$\+locally contraherent
$\cA$\+modules on~$X$.

 A complex in the exact category $\cA\Lcth_\bW^{\cA\dlct}$ is said to be
\emph{Positselski-semi\-con\-tra\-acyclic} (relative to
$X\Lcth_\bW^\lct$) if, \emph{viewed as a complex in $X\Lcth_\bW^\lct$},
it is Positselski-contraacyclic in the exact category $X\Lcth_\bW^\lct$.
 The full subcategory of Positselski-semicontraacyclic complexes is
denoted by $\Ac^\si(\cA\Lcth_\bW^{\cA\dlct})\subset
\Hot(\cA\Lcth_\bW^{\cA\dlct})$.
 By Lemma~\ref{forg-functor-preserves-positselski-contraacyclicity}(c),
all Positselski-contraacyclic complexes in $\cA\Lcth_\bW^{\cA\dlct}$ are
Positselski-semicontraacyclic.
 By Lemma~\ref{positselski-contraacyclic-complexes-are-acyclic}(b)
(for the quasi-coherent quasi-algebra $\cO_X$ over~$X$), all
Positselski-contraacyclic complexes in $X\Lcth_\bW^\lct$ are acyclic
in $X\Lcth_\bW^\lct$.
 Therefore, it follows from
Lemma~\ref{lcth-forgetful-functor-reflects-acyclicity}(c) that all
Positselski-semicontraacyclic complexes in $\cA\Lcth_\bW^{\cA\dlct}$
are acyclic in $\cA\Lcth_\bW^{\cA\dlct}$.
 So we have
$$
 \Ac^\ctr(\cA\Lcth_\bW^{\cA\dlct})\subset
 \Ac^\si(\cA\Lcth_\bW^{\cA\dlct})\subset
 \Ac(\cA\Lcth_\bW^{\cA\dlct}).
$$
 The triangulated Verdier quotient category
$$
 \sD^\si(\cA\Lcth_\bW^{\cA\dlct})=
 \Hot(\cA\Lcth_\bW^{\cA\dlct})/\Ac^\si(\cA\Lcth_\bW^{\cA\dlct})
$$
is called the Positselski \emph{semi}(\emph{contra})\emph{derived
category} of $\cA$\+locally cotorsion $\bW$\+locally contraherent
$\cA$\+modules on~$X$.

 In the case of the open covering $\bW=\{X\}$, we will write
\begin{align*}
 \sD^\si(\cA\Ctrh)&=\sD^\si(\cA\Lcth_{\{X\}}), \\
 \sD^\si(\cA\Ctrh^{X\dlct})&=\sD^\si(\cA\Lcth_{\{X\}}^{X\dlct}), \\
 \sD^\si(\cA\Ctrh^{\cA\dlct})&=\sD^\si(\cA\Lcth_{\{X\}}^{\cA\dlct}),
\end{align*}
etc.

\begin{cor} \label{positselski-semicontrader-independ-of-covering}
 Let $X$ be a quasi-compact semi-separated scheme with an open
covering\/ $\bW$ and $\cA$ be a quasi-coherent quasi-algebra over~$X$.
 In this context: \par
\textup{(a)} The inclusions of exact categories $\cA\Ctrh_\bW
\rarrow\cA\Lcth_\bW$ and $X\Ctrh_\bW\rarrow X\Lcth_\bW$ induce
a triangulated equivalence of the semiderived categories
$$
 \sD^\si(\cA\Ctrh)\simeq\sD^\si(\cA\Lcth_\bW).
$$ \par
\textup{(b)} The inclusions of exact categories $\cA\Ctrh_\bW^{X\dlct}
\rarrow\cA\Lcth_\bW^{X\dlct}$ and $X\Ctrh_\bW^\lct\allowbreak\rarrow
X\Lcth_\bW^\lct$ induce a triangulated equivalence of the semiderived
categories
$$
 \sD^\si(\cA\Ctrh^{X\dlct})\simeq\sD^\si(\cA\Lcth_\bW^{X\dlct}).
$$ \par
\textup{(c)} The inclusions of exact categories $\cA\Ctrh_\bW^{\cA\dlct}
\rarrow\cA\Lcth_\bW^{\cA\dlct}$ and $X\Ctrh_\bW^\lct\allowbreak\rarrow
X\Lcth_\bW^\lct$ induce a triangulated equivalence of the semiderived
categories
$$
 \sD^\si(\cA\Ctrh^{\cA\dlct})\simeq\sD^\si(\cA\Lcth_\bW^{\cA\dlct}).
$$
\end{cor}

\begin{proof}
 Part~(a): firstly, by
Corollary~\ref{contraderived-indep-of-covering-or-antilocal}(a)
(for the quasi-coherent CDG\+quasi-algebra $\cB^\cu=\cO_X$) or
by~\cite[Corollary~4.7.4(a)]{Pcosh}, a complex in $\cA\Ctrh$
is Positselski-semicontraacyclic in $\cA\Ctrh$ if and only if it is
Positselski-semicontraacyclic in $\cA\Lcth_\bW$.
 Secondly, by
Corollary~\ref{contraderived-indep-of-covering-or-antilocal}(a)
and the proof of Proposition~\ref{second-kind-finite-resolutions}(d)
(for the quasi-coherent CDG\+quasi-algebra $\cB^\cu=\cA$), for
every complex $\Q^\bu$ in $\cA\Lcth_\bW$ there exists a complex
$\P^\bu$ in $\cA\Ctrh$ together with a closed morphism of complexes
$\P^\bu\rarrow\Q^\bu$ whose cone is absolutely acyclic in
$\cA\Lcth_\bW$.
 Finally, any complex that is absolutely acyclic in $\cA\Lcth_\bW$
is also Positselski-contraacyclic (in fact, absolutely acyclic) in
$X\Lcth_\bW$ by
Lemma~\ref{forg-functor-preserves-positselski-contraacyclicity}(a).
 The proofs of parts~(b\+-c) are similar and based on
Corollary~\ref{contraderived-indep-of-covering-or-antilocal}(b\+-c),
Lemma~\ref{forg-functor-preserves-positselski-contraacyclicity}(b\+-c),
and~\cite[Corollary~4.7.5(a)]{Pcosh}.
\end{proof}

 For a triangulated equivalence $\sD^\si(\cA\Lcth_\bW^{X\dlct})
\simeq\sD^\si(\cA\Lcth_\bW)$ proved under more restrictive assumptions
on the scheme $X$, see
Corollary~\ref{semicontraderived-X-lcta-X-lct-A-mod-equivalence} below.

\Section{Becker Semicontraderived Category}
\label{becker-semicontraderived-secn}

 In this section we prove the basic properties of Becker-contraacyclic
complexes of $\bW$\+locally contraherent $\cA$\+modules.
 The aim is to make sense of the definition of the Becker
semicontraderived category $\sD^\bsi(\cA\Lcth_\bW)$ for a quasi-coherent
quasi-algebra $\cA$ over a quasi-compact semi-separated scheme $X$ with
an open covering~$\bW$.
 Specifically, we need to show that the forgetful functor
$\cA\Lcth_\bW\rarrow X\Lcth_\bW$ preserves the Becker contraacyclicity
of complexes.

 This section uses more knowledge of the machinery of \emph{cotorsion
pairs} then the rest of this paper.
 We refer to the introduction to the paper~\cite{Pctrl} for
an introductory discussion and to~\cite[Appendix~B]{Pcosh} for
the technical background.
 The argument in this section follows the ``alternative approach''
outlined in~\cite[Corollary~5.1.10, Remark~5.4.7, and proof of
Corollary~5.4.10]{Pcosh}.

\subsection{Locality of Becker coacyclicity}
 For the sake of context and completeness of the exposition, we start
with proving locality of the Becker coacyclicity property of complexes
of quasi-coherent $\cA$\+modules on a quasi-compact semi-separated
scheme before proceeding to present the more complicated proof of
the locality of the Becker contraacyclicity of complexes of
$\bW$\+locally contraherent $\cA$\+modules.
 The proof of the locality of the Positselski coacyclicity and
contraacyclicity is easy (as the preservation of the Positselski
co/contraacyclicity by the restrictions to open subschemes and
the direct images from affine open immersions is obvious),
so we skip it.
 The basic argument goes back to~\cite[Remark~1.3]{EP}.

\begin{lem} \label{qcoh-becker-coacyclic-inverse-image}
 Let $X$ be a scheme, $Y\subset X$ be an open subscheme with an open
immersion morphism $f\:Y\rarrow X$, and $\cA$ be a quasi-coherent
quasi-algebra over~$X$.
 Then the inverse image functor $f^*\:\cA\Qcoh\rarrow\cA|_Y\Qcoh$
defined in Sections~\ref{direct-images-of-A-co-sheaves-subsecn}--%
\ref{inverse-images-of-A-co-sheaves-subsecn} takes Becker-coacyclic
complexes in $\cA\Qcoh$ to Becker-coacyclic complexes in $\cA|_Y\Qcoh$.
\end{lem}

\begin{proof}
 Notice that the inverse image functor $f^*\:\cA\Qcoh\rarrow\cA|_Y\Qcoh$
is exact (as a functor between abelian categories) and preserves
infinite direct sums.
 Under the simplifying assumption that the open immersion morphism~$f$
is quasi-compact, the functor~$f^*$ is left adjoint to the direct image
functor $f_*\:\cA|_Y\Qcoh\rarrow\cA\Qcoh$;
so~\cite[Lemma~B.7.5(a)]{Pcosh} or
Lemma~\ref{DG-functors-preserve-Becker-co-contra-acyclicity}(a) above
is applicable.
 In the general case, one can refer to~\cite[Lemma~A.5]{Psemten}
or Lemma~\ref{Grothendick-abelian-DG-functor-preserves-Becker-co}
above; cf.~\cite[proof of Lemma~A.2.3]{Pcosh}.
\end{proof}

\begin{lem} \label{qcoh-becker-coacyclic-affine-direct-image}
 Let $X$ be a scheme, $Y\subset X$ be an open subscheme such tha
the open immersion morphism $f\:Y\rarrow X$ is affine, and $\cA$ be
a quasi-coherent quasi-algebra over~$X$.
 Then the direct image functor $f_*\:\cA|_Y\Qcoh\rarrow\cA\Qcoh$
defined in Section~\ref{direct-images-of-A-co-sheaves-subsecn}
takes Becker-coacyclic complexes in $\cA|_Y\Qcoh$ to Becker-coacyclic
complexes in $\cA\Qcoh$.
\end{lem}

\begin{proof}
 This is also a particular case of~\cite[Lemma~A.5]{Psemten} or
Lemma~\ref{Grothendick-abelian-DG-functor-preserves-Becker-co} above.
 The point is that the direct image functor~$f_*$ is exact and
preserves infinite direct sums (cf.~\cite[proof of Lemma~5.2.4]{Pcosh}).
\end{proof}

 The following corollary is a generalization
of~\cite[Corollary~5.2.5]{Pcosh} (see
also~\cite[Corollary~6.4.12(a)]{Pcosh}).

\begin{cor} \label{locality-of-becker-coacyclicity}
 Let $X$ be a quasi-compact semi-separated scheme and
$X=\bigcup_{\alpha=1}^N U_\alpha$ be its finite affine open covering.
 Let $\cA$ be a quasi-coherent quasi-algebra over~$X$.
 Then a complex $\cB^\bu$ in $\cA\Qcoh$ is Becker-coacyclic if and
only if, for every index~$\alpha$, the complex of
$\cA(U_\alpha)$\+modules $\cB^\bu(U_\alpha)$ is Becker-coacyclic
in $\cA(U_\alpha)\Modl$.
\end{cor}

\begin{proof}
 The ``only if'' implications holds by
Lemma~\ref{qcoh-becker-coacyclic-inverse-image}.
 The ``if'' implication is nontrivial.
 The following argument goes back to~\cite[Remark~1.3]{EP}.
 Consider the \v Cech complex of complexes of quasi-coherent
$\cA$\+modules (cf.~\eqref{qcoh-sheaf-of-rings-cech-coresolution})
\begin{multline} \label{qcoh-sheaf-of-rings-cech-coresol-of-complex}
 0\lrarrow\cB^\bu\lrarrow\bigoplus\nolimits_{\alpha=1}^N
 j_\alpha{}_*j_\alpha^*\cB^\bu\lrarrow
 \bigoplus\nolimits_{1\le\alpha<\beta\le N}
 j_{\alpha,\beta}{}_*j_{\alpha,\beta}^*\cB^\bu \\
 \lrarrow\dotsb\lrarrow j_{1,\dotsc,N}{}_*j_{1,\dotsc,N}^*\cB^\bu
 \lrarrow0,
\end{multline}
where, for every subsequence of indices $1\le\alpha_1<\dotsb<\alpha_k
\le N$, we denote by $j_{\alpha_1,\dotsc,\alpha_k}$ the open
immersion morphism $\bigcap_{s=1}^k U_{\alpha_s}\rarrow X$.
 Then, by Lemmas~\ref{qcoh-becker-coacyclic-inverse-image} and
\ref{qcoh-becker-coacyclic-affine-direct-image}, all
the terms of~\eqref{qcoh-sheaf-of-rings-cech-coresol-of-complex}
except perhaps the leftmost one are Becker-coacyclic complexes
in $\cA\Qcoh$.
 Since the finite complex of
complexes~\eqref{qcoh-sheaf-of-rings-cech-coresol-of-complex} is
acyclic in $\cA\Qcoh$ and all absolutely acyclic complexes are
Becker-coacyclic, it follows that the complex $\cB^\bu$ is
Becker-coacyclic in $\cA\Qcoh$ as well.
\end{proof}

\subsection{Direct images of Becker-contraacyclic complexes of
antilocal contraherent modules}
 Given a quasi-coherent quasi-algebra $\cA$ over a scheme $X$, we
denote by $\cA\Qcoh_\vflp\subset\cA\Qcoh_\flp\subset\cA\Qcoh_\fl
\subset\cA\Qcoh$ the full subcategories of very flaprojective,
flaprojective, and flat quasi-coherent $\cA$\+modules in the abelian
category $\cA\Qcoh$ of quasi-coherent left $\cA$\+modules on~$X$
(see Sections~\ref{A-loc-cotors-loc-inj-cosheaves-subsecn}
and~\ref{antilocality-of-X-contraadjusted-subsecn}\+-%
\ref{antilocality-of-A-cotorsion-subsecn} for the definitions).
 The notation $\cA\Qcoh^{\cA\dcot}\subset\cA\Qcoh^{X\dcot}\subset
\cA\Qcoh^{X\dcta}\subset\cA\Qcoh$ for the full exact subcategories
of $\cA$\+cotorsion, $X$\+cotorsion, and $X$\+contraadjusted
quasi-coherent $\cA$\+modules in the abelian category $\cA\Qcoh$
was introduced in Section~\ref{background-derived-quasi-coherent}.

 The following lemma is our version of~\cite[Lemma~5.1.9]{Pcosh}.
 The formulation of the lemma involves Becker-contraacyclic complexes
in certain exact categories of quasi-coherent $\cA$\+modules, which
may appear to be somewhat unusual.
 Notice, however, that (\emph{unlike} the abelian category $\cA\Qcoh$)
the exact categories $\cA\Qcoh^{X\dcta}$, $\cA\Qcoh^{X\dcot}$, and
$\cA\Qcoh^{\cA\dcot}$ have enough projective objects.
 In fact, the intersection $\cA\Qcoh_\vflp^{X\dcta}=\cA\Qcoh_\vflp
\cap\cA\Qcoh^{X\dcta}$ is the full subcategory of projective objects
in $\cA\Qcoh^{X\dcta}$, the intersection $\cA\Qcoh_\flp^{X\dcot}=
\cA\Qcoh_\flp\cap\cA\Qcoh^{X\dcot}$ is the full subcategory of
projective objects in $\cA\Qcoh^{X\dcot}$, and the intersection
$\cA\Qcoh_\fl^{\cA\dcot}=\cA\Qcoh_\fl\cap\cA\Qcoh^{\cA\dcot}$ is
the full subcategory of projective objects in $\cA\Qcoh^{\cA\dcot}$;
and in each case there are enough of such objects, as one can see from
Theorems~\ref{qcomp-qsep-very-flaproj-complete-cotorsion-pair-thm}(a),
\ref{qcomp-qsep-flaproj-complete-cotorsion-pair-thm}(a),
and~\ref{qcomp-qsep-A-flat-complete-cotorsion-pair-thm}(a)
(cf.~\cite[Lemma~1.9(b)]{Pfltp}).

 Recall the notation $\Com(\sA)=\sZ^0(\bCom(\sA))$ for the additive
category of complexes (and closed morphisms of degree~$0$ between
them) in an additive category~$\sA$, as in
Section~\ref{semiderived-quasi-coherent-subsecn}.
 Given an exact category $\sE$, we will use the notation
$\Ac^\bctr(\sE)\subset\Com(\sE)$ for the full subcategory of
Becker-contraacyclic complexes in $\Com(\sE)$ (cf.\
Section~\ref{semiderived-contraherent-subsecn}).

\begin{lem} \label{all-flapr-becker-contraacyclic-cot-cotorsion-pair}
 Let $X$ be a quasi-compact semi-separated scheme and $\cA$ be 
a quasi-coherent quasi-algebra over~$X$.
 Then \par
\textup{(a)} The pair of classes of all complexes of very flaprojective
quasi-coherent $\cA$\+modules\/ $\sF=\Com(\cA\Qcoh_\vflp)$ and
Becker-contraacyclic complexes of $X$\+con\-tra\-ad\-justed
quasi-coherent $\cA$\+modules\/ $\sC=\Ac^\bctr(\cA\Qcoh^{X\dcta})$ is
a hereditary complete cotorsion pair $(\sF,\sC)$ in the abelian
category\/ $\Com(\cA\Qcoh)$. \par
\textup{(b)}  The pair of classes of all complexes of flaprojective
quasi-coherent $\cA$\+modules\/ $\sF=\Com(\cA\Qcoh_\flp)$ and
Becker-contraacyclic complexes of $X$\+cotorsion quasi-coherent
$\cA$\+modules\/ $\sC=\Ac^\bctr(\cA\Qcoh^{X\dcot})$ is a hereditary
complete cotorsion pair $(\sF,\sC)$ in the abelian category\/
$\Com(\cA\Qcoh)$. \par
\textup{(c)} The pair of classes of all complexes of flat
quasi-coherent $\cA$\+modules\/ $\sF=\Com(\cA\Qcoh_\fl)$ and
Becker-contraacyclic complexes of $\cA$\+cotorsion quasi-coherent
$\cA$\+modules\/ $\sC=\Ac^\bctr(\cA\Qcoh^{\cA\dcot})$ is a hereditary
complete cotorsion pair $(\sF,\sC)$ in the abelian category\/
$\Com(\cA\Qcoh)$.
\end{lem}

\begin{proof}
 All the three assertions are special cases
of~\cite[Proposition~2.2]{Pfltp}, which is a general result about
cotorsion pairs in Grothendieck categories (see
also~\cite[Proposition~3.2, Theorem~5.5, and Section~5.3]{Gil2}
and~\cite[Lemma~4.9]{Gil3}).
 For applicability of~\cite[Proposition~2.2]{Pfltp}, one only needs
to show that each of the original pairs of classes ($\cA\Qcoh_\vflp$,
$\cA\Qcoh^{X\dcta}$), \ ($\cA\Qcoh_\flp$, $\cA\Qcoh^{X\dcot}$),
and ($\cA\Qcoh_\fl$, $\cA\Qcoh^{\cA\dcot}$) is a hereditary complete
cotorsion pair generated by a set of objects in $\cA\Qcoh$.
 In the context of part~(a), this is explained
in~\cite[Example~6.12(2)]{Pfltp}, and in the context of part~(c),
in~\cite[Example~8.10]{Pfltp}.
 Part~(b) is similar; let us spell out some details.

 The pair of classes ($\cA\Qcoh_\flp$, $\cA\Qcoh^{X\dcot}$) is
a hereditary complete cotorsion pair in $\cA\Qcoh$ by
Theorem~\ref{qcomp-qsep-flaproj-complete-cotorsion-pair-thm}
and Remark~\ref{X-cotorsion-A-modules-concluding-remark}.
 Using the Eklof lemma~\cite[Lemma~7.5]{PS6}, in order to prove that
this cotorsion pair is generated by a set of objects in $\cA\Qcoh$,
it suffices to check that the class of flaprojective quasi-coherent
left $\cA$\+modules $\cA\Qcoh_\flp$ is deconstructible in $\cA\Qcoh$.

 Let $X=\bigcup_\alpha U_\alpha$ be a finite affine open
covering of~$X$.
 Then a quasi-coherent $\cA$\+module $\F$ is flaprojective if and only
if, for every index~$\alpha$, the $\cA(U_\alpha)$\+module $\F(U_\alpha)$
is $\cA(U_\alpha)/\cO_X(U_\alpha)$\+flaprojective.
 As explained in
Remark~\ref{very-flaprojective-cotorsion-pair-generated-by} (using
also the result of~\cite[Theorem~7.13]{GT}), the class of all
$\cA(U_\alpha)/\cO_X(U_\alpha)$\+flaprojective left
$\cA(U_\alpha)$\+modules is deconstructible in $\cA(U_\alpha)\Modl$.
 Now an argument based on the Hill lemma for filtered
modules~\cite[Theorem~7.10]{GT} (see~\cite[Theorem~2.1]{Sto-hill}
for a Grothendieck category version) shows that the class of
flaprojective quasi-coherent left $\cA$\+modules $\cA\Qcoh_\flp$
is deconstructible in $\cA\Qcoh$, as desired (this is similar
to~\cite[Lemma~B.3.5]{Pcosh}).
\end{proof}

\begin{cor} \label{qcoh-cta-cot-becker-contraacycl-direct-image}
 Let $X$ be a quasi-compact semi-separated scheme, and let $Y$ be
a quasi-compact open subscheme in $X$ with the open immersion morphism
$f\:Y\rarrow X$.
 Let $\cA$ be a quasi-coherent quasi-algebra over~$X$.
 Consider the direct image functor $f_*\:\cA|_Y\Qcoh\rarrow\cA\Qcoh$,
as defined in Section~\ref{direct-images-of-A-co-sheaves-subsecn}.
 Then \par
\textup{(a)} the functor~$f_*$ takes Becker-contraacyclic complexes in
$\cA|_Y\Qcoh^{Y\dcta}$ to Becker-contraacyclic complexes in
$\cA\Qcoh^{X\dcta}$; \par
\textup{(b)} the functor~$f_*$ takes Becker-contraacyclic complexes in
$\cA|_Y\Qcoh^{Y\dcot}$ to Becker-contraacyclic complexes in
$\cA\Qcoh^{X\dcot}$; \par
\textup{(c)} the functor~$f_*$ takes Becker-contraacyclic complexes in
$\cA|_Y\Qcoh^{\cA|_Y\dcot}$ to Becker-contraacyclic complexes in
$\cA\Qcoh^{\cA\dcot}$.
\end{cor}

\begin{proof}
 This is our version of~\cite[Corollary~5.1.10]{Pcosh}.
 In all the three contexts of parts~(a\+-c), the assertion follows from
the description of the class of Becker-contraacyclic complexes as
the right $\Ext^1$\+orthogonal class $\sC=\sF^{\perp_1}$ provided by
Lemma~\ref{all-flapr-becker-contraacyclic-cot-cotorsion-pair}.
 The argument is based on the $\Ext^1$\+adjunction result
of~\cite[Lemma~1.7(b)]{Pal} applied to the pair of adjoint functors
$f^*\:\cA\Qcoh\rarrow\cA|_Y\Qcoh$ and $f_*\:\cA|_Y\Qcoh\rarrow\cA\Qcoh$.
 The point is that the exact functor~$f^*$ takes the class of complexes
$\sF$ on $X$ into the respective class of complexes $\sF$ on~$Y$.
\end{proof}

\begin{cor} \label{ctrh-antiloc-becker-contraacycl-direct-image}
 Let $X$ be a quasi-compact semi-separated scheme, and let $Y$ be
a quasi-compact open subscheme in $X$ with the open immersion morphism
$f\:Y\rarrow X$.
 Let $\cA$ be a quasi-coherent quasi-algebra over~$X$.
 Then \par
\textup{(a)} the direct image functor
$f_!\:\cA|_Y\Ctrh_\al\rarrow\cA\Ctrh_\al$ from
Lemma~\textup{\ref{underived-naive-co-contra-direct-image-lemma}(a)}
takes Becker-contraacyclic complexes in $\cA|_Y\Ctrh_\al$ to
Becker-contraacyclic complexes in $\cA\Ctrh_\al$; \par
\textup{(b)} the direct image functor $f_!\:\cA|_Y\Ctrh^{Y\dlct}_\al
\rarrow\cA\Ctrh^{X\dlct}_\al$ from
Lemma~\textup{\ref{underived-naive-co-contra-direct-image-lemma}(b)}
takes Becker-contraacyclic complexes in $\cA|_Y\Ctrh^{Y\dlct}_\al$ to
Becker-contraacyclic complexes in $\cA\Ctrh^{X\dlct}_\al$; \par
\textup{(c)} the direct image functor
$f_!\:\cA|_Y\Ctrh^{\cA|_Y\dlct}_\al\rarrow\cA\Ctrh^{\cA\dlct}_\al$ from
Lemma~\textup{\ref{underived-naive-co-contra-direct-image-lemma}(c)}
takes Becker-contraacyclic complexes in $\cA|_Y\Ctrh^{\cA|_Y\dlct}_\al$
to Becker-contraacyclic complexes in $\cA\Ctrh^{\cA\dlct}_\al$.
\hfuzz=1.7pt
\end{cor}

\begin{proof}
 Compare Lemmas~\ref{quasi-algebra-underived-naive-co-contra}
and~\ref{underived-naive-co-contra-direct-image-lemma} with
Corollary~\ref{qcoh-cta-cot-becker-contraacycl-direct-image}.
\end{proof}

\subsection{Locality of Becker contraacyclicity}
 The following version of
Corollary~\ref{ctrh-antiloc-becker-contraacycl-direct-image} drops
the antilocality assumption at the expense of imposing the assumption
of affineness of the open immersion morphism.

\begin{cor} \label{lcth-becker-contraacycl-affine-direct-image}
 Let $X$ be a quasi-compact semi-separated scheme, and let $Y$ be
an open subscheme in $X$ such that the open immersion morphism
$f\:Y\rarrow X$ is affine.
 Let\/ $\bW$ be an open covering of $X$, and let $\cA$ be
a quasi-coherent quasi-algebra over~$X$.
 Then \par
\textup{(a)} the direct image functor
$f_!\:\cA|_Y\Lcth_{\bW|_Y}\rarrow\cA\Lcth_\bW$ defined in
Section~\ref{direct-images-of-A-co-sheaves-subsecn} takes
Becker-contraacyclic complexes in $\cA|_Y\Lcth_{\bW|_Y}$ to
Becker-contraacyclic complexes in $\cA\Lcth_\bW$; \par
\textup{(b)} the direct image functor
$f_!\:\cA|_Y\Lcth_{\bW|_Y}^{Y\dlct}\rarrow\cA\Lcth_{\bW|_Y}^{X\dlct}$
defined in Section~\ref{direct-images-of-A-co-sheaves-subsecn}
takes Becker-contraacyclic complexes in $\cA|_Y\Lcth_{\bW|_Y}^{Y\dlct}$
to Becker-contraacyclic complexes in $\cA\Lcth_\bW^{X\dlct}$; \par
\textup{(c)} the direct image functor
$f_!\:\cA|_Y\Lcth_{\bW|_Y}^{\cA|_Y\dlct}\rarrow\cA\Lcth_\bW^{\cA\dlct}$
defined in Section~\ref{direct-images-of-A-co-sheaves-subsecn} takes
Becker-contraacyclic complexes in $\cA|_Y\Lcth_{\bW|_Y}^{\cA|_Y\dlct}$
to Becker-contraacyclic complexes in $\cA\Lcth_\bW^{\cA\dlct}$.
\hbadness=1200
\end{cor}

\begin{proof}
 Let us prove part~(a).
 Let $\gB^\bu$ be a Becker-contraacyclic complex in the exact category
$\cA|_Y\Lcth_{\bW|_Y}$.
 Using Lemma~\ref{antilocal-resolving-resolution-dimension}(a),
the construction of~\cite[Lemma~A.3.4]{Pcosh}, and the fact that
the resolution dimension does not depend on the choice of
a resolution, one can construct a finite acyclic complex of complexes
$0\rarrow\gA_{N-1}^\bu\rarrow\dotsb\rarrow\gA_0^\bu\rarrow\gB^\bu
\rarrow0$ in $\cA|_Y\Lcth_{\bW|_Y}$ such that all the terms of
the complexes $\gA_i^\bu$, \,$0\le i\le N-1$, are antilocal contraherent
$\cA|_Y$\+modules on~$Y$.
 Here $N$ is the number of affine open subschemes in a finite affine
open covering of $Y$ subordinate to~$\bW|_Y$.
 Moreover, following the construction of~\cite[Lemma~A.3.4]{Pcosh}
more closely, one can make the complexes $\gA_i^\bu$ contractible for
all $0\le i<N-1$.
 Then, since all absolutely acyclic complexes are Becker-contraacyclic,
it follows that the complex $\gA_{N-1}^\bu$ is Becker-contraacyclic
in $\cA|_Y\Lcth_{\bW|_Y}$, too.

 The full subcategory $\cA|_Y\Ctrh_\al$ is resolving and closed under
direct summands in $\cA|_Y\Lcth_{\bW|_Y}$, so the classes of projective
objects in the exact categories $\cA|_Y\Ctrh_\al$ and
$\cA|_Y\Lcth_{\bW|_Y}$ coincide.
 Hence a complex in $\cA|_Y\Ctrh_\al$ is Becker-contraacyclic in
$\cA|_Y\Ctrh_\al$ if and only if it is Becker-contraacyclic in
$\cA|_Y\Lcth_{\bW|_Y}$ (cf.\
Corollary~\ref{contraderived-indep-of-covering-or-antilocal}(a)).
 So we have proved that the complex $\gA_{N-1}^\bu$ is
Becker-contraacyclic in $\cA|_Y\Ctrh_\al$.
 By Corollary~\ref{ctrh-antiloc-becker-contraacycl-direct-image}(a),
it follows that the complex $f_!\gA_{N-1}^\bu$ is
Becker-contraacyclic in $\cA\Ctrh_\al$, hence also in $\cA\Lcth_\bW$.

 Finally, the functor $f_!\:\cA|_Y\Lcth_{\bW|_Y}\rarrow\cA\Lcth_\bW$
is exact for an affine open immersion morphism $f\:Y\rarrow X$.
 So the finite complex of complexes $0\rarrow f_!\gA_{N-1}^\bu\rarrow
\dotsb\rarrow f_!\gA_0^\bu\rarrow f_!\gB^\bu\rarrow0$ is acyclic in
$\cA\Lcth_\bW$.
 Furthermore, the complexes $f_!\gA_i^\bu$ are contractible for
all $0\le i<N-1$.
 Once again, since all absolutely acyclic complexes are
Becker-contraacyclic, it follows that the complex $f_!\gB^\bu$ is
Becker-contraacyclic in $\cA\Lcth_\bW$, too.

 The proofs of parts~(b\+-c) are similar and based on
Lemma~\ref{antilocal-resolving-resolution-dimension}(b\+-c),
Corollary~\ref{contraderived-indep-of-covering-or-antilocal}(b\+-c),
and Corollary~\ref{ctrh-antiloc-becker-contraacycl-direct-image}(b\+-c).
\end{proof}

\begin{lem} \label{lcth-becker-contraacycl-inverse-image}
 Let $X$ be a quasi-compact semi-separated scheme, and let $Y$ be
a quasi-compact open subscheme in $X$ with the open immersion morphism
$f\:Y\rarrow X$.
 Let\/ $\bW$ be an open covering of $X$, and let $\cA$ be
a quasi-coherent quasi-algebra over~$X$.
 Then \par
\textup{(a)} the inverse image functor $f^!\:\cA\Lcth_\bW\rarrow
\cA|_Y\Lcth_{\bW|_Y}$ defined in
Section~\ref{direct-images-of-A-co-sheaves-subsecn} and at the end of
Section~\ref{inverse-images-of-A-co-sheaves-subsecn} takes
Becker-contraacyclic complexes in $\cA\Lcth_\bW$ to
Becker-contraacyclic complexes in $\cA|_Y\Lcth_{\bW|_Y}$; \par
\textup{(b)} the inverse image functor $f^!\:\cA\Lcth_\bW^{X\dlct}
\rarrow\cA|_Y\Lcth_{\bW|_Y}^{Y\dlct}$ defined in
Section~\ref{direct-images-of-A-co-sheaves-subsecn} and at the end of
Section~\ref{inverse-images-of-A-co-sheaves-subsecn} takes
Becker-contraacyclic complexes in $\cA\Lcth_\bW^{X\dlct}$ to
Becker-contraacyclic complexes in $\cA|_Y\Lcth_{\bW|_Y}^{Y\dlct}$; \par
\textup{(c)} the inverse image functor $f^!\:\cA\Lcth_\bW^{\cA\dlct}
\rarrow\cA|_Y\Lcth_{\bW|_Y}^{\cA|_Y\dlct}$ defined in
Section~\ref{direct-images-of-A-co-sheaves-subsecn} and at the end of
Section~\ref{inverse-images-of-A-co-sheaves-subsecn} takes
Becker-contraacyclic complexes in $\cA\Lcth_\bW^{\cA\dlct}$ to
Becker-contraacyclic complexes in $\cA|_Y\Lcth_{\bW|_Y}^{\cA|_Y\dlct}$.
\end{lem}

\begin{proof}
 Let us spell out part~(a).
 Notice, first of all, that the inverse image functor
$f^!\:\cA\Lcth_\bW\rarrow\cA|_Y\Lcth_{\bW|_Y}$ is exact (as a functor
between exact categories).
 Under the simplifying assumption that the open immersion morphism~$f$
is affine, the functor~$f^!$ has a left adjoint functor
$f_!\:\cA|_Y\Lcth_{\bW|_Y}\rarrow\cA\Lcth_\bW$;
so~\cite[Lemma~B.7.5(b)]{Pcosh} or
Lemma~\ref{DG-functors-preserve-Becker-co-contra-acyclicity}(b) above
is applicable.
 In the general case, the argument is essentially the same.
 One needs to observe that all the projective $\bW|_Y$\+locally
contraherent $\cA|_Y$\+modules are antilocal (as mentioned in the proof
of Corollary~\ref{lcth-becker-contraacycl-affine-direct-image}),
the functor of direct image of antilocal contraherent modules
$f_!\:\cA|_Y\Ctrh_\al\rarrow\cA\Ctrh_\al$ is well-defined as per
Lemma~\ref{underived-naive-co-contra-direct-image-lemma}(a),
and the functor~$f_!$ is a partial left adjoint to~$f^!$ according
to formula~\eqref{open-immers-direct-inverse-cosheaf-adjunction}
for the ringed spaces $(X,\cA)$ and $(Y,\cA|_Y)$.
 The proofs of parts~(b\+-c) are similar and based on
Lemma~\ref{underived-naive-co-contra-direct-image-lemma}(b)
and Corollary~\ref{contrah-al-A-lct-direct-image}.
\end{proof}

 The following proposition, establishing the locality of
Becker-contraacyclicity, is our generalization
of~\cite[Theorem~5.4.2]{Pcosh} (see
also~\cite[Corollary~6.4.12(b\+-c)]{Pcosh}).

\begin{prop} \label{locality-of-becker-contraacyclicity}
 Let $X$ be a quasi-compact semi-separated scheme, $\bW$ be an open
covering of $X$, and $X=\bigcup_{\alpha=1}^N U_\alpha$ be a finite
affine open covering of $X$ subordinate to\/~$\bW$.
 Let $\cA$ be a quasi-coherent quasi-algebra over~$X$.
 In this context: \hbadness=1075 \par
\textup{(a)} a complex\/ $\gB^\bu$ in $\cA\Lcth_\bW$ is
Becker-contraacyclic if and only if, for every index~$\alpha$,
the complex of $\cO_X(U_\alpha)$\+contraadjusted
$\cA(U_\alpha)$\+modules\/ $\gB^\bu[U_\alpha]$ is Becker-contraacyclic
in the exact category of $\cO_X(U_\alpha)$\+contraadjusted
$\cA(U_\alpha)$\+modules $\cA(U_\alpha)\Modl^{\cO_X(U_\alpha)\dcta}$;
\par
\textup{(b)} a complex\/ $\gB^\bu$ in $\cA\Lcth_\bW^{X\dlct}$ is
Becker-contraacyclic if and only if, for every index~$\alpha$,
the complex of $\cO_X(U_\alpha)$\+cotorsion $\cA(U_\alpha)$\+modules\/
$\gB^\bu[U_\alpha]$ is Becker-contraacyclic in the exact category
of $\cO_X(U_\alpha)$\+cotorsion $\cA(U_\alpha)$\+modules
$\cA(U_\alpha)\Modl^{\cO_X(U_\alpha)\dcot}$; {\hbadness=2825\par}
\textup{(c)} a complex\/ $\gB^\bu$ in $\cA\Lcth_\bW^{\cA\dlct}$ is
Becker-contraacyclic if and only if, for every index~$\alpha$,
the complex of cotorsion $\cA(U_\alpha)$\+modules\/ $\gB^\bu[U_\alpha]$
is Becker-contraacyclic in the exact category of cotorsion
$\cA(U_\alpha)$\+modules $\cA(U_\alpha)\Modl^\cot$.
\end{prop}

\begin{proof}
 The argument is dual-analogous to the proof of
Corollary~\ref{locality-of-becker-coacyclicity}.
 Let us prove part~(a).
 The ``only if'' implication holds by
Lemma~\ref{lcth-becker-contraacycl-inverse-image}(a).
 The ``if'' implication is nontrivial.
 Consider the \v Cech complex of complexes of $\bW$\+locally
contraherent $\cA$\+modules
(cf.~\eqref{lcth-sheaf-of-rings-cech-resolution})
\begin{multline} \label{lcth-sheaf-of-rings-cech-resolution-of-complex}
 0\lrarrow j_{1,\dotsc,N}{}_!j_{1,\dotsc,N}^!\gB^\bu \lrarrow\dotsb \\
 \lrarrow\bigoplus\nolimits_{1\le\alpha<\beta\le N}
 j_{\alpha,\beta}{}_!j_{\alpha,\beta}^!\gB^\bu\lrarrow
 \bigoplus\nolimits_{\alpha=1}^N j_\alpha{}_!j_\alpha^!\gB^\bu
 \lrarrow\gB^\bu\lrarrow0,
\end{multline}
where, for every subsequence of indices $1\le\alpha_1<\dotsb<\alpha_k
\le N$, we denote by $j_{\alpha_1,\dotsc,\alpha_k}$ the open
immersion morphism $\bigcap_{s=1}^k U_{\alpha_s}\rarrow X$.
 Then, by Lemma~\ref{lcth-becker-contraacycl-inverse-image}(a) and
Corollary~\ref{lcth-becker-contraacycl-affine-direct-image}(a), all
the terms of~\eqref{lcth-sheaf-of-rings-cech-resolution-of-complex}
except perhaps the rightmost one are Becker-contraacyclic complexes
in $\cA\Lcth_\bW$.
 Since the finite complex of
complexes~\eqref{lcth-sheaf-of-rings-cech-resolution-of-complex} is
acyclic in $\cA\Lcth_\bW$ and all absolutely acyclic complexes are
Becker-contraacyclic, it follows that the complex $\gB^\bu$ is
Becker-contraacyclic in $\cA\Lcth_\bW$ as well.
 The proofs of parts~(b\+-c) are similar and based on
Corollary~\ref{lcth-becker-contraacycl-affine-direct-image}(b\+-c)
and Lemma~\ref{lcth-becker-contraacycl-inverse-image}(b\+-c).
\end{proof}

\subsection{Becker contraacyclicity in module categories}
 We start with restating the result of~\cite[Proposition~B.7.12]{Pcosh}
and its obvious generalization to $R$\+contraadjusted $A$\+modules
before passing to a partial generalization to $R$\+cotorsion
$A$\+modules.

 Given a homomorphism of associative rings $R\rarrow A$, we denote by
$A\Modl^{R\dcot}$ and $A\Modl_{A/R\dflp}\subset A\Modl$ the full exact
subcategories of $R$\+cotorsion and $A/R$\+flaprojective left
$A$\+modules, respectively.
 When $R$ is a commutative ring, the similar notation $A\Modl^{R\dcta}$
and $A\Modl_{A/R\dvflp}\subset A\Modl$ stands for the full exact
subcategories of $R$\+contraadjusted and $A/R$\+very flaprojective
left $A$\+modules.
 The exact categories of cotorsion and flat left $A$\+modules are
denoted by $A\Modl^\cot$ and $A\Modl_\fl$, respectively.
 For any exact category $\sE$, the notation $\Ac(\sE)\subset\Com(\sE)$
stands for the full subcategory of acyclic complexes.

\begin{prop} \label{becker-contraacycl-in-module-categories-agrees}
\textup{(a)} Let $R$ be a commutative ring, $A$ be an associative ring,
and $R\rarrow A$ be a ring homomorphism.
 Let $B^\bu$ be a complex of $R$\+contraadjusted left $A$\+modules.
 Then $B^\bu$ is Becker-contraacyclic as a complex in $A\Modl^{R\dcta}$
if and only if $B^\bu$ is Becker-contraacyclic as a complex in $A\Modl$.
 Moreover, the inclusion of exact/abelian categories $A\Modl^{R\dcta}
\rarrow A\Modl$ induces a triangulated equivalence of the Becker 
contraderived categories\/ $\sD^\bctr(A\Modl^{R\dcta})\rarrow
\sD^\bctr(A\Modl)$. \par
\textup{(b)} Let $A$ be an associative ring and $B^\bu$ be a complex
of cotorsion left $A$\+modules.
 Then $B^\bu$ is Becker-contraacyclic as a complex in $A\Modl^\cot$
if and only if $B^\bu$ is Becker-contraacyclic as a complex in $A\Modl$.
 Moreover, the inclusion of exact/abelian categories $A\Modl^\cot\rarrow
A\Modl$ induces a triangulated equivalence of the Becker contraderived
categories\/ $\sD^\bctr(A\Modl^\cot)\rarrow\sD^\bctr(A\Modl)$.
\end{prop}

\begin{proof}
 Part~(b) is~\cite[Proposition~B.7.12(b)]{Pcosh}
or~\cite[Theorem~7.19 and Corollary~7.21]{Pphil}.
 Part~(a) is a generalization of~\cite[Proposition~B.7.12(a)]{Pcosh}
provable in the same way.
 Specifically, the pair of classes of $A/R$\+very flaprojective left
$A$\+modules $\sF=A\Modl_{A/R\dvflp}$ and $R$\+contraadjusted left
$A$\+modules $\sC=A\Modl^{R\dcta}$ is a hereditary complete cotorsion
pair in the abelian category $\sE=A\Modl$ (by
Lemmas~\ref{very-flaprojective-cotorsion-pair-hereditary},
\ref{very-flaprojective-is-a-cotorsion-pair}(b), and
Theorem~\ref{very-flaprojective-cotorsion-pair-complete}) and
the $\sC$\+coresolution dimensions of all objects of $\sE$ do not
exceed a fixed constant $n=1$.
 So~\cite[Proposition~B.7.10]{Pcosh} is applicable.
\end{proof}

\begin{lem} \label{becker-contraacycl-in-module-category-implies}
 Let $R\rarrow A$ be a homomorphism of associative rings.
 Then any Becker-contraacyclic complex in $A\Modl^{R\dcot}$ is also
Becker-contraacyclic in $A\Modl$.
 So the inclusion of exact/abelian categories $A\Modl^{R\dcot}\rarrow
A\Modl$ induces a well-defined triangulated functor\/
$\sD^\bctr(A\Modl^\cot)\rarrow\sD^\bctr(A\Modl)$.
\end{lem}

\begin{proof}
 This is a special case of the result of~\cite[Corollary~5.2]{Pfltp}
for the nested pair of cotorsion pairs ($A\Modl_\proj$, $A\Modl$)
and ($A\Modl_{A/R\dvflp}$, $A\Modl^{R\dcot}$) in $A\Modl$ (where
$A\Modl_\proj$ denotes the class of projective left $A$\+modules).
\end{proof}

\begin{prop} \label{flaprojective-conjecture-prop}
 Let $R\rarrow A$ be a homomorphism of associative rings.
 Then the following conditions are equivalent:
\begin{enumerate}
\item any complex of $R$\+cotorsion left $A$\+modules that is
Becker-contraacyclic as a complex in $A\Modl$ is also
Becker-contraacyclic as a complex in $A\Modl^{R\dcot}$;
\item the triangulated functor\/ $\sD^\bctr(A\Modl^{R\dcot})\rarrow
\sD^\bctr(A\Modl)$ induced by the inclusion of exact/abelian categories
$A\Modl^{R\dcot}\rarrow A\Modl$ is a triangulated equivalence;
\item in any acyclic complex of $A/R$\+flaprojective left $A$\+modules
with flat $A$\+modules of cocycles, the $A$\+modules of cocycles are
actually $A/R$\+flaprojective;
\item the triangulated functor\/ $\sD(A\Modl_{A/R\dflp})\rarrow
\sD(A\Modl_\fl)$ induced by the inclusion of exact categories
$A\Modl_{A/R\dflp}\rarrow A\Modl_\fl$ is a triangulated equivalence.
\item the triangulated functor\/ $\sD^\bco(A\Modl_{A/R\dflp})\rarrow
\sD^\bco(A\Modl_\fl)$ induced by the inclusion of exact categories
$A\Modl_{A/R\dflp}\rarrow A\Modl_\fl$ is a triangulated equivalence.
\end{enumerate}
\end{prop}

\begin{proof}
 This is a part of~\cite[Theorem~9.9 and Corollary~10.5]{Pfltp}.
\end{proof}

 The \emph{Flaprojective Conjecture} predicts that the equivalent
conditions of Proposition~\ref{flaprojective-conjecture-prop} hold for
all homomorphisms of associative rings $R\rarrow A$.
 See~\cite[Conjecture~0.1 or Conjecture~10.6]{Pfltp}.

 It is clear from~\cite[end of Section~7 and Example~8.7]{Pfltp} that
the equivalent conditions of
Proposition~\ref{flaprojective-conjecture-prop} are satisfied for
a ring homorphism $R\rarrow A$ whenever all flat left $R$\+modules
have finite projective dimensions.
 In particular, this holds whenever the ring $R$ is commutative and
Noetherian of finite Krull dimension~\cite[Corollaire~II.3.2.7]{RG}.
 See Theorem~\ref{lcth-lcta-lct-cdg-abs-contraderived-equiv-thm} below
for a nonaffine version.

\subsection{Forgetful functor preserves Becker contraacyclicity}
 We start with another module-theoretic lemma.

\begin{lem} \label{module-restr-of-scal-preserves-becker-contraacycl}
 Let $R\rarrow A$ be a homomorphism of associative rings.
 Then any Becker-contraacyclic complex in $A\Modl$ is also
Becker-contraacyclic as a complex in $R\Modl$.
\end{lem}

\begin{proof}
 The forgetful (restriction-of-scalars) functor $A\Modl\rarrow R\Modl$
is exact and has a left adjoint functor
$A\ot_R{-}\,\:R\Modl\rarrow A\Modl$.
 By~\cite[Lemma~B.7.5(b)]{Pcosh} or
Lemma~\ref{DG-functors-preserve-Becker-co-contra-acyclicity}(b) above,
it follows that the forgetful functor takes Becker-contraacyclic
complexes in $A\Modl$ to Becker-contraacyclic complexes in $R\Modl$.
\end{proof}

 Now we can prove the main result of this
Section~\ref{becker-semicontraderived-secn}.

\begin{cor} \label{forg-functor-preserves-becker-contraacyclicity}
  Let $X$ be a quasi-compact semi-separated scheme with an open
covering\/ $\bW$ and $\cA$ be a quasi-coherent quasi-algebra over~$X$.
 In this context: \par
\textup{(a)} the forgetful functor $\cA\Lcth_\bW\rarrow X\Lcth_\bW$
takes Becker-contraacyclic complexes in $\cA\Lcth_\bW$ to
Becker-contraacyclic complexes in $X\Lcth_\bW$; \par
\textup{(b)} the forgetful functor $\cA\Lcth_\bW^{X\dlct}\rarrow
X\Lcth_\bW^\lct$ takes Becker-contraacyclic complexes in
$\cA\Lcth_\bW^{X\dlct}$ to Becker-contraacyclic complexes
in $X\Lcth_\bW^\lct$; \par
\textup{(c)} the forgetful functor $\cA\Lcth_\bW^{\cA\dlct}\rarrow
X\Lcth_\bW^\lct$ takes Becker-contraacyclic complexes in
$\cA\Lcth_\bW^{\cA\dlct}$ to Becker-contraacyclic complexes
in $X\Lcth_\bW^\lct$.
\end{cor}

\begin{proof}
 This is actually based on the easy ``only if'' implications in
Proposition~\ref{locality-of-becker-contraacyclicity} and the difficult
``if'' implications in~\cite[Theorem~5.4.2]{Pcosh}.
 Let us prove part~(b).
 Let $X=\bigcup_\alpha U_\alpha$ be an affine open covering of $X$ and
$\gB^\bu$ be a Becker-contraacyclic complex in $\cA\Lcth_\bW^{X\dlct}$.
 By Proposition~\ref{locality-of-becker-contraacyclicity}(b),
the complex of $\cO_X(U_\alpha)$\+cotorsion $\cA(U_\alpha)$\+modules
$\gB^\bu[U_\alpha]$ is Becker-contraacyclic in
$\cA(U_\alpha)\Modl^{\cO_X(U_\alpha)\dcot}$ for every index~$\alpha$.
 By Lemma~\ref{becker-contraacycl-in-module-category-implies},
it follows that the complex of $\cA(U_\alpha)$\+modules
$\gB^\bu[U_\alpha]$ is also Becker-contraacyclic in the abelian
category $\cA(U_\alpha)\Modl$.
 Applying
Lemma~\ref{module-restr-of-scal-preserves-becker-contraacycl},
we see that the complex $\gB^\bu[U_\alpha]$ is Becker-contraacyclic
in the abelian category $\cO_X(U_\alpha)\Modl$.
 By~\cite[Proposition~B.7.12(b)]{Pcosh} or
Proposition~\ref{becker-contraacycl-in-module-categories-agrees}(b)
above (for the ring $\cO_X(U_\alpha)$), it follows that the complex
$\gB^\bu[U_\alpha]$ is Becker-contraacyclic in the exact category
$\cO_X(U_\alpha)\Modl^\cot$.
 Finally, using~\cite[Theorem~5.4.2(b)]{Pcosh} or
Proposition~\ref{locality-of-becker-contraacyclicity}(b)
(for the quasi-coherent quasi-algebra $\cO_X$ on~$X$), we can
conclude that the complex $\gB^\bu$ is Becker-contraacyclic in
the exact category $X\Lcth_\bW^\lct$.
 The proofs of parts~(a) and~(c) are similar and use
Proposition~\ref{becker-contraacycl-in-module-categories-agrees}(a\+-b)
instead of Lemma~\ref{becker-contraacycl-in-module-category-implies}.
\end{proof}

\begin{cor} \label{becker-contraacyclic-complexes-are-acyclic}
 Let $X$ be a quasi-compact semi-separated scheme with an open
covering\/ $\bW$ and $\cA$ be a quasi-coherent quasi-algebra over~$X$.
 In this context: \par
\textup{(a)} any Becker-contraacyclic complex in the exact category
$\cA\Lcth_\bW$ is acyclic in $\cA\Lcth_\bW$; \par
\textup{(b)} any Becker-contraacyclic complex in the exact category
$\cA\Lcth_\bW^{X\dlct}$ is acyclic in $\cA\Lcth_\bW^{X\dlct}$; \par
\textup{(c)} any Becker-contraacyclic complex in the exact
category $\cA\Lcth_\bW^{\cA\dlct}$ is acyclic in
$\cA\Lcth_\bW^{\cA\dlct}$.
\end{cor}

\begin{proof}
 In view of
Corollary~\ref{forg-functor-preserves-becker-contraacyclicity}
and Lemma~\ref{lcth-forgetful-functor-reflects-acyclicity},
it suffices to show that any Becker-contraacyclic complex in
$X\Lcth_\bW$ is acyclic in $X\Lcth_\bW$ and any Becker-contraacyclic
complex in $X\Lcth_\bW^\lct$ is acyclic in $X\Lcth_\bW^\lct$.
 This is the result of~\cite[Corollary~5.4.3]{Pcosh}.
\end{proof}

\subsection{The Becker semicontraderived category}
\label{becker-semicontraderived-defined-subsecn}
 Now we can spell out the definitions of the Becker semicontraderived
categories, similarly to the definitions of the Becker semicoderived
categories in Section~\ref{semiderived-quasi-coherent-subsecn}
and the definitions of the Positselski semicontraderived categories
in Section~\ref{semiderived-contraherent-subsecn}.
 Let $X$ be a quasi-compact semi-separated scheme with an open
covering $\bW$ and $\cA$ be a quasi-coherent quasi-algebra over~$X$.

 A complex in the exact category $\cA\Lcth_\bW$ is said to be
\emph{Becker-semicontraacyclic} (relative to $X\Lcth_\bW$) if,
\emph{viewed as a complex in $X\Lcth_\bW$}, it is
Becker-contraacyclic in the exact category $X\Lcth_\bW$.
 The full subcategory of Becker-semicontraacyclic complexes is
denoted by $\Ac^\bsi(\cA\Lcth_\bW)\subset\Hot(\cA\Lcth_\bW)$.
 By Corollary~\ref{forg-functor-preserves-becker-contraacyclicity}(a),
all Becker-contraacyclic complexes in $\cA\Lcth_\bW$ are
Becker-semicontraacyclic.
 By~\cite[Corollary~5.4.3(a)]{Pcosh} or
Corollary~\ref{becker-contraacyclic-complexes-are-acyclic}(a) above
(for the quasi-coherent quasi-algebra $\cO_X$ over~$X$),
all Becker-contraacyclic complexes in $X\Lcth_\bW$ are acyclic
in $X\Lcth_\bW$.
 Therefore, it follows from
Lemma~\ref{lcth-forgetful-functor-reflects-acyclicity}(a) that all
Becker-semicontraacyclic complexes in $\cA\Lcth_\bW$ are acyclic
in $\cA\Lcth_\bW$.
 So we have
$$
 \Ac^\bctr(\cA\Lcth_\bW)\subset\Ac^\bsi(\cA\Lcth_\bW)\subset
 \Ac(\cA\Lcth_\bW).
$$
 The triangulated Verdier quotient category
$$
 \sD^\bsi(\cA\Lcth_\bW)=\Hot(\cA\Lcth_\bW)/\Ac^\bsi(\cA\Lcth_\bW)
$$
is called the Becker \emph{semi}(\emph{contra})\emph{derived category}
of $\bW$\+locally contraherent $\cA$\+mod\-ules on~$X$.

 A complex in the exact category $\cA\Lcth_\bW^{X\dlct}$ is said to be
\emph{Becker-semi\-con\-tra\-acyclic} (relative to
$X\Lcth_\bW^\lct$) if, \emph{viewed as a complex in $X\Lcth_\bW^\lct$},
it is Becker-contraacyclic in the exact category $X\Lcth_\bW^\lct$.
 The full subcategory of Becker-semicontraacyclic complexes is denoted
by $\Ac^\bsi(\cA\Lcth_\bW^{X\dlct})\subset\Hot(\cA\Lcth_\bW^{X\dlct})$.
 By Corollary~\ref{forg-functor-preserves-becker-contraacyclicity}(b),
all Becker-contraacyclic complexes in $\cA\Lcth_\bW^{X\dlct}$ are
Becker-semicontraacyclic.
 By~\cite[Corollary~5.4.3(b)]{Pcosh} or
Corollary~\ref{becker-contraacyclic-complexes-are-acyclic}(b) above
(for the quasi-coherent quasi-algebra $\cO_X$ over~$X$), all
Becker-contraacyclic complexes in $X\Lcth_\bW^\lct$ are acyclic
in $X\Lcth_\bW^\lct$.
 Therefore, it follows from
Lemma~\ref{lcth-forgetful-functor-reflects-acyclicity}(b) that all
Becker-semicontraacyclic complexes in $\cA\Lcth_\bW^{X\dlct}$ are
acyclic in $\cA\Lcth_\bW^{X\dlct}$.
 So we have {\hbadness=1150
$$
 \Ac^\bctr(\cA\Lcth_\bW^{X\dlct})\subset\Ac^\bsi(\cA\Lcth_\bW^{X\dlct})
 \subset\Ac(\cA\Lcth_\bW^{X\dlct}).
$$
 The} triangulated Verdier quotient category
$$
 \sD^\bsi(\cA\Lcth_\bW^{X\dlct})=
 \Hot(\cA\Lcth_\bW^{X\dlct})/\Ac^\bsi(\cA\Lcth_\bW^{X\dlct})
$$
is called the Becker \emph{semi}(\emph{contra})\emph{derived category}
of $X$\+locally cotorsion $\bW$\+locally contraherent $\cA$\+modules
on~$X$.

 A complex in the exact category $\cA\Lcth_\bW^{\cA\dlct}$ is said to be
\emph{Becker-semi\-con\-tra\-acyclic} (relative to
$X\Lcth_\bW^\lct$) if, \emph{viewed as a complex in $X\Lcth_\bW^\lct$},
it is Becker-contraacyclic in the exact category $X\Lcth_\bW^\lct$.
 The full subcategory of Becker-semicontraacyclic complexes is
denoted by $\Ac^\bsi(\cA\Lcth_\bW^{\cA\dlct})\subset
\Hot(\cA\Lcth_\bW^{\cA\dlct})$.
 By Corollary~\ref{forg-functor-preserves-becker-contraacyclicity}(c),
all Becker-contraacyclic complexes in $\cA\Lcth_\bW^{\cA\dlct}$ are
Becker-semicontraacyclic.
 As mentioned in the previous paragraph, all Becker-contraacyclic
complexes in $X\Lcth_\bW^\lct$ are acyclic in $X\Lcth_\bW^\lct$.
 Therefore, it follows from
Lemma~\ref{lcth-forgetful-functor-reflects-acyclicity}(c) that all
Becker-semicontraacyclic complexes in $\cA\Lcth_\bW^{\cA\dlct}$
are acyclic in $\cA\Lcth_\bW^{\cA\dlct}$.
 So we have
$$
 \Ac^\bctr(\cA\Lcth_\bW^{\cA\dlct})\subset
 \Ac^\bsi(\cA\Lcth_\bW^{\cA\dlct})\subset
 \Ac(\cA\Lcth_\bW^{\cA\dlct}).
$$
 The triangulated Verdier quotient category
$$
 \sD^\bsi(\cA\Lcth_\bW^{\cA\dlct})=
 \Hot(\cA\Lcth_\bW^{\cA\dlct})/\Ac^\bsi(\cA\Lcth_\bW^{\cA\dlct})
$$
is called the Becker \emph{semi}(\emph{contra})\emph{derived category}
of $\cA$\+locally cotorsion $\bW$\+locally contraherent $\cA$\+modules
on~$X$.

 In the case of the open covering $\bW=\{X\}$, we will write
\begin{align*}
 \sD^\bsi(\cA\Ctrh)&=\sD^\bsi(\cA\Lcth_{\{X\}}), \\
 \sD^\bsi(\cA\Ctrh^{X\dlct})&=\sD^\bsi(\cA\Lcth_{\{X\}}^{X\dlct}), \\
 \sD^\bsi(\cA\Ctrh^{\cA\dlct})&=\sD^\bsi(\cA\Lcth_{\{X\}}^{\cA\dlct}),
\end{align*}
etc.

\begin{cor} \label{becker-semicontrader-independ-of-covering}
 Let $X$ be a quasi-compact semi-separated scheme with an open
covering\/ $\bW$ and $\cA$ be a quasi-coherent quasi-algebra over~$X$.
 In this context: \par
\textup{(a)} The inclusions of exact categories $\cA\Ctrh_\bW
\rarrow\cA\Lcth_\bW$ and $X\Ctrh_\bW\rarrow X\Lcth_\bW$ induce
a triangulated equivalence of the semiderived categories
$$
 \sD^\bsi(\cA\Ctrh)\simeq\sD^\bsi(\cA\Lcth_\bW).
$$ \par
\textup{(b)} The inclusions of exact categories $\cA\Ctrh_\bW^{X\dlct}
\rarrow\cA\Lcth_\bW^{X\dlct}$ and $X\Ctrh_\bW^\lct\allowbreak\rarrow
X\Lcth_\bW^\lct$ induce a triangulated equivalence of the semiderived
categories
$$
 \sD^\bsi(\cA\Ctrh^{X\dlct})\simeq\sD^\bsi(\cA\Lcth_\bW^{X\dlct}).
$$ \par
\textup{(c)} The inclusions of exact categories $\cA\Ctrh_\bW^{\cA\dlct}
\rarrow\cA\Lcth_\bW^{\cA\dlct}$ and $X\Ctrh_\bW^\lct\allowbreak\rarrow
X\Lcth_\bW^\lct$ induce a triangulated equivalence of the semiderived
categories
$$
 \sD^\bsi(\cA\Ctrh^{\cA\dlct})\simeq\sD^\bsi(\cA\Lcth_\bW^{\cA\dlct}).
$$
\end{cor}

\begin{proof}
 The proof is similar to that of
Corollary~\ref{positselski-semicontrader-independ-of-covering}.
 Instead of
Lemma~\ref{forg-functor-preserves-positselski-contraacyclicity},
one can refer to
Corollary~\ref{forg-functor-preserves-becker-contraacyclicity}
or~\cite[Lemma~B.7.1(a)]{Pcosh}.
\end{proof}

 For a triangulated equivalence $\sD^\bsi(\cA\Lcth_\bW^{X\dlct})
\simeq\sD^\bsi(\cA\Lcth_\bW)$ proved under more restrictive assumptions
on the scheme $X$, see
Corollary~\ref{semicontraderived-X-lcta-X-lct-A-mod-equivalence} below.

 The following corollary is only stated for quasi-compact
semi-separated schemes, as the Becker semicontraderived categories
were only defined above under these assumptions.

\begin{cor} \label{noetherian-Positselski=Becker-semicontraderived}
\textup{(a)} Let $X$ be a semi-separated Noetherian scheme of finite
Krull dimension with an open covering\/ $\bW$ and $\cA$ be
a quasi-coherent quasi-algebra over~$X$.
 Then the classes of Positselski-semicontraacyclic and
Becker-semicontraacyclic complexes of\/ $\bW$\+locally contraherent
$\cA$\+modules on $X$ coincide.
 So one has\/ $\Ac^\si(\cA\Lcth_\bW)=\Ac^\bsi(\cA\Lcth_\bW)$ and\/
$\sD^\si(\cA\Lcth_\bW)=\sD^\bsi(\cA\Lcth_\bW)$. {\hbadness=1800\par}
\textup{(b)} Let $X$ be a semi-separated Noetherian scheme with
an open covering\/ $\bW$ and $\cA$ be a quasi-coherent quasi-algebra
over~$X$.
 Then the classes of Positselski-semicontraacyclic and
Becker-semicontraacyclic complexes of $X$\+locally cotorsion\/
$\bW$\+locally contraherent $\cA$\+modules on $X$ coincide.
 So one has\/
$\Ac^\si(\cA\Lcth_\bW^{X\dlct})=\Ac^\bsi(\cA\Lcth_\bW^{X\dlct})$ and\/
$\sD^\si(\cA\Lcth_\bW^{X\dlct})=\sD^\bsi(\cA\Lcth_\bW^{X\dlct})$. \par
\textup{(c)} Let $X$ be a semi-separated Noetherian scheme with
an open covering\/ $\bW$ and $\cA$ be a quasi-coherent quasi-algebra
over~$X$.
 Then the classes of Positselski-semicontraacyclic and
Becker-semicontraacyclic complexes of $\cA$\+locally cotorsion\/
$\bW$\+locally contraherent $\cA$\+modules on $X$ coincide.
 So one has\/ $\Ac^\si(\cA\Lcth_\bW^{\cA\dlct})=
\Ac^\bsi(\cA\Lcth_\bW^{\cA\dlct})$ and\/
$\sD^\si(\cA\Lcth_\bW^{\cA\dlct})=\sD^\bsi(\cA\Lcth_\bW^{\cA\dlct})$.
\end{cor}

\begin{proof}
 Parts~(b\+-c) hold because, for any locally Noetherian scheme $X$,
the classes of Positselski-contraacyclic and Becker-contraacyclic
complexes in the exact category $X\Lcth_\bW^\lct$
coincide~\cite[Theorem~6.4.10(b)]{Pcosh}.
 Part~(a) holds because, for any Noetherian scheme $X$ of finite Krull
dimension, the classes of Positselski-contraacyclic and
Becker-contraacyclic complexes in the exact category $X\Lcth_\bW$
coincide~\cite[Theorem~6.4.10(d)]{Pcosh}.
\end{proof}

\Section{Twisted Lie Algebroids}  \label{twisted-lie-algebroids-secn}

\subsection{Finite locally free and locally finitely presented sheaves}
\label{loc-free-sheaves-subsecn}
 Let $X$ be a scheme and $\F$ be a quasi-coherent sheaf on~$X$.
 Following~\cite[Definitions Tags~00NW and~01C6]{SP}, we will say
that $\F$ is a \emph{finite locally free sheaf} if, for every
affine open subscheme $U\subset X$, the $\cO_X(U)$\+module $\F(U)$
is finitely generated and projective.
 It suffices to check this condition for affine open subschemes $U$
belonging to any given affine open covering of the scheme~$X$.
 Equivalently, $\F$ is finite locally free if and only if every
point $x\in X$ has an affine open neighborhood $x\in U\subset X$
such that the $\cO_X(U)$\+module $\F(U)$ is finitely generated and
free~\cite[Lemma Tag~00NX]{SP}.

 For any finite locally free sheaf $\F$ on $X$, the rank of
the free $\cO_{X,x}$\+module $\F_x$ is a locally constant function
of a point $x\in X$ \,\cite[Lemma Tag~01C9]{SP}.
 A quasi-coherent sheaf $\F$ on $X$ is said to be \emph{finite locally
free of} (\emph{constant}) \emph{rank~$r$} (where $r\ge0$ is an integer)
if every point $x\in X$ has an affine open neighborhood
$x\in U\subset X$ such that the $\cO_X(U)$\+module $\F(U)$ is free
of rank~$r$.
 A finite locally free sheaf $\F$ on $X$ has (constrant) rank~$r$ if
and only if, for every point $x\in X$, the free $\cO_{X,x}$\+module
$\F_x$ has rank~$r$.
 A finitely generated projective module $F$ over a commutative ring $R$
is said to be \emph{locally free of} (\emph{constant}) \emph{rank~$r$}
if the related quasi-coherent sheaf on $\Spec R$ is finite locally free
of constant rank~$r$.

 A finite locally free sheaf $\F$ on $X$ is said to have \emph{bounded
rank} if the locally constant function assigning to a point $x\in X$
the rank of the free $\cO_{X,x}$\+module $\F_x$ is bounded on~$X$.
 A finite locally free sheaf $\F$ has bounded rank if and only if
there exists a finite integer~$r$ such that every point $x\in X$
has an affine open neighborhood $x\in U\subset X$ for which
the $\cO_X(U)$\+module $\F(U)$ is free of rank~$\le r$.

 Given two quasi-coherent sheaves $\M$ and $\N$ on $X$, the sheaf
of $\cO_X$\+modules $\cHom_{\cO_X}(\M,\N)$ on $X$ is defined by
the rule
$$
 \cHom_{\cO_X}(\M,\N)(U)=\Hom_{\cO_X(U)}(\M(U),\N(U))
$$
for all affine open subschemes $U\subset X$, with the obvious
construction of the restriction maps $\Hom_{\cO_X(U)}(\M(U),\N(U))
\rarrow\Hom_{\cO_X(V)}(\M(V),\N(V))$ for pairs of affine open
subschemes $V\subset U\subset X$.
 The sheaf of $\cO_X$\+modules $\cHom_{\cO_X}(\M,\N)$ is \emph{not}
quasi-coherent for quasi-coherent sheaves $\M$ and $\N$ on $X$ in
general.
 However, for any finite locally free sheaf $\F$ and any
quasi-coherent sheaf $\N$ on $X$, the sheaf of $\cO_X$\+modules
$\cHom_{\cO_X}(\F,\N)$ \emph{is} quasi-coherent.

 For any two finite locally free sheaves $\F$ and $\G$ on $X$,
the quasi-coherent sheaves $\F\ot_{\cO_X}\G$ and
$\Hom_{\cO_X}(\F,\G)$ are finite locally free sheaves again.
 Moreover, there are natural isomorphisms of quasi-coherent sheaves
\begin{gather}
 \cHom_{\cO_X}(\F,\N)\simeq\cHom_{\cO_X}(\F,\cO_X)\ot_{\cO_X}\N; \\
 \cHom_{\cO_X}(\cHom_{\cO_X}(\F,\cO_X),\cO_X)\simeq\F
\end{gather}
for any finite locally free sheaf $\F$ and any quasi-coherent
sheaf $\N$ on~$X$.

 More generally, let us say that a quasi-coherent sheaf $\E$ on $X$ is
\emph{locally finitely presented} if, for every affine open subscheme
$U\subset X$, the $\cO_X(U)$\+module $\E(U)$ is finitely presented.
 By~\cite[Lemma~3.1(b)]{PS6}, it suffices to check this condition for
affine open subschemes $U$ belonging to any chosen affine open covering
of the scheme~$X$.

 For any locally finitely presented quasi-coherent sheaf $\E$ and
any quasi-coherent sheaf $\M$ on $X$, the sheaf of $\cO_X$\+modules
$\cHom_{\cO_X}(\E,\M)$ is quasi-coherent.
 This assertion follows immediately from the next lemma.

\begin{lem}
 Let $R\rarrow S$ be a homomorphism of commutative rings such that
$S$ is flat $R$\+module.
 Let $M$ be an $R$\+module and $E$ be a finitely presented $R$\+module.
 Then the natural $S$\+module map
$$
 S\ot_R\Hom_R(E,M)\lrarrow\Hom_S(S\ot_RE,\>S\ot_RM)
$$ 
is an isomorphism.
\end{lem}

\begin{proof}
 Represent $E$ as the cokernel of a morphism of finitely generated
free $R$\+modules and use the assumption that the tensor product
functor $S\ot_R{-}$ preserves kernels.
\end{proof}

\subsection{Lie algebroids} \label{lie-algebroids-subsecn}
 A \emph{Lie algebroid}~\cite[Section~1.2]{BB}, \cite[Section~1]{MM},
\cite[Section~10.9]{Prel}
(also known as a \emph{Lie--Rinehart algebra}~\cite[Section~2]{Rin})
is a set of data consisting of
\begin{itemize}
\item a commutative ring~$R$;
\item a Lie algebra~$\g$ over the ring of integers~$\boZ$;
\item an $R$\+module structure on~$\g$;
\item a $\g$\+module structure on~$R$.
\end{itemize}
 The following axioms are imposed:
\begin{enumerate}
\renewcommand{\theenumi}{\roman{enumi}}
\item $\g$~acts on $R$ by derivations of $R$, i.~e., one has
$v(rs)=v(r)s+rv(s)\in R$ for all $v\in\g$ and $r$, $s\in R$;
\item the action of $\g$ on $R$ is compatible with the $R$\+module
structures on $\g$ and~$R$, i.~e., one has $(rv)(s)=r\cdot v(s)\in R$
for all $r$, $s\in R$ and $v\in\g$;
\item $\g$~acts on $R$ and itself by derivations of the action of
$R$ on~$\g$, i.~e., one has $[v,rw]=v(r)w+r[v,w]\in\g$ for all
$v$, $w\in\g$ and $r\in R$.
\end{enumerate}

\begin{lem} \label{localizing-derivation-of-commutative-ring}
 Let $f\:R\rarrow S$ be a homomorphism of commutative rings such that
the related morphism of affine schemes\/ $\Spec S\rarrow\Spec R$ is
an open immersion.
 Let $D_R\:R\rarrow R$ be a derivation of the ring~$R$.
 Then there exists a unique derivation $D_S\:S\rarrow S$ of 
the ring $S$ such that $f\circ D_R=D_S\circ f$.
\end{lem}

\begin{proof}
 This is provable similarly to (but simpler than)
Lemma~\ref{localizing-derivation-by-affine-open-immersion}.
 The alternative argument based on~\cite[Proposition~6.6]{Ptd} is
also applicable, and proves the assertion of the lemma for any
flat epimorphism of commutative rings $f\:R\rarrow S$.
\end{proof}

\begin{lem} \label{localizing-lie-algebroid-lemma}
 Let $f\:R\rarrow S$ be a homomorphism of commutative rings such that
the related morphism of affine schemes\/ $\Spec S\rarrow\Spec R$ is
an open immersion, and let $(R,\g)$ be a Lie algebroid.
 Then there exists a unique Lie algebroid structure on the pair $(S,\h)$
with the $S$\+module\/ $\h=S\ot_R\g$ such that the natural map\/
$\g\rarrow\h$ is a Lie algebra homomorphism and the map $f\:R\rarrow S$
is homomorphism of\/ $\g$\+modules (where the\/ $\g$\+module structure
on $S$ comes from the\/ $\h$\+module structure and the Lie algebra
homomorphism\/ $\g\rarrow\h$).
\end{lem}

\begin{proof}
 For every element $v\in\g$, consider the derivation
$D_{R,v}\:R\rarrow R$ of the ring $R$ given by the action of~$v$
on $R$, that is $D_{R,v}(r)=v(r)$ for all $r\in R$.
 By Lemma~\ref{localizing-derivation-of-commutative-ring},
there exists a unique derivation $D_{S,v}\:S\rarrow S$ of
the ring $S$ such that $f\circ D_{R,v}=D_{S,v}\circ f$.
 The uniqueness assertion of
Lemma~\ref{localizing-derivation-of-commutative-ring} and the identities
$D_{R,rv}=rD_{R,v}$ and $D_{R,[v',v'']}=[D_{R,v'},D_{R,v''}]$
for all $r\in R$ and $v$, $v'$, $v''\in\g$ imply the identities
$D_{S,rv}=rD_{S,v}$ and $D_{S,[v',v'']}=[D_{S,v'},D_{S,v''}]$.

 Now we can define the Lie algebra structure on~$\h$ by the formula
$$
 [s'\ot v',\>s''\ot v''] = 
 s'D_{S,v'}(s'')\ot v'' - s''D_{S,v''}(s')\ot v'
 + s's''\ot[v',v'']
$$
for all $s'$, $s''\in S$ and $v'$, $v''\in\g$, and the action of $\h$
on $S$ by the rule $(s\ot v)(t)=sD_{S,v}(t)$ for all $s$, $t\in S$
and $v\in\g$.
 The uniqueness assertion of the lemma is obvious from this
construction.
\end{proof}

 Let $X$ be a scheme.
 A \emph{quasi-coherent Lie algebroid over~$X$} is a set of data
consisting of
\begin{itemize}
\item a quasi-coherent sheaf~$\g$ on~$X$;
\item a structure of a sheaf of Lie algebras (over the ring~$\boZ$)
on the sheaf~$\g$ on~$X$;
\item a structure of a sheaf of $\g$\+modules on the structure
sheaf $\cO_X$ of the scheme~$X$.
\end{itemize}
 The following axiom is imposed:
\begin{enumerate}
\renewcommand{\theenumi}{\roman{enumi}}
\setcounter{enumi}{3}
\item For every affine open subscheme $U\subset X$, the commutative
ring $\cO_X(U)$ and the Lie algebra $\g(U)$, together with
the $\cO_X(U)$\+module structure on $\g(U)$ that is a part of
the structure of a (quasi-coherent) sheaf of $\cO_X$\+modules on~$\g$,
and together with the $\g(U)$\+module structure on $\cO_X(U)$
that is a part of the structure of a sheaf of $\g$\+modules on $\cO_X$,
form a Lie algebroid.
 In other words, conditions~(i\+-iii) above are satisfied for $\cO_X(U)$
and~$\g(U)$.
\end{enumerate}

 It suffices to check condition~(iv) for affine open subschemes $U$
subordinate to any given open covering $\bW$ of the scheme~$X$.
 Furthermore, for an affine scheme $U$, the datum of a quasi-coherent
Lie algebroid~$\g$ over $U$ is equivalent to the datum of a Lie
algebroid $(R,\g)$ with the commutative ring $R=\cO(U)$.
 The latter assertion is clear from
Lemma~\ref{localizing-lie-algebroid-lemma}.

\subsection{Twisted Lie algebroids}
\label{twisted-lie-algebroids-subsecn}
 By a \emph{twisted Lie algebroid} we mean what otherwise might be
called a central extension of a Lie algebroid.
 Specifically, a twisted Lie algebroid is a set of data consiting of
\begin{itemize}
\item a Lie algebroid $(R,\g)$;
\item a Lie algebroid $(R,\widetilde\g)$;
\item a short exact sequence of $R$\+modules $0\rarrow R\overset\iota
\rarrow\widetilde\g\overset\pi\rarrow\g\rarrow0$.
\end{itemize}
 The following axioms are imposed:
\begin{enumerate}
\renewcommand{\theenumi}{\roman{enumi}}
\setcounter{enumi}{4}
\item $\pi\:\widetilde\g\rarrow\g$ is a Lie algebra homomorphism;
\item the elements of $\iota(R)\subset\widetilde\g$ are central in
the Lie algebra~$\widetilde\g$, i.~e., $[\widetilde v,\iota(r)]=0$
in~$\widetilde\g$ for all $\widetilde v\in\widetilde\g$ and $r\in R$;
\item the elements of $\iota(R)\subset\widetilde\g$ act by zero on~$R$,
i.~e., $\widetilde r(s)=0$ in $R$ for all $\widetilde r=\iota(r)\in
\widetilde\g$ and $r$, $s\in R$.
\end{enumerate}

 We will say that the Lie algebroid $(R,\g)$ is the \emph{underlying
Lie algebroid} of a twisted Lie algebroid $(R,\g,\widetilde\g)$.
 A twisted Lie algebroid $(R,\g,\widetilde\g)$ is said to be
\emph{split} if a \emph{splitting map} $\sigma\:\g\rarrow\widetilde\g$
is chosen such that $\sigma$~is an $R$\+module map and a Lie
algebra homomorphism, and $\pi\circ\sigma=\id_\g$.
 The datum of a split twisted Lie algebroid $(R,\g,\widetilde\g)$ is
equivalent to the datum of the underlying Lie algebroid $(R,\g)$:
the Lie algebroid $(R,\widetilde\g)$ is uniquely recovered by setting
$\widetilde\g=R\oplus\g$.

 Let $X$ be a scheme.
 A \emph{quasi-coherent twisted Lie algebroid over~$X$} is a set of
data consisting of
\begin{itemize}
\item a quasi-coherent Lie algebroid~$\g$ over~$X$;
\item a quasi-coherent Lie algebroid~$\widetilde\g$ over~$X$;
\item a short exact sequence of quasi-coherent sheaves
$0\rarrow\cO_X\overset\iota\rarrow\widetilde\g\overset\pi\rarrow\g
\rarrow0$ on~$X$.
\end{itemize}
 The following axioms are imposed:
\begin{enumerate}
\renewcommand{\theenumi}{\roman{enumi}}
\setcounter{enumi}{7}
\item $\pi\:\widetilde\g\rarrow\g$ is a morphism of sheaves of
Lie algebras;
\item the subsheaf $\iota(\cO_X)\subset\widetilde\g$ is central in
the sheaf of Lie algebras~$\widetilde\g$, i.~e.,
$[\widetilde v,\iota(r)]=0$ in $\widetilde\g(U)$ for all
$\widetilde v\in\widetilde\g(U)$ and $r\in\cO_X(U)$, where $U\subset X$
is any open subscheme;
\item the subsheaf $\iota(R)\subset\widetilde\g$ acts by zero
on $\cO_X$, i.~e., $\widetilde r(s)=0$ in $\cO_X(U$) for all
$\widetilde r=\iota(r)\in\widetilde\g(U)$ and $r$, $s\in\cO_X(U)$,
where $U\subset X$ is any open subscheme.
\end{enumerate}

 It suffices to check conditions~(ix\+-x) for affine open subschemes
$U$ subordinate to any given open covering $\bW$ of the scheme~$X$.
 For an affine scheme $U$, the datum of a quasi-coherent twisted Lie
algebroid $(\g,\widetilde\g)$ over $U$ is equivalent to the datum of
a twisted Lie algebroid $(R,\g,\widetilde\g)$ with the commutative
ring $R=\cO(U)$.

 We will say that the quasi-coherent Lie algebroid~$\g$ is
the \emph{underlying} (\emph{quasi-coherent}) \emph{Lie algebroid}
of a quasi-coherent twisted Lie algebroid~$(\g,\widetilde\g)$.
 A quasi-coherent twisted Lie algebroid $(\g,\widetilde\g)$ is said
to be \emph{split} if a \emph{splitting morphism} of quasi-coherent
sheaves $\sigma\:\g\rarrow\widetilde\g$ is chosen such that
$\sigma$~is also a morphism of sheaves of Lie algebras and
$\pi\circ\sigma=\id_\g$.
 The datum of a split quasi-coherent twisted Lie algebroid
$(\g,\widetilde\g)$ over $X$ is equivalent to the datum of
the underlying Lie algebroid~$\widetilde\g$: the quasi-coherent
Lie algebroid~$\widetilde\g$ is uniquely recovered by setting
$\widetilde\g=\cO_X\oplus\g$.

\subsection{Enveloping algebra} \label{enveloping-algebra-subsecn}
 The (\emph{twisted, universal}) \emph{enveloping ring}
$A_R(\g,\widetilde\g)$ of a twisted Lie algebroid
$(R,\g,\widetilde\g)$ (cf.~\cite[Section~2]{Rin},
\cite[Section~1.2.5]{BB}, \cite[Section~1]{MM},
\cite[Section~10.9]{Prel}) is constructed as  the associative ring
with the following generators and relations.

 The generators of $A_R(\g,\widetilde\g)$ are the symbols
$\lambda(\widetilde v)$, where $\widetilde v\in\widetilde\g$.
 To write down the relations, let us adopt the notational convention
that the symbol~$*$ denotes the multiplication in
$A_R(\g,\widetilde\g)$.
 Let $1_R\in R$ and $1_{A_R(\g,\widetilde\g)}\in A_R(\g,\widetilde\g)$
denote the unit elements of the respective rings.
 The relations are:
\begin{enumerate}
\renewcommand{\theenumi}{\alph{enumi}}
\item $\lambda(\widetilde v)+\lambda(\widetilde w)=
\lambda(\widetilde v+\widetilde w)$ for all $\widetilde v$,
$\widetilde w\in\widetilde\g$;
\item $\lambda(\iota(1_R))=1_{A_R(\g,\widetilde\g)}$;
\item $\lambda(\iota(r))*\lambda(\widetilde v)=
\lambda(r\widetilde v)$ for all $r\in R$ and
$\widetilde v\in\widetilde\g$;
\item $\lambda(\widetilde v)*\lambda(\iota(r))=
\lambda(r\widetilde v)+\lambda(\iota(\widetilde v(r)))$ for all
$r\in R$ and $\widetilde v\in\widetilde\g$;
\item $\lambda(\widetilde v)*\lambda(\widetilde w)-
\lambda(\widetilde w)*\lambda(\widetilde v)=
\lambda([\widetilde v,\widetilde w])$ for all $\widetilde v$,
$\widetilde w\in\widetilde\g$.
\end{enumerate}
 It is clear from relations~(b\+-c) that the map $\lambda\circ\iota\:
R\rarrow A_R(\g,\widetilde\g)$ is a ring homomorphism, while
relation~(e) implies that the map $\lambda\:\widetilde v\rarrow
A_R(\g,\widetilde\g)$ is a Lie algebra homomorphism from~$\g$ to
the associative ring $A_R(\g,\widetilde\g)$ viewed as a Lie algebra
(over~$\boZ$) with the commutator bracket.

 The natural multiplicative increasing filtration $F$ on the ring
$A_R(\g,\widetilde\g)$ is defined by the following rules.
 Firstly, $F_0A_R(\g,\widetilde\g)$ is the image of the ring
homomorphism $\lambda\circ\iota\:R\rarrow A_R(\g,\widetilde\g)$,
while $F_1A_R(\g,\widetilde\g)$ is the image of the map
$\lambda\:\widetilde\g\rarrow A_R(\g,\widetilde\g)$.
 Notice that $F_1A_R(\g,\widetilde\g)$ is an $R$\+$R$\+subbimodule
in $A_R(\g,\widetilde\g)$, or in other words, one has
$F_0A_R(\g,\widetilde\g)*F_1A_R(\g,\widetilde\g)\subset
F_1A_R(\g,\widetilde\g)$ and $F_1A_R(\g,\widetilde\g)*
F_0A_R(\g,\widetilde\g)\subset F_1A_R(\g,\widetilde\g)$,
as one can see from relations~(c\+-d).
 Secondly, the multiplicative increasing filtration $F$ on
$A_R(\g,\widetilde\g)$ is generated by $F_1$ over~$F_0$.
 This means that, for every $n\ge1$, the subgroup (in fact,
$R$\+$R$\+subbimodule) $F_nA_R(\g,\widetilde\g)\subset
A_R(\g,\widetilde\g)$ is spanned by the products
$\lambda(\widetilde v_1)\dotsm\lambda(\widetilde v_n)$
of at most~$n$ elements from $F_1A_R(\g,\widetilde\g)$
in $A_R(\g,\widetilde\g)$.
 Obviously, we have $A_R(\g,\widetilde\g)=\bigcup_{n\ge0}
F_nA_R(\g,\widetilde\g)$; so the filtration $F$ on
$A_R(\g,\widetilde\g)$ is exhaustive.

 For convenience, we put $F_{-1}A_R(\g,\widetilde\g)=0$.
 For any associative ring $A$ with a multiplicative increasing
filtration $0=F_{-1}A\subset F_0A\subset F_1A\subset F_2A\subset
\dotsb\subset A$, we denote by $\gr^F_*A=\bigoplus_{n=0}^\infty
F_nA/F_{n-1}A$ the associated graded ring of the filtered ring~$A$.

 It is clear from relations~(c\+-d) that the left and right
actions of the ring $R$ in the quotient $R$\+$R$\+bimodule
$F_1A_R(\g,\widetilde\g)/F_0A_R(\g,\widetilde\g)$ agree.
 It follows that the left and right actions of $R$ in the quotient
bimodule $F_nA_R(\g,\widetilde\g)/F_{n-1}A_R(\g,\widetilde\g)$ for
every integer $n\ge1$ agree as well.
 Therefore, $A_R(\g,\widetilde\g)$ is a quasi-module over $R$
in the sense of Section~\ref{prelim-quasi-modules-subsecn} (in fact,
a strong quasi-module in the sense of~\cite[Section~2.2]{Ptd}).
 Hence $A_R(\g,\widetilde\g)$ is a quasi-algebra over $R$ in the sense
of Section~\ref{prelim-quasi-algebras-subsecn}.

\begin{lem} \label{enveloping-algebra-commutes-with-localization}
 Let $f\:R\rarrow S$ be a homomorphism of commutative rings such that
the related morphism of affine schemes\/ $\Spec S\rarrow\Spec R$ is
an open immersion, and let $(R,\g,\widetilde\g)$ be a twisted
Lie algebroid.
 Consider the $S$\+modules\/ $\h=S\ot_R\g$ and\/ $\widetilde\h=
S\ot_R\widetilde\g$, and let $(S,\h,\widetilde\h)$ be the twisted Lie
algebroid constructed using Lemma~\ref{localizing-lie-algebroid-lemma}.
 Put $A=A_R(\g,\widetilde\g)$ and $B=A_S(\h,\widetilde\h)$; so
$A$ is a quasi-algebra over $R$ and $B$ is a quasi-algebra over~$S$.
 Then the ring homomorphism $A\rarrow B$ induced by the natural map\/
$\widetilde\g\rarrow\widetilde\h$ induces an isomorphism of
quasi-algebras $S\ot_RA\simeq B$ over~$S$
(cf.\ Lemma~\ref{quasi-algebra-co-extension-of-scalars}(a)).
 Moreover, the isomorphism $S\ot_RA\simeq B$ restricts to isomorphisms
of the filtration components $S\ot_RF_nA\simeq F_nB$ for all $n\ge0$.
\end{lem}

\begin{proof}
 Clear from the definitions.
\end{proof}

 In the case of a split twisted Lie algebroid $\widetilde\g=R\oplus\g$
(as defined in Section~\ref{twisted-lie-algebroids-subsecn}), we put
$A_R(\g)=A_R(\g,\widetilde\g)$.
 This is the usual notion of the enveloping algebra of a Lie
algebroid $(R,\g)$.

 Let $X$ be a scheme and $(\g,\widetilde\g)$ be a quasi-coherent
twisted Lie algebroid over~$X$.
 Then, in view of
Lemma~\ref{enveloping-algebra-commutes-with-localization}, the rule
$\cA_X(\g,\widetilde\g)(U)=A_{\cO_X(U)}(\g(U),\widetilde\g(U))$ for
all affine open subschemes $U\subset X$ defines a quasi-coherent
quasi-algebra $\cA_X(\g,\widetilde\g)$ over~$X$.
 The rule $F_n\cA_X(\g,\widetilde\g)(U)=
F_nA_{\cO_X(U)}(\g(U),\widetilde\g(U))$ for affine open subschemes
$U\subset X$ defines a natural exhaustive multiplicative increasing
filtration by quasi-coherent sub-quasi-modules
$F_n\cA_X(\g,\widetilde\g)\subset\cA_X(\g,\widetilde\g)$, \,$n\ge0$,
on the quasi-coherent quasi-algebra $\cA_X(\g,\widetilde\g)$.
 The quasi-coherent quasi-algebra $\cA_X(\g,\widetilde\g)$ is called
the (\emph{twisted, universal}) \emph{enveloping quasi-algebra} of
the quasi-coherent twisted Lie algebroid $(\g,\widetilde\g)$ over~$X$.

 For any quasi-coherent quasi-algebra $\cA$ over $X$ endowed with
a multiplicative increasing filtration $0=F_{-1}\cA\subset F_0\cA
\subset F_1\cA\subset F_2\cA\subset\dotsb\subset\cA$ by quasi-coherent
sub-quasi-modules, we denote by $\gr^F_*\cA=\bigoplus_{n=0}^\infty
F_n\cA/F_{n-1}\cA$ the associated graded quasi-coherent quasi-algebra
of the filtered quasi-coherent quasi-algebra~$\cA$.

 In the case of a split quasi-coherent twisted Lie algebroid
$\widetilde\g=\cO_X\oplus\g$ (as defined in
Section~\ref{twisted-lie-algebroids-subsecn}), we put
$\cA_X(\g)=\cA_X(\g,\widetilde\g)$.
 This is the usual notion of the enveloping algebra of a quasi-coherent
Lie algebroid~$\g$ over a scheme~$X$.

\subsection{Poincar\'e--Birkhoff--Witt theorem}
\label{lie-algebroid-pbw-theorem-subsecn}
 Let $R$ be a commutative ring and $M$ be an $R$\+module.
 The \emph{symmetric algebra} $\Sym_R^*(M)$ of the $R$\+module $M$ is
the nonnegatively graded commutative, associative $R$\+algebra freely
generated by the $R$\+module $M$ sitting in degree~$1$.
 In other words, $\Sym_R^*(M)$ is the associative $R$\+algebra generated
by the $R$\+module $M$ subject to the commutativity relations
$m'*m''=m''*m'$ for all $m'$, $m''\in M$.
 Here we denote by~$*$ the multiplication in the symmetric algebra.
 So the low-degree grading components of $\Sym_R^*(M)$ are
$\Sym_R^0(M)=R$ and $\Sym_R^1(M)=M$, while $\Sym_R^2(M)$ is
the quotient $R$\+module of $M\ot_RM$ by the submodule spanned
by the elements $m'\ot m''-m''\ot m'$, etc.

\begin{thm} \label{twisted-lie-algebroid-pbw-theorem}
 Let $(R,\g,\widetilde\g)$ be a twisted Lie algebroid.
 Then \par
\textup{(a)} the map $\lambda\:\widetilde\g\rarrow
F_1A_R(\g,\widetilde\g)$ induces a surjective homomorphism of
nonnegatively graded $R$\+algebras\/ $\Sym_R^*(\g)\rarrow
\gr^F_*A_R(\g,\widetilde\g)=\bigoplus_{n=0}^\infty
F_nA_R(\g,\widetilde\g)/F_{n-1}A_R(\g,\widetilde\g)$; \par
\textup{(b)} if the $R$\+module\/~$\g$ is flat, then the natural map\/
$\Sym_R^*(\g)\rarrow\gr^F_*A_R(\g,\widetilde\g)$ is an isomorphism
of graded $R$\+algebras.
\end{thm}

\begin{proof}
 In part~(a), the point is that the associated graded algebra
$\gr^FA_R(\g,\widetilde\g)$ is commutative, as one can see from
the relations~(d\+-e) from Section~\ref{enveloping-algebra-subsecn}.
 Furthermore, the image of the composition $\lambda\circ\iota\:R
\rarrow F_1A_R(\g,\widetilde\g)$ is contained in (in fact, coincides
with) $F_0A_R(\g,\widetilde\g)$.
 Therefore, the composition $\widetilde\g\rarrow
F_1A_R(\g,\widetilde\g)/F_0A_R(\g,\widetilde\g)$ annihilates
$\iota(R)\subset\widetilde\g$, so the map $\widetilde\g\rarrow
F_1A_R(\g,\widetilde\g)/F_0A_R(\g,\widetilde\g)$ factorizes through
the map $\pi\:\widetilde\g\rarrow\g$ and induces a map $\bar\lambda\:
\g\rarrow F_1A_R(\g,\widetilde\g)/F_0A_R(\g,\widetilde\g)$.
 The map~$\bar\lambda$ is an $R$\+module morphism by relations~(c\+-d).
 The map~$\bar\lambda$ is surjective, because the map $\lambda\:
\g\rarrow F_1A_R(\g,\widetilde\g)$ is.
 The ring homomorphism $\lambda\circ\iota\:R\rarrow
F_0A_R(\g,\widetilde\g)$ is surjective, too.
 Finally the graded ring $\gr^FA_R(\g,\widetilde\g)$ is generated
by $F_1A_R(\g,\widetilde\g)/F_0A_R(\g,\widetilde\g)$ over
$F_0A_R(\g,\widetilde\g)$ by the definition of the filtration $F$
on $A_R(\g,\widetilde\g)$.
 Part~(a) follows immediately from these observations.

 In the case when the $R$\+module~$\g$ is finitely generated and
projective (which is the only case that will be relevant for us in
the rest of this paper), part~(b) can be obtained as a special case
of~\cite[Theorem~4.25]{Prel}.
 The construction of the Chevalley--Eilenberg CDG\+ring in
Section~\ref{chevalley-eilenberg-cdg-ring-subsecn} below is
important in this context.
 The Chevalley--Eilenberg CDG\+ring $C^\cu_R(\g,\widetilde\g)$
is the right finitely projective Koszul CDG\+ring nonhomogeneous
Koszul dual to the left finitely projective nonhomogeneous Koszul
filtered ring $(A_R(\g,\widetilde\g),F)$ in the sense of~\cite{Prel};
so one needs to apply~\cite[Theorem~4.25]{Prel} to the CDG\+ring
$C^\cu_R(\g,\widetilde\g)$ in order to deduce part~(b) when $R$ is
a finitely generated projective $R$\+module.
 The general case of a flat $R$\+module~$\g$ is provable similarly
(one only needs to take care to avoid the double passage to the dual
$R$\+modules, and then the same argument works), or it can
be obtained as a special case of the complicated general result
of~\cite[Theorem~11.6]{Psemi}.
 For another version of the PBW theorem for Lie algebroids, see
also~\cite{Cal}.
\end{proof}

 Let $R\rarrow S$ be a homomorphism of commutative rings and $M$
be an $R$\+module.
 Consider the $S$\+module $N=S\ot_RM$.
 Then the map of rings $R\rarrow S$ and the natural $R$\+module map
$M\rarrow N$ induce a homomorphism of graded $R$\+algebras
$\Sym_R^*(M)\rarrow\Sym_S^*(N)$.
 One can easily see that the induced homomorphism of graded
$S$\+algebras $S\ot_R\Sym_R^*(M)\rarrow\Sym_S^*(N)$ is an isomorphism.

 Now let $\M$ be a quasi-coherent sheaf over a scheme~$X$.
 Then it follows from the previous paragraph that the rule
$\Sym_X^*(\M)(U)=\Sym_{\cO_X(U)}^*(\M(U))$ for all affine open
subschemes $U\subset X$ defines a quasi-coherent graded algebra
$\Sym_X^*(\M)$ over~$X$.
 The quasi-coherent graded algebra $\Sym_X^*(\M)$ is called
the \emph{symmetric algebra} of a quasi-coherent sheaf $\M$ on~$X$.

\begin{cor} \label{twisted-lie-algebroid-pbw-qcoh-cor}
 Let $(\g,\widetilde\g)$ be a twisted Lie algebroid over a scheme~$X$.
 Then \par
\textup{(a)} the morphism of quasi-coherent sheaves
$\lambda\:\widetilde\g\rarrow F_1\cA_X(\g,\widetilde\g)$ on $X$
(where the quasi-coherent quasi-module $F_1\cA_X(\g,\widetilde\g)$ is
viewed as a quasi-coherent sheaf with respect to its left
$\cO_X$\+module structure) induces a surjective homomorphism of
quasi-coherent nonnegatively graded algebras\/ $\Sym_X^*(\g)\rarrow
\gr^F_*\cA_X(\g,\widetilde\g)=\bigoplus_{n=0}^\infty
F_n\cA_X(\g,\widetilde\g)/F_{n-1}\cA_X(\g,\widetilde\g)$ over~$X$; \par
\textup{(b)} if the quasi-coherent sheaf\/~$\g$ over $X$ is flat, then
the natural map\/ $\Sym_X^*(\g)\rarrow\gr^F_*\cA_X(\g,\widetilde\g)$
is an isomorphism of quasi-coherent graded algebras over~$X$.
\end{cor}

\begin{proof}
 Follows from Theorem~\ref{twisted-lie-algebroid-pbw-theorem}
by passing to the affine open subschemes $U\subset X$.
 As a terminological observation, recall that the $\cO_X(U)$\+module
$\g(U)$ is flat for all such $U$ if and only if the quasi-coherent
sheaf~$\g$ on $X$ is flat, and $\cO_X(U)$\+module $\g(U)$ is finitely
generated and projective for all $U$ if and only if the quasi-coherent
sheaf~$\g$ on $X$ is finitely locally free in the sense of
Section~\ref{loc-free-sheaves-subsecn}.
\end{proof}

\subsection{Skew-symmetric tensors}
\label{skew-symmetric-tensors-subsecn}
 Let $R$ be a commutative ring and $M$ be an $R$\+module.
 The \emph{exterior algebra} $\bigwedge^*_R(M)$ of the $R$\+module $M$
is the nonnegatively graded associative $R$\+algebra generated by
the $R$\+module $M$ sitting in degree~$1$, subject to
the graded commutativity relations $m'\wedge m''+m''\wedge m'=0$
and $m\wedge m=0$ for all $m$, $m'$, $m''\in M$.
 Here we denote by~$\wedge$ the multiplication in the exterior algebra.
 So the low-degree grading components of $\bigwedge^*_R(M)$ are
$\bigwedge_R^0(M)=R$ and $\bigwedge_R^1(M)=M$, while
$\bigwedge_R^2(M)$ is the quotient $R$\+module of $M\ot_RM$ by
the submodule spanned by the elements $m\ot m$, etc.

 Let $n\ge0$ be an integer.
 An element $t\in M^{\ot_R\,n}=M\ot_RM\ot_R\dotsb\ot_RM$ ($n$~factors)
is called \emph{skew-symmetric} if, for every permutation~$\sigma$
of the set $\{1,2,\dotsc,n\}$ indexing the tensor factors, one has
$\sigma(t)=(-1)^{|\sigma|}t$, where $|\sigma|$~is the parity of
the permutation~$\sigma$.
 The natural \emph{skew-symmetrization map}
$$
 \sk_n\:\bigwedge\nolimits_R^n(M)\lrarrow M^{\ot_R\,n}
$$
is an $R$\+module map given by the formula
$$
 \sk_n(m_1\wedge\dotsb\wedge m_n)=\sum\nolimits_{\sigma\in\Sigma_n}
 (-1)^{|\sigma|}m_{\sigma(1)}\ot\dotsb m_{\sigma(n)},
$$
where $\Sigma_n$ denotes the group of all permutations of
$\{1,2,\dotsc,n\}$ and $m_1$,~\dots, $m_n\in M$.
 One can easily check that the map $\sk_n$ is well-defined.
 Clearly, all the tensors in the image of the skew-symmetrization map
are skew-symmetric.

\begin{lem} \label{skew-symmetric-tensors-lemma}
 Let $R$ be a commutative ring and $M$ be a flat $R$\+module.
 Then, for every integer $n\ge0$, the map\/ $\sk_n$ is an isomorphism
between the $R$\+module $\bigwedge_R^n(M)$ and the $R$\+submodule
of skew-symmetric tensors in $M^{\ot_R\,n}$.
\end{lem}

\begin{proof}
 It is worth pointing out that the analogue of this lemma for
the symmetric powers $\Sym^n_R(M)$, the submodule of symmetric tensors
in $M^{\ot_R\,n}$, and the symmetrization map is \emph{only} true in
characteristic~$0$ (i.~e., for the rings $R$ containing the field of
rational numbers~$\boQ$).
 Cf.\ Lemma~\ref{symmetric-tensors-lemma} below.
 However, the skew-symmetric version of the assertion holds true
irrespectively of the characteristic.
 The proof is straightforward: representing $M$ as a filtered direct
limit of free $R$\+modules, one reduces the question to the case when
$M$ is a free $R$\+module.
 Then all the $R$\+modules involved are free, and the assertion is
can be checked in explicit natural bases of all the $R$\+modules
involved, induced by the choice of a free $R$\+module basis of~$M$.
\end{proof}

\begin{lem} \label{exterior-dual-lemma}
 Let $R$ be a commutative ring and $M$ be a finitely generated
projective $R$\+module.
 Then, for every integer $n\ge0$, there is a natural isomorphism
$$
 \bigwedge\nolimits_R^n(\Hom_R(M,R))\simeq
 \Hom_R\Bigl(\bigwedge\nolimits_R^n(M),R\Bigr)
$$
provided by the natural $R$\+bilinear pairing
$$
 \bigwedge\nolimits_R^n(M)\times\bigwedge\nolimits_R^n(\Hom_R(M,R))
 \lrarrow R
$$
given by the formula
$$
 \langle m_1\wedge\dotsb\wedge m_n,\>b_n\wedge\dotsb\wedge b_1\rangle
 = \sum\nolimits_{\sigma\in\Sigma_n}(-1)^{|\sigma|}
 \prod\nolimits_{i=1}^n\langle m_i,b_{\sigma(i)}\rangle,
$$
where $m$, $m_i\in M$, \ $b$, $b_i\in\Hom_R(M,R)$, and
$(m,b)\longmapsto\langle m,b\rangle$ denotes the natural
$R$\+bilinear pairing $M\times\Hom_R(M,R)\rarrow R$.
\end{lem}

\begin{proof}
 Passing to retracts, one reduces the question to the case of
a finitely generated free $R$\+module $M$, when the assertion is
provable with explicit free $R$\+module bases (similarly to
the proof of Lemma~\ref{skew-symmetric-tensors-lemma}).
 Alternatively, one can deduce the desired assertion from
Lemma~\ref{skew-symmetric-tensors-lemma} by observing that the natural
isomorphism $\Hom_R(M^{\ot_R\,n},R)\simeq\Hom_R(M,R)^{\ot_R\,n}$
obviously induces a natural isomorphism between the $R$\+module
$\Hom_R(\bigwedge^n_R(M),R)$ and the submodule of all
skew-symmetric tensors in $\Hom_R(M,R)^{\ot_R\,n}$.
\end{proof}

 Let $R\rarrow S$ be a homomorphism of commutative rings and $M$
be an $R$\+module.
 Consider the $S$\+module $N=S\ot_RM$.
 Then, similarly to the discussion of the symmetric algebras in
Section~\ref{lie-algebroid-pbw-theorem-subsecn}, the map of rings
$R\rarrow S$ and the natural $R$\+module map $M\rarrow N$ induce
a homomorphism of graded $R$\+algebras $\bigwedge_R^*(M)\rarrow
\bigwedge_S^*(N)$.
 One can easily see that the induced homomorphism of graded
$S$\+algebras $S\ot_R\bigwedge_R^*(M)\rarrow\bigwedge_S^*(N)$
is an isomorphism.

 Let $\M$ be a quasi-coherent sheaf over a scheme~$X$.
 Then it follows from the previous paragraph that the rule
$\bigwedge_X^*(\M)(U)=\bigwedge_{\cO_X(U)}^*(\M(U))$ for all affine
open subschemes $U\subset X$ defines a quasi-coherent graded algebra
$\bigwedge_X^*(\M)$ over~$X$.
 The quasi-coherent graded algebra $\bigwedge_X^*(\M)$ is called
the \emph{exterior algebra} of a quasi-coherent sheaf $\M$ on~$X$.

\subsection{Chevalley--Eilenberg CDG-ring}
\label{chevalley-eilenberg-cdg-ring-subsecn}
 Let $(R,\g)$ be a Lie algebroid such that $\g$~is a finitely
generated projective $R$\+module.
 The \emph{Chevalley--Eilenberg complex} $C^\bu_R(\g)$ of
the Lie algebroid $(R,\g)$ is constructed as
follows~\cite[Section~III.14]{ChE}, \cite[Section~4]{Rin},
\cite[Section~2.2]{Cal}, \cite[Section~10.9]{Prel}.
 Consider the $R$\+module $\Hom_R(\g,R)$, and put
$C^*_R(\g)=\bigwedge_R^*(\Hom_R(\g,R))$.
 Let us define the differential $d\:C^n_R(\g)\rarrow C^{n+1}_R(\g)$,
for all $n\ge0$.

 The map $d_0\:C^0_R(\g)\rarrow C^1_R(\g)$ is defined by the rule
$$
 \langle v,d_0(r)\rangle = v(r)
 \quad\text{for all $v\in\g$ and $r\in R$}
$$
in the pairing notation $\langle{-},{-}\rangle$ from
Lemma~\ref{exterior-dual-lemma}.
 In the same notation, the map $d_1\:C^1_R(\g)\rarrow C^2_R(\g)$ is
defined by the formula
$$
 \langle v\wedge w,\>d_1(b)\rangle=
 \langle[v,w],b\rangle-v(\langle w,b\rangle)+w(\langle v,b\rangle)
 \quad\text{for all $v$, $w\in\g$ and $b\in\Hom_R(\g,R)$}.
$$
 The maps~$d_0$ and~$d_1$ can be (obviously uniquely) extended to
an odd derivation $d\:C^*_R(\g)\rarrow C^*_R(\g)$ of degree~$1$,
as one can compute.
 Moreover, the differential~$d$ on $C^*_R(\g)$ has zero square,
$d\circ d=0$.
 See~\cite[Proposition~3.16]{Prel} for a much more general assertion
(based on the presumption that the nonhomogeneous quadratic relations
defining the enveloping ring $A_R(\g)$ are self-consistent in
the sense of~\cite[Section~3.3]{Prel}, which still needs
to be checked).

 Notice that the submodule that in the notation
of~\cite[Section~3]{Prel} is denoted by $I\subset\g\ot_R\g$ is
precisely the submodule of all skew-symmetric tensors
in $\g^{\ot_R\,2}$.
 Similary, the submodule $I^{(3)}=I\ot_R\g\cap\g\ot_RI\subset
\g\ot_R\g\ot_R\g$ is precisely the submodule of all skew-symmetric
tensors in $\g^{\ot_R\,3}$.
 By Lemma~\ref{skew-symmetric-tensors-lemma}, these two $R$\+modules
are naturally isomorphic to $\bigwedge^2_R(\g)$ and
$\bigwedge^3_R(\g)$, respectively.

 We have constructed the Chevalley--Eilenberg complex $C^\bu_R(\g)$.
 The complex $C^\bu_R(\g)$ is a DG\+ring whose underlying graded ring
$C^*_R(\g)$ is a graded algebra over the commutative ring $R$, but
the differential~$d$ is \emph{not} $R$\+linear.

 Now let $(R,\g,\widetilde\g)$ be a twisted Lie algebroid such that
$\g$~is a finitely generated projective $R$\+module.
 Let us construct the \emph{Chevalley--Eilenberg CDG\+ring}
$C^\cu_R(\g,\widetilde\g)$.
 We refer to Section~\ref{cdg-rings-cdg-modules-subsecn} for
the background definitions of a CDG\+ring and a morphism of CDG\+rings.

 The underlying graded ring of the CDG\+ring $C^\cu_R(\g,\widetilde\g)$
is the graded $R$\+algebra $C^*_R(\g,\widetilde\g)=C^*_R(\g)$.
 The differential~$d$ on $C^*_R(\g,\widetilde\g)$ is just
the differential~$d $ on $C^*_R(\g)$ defined above; so the differential
on $C^*_R(\g,\widetilde\g)$ does not depend on~$\widetilde\g$ but
only on~$\g$.
 The curvature element $h\in C^2_R(\g,\widetilde\g)$, which depends
on~$\widetilde\g$ and is \emph{not} determined by $R$ and~$\g$, is
constructed as follows.

 Notice that the graded algebra $C^*_R(\g,\widetilde\g)=C^*_R(\g)$ is
graded commutative, so \emph{all} CDG\+ring structures on
$C^*_R(\g,\widetilde\g)$ have differentials with zero squares.
 \emph{Any} choice of a curvature element $h\in C^2_R(\g)$ on top of
a given DG\+ring structure on $C^*_R(\g)$ defines a CDG\+ring
structure on $C^*_R(\g)$, as the commutator with~$h$ is a zero map.
 In the construction below, we use the datum of a twisted Lie
algebroid $(R,\g,\widetilde\g)$ in order to produce a CDG\+ring
structure by adding the curvature element to what is otherwise
a DG\+ring structure depending only on the Lie algebroid $(R,\g)$.

 Since $\g$~is a projective $R$\+module, the short exact sequence of
$R$\+modules $0\rarrow R\overset\iota\rarrow\widetilde\g\overset\pi
\rarrow\g\rarrow0$ splits.
 Choose an arbitrary splitting $\sigma\:\g\rarrow\widetilde\g$; so
$\sigma$~is an $R$\+module map and $\pi\circ\sigma=\id_\g$.
 Notice that \emph{no} compatibility of~$\sigma$ with the Lie algebra
structures on $\g$ and~$\widetilde\g$ is assumed here.

 The element $h\in C^2_R(\g,\widetilde\g)$ is defined by the formula
$$
 \iota(\langle v\wedge w,\>h\rangle)=
 \sigma([v,w])-[\sigma(v),\sigma(w)]
 \quad\text{for all $v$, $w\in\g$}
$$
in the pairing notation $\langle{-},{-}\rangle$ from
Lemma~\ref{exterior-dual-lemma}.
 Here $[v,w]$ is the commutator computed in the Lie algebra~$\g$,
while $[\sigma(v),\sigma(w)]$ is the commutator in the Lie
algebra~$\widetilde\g$ and $\iota\:R\rarrow\widetilde\g$ is
the injective $R$\+module map from the definition of a twisted
Lie algebroid in Section~\ref{twisted-lie-algebroids-subsecn}.

 The CDG\+ring $C^\cu_R(\g,\widetilde\g)$ is defined uniquely up to
a natural isomorphism.
 Given two $R$\+linear splittings $\sigma'\:\g\rarrow\widetilde\g$
and $\sigma''\:\g\rarrow\widetilde\g$ of the surjective map
$\pi\:\widetilde\g\rarrow\g$, the difference $\sigma''-\sigma'$
provides an $R$\+linear map $a\:\g\rarrow R$ such that
$\sigma''(v)-\sigma'(v)=\iota(a(v))$ for all $v\in\g$.
 So we have an element $a\in C^1_R(\g,\widetilde\g)$.
 Denote by $(C^*_R(\g,\widetilde\g),d,h')$ and
$(C^*_R(\g,\widetilde\g),d,h'')$ the two CDG\+ring structures
constructed as above using the two splittings $\sigma'$ and~$\sigma''$.
 Then one has a natural change-of-connection isomorphism of
CDG\+rings $(\id,a)\:(C^*_R(\g,\widetilde\g),d,h'')\rarrow
(C^*_R(\g,\widetilde\g),d,h')$ \,\cite[Proposition~3.20]{Prel}.
 This simply means that $h''=h'+d(a)$ in $C^2(\g,\widetilde\g)$ (as
$a^2=0$ in $C^2(\g,\widetilde\g)$ for all $a\in C^1(\g,\widetilde\g)$).

 In particular, if $\widetilde\g=R\oplus\g$ is a split twisted Lie
algebroid, then we are given a splitting $\sigma\:\g\rarrow
\widetilde\g$ that is both an $R$\+linear map and a homomorphism
of Lie algebras.
 Using such a splitting~$\sigma$ for the construction of the curvature
element~$h$ above, we get $h=0$.
 Hence, in the case a split twisted Lie algebroid $(R,\g,\widetilde\g)$,
the Chevalley--Eilenberg CDG\+ring $C^\cu_R(\g,\widetilde\g)$ coincides
with the Chevalley--Eilenberg DG\+ring $C^\bu_R(\g)$ of the Lie
algebroid~$\g$, that is $C^\cu_R(\g,\widetilde\g)=C^\bu_R(\g)$.

\subsection{Chevalley--Eilenberg quasi-coherent CDG\+quasi-algebra}
\label{chevalley-eilenberg-qcoh-cdg-quasi-algebra-subsecn}
 Let $X$ be a scheme and $(\g,\widetilde\g)$ be a quasi-coherent
twisted Lie algebroid over~$X$.
 We will say that the quasi-coherent twisted Lie algebroid
$(\g,\widetilde\g)$ over $X$ is \emph{finite locally free} if
the quasi-coherent sheaf~$\g$ on $X$ is finite locally free (in
the sense of Section~\ref{loc-free-sheaves-subsecn}).

 Let $(\g,\widetilde\g)$ be a finite locally free quasi-coherent
twisted Lie algebroid over~$X$.
 Let us construct the \emph{Chevalley--Eilenberg quasi-coherent
CDG\+quasi-algebra} $\cC^\cu_X(\g,\widetilde\g)$ over~$X$.
 We refer to Section~\ref{qcoh-cdg-quasi-algebras-subsecn} for
the background definitions of a quasi-coherent CDG\+quasi-algebra
and a morphism of quasi-coherent CDG\+quasi-algebras.

 The underlying quasi-coherent graded quasi-algebra
$\cC^*_X(\g,\widetilde\g)$ of $\cC^\cu_X(\g,\widetilde\g)$ is,
actually, a quasi-coherent graded algebra over~$X$.
 The quasi-coherent graded algebra $\cC^*_X(\g,\widetilde\g)$ only
depends on the finite locally free sheaf~$\g$ on $X$, and does not
depend on~$\widetilde\g$.
 Specifically, we put $\cC^*_X(\g,\widetilde\g)=
\bigwedge_X^*(\cHom_{\cO_X}(\g,\cO_X))$, where the finite locally
free sheaf $\cHom_{\cO_X}(\F,\cO_X)$ was defined in
Section~\ref{loc-free-sheaves-subsecn} and the quasi-coherent
graded algebra $\bigwedge_X^*(\M)$ was defined in
Section~\ref{skew-symmetric-tensors-subsecn}.

 For any affine open subscheme $U\subset X$, the differential~$d$
on the graded ring $\cC^*_X(\g,\widetilde\g)(U)$ is defined without
any arbitrary choices, and also does not depend on~$\widetilde\g$
but only on the quasi-coherent Lie algebroid~$\g$ on~$X$.
 Specifically, we have $\cC^*_X(\g,\widetilde\g)(U)=
\bigwedge^*_{\cO_X(U)}(\Hom_{\cO_X(U)}(\g(U),\cO_X(U)))
=C^*_{\cO_X(U)}(\g(U))$, and the differential~$d$ on
$\cC^*_X(\g,\widetilde\g)(U)$ is simply the differential on
the Chevalley--Eilenberg complex $C_{\cO_X(U)}^\bu(\g(U))$ of
the Lie algebroid $(\cO_X(U),\g(U))$, as constructed in
Section~\ref{chevalley-eilenberg-cdg-ring-subsecn}.

 In order to construct the curvature elements of
the CDG\+quasi-algebra $\cC^\cu_X(\g,\widetilde\g)$, we need to make
an arbitrary choice of a family of spllittings.
 For every affine open subscheme $U\subset X$, the short exact
sequence of $\cO_X(U)$\+modules $0\rarrow\cO_X(U)\overset\iota
\rarrow\widetilde\g(U)\overset\pi\rarrow\g(U)\rarrow0$ splits, since
$\g(U)$ is a (finitely generated) projective $\cO_X(U)$\+module.
 Choose an $\cO_X(U)$\+linear splitting $\sigma_U\:\g(U)\rarrow
\widetilde\g(U)$ of this short exact sequence of $\cO_X(U)$\+modules,
one per every affine open subscheme $U\subset X$.
 Let us emphasize that \emph{neither} any compatibility of
the splittings~$\sigma_U$ with the Lie algebra structures on
$\g(U)$ and $\widetilde\g(U)$, \emph{nor} any compatibility between
$\sigma_U$ and $\sigma_V$ for varying affine open subschemes $U$
and $V\subset X$, is presumed.

 Once the splitting~$\sigma_U$ is chosen, the construction from 
Section~\ref{chevalley-eilenberg-cdg-ring-subsecn} applied to
the twisted Lie algebroid $(\cO_X(U),\g(U),\widetilde\g(U))$ and
the splitting~$\sigma_U$ provides the desired curvature element
$h_U\in\cC^2_X(\g,\widetilde\g)(U)=
C^2_{\cO_X(U)}(\g(U),\widetilde\g(U))$.
 It remains to construct the change-of-connection elements
$a_{VU}\in\cC^1_X(\g,\widetilde\g)(V)$ for all pairs of affine
open subschemes $V\subset U\subset X$.

 We have two $\cO_X(V)$\+linear splittings of the surjective
$\cO_X(V)$\+module map $\pi\:\widetilde\g(V)\rarrow\g(V)$: viz.,
$\cO_X(V)\ot_{\cO_X(U)}\sigma_U\:\g(V)\rarrow\widetilde\g(V)$
and $\sigma_V\:\g(V)\rarrow\widetilde\g(V)$.
 Let $a_{VU}\in\Hom_{\cO_X(V)}(\g(V),\cO_X(V))$ be the map
defined by the rule $\iota(a_{VU}(v))=
(\cO_X(V)\ot_{\cO_X(U)}\sigma_U)(v)-\sigma_V(v)$ for all $v\in\g(V)$.
 Then $a_{VU}\in\cC^1_X(\g,\widetilde\g)(V)=
\Hom_{\cO_X(V)}(\g(V),\cO_X(V))$ is the desired change-of-connection
element.

 The quasi-coherent CDG\+quasi-algebra $\cC^\cu_X(\g,\widetilde\g)$
is defined uniquely up to a natural isomorphism.
 Given two families of $\cO_X(U)$\+linear splittings
$\sigma'_U\:\g(U)\rarrow\widetilde\g(U)$ and $\sigma''_U\:\g(U)\rarrow
\widetilde\g(U)$, a natural change-of-connection isomorphism between
the two related quasi-coherent CDG\+quasi-algebras over $X$ is
constructed similarly to the construction in
Section~\ref{chevalley-eilenberg-cdg-ring-subsecn}.

 Let $\widetilde\g=\cO_X\oplus\g$ be a split (finite locally free)
quasi-coherent twisted Lie algebroid over $X$ with the global
$\cO_X$\+linear splitting $\sigma\:\g\rarrow\widetilde\g$ compatible
with the sheaf of Lie algebra structures on $\g$ and~$\widetilde\g$.
 In this context, one can choose $\sigma_U=\sigma(U)\:\g(U)\rarrow
\widetilde\g(U)$ to be the splitting maps provided by the given
splitting of the whole quasi-coherent twisted Lie algebroid
$(\g,\widetilde\g)$ over~$X$.
 Then the construction above produces a quasi-coherent
DG\+quasi-algebra, which we denote by
$\cC^\bu_X(\g)=\cC^\cu_X(\g,\widetilde\g)$ and call
the \emph{Chevalley--Eilenberg quasi-coherent DG\+quasi-algebra}
of the quasi-coherent Lie algebroid~$\g$.
 This is the usual construction of the Chevalley--Eilenberg
DG\+(quasi)-algebra of a quasi-coherent Lie algebroid over a scheme.
 Notice that we use the term ``DG\+quasi-algebra'' rather than
``DG\+algebra'' in application to $\cC^\bu_X(\g)$ for the only reason
that the differential in $\cC^\bu_X(\g)$ is not $\cO_X$\+linear.

\subsection{Chevalley--Eilenberg CDG-modules}
\label{chevalley-eilenberg-cdg-modules-subsecn}
 Let $(R,\g,\widetilde\g)$ be a Lie algebroid such that $\g$~is
a finitely generated projective $R$\+module.
 The aim of this section is to construct two CDG\+modules over
the CDG\+ring $C^\cu_R(\g,\widetilde\g)$, a left and a right one.
 In the nontwisted (or split twisted) case, this construction goes
back to~\cite[Section~4]{Rin}.

 The left CDG\+module, which is actually a CDG\+bimodule over
the CDG\+ring $C^\cu_R(\g,\widetilde\g)$ and the enveloping ring
$A_R(\g,\widetilde\g)$ viewed as a CDG\+ring, is denoted by
$C^\cu_R(\g,\widetilde\g,A_R)$.
 (See Section~\ref{cdg-rings-cdg-modules-subsecn} for
the definition of a CDG\+bimodule.)
 The underlying graded bimodule $C^*_R(\g,\widetilde\g,A_R)$
of $C^\cu_R(\g,\widetilde\g,A_R)$ is the graded bimodule
$C^*_R(\g,\widetilde\g,A_R)=C^*_R(\g,\widetilde\g)\ot_R
A_R(\g,\widetilde\g)$ over the graded ring $C^*_R(\g,\widetilde\g)=
C^*_R(\g)$ and the ring $A_R(\g,\widetilde\g)$ viewed as a graded ring
concentrated in degree~$0$.
 The graded bimodule $C^*_R(\g,\widetilde\g,A_R)$ is concentrated
in a finite interval of nonnegative cohomological degrees.

 The right CDG\+module, which is actually a CDG\+bimodule over
the enveloping ring $A_R(\g,\widetilde\g)$ viewed as a CDG\+ring
and the CDG\+ring $C^\cu_R(\g,\widetilde\g)$, is denoted by
$C_\cu^R(A_R,\g,\widetilde\g)$.
 The underlying graded bimodule $C_*^R(A_R,\g,\widetilde\g)$
of $C_\cu^R(A_R,\g,\widetilde\g)$ is the graded bimodule
$C_*^R(A_R,\g,\widetilde\g)=A_R(\g,\widetilde\g)\ot_R
\bigwedge_R^*(\g)$ over the ring $A_R(\g,\widetilde\g)$ and
the graded ring $C^*_R(\g,\widetilde\g)=C^*_R(\g)$.
 Here the graded right $C^*_R(\g)$\+module structure on the graded
$R$\+module $\bigwedge_R^*(\g)$ comes from the free graded left
$C^*_R(\g)$\+module structure on the graded $R$\+module
$C^*_R(\g)=\bigwedge^*_R(\Hom_R(\g,R))$ via the isomorphism
$\bigwedge_R^*(\g)=\Hom_R(\bigwedge^*_R(\Hom_R(\g,R)),R)$
provided by Lemma~\ref{exterior-dual-lemma}.
 The graded bimodule $C_*^R(A_R,\g,\widetilde\g)$ is concentrated
in a finite interval of nonnegative homological degrees or in
a finite interval of nonpositive cohomological degrees.

 We start with constructing the left CDG\+module
$C^\cu_R(\g,\widetilde\g,A_R)$.
 The differential~$d_\sigma$ on $C^*_R(\g,\widetilde\g,A_R)$ depends
on the choice of an $R$\+linear splitting $\sigma\:\g\rarrow
\widetilde\g$, as in Section~\ref{chevalley-eilenberg-cdg-ring-subsecn}.
 Let $e_\sigma\in\Hom_R(\g,R)\ot_R\widetilde\g$ be the element
corresponding to the $R$\+linear map~$\sigma$ under the natural
$R$\+module isomorphism $\Hom_R(\g,\widetilde\g)\simeq\Hom_R(\g,R)
\ot_R\widetilde\g$ (which holds since $\g$~is a finitely generated
projective $R$\+module).
 Recall that $\Hom_R(\g,R)=C^1_R(\g)=C^1_R(\g,\widetilde\g)$ by
the definition of $C^*_R(\g,\widetilde\g)$, and the map
$\lambda\:\widetilde\g\rarrow A_R(\g,\widetilde\g)$ is a left
$R$\+module map by relation~(c) from
Section~\ref{enveloping-algebra-subsecn}.
 Denote by $l_\sigma\in C^1_R(\g,\widetilde\g)\ot_RA_R(\g,\widetilde\g)$
the image of~$e_\sigma$ under the map $\id\ot_R\lambda\:
\Hom_R(\g,R)\ot_R\widetilde\g\rarrow C^1_R(\g,\widetilde\g)\ot_R
A_R(\g,\widetilde\g)$.
 Then we put
$$
 d_\sigma(c\ot a)=d(c)\ot a+(-1)^{|c|}cl_\sigma a
 \quad\text{for all $c\in C^{|c|}_R(\g,\widetilde\g)$ and
 $a\in A_R(\g,\widetilde\g)$}.
$$

 Clearly, the differential~$d_\sigma$ is a right
$A_R(\g,\widetilde\g)$\+linear map.
 As a special case of~\cite[Proposition~6.4]{Prel}, one can obtain
the assertion that the differential~$d_\sigma$ defines a structure
of left CDG\+module over $C^\cu_R(\g,\widetilde\g,A_R)$ on the graded
left module $C^*_R(\g,\widetilde\g,A_R)$ over $C^*_R(\g,\widetilde\g)$.
 Finally, according to~\cite[Lemma~6.5]{Prel}, the resulting
CDG\+module $C^\cu_R(\g,\widetilde\g,A_R)$ over the CDG\+ring
$C^\cu_R(\g,\widetilde\g)$ does not depend on the choice of
the splitting $\sigma\:\g\rarrow\widetilde\g$.
 The latter assertions means that, for any two splittings $\sigma'$
and~$\sigma''$, the differentials $d_{\sigma'}$ and~$d_{\sigma''}$ are
transformed into each other, as per formula~(vi) from
Section~\ref{cdg-rings-cdg-modules-subsecn}, by the change of scalars
with respect to the natural change-of-connection isomorphism of
CDG\+rings constructed in
Section~\ref{chevalley-eilenberg-cdg-ring-subsecn}.

 The CDG\+bimodule $C_\cu^R(A_R,\g,\widetilde\g)$ over
$A_R(\g,\widetilde\g)$ and $C^\cu_R(\g,\widetilde\g)$ can be now
produced as the Hom CDG\+bimodule
$$
 C_\cu^R(A_R,\g,\widetilde\g)=\Hom^\cu_{A_R(\g,\widetilde\g)^\rop}
 (C^\cu_R(\g,\widetilde\g,A_R),A_R(\g,\widetilde\g)).
$$
 See the end of Section~\ref{cdg-rings-cdg-modules-subsecn} for
the constructions of the Hom CDG\+bimodules.
 Cf.~\cite[Remark~6.6]{Prel}.

 In the case of a split twisted Lie algebroid $\widetilde\g=R\oplus\g$,
the Chevalley--Eilenberg CDG\+ring is a DG\+ring,
$C^\cu_R(\g,\widetilde\g)=C^\bu(\g)$, as mentioned at the end of
Section~\ref{chevalley-eilenberg-cdg-ring-subsecn}.
 In this case, the constructions above produce two DG\+modules,
a left DG\+module $C^\bu_R(\g,A_R)=C^\cu_R(\g,\widetilde\g,A_R)$
and a right DG\+module $C_\bu^R(A_R,\g)=C_\cu^R(A_R,\g,\widetilde\g)$,
over the DG\+ring $C^\bu(\g)$.
 More precisely, $C^\bu_R(\g,A_R)$ is a DG\+bimodule over
$C^\bu_R(\g)$ and $A_R$, while $C_\bu^R(A_R,\g)$ is a DG\+bimodule
over $A_R$ and $C^\bu_R(\g)$.

 Now let $X$ be a scheme and $(\g,\widetilde\g)$ be a finite locally
free quasi-coherent twisted Lie algebroid over~$X$.
 Denote by $\cC^*_X(\g,\widetilde\g,\cA_X)$ the quasi-coherent graded
bimodule $\cC^*_X(\g,\widetilde\g,\cA_X)\allowbreak=
\cC^*_X(\g,\widetilde\g)\ot_{\cO_X}\cA_X(\g,\widetilde\g)$ over
the quasi-coherent graded algebra $\cC^*_X(\g,\widetilde\g)=\cC^*_X(\g)$
and the quasi-coherent quasi-algebra $\cA_X(\g,\widetilde\g)$ over~$X$.
 The construction of the differential~$d_\sigma$ above in this section
produces a quasi-coherent left CDG\+module
$\cC^\cu_X(\g,\widetilde\g,\cA_X)$ over the quasi-coherent
CDG\+quasi-algebra $\cC^\cu_X(\g,\widetilde\g)$ over~$X$,
with the underlying quasi-coherent graded module
$\cC^*_X(\g,\widetilde\g,\cA_X)$.
 The differentials on the CDG\+modules
$\cC^\cu_X(\g,\widetilde\g,\cA_X)(U)=
C^\cu_{\cO_X(U)}(\g(U),\widetilde\g(U),A_{\cO_X(U)})$ are right
$\cA_X(\g,\widetilde\g)(U)$\+linear for all affine open subschemes
$U\subset X$.
 In other words, this means that $\cC^\cu_X(\g,\widetilde\g,\cA_X)$ is
a quasi-coherent CDG\+bimodule over $\cC^\cu_X(\g,\widetilde\g)$ and
$\cA_X(\g,\widetilde\g)$ in the sense of the definition in
Section~\ref{fHom-and-contratensor-of-cdg-modules} below.

 Similarly, denote by $\cC_*^X(\cA_X,\g,\widetilde\g)$
the quasi-coherent graded bimodule
$\cC_*^X(\cA_X,\g,\widetilde\g)\allowbreak=
\cA_X(\g,\widetilde\g)\ot_{\cO_X}\bigwedge_X^*(\g)$ over
the quasi-coherent quasi-algebra $\cA_X(\g,\widetilde\g)$ and
the quasi-coherent graded algebra $\cC^*_X(\g,\widetilde\g)=
\cC^*_X(\g)$ over~$X$.
 Here the right module structures over the graded rings
$\cC^*_X(\g)(U)=C^*_{\cO_X(U)}(\g(U))$ on the graded $\cO_X(U)$\+modules
$\bigwedge_X^*(\g)(U)=\bigwedge_{\cO_X(U)}^*(\g(U))$ for all affine
open subschemes $U\subset X$ are constructed as explained above
in this section.
 The construction of the differential on the Hom CDG\+module as above
and in end of Section~\ref{cdg-rings-cdg-modules-subsecn} produces
a quasi-coherent right CDG\+module
$\cC_\cu^X(\cA_X,\g,\widetilde\g)$ over the quasi-coherent
CDG\+quasi-algebra $\cC^\cu_X(\g,\widetilde\g)$ over~$X$,
with the underlying quasi-coherent graded module
$\cC_*^X(\cA_X,\g,\widetilde\g)$.
 The differentials on the CDG\+modules
$\cC_\cu^X(\cA_X,\g,\widetilde\g)(U)=
C_\cu^{\cO_X(U)}(A_{\cO_X(U)},\g(U),\widetilde\g(U))$ are left
$\cA_X(\g,\widetilde\g)(U)$\+linear for all affine open subschemes
$U\subset X$.
 In other words, this means that $\cC_\cu^X(\cA_X,\g,\widetilde\g)$ is
a quasi-coherent CDG\+bimodule over $\cA_X(\g,\widetilde\g)$ and
$\cC^\cu_X(\g,\widetilde\g)$ in the sense of the definition in
Section~\ref{fHom-and-contratensor-of-cdg-modules} below.

 In the case of a split quasi-coherent twisted Lie algebroid
$\widetilde\g=\cO_X\oplus\g$, the Chevalley--Eilenberg quasi-coherent
CDG\+quasi-algebra is a quasi-coherent DG\+quasi-algebra,
$\cC^\cu_X(\g,\widetilde\g)=\cC^\bu_X(\g)$, as explained at the end of
Section~\ref{chevalley-eilenberg-qcoh-cdg-quasi-algebra-subsecn}.
 In this case, the constructions above produce two quasi-coherent
DG\+modules, a quasi-coherent left DG\+module
$\cC^\bu_X(\g,\cA_X)=\cC^\cu_X(\g,\widetilde\g,\cA_X)$ and
a quasi-coherent right DG\+module
$\cC_\bu^X(\cA_X,\g)=\cC_\cu^X(\cA_X,\g,\widetilde\g)$, over
the quasi-coherent DG\+quasi-algebra $\cC^\bu_X(\g)$.

\subsection{Opposite CDG-rings and opposite twisted Lie algebroids}
\label{opposite-cdg-rings-and-twisted-lie-subsecn}
 Let $B^\cu=(B^*,d,h)$ be a CDG\+ring (as defined in
Section~\ref{cdg-rings-cdg-modules-subsecn}).
 The \emph{opposite CDG\+ring} $B^\rop{}^\cu=
(B^\rop{}^*,d^\rop,-h^\rop)$
is defined by the following rules~\cite[Remark~9.10]{Prel}:
\begin{itemize}
\item the graded ring $B^\rop{}^*$ is the opposite graded ring
to~$B^*$; this means that the elements of $B^\rop{}^*$ are
the symbols $b^\rop$, where $b\in B$, the grading on $B^\rop{}^*$
agrees with the grading on $B^*$, and the multiplication
in $B^\rop{}^*$ is given by the formula
$b^\rop c^\rop=(-1)^{|b||c|}(cb)^\rop$ for all $b\in B^{|b|}$
and $c\in B^{|c|}$;
\item the differential~$d^\rop$ on $B^\rop{}^*$ is equal to
the differential~$d$ on~$B$; more precisely,
$d^\rop(b^\rop)=(d(b))^\rop$ for all $b\in B$;
\item the curvature element in $B^\rop{}^2$ only differs from
the curvature element in $B^2$ by the minus sign; so the curvature
element of $B^\rop{}^\cu$ is the element $-h^\rop\in B^\rop{}^2$.
\end{itemize}

 Given a morphism of CDG\+rings $(f,a)\:B^\cu\rarrow A^\cu$,
the related morphism of the opposite CDG\+rings is
$(f^\rop,-a^\rop)\:B^\rop{}^\cu\rarrow A^\rop{}^\cu$.
 Here the pair $(f^\rop,-a^\rop)$ is constructed by the rules:
\begin{itemize}
\item the morphism of graded rings $f^\rop\:B^\rop{}^*\rarrow
A^\rop{}^*$ is the same map as the morphism of graded rings
$f\:B^*\rarrow A^*$; so $f^\rop(b^\rop)=f(b)^\rop$ for all
elements $b\in B$;
\item the change-of-connection element $-a^\rop\in A^\rop{}^1$
in the morphism of CDG\+rings $(f^\rop,-a^\rop)$ only differs from
the change-of-connection element $a\in A^1$ in the morphism of
CDG\+rings $(f,a)$ by the minus sign.
\end{itemize}

 For any CDG\+ring $B^\cu$, there are natural equivalences (in fact,
isomorphisms) of DG\+categories $B^\cu\bModl\simeq\bModr B^\rop{}^\cu$
and $\bModr B^\cu\simeq B^\rop{}^\cu\bModl$.
 Specifically, given a left CDG\+module $M^\cu=(M^*,d_M)$ over $B^\cu$,
the related right CDG\+module $M^\rop{}^\cu=(M^\rop{}^*,d_M^\rop)$
over $B^\rop{}^\cu$ is defined by the rules
\begin{itemize}
\item the elements of $M^\rop{}^*$ are the symbols $x^\rop$, where
$x\in M^*$, the grading on $M^\rop{}^*$ agrees with the grading on
$M^*$, and the action of $B^\rop{}^*$ on $M^\rop{}^*$ is given
by the formula $x^\rop b^\rop=(-1)^{|x||b|}(bx)^\rop$ for all
$x\in M^{|x|}$ and $b\in B^{|b|}$;
\item the differential~$d_M^\rop$ on $M^\rop{}^*$ is equal to
the differential~$d_M$ on~$M$; more precisely, $d_M^\rop(x^\rop)=
(d_M(x))^\rop$ for all $x\in M$.
\end{itemize}
 The rules for producing a left CDG\+module $N^\rop{}^\cu$ over
$B^\rop{}^\cu$ from a right CDG\+module $N^\cu$ over $B^\cu$
are similar.

 Let $R$ be a commutative ring and $B$ be a quasi-module over~$R$.
 Then the opposite quasi-module $B^\rop$ over $R$ is defined by
the rules $B^\rop=B$ and $r(b^\rop)s=(sbr)^\rop$ for all $b\in B$
and $r$, $s\in R$.
 So the left action of $R$ on $B$ corresponds to the right action of
$R$ on $B^\rop$ and vice versa.
 Similarly one defines the opposite quasi-coherent quasi-module
$\cB^\rop$ to a quasi-coherent quasi-module $\cB$ over a scheme~$X$:
the passage from $\cB$ to $\cB^\rop$ switches the left and right
actions of~$\cO_X$.

 Let $R$ be a commutative ring and $A$ be a quasi-algebra over~$R$.
 Then the opposite ring $A^\rop$ to $A$ is also naturally
a quasi-algebra over~$R$.
 The passage from $A$ to $A^\rop$ switches the left and right actions
of~$R$; so the left action of $R$ on $A$ corresponds to the right action
of $R$ on $A^\rop$ and vice versa.

 Let $X$ be a scheme and $\cA$ be a quasi-coherent quasi-algebra
over~$X$.
 Then the rule $\cA^\rop(U)=\cA(U)^\rop$ for all affine open subschemes
$U\subset X$ defines a quasi-coherent quasi-algebra $\cA^\rop$ over~$X$.
 The left action of the structure sheaf $\cO_X$ on $\cA$ corresponds to
the right action of $\cO_X$ on $\cA^\rop$ and vice versa.
 There are natural equivalences (in fact, isomorphisms) of the abelian
categories of quasi-coherent modules $\cA\Qcoh\simeq\Qcohr\cA^\rop$ and
$\Qcohr\cA\simeq\cA^\rop\Qcoh$.

 Now let $\cB^\cu$ be a quasi-coherent CDG\+quasi-algebra over~$X$.
 Applying the rules listed above, one can produce the opposite
quasi-coherent CDG\+quasi-algebra $\cB^\rop{}^\cu$ over $X$ together
with equivalences (in fact, isomorphisms) of the DG\+categories of
quasi-coherent CDG\+modules $\cB^\cu\bQcoh\simeq\bQcohr\cB^\rop{}^\cu$
and $\bQcohr\cB^\cu\simeq\cB^\rop{}^\cu\bQcoh$.

 Finally, let $(\g,\widetilde\g)$ be a quasi-coherent twisted Lie
algebroid over~$X$ (as defined in
Section~\ref{twisted-lie-algebroids-subsecn}).
 Define the \emph{opposite} quasi-coherent twisted Lie algebroid
$(\g,\widetilde\g^\circ)$ over $X$ by the following rules.
 All the pieces of data consituting the quasi-coherent Lie algebroid
$(\g,\widetilde\g)$ remain unchanged by the passage to
$(\g,\widetilde\g^\circ)$, with the only exception of the morphism
of quasi-coherent sheaves $\iota\:\cO_X\rarrow\widetilde\g$,
which changes the sign.
 So we have $\widetilde\g^\circ=\widetilde\g$ and $\iota^\circ=-\iota$.

 Assume that $\g$~is a finite locally free sheaf on $X$, and let
$\cC^\cu_X(\g,\widetilde\g)$ be the Chevalley--Eilenberg quasi-coherent
CDG\+quasi-algebra of $(\g,\widetilde\g)$ (as constructed in
Section~\ref{chevalley-eilenberg-qcoh-cdg-quasi-algebra-subsecn}).
 Recall that the construction of the Chevalley--Eilenberg
quasi-coherent CDG\+quasi-algebra used an arbitrary choice of a family
of $\cO_X(U)$\+linear splittings $\sigma_U\:\g(U)\rarrow
\widetilde\g(U)$ for all affine open subschemes $U\subset X$.
 Choose the same $\cO_X(U)$\+linear splittings $\sigma^\circ_U\:\g(U)
\rarrow\widetilde\g^\circ(U)$ for the opposite quasi-coherent twisted
Lie algebroid $(\g,\widetilde\g^\circ)$; so $\sigma^\circ_U=\sigma_U$
for all affine open subschemes $U\subset X$.
 Then the identity endomorphism of the exterior algebra
$\bigwedge_X^*(\cHom_{\cO_X}(\g,\cO_X))$ is a strict isomorphism
of quasi-coherent CDG\+quasi-algebras $\cC^\cu_X(\g,\widetilde\g)^\rop
\simeq\cC^\cu_X(\g,\widetilde\g^\circ)$ over~$X$.
 Here $\cC^\cu_X(\g,\widetilde\g)^\rop$ denotes the opposite
quasi-coherent CDG\+quasi-algebra to the Chevalley--Eilenberg
quasi-coherent CDG\+quasi-algebra $\cC^\cu_X(\g,\widetilde\g)$
over~$X$.

 For a split quasi-coherent twisted Lie algebroid $\widetilde\g=
\cO_X\oplus\g$, the map $(-\id_{\cO_X},\id_\g)\:\widetilde\g
\rarrow\widetilde\g$ provides a natural isomorphism of quasi-coherent
twisted Lie algebroids $(\g,\widetilde\g^\circ)\simeq(\g,\widetilde\g)$.
 The quasi-coherent DG\+quasi-algebra $\cC^\bu_X(\g)$ is
graded commutative, so the identity endomorphism is an isomorphism
of quasi-coherent DG\+quasi-algebras $\cC^\bu_X(\g)^\rop\simeq
\cC^\bu_X(\g)$ (cf.~\cite[Remark~9.10]{Prel}). {\hbadness=1325\par}

 Let us emphasize, however, that the quasi-coherent quasi-algebra
$\cA_X(\g,\widetilde\g^\circ)$ is usually \emph{not} isomorphic to
the opposite quasi-coherent quasi-algebra $\cA_X(\g,\widetilde\g)^\rop$
to the quasi-coherent quasi-algebra $\cA_X(\g,\widetilde\g)$.
 In particular, in the case of a split quasi-coherent twisted Lie
algebroid $\widetilde\g=\cO_X\oplus\g$, the quasi-coherent
quasi-algebra $\cA_X(\g)$ is usually \emph{not} isomorphic to
its opposite quasi-coherent quasi-algebra $\cA_X(\g)^\rop$.
 We will continue this discussion in
Section~\ref{dual-koszul-subsecn}.

\subsection{Vector fields} \label{vector-fields-subsecn}
 Let $\varkappa\:K\rarrow R$ be a homomorphism of commutative rings.
 The $R$\+module of \emph{K\"ahler differentials} $\Omega_{R/K}$
is spanned by the symbols $d(r)$, \,$r\in R$, subject to
the relations~\cite[Section Tag~00RM]{SP}
\begin{itemize}
\item $d(r+s)=d(r)+d(s)$ for all $r$, $s\in R$;
\item $d(rs)=rd(s)+sd(r)$ for all $r$, $s\in R$;
\item $d(\varkappa(k))=0$ for all $k\in K$.
\end{itemize}

 Given an $R$\+module $M$, one can say that a $K$\+linear map
$D\:R\rarrow M$ is an \emph{$R$\+derivation} if $D(rs)=rD(s)+sD(r)$
in $M$ for all $r$, $s\in R$.
 For any $K$\+linear $R$\+derivation $D\:R\rarrow M$ and any element
$t\in R$, the map $tD\:R\rarrow M$ defined by the rule $(tD)(r)=t(D(r))$
for all $r\in R$ is also a $K$\+linear $R$\+derivation.
 So $K$\+linear $R$\+derivations $R\rarrow M$ form an $R$\+module.
 One can immediately see from the definitions above that
the $R$\+module of $K$\+linear $R$\+derivations $R\rarrow M$ is
naturally isomorphic to the $R$\+module $\Hom_R(\Omega_{R/K},M)$;
this is a universal property defining the $R$\+module of K\"ahler
differentials~$\Omega_{R/K}$.

\begin{lem} \label{kaehler-differentials-ring-epimorphisms}
\textup{(a)} Let $K\rarrow L\rarrow S$ be homomorphisms of
commutative rings such that $K\rarrow L$ is a ring epimorphism.
 Then the natural $S$\+module map\/ $\Omega_{S/K}\rarrow\Omega_{S/L}$
is an isomorphism. \par
\textup{(b)} Let $K\rarrow R\rarrow S$ be homomorphisms of
commutative rings such that $R\rarrow S$ is a flat ring epimorphism.
 Then the natural $S$\+module map $S\ot_R\Omega_{R/K}\rarrow
\Omega_{S/K}$ is an isomorphism.
\end{lem}

\begin{proof}
 Part~(a): by~\cite[Lemma Tag~00RS]{SP}, for any commutative ring
homomorphisms $K\rarrow L\rarrow S$ there is a natural exact sequence
of $S$\+modules $S\ot_L\Omega_{L/K}\rarrow\Omega_{S/K}\rarrow
\Omega_{S/L}\rarrow0$.
 It remains to show that $\Omega_{L/K}=0$ for any commutative ring
epimorphism $K\rarrow L$.
 Indeed, let $M$ be an $L$\+module and $D\:L\rarrow M$ be
a $K$\+linear $L$\+derivation.
 Since $K\rarrow L$ is a ring epimorphism, the forgetful functor
$L\Modl\rarrow K\Modl$ is fully faithful, i.~e., any $K$\+linear
map acting between $L$\+modules is an $L$\+linear map.
 So $D$ is an $L$\+linear map; and since $D(1)=0$ for the unit element
$1\in L$, we can conclude that $D=0$.
 Thus $\Hom_L(\Omega_{L/K},M)=0$.
 As this holds for all $L$\+modules $M$, we have $\Omega_{L/K}=0$.

 Part~(b): in view of the proof of part~(a), we have $\Omega_{S/R}=0$,
hence the map $h\:S\ot_R\Omega_{R/K}\rarrow\Omega_{S/K}$ is surjective.
 In order to prove that $h$~is injective, suppose given an $S$\+module
$N$ and an $S$\+module map $g\:S\ot_R\Omega_{R/K}\rarrow N$.
 Equivalently, this means the datum of an $R$\+module map
$\Omega_{R/K}\rarrow N$; so we have the corresponding $K$\+linear
$R$\+derivation $D_R\:R\rarrow N$.
 The map $D_R$ is a $K$\+linear strongly $R$\+differential operator
of order~$\le1$ in the sense of~\cite[Section~3.1]{Ptd}.
 By~\cite[Proposition~6.6]{Ptd}, it follows that the map $D_R$ can
be uniquely extended to a $K$\+linear strongly $S$\+differential
operator $D_S\:S\rarrow S\ot_RN=N$ of order~$\le1$.
 Using the uniqueness assertion of~\cite[Proposition~6.6]{Ptd}, one
easily proves that the equation $D_S(rs)=rD_S(s)+sD_S(r)$ holds for
all $r\in R$ and $s\in S$ (since it holds for all $r$, $s\in R$),
and then that it holds for all $r$, $s\in S$.
 So $D_S$ is a $K$\+linear $S$\+derivation, and we have
the corresponding $S$\+module map $\Omega_{S/K}\rarrow N$.
 It follows that the map~$g$ factorizes through the surjective map~$h$.
 As this holds for all $N$ and~$g$, we can conclude that $h$~is
an isomorphism.
\end{proof}

 Given a homomorphism of commutative rings $K\rarrow R$, denote by
$\vect_{R/K}=\Hom_R(\Omega_{R/K},R)$ the $R$\+module of
$K$\+linear $R$\+derivations $R\rarrow R$.
 Such derivations are called \emph{vector fields} on $R$ over~$K$.
 The \emph{commutator of vector fields} is defined by the obvious
formula $[v,w](r)=v(w(r))-w(v(r))\in R$ for all $r\in R$ and
$v$, $w\in\vect_{R/K}$.
 The Lie algebra structure on $\vect_{R/K}$ given by the commutator
of vector fields, together with the $R$\+module structure on
$\vect_{R/K}$ defined above and with the obvious action of
the Lie algebra $\vect_{R/K}$ on the ring $R$ form a natural structure
of \emph{Lie algebroid} (in the sense of
Section~\ref{lie-algebroids-subsecn}) on the pair $(R,\vect_{R/K})$.
 We will call the Lie algebroid $(R,\vect_{R/K})$ the \emph{Lie
algebroid of vector fields on $R$ over~$K$}.

 Now let $\tau\:X\rarrow T$ be a morphism of schemes.
 The \emph{quasi-coherent sheaf of K\"ahler differentials}
$\Omega_{X/T}$ is defined as follows.
 For any pair of affine open subschemes $V\subset X$ and $U\subset T$
such that $\tau(V)\subset U$, the $\cO_X(U)$\+module
$\Omega_{X/T}(V)$ is $\Omega_{X/T}(V)=\Omega_{\cO_X(V)/\cO_T(U)}$.
 If $U'\subset U$ is another affine open subscheme such that
$\tau(V)\subset U'$, then there is a natural isomorphism
$\Omega_{\cO_X(V)/\cO_T(U)}\simeq\Omega_{\cO_X(V)/\cO_T(U')}$
by Lemma~\ref{kaehler-differentials-ring-epimorphisms}(a); so
the $\cO_X(V)$\+module $\Omega_{X/T}(V)$ is well-defined.
 (To be more precise, some additional care is required in the case of
a non-semi-separated scheme~$T$; see~\cite[Section~3.3]{Pcosh} for
a similar kind of argument in the context of inverse images of
contraherent cosheaves, mentioned in
Section~\ref{inverse-images-of-O-co-sheaves-subsecn} above.
 A version of this argument for quasi-coherent sheaves can be also
found in Section~\ref{inverse-images-of-A-co-sheaves-subsecn}.)
 For any affine open subscheme $V'\subset V$, we have a natural
isomorphism $\cO_X(V')\ot_{\cO_X(V)}\Omega_{X/T}(V)\simeq
\Omega_{X/T}(V')$ by
Lemma~\ref{kaehler-differentials-ring-epimorphisms}(b).
 In view of Corollary~\ref{enochs-estrada}, the quasi-coherent sheaf
$\Omega_{X/T}$ on $X$ is well-defined by this construction.

 The sheaf of vector fields $\vect_{X/T}$ can be defined as the sheaf
of $\cO_X$\+modules $\Hom_{\cO_X}(\Omega_{X/T},\cO_X)$.
 By the definitions, this means that, for every pair of affine open
subschemes $V\subset X$ and $U\subset T$ such that $\tau(V)\subset U$,
one has $\vect_{X/T}(V)=\vect_{\cO_X(V)/\cO_T(U)}$.
 So $\vect_{X/T}$ is a sheaf of $\cO_X$\+modules and a sheaf of Lie
algebras on $X$ acting naturally on the sheaf of commutative
rings~$\cO_X$; but $\vect_{X/T}$ \emph{need not} be a quasi-coherent
sheaf on $X$ is general.

 The definition of a \emph{morphism of schemes locally of finite
presentation} was given in
Example~\ref{fiberwise-differential-operators-examples}(4).
 By~\cite[Corollaire~IV.16.4.22]{EGA4} or~\cite[Lemma Tag~00RY]{SP},
the $R$\+module of K\"ahler differentials $\Omega_{R/K}$ is finitely
presented whenever $R$ is a finitely presented $K$\+algebra.
 Therefore, the quasi-coherent sheaf of K\"ahler differentials
$\Omega_{X/T}$ is locally finitely presented (in the sense of
Section~\ref{loc-free-sheaves-subsecn}) whenever $\tau\:X\rarrow T$
is a morphism locally of finite presentation.
 According to the discussion in
Section~\ref{loc-free-sheaves-subsecn}, it follows that
$\vect_{X/T}$ is a quasi-coherent sheaf on $X$ whenever
$\tau\:X\rarrow T$ is a morphism locally of finite presentation.

 More generally, $\vect_{X/T}$ is a quasi-coherent sheaf on $X$ whenever
$\Omega_{X/T}$ is a locally finitely presented quasi-coherent sheaf.
 If this is the case, then the quasi-coherent sheaf $\vect_{X/T}$ has
a natural structure of a quasi-coherent Lie algebroid over $X$, as one
can see from the discussion above.

\subsection{Crystalline differential operators}
\label{crystalline-diffoperators-subsecn}
 An introductory discussion of crystalline differential operators vs.\
differential in the sense of Grothendieck can be found in
the book~\cite[Section~10.3]{Prel}; see in
particular~\cite[Examples~10.9]{Prel}.

 Let $\varkappa\:K\rarrow R$ be a homomorphism of commutative rings.
 We refer to~\cite[Section Tag~00T1]{SP} for a discussion of
smooth morphisms of commutative rings.
 Following~\cite[Section~10.8]{Prel}, we will say that the ring
homomorphism~$\varkappa$ is \emph{weakly smooth} if the $R$\+module
of K\"ahler differentials $\Omega_{R/K}$ is finitely generated
and projective.
 According to~\cite[Section Tag~00T1]{SP}, every smooth morphism is
weakly smooth.
 We refer to~\cite[Remark~10.53]{Prel} for examples of weakly smooth
morphisms that are not smooth.

 Following the discussion in Section~\ref{vector-fields-subsecn},
the $R$\+module of vector fields $\vect_{R/K}$ is computed as
$\vect_{R/K}=\Hom_R(\Omega_{R/K},R)$.
 So the $R$\+module $\vect_{R/K}$ is finitely generated and projective
whenever the ring homomorphism~$\varkappa$ is weakly smooth.

 Assuming that the ring homomorphism~$\varkappa$ is weakly smooth, we
define the ring $\cD^\cry_{R/K}$ of \emph{crystalline differential
operators on $R$ over~$K$} as the universal enveloping ring of
the Lie algebroid $(R,\vect_{R/K})$, in the sense of
Section~\ref{enveloping-algebra-subsecn}; that is
$\cD^\cry_{R/K}=A_R(\vect_{R/K})$.
 Following the discussion in Section~\ref{enveloping-algebra-subsecn},
the ring $\cD^\cry_{R/K}$ is a quasi-algebra over~$R$.
 We also define the \emph{de~Rham DG\+ring} $\Omega^\bu_{R/K}$
as the Chevalley--Eilenberg complex (DG\+ring) of the Lie algebroid
$(R,\vect_{R/K})$, in the sense of
Section~\ref{chevalley-eilenberg-cdg-ring-subsecn}; that is
$\Omega^\bu_{R/K}=C^\bu_R(\vect_{R/K})$.
 This definition agrees with the one in~\cite[Section~0FKF]{SP}
(where the setting is much more general).

 By the definition, we have $\Omega^0_{R/K}=R$.
 Due to the assumption that the $R$\+module $\Omega_{R/K}$ is finitely
generated and projective, we have $\Hom_R(\vect_{R/K},R)=\Omega_{R/K}$.
 Hence the degree~$1$ component of the de~Rham DG\+ring
$\Omega^\bu_{R/K}$ is the $R$\+module of K\"ahler differentials,
$\Omega^1_{R/K}=\Omega_{R/K}$, just as the notation suggests.
 Under this identification, the component $d_0\:\Omega^0_{R/K}
\rarrow\Omega^1_{R/K}$ of the differential on $\Omega^\bu_{R/K}$
corresponds to the structure map $d\:R\rarrow\Omega_{R/K}$ from
the construction of the $R$\+module of K\"ahler
differentials~$\Omega_{R/K}$.
 
 The constructions from
Section~\ref{chevalley-eilenberg-cdg-modules-subsecn} also provide
two \emph{de~Rham DG\+modules} over $\Omega^\bu_{R/K}$, the left
DG\+module $C^\bu_R(\vect_{R/K},A_R)=
\Omega^\bu_{R/K}\ot_R\cD^\cry_{R/K}$ and the right DG\+module
$C_\bu^R(A_R,\vect_{R/K})=
\cD^\cry_{R/K}\ot_R\Hom_R(\Omega^\bu_{R/K},R)$.
 Let us emphasize that the subformula $\Hom_R(\Omega^\bu_{R/K},R)$ of
the latter formula does not denote any complex, as the differential on
$\Omega^\bu_{R/K}$ is \emph{not} $R$\+linear and the functor
$\Hom_R({-},R)$ \emph{cannot} be applied to it.
 There is \emph{no} natural differential on the graded $R$\+module
of polyvector fields $\Hom_R(\Omega^*_{R/K},R)\simeq
\bigwedge_R^*(\vect_{R/K})$.
 It is only the tensor product
$\cD^\cry_{R/K}\ot_R\Hom_R(\Omega^\bu_{R/K},R)$ that is naturally
endowed with the de~Rham differential.

 Now let $\tau\:X\rarrow T$ be a morphism of schemes.
 A discussion of smooth morphisms of schemes can be found
in~\cite[Section Tag~01V4]{SP}.
 Let us say that the morphism~$\tau$ is \emph{weakly smooth} if,
for every pair of affine open subschemes $V\subset X$ and $U\subset T$
such that $\tau(V)\subset U$, the related morphism of commutative rings
$\cO_T(U)\rarrow\cO_X(V)$ is weakly smooth.
 Equivalently, in view of the discussion in
Section~\ref{vector-fields-subsecn}, a morphism of schemes
$\tau\:X\rarrow T$ is weakly smooth if and only if the quasi-coherent
sheaf of K\"ahler differentials $\Omega_{X/T}$ is finite locally free
on $X$ is the sense of Section~\ref{loc-free-sheaves-subsecn}.

 A weakly smooth morphism of schemes $\tau\:X\rarrow T$ is said to
have \emph{bounded relative dimension} if the finite locally free
sheaf $\Omega_{X/T}$ on $X$ has bounded rank.
 A weakly smooth morphism~$\tau$ is said to have \emph{relative
dimension~$d$} if the finite locally free sheaf $\Omega_{X/T}$ on $X$
has constant rank~$d$.

 Once again, following
Section~\ref{vector-fields-subsecn}, the sheaf of $\cO_X$\+modules
$\vect_{X/T}$ on $X$ is computed as
$\vect_{X/T}=\cHom_{\cO_X}(\Omega_{X/T},\cO_X)$.
 So the sheaf of $\cO_X$\+modules $\vect_{X/T}$ is finite locally free
whenever the morphism of schemes~$\tau$ is weakly smooth.
 In particular, the sheaf $\vect_{X/T}$ is quasi-coherent in this case;
so $\vect_{X/T}$ is a quasi-coherent Lie algebroid over~$X$.

 Assuming that the morphism of schemes~$\tau$ is weakly smooth,
we define the \emph{quasi-coherent quasi-algebra of crystalline
differential operators on $X$ over~$T$} as the universal enveloping
quasi-algebra of the quasi-coherent Lie algebroid $\vect_{X/T}$,
in the sense of Section~\ref{enveloping-algebra-subsecn}; that is
$\cD^\cry_{X/T}=\cA_X(\vect_{X/T})$.
 We also define the \emph{de~Rham quasi-coherent DG\+quasi-algebra}
$\Omega^\bu_{X/T}$ as the Chevalley--Eilenberg quasi-coherent
DG\+quasi-algebra of the quasi-coherent Lie algebroid
$\vect_{X/T}$ over $X$, in the sense of
Section~\ref{chevalley-eilenberg-qcoh-cdg-quasi-algebra-subsecn};
that is $\Omega^\bu_{X/T}=\cC^\bu_X(\vect_{X/T})$.
 This definition agrees with the ones in~\cite[Section~IV.16.6]{EGA4}
and~\cite[Section Tag~0FKL]{SP} (where the settings are much
more general).
 In particular, the quasi-coherent sheaf of K\"ahler differentials
$\Omega_{X/T}$ is the degree~$1$ component of the de~Rham
quasi-coherent DG\+quasi-algebra $\Omega^\bu_{X/T}$, that is
$\Omega^1_{X/T}=\Omega_{X/T}$.
 The constructions from
Section~\ref{chevalley-eilenberg-cdg-modules-subsecn} also provide
two \emph{de~Rham quasi-coherent DG\+modules} over $\Omega^\bu_{X/T}$,
the left DG\+module $\cC^\bu_X(\vect_{X/T},\cA_X)=
\Omega^\bu_{X/T}\ot_{\cO_X}\cD^\cry_{X/T}$ and the right DG\+module
$\cC_\bu^X(\cA_X,\vect_{X/T})=
\cD^\cry_{X/T}\ot_{\cO_X}\cHom_{\cO_X}(\Omega^\bu_{X/T},\cO_X)$.

\subsection{Symmetric tensors}
 Let $R$ be a commutative ring, $M$ be an $R$\+module, and $n\ge0$ be
an integer.
 The definition of the symmetric algebra $\Sym^*_R(M)$ was given in
Section~\ref{lie-algebroid-pbw-theorem-subsecn}.

 An element $t\in M^{\ot_R\,n}=M\ot_RM\ot_R\dotsb\ot_RM$ ($n$~factors)
is called \emph{symmetric} if one has $\sigma(t)=t$ for every
permutation~$\sigma$ of the set $\{1,2,\dotsc,n\}$ indexing the tensor
factors.
 The natural \emph{symmetrization map}
$$
 \sy_n\:\Sym_R^n(M)\lrarrow M^{\ot_R\,n}
$$
is an $R$\+module map given by the formula
$$
 \sy_n(m_1*\dotsb*m_n)=\sum\nolimits_{\sigma\in\Sigma_n}
 m_{\sigma(1)}\ot\dotsb m_{\sigma(n)},
$$
where $\Sigma_n$ denotes the group of all permutations of
$\{1,2,\dotsc,n\}$ and $m_1$,~\dots, $m_n\in M$.
 One can easily check that the map $\sy_n$ is well-defined.
 Clearly, all the tensors in the image of the symmetrization map
are symmetric.

\begin{lem} \label{symmetric-tensors-lemma}
 Let $R$ be a commutative ring contaning the field\/ $\boQ$ as
a subring, and let $M$ be a flat $R$\+module.
 Then, for every integer $n\ge0$, the map\/ $\sy_n$ is an isomorphism
between the $R$\+module $\Sym_R^n(M)$ and the $R$\+submodule
of symmetric tensors in $M^{\ot_R\,n}$.
\end{lem}

\begin{proof}
 The proof is similar to that of
Lemma~\ref{skew-symmetric-tensors-lemma}.
 However, the assertion of the present lemma does \emph{not} hold
without the assumption that all the positive integers are invertible
in~$R$.
 It suffices to take $M=R$ for a counterexample.
\end{proof}

\begin{lem} \label{symmetric-dual-lemma}
 Let $R$ be a commutative ring contaning the field\/ $\boQ$ as
a subring, and let $M$ be a finitely generated projective $R$\+module.
 Then, for every integer $n\ge0$, there is a natural isomorphism
$$
 \Sym\nolimits_R^n(\Hom_R(M,R))\simeq
 \Hom_R(\Sym\nolimits_R^n(M),R)
$$
provided by the natural $R$\+bilinear pairing
$$
 \Sym\nolimits_R^n(M)\times\Sym\nolimits_R^n(\Hom_R(M,R))
 \lrarrow R
$$
given by the formula
$$
 \langle m_1*\dotsb*m_n,\>b_n*\dotsb*b_1\rangle
 = \sum\nolimits_{\sigma\in\Sigma_n}
 \prod\nolimits_{i=1}^n\langle m_i,b_{\sigma(i)}\rangle,
$$
where $m$, $m_i\in M$, \ $b$, $b_i\in\Hom_R(M,R)$, and
$(m,b)\longmapsto\langle m,b\rangle$ denotes the natural
$R$\+bilinear pairing $M\times\Hom_R(M,R)\rarrow R$.
\end{lem}

\begin{proof}
 Once again, the assertion does \emph{not} hold without the assumption
that all the positive integers are invertible in~$R$ (e.~g., $M=R$ is
a counterexample).
 Nevertheless, the proof is similar to that of
Lemma~\ref{exterior-dual-lemma}, which holds over any commutative ring.
 Both a direct argument and a reduction to
Lemma~\ref{symmetric-tensors-lemma} are possible.
\end{proof}

\subsection{Differential operators in characteristic zero}
\label{diffoperators-characteristic-zero-subsecn}
 Let $(R,\g)$~be a Lie algebroid.
 Then the rules
\begin{itemize}
\item $\lambda(\iota(r))(s)=rs$ for all $r\in R$ and $s\in S$;
\item $\lambda(v)(s)=v(s)$ for all $v\in\g$ and $s\in S$
\end{itemize}
define a left action of the enveloping ring $A_R(\g)$ of the Lie
algebroid $(R,\g)$ on the underlying abelian group of the ring~$R$.
 So $R$ is naturally a left $A_R(\g)$\+module.
 
 Recall that, according to Section~\ref{enveloping-algebra-subsecn},
the enveloping ring $A_R(\g)$ is a quasi-algebra over~$R$.
 In fact, the $R$\+$R$\+bimodule $A_R(\g)$ is a strong quasi-module
over $R$ in the sense of~\cite[Section~2.2]{Ptd}.

 It follows that the action of $A_R(\g)$ in $R$ is an action by
strongly $R$\+differential operators (in the sense of
Examples~\ref{diffoperators-quasi-algebras-quasi-modules-exs}).
 So the action of $A_R(\g)$ in $R$ induces a ring homomorphism
$\theta\:A_R(\g)\rarrow\cD^\str_{R/\boZ}(R,R)$.

 For every integer $n\ge0$, the elements of $F_nA_R(\g)$ act on $R$
by strongly $R$\+differential operators of order at most~$n$;
so $\theta\:A_R(\g)\rarrow\cD^\str_{R/\boZ}(R,R)$ is a homomorphism of
filtered rings with their natural filtrations~$F$.
 In particular, for $n=0$, the ring homomorphism $\theta\:A_R(\g)
\rarrow\cD^\str_{R/\boZ}(R,R)$ restricts to the identity isomorphism
$F_0A_R(\g)=R\rarrow R=F_0\cD^\str_{R/\boZ}(R,R)$.

 Let $\varkappa\:K\rarrow R$ be a weakly smooth homomorphism of
commutative rings.
 We are interested in the Lie algebroid $(R,\g)=(R,\vect_{R/K})$.
 By the definition, we have $A_R(\vect_{R/K})=\cD^\cry_{R/K}$.
 Clearly, the action of $\cD^\cry_{R/K}$ in $R$ is an action by
$K$\+linear maps.
 For simplicity of notation, let us write $\cD^\str_{R/K}=
\cD^\str_{R/K}(R,R)$.
 We have constructed a natural comparison homomorphism between two
filtered rings of differential operators
\begin{equation} \label{Dcr-Dst-comparison-affine}
 \theta=\theta_{R/K}\:(\cD^\cry_{R/K},F)\lrarrow(\cD^\str_{R/K},F).
\end{equation}

 The following proposition is a well-known result.

\begin{prop} \label{crystalline-strong-diffops-comparison-prop}
 Let $\varkappa\:K\rarrow R$ be a weakly smooth homomorphism of
commutative rings.
 Assume that the ring $R$ contains the field of rational numbers\/
$\boQ$ as a subring.
 Then the comparison map~\eqref{Dcr-Dst-comparison-affine} is
an isomorphism of filtered rings.
\end{prop}

 Before proving
Proposition~\ref{crystalline-strong-diffops-comparison-prop},
let us state a simple lemma.

\begin{lem} \label{commutator-derivation-lemma}
 Let $K\rarrow R$ be a homomorphism of commutative rings and
$D\:R\rarrow R$ be a $K$\+linear strongly $R$\+differential operator
of order~$\le1$.
 Consider the map $[D,{-}]\:R\rarrow R$ assigning to every element
$r\in R$ the element $s=[D,r]\in R$ such that the operator
$t\longmapsto D(rt)-rD(t)\:R\rarrow R$ (where $t\in R$ is an arbitrary
element) is equal to the multiplication with~$s$.
 Then the map $[D,{-}]\:R\rarrow R$ is a $K$\+linear derivation of
the ring $R$, i.~e., $[D,{-}]\in\vect_{R/K}$.
\end{lem}

\begin{proof}
 Straightforward.
\end{proof}

\begin{proof}[Proof of
Proposition~\ref{crystalline-strong-diffops-comparison-prop}]
 For every integer~$n\ge0$, the map~$\theta$ induces a map of
the successive quotient $R$\+modules
\begin{equation} \label{gr-F-theta-map}
 \gr^F_n\theta\:F_n\cD^\cry_{R/K}/F_{n-1}\cD^\cry_{R/K}
 \lrarrow F_n\cD^\str_{R/K}/F_{n-1}\cD^\str_{R/K}.
\end{equation}
 According to Theorem~\ref{twisted-lie-algebroid-pbw-theorem}(a),
we have a natural surjective $R$\+module map
\begin{equation} \label{surjection-from-symmetric-power-map}
 \Sym_R^n(\vect_{R/K})\lrarrow F_n\cD^\cry_{R/K}/F_{n-1}\cD^\cry_{R/K}.
\end{equation}
 Let us construct a natural injective $R$\+module map
\begin{equation} \label{injection-into-symmetric-tensors-map}
 F_n\cD^\str_{R/K}/F_{n-1}\cD^\str_{R/K}\lrarrow
 \Hom_R(\Sym^n_R(\Omega_{R/K}),R).
\end{equation}

 The map~\eqref{injection-into-symmetric-tensors-map} is the map of
symbols of differential operators.
 It is constructed as follows.
 Given a differential operator $D\in F_n\cD^\str_{R/K}$ and
an $n$\+tuple of elements $(s_1,\dotsc,s_n)$ in $R$, consider
the iterated commutator $[\dotsm[[D,s_1],s_2]\dotsm s_n]$
($n$~nested brackets).
 By the definition of the filtration $F$ on $\cD^\str_{R/K}$,
the map $[\dotsm[[D,s_1],s_2]\dotsm s_n]\:R\rarrow R$ is
a differential operator of order~$0$, i.~e., the operator of
multiplication with an element $\delta_D(s_1,\dotsc,s_n)\in R$.

 One can easily check that the map $\delta_D$ is symmetric,
i.~e., the element $\delta_D(s_1,\dotsc,s_n)\in R$ remains unchanged
by permutations of the elements $s_1$,~\dots, $s_n\in R$.
 By Lemma~\ref{commutator-derivation-lemma}, for every collection
of elements $s_1$,~\dots, $s_{n-1}\in R$ there exists a unique
$R$\+linear map $v_D(s_1,\dotsc,s_{n-1})\:\Omega_{R/K}\rarrow R$
such that $\delta_D(s_1,\dotsc,s_n)=v_D(s_1,\dotsc,s_{n-1})(ds_n)$
for all $s_n\in R$.  {\hbadness=1425\par}

 Set the image of the differential operator $D$ under
the map~\eqref{injection-into-symmetric-tensors-map} to be
the $R$\+linear map $\Sym^n_R(\Omega_R/K)\rarrow R$ given by
the rule
\begin{equation} \label{symbol-defined}
 r_1d(s_1)*r_2d(s_2)*\dotsb*r_nd(s_n)=
 r_1\dotsm r_n\delta_D(s_1,\dotsc,s_n).
\end{equation}
 It follows from the previous paragraph that the desired map
$\Sym^n_R(\Omega_R/K)\rarrow R$ is well-defined by
the formula~\eqref{symbol-defined}.
 Thus the map~\eqref{injection-into-symmetric-tensors-map} is
well-defined, too.

 Any differential operator $D\in F_n\cD^\str_{R/K}$ annihilated
by the map~\eqref{injection-into-symmetric-tensors-map} belongs,
by the definition, to $F_{n-1}\cD^\str_{R/K}$.
 Thus the map~\eqref{injection-into-symmetric-tensors-map} is injective.

 It remains to notice that the composition of the three
maps~\eqref{surjection-from-symmetric-power-map},
\eqref{gr-F-theta-map}, and~\eqref{injection-into-symmetric-tensors-map}
is equal to the map $\Sym^n_R(\Hom_R(\Omega_{R/K},R))\rarrow
\Hom_R(\Sym^n_R(\Omega_{R/K}),R)$ from
Lemma~\ref{symmetric-dual-lemma} for the $R$\+module $M=\Omega_{R/K}$.
 As the $R$\+module $\Omega_{R/K}$ is finitely generated and projective
by assumption, Lemma~\ref{symmetric-dual-lemma} tells us that
the composition of the three
maps~\eqref{surjection-from-symmetric-power-map},
\eqref{gr-F-theta-map}, and~\eqref{injection-into-symmetric-tensors-map}
is an isomorphism.
 Since the map~\eqref{surjection-from-symmetric-power-map} is
surjective and the map~\eqref{injection-into-symmetric-tensors-map}
is injective, we can conclude that all the three
maps~\eqref{surjection-from-symmetric-power-map},
\eqref{gr-F-theta-map}, and~\eqref{injection-into-symmetric-tensors-map}
are isomorphisms.
 In particular, the map~\eqref{gr-F-theta-map} is an isomorphism,
and it follows that the map~\eqref{Dcr-Dst-comparison-affine} is
a filtered isomorphism, as desired.

 In addition to a proof of
Proposition~\ref{crystalline-strong-diffops-comparison-prop},
we have also obtained an alternative proof of
the Poincar\'e--Birkhoff--Witt theorem
(Theorem~\ref{twisted-lie-algebroid-pbw-theorem}(b)) for
crystalline differential operators in characteristic zero.
\end{proof}

 Let $X$ be a scheme and $\g$~be a Lie algebroid over~$X$.
 Then the collection of natural left actions of the enveloping rings
$\cA_X(\g)(V)=A_{\cO_X(V)}(\g(V))$ on the rings $\cO_X(V)$ for
affine open subschemes $V\subset X$ defines a left action
of the quasi-coherent quasi-algebra $\cA_X(\g)$ on the structure
sheaf $\cO_X$ of the scheme~$X$.
 So $\cO_X$ is naturally a quasi-coherent left module over
the quasi-coherent quasi-algebra~$\cA_X(\g)$.

 Now let $\tau\:X\rarrow T$ be a weakly smooth morphism of schemes.
 Then the collection of homomorphisms~$\theta_{R/K}$ for the rings
$R=\cO_X(V)$ and $K=\cO_T(U)$, where $V\subset X$ and $U\subset T$
are affine open subschemes such that $\tau(V)\subset U$, defines
a morphism of sheaves of (filtered) quasi-algebras
$$
 \theta_{X/T}\:\cD^\cry_{X/T}\lrarrow\cD^\str_{X/T}.
$$
 Here we use the notation $\cD^\str_{X/T}=\cD^\str_{X/T}(\cO_X,\cO_X)$
for the sheaf of quasi-algebras of strongly $X/T$\+differential
operators acting on the structure sheaf $\cO_X$ of the scheme~$X$ (see
Examples~\ref{fiberwise-differential-operators-examples}(2\+-3)).
 The quasi-coherent quasi-algebra $\cD^\cry_{X/T}$ over $X$ was
defined in Section~\ref{crystalline-diffoperators-subsecn}.

\begin{cor} \label{crystalline-strong-diffops-comparison-cor}
 Let $\tau\:X\rarrow T$ be a weakly smooth morphism of schemes such
that $X$ is a scheme over\/ $\Spec\boQ$.
 Then the morphism~$\theta_{X/T}$ is an isomorphism of sheaves of 
filtered quasi-algebras $\cD^\cry_{X/T}\lrarrow\cD^\str_{X/T}$ over~$X$.
\end{cor}

\begin{proof}
 Follows from
Proposition~\ref{crystalline-strong-diffops-comparison-prop}.
\end{proof}

 In particular, since $\cD^\cry_{X/T}$ is a quasi-coherent
quasi-algebra over $X$ according to
Section~\ref{crystalline-diffoperators-subsecn}, it follows from
Corollary~\ref{crystalline-strong-diffops-comparison-cor} that
$\cD^\str_{X/T}$ is a quasi-coherent quasi-algebra over $X$, too.
 (See Example~\ref{fiberwise-differential-operators-examples}(4) for
an essentially much more general version of this assertion.)

\subsection{Twisted differential operators}
\label{twisted-diffoperators-subsecn}
 The notion of \emph{twisted differential operators} seems to go back
to~\cite[Section~2]{BB}.
 See~\cite[Section~10.5]{Prel} for an introductory discussion.

 Let $\varkappa\:K\rarrow R$ be a weakly smooth homomorphism of
commutative rings and $(R,\vect_{R/K})$ be the Lie algebroid of
vector fields on $R$ over~$K$.
 Suppose given a twisted Lie algebroid $(R,\g,\widetilde\g)$ such that
$\g=\vect_{R/K}$ and the underlying Lie algebroid $(R,\g)$ of
the twisted Lie algebroid $(R,\g,\widetilde\g)$ is $(R,\g)=
(R,\vect_{R/K})$.
 Put $\widetilde\vect_{R/K}=\widetilde\g$.
 Then the twisted universal enveloping ring
$A_R(\vect_{R/K},\widetilde\vect_{R/K})$ is
called a ring of \emph{twisted} (\emph{crystalline})
\emph{differential operators} on $R$ over~$K$.

 The related Chevalley--Eilenberg CDG\+ring
$C^\cu_R(\g,\widetilde\g)=C^\cu(\vect_{R/K},\widetilde\vect_{R/K})$
is a CDG\+ring structure on the graded commutative de~Rham DG\+ring
$\Omega^\bu_{R/K}$.
 So one has $C^*_R(\g,\widetilde\g)=\Omega^*_{R/K}$, and
the differential in $C^\cu_R(\g,\widetilde\g)$ coincides with
the de~Rham differential in $\Omega^\bu_{R/K}$; while the curvature
element in $C^2_R(\g,\widetilde\g)$ is just some element encoding
the datum of the central extension of Lie algebras $0\rarrow R
\overset\iota\rarrow\widetilde\vect_{R/K}\overset\pi\rarrow\vect_{R/K}
\rarrow0$ over $\boZ$ and the action of $R$
on~$\widetilde\vect_{R/K}$.

 Now let $\tau\:X\rarrow T$ be a weakly smooth morphism of schemes
and $\vect_{X/T}$ be the quasi-coherent Lie algebroid of vector fields
on $X$ over~$T$.
 Suppose given a quasi-coherent twisted Lie algebroid
$(\g,\widetilde\g)$ over $X$ such that the underlying quasi-coherent
Lie algebroid~$\g$ of $(\g,\widetilde\g)$ is $\g=\vect_{X/T}$.
 Put $\widetilde\vect_{X/T}=\widetilde\g$.
 Then the twisted universal enveloping quasi-algebra
$\cA_X(\vect_{X/T},\widetilde\vect_{X/T})$ is called a quasi-algebra
of \emph{twisted} (\emph{crystalline}) \emph{differential operators}
on $X$ over~$T$.

 The related Chevalley--Eilenberg quasi-coherent CDG\+quasi-algebra
$\cC^\cu_X(\g,\widetilde\g)=
\cC^\cu_X(\vect_{X/T},\widetilde\vect_{X/T})$ is a quasi-coherent
CDG\+quasi-algebra structure on the graded commutative de~Rham
DG\+quasi-algebra $\Omega^\bu_{X/T}$.
 So one has $\cC^*_X(\g,\widetilde\g)=\Omega^*_{X/T}$, and, for every
affine open subscheme $U\subset X$, the differential on
$\cC^*_X(\g,\widetilde\g)(U)$ coincides with the de~Rham differential
on $\Omega^*_{X/T}(U)$.
 There are, however, nontrivial curvature elements
$h_U\in\cC^2_X(\g,\widetilde\g)(U)$ and nontrivial change-of-connection
elements $a_{VU}\in\cC^1_X(\g,\widetilde\g)(V)$ for pairs of affine
open subschemes $V\subset U\subset X$, forming the CDG\+quasi-algebra
structure $\cC^\cu_X(\g,\widetilde\g)$ on the quasi-coherent
graded quasi-algebra $\cC^*_X(\g,\widetilde\g)$ over~$X$.

 One can say informally that the quasi-coherent twisted Lie algebroid
$(\g,\widetilde\g)=(\vect_{X/T},\widetilde\vect_{X/T})$ and
the quasi-algebra of twisted differential operators
$\cA_X(\g,\widetilde\g)$ over $X$ are ``twisted in two ways''.
 There is a local twisting given by the central extension of Lie
algebras $0\rarrow \cO_X(U)\overset\iota\rarrow\widetilde\vect_{X/T}(U)
\overset\pi\rarrow\vect_{X/T}(U)\rarrow0$ for every affine open
subscheme $U\subset X$.
 This central extension of Lie algebras splits as a short exact
sequence of $\cO_X(U)$\+modules, but not as a short exact sequence
of Lie algebras.
 There is also a global twisting given by the extension of finite
locally free sheaves $0\rarrow\cO_X\overset\iota\rarrow
\widetilde\vect_{X/T}\overset\pi\rarrow\vect_{X/T}\rarrow0$ on $X$,
which does not even have an $\cO_X$\+linear splitting.
 The local twisting gives rise to the curvature elements
$h_U\in\cC^2_X(\g,\widetilde\g)(U)$, while the global twisting produces
the change-of-connection elements $a_{VU}\in\cC^1_X(\g,\widetilde\g)(V)$
in the quasi-coherent CDG\+quasi-algebra $\cC^\cu_X(\g,\widetilde\g)$.

\Section{Thick Graded Modules and Thick CDG-Modules}
\label{thick-graded-modules-secn}

 The terminology ``thick graded modules'' (over the exterior algebra
of a locally free sheaf over a smooth algebraic variety) comes from
the paper~\cite[Section~3.5]{Ryb}.
 In this section, we extend this concept to not necessarily regular
schemes.

 Basically, let $k$~be a field, $R$ be a commutative $k$\+algebra,
$d\ge0$ be an integer, $V=k^d$ be the $d$\+dimensional vector space
over~$k$, and $E=R^d=R\ot_k V$ be the free $R$\+module of rank~$d$.
 Then a graded module $M^*$ over the exterior algebra
$\bigwedge^*_R(E)$ over $R$ is thick if and only if $M^*$ is
a projective/injective module over the exterior algebra
$\bigwedge^*_k(V)$ over~$k$.
 This is a particular case of
Lemma~\ref{thick-graded-modules-restriction-of-scalars} below.

 So a thick graded $\bigwedge^*_R(E)$\+module $M^*$ is
projective/injective along $\bigwedge^*_k(V)$, while no restriction
is imposed on the underlying $R$\+module of~$M^*$.
 In the case of a finitely generated projective $R$\+module $E$ that
is not free or does not have a chosen basis, the more complicated
definition spelled out in this section is needed.

\subsection{Finitely filtered rings and modules}
\label{finitely-filtered-rings-and-modules-subsecn}
 In this section we consider an associative ring $B$ with a finite
decreasing filtration $B=F^0B\supset F^1B\supset F^2B\supset\dotsb\supset
F^nB\supset F^{n+1}B=0$, where $n\ge1$ is an integer.
 The filtration is presumed to be multiplicative, i.~e.,
$F^iB\cdot F^jB\subset F^{i+j}B$ for all $i$, $j\ge0$.

 A \emph{finitely filtered} (\emph{left}) \emph{$B$\+module} $M$ is
a $B$\+module endowed with a decreasing filtration indexed by
the integers, $M=F^{-m}M\supset F^{-m+1}M\supset\dotsb\supset F^mM
\supset F^{m+1}M=0$, where $m\ge1$ is an integer that may depend on~$M$.
 The filtration on filtered $B$\+modules is presumed to be
multiplicative, i.~e., $F^iB\cdot F^jM\subset F^{i+j}M$ for all $i\ge0$,
\,$j\in\boZ$ (and similarly for right $B$\+modules).

 The associated graded ring to $(B,F)$ is
$\gr_F^*B=\bigoplus_{i=0}^n F^iB/F^{i+1}B$.
 The associated graded module to $(M,F)$ is
$\gr_F^*M=\bigoplus_{j=-m}^m F^jM/F^{j+1}M$.
 So $\gr_F^*B$ is naturally a graded ring, and $\gr_F^*M$ is
naturally a graded (left) module over $\gr_F^*B$.

\begin{lem} \label{filtered-modules-ext-spectral-sequence}
 Let $(L,F)$ and $(M,F)$ be two finitely filtered left $B$\+modules.
 Then there is a natural spectral sequence of abelian groups
$$
 E_1^{p,i}=\Ext_{\gr_F^*B}^{i,p}(\gr_F^*L,\gr_F^*M)
 \Longrightarrow E_\infty^{p,i}=\gr_F^p\Ext^i_B(L,M),
 \qquad i\ge0, \ p\in\boZ
$$
with the differentials $d_r^{p,i}\:E_r^{p,i}\rarrow E_r^{p+r,i+1}$.
 Here the first grading~$i$ on the abelian groups\/
$\Ext_{\gr_F^*B}^*(\gr_F^*L,\gr_F^*M)$ is the usual cohomological
grading, while the second grading~$p$ is the \emph{internal} grading
induced by the grading of the graded ring $\gr_F^*B$ and the graded
modules $\gr_F^*L$ and $\gr_F^*M$.
 The natural decreasing filtrations $F$ on the abelian groups
$\Ext^i_B(L,M)$ are finite, i.~e., for every $i\ge0$ there exists
an integer $e_i\ge1$ such that $F^{-e_i}\Ext^i_B(L,M)=\Ext^i_B(L,M)$
and $F^{e_i+1}\Ext^i_B(L,M)=0$.
\end{lem}

\begin{proof}
 The construction uses a filtered projective resolution of
the filtered $B$\+module $(L,F)$.
 A complex of finitely filtered left $B$\+modules
$$
 \dotsb\lrarrow(P_2,F)\lrarrow(P_1,F)\lrarrow(P_0,F)\lrarrow(L,F)
 \lrarrow0
$$
is said to be a \emph{filtered projective resolution} of $(L,F)$ if
the left $B$\+modules $P_i$ are projective, the graded left
modules $\gr_F^*P_i$ over the graded ring $\gr_F^*B$ are projective,
and the complex $\gr_F^*P_\bu$ is exact, i.~e., it is a resolution
of the graded module $\gr_F^*L$.
 Equivalently, one can say that $(P_\bu,F)$ must be a projective
resolution of $(L,F)$ in the exact category of finitely filtered
left $B$\+modules.

 To construct a filtered projective resolution of a filtered $B$\+module
$(L,F)$, one proceeds by induction.
 Suppose that $F^{-l}L=L$ and $F^{l+1}L=0$.
 For every filtration degree $-l\le j\le l$, pick a set of elements
$X_j\subset F^jL$ such that one has
$$
 F^kL=\sum\nolimits_{j=-l}^k\sum\nolimits_{x\in X_j}F^{k-j}B\cdot x
$$
for every $-l\le k\le l$.
 So the submodule $F^lL\subset L$ must be generated by its subset
$X_l\subset F^lL$, the submodule $F^{l-1}L\subset L$ must be
generated by its subset $X_l\cup X_{l-1}\subset F^{l-1}L$, etc.,
and the submodule $F^{-l}L=L$ must be generated by its subset
$\bigcup_{j=-l}^l X_j\subset L$.
 One can say that the collection of subsets $X_j$ must be a set of
filtered generators of $L$, or that the family of subsets $X_j$ must
generate the filtered $B$\+module $(L,F)$ with its filtration.

 For convenience of notation, put $Z_j=X_j$ for every $-l\le j\le l$.
 The difference between $X_j$ and $Z_j$ is that $X_j$ is a subset
in $L$, while $Z_j$ will be a set of free generators of the $B$\+module
we are about to construct.
 Put $Z=\bigcup_{j=-l}^lZ_j$ and $Z_j=\varnothing$ for $j>l$, and let
$P_0=\bigoplus_{x\in X}Bx$ be the free $B$\+module generated by
the set~$Z$.
 Endow the $B$\+module $P_0$ with the filtration $F$ defined by the rule
$F^kP_0=\sum_{j=-l}^k\sum_{z\in Z_j}F^{k-j}B\cdot z\allowbreak\subset
P_0$ for every $-l\le k\le l+n$.
 Then the collection of identity maps of the set of generators
$Z_j\rarrow X_j$ defines a surjective $B$\+module map $P_0\rarrow L$
respecting the filtrations.
 In fact, by construction, the map of the filtration components
$F^kP_0\rarrow F^kL$ is surjective for every $k\in\boZ$.

 Denote by $L_1$ the kernel of the map $P_0\rarrow L_0$, endowed with
the filtration $F$ induced by the filtration $F$ on~$P_0$.
 Applying the same construction to the filtered $B$\+module $(L_1,F)$
instead of $(L,F)$, we obtain a free $B$\+module $P_1$ endowed with
a filtration $F$ and a surjective $B$\+module map $P_1\rarrow L_1$, etc.
 Proceeding in this way, we obtain the desired filtered projective
resolution $(P_\bu,F)$ of the filtered $B$\+module~$L$.

 Finally, for any two finitely filtered $B$\+modules $(Q,F)$ and
$(M,F)$, the abelian group $\Hom_B(Q,M)$ is naturally finitely
filtered, with the filtration component $F^k\Hom_B(Q,M)$ consisting
of all $B$\+linear maps $Q\rarrow M$ taking $F^jQ$ into $F^{j+k}M$
for every $j\in\boZ$.
 Therefore, the complex of abelian groups $\Hom_B(P_\bu,M)$ carries
a natural filtration~$F$.
 The desired spectral sequence $E_r^{p,i}$ is the spectral sequence
associated with the filtered complex $(\Hom_B(P_\bu,M),F)$.
\end{proof}

\subsection{Very weakly relatively adjusted modules}
\label{very-weak-relatively-adjusted-subsecn}
 Given a ring $B$ with a subring $A\subset B$, a left $B$\+module $M$
is usually called \emph{projective relative to~$A$} if the functor
$\Hom_B(M,{-})$ takes short exact sequences of left $B$\+modules that
are \emph{split as short exact sequences of left $A$\+modules} to
short exact sequences of abelian groups.
 A left $B$\+module $M$ is projective relative to $A$ if and only if
$M$ is isomorphic to a direct summand of the $B$\+module $B\ot_AL$
for some left $A$\+module~$L$.
 Left $B$\+modules \emph{injective relative to~$A$} are defined and
described in the dual way.

 A left $B$\+module $M$ is called \emph{weakly projective relative
to~$A$} if the functor $\Hom_B(M,{-})$ takes short exact sequences of
left $B$\+modules in which all the three terms are \emph{injective
as left $A$\+modules} to short exact sequences of abelian
groups~\cite[Sections~4.1 and~4.3]{BP}, \cite[Section~5]{Pfp},
\cite[Example~3.9]{Pctrl}.
 Under certain assumptions on $A$ and $B$, one can prove that
a $B$\+module is weakly projective relative to $A$ if and only if it
is a direct summand of a $B$\+module finitely filtered by $B$\+modules
of the form $B\ot_AL$ \,\cite[Example~3.9]{Pctrl}.
 Left $B$\+modules \emph{weakly injective relative to~$A$} are
defined and described in the dual way~\cite[Example~2.10]{Pctrl}.

 In this section we discuss certain properties that appear to be
even weaker than the weak relative projectivity/injectivity.
 The natural setting for these properties, which we call \emph{very
weak relative projectivity/injectivity}, is that of a finitely
filtered ring rather than a ring with a subring.

\begin{lem} \label{very-weakly-relatively-projective-lemma}
 Let $B=F^0B\supset F^1B\supset F^2B\supset\dotsb\supset
F^nB\supset F^{n+1}B=0$ be a finitely filtered associative ring.
 Assume that the successive quotients $F^iB/F^{i+1}B$, \, $0\le i\le n$,
are flat as right modules over the ring $B/F^1B$.
 Let $M$ be a left $B$\+module.
 Then the following conditions are equivalent:
\begin{enumerate}
\item $\Tor^B_1(B/F^1B,M)=0$;
\item $\Tor^B_i(G,M)=0$ for all $i>0$ and any flat right
$B/F^1B$\+module $G$ viewed as a right $B$\+module with a zero action
of~$F^1B$;
\item $\Ext_B^1(M,J)=0$ for any injective left $B/F^1B$\+module $J$
viewed as a left $B$\+module with a zero action of~$F^1B$;
\item $\Ext_B^i(M,J)=0$ for all $i>0$ and any injective left
$B/F^1B$\+module $J$ viewed as a left $B$\+module with a zero action
of~$F^1B$;
\item there exists a finite decreasing filtration $F$ on
the $B$\+module $M$ compatible with the filtration $F$ on the ring $B$
such that the associated graded module\/ $\gr_F^*M$ over the graded
ring\/ $\gr_F^*B$ has the form
$$
 \gr_F^*M\simeq\gr_F^*B\ot_{B/F^1B}N
$$
(an isomorphism of graded left modules over\/ $\gr_F^*B$) for some
left $B/F^1B$\+mod\-ule~$N$.
\end{enumerate}
\end{lem}

\begin{proof}
 (2)\,$\Longrightarrow$\,(1) and (4)\,$\Longrightarrow$\,(3) Obvious.

 (1)\,$\Longleftrightarrow$\,(3) and (2)\,$\Longleftrightarrow$\,(4)
 Use the fact that a left $B/F^1B$\+module $J$ is injective if and only
if it is a direct summand of the character $B/F^1B$\+module
$\Hom_\boZ(G,\boQ/\boZ)$, where $G$ is some flat right $B/F^1B$\+module.

 (1)\,$\Longrightarrow$\,(5) Since the functor $\Tor$ preserves filtered
inductive limits of modules, and all flat modules are filtered inductive
limits of projective ones, it follows from~(1) that
$\Tor^B_1(G,M)=0$ for all flat right $B/F^1B$\+modules $G$ viewed
as right $B$\+modules with a zero action of~$F^1B$.
 In particular, since the right $B/F^1B$\+modules $F^iB/F^{i+1}B$ are
flat by assumption, condition~(1) implies the Tor vanishing
$\Tor^B_1(F^iB/F^{i+1}B,M)=0$ for all $0\le i\le n$.
 Consequently, we have $\Tor^B_1(B/F^iB,M)=0$ for all $1\le i\le n$.

 It follows that the map $F^iB\ot_BM\rarrow M$ induced by the inclusion
of right $B$\+modules $F^iB\rarrow B$ is injective.
 In other words, the left $B$\+module map
$$
 F^iB\ot_BM\lrarrow F^iB\cdot M\subset M
$$
is an isomorphism.
 Set $F^iM=F^iB\cdot M\subset M$ for all $i\ge0$.
 Then we have
\begin{multline*}
 F^iM/F^{i+1}M\simeq (F^iB\ot_BM)/(F^{i+1}M\ot_BM) \\
 \simeq (F^iB/F^{i+1}B)\ot_BM\simeq
 (F^iB/F^{i+1}B)\ot_{B/F^1B}(M/F^1M),
\end{multline*}
so the associated graded module $\gr_F^*M$ over $\gr_F^*B$ has
the form
$$
 \gr_F^*M\simeq\gr_F^*B\ot_{B/F^1B}N,
 \quad\text{where } N=M/F^1M=B/F^1B\ot_BM.
$$

 (5)\,$\Longrightarrow$\,(4) The argument is based on
Lemma~\ref{filtered-modules-ext-spectral-sequence}.
 Let us view $J$ as a filtered $B$\+module with the trivial filtration
$F^0J=J$, \ $F^1J=0$; so $\gr_F^*J=J$.
 Then we have
$$
 \Ext^i_{\gr_F^*B}(\gr_F^*M,J)\simeq
 \Ext^i_{\gr_F^*B}(\gr_F^*B\ot_{B/F^1B}N,\>J)\simeq
 \Ext^i_{B/F^1B}(N,J)=0
$$
for all $i>0$.
 Here the second isomorphism holds by
Lemma~\ref{Ext-homological-formula} and uses the assumption that
the graded right $B/F^1B$\+module $\gr_F^*B$ is flat.
 By Lemma~\ref{filtered-modules-ext-spectral-sequence}, it follows
that $\Ext^i_B(M,J)=0$.
\end{proof}

\begin{lem} \label{very-weakly-relatively-injective-lemma}
 Let $B=F^0B\supset F^1B\supset F^2B\supset\dotsb\supset
F^nB\supset F^{n+1}B=0$ be a finitely filtered associative ring.
 Assume that the successive quotients $F^iB/F^{i+1}B$, \, $0\le i\le n$,
are projective as left modules over the ring $B/F^1B$.
 Let $M$ be a left $B$\+module.
 Then the following conditions are equivalent:
\begin{enumerate}
\item $\Ext^1_B(B/F^1B,M)=0$;
\item $\Ext^i_B(P,M)=0$ for all $i>0$ and any projective left
$B/F^1B$\+module $P$ viewed as a left $B$\+module with a zero action
of~$F^1B$;
\item there exists a finite decreasing filtration $F$ on
the $B$\+module $M$ compatible with the filtration $F$ on the ring $B$
such that the associated graded module\/ $\gr_F^*M$ over the graded
ring\/ $\gr_F^*B$ has the form
$$
 \gr_F^*M\simeq\Hom_{B/F^1B}(\gr_F^*B,N)
$$
(an isomorphism of graded left modules over\/ $\gr_F^*B$) for some
left $B/F^1B$\+mod\-ule~$N$.
\end{enumerate}
\end{lem}

\begin{proof}
 (2)\,$\Longrightarrow$\,(1) Obvious.

 (1)\,$\Longrightarrow$\,(3) Since the functor $\Ext$ takes infinite
direct sums of modules in its first argument to direct products of
abelian groups, it follows from~(1) that $\Ext_B^1(P,M)=0$ for all
projective left $B/F^1B$\+modules $P$ viewed as left $B$\+modules with
a zero action of $F^1B$.
 In particular, since the left $B/F^1B$\+modules $F^iB/F^{i+1}B$ are
projective by assumption, condition~(1) implies the Ext vanishing
$\Ext^1_B(F^iB/F^{i+1}B,M)=0$ for all $0\le i\le n$.
 Consequently, we have $\Ext^1_B(B/F^iB,M)=0$ for all $1\le i\le n$.

 It follows that the map $M\rarrow\Hom_B(F^iB,M)$ induced by
the inclusion of left $B$\+modules $F^iB\rarrow B$ is surjective.
 In other words, we have a short exact sequence of left $B$\+modules
$$
 0\lrarrow\Hom_B(B/F^iB,M)\lrarrow M\lrarrow\Hom_B(F^iB,M)\lrarrow0.
$$
 For every $i\ge0$, let $F^{-i}M=\Hom_B(B/F^{i+1}B,M)$ be the submodule
of all elements annihilated by $F^iB$ in~$M$.
 So the filtration $F$ on $M$ is actually concentrated in the finite
interval of nonpositive integers $-n\le j\le0$.

 Now we have
\begin{multline*}
 F^jM/F^{j+1}M\simeq\ker(\Hom_B(F^{-j}B,M)\to\Hom_B(F^{-j+1}B,M)) \\
 \simeq\Hom_B(F^{-j}B/F^{-j+1}B),M)\simeq
 \Hom_{B/F^1B}(F^{-j}B/F^{-j+1}B,F^0M).
\end{multline*}
 So the associated graded module $\gr_F^*M$ over $\gr_F^*B$ has the form
$$
 \gr_F^*M\simeq\Hom_{B/F^1B}(\gr_F^*B,N),
 \quad\text{where } N=F^0M=\Hom_B(B/F^1B,M).
$$

 (3)\,$\Longrightarrow$\,(2) The argument is based on
Lemma~\ref{filtered-modules-ext-spectral-sequence}.
 Let us view $P$ as a filtered $B$\+module with the trivial filtration
$F^0P=J$, \ $F^1P=0$; so $\gr_F^*P=P$.
 Then we have
$$
 \Ext^i_{\gr_F^*B}(P,\gr_F^*M)\simeq
 \Ext^i_{\gr_F^*B}(P,\,\Hom_{B/F^1B}(\gr_F^*B,N))\simeq
 \Ext^i_{B/F^1B}(P,N)=0
$$
for all $i>0$.
 Here the second isomorphism holds by
Lemma~\ref{Ext-homological-formula} and uses the assumption that
the graded left $B/F^1B$\+module $\gr_F^*B$ is projective.
 By Lemma~\ref{filtered-modules-ext-spectral-sequence}, it follows
that $\Ext^i_B(P,M)=0$.
\end{proof}

 We will say that a $B$\+module $M$ is \emph{very weakly relatively
projective} (relative to the finite decreasing filtration $F$ on~$B$)
if it satisfies any one of the equivalent conditions of
Lemma~\ref{very-weakly-relatively-projective-lemma}.
 Similarly, a $B$\+module $M$ is said to be \emph{very weakly relatively
injective} (relative to the filtration $F$ on~$B$) if it satisfies
any one of the equivalent conditions of
Lemma~\ref{very-weakly-relatively-injective-lemma}.

\subsection{Canonical filtrations on very weakly relatively
adjusted modules} \label{canonical-filtrations-on-vwra-modules-subsecn}
 Every module over a finitely filtered ring has two natural filtrations,
which were constructed in the proofs of
Lemmas~\ref{very-weakly-relatively-projective-lemma}
and~\ref{very-weakly-relatively-injective-lemma}.
 For very weakly relatively projective/injective modules, such
canonical filtrations are particularly well-behaved, in that they have
good exactness properties.

\begin{lem} \label{vwr-projective-canonical-filtration-lemma}
 Let $B=F^0B\supset F^1B\supset F^2B\supset\dotsb\supset
F^nB\supset F^{n+1}B=0$ be a finitely filtered associative ring.
 Assume that the successive quotients $F^iB/F^{i+1}B$, \, $0\le i\le n$,
are flat as right modules over the ring $B/F^1B$.
 Let $M$ be a very weakly relatively projective left $B$\+module.
 Then a filtration $F$ on $M$ satisfying the conditions of
Lemma~\textup{\ref{very-weakly-relatively-projective-lemma}(5)} is
unique, and preserved by all morphisms of very weakly relatively
projective left $B$\+modules.
 For any short exact sequence of very weakly relatively projective
left $B$\+modules\/ $0\rarrow K\rarrow L\rarrow N\rarrow0$ and every
integer $i\in\boZ$, the short sequence of abelian groups (in fact,
$B$\+modules)\/ $0\rarrow F^iK\rarrow F^iL\rarrow F^iN\rarrow0$
is exact.
\end{lem}

\begin{proof}
 The isomorphism of graded $\gr_F^*B$\+modules $\gr_F^*M\simeq
\gr_F^*B\ot_{B/F^1B}N$ implies that the graded $\gr_F^*B$\+module
$\gr_F^*M$ is generated by its degree-zero component
$\gr_F^0M=N$.
 It follows that $\gr_F^iM=0$ for $i<0$, so $F^0M=M$, and
the filtration $F$ on $M$ is generated by $F^0M$, i.~e.,
$F^iM=F^iB\cdot F^0M=F^iB\cdot M\subset M$ for all $i>0$.
 This proves uniqueness of the filtration $F$ on~$M$.
 The preservation by morphisms in now obvious.
 Finally, according to the proof of
Lemma~\ref{very-weakly-relatively-projective-lemma}\,%
(1)\,$\Rightarrow$\,(5), the natural map
$F^iB\ot_BM\rarrow F^iB\cdot M$ is an isomorphism.
 It remains to observe that the short sequence
$0\rarrow F^iB\ot_RK\rarrow F^iB\ot_RL\rarrow F^iB\ot_RN\rarrow0$
is exact because $\Tor_1^B(F^iB,N)=0$.
 The latter $\Tor$ vanishing holds since the filtration $F$ on $B$ is
finite and $\Tor_1^B(F^jB/F^{j+1}B,M)=0$ for all $j\ge0$.
\end{proof}

\begin{lem} \label{vwr-injective-canonical-filtration-lemma}
 Let $B=F^0B\supset F^1B\supset F^2B\supset\dotsb\supset
F^nB\supset F^{n+1}B=0$ be a finitely filtered associative ring.
 Assume that the successive quotients $F^iB/F^{i+1}B$, \, $0\le i\le n$,
are projective as left modules over the ring $B/F^1B$.
 Let $M$ be a very weakly relatively injective left $B$\+module.
 Then a filtration $F$ on $M$ satisfying the conditions of
Lemma~\textup{\ref{very-weakly-relatively-injective-lemma}(3)} is
unique, and preserved by all morphisms of very weakly relatively
injective left $B$\+modules.
 For any short exact sequence of very weakly relatively injective
left $B$\+modules\/ $0\rarrow K\rarrow L\rarrow N\rarrow0$ and every
integer $i\in\boZ$, the short sequence of abelian groups (in fact,
$B$\+modules)\/ $0\rarrow F^iK\rarrow F^iL\rarrow F^iN\rarrow0$
is exact.
\end{lem}

\begin{proof}
 The isomorphism of graded $\gr_F^*B$\+modules $\gr_F^*M\simeq
\Hom_{B/F^1B}(\gr_F^*B,N)$ implies that the graded $\gr^F_*B$\+module
$\gr_F^*M$ is cogenerated by its degree-zero component $\gr_F^0M=N$,
i.~e., the natural map $\gr_F^*M\rarrow\Hom_\boZ(\gr_F^*B,\gr_F^0M)$
is injective.
 It follows that $\gr_F^iM=0$ for $i>0$, so $F^1M=0$, and
the filtration $F$ on $M$ is cogenerated by $M/F^1M$, i.~e.,
$F^{i+1}M$ is the kernel of the map $M=\Hom_B(B,M)\rarrow
\Hom_B(F^{-i}B,M/F^1M)=\Hom_B(F^{-i}B,M)$ for all $i<0$.
 This proves uniqueness of the filtration $F$ on $M$ and its
preservation by morphisms.
 Finally, according to the proof of
Lemma~\ref{very-weakly-relatively-injective-lemma}\,%
(1)\,$\Rightarrow$\,(3), the natural map
$F^{i+1}M\rarrow\Hom_B(B/F^{-i}B,M)$ is an isomorphism.
 It remains to observe that the short sequence
$0\rarrow\Hom_B(B/F^{-i}B,K)\rarrow\Hom_B(B/F^{-i}B,L)\rarrow
\Hom_B(B/F^{-i}B,N)\rarrow0$ is exact because $\Ext^1_B(B/F^{-i}B,K)=0$
(according to the same proof).
\end{proof}

\subsection{Thick graded modules}  \label{thick-graded-modules-subsecn}
 Let $R$ be a commutative ring and $E$ be a finitely generated
projective $R$\+module.
 We are interested in $\boZ$\+graded modules over the graded ring
$B^*=\bigwedge_R^*(E)$, i.~e., the exterior algebra of $R$\+module~$E$
(see Section~\ref{skew-symmetric-tensors-subsecn}).
 
 We endow the graded ring $B^*$ with the decreasing filtration $F$
induced by the grading, i.~e., $F^iB^*=\bigoplus_{j=i}^\infty B^j$.
 Notice that $\bigwedge_R^{n+1}(E)=0$ if $E$ is generated by
$n$~elements; so the decreasing filtration $F$ on $B^*$ is finite.
 Obviously, all the components $F^iB^*$ of the filtration $F$ on $B^*$
are homogeneous subgroups (in fact, homogeneous ideals) in
the graded ring~$B^*$.

 We will apply the theory developed in
Sections~\ref{finitely-filtered-rings-and-modules-subsecn}\+-%
\ref{very-weak-relatively-adjusted-subsecn} to $\boZ$\+graded
(rings and) modules endowed with filtrations by homogenenous
ideals/submodules.
 So our modules $M^*$ are $\boZ$\+graded, our filtered modules $(M^*,F)$
carry \emph{both} a grading and a filtration, and our associated graded
modules $\gr_F^*M^*$ are bigraded.

 As usual (per the discussion in Section~\ref{graded-modules-subsecn}),
the results concerning ungraded rings and modules remain true in
the graded world.
 In particular, this applies to
Lemmas~\ref{very-weakly-relatively-projective-lemma}
and~\ref{very-weakly-relatively-injective-lemma}, which remain valid
in the context of $\boZ$\+graded $B^*$\+modules~$M^*$.
 Notice that we have $B^*/F^1B^*=R$, and the graded $R$\+modules
$F^iB^*/F^{i+1}B^*=B^i$ are projective (and finitely generated)
for all $i\ge0$.
 So one can speak of \emph{very weakly relative projective} and
\emph{very weakly relatively injective} graded $B^*$\+modules.

\begin{lem} \label{thick-graded-modules-lemma}
 Let $R$ be a commutative ring and $E$ be a finitely generated
projective $R$\+module.
 Consider the graded ring $B^*=\bigwedge_R^*(E)$, endowed with
with the finite decreasing filtration $F$ defined above.
 Then a graded $B^*$\+module $M^*$ is very weakly relatively projective
if and only if it is very weakly relatively injective.
\end{lem}

\begin{proof}
 Essentially, the assertion of the lemma holds because the exterior
algebra over a field is Frobenius, and the exterior algebra of a locally
free module of constant rank over a commutative ring is relatively
Frobenius in a suitable sense (cf.~\cite[Section~9.4]{Prel}).
 To be more precise, consider the affine scheme $\Spec R$, and let us
assign to every point $\p\in\Spec R$ the rank~$d$ of the free module
$E_\p$ over the local ring~$R_\p$.
 The resulting function (taking a finite number of integer values)
is locally constant on $\Spec R$ by~\cite[Lemma Tag~00NX]{SP}; so
it induces a decomposition of $R$ into a finite product of
commutative rings, $R=\prod_{d=0}^n R_d$, such that, in the related
direct sum decomposition $E=\bigoplus_{d=0}^n E_d$ of the $R$\+module
$E$, the $R_d$\+module $E_d$ is locally free of constant rank~$d$
(in the sense of~\cite[Definition Tag~00NW]{SP} and
Section~\ref{loc-free-sheaves-subsecn}).

 The graded ring $B^*$ decomposes accordingly into a finite direct
product $B^*=\prod_{d=0}^n B^*_d$, and the graded $B^*$\+module
$M^*$ decomposes into a finite direct sum $M^*=\bigoplus_{d=0}^n M^*_d$
of graded $B^*_d$\+modules~$M^*_d$.
 The filtration $F$ on the graded ring $B^*$ is the direct product of
the similar filtrations $F$ on the graded rings~$B^*_d$.
 Now one can easily see that the graded $B^*$\+module $M^*$ is very
weakly relatively projective if and only if the graded $B^*_d$\+module
$M^*_r$ is very weakly relatively projective  for every $0\le d\le n$;
and similarly for the very weak relative injectivity.
 This argument reduces the assertion of the lemma to the case when
the $R$\+module $E$ is finite locally free of constant rank~$d$.

 Now we have $\bigwedge^{d+1}_R(E)=0$, while the $R$\+module
$\bigwedge^d_R(E)$ is finite locally free of constant rank~$1$.
 So the pairing map $\bigwedge^d_R(E)\ot_R\bigwedge^d_R(\Hom_R(E,R))
\rarrow R$ from Lemma~\ref{exterior-dual-lemma} is an isomorphism
of $R$\+modules.
 Finally, we can conclude that, for any fixed finite decreasing
filtration $F$ on $M^*$ compatible with the filtration $F$ on $B^*$,
property~(5) from Lemma~\ref{very-weakly-relatively-projective-lemma}
is equivalent to property~(3) from
Lemma~\ref{very-weakly-relatively-injective-lemma}.
 Indeed, we have $\gr_F^*B^*=B^*$ and $B/F^1B=R$, and
$B^*\ot_RN^*\simeq\Hom_R(B^*,K^*)$ as graded $B^*$\+modules whenever
graded $R$\+modules $N^*$ and $K^*$ are connected by the rules
$K^*=\bigwedge^d_R(E)\ot_RN^*(-d)$ or $N^*=\bigwedge^d_R(\Hom_R(E,R))
\ot_RK^*(d)$.
 Here the notation $N^*\longmapsto N^*(-d)$ and $K^*\longmapsto K^*(d)$
stands for the grading shift.
\end{proof}

 The terminology ``thick modules'' in application to graded modules
over the exterior algebra was suggested in~\cite[Section~3.5]{Ryb}.
 The context in~\cite{Ryb} was that of the quasi-coherent graded algebra
of differential forms $\Omega^*_{X/k}$ for a smooth algebraic variety
$X$ over an algebraically closed field~$k$.
 Hence, in our notation, the commutative rings $R$ in the context
of~\cite{Ryb} were regular of finite Krull dimension, so they had
finite global dimension.
 Our context in the present paper is more general in that the global
dimension of the commutative ring $R$ may well be infinite.

 Nevertheless, the terminology of~\cite{Ryb} is convenient, and we
adopt it to our more general context.
 So we will say that a $\boZ$\+graded module $M^*$ over
the graded ring $B^*=\bigwedge_R^*(E)$ is \emph{thick} if it satisfies
the equivalent conditions of Lemma~\ref{thick-graded-modules-lemma}.

 The reader can notice the difference between our lists of equivalent
conditions in Lemmas~\ref{very-weakly-relatively-projective-lemma}
and~\ref{very-weakly-relatively-injective-lemma}, on the one hand,
and the list of equivelent conditions in~\cite[Theorem~3.5]{Ryb} on
the other hand, in that in our two lemmas \emph{all the $B$\+modules
$G$, $J$, and $P$ are annihilated by~$F^1B$}.
 In this sense, the conditions in~\cite[Theorem~3.5]{Ryb} are closer
to the context of \emph{weak} relative projectivity and \emph{weak} 
relative injectivity as in~\cite[Examples~2.10 and~3.9]{Pctrl}, rather
than our \emph{very weak} relative projectivity/injectivity.
 In the case of a commutative ring $R$ of finite global dimension,
our definition of a thick graded module is equivalent to the one
in~\cite{Ryb}.

\begin{lem} \label{thick-graded-modules-closed-under-co-kernels}
 Let $R$ be a commutative ring, $E$ be a finitely generated
projective $R$\+module, and $B^*$ be the graded ring
$B^*=\bigwedge_R^*(E)$.
 Then the class of thick graded $B^*$\+modules is closed under
extensions, kernels of epimorphisms, cokernels of monomorphisms,
direct summands, infinite direct sums, and infinite products in
the abelian category of graded $B^*$\+modules $B^*\Modl$.
\end{lem}

\begin{proof}
 More generally, for any finitely filtered graded ring $(B^*,F)$,
the class of very weakly relatively projective graded $B^*$\+modules
is closed under extensions and infinite direct sums in $B^*\Modl$,
as one can see from (the graded version of) any one of
conditions~(1\+-4) of
Lemma~\ref{very-weakly-relatively-projective-lemma}.
 Any one of conditions~(2) or~(4) of the same lemma implies
that the class of very weakly relatively projective graded
$B^*$\+modules is also closed under kernels of epimorphisms.

 Dual-analogously, for any finitely filtered graded ring $(B^*,F)$,
the class of very weakly relatively injective graded $B^*$\+modules
is closed under extensions and infinite products in $B^*\Modl$,
as one can see from any one of conditions~(1\+-2) of
Lemma~\ref{very-weakly-relatively-injective-lemma}.
 Condition~(2) of the same lemma implies that the class of very
weakly relatively injective graded $B^*$\+modules is also closed
under cokernels of monomorphisms.
\end{proof}

 We will denote the full subcategory of thick graded modules in
the abelian category of graded $B^*$\+modules $B^*\Modl$ by
$B^*\Modl_\thk=B^*\Modl^\thk\subset B^*\Modl$.
 It is clear from
Lemma~\ref{thick-graded-modules-closed-under-co-kernels} that
the full subcategory $B^*\Modl_\thk$ inherits an exact category
structure from the abelian exact structure of $B^*\Modl$.

\begin{exs} \label{thick-graded-modules-examples}
 The following examples illustrate the concept of thick graded module
over the exterior algebra and its canonical filtration.
 Let $R$ be a commutative ring and $E=Re$ be the free $R$\+module with
one generator~$e$.
 Then a graded module $M^*$ over the graded $R$\+algebra
$B^*=\bigwedge_R^*(E)$ is the same thing as a complex of $R$\+modules,
with the differential provided by the action of the exterior algebra
generator $e\in B^1$.
 Now the filtration $F$ on $M^*$ constructed as in the proof of
Lemma~\ref{very-weakly-relatively-projective-lemma}\,%
(1)\,$\Rightarrow$\,(5) has the form $F^0M^*=M^*$, \ $F^1M^*=eM^*$,
and $F^2M^*=0$.
 The filtration $F$ on $M^*$ constructed as in the proof of
Lemma~\ref{very-weakly-relatively-injective-lemma}\,%
(1)\,$\Rightarrow$\,(3) has the form $F^{-1}M^*=M^*$, \
$F^0M^*=\ker(e\:M^*\to M^*)$, and $F^1M^*=0$.
 A graded $B^*$\+module $M^*$ is thick if and only if the corresponding
complex of graded $R$\+modules is acyclic.
 In this case, the two filtrations on $M^*$ agree up to a shift of
the filtration degree, and the associated graded module $\gr_F^*M^*$
has the form $\gr_F^*M^*=M^*/eM^*\oplus eM^*$.

 For a more complicated example, one can consider the free $R$\+module
$E$ with two generators, $E=Re_1+Re_2$.
 Then a graded module $M^*$ over the graded ring
$B^*=\bigwedge_R^*(E)$ is a graded $R$\+module with two anti-commuting
differentials of degree~$1$ with zero squares, etc.
 The notion of a thick graded module over $B^*$ captures the idea of
this kind of polycomplex which is ``maximally nondegenerate'' with
respect to the differentials (while its underlying graded $R$\+module
may be rather arbitrary).
\end{exs}

\subsection{Change-of-scalar properties of thick graded modules}
 We discuss the restriction, extension, and coextension of scalars
before passing to the locality and colocality properties.

\begin{lem} \label{thick-graded-modules-restriction-of-scalars}
 Let $R\rarrow S$ be a homomorphism of commutative rings, $E$ be
a finitely generated projective $R$\+module, $A^*=\bigwedge_R^*(E)$
be the related exterior algebra over $R$, and $B^*=S\ot_RA^*=
\bigwedge_S^*(S\ot_RE)$ be the related exterior algebra over~$S$.
 Then a graded $B^*$\+module $N^*$ is thick if and only if $N^*$ is
thick as a graded $A^*$\+module.
\end{lem}

\begin{proof}
 It is convenient to use condition~(1) from
Lemma~\ref{very-weakly-relatively-projective-lemma} or condition~(1)
from Lemma~\ref{very-weakly-relatively-injective-lemma}.
 Consider the graded $A^*$\+module $R=A^*/F^1A^*$ and the graded
$B^*$\+module $S=B^*/F^1B^*$.
 Then any graded free resolution $P^*_\bu$ of the graded $A^*$\+module
$R$ is exact in the exact category of $R$\+projective graded
$A^*$\+modules, so $S\ot_RP^*_\bu$ is a graded free resolution of
the graded $B^*$\+module~$S$.
 Furthermore, there is a natural isomorphism of graded $B^*$\+modules
$S\ot_RP^*\simeq B^*\ot_{A^*}P^*$ for any graded $A^*$\+module~$P^*$.
 Hence we have $\Tor_i^{A^*}(R,N^*)\simeq\Tor_i^{B^*}(S,N^*)$ and
$\Ext^i_{A^*}(R,N^*)\simeq\Ext^i_{B^*}(S,N^*)$ for all $i\ge0$, and
the desired assertion follows.
\end{proof}

\begin{lem} \label{thick-graded-modules-restr-of-scal-filtration}
 Let $R\rarrow S$ be a homomorphism of commutative rings, $E$ be
a finitely generated projective $R$\+module, $A^*=\bigwedge_R^*(E)$
be the related exterior algebra over $R$, and $B^*=S\ot_RA^*=
\bigwedge_S^*(S\ot_RE)$ be the related exterior algebra over~$S$.
 Let $N^*$ be a thick graded $B^*$\+module.
 In this setting: \par
\textup{(a)} the canonical filtrations from
Lemma~\ref{vwr-projective-canonical-filtration-lemma} on the very weakly
relatively projective graded $B^*$\+module $N^*$ and the very weakly
relatively projective graded $A^*$\+module $N^*$ coincide; \par
\textup{(b)} the canonical filtrations from
Lemma~\ref{vwr-injective-canonical-filtration-lemma} on the very weakly
relatively injective graded $B^*$\+module $N^*$ and the very weakly
relatively injective graded $A^*$\+module $N^*$ coincide.
\end{lem}

\begin{proof}
 In part~(a), the point is that an isomorphism of
bigraded $\gr_F^*B^*$\+modules $\gr_F^*N^*\simeq\gr_F^*B^*\ot_SL^*$
for a graded $S$\+module $L^*$ implies an isomorphism of
bigraded $\gr_F^*A^*$\+modules $\gr_F^*N^*\simeq\gr_F^*A^*\ot_RL^*$.
 Similarly, in part~(b), an isomorphism of graded $\gr_F^*B^*$\+modules
$\gr_F^*N^*\simeq\Hom_S(\gr_F^*B^*,L^*)$ for a graded $S$\+module $L^*$
implies an isomorphism of bigraded $\gr_F^*A^*$\+modules
$\gr_F^*N^*\simeq\Hom_R(\gr_F^*A^*,L^*)$.
\end{proof}

\begin{lem} \label{thick-graded-modules-co-extension-of-scalars}
 Let $R\rarrow S$ be a homomorphism of commutative rings, $E$ be
a finitely generated projective $R$\+module, $A^*=\bigwedge_R^*(E)$
be the related exterior algebra over $R$, and $B^*=S\ot_RA^*=
\bigwedge_S^*(S\ot_RE)$ be the related exterior algebra over~$S$.
 Let $M^*$ be a thick graded $A^*$\+module.
 In this setting: \par
\textup{(a)} if the flat dimension of the $R$\+module $S$ does not
exceed\/~$1$ and\/ $\Tor_1^R(S,M^n)=0$ for all $n\in\boZ$, then
$S\ot_RM^*=B^*\ot_{A^*}M^*$ is a thick graded $B^*$\+module; \par
\textup{(b)} if the projective dimension of the $R$\+module $S$ does
not exceed\/~$1$ and\/ $\Ext^1_R(S,M^n)=0$ for all $n\in\boZ$, then\/
$\Hom_R^*(S,M^*)=\Hom_{A^*}^*(B^*,M^*)$ is a thick graded $B^*$\+module.
\end{lem}

\begin{proof}
 It is convenient to use condition~(5) from
Lemma~\ref{very-weakly-relatively-projective-lemma} or condition~(3)
from Lemma~\ref{very-weakly-relatively-injective-lemma}.
 Let us prove part~(b).
 Arguing as in the proof of Lemma~\ref{thick-graded-modules-lemma},
we can assume that the $R$\+module $E$ is locally free of constant
rank~$d$.
 By the graded version of
Lemma~\ref{very-weakly-relatively-injective-lemma}(3), we have
a finite decreasing filtration $F$ on the graded $A^*$\+module $M^*$
compatible with the natural decreasing filtration $F$ on $A^*$ such
that $F^{-d}M^*=M^*$, \ $F^1M^*=0$, and $F^{-i}M^n/F^{-i+1}M^n\simeq
\Hom_R(A^i,F^0M^{i+n})$ for all $0\le i\le d$ and $n\in\boZ$.

 Now the assumptions of~(b) imply that $\Ext^1_R(S,M^n/F^{-d+1}M^n)=0$
for all $n\in\boZ$, as the class of all $R$\+modules $C$ satisfying
$\Ext^1_R(S,C)=0$ is closed under quotients.
 So we have $\Ext^1_R(S,\Hom_R(A^d,F^0M^{d+n}))=0$ for all $n\in\boZ$.
 As $A^d$ is a locally free $R$\+module of constant rank~$1$, it
follows that $\Ext^1_R(S,F^0M^n)=0$ for all $n\in\boZ$.
 As $A^i$ is a finitely generated projective $R$\+module for every~$i$,
we can conclude that $\Ext^1_R(S,F^{-i}M^n/F^{-i+1}M^n)=0$ for
all $i$, $n\in\boZ$.
 Therefore, the rule $F^{-i}\Hom_R(S,M^n)=\Hom_R(S,F^{-i}M^n)$ defines
a finite decreasing filtration $F$ satisfying the condition
of Lemma~\ref{very-weakly-relatively-injective-lemma}(3) for
the graded $B^*$\+module $\Hom_R^*(S,M^*)$.

 The proof of part~(a) is dual-analogous.
\end{proof}

\begin{cor} \label{canonical-filtration-qcoh-lcth-compatible}
 Let $R\rarrow S$ be a homomorphism of commutative rings such that
the related morphism of affine schemes\/ $\Spec S\rarrow\Spec R$ is
an open immersion.
 Let $E$ be a finitely generated projective $R$\+module,
$A^*=\bigwedge_R^*(E)$ be the related exterior algebra over $R$, and
$B^*=S\ot_RA^*=\bigwedge_S^*(S\ot_RE)$ be the related exterior algebra
over~$S$.
 Let $M^*$ be a thick graded $A^*$\+module.
 In this setting: \par
\textup{(a)} Denoting by $F$ the canonical filtrations from
Lemma~\ref{vwr-projective-canonical-filtration-lemma} on the very weakly
relatively projective graded $A^*$\+module $M^*$ and the very weakly
relatively projective graded $B^*$\+module $N^*=S\ot_RM^*$, one has
$F^iN^*=S\ot_RF^iM^*$ and $\gr_F^iN^*=S\ot_R\gr_F^iM^*$ for
all $i\in\boZ$. \par
\textup{(b)} Assume that\/ $\Ext^1_R(S,M^n)=0$ for all $n\in\boZ$.
 Then, denoting by $F$ the canonical filtrations from
Lemma~\ref{vwr-injective-canonical-filtration-lemma} on the very weakly
relatively injective graded $A^*$\+module $M^*$ and the very weakly
relatively injective graded $B^*$\+module $N^*=\Hom_R^*(S,M^*)$, one
has $F^iN^*=\Hom_R^*(S,F^iM^*)$ and $\gr_F^iN^*=\Hom_R^*(S,\gr_F^iM^*)$
for all $i\in\boZ$.
\end{cor}

\begin{proof}
 This is a straightforward corollary of the proof of
Lemma~\ref{thick-graded-modules-co-extension-of-scalars}.
\end{proof}

\begin{lem} \label{thick-graded-modules-locality}
 Let $R\rarrow S_\alpha$, \ $1\le\alpha\le N$, be a finite collection
of homomorphisms of commutative rings such that the collection of
induced maps of the spectra\/ $\Spec S_\alpha\rarrow\Spec R$ is
an affine open covering of an affine scheme.
 Let $E$ be a finitely generated projective $R$\+module,
$A^*=\bigwedge_R^*(E)$ be the related exterior algebra over $R$,
and $B^*_\alpha=S_\alpha\ot_RA^*=\bigwedge_{S_\alpha}^*(S_\alpha\ot_RE)$
be the related exterior algebras over~$S_{\alpha}$.
 Then a graded $A^*$\+module $M^*$ is thick if and only if the graded
$B^*_\alpha$\+module $S_\alpha\ot_RM^*$ is thick for every index\/
$1\le\alpha\le N$.
\end{lem}

\begin{proof}
 The ``only if'' implication holds by
Lemma~\ref{thick-graded-modules-co-extension-of-scalars}(a).
 To prove the ``if'', consider the \v Cech
coresolution~\eqref{module-over-quasi-algebra-cech-coresolution}
\begin{multline} \label{module-over-exterior-algebra-cech-coresolution}
 0\lrarrow M^*\lrarrow\bigoplus\nolimits_{\alpha=1}^N S_\alpha\ot_RM^*
 \lrarrow\bigoplus\nolimits_{1\le\alpha<\beta\le N}
 S_\alpha\ot_RS_\beta\ot_RM^* \\
 \lrarrow\dotsb\lrarrow S_1\ot_R S_2\ot_R\dotsb\ot_R S_N\ot_RM^*
 \lrarrow0.
\end{multline}
 By Lemmas~\ref{thick-graded-modules-restriction-of-scalars}
and~\ref{thick-graded-modules-co-extension-of-scalars}(a), all
the terms of~\eqref{module-over-exterior-algebra-cech-coresolution},
except perhaps the leftmost one, are thick graded $A^*$\+modules.
 By Lemma~\ref{thick-graded-modules-closed-under-co-kernels}, it follows
that the leftmost term $M^*$ is a thick graded $A^*$\+module, too.
\end{proof}

\begin{lem} \label{thick-graded-modules-colocality}
 Let $R\rarrow S_\alpha$, \ $1\le\alpha\le N$, be a finite collection
of homomorphisms of commutative rings such that the collection of
induced maps of the spectra\/ $\Spec S_\alpha\rarrow\Spec R$ is
an affine open covering of an affine scheme.
 Let $E$ be a finitely generated projective $R$\+module,
$A^*=\bigwedge_R^*(E)$ be the related exterior algebra over $R$,
and $B^*_\alpha=S_\alpha\ot_RA^*=\bigwedge_{S_\alpha}^*(S_\alpha\ot_RE)$
be the related exterior algebras over~$S_{\alpha}$.
 Let $P^*$ be a graded $A^*$\+module such that the $R$\+module $P^n$
is contraadjusted for every $n\in\boZ$.
 Then the graded $A^*$\+module $P^*$ is thick if and only if the graded
$B^*_\alpha$\+module\/ $\Hom_R(S_\alpha,P^*)$ is thick for
every index\/ $1\le\alpha\le N$.
\end{lem}

\begin{proof}
 The ``only if'' implication holds by
Lemma~\ref{thick-graded-modules-co-extension-of-scalars}(b).
 To prove the ``if'', consider the \v Cech
resolution~\eqref{module-over-quasi-algebra-cech-resolution}
\begin{multline} \label{module-over-exterior-algebra-cech-resolution}
 0\lrarrow\Hom_R(S_1\ot_R\dotsb\ot_RS_N,\>P^*)\lrarrow\dotsb \\
 \lrarrow\bigoplus\nolimits_{\alpha<\beta}\Hom_R(S_\alpha\ot_R
 S_\beta,\>P^*)\lrarrow\bigoplus\nolimits_\alpha\Hom_R(S_\alpha,P^*)
 \lrarrow P^*\lrarrow0.
\end{multline}
 The complex~\eqref{module-over-exterior-algebra-cech-resolution}
is exact since all the grading components of $P^*$ are contraadjusted
$R$\+modules by assumption~\cite[formula~(1.3) in
Lemma~1.2.6(b)]{Pcosh}.
 By Lemmas~\ref{thick-graded-modules-restriction-of-scalars}
and~\ref{thick-graded-modules-co-extension-of-scalars}(b), all
the terms of~\eqref{module-over-exterior-algebra-cech-resolution},
except perhaps the rightmost one, are thick graded $A^*$\+modules.
 By Lemma~\ref{thick-graded-modules-closed-under-co-kernels}, it follows
that the rightmost term $P^*$ is a thick graded $A^*$\+module, too.
\end{proof}

\subsection{Thick quasi-coherent graded modules}
\label{thick-quasi-coherent-modules-subsecn}
 Let $X$ be a scheme and $\E$ be a finite locally free sheaf on~$X$.
 Consider the exterior algebra $\cB^*=\bigwedge^*_X(\E)$ of
the quasi-coherent sheaf $\E$, as defined at the end of
Section~\ref{skew-symmetric-tensors-subsecn}.
 So $\cB^*$ is a quasi-coherent graded algebra over~$X$.

 A quasi-coherent graded $\cB^*$\+module $\M^*$ is said to be
\emph{thick} if, for every affine open subscheme $U\subset X$,
the graded $\cB^*(U)$\+module $\M^*(U)$ is thick (in the sense of
Section~\ref{thick-graded-modules-subsecn}).
 In view of Lemma~\ref{thick-graded-modules-locality}, it suffices to
check this condition for affine open subschemes $U$ belonging to any
given affine open covering of the scheme~$X$.

 It is clear from
Lemma~\ref{thick-graded-modules-closed-under-co-kernels} that
the full subcategory of thick quasi-coherent graded modules
$\cB^*\Qcoh_\thk\subset\cB^*\Qcoh$ is closed under extensions,
kernels of epimorphisms, cokernels of monomorphisms, direct summands,
and infinite direct sums in the abelian category $\cB^*\Qcoh$.
 So the full subcategory $\cB^*\Qcoh_\thk$ inherits an exact category
structure from the abelian exact structure of $\cB^*\Qcoh$.

 The following two lemmas are concerned with the direct and inverse
images of thick quasi-coherent graded modules.
 The discussion of the direct and inverse images of quasi-coherent
$\cB$\+modules and $\cA$\+modules in
Sections~\ref{direct-images-of-A-co-sheaves-subsecn}\+-%
\ref{inverse-images-of-A-co-sheaves-subsecn} presumed the context of
a morphism of quasi-ringed schemes $f\:(Y,\cB)\rarrow(X,\cA)$.
 Recall from that discussion that the direct image functor
$f_*\:\cB\Qcoh\rarrow\cA\Qcoh$, by construction, always agrees with
the direct image functor $f_*\:Y\Qcoh\rarrow X\Qcoh$, while
the inverse image functor $f^*\:\cA\Qcoh\rarrow\cB\Qcoh$, generally
speaking, does \emph{not} agree with the inverse image functor
$f^*\:X\Qcoh\rarrow Y\Qcoh$.

 In the context of the next two lemmas, however, the inverse image 
functor $f^*\:\cA^*\Qcoh\rarrow\cB^*\Qcoh$ \emph{agrees} with
the inverse image functor $f^*\:X\Qcoh\rarrow Y\Qcoh$.
 This happens for two reasons.
 Firstly, $\cA^*$ is a quasi-coherent graded algebra (rather than
merely a quasi-coherent graded quasi-algebra) over~$X$.
 Consequently, there is a natural structure of quasi-coherent graded
algebra on the inverse image $f^*\cA^*$ of the graded quasi-coherent
sheaf $\cA^*$ on $X$ to~$Y$.
 Secondly, the morphism of quasi-coherent graded algebras
$f^*\cA^*\rarrow\cB^*$ on~$Y$ (corresponding by adjunction to
the morphism of sheaves of graded rings $\cA^*\rarrow f_*\cB^*$
on~$X$) is an isomorphism.

\begin{lem} \label{qcoh-thick-direct-image}
 Let $f\:Y\rarrow X$ be an affine morphism of schemes, $\E$ be a finite
locally free sheaf on $X$, and $f^*\E$ be the inverse image of $\E$
on~$Y$.
 Consider the related exterior algebras $\cA^*=\bigwedge^*_X(\E)$
on $X$ and $\cB^*=\bigwedge^*_Y(f^*\E)=f^*\cA^*$ on~$Y$.
 Then a quasi-coherent graded $\cB^*$\+module $\N^*$ on $Y$ is thick
if and only if the quasi-coherent graded $\cA^*$\+module $f_*\N^*$
on $X$ is thick.
\end{lem}

\begin{proof}
 Follows from
Lemma~\ref{thick-graded-modules-restriction-of-scalars}.
\end{proof}

\begin{lem} \label{qcoh-thick-inverse-image}
 Let $f\:Y\rarrow X$ be a flat morphism of schemes, $\E$ be a finite
locally free sheaf on $X$, and $f^*\E$ be the inverse image of $\E$
on~$Y$.
 Consider the related exterior algebras $\cA^*=\bigwedge^*_X(\E)$
on $X$ and $\cB^*=\bigwedge^*_Y(f^*\E)=f^*\cA^*$ on~$Y$.
 Then, for any thick quasi-coherent graded $\cA^*$\+module $\M^*$ on
$X$, the quasi-coherent graded $\cB^*$\+module $f^*\M$ on $Y$ is thick.
\end{lem}

\begin{proof}
 Follows from
Lemma~\ref{thick-graded-modules-co-extension-of-scalars}(a).
\end{proof}

\subsection{Thick locally contraherent graded modules}
\label{thick-loc-contraherent-modules-subsecn}
 Let $X$ be a scheme with an open covering $\bW$ and $\E$ be a finite
locally free sheaf on~$X$.
 Put $\cB^*=\bigwedge_X^*(\E)$.

 Let $\P^*$ be a $\bW$\+locally contraherent graded $\cB^*$\+module
on~$X$ (as defined in Sections~\ref{cosheaves-of-A-modules-subsecn}
and~\ref{graded-modules-subsecn}).
 A $\bW$\+locally contraherent graded $\cB^*$\+module $\P^*$ is said
to be \emph{thick} if, for every affine open subscheme $U\subset X$
subordinate to $\bW$, the graded $\cB^*(U)$\+module $\P^*[U]$ is thick
(in the sense of Section~\ref{thick-graded-modules-subsecn}).
 In view of Lemma~\ref{thick-graded-modules-colocality}, it suffices to
check this condition for affine open subschemes $U$ belonging to any
chosen affine open covering of the scheme $X$ subordinate to~$\bW$.
 In particular, the property of a $\bW$\+locally contraherent graded
$\cB^*$\+module to be thick is not changed by refinements of
the open covering~$\bW$.

 It is clear from
Lemma~\ref{thick-graded-modules-closed-under-co-kernels} that
the full subcategory of thick $\bW$\+locally contraherent graded
modules $\cB^*\Lcth_\bW^\thk\subset\cB^*\Lcth_\bW$ is closed under 
extensions, kernels of admissible epimorphisms, cokernels of admissible 
monomorphisms, direct summands, and infinite products in the exact
category $\cB^*\Lcth_\bW$.
 So the full subcategory $\cB^*\Lcth_\bW^\thk$ inherits an exact
category structure from $\cB^*\Lcth_\bW$.

 The following two lemmas are concerned with the direct and inverse
images of thick locally contraherent graded modules.
 The discussion of the direct and inverse images of locally
contraherent $\cB$\+modules and $\cA$\+modules in
Sections~\ref{direct-images-of-A-co-sheaves-subsecn}\+-%
\ref{inverse-images-of-A-co-sheaves-subsecn} presumed the context of
a morphism of quasi-ringed schemes $f\:(Y,\cB)\rarrow(X,\cA)$.
 Recall from that discussion that the direct image functor
$f_!\:\cB\Lcth_\bT\rarrow\cA\Lcth_\bW$, by construction, always agrees
with the direct image functor $f_!\:Y\Lcth_\bT\rarrow X\Lcth_\bW$.
 However, the inverse image functor~$f^!$ acting from a suitable full
subcategory of $\cA\Lcth_\bW$ to the similar full subcategory of
$\cB\Lcth_\bT$, generally speaking, does \emph{not} agree with
the inverse image functor~$f^!$ acting from a full subcategory
of $X\Lcth_\bW$ to a full subcategory of $Y\Lcth_\bT$.
 Moreover, the relevant full subcategories of $\cA\Lcth_\bW$ mentioned
in Section~\ref{inverse-images-of-A-co-sheaves-subsecn} were
rather narrow, as inverse image $f^!\P$ of a locally contraherent
$\cA$\+module $\P$ on $X$ was only defined in
formulas~\eqref{ctrh-lct-A-modules-inverse-image}
or~\eqref{ctrh-lin-A-modules-inverse-image} under somewhat
restrictive assumptions on~$\P$.

 Nevertheless, in the context of Lemma~\ref{lcth-thick-inverse-image}
below in this section, the inverse image functor $f^!\:\cA^*\Lcth_\bW
\rarrow\cB^*\Lcth_\bT$ is well-defined on the whole category of
$\bW$\+locally contraherent graded $\cA^*$\+modules on $X$ and agrees
with the inverse image functor $f^!\:X\Lcth_\bW\rarrow Y\Lcth_\bT$
from Section~\ref{inverse-images-of-O-co-sheaves-subsecn}.
 This happens for the same reasons that were explained in
Section~\ref{thick-quasi-coherent-modules-subsecn}.

\begin{lem} \label{lcth-thick-direct-image}
 Let $X$ be a scheme with an open covering\/ $\bW$ and $Y$ be
a scheme with an open covering\/~$\bT$.
 Let $f\:Y\rarrow X$ be a $(\bW,\bT)$\+affine morphism of schemes,
$\E$ be a finite locally free sheaf on $X$, and $f^*\E$ be
the inverse image of $\E$ on~$Y$.
 Consider the related exterior algebras $\cA^*=\bigwedge^*_X(\E)$
on $X$ and $\cB^*=\bigwedge^*_Y(f^*\E)=f^*\cA^*$ on~$Y$.
 Then a\/ $\bT$\+locally contraherent graded $\cB^*$\+module\/ $\Q^*$
on $Y$ is thick if and only if the\/ $\bW$\+locally contraherent
graded $\cA^*$\+module $f_!\Q^*$ on $X$ is thick.
\end{lem}

\begin{proof}
 Follows from
Lemma~\ref{thick-graded-modules-restriction-of-scalars}.
\end{proof}

\begin{lem} \label{lcth-thick-inverse-image}
 Let $X$ be a scheme with an open covering\/ $\bW$ and $Y$ be
a scheme with an open covering\/~$\bT$.
 Let $f\:Y\rarrow X$ be a very flat $(\bW,\bT)$\+coaffine
morphism of schemes, $\E$ be a finite locally free sheaf on $X$,
and $f^*\E$ be the inverse image of $\E$ on~$Y$.
 Consider the related exterior algebras $\cA^*=\bigwedge^*_X(\E)$
on $X$ and $\cB^*=\bigwedge^*_Y(f^*\E)=f^*\cA^*$ on~$Y$.
 Then, for any thick\/ $\bW$\+locally contraherent graded
$\cA^*$\+module\/ $\P^*$ on $X$, the\/ $\bT$\+locally contraherent
graded $\cB^*$\+module $f^!\P^*$ on $Y$ is thick.
\end{lem}

\begin{proof}
 Follows from
Lemma~\ref{thick-graded-modules-co-extension-of-scalars}(b).
\end{proof}

\subsection{Cotorsion modules over a finite locally free quasi-algebra}
 Since Sections~\ref{cosheaves-of-A-modules-subsecn}\+-%
\ref{A-loc-cotors-loc-inj-cosheaves-subsecn}
and~\ref{graded-modules-subsecn}, we had to deal with the distinction
between $X$\+locally cotorsion and $\cA^*$\+locally cotorsion locally
contraherent graded modules for a quasi-coherent graded quasi-algebra
$\cA^*$ over a scheme~$X$. 
 Since Sections~\ref{antilocality-of-X-cotorsion-subsecn}\+-%
\ref{antilocality-of-A-cotorsion-subsecn}
and~\ref{exact-dg-of-cta-cot-qcoh-cdg}, we also dealt with
the distinction between $X$\+cotorsion and $\cA^*$\+cotorsion
quasi-coherent graded $\cA^*$\+modules.
 In this section we will see that, for quasi-coherent graded
quasi-algebra $\cA^*$ that is a finite locally free sheaf on $X$
with respect to its right $\cO_X$\+module structure,
there is \emph{no difference} between the two versions of
the cotorsion or locally cotorsion concept.
 The discussion is based on
Proposition~\ref{relative-bass-theorem-cotorsion-prop}.

 We will say that a graded quasi-algebra $A^*$ over a commutative ring
$R$ is \emph{finitely generated projective on the right} if $A^*$ has
only a finite number of nonzero grading components, i.~e., $A^n=0$ for
all but a finite set of integers~$n$, and the grading component $A^n$
is a finitely generated projective $R$\+module with respect to its
right $R$\+module structure for all $n\in\boZ$.
 A quasi-coherent graded quasi-algebra $\cA^*$ over a scheme $X$ is
said to be \emph{finite locally free on the right} if the graded
quasi-algebra $\cA^*(U)$ over the commutative ring $U$ is finitely
generated projective on the right for all affine open subschemes
$U\subset X$.
 It suffices to check this condition for affine open subschemes $U$
belonging to any given affine open covering of the scheme~$X$.

\begin{prop} \label{cotorsion-modules-over-fgpgr-ring-prop}
 Let $R$ be an associative ring and $R\rarrow A^*$ be a homomorphism
of graded associative rings such that $A^n=0$ for all but a finite
set of degrees~$n$ and $A^n$ is a finitely generated projective right
$R$\+module for all $n\in\boZ$.
 Then a graded left $A^*$\+module $C^*$ is cotorsion if and only if its
grading components $C^n$ are cotorsion $R$\+modules for all $n\in\boZ$.
\end{prop}

\begin{proof}
 This is a particular case of the graded version of
Proposition~\ref{relative-bass-theorem-cotorsion-prop}.
\end{proof}

\begin{cor} \label{lct-lcth-modules-over-flfrqa-cor}
 Let $X$ be a scheme with an open covering\/ $\bW$ and $\cA^*$ be
a quasi-coherent graded quasi-algebra over $X$ that that $\cA^*$ is
finite locally free on the right.
 Then the classes of $X$\+locally cotorsion and $\cA^*$\+locally
cotorsion $\bW$\+locally contraherent graded $\cA^*$\+modules on $X$
coincide.
\end{cor}

\begin{proof}
 The assertion follows immediately from the definitions
(see Sections~\ref{cosheaves-of-A-modules-subsecn}\+-%
\ref{A-loc-cotors-loc-inj-cosheaves-subsecn}
and~\ref{graded-modules-subsecn}) and
Proposition~\ref{cotorsion-modules-over-fgpgr-ring-prop}.
\end{proof}

\begin{cor} \label{cot-qcoh-modules-over-flfrqa-cor}
 Let $X$ be a quasi-compact semi-separated scheme and $\cA^*$ be
a quasi-coherent graded quasi-algebra over $X$ that that $\cA^*$ is
finite locally free on the right.
 Then the classes of $X$\+cotorsion and $\cA^*$\+cotorsion
quasi-coherent graded $\cA^*$\+modules on $X$ coincide.
\end{cor}

\begin{proof}
 This is provable by a suitable version of the argument
from~\cite[Section~2]{PBas}, extended from the realm of rings and
modules to the realm of schemes, quasi-coherent quasi-algebras and
quasi-coherent modules, and based on
Theorem~\ref{qcoh-cotorsion-periodicity}.
 Alternatively, the following approach allows to reduce
the quasi-compact semi-separated scheme case to the case of
an affine scheme.

 We use the idea of the argument from~\cite[Section~8]{Pphil} involving
the passage from quasi-coherent sheaves (for which the cotorsion
properties are \emph{not} local) to contraherent cosheaves (for which
the related locally cotorsion properties \emph{are} local).
 By the graded version of
Lemma~\ref{quasi-algebra-underived-naive-co-contra}(b\+-c)
(as per the discussion in Section~\ref{graded-modules-subsecn}),
for any quasi-coherent graded quasi-algebra $\cA^*$ over
a quasi-compact semi-separated scheme $X$ there is a commutative
square diagram of exact category equivalences and exact, fully faithful
inclusion functors
\begin{equation} \label{A-cot-A-lct-X-cot-X-lct-comparison}
\begin{gathered}
 \qquad\xymatrix{
  \text{\llap{$\fHom_{\cA^*}^*(\cA^*,{-})\:$}}\cA^*\Qcoh^{X\dcot}
  \ar@{=}[rr] 
  && \cA^*\Ctrh_\al^{X\dlct}
  \text{\rlap{$\,\,:\!\cA^*\ocn_{\cA^*}{-}$}} \\
  \text{\llap{$\fHom_{\cA^*}^*(\cA^*,{-})\:$}} \cA^*\Qcoh^{\cA^*\dcot}
  \ar@{=}[rr] \ar@{>->}[u]
  && \cA^*\Ctrh_\al^{\cA^*\dlct}
  \text{\rlap{$\,\,:\!\cA^*\ocn_{\cA^*}{-}$}} \ar@{>->}[u]
 }
\end{gathered}
\end{equation}
 In the situation at hand when the quasi-coherent graded quasi-algebra
$\cA^*$ is finite locally free on the right, the rightmost vertical
fully faithful functor on
the diagram~\eqref{A-cot-A-lct-X-cot-X-lct-comparison} is a category
equivalence by Corollary~\ref{lct-lcth-modules-over-flfrqa-cor},
and it follows that the leftmost vertical functor is a category
equivalence, too.
\end{proof}

\subsection{Filtered quasi-coherent and locally contraherent modules}
\label{filtered-qcoh-lcth-subsecn}
 The context of the exposition in this section is much more general
than in the preceding and subsequent ones.

 Let $\dotsb\supset F^{-i}\cA\supset\dotsb\supset F^0\cA\supset\dotsb
\supset F^i\cA\supset\dotsb$ be a filtered quasi-coherent
quasi-algebra over a scheme~$X$.
 We always presume the compatibility of the filtration with
the multiplicative structure; so $F^i\cA(U)\cdot F^j\cA(U)\subset
F^{i+j}\cA(U)$ and $1_{\cA(U)}\in F^0\cA(U)$ for all affine open
subschemes $U\subset X$ and all $i$, $j\in\boZ$.
 Here $1_{\cA(U)}$ denotes the unit element of the ring $\cA(U)$.
 In the context of this section, we are mostly interested in finite
filtrations ($F^{-n}\cA=\cA$ and $F^n\cA=0$ for $n$~large enough),
so we refrain from a discussion of the (co)limit conditions at
$i\to\pm\infty$.

 To a filtered quasi-coherent quasi-algebra $(\cA,F)$, one can assign
its quasi-coherent graded \emph{Rees quasi-algebra}
$\bigoplus_{i=-\infty}^\infty F^i\cA$.
 A \emph{filtered quasi-coherent module} $(\M,F)$ over $(\cA,F)$ is
the same thing as a quasi-coherent graded module
$\bigoplus_{i=-\infty}^\infty F^i\M$ over
$\bigoplus_{i=-\infty}^\infty F^i\cA$ satisfying the additional
condition that the natural morphisms of quasi-coherent sheaves
$F^{i+1}\M\rarrow F^i\M$ are injective.

 Here the morphisms $F^{i+1}\M\rarrow F^i\M$ are determined by
the structure of quasi-coherent graded module over
$\bigoplus_{i=-\infty}^\infty F^i\cA$ in the following way.
 For every affine open subscheme $U\subset X$, the unit element
$1_{\cA(U)}\in F^0\cA(U)$ can be also viewed as an element of
$F^{-1}\cA(U)\supset F^0\cA(U)$.
 We denote this element by $t(U)\in F^{-1}\cA(U)$.
 Clearly, $t(U)$ is a central element in the Rees ring
$\bigoplus_{i=-\infty}^\infty F^i\cA(U)$.
 The element~$t$ is a natural global section of the quasi-coherent
sheaf $F^{-1}\cA$ on $X$, i.~e., $t\in F^{-1}\cA(X)$.
 The natural morphisms of quasi-coherent sheaves
$F^{i+1}\M\rarrow F^i\M$ are provided by the action of
the element $t\in F^{-1}\cA(X)$.

 Now let $\bW$ be an open covering of the scheme~$X$.
 A \emph{filtered\/ $\bW$\+locally contraherent module} $(\P,F)$ over
$(\cA,F)$ can be defined as a $\bW$\+locally contraherent graded module
$\bigoplus_{i=-\infty}^\infty F^i\P$ over the quasi-coherent graded
quasi-algebra $\bigoplus_{i=-\infty}^\infty F^i\cA$ over~$X$ satisfying
the following condition.
 The action of the element $t\in F^{-1}\cA$ on $\P$ provides
natural morphisms of $\bW$\+locally contraherent cosheaves
$F^{i+1}\P\rarrow F^i\P$.
 These morphisms must be admissible monomorphisms in $X\Lcth_\bW$
for all $i\in\boZ$.

 \emph{$X$\+locally cotorsion filtered\/ $\bW$\+locally contraherent
modules} $(\P,F)$ over $(\cA,F)$ are defined similarly; so,
the $\bW$\+locally contraherent cosheaves $F^i\P$ on $X$ are
assumed to be locally cotorsion in this case.
 By a \emph{filtered contraherent module} over $(\cA,F)$ we mean
a filtered $\bW$\+locally contraherent module for the filtration
$\bW=\{X\}$.

 Notice that a morphism in $X\Lcth_\bW^\lct$ is an admissible
monomorphism if and only if it is an admissible monomorphism
in $X\Lcth_\bW$.
 However, given a refinement $\bW'$ of the open covering $\bW$,
a morphism in $X\Lcth_\bW$ can be an admissible monomorphism in
$X\Lcth_{\bW'}$ without being an admissible monomorphism in
$X\Lcth_\bW$ (see Section~\ref{exact-categories-of-contrah-subsecn}).
 In fact, a morphism $\P\rarrow\Q$ is an admissible monomorphism
in $X\Lcth_\bW$ or $X\Lcth_\bW^\lct$ if and only if the map
$\P[U]\rarrow\Q[U]$ is injective for all affine open subschemes
$U\subset X$ subordinate to~$\bW$.
 So one has to be careful.

 Assume that the scheme $X$ is quasi-compact and semi-separated.
 A filtered $\bW$\+locally contraherent module $(\P,F)$ over
$(\cA,F)$ is said to be \emph{antilocal} if the $\bW$\+locally
contraherent cosheaves $F^i\P$ are antilocal for all $i\in\boZ$
\emph{and} the natural morphisms $F^{i+1}\P\rarrow F^i\P$ are
admissible monomorphisms in the exact category of antilocal
($\bW$\+locally) contraherent cosheaves on~$X$.
 In other words, this means that the successive quotient
$\bW$\+locally contraherent cosheaves $F^i\P/F^{i+1}\P$ must be
antilocal.
 (Notice that the full subcategory of antilocal contraherent cosheaves
$X\Ctrh_\al$ is usually \emph{not} closed under cokernels of admissible
monomorphisms in $X\Ctrh$ or $X\Lcth_\bW$.)

 Similarly, a filtered quasi-coherent module $(\M,F)$ over $(\cA,F)$
is said to be \emph{$X$\+contraadjusted} if the quasi-coherent sheaf
$F^i\M$ on $X$ is contraadjusted for all $i\in\boZ$.
 A filtered quasi-coherent module $(\M,F)$ over $(\cA,F)$ is said to be
\emph{$X$\+cotorsion} if the quasi-coherent sheaf $F^i\M$ is
cotorsion for all $i\in\boZ$.
 Notice that the full subcategories of contraadjusted/cotorsion
quasi-coherent sheaves $X\Qcoh^\cta$ and $X\Qcoh^\cot$ are closed
under cokernels of admissible monomorphisms in $X\Qcoh$ (see
Sections~\ref{antilocality-of-X-contraadjusted-subsecn}\+-%
\ref{antilocality-of-X-cotorsion-subsecn}), so under the respective
condition above the quasi-coherent sheaves $F^i\M/F^{i+1}\M$ on $X$
are contraadjusted/cotorsion as well.

 For any filtered quasi-coherent quasi-algebra $(\cA,F)$, one can
consider the associated graded quasi-coherent quasi-algebra
$\gr_F^*\cA$ with the grading components $\gr_F^i\cA=F^i\cA/F^{i+1}\cA$.
 For any filtered quasi-coherent module $(\M,F)$ over $(\cA,F)$,
the associated graded quasi-coherent sheaf $\gr_F^*\M$ with
the grading components $\gr_F^i\M=F^i\M/F^{i+1}\M$ is a quasi-coherent
graded module over $\gr_F^*\cA$.
 For any filtered $\bW$\+locally contraherent module $(\P,F)$ over
$(\cA,F)$, the associated graded $\bW$\+locally contraherent cosheaf
$\gr_F^*\P$ with the grading components $\gr_F^i\P=F^i\P/F^{i+1}\P$
is a $\bW$\+locally contraherent graded module over $\gr_F^*\cA$.

 From this point on, we enforce the assumption that the quasi-coherent
quasi-algebra $\cA$ with the filtration $F$ is \emph{finitely filtered},
i.~e., there exist an integer $n>0$ for which $F^{-n}\cA=\cA$ and
$F^n\cA=0$.
 Similarly, a filtered quasi-coherent module $(\M,F)$ over $(\cA,F)$ is
said to be \emph{finitely filtered} if there exists an integer $m>0$
for which $F^{-m}\M=\M$ and $F^m\M=0$.
 A filtered locally contraherent module $(\P,F)$ over $(\cA,F)$ is
said to be \emph{finitely filtered} if there exists an integer $p>0$
for which $F^{-p}\M=\M$ and $F^p\M=0$.

 As usual, all \emph{mophisms} of filtered modules are presumed to
preserve the filtrations.
 A short sequence of filtered quasi-coherent modules $0\rarrow\cL
\rarrow\M\rarrow\N\rarrow0$ is said to be (\emph{admissible})
\emph{exact} if the short sequence $0\rarrow F^i\cL\rarrow F^i\M
\rarrow F^i\N\rarrow0$ is exact in $X\Qcoh$ for every $i\in\boZ$.
 A short sequence of filtered $\bW$\+locally contraherent modules
$0\rarrow\P\rarrow\Q\rarrow\R\rarrow0$ is said to be \emph{exact} if
the short sequence $0\rarrow F^i\P\rarrow F^i\Q\rarrow F^i\R\rarrow0$
is exact in $X\Lcth_\bW$ for every $i\in\boZ$.

 We denote the exact category of finitely filtered quasi-coherent
modules over $(\cA,F)$ by $(\cA,F)\Qcoh^\ff$ and its full exact
subcategories of $X$\+cotorsion and $X$\+contraadjusted finitely
filtered quasi-coherent modules by $(\cA,F)\Qcoh^{\ff,X\dcot}
\subset(\cA,F)\Qcoh^{\ff,X\dcta}\subset(\cA,F)\Qcoh^\ff$.
 Similarly, the exact category of finitely filtered $\bW$\+locally
contraherent modules over $(\cA,F)$ is denoted by
$(\cA,F)\Lcth_\bW^\ff$, and its full subcategory of $X$\+locally
cotorsion finitely filtered $\bW$\+locally contraherent modules
is denoted by $(\cA,F)\Lcth_\bW^{\ff,X\dlct}\subset
(\cA,F)\Lcth_\bW^\ff$.
 Finally, the full exact subcategories of antilocal ($X$\+locally
contraadjusted or $X$\+locally cotorsion) finitely filtered
contraherent modules over $(\cA,F)$ are denoted by
$(\cA,F)\Ctrh_\al^\ff\subset(\cA,F)\Lcth_\bW^\ff$ and
$(\cA,F)\Ctrh_\al^{\ff,X\dlct}\subset(\cA,F)\Lcth_\bW^{\ff,X\dlct}$.
{\hbadness=1150\par}

\begin{thm} \label{quasi-algebra-filtered-naive-co-contra}
 Let $X$ be a quasi-compact semi-separated scheme and $(\cA,F)$ be
a finitely filtered quasi-coherent quasi-algebra over~$X$.
 Then there exist commutative square diagrams of exact functors and
exact category equivalences
\begin{equation} \label{X-cta-filtered-naive-co-contra}
\begin{gathered}
 \qquad\xymatrix{
  \text{\llap{$\fHom_\cA(\cA,{-})\:$}}(\cA,F)\Qcoh^{\ff,X\dcta}
  \ar@{=}[rr] \ar[d]_{\gr_F^*}
  && (\cA,F)\Ctrh_\al^\ff
  \text{\rlap{$\,\,:\!\cA\ocn_\cA{-}$}} \ar[d]^{\gr_F^*} \\
  \text{\llap{$\fHom_{\gr_F^*\cA}(\gr_F^*\cA,{-})\:$}}
  \gr_F^*\cA\Qcoh^{X\dcta} \ar@{=}[rr] 
  && \gr_F^*\cA\Ctrh_\al
  \text{\rlap{$\,\,:\!\gr_F^*\cA\ocn_{\gr_F^*\cA}{-}$}}
 }
\end{gathered}
\end{equation}
\begin{equation} \label{X-cot-filtered-naive-co-contra}
\begin{gathered}
 \qquad\xymatrix{
  \text{\llap{$\fHom_\cA(\cA,{-})\:$}}(\cA,F)\Qcoh^{\ff,X\dcot}
  \ar@{=}[rr] \ar[d]_{\gr_F^*}
  && (\cA,F)\Ctrh_\al^{\ff,X\dlct}
  \text{\rlap{$\,\,:\!\cA\ocn_\cA{-}$}} \ar[d]^{\gr_F^*} \\
  \text{\llap{$\fHom_{\gr_F^*\cA}(\gr_F^*\cA,{-})\:$}}
  \gr_F^*\cA\Qcoh^{X\dcot} \ar@{=}[rr] 
  && \gr_F^*\cA\Ctrh_\al^{X\dlct}
  \text{\rlap{$\,\,:\!\gr_F^*\cA\ocn_{\gr_F^*\cA}{-}$}}
 }
\end{gathered}
\end{equation}
 Here the exact categories in the lower lines of the two diagrams are
the categories of quasi-coherent and contraherent graded modules
over the quasi-coherent graded algebra\/ $\gr_F^*\cA$, as defined
in Sections~\ref{graded-modules-subsecn}
and~\ref{exact-dg-of-cta-cot-qcoh-cdg}\+-%
\ref{exact-dg-categories-of-antilocal-ctrh-subsecn}.
 The exact category equivalences in the lower lines are the graded
versions of Lemma~\ref{quasi-algebra-underived-naive-co-contra}(a\+-b).
 The vertical exact functors are $(\M,F)\longmapsto\gr_F^*\M$
and $(\P,F)\longmapsto\gr_F^*\P$.
\end{thm}

\begin{proof}
 The simplest way to frame the argument is to apply the graded version
of Lemma~\ref{quasi-algebra-underived-naive-co-contra}(a\+-b) to
quasi-coherent and contraherent graded modules over the quasi-coherent
graded Rees quasi-algebra $\bigoplus_{i=-\infty}^\infty F^i\cA$.
 Let $(F^i\M)_{i\in\boZ}$ and $(F^i\P)_{i\in\boZ}$ be a quasi-coherent
and a contraherent graded module over
$\bigoplus_{i=-\infty}^\infty F^i\cA$ corresponding to each other under
the graded version of
Lemma~\ref{quasi-algebra-underived-naive-co-contra}.
 Then, for every $i\in\boZ$, the quasi-coherent sheaf $F^i\M$ and
the contraherent cosheaf $F^i\P$ correspond to each other under
the na\"\i ve co-contra correspondence
of~\cite[Lemma~4.8.2 or~4.8.4(a)]{Pcosh}.
 One needs to observe that the morphisms $F^{i+1}\M\rarrow F^i\M$
and $F^{i+1}\P\rarrow F^i\P$ provided by the action of the canonical
central element $t\in F^{-1}\cA(X)$ correspond to each other under
the na\"\i ve co-contra correspondence from~\cite[Section~4.8]{Pcosh}.
 Then it follows that the morphisms $F^{i+1}\M\rarrow F^i\M$ are
admissible monomorphisms in $X\Qcoh^\cta$ (respectively,
$X\Qcoh^\cot$) if and only if the morphisms $F^{i+1}\P\rarrow F^i\P$
are admissible monomorphisms in $X\Ctrh_\al$ (resp.,
$X\Ctrh_\al^{X\dlct}$).

 To prove commutativity of
the diagrams~(\ref{X-cta-filtered-naive-co-contra}\+-%
\ref{X-cot-filtered-naive-co-contra}), one can point out that,
over a semi-separated scheme $X$, the contratensor product
functor~$\ocn$ from
Section~\ref{contratensor-over-qcoh-quasi-algebra-subsecn} preserves
cokernels in all of its arguments, essentially by construction.
 So, in particular, the functor~$\ocn$ commutes with
the passage to the quotients by the action of the canonical central
element $t\in F^{-1}\cA(X)$.
\end{proof}

\subsection{Associated graded modules are contraadjusted/cotorsion}
\label{assoc-graded-modules-are-cta-cot}
 The aim of this section is to show that when a thick graded module
has contraadjusted/cotorsion grading components, then its associated
graded module with respect to the canonical filtration has
contraadjusted/cotorsion bigrading components.

\begin{lem} \label{thick-cta-module-associated-graded-lemma}
 Let $R$ be a commutative ring, $E$ be a finitely generated projective
$R$\+module, and $B^*$ be the graded algebra $B^*=\bigwedge_R^*(E)$
over~$R$.
 Let $C^*$ be a thick graded $B^*$\+module whose grading components
$C^n$, \,$n\in\boZ$, are contraadjusted $R$\+modules.
 In this context: \par
\textup{(a)} Denote by $F$ the canonical filtration from
Lemma~\ref{vwr-projective-canonical-filtration-lemma} on the very weakly
relatively projective graded $B^*$\+module~$C^*$.
 Then the bigrading components $\gr_F^iC^n$, \,$i$,~$n\in\boZ$, of
the associated graded module $\gr_F^*C^*$ are contraadjusted
$R$\+modules. \par
\textup{(b)} Denote by $F$ the canonical filtration from
Lemma~\ref{vwr-injective-canonical-filtration-lemma} on the very weakly
relatively injective graded $B^*$\+module~$C^*$.
 Then the bigrading components $\gr_F^iC^n$, \,$i$,~$n\in\boZ$, of
the associated graded module $\gr_F^*C^*$ are contraadjusted
$R$\+modules.
\end{lem}

\begin{proof}
 Part~(b) follows from the proof of
Lemma~\ref{thick-graded-modules-co-extension-of-scalars}.
 Part~(a) is deduced from part~(b) by the argument from the proof
of Lemma~\ref{thick-graded-modules-lemma}.
 Notice that, for any contraadjusted $R$\+module $D$ and any finitely
generated projective $R$\+module $P$, the $R$\+module
$P\ot_RD\simeq\Hom_R(\Hom_R(P,R),D)$ is contraadjusted by
Lemma~\ref{vfl-cta-tensor-hom-lemma}(b).
\end{proof}

\begin{lem} \label{thick-cot-module-associated-graded-lemma}
 Let $R$ be a commutative ring, $E$ be a finitely generated projective
$R$\+module, and $A^*$ be the graded algebra $B^*=\bigwedge_R^*(E)$
over~$R$.
 Let $C^*$ be a thick graded $B^*$\+module whose grading components
$C^n$, \,$n\in\boZ$, are cotorsion $R$\+modules.
 In this context: \par
\textup{(a)} Denote by $F$ the canonical filtration from
Lemma~\ref{vwr-projective-canonical-filtration-lemma} on the very weakly
relatively projective graded $B^*$\+module~$C^*$.
 Then the bigrading components $\gr_F^iC^n$, \,$i$,~$n\in\boZ$, of
the associated graded module $\gr_F^*C^*$ are cotorsion $R$\+modules.
\par
\textup{(b)} Denote by $F$ the canonical filtration from
Lemma~\ref{vwr-injective-canonical-filtration-lemma} on the very weakly
relatively injective graded $B^*$\+module~$C^*$.
 Then the bigrading components $\gr_F^iC^n$, \,$i$,~$n\in\boZ$, of
the associated graded module $\gr_F^*C^*$ are cotorsion $R$\+modules.
\end{lem}

\begin{proof}
 Part~(b): let us reduce the assertion in question to the case when
$E$ is a free $R$\+module.
 In view of the discussion in
Section~\ref{loc-free-sheaves-subsecn}, there exists a finite affine
open covering $\Spec R=\bigcup_{\alpha=1}^N\Spec S_\alpha$ of
the affine scheme $\Spec R$ such that the $S_\alpha$\+module
$S_\alpha\ot_RE$ is free for every index~$\alpha$.
 Put $A_\alpha^*=S_\alpha\ot_RB^*=
\bigwedge_{S_\alpha}^*(S_\alpha\ot_RE)$ and $D_\alpha^*=\Hom_R(S,C^*)$.
 Denote by $F$ the canonical filtration from
Lemma~\ref{vwr-injective-canonical-filtration-lemma} on the very weakly
relatively injective graded $A_\alpha^*$\+modules~$D_\alpha^*$.
 By Lemma~\ref{thick-cta-module-associated-graded-lemma}(b),
the $R$\+modules $\gr_F^iC^n$ are contraadjusted.
 By Corollary~\ref{canonical-filtration-qcoh-lcth-compatible}(b),
we have $\gr_F^iD_\alpha^n=\Hom_R(S,\gr_F^iC^n)$.
 If we manage to show that the $S_\alpha$\+modules
$\gr_F^iD_\alpha^n$ are cotorsion, then it would follow that
the $R$\+modules $\gr_F^iC^n$ are cotorsion
by~\cite[Lemma~1.3.6(a)]{Pcosh}.

 So we can assume that $E$ is a free $R$\+module with $m$~free
generators.
 We have a natural isomorphism of $R$\+modules
$\gr_F^iC^n\simeq\Hom_R(B^{-i},\gr_F^0C^{-i+n})$ for all
$i$, $n\in\boZ$.
 The $R$\+modules $B^{-i}$ are finitely generated projective (in fact,
free), so it suffices to show that the grading components of
$\gr_F^0C^*=\Hom_{B^*}^*(B^*/F^1B^*,C^*)$ are cotorsion $R$\+modules.
 Proceeding by induction on~$m$, in the case of $m=0$ we have
$C^*=\gr_F^0C^*$ and there is nothing to prove.

 For $m\ge1$, pick a direct sum decomposition $E=Re\oplus D$, where $D$
is a free $R$\+module with $m-1$ free generators and $Re$ is a free
$R$\+module with one generator~$e$.
 Put $A^*=\bigwedge_R^*(D)=B^*/B^*e$; so the graded commutative
$R$\+algebra $B^*$ is obtained by adjoining to $A^*$ one free (graded
commutative) generator~$e$ of degree~$1$.
 In order to make the induction step, we need to check that the graded
$A^*$\+module $K^*=\ker(e\:C^*\to C^*)$ is thick and its grading
components are cotorsion $R$\+modules.

 Consider the complex of graded $B^*$\+modules
\begin{equation} \label{multiplication-with-e-complex}
 \dotsb\lrarrow C^*(-1)\overset e\lrarrow C^*\overset e\lrarrow
 C^*(1)\lrarrow\dotsb
\end{equation}
where all the differentials are the operators of multiplication
with~$e$ and $M^*\longmapsto M^*(n)$, \,$n\in\boZ$, denotes
the grading shift.
 The passage to associated (bi)graded $\gr_F^*B^*$\+modules
takes~\eqref{multiplication-with-e-complex} to the complex
\begin{equation} \label{multiplication-with-e-assoc-graded-complex}
 \dotsb\lrarrow \gr_F^*C^*(-1)\overset e\lrarrow \gr_F^*C^*
 \overset e\lrarrow \gr_F^*C^*(1)\lrarrow\dotsb
\end{equation}
 We have $\gr_F^*C^*\simeq\Hom_R^*(\gr_F^*B^*,N^*)$ for the graded
$R$\+module $N^*=\gr_F^0C^*$.
 Furthermore, the bigraded $R$\+algebra $\gr_F^*B^*$ is naturally
isomorphic to the graded $R$\+algebra $B^*$ endowed with the obvious
(diagonal) second grading.
 One can easily see that
the complex~\eqref{multiplication-with-e-assoc-graded-complex}
is acyclic.
 Thus the complex~\eqref{multiplication-with-e-complex} is acyclic
as well.

 Consider the filtration $F$ on $K^*$ induced by the filtration
$F$ on $C^*$; then exactness of
the complex~\eqref{multiplication-with-e-assoc-graded-complex} implies
the isomorphism $\gr_F^*K^*\simeq
\ker(e\:\gr_F^*C^*\to\nobreak\gr_F^*C^*)\allowbreak\simeq
\Hom_R^*(\gr_F^*A^*,N^*)$.
 So the graded $A^*$\+module $K^*$ satisfies the conditions of
Lemma~\ref{very-weakly-relatively-injective-lemma}(3).
 Applying Theorem~\ref{module-cotorsion-periodicity} to
the exact complex~\eqref{multiplication-with-e-complex} for
the graded $B^*$\+module $C^*$, we see that the grading components
of $K^*$ are cotorsion $R$\+modules.
 Finally, we can say that the grading components of
$\gr_F^0C^*=\gr_F^0K^*$ are cotorsion $R$\+modules by the assumption
of induction on~$m$.

 Part~(a) is deduced from part~(b) similarly to the proof of
Lemma~\ref{thick-cta-module-associated-graded-lemma}(a).
\end{proof}

 Let $X$ be a scheme with an open covering~$\bW$ and $\E$ be
a finite locally free sheaf on~$X$.
 Consider the quasi-coherent graded algebra $\cB^*=\bigwedge_X^*(\E)$
over~$X$.
 The definition of a thick $\bW$\+locally contraherent graded
$\cB^*$\+module on $X$ was given in
Section~\ref{thick-loc-contraherent-modules-subsecn}.
 We will say that an $X$\+locally cotorsion $\bW$\+locally contraherent
graded $\cB^*$\+module on $X$ is \emph{thick} if it is thick as
a $\bW$\+locally contraherent graded $\cB^*$\+module.
 Let us introduce notation for the corresponding full subcategory
in $\cB^*\Lcth_\bW$:
$$
 \cB^*\Lcth_\bW^{X\dlct,\thk} =
 \cB^*\Lcth_\bW^\thk\cap\cB^*\Lcth_\bW^{X\dlct}.
$$
 The full subcategory $\cB^*\Lcth_\bW^{X\dlct,\thk}\subset
\cB^*\Lcth_\bW^{X\dlct}$ is closed under extensions, kernels of
admissible epimorphisms, cokernels of admissible monomorphisms, and
infinite products in the exact category $\cB^*\Lcth_\bW^{X\dlct}$.
 So the full subcategory $\cB^*\Lcth_\bW^{X\dlct,\thk}$ inherits
an exact category structure from $\cB^*\Lcth_\bW^{X\dlct}$.

 In the following exposition, we apply the machinery of filtered
quasi-coherent and locally contraherent graded modules developed in
Section~\ref{filtered-qcoh-lcth-subsecn} to quasi-coherent and
contraherent graded module over the graded algebra $\cB^*$ endowed
with its natural decreasing filtration associated with the grading.
 So we consider quasi-coherent and locally contraherent
$\cB^*$\+modules endowed with \emph{both} a grading and a finite
filtration by graded submodules.
 As usual, the (obvious) fact that the results of
Section~\ref{filtered-qcoh-lcth-subsecn} remain valid in the setting
with an additional grading is used without much explanation.

 The definition of a thick quasi-coherent graded $\cB^*$\+module on $X$
was given in Section~\ref{thick-quasi-coherent-modules-subsecn}.

\begin{cor} \label{thick-qcoh-sheaves-are-filtered}
 Let $X$ be a scheme and $\E$ be a finite locally free sheaf on~$X$.
 Let $\cB^*$ be the quasi-coherent graded algebra
$\cB^*=\bigwedge_X^*(\E)$ over~$X$.
 Assume for simplicity that the rank of $\E$ is bounded; so $\cB^n=0$
for $n$~large enough.
 Endow the quasi-coherent graded algebra $\cB^*$ with the decreasing
filtration induced by the grading,
$F^i\cB^*=\bigoplus_{j=i}^\infty\cB^j$ for all $i\in\boZ$.
 Let $\M^*$ be a thick quasi-coherent graded $\cB^*$\+module on~$X$.
 In this context: \par
\textup{(a)} Define a finite decreasing filtration $F$ on $\M^*$ by
the rule that, for every affine open subscheme $U\subset X$,
the filtration $F$ on the thick graded $\cB^*(U)$\+module $\M^*(U)$ is
the canonical filtration from
Lemma~\ref{vwr-projective-canonical-filtration-lemma}.
 Then $(\M^*,F)$ is a finitely filtered quasi-coherent graded module
over $(\cB^*,F)$. \par
\textup{(b)} Define a finite decreasing filtration $F$ on $\M^*$ by
the rule that, for every affine open subscheme $U\subset X$,
the filtration $F$ on the thick graded $\cB^*(U)$\+module $\M^*(U)$ is
the canonical filtration from
Lemma~\ref{vwr-injective-canonical-filtration-lemma}.
 Then $(\M^*,F)$ is a finitely filtered quasi-coherent graded module
over $(\cB^*,F)$.
\end{cor}

\begin{proof}
 Part~(a) follows immediately from
Corollary~\ref{canonical-filtration-qcoh-lcth-compatible}(a).
 To prove part~(b), one can represent the scheme $X$ as a finite
disjoint union $X=\coprod_d X_d$ of open subschemes $X_d\subset X$ on
which the finite locally free sheaf $\E$ has constant rank~$d$, and
notice that over every scheme $X_d$ the two filtrations only differ
by a shift of the filtration degree.
\end{proof}

 Clearly, if $\M^*$ is a thick quasi-coherent graded $\cB^*$\+module
and $F$ is the canonical filtration on $M^*$ described in
Corollary~\ref{thick-qcoh-sheaves-are-filtered}(a), then the associated
bigraded quasi-coherent module $\gr_F^*\M^*$ over $\gr_F^*\cB^*=\cB^*$
is computed as $\gr_F^*\M^*\simeq\gr_F^*\cB^*\ot_{\cO_X}\gr_F^0\M^*$.
 Conversely, if $\M^*$ is a quasi-coherent graded $\cB^*$\+module
admittting a finite decreasing filtration $F$ such that $(\M^*,F)$
is a filtered quasi-coherent module over $(\cB^*,F)$, and
the associated bigraded quasi-coherent module $\gr_F^*\M^*$
over $\gr_F^*\cB^*$ has the form $\gr_F^*\M^*\simeq\gr_F^*\cB^*
\ot_{\cO_X}\N^*$ for some graded quasi-coherent sheaf $\N^*$ on $X$,
then $\M^*$ is a thick quasi-coherent graded $\cB^*$\+module.
 The filtration $F$ on $\M^*$ with the described properties is
unique by Lemma~\ref{vwr-projective-canonical-filtration-lemma}.

 Similarly, if $\M^*$ is a thick quasi-coherent graded $\cB^*$\+module
and $F$ is the canonical filtration on $M^*$ described in
Corollary~\ref{thick-qcoh-sheaves-are-filtered}(b), then the associated
bigraded quasi-coherent module $\gr_F^*\M^*$ over $\gr_F^*\cB^*=\cB^*$
is computed as $\gr_F^*\M^*\simeq
\cHom_{\cO_X}(\gr_F^*\cB^*,\gr_F^0\M^*)$.
 Conversely, if $\M^*$ is a quasi-coherent graded $\cB^*$\+module
admittting a finite decreasing filtration $F$ such that $(\M^*,F)$
is a filtered quasi-coherent module over $(\cB^*,F)$, and
the associated bigraded quasi-coherent module $\gr_F^*\M^*$
over $\gr_F^*\cB^*$ has the form $\gr_F^*\M^*\simeq
\cHom_{\cO_X}(\gr_F^*\cB^*,\N^*)$ for some graded quasi-coherent
sheaf $\N^*$ on $X$, then $\M^*$ is a thick quasi-coherent graded
$\cB^*$\+module.
 The filtration $F$ on $\M^*$ with the described properties is
unique by Lemma~\ref{vwr-injective-canonical-filtration-lemma}.
{\emergencystretch=1em\par}

\begin{cor} \label{thick-lcta-lct-cosheaves-are-filtered}
 Let $X$ be a scheme with an open covering\/ $\bW$ and $\E$ be a finite
locally free sheaf on~$X$.
 Let $\cB^*$ be the quasi-coherent graded algebra
$\cB^*=\bigwedge_X^*(\E)$ over~$X$.
 Assume for simplicity that the rank of $\E$ is bounded; so $\cB^n=0$
for $n$~large enough.
 Endow the quasi-coherent graded algebra $\cB^*$ with the decreasing
filtration induced by the grading,
$F^i\cB^*=\bigoplus_{j=i}^\infty\cB^j$ for all $i\in\boZ$.
 In this setting: \par
\textup{(a)} Let\/ $\P^*$ be a thick\/ $\bW$\+locally contraherent
graded $\cB^*$\+module on~$X$.
 Define a finite decreasing filtration $F$ on\/ $\P^*$ by the rule
that, for every affine open subscheme $U\subset X$ subordinate
to\/~$\bW$, the filtration $F$ on the thick graded $\cB^*(U)$\+module\/
$\P^*[U]$ is the canonical filtration from
Lemma~\ref{vwr-injective-canonical-filtration-lemma}.
 Then $(\P^*,F)$ is a finitely filtered\/ $\bW$\+locally contraherent
graded module over $(\cB^*,F)$. \par
\textup{(b)} Let\/ $\P^*$ be a thick $X$\+locally cotorsion\/
$\bW$\+locally contraherent graded $\cB^*$\+module on~$X$.
 Define a finite decreasing filtration $F$ on\/ $\P^*$ by the rule
that, for every affine open subscheme $U\subset X$ subordinate to\/
$\bW$, the filtration $F$ on the thick graded $\cB^*(U)$\+module\/
$\P^*[U]$ is the canonical filtration from
Lemma~\ref{vwr-injective-canonical-filtration-lemma}.
 Then $(\P^*,F)$ is a finitely filtered $X$\+locally cotorsion\/
$\bW$\+locally contraherent graded module over $(\cB^*,F)$.
\end{cor}

\begin{proof}
 Part~(a) follows from
Corollary~\ref{canonical-filtration-qcoh-lcth-compatible}(b)
and Lemma~\ref{thick-cta-module-associated-graded-lemma}(b).
 Part~(b) follows from
Corollary~\ref{canonical-filtration-qcoh-lcth-compatible}(b)
and Lemma~\ref{thick-cot-module-associated-graded-lemma}(b).
\end{proof}

 Let $\P^*$ be a thick $\bW$\+locally contraherent graded
$\cB^*$\+module, and let $F$ be the canonical filtration on $\P^*$
described in Corollary~\ref{thick-lcta-lct-cosheaves-are-filtered}(a).
 Then the associated bigraded $\bW$\+locally contraherent module
$\gr_F^*\P^*$ over $\gr_F^*\cB^*=\cB^*$ is computed as
$\gr_F^*\P^*\simeq\Cohom_X(\gr_F^*\cB^*,\gr_F^0\P^*)$.
 (We refer to~\cite[Sections~2.4 and~3.6]{Pcosh} and
Section~\ref{cohom-from-quasi-modules-subsecn} above for the definition
of the functor $\Cohom_X$.)
 Conversely, if $\P^*$ is a $\bW$\+locally contraherent graded
$\cB^*$\+module admitting a finite decreasing filtration $F$ such
that $(\P^*,F)$ is a filtered $\bW$\+locally contraherent graded
modules over $(\cB^*,F)$, and the associated bigraded $\bW$\+locally
contraherent module $\gr_F^*\P^*$ over $\gr_F^*\cB^*$ has the form
$\gr_F^*\P^*\simeq\Cohom_X(\gr_F^*\cB^*,\Q^*)$ for some graded
$\bW$\+locally contraherent cosheaf $\Q^*$ on $X$, then $\P^*$ is
a thick $\bW$\+locally contraherent graded $\cB^*$\+module.
 The filtration $F$ on $\P^*$ with the described properties is unique
by Lemma~\ref{vwr-injective-canonical-filtration-lemma}.

\subsection{Associated graded quasi-coherent modules are
contraadjusted} \label{qcoh-assoc-graded-are-cta-subsecn}
 In this section and the next one we extend the results of the previous
Section~\ref{assoc-graded-modules-are-cta-cot} to thick quasi-coherent
graded modules over a quasi-compact semi-separated scheme.
 This is nontrivial, as the notions of a contraadjusted and cotorsion
quasi-coherent sheaves over a scheme are \emph{not local};
so global arguments are needed.
 We also prove a similar result for thick antilocal contraherent
graded modules.

 Let $X$ be a scheme and $\E$ be a finite locally free sheaf on~$X$.
 Consider the quasi-coherent graded algebra $\cB^*=\bigwedge_X^*(\E)$
over~$X$.
 We will say that an $X$\+contraadjusted or $X$\+cotorsion
quasi-coherent graded $\cB^*$\+module is \emph{thick} if it is thick
as a quasi-coherent graded $\cB^*$\+module (in the sense of
the definition in
Section~\ref{thick-quasi-coherent-modules-subsecn}).
 Let us introduce notation for the corresponding full subcategories
in $\cB^*\Qcoh$:
\begin{align*}
 \cB^*\Qcoh^{X\dcta}_\thk &= \cB^*\Qcoh_\thk\cap\cB^*\Qcoh^{X\dcta}, \\
 \cB^*\Qcoh^{X\dcot}_\thk &= \cB^*\Qcoh_\thk\cap\cB^*\Qcoh^{X\dcot}.
\end{align*}
 The full subcategories $\cB^*\Qcoh^{X\dcot}_\thk\subset
\cB^*\Qcoh^{X\dcta}_\thk\subset\cB^*\Qcoh$ are closed under extensions
and (for a quasi-compact semi-separated scheme~$X$) under cokernels of 
monomorphisms in $\cB^*\Qcoh$.
 So both the full subcategories inherit exact category structures
from the abelian exact structure of $\cB^*\Qcoh$.

 We start with an auxiliary lemma.
 Given a scheme $X$ and two quasi-coherent sheaves $\F$ and $\C$ on
$X$, we denote by $\cHom_{X\dqc}(\F,\C)$ the \emph{quasi-coherent
internal $\cHom$ sheaf}, defined by the rules that
$\cHom_{X\dqc}(\F,\C)\in X\Qcoh$ and $\Hom_X(\M,\cHom_{X\dqc}(\F,\C))
\simeq\Hom_X(\F\ot_{\cO_X}\M,\>\F)$ for all $\M\in X\Qcoh$.
 The quasi-coherent sheaf $\cHom_{X\dqc}(\F,\C)$ can be obtained by
applying the coherator construction~\cite[Appendix~B]{TT},
\cite[Proposition Tag~077P]{SP} to the sheaf of $\cO_X$\+modules
$\cHom_{\cO_X}(\F,\C)$ discussed in
Section~\ref{loc-free-sheaves-subsecn}.

\begin{lem} \label{qcoh-internal-Hom-sheaf-cta-cot}
\textup{(a)} Let $X$ be a scheme, $\F$ be a very flat quasi-coherent
sheaf on $X$, and $\C$ be a contraadjusted quasi-coherent sheaf on~$X$.
 Then the quasi-coherent sheaf $\cHom_{X\dqc}(\F,\C)$ on $X$ is 
contraadjusted. \par
\textup{(b)} Let $X$ be a scheme, $\F$ be a flat quasi-coherent sheaf
on $X$, and $\C$ be a cotorsion quasi-coherent sheaf on~$X$.
 Then the quasi-coherent sheaf $\cHom_{X\dqc}(\F,\C)$ on $X$
is cotorsion.
\end{lem}

\begin{proof}
 Part~(a) is a quasi-coherent version of
Lemma~\ref{hom-injective-cotorsion}(a), while part~(b) is
a quasi-coherent generalization of
Lemma~\ref{vfl-cta-tensor-hom-lemma}(b).
 To prove both parts~(a) and~(b), consider the pair of adjoint
functors $\G\longmapsto\F\ot_{\cO_X}\G\:X\Qcoh\rarrow X\Qcoh$ and
$\C\longmapsto\cHom_{X\dqc}(\F,\C)\:X\Qcoh\rarrow X\Qcoh$.
 By~\cite[Lemma~1.7(b)]{Pal}, for any flat quasi-coherent sheaf $\F$
and any quasi-coherent sheaves $\G$ and $\C$ on $X$, there is
an injective morphism of abelian groups
$\Ext_X^1(\G,\cHom_{X\dqc}(\F,\C))\rarrow
\Ext_X^1(\F\ot_{\cO_X}\G,\>\C)$; so
$\Ext_X^1(\F\ot_{\cO_X}\G,\>\C)=0$ implies
$\Ext_X^1(\G,\cHom_{X\dqc}(\F,\C))=0$.
 It remains to recall that the tensor product of any two flat
(respectively, very flat) quasi-coherent sheaves on $X$ is a flat
(respectively, very flat) quasi-coherent sheaf again.
\end{proof}

 Let $X$ be a quasi-compact semi-separated scheme.
 By the \emph{contraadjusted dimension} of a quasi-coherent sheaf
$\M$ on $X$ we mean its coresolution dimension with respect to
the coresolving subcategory $X\Qcoh^\cta\subset X\Qcoh$ of
contraadjusted quasi-coherent sheaves (in the sense of
Section~\ref{co-resolution-dimensions-subsecn}).
 By~\cite[Lemma~4.7.1(b)]{Pcosh}, the contraadjusted dimension of
any quasi-coherent sheaf on $X$ does not exceed the minimal number
of affine open subschemes in an affine open covering of~$X$.

\begin{cor} \label{qcoh-internal-Hom-sheaf-cta-dimension}
 Let $X$ be a quasi-compact semi-separated scheme, $\F$ be a finite
locally free sheaf on $X$, and $\M$ be a quasi-coherent sheaf on~$X$.
 Then the contraadjusted dimension of the quasi-coherent sheaf
$\cHom_{\cO_X}(\F,\M)$ on $X$ does not exceed the contraadjusted
dimension of the quasi-coherent sheaf~$\M$.
\end{cor}

\begin{proof}
 For any quasi-coherent sheaf $\M$, the sheaf of $\cO_X$\+modules
$\cHom_{\cO_X}(\F,\M)$ on $X$ is quasi-coherent; so we have
$\cHom_{X\dqc}(\F,X)=\cHom_{\cO_X}(\F,\M)$.
 The functor $\cHom_{\cO_X}(\F,{-})\:X\Qcoh\rarrow X\Qcoh$ is obviously
exact; by Lemma~\ref{qcoh-internal-Hom-sheaf-cta-cot}(a), it also takes
contraadjusted quasi-coherent sheaves to contraadjusted quasi-coherent
sheaves.
 Hence this functor takes any contraadjusted coresolution of $\M$ to
a contraadjusted coresolution of $\cHom_{\cO_X}(\F,\M)$.
\end{proof}

\begin{cor} \label{thick-cta-qcoh-sheaves-are-filtered-by-cta}
 Let $X$ be a quasi-compact semi-separated scheme and $\E$ be a finite
locally free sheaf on~$X$.
 Let $\cB^*$ be the quasi-coherent graded algebra
$\cB^*=\bigwedge_X^*(\E)$ over~$X$.
 Endow the quasi-coherent graded algebra $\cB^*$ with the decreasing
filtration induced by the grading,
$F^i\cB^*=\bigoplus_{j=i}^\infty\cB^j$ for all $i\in\boZ$.
 Let $\M^*$ be a thick $X$\+contraadjusted quasi-coherent graded
$\cB^*$\+module on~$X$.
 Consider the finite decreasing filtration $F$ on $\M^*$ defined in
Corollary~\ref{thick-qcoh-sheaves-are-filtered}(a) or~(b).
 Then the quasi-coherent sheaves\/ $\gr_F^i\M^n$ on $X$ are
contraadjusted for all $i$, $n\in\boZ$.
 In other words, $(\M^*,F)$ is an $X$\+contraadjusted finitely filtered
quasi-coherent graded module over $(\cB^*,F)$.
\end{cor}

\begin{proof}
 In view of the proof of
Corollary~\ref{thick-qcoh-sheaves-are-filtered}(b),
it suffices to consider the filtration $F$ on $\M^*$ from
Corollary~\ref{thick-qcoh-sheaves-are-filtered}(a).
 Then the quasi-coherent bigraded module $\gr_F^*\M^*$ over
the quasi-coherent bigraded algebra $\gr_F^*\cB^*=\cB^*$ can be
computed as $\gr_F^*\M^*\simeq\gr_F^*\cB^*\ot_{\cO_X}\gr_F^0\M^*$.
 Let $N$ be the minimal number of affine open subschemes in an affine
open covering of~$X$.
 Then the contraadjusted dimensions of the grading components of
the graded quasi-coherent sheaf $F^1\M^*$ on $X$ cannot exceed $N$
by~\cite[Lemma~4.7.1(b)]{Pcosh}.
 Since the coresolution dimension does not depend on the choice of
a coresolution and the grading components of the graded quasi-coherent
sheaf $\M^*$ are contraadjusted by assumption, it follows that
the contraadjusted dimensions of the grading components of the graded
quasi-coherent sheaf $\gr_F^0\M^*=\M^*/F^1\M^*$ do not exceed $N-1$
if $N\ge1$.

 Now we have $\gr_F^i\M^*\simeq\cB^i\ot_{\cO_X}\gr_F^0\M^*\simeq
\cHom_{\cO_X}(\cHom_{\cO_X}(\cB^i,\cO_X),\gr_F^0\M^*)$ for all $i>0$.
 By Corollary~\ref{qcoh-internal-Hom-sheaf-cta-dimension}, it follows
that the contraadjusted dimensions of the grading components of
the graded quasi-coherent sheaf $\gr_F^i\M^*$ do not exceed $N-1$.
 As the coresolution dimension is not increased by finitely iterated
extensions~\cite[Proposition~2.3(2)]{Sto},
\cite[Lemma~A.5.6(a)]{Pcosh}, it follows that the contraadjusted
dimensions of the grading components of the graded quasi-coherent
sheaf $F^1\M^*$ do not exceed~$N-1$.
 Thus the contraadjusted dimensions of the grading components of
the graded quasi-coherent sheaf $\gr_F^0\M^*=\M^*/F^1\M^*$
do not exceed $N-2$ if $N\ge2$.
 Proceeding in this way, we bring the upper bound on the contraadjusted
dimensions of the quasi-coherent sheaves $\gr_F^i\M^n$ on $X$ down
to~$0$.
\end{proof}

 Assume that the scheme $X$ is quasi-compact and semi-separated, and
let $\bW$ be an open covering of~$X$.
 We will say that an antilocal contraherent graded $\cB^*$\+module or
an antilocal $X$\+locally cotorsion contraherent graded $\cB^*$\+module
is \emph{thick} if it is thick as a locally contraherent graded
$\cB^*$\+module (in the sense of the definitions in
Sections~\ref{thick-loc-contraherent-modules-subsecn}
and~\ref{assoc-graded-modules-are-cta-cot}).
 Let us introduce notation for the corresponding full subcategories
in $\cB^*\Lcth_\bW$:
\begin{align*}
 \cB^*\Ctrh_\al^\thk &=
 \cB^*\Lcth_\bW^\thk\cap\cB^*\Ctrh_\al, \\
 \cB^*\Ctrh_\al^{X\dlct,\thk} &=
 \cB^*\Lcth_\bW^\thk\cap\cB^*\Ctrh_\al^{X\dlct}.
\end{align*}
 The full subcategories $\cB^*\Ctrh_\al^\thk$ and 
$\cB^*\Ctrh_\al^{X\dlct,\thk}$
are closed under extensions, kernels of admissible epimorphisms,
cokernels of admissible monomorphisms, and infinite products in
their ambient exact categories $\cB^*\Ctrh_\al$ and
$\cB^*\Ctrh_\al^{X\dlct}$.
 So both the full subcategories inherit exact category structures.

 Once again, we start with an auxiliary lemma.
 For a more general version, see
Lemma~\ref{Cohom-from-quasi-mod-into-antilocal-is-antilocal} below.

\begin{lem} \label{Cohom-into-antilocal-cosheaf-is-antilocal}
 Let $X$ be a quasi-compact semi-separated scheme.
 Then \par
\textup{(a)} for any very flat quasi-coherent sheaf $\F$ on $X$ and
any antilocal contraherent cosheaf\/ $\P$ on $X$, the contraherent
cosheaf\/ $\Cohom_X(\F,\P)$ on $X$ is antilocal; \par
\textup{(b)} for any flat quasi-coherent sheaf $\F$ on $X$ and any
antilocal localy cotorsion contraherent cosheaf\/ $\P$ on $X$,
the locally cotorsion contraherent cosheaf\/ $\Cohom_X(\F,\P)$ on $X$
is antilocal.
\end{lem}

\begin{proof}
 Let us prove part~(a).
 By construction (see~\cite[Sections~2.4 and~3.6]{Pcosh}), the functor
$\Cohom_X(\F,{-})\:X\Lcth_\bW\rarrow X\Lcth_\bW$ is exact.
 In view of~\cite[Corollary~4.3.5(c)]{Pcosh}, it suffices to check
that, for every affine open subscheme $U\subset X$ with the open
immersion morphism $j\:U\rarrow X$ and any contraherent cosheaf $\Q$
on $U$, the contraherent cosheaf $\Cohom_X(\F,j_!\Q)$ on $X$
is antilocal.
 It remains to point out that the natural isomorphism
$\Cohom_X(\F,j_!\Q)\simeq f_!\Cohom_U(j^*\F,\Q)$ of contraherent
cosheaves on~$X$ holds by~\cite[formula~(3.21) in Section~3.8]{Pcosh}.
 The proof of part~(b) is similar and uses exactness of the functor
$\Cohom_X(\F,{-})\:X\Lcth_\bW^\lct\rarrow X\Lcth_\bW^\lct$
together with~\cite[Corollary~4.3.6(c)]{Pcosh}.
\end{proof}

 Let $X$ be a quasi-compact semi-separated scheme with an open
covering~$\bW$.
 By the \emph{antilocal dimension} of a $\bW$\+locally contraherent
cosheaf $\P$ on $X$ we mean its resolution dimension with respect to
the resolving subcategory $X\Ctrh_\al\subset X\Lcth_\bW$ of antilocal
contraherent cosheaves (in the sense of
Section~\ref{co-resolution-dimensions-subsecn}).
 Since, as the open covering $\bW$ varies, the full exact subcategories
$X\Lcth_\bW$ are closed under kernels of admissible epimorphisms in
each other, the antilocal dimension of $\P$ remains unchanged when
the open covering $\bW$ is refined.
 By~\cite[Lemma~4.7.1(c)]{Pcosh}, the antilocal dimension of
any $\bW$\+locally contraherent cosheaf on $X$ does not exceed $N-1$,
where $N$ is the minimal number of affine open subschemes in an affine
open covering of $X$ subordinate to~$\bW$.

\begin{cor} \label{Cohom-does-not-increase-antilocal-dimension}
 Let $X$ be a quasi-compact semi-separated scheme with an open
covering~$\bW$.
 Let $\F$ be a very flat quasi-coherent sheaf on $X$ and\/ $\P$ be
a\/ $\bW$\+locally contraherent cosheaf on~$X$.
 Then the antilocal dimension of the\/ $\bW$\+locally contraherent
cosheaf\/ $\Cohom_X(\F,\P)$ on $X$ does not exceed the antilocal
dimension of the\/ $\bW$\+locally contraherent cosheaf\/~$\P$.
\end{cor}

\begin{proof}
 The functor $\Cohom_X(\F,{-})\:X\Lcth_\bW\rarrow X\Lcth_\bW$ is
exact by construction; by
Lemma~\ref{Cohom-into-antilocal-cosheaf-is-antilocal}(a), it also takes
antilocal contraherent cosheaves to antilocal contraherent cosheaves.
 Hence this functor takes any antilocal resolution of $\P$ to
an antilocal resolution of $\Cohom_X(\F,\P)$.
\end{proof}

\begin{cor} \label{thick-antilocal-cosheaves-are-filtered-by-antilocal}
 Let $X$ be a quasi-compact semi-separated scheme and $\E$ be a finite
locally free sheaf on~$X$.
 Let $\cB^*$ be the quasi-coherent graded algebra
$\cB^*=\bigwedge_X^*(\E)$ over~$X$.
 Endow the quasi-coherent graded algebra $\cB^*$ with the decreasing
filtration induced by the grading,
$F^i\cB^*=\bigoplus_{j=i}^\infty\cB^j$ for all $i\in\boZ$.
 Let\/ $\P^*$ be a thick antilocal contraherent graded $\cB^*$\+module
on~$X$.
 Consider the finite decreasing filtration $F$ on\/ $\P^*$ defined in
Corollary~\ref{thick-lcta-lct-cosheaves-are-filtered}(a).
 Then the contraherent cosheaves\/ $\gr_F^i\P^n$ on $X$ are antilocal
for all $i$, $n\in\boZ$.
 In other words, $(\P^*,F)$ is an antilocal finitely filtered
contraherent graded module over $(\cB^*,F)$.
\end{cor}

\begin{proof}
 The argument is dual-analogous to the proof of
Corollary~\ref{thick-cta-qcoh-sheaves-are-filtered-by-cta}.
 We use Corollary~\ref{thick-lcta-lct-cosheaves-are-filtered}(a)
for the open covering $\bW=\{X\}$ of the scheme~$X$.
 The contraherent bigraded module $\gr_F^*\P^*$ over the quasi-coherent
bigraded algebra $\gr_F^*\cB^*$ over $X$ can be computed as
$\gr_F^*\P^*\simeq\Cohom_X(\gr_F^*\cB^*,\gr_F^0\P^*)$.

 Let $N$ be the minimal number of affine open subschemes in an affine
open covering of~$X$.
 Then the antilocal dimensions of the grading components of
the graded contraherent cosheaf $\P^*/F^0\P^*$ on $X$ cannot
exceed $N-1$ by~\cite[Lemma~4.7.1(c)]{Pcosh}.
 Since the resolution dimension does not depend on the choice of
a resolution and the grading components of the graded quasi-coherent
sheaf $\P^*$ are antilocal by assumption, it follows that
the grading components of the graded quasi-coherent sheaf
$\gr_F^0\P^*=F^0\P^*$ do not exceed $N-2$ if $N\ge2$.

 Now we have $\gr_F^i\P^*=\Cohom_X(\cB^{-i},\gr_F^0\P^*)$ for all $i<0$.
 By Corollary~\ref{Cohom-does-not-increase-antilocal-dimension}, it
follows that the grading components of the graded contraherent cosheaf
$\gr_F^i\P^*$ do not exceed $N-2$.
 As the resolution dimension is not increased by finitely iterated
extensions~\cite[Proposition~2.3(2)]{Sto},
\cite[Lemma~A.5.6(a)]{Pcosh}, it follows that the antilocal dimensions
of the grading components of the graded contraherent cosheaf
$\P^*/F^0\P^*$ do not exceed~$N-2$.
 Thus the antilocal dimensions of the grading components of the graded
contraherent cosheaf $\gr_F^0\P^*=F^0\P^*$ do not exceed $N-3$
if $N\ge3$, etc.
\end{proof}

\subsection{Associated graded quasi-coherent modules are cotorsion}
 The proof of Corollary~\ref{thick-cta-qcoh-sheaves-are-filtered-by-cta}
in the previous section was based on finiteness of the contraadjusted
dimensions of quasi-coherent sheaves over quasi-compact semi-separated
schemes.
 The cotorsion dimensions need not be finite; so a different argument
is needed for thick $X$\+cotorsion graded modules.

\begin{lem} \label{internal-Hom-and-Cohom-correspond-under-co-contra}
 Let $X$ be a quasi-compact semi-separated scheme. \par
\textup{(a)} Let $\C$ be a contraadjusted quasi-coherent sheaf and\/
$\P$ be an antilocal contraherent cosheaf corresponding to each other
under the equivalence of categories $X\Qcoh^\cta\simeq X\Ctrh_\al$
from~\cite[Lemma~4.8.2]{Pcosh}.
 Let $\F$ be a very flat quasi-coherent sheaf on~$X$.
 Then the contraadjusted quasi-coherent sheaf $\cHom_{X\dqc}(\F,\C)$
and the antilocal contraherent cosheaf\/ $\Cohom_X(\F,\P)$ correspond
to each other under the same equivalence of categories. \par
\textup{(b)} Let $\C$ be a cotorsion quasi-coherent sheaf and\/
$\P$ be an antilocal locally cotorsion contraherent cosheaf
corresponding to each other under the equivalence of categories
$X\Qcoh^\cot\simeq X\Ctrh^\lct_\al$ from~\cite[Lemma~4.8.4(a)]{Pcosh}.
 Let $\F$ be a flat quasi-coherent sheaf on~$X$.
 Then the cotorsion quasi-coherent sheaf $\cHom_{X\dqc}(\F,\C)$
and the antilocal locally cotorsion contraherent cosheaf\/
$\Cohom_X(\F,\P)$ correspond to each other under the same equivalence
of categories.
\end{lem}

\begin{proof}
 Let us prove part~(a).
 The quasi-coherent sheaf $\cHom_{X\dqc}(\F,\C)$ on $X$ is
contraadjusted by Lemma~\ref{qcoh-internal-Hom-sheaf-cta-cot}(a),
while the contraherent cosheaf $\Cohom_X(\F,\P)$ on $X$ is
antilocal by Lemma~\ref{Cohom-into-antilocal-cosheaf-is-antilocal}(a).
 According to the construction of the category equivalence
$X\Qcoh^\cta\simeq X\Ctrh_\al$, we have $\P=\fHom_X(\cO_X,\C)$
and $\C=\cO_X\ocn_X\P$.
 It remains to point out the natural isomorphisms of contraherent 
cosheaves $\fHom_X(\G,\cHom_{X\dqc}(\F,\C))\simeq
\fHom_X(\F\ot_{\cO_X}\G,\>\C)\simeq\Cohom_X(\F,\fHom_X(\G,\C))$,
which hold for any very flat quasi-coherent sheaves $\F$ and $\G$
and any contraadjusted quasi-coherent sheaf $\C$ on~$X$
(and in particular for $\G=\cO_X$); see~\cite[formulas~(2.15)
and~(2.17) in Section~2.5]{Pcosh} (cf.\
Sections~\ref{Cohom-into-fHom-subsecn}\+-\ref{fHom-into-cHom-subsecn}
below).
 Alternatively, in the case of a finite locally free sheaf $\F$ on $X$,
one can use the natural isomorphisms $\cHom_{X\dqc}(\F,\C)\allowbreak
\simeq\Hom_{\cO_X}(\F,\cO_X)\ot_{\cO_X}\C$ and $\Cohom_X(\F,\P)\simeq
\Hom_{\cO_X}(\F,\cO_X)\ot_X\P$ (see~\cite[Section~3.7]{Pcosh} for
the latter isomorphism, including the definition of its right-hand
side).
 Then one can refer to the natural isomorphisms of quasi-coherent
sheaves $\E\ot_{\cO_X}(\G\ocn_X\P)\simeq(\E\ot_{\cO_X}\G)\ocn_X\P
\simeq\G\ocn_X(\E\ot_X\P)$, which hold for any quasi-coherent sheaves
$\E$ and $\G$ and any cosheaf of $\cO_X$\+modules $\P$ on~$X$
(and in particular, for $\E=\Hom_{\cO_X}(\F,\cO_X)$ and $\G=\cO_X$);
see~\cite[formula~(2.20) in Section~2.6 and formula~(3.15) in
Section~3.7]{Pcosh}.
\end{proof}

\begin{thm} \label{thick-corresponds-to-thick-under-co-contra}
 Let $X$ be a quasi-compact semi-separated scheme and $\E$ be a finite
locally free sheaf on~$X$.
 Let $\cB^*$ be the quasi-coherent graded algebra
$\cB^*=\bigwedge_X^*(\E)$ over~$X$.
 Endow the quasi-coherent graded algebra $\cB^*$ with the decreasing
filtration induced by the grading,
$F^i\cB^*=\bigoplus_{j=i}^\infty\cB^j$ for all $i\in\boZ$.
 Let $\M^*$ be an $X$\+contraadjusted quasi-coherent graded
$\cB^*$\+module and\/ $\P^*$ be an antilocal contraherent graded
$\cB^*$\+module on $X$ such that $\M^*$ and\/ $\P^*$ correspond to
each other under the equivalence of categories from the graded version
of Lemma~\textup{\ref{quasi-algebra-underived-naive-co-contra}(a)}.
 Then the $X$\+contraadjusted quasi-coherent graded $\cB^*$\+module
$\M^*$ is thick if and only if the antilocal contraherent graded
$\cB^*$\+module\/ $\P^*$ is thick.

 If this is the case, consider the canonical filtration $F$ on $\M^*$
from Corollary~\textup{\ref{thick-qcoh-sheaves-are-filtered}(b)}
and the canonical filtration $F$ on\/ $\P^*$ from
Corollary~\textup{\ref{thick-lcta-lct-cosheaves-are-filtered}(a)}.
 Then the $X$\+contraadjusted filtered quasi-coherent graded module
$(\M^*,F)$ over $(\cB^*,F)$ and the antilocal filtered contraherent
graded module $(\P^*,F)$ over $(\cB^*,F)$ correspond to each other
under the equivalence of categories from the graded version of
the upper line of diagram~\eqref{X-cta-filtered-naive-co-contra} in
Theorem~\textup{\ref{quasi-algebra-filtered-naive-co-contra}}.
\end{thm}

\begin{proof}
 If the $X$\+contraadjusted quasi-coherent graded $\cB^*$\+module
$\M^*$ is thick, then the filtered quasi-coherent graded module
$(\M^*,F)$ over $(\cB^*,F)$ with the canonical filtration $F$
from Corollary~\ref{thick-qcoh-sheaves-are-filtered}(b) is
$X$\+contraadjusted by
Corollary~\ref{thick-cta-qcoh-sheaves-are-filtered-by-cta}.
 According to (the graded version of) the upper line of
diagram~\eqref{X-cta-filtered-naive-co-contra} in
Theorem~\ref{quasi-algebra-filtered-naive-co-contra},
there is a corresponding filtration $F$ on the antilocal contraherent
graded $\cB^*$\+module $\P^*$ such that $(\P,F)$ is an antilocal
finitely filtered contraherent graded $\cB^*$\+module.

 Furthermore, we have $\gr_F^*\M^*\simeq
\cHom_{\cO_X}(\gr_F^*\cB^*,\N^*)$ for the graded quasi-coherent
sheaf $\N^*=\gr_F^0\M^*$ with contraadjusted grading components.
 The commutative diagram~\eqref{X-cta-filtered-naive-co-contra}
now tells us that the associated bigraded antilocal contraherent
$\gr_F^*\cB^*$\+module $\gr_F^*\P^*$ to $(\P^*,F)$ corresponds to
the $X$\+contraadjusted quasi-coherent bigraded $\gr_F^*\cB^*$\+module
$\cHom_{\cO_X}(\gr_F^*\cB^*,\N^*)$ under the bigraded version
of the category equivalence from
Lemma~\ref{quasi-algebra-underived-naive-co-contra}(a).

 In view of commutativity of
the diagram~\eqref{X-cta-underived-naive-co-contra} in
Lemma~\ref{quasi-algebra-underived-naive-co-contra}(a), we can conclude
from Lemma~\ref{internal-Hom-and-Cohom-correspond-under-co-contra}(a)
that the antilocal contraherent bigraded $\gr_F^*\cB^*$\+module
$\gr_F^*\P^*$ has the form $\gr_F^*\P^*\simeq
\Cohom_X(\gr_F^*\cB^*,\Q^*)$.
 Here $\Q^*=\fHom_X(\cO_X,\N^*)$ is the graded contraherent cosheaf
with antilocal grading components correponding to the graded
quasi-coherent sheaf $\N^*$ under the equivalence of categories
from~\cite[Lemma~4.8.2]{Pcosh}.
 It follows that the antilocal contraherent graded $\cB^*$\+module
$\P^*$ is thick and $F$ is its canonical filtration from
Corollary~\ref{thick-lcta-lct-cosheaves-are-filtered}(a).

 Conversely, if the antilocal contraherent graded $\cB^*$\+module
$\P^*$ is thick, then the filtered contraherent graded module
$(\P^*,F)$ over $(\cB^*,F)$ with the canonical filtration from
Corollary~\ref{thick-lcta-lct-cosheaves-are-filtered}(a) is antilocal
by Corollary~\ref{thick-antilocal-cosheaves-are-filtered-by-antilocal}.
 According to (the graded version of) the upper line of
diagram~\eqref{X-cta-filtered-naive-co-contra} in
Theorem~\ref{quasi-algebra-filtered-naive-co-contra},
there is a corresponding filtration $F$ on the $X$\+contraadjusted
quasi-coherent graded $\cB^*$\+module $\M^*$ such that $(\M,F)$ is
an $X$\+contraadjusted finitely filtered quasi-coherent graded
$\cB^*$\+module.

 Furthermore, we have $\gr_F^*\P^*\simeq
\Cohom_X(\gr_F^*\cB^*,\Q^*)$ for the graded contraherent cosheaf
$\Q^*=\gr_F^0\P^*$ with antilocal grading components.
 The commutative diagram~\eqref{X-cta-filtered-naive-co-contra}
now tells us that the associated bigraded $X$\+contraadjusted
quasi-coherent $\gr_F^*\cB^*$\+module $\gr_F^*\M^*$ to $(\M^*,F)$ 
corresponds to the antilocal contraherent bigraded
$\gr_F^*\cB^*$\+module $\Cohom_X(\gr_F^*\cB^*,\Q^*)$ under
the bigraded version of the category equivalence from
Lemma~\ref{quasi-algebra-underived-naive-co-contra}(a).

 In view of commutativity of
the diagram~\eqref{X-cta-underived-naive-co-contra} in
Lemma~\ref{quasi-algebra-underived-naive-co-contra}(a), we can conclude
from Lemma~\ref{internal-Hom-and-Cohom-correspond-under-co-contra}(a)
that the $X$\+contraadjusted quasi-coherent bigraded
$\gr_F^*\cB^*$\+mod\-ule $\gr_F^*\M^*$ has the form $\gr_F^*\M^*\simeq
\cHom_{\cO_X}(\gr_F^*\cB^*,\N^*)$.
 Here $\N^*=\cO_X\ocn_X\Q^*$ is the graded quasi-coherent sheaf
with contraadjusted grading components correponding to the graded
contraherent cosheaf $\Q^*$ under the equivalence of categories
from~\cite[Lemma~4.8.2]{Pcosh}.
 It follows that the $X$\+contraadjusted quasi-coherent graded
$\cB^*$\+module $\M^*$ is thick and $F$ is its canonical filtration
from Corollary~\ref{thick-qcoh-sheaves-are-filtered}(b).
\end{proof}

 Now we can prove the promised cotorsion version of
Corollary~\ref{thick-cta-qcoh-sheaves-are-filtered-by-cta}.

\begin{cor} \label{thick-cot-qcoh-sheaves-are-filtered-by-cot}
 Let $X$ be a quasi-compact semi-separated scheme and $\E$ be a finite
locally free sheaf on~$X$.
 Let $\cB^*$ be the quasi-coherent graded algebra
$\cB^*=\bigwedge_X^*(\E)$ over~$X$.
 Endow the quasi-coherent graded algebra $\cB^*$ with the decreasing
filtration induced by the grading,
$F^i\cB^*=\bigoplus_{j=i}^\infty\cB^j$ for all $i\in\boZ$.
 Let $\M^*$ be a thick $X$\+cotorsion quasi-coherent graded
$\cB^*$\+module on~$X$.
 Consider the finite decreasing filtration $F$ on $\M^*$ defined in
Corollary~\ref{thick-qcoh-sheaves-are-filtered}(a) or~(b).
 Then the quasi-coherent sheaves\/ $\gr_F^i\M^n$ on $X$ are
cotorsion for all $i$, $n\in\boZ$.
 In other words, $(\M^*,F)$ is an $X$\+cotorsion finitely filtered
quasi-coherent graded module over $(\cB^*,F)$.
\end{cor}

\begin{proof}
 Similarly to the proof of
Corollary~\ref{cot-qcoh-modules-over-flfrqa-cor}, we use the idea
of the argument from~\cite[Section~8]{Pphil} involving the reduction
of a problem about quasi-coherent sheaves to a related question
concerning contraherent cosheaves, which is local and therefore
easier to resolve.
 In view of the proof of
Corollary~\ref{thick-qcoh-sheaves-are-filtered}(b),
it suffices to consider the filtration $F$ on $\M^*$ from
Corollary~\ref{thick-qcoh-sheaves-are-filtered}(b).
 Let $\P^*$ be the antilocal $X$\+locally cotorsion contraherent graded
$\cB^*$\+module on $X$ corresponding to $\M^*$ under the equivalence
of categories from the graded version of
Lemma~\ref{quasi-algebra-underived-naive-co-contra}(b).

 According to
Theorem~\ref{thick-corresponds-to-thick-under-co-contra},
the contraherent graded $\cB^*$\+module $\P^*$ is thick, and
the filtration $F$ on $\M^*$ corresponds to the canonical filtration $F$
on $\P^*$ from Corollary~\ref{thick-lcta-lct-cosheaves-are-filtered}(a)
under the equivalence of categories from the graded version of
the upper line of diagram~\eqref{X-cta-filtered-naive-co-contra} in
Theorem~\ref{quasi-algebra-filtered-naive-co-contra}.
 By Corollary~\ref{thick-lcta-lct-cosheaves-are-filtered}(b),
$(\P^*,F)$ is a finitely filtered $X$\+locally cotorsion contraherent
graded module over $(\cB^*,F)$; in other words, the contraherent
cosheaves $\gr_F^i\P^n$ on $X$ are locally cotorsion.

 The graded versions of the upper lines of the two
diagrams~\eqref{X-cta-filtered-naive-co-contra}
and~\eqref{X-cot-filtered-naive-co-contra} in
Theorem~\ref{quasi-algebra-filtered-naive-co-contra} form
a commutative diagram of exact category equivalences and exact,
fully faithful inclusion functors
\begin{equation} \label{X-cta-X-cot-filtered-naive-co-contra}
\begin{gathered}
 \qquad\xymatrix{
  \text{\llap{$\fHom_{\cB^*}^*(\cB^*,{-})\:$}}
  (\cB^*,F)\Qcoh^{\ff,X\dcta} \ar@{=}[rr]
  && (\cB^*,F)\Ctrh_\al^\ff
  \text{\rlap{$\,\,:\!\cB^*\ocn_{\cB^*}{-}$}} \\
  \text{\llap{$\fHom_{\cB^*}^*(\cB^*,{-})\:$}}
  (\cB^*,F)\Qcoh^{\ff,X\dcot} \ar@{=}[rr] \ar@{>->}[u]
  && (\cB^*,F)\Ctrh_\al^{\ff,X\dlct}
  \text{\rlap{$\,\,:\!\cB^*\ocn_{\cB^*}{-}$}} \ar@{>->}[u]
 }
\end{gathered}
\end{equation}
 Since $(\M^*,F)$ and $(\P^*,F)$ correspond to each other under
the category equivalence in the upper line
of~\eqref{X-cta-X-cot-filtered-naive-co-contra}, the lower line
of~\eqref{X-cta-X-cot-filtered-naive-co-contra} is also known to be
a category equivalence, and $(\P^*,F)$ belongs to the essential image of
the rightmost vertical functor, we can conclude that $(\M^*,F)$ belongs
to the essential image of the leftmost vertical functor.
 Thus $(\M^*,F)$ is an $X$\+cotorsion finitely filtered quasi-coherent
graded module over $(\cB^*,F)$, as desired.
\end{proof}

\subsection{Exact DG-category of thick CDG-modules}
\label{exact-dg-categ-of-thick-cdg-mods}
 We start with a preliminary lemma.

\begin{lem} \label{thick-submodule-and-quotient-lemma}
 Let $R$ be a commutative ring, $E$ be a finitely generated projective
$R$\+module, and $B^*$ be the graded ring $B^*=\bigwedge_R^*(E)$.
 Let\/ $0\rarrow M^*\rarrow L^*\rarrow K^*\rarrow0$ be a short exact
sequence of graded $B^*$\+modules.
 Assume that the graded $B^*$\+module $L^*$ is thick (in the sense of
Section~\ref{thick-graded-modules-subsecn}) and the graded
$R$\+module map $M^*/(F^1B^*\cdot M^*)\rarrow L^*/(F^1B^*\cdot L^*)$
induced by the graded $B^*$\+module map $M^*\rarrow L^*$ is injective.
 Then the graded $B^*$\+modules $M^*$ and $K^*$ are thick, too.
\end{lem}

\begin{proof}
 We have $M^*/(F^1B^*\cdot M^*)=B^*/F^1B^*\ot_{B^*}M^*$ and
$L^*/(F^1B^*\cdot L^*)=B^*/F^1B^*\ot_{B^*}L^*$.
 Since the graded $B^*$\+module $L^*$ is thick, we also have
$\Tor_1^{B^*}(B^*/F^1B^*,L^*)=0$ by (the graded version of)
Lemma~\ref{very-weakly-relatively-projective-lemma}(1).
 So, in view of the long exact sequence of
$\Tor_*^{B^*}(B^*/F^1B^*,{-})$, it follows from the assumptions of
the lemma that $\Tor_1^{B^*}(B^*/F^1B^*,K^*)=0$.
 Hence the graded $B^*$\+module $K^*$ is thick by
Lemma~\ref{very-weakly-relatively-projective-lemma}(1).
 Applying Lemma~\ref{very-weakly-relatively-projective-lemma}(2)
or Lemma~\ref{thick-graded-modules-closed-under-co-kernels}, we
conclude that the graded $B^*$\+module $M^*$ is thick, too.
\hbadness=1325
\end{proof}

 Let $(R,\g,\widetilde\g)$ be a twisted Lie algebroid (as defined in
Section~\ref{twisted-lie-algebroids-subsecn}).
 Assume that $\g$~is a finitely generated projective $R$\+module,
and let $B^\cu=C^\cu_R(\g,\widetilde\g)$ be the Chevalley--Eilenberg
CDG\+ring of the twisted Lie algebroid $(R,\g,\widetilde\g)$, as
constructed in Section~\ref{chevalley-eilenberg-cdg-ring-subsecn}.
 So the underlying graded ring of $B^\cu$ is $B^*=
\Lambda_R^*(\Hom_R(\g,R))$, and the definition of a thick graded module
from Section~\ref{thick-graded-modules-subsecn} is applicable to
graded $B^*$\+modules.

\begin{lem} \label{thick-G-plus-lemma}
 Let $M^*$ be a graded $B^*$\+module.
 Consider the CDG\+module $G^+(M^*)$ over $B^\cu$, as constructed
in Lemma~\ref{G-plus-G-minus-functors-for-modules}.
 Assume that the graded $B^*$\+module $G^+(M^*)^\#$ is thick.
 Then the graded $B^*$\+module $M^*$ is thick, too.
\end{lem}

\begin{proof}
 We apply Lemma~\ref{thick-submodule-and-quotient-lemma} to the short
exact sequence of graded $B^*$\+modules $0\rarrow M^*\rarrow
G^+(M^*)^\#\rarrow M^*[-1]\rarrow0$ from
Lemma~\ref{G-plus-G-minus-functors-for-modules}.
 In order to check the injectivity assumption of
Lemma~\ref{thick-submodule-and-quotient-lemma}, we use the explicit
construction of the CDG\+module $G^+(M^*)$ spelled out in~\cite[proof
of Theorem~3.6]{Pkoszul} and mentioned in the proof of
Lemma~\ref{G-plus-G-minus-functors-for-modules}.

 According to this construction, the graded $B^*$\+module $G^+(M^*)^\#$
is the set of all formal expressions $x'+dx''$ with $x'$, $x''\in M$.
 The graded $B^*$\+submodule $M^*\subset G^+(M^*)$ is the set of
all expressions $x'+d(0)\in G^+(M^*)^\#$.
 Denote by $K^*\subset G^+(M^*)$ the set of all expresstions
$0+dx''\in G^+(M^*)$.
 Then $K^*$ is \emph{not} a graded $B^*$\+submodule but only
a homogeneous abelian subgroup in $G^+(M^*)^\#$, and we have a direct
sum decomposition of graded abelian groups $G^+(M^*)^\#=M^*\oplus K^*$.

 Following the construction, the action of the graded ring $B^*$ on
$G^+(M^*)$ is given by the formula $b(x'+dx'')=(bx'-(-1)^{|b|}d_B(b)x'')
+(-1)^{|b|}d(bx'')$ for all $b\in B^{|b|}$ and $x'$, $x''\in M^*$.
 Here $d_B$~denotes the differential in the CDG\+ring $B^\cu=
(B^*,d_B,h_B)$.
 It is clear from the formula that $F^iB^*\cdot K^*\subset
(F^{i+1}B^*\cdot M^*)\oplus K^*\subset G^+(M^*)^\#$ for all $i\ge0$.
 In particular, we have $F^1B^*\cdot K^*\subset
(F^2B^*\cdot M^*)\oplus K^*\subset (F^1B^*\cdot M^*)\oplus K^*$.
 Hence $M^*\cap (F^1B^*\cdot G^+(M^*)^\#)=F^1B^*\cdot M^*$ and
the map $M^*/(F^1B^*\cdot M^*)\rarrow G^+(M^*)^\#/
(F^1B^*\cdot G^+(M^*)^\#)$ is injective, as desired.
\end{proof}

 Consider the abelian DG\+category $B^\cu\bModl$ of left
CDG\+modules over $B^\cu$, as per the discussion in
Section~\ref{abelian-and-exact-dg-categs-of-cdg-modules-subsecn}.
 We will say that a CDG\+module $M^\cu$ over $B^\cu$ is \emph{thick}
if its underlying graded $B^*$\+module $M^*$ is thick.
 Denote the full DG\+subcategory of thick CDG\+modules by
$B^\cu\bModl_\bth\subset B^\cu\bModl$.

\begin{cor} \label{exact-dg-category-of-thick-cdg-modules}
 Let $(R,\g,\widetilde\g)$ be a twisted Lie algebroid with a finitely
generated projective $R$\+module\/~$\g$, and let
$B^\cu=C^\cu_R(\g,\widetilde\g)$ be its Chevalley--Eilenberg CDG\+ring.
 Then there is a commutative diagram of additive category equivalences
and fully faithful identity inclusion functors
\begin{equation} \label{Upsilon-thick-CDG-modules-diagram}
\begin{gathered}
 \xymatrix{
  B^*\Modl_\thk
  \ar@<0.4ex>[rr]^-{\Upsilon_{B^\subcu}}
  \ar@{>->}[d]
  && \sZ^0((B^\cu\bModl_\bth)^\bec) \ar@<0.4ex>@{-}[ll]
  \ar@{>->}[d] \\
  B^*\Modl \ar@<0.4ex>[rr]^-{\Upsilon_{B^\subcu}}
  && \sZ^0((B^\cu\bModl)^\bec) \ar@<0.4ex>@{-}[ll]
 }
\end{gathered}
\end{equation}
with horizontal double lines showing category equivalences and vertical
arrows with tails showing fully faithful inclusions.
 The category equivalence in the lower horizontal line in the one
from formula~\eqref{CDG-modules-over-CDG-ring-Upsilon-equivalence}
in Section~\ref{abelian-and-exact-dg-categs-of-cdg-modules-subsecn}.
 The DG\+category $B^\cu\bModl_\bth$ has a natural structure of exact
DG\+category such that the related exact category structure on
the additive category $B^*\Modl_\thk$ is inherited from the abelian
exact structure of the abelian category $B^*\Modl$.
\end{cor}

\begin{proof}
 One needs to use Lemma~\ref{thick-G-plus-lemma} in order to show that
upper horizontal functor is essentially surjective.
 In the exact category structure on the additive category
$\sZ^0(B^\cu\bModl_\bth)$, a short sequence is exact if and only if
the forgetful functor~$\#$ takes it to a short exact sequence in
$B^*\Modl_\thk$, or equivalently, in $B^*\Modl$.
 One easily checks that this is a DG\+compatible exact structure in
the sense of~\cite[Section~4.2]{Pedg} inducing the original exact
structure on the underlying category of graded objects $B^*\Modl_\thk$
(Lemma~\ref{G-plus-reflects-exactness} is helpful here).
 So~\cite[Theorem~4.17]{Pedg} is applicable.
\end{proof}

\subsection{Exact DG-categories of thick quasi-coherent CDG-modules}
\label{exact-dg-categories-of-thick-qcoh-subsecn}
 Let $X$ be a scheme and $\E$ be a finite locally free sheaf on~$X$.
 As in Section~\ref{thick-quasi-coherent-modules-subsecn}, we consider
the quasi-coherent graded algebra $\cB^*=\bigwedge^*_X(\E)$ over~$X$.

 The definitions of the exact category of thick $X$\+contraadjusted
quasi-coherent graded modules $\cB^*\Qcoh^{X\dcta}_\thk$ and the exact
category of thick $X$\+cotorsion quasi-coherent graded modules
$\cB^*\Qcoh^{X\dcot}_\thk$ over $\cB^*$ were given in
Section~\ref{qcoh-assoc-graded-are-cta-subsecn}.
 Notice that we do \emph{not} consider $\cB^*$\+cotorsion
quasi-coherent graded $\cB^*$\+modules, as these coincide with
the $X$\+cotorsion quasi-coherent graded $\cB^*$\+modules (for
a quasi-compact semi-separated scheme~$X$) by
Corollary~\ref{cot-qcoh-modules-over-flfrqa-cor}.

 Let $(\g,\widetilde\g)$ be a quasi-coherent twisted Lie algebroid
over a scheme~$X$ (as defined in
Section~\ref{twisted-lie-algebroids-subsecn}).
 Assume that $\g$~is a finite locally free sheaf on $X$, and let
$\cB^\cu=\cC^\cu_X(\g,\widetilde\g)$ be the Chevalley--Eilenberg
quasi-coherent CDG\+quasi-algebra of the quasi-coherent twisted Lie
algebroid $(\g,\widetilde\g)$, as constructed in
Section~\ref{chevalley-eilenberg-qcoh-cdg-quasi-algebra-subsecn}.
 So the underlying quasi-coherent graded ring of $\cB^\cu$ is
$\cB^*=\Lambda_X^*(\cHom_X(\g,\cO_X))$, and the definition of
a thick quasi-coherent graded module from
Section~\ref{thick-quasi-coherent-modules-subsecn} is applicable to
quasi-coherent graded $\cB^*$\+modules on~$X$.

 Consider the abelian DG\+category $\cB^\cu\bQcoh$ of quasi-coherent
left CDG\+modules over $\cB^\cu$, as per
Corollary~\ref{abelian-dg-category-of-qcoh-cdg-modules}.
 We will say that a quasi-coherent CDG\+module $\M^\cu$ over $\cB^\cu$
is \emph{thick} if its underlying quasi-coherent graded $\cB^*$\+module
$\M^*$ is thick.
 Denote the full DG\+subcategory of thick quasi-coherent CDG\+modules
by $\cB^\cu\bQcoh_\bth\subset\cB^\cu\bQcoh$.

 See Section~\ref{exact-dg-of-cta-cot-qcoh-cdg} for a discussion of
$X$\+contraadjusted and $X$\+cotorsion quasi-coherent CDG\+modules
over~$\cB^\cu$.
 We will say that an $X$\+contraadjusted or $X$\+cotorsion
quasi-coherent CDG\+module over~$\cB^\cu$ is \emph{thick} if it is
thick as a quasi-coherent CDG\+module over~$\cB^\cu$.
 Let us introduce notation for the corresponding full DG\+subcategories
in $\cB^\cu\bQcoh$, similar to the notation for the underlying
categories of graded objects:
\begin{align*}
 \cB^\cu\bQcoh^{X\dcta}_\bth &=
 \cB^\cu\bQcoh_\bth\cap\cB^\cu\bQcoh^{X\dcta}, \\
 \cB^\cu\bQcoh^{X\dcot}_\bth &=
 \cB^\cu\bQcoh_\bth\cap\cB^\cu\bQcoh^{X\dcot}.
\end{align*}

\begin{cor} \label{exact-dg-categories-of-qcoh-thick-cdg-modules}
 Let $X$ be a quasi-compact semi-separated scheme and
$(\g,\widetilde\g)$ be a quasi-coherent twisted Lie algebroid over~$X$.
 Assume that\/ $\g$~is a finite locally free sheaf on $X$, and let
$\cB^\cu=\cC^\cu_X(\g,\widetilde\g)$ be the related Chevalley--Eilenberg
quasi-coherent CDG\+quasi-algebra over~$X$.
 Then there are commutative diagrams of additive category equivalences
and fully faithful identity inclusion functors
\begin{equation} \label{Upsilon-qcoh-X-cot-thick-diagram}
\begin{gathered}
 \xymatrix{
  \cB^*\Qcoh^{X\dcot}_\thk
  \ar@<0.4ex>[rr]^-{\Upsilon_{\cB^\subcu}^\qc}
  \ar@{>->}[d]
  && \sZ^0((\cB^\cu\bQcoh^{X\dcot}_\bth)^\bec) \ar@<0.4ex>@{-}[ll]
  \ar@{>->}[d] \\
  \cB^*\Qcoh^{X\dcot}
  \ar@<0.4ex>[rr]^-{\Upsilon_{\cB^\subcu}^\qc}
  && \sZ^0((\cB^\cu\bQcoh^{X\dcot})^\bec) \ar@<0.4ex>@{-}[ll]
 }
\end{gathered}
\end{equation}
\begin{equation} \label{Upsilon-qcoh-X-cta-thick-diagram}
\begin{gathered}
 \xymatrix{
  \cB^*\Qcoh^{X\dcta}_\thk
  \ar@<0.4ex>[rr]^-{\Upsilon_{\cB^\subcu}^\qc}
  \ar@{>->}[d]
  && \sZ^0((\cB^\cu\bQcoh^{X\dcta}_\bth)^\bec) \ar@<0.4ex>@{-}[ll]
  \ar@{>->}[d] \\
  \cB^*\Qcoh^{X\dcta}
  \ar@<0.4ex>[rr]^-{\Upsilon_{\cB^\subcu}^\qc}
  && \sZ^0((\cB^\cu\bQcoh^{X\dcta})^\bec) \ar@<0.4ex>@{-}[ll]
 }
\end{gathered}
\end{equation}
\begin{equation} \label{Upsilon-qcoh-X-thick-diagram}
\begin{gathered}
 \xymatrix{
  \cB^*\Qcoh_\thk \ar@<0.4ex>[rr]^-{\Upsilon_{\cB^\subcu}^\qc}
  \ar@{>->}[d]
  && \sZ^0((\cB^\cu\bQcoh_\bth)^\bec) \ar@<0.4ex>@{-}[ll]
  \ar@{>->}[d] \\
  \cB^*\Qcoh \ar@<0.4ex>[rr]^-{\Upsilon_{\cB^\subcu}^\qc}
  && \sZ^0((\cB^\cu\bQcoh)^\bec) \ar@<0.4ex>@{-}[ll]
 }
\end{gathered}
\end{equation}
with horizontal double lines showing category equivalences and
vertical arrows with tails showing fully faithful inclusions.
 The category equivalence in the lower horizontal line of the latter
diagram~\eqref{Upsilon-qcoh-X-thick-diagram} is
the equivalence~\eqref{Upsilon-qcoh-equivalence}
from Corollary~\textup{\ref{abelian-dg-category-of-qcoh-cdg-modules}},
while the category equivalences in the lower horizontal lines of
the former two
diagrams~\textup{(\ref{Upsilon-qcoh-X-cot-thick-diagram}\+-%
\ref{Upsilon-qcoh-X-cta-thick-diagram})} are the ones from
diagram~\eqref{Upsilon-qcoh-diagram} in
Corollary~\textup{\ref{exact-dg-categories-of-qcoh-cta-cot-cdg-modules}}.
 The DG\+categories $\cB^\cu\bQcoh_\bth$, \
$\cB^\cu\bQcoh^{X\dcta}_\bth$, and $\cB^\cu\bQcoh^{X\dcot}_\bth$ have
natural structures of exact DG\+categories such that the related exact
category structures on the additive categories in the rightmost column
of the diagram agree with the exact category structures on the additive
categories in the leftmost column of the diagram inherited from
the abelian exact structure of $\cB^*\Qcoh$.
\end{cor}

\begin{proof}
 The proof is similar to those of
Corollaries~\ref{exact-dg-categories-of-qcoh-cta-cot-cdg-modules}
and~\ref{exact-dg-category-of-thick-cdg-modules}.
 One needs to use Lemmas~\ref{qcoh-G-plus-periodicity-lemma}(a\+-b)
and~\ref{thick-G-plus-lemma} in order to show that the upper horizontal
functors are essentially surjective.
\end{proof}

\subsection{Exact DG-categories of thick contraherent CDG-modules}
\label{exact-dg-categories-of-thick-lcth-subsecn}
 Let $X$ be a scheme with an open covering $\bW$ and $\E$ be a finite
locally free sheaf on~$X$.
 As above, we consider the quasi-coherent graded algebra
$\cB^*=\bigwedge^*_X(\E)$ over~$X$.

 The definitions of the exact category of thick $\bW$\+locally
contraherent graded $\cB^*$\+modules $\cB^*\Lcth_\bW^\thk$ and
the exact category of thick $X$\+locally cotorsion $\bW$\+locally
contraherent graded $\cB^*$\+modules $\cB^*\Lcth_\bW^{X\dlct,\thk}$
on $X$ were given in
Sections~\ref{thick-loc-contraherent-modules-subsecn}
and~\ref{assoc-graded-modules-are-cta-cot}.
 Notice that we do \emph{not} consider $\cB^*$\+locally cotorsion
$\bW$\+locally contraherent graded $\cB^*$\+modules, as these coincide
with the $X$\+locally cotorsion $\bW$\+locally contraherent graded
$\cB^*$\+modules by Corollary~\ref{lct-lcth-modules-over-flfrqa-cor}.

 Let $(\g,\widetilde\g)$ be a quasi-coherent twisted Lie algebroid
over a scheme~$X$ (as defined in
Section~\ref{twisted-lie-algebroids-subsecn}).
 Assume that $\g$~is a finite locally free sheaf on $X$, and let
$\cB^\cu=\cC^\cu_X(\g,\widetilde\g)$ be the Chevalley--Eilenberg
quasi-coherent CDG\+quasi-algebra of the quasi-coherent twisted Lie
algebroid $(\g,\widetilde\g)$, as constructed in
Section~\ref{chevalley-eilenberg-qcoh-cdg-quasi-algebra-subsecn}.
 Then the definition of a thick $\bW$\+locally contraherent graded
module is applicable to $\bW$\+locally contraherent graded
$\cB^*$\+modules on~$X$ (cf.\
Section~\ref{exact-dg-categories-of-thick-qcoh-subsecn}).

 Consider the exact DG\+categories $\bLcth_\bW^{X\dlct}\subset
\cB^\cu\bLcth_\bW$ of $\bW$\+locally contraherent CDG\+modules over
$\cB^\cu$, as per
Corollary~\ref{exact-dg-categories-of-lcth-cdg-modules}.
 We will say that a $\bW$\+locally contraherent CDG\+module
(respectively, an $X$\+locally cotorsion $\bW$\+locally contraherent
CDG\+module) $\P^\cu$ over $\cB^\cu$ is \emph{thick} if its underlying
$\bW$\+locally contraherent graded $\cB^*$\+module is thick.
 Denote the full DG\+subcategories of thick $\bW$\+locally
contraherent CDG\+modules by
\begin{align*}
 \cB^\cu\bLcth_\bW^\bth &\subset \cB^\cu\bLcth_\bW, \\
 \cB^\cu\bLcth_\bW^{X\dlct,\bth} &\subset
 \cB^\cu\bLcth_\bW^{X\dlct}.
\end{align*}

\begin{cor} \label{exact-dg-categories-of-lcth-thick-cdg-modules}
 Let $X$ be a scheme with an open covering\/ $\bW$ and
$(\g,\widetilde\g)$ be a quasi-coherent twisted Lie algebroid over~$X$.
 Assume that\/ $\g$~is a finite locally free sheaf on $X$, and let
$\cB^\cu=\cC^\cu_X(\g,\widetilde\g)$ be the related Chevalley--Eilenberg
quasi-coherent CDG\+quasi-algebra over~$X$.
 Then there are commutative diagrams of additive category equivalences
and fully faithful identity inclusion functors
\begin{equation} \label{Upsilon-lcth-X-lct-thick-diagram}
\begin{gathered}
 \xymatrix{
  \cB^*\Lcth_\bW^{X\dlct,\thk}
  \ar@<0.4ex>[rr]^-{\Upsilon_{\cB^\subcu}^\ct}
  \ar@{>->}[d]
  && \sZ^0((\cB^\cu\bLcth_\bW^{X\dlct,\bth})^\bec) \ar@<0.4ex>@{-}[ll]
  \ar@{>->}[d] \\
  \cB^*\Lcth_\bW^{X\dlct}
  \ar@<0.4ex>[rr]^-{\Upsilon_{\cB^\subcu}^\ct}
  && \sZ^0((\cB^\cu\bLcth_\bW^{X\dlct})^\bec) \ar@<0.4ex>@{-}[ll]
 }
\end{gathered}
\end{equation}
\begin{equation} \label{Upsilon-lcth-X-thick-diagram}
\begin{gathered}
 \xymatrix{
  \cB^*\Lcth_\bW^\thk \ar@<0.4ex>[rr]^-{\Upsilon_{\cB^\subcu}^\ct}
  \ar@{>->}[d]
  && \sZ^0((\cB^\cu\bLcth_\bW^\bth)^\bec) \ar@<0.4ex>@{-}[ll]
  \ar@{>->}[d] \\
  \cB^*\Lcth_\bW \ar@<0.4ex>[rr]^-{\Upsilon_{\cB^\subcu}^\ct}
  && \sZ^0((\cB^\cu\bLcth_\bW)^\bec) \ar@<0.4ex>@{-}[ll]
 }
\end{gathered}
\end{equation}
with horizontal double lines showing category equivalences and
vertical arrows with tails showing fully faithful inclusions.
 The category equivalences in the lower horizontal lines are the ones
from diagram~\eqref{Upsilon-lcth-diagram}
in Corollary~\textup{\ref{exact-dg-categories-of-lcth-cdg-modules}}.
 The DG\+categories $\cB^\cu\bLcth_\bW^\bth$ and
$\cB^\cu\bLcth_\bW^{X\dlct,\bth}$ have natural structures of exact
DG\+categories such that the related exact category structures on
the additive categories in the rightmost column of the diagram agree
with the exact category structures on the additive categories in
the leftmost column of the diagram defined in
Sections~\textup{\ref{thick-loc-contraherent-modules-subsecn}}
and~\textup{\ref{assoc-graded-modules-are-cta-cot}}.
\end{cor}

\begin{proof}
 The proof is similar to those of
Corollaries~\ref{exact-dg-categories-of-lcth-cdg-modules}
and~\ref{exact-dg-category-of-thick-cdg-modules}.
 One needs to use Lemmas~\ref{G-plus-periodicity-lemma}(a\+-b)
and~\ref{thick-G-plus-lemma} in order to show that the upper horizontal
functors are essentially surjective.
\end{proof}

\subsection{Exact DG-categories of thick antilocal contraherent
CDG-modules} \label{exact-dg-categories-of-thick-antilocal-subsecn}
 Let $X$ be a quasi-compact semi-separated scheme with an open covering
$\bW$ and $\E$ be a finite locally free sheaf on~$X$.
 As above, we consider the quasi-coherent graded algebra
$\cB^*=\bigwedge^*_X(\E)$ over~$X$.

 The definitions of the exact category of thick antilocal contraherent
graded $\cB^*$\+modules $\cB^*\Ctrh_\al^\thk$ and the exact category
of thick antilocal $X$\+locally cotorsion contraherent graded
$\cB^*$\+modules $\cB^*\Ctrh_\al^{X\dlct,\thk}$
were given in Section~\ref{qcoh-assoc-graded-are-cta-subsecn}.

 Let $(\g,\widetilde\g)$ be a quasi-coherent twisted Lie algebroid
over a scheme~$X$.
 Assume that $\g$~is a finite locally free sheaf on $X$, and let
$\cB^\cu=\cC^\cu_X(\g,\widetilde\g)$ be the related Chevalley--Eilenberg
quasi-coherent CDG\+quasi-algebra over~$X$, as in
Sections~\ref{exact-dg-categories-of-thick-qcoh-subsecn}\+-%
\ref{exact-dg-categories-of-thick-lcth-subsecn}.

 Consider the exact DG\+categories $\cB^\cu\bCtrh_\al^{X\dlct}\subset
\cB^\cu\bCtrh_\al$ of antilocal contraherent CDG\+modules over
$\cB^\cu$, as per
Corollary~\ref{exact-dg-categories-of-antilocal-ctrh-cdg-modules}.
 We will say that an antilocal contraherent CDG\+module
(respectively, an antilocal $X$\+locally cotorsion contraherent
CDG\+module) $\P^\cu$ over $\cB^\cu$ is \emph{thick} if its underlying
locally contraherent graded $\cB^*$\+module is thick.
 Denote the full DG\+subcategories of thick antilocal contraherent
CDG\+modules by
\begin{align*}
 \cB^\cu\bCtrh_\al^\bth &\subset \cB^\cu\bCtrh_\al, \\
 \cB^\cu\bCtrh_\al^{X\dlct,\bth} &\subset
 \cB^\cu\bCtrh_\al^{X\dlct}.
\end{align*}

\begin{cor} \label{exact-dg-categories-of-antilocal-thick-cdg-modules}
 Let $X$ be a quasi-compact semi-separated scheme and
$(\g,\widetilde\g)$ be a quasi-coherent twisted Lie algebroid over~$X$.
 Assume that\/ $\g$~is a finite locally free sheaf on $X$, and let
$\cB^\cu=\cC^\cu_X(\g,\widetilde\g)$ be the related Chevalley--Eilenberg
quasi-coherent CDG\+quasi-algebra over~$X$.
 Then there are commutative diagrams of additive category equivalences
and fully faithful identity inclusion functors
\begin{equation} \label{Upsilon-ctrh-al-X-lct-thick-diagram}
\begin{gathered}
 \xymatrix{
  \cB^*\Ctrh_\al^{X\dlct,\thk}
  \ar@<0.4ex>[rr]^-{\Upsilon_{\cB^\subcu}^\ct}
  \ar@{>->}[d]
  && \sZ^0((\cB^\cu\bCtrh_\al^{X\dlct,\bth})^\bec) \ar@<0.4ex>@{-}[ll]
  \ar@{>->}[d] \\
  \cB^*\Ctrh_\al^{X\dlct}
  \ar@<0.4ex>[rr]^-{\Upsilon_{\cB^\subcu}^\ct}
  && \sZ^0((\cB^\cu\bCtrh_\al^{X\dlct})^\bec) \ar@<0.4ex>@{-}[ll]
 }
\end{gathered}
\end{equation}
\begin{equation} \label{Upsilon-ctrh-al-X-thick-diagram}
\begin{gathered}
 \xymatrix{
  \cB^*\Ctrh_\al^\thk \ar@<0.4ex>[rr]^-{\Upsilon_{\cB^\subcu}^\ct}
  \ar@{>->}[d]
  && \sZ^0((\cB^\cu\bCtrh_\al^\bth)^\bec) \ar@<0.4ex>@{-}[ll]
  \ar@{>->}[d] \\
  \cB^*\Ctrh_\al \ar@<0.4ex>[rr]^-{\Upsilon_{\cB^\subcu}^\ct}
  && \sZ^0((\cB^\cu\bCtrh_\al)^\bec) \ar@<0.4ex>@{-}[ll]
 }
\end{gathered}
\end{equation}
with horizontal double lines showing category equivalences and
vertical arrows with tails showing fully faithful inclusions.
 The category equivalences in the lower horizontal lines are
the ones from diagrams~\textup{(\ref{Upsilon-A-lct-antiloc-diagram}\+-%
\ref{Upsilon-X-lcta-antiloc-diagram})} in Corollary~%
\textup{\ref{exact-dg-categories-of-antilocal-ctrh-cdg-modules}}.
 The DG\+categories $\cB^\cu\bCtrh_\al^\bth$ and
$\cB^\cu\bCtrh_\al^{X\dlct,\bth}$ have natural structures of exact
DG\+categories such that the related exact category structures on
the additive categories in the rightmost column of the diagram agree
with the exact category structures on the additive categories in
the leftmost column of the diagram defined in
Section~\textup{\ref{qcoh-assoc-graded-are-cta-subsecn}}.
\end{cor}

\begin{proof}
 The proof is similar to those of
Corollaries~\ref{exact-dg-categories-of-antilocal-ctrh-cdg-modules},
\ref{exact-dg-category-of-thick-cdg-modules},
and~\ref{exact-dg-categories-of-lcth-thick-cdg-modules}.
 One needs to use Lemmas~\ref{G-plus-periodicity-lemma}(a\+-b),
\ref{antilocal-G-plus-periodicity-lemma}, and~\ref{thick-G-plus-lemma}
in order to show that the upper horizontal functors are
essentially surjective.
\end{proof}

\Section{Reduced Contraderived Categories of CDG-Modules}
\label{cdg-reduced-contraderived-secn}

 In this section, unless otherwise mentioned, we consider
a quasi-compact semi-separated scheme $X$ and a quasi-coherent twisted
Lie algebroid $(\g,\widetilde\g)$ over~$X$ (in the sense of
Section~\ref{twisted-lie-algebroids-subsecn}) such that
the quasi-coherent sheaf~$\g$ on $X$ is finite locally free.
 We denote by $\cB^\cu=\cC^\cu_X(\g,\widetilde\g)$ the related
Chevalley--Eilenberg quasi-coherent CDG\+quasi-algebra over $X$,
as constructed in
Section~\ref{chevalley-eilenberg-qcoh-cdg-quasi-algebra-subsecn}.

 We also denote by $\E$ a finite locally free sheaf on $X$, which is
eventually supposed to be $\E=\cHom_{\cO_X}(\g,\cO_X)$.
 When we are not interested in the whole quasi-coherent
CDG\+quasi-algebra $\cB^\cu$, but only in its underlying quasi-coherent
graded algebra $\cB^*$ over $X$, we put $\cB^*=\bigwedge_X^*(\E)$,
in the notation of Section~\ref{skew-symmetric-tensors-subsecn}.

\subsection{Coderived categories of thick CDG-modules}
\label{coderived-of-thick-cdg-modules-subsecn}
 We start with the following couple of simple lemmas.
 The construction of the direct image functor for quasi-coherent modules
over quasi-coherent quasi-algebras was explained in
Section~\ref{direct-images-of-A-co-sheaves-subsecn};
see formula~\eqref{qcoh-A-modules-direct-image}.
 The inverse images with respect to open immersions of schemes are
also mentioned in the same section.

\begin{lem} \label{injective-qcoh-modules-strongly-antilocal}
 Let $X$ be a quasi-compact quasi-separated scheme with a finite affine
open covering $X=\bigcup_{\alpha=1}^N U_\alpha$ and $\cA$ be
a quasi-coherent quasi-algebra over~$X$.
 Denote by $j_\alpha\:U_\alpha\rarrow X$ the open immersion morphisms.
 Then a quasi-coherent $\cA$\+module $\J$ on $X$ is injective if and
only if $\J$ is a direct summand of a quasi-coherent $\cA$\+module of
the form\/ $\bigoplus_{\alpha=1}^N j_\alpha{}_*\J_\alpha$, where
$\J_\alpha$ are injective quasi-coherent graded modules over
the quasi-coherent quasi-algebras $j_\alpha^*\cA$ over~$U_\alpha$.
\end{lem}

\begin{proof}
 This is an easy version of the results of
Sections~\ref{antilocality-of-X-contraadjusted-subsecn}\+-%
\ref{antilocality-of-A-cotorsion-subsecn}, or more specifically of
Theorems~\ref{qcomp-qsep-very-flaproj-complete-cotorsion-pair-thm}(c),
\ref{qcomp-qsep-flaproj-complete-cotorsion-pair-thm}(c),
and~\ref{qcomp-qsep-A-flat-complete-cotorsion-pair-thm}(c).
 The ``if'' assertion holds, since the direct image functor
$j_\alpha{}_*\:j_\alpha^*\cA\Qcoh\rarrow\cA\Qcoh$ is right adjoint
to the inverse image functor $j_\alpha^*\:\cA\Qcoh\rarrow
j_\alpha^*\cA\Qcoh$, and the latter functor is exact.
 The proof of the ``if'' is based on the construction from the proof
of Lemma~\ref{X-contraadjusted-modules-cogenerating-class}.
 For every index~$\alpha$, consider the quasi-coherent
$\cA|_{U_\alpha}$\+module $j_\alpha^*\J$ on $U_\alpha$, and pick
an injective quasi-coherent $\cA|_{U_\alpha}$\+module $\J_\alpha$
together with an injective morphism of quasi-coherent
$\cA|_{U_\alpha}$\+modules $j_\alpha^*\J\rarrow\J_\alpha$ on~$U_\alpha$.
 This simply means that $\J_\alpha(U_\alpha)$ is an injective
$\cA(U_\alpha)$ module such that the $\cA(U_\alpha)$\+module
$\J(U_\alpha)$ is a submodule in $\J_\alpha(U_\alpha)$.
 By adjunction, we obtain a morphism of quasi-coherent $\cA$\+modules
$\J\rarrow\bigoplus_\alpha j_\alpha{}_*\J_\alpha$, which is injective
because it is injective in restriction to every~$U_\alpha$.
 As $\J$ is an injective object of $\cA\Qcoh$, it follows that $\J$ is
a direct summand of $\bigoplus_\alpha j_\alpha{}_*\J_\alpha$.
\end{proof}

 The definition of a thick quasi-coherent graded module over $\cB^*$
was given in Section~\ref{thick-quasi-coherent-modules-subsecn}.

\begin{lem} \label{injective-qcoh-graded-modules-are-thick}
 Let $X$ be a quasi-compact semi-separated scheme and $\E$ be a finite
locally free sheaf on~$X$.
 Consider the quasi-coherent graded algebra $\cB^*=\bigwedge_X^*(\E)$
on~$X$.
 Then every injective quasi-coherent graded $\cB^*$\+module on $X$
is thick.
\end{lem}

\begin{proof}
 The affine case is obvious: for any commutative ring $R$ and any
finitely generated projective $R$\+module $E$, any injective graded
module over the graded ring $B^*=\bigwedge_R^*(E)$ is thick by
the definition (see
Lemma~\ref{very-weakly-relatively-injective-lemma}(1\+-2)
and Section~\ref{thick-graded-modules-subsecn}).
 To deduce the similar assertion for a quasi-compact semi-separated
scheme $X$, pick a finite affine open covering
$X=\bigcup_{\alpha=1}^N U_\alpha$.
 Notice that any injective quasi-coherent graded $\cB^*$\+module $\J^*$
is a direct summand of a quasi-coherent graded $\cB^*$\+module
$\bigoplus_\alpha j_\alpha{}_*\J_\alpha^*$ for some injective
quasi-coherent graded modules $\J_\alpha$ over $\cB^*|_{U_\alpha}$,
by the graded version of
Lemma~\ref{injective-qcoh-modules-strongly-antilocal}.
 It remains to refer to Lemma~\ref{qcoh-thick-direct-image}.
\end{proof}

 We refer to Section~\ref{derived-second-kind-subsecn}
for the definitions of the coderived category in the sense of
Positselski $\sD^\co(\bE)$ and the coderived category in the sense
of Becker $\sD^\bco(\bE)$ for an exact DG\+category~$\bE$.

\begin{cor} \label{thick-cdg-becker-coderived-equivs-cor}
 Let $X$ be a quasi-compact semi-separated scheme and
$(\g,\widetilde\g)$ be a quasi-coherent twisted Lie algebroid over~$X$.
 Assume that\/ $\g$~is a finite locally free sheaf on $X$, and let
$\cB^\cu=\cC^\cu_X(\g,\widetilde\g)$ be the related Chevalley--Eilenberg
quasi-coherent CDG\+quasi-algebra over~$X$.
 Then there is a commutative diagram of triangulated category
equivalences of the Becker coderived categories induced by the identity
inclusions of exact/abelian DG\+categories
\begin{equation} \label{thick-cdg-becker-coderived-equivs-diagram}
\begin{gathered}
 \xymatrix{
 \sD^\bco(\cB^\cu\bQcoh_\bth) \ar@<-2pt>[r] \ar@{-}@<-2pt>[d]
  & \sD^\bco(\cB^\cu\bQcoh) \ar@{-}@<-2pt>[l]
  \ar@{-}@<-2pt>[d] \\
 \sD^\bco(\cB^\cu\bQcoh_\bth^{X\dcta}) \ar@<-2pt>[r]
  \ar@<-2pt>[u] \ar@{-}@<-2pt>[d]
  & \sD^\bco(\cB^\cu\bQcoh^{X\dcta}) \ar@{-}@<-2pt>[l]
  \ar@<-2pt>[u] \ar@{-}@<-2pt>[d] \\
  \sD^\bco(\cB^\cu\bQcoh_\bth^{X\dcot}) \ar@<-2pt>[r]
  \ar@<-2pt>[u]
  & \sD^\bco(\cB^\cu\bQcoh^{X\dcot}) \ar@{-}@<-2pt>[l]
  \ar@<-2pt>[u]
 }
\end{gathered}
\end{equation}
\end{cor}

\begin{proof}
 The abelian DG\+category $\cB^\cu\bQcoh$ is described in
Corollary~\ref{abelian-dg-category-of-qcoh-cdg-modules}, while
the two nonabelian exact DG\+categories appearing in the rightmost
column of the diagram are described in
Corollary~\ref{exact-dg-categories-of-qcoh-cta-cot-cdg-modules}.
 We do not consider $\cB^*$\+cotorsion CDG\+modules over $\cB^\cu$,
as these coincide with the $X$\+cotorsion CDG\+modules by
Corollary~\ref{cot-qcoh-modules-over-flfrqa-cor}.
 The three exact DG\+categories appearing in the leftmost column of
the diagram are described in
Corollary~\ref{exact-dg-categories-of-qcoh-thick-cdg-modules}.

 The vertical equivalences in the rightmost column of the diagram were
proved in Corollary~\ref{qcoh-cta-cot-cdg-becker-coderived-equivalence}.
 According to the proof of that corollary, the full subcategories
$\cB^*\Qcoh^{X\dcot}\subset\cB^*\Qcoh^{X\dcta}$ are coresolving
subcategories closed under direct summands in the abelian category
$\cB^*\Qcoh$; in particular, these full subcategories contain all
the injective quasi-coherent graded $\cB^*$\+modules.
 In view of Lemma~\ref{injective-qcoh-graded-modules-are-thick},
the full subcategories $\cB^*\Qcoh^{X\dcot}_\thk\subset
\cB^*\Qcoh^{X\dcta}_\thk\subset\cB^*\Qcoh_\thk$ also contain all
the injective objects of $\cB^*\Qcoh$.
 By Lemma~\ref{thick-graded-modules-closed-under-co-kernels} or
the discussions in Sections~\ref{thick-quasi-coherent-modules-subsecn}
and~\ref{exact-dg-categories-of-thick-qcoh-subsecn}, all the mentioned
full subcategories are coresolving in $\cB^*\Qcoh$.
 So the class of injective objects in every one of these exact
subcategories coincides with the class of injective objects
in $\cB^*\Qcoh$.

 Thus an object of every one of the exact DG\+subcategories involved
is graded-injective if and only if it is graded-injective in
the abelian DG\+category $\cB^\cu\bQcoh$.
 It follows that an object of every one of these exact DG\+categories
is Becker-coacyclic if and only if it is Becker-coacyclic as an object
of $\cB^\cu\bQcoh$.
 The argument finishes similarly to the proof of
Corollary~\ref{qcoh-cta-cot-cdg-becker-coderived-equivalence}.
\end{proof}

\begin{cor} \label{thick-cdg-positselski-coderived-equiv-cor}
 Let $X$ be a quasi-compact semi-separated scheme and
$(\g,\widetilde\g)$ be a quasi-coherent twisted Lie algebroid over~$X$.
 Assume that\/ $\g$~is a finite locally free sheaf on $X$, and let
$\cB^\cu=\cC^\cu_X(\g,\widetilde\g)$ be the related Chevalley--Eilenberg
quasi-coherent CDG\+quasi-algebra over~$X$.
 Then the inclusion of exact/abelian DG\+categories
$\cB^\cu\bQcoh_\bth\rarrow\cB^\cu\bQcoh$ induces an equivalence of
the Positselski coderived categories
$$
 \sD^\co(\cB^\cu\bQcoh_\bth)\simeq\sD^\co(\cB^\cu\bQcoh).
$$
\end{cor}

\begin{proof}
 For the assertion of this corollary to make sense, it is important
that the full subcategory $\cB^*\Qcoh_\thk$ is closed under infinite
direct sums in $\cB^*\Qcoh$ (while the other coresolving subcategories
in $\cB^*\Qcoh$ appearing in the context of
Corollary~\ref{thick-cdg-becker-coderived-equivs-cor} do not have
infinite direct sums, so the construction of the Positselski coderived
category does not make sense for the respective exact DG\+categories).

 In order to prove the corollary, it suffices to point out that
the full subcategory $\cB^*\Qcoh_\thk$ is coresolving in $\cB^*\Qcoh$,
while the full DG\+subcategory $\cB^\cu\bQcoh_\bth$ is closed under
twists and infinite direct sums in $\cB^\cu\bQcoh$.
 So the dual version of
Proposition~\ref{second-kind-infinite-resolutions}(a) is applicable.
\end{proof}

 The next corollary is similar to
Corollary~\ref{qcoh-cta-cdg-absolute-derived-equivalence}.

\begin{cor} \label{thick-cdg-qcoh-cta-abs-derived-equiv-cor}
 Let $X$ be a quasi-compact semi-separated scheme and
$(\g,\widetilde\g)$ be a quasi-coherent twisted Lie algebroid over~$X$.
 Assume that\/ $\g$~is a finite locally free sheaf on $X$, and let
$\cB^\cu=\cC^\cu_X(\g,\widetilde\g)$ be the related Chevalley--Eilenberg
quasi-coherent CDG\+quasi-algebra over~$X$.
 Then the inclusion of exact DG\+categories
$\cB^\cu\bQcoh^{X\dcta}_\bth\rarrow\cB^\cu\bQcoh_\bth$ induces
an equivalence of the absolute derived categories
$$
 \sD^\abs(\cB^\cu\bQcoh^{X\dcta}_\bth)\simeq
 \sD^\abs(\cB^\cu\bQcoh_\bth).
$$
\end{cor}

\begin{proof}
 It is clear from the upper line of
diagram~\eqref{Upsilon-qcoh-X-cta-thick-diagram}
in Corollary~\ref{exact-dg-categories-of-qcoh-thick-cdg-modules}
that $\cB^\cu\bQcoh^{X\dcta}_\bth$ is a strict exact DG\+subcategory
in $\cB^\cu\bQcoh_\bth$.
 It was mentioned in the proof of
Corollary~\ref{thick-cdg-becker-coderived-equivs-cor} that the full
subcategory $\cB^*\Qcoh^{X\dcta}_\thk$ is coresolving in
$\cB^*\Qcoh$, hence also in $\cB^*\Qcoh_\thk$.
 In view of the graded version of
Lemma~\ref{X-A-cta-cot-coresolving-coresolution-dimension}(b) and
the fact that the coresolution dimension does not depend on the choice
of a coresolution, it follows that the coresolution dimension of any
object of $\cB^*\Qcoh_\thk$ with respect to $\cB^*\Qcoh^{X\dcta}_\thk$ 
does not exceed the number $N$ of affine open subschemes in any given
affine open covering of~$X$ (cf.\ the dual version
of~\cite[Corollary~A.5.5]{Pcosh}).
 Thus the dual version of
Proposition~\ref{second-kind-finite-resolutions}(a) is applicable.
\end{proof}

 For a similar equivalence involving also the absolute derived category
of thick $X$\+cotorsion quasi-coherent CDG\+modules over $\cB^\cu$
(which holds under more restrictive assumptions on the scheme~$X$), see
Corollary~\ref{thick-cdg-qcoh-cot-abs-derived-equiv-cor} below.

\subsection{Contraderived categories of thick CDG-modules}
\label{contraderived-of-thick-cdg-modules-subsecn}
 We start with a simple module-theoretic lemma.
 The lemma does not mention cotorsion graded $B^*$\+modules, as these
coincide with graded $B^*$\+modules with $R$\+cotorsion grading
components by Proposition~\ref{cotorsion-modules-over-fgpgr-ring-prop}.

\begin{lem} \label{enough-thick-contraadjusted-or-cotorsion}
 Let $R$ be a commutative ring and $E$ be a finitely generated
projective $R$\+module.
 Consider the graded ring $B^*=\bigwedge_R^*(E)$.
 In this setting: \par
\textup{(a)} Let $P^*$ be a graded $B^*$\+module such that
the $R$\+module $P^i$ is contraadjusted for all $i\in\boZ$.
 Then there exists a short exact sequence of graded $B^*$\+modules\/
$0\rarrow C^*\rarrow Q^*\rarrow P^*\rarrow0$ such that the $R$\+modules
$C^i$ and $Q^i$ are contraadjusted for all $i\in\boZ$, while
the graded $B^*$\+module $Q^*$ is thick. \par
\textup{(b)} Let $P^*$ be a graded $B^*$\+module such that
the $R$\+module $P^i$ is cotorsion for all $i\in\boZ$.
 Then there exists a short exact sequence of graded $B^*$\+modules\/
$0\rarrow C^*\rarrow Q^*\rarrow P^*\rarrow0$ such that the $R$\+modules
$C^i$ and $Q^i$ are cotorsion for all $i\in\boZ$, while
the graded $B^*$\+module $Q^*$ is thick.
\end{lem}

\begin{proof}
 The point is that every flat graded $B^*$\+module is thick,
because it satisfies condition~(1) or~(2) of
Lemma~\ref{very-weakly-relatively-projective-lemma}.
 By the graded version of Theorem~\ref{flat-cotorsion-pair-complete}(a),
for any graded $B^*$\+module $P^*$ there exists a short exact sequence
of graded $B^*$\+modules $0\rarrow C^*\rarrow Q^*\rarrow P^*\rarrow0$
with a flat graded $B^*$\+module $Q^*$ and a cotorsion graded
$B^*$\+module~$C^*$.
 It remains to recall that any cotorsion graded $B^*$\+module is
a cotorsion graded $R$\+module by the graded version of
Lemma~\ref{restriction-coextension-injective-cotorsion}(a); so
the graded $R$\+modules $C^i$ are cotorsion, hence contraadjusted.
 Finally, one needs to use the fact that the classes of contraadjusted
$R$\+modules and cotorsion $R$\+modules are closed under extensions
in the respective module categories.
 Therefore, if the $R$\+modules $P^i$ are cotorsion or contraadjusted,
then so are the $R$\+modules~$Q^i$.
\end{proof}

\begin{lem} \label{enough-thick-antilocal-in-lcth}
 Let $X$ be a quasi-compact semi-separated scheme with an open
covering\/ $\bW$ and $\E$ be a finite locally free sheaf on~$X$.
 Let $\cB^*$ be the quasi-coherent graded algebra
$\cB^*=\bigwedge_X^*(\E)$ over~$X$.
 Then \par
\textup{(a)} for every\/ $\bW$\+locally contraherent graded
$\cB^*$\+module\/ $\P^*$ over $X$, there exists a thick antilocal
contraherent graded $\cB^*$\+module\/ $\Q^*$ over $X$ together with
an admissible epimorphism\/ $\Q^*\rarrow\P^*$ in $\cB^*\Lcth_\bW$; \par
\textup{(b)} for every $X$\+locally cotorsion\/ $\bW$\+locally 
contraherent graded $\cB^*$\+module\/ $\P^*$ over $X$, there exists
a thick antilocal $X$\+locally cotorsion contraherent graded
$\cB^*$\+module\/ $\Q^*$ over $X$ together with an admissible
epimorphism\/ $\Q^*\rarrow\P^*$ in $\cB^*\Lcth_\bW^{X\dlct}$.
\end{lem}

\begin{proof}
 The argument is based on
Lemma~\ref{enough-thick-contraadjusted-or-cotorsion} and
the graded version of the \v Cech
resolution~\eqref{lcth-sheaf-of-rings-cech-resolution} from
Section~\ref{cech-subsecn}.
 Let us prove part~(a).
 Similarly to the second paragraph of the proof of
Lemma~\ref{direct-image-of-A-lin-from-affine-Ext-orthogonality}(b),
choose a finite affine open covering $X=\bigcup_\alpha U_\alpha$
of the scheme $X$ subordinate to $\bW$, and denote by
$j_\alpha\:U_\alpha\rarrow X$ the open immersion morphisms.
 Then the natural morphism
$\bigoplus_\alpha j_\alpha{}_!j_\alpha^!\P^*\rarrow\P^*$ is
an admissible epimorphism in $\cB^*\Lcth_\bW$.
 For every index~$\alpha$, choose a short exact sequence of graded
$\cB^*(U_\alpha)$\+modules $0\rarrow C_\alpha^*\rarrow Q_\alpha^*\rarrow
\P^*[U_\alpha]\rarrow0$ such that the graded $B^*(U_\alpha)$\+module
$Q_\alpha^*$ is thick and the $\cO_X(U_\alpha)$\+modules $C_\alpha^i$
and $Q_\alpha^i$ are contraadjusted for all $i\in\boZ$, as per
Lemma~\ref{enough-thick-contraadjusted-or-cotorsion}(a).
 Let $\Q_\alpha^*$ be the thick contraherent graded
$\cB^*|_{U_\alpha}$\+module on $U_\alpha$ corresponding to
the thick graded $\cB^*(U_\alpha)$\+module~$Q_\alpha^*$ with
$\cO_X(U_\alpha)$\+contraadjusted grading components (cf.\
Lemma~\ref{thick-graded-modules-co-extension-of-scalars}(b)).
 Then we have an admissible epimorphism $\Q_\alpha^*\rarrow
j_\alpha^!\P_\alpha^*$ in $\cB^*|_{U_\alpha}\Ctrh$ and an admissible
epimorphism $\bigoplus_\alpha j_\alpha{}_!\Q_\alpha^*\rarrow
\bigoplus_\alpha j_\alpha{}_!j_\alpha^!\P^*$ in $\cB^*\Lcth_\bW$.
 The composition $\bigoplus_\alpha j_\alpha{}_!\Q_\alpha^*\rarrow
\bigoplus_\alpha j_\alpha{}_!j_\alpha^!\P^*\rarrow\P^*$ is
an admissible epimorphism in $\cB^*\Lcth_\bW$, and
$\Q^*=\bigoplus_\alpha j_\alpha{}_!\Q_\alpha^*$ is an antilocal
contraherent graded $\cB^*$\+module on~$X$.
 By Lemma~\ref{lcth-thick-direct-image}, the contraherent graded
$\cB^*$\+module $\Q^*$ on $X$ is thick.
 The proof of part~(b) is similar.
\end{proof}

 The following corollary is similar to
Corollary~\ref{contraderived-indep-of-covering-or-antilocal}.

\begin{cor} \label{thick-abs-ctr-derived-indep-of-covering-or-al}
  Let $X$ be a quasi-compact semi-separated scheme with an open
covering\/ $\bW$ and $(\g,\widetilde\g)$ be a quasi-coherent twisted
Lie algebroid over~$X$.
 Assume that\/ $\g$~is a finite locally free sheaf on $X$, and let
$\cB^\cu=\cC^\cu_X(\g,\widetilde\g)$ be the related Chevalley--Eilenberg
quasi-coherent CDG\+quasi-algebra over~$X$.
 Let\/ $\st=\abs$, $\ctr$, or\/~$\bctr$ be an absolute derived or
contraderived category symbol.
 In this context: \par
\textup{(a)} The inclusions of exact DG\+categories
$\cB^\cu\bCtrh_\al^\bth\rarrow\cB^\cu\bCtrh^\bth\rarrow
\cB^\cu\bLcth_\bW^\bth$ induce equivalences of triangulated categories
$$
 \sD^\st(\cB^\cu\bCtrh_\al^\bth)\simeq\sD^\st(\cB^\cu\bCtrh^\bth)
 \simeq\sD^\st(\cB^\cu\bLcth_\bW^\bth).
$$ \par
\textup{(b)} The inclusions of exact DG\+categories
$\cB^\cu\bCtrh^{X\dlct,\bth}_\al\rarrow\cB^\cu\bCtrh^{X\dlct,\bth}
\rarrow\cB^\cu\bLcth_\bW^{X\dlct,\bth}$ induce equivalences of
triangulated categories
$$
 \sD^\st(\cB^\cu\bCtrh^{X\dlct,\bth}_\al)\simeq
 \sD^\st(\cB^\cu\bCtrh^{X\dlct,\bth})\simeq
 \sD^\st(\cB^\cu\bLcth_\bW^{X\dlct,\bth}).
$$
\end{cor}

\begin{proof}
 We do not mention the $\cB^*$\+locally cotorsion $\bW$\+locally
contraherent CDG\+mod\-ules over $\cB^\cu$, as these coincide with
the $X$\+locally cotorsion $\bW$\+locally contraherent CDG\+modules
by Corollary~\ref{lct-lcth-modules-over-flfrqa-cor}.
 It is clear from
Corollaries~\ref{exact-dg-categories-of-lcth-thick-cdg-modules}\+-%
\ref{exact-dg-categories-of-antilocal-thick-cdg-modules} that
all the full DG\+subcategories in question are strict exact
DG\+subcategories.
 All the related full subcategories of underlying graded objects are
resolving in each other by Lemma~\ref{enough-thick-antilocal-in-lcth}
and the graded version of
Lemma~\ref{antilocal-resolving-resolution-dimension},
and also closed under direct summands in each other.
 Moreover, in view of the resolution dimension assertions in the graded
version of Lemma~\ref{antilocal-resolving-resolution-dimension} and
the fact that the resolution dimension does not depend on the choice of
a resolution, it follows that, in each of the cases~(a\+-b),
the respective resolution dimension does not exceed $N-1$, where $N$ is
the number of affine open subschemes in any given affine open covering
of~$X$.
 Therefore, Proposition~\ref{second-kind-finite-resolutions}(a,c,d)
is applicable. \emergencystretch=1em
\end{proof}

\begin{cor} \label{thick-cdg-contraderived-equiv-cor}
 Let $X$ be a quasi-compact semi-separated scheme with an open
covering\/ $\bW$ and $(\g,\widetilde\g)$ be a quasi-coherent twisted
Lie algebroid over~$X$.
 Assume that\/ $\g$~is a finite locally free sheaf on $X$, and let
$\cB^\cu=\cC^\cu_X(\g,\widetilde\g)$ be the related Chevalley--Eilenberg
quasi-coherent CDG\+quasi-algebra over~$X$.
 Let\/ $\st=\ctr$ or\/~$\bctr$ be a (Positselski or Becker)
contraderived category symbol.
 Then there are commutative diagrams of triangulated equivalences of
the contraderived categories induced by the identity inclusions of
exact DG\+categories
\begin{equation} \label{thick-lcta-cdg-contrader-equivs-diagram}
\begin{gathered}
 \xymatrix{
  \sD^\st(\cB^\cu\bLcth_\bW^\bth) \ar@<-2pt>[r] \ar@{-}@<-2pt>[d]
  & \sD^\st(\cB^\cu\bLcth_\bW) \ar@{-}@<-2pt>[l]
  \ar@{-}@<-2pt>[d] \\
  \sD^\st(\cB^\cu\bCtrh^\bth) \ar@<-2pt>[r]
  \ar@<-2pt>[u] \ar@{-}@<-2pt>[d]
  & \sD^\st(\cB^\cu\bCtrh) \ar@{-}@<-2pt>[l]
  \ar@<-2pt>[u] \ar@{-}@<-2pt>[d] \\
  \sD^\st(\cB^\cu\bCtrh_\al^\bth) \ar@<-2pt>[r]
  \ar@<-2pt>[u]
  & \sD^\st(\cB^\cu\bCtrh_\al) \ar@{-}@<-2pt>[l]
  \ar@<-2pt>[u]
 }
\end{gathered}
\end{equation}
\begin{equation} \label{thick-X-lct-cdg-contrader-equivs-diagram}
\begin{gathered}
 \xymatrix{
  \sD^\st(\cB^\cu\bLcth_\bW^{X\dlct,\bth})
  \ar@<-2pt>[r] \ar@{-}@<-2pt>[d]
  & \sD^\st(\cB^\cu\bLcth_\bW^{X\dlct}) \ar@{-}@<-2pt>[l]
  \ar@{-}@<-2pt>[d] \\
  \sD^\st(\cB^\cu\bCtrh^{X\dlct,\bth}) \ar@<-2pt>[r]
  \ar@<-2pt>[u] \ar@{-}@<-2pt>[d]
  & \sD^\st(\cB^\cu\bCtrh^{X\dlct}) \ar@{-}@<-2pt>[l]
  \ar@<-2pt>[u] \ar@{-}@<-2pt>[d] \\
  \sD^\st(\cB^\cu\bCtrh_\al^{X\dlct,\bth}) \ar@<-2pt>[r]
  \ar@<-2pt>[u]
  & \sD^\st(\cB^\cu\bCtrh_\al^{X\dlct}) \ar@{-}@<-2pt>[l]
  \ar@<-2pt>[u]
 }
\end{gathered}
\end{equation}
\end{cor}

\begin{proof}
 The vertical equivalences in the rightmost columns of both
the diagrams were proved in
Corollary~\ref{contraderived-indep-of-covering-or-antilocal},
while the vertical equivalences in the leftmost columns are provided by
Corollary~\ref{thick-abs-ctr-derived-indep-of-covering-or-al}.
 All the triangulated equivalences on the two diagrams can be obtained
by applying Proposition~\ref{second-kind-infinite-resolutions}.
 The point is that all the full exact DG\+subcategories in question are
strict, have twists and exact infinite products, and are closed under
infinite products and direct summands in each other.
 Furthermore, the related full subcategories of underlying graded
objects are resolving in each other by
Lemma~\ref{enough-thick-antilocal-in-lcth} together with
Lemma~\ref{thick-graded-modules-closed-under-co-kernels} or
with the discussions of the closedness properties in
Sections~\ref{thick-loc-contraherent-modules-subsecn}
and~\ref{assoc-graded-modules-are-cta-cot}\+-%
\ref{qcoh-assoc-graded-are-cta-subsecn}.
 So Proposition~\ref{second-kind-infinite-resolutions}(a\+-b)
is applicable.
\end{proof}

 Under more restrictive assumptions on the scheme $X$, an equivalence
between the triangulated categories appearing in parts~(a) and~(b)
of Corollary~\ref{thick-abs-ctr-derived-indep-of-covering-or-al} is
provided by
Corollary~\ref{thick-lcth-al-abs-ctr-derived-lcta-lct-equiv-cor} below.
 In fact, there is a three-dimensional commutative diagram of
triangulated equivalences uniting the two
diagrams~(\ref{thick-lcta-cdg-contrader-equivs-diagram}\+-
\ref{thick-X-lct-cdg-contrader-equivs-diagram}) from
Corollary~\ref{thick-cdg-contraderived-equiv-cor} in this case.

\subsection{Koszul coresolutions of trivial quasi-coherent CDG-modules}
\label{koszul-coresolutions-of-trivial-qcoh-cdg-mods}
 We start with a category-theoretic lemma.

\begin{lem} \label{filtered-by-coacyclic-is-coacyclic-lemma}
 Let\/ $\bE$ be an exact DG\+category.
 Assume that every countable sequence of admissible monomorphisms
$X_1\rarrow X_2\rarrow X_3\rarrow\dotsb$ has an inductive limit in
the exact category\/~$\sZ^0(\bE)$.
 Let\/ $0=Z_0\rarrow Z_1\rarrow Z_2\rarrow\dotsb$ be a countable sequence
of admissible monomorphisms in\/~$\sZ^0(\bE)$.
 Consider the inductive limit $Z=\varinjlim_n Z_n\in\sZ^0(\bE)$ and
the sequence of successive quotients $Z_n/Z_{n-1}=
\coker(Z_{n-1}\to Z_n)\in\sZ^0(\bE)$.
 In this setting: \par
\textup{(a)} Assume further that, for any sequence of (admissible) short
exact sequences\/ $0\rarrow C_n\rarrow D_n\rarrow E_n\rarrow0$, indexed
by the nonnegative integers $n\ge1$, and any sequence of termwise
admissible monomorphisms of short exact sequences\/ $(0\to C_n\to D_n
\to E_n\to0)\rarrow(0\to C_{n+1}\to D_{n+1} \to E_{n+1}\to0)$ in
the exact category\/ $\sZ^0(\bE^\bec)$, the short sequence of inductive
limits\/ $0\rarrow\varinjlim_n C_n\rarrow\varinjlim_n D_n\rarrow
\varinjlim_n E_n\rarrow0$ is exact in\/~$\sZ^0(\bE^\bec)$.
 Assuming additionally that infinite direct sums exist and are exact
in\/ $\bE$, if all the successive quotients $Z_n/Z_{n-1}$, \,$n\ge1$,
are Positselski-coacyclic objects of\/ $\bE$, then $Z$ is
a Positselski-coacyclic object of\/~$\bE$. \par
\textup{(b)} If all the successive quotients $Z_n/Z_{n-1}$, \,$n\ge1$,
are Becker-coacyclic objects of\/ $\bE$, then $Z$ is
a Becker-coacyclic object of\/~$\bE$.
\end{lem}

\begin{proof}
 One can easily see from the definition of an exact DG\+category and
the construction of the DG\+category\/ $\bE^\bec$ that existence of
the inductive limits of countable sequences of admissible monomorphisms
in $\sZ^0(\bE)$ implies existence of such inductive limits
in~$\sZ^0(\bE^\bec)$.
 Now, to prove part~(a), consider the telescope short sequence
\begin{equation} \label{telescope-sequence}
 0\lrarrow\bigoplus\nolimits_{n=1}^\infty Z_n
 \xrightarrow{\id-\textit{shift}}\bigoplus\nolimits_{n=1}^\infty Z_n
 \lrarrow\varinjlim\nolimits_n Z_n\lrarrow0
\end{equation}
in the exact category~$\sZ^0(\bE)$.
 To prove that the sequence~\eqref{telescope-sequence} is exact
in $\sZ^0(\bE)$, it suffices to check that the functor $\Phi_\bE$
takes it to an exact sequence in $\sZ^0(\bE^\bec)$.
 For this purpose, one observes that the functor $\Phi_\bE$, being
both a left and a right adjoint, preserves all inductive and
projective limits.
 Furthermore, the sequence
\begin{equation} \label{Phi-of-telescope-sequence}
 0\lrarrow\bigoplus\nolimits_{n=1}^\infty\Phi_\bE(Z_n)
 \xrightarrow{\id-\textit{shift}}
 \bigoplus\nolimits_{n=1}^\infty\Phi_\bE(Z_n)
 \lrarrow\varinjlim\nolimits_n\Phi_\bE(Z_n)\lrarrow0
\end{equation}
is the inductive limit of the sequence of split short exact sequences
$$
 0\lrarrow\bigoplus\nolimits_{i=1}^{n-1}\Phi_\bE(Z_i)
 \xrightarrow{\id-\textit{shift}}
 \bigoplus\nolimits_{i=1}^n\Phi_\bE(Z_i)\lrarrow\Phi_\bE(Z_n)\lrarrow0
$$
and termwise admissible monomorphisms between them
in~$\sZ^0(\bE^\bec)$.
 By assumption, it follows that the short
sequence~\eqref{Phi-of-telescope-sequence} is (admissible) exact in
$\sZ^0(\bE^\bec)$, hence the short sequence~\eqref{telescope-sequence}
is admissible exact in~$\sZ^0(\bE)$.

 Now if $Z_0=0$ and the successive quotients $Z_n/Z_{n-1}$ are
Positselski-coacyclic, then it follows easily by induction on~$n$
that all the objects $Z_n$ are Positselski-coacyclic (since
the totalizations of the short exact sequences $0\rarrow Z_{n-1}
\rarrow Z_n\rarrow Z_n/Z_{n-1}\rarrow0$ are absolutely acyclic, hence
Positselski-coacyclic, by the definition).
 Therefore, the countable direct sum $\bigoplus_{n=1}^\infty Z_n$
is Positselski-coacyclic, too.
 Since the totalization of~\eqref{telescope-sequence} is absolutely
acyclic, hence Positselski-coacyclic, by the definition, it follows
that the object $Z$ is Positselski-coacyclic as well.

 Under the exactness assumptions of part~(a), part~(b) can be proved by
the same argument as part~(a).
 Without such assumptions, the point is that any injective object of
the exact category $\sZ^0(\bE^\bec)$ sees the relevant short sequences
in $\sZ^0(\bE^\bec)$ as exact.
 One starts with showing that $\Hom_\bE^\bu(\varinjlim_n Z_n,\>J)
=\varprojlim_n\Hom_\bE^\bu(Z_n,J)$ in the category of complexes of
abelian groups for every object $J\in\bE$.
 The fact that the functor~$\Phi_\bE$ preserves inductive limits
can be used here together with the result of~\cite[Lemma~3.9]{Pedg}
(cf.~\cite[proof of Lemma~9.2]{Pedg}).

 Now let $J\in\bE$ be a graded-injective object.
 Then the morphisms of complexes of abelian groups $\Hom_\bE^\bu(Z_n,J)
\rarrow\Hom_\bE^\bu(Z_{n-1},J)$ are termwise surjective.
 This is also provable using~\cite[Lemma~3.9]{Pedg} together with
the assumption that the functor $\Phi_\bE$ preserves admissible
monomorphisms, which is a part of the definition of
an exact DG\+category.  
 It follows that the short sequence of complexes of abelian groups
\begin{multline} \label{telescope-sequence-of-Hom-complexes}
 0\lrarrow\Hom_\bE^\bu(\varinjlim\nolimits_n Z_n,\>J) \\
 \lrarrow\prod\nolimits_{n=1}^\infty\Hom_\bE^\bu(Z_n,J)
 \lrarrow\prod\nolimits_{n=1}^\infty\Hom_\bE^\bu(Z_n,J)
 \lrarrow0
\end{multline}
is exact.

 In view of the short exact sequences of complexes of abelian groups
$0\rarrow\Hom_\bE^\bu(Z_n/Z_{n-1},J) \rarrow\Hom_\bE^\bu(Z_n,J)\rarrow
\Hom_\bE^\bu(Z_{n-1},J)\rarrow0$, acyclicity of the complexes
$\Hom_\bE^\bu(Z_n/Z_{n-1},J)$ implies acyclicity
of the complexes $\Hom_\bE^\bu(Z_n,J)$, as one can easily prove
by induction on~$n$.
 Infinite products of acyclic complexes of abelian groups are acyclic,
and it follows from exactness
of~\eqref{telescope-sequence-of-Hom-complexes} that the complex
$\Hom_\bE^\bu(\varinjlim\nolimits_n Z_n,\>J)$ is acyclic, too.
\end{proof}

 Let $X$ be a scheme and $(\g,\widetilde\g)$ be a quasi-coherent
twisted Lie algebroid over~$X$ such that the quasi-coherent sheaf~$\g$
on $X$ is finite locally free.
 Denote by $\cB^\cu=\cC_X^\cu(\g,\widetilde\g)$ the related
Chevalley--Eilenberg quasi-coherent CDG\+quasi-algebra over~$X$.

 A quasi-coherent graded $\cB^*$\+module $\N^*$ on $X$ is said to be
\emph{trivial} (or more precisely, having a \emph{trivial action
of~$\cB^*$}) if $\N^*$ is annihilated by the action of $\cB^i$ for
all $i>0$.
 This means that, for every affine open subscheme $U\subset X$,
the graded module $\N^*(U)$ is annihilated by the action of
the grading components $\cB^i(U)$ of the graded ring $\cB^*(U)$
for all $i>0$.
 Clearly, it suffices to check this condition for affine open
subschemes $U$ belonging to any given affine open covering of
the scheme~$X$.

 A quasi-coherent CDG\+module $\N^\cu$ over $\cB^\cu$ is said to
be \emph{trivial} if its underlying quasi-coherent graded
$\cB^*$\+module $\N^*$ is trivial.
 The differentials in a trivial quasi-coherent CDG\+module over
$\cB^\cu$ have zero square, as the curvature elements
$h_U\in\cB^2(U)$ act by zero on~$\N^*(U)$.
 The change-of-connection elements $a_{VU}\in\cB^1(V)$ also act by
zero on $\N^*(V)$; and it follows that the datum of a structure of
a trivial quasi-coherent CDG\+module over $\cB^\cu$ on a given
graded quasi-coherent sheaf $\N^*$ on $X$ is equivalent to
the datum of a differential making $\N^*$ a complex of quasi-coherent
sheaves on~$X$.
 So the full DG\+subcategory in $\cB^\cu\bQcoh$ formed by the trivial
quasi-coherent CDG\+modules is naturally equivalent (in fact,
isomorphic) to the DG\+category of complexes of quasi-coherent
sheaves on~$X$.
 The same applies to the full DG\+subcategory in $\bQcohr\cB^\cu$
formed by the trivial quasi-coherent right CDG\+modules over~$\cB^\cu$.

 The quasi-coherent right CDG\+module
$\cC^X_\cu(\cA_X,\g,\widetilde\g)$ over the quasi-coherent
CDG\+quasi-algebra $\cB^\cu=\cC^\cu_X(\g,\widetilde\g)$ was
constructed in Section~\ref{chevalley-eilenberg-cdg-modules-subsecn}.
{\hbadness=1125\par}

\begin{lem} \label{trivial-qcoh-koszul-coresolution-lemma}
 Let $\N^\bu$ be a trivial quasi-coherent right CDG\+module
over~$\cB^\cu$.
 Then the tensor product
$\N^\bu\ot_{\cO_X}\cC^X_\cu(\cA_X,\g,\widetilde\g)$ is a quasi-coherent
right CDG\+module over $\cB^\cu$ endowed with a natural closed morphism
of quasi-coherent right CDG\+modules
$\N^\bu\rarrow\N^\bu\ot_{\cO_X}\cC^X_\cu(\cA_X,\g,\widetilde\g)$
whose cone is Positselski-coacyclic in the abelian DG\+category\/
$\bQcohr\cB^\cu$.
\end{lem}

\begin{proof}
 The quasi-coherent right CDG\+module structure on the tensor product
$\N^\bu\ot_{\cO_X}\cC^X_\cu(\cA_X,\g,\widetilde\g)$ is induced by
the quasi-coherent right CDG\+module structure on
the Chevalley--Eilenberg CDG\+module $\cC^X_\cu(\cA_X,\g,\widetilde\g)$.
 (Recall that the differentials on $\cC^X_\cu(\cA_X,\g,\widetilde\g)$
are left $\cA_X(\g,\widetilde\g)$\+linear, hence also
left $\cO_X$\+linear.)
 The map of quasi-coherent graded right modules
$\N^*\rarrow\N^*\ot_{\cO_X}\cC^X_*(\cA_X,\g,\widetilde\g)$ over
the quasi-coherent graded algebra $\cB^*$ on $X$ is obtained by taking
the tensor product of the graded quasi-coherent sheaf $\N^*$ with
the natural map of quasi-coherent graded quasi-modules
$\cO_X\rarrow\cC_*^X(\cA,\g,\widetilde\g)$ over~$X$.
 The latter map is the composition of the natural injective morphism of
quasi-coherent quasi-algebras $\cO_X\rarrow\cA_X(\g,\widetilde\g)$
and the direct summand/degree~$0$ component inclusion
$\cA_X(\g,\widetilde\g)=
\cA_X(\g,\widetilde\g)\ot_{\cO_X}\bigwedge^0_X(\g)\rarrow
\cA_X(\g,\widetilde\g)\ot_{\cO_X}\bigwedge^*_X(\g)=
\cC_*^X(\cA,\g,\widetilde\g)$.
 One can easily check that the resulting map is a closed morphism
$\N^\bu\rarrow\N^\bu\ot_{\cO_X}\cC^X_\cu(\cA_X,\g,\widetilde\g)$ is
of quasi-coherent right CDG\+modules over~$\cB^\cu$.

 The closed morphism of quasi-coherent CDG\+modules $\N^\bu\rarrow
\N^\bu\ot_{\cO_X}\cC^X_\cu(\cA_X,\g,\widetilde\g)$ is injective,
so showing that its cone is coacyclic in $\bQcohr\cB^\cu$ is equivalent
to checking that its cokernel is coacyclic in $\bQcohr\cB^\cu$.
 For this purpose, we use
Lemma~\ref{filtered-by-coacyclic-is-coacyclic-lemma}.
 Let $F$ be the natural increasing filtration on the quasi-coherent
quasi-algebra $\cA_X(\g,\widetilde\g)$ defined in
Section~\ref{enveloping-algebra-subsecn}; so we have
$\gr^F_*\cA_X(\g,\widetilde\g)\simeq\Sym_X^*(\g)$ by
Corollary~\ref{twisted-lie-algebroid-pbw-qcoh-cor}(b).
 Define a natural increasing filtration $F$ on the graded
quasi-coherent sheaf $\bigwedge_X^*(\g)$ by the rule
$F_n\bigl(\bigwedge_X^*(\g)\bigr)=\bigoplus_{i=0}^n\bigwedge_X^i(\g)$.
 Let the increasing filtration $F$ on the tensor product
$\cC_*^X(\cA,\g,\widetilde\g)=
\cA_X(\g,\widetilde\g)\ot_{\cO_X}\bigwedge^*_X(\g)$ be obtained
as the tensor product of the filtrations $F$ on
$\cA_X(\g,\widetilde\g)$ and $\bigwedge_X^*(\g)$.
 Finally, denote also by $F$ the increasing filtration on the tensor
product $\N^\bu\ot_{\cO_X}\cC^X_\cu(\cA_X,\g,\widetilde\g)$ induced
by the increasing filtration $F$ on the Chevalley--Eilenberg
quasi-coherent CDG\+module $\cC_\cu^X(\cA,\g,\widetilde\g)$.
{\hbadness=2100\par}

 One can easily check that $F$ is a filtration of
$\N^\bu\ot_{\cO_X}\cC^X_\cu(\cA_X,\g,\widetilde\g)$ by
quasi-coherent CDG\+submodules over~$\cB^\cu$.
 The associated graded quasi-coherent CDG\+module to the filtration $F$
on $\N^\bu\ot_{\cO_X}\cC^X_\cu(\cA_X,\g,\widetilde\g)$ is a trivial
right CDG\+module over~$\cB^\cu$ corresponding to the complex of
quasi-coherent sheaves
$$
 \N^\bu\ot_{\cO_X}\gr^F_*\cC^X_\cu(\cA_X,\g,\widetilde\g)
 \simeq\N^\bu\ot_{\cO_X}\Sym_X^*(\g)\ot_{\cO_X}
 \bigwedge\nolimits_X^*(\g).
$$
 Here $\N^\bu\ot_{\cO_X}\Sym_X^*(\g)\ot_{\cO_X}\bigwedge_X^*(\g)$ is
the tensor product of two complexes of quasi-coherent sheaves on $X$,
viz., the complex $\N^\bu$ and the Koszul complex
$$
 \Sym_X^*(\g)\ot_{\cO_X}\bigwedge\nolimits_X^*(\g)
$$
with its usual Koszul differential.

 Now the Koszul complex is bigraded; it has the homological grading,
corresponding to the homological grading on the Chevalley--Eilenberg
CDG\+module $\cC^X_\cu(\cA_X,\g,\widetilde\g)$, and the internal
grading, corresponding to the degrees of the filtration~$F$.
 The Koszul complex is the direct sum of its internal grading
components; the internal grading component of degree~$0$ is
isomorphic to $\cO_X$, while the internal grading components of
positive degrees are finite acyclic complexes of finite locally
free sheaves on~$X$.
 Any finite acyclic complex of finite locally free sheaves is
absolutely acyclic in the exact category of finite locally free sheaves,
hence also in the exact category of flat quasi-coherent sheaves on~$X$;
the tensor product of such a complex with any complex in $X\Qcoh$ is
absolutely acyclic in $X\Qcoh$.
 Any Positselski-coacyclic complex in $X\Qcoh$, viewed as
a trivial quasi-coherent CDG\+module over $\cB^\cu$, is
a Positselski-coacyclic quasi-coherent CDG\+module over~$\cB^\cu$.

 Thus the map
$\N^\bu\rarrow\N^\bu\ot_{\cO_X}\cC^X_\cu(\cA_X,\g,\widetilde\g)$ is
a closed isomorphism of quasi-coherent CDG\+modules $\N^\bu\simeq
F_0\bigl(\N^\bu\ot_{\cO_X}\cC^X_\cu(\cA_X,\g,\widetilde\g)\bigr)$,
while the successive quotient CDG\+module
$$
 F_n\bigl(\N^\bu\ot_{\cO_X}\cC^X_\cu(\cA_X,\g,\widetilde\g)\bigr)\big/
 F_{n-1}\bigl(\N^\bu\ot_{\cO_X}\cC^X_\cu(\cA_X,\g,\widetilde\g)\bigr)
$$
is Positselski-coacyclic (in fact, absolutely acyclic) in
$\bQcohr\cB^\cu$ for every $n\ge1$.
 By Lemma~\ref{filtered-by-coacyclic-is-coacyclic-lemma}, it follows
that the cone of the closed morphism
$\N^\bu\rarrow\N^\bu\ot_{\cO_X}\cC^X_\cu(\cA_X,\g,\widetilde\g)$
is Positselski-coacyclic in $\bQcohr\cB^\cu$.
\end{proof}

\subsection{Koszul resolutions of trivial contraherent CDG-modules}
\label{koszul-resolutions-of-trivial-ctrh-cdg-mods}
 Let $X$ be a scheme with an open covering $\bW$ and $(\g,\widetilde\g)$
be a quasi-coherent twisted Lie algebroid over~$X$ such that
the quasi-coherent sheaf~$\g$ on $X$ is finite locally free.
 Denote by $\cB^\cu=\cC_X^\cu(\g,\widetilde\g)$ the related
Chevalley--Eilenberg quasi-coherent CDG\+quasi-algebra over~$X$.

 A $\bW$\+locally contraherent graded $\cB^*$\+module $\Q^*$ on $X$
is said to be \emph{trivial} (or more precisely, having a \emph{trivial
action of~$\cB^*$}) if $\Q^*$ is annihilated by the action of $\cB^i$
for all $i>0$.
 This means that, for every affine open subscheme $U\subset X$,
the graded module $\Q^*[U]$ is annihilated by the action of
the grading components $\cB^i(U)$ of the graded ring $\cB^*(U)$
for all $i>0$.
 Clearly, it suffices to check this condition for affine open
subschemes $U$ belonging to any given affine open covering of
the scheme $X$ subordinate to~$\bW$.

 A $\bW$\+locally contraherent CDG\+module $\Q^\cu$ over $\cB^\cu$ is
said to be \emph{trivial} if its underlying $\bW$\+locally contraherent
graded $\cB^*$\+module $\Q^*$ is trivial.
 The differentials in a trivial $\bW$\+locally contraherent CDG\+module
over $\cB^\cu$ have zero square, as the curvature elements
$h_U\in\cB^2(U)$ act by zero on~$\Q^*[U]$.
 The change-of-connection elements $a_{VU}\in\cB^1(V)$ also act by
zero on $\Q^*[V]$; and it follows that the datum of a structure of
a trivial $\bW$\+locally contraherent CDG\+module over $\cB^\cu$ on
a given graded $\bW$\+locally contraherent cosheaf $\Q^*$ on $X$ is
equivalent to the datum of a differential making $\Q^*$ a complex of
$\bW$\+locally contraherent cosheaves on~$X$.
 So the full DG\+subcategory in $\cB^\cu\bLcth_\bW$ formed by
the trivial $\bW$\+locally contraherent CDG\+modules is naturally
equivalent (in fact, isomorphic) to the DG\+category of complexes
of $\bW$\+locally contraherent cosheaves on~$X$.

 An $X$\+locally cotorsion $\bW$\+locally contraherent graded
$\cB^*$\+module $\Q^*$ on $X$ is said to be \emph{trivial} if it is
trivial as a $\bW$\+locally contraherent graded $\cB^*$\+module.
 An $X$\+locally cotorsion $\bW$\+locally contraherent CDG\+module
$\Q^\cu$ over $\cB^\cu$ is said to be \emph{trivial} if its underlying
$X$\+locally cotorsion $\bW$\+locally contraherent graded
$\cB^*$\+module $\Q^*$ is trivial.
 The full DG\+subcategory in $\cB^\cu\bLcth_\bW^{X\dlct}$ formed by
the trivial $X$\+locally cotorsion $\bW$\+locally contaherent
CDG\+modules is naturally equivalent (in fact, isomorphic) to
the DG\+category of complexes of locally cotorsion $\bW$\+locally
contraherent cosheaves on~$X$.

 The construction in the following lemma uses the functor $\Cohom_X$
from Section~\ref{cohom-from-quasi-modules-subsecn}.
 Notice that the quasi-coherent quasi-algebra $\cA_X(\g,\widetilde\g)$
is very flat, and in fact, locally projective as a quasi-coherent sheaf
on $X$ with respect to its left (as well as right) $\cO_X$\+module
structure, as one can see from
Corollary~\ref{twisted-lie-algebroid-pbw-qcoh-cor}(b).

\begin{lem} \label{trivial-lcth-koszul-resolution-lemma}
\textup{(a)} Let\/ $\Q^\bu$ be a trivial\/ $\bW$\+locally contraherent
left CDG\+module over~$\cB^\cu$.
 Then the\/ $\Cohom$ cosheaf\/
$\Cohom_X(\cC^X_\cu(\cA_X,\g,\widetilde\g),\Q^\bu)$ is
a\/ $\bW$\+locally contraherent left CDG\+module over $\cB^\cu$ endowed
with a natural closed morphism of\/ $\bW$\+locally contraherent
CDG\+modules\/ $\Cohom_X(\cC^X_\cu(\cA_X,\g,\widetilde\g),\Q^\bu)
\rarrow\Q^\bu$ whose cone is Positselski-contraacyclic in the exact
DG\+category $\cB^\cu\bLcth_\bW$. \par
\textup{(b)} In the context of part~\textup{(a)}, if\/ $\Q^\bu$ is
a trivial $X$\+locally cotorsion\/ $\bW$\+locally contraherent
CDG\+module over $\cB^\cu$, then\/
$\Cohom_X(\cC^X_\cu(\cA_X,\g,\widetilde\g),\Q^\bu)$ is
an $X$\+locally cotorsion\/ $\bW$\+locally contraherent CDG\+module
over $\cB^\cu$ and the cone of the closed morphism\/
$\Cohom_X(\cC^X_\cu(\cA_X,\g,\widetilde\g),\Q^\bu)\rarrow\Q^\bu$ is
Positselski-contraacyclic in the exact DG\+category
$\cB^\cu\bLcth_\bW^{X\dlct}$.
\end{lem}

\begin{proof}
 This is dual-analogous to
Lemma~\ref{trivial-qcoh-koszul-coresolution-lemma}.
 Let us prove part~(a).

 First of all, the graded quasi-coherent sheaf
$\cC^X_*(\cA_X,\g,\widetilde\g)$ in its left (and well as right)
$\cO_X$\+module structure is very flat (in fact, locally projective)
by Corollary~\ref{twisted-lie-algebroid-pbw-qcoh-cor}(b); so the graded
$\Cohom$ cosheaf $\Cohom_X(\cC^X_*(\cA_X,\g,\widetilde\g),\Q^\bu)$ is
well-defined and $\bW$\+locally contraherent on~$X$.
 The $\bW$\+locally contraherent left CDG\+module structure on
the $\Cohom$ cosheaf $\Cohom_X(\cC^X_\cu(\cA_X,\g,\widetilde\g),\Q^\bu)$
is induced by the quasi-coherent right CDG\+module structure on
the Chevalley--Eilenberg CDG\+module $\cC^X_\cu(\cA_X,\g,\widetilde\g)$
(see Section~\ref{cdg-rings-cdg-modules-subsecn} for the discussion
of Hom CDG\+modules).
 The map of $\bW$\+locally contraherent graded $\cB^*$\+modules
$\Cohom_X(\cC^X_*(\cA_X,\g,\widetilde\g),\Q^*)\rarrow\Q^*$ on $X$
is obtained by applying the $\Cohom_X({-},\Q^*)$ functor to
the map of quasi-coherent graded quasi-modules
$\cO_X\rarrow\cC_*^X(\cA,\g,\widetilde\g)$ constructed in the proof
of Lemma~\ref{trivial-qcoh-koszul-coresolution-lemma}.
 The resulting map is a closed morphism
$\Cohom_X(\cC^X_\cu(\cA_X,\g,\widetilde\g),\Q^\bu)\rarrow\Q^\bu$
of $\bW$\+locally contraherent CDG\+modules over $\cB^\cu$, since
the map $\cO_X\rarrow\cC_*^X(\cA,\g,\widetilde\g)$ is a closed
morphism of quasi-coherent right CDG\+modules over~$\cB^\cu$.

 The closed morphism of $\bW$\+locally contraherent CDG\+modules
$$
 \Cohom_X(\cC^X_\cu(\cA_X,\g,\widetilde\g),\Q^\bu)\lrarrow\Q^\bu
$$
is an admissible epimorphism in $\sZ^0(\bW\bLcth_\bW)$, since
$\cO_X\rarrow\cC_*^X(\cA,\g,\widetilde\g)$ is an injective map
whose cokernel is a very flat (and in fact, locally projective)
graded quasi-coherent sheaf on~$X$.
 So checking that the cone of this closed morphism is contraacyclic
in $\cB^\cu\bLcth_\bW$ is equivalent to checking that its kernel is
contraacyclic in $\cB^\cu\bLcth_\bW$.
 For this purpose, we use the dual version of
Lemma~\ref{filtered-by-coacyclic-is-coacyclic-lemma}.

 Let $F$ be the increasing filtration on the quasi-coherent
CDG\+module $\cC^X_\cu(\cA_X,\g,\widetilde\g)$ constructed in
the proof of Lemma~\ref{trivial-qcoh-koszul-coresolution-lemma}.
 Denote also by $F$ the decreasing filtration on the $\Cohom$
cosheaf $\Cohom_X(\cC^X_\cu(\cA_X,\g,\widetilde\g),\Q^\bu)$
induced by the increasing filtration $F$ on the Chevalley--Eilenberg
quasi-coherent CDG\+module $\cC^X_\cu(\cA_X,\g,\widetilde\g)$; so
$$
 F^n\Cohom_X(\cC^X_\cu(\cA_X,\g,\widetilde\g),\Q^\bu)
 =\Cohom_X\bigl(\cC^X_\cu(\cA_X,\g,\widetilde\g)/
 F_{n-1}\cC^X_\cu(\cA_X,\g,\widetilde\g),\>\Q^\bu\bigr)
$$
for all $n\in\boZ$.
 
 Then one can easily check that the graded $\bW$\+locally contraherent
cosheaves $F^n\Cohom(\cC^X_\cu(\cA_X,\g,\widetilde\g),\Q^\bu)$ on $X$
have  unique structures of $\bW$\+locally contraherent CDG\+modules
over $\cB^\cu$ such that the natural maps
$F^{n+1}\Cohom(\cC^X_\cu(\cA_X,\g,\widetilde\g),\Q^\bu)\rarrow
\Cohom(\cC^X_\cu(\cA_X,\g,\widetilde\g),\Q^\bu)$ are closed morphisms
in the DG\+category $\cB^\cu\bLcth_\bW$, and in fact, admissible
monomorphisms in the exact category $\sZ^0(\cB^\cu\bLcth_\bW)$.
 Furthermore, the decreasing filtration $F$ on
$\Cohom(\cC^X_\cu(\cA_X,\g,\widetilde\g),\Q^\bu)$ is separated
and complete, that is
\begin{multline*}
 \Cohom_X(\cC^X_\cu(\cA_X,\g,\widetilde\g),\Q^\bu) \\
 = \varprojlim\nolimits_n
 \Cohom_X(\cC^X_\cu(\cA_X,\g,\widetilde\g),\Q^\bu)/
 F^n\Cohom_X(\cC^X_\cu(\cA_X,\g,\widetilde\g),\Q^\bu)
\end{multline*}
in $\sZ^0(\cB^\cu\bLcth_\bW)$, where
\begin{multline*}
 \Cohom_X(\cC^X_\cu(\cA_X,\g,\widetilde\g),\Q^\bu)/
 F^n\Cohom_X(\cC^X_\cu(\cA_X,\g,\widetilde\g),\Q^\bu) \\
 = \Cohom_X(F_{n-1}\cC^X_\cu(\cA_X,\g,\widetilde\g),\Q^\bu).
\end{multline*}
 The associated graded $\bW$\+locally contraherent CDG\+module to
the decreasing filtration $F$ on
$\Cohom_X(\cC^X_\cu(\cA_X,\g,\widetilde\g),\Q^\bu)$ is a trivial
CDG\+module over $\cB^\cu$ corresponding to the complex of
$\bW$\+locally contraherent cosheaves
$$
 \Cohom_X\Bigl(\Sym_X^*(\g)\ot_{\co_X}\bigwedge\nolimits_X^*(\g),
 \>\Q^\bu\Bigr).
$$
 Here the latter complex is the $\Cohom_X$ from the Koszul complex
of quasi-coherent sheaves $\Sym_X^*(\g)\ot_{\cO_X}\bigwedge_X^*(\g)$
(with its usual Koszul differential) into the complex of
$\bW$\+locally contraherent cosheaves $\Q^\bu$ on~$X$.

 Continuing to argue similarly to the proof of
Lemma~\ref{trivial-qcoh-koszul-coresolution-lemma}, it remains to
point out that any finite acyclic complex of finite locally free
sheaves on $X$ is absolutely acyclic in the exact category of very
flat quasi-coherent sheaves on~$X$; the $\Cohom$ from such a complex
into any complex in $X\Lcth_\bW$ is absolutely acyclic in $X\Lcth_\bW$.
 Any Positselski-contraacyclic complex in $X\Lcth_\bW$, viewed as
a trivial $\bW$\+locally contraherent CDG\+module over $\cB^\cu$,
is a Positselski-contraacyclic $\bW$\+locally contraherent CDG\+module
over~$\cB^\cu$.

 So the map
$\Cohom_X(\cC^X_\cu(\cA_X,\g,\widetilde\g),\Q^\bu)\rarrow\Q^\bu$
induces a closed isomorphism of $\bW$\+locally contraherent
CDG\+modules
$$
 \Cohom_X(\cC^X_\cu(\cA_X,\g,\widetilde\g),\Q^\bu)/
 F^1\Cohom_X(\cC^X_\cu(\cA_X,\g,\widetilde\g),\Q^\bu)
 \,\simeq\,\Q^\bu,
$$
while the successive quotients
$$
 F^n\Cohom_X(\cC^X_\cu(\cA_X,\g,\widetilde\g),\Q^\bu)
 /F^{n+1}\Cohom_X(\cC^X_\cu(\cA_X,\g,\widetilde\g),\Q^\bu)
$$
are Positselski-contraacyclic (in fact, absolutely acyclic) in
$\cB^\cu\bLcth_\bW$ for all $n\ge1$.
 By the dual version of
Lemma~\ref{filtered-by-coacyclic-is-coacyclic-lemma}, it follows
that the cone of the closed morphism
$\Cohom_X(\cC^X_\cu(\cA_X,\g,\widetilde\g),\Q^\bu)\rarrow\Q^\bu$
is Positselski-contraacyclic in $\cB^\cu\bLcth_\bW$.

 The proof of part~(b) is similar.
\end{proof}

 An antilocal contraherent graded graded $\cB^*$\+module $\Q^*$ on $X$
is said to be \emph{trivial} if it is trivial as a ($\bW$\+locally)
contraherent graded $\cB^*$\+module.
 Similarly, an  antilocal $X$\+locally cotorsion graded $\cB^*$\+module
$\Q^*$ on $X$ is said to be \emph{trivial} if it is trivial as
an ($X$\+locally cotorsion $\bW$\+locally) contraherent graded
$\cB^*$\+module.
 An antilocal ($X$\+locally contraadjusted or $X$\+locally cotorsion) 
contraherent CDG\+module $\Q^\cu$ over $\cB^\cu$ is said to be
\emph{trivial} if its underlying antilocal ($X$\+locally contraadjusted
or $X$\+locally cotorsion) graded $\cB^*$\+module $\Q^*$ is trivial.
 The full subcategory in $\cB^\cu\bCtrh_\al$ formed by the trivial
antilocal contraherent CDG\+modules is isomorphic to the DG\+category
of complexes of antilocal contraherent cosheaves on~$X$.
 The full subcategory in $\cB^\cu\bCtrh_\al^{X\dlct}$ formed by
the trivial antilocal $X$\+locally cotorsion contraherent CDG\+modules
is isomorphic to the DG\+category of complexes of antilocal
locally cotorsion contraherent cosheaves on~$X$.

 We need the following auxiliary lemma.
 
\begin{lem} \label{Cohom-from-quasi-mod-into-antilocal-is-antilocal}
 Let $X$ be a quasi-compact semi-separated scheme.
 Then \par
\textup{(a)} for any quasi-coherent quasi-module $\cB$ on $X$ that is
very flat as a quasi-coherent sheaf with respect to its left
$\cO_X$\+module structure, and for any antilocal contraherent cosheaf\/
$\P$ on $X$, the contraherent cosheaf\/ $\Cohom_X(\cB,\P)$ on $X$
is antilocal; \par
\textup{(b)} for any quasi-coherent quasi-module $\cB$ on $X$ that is
flat as a quasi-coherent sheaf with respect to its left $\cO_X$\+module
structure, and for any antilocal localy cotorsion contraherent cosheaf\/
$\P$ on $X$, the locally cotorsion contraherent cosheaf\/
$\Cohom_X(\cB,\P)$ on $X$ is antilocal.
\end{lem}

\begin{proof}
 The argument from the proof of
Lemma~\ref{Cohom-into-antilocal-cosheaf-is-antilocal} works in the more
general context of the present lemma as well.
 The only changes are that one needs to use exactness of the functor
$\Cohom_X(\cB,{-})\:X\Lcth_\bW\rarrow X\Lcth_\bW$ (in part~(a)) or
the functor $\Cohom_X(\cB,{-})\:X\Lcth_\bW^\lct\rarrow X\Lcth_\bW^\lct$
(in part~(b)), which is clear from the construction of the respective
functor $\Cohom_X$ in Section~\ref{cohom-from-quasi-modules-subsecn}.
 One also needs to use
formula~\eqref{qcoh-quasi-mod-cohom-projection-formula} from
Section~\ref{cohom-from-quasi-modules-subsecn} in lieu
of~\cite[formula~(3.21)]{Pcosh}.
\end{proof}

\begin{lem} \label{trivial-antiloc-ctrh-koszul-resolution-lemma}
 Assume that the scheme $X$ is quasi-compact and semi-separated. \par
\textup{(a)} Let\/ $\Q^\bu$ be a trivial antilocal contraherent
left CDG\+module over~$\cB^\cu$.
 Then the\/ $\Cohom$ cosheaf\/
$\Cohom_X(\cC^X_\cu(\cA_X,\g,\widetilde\g),\Q^\bu)$ is
an antilocal contraherent left CDG\+mod\-ule over $\cB^\cu$, and
the cone of the natural closed morphism of antilocal contraherent
CDG\+modules\/ $\Cohom_X(\cC^X_\cu(\cA_X,\g,\widetilde\g),\Q^\bu)
\rarrow\Q^\bu$ is Positselski-contraacyclic in the exact DG\+category
$\cB^\cu\bCtrh_\al$. \par
\textup{(b)} In the context of part~\textup{(a)}, if\/ $\Q^\bu$ is
a trivial antilocal $X$\+locally cotorsion contraherent CDG\+module
over $\cB^\cu$, then\/
$\Cohom_X(\cC^X_\cu(\cA_X,\g,\widetilde\g),\Q^\bu)$ is
an antilocal $X$\+locally cotorsion contraherent CDG\+module over
$\cB^\cu$ and the cone of the closed morphism\/
$\Cohom_X(\cC^X_\cu(\cA_X,\g,\widetilde\g),\Q^\bu)\rarrow\Q^\bu$ is
Positselski-contraacyclic in the exact DG\+category
$\cB^\cu\bCtrh_\al^{X\dlct}$.
\end{lem}

\begin{proof}
 Similar to the proof of
Lemma~\ref{trivial-lcth-koszul-resolution-lemma}.
 The result of
Lemma~\ref{Cohom-from-quasi-mod-into-antilocal-is-antilocal}
needs to be used.
\end{proof}

\subsection{Reduced coderived category of quasi-coherent CDG-modules}
\label{reduced-coderived-of-qcoh-subsecn}
 Let $X$ be a scheme and $(\g,\widetilde\g)$ be a quasi-coherent
twisted Lie algebroid over~$X$ such that the quasi-coherent sheaf~$\g$
on $X$ is finite locally free.
 Denote by $\cB^\cu=\cC_X^\cu(\g,\widetilde\g)$ the related
Chevalley--Eilenberg quasi-coherent CDG\+quasi-algebra over~$X$.
 Assume for simplicity that the rank of the quasi-coherent sheaf~$\g$
on $X$ is bounded; so $\cB^n=0$ for $n$~large enough.

 Following the discussion in
Section~\ref{koszul-coresolutions-of-trivial-qcoh-cdg-mods},
any complex of quasi-coherent sheaves $\N^\bu$ on $X$ can be viewed
as a trivial quasi-coherent CDG\+module over~$\cB^\cu$.
 The \emph{reduced coderived category of quasi-coherent CDG\+modules
over~$\cB^\cu$} is obtained from the coderived category
$\sD^\co(\cB^\cu\bQcoh)$ by annihilating the trivial quasi-coherent
CDG\+modules corresponding to \emph{acyclic} complexes of
quasi-coherent sheaves~$\N^\bu$.

 Let us say that a quasi-coherent CDG\+module $\M^\cu$ over $\cB^\cu$
is \emph{Positselski reduced-coacyclic} if $\M^\cu$ belongs to
the minimal thick subcategory of the homotopy category
$\sH^0(\cB^\cu\bQcoh)$ containing the Positselski-coacyclic
quasi-coherent CDG\+modules over $\cB^\cu$ and the trivial
quasi-coherent CDG\+modules corresponding to acyclic complexes of
quasi-coherent sheaves on~$X$.
 The reduced (Positselski) coderived category
$$
 \sD^\co_{X\red}(\cB^\cu\bQcoh)=
 \sH^0(\cB^\cu\bQcoh)/\Ac^\co_{X\red}(\cB^\cu\bQcoh)
$$
is constructed as the triangulated Verdier quotient category of
the homotopy category $\sH^0(\cB^\cu\bQcoh)$ by the thick subcategory
$\Ac^\co_{X\red}(\cB^\cu\bQcoh)\subset\sH^0(\cB^\cu\bQcoh)$ of
Positselski reduced-coacyclic CDG\+modules.
 We will see below in
Remark~\ref{reduced-coacyclic-direct-sum-closed-remark} that, for
a quasi-compact semi-separated scheme $X$, the full triangulated
subcategory $\Ac^\co_{X\red}(\cB^\cu\bQcoh)$ is closed under infinite
direct sums in $\sH^0(\cB^\cu\bQcoh)$.

 Notice that the trivial quasi-coherent CDG\+modules over $\cB^\cu$
are essentially \emph{never} thick (if $\g\ne0$).
 For this reason, the definition of the reduced coderived category of
thick CDG\+modules is different.
 We start with the following couple of lemmas.

\begin{lem} \label{canonical-filtrations-compatible-with-differential}
 Let $(R,\g,\widetilde\g)$ be a twisted Lie algebroid with a finitely
generated projective $R$\+module\/~$\g$, and let
$B^\cu=C^\cu_R(\g,\widetilde\g)$ be its Chevalley--Eilenberg CDG\+ring.
 Let $M^\cu=(M^*,d_M)$ be a thick CDG\+module over $B^\cu$ (i.~e.,
$M^\cu$ is a CDG\+module over $B^\cu$ such that the graded $B^*$\+module
$M^*$ is thick; see Section~\ref{exact-dg-categ-of-thick-cdg-mods}).
 Denote by $F$ the canonical decreasing filtration on $M^*$ as per
Lemma~\ref{vwr-projective-canonical-filtration-lemma}
or~\ref{vwr-injective-canonical-filtration-lemma}.
 Then $F$ is a filtration of the CDG\+module $M^\cu$ by
CDG\+submodules.
\end{lem}

\begin{proof}
 The claim is that, for either one of the two filtrations from the two
lemmas, one has $d_M(F^iM^*)\subset F^i(M^*)$ for all $i\in\boZ$.
 In both cases, the assertion follows from the construction of
the respective filtration $F$ on $M^*$ together with the fact that
$d(F^iB^*)\subset F^iB^*$ for all $i\in\boZ$ (where $d$~is
the Chevalley--Eilenberg differential in~$B^\cu$) and the assumption
that $d_M$~is an odd derivation on the graded $B^*$\+module $M^*$
compatible with the derivation~$d$ on~$B^*$.
\end{proof}

\begin{lem} \label{reduced-acyclic-thick-qcoh-cdg-modules}
 Let $\M^\cu$ be a thick quasi-coherent CDG\+module over
the Chevalley--Eilenberg quasi-coherent CDG\+quasi-algebra $\cB^\cu$,
and $F$ be the canonical decreasing filtration on the quasi-coherent
graded module $\M^*$ over $\cB^*$ as per
Corollary~\ref{thick-qcoh-sheaves-are-filtered}(a).
 Then $F$ is a filtration of the quasi-coherent CDG\+module $\M^\cu$
by quasi-co\-her\-ent CDG\+submodules, and the successive quotient
quasi-coherent CDG\+modules $F^i\M^\cu/F^{i+1}\M^\cu$ are trivial
quasi-coherent CDG\+modules over~$\cB^\cu$.
 Furthermore, the following two conditions are equivalent:
\begin{enumerate}
\item the complex of quasi-coherent sheaves on $X$ corresponding to
the trivial quasi-coherent CDG\+module\/
$\gr_F^0\M^\cu=\M^\cu/F^1\M^\cu$ is acyclic; \par
\item the complexes of quasi-coherent sheaves on $X$ corresponding to
the trivial quasi-coherent CDG\+modules\/ $\gr_F^i\M^\cu$ are acyclic
for all $i\in\boZ$.
\end{enumerate}
\end{lem}

\begin{proof}
 The first assertion follows from
Lemma~\ref{canonical-filtrations-compatible-with-differential},
and the second assertion is obvious.
 The prove the last assertion, put $\N^\bu=\M^\cu/F^1\M^\cu$,
and recall the isomorphisms of graded quasi-coherent sheaves
$\gr_F^i\M^*\simeq\gr_F^i\cB^*\ot_{\cO_X}\N^*=\cB^i\ot_{\cO_X}\N^*$
from Section~\ref{assoc-graded-modules-are-cta-cot}.
 For every affine open subscheme $U\subset X$,
the differential~$d_{\cB,U}$ on the graded ring $\cB^*(U)$ has
the property that $d_{\cB,U}(F^i\cB^*(U))\subset F^{i+1}\cB^*(U)$,
and it follows that $\gr_F^i\M^\cu\simeq\cB^i\ot_{\cO_X}\N^\bu$ is
an isomorphism of complexes of quasi-coherent sheaves on~$X$.
 Since $\cB^i$ is a flat (in fact, finite locally free) quasi-coherent
sheaf on $X$, it is clear that acyclicity of $\N^\bu$ implies
acyclicity of $\cB^i\ot_{\cO_X}\N^\bu$.
\end{proof}

 Let us say that a thick quasi-coherent CDG\+module $\M^\cu$ over
$\cB^\cu$ is \emph{reduced-acyclic} if it satisfies any one of
the equivalent conditions~(1) or~(2) from
Lemma~\ref{reduced-acyclic-thick-qcoh-cdg-modules}.
 
\begin{lem} \label{reduced-acyclic-thick-qcoh-localizing-subcategory}
 If two thick quasi-coherent CDG\+modules over $\cB^\cu$ are isomorphic
as objects of the Positselski coderived category\/
$\sD^\co(\cB^\cu\bQcoh_\bth)$, then one of them is reduced-acyclic if
and only if the other one is.
 In particular, all coacyclic thick quasi-coherent CDG\+modules over
$\cB^\cu$ are reduced-acyclic.
 Moreover, the reduced-acyclic thick quasi-coherent CDG\+modules form
a full triangulated subcategory closed under infinite direct sums
in\/ $\sD^\co(\cB^\cu\bQcoh_\bth)$.
\end{lem}

\begin{proof}
 Consider the DG\+functor
$$
 \gr_F^0\:\cB^\cu\bQcoh_\bth\lrarrow\bCom(X\Qcoh)
$$
assigning to every thick quasi-coherent CDG\+module $\M^\cu$
the complex of quasi-coherent sheaves $\gr_F^0\M^\cu=\M^\cu/F^1\M^\cu$
on~$X$.
 By (the graded version of) the last assertion of
Lemma~\ref{vwr-projective-canonical-filtration-lemma},
the functor $\gr_F^0$ is exact as a functor between exact
DG\+categories, i.~e., the functor $\gr_F^0\:\sZ^0(\cB^\cu\bQcoh)
\rarrow\Com(X\Qcoh)$ takes (admissible) short exact sequences to
short exact sequences.
 It is also clear that the functor $\gr_F^0$ preserves infinite
direct sums.
 Therefore, the functor $\gr_F^0$ takes Positselski-coacyclic objects
in $\cB^\cu\bQcoh_\bth$ to Positselski-coacyclic objects in
$\bCom(X\Qcoh)$, and induces a triangulated functor between
the Positselski coderived categories
$$
 \gr_F^0\:\sD^\co(\cB^\cu\bQcoh_\bth)\lrarrow
 \sD^\co(X\Qcoh).
$$
 It also follows that the latter triangulated functor preserves
infinite direct sums.
 Now the full subcategory of reduced-acyclic objects in
$\sD^\co(\cB^\cu\bQcoh_\bth)$ is the preimage of the full subcategory
of acyclic complexes in $\sD^\co(X\Qcoh)$ under the functor $\gr_F^0$.
 It remains to point out that the full subcategory of acyclic complexes
in $\sD^\co(X\Qcoh)$ is a strictly full triangulated subcategory closed
under infinite direct sums.
\end{proof}

\begin{rem} \label{proving-becker-coacyclicity-preserved-remark}
 The version of
Lemma~\ref{reduced-acyclic-thick-qcoh-localizing-subcategory} for
the Becker coderived category $\sD^\bco(\cB^\cu\bQcoh_\bth)$ would be
harder to prove, because it is generally more difficult to check
that a functor preserves the Becker-coacyclicity than
the Positselski-coacyclicity.
 The result of
Lemma~\ref{DG-functors-preserve-Becker-co-contra-acyclicity}(a) is not
applicable, as the exact DG\+functor $\gr_F^0\:\cB^\cu\bQcoh_\bth
\rarrow\bCom(X\Qcoh)$ has neither a left nor a right adjoint.
 An approach based on the description of Becker-coacyclic objects
in the spirit of~\cite[Corollary~7.17]{PS5} and the proof of
Lemma~\ref{Grothendick-abelian-DG-functor-preserves-Becker-co} above
might be possible, but it would require extending this result
from~\cite{PS5} from Grothendieck abelian DG\+categories to
a suitable class of exact DG\+categories.
 This has not been worked out yet.
 This is the reason why we do not consider the Becker coderived 
categories in this Section~\ref{reduced-coderived-of-qcoh-subsecn}.
 This discussion will be continued in
Remark~\ref{becker=positselski-at-the-end-remark} below.
\end{rem}

 We define the \emph{reduced} (\emph{Positselski}) \emph{coderived
category of thick quasi-coherent CDG\+modules over~$\cB^\cu$} as
the triangulated Verdier quotient category
$$
 \sD^\co_{X\red}(\cB^\cu\bQcoh_\bth)=
 \sD^\co(\cB^\cu\bQcoh_\bth)/\Ac_{X\red}\sD^\co(\cB^\cu\bQcoh_\bth)
$$
of the coderived category of thick quasi-coherent CDG\+modules
$\sD^\co(\cB^\cu\bQcoh_\bth)$ by the thick subcategory
$\Ac_{X\red}\sD^\co(\cB^\cu\bQcoh_\bth)\subset
\sD^\co(\cB^\cu\bQcoh_\bth)$ of reduced-acyclic thick
quasi-coherent CDG\+modules.

 It follows from
Lemma~\ref{reduced-acyclic-thick-qcoh-localizing-subcategory}
that the reduced-acyclic quasi-coherent CDG\+modules over $\cB^\cu$
also form a strictly full and thick subcategory in the absolute
derived category $\sD^\abs(\cB^\cu\bQcoh_\bth)$.
 We define the \emph{reduced absolute derived category of thick
quasi-coherent CDG\+modules over~$\cB^\cu$} as the triangulated
quotient category
$$
 \sD^\abs_{X\red}(\cB^\cu\bQcoh_\bth)=
 \sD^\abs(\cB^\cu\bQcoh_\bth)/\Ac_{X\red}\sD^\abs(\cB^\cu\bQcoh_\bth)
$$
of the absolute derived category $\sD^\abs(\cB^\cu\bQcoh_\bth)$ of
thick quasi-coherent CDG\+mod\-ules by the thick subcategory
$\Ac_{X\red}\sD^\abs(\cB^\cu\bQcoh_\bth)\subset
\sD^\abs(\cB^\cu\bQcoh_\bth)$ of re\-duced-acyclic thick
quasi-coherent CDG\+modules.

 It is clear from the two definitions above that the identity functor
is an equivalence (in fact, isomorphism) of triangulated categories
\begin{equation} \label{reduced-coderived=absolute-derived}
 \sD^\abs_{X\red}(\cB^\cu\bQcoh_\bth)\simeq
 \sD^\co_{X\red}(\cB^\cu\bQcoh_\bth).
\end{equation}
 So the reduced versions of the coderived and absolute derived
categories of thick quasi-coherent CDG\+modules agree with each other.

\begin{thm} \label{qcoh-reduced-coderived-category-equivalence-thm}
 Let $X$ be a quasi-compact semi-separated scheme and
$(\g,\widetilde\g)$ be a quasi-coherent twisted Lie algebroid over~$X$.
 Assume that\/ $\g$~is a finite locally free sheaf on $X$, and let
$\cB^\cu=\cC^\cu_X(\g,\widetilde\g)$ be the related Chevalley--Eilenberg
quasi-coherent CDG\+quasi-algebra over~$X$.
 Then the inclusion of exact/abelian DG\+categories\/
$\cB^\cu\bQcoh_\bth\rarrow\cB^\cu\bQcoh$ induces an equivalence of
the reduced coderived categories
$$
 \sD^\co_{X\red}(\cB^\cu\bQcoh_\bth)\simeq
 \sD^\co_{X\red}(\cB^\cu\bQcoh).
$$
\end{thm}

\begin{proof}
 By Corollary~\ref{thick-cdg-positselski-coderived-equiv-cor},
the inclusion of exact/abelian DG\+categories $\cB^\cu\bQcoh_\bth
\allowbreak\rarrow\cB^\cu\bQcoh$ induces an equivalence of
the coderived categories
\begin{equation} \label{unreduced-coderived-equivalence}
 \sD^\co(\cB^\cu\bQcoh_\bth)\simeq\sD^\co(\cB^\cu\bQcoh).
\end{equation}
 By the definition, the reduced coderived category
$\sD^\co_{X\red}(\cB^\cu\bQcoh_\bth)$ is the quotient category of
the coderived category $\sD^\co(\cB^\cu\bQcoh_\bth)$ by the thick
subcategory of reduced-acyclic CDG\+modules.
 It is clear from the definition that the reduced coderived category
$\sD^\co_{X\red}(\cB^\cu\bQcoh)$ is the quotient category of
the coderived category $\sD^\co(\cB^\cu\bQcoh)$ by the thick subcategory
spanned by the trivial CDG\+modules corresponding to acyclic complexes
of quasi-coherent sheaves on~$X$.
 In order to prove the theorem, it remains to check that
the triangulated equivalence~\eqref{unreduced-coderived-equivalence}
identifies the respective thick subcategories in
$\sD^\co(\cB^\cu\bQcoh_\bth)$ and $\sD^\co(\cB^\cu\bQcoh)$.

 Indeed, let $\M^\cu$ be a reduced-acyclic quasi-coherent CDG\+module
over~$\cB^\cu$.
 Then $\M^\cu$ has a finite filtration $F$ whose successive quotients
are trivial CDG\+modules corresponding to acyclic complexes of
quasi-coherent sheaves.
 Therefore, viewed as an object of the coderived category
$\sD^\co(\cB^\cu\bQcoh)$ (or even of the absolute derived category
$\sD^\abs(\cB^\cu\bQcoh)$), the quasi-coherent CDG\+module $\M^\cu$
belongs to the thick subcategory spanned by trivial CDG\+modules
corresponding to acyclic complexes of quasi-coherent sheaves.

 Conversely, let $\N^\bu$ be an acyclic complex of quasi-coherent
sheaves on $X$, viewed as a trivial quasi-coherent CDG\+module.
 We need to show that the corresponding object of the coderived
category of thick quasi-coherent CDG\+modules
$\sD^\co(\cB^\cu\bQcoh_\bth)$ is reduced-acyclic.
 At this point, it is convenient to pass from the left to the right
quasi-coherent CDG\+modules, which we can do in view of the discussion
in Section~\ref{opposite-cdg-rings-and-twisted-lie-subsecn}.
 Without loss of generality, we can also assume that $\g$~is
a finite locally free sheaf of constant rank~$d$ on $X$ (cf.\
the proof of Corollary~\ref{thick-qcoh-sheaves-are-filtered}(b)).

 So let us view $\N^\bu$ as a trivial quasi-coherent right
CDG\+module over~$\cB^\cu$.
 Then, by Lemma~\ref{trivial-qcoh-koszul-coresolution-lemma}, we
have a closed morphism of quasi-coherent right CDG\+modules
$\N^\bu\rarrow\N^\bu\ot_{\cO_X}\cC^X_\cu(\cA_X,\g,\widetilde\g)$
that is an isomorpism in the coderived category
$\sD^\co(\bQcohr\cB^\cu)$.
 Now it is clear that the quasi-coherent right CDG\+module
$\N^\bu\ot_{\cO_X}\cC^X_\cu(\cA_X,\g,\widetilde\g)$ over $\cB^\cu$
is thick.

 Let us \emph{warn} the reader that the canonical decreasing filtration
$F$ on the thick quasi-coherent CDG\+module
$\N^\bu\ot_{\cO_X}\cC^X_\cu(\cA_X,\g,\widetilde\g)$, as per the opposite
version of Corollary~\ref{thick-qcoh-sheaves-are-filtered}(a), is
very different from the increasing filtration $F$ used in the proof of
Lemma~\ref{trivial-qcoh-koszul-coresolution-lemma}.
 The essential difference is that the filtration in the proof of
Lemma~\ref{trivial-qcoh-koszul-coresolution-lemma} incorporated
the Poincar\'e--Birkhoff--Witt filtration on the enveloping
quasi-algebra $\cA_X(\g,\widetilde\g)$, while the canonical filtration
from Corollary~\ref{thick-qcoh-sheaves-are-filtered}(a) ignores
the Poincar\'e--Birkhoff--Witt filtration on $\cA_X(\g,\widetilde\g)$
completely.

 Recall that the underlying graded quasi-coherent right $\cB^*$\+module
of the quasi-coherent right CDG\+module
$\N^\bu\ot_{\cO_X}\cC^X_\cu(\cA_X,\g,\widetilde\g)$ is
$\N^*\ot_{\cO_X}\cC^X_*(\cA_X,\g,\widetilde\g)=
\N^*\ot_{\cO_X}\cA_X(\g,\widetilde\g)\ot_{\cO_X}\bigwedge_X^*(\g)$.
 The decreasing filtration $F$ from
Corollary~\ref{thick-qcoh-sheaves-are-filtered}(a) is given by
the rule $F^i\bigl(\N^*\ot_{\cO_X}\cA_X(\g,\widetilde\g)\ot_{\cO_X}\bigwedge_X^*(\g)\bigr)=\N^*\ot_{\cO_X}\cA_X(\g,\widetilde\g)
\ot_{\cO_X}\bigoplus_{j=0}^{d-i}\bigwedge_X^j(\g)$.
 The complex $\gr_F^0\bigl(\N^\bu\ot_{\cO_X}
\cC^X_\cu(\cA_X,\g,\widetilde\g)\bigr)$ is isomorphic to the complex
$\N^\bu\ot_{\cO_X}\cA_X(\g,\widetilde\g)\ot_{\cO_X}\bigwedge^d_X(\g)$.
 Since the complex of quasi-coherent sheaves $\N^\bu$ is acyclic by
assumption, so is the complex of quasi-coherent sheaves
$\gr_F^0\bigl(\N^\bu\ot_{\cO_X}
\cC^X_\cu(\cA_X,\g,\widetilde\g)\bigr)$ on~$X$.
 Thus the thick quasi-coherent right CDG\+module $\N^\bu\ot_{\cO_X}
\cC^X_\cu(\cA_X,\g,\widetilde\g)$ over $\cB^\cu$ is reduced-acyclic.
\emergencystretch=1em
\end{proof}

\begin{rem} \label{reduced-coacyclic-direct-sum-closed-remark}
 Now we can show that the full triangulated subcategory of
reduced-coacyclic CDG\+modules $\Ac^\co_{X\red}(\cB^\cu\bQcoh)$ is
closed under infinite direct sums in the homotopy category
$\sH^0(\cB^\cu\bQcoh)$, as promised in the beginning of this section.
 Indeed, the triangulated Verdier quotient functor
$\sH^0(\cB^\cu\bQcoh)\rarrow\sD^\co(\cB^\cu\bQcoh)$ preserves
infinite direct sums, and the kernel $\Ac^\co(\cB^\cu\bQcoh)$
of this functor is contained in $\Ac^\co_{X\red}(\cB^\cu\bQcoh)$.
 Therefore, it suffices to check that the triangulated Verdier quotient
category $\Ac^\co_{X\red}(\cB^\cu\bQcoh)/\Ac^\co(\cB^\cu\bQcoh)$ is
closed under infinite direct sums as a full subcategory
in $\sD^\co(\cB^\cu\bQcoh)$.
 Viewed as a full subcategory in $\sD^\co(\cB^\cu\bQcoh)$, the quotient
category $\Ac^\co_{X\red}(\cB^\cu\bQcoh)/\Ac^\co(\cB^\cu\bQcoh)$ is
just the thick subcategory in $\sD^\co(\cB^\cu\bQcoh)$ spanned by
the trivial CDG\+modules corresponding to acyclic complexes of
quasi-coherent sheaves.
 Following the proof of
Theorem~\ref{qcoh-reduced-coderived-category-equivalence-thm},
this thick subcategory in $\sD^\co(\cB^\cu\bQcoh)$ corresponds to
the full subcategory $\Ac_{X\red}\sD^\co(\cB^\cu\bQcoh_\bth)\subset
\sD^\co(\cB^\cu\bQcoh_\bth)$ under the triangulated equivalence
$\sD^\co(\cB^\cu\bQcoh_\bth)\simeq\sD^\co(\cB^\cu\bQcoh)$ from
Corollary~\ref{thick-cdg-positselski-coderived-equiv-cor}.
 Finally, the full triangulated subcategory
$\Ac_{X\red}\sD^\co(\cB^\cu\bQcoh_\bth)$ is closed under infinite
direct sums in $\sD^\co(\cB^\cu\bQcoh_\bth)$ by
Lemma~\ref{reduced-acyclic-thick-qcoh-localizing-subcategory}.
\end{rem}

\subsection{Reduced absolute derived category of thick
$X$-contraadjusted quasi-coherent CDG-modules}
\label{reduced-coderived-of-qcoh-X-cta-subsecn}
 Let $X$ be a quasi-compact semi-separated scheme and
$(\g,\widetilde\g)$ be a quasi-coherent twisted Lie algebroid over~$X$
such that the quasi-coherent sheaf~$\g$ on $X$ is finite locally free.
 Denote by $\cB^\cu=\cC_X^\cu(\g,\widetilde\g)$ the related
Chevalley--Eilenberg quasi-coherent CDG\+quasi-algebra over~$X$.

 The following lemma is an $X$\+contraadjusted version of
Lemma~\ref{reduced-acyclic-thick-qcoh-cdg-modules}.

\begin{lem} \label{reduced-acyclic-thick-X-cta-qcoh-cdg-modules}
 Let $\M^\cu$ be a thick $X$\+contraadjusted quasi-coherent CDG\+module
over $\cB^\cu$, and $F$ be the canonical decreasing filtration on
the $X$\+contraadjusted quasi-coherent graded module $\M^*$ over $\cB^*$
as per Corollaries~\ref{thick-qcoh-sheaves-are-filtered}(a)
and~\ref{thick-cta-qcoh-sheaves-are-filtered-by-cta}.
 Then the following two conditions are equivalent:
\begin{enumerate}
\item the complex of contraadjusted quasi-coherent sheaves on $X$
corresponding to the trivial quasi-coherent CDG\+module\/
$\gr_F^0\M^\cu=\M^\cu/F^1\M^\cu$ is acyclic in $X\Qcoh^\cta$; \par
\item the complexes of contraadjusted quasi-coherent sheaves on $X$
corresponding to the trivial quasi-coherent CDG\+modules\/
$\gr_F^i\M^\cu$ are acyclic in $X\Qcoh^\cta$ for all $i\in\boZ$.
\end{enumerate}
\end{lem}

\begin{proof}
 By~\cite[Corollary~4.7.3(b)]{Pcosh} or
Corollary~\ref{qcoh-X-cta-derived-equivalence} above
(see also~\cite[Corollary~A.5.4]{Pcosh}), a complex in $X\Qcoh^\cta$ is
acyclic in $X\Qcoh^\cta$ if and only if it is acyclic in $X\Qcoh$.
 So the assertion follows from
Lemma~\ref{reduced-acyclic-thick-qcoh-cdg-modules}.
 Alternatively, the argument from the proof of
Lemma~\ref{reduced-acyclic-thick-qcoh-cdg-modules} is applicable,
together with Lemma~\ref{qcoh-internal-Hom-sheaf-cta-cot}(a).
\emergencystretch=1em
\end{proof}

 Let us also state the $X$\+cotorsion version.

\begin{lem} \label{reduced-acyclic-thick-X-cot-qcoh-cdg-modules}
 Let $\M^\cu$ be a thick $X$\+cotorsion quasi-coherent CDG\+module
over $\cB^\cu$, and $F$ be the canonical decreasing filtration on
the $X$\+cotorsion quasi-coherent graded module $\M^*$ over $\cB^*$
as per Corollaries~\ref{thick-qcoh-sheaves-are-filtered}(a)
and~\ref{thick-cot-qcoh-sheaves-are-filtered-by-cot}.
 Then the following two conditions are equivalent:
\begin{enumerate}
\item the complex of cotorsion quasi-coherent sheaves on $X$
corresponding to the trivial quasi-coherent CDG\+module\/
$\gr_F^0\M^\cu=\M^\cu/F^1\M^\cu$ is acyclic in $X\Qcoh^\cot$; \par
\item the complexes of cotorsion quasi-coherent sheaves on $X$
corresponding to the trivial quasi-coherent CDG\+modules\/
$\gr_F^i\M^\cu$ are acyclic in $X\Qcoh^\cot$ for all $i\in\boZ$.
\end{enumerate}
\end{lem}

\begin{proof}
 By the cotorsion periodicity theorem for quasi-coherent sheaves
(Theorem~\ref{qcoh-cotorsion-periodicity}), a complex in $X\Qcoh^\cot$
is acyclic in $X\Qcoh^\cot$ if and only if it is acyclic in $X\Qcoh$.
 So the assertion follows from
Lemma~\ref{reduced-acyclic-thick-qcoh-cdg-modules}.
 Alternatively, the argument from the proof of
Lemma~\ref{reduced-acyclic-thick-qcoh-cdg-modules} is applicable,
together with Lemma~\ref{qcoh-internal-Hom-sheaf-cta-cot}(b).
\end{proof}

 We will say that a thick $X$\+contraadjusted quasi-coherent
CDG\+module $\M^\cu$ over $\cB^\cu$ is \emph{reduced-acyclic} if
it satisfies any one of the equivalent conditions~(1) or~(2) from
Lemma~\ref{reduced-acyclic-thick-X-cta-qcoh-cdg-modules}.
 Similarly, a thick $X$\+cotorsion quasi-coherent CDG\+module $\M^\cu$
over $\cB^\cu$ is said to be \emph{reduced-acyclic} if it satisfies
any one of the equivalent conditions~(1) or~(2) from
Lemma~\ref{reduced-acyclic-thick-X-cot-qcoh-cdg-modules}.

\begin{lem} \label{reduced-acyclic-thick-X-cta-cot-qcoh-thick-subcat}
\textup{(a)} If two thick $X$\+contraadjusted quasi-coherent
CDG\+modules over $\cB^\cu$ are isomorphic as objects of the absolute
derived category\/ $\sD^\abs(\cB^\cu\bQcoh^{X\dcta}_\bth)$, then one
of them is reduced-acyclic if and only if the other one is.
 In particular, all absolutely acyclic thick $X$\+contraadjusted
quasi-coherent CDG\+modules over $\cB^\cu$ are reduced-acyclic.
 Moreover, the reduced-acyclic thick $X$\+contraadjusted quasi-coherent
CDG\+modules form a thick subcategory in the triangulated category\/
$\sD^\abs(\cB^\cu\bQcoh^{X\dcta}_\bth)$. \par
\textup{(b)} If two thick $X$\+cotorsion quasi-coherent CDG\+modules
over $\cB^\cu$ are isomorphic as objects of the absolute derived
category\/ $\sD^\abs(\cB^\cu\bQcoh^{X\dcot}_\bth)$, then one of them is
reduced-acyclic if and only if the other one is.
 In particular, all absolutely acyclic thick $X$\+cotorsion
quasi-coherent CDG\+modules over $\cB^\cu$ are reduced-acyclic.
 Moreover, the reduced-acyclic thick $X$\+cotorsion quasi-coherent
CDG\+modules form a thick subcategory in the triangulated category\/
$\sD^\abs(\cB^\cu\bQcoh^{X\dcot}_\bth)$.
\end{lem}

\begin{proof}
 The proof is similar to that of
Lemma~\ref{reduced-acyclic-thick-qcoh-localizing-subcategory}.
 Let us spell out part~(a).
 One observes that the functor
$$
 \gr_F^0\:\cB^\cu\bQcoh^{X\dcta}_\bth\lrarrow
 \bCom(X\Qcoh^\cta)
$$
is well-defined by
Corollary~\ref{thick-cta-qcoh-sheaves-are-filtered-by-cta}
and exact as a functor between exact DG\+categories by the last
assertion of Lemma~\ref{vwr-projective-canonical-filtration-lemma}.
 Therefore, the functor $\gr_F^0$ takes absolutely acyclic objects
in $\cB^\cu\bQcoh^{X\dcta}_\bth$ to absolutely acyclic objects
in $\bCom(X\Qcoh^\cta)$, and induces a triangulated functor between
the absolute derived categories
$$
 \gr_F^0\:\sD^\abs(\cB^\cu\bQcoh^{X\dcta}_\bth)
 \lrarrow\sD^\abs(X\Qcoh^\cta).
$$
 The thick subcategory of reduced-acyclic CDG\+modules in
$\sD^\abs(\cB^\cu\bQcoh^{X\dcta}_\bth)$ is the full preimage of
the thick subcategory of acyclic complexes in
$\sD^\abs(X\Qcoh^\cta)$ under this functor.
 The proof of part~(b) is similar (the result of
Corollary~\ref{thick-cot-qcoh-sheaves-are-filtered-by-cot} is presumed).
\end{proof}

 We define the \emph{reduced absolute derived category of thick
$X$\+contraadjusted quasi-coherent CDG\+modules over~$\cB^\cu$} as
the triangulated Verdier quotient category
$$
 \sD^\abs_{X\red}(\cB^\cu\bQcoh^{X\dcta}_\bth)=
 \sD^\abs(\cB^\cu\bQcoh^{X\dcta}_\bth)/
 \Ac_{X\red}\sD^\abs(\cB^\cu\bQcoh^{X\dcta}_\bth)
$$
of the absolute derived category of thick $X$\+contraadjusted
quasi-coherent CDG\+mod\-ules $\sD^\abs(\cB^\cu\bQcoh^{X\dcta}_\bth)$
by the thick subcategory
$\Ac_{X\red}\sD^\abs(\cB^\cu\bQcoh^{X\dcta}_\bth)\subset
\sD^\abs(\cB^\cu\bQcoh^{X\dcta}_\bth)$ of reduced-acyclic thick
$X$\+contraadjusted quasi-coherent CDG\+modules.
{\hbadness=4300\par}

 Similarly, the \emph{reduced absolute derived category of thick
$X$\+cotorsion quasi-coherent CDG\+modules over~$\cB^\cu$} is defined
as the triangulated quotient category
$$
 \sD^\abs_{X\red}(\cB^\cu\bQcoh^{X\dcot}_\bth)=
 \sD^\abs(\cB^\cu\bQcoh^{X\dcot}_\bth)/
 \Ac_{X\red}\sD^\abs(\cB^\cu\bQcoh^{X\dcot}_\bth)
$$
of the absolute derived category of thick $X$\+cotorsion quasi-coherent
CDG\+mod\-ules $\sD^\abs(\cB^\cu\bQcoh^{X\dcot}_\bth)$ by the thick
subcategory $\Ac_{X\red}\sD^\abs(\cB^\cu\bQcoh^{X\dcot}_\bth)\subset
\sD^\abs(\cB^\cu\bQcoh^{X\dcot}_\bth)$ of reduced-acyclic thick
$X$\+cotorsion quasi-coherent CDG\+mod\-ules.

\begin{cor} \label{thick-cdg-qcoh-cta-reduced-abs-derived-equiv-cor}
 Let $X$ be a quasi-compact semi-separated scheme and
$(\g,\widetilde\g)$ be a quasi-coherent twisted Lie algebroid over~$X$.
 Assume that\/ $\g$~is a finite locally free sheaf on $X$, and let
$\cB^\cu=\cC^\cu_X(\g,\widetilde\g)$ be the related Chevalley--Eilenberg
quasi-coherent CDG\+quasi-algebra over~$X$.
 Then there is a commutative diagram of triangulated equivalences and
triangulated Verdier quotient functors
\begin{equation} \label{thick-cdg-qcoh-cta-reduced-abs-derived-diagram}
\begin{gathered}
 \xymatrix{
  \sD^\abs(\cB^\cu\bQcoh^{X\dcta}_\bth) \ar@{=}[rr] \ar@{->>}[d]
  && \sD^\abs(\cB^\cu\bQcoh_\bth) \ar@{->>}[d] \\
  \sD^\abs_{X\red}(\cB^\cu\bQcoh^{X\dcta}_\bth) \ar@{=}[rr]
  && \sD^{\co=\abs}_{X\red}(\cB^\cu\bQcoh_\bth)
 }
\end{gathered}
\end{equation}
where the vertical arrows with double heads show the natural
triangulated Verdier quotient functors, while the horizontal
triangulated equivalences are induced by the inclusion of exact
DG\+categories $\cB^\cu\bQcoh^{X\dcta}_\bth\rarrow\cB^\cu\bQcoh_\bth$.
\end{cor}

\begin{proof}
 The upper horizontal equivalence is provided by
Corollary~\ref{thick-cdg-qcoh-cta-abs-derived-equiv-cor}.
 One has $\sD^\co_{X\red}(\cB^\cu\bQcoh_\bth)=
\sD^\abs_{X\red}(\cB^\cu\bQcoh_\bth)$ by
formula~\eqref{reduced-coderived=absolute-derived}.

 To construct the lower horizontal equivalence and the whole diagram,
it remains to check that the upper horizontal equivalence identifies
the thick subcategory $\Ac_{X\red}\sD^\abs(\cB^\cu\bQcoh^{X\dcta}_\bth)
\subset\sD^\abs(\cB^\cu\bQcoh^{X\dcta}_\bth)$ with the thick
subcategory $\Ac_{X\red}\sD^\abs(\cB^\cu\bQcoh_\bth)\subset
\sD^\abs(\cB^\cu\bQcoh_\bth)$.
 This follows from commutativity of the diagram formed by
the functors $\gr_F^0$ from the proofs of
Lemmas~\ref{reduced-acyclic-thick-qcoh-localizing-subcategory}
and~\ref{reduced-acyclic-thick-X-cta-cot-qcoh-thick-subcat}(a),
{\hbadness=1725
$$
 \xymatrix{
  \sD^\abs(\cB^\cu\bQcoh^{X\dcta}_\bth) \ar@{=}[rr] \ar[d]_{\gr_F^0}
  && \sD^\abs(\cB^\cu\bQcoh_\bth) \ar[d]^{\gr_F^0} \\
  \sD^\abs(X\Qcoh^\cta) \ar@{=}[rr]
  && \sD^\abs(X\Qcoh)
 }
$$
together} with the fact that the functor
$\sD^\abs(X\Qcoh^\cta)\rarrow\sD^\abs(X\Qcoh)$ is a triangulated
equivalence (as shown in the lower line of the diagram).
 The latter fact is the result of~\cite[Corollary~4.7.3(b)]{Pcosh}
and a particular case of the dual version of
Proposition~\ref{second-kind-finite-resolutions}(a) above.
 Furthermore, one needs to use the fact that the triangulated
equivalence $\sD^\abs(X\Qcoh^\cta)\rarrow\sD^\abs(X\Qcoh)$ identifies
the full subcategory of acyclic complexes in $\sD^\abs(X\Qcoh^\cta)$
with the full subcategory of acyclic complexes in
$\sD^\abs(X\Qcoh)$, which follows from commutativity of the square
diagram of triangulated equivalences and triangulated Verdier
quotient functors
$$
 \xymatrix{
  \sD^\abs(X\Qcoh^\cta) \ar@{=}[rr] \ar@{->>}[d]
  && \sD^\abs(X\Qcoh) \ar@{->>}[d] \\
  \sD(X\Qcoh^\cta) \ar@{=}[rr] && \sD(X\Qcoh)
 }
$$
 The fact that the functor $\sD(X\Qcoh^\cta)\rarrow\sD(X\Qcoh)$ is
a triangulated equivalence, which is another assertion
of~\cite[Corollary~4.7.3(b)]{Pcosh} and a particular case of
Corollary~\ref{qcoh-X-cta-derived-equivalence}, also needs to be used
here.
\end{proof}

 For a version of
Corollary~\ref{thick-cdg-qcoh-cta-reduced-abs-derived-equiv-cor}
involving also $X$\+cotorsion quasi-coherent CDG\+modules, see
Corollary~\ref{thick-reduced-co-contra-cta-cot-equivs-cor} below.
 A more direct argument, which works under more restrictive
assumptions on the scheme $X$, can be found in
Corollary~\ref{thick-cdg-qcoh-cot-reduced-abs-derived-equiv-cor}.
{\hbadness=1325\par}

\begin{rem} \label{reduced-acyclicity-of-qcoh-with-other-filtration}
 All the definitions of reduced-acyclic thick quasi-coherent
CDG\+modules in this section, as well as in the previous
Section~\ref{reduced-coderived-of-qcoh-subsecn}, can be equivalently
stated in terms of the filtration $F$ from
Corollary~\ref{thick-qcoh-sheaves-are-filtered}(b) instead of the one
from Corollary~\ref{thick-qcoh-sheaves-are-filtered}(a).
 In particular, if $\M^\cu$ is a thick $X$\+contraadjusted
quasi-coherent CDG\+module over $\cB^\cu$ and $F$ is the canonical
decreasing filtration on $\M^\cu$ as per
Corollary~\ref{thick-qcoh-sheaves-are-filtered}(b), then
the quasi-coherent sheaves $\gr_F^i\M^n$ on $X$ are contraadjusted
by Corollary~\ref{thick-cta-qcoh-sheaves-are-filtered-by-cta}.
 Using the discussion in
Section~\ref{assoc-graded-modules-are-cta-cot} and similarly to
the proof of Lemma~\ref{reduced-acyclic-thick-qcoh-cdg-modules},
there are natural isomorphisms of complexes of quasi-coherent sheaves
$\gr_F^i\M^\cu\simeq\cHom_{\cO_X}(\cB^{-i},\N^\bu)$, where
$\N^\bu=\gr_F^0\M^\cu=F^0\M^\cu$.
 By Lemma~\ref{qcoh-internal-Hom-sheaf-cta-cot}(a), the complexes of
contraadjusted quasi-coherent sheaves $\gr_F^i\M^\cu$ are acyclic in
$X\Qcoh^\cta$ for all $i\in\boZ$ if and only if the complex
$\gr_F^0\M^\cu=F^0\M^\cu$ is acyclic in $X\Qcoh^\cta$.
 In view of the proof of
Corollary~\ref{thick-qcoh-sheaves-are-filtered}(b), any one of these
two equivalent conditions holds if and only if the thick
$X$\+contraadjusted quasi-coherent CDG\+module $\M^\cu$ over $\cB^\cu$
is reduced-acyclic in the sense of the definition above in this section.
 The same applies to thick $X$\+cotorsion quasi-coherent CDG\+modules
over $\cB^\cu$, etc.  \hbadness=3100
\end{rem}

\subsection{Reduced contraderived category of contraherent
CDG-modules} \label{reduced-contraderived-of-lcth-subsecn}
 Let $X$ be a scheme with an open covering $\bW$ and $(\g,\widetilde\g)$
be a quasi-coherent twisted Lie algebroid over~$X$ such that
the quasi-coherent sheaf~$\g$ on $X$ is finite locally free.
 Denote by $\cB^\cu=\cC_X^\cu(\g,\widetilde\g)$ the related
Chevalley--Eilenberg quasi-coherent CDG\+quasi-algebra over~$X$.
 Assume for simplicity that the rank of the quasi-coherent sheaf~$\g$
on $X$ is bounded; so $\cB^n=0$ for $n$~large enough.

 Following the discussion in
Section~\ref{koszul-resolutions-of-trivial-ctrh-cdg-mods},
any complex of $\bW$\+locally contraherent cosheaves $\Q^\bu$ on $X$
can be viewed as a trivial $\bW$\+locally contraherent CDG\+module
over~$\cB^\cu$.
 The \emph{reduced contraderived category of\/ $\bW$\+locally
contraherent CDG\+modules over~$\cB^\cu$} is obtained from
the contraderived category $\sD^\ctr(\cB^\cu\bLcth_\bW)$ by annihilating
the trivial $\bW$\+locally contraherent CDG\+modules corresponding to
\emph{acyclic} complexes of $\bW$\+locally contraherent
cosheaves~$\Q^\bu$.

 Let us say that a $\bW$\+locally contraherent CDG\+module $\P^\cu$
over $\cB^\cu$ is \emph{Positselski reduced-contraacyclic} if $\P^\cu$ 
belongs to the minimal thick subcategory of the homotopy category
$\sH^0(\cB^\cu\bLcth_\bW)$ containing the Positselski-contraacyclic
$\bW$\+locally contraherent CDG\+modules over $\cB^\cu$ and the trivial
$\bW$\+locally contraherent CDG\+modules corresponding to acyclic
complexes of $\bW$\+locally contraherent cosheaves on~$X$ (i.~e.,
the acyclic complexes in the exact category $X\Lcth_\bW$).
 The reduced (Positselski) contraderived category
$$
 \sD^\ctr_{X\red}(\cB^\cu\bLcth_\bW)=
 \sH^0(\cB^\cu\bLcth_\bW)/\Ac^\ctr_{X\red}(\cB^\cu\bLcth_\bW)
$$
is constructed as the triangulated Verdier quotient category of
the homotopy category $\sH^0(\cB^\cu\bLcth_\bW)$ by the thick
subcategory $\Ac^\ctr_{X\red}(\cB^\cu\bLcth_\bW)\subset
\sH^0(\cB^\cu\bLcth_\bW)$ of Positselski reduced-contraacyclic
$\bW$\+locally contraherent CDG\+modules.
 We will see below in
Remark~\ref{reduced-contraacyclic-product-closed-remark} that, for
a quasi-compact semi-separated scheme $X$, the full triangulated
subcategory $\Ac^\ctr_{X\red}(\cB^\cu\bLcth_\bW)$ is closed under
infinite products in $\sH^0(\cB^\cu\bLcth_\bW)$.

 Similarly, we will say that an $X$\+locally cotorsion $\bW$\+locally
contraherent CDG\+module $\P^\cu$ over $\cB^\cu$ is \emph{Positselski
reduced-contraacyclic} if $\P^\cu$  belongs to the minimal thick
subcategory of the homotopy category $\sH^0(\cB^\cu\bLcth_\bW^{X\dlct})$
containing the Positselski-contraacyclic $X$\+locally cotorsion
$\bW$\+locally contraherent CDG\+modules over $\cB^\cu$ and the trivial
$X$\+locally cotorsion $\bW$\+locally contraherent CDG\+modules
corresponding to acyclic complexes of locally cotorsion $\bW$\+locally
contraherent cosheaves on~$X$ (i.~e., the acyclic complexes in
the exact category $X\Lcth_\bW^\lct$).
 The reduced (Positselski) contraderived category {\hbadness=1550
$$
 \sD^\ctr_{X\red}(\cB^\cu\bLcth_\bW^{X\dlct})=
 \sH^0(\cB^\cu\bLcth_\bW^{X\dlct})/
 \Ac^\ctr_{X\red}(\cB^\cu\bLcth_\bW^{X\dlct})
$$
is} constructed as the triangulated Verdier quotient category of
the homotopy category $\sH^0(\cB^\cu\bLcth_\bW^{X\dlct})$ by the thick
subcategory $\Ac^\ctr_{X\red}(\cB^\cu\bLcth_\bW^{X\dlct})\subset
\sH^0(\cB^\cu\bLcth_\bW^{X\dlct})$ of Positselski reduced-contraacyclic
$X$\+locally cotorsion $\bW$\+locally contraherent CDG\+modules.
 For a quasi-compact semi-separated scheme $X$, the full triangulated
subcategory $\Ac^\ctr_{X\red}(\cB^\cu\bLcth_\bW^{X\dlct})$ is closed
under infinite products in $\sH^0(\cB^\cu\bLcth_\bW^{X\dlct})$, as
we will see below.

 Notice that the trivial $\bW$\+locally contraherent CDG\+modules over
$\cB^\cu$ are essentially \emph{never} thick (if $\g\ne0$).
 For this reason, the definitions of the reduced contraderived
categories of thick CDG\+modules are different.
 We start with the following couple of lemmas, very similar to
each other and dual-analogous to
Lemma~\ref{reduced-acyclic-thick-qcoh-cdg-modules}.

\begin{lem} \label{reduced-acyclic-thick-lcth-cdg-modules}
 Let\/ $\P^\cu$ be a thick\/ $\bW$\+locally contraherent CDG\+module
over $\cB^\cu$ and $F$ be the canonical decreasing filtration on
the\/ $\bW$\+locally contraherent graded module\/ $\P^*$ over $\cB^*$
as per Corollary~\ref{thick-lcta-lct-cosheaves-are-filtered}(a).
 Then there are unique structures of\/ $\bW$\+locally contraherent
CDG\+modules $F^i\P^\cu$ over $\cB^\cu$ on the\/ $\bW$\+locally
contraherent graded modules $F^i\P^*$ over $\cB^*$ such that
the inclusion maps $F^i\P^\cu\rarrow\P^\cu$ are closed morphisms
in the DG\+category $\cB^\cu\bLcth_\bW$, and in fact, admissible
monomorphisms in the exact category\/ $\sZ^0(\cB^\cu\bLcth_\bW)$.
 The successive quotient\/ $\bW$\+locally contraherent CDG\+modules
$F^i\P^\cu/F^{i+1}\P^\cu$ are trivial\/ $\bW$\+locally contraherent
CDG\+modules over~$\cB^\cu$.
 Furthermore, the following two conditions are equivalent:
\begin{enumerate}
\item the complex of\/ $\bW$\+locally contraherent cosheaves on $X$
corresponding to the trivial\/ $\bW$\+locally contraherent
CDG\+module\/ $\gr_F^0\P^\cu=F^0\P^\cu$ is acyclic in $X\Lcth_\bW$;
\item the complexes of\/ $\bW$\+locally contraherent cosheaves on $X$
corresponding to the trivial\/ $\bW$\+locally contraherent
CDG\+modules\/ $\gr_F^i\P^\cu$ are acyclic in $X\Lcth_\bW$ for
all $i\in\boZ$.
\end{enumerate}
\end{lem}

\begin{proof}
 The first assertion follows from
Lemma~\ref{canonical-filtrations-compatible-with-differential},
and the second assertion is obvious.
 The prove the last assertion, put $\Q^\bu=F^0\P^\cu$, and recall
the isomorphisms of graded $\bW$\+locally contraherent cosheaves
$\gr_F^i\P^*\simeq\Cohom_X(\gr_F^{-i}\cB^*,\Q^*)=
\Cohom_X(\cB^{-i},\Q^*)$
from Section~\ref{assoc-graded-modules-are-cta-cot}.
 For every affine open subscheme $U\subset X$ subordinate to $\bW$,
the differential~$d_{\cB,U}$ on the graded ring $\cB^*(U)$ has
the property that $d_{\cB,U}(F^i\cB^*(U))\subset F^{i+1}\cB^*(U)$,
and it follows that $\gr_F^i\P^\cu\simeq\Cohom_X(\cB^{-i},\Q^\bu)$ is
an isomorphism of complexes of $\bW$\+locally contraherent cosheaves
on~$X$.
 Since $\cB^{-i}$ is a very flat (in fact, finite locally free)
quasi-coherent sheaf on $X$, acyclicity of $\Q^\bu$ implies
acyclicity of $\Cohom_X(\cB^{-i},\Q^\bu)$ in $X\Lcth_\bW$.
\end{proof}

\begin{lem}  \label{reduced-acyclic-thick-X-lct-lcth-cdg-modules}
 Let\/ $\P^\cu$ be a thick $X$\+locally cotorsion\/ $\bW$\+locally
contraherent CDG\+module over $\cB^\cu$ and $F$ be the canonical
decreasing filtration on the $X$\+locally cotorsion\/ $\bW$\+locally
contraherent graded module\/ $\P^*$ over $\cB^*$ as per
Corollary~\ref{thick-lcta-lct-cosheaves-are-filtered}(b).
 Then there are unique structures of $X$\+locally cotorsion\/
$\bW$\+locally contraherent CDG\+modules $F^i\P^\cu$ over $\cB^\cu$
on the\/ $X$\+locally cotorsion $\bW$\+locally contraherent graded
modules $F^i\P^*$ over $\cB^*$ such that the inclusion maps
$F^i\P^\cu\rarrow\P^\cu$ are closed morphisms in the DG\+category
$\cB^\cu\bLcth_\bW^{X\dlct}$, and in fact, admissible monomorphisms
in the exact category\/ $\sZ^0(\cB^\cu\bLcth_\bW^{X\dlct})$.
 The successive quotient $X$\+locally cotorsion\/ $\bW$\+locally
contraherent CDG\+modules $F^i\P^\cu/F^{i+1}\P^\cu$ are trivial\/
$X$\+locally cotorsion $\bW$\+locally contraherent CDG\+modules
over~$\cB^\cu$.
 Furthermore, the following two conditions are equivalent:
\begin{enumerate}
\item the complex of locally cotorsion\/ $\bW$\+locally contraherent
cosheaves on $X$ corresponding to the trivial $X$\+locally cotorsion\/
$\bW$\+locally contraherent CDG\+mod\-ule\/ $\gr_F^0\P^\cu=F^0\P^\cu$
is acyclic in $X\Lcth_\bW^\lct$;
\item the complexes of locally cotorsion\/ $\bW$\+locally contraherent
cosheaves on $X$ corresponding to the trivial $X$\+locally cotorsion\/
$\bW$\+locally contraherent CDG\+mod\-ules\/ $\gr_F^i\P^\cu$ are
acyclic in $X\Lcth_\bW^\lct$ for all $i\in\boZ$.
\hbadness=1225
\end{enumerate}
\end{lem}

\begin{proof}
 The first assertion follows from
Lemma~\ref{canonical-filtrations-compatible-with-differential},
and the second assertion is obvious.
 Since a complex in $X\Lcth_\bW^\lct$ is acyclic if and only if
it is acyclic in $X\Lcth_\bW$ (by
Theorem~\ref{module-cotorsion-periodicity}
or~\cite[Lemma~3.1.2]{Pcosh}), the last assertion actually follows
from the last assertion of
Lemma~\ref{reduced-acyclic-thick-lcth-cdg-modules}.
 Alternatively, without referring to the difficult cotorsion
periodicity theorem, the last assertion is provable by the same
argument as the last assertion of
Lemma~\ref{reduced-acyclic-thick-lcth-cdg-modules}.
 It is only important here that the grading components $\cB^{-i}$ of
the quasi-coherent graded algebra $\cB$ are flat quasi-coherent sheaves
on~$X$.
\end{proof}

 Let us say that a thick $\bW$\+locally contraherent CDG\+module
$\P^\cu$ over $\cB^\cu$ is \emph{reduced-acyclic} if it satisfies any
one of the equivalent conditions~(1) or~(2) from
Lemma~\ref{reduced-acyclic-thick-lcth-cdg-modules}.
 Similarly, we will say that a thick $X$\+locally cotorsion
$\bW$\+locally contraherent CDG\+module $\P^\cu$ over $\cB^\cu$ is
\emph{reduced-acyclic} if it satisfies any one of the equivalent
conditions~(1) or~(2) from
Lemma~\ref{reduced-acyclic-thick-X-lct-lcth-cdg-modules}.

\begin{lem} \label{reduced-acyclic-thick-lcth-colocalizing-subcategory}
\textup{(a)} If two thick\/ $\bW$\+locally contraherent CDG\+modules
over $\cB^\cu$ are isomorphic as objects of the Positselski
contraderived category\/ $\sD^\ctr(\cB^\cu\bLcth_\bW^\bth)$, then one
of them is reduced-acyclic if and only if the other one is.
 In particular, all contraacyclic thick\/ $\bW$\+locally contraherent
CDG\+modules over $\cB^\cu$ are reduced-acyclic.
 Moreover, the reduced-acyclic thick\/ $\bW$\+locally contraherent
CDG\+modules form a full triangulated subcategory closed under infinite
products in\/ $\sD^\ctr(\cB^\cu\bLcth_\bW^\bth)$. \par
\textup{(b)} If two thick $X$\+locally cotorsion\/ $\bW$\+locally
contraherent CDG\+modules over $\cB^\cu$ are isomorphic as objects of
the Positselski contraderived category\/
$\sD^\ctr(\cB^\cu\bLcth_\bW^{X\dlct,\bth})$, then one of them is
reduced-acyclic if and only if the other one is.
 In particular, all contraacyclic thick $X$\+locally cotorsion\/
$\bW$\+locally contraherent CDG\+modules over $\cB^\cu$ are
reduced-acyclic.
 Moreover, the reduced-acyclic thick $X$\+locally cotorsion\/
$\bW$\+locally contraherent CDG\+modules form a full triangulated
subcategory closed under infinite products in\/
$\sD^\ctr(\cB^\cu\bLcth_\bW^{X\dlct,\bth})$.  \hbadness=2800
\end{lem}

\begin{proof}
 This is dual-analogous to
Lemma~\ref{reduced-acyclic-thick-qcoh-localizing-subcategory}.
 Let us prove part~(a).
 Consider the DG\+functor
$$
 \gr_F^0\:\cB^\cu\bLcth_\bW^\bth\lrarrow
 \bCom(X\Lcth_\bW)
$$
assigning to every thick $\bW$\+locally contraherent CDG\+module
$\P^\cu$ the complex of $\bW$\+locally contraherent cosheaves
$\gr_F^0\P^\cu=F^0\P^\cu$ on~$X$.
 By (the graded version of) the last assertion of
Lemma~\ref{vwr-injective-canonical-filtration-lemma}, the functor
$\gr_F^0$ is exact as a functor between exact DG\+categories, i.~e.,
the functor $\gr_F^0\:\sZ^0(\cB^\cu\bLcth_\bW)\rarrow\Com(X\Lcth_\bW)$
takes (admissible) short exact sequences to short exact sequences.
 It is also clear that the functor $\gr_F^0$ preserves infinite
products.
 Therefore, the functor $\gr_F^0$ takes Positselski-contraacyclic
objects in $\cB^\cu\bLcth_\bW^\bth$ to Positselski-contraacyclic
objects in $\bCom(X\Lcth_\bW)$, and induces a triangulated functor
between the Positselski contraderived categories
$$
 \gr_F^0\:\sD^\ctr(\cB^\cu\bLcth_\bW^\bth)\lrarrow
 \sD^\ctr(X\Lcth_\bW).
$$
 It also follows that the latter triangulated functor preserves
infinite products.
 Now the full subcategory of reduced-acyclic objects in
$\sD^\ctr(\cB^\cu\bLcth_\bW^\bth)$ is the preimage of the full
subcategory of acyclic complexes in $\sD^\ctr(X\Lcth_\bW)$ under
the functor~$\gr_F^0$.
 It remains to point out that the full subcategory of acyclic
complexes in $\sD^\ctr(X\Lcth_\bW)$ is a strictly full triangulated
subcategory closed under infinite products.
 The proof of part~(b) is similar.
\end{proof}

\begin{rem} \label{proving-becker-contraacyclicity-preserved-remark}
 We do not consider the Becker contraderived categories in this
Section~\ref{reduced-contraderived-of-lcth-subsecn} because we do
\emph{not know} whether the assertions of
Lemma~\ref{reduced-acyclic-thick-lcth-colocalizing-subcategory} are
true for the Becker contraderived categories.
 Proving that a functor preserves the Becker-contraacyclicity tends
to be difficult and much harder than checking the preservation of
the Positselski-contraacyclicity.
 The result of
Lemma~\ref{DG-functors-preserve-Becker-co-contra-acyclicity}(b) is not
applicable, as the exact DG\+functor
$\gr_F^0\:\cB^\cu\bLcth_\bW^\bth\rarrow\bCom(X\Lcth_\bW)$ has
neither a left nor a right adjoint.
 It is compeletely unclear what a dual-analogous version
of~\cite[Corollary~7.17]{PS5} for the Becker-contraderived categories
might consist in (cf.\
Remark~\ref{proving-becker-coacyclicity-preserved-remark}).
 This discussion will be continued in
Remark~\ref{becker=positselski-at-the-end-remark} below.
\end{rem}

 We define the \emph{reduced} (\emph{Positselski}) \emph{contraderived
category of thick $\bW$\+locally contraherent CDG\+modules
over~$\cB^\cu$} as the triangulated Verdier quotient category
$$
 \sD^\ctr_{X\red}(\cB^\cu\bLcth_\bW^\bth)=
 \sD^\ctr(\cB^\cu\bLcth_\bW^\bth)/
 \Ac_{X\red}\sD^\ctr(\cB^\cu\bLcth_\bW^\bth)
$$
of the contraderived category of thick $\bW$\+locally contraherent
CDG\+modules $\sD^\ctr(\cB^\cu\bLcth_\bW^\bth)$ by the thick subcategory
$\Ac_{X\red}\sD^\ctr(\cB^\cu\bLcth_\bW^\bth)\subset
\sD^\ctr(\cB^\cu\bLcth_\bW^\bth)$ of reduced-acyclic thick
$\bW$\+locally contraherent CDG\+modules.
{\hbadness=1600\hfuzz=2.2pt\par}

 It follows from
Lemma~\ref{reduced-acyclic-thick-lcth-colocalizing-subcategory}(a)
that the reduced-acyclic $\bW$\+locally contraherent CDG\+modules over
$\cB^\cu$ also form a strictly full and thick subcategory in
the absolute derived category $\sD^\abs(\cB^\cu\bLcth_\bW^\bth)$.
 We define the \emph{reduced absolute derived category of thick
$\bW$\+locally contraherent CDG\+modules over~$\cB^\cu$} as
the triangulated quotient category
$$
 \sD^\abs_{X\red}(\cB^\cu\bLcth_\bW^\bth)=
 \sD^\abs(\cB^\cu\bLcth_\bW^\bth)/
 \Ac_{X\red}\sD^\abs(\cB^\cu\bLcth_\bW^\bth)
$$
of the absolute derived category $\sD^\abs(\cB^\cu\bLcth_\bW^\bth)$ of
thick $\bW$\+locally contraherent CDG\+modules by the thick subcategory
$\Ac_{X\red}\sD^\abs(\cB^\cu\bLcth_\bW^\bth)\subset
\sD^\abs(\cB^\cu\bLcth_\bW^\bth)$ of reduced-acyclic thick
$\bW$\+locally contraherent CDG\+modules.

 It is clear from the two definitions above that the identity functor
is an equivalence (in fact, isomorphism) of triangulated categories
\begin{equation} \label{lcth-reduced-contraderived=absolute-derived}
 \sD^\abs_{X\red}(\cB^\cu\bLcth_\bW^\bth)\simeq
 \sD^\ctr_{X\red}(\cB^\cu\bLcth_\bW^\bth).
\end{equation}
 So the reduced versions of the contraderived and absolute derived
categories of thick $\bW$\+locally contraherent CDG\+modules agree
with each other.

 Similarly, the \emph{reduced} (\emph{Positselski}) \emph{contraderived
category of thick $X$\+locally cotorsion $\bW$\+locally contraherent
CDG\+modules over~$\cB^\cu$} is defined as the triangulated Verdier
quotient category
$$
 \sD^\ctr_{X\red}(\cB^\cu\bLcth_\bW^{X\dlct,\bth})=
 \sD^\ctr(\cB^\cu\bLcth_\bW^{X\dlct,\bth})/
 \Ac_{X\red}\sD^\ctr(\cB^\cu\bLcth_\bW^{X\dlct,\bth})
$$
of the contraderived category of thick $X$\+locally cotorsion
$\bW$\+locally contraherent CDG\+modules
$\sD^\ctr(\cB^\cu\bLcth_\bW^{X\dlct,\bth})$ by the thick subcategory
$$
 \Ac_{X\red}\sD^\ctr(\cB^\cu\bLcth_\bW^{X\dlct,\bth})\subset
\sD^\ctr(\cB^\cu\bLcth_\bW^{X\dlct,\bth})
$$
of reduced-acyclic thick $X$\+locally cotorsion $\bW$\+locally 
contraherent CDG\+modules.

 It follows from
Lemma~\ref{reduced-acyclic-thick-lcth-colocalizing-subcategory}(b)
that the reduced-acyclic $X$\+locally cotorsion $\bW$\+locally
contraherent CDG\+modules over $\cB^\cu$ also form a strictly full
and thick subcategory in the absolute derived category
$\sD^\abs(\cB^\cu\bLcth_\bW^{X\dlct,\bth})$.
 We define the \emph{reduced absolute derived category of thick
$X$\+locally cotorsion $\bW$\+locally contraherent CDG\+modules
over~$\cB^\cu$} as the triangulated quotient category
$$
 \sD^\abs_{X\red}(\cB^\cu\bLcth_\bW^{X\dlct,\bth})=
 \sD^\abs(\cB^\cu\bLcth_\bW^{X\dlct,\bth})/
 \Ac_{X\red}\sD^\abs(\cB^\cu\bLcth_\bW^{X\dlct,\bth})
$$
of the absolute derived category
$\sD^\abs(\cB^\cu\bLcth_\bW^{X\dlct,\bth})$ of thick $X$\+locally
cotorsion $\bW$\+locally contraherent CDG\+modules by the thick
subcategory
$$
 \Ac_{X\red}\sD^\abs(\cB^\cu\bLcth_\bW^{X\dlct,\bth})\subset
 \sD^\abs(\cB^\cu\bLcth_\bW^{X\dlct,\bth})
$$
of reduced-acyclic thick $X$\+locally cotorsion $\bW$\+locally
contraherent CDG\+modules.

 It is clear from the two definitions above that the identity functor
is an equivalence (in fact, isomorphism) of triangulated categories
\begin{equation} \label{lcth-lct-reduced-contraderived=absolute-derived}
 \sD^\abs_{X\red}(\cB^\cu\bLcth_\bW^{X\dlct,\bth})\simeq
 \sD^\ctr_{X\red}(\cB^\cu\bLcth_\bW^{X\dlct,\bth}).
\end{equation}
 So the reduced versions of the contraderived and absolute derived
categories of thick $X$\+locally cotorsion $\bW$\+locally contraherent
CDG\+modules agree with each other.

\begin{thm} \label{lcth-reduced-contrader-category-equivalences-thm}
 Let $X$ be a quasi-compact semi-separated scheme with an open
covering\/ $\bW$ and $(\g,\widetilde\g)$ be a quasi-coherent twisted
Lie algebroid over~$X$.
 Assume that\/ $\g$~is a finite locally free sheaf on $X$, and let
$\cB^\cu=\cC^\cu_X(\g,\widetilde\g)$ be the related Chevalley--Eilenberg
quasi-coherent CDG\+quasi-algebra over~$X$.
 Then \par
\textup{(a)} the inclusion of exact DG\+categories\/
$\cB^\cu\bLcth_\bW^\bth\rarrow\cB^\cu\bLcth_\bW$ induces an equivalence
of the reduced contraderived categories
$$
 \sD^\ctr_{X\red}(\cB^\cu\bLcth_\bW^\bth)\simeq
 \sD^\ctr_{X\red}(\cB^\cu\bLcth_\bW);
$$ \par
\textup{(b)} the inclusion of exact DG\+categories\/
$\cB^\cu\bLcth_\bW^{X\dlct,\bth}\rarrow\cB^\cu\bLcth_\bW^{X\dlct}$
induces an equivalence of the reduced contraderived categories
$$
 \sD^\ctr_{X\red}(\cB^\cu\bLcth_\bW^{X\dlct,\bth})\simeq
 \sD^\ctr_{X\red}(\cB^\cu\bLcth_\bW^{X\dlct}).
$$
\end{thm}

\begin{proof}
 This is the dual-analogous version of
Theorem~\ref{qcoh-reduced-coderived-category-equivalence-thm}.
 Let us prove part~(a).
 According to the upper line of
diagram~\eqref{thick-lcta-cdg-contrader-equivs-diagram} in
Corollary~\ref{thick-cdg-contraderived-equiv-cor}, the inclusion of
exact DG\+categories $\cB^\cu\bLcth_\bW^\bth\rarrow\cB^\cu\bLcth_\bW$
induces an equivalence of the contraderived categories
\begin{equation} \label{unreduced-contraderived-equivalence}
 \sD^\ctr(\cB^\cu\bLcth_\bW^\bth)\simeq\sD^\ctr(\cB^\cu\bLcth_\bW).
\end{equation}
 By the definition, the reduced contraderived category
$\sD^\ctr_{X\red}(\cB^\cu\bLcth_\bW^\bth)$ is the quotient category of
the contraderived category $\sD^\ctr(\cB^\cu\bLcth_\bW^\bth)$ by
the thick subcategory of reduced-acyclic CDG\+modules.
 It is clear from the definition that the reduced contraderived category
$\sD^\ctr_{X\red}(\cB^\cu\bLcth_\bW)$ is the quotient category of
the contraderived category $\sD^\ctr(\cB^\cu\bLcth_\bW)$ by the thick
subcategory spanned by the trivial CDG\+modules corresponding to acyclic
complexes of $\bW$\+locally contraherent cosheaves on~$X$.
 In order to prove the theorem, it remains to check that
the triangulated equivalence~\eqref{unreduced-contraderived-equivalence}
identifies the respective thick subcategories in
$\sD^\ctr(\cB^\cu\bLcth_\bW^\bth)$ and $\sD^\ctr(\cB^\cu\bLcth_\bW)$.

 Indeed, let $\P^\cu$ be a reduced-acyclic $\bW$\+locally contraherent
CDG\+module over~$\cB^\cu$.
 Then $\P^\cu$ has a finite filtration $F$ whose successive quotients
are trivial CDG\+modules corresponding to acyclic complexes in
$X\Lcth_\bW$.
 Therefore, viewed as an object of the contraderived category
$\sD^\ctr(\cB^\cu\bLcth_\bW)$ (or even of the absolute derived category
$\sD^\abs(\cB^\cu\bLcth_\bW)$), the $\bW$\+locally contraherent
CDG\+module $\P^\cu$ belongs to the thick subcategory spanned by
trivial CDG\+modules corresponding to acyclic complexes of
$\bW$\+locally contraherent cosheaves.

 Conversely, let $\Q^\bu$ be an acyclic complex of $\bW$\+locally
contraherent cosheaves on $X$, viewed as a trivial $\bW$\+locally
contraherent CDG\+module.
 We need to show that the corresponding object of the contraderived
category of thick $\bW$\+locally contraherent CDG\+modules
$\sD^\ctr(\cB^\cu\bLcth_\bW^\bth)$ is reduced-acyclic.
 Without loss of generality, we can also assume that $\g$~is
a finite locally free sheaf of constant rank~$d$ on $X$ (cf.\
the proof of Corollary~\ref{thick-qcoh-sheaves-are-filtered}(b)).

 By Lemma~\ref{trivial-lcth-koszul-resolution-lemma}(a), we have
a closed morphism of\/ $\bW$\+locally contraherent
CDG\+modules\/ $\Cohom_X(\cC^X_\cu(\cA_X,\g,\widetilde\g),\Q^\bu)
\rarrow\Q^\bu$ that is an isomorphism in the contraderived category
$\sD^\ctr(\cB^\cu\bLcth_\bW)$.
 Now it is clear that the $\bW$\+locally contraherent CDG\+module
$\Cohom_X(\cC^X_\cu(\cA_X,\g,\widetilde\g),\Q^\bu)$ over $\cB^\cu$
is thick.

 As in the proof of
Theorem~\ref{qcoh-reduced-coderived-category-equivalence-thm},
we \emph{warn} the reader that the canonical decreasing filtration $F$
on the thick $\bW$\+locally contraherent CDG\+module
$\Cohom_X(\cC^X_\cu(\cA_X,\g,\widetilde\g),\Q^\bu)$, as per
Corollary~\ref{thick-lcta-lct-cosheaves-are-filtered}(a),
is very different from the decreasing filtration $F$ used in
the proof of Lemma~\ref{trivial-lcth-koszul-resolution-lemma}(a).
 The essential difference is that the filtration in the proof of
Lemma~\ref{trivial-lcth-koszul-resolution-lemma} incorporated
the Poincar\'e--Birkhoff--Witt increasing filtration on the enveloping
quasi-algebra $\cA_X(\g,\widetilde\g)$, while the canonical filtration
from Corollary~\ref{thick-lcta-lct-cosheaves-are-filtered}(a) ignores
the Poincar\'e--Birkhoff--Witt filtration on $\cA_X(\g,\widetilde\g)$.

 The underlying graded $\bW$\+locally contraherent $\cB^*$\+module
of the $\bW$\+locally contraherent CDG\+module
$\Cohom_X(\cC^X_\cu(\cA_X,\g,\widetilde\g),\Q^\bu)$ is
$\Cohom_X(\cC^X_*(\cA_X,\g,\widetilde\g),\Q^*)=
\Cohom_X\bigl(\cA_X(\g,\widetilde\g)\ot_{\cO_X}\bigwedge_X^*(\g),\>
\Q^*\bigr)$.
 The decreasing filtration $F$ from
Corollary~\ref{thick-lcta-lct-cosheaves-are-filtered}(a) is given by
the rule
\begin{multline*}
 F^i\Cohom_X\Bigl(\cA_X(\g,\widetilde\g)\ot_{\cO_X}
 \bigwedge\nolimits_X^*(\g),\>\Q^*\Bigr) \\ \,=\,
 \Cohom_X\Bigl(\cA_X(\g,\widetilde\g)\ot_{\cO_X}
 \bigoplus\nolimits_{j=d+i}^d\bigwedge\nolimits_X^j(\g),\>\Q^*\Bigr).
\end{multline*}
 The complex
$\gr_F^0\Cohom_X(\cC^X_\cu(\cA_X,\g,\widetilde\g),\Q^\bu)$
is isomorphic to the complex
$$
 \Cohom_X\Bigl(\cA_X(\g,\widetilde\g)
 \ot_{\cO_X}\bigwedge\nolimits_X^d(\g),\>\Q^\bu\Bigr).
$$
 Since the complex of $\bW$\+locally contraherent cosheaves $\Q^\bu$
is acyclic by assumption, so is the complex of $\bW$\+locally
contraherent cosheaves
$\gr_F^0\Cohom_X(\cC^X_\cu(\cA_X,\g,\widetilde\g),\Q^\bu)$ on~$X$.
 Thus the thick $\bW$\+locally contraherent CDG\+module
$\Cohom_X(\cC^X_\cu(\cA_X,\g,\widetilde\g),\allowbreak\Q^\bu)$
over $\cB^\cu$ is reduced-acyclic.
 The proof of part~(b) is similar.
\end{proof}

\begin{rem} \label{reduced-contraacyclic-product-closed-remark}
 Now we can show that the full triangulated subcategory of
reduced-contraacyclic CDG\+modules
$\Ac^\ctr_{X\red}(\cB^\cu\bLcth_\bW)$ is closed under infinite products
in the homotopy category $\sH^0(\cB^\cu\bLcth_\bW)$, as promised in
the beginning of this section.
 Indeed, the triangulated Verdier quotient functor
$\sH^0(\cB^\cu\bLcth_\bW)\rarrow\sD^\ctr(\cB^\cu\bLcth_\bW)$ preserves
infinite products, and the kernel $\Ac^\ctr(\cB^\cu\bLcth_\bW)$
of this functor is contained in $\Ac^\ctr_{X\red}(\cB^\cu\bLcth_\bW)$.
 Therefore, it suffices to check that the triangulated Verdier quotient
category
$\Ac^\ctr_{X\red}(\cB^\cu\bLcth_\bW)/\Ac^\ctr(\cB^\cu\bLcth_\bW)$ is
closed under infinite direct sums as a full subcategory
in $\sD^\ctr(\cB^\cu\bLcth_\bW)$.
 Viewed as a full subcategory in $\sD^\ctr(\cB^\cu\bLcth_\bW)$,
the quotient category \emergencystretch=1em
$$
 \Ac^\ctr_{X\red}(\cB^\cu\bLcth_\bW)/\Ac^\ctr(\cB^\cu\bLcth_\bW)
$$
is just the thick subcategory in $\sD^\ctr(\cB^\cu\bLcth_\bW)$ spanned
by the trivial CDG\+modules corresponding to acyclic complexes of
$\bW$\+locally contraherent cosheaves.
 Following the proof of
Theorem~\ref{lcth-reduced-contrader-category-equivalences-thm},
this thick subcategory in $\sD^\ctr(\cB^\cu\bLcth_\bW)$ corresponds to
the full subcategory $\Ac_{X\red}\sD^\ctr(\cB^\cu\bLcth_\bW^\bth)
\subset\sD^\ctr(\cB^\cu\bLcth_\bW^\bth)$ under the triangulated
equivalence $\sD^\ctr(\cB^\cu\bLcth_\bW^\bth)\simeq
\sD^\ctr(\cB^\cu\bLcth_\bW)$ from
diagram~~\eqref{thick-lcta-cdg-contrader-equivs-diagram} in
Corollary~\ref{thick-cdg-contraderived-equiv-cor}.
 Finally, the full triangulated subcategory
$\Ac_{X\red}\sD^\ctr(\cB^\cu\bLcth_\bW^\bth)$ is closed under infinite
products in $\sD^\ctr(\cB^\cu\bLcth_\bW^\bth)$ by
Lemma~\ref{reduced-acyclic-thick-lcth-colocalizing-subcategory}.

 Similarly one shows that the full triangulated subcategory of
reduced-contraacyclic CDG\+modules
$\Ac^\ctr_{X\red}(\cB^\cu\bLcth_\bW^{X\dlct})$ is closed under infinite
products in the homotopy category $\sH^0(\cB^\cu\bLcth_\bW^{X\dlct})$.
\end{rem}

 The following corollary claiming a triangulated equivalence between
the categories appearing in parts~(a) and~(b) of
Theorem~\ref{lcth-reduced-contrader-category-equivalences-thm}
is based on the reduced Koszul duality theorem from
Section~\ref{reduced-koszul-duality-contra-side-subsecn} below.
 For a more direct proof requiring more restrictive assumptions
on the scheme $X$, see
Corollary~\ref{lcth-reduced-X-lcta-X-lct-equivalences-cor}.

\begin{cor} \label{lcth-reduced-X-lcta-X-lct-equivs-cor}
 Let $X$ be a quasi-compact semi-separated scheme with an open
covering\/ $\bW$ and $(\g,\widetilde\g)$ be a quasi-coherent twisted
Lie algebroid over~$X$.
 Assume that\/ $\g$~is a finite locally free sheaf on $X$, and let
$\cB^\cu=\cC^\cu_X(\g,\widetilde\g)$ be the related Chevalley--Eilenberg
quasi-coherent CDG\+quasi-algebra over~$X$.
 Then there is a commutative diagram of triangulated equivalences
induced by the inclusions of exact DG\+categories
\begin{equation} \label{thick-all-X-lcta-X-lct-reduced-contrader-diag}
\begin{gathered}
 \xymatrix{
  \sD^{\ctr=\abs}_{X\red}(\cB^\cu\bLcth_\bW^\bth)
  \ar@<-2pt>[rr] \ar@{-}@<-2pt>[d]
  && \sD^\ctr_{X\red}(\cB^\cu\bLcth_\bW)
  \ar@{-}@<-2pt>[ll] \ar@{-}@<-2pt>[d] \\
  \sD^{\ctr=\abs}_{X\red}(\cB^\cu\bLcth_\bW^{X\dlct,\bth})
  \ar@<-2pt>[rr] \ar@<-2pt>[u]
  && \sD^\ctr_{X\red}(\cB^\cu\bLcth_\bW^{X\dlct})
  \ar@{-}@<-2pt>[ll] \ar@<-2pt>[u]
 }
\end{gathered}
\end{equation}
\end{cor}

\begin{proof}
 First of all, we have $\sD^\ctr_{X\red}(\cB^\cu\bLcth_\bW^\bth)
=\sD^\abs_{X\red}(\cB^\cu\bLcth_\bW^\bth)$ by
formula~\eqref{lcth-reduced-contraderived=absolute-derived}
$\sD^\ctr_{X\red}(\cB^\cu\bLcth_\bW^{X\dlct,\bth})=
\sD^\abs_{X\red}(\cB^\cu\bLcth_\bW^{X\dlct,\bth})$ by
formula~\eqref{lcth-lct-reduced-contraderived=absolute-derived}.
 Now the horizontal functors
in~\eqref{thick-all-X-lcta-X-lct-reduced-contrader-diag} are
triangulated equivalences by
Theorem~\ref{lcth-reduced-contrader-category-equivalences-thm}(a\+-b).
 The rightmost vertical functor is a triangulated equivalence by
Corollary~\ref{reduced-koszul-duality-contra-side-all-equiv} below.
 The commutativity of the diagram is obvious, so it follows that
the leftmost vertical functor is a triangulated equivalence, too.
\end{proof}

 As usual in our notation system, in the case of the open covering
$\bW=\{X\}$ of a scheme $X$, we put
\begin{align*}
 \sD^\ctr_{X\red}(\cB^\cu\bCtrh) &=
 \sD^\ctr_{X\red}(\cB^\cu\bLcth_{\{X\}}), \\
  \sD^\ctr_{X\red}(\cB^\cu\bCtrh^{X\dlct}) &=
 \sD^\ctr_{X\red}(\cB^\cu\bLcth_{\{X\}}^{X\dlct}),
\end{align*}
and
\begin{align*}
 \sD^\ctr_{X\red}(\cB^\cu\bCtrh^\bth) &=
 \sD^\ctr_{X\red}(\cB^\cu\bLcth_{\{X\}}^\bth), \\
  \sD^\ctr_{X\red}(\cB^\cu\bCtrh^{X\dlct,\bth}) &=
 \sD^\ctr_{X\red}(\cB^\cu\bLcth_{\{X\}}^{X\dlct,\bth}).
\end{align*}

\subsection{Reduced contraderived category of antilocal
contraherent CDG-mod\-ules}
\label{reduced-contraderived-of-antilocal-subsecn}
 Let $X$ be a quasi-compact semi-separated scheme and
$(\g,\widetilde\g)$ be a quasi-coherent twisted Lie algebroid over~$X$
such that the quasi-coherent sheaf~$\g$ on $X$ is finite locally free.
 Denote by $\cB^\cu=\cC_X^\cu(\g,\widetilde\g)$ the related
Chevalley--Eilenberg quasi-coherent CDG\+quasi-algebra over~$X$.

 Following the discussion in
Section~\ref{koszul-resolutions-of-trivial-ctrh-cdg-mods},
any complex of antilocal contraherent cosheaves $\Q^\bu$ on $X$ can be
viewed as a trivial antilocal contraherent CDG\+module over~$\cB^\cu$.
 The \emph{reduced contraderived category of antilocal contraherent
CDG\+modules over~$\cB^\cu$} is obtained from the contraderived category
$\sD^\ctr(\cB^\cu\bCtrh_\al)$ by annihilating the trivial antilocal 
contraherent CDG\+modules corresponding to \emph{acyclic} complexes of
antilocal contraherent cosheaves~$\Q^\bu$.

 Let us say that an antilocal contraherent CDG\+module $\P^\cu$
over $\cB^\cu$ is \emph{Positselski reduced-contraacyclic} if $\P^\cu$ 
belongs to the minimal thick subcategory of the homotopy category
$\sH^0(\cB^\cu\bCtrh_\al)$ containing the Positselski-contraacyclic
antilocal contraherent CDG\+modules over $\cB^\cu$ and the trivial
antilocal contraherent CDG\+modules corresponding to acyclic complexes
of antilocal contraherent cosheaves on~$X$ (i.~e., the acyclic complexes
in the exact category $X\Ctrh_\al$).
 The reduced (Positselski) contraderived category
$$
 \sD^\ctr_{X\red}(\cB^\cu\bCtrh_\al)=
 \sH^0(\cB^\cu\bCtrh_\al)/\Ac^\ctr_{X\red}(\cB^\cu\bCtrh_\al)
$$
is constructed as the triangulated Verdier quotient category of
the homotopy category $\sH^0(\cB^\cu\bCtrh_\al)$ by the thick
subcategory $\Ac^\ctr_{X\red}(\cB^\cu\bCtrh_\al)\subset
\sH^0(\cB^\cu\bCtrh_\al)$ of Positselski reduced-contraacyclic
antilocal contraherent CDG\+modules.
 We will see below in
Remark~\ref{antiloc-reduced-contraacyclic-product-closed-remark} that
the full triangulated subcategory $\Ac^\ctr_{X\red}(\cB^\cu\bCtrh_\al)$
is closed under infinite products in $\sH^0(\cB^\cu\bCtrh_\al$.

 Similarly, we will say that an antilocal $X$\+locally cotorsion
contraherent CDG\+mod\-ule $\P^\cu$ over $\cB^\cu$ is \emph{Positselski
reduced-contraacyclic} if $\P^\cu$  belongs to the minimal thick
subcategory of the homotopy category $\sH^0(\cB^\cu\bCtrh_\al^{X\dlct})$
containing the Positselski-contraacyclic antilocal $X$\+locally
cotorsion contraherent CDG\+modules over $\cB^\cu$ and the trivial
antilocal $X$\+locally cotorsion contraherent CDG\+modules
corresponding to acyclic complexes of antilocal locally cotorsion
contraherent cosheaves on~$X$ (i.~e., the acyclic complexes in
the exact category $X\Ctrh_\al^\lct$).
 The reduced (Positselski) contraderived category
$$
 \sD^\ctr_{X\red}(\cB^\cu\bCtrh_\al^{X\dlct})=
 \sH^0(\cB^\cu\bCtrh_\al^{X\dlct})/
 \Ac^\ctr_{X\red}(\cB^\cu\bCtrh_\al^{X\dlct})
$$
is constructed as the triangulated Verdier quotient category of
the homotopy category $\sH^0(\cB^\cu\bCtrh_\al^{X\dlct})$ by the thick
subcategory $\Ac^\ctr_{X\red}(\cB^\cu\bCtrh_\al^{X\dlct})\subset
\sH^0(\cB^\cu\bCtrh_\al^{X\dlct})$ of Positselski reduced-contraacyclic
antilocal $X$\+locally cotorsion contraherent CDG\+modules.
 The full triangulated subcategory
$\Ac^\ctr_{X\red}(\cB^\cu\bCtrh_\al^{X\dlct})$ is closed
under infinite products in $\sH^0(\cB^\cu\bCtrh_\al^{X\dlct})$, as
we will see below. {\hbadness=1500\par}

 Similarly to the discussion in
Section~\ref{reduced-contraderived-of-lcth-subsecn}, the definitions
of the reduced contraderived categories of thick antilocal CDG\+modules
require a different approach.
 The following two lemmas are relevant.

\begin{lem} \label{reduced-acyclic-thick-antiloc-ctrh-cdg-modules}
 Let\/ $\P^\cu$ be a thick antilocal contraherent CDG\+module over
$\cB^\cu$ and $F$ be the canonical decreasing filtration on
the antilocal contraherent graded module\/ $\P^*$ over $\cB^*$
as per Corollaries~\ref{thick-lcta-lct-cosheaves-are-filtered}(a)
and~\ref{thick-antilocal-cosheaves-are-filtered-by-antilocal}.
 Then the following two conditions are equivalent:
\begin{enumerate}
\item the complex of antilocal contraherent cosheaves on $X$
corresponding to the trivial antilocal contraherent CDG\+module\/
$\gr_F^0\P^\cu=F^0\P^\cu$ is acyclic in $X\Ctrh_\al$;
\item the complexes of antilocal contraherent cosheaves on $X$
corresponding to the trivial antilocal contraherent
CDG\+modules\/ $\gr_F^i\P^\cu$ are acyclic in $X\Ctrh_\al$ for
all $i\in\boZ$.
\end{enumerate}
\end{lem}

\begin{proof}
 By~\cite[Corollary~4.7.4(a)]{Pcosh} or
Corollary~\ref{A-lcth-ctrh-al-derived-equivalences}(a) above
(see also~\cite[Corollary~A.5.4]{Pcosh}), a complex in $X\Ctrh_\al$ is
acyclic in $X\Ctrh_\al$ if and only if it is acyclic in $X\Ctrh$.
 So the assertion follows from
Lemma~\ref{reduced-acyclic-thick-lcth-cdg-modules}.
 Alternatively, the argument from the proof of
Lemma~\ref{reduced-acyclic-thick-lcth-cdg-modules} is applicable,
together with Lemma~\ref{Cohom-into-antilocal-cosheaf-is-antilocal}(a).
\end{proof}

\begin{lem}  \label{reduced-acyclic-thick-X-lct-antiloc-cdg-modules}
 Let\/ $\P^\cu$ be a thick antilocal $X$\+locally cotorsion
contraherent CDG\+module over $\cB^\cu$ and $F$ be the canonical
decreasing filtration on the antilocal $X$\+locally cotorsion
contraherent graded module\/ $\P^*$ over $\cB^*$ as per
Corollaries~\ref{thick-lcta-lct-cosheaves-are-filtered}(b)
and~\ref{thick-antilocal-cosheaves-are-filtered-by-antilocal}.
 Then the following two conditions are equivalent:  \hbadness=1250
\begin{enumerate}
\item the complex of antilocal locally cotorsion contraherent cosheaves
on $X$ corresponding to the trivial antilocal $X$\+locally cotorsion
contraherent CDG\+module\/ $\gr_F^0\P^\cu=F^0\P^\cu$ is acyclic in
$X\Ctrh_\al^\lct$;
\item the complexes of antilocal locally cotorsion contraherent
cosheaves on $X$ corresponding to the trivial antilocal $X$\+locally
cotorsion contraherent CDG\+modules\/ $\gr_F^i\P^\cu$ are
acyclic in $X\Ctrh_\al^\lct$ for all $i\in\boZ$.
\end{enumerate}
\end{lem}

\begin{proof}
 Similar to the proof of
Lemma~\ref{reduced-acyclic-thick-antiloc-ctrh-cdg-modules}.
 One can refer to~\cite[Corollary~4.7.5(a)]{Pcosh} or
Corollary~\ref{A-lcth-ctrh-al-derived-equivalences}(b) above,
as well as to (the proof of)
Lemma~\ref{reduced-acyclic-thick-X-lct-lcth-cdg-modules}.
\end{proof}

 Let us say that a thick antilocal contraherent CDG\+module $\P^\cu$
over $\cB^\cu$ is \emph{reduced-acyclic} if it satisfies any
one of the equivalent conditions~(1) or~(2) from
Lemma~\ref{reduced-acyclic-thick-antiloc-ctrh-cdg-modules}.
 Similarly, we will say that a thick antilocal $X$\+locally cotorsion
contraherent CDG\+module $\P^\cu$ over $\cB^\cu$ is
\emph{reduced-acyclic} if it satisfies any one of the equivalent
conditions~(1) or~(2) from
Lemma~\ref{reduced-acyclic-thick-X-lct-antiloc-cdg-modules}.

\begin{lem} \label{reduced-acycl-thick-antiloc-ctrh-colocal-subcat}
\textup{(a)} If two thick antilocal contraherent CDG\+modules over
$\cB^\cu$ are isomorphic as objects of the Positselski contraderived
category\/ $\sD^\ctr(\cB^\cu\bCtrh_\al^\bth)$, then one of them is
reduced-acyclic if and only if the other one is.
 In particular, all contraacyclic thick antilocal contraherent
CDG\+modules over $\cB^\cu$ are reduced-acyclic.
 Moreover, the reduced-acyclic thick antilocal contraherent
CDG\+modules form a full triangulated subcategory closed under infinite
products in\/ $\sD^\ctr(\cB^\cu\bCtrh_\al^\bth)$. \par
\textup{(b)} If two thick antilocal $X$\+locally cotorsion
contraherent CDG\+modules over $\cB^\cu$ are isomorphic as objects of
the Positselski contraderived category\/
$\sD^\ctr(\cB^\cu\bCtrh_\al^{X\dlct,\bth})$, then one of them is
reduced-acyclic if and only if the other one is.
 In particular, all contraacyclic thick a antilocal $X$\+locally
cotorsion locally contraherent CDG\+modules over $\cB^\cu$ are
reduced-acyclic.
 Moreover, the reduced-acyclic thick antilocal $X$\+locally cotorsion
contraherent CDG\+modules form a full triangulated subcategory closed
under infinite products in\/
$\sD^\ctr(\cB^\cu\bCtrh_\al^{X\dlct,\bth})$.
\end{lem}

\begin{proof}
 Similar to
Lemma~\ref{reduced-acyclic-thick-lcth-colocalizing-subcategory}.
\end{proof}

 We define the \emph{reduced} (\emph{Positselski}) \emph{contraderived
category of thick antilocal contraherent CDG\+modules over~$\cB^\cu$}
as the triangulated Verdier quotient category
$$
 \sD^\ctr_{X\red}(\cB^\cu\bCtrh_\al^\bth)=
 \sD^\ctr(\cB^\cu\bCtrh_\al^\bth)/
 \Ac_{X\red}\sD^\ctr(\cB^\cu\bCtrh_\al^\bth)
$$
of the contraderived category of thick antilocal contraherent
CDG\+modules $\sD^\ctr(\cB^\cu\bCtrh_\al^\bth)$ by the thick subcategory
$\Ac_{X\red}\sD^\ctr(\cB^\cu\bCtrh_\al^\bth)\subset
\sD^\ctr(\cB^\cu\bCtrh_\al^\bth)$ of reduced-acyclic thick antilocal 
contraherent CDG\+modules. {\hbadness=2650\hfuzz=2.7pt\par}

 It follows from
Lemma~\ref{reduced-acycl-thick-antiloc-ctrh-colocal-subcat}(a)
that the reduced-acyclic antilocal contraherent CDG\+modules over
$\cB^\cu$ also form a strictly full and thick subcategory in
the absolute derived category $\sD^\abs(\cB^\cu\bCtrh_\al^\bth)$.
 We define the \emph{reduced absolute derived category of thick
antilocal contraherent CDG\+modules over~$\cB^\cu$} as the triangulated
quotient category
$$
 \sD^\abs_{X\red}(\cB^\cu\bCtrh_\al^\bth)=
 \sD^\abs(\cB^\cu\bCtrh_\al^\bth)/
 \Ac_{X\red}\sD^\abs(\cB^\cu\bCtrh_\al^\bth)
$$
of the absolute derived category $\sD^\abs(\cB^\cu\bCtrh_\al^\bth)$ of
thick antilocal contraherent CDG\+modules by the thick subcategory
$\Ac_{X\red}\sD^\abs(\cB^\cu\bCtrh_\al^\bth)\subset
\sD^\abs(\cB^\cu\bCtrh_\al^\bth)$ of reduced-acyclic thick
antilocal contraherent CDG\+modules.

 It is clear from the two definitions above that the identity functor
is an equivalence (in fact, isomorphism) of triangulated categories
\begin{equation} \label{al-reduced-contraderived=absolute-derived}
 \sD^\abs_{X\red}(\cB^\cu\bCtrh_\al^\bth)\simeq
 \sD^\ctr_{X\red}(\cB^\cu\bCtrh_\al^\bth).
\end{equation}
 So the reduced versions of the contraderived and absolute derived
categories of thick antilocal contraherent CDG\+modules agree with
each other.

 Similarly, the \emph{reduced} (\emph{Positselski}) \emph{contraderived
category of thick antilocal $X$\+lo\-cally cotorsion contraherent
CDG\+modules over~$\cB^\cu$} is defined as the triangulated Verdier
quotient category
$$
 \sD^\ctr_{X\red}(\cB^\cu\bCtrh_\al^{X\dlct,\bth})=
 \sD^\ctr(\cB^\cu\bCtrh_\al^{X\dlct,\bth})/
 \Ac_{X\red}\sD^\ctr(\cB^\cu\bCtrh_\al^{X\dlct,\bth})
$$
of the contraderived category of thick antilocal $X$\+locally cotorsion
contraherent CDG\+modules $\sD^\ctr(\cB^\cu\bCtrh_\al^{X\dlct,\bth})$
by the thick subcategory
$$
 \Ac_{X\red}\sD^\ctr(\cB^\cu\bCtrh_\al^{X\dlct,\bth})\subset
 \sD^\ctr(\cB^\cu\bCtrh_\al^{X\dlct,\bth})
$$
of reduced-acyclic thick antilocal $X$\+locally cotorsion
contraherent CDG\+modules.

 It follows from
Lemma~\ref{reduced-acycl-thick-antiloc-ctrh-colocal-subcat}(b)
that the reduced-acyclic antilocal $X$\+locally cotorsion contraherent
CDG\+modules over $\cB^\cu$ also form a strictly full and thick
subcategory in the absolute derived category
$\sD^\abs(\cB^\cu\bCtrh_\al^{X\dlct,\bth})$.
 We define the \emph{reduced absolute derived category of thick
antilocal $X$\+locally cotorsion contraherent CDG\+modules
over~$\cB^\cu$} as the triangulated quotient category
$$
 \sD^\abs_{X\red}(\cB^\cu\bCtrh_\al^{X\dlct,\bth})=
 \sD^\abs(\cB^\cu\bCtrh_\al^{X\dlct,\bth})/
 \Ac_{X\red}\sD^\abs(\cB^\cu\bCtrh_\al^{X\dlct,\bth})
$$
of the absolute derived category
$\sD^\abs(\cB^\cu\bCtrh_\al^{X\dlct,\bth})$ of thick antilocal
$X$\+locally cotorsion contraherent CDG\+modules by the thick
subcategory
$$
 \Ac_{X\red}\sD^\abs(\cB^\cu\bCtrh_\al^{X\dlct,\bth})\subset
 \sD^\abs(\cB^\cu\bCtrh_\al^{X\dlct,\bth})
$$
of reduced-acyclic thick antilocal $X$\+locally cotorsion contraherent
CDG\+modules.

 It is clear from the two definitions above that the identity functor
is an equivalence (in fact, isomorphism) of triangulated categories
\begin{equation} \label{al-lct-reduced-contraderived=absolute-derived}
 \sD^\abs_{X\red}(\cB^\cu\bCtrh_\al^{X\dlct,\bth})\simeq
 \sD^\ctr_{X\red}(\cB^\cu\bCtrh_\al^{X\dlct,\bth}).
\end{equation}
 So the reduced versions of the contraderived and absolute derived
categories of thick antilocal $X$\+locally cotorsion contraherent
CDG\+modules agree with each other.

\begin{thm} \label{antiloc-reduced-contrader-categ-equivs-thm}
 Let $X$ be a quasi-compact semi-separated scheme and
$(\g,\widetilde\g)$ be a quasi-coherent twisted Lie algebroid over~$X$.
 Assume that\/ $\g$~is a finite locally free sheaf on $X$, and let
$\cB^\cu=\cC^\cu_X(\g,\widetilde\g)$ be the related Chevalley--Eilenberg
quasi-coherent CDG\+quasi-algebra over~$X$.
 Then \par
\textup{(a)} the inclusion of exact DG\+categories\/
$\cB^\cu\bCtrh_\al^\bth\rarrow\cB^\cu\bCtrh_\al$ induces an equivalence
of the reduced contraderived categories
$$
 \sD^\ctr_{X\red}(\cB^\cu\bCtrh_\al^\bth)\simeq
 \sD^\ctr_{X\red}(\cB^\cu\bCtrh_\al);
$$ \par
\textup{(b)} the inclusion of exact DG\+categories\/
$\cB^\cu\bCtrh_\al^{X\dlct,\bth}\rarrow\cB^\cu\bCtrh_\al^{X\dlct}$
induces an equivalence of the reduced contraderived categories
$$
 \sD^\ctr_{X\red}(\cB^\cu\bCtrh_\al^{X\dlct,\bth})\simeq
 \sD^\ctr_{X\red}(\cB^\cu\bCtrh_\al^{X\dlct}).
$$
\end{thm}

\begin{proof}
 Similar to the proof of
Theorem~\ref{lcth-reduced-contrader-category-equivalences-thm},
and based on the triangulated equivalences in the lower lines of
the two diagrams in Corollary~\ref{thick-cdg-contraderived-equiv-cor}
together with Lemma~\ref{trivial-antiloc-ctrh-koszul-resolution-lemma}.
 One also needs to use
Lemma~\ref{Cohom-from-quasi-mod-into-antilocal-is-antilocal}.
\end{proof}

\begin{rem} \label{antiloc-reduced-contraacyclic-product-closed-remark}
 Similarly to Remark~\ref{reduced-contraacyclic-product-closed-remark},
one shows that the full triangulated subcategories of
reduced-contraacyclic antilocal CDG\+modules
$\Ac^\ctr_{X\red}(\cB^\cu\bCtrh_\al)$ and
$\Ac^\ctr_{X\red}(\cB^\cu\bCtrh_\al^{X\dlct})$ are closed under
infinite products in their respective homotopy categories
$\sH^0(\cB^\cu\bCtrh_\al)$ and $\sH^0(\cB^\cu\bCtrh_\al^{X\dlct})$, as
promised in the beginning of this section.
 The argument is based on
Lemma~\ref{reduced-acycl-thick-antiloc-ctrh-colocal-subcat} and
the proof of Theorem~\ref{antiloc-reduced-contrader-categ-equivs-thm}.
\end{rem}

 The following corollary is a reduced version of
Corollary~\ref{thick-cdg-contraderived-equiv-cor}.

\begin{cor} \label{reduced-cdg-contraderived-equiv-cor}
 Let $X$ be a quasi-compact semi-separated scheme with an open
covering\/ $\bW$ and $(\g,\widetilde\g)$ be a quasi-coherent twisted
Lie algebroid over~$X$.
 Assume that\/ $\g$~is a finite locally free sheaf on $X$, and let
$\cB^\cu=\cC^\cu_X(\g,\widetilde\g)$ be the related Chevalley--Eilenberg
quasi-coherent CDG\+quasi-algebra over~$X$.
 Then there are commutative diagrams of triangulated equivalences of
the reduced contraderived categories induced by the identity inclusions
of exact DG\+categories
\begin{equation} \label{reduced-lcta-cdg-contrader-equivs-diagram}
\begin{gathered}
 \xymatrix{
  \sD^{\ctr=\abs}_{X\red}(\cB^\cu\bLcth_\bW^\bth)
  \ar@<-2pt>[r] \ar@{-}@<-2pt>[d]
  & \sD^\ctr_{X\red}(\cB^\cu\bLcth_\bW) \ar@{-}@<-2pt>[l]
  \ar@{-}@<-2pt>[d] \\
  \sD^{\ctr=\abs}_{X\red}(\cB^\cu\bCtrh^\bth) \ar@<-2pt>[r]
  \ar@<-2pt>[u] \ar@{-}@<-2pt>[d]
  & \sD^\ctr_{X\red}(\cB^\cu\bCtrh) \ar@{-}@<-2pt>[l]
  \ar@<-2pt>[u] \ar@{-}@<-2pt>[d] \\
  \sD^{\ctr=\abs}_{X\red}(\cB^\cu\bCtrh_\al^\bth) \ar@<-2pt>[r]
  \ar@<-2pt>[u]
  & \sD^\ctr_{X\red}(\cB^\cu\bCtrh_\al) \ar@{-}@<-2pt>[l]
  \ar@<-2pt>[u]
 }
\end{gathered}
\end{equation}
\begin{equation} \label{reduced-X-lct-cdg-contrader-equivs-diagram}
\begin{gathered}
 \xymatrix{
  \sD^{\ctr=\abs}_{X\red}(\cB^\cu\bLcth_\bW^{X\dlct,\bth})
  \ar@<-2pt>[r] \ar@{-}@<-2pt>[d]
  & \sD^\ctr_{X\red}(\cB^\cu\bLcth_\bW^{X\dlct}) \ar@{-}@<-2pt>[l]
  \ar@{-}@<-2pt>[d] \\
  \sD^{\ctr=\abs}_{X\red}(\cB^\cu\bCtrh^{X\dlct,\bth}) \ar@<-2pt>[r]
  \ar@<-2pt>[u] \ar@{-}@<-2pt>[d]
  & \sD^\ctr_{X\red}(\cB^\cu\bCtrh^{X\dlct}) \ar@{-}@<-2pt>[l]
  \ar@<-2pt>[u] \ar@{-}@<-2pt>[d] \\
  \sD^{\ctr=\abs}_{X\red}(\cB^\cu\bCtrh_\al^{X\dlct,\bth})
  \ar@<-2pt>[r] \ar@<-2pt>[u]
  & \sD^\ctr_{X\red}(\cB^\cu\bCtrh_\al^{X\dlct}) \ar@{-}@<-2pt>[l]
  \ar@<-2pt>[u]
 }
\end{gathered}
\end{equation}
\end{cor}

\begin{proof}
 The isomorphisms between the reduced contraderived and reduced absolute
derived category in every node of the leftmost columns of the diagrams
hold by formulas~(\ref{lcth-reduced-contraderived=absolute-derived}\+-%
\ref{lcth-lct-reduced-contraderived=absolute-derived})
and~(\ref{al-reduced-contraderived=absolute-derived}\+-%
\ref{al-lct-reduced-contraderived=absolute-derived}).
 The upper and middle horizontal equivalences on
diagram~\eqref{reduced-lcta-cdg-contrader-equivs-diagram}
are provided by
Theorem~\ref{lcth-reduced-contrader-category-equivalences-thm}(a).
 The upper and middle horizontal equivalences on
diagram~\eqref{reduced-X-lct-cdg-contrader-equivs-diagram}
are provided by
Theorem~\ref{lcth-reduced-contrader-category-equivalences-thm}(b).
 The lower horizontal equivalence on
diagram~\eqref{reduced-lcta-cdg-contrader-equivs-diagram}
is provided by
Theorem~\ref{antiloc-reduced-contrader-categ-equivs-thm}(a).
 The lower horizontal equivalence on
diagram~\eqref{reduced-X-lct-cdg-contrader-equivs-diagram}
is provided by
Theorem~\ref{antiloc-reduced-contrader-categ-equivs-thm}(b).

 Let us prove the vertical equivalences on
diagram~\eqref{reduced-lcta-cdg-contrader-equivs-diagram}.
 To show that the rightmost vertical functor
$$
 \sD^\ctr_{X\red}(\cB^\cu\bCtrh_\al)\lrarrow
 \sD^\ctr_{X\red}(\cB^\cu\bLcth_\bW)
$$
is a triangulated equivalence, one has to make the following
observations.
 First of all, the functor
$$
 \sD^\ctr(\cB^\cu\bCtrh_\al)\lrarrow\sD^\ctr(\cB^\cu\bLcth_\bW)
$$
is a triangulated equivalence according to the rightmost column
of diagram~\eqref{thick-lcta-cdg-contrader-equivs-diagram} in
Corollary~\ref{thick-cdg-contraderived-equiv-cor}.

 Secondly, the reduced contraderived category
$\sD^\ctr_{X\red}(\cB^\cu\bLcth_\bW)$ is the quotient category of
the contraderived category $\sD^\ctr(\cB^\cu\bLcth_\bW)$ by the thick
subcategory spanned by acyclic complexes in $X\Lcth_\bW$, while
the reduced contraderived category $\sD^\ctr_{X\red}(\cB^\cu\bCtrh_\al)$
is the quotient category of the contraderived category
$\sD^\ctr(\cB^\cu\bCtrh_\al)$ by the thick subcategory spanned by
acyclic complexes in $X\Ctrh_\al$ viewed as trivial CDG\+modules.
 Now the construction of the trivial CDG\+module arizing from
a complex provides a commutative diagram of triangulated functors
and triangulated equivalences
$$
 \xymatrix{
  \sD^\ctr(X\Lcth_\bW) \ar[r] \ar@{-}@<-2pt>[d]
  & \sD^\ctr(\cB^\cu\bLcth_\bW) \ar@{-}@<-2pt>[d] \\
  \sD^\ctr(X\Ctrh_\al) \ar[r] \ar@<-2pt>[u]
  & \sD^\ctr(\cB^\cu\bCtrh_\al) \ar@<-2pt>[u]
 }
$$
where the leftmost vertical equivalence
$\sD^\ctr(X\Ctrh_\al)\simeq\sD^\ctr(X\Lcth_\bW)$ is the result
of~\cite[Corollary~4.7.4(a)]{Pcosh}.

 Finally, we have a commutative diagram of triangulated Verdier
quotient functors and triangulated equivalences
$$
 \xymatrix{
  \sD^\ctr(X\Lcth_\bW) \ar@{->>}[r] \ar@{-}@<-2pt>[d]
  & \sD(X\Lcth_\bW) \ar@{-}@<-2pt>[d] \\
  \sD^\ctr(X\Ctrh_\al) \ar@{->>}[r] \ar@<-2pt>[u]
  & \sD(X\Ctrh_\al) \ar@<-2pt>[u]
 }
$$
(where the rightmost vertical equivalence is also provided
by~\cite[Corollary~4.7.4(a)]{Pcosh}), implying that the thick
subcategories of acyclic complexes in $\sD^\ctr(X\Ctrh_\al)$ and
$\sD^\ctr(X\Lcth_\bW)$ are identified by the triangulated
equivalence $\sD^\ctr(X\Ctrh_\al)\simeq\sD^\ctr(X\Lcth_\bW)$.

 To show that the leftmost vertical functor
$$
 \sD^\ctr_{X\red}(\cB^\cu\bCtrh_\al^\bth)\lrarrow
 \sD^\ctr_{X\red}(\cB^\cu\bLcth_\bW^\bth)
$$
is a triangulated equivalence, we make the following observations.
 First of all, the functor
$$
 \sD^\ctr(\cB^\cu\bCtrh_\al^\bth)\lrarrow
 \sD^\ctr(\cB^\cu\bLcth_\bW^\bth)
$$
is a triangulated equivalence according to the leftmost column
of diagram~\eqref{thick-lcta-cdg-contrader-equivs-diagram} in
Corollary~\ref{thick-cdg-contraderived-equiv-cor}.

 Secondly, the reduced contraderived categories of thick CDG\+modules
$\sD^\ctr_{X\red}(\cB^\cu\allowbreak\bLcth_\bW^\bth)$ and
$\sD^\ctr_{X\red}(\cB^\cu\bCtrh_\al^\bth)$ are the quotient categories
of the contraderived categories $\sD^\ctr(\cB^\cu\bLcth_\bW^\bth)$ and
$\sD^\ctr(\cB^\cu\bCtrh_\al^\bth)$ by their respective thick
subcategories of reduced-acyclic thick CDG\+modules,
\begin{align*}
 \sD^\ctr_{X\red}(\cB^\cu\bLcth_\bW^\bth) &=
 \sD^\ctr(\cB^\cu\bLcth_\bW^\bth)/
 \Ac_{X\red}\sD^\ctr(\cB^\cu\bLcth_\bW^\bth), \\
  \sD^\ctr_{X\red}(\cB^\cu\bCtrh_\al^\bth) &=
 \sD^\ctr(\cB^\cu\bCtrh_\al^\bth)/
 \Ac_{X\red}\sD^\ctr(\cB^\cu\bCtrh_\al^\bth).
\end{align*}
 The latter thick subcategories were described in
Lemma~\ref{reduced-acyclic-thick-lcth-colocalizing-subcategory}(a)
and~\ref{reduced-acycl-thick-antiloc-ctrh-colocal-subcat}(a).
 Following these descriptions, the following commutative diagram
of triangulated functors and triangulated equivalences is relevant:
$$
 \xymatrix{
  \sD^\ctr(\cB^\cu\bLcth_\bW^\bth) \ar[rr]^-{\gr_F^0} \ar@{-}@<-2pt>[d]
  && \sD^\ctr(X\Lcth_\bW) \ar@{-}@<-2pt>[d] \\
  \sD^\ctr(\cB^\cu\bCtrh_\al^\bth) \ar[rr]^-{\gr_F^0} \ar@<-2pt>[u]
  && \sD^\ctr(X\Ctrh_\al) \ar@<-2pt>[u]
 }
$$

 To finish the argument, one needs to use the observation that
the thick subcategories of acyclic complexes in the contraderived
categories $\sD^\ctr(X\Ctrh_\al)$ and $\sD^\ctr(X\Lcth_\bW)$ are
identified by the triangulated equivalence
$\sD^\ctr(X\Ctrh_\al)\simeq\sD^\ctr(X\Lcth_\bW)$.
 That was explained already above in this proof.

 The proofs of the vertical equivalences on
diagram~\eqref{reduced-X-lct-cdg-contrader-equivs-diagram} are similar.
 One needs to use the vertical eqiuvalences on
diagram~\eqref{thick-X-lct-cdg-contrader-equivs-diagram} in
Corollary~\ref{thick-cdg-contraderived-equiv-cor}, the constructions of
the triangulated functors $\gr_F^0$ from the proofs of
Lemmas~\ref{reduced-acyclic-thick-lcth-colocalizing-subcategory}(b)
and~\ref{reduced-acycl-thick-antiloc-ctrh-colocal-subcat}(b),
and the results of~\cite[Corollary~4.7.5(a)]{Pcosh}.
\end{proof}

 Now we can prove a triangulated equivalence between the triangulated
categories appearing in parts~(a) and~(b) of
Theorem~\ref{antiloc-reduced-contrader-categ-equivs-thm}.
 In fact, there is a three-dimensional commutative diagram of
triangulated equivalences uniting the two
diagrams~(\ref{reduced-lcta-cdg-contrader-equivs-diagram}\+-%
\ref{reduced-X-lct-cdg-contrader-equivs-diagram}).
 The following corollary describes this three-dimensional diagram in
terms of its two-dimensional slices.
 The proof is essentially based on the reduced Koszul duality theorem
from Section~\ref{reduced-koszul-duality-contra-side-subsecn} below.
 For a more direct argument applicable under more restrictive
assumptions on the scheme $X$, see
Corollary~\ref{antilocal-reduced-X-lcta-X-lct-equivalences-cor}.

\begin{cor} \label{reduced-thick-lcth-al-ctrder-lcta-lct-equiv-cor}
 Let $X$ be a quasi-compact semi-separated scheme with an open
covering\/ $\bW$ and $(\g,\widetilde\g)$ be a quasi-coherent twisted
Lie algebroid over~$X$.
 Assume that\/ $\g$~is a finite locally free sheaf on $X$, and let
$\cB^\cu=\cC^\cu_X(\g,\widetilde\g)$ be the related Chevalley--Eilenberg
quasi-coherent CDG\+quasi-algebra over~$X$.
 Then \par
\textup{(a)} there is a commutative diagram of triangulated equivalences
induced by the inclusions of exact DG\+categories
\begin{equation} \label{reduced-lcth-al-X-lcta-X-lct-cdg-ctrder-diag}
\begin{gathered}
 \xymatrix{
  \sD^\ctr_{X\red}(\cB^\cu\bCtrh_\al)
  \ar@<-2pt>[r] \ar@{-}@<-2pt>[d]
  & \sD^\ctr_{X\red}(\cB^\cu\bCtrh)
  \ar@<-2pt>[r] \ar@{-}@<-2pt>[d] \ar@{-}@<-2pt>[l]
  & \sD^\ctr_{X\red}(\cB^\cu\bLcth_\bW)
  \ar@{-}@<-2pt>[l] \ar@{-}@<-2pt>[d] \\
  \sD^\ctr_{X\red}(\cB^\cu\bCtrh^{X\dlct}_\al)
  \ar@<-2pt>[r] \ar@<-2pt>[u]
  & \sD^\ctr_{X\red}(\cB^\cu\bCtrh^{X\dlct})
  \ar@<-2pt>[r] \ar@<-2pt>[u] \ar@{-}@<-2pt>[l]
  & \sD^\ctr_{X\red}(\cB^\cu\bLcth_\bW^{X\dlct})
  \ar@{-}@<-2pt>[l] \ar@<-2pt>[u]
 }
\end{gathered}
\end{equation} \par
\textup{(b)} there is a commutative diagram of triangulated equivalences
induced by the inclusions of exact DG\+categories
\begin{equation} \label{reduced-thick-lcth-al-lcta-lct-cdg-ctrder-diag}
\begin{gathered}
 \xymatrix{
  \sD^\ctr_{X\red}(\cB^\cu\bCtrh_\al^\bth)
  \ar@<-2pt>[r] \ar@{-}@<-2pt>[d]
  & \sD^\ctr_{X\red}(\cB^\cu\bCtrh^\bth)
  \ar@<-2pt>[r] \ar@{-}@<-2pt>[d] \ar@{-}@<-2pt>[l]
  & \sD^\ctr_{X\red}(\cB^\cu\bLcth_\bW^\bth)
  \ar@{-}@<-2pt>[l] \ar@{-}@<-2pt>[d] \\
  \sD^\ctr_{X\red}(\cB^\cu\bCtrh^{X\dlct,\bth}_\al)
  \ar@<-2pt>[r] \ar@<-2pt>[u]
  & \sD^\ctr_{X\red}(\cB^\cu\bCtrh^{X\dlct,\bth})
  \ar@<-2pt>[r] \ar@<-2pt>[u] \ar@{-}@<-2pt>[l]
  & \sD^\ctr_{X\red}(\cB^\cu\bLcth_\bW^{X\dlct,\bth})
  \ar@{-}@<-2pt>[l] \ar@<-2pt>[u]
 }
\end{gathered}
\end{equation} \par
\textup{(c)} There is a commutative diagram of triangulated equivalences
induced by the inclusions of exact DG\+categories
\begin{equation} \label{antil-thick-all-lcta-lct-reduced-ctrder-diag}
\begin{gathered}
 \xymatrix{
  \sD^{\ctr=\abs}_{X\red}(\cB^\cu\bCtrh_\al^\bth)
  \ar@<-2pt>[rr] \ar@{-}@<-2pt>[d]
  && \sD^\ctr_{X\red}(\cB^\cu\bCtrh_\al)
  \ar@{-}@<-2pt>[ll] \ar@{-}@<-2pt>[d] \\
  \sD^{\ctr=\abs}_{X\red}(\cB^\cu\bCtrh_\al^{X\dlct,\bth})
  \ar@<-2pt>[rr] \ar@<-2pt>[u]
  && \sD^\ctr_{X\red}(\cB^\cu\bCtrh_\al^{X\dlct})
  \ar@{-}@<-2pt>[ll] \ar@<-2pt>[u]
 }
\end{gathered}
\end{equation}
\end{cor}

\begin{proof}
 Part~(a): the horizontal triangulated equivalences
in~\eqref{reduced-lcth-al-X-lcta-X-lct-cdg-ctrder-diag} are provided
by the rightmost columns of the two
diagrams~(\ref{reduced-lcta-cdg-contrader-equivs-diagram}\+-%
\ref{reduced-X-lct-cdg-contrader-equivs-diagram}) from
Corollary~\ref{reduced-cdg-contraderived-equiv-cor}.
 The triangulated equivalences in the middle and rightmost columns
are the result of
Corollary~\ref{reduced-koszul-duality-contra-side-all-equiv} below.
 The commutativity of the diagram is obvious, so it follows that
the vertical functor in the leftmost column is a triangulated
equivalence, too.
 
 Part~(b): the horizontal triangulated equivalences
in~\eqref{reduced-thick-lcth-al-lcta-lct-cdg-ctrder-diag} are provided
by the leftmost columns of the two
diagrams~(\ref{reduced-lcta-cdg-contrader-equivs-diagram}\+-%
\ref{reduced-X-lct-cdg-contrader-equivs-diagram}) from
Corollary~\ref{reduced-cdg-contraderived-equiv-cor}.
 The triangulated equivalences in the middle and rightmost columns
are the result of Corollary~\ref{lcth-reduced-X-lcta-X-lct-equivs-cor}.
 The commutativity of the diagram is obvious, so it follows that
the vertical functor in the leftmost column is a triangulated
equivalence, too.

 Part~(c): the horizontal triangulated equivalences
in~\eqref{antil-thick-all-lcta-lct-reduced-ctrder-diag} are provided
by the lower horizontal lines of commutative
diagrams~(\ref{reduced-lcta-cdg-contrader-equivs-diagram}\+-%
\ref{reduced-X-lct-cdg-contrader-equivs-diagram})
from Corollary~\ref{reduced-cdg-contraderived-equiv-cor}.
 The vertical triangulated equivalences are provided by the leftmost
columns of
diagrams~(\ref{reduced-lcth-al-X-lcta-X-lct-cdg-ctrder-diag}\+-%
\ref{reduced-thick-lcth-al-lcta-lct-cdg-ctrder-diag})
from parts~(a\+-b).
 The commutativity of the diagram is obvious.
\end{proof}

\Section{Co-Contra Correspondence for CDG-Modules}
\label{cdg-co-contra-secn}

 We keep the notation of Section~\ref{cdg-reduced-contraderived-secn};
so $\cB^\cu=\cC^\cu_X(\g,\widetilde\g)$ is the Chevalley--Eilenberg
quasi-coherent CDG\+quasi-algebra of a quasi-coherent twisted Lie
algebroid $(\g,\widetilde\g)$ over a scheme $X$ such that
the quasi-coherent sheaf~$\g$ on $X$ is finite locally free.
 The scheme $X$ will be eventually assumed to be quasi-compact and
semi-separated.

 The aim of this section is to establish a natural triangulated
equivalence between the reduced coderived category of quasi-coherent
CDG\+modules and the reduced contraderived category of $\bW$\+locally
contraherent CDG\+modules over~$\cB^\cu$,
$$
 \sD^\co_{X\red}(\cB^\cu\bQcoh)\simeq
 \sD^\ctr_{X\red}(\cB^\cu\bLcth_\bW).
$$
 This is the result of
Corollary~\ref{cdg-module-co-contra-correspondence-cor}.

\subsection{$\fHom$ and contratensor product of CDG-modules}
\label{fHom-and-contratensor-of-cdg-modules}
 The constructions of the functors $\fHom_\cA$ and~$\ocn_\cB$
from Sections~\ref{fHom-over-qcoh-quasi-algebra-subsecn}
and~\ref{contratensor-over-qcoh-quasi-algebra-subsecn} are applicable
to graded modules over graded quasi-algebras $\cA^*$ and $\cB^*$
according to the discussion in
Section~\ref{graded-modules-subsecn}.
 The aim of this section is to extend these constructions to the realm
of CDG\+modules over CDG\+quasi-algebras $\cA^\cu$ and~$\cB^\cu$.

 The definition of a quasi-coherent CDG\+module over a quasi-coherent
CDG\+quasi-algebra $\cA^\cu$ over a scheme $X$ was given in
Section~\ref{qcoh-lcth-cdg-modules-subsecn}.
 Let us define \emph{quasi-coherent CDG\+bimodules} over a pair of
quasi-coherent CDG\+quasi-algebras $\cA^\cu$ and $\cB^\cu$ over~$X$.

 Let $(f,a)\:A^\cu{}'\rarrow A^\cu$ and $(g,b)\:B^\cu{}'\rarrow B^\cu$
be morphisms of CDG\+rings.
 Let $M^\cu{}'=(M^*{}',d_{M'})$ be a CDG\+bimodule over $A^\cu{}'$
and $B^\cu{}'$, and let $M^\cu=(M^*,d_M)$ be a CDG\+bimodule over
$A^\cu$ and $B^\cu$.
 We will say that a homogeneous additive map $k\:M^\cu{}'\rarrow M^\cu$
of degree~$0$ is a \emph{map of CDG\+bimodules compatible with
the morphisms of CDG\+rings $(f,a)$ and~$(g,b)$} if $k$~is a morphism
of graded $A^*{}'$\+$B^*{}'$\+bimodules and the map~$k$ forms
a commutative square diagram with the differential $d_{M'}$ on $M^*{}'$
and the differential $d'_M$ on $M^*$ given by formula~(x) from
Section~\ref{cdg-rings-cdg-modules-subsecn}.

 Let $X$ be a scheme with an open covering $\bW$, and let $\cA^\cu$
and $\cB^\cu$ be quasi-coherent CDG\+quasi-algebras over~$X$.
 A \emph{quasi-coherent CDG\+bimodule $\cE^\cu$ over $\cA^\cu$
and~$\cB^\cu$} is a set of data consisting of
\begin{itemize}
\item a quasi-coherent graded bimodule $\cE^*=\bigoplus_{i\in\boZ}\cE^i$
over the quasi-coherent graded quasi-algebras $\cA^*$ and $\cB^*$
over~$X$ (see Section~\ref{fHom-over-qcoh-quasi-algebra-subsecn} for
the definition in the ungraded case and
Section~\ref{graded-modules-subsecn} for a discussion of the graded
context);
\item an odd derivation $d_{\cE,U}\:\cE^*(U)\rarrow\cE^*(U)$ of
degree~$1$ on the graded $\cA^*(U)$\+$\cB^*(U)$\+bimodule $\cE^*(U)$
compatible with the odd derivation $d_{\cA,U}\:\cA^*(U)\allowbreak
\rarrow\cA^*(U)$ on the graded ring $\cA^*(U)$ as well as with the odd
derivation $d_{\cB,U}\:\cB^*(U)\rarrow\cB^*(U)$ on the graded ring
$\cB^*(U)$, given for each affine open subscheme $U\subset X$
subordinate to~$\bW$.  \hfuzz=3.5pt
\end{itemize}
 The following two axioms must be satisfied:
\begin{enumerate}
\renewcommand{\theenumi}{\roman{enumi}}
\setcounter{enumi}{19}
\item for each affine open subscheme $U\subset X$ subordinate to $\bW$,
the pair $\cE^\cu(U)=(\cE^*(U),d_{\cE,U})$ is a CDG\+bimodule over
the CDG\+rings $\cA^\cu(U)=(\cA^*(U),d_{\cA,U},h_{\cA,U})$ and
$\cB^\cu(U)=(\cB^*(U),d_{\cB,U},h_{\cB,U})$;
\item for each pair of affine open subschemes $V\subset U\subset X$
subordinate to $\bW$, the map of restriction of sections
$\cE^*(U)\rarrow\cE^*(V)$ in the quasi-coherent graded quasi-module
$\cE^*$ is a map of CDG\+bimodules compatible with the morphisms of
CDG\+rings $(\rho_{\cA,VU},a_{\cA,VU})\:(\cA^*(U),d_{\cA,U},h_{\cA,U})
\rarrow(\cA^*(V),d_{\cA,V},h_{\cA,V})$ and
$(\rho_{\cB,VU},a_{\cB,VU})\:(\cB^*(U),d_{\cB,U},h_{\cB,U})
\rarrow(\cB^*(V),d_{\cB,V},h_{\cB,V})$.
\end{enumerate}

 Similarly to Lemma~\ref{strict-gluing-module-derivation-on-affine}(a)
and Corollary~\ref{qcoh-cdg-module-not-dependent-on-covering}, one
shows that the notion of a quasi-coherent CDG\+bimodule over two
quasi-coherent CDG\+quasi-algebras $\cA^\cu$ and $\cB^\cu$ over
a scheme~$X$ does not depend on an open covering~$\bW$.
 (The result of~\cite[Proposition~6.6]{Ptd} about uniqueness of
localizations of differential operators may be helpful for proving
that the left and right versions of the construction of the differential
in Corollary~\ref{cdg-quasi-algebra-co-extension-of-scalars}(a) agree
in the context of CDG\+bimodules.)

 In particular, any quasi-coherent CDG\+quasi-algebra $\cA^\cu$ over $X$
has an obvious natural structure of a CDG\+bimodule over itself.
 The differentials $d_{\cA,U}$ in the CDG\+bimodule $\cA^\cu$ over
$\cA^\cu$ and $\cA^\cu$ coincide with the differentials $d_{\cA,U}$
in the CDG\+quasi-algebra $\cA^\cu$ for all affine open subschemes
$U\subset X$.

 We will also need the following version of the previous definition.
 Let $X$ be a scheme with an open covering $\bW$, let $\cA^\cu$ be
a quasi-coherent CDG\+quasi-algebra over $X$, and let
$R^\cu=(R^*,d_R,h_R)$ be a CDG\+ring.
 A \emph{quasi-coherent CDG\+bimodule $\F^\cu$ over $\cA^\cu$
and~$R^\cu$} is a set of data consisting of
\begin{itemize}
\item a quasi-coherent graded left module
$\F^*=\bigoplus_{i\in\boZ}\F^i$ over the quasi-coherent graded
quasi-algebra $\cA^*$ over $X$, endowed with a right action of
the graded ring $R^*$ by homogenenous quasi-coherent left
$\cA^*$\+module endomorphisms, given for each affine open subscheme
$U\subset X$ subordinate to~$\bW$;
\item ad odd derivation $d_{\F,U}\:\F^*(U)\rarrow\F^*(U)$ of degree~$1$
on the graded $\cA^*(U)$\+$R^*$\+bimodule $\F^*(U)$ compatible with
the odd derivation $d_{\cA,U}\:\cA^*(U)\allowbreak\rarrow\cA^*(U)$ on
the graded ring $\cA^*(U)$ as well as with the odd derivation
$d_R\:R^*\rarrow R^*$ on the graded ring $R^*$, given for each affine
open subscheme $U\subset X$ subordinate to~$\bW$.
\end{itemize}
 The following two axioms must be satisfied:
\begin{enumerate}
\renewcommand{\theenumi}{\roman{enumi}}
\setcounter{enumi}{21}
\item for each affine open subscheme $U\subset X$ subordinate to $\bW$,
the pair $\F^\cu(U)=(\F^*(U),d_{\F,U})$ is a CDG\+bimodule over
the CDG\+rings $\cA^\cu(U)=(\cA^*(U),d_{\cA,U},h_{\cA,U})$ and
$R^\cu=(R,d_R,h_R)$;
\item for each pair of affine open subschemes $V\subset U\subset X$
subordinate to $\bW$, the map of restriction of sections
$\F^*(U)\rarrow\F^*(V)$ in the quasi-coherent sheaf $\F^*$ is a map of
CDG\+bimodules compatible with the morphisms of CDG\+rings
$(\rho_{\cA,VU},a_{\cA,VU})\:(\cA^*(U),d_{\cA,U},h_{\cA,U})
\rarrow(\cA^*(V),d_{\cA,V},h_{\cA,V})$ and
$(\id_{R^*},0)\:R^\cu\rarrow R^\cu$.
\end{enumerate}

 Similarly to Lemma~\ref{strict-gluing-module-derivation-on-affine}(a)
and Corollary~\ref{qcoh-cdg-module-not-dependent-on-covering}, one
shows that the notion of a quasi-coherent CDG\+bimodule over
a quasi-coherent CDG\+quasi-algebra $\cA^\cu$ and a CDG\+ring $R^\cu$
does not depend on an open covering $\bW$ of the scheme~$X$.

 Let $X$ be a scheme, $\cA^\cu$ be a quasi-coherent CDG\+quasi-algebra
over $X$, and $Y\subset X$ be an open subscheme with the open immersion
morphism $j\:Y\rarrow X$.
 Then the quasi-coherent graded quasi-algebra $j^*\cA^*=\cA^*|_Y$
has a natural structure of a CDG\+quasi-algebra $j^*\cA^\cu$ over $Y$,
with the same differentials $d_{j^*\cA,U}=d_{\cA,U}$, the same
curvature elements $h_{j^*\cA,U}=h_{\cA,U}$, and the same
change-of-connection elements $a_{j^*\cA,VU}=a_{\cA,VU}$ for all
affine open subschemes $V\subset U\subset Y$.

 Let $\cA^\cu$ and $\cB^\cu$ be two quasi-coherent CDG\+quasi-algebras
over~$X$.
 Let $\cE^\cu$ be a quasi-coherent CDG\+bimodule over $\cA^\cu$ and
$\cB^\cu$, and let $U\subset X$ be an affine open subscheme with
the open immersion morphism $j\:U\rarrow X$.
 Then the quasi-coherent graded quasi-module $j^*\cE^*$ on $U$ acquires
a natural structure of a quasi-coherent CDG\+bimodule $j^*\cE^\cu$ over
the quasi-coherent CDG\+quasi-algebra $j^*\cA^\cu$ and the CDG\+ring
$R^\cu=(\cB^*(U),d_{\cB,U},h_{\cB,U})$.
 There is one delicate aspect to this construction where one has to
be careful.
 In order to define the differential~$d_{j^*\cE,V}$ involved in
the quasi-coherent CDG\+bimodule $j^*\cE^\cu$ over $j^*\cA^\cu$ and
$R^\cu$ for an affine open subscheme $V\subset U$, one needs to twist
the differential~$d_{\cE,V}$ with the change-of-connection element
$a_{\cB,VU}\in\cB^1(V)$ as per formulas~(viii) and~(x)
from Section~\ref{cdg-rings-cdg-modules-subsecn}; so
$$
 d_{j^*\cE,V}(e)=d_{\cE,V}(e)-(-1)^{|e|}ea_{\cB,VU}
 \quad\textup{for all $e\in\cE^{|e|}(V)$}.
$$

 Let $Y\subset X$ be an open subscheme such that the open immersion
morphism $j\:Y\rarrow X$ is quasi-compact.
 Let $\cA^\cu$ be a quasi-coherent CDG\+quasi-algebra over $X$,
let $R^\cu$ be a CDG\+ring, and let $\G^\cu$ be a quasi-coherent
CDG\+bimodule over $j^*\cA^\cu$ and $R^\cu$ on~$Y$.
 Then the quasi-coherent CDG\+bimodule $j_*\G^\cu$ over $\cA^\cu$
and $R^\cu$ on $X$ is constructed in the following way.
 The underlying quasi-coherent graded $\cA^*$\+module of $j_*\G^\cu$
is the quasi-coherent graded $\cA^*$\+module $j_*\G^*$ provided
by the graded version of the construction of
the functor~\eqref{qcoh-A-modules-direct-image} in
Section~\ref{direct-images-of-A-co-sheaves-subsecn}.
 The right action of $R^*$ on $j_*\G^*$ is induced by the right
action of $R^*$ on $\G^*$ in the obvious way.

 Let us explain how to construct the differential $d_{j_*\G,U}\:
(j_*\G^*)(U)\rarrow(j_*\G^*)(U)$ for every affine open subscheme
$U\subset X$.
 Let $Y\cap U=\bigcup_\alpha V_\alpha$ be a finite affine open
covering of the intersection of open subschemes $Y\cap U\subset X$.
 Let $d'_{\G,V_\alpha}\:\G^*(V_\alpha)\rarrow\G^*(V_\alpha)$ be
the differential obtained from the differential~$d_{\G,V_\alpha}$
by twisting with the change-of-connection element
$a_{\cA,V_\alpha U}\in\cA^1(V_\alpha)$ as per formulas~(vi) and~(x)
from Section~\ref{cdg-rings-cdg-modules-subsecn}; so
$d'_{\G,V_\alpha}(z)=d_{\G,V_\alpha}(z)+a_{\cA,V_\alpha U}z$
for all $z\in\G(V_\alpha$).
 Similarly, we notice that, for every pair of indices $\alpha$
and~$\beta$, the intersection $V_\alpha\cap V_\beta$ is an affine
open subscheme in $X$, and define the differential
$d'_{\G,V_\alpha\cap V_\beta}\:\G^*(V_\alpha\cap V_\beta)\rarrow
\G^*(V_\alpha\cap V_\beta)$ by twisting the differential
$d_{\G,V_\alpha\cap V_\beta}$ with the change-of-connection element
$a_{\cA,(V_\alpha\cap V_\beta)U}\in\cA^1(V_\alpha\cap V_\beta)$.
 Finally, consider the \v Cech right exact sequence
\begin{equation} \label{cech-right-exact-sequence}
 \bigoplus\nolimits_{\alpha<\beta}\G^*(V_\alpha\cap V_\beta)
 \lrarrow\bigoplus\nolimits_\alpha\G^*(V_\alpha)\lrarrow
 (j_*\G)(U)\lrarrow0.
\end{equation}
 Define the differential~$d$ on
$\bigoplus_{\alpha<\beta}\G^*(V_\alpha\cap V_\beta)$ as the direct
sum of the differentials $d'_{\G,V_\alpha\cap V_\beta}$ and
the differential~$d$ on $\bigoplus_\alpha\G^*(V_\alpha)$ as
the direct sum of the differentials~$d'_{\G,V_\alpha}$.
 Then the leftmost arrow in~\eqref{cech-right-exact-sequence} forms
a commutative square diagram with the differentials~$d$.
 Passing to the cokernel, we obtain the desired differential
$d_{j_*\G,U}\:(j_*\G^*)(U)\rarrow(j_*\G^*)(U)$.

 Let $\F^\cu$ be a quasi-coherent CDG\+bimodule over $\cA^\cu$ and
$R^\cu$ on $X$, and let $\C^\cu$ be a quasi-coherent left
CDG\+module over~$\cA^\cu$.
 Then the left CDG\+module $\Hom^\cu_{\cA^*}(\F^\cu,\C^\cu)$ over
the CDG\+ring $R^\cu$ is constructed as follows.
 The underlying graded abelian group of $\Hom^*_{\cA^*}(\F^*,\C^*)$
of the CDG\+module $\Hom^\cu_{\cA^*}(\F^\cu,\C^\cu)$ is constructed
as explained in Section~\ref{dg-categories-of-cdg-modules-subsecn}.
 The left action of the graded ring $R^*$ on
$\Hom^*_{\cA^*}(\F^*,\C^*)$ is induced by the right action of $R^*$
on $\F^*$, with the sign rule as mentioned in
Section~\ref{cdg-rings-cdg-modules-subsecn} and spelled out
in~\cite[Section~6.1]{Prel}.
 Finally, given a homogeneous map of quasi-coherent graded
$\cA^*$\+modules $f\in\Hom_{\cA^*}^{|f|}(\F^*,\C^*)$, its differential
$d(f)\in\Hom_{\cA^*}^{|f|+1}(\F^*,\C^*)$ is defined by the rule
$d(f)(U)=d(f(U))\:\F^*(U)\rarrow\M^*(U)$ for all affine open subschemes
$U\subset X$, where $d(f(U))$ is the differential of the homogeneous
$\cA^*(U)$\+linear map $f(U)\:\F^*(U)\rarrow\C^*(U)$ between
the underlying graded $\cA^*(U)$\+modules of the CDG\+bimodule
$\F^\cu(U)$ over $\cA^\cu(U)$ and $R^\cu$ and the left CDG\+module
$\C^\cu(U)$ over $\cA^\cu(U)$, as per the formula in
Section~\ref{cdg-rings-cdg-modules-subsecn} (see also
Section~\ref{dg-categories-of-cdg-modules-subsecn}).

 We follow the exposition in
Section~\ref{fHom-over-qcoh-quasi-algebra-subsecn}.
 Let $X$ be a quasi-separated scheme, and let $\cA^\cu$  and $\cB^\cu$
be two quasi-coherent CDG\+quasi-algebras over~$X$.
 Let $\cE^\cu$ be a quasi-coherent CDG\+bimodule over $\cA^\cu$ and
$\cB^\cu$, and let $\C^\cu$ be a quasi-coherent left CDG\+module
over~$\cA^\cu$.
 Assume that the quasi-coherent graded bimodule $\cE^*$ over $\cA^*$
and $\cB^*$ and the quasi-coherent graded $\cA^*$\+module $\C^*$
satisfy one of the sufficient conditions for the contraherent graded
$\cB^*$\+module $\fHom_{\cA^*}^*(\cE^*,\C^*)$ to be well-defined, as per
the graded version of Lemma~\ref{fHom-contraadjusted-and-more-lemma}
and the subsequent discussion in
Section~\ref{fHom-over-qcoh-quasi-algebra-subsecn}.

 Let $U\subset X$ be an affine open subscheme with the open immersion
morphism $j\:U\rarrow X$.
 Following the constructions above, we have a quasi-coherent
CDG\+bi\-mod\-ule $\G^\cu=j^*\cE^\cu$ over the quasi-coherent
CDG\+quasi-algebra $j^*\cA^\cu$ and the CDG\+ring
$R^\cu=(\cB^*(U),d_{\cB,U},h_{\cB,U})$.
 Furthermore, we have the quasi-coherent CDG\+bimodule
$\F^\cu=j_*\G^\cu=j_*j^*\cE^\cu$ over the CDG\+quasi-algebra $\cA^\cu$
and the CDG\+ring~$R^\cu$.
 Finally, the construction above provides a left CDG\+module
$\Hom^\cu_{\cA^*}(\F^\cu,\C^\cu)=\Hom^\cu_{\cA^*}(j_*j^*\cE^\cu,\C^\cu)$
over the CDG\+ring $R^\cu=(\cB^*(U),d_{\cB,U},h_{\cB,U})$.

 The underlying graded left
$R^*$\+module $\Hom^\cu_{\cA^*}(\F^\cu,\C^\cu)$ of the CDG\+module
$\Hom^\cu_{\cA^*}(\F^\cu,\C^\cu)$ is the graded left module
$\fHom_{\cA^*}^*(\cE^*,\C^*)[U]$ over the graded ring $\cB^*(U)$,
as per the graded version of the construction of
Section~\ref{fHom-over-qcoh-quasi-algebra-subsecn}.
 So we have produced a differential on the graded left module
$\fHom_{\cA^*}^*(\cE^*,\C^*)[U]$ over the graded ring~$\cB^*(U)$.
 The collection of all such differentials, indexed over the affine
open subschemes $U\subset X$, defines a structure of contraherent
left CDG\+module $\fHom^\cu_{\cA^*}(\cE^\cu,\C^\cu)$ over
the quasi-coherent CDG\+quasi-algebra $\cB^\cu$ on the contraherent
graded left module $\fHom_{\cA^*}^*(\cE^*,\C^*)$ over the quasi-coherent
graded quasi-algebra~$\cB^*$. {\hbadness=1375\par}

 In particular, if the quasi-coherent graded $\cA^*$\+$\cB^*$\+bimodule
$\cE^*$ is $\cA^*$\+very flaprojective (in the sense of the graded
versions of the definitions from
Sections~\ref{prelim-flaprojective-subsecn},
\ref{antilocality-of-X-contraadjusted-subsecn},
and~\ref{fHom-over-qcoh-quasi-algebra-subsecn}), then the construction
above provides a DG\+functor
\begin{equation} \label{fHom-DG-functor-for-X-contraadjusted}
 \fHom^\cu_{\cA^*}(\cE^\cu,{-})\:\cA^\cu\bQcoh^{X\dcta}
 \lrarrow\cB^\cu\bCtrh.
\end{equation}
 Similarly, if the quasi-coherent graded $\cA^*$\+$\cB^*$\+bimodule
$\cE^*$ is $\cA^*$\+robustly flaprojective (in the sense of the graded
versions of the definitions from
Sections~\ref{prelim-robustly-flaprojective-subsecn}
and~\ref{fHom-over-qcoh-quasi-algebra-subsecn}), then the construction
above provides a DG\+functor
\begin{equation} \label{fHom-DG-functor-for-X-cotorsion}
 \fHom^\cu_{\cA^*}(\cE^\cu,{-})\:\cA^\cu\bQcoh^{X\dcot}
 \lrarrow\cB^\cu\bCtrh^{X\dlct}.
\end{equation}
 Following~\cite[Section~6.1]{Prel}, there is \emph{no} sign rule
involved in the construction of the action of the DG\+functors
$\fHom^\cu_{\cA^*}(\cE^\cu,{-})$ on homogeneous morphisms of odd
cohomological degrees in the DG\+categories of quasi-coherent
CDG\+modules over~~$\cA^\cu$.

 Let $X$ be a scheme and $\cA^\cu$ be a quasi-coherent
CDG\+quasi-algebra over~$X$.
 Let $R^\cu$ be a CDG\+ring, $\F^\cu$ be a quasi-coherent CDG\+bimodule
over $\cA^\cu$ and $R^\cu$, and $M^\cu$ be a left CDG\+module
over~$R^\cu$.
 Then the tensor product $\F^\cu\ot_{R^*}M^\cu$ is a quasi-coherent
left CDG\+module over $\cA^\cu$ defined by the rule
$$
 (\F^\cu\ot_{R^*}M^\cu)(U)=\F^\cu(U)\ot_{R^*}M^\cu
$$
for all affine open subschemes $U\subset X$.
 Here $\F^\cu(U)\ot_{R^*}M^\cu$ is the tensor product of
the CDG\+bimodule $\F^\cu(U)$ over the CDG\+rings $\cA^\cu(U)$ and
$R^\cu$ with the left CDG\+module $M^\cu$ over~$R^\cu$.
 So $\F^\cu(U)\ot_{R^*}M^\cu$ is a left CDG\+module over $\cA^\cu(U)$,
as per the discussion in Section~\ref{cdg-rings-cdg-modules-subsecn}.
 The collection of the differentials of the CDG\+modules
$(\F^\cu\ot_{R^*}M^\cu)(U)=\F^\cu(U)\ot_{R^*}M^\cu$ over $\cA^\cu(U)$,
indexed over the affine open subschemes $U\subset X$, defines
the desired structure of a quasi-coherent left CDG\+module
$\F^\cu\ot_{R^*}M^\cu$ over $\cA^\cu$ on the quasi-coherent graded
left $\cA^*$\+module $\F^*\ot_{R^*}M^*$ produced by the graded version
of the construction from
Section~\ref{fHom-over-qcoh-quasi-algebra-subsecn}.

 We roughly follow the exposition in
Section~\ref{contratensor-over-qcoh-quasi-algebra-subsecn}.
 Let $X$ be a quasi-separated scheme with an open covering $\bW$, and
let $\cA^\cu$ and $\cB^\cu$ be two quasi-coherent CDG\+quasi-algebras
over~$X$.
 Let $\cE^\cu$ be a quasi-coherent CDG\+bimodule over $\cA^\cu$ and
$\cB^\cu$, and let $\P^\cu$ be a $\bW$\+locally contraherent
CDG\+module over~$\cB^\cu$.
 
 For every affine open subscheme $U\subset X$ subordinate to $\bW$,
we denote by $j\:U\rarrow X$ the open immersion morphism, and
consider the quasi-coherent CDG\+bimodule $\F^\cu=j_*j^*\cE^\cu$
over the quasi-coherent CDG\+bimodule $\F^\cu=j_*\G^\cu=j_*j^*\cE^\cu$
over the CDG\+quasi-algebra $\cA^\cu$ and the CDG\+ring
$R^\cu=\cB^\cu(U)=(\cB^*(U),d_{\cB,U},h_{\cB,U})$, as above.
 We also have the left CDG\+module $\P^\cu(U)=(\P^*(U),d_{\P^*,U})$
over $R^\cu=\cB^\cu(U)$.
 The construction of the tensor product above provides a quasi-coherent
left CDG\+module $\F^\cu\ot_{R^*}\P^\cu(U)=
j_*j^*\cE^\cu\ot_{R^*}\P^\cu(U)$ over the CDG\+quasi-algebra~$\cA^\cu$.

 Let $\bB$ denote the topology base of affine open subschemes in $X$
subordinate to~$\bW$.
 The quasi-coherent left CDG\+module $\cA^\cu\ocn_{\cB^*}\P^\cu$ over
the quasi-coherent CDG\+quasi-algebra $\cA^\cu$ is defined as
the inductive limit
$$
 \cE^\cu\ocn_{\cB^*}\P^\cu=
 \varinjlim\nolimits_{U\in\bB}((j_*j^*\cE^\cu)\ot_{\cB^*(U)}\P^\cu[U]).
$$
 Here the (nonfiltered) inductive limit is taken in the abelian
category of quasi-coherent CDG\+modules $\sZ^0(\cA^\cu\bQcoh)$,
while the transition maps in the diagram are provided by
the construction from
Section~\ref{contratensor-over-qcoh-quasi-algebra-subsecn}.
 Recall that the forgetful functor $\sZ^0(\cA^\cu\bQcoh)\rarrow
\cA^*\Qcoh$ preserves all limits and colimits (and in fact, has adjoints
on both sides by Lemma~\ref{G-plus-G-minus-functors-for-qcoh-lcth}(a)).

 We have defined a structure of quasi-coherent left CDG\+module
$\cE^\cu\ocn_{\cB^*}\P^\cu$ over the quasi-coherent CDG\+quasi-algebra
$\cA^\cu$ on the quasi-coherent graded left module
$\cE^*\ocn_{\cB^*}\P^*$ over $\cA^*$ provided by the graded version
of the construction from
Section~\ref{contratensor-over-qcoh-quasi-algebra-subsecn}.
 So we obtain a DG\+functor
\begin{equation} \label{contratensor-DG-functor}
 \cE^\cu\ocn_{\cB^*}{-}\,\:\cB^\cu\bLcth_\bW\lrarrow\cA^\cu\bQcoh.
\end{equation}
 There is a sign rule involved in the construction of the action of
the DG\+functor $\cE^\cu\ocn_{\cB^*}\nobreak{-}$ on homogeneous
morphisms of odd cohomological degrees in the DG\+cat\-e\-gories of
locally contraherent CDG\+modules over~$\cB^\cu$;
see~\cite[Section~6.1]{Prel}.

 The following lemma establishes a partial adjunction between
the DG\+func\-tors~\eqref{fHom-DG-functor-for-X-contraadjusted}
or~\eqref{fHom-DG-functor-for-X-cotorsion}, on the one hand, and
the DG\+functor~\eqref{contratensor-DG-functor}, on the other hand.

\begin{lem} \label{fHom-contratensor-DG-adjunction}
 Let $X$ be a quasi-separated scheme with an open covering $\bW$, and
let $\cA^\cu$ and $\cB^\cu$ be two quasi-coherent CDG\+quasi-algebras
over~$X$.
 Let $\cE^\cu$ be a quasi-coherent CDG\+bimodule over $\cA^\cu$ and
$\cB^\cu$, left $\C^\cu$ be a quasi-coherent left CDG\+module over
$\cA^\cu$, and let\/ $\P^\cu$ be a $\bW$\+locally contraherent
CDG\+module over~$\cB^\cu$.
  Assume that the quasi-coherent graded bimodule $\cE^*$ over $\cA^*$
and $\cB^*$ and the quasi-coherent graded $\cA^*$\+module $\C^*$
satisfy one of the sufficient conditions for the contraherent graded
$\cB^*$\+module\/ $\fHom_{\cA^*}^*(\cE^*,\C^*)$ to be well-defined,
as per the graded version of the results of
Section~\textup{\ref{fHom-over-qcoh-quasi-algebra-subsecn}}.
 Then the graded version of the adjunction
isomorphism~\eqref{fHom-contratensor-adjunction} from
Section~\textup{\ref{contratensor-over-qcoh-quasi-algebra-subsecn}},
$$
 \Hom^{\cB^*,*}(\P^*,\fHom_{\cA^*}^*(\cE^*,\C^*))\simeq
 \Hom_{\cA^*}^*(\cE^*\ocn_{\cB^*}\P^*,\>\C^*)
$$
is actually an isomorphism of complexes of abelian groups
\begin{equation} \label{fHom-contratensor-adjunction-of-complexes}
 \Hom^{\cB^*,\bu}(\P^\cu,\fHom_{\cA^*}^\cu(\cE^\cu,\C^\cu))\simeq
 \Hom_{\cA^*}^\bu(\cE^\cu\ocn_{\cB^*}\P^\cu,\>\C^\cu).
\end{equation}
\end{lem}

\begin{proof}
 We leave this straightforward verification to the reader.
\end{proof}

\subsection{DG-category equivalences}
 The next lemma follows directly from the results of
Section~\ref{naive-co-contra-subsecn} and the discussion in
Section~\ref{fHom-and-contratensor-of-cdg-modules}.
 The definition of an exact DG\+functor between exact DG\+categories
can be found in~\cite[Section~4.4]{Pedg} or in
Section~\ref{abelian-and-exact-dg-categs-of-cdg-modules-subsecn} above.

\begin{lem} \label{cta-cot-al-co-contra-dg-equivalence}
 Let $X$ be a quasi-compact semi-separated scheme with an open
covering\/ $\bW$ and $(\g,\widetilde\g)$ be a quasi-coherent twisted
Lie algebroid over~$X$.
 Assume that\/ $\g$~is a finite locally free sheaf on $X$, and let
$\cB^\cu=\cC^\cu_X(\g,\widetilde\g)$ be the related Chevalley--Eilenberg
quasi-coherent CDG\+quasi-algebra over~$X$.
 In this setting: \par
\textup{(a)} The adjoint DG\+functors\/ $\cB^\cu\ocn_{\cB^*}{-}$ and\/
$\fHom^\cu_{\cB^*}(\cB^\cu,{-})$ are mutually inverse equivalences
of exact DG\+categories
$$
 \fHom^\cu_{\cB^*}(\cB^\cu,{-})\:
 \cB^\cu\bQcoh^{X\dcta}\simeq
 \cB^\cu\bCtrh_\al\,:\!\cB^\cu\ocn_{\cB^*}{-}.
$$ \par
\textup{(b)} The adjoint DG\+functors\/ $\cB^\cu\ocn_{\cB^*}{-}$ and\/
$\fHom^\cu_{\cB^*}(\cB^\cu,{-})$ are mutually inverse equivalences
of exact DG\+categories
$$
 \fHom^\cu_{\cB^*}(\cB^\cu,{-})\:
 \cB^\cu\bQcoh^{X\dcot}\simeq
 \cB^\cu\bCtrh_\al^{X\dlct}\,:\!\cB^\cu\ocn_{\cB^*}{-}.
$$ \par
\textup{(c)} The DG\+category equivalences from parts~(a) and~(b) form
a commutative diagram of exact DG\+functors and exact DG\+category
equivalences
$$
 \qquad\quad\xymatrix{
  \text{\llap{$\fHom^\cu_{\cB^*}(\cB^\cu,{-})\:$}}
  \cB^\cu\bQcoh^{X\dcta} \ar@{=}[r]
  & \cB^\cu\bCtrh_\al
  \text{\rlap{$\,\,:\!\cB^\cu\ocn_{\cB^*}{-}$}}  \\
  \text{\llap{$\fHom^\cu_{\cB^*}(\cB^\cu,{-})\:$}}
  \cB^\cu\bQcoh^{X\dcot} \ar@{=}[r] \ar@{>->}[u]
  & \cB^\cu\bCtrh_\al^{X\dlct}
  \text{\rlap{$\,\,:\!\cB^\cu\ocn_{\cB^*}{-}$}} \ar@{>->}[u]
 }
$$
where the vertical functors, shown by arrows with tails, are
the identity inclusions of DG\+categories.
\end{lem}

\begin{proof}
 Parts~(a) and~(b): first of all, by the graded version of
Lemma~\ref{quasi-algebra-underived-naive-co-contra}(a\+-b),
the functor $\fHom_{\cB^*}^*(\cB^*,{-})$ takes $X$\+contraadjusted
quasi-coherent graded $\cB^*$\+modules to antilocal contraherent
graded $\cB^*$\+modules.
 By the same lemma, the functor $\cB^*\ocn_{\cB^*}{-}$ takes antilocal
contraherent graded $\cB^*$\+modules to $X$\+contraadjusted
quasi-coherent graded $\cB^*$\+modules, and antilocal $X$\+locally
cotorsion contraherent graded $\cB^*$\+modules to $X$\+cotorsion
quasi-coherent graded $\cB^*$\+modules.
 Therefore, the DG\+functors in part~(a) and~(b) are well-defined.
 The respective DG\+functors are adjoint to each other by
Lemma~\ref{fHom-contratensor-DG-adjunction}.
 It remains to point out that, according to the graded version of
Lemma~\ref{quasi-algebra-underived-naive-co-contra}(a\+-b),
the DG\+functors in question induce equivalences of the underlying
exact categories of graded objects,
\begin{align*}
 \fHom^*_{\cB^*}(\cB^*,{-})\:
 \cB^*\Qcoh^{X\dcta} &\simeq
 \cB^*\Ctrh_\al\,:\!\cB^*\ocn_{\cB^*}{-}, \\
 \fHom^*_{\cB^*}(\cB^*,{-})\:
 \cB^*\Qcoh^{X\dcot} &\simeq
 \cB^*\Ctrh_\al^{X\dlct}\,:\!\cB^*\ocn_{\cB^*}{-}
\end{align*}
(cf.\ Corollaries~\ref{exact-dg-categories-of-qcoh-cta-cot-cdg-modules}
and~\ref{exact-dg-categories-of-antilocal-ctrh-cdg-modules}).
 Part~(c) is clear from the constructions of the functors involved.
\end{proof}

\begin{lem} \label{cta-cot-al-thick-co-contra-dg-equivalence}
 Let $X$ be a quasi-compact semi-separated scheme with an open
covering\/ $\bW$ and $(\g,\widetilde\g)$ be a quasi-coherent twisted
Lie algebroid over~$X$.
 Assume that\/ $\g$~is a finite locally free sheaf on $X$, and let
$\cB^\cu=\cC^\cu_X(\g,\widetilde\g)$ be the related Chevalley--Eilenberg
quasi-coherent CDG\+quasi-algebra over~$X$.
 In this setting: \par
\textup{(a)} The adjoint DG\+functors\/ $\cB^\cu\ocn_{\cB^*}{-}$ and\/
$\fHom^\cu_{\cB^*}(\cB^\cu,{-})$ restrict to mutually inverse
equivalences of exact DG\+categories
$$
 \fHom^\cu_{\cB^*}(\cB^\cu,{-})\:
 \cB^\cu\bQcoh^{X\dcta}_\bth\simeq
 \cB^\cu\bCtrh_\al^\bth\,:\!\cB^\cu\ocn_{\cB^*}{-}.
$$ \par
\textup{(b)} The adjoint DG\+functors\/ $\cB^\cu\ocn_{\cB^*}{-}$ and\/
$\fHom^\cu_{\cB^*}(\cB^\cu,{-})$ restrict to mutually inverse
equivalences of exact DG\+categories
$$
 \fHom^\cu_{\cB^*}(\cB^\cu,{-})\:
 \cB^\cu\bQcoh^{X\dcot}_\bth\simeq
 \cB^\cu\bCtrh_\al^{X\dlct,\bth}\,:\!\cB^\cu\ocn_{\cB^*}{-}.
$$
\end{lem}

\begin{proof}
 We only need to point out that the DG\+category equivalences from
Lemma~\ref{cta-cot-al-co-contra-dg-equivalence} take thick
quasi-coherent CDG\+modules to thick contraherent CDG\+modules and back.
 This is essentially the result of
Theorem~\ref{thick-corresponds-to-thick-under-co-contra}.
\end{proof}

\subsection{Co-contra correspondence}  \label{co-contra-subsecn}
 Once again, we refer to the introduction to the paper~\cite{Pmgm} for
a general discussion of the co-contra correspondence phenomenon.
 This section is a nonaffine, nonsmooth generalization
of~\cite[Section~B.3]{Pkoszul}.
 As compared to~\cite[Section~8.5]{Pkoszul}, the results of this section
are more general in that they make no assumptions of finiteness of
homological dimension and apply to nonaffine schemes, but also less
general in that they are only applicable to Chevalley--Eilenberg
CDG\+rings of twisted Lie algebroids.

\begin{thm} \label{thick-reduced-co-contra-correspond-theorem}
 Let $X$ be a quasi-compact semi-separated scheme and
$(\g,\widetilde\g)$ be a quasi-coherent twisted Lie algebroid over~$X$.
 Assume that\/ $\g$~is a finite locally free sheaf on $X$, and let
$\cB^\cu=\cC^\cu_X(\g,\widetilde\g)$ be the related Chevalley--Eilenberg
quasi-coherent CDG\+quasi-algebra over~$X$.
 In this setting: \par
\textup{(a)} The adjoint DG\+functors\/ $\cB^\cu\ocn_{\cB^*}{-}$ and\/
$\fHom^\cu_{\cB^*}(\cB^\cu,{-})$ induce mutually inverse triangulated
equivalences between the reduced absolute derived categories
$$
 \fHom^\cu_{\cB^*}(\cB^\cu,{-})\:
 \sD^\abs_{X\red}(\cB^\cu\bQcoh^{X\dcta}_\bth)\,\simeq\,
 \sD^\abs_{X\red}(\cB^\cu\bCtrh_\al^\bth)
 \,:\!\cB^\cu\ocn_{\cB^*}{-}.
$$ \par
\textup{(b)} The adjoint DG\+functors\/ $\cB^\cu\ocn_{\cB^*}{-}$ and\/
$\fHom^\cu_{\cB^*}(\cB^\cu,{-})$ induce mutually inverse triangulated
equivalences between the reduced absolute derived categories
$$
 \fHom^\cu_{\cB^*}(\cB^\cu,{-})\:
 \sD^\abs_{X\red}(\cB^\cu\bQcoh^{X\dcot}_\bth)\,\simeq\,
 \sD^\abs_{X\red}(\cB^\cu\bCtrh_\al^{X\dlct,\bth})
 \,:\!\cB^\cu\ocn_{\cB^*}{-}.
$$ \par
\textup{(c)} The category equivalences from parts~(a) and~(b) form
a commutative diagram of triangulated functors and triangulated
equivalences
\begin{equation} \label{thick-reduced-co-contra-X-cta-X-cot-diagram}
\begin{gathered}
 \qquad\quad\xymatrix{
  \text{\llap{$\fHom^\cu_{\cB^*}(\cB^\cu,{-})\:$}}
  \sD^\abs_{X\red}(\cB^\cu\bQcoh^{X\dcta}_\bth) \ar@{=}[r]
  & \sD^\abs_{X\red}(\cB^\cu\bCtrh_\al^\bth)
  \text{\rlap{$\,\,:\!\cB^\cu\ocn_{\cB^*}{-}$}} \\
  \text{\llap{$\fHom^\cu_{\cB^*}(\cB^\cu,{-})\:$}}
  \sD^\abs_{X\red}(\cB^\cu\bQcoh^{X\dcot}_\bth)
  \ar@{=}[r] \ar[u]
  & \sD^\abs_{X\red}(\cB^\cu\bCtrh_\al^{X\dlct,\bth})
  \text{\rlap{$\,\,:\!\cB^\cu\ocn_{\cB^*}{-}$}} \ar[u]
 }
\end{gathered}
\end{equation}
where the vertical functors are induced by the inclusions
of exact DG\+categories $\cB^\cu\bQcoh^{X\dcot}_\bth\rarrow
\cB^\cu\bQcoh^{X\dcta}_\bth$ and $\cB^\cu\bCtrh_\al^{X\dlct,\bth}
\rarrow\cB^\cu\bCtrh_\al^\bth$.
\end{thm}

\begin{proof}
 Parts~(a) and~(b): the equivalences of exact DG\+categories from
Lemma~\ref{cta-cot-al-thick-co-contra-dg-equivalence} induce
equivalences of absolute derived categories
\begin{align*}
 \fHom^\cu_{\cB^*}(\cB^\cu,{-})\:
 \sD^\abs(\cB^\cu\bQcoh^{X\dcta}_\bth) &\,\simeq\,
 \sD^\abs(\cB^\cu\bCtrh_\al^\bth)
 \,:\!\cB^\cu\ocn_{\cB^*}{-}, \\
 \fHom^\cu_{\cB^*}(\cB^\cu,{-})\:
 \sD^\abs(\cB^\cu\bQcoh^{X\dcot}_\bth) &\,\simeq\,
 \sD^\abs(\cB^\cu\bCtrh_\al^{X\dlct,\bth})
 \,:\!\cB^\cu\ocn_{\cB^*}{-}.
\end{align*}
 It remains to show that these triangulated equivalences take
reduced-acyclic CDG\+modules to reduced-acyclic CDG\+modules.
 For this purpose, it is convenient to use the definition of
reduced-acyclicity of thick ($X$\+contraadjusted or $X$\+cotorsion)
CDG\+modules based on the filtration $F$ from
Corollary~\ref{thick-qcoh-sheaves-are-filtered}(b), as per
Remark~\ref{reduced-acyclicity-of-qcoh-with-other-filtration}.
{\hbadness=1225\par}

 Let $\M^\cu$ and $\P^\cu$ be a quasi-coherent and a contraherent
CDG\+module corresponding to each other under the equivalence of
exact DG\+categories from
Lemma~\ref{cta-cot-al-thick-co-contra-dg-equivalence}(a) or~(b).
 Consider the canonical filtration $F$ on $\M^\cu$ as per
Corollary~\ref{thick-qcoh-sheaves-are-filtered}(b)
and the canonical filtration $F$ on $\P^*$ as per
Corollary~\ref{thick-lcta-lct-cosheaves-are-filtered}.
 According to Theorem~\ref{thick-corresponds-to-thick-under-co-contra},
the filtered CDG\+modules $(\M^\cu,F)$ and $(\P^\cu,F)$ correspond
to each other under the equivalence of categories from the upper line
of the commutative diagram~\eqref{X-cta-filtered-naive-co-contra}
or~\eqref{X-cot-filtered-naive-co-contra}
in Theorem~\ref{quasi-algebra-filtered-naive-co-contra}.

 In view of commutativity of the diagrams in
Theorem~\ref{quasi-algebra-filtered-naive-co-contra} and commutativity
of the diagrams in
Lemma~\ref{quasi-algebra-underived-naive-co-contra}(a\+-b), it follows
that the complex of quasi-coherent sheaves $\gr_F^0\M^\cu$ and
the complex of contraherent cosheaves $\gr_F^0\P^\cu$ correspond to
each other under the equivalence of exact categories
$X\Qcoh^\cta\simeq X\Ctrh_\al$ or $X\Qcoh^\cot\simeq X\Ctrh_\al^\lct$
from~\cite[Lemma~4.8.2 or~4.8.4(a)]{Pcosh}.
 Thus the complex $\gr_F^0\M^\cu$ is acyclic in $X\Qcoh^\cta$
(respectively, in $X\Qcoh^\cot$) if and only if the complex
$\gr_F^0\P^\cu$ is acyclic in $X\Ctrh_\al$ (respectively, in
$X\Ctrh_\al^\lct$).

 Part~(c) follows from
Lemma~\ref{cta-cot-al-co-contra-dg-equivalence}(c).
\end{proof}

 We will see below in
Corollary~\ref{thick-reduced-co-contra-cta-cot-equivs-cor} that all
the functors in commutative
diagram~\eqref{thick-reduced-co-contra-X-cta-X-cot-diagram}
of Theorem~\ref{thick-reduced-co-contra-correspond-theorem}
are triangulated equivalences.

 The following corollary is the main result of
Section~\ref{cdg-co-contra-secn}.

\begin{cor} \label{cdg-module-co-contra-correspondence-cor}
 Let $X$ be a quasi-compact semi-separated scheme with an open
covering\/ $\bW$ and $(\g,\widetilde\g)$ be a quasi-coherent twisted
Lie algebroid over~$X$.
 Assume that\/ $\g$~is a finite locally free sheaf on $X$, and let
$\cB^\cu=\cC^\cu_X(\g,\widetilde\g)$ be the related Chevalley--Eilenberg
quasi-coherent CDG\+quasi-algebra over~$X$.
 Then there is a natural triangulated equivalence between the reduced
coderived and the reduced contraderived category
\begin{equation} \label{cdg-module-co-contra-formula}
 \boR\fHom_{\cB^*}^\cu(\cB^\cu,{-})\:
 \sD^\co_{X\red}(\cB^\cu\bQcoh)\,\simeq\,
 \sD^\ctr_{X\red}(\cB^\cu\bLcth_\bW)\,:\!
 \cB^\cu\ocn_{\cB^*}^\boL{-}
\end{equation}
provided by the right derived functor of contraherent\/ $\fHom$ from
$\cB^\cu$ and the left derived functor of contratensor product with
$\cB^\cu$ over~$\cB^*$.
 The desired derived functors are constructed as the compositions of
triangulated equivalences on the diagram
\begin{equation} \label{cdg-module-co-contra-construction-diagram}
\begin{gathered}
 \xymatrix{
  \sD^\co_{X\red}(\cB^\cu\bQcoh) \ar@{-}@<2pt>[d]
  & \sD^\ctr_{X\red}(\cB^\cu\bLcth_\bW) \ar@{-}@<-2pt>[d] \\
  \sD^\co_{X\red}(\cB^\cu\bQcoh_\bth) \ar@<2pt>[u] \ar@{=}[d]
  & \sD^\ctr_{X\red}(\cB^\cu\bLcth_\bW^\bth)
  \ar@<-2pt>[u] \ar@{=}[d] \\
  \sD^\abs_{X\red}(\cB^\cu\bQcoh_\bth) \ar@{-}@<2pt>[d]
  & \sD^\abs_{X\red}(\cB^\cu\bLcth_\bW^\bth) \ar@{-}@<-2pt>[d] \\
  \text{\llap{$\fHom^\cu_{\cB^*}(\cB^\cu,{-})\:$}}
  \sD^\abs_{X\red}(\cB^\cu\bQcoh^{X\dcta}_\bth)
  \ar@{=}[r] \ar@<2pt>[u]
  & \sD^\abs_{X\red}(\cB^\cu\bCtrh_\al^\bth)
  \text{\rlap{$\,\,:\!\cB^\cu\ocn_{\cB^*}{-}$}} \ar@<-2pt>[u]
 }
\end{gathered}
\end{equation}
\end{cor}

\begin{proof}
 The triangulated equivalence in the lower line of the diagram
is the result of
Theorem~\ref{thick-reduced-co-contra-correspond-theorem}(a).
 The triangulated equivalence
$\sD^\abs_{X\red}(\cB^\cu\bQcoh^{X\dcta}_\bth)\simeq
\sD^\abs_{X\red}(\cB^\cu\bQcoh_\bth)$ in the leftmost column
is provided by
Corollary~\ref{thick-cdg-qcoh-cta-reduced-abs-derived-equiv-cor}.
 The isomorphism of triangulated categories
$\sD^\abs_{X\red}(\cB^\cu\bQcoh_\bth)=
\sD^\co_{X\red}(\cB^\cu\bQcoh_\bth)$ in the leftmost column
was mentioned in formula~\eqref{reduced-coderived=absolute-derived}.
 The triangulated equivalence $\sD^\co_{X\red}(\cB^\cu\bQcoh_\bth)
\simeq\sD^\co_{X\red}(\cB^\cu\bQcoh)$ in the leftmost column
of~\eqref{cdg-module-co-contra-construction-diagram} is the result of
Theorem~\ref{qcoh-reduced-coderived-category-equivalence-thm}.

 The triangulated equivalence
$\sD^\abs_{X\red}(\cB^\cu\bCtrh_\al^\bth)\simeq
\sD^\abs_{X\red}(\cB^\cu\bLcth_\bW^\bth)$ in the rightmost column
of~\eqref{cdg-module-co-contra-construction-diagram} is a part of
the leftmost column of
diagram~\eqref{reduced-lcta-cdg-contrader-equivs-diagram}
in Corollary~\ref{reduced-cdg-contraderived-equiv-cor}.
 The isomorphism of triangulated categories
$\sD^\abs_{X\red}(\cB^\cu\bLcth_\bW^\bth)=
\sD^\ctr_{X\red}(\cB^\cu\bLcth_\bW^\bth)$ in the rightmost column
of~\eqref{cdg-module-co-contra-construction-diagram} was mentioned in
formula~\eqref{lcth-reduced-contraderived=absolute-derived}.
 The triangulated equivalence $\sD^\ctr_{X\red}(\cB^\cu\bLcth_\bW^\bth)
\simeq\sD^\ctr_{X\red}(\cB^\cu\bLcth_\bW)$ in the rightmost column
of~\eqref{cdg-module-co-contra-construction-diagram} is the result of
Theorem~\ref{lcth-reduced-contrader-category-equivalences-thm}(a).
\end{proof}

\subsection{Co-contra correspondence via $X$-cotorsion/$X$-locally
cotorsion CDG-modules} \label{co-contra-cot-lct-subsecn}
 This section is essentially based on the reduced Koszul duality theorem
from Section~\ref{reduced-koszul-duality-contra-side-subsecn} below.
 For a more direct approach applicable under more restrictive
assumptions on the scheme $X$, see
Section~\ref{reduced-co-contrader-X-cot-X-lct-cdg-modules}.

\begin{cor} \label{thick-reduced-co-contra-cta-cot-equivs-cor}
 Let $X$ be a quasi-compact semi-separated scheme and
$(\g,\widetilde\g)$ be a quasi-coherent twisted Lie algebroid over~$X$.
 Assume that\/ $\g$~is a finite locally free sheaf on $X$, and let
$\cB^\cu=\cC^\cu_X(\g,\widetilde\g)$ be the related Chevalley--Eilenberg
quasi-coherent CDG\+quasi-algebra over~$X$.
 Then all the functors in the commutative
diagram~\eqref{thick-reduced-co-contra-X-cta-X-cot-diagram} from
Theorem~\textup{\ref{thick-reduced-co-contra-correspond-theorem}(c)}
are triangulated equivalences:
\begin{equation} \label{thick-reduced-co-contra-cta-cot-equivs-diagram}
\begin{gathered}
 \qquad\quad\xymatrix{
  \text{\llap{$\fHom^\cu_{\cB^*}(\cB^\cu,{-})\:$}}
  \sD^\abs_{X\red}(\cB^\cu\bQcoh^{X\dcta}_\bth)
  \ar@{=}[r] \ar@{-}@<2pt>[d]
  & \sD^\abs_{X\red}(\cB^\cu\bCtrh_\al^\bth)
  \text{\rlap{$\,\,:\!\cB^\cu\ocn_{\cB^*}{-}$}} \ar@{-}@<-2pt>[d] \\
  \text{\llap{$\fHom^\cu_{\cB^*}(\cB^\cu,{-})\:$}}
  \sD^\abs_{X\red}(\cB^\cu\bQcoh^{X\dcot}_\bth)
  \ar@{=}[r] \ar@<2pt>[u]
  & \sD^\abs_{X\red}(\cB^\cu\bCtrh_\al^{X\dlct,\bth})
  \text{\rlap{$\,\,:\!\cB^\cu\ocn_{\cB^*}{-}$}} \ar@<-2pt>[u]
 }
\end{gathered}
\end{equation}
 The inclusions of exact DG\+categories
$\cB^\cu\bQcoh^{X\dcot}_\bth\rarrow\cB^\cu\bQcoh^{X\dcta}_\bth
\rarrow\cB^\cu\bQcoh_\bth$ induce equivalences of the reduced
absolute derived categories  \hbadness=1100
\begin{equation} \label{qcoh-thick-cot-cta-all-reduced-abs-derived}
 \sD^\abs_{X\red}(\cB^\cu\bQcoh^{X\dcot}_\bth)\simeq
 \sD^\abs_{X\red}(\cB^\cu\bQcoh^{X\dcta}_\bth)\simeq
 \sD^\abs_{X\red}(\cB^\cu\bQcoh_\bth).
\end{equation}
\end{cor}

\begin{proof}
 The horizontal triangulated equivalences
in~\eqref{thick-reduced-co-contra-cta-cot-equivs-diagram} are provided
by Theorem~\ref{thick-reduced-co-contra-correspond-theorem}(a\+-b),
and the diagram is commutative by
Theorem~\ref{thick-reduced-co-contra-correspond-theorem}(c).
 The rightmost vertical functor is a triangulated equivalence by
Corollary~\ref{reduced-thick-lcth-al-ctrder-lcta-lct-equiv-cor}(b)
or~(c).
 It follows that the leftmost vertical functor is a triangulated
equivalence, too.
 The second assertion of the corollary follows from the first one
together with
Corollary~\ref{thick-cdg-qcoh-cta-reduced-abs-derived-equiv-cor}.
\end{proof}

\begin{cor} \label{X-cot-X-lct-cdg-module-co-contra-corresp-cor}
 Let $X$ be a quasi-compact semi-separated scheme with an open
covering\/ $\bW$ and $(\g,\widetilde\g)$ be a quasi-coherent twisted
Lie algebroid over~$X$.
 Assume that\/ $\g$~is a finite locally free sheaf on $X$, and let
$\cB^\cu=\cC^\cu_X(\g,\widetilde\g)$ be the related Chevalley--Eilenberg
quasi-coherent CDG\+quasi-algebra over~$X$.
 Then there is a natural triangulated equivalence between the reduced
coderived and the reduced contraderived category
\begin{equation} \label{X-cot-X-lct-cdg-module-co-contra-formula}
 \boR\fHom_{\cB^*}^\cu(\cB^\cu,{-})\:
 \sD^\co_{X\red}(\cB^\cu\bQcoh)\,\simeq\,
 \sD^\ctr_{X\red}(\cB^\cu\bLcth_\bW^{X\dlct})\,:\!
 \cB^\cu\ocn_{\cB^*}^\boL{-}
\end{equation}
provided by the right derived functor of contraherent\/ $\fHom$ from
$\cB^\cu$ and the left derived functor of contratensor product with
$\cB^\cu$ over~$\cB^*$.
 The desired derived functors are constructed as the compositions of
triangulated equivalences on the diagram
\begin{equation} \label{X-cot-X-lct-cdg-module-co-contra-diagram}
\begin{gathered}
 \xymatrix{
  \sD^\co_{X\red}(\cB^\cu\bQcoh) \ar@{-}@<2pt>[d]
  & \sD^\ctr_{X\red}(\cB^\cu\bLcth_\bW^{X\dlct}) \ar@{-}@<-2pt>[d] \\
  \sD^\co_{X\red}(\cB^\cu\bQcoh_\bth) \ar@<2pt>[u] \ar@{=}[d]
  & \sD^\ctr_{X\red}(\cB^\cu\bLcth_\bW^{X\dlct,\bth})
  \ar@<-2pt>[u] \ar@{=}[d] \\
  \sD^\abs_{X\red}(\cB^\cu\bQcoh_\bth) \ar@{-}@<2pt>[d]
  & \sD^\abs_{X\red}(\cB^\cu\bLcth_\bW^{X\dlct,\bth})
  \ar@{-}@<-2pt>[d] \\
  \text{\llap{$\fHom^\cu_{\cB^*}(\cB^\cu,{-})\:$}}
  \sD^\abs_{X\red}(\cB^\cu\bQcoh^{X\dcot}_\bth)
  \ar@{=}[r] \ar@<2pt>[u]
  & \sD^\abs_{X\red}(\cB^\cu\bCtrh_\al^{X\dlct,\bth})
  \text{\rlap{$\,\,:\!\cB^\cu\ocn_{\cB^*}{-}$}} \ar@<-2pt>[u]
 }
\end{gathered}
\end{equation}
\end{cor}

\begin{proof}
 The triangulated equivalence in the lower line of the diagram
is the result of
Theorem~\ref{thick-reduced-co-contra-correspond-theorem}(b).
 The triangulated equivalence
$\sD^\abs_{X\red}(\cB^\cu\bQcoh^{X\dcot}_\bth)\simeq
\sD^\abs_{X\red}(\cB^\cu\bQcoh_\bth)$ in the leftmost column is
provided by formula~\eqref{qcoh-thick-cot-cta-all-reduced-abs-derived}
in Corollary~\ref{thick-reduced-co-contra-cta-cot-equivs-cor}.
 The remaining triangulated equivalences in the leftmost column
are a part of Corollary~\ref{cdg-module-co-contra-correspondence-cor}.

 The triangulated equivalence
$\sD^\abs_{X\red}(\cB^\cu\bCtrh_\al^{X\dlct,\bth})\simeq
\sD^\abs_{X\red}(\cB^\cu\bLcth_\bW^{X\dlct,\bth})$ in the rightmost
column of~\eqref{X-cot-X-lct-cdg-module-co-contra-diagram} is a part of
the leftmost column of
diagram~\eqref{reduced-X-lct-cdg-contrader-equivs-diagram}
in Corollary~\ref{reduced-cdg-contraderived-equiv-cor}.
 The isomorphism of triangulated categories
$\sD^\abs_{X\red}(\cB^\cu\bLcth_\bW^{X\dlct,\bth})=
\sD^\ctr_{X\red}(\cB^\cu\bLcth_\bW^{X\dlct,\bth})$ in the rightmost
column of~\eqref{X-cot-X-lct-cdg-module-co-contra-diagram} was mentioned
in formula~\eqref{lcth-lct-reduced-contraderived=absolute-derived}.
 The triangulated equivalence
$\sD^\ctr_{X\red}(\cB^\cu\bLcth_\bW^{X\dlct,\bth})\allowbreak\simeq
\sD^\ctr_{X\red}(\cB^\cu\bLcth_\bW^{X\dlct})$ in the rightmost column
of~\eqref{X-cot-X-lct-cdg-module-co-contra-diagram} is the result of
Theorem~\ref{lcth-reduced-contrader-category-equivalences-thm}(b).
\hbadness=1150
\end{proof}

\begin{cor} \label{cdg-co-contra-via-lcta-and-lct-agree-cor}
 Let $X$ be a quasi-compact semi-separated scheme with an open
covering\/ $\bW$ and $(\g,\widetilde\g)$ be a quasi-coherent twisted
Lie algebroid over~$X$.
 Assume that\/ $\g$~is a finite locally free sheaf on $X$, and let
$\cB^\cu=\cC^\cu_X(\g,\widetilde\g)$ be the related Chevalley--Eilenberg
quasi-coherent CDG\+quasi-algebra over~$X$.
 Then there is a commutative diagram of triangulated equivalences
\begin{equation} \label{cdg-co-contra-via-lcta-and-lct-diagram}
\begin{gathered}
 \xymatrix{
  & \sD^\ctr_{X\red}(\cB^\cu\bLcth_\bW)
  \ar@{=}[ld] \ar@{-}@<-2pt>[dd] \\
  \sD^\co_{X\red}(\cB^\cu\bQcoh) \ar@{=}[rd] \\
  & \sD^\ctr_{X\red}(\cB^\cu\bLcth_\bW^{X\dlct}) \ar@<-2pt>[uu]
 }
\end{gathered}
\end{equation}
 Here the upper diagonal double line is the triangulated
equivalence~\eqref{cdg-module-co-contra-formula}
from Corollary~\textup{\ref{cdg-module-co-contra-correspondence-cor}},
the lower diagonal double line is the triangulated
equivalence~\eqref{X-cot-X-lct-cdg-module-co-contra-formula} from
Corollary~\textup{\ref{X-cot-X-lct-cdg-module-co-contra-corresp-cor}},
and the vertical triangulated equivalence is induced by the inclusion
of exact DG\+categories $\cB^\cu\bLcth_\bW^{X\dlct}\rarrow
\cB^\cu\bLcth_\bW$, as per
Corollary~\textup{\ref{reduced-koszul-duality-contra-side-all-equiv}}
below.
\end{cor}

\begin{proof}
 We have to check commutativity of the natural diagram of triangulated
equivalences uniting
the diagrams~\eqref{cdg-module-co-contra-construction-diagram}
and~\eqref{X-cot-X-lct-cdg-module-co-contra-diagram}.
 Indeed, the constructions in the bottom lines
of~\eqref{cdg-module-co-contra-construction-diagram}
and~\eqref{X-cot-X-lct-cdg-module-co-contra-diagram} agree with each
other by Theorem~\ref{thick-reduced-co-contra-correspond-theorem}(c)
or Corollary~\ref{thick-reduced-co-contra-cta-cot-equivs-cor}.
 The leftmost columns
of~\eqref{cdg-module-co-contra-construction-diagram}
and~\eqref{X-cot-X-lct-cdg-module-co-contra-diagram} agree with each
other in view of the commutativity of the triangular diagram of
three triangulated equivalences in
formula~\eqref{qcoh-thick-cot-cta-all-reduced-abs-derived}
from Corollary~\ref{thick-reduced-co-contra-cta-cot-equivs-cor}.

 Finally, the rightmost columns
of~\eqref{cdg-module-co-contra-construction-diagram}
and~\eqref{X-cot-X-lct-cdg-module-co-contra-diagram} agree with each
other in view of commutativity of the diagram
\begin{equation} \label{co-contra-X-lct-X-lcta-right-columns-diag}
\begin{gathered}
  \xymatrix{
  \sD^\ctr_{X\red}(\cB^\cu\bLcth_\bW^{X\dlct})
  \ar@<-2pt>[r] \ar@{-}@<-2pt>[d]
  & \sD^\ctr_{X\red}(\cB^\cu\bLcth_\bW)
  \ar@{-}@<-2pt>[l] \ar@{-}@<-2pt>[d] \\
  \sD^\ctr_{X\red}(\cB^\cu\bLcth_\bW^{X\dlct,\bth})
  \ar@<-2pt>[r] \ar@<-2pt>[u] \ar@{=}[d]
  & \sD^\ctr_{X\red}(\cB^\cu\bLcth_\bW^\bth)
  \ar@{-}@<-2pt>[l] \ar@<-2pt>[u] \ar@{=}[d] \\
  \sD^\abs_{X\red}(\cB^\cu\bLcth_\bW^{X\dlct,\bth})
  \ar@<-2pt>[r] \ar@{-}@<-2pt>[d]
  & \sD^\abs_{X\red}(\cB^\cu\bLcth_\bW^\bth)
  \ar@{-}@<-2pt>[l] \ar@{-}@<-2pt>[d] \\
  \sD^\abs_{X\red}(\cB^\cu\bCtrh_\al^{X\dlct,\bth})
  \ar@<-2pt>[r] \ar@<-2pt>[u]
  & \sD^\abs_{X\red}(\cB^\cu\bCtrh_\al^\bth)
  \ar@{-}@<-2pt>[l] \ar@<-2pt>[u]
 }
\end{gathered}
\end{equation}
where the horizontal triangulated equivalences are induced by
the inclusions of the respective exact DG\+categories.
 Here the upper square
in~\eqref{co-contra-X-lct-X-lcta-right-columns-diag} can be found in
Corollary~\ref{lcth-reduced-X-lcta-X-lct-equivs-cor}, the commutativity
of the middle square is obvious, and the lower square appears in
Corollary~\ref{reduced-thick-lcth-al-ctrder-lcta-lct-equiv-cor}(b)
(see also formulas~(\ref{al-reduced-contraderived=absolute-derived}\+-%
\ref{al-lct-reduced-contraderived=absolute-derived})).
\end{proof}

\begin{rem} \label{weakly-smooth-de-Rham-remark}
 Let $X\rarrow T$ be a weakly smooth morphism of schemes.
 Then the de~Rham quasi-coherent DG\+quasi-algebra $\Omega^\bu_{X/T}$
is the Chevalley--Eilenberg quasi-coherent DG\+quasi-algebra of
the quasi-coherent Lie algebroid of fiberwise vector fields
$\vect_{X/T}$, which is a finite locally free sheaf on~$X$.
 So we have $\Omega^\bu_{X/T}=\cC^\bu_X(\vect_{X/T})$ (see
Section~\ref{crystalline-diffoperators-subsecn}).
 Therefore, assuming that the scheme $X$ is quasi-compact and
semi-separated, the results of
Corollaries~\ref{cdg-module-co-contra-correspondence-cor}
and~\ref{X-cot-X-lct-cdg-module-co-contra-corresp-cor} are applicable
to the de~Rham quasi-coherent DG\+quasi-algebra
$\cB^\cu=\Omega^\bu_{X/T}$ on~$X$.
\end{rem}

\Section{$\cD$-$\Omega$ Duality}  \label{D-Omega-duality-secn}

 In this section, as in
Sections~\ref{cdg-reduced-contraderived-secn}\+-%
\ref{cdg-co-contra-secn}, we mostly work with a quasi-coherent twisted
Lie algebroid $(\g,\widetilde\g)$ over a scheme $X$ such that
the quasi-coherent sheaf~$\g$ on $X$ is finite locally free.
 We consider the twisted universal enveloping quasi-coherent
quasi-algebra $\cA_X(\g,\widetilde\g)$ of $(\g,\widetilde\g)$
and the Chevalley--Eilenberg quasi-coherent CDG\+quasi-algebra
$\cB^\cu=\cC^\cu_X(\g,\widetilde\g)$ of $(\g,\widetilde\g)$,
as constructed in Sections~\ref{enveloping-algebra-subsecn}
and~\ref{chevalley-eilenberg-qcoh-cdg-quasi-algebra-subsecn}.

\subsection{Tensor product of quasi-coherent CDG-modules}
\label{tensor-product-of-qcoh-cdg-bimods}
 We start with an easy module-theoretic lemma.

\begin{lem} \label{tensor-product-over-quasi-algebra-lemma}
 Let $R$ be a commutative ring, and let $B$ be a quasi-algebra over~$R$.
 Let $N$ be a right $B$\+module and $M$ be a left $B$\+module.
 Let $R\rarrow S$ be flat epimorphism of commutative rings.
 In this context: \par
\textup{(a)} Assume that $N$ is endowed with an additional left
$R$\+module structure making $N$ an $R$\+$B$\+bimodule whose underlying
$R$\+$R$\+bimodule is a quasi-module over~$R$.
 Then the ring homomorphism $R\rarrow S$ induces an isomorphism of
$S$\+modules $S\ot_R(N\ot_BM)\simeq(S\ot_RN)\ot_{S\ot_RB}(S\ot_RM)$.
\par
\textup{(b)} Assume that $M$ is endowed with an additional right
$R$\+module structure making $M$ an $B$\+$R$\+bimodule whose underlying
$R$\+$R$\+bimodule is a quasi-module over~$R$.
 Then the ring homomorphism $R\rarrow S$ induces an isomorphism
of $S$\+modules $S\ot_R(N\ot_B\nobreak M)\allowbreak\simeq
(N\ot_RS)\ot_{S\ot_RB}(S\ot_RM)$. \par
\textup{(c)} Under the combined assumptions of both parts~\textup{(a)}
and~\textup{(b)}, the tensor product $N\ot_BM$ with its induced
$R$\+$R$\+bimodule structure is a quasi-module over~$R$.
\end{lem}

\begin{proof}
 Part~(a): mostly, it just needs to be explained what the assertion
means.
 We have $S\ot_RN\simeq N\ot_RS$ by
Lemma~\ref{quasi-module-localization-lemma}(a).
 The tensor product $S\ot_RB\simeq B\ot_RS$ is an associative ring (in
fact, a quasi-algebra over~$S$) by
Lemma~\ref{quasi-algebra-co-extension-of-scalars}(a).
 The tensor product $N\ot_RS$ has a natural structure of right module
over $B\ot_RS$ by the left-right opposite version of
Lemma~\ref{quasi-algebra-co-extension-of-scalars}(b).
 One can easily see that this module structure agrees with the right
$B$\+module structure on $S\ot_RN$.
 The tensor product $S\ot_RM$ has a natural structure of left module
over $S\ot_RB$ by Lemma~\ref{quasi-algebra-co-extension-of-scalars}(b);
so we have $S\ot_RM\simeq(S\ot_RB)\ot_BM$ according to the alternative
proof of the latter lemma.
 Hence $(S\ot_RN)\ot_{S\ot_RB}(S\ot_RM)\simeq(S\ot_RN)\ot_BM\simeq
S\ot_R(N\ot_BM)$.

 Part~(b) is left-right opposite to part~(a).
 In part~(c), the $R$\+$R$\+bimodule $N\ot_BM$ is a quotient bimodule
of the $R$\+$R$\+bimodule $N\ot_RM$, which is a quasi-module over $R$
by~\cite[Proposition~2.4]{Ptd} (see
Section~\ref{prelim-quasi-modules-subsecn}).
\end{proof}

 Let $X$ be a scheme, and let $\cA$, $\cB$, and $\cC$ be three
quasi-coherent quasi-algebras over~$X$.
 Let $\N$ be a quasi-coherent $\cA$\+$\cB$\+bimodule and $\M$ be
a quasi-coherent $\cB$\+$\cC$\+bimodule on~$X$ (in the sense of
the definition in Section~\ref{fHom-over-qcoh-quasi-algebra-subsecn}).
 Then the quasi-coherent $\cA$\+$\cC$\+bimodule $\N\ot_\cB\M$ is
defined by the obvious rule
$$
 (\N\ot_\cB\M)(U)=\N(U)\ot_{\cB(U)}\M(U)
$$
for all affine open subschemes $U\subset X$.
 It is clear from
Lemma~\ref{tensor-product-over-quasi-algebra-lemma} that
the $\cA(U)$\+$\cC(U)$\+bimodules $\N(U)\ot_{\cB(U)}\M(U)$ glue together
to form a quasi-coherent $\cA$\+$\cC$\+bimodule $\N\ot_\cB\M$ on~$X$.

 Similarly, let $\N$ be a quasi-coherent $\cA$\+$\cB$\+bimodule and
$\M$ be a quasi-coherent left $\cB$\+module on~$X$ (in the sense of
the definition in Section~\ref{cosheaves-of-A-modules-subsecn}).
 Then the left $\cA(U)$\+modules $(\N\ot_\cB\M)(U)=
\N(U)\ot_{\cB(U)}\M(U)$ glue together to form a quasi-coherent left
$\cA$\+module $\N\ot_\cB\M$ on~$X$.
 This is clear from
Lemma~\ref{tensor-product-over-quasi-algebra-lemma}(a).

 Finally, let $\N$ be a quasi-coherent right $\cB$\+module and
$\M$ be a quasi-coherent $\cB$\+$\cC$\+bimodule on~$X$.
 Then the right $\cC(U)$\+modules $(\N\ot_\cB\M)(U)=
\N(U)\ot_{\cB(U)}\M(U)$ glue together to form a quasi-coherent right
$\cC$\+module $\N\ot_\cB\M$ on~$X$.
 This is clear from
Lemma~\ref{tensor-product-over-quasi-algebra-lemma}(b).

 Now let $\cA^\cu$, $\cB^\cu$, and $\cC^\cu$ be three quasi-coherent
CDG\+quasi-algebras over~$X$.
 Let $\N^\cu$ be a quasi-coherent CDG\+bimodule over $\cA^\cu$ and
$\cB^\cu$, and let $\M^\cu$ be a quasi-coherent CDG\+bimodule over
$\cB^\cu$ and~$\cC^\cu$ (in the sense of the definition in
Section~\ref{fHom-and-contratensor-of-cdg-modules}).
 Then the graded version of the construction above produces
a quasi-coherent graded $\cA^*$\+$\cC^*$\+bimodule
$\N^*\ot_{\cB^*}\M^*$ on~$X$.
 The rule
$$
 (\N^\cu\ot_{\cB^*}\M^\cu)(U)=\N^\cu(U)\ot_{\cB^*(U)}\M^\cu(U),
$$
with the reference to the definition of the tensor product of
CDG\+bimodules in Section~\ref{cdg-rings-cdg-modules-subsecn},
produces a differential on the graded $\cA^*(U)$\+$\cC^*(U)$\+bimodule
$(\N^*\ot_{\cB^*}\M^*)(U)$ for every affine open subscheme $U\subset X$.
 The collection of all such differentials defines a structure of
quasi-coherent CDG\+bimodule over $\cA^\cu$ and $\cC^\cu$ on the tensor
product of quasi-coherent graded bimodules $\M^*\ot_{\cB^*}\M^*$.
 We denote the resulting quasi-coherent CDG\+bimodule over $\cA^\cu$
and $\cC^\cu$ by $\N^\cu\ot_{\cB^*}\M^\cu$.

 Similarly, let $\N^\cu$ be a quasi-coherent CDG\+bimodule over
$\cA^\cu$ and $\cB^\cu$, and let $\M^\cu$ be a quasi-coherent left
CDG\+module over $\cB^\cu$ (in the sense of the definition in
Section~\ref{qcoh-lcth-cdg-modules-subsecn}).
 Then the graded version of the construction above produces
a quasi-coherent graded left $\cA^*$\+module $\N^*\ot_{\cB^*}\M^*$
on~$X$.
 The same rule
$$
 (\N^\cu\ot_{\cB^*}\M^\cu)(U)=\N^\cu(U)\ot_{\cB^*(U)}\M^\cu(U)
$$
endows the graded left $\cA^*(U)$\+module $(\N^*\ot_{\cB^*}\M^*)(U)$
with a differential.
 The collection of such differentials for all affine open subschemes
$U\subset X$ defines a structure of quasi-coherent left CDG\+module over
$\cA^\cu$ on the tensor product of quasi-coherent graded (bi)modules
$\M^*\ot_{\cB^*}\M^*$.
 We denote the resulting quasi-coherent left CDG\+module over $\cA^\cu$
by $\N^\cu\ot_{\cB^*}\M^\cu$.

 Finally, let $\N^\cu$ be a quasi-coherent right CDG\+module over
$\cB^\cu$, and let $\M^\cu$ be a quasi-coherent CDG\+bimodule over
$\cB^\cu$ and~$\cC^\cu$.
 Then the graded version of the construction above produces
a quasi-coherent graded right $\cC^*$\+module $\N^*\ot_{\cB^*}\M^*$
on~$X$.
 The very same rule
$$
 (\N^\cu\ot_{\cB^*}\M^\cu)(U)=\N^\cu(U)\ot_{\cB^*(U)}\M^\cu(U)
$$
endows the graded right $\cC^*(U)$\+module $(\N^*\ot_{\cB^*}\M^*)(U)$
with a differential.
 The collection of such differentials for all affine open subschemes
$U\subset X$ defines a structure of quasi-coherent right CDG\+module
over $\cC^\cu$ on the tensor product of quasi-coherent graded
(bi)modules $\M^*\ot_{\cB^*}\M^*$.
 We denote the resulting quasi-coherent right CDG\+module over $\cC^\cu$
by $\N^\cu\ot_{\cB^*}\M^\cu$.

\subsection{$\Cohom$ into contraherent CDG-module}
\label{Cohom-into-lcth-cdg-module}
 Once again, we start with a couple of module-theoretic lemmas.

\begin{lem} \label{Hom-over-quasi-algebra-lemma}
 Let $R$ be a commutative ring and $B$ be a quasi-algebra over~$R$.
 Let $E$ be a $B$\+$R$\+bimodule whose underlying $R$\+$R$\+bimodule is
a quasi-module over~$R$, and let $P$ be a left $B$\+module.
 Let $R\rarrow S$ be a flat epimorphism of commutative rings.
 Then the ring homomorphism $R\rarrow S$ induces an isomorphism of
$S$\+modules\/ $\Hom_{S\ot_RB}(S\ot_RE,\>\Hom_R(S,P))\simeq
\Hom_R(S,\Hom_B(E,P))$.
\end{lem}

\begin{proof}
 This is dual-analogous to
Lemma~\ref{tensor-product-over-quasi-algebra-lemma}(a\+-b).
 We have $S\ot_RE\simeq E\ot_RS$ by
Lemma~\ref{quasi-module-localization-lemma}(a).
 The tensor product $S\ot_RE$ and the Hom module $\Hom_R(S,P)$ have
natural structures of left modules over $S\ot_RB$ by
Lemma~\ref{quasi-algebra-co-extension-of-scalars}(b\+-c); furthermore,
we have $\Hom_R(S,P)\simeq\Hom_B(B\ot_RS,\>P)$ according to
the alternative proof of 
Lemma~\ref{quasi-algebra-co-extension-of-scalars}(c).
 One can easily see that the left $(S\ot_RB)$\+module structure on
$S\ot_RE$ agrees with the left $B$\+module structure on $E\ot_RS$.
 Hence $\Hom_{B\ot_RS}(E\ot_R\nobreak S,\>\Hom_R(S,P))\simeq
\Hom_B(E\ot_RS,\>P)\simeq\Hom_R(S,\Hom_B(E,P))$.
\end{proof}

 The next lemma is a simpler affine version of
Lemma~\ref{fHom-contraadjusted-and-more-lemma}.

\begin{lem} \label{module-Hom-contraadjusted-and-more-lemma}
 Let $R$ be a commutative ring, let $A$ and $B$ be quasi-coherent
quasi-algebras over $R$, and let $E$ be an $A$\+$B$ bimodule whose
underlying $R$\+$R$\+bimodule is a quasi-module over~$R$.
 Let $P$ and $J$ be left $A$\+modules.
 In this context: \par
\textup{(a)} if $E$ is $A/R$\+very flaprojective as a left $A$\+module
and $P$ is contraadjusted as an $R$\+module, then\/ $\Hom_A(E,P)$ is
contraadjusted as an $R$\+module; \par
\textup{(b)} if $E$ is $(A/R,R)$\+robustly flaprojective as
an $A$\+$R$\+bimodule and $P$ is cotorsion as an $R$\+module, then\/
$\Hom_A(E,P)$ is cotorsion as an $R$\+module; \par
\textup{(c)} if $E$ is flat as a left $A$\+module and $P$ is
a cotorsion left $A$\+module, then\/ $\Hom_A(E,P)$ is a cotorsion
left $B$\+module; \par
\textup{(d)} if $E$ is flat as a right $B$\+module and $J$ is
an injective left $A$\+module, then\/ $\Hom_A(E,J)$ is an injective
left $B$\+module; \par
\textup{(e)} if $J$ is an injective left $A$\+module, then\/
$\Hom_A(E,J)$ is a cotorsion left $B$\+module.
\end{lem}

\begin{proof}
 Part~(a) is a restatement of
Lemma~\ref{very-flaprojective-Hom-tensor-lemma}(b).
 Part~(b) is a particular case of
Lemma~\ref{robustly-flaprojective-Hom-tensor-lemma}(b).
 Parts~(c\+-e) are a restatement of Lemma~\ref{hom-injective-cotorsion}.
\end{proof}

 Let $X$ be a scheme with an open covering~$\bW$.
 Denote by $\bB$ the base of open subsets in $X$ consisting of all
the affine open subschemes $U\subset X$ subordinate to~$\bW$.

 Let $\cA$ and $\cB$ be quasi-coherent quasi-algebras over~$X$.
 Let $\cE$ be a quasi-coherent $\cA$\+$\cB$\+bimodule and $\P$ be
a $\bW$\+locally contraherent $\cA$\+module on~$X$.
 Similarly to Section~\ref{cohom-from-quasi-modules-subsecn},
the copresheaf of $\cB$\+modules $\Cohom_\cA(\cE,\P)$ on $\bB$ is
defined by the rule
$$
 \Cohom_\cA(\cE,\P)[U]=\Hom_{\cA(U)}(\cE(U),\P[U])
$$
for all affine open subschemes $U\subset X$ subordinate to~$\bW$.

\begin{cor} \label{Cohom-contraadjusted-and-more-cor}
\textup{(a)} If a quasi-coherent $\cA$\+$\cB$\+bimodule $\cE$ on $X$
is $\cA$\+very flaprojective (in the sense of
Section~\ref{fHom-over-qcoh-quasi-algebra-subsecn}), and\/ $\P$ is
a\/ $\bW$\+locally contraherent $\cA$\+module, then
the $\cO_X(U)$\+module\/ $\Cohom_\cA(\cE,\P)[U]$ is contraadjusted; \par
\textup{(b)} If a quasi-coherent $\cA$\+$\cB$\+bimodule $\cE$ on $X$
is $\cA$\+robustly flaprojective (in the sense of
Section~\ref{fHom-over-qcoh-quasi-algebra-subsecn}), and\/ $\P$ is
an $X$\+locally cotorsion\/ $\bW$\+locally contraherent $\cA$\+module,
then the $\cO_X(U)$\+module\/ $\Cohom_\cA(\cE,\P)[U]$ is cotorsion; \par
\textup{(c)} If a quasi-coherent $\cA$\+$\cB$\+bimodule $\cE$ on $X$
is $\cA$\+flat and\/ $\P$ is an $\cA$\+locally cotorsion\/
$\bW$\+locally contraherent $\cA$\+module, then the $\cB(U)$\+module\/
$\Cohom_\cA(\cE,\P)[U]$ is cotorsion; \par
\textup{(d)} If a quasi-coherent $\cA$\+$\cB$\+bimodule $\cE$ on $X$
is $\cB$\+flat and\/ $\gJ$ is an $\cA$\+locally injective\/
$\bW$\+locally contraherent $\cA$\+module, then the $\cB(U)$\+module\/
$\Cohom_\cA(\cE,\gJ)[U]$ is injective; \par
\textup{(e)} If $\cE$ is a quasi-coherent $\cA$\+$\cB$\+bimodule and\/
$\gJ$ is an $\cA$\+locally injective\/ $\bW$\+locally contraherent
$\cA$\+module on $X$, then the $\cB(U)$\+module\/
$\Cohom_\cA(\cE,\gJ)[U]$ is cotorsion.
\end{cor}

\begin{proof}
 All the assertions follow immediately from the definitions and
Lemma~\ref{module-Hom-contraadjusted-and-more-lemma}.
\end{proof}

 Under any one of the assumptions~(a\+-e) of
Corollary~\ref{Cohom-contraadjusted-and-more-cor}, the corollary
implies, in particular, that the copresheaf $\Cohom_\cA(\cE,\P)$ or
$\Cohom_\cA(\cE,\gJ)$ on $\bB$, viewed as a copresheaf of
$\cO_X$\+modules, satisfies the contraadjustedness axiom~(ii)
from Section~\ref{locally-contraherent-cosheaves-subsecn}.
 The contraherence axiom~(i) follows from
Lemma~\ref{Hom-over-quasi-algebra-lemma}.
 By Lemma~\ref{contraherence+contraadjustedness-imply-cosheaf}
and Theorem~\ref{extension-of-co-sheaves-from-topology-base}(b),
the copresheaf $\Cohom_\cA(\cE,\P)$ or $\Cohom_\cA(\cE,\gJ)$ on $\bB$
extends uniquely to a $\bW$\+locally contraherent cosheaf of
$\cB$\+modules on $X$, which we will denote also by
$\Cohom_\cA(\cE,\P)$ or $\Cohom_\cA(\cE,\gJ)$.

 To summarize:
\begin{itemize}
\item $\Cohom_\cA(\cE,\P)$ is a $\bW$\+locally contraherent
$\cB$\+module whenever $\cE$ is $\cA$\+very flaprojective and $\P$
is $\bW$\+locally contraherent;
\item $\Cohom_\cA(\cE,\P)$ is an $X$\+locally cotorsion $\bW$\+locally
contraherent $\cB$\+module whenever $\cE$ is $\cA$\+robustly
flaprojective and $\P$ is $X$\+locally cotorsion $\bW$\+locally
contraherent;
\item $\Cohom_\cA(\cE,\P)$ is a $\cB$\+locally cotorsion $\bW$\+locally
contraherent $\cB$\+module whenever $\cE$ is $\cA$\+flat and $\P$ is
$\cA$\+locally cotorsion $\bW$\+locally contraherent;
\item $\Cohom_\cA(\cE,\gJ)$ is a $\cB$\+locally injective $\bW$\+locally
contraherent $\cB$\+module whenever $\cE$ is $\cB$\+flat and $\gJ$ is
locally injective $\bW$\+locally contraherent;
\item $\Cohom_\cA(\cE,\gJ)$ is a $\cB$\+locally cotorsion $\bW$\+locally
contraherent $\cB$\+module whenever  $\gJ$ is locally injective
$\bW$\+locally contraherent.
\end{itemize}

 Now let $\cA^\cu$ and $\cB^\cu$ be two quasi-coherent
CDG\+quasi-algebras over~$X$.
 Let $\cE^\cu$ be a quasi-coherent CDG\+bimodule over $\cA^\cu$ and
$\cB^\cu$ (in the sense of
Section~\ref{fHom-and-contratensor-of-cdg-modules}) and $\P^\cu$ be
a $\bW$\+locally contraherent CDG\+module over $\cA^\cu$ (in the sense
of Section~\ref{qcoh-lcth-cdg-modules-subsecn}).
 Then, under the graded version of any one of the assumptions~(a\+-e)
of Corollary~\ref{Cohom-contraadjusted-and-more-cor}, the graded
version of the construction above produces a $\bW$\+locally
contraherent graded module $\Cohom_{\cA^*}^*(\cE^*,\P^*)$ or
over the quasi-coherent graded quasi-algebra $\cB^*$ over~$X$.
 Depending on which of the assumptions~(a\+-e) is satisfied,
the $\bW$\+locally contraherent graded module
$\Cohom_{\cA^*}^*(\cE^*,\P^*)$ over $\cB^*$ has further adjustedness
properties as per the discussion above.

 The rule
$$
 \Cohom^\cu_{\cA^*}(\cE^\cu,\P^\cu)[U]=
 \Hom_{\cA^*(U)}^\cu(\cE^\cu(U),\P^\cu[U]),
$$
with the reference to the definition of the Hom CDG\+module in
Section~\ref{cdg-rings-cdg-modules-subsecn},
produces a differential on the graded $\cA^*(U)$\+module
$\Cohom_{\cA^*}^*(\cE^*,\P^*)[U]$ for every affine open subscheme
$U\subset X$ subordinate to~$\bW$.
 The collection of all such differentials defines a structure of
$\bW$\+locally contraherent CDG\+module over $\cB^\cu$ on
the $\bW$\+locally contraherent graded module
$\Cohom_{\cA^*}^*(\cE^*,\P^*)$ over~$\cB^*$.
 We denote the resulting $\bW$\+locally contraherent CDG\+module
over $\cB^\cu$ by $\Cohom_{\cA^*}^\cu(\cE^\cu,\P^\cu)$.

\subsection{Koszul duality functors on the co side}
\label{koszul-duality-functors-co-side}
 Let $X$ be a scheme and $(\g,\widetilde\g)$ be a quasi-coherent
twisted Lie algebroid over $X$ such that the quasi-coherent sheaf~$\g$
on $X$ is finite locally free.
 Let $\cA=\cA_X(\g,\widetilde\g)$ be the twisted universal enveloping
quasi-coherent quasi-algebra and $\cB^\cu=\cC^\cu_X(\g,\widetilde\g)$
be the Chevalley--Eilenberg quasi-coherent CDG\+quasi-algebra
of $(\g,\widetilde\g)$.

 We recall that the notation $\bCom(\Qcohr\cA)$ stands for
the abelian DG\+category of complexes of quasi-coherent right
$\cA$\+modules on $X$, while $\bQcohr\cB^\cu$ is the abelian
DG\+category of quasi-coherent right CDG\+modules over~$\cB^\cu$.

\begin{lem} \label{koszul-duality-dg-functors-right-co-side}
 There is a natural pair of adjoint exact DG\+functors
$$
 \xymatrix{
  {-}\ot_{\cO_X}\cHom_{\cO_X}(\cB^*,\cO_X)\:
  \bCom(\Qcohr\cA) \ar@<-2pt>[r]
  & \bQcohr\cB^\cu \,:{-}\ot_{\cO_X}\cA, \ar@<-2pt>[l]
 }
$$
where the DG\+functor\/ ${-}\ot_{\cO_X}\cHom_{\cO_X}(\cB^*,\cO_X)$ is
the right adjoint and the DG\+functor\/ ${-}\ot_{\cO_X}\cA$ is
the left adjoint.
 To be more precise, the differentials on the respective tensor
products over $\cO_X$ are defined in terms of the identifications
$$
 {-}\ot_{\cO_X}\cHom_{\cO_X}(\cB^*,\cO_X) =
 \,{-}\ot_\cA\cC_\cu^X(\cA_X,\g,\widetilde\g)
 \:\bCom(\Qcohr\cA)\lrarrow\bQcohr\cB^\cu
$$
and
$$
 {-}\ot_{\cO_X}\cA = \,{-}\ot_{\cB^*}\cC^\cu_X(\g,\widetilde\g,\cA_X)
 \:\bQcohr\cB^\cu\lrarrow\bCom(\Qcohr\cA).
$$
 Here $\cC^\cu_X(\g,\widetilde\g,\cA_X)$ and
$\cC_\cu^X(\cA_X,\g,\widetilde\g)$ are the Chevalley--Eilenberg
quasi-coherent CDG\+(bi)modules constructed in
Section~\ref{chevalley-eilenberg-cdg-modules-subsecn}.
\end{lem}

\begin{proof}
 We recall that $\cC^\cu_X(\g,\widetilde\g,\cA_X)$ is a quasi-coherent
CDG\+bimodule over $\cB^\cu$ and $\cA$ with the underlying
quasi-coherent graded bimodule
$$
 \cB^*\ot_{\cO_X}\cA\simeq
 \bigwedge\nolimits_X^*\bigl(\cHom_{\cO_X}(\g,\cO_X)\bigr)
 \ot_{\cO_X}\cA,
$$
while $\cC_\cu^X(\cA_X,\g,\widetilde\g)$ is a quasi-coherent
CDG\+bimodule over $\cA$ and $\cB^\cu$ with the underlying
quasi-coherent graded bimodule
$$
 \cA\ot_{\cO_X}\cHom_{\cO_X}(\cB^*,\cO_X)\simeq
 \cA\ot_{\cO_X}\bigwedge\nolimits_X^*(\g).
$$
 The construction of the tensor product of quasi-coherent
CDG\+(bi)modules from Section~\ref{tensor-product-of-qcoh-cdg-bimods}
defines the quasi-coherent CDG\+module structures on the tensor
products with these quasi-coherent CDG\+bimodules, providing
the desired DG\+functors in both directions.
 So the notation ${-}\ot_{\cO_X}\bigwedge^*_X(\g)$ can be used as
a shorthand for ${-}\ot_{\cO_X}\cHom_{\cO_X}(\cB^*,\cO_X)$;
but one has to keep in mind that the cohomological grading on
$\bigwedge^*_X(\g)=\cHom_{\cO_X}(\cB^*,\cO_X)$ is nonpositive.

 Now let $\M^\bu$ be a complex of quasi-coherent right $\cA$\+modules
and $\N^\cu$ be a quasi-coherent right CDG\+module over~$\cB^\cu$.
 In order to establish the adjunction of the two DG\+functors, one
needs to compute that the two complexes of abelian groups
$$ \textstyle
 \Hom^\bu_{\cB^\rop{}^*}(\N^\cu,\>\M^\bu\ot_{\cO_X}\bigwedge^*_X(\g))
 \quad\text{and}\quad
 \Hom^\bu_{\cA^\rop}(\N^\cu\ot_{\cO_X}\cA,\>\M^\bu)
$$
have one and the same (up to natural isomorphism) underlying graded
abelian group $\Hom_X^*(\N^*,\M^*)$, and also one and the same
natural differential.
 The differential has three summands: one coming from the differential
on $\M^\bu$, the second one coming from the differential on $\N^\cu$,
and the third one determined by the action of $\cA$ on $\M^*$
together with the action of $\cB^*$ on~$\N^*$.
 The latter one is a kind of Koszul differential.

 Finally, both the DG\+functors are exact since the graded
quasi-coherent sheaf $\cHom_{\cO_X}(\cB^*,\cO_X)$ on $X$ is finite
locally free (hence flat), while the quasi-coherent quasi-module $\cA$
is locally projective (hence also flat as a quasi-coherent sheaf) in
its left (as well as right) $\cO_X$\+module structure.
 See Corollary~\ref{twisted-lie-algebroid-pbw-qcoh-cor}(b).
\end{proof}

\subsection{Koszul duality functors on the contra side}
\label{koszul-duality-functors-contra-side}
 We keep the notation of Section~\ref{koszul-duality-functors-co-side}.
 Let $\bW$ be an open covering of~$X$.
 
 Recall that the notation $\bCom(\cA\Lcth_\bW)$ stands for the exact
DG\+category of complexes of $\bW$\+locally contraherent $\cA$\+modules
on $X$, while $\cB^\cu\bLcth_\bW$ is the exact DG\+category of
$\bW$\+locally contraherent CDG\+modules over~$\cB^\cu$.
 Furthermore, $\bCom(\cA\Lcth_\bW^{\cA\dlct})\subset
\bCom(\cA\Lcth_\bW^{X\dlct})\subset\bCom(\cA\Lcth_\bW)$ are the full
exact DG\+subcategories of complexes of $\cA$\+locally cotorsion and
$X$\+locally cotorsion $\bW$\+locally contraherent $\cA$\+modules,
while $\cB^\cu\bLcth_\bW^{X\dlct}\subset\cB^\cu\bLcth_\bW$ is the full
exact DG\+subcategory of $X$\+locally cotorsion $\bW$\+locally
contraherent CDG\+modules over~$\cB^\cu$.
 In view of Corollary~\ref{lct-lcth-modules-over-flfrqa-cor}, we do not
need to consider $\cB^*$\+locally cotorsion locally contraherent
CDG\+modules over~$\cB^\cu$.

\begin{lem} \label{koszul-duality-dg-functors-contra-side}
\textup{(a)} There is a natural pair of adjoint exact DG\+functors
$$
 \xymatrix{
  \Cohom_X(\cHom_{\cO_X}(\cB^*,\cO_X),{-})\:
  \bCom(\cA\Lcth_\bW) \ar@<2pt>[r]
  & \cB^\cu\bLcth_\bW \,:\!\Cohom_X(\cA,{-}), \ar@<2pt>[l]
 }
$$
where the DG\+functor\/ $\Cohom_X(\cHom_{\cO_X}(\cB^*,\cO_X),{-})$ is
the left adjoint and the DG\+functor\/ $\Cohom_X(\cA,{-})$ is
the right adjoint. \par
\textup{(b)} There is a natural pair of adjoint exact DG\+functors
\begin{multline*}
 \Cohom_X(\cHom_{\cO_X}(\cB^*,\cO_X),{-})\:
 \bCom(\cA\Lcth_\bW^{X\dlct}) \\
 \xymatrix{
  {} \ar@<2pt>[r]
  & \cB^\cu\bLcth_\bW^{X\dlct} \,:\!\Cohom_X(\cA,{-}), \ar@<2pt>[l]
 }
\end{multline*}
where the DG\+functor\/ $\Cohom_X(\cHom_{\cO_X}(\cB^*,\cO_X),{-})$ is
the left adjoint and the DG\+functor\/ $\Cohom_X(\cA,{-})$ is
the right adjoint. \par
\textup{(c)} There is a natural pair of adjoint exact DG\+functors
\begin{multline*}
 \Cohom_X(\cHom_{\cO_X}(\cB^*,\cO_X),{-})\:
 \bCom(\cA\Lcth_\bW^{\cA\dlct}) \\
 \xymatrix{
  {} \ar@<2pt>[r]
  & \cB^\cu\bLcth_\bW^{X\dlct} \,:\!\Cohom_X(\cA,{-}), \ar@<2pt>[l]
 }
\end{multline*}
where the DG\+functor\/ $\Cohom_X(\cHom_{\cO_X}(\cB^*,\cO_X),{-})$ is
the left adjoint and the DG\+functor\/ $\Cohom_X(\cA,{-})$ is
the right adjoint.

 To be more precise, in all the three cases~(a\+-c), the differentials
on the respective\/ $\bW$\+locally contraherent modules\/ $\Cohom_X$
are defined in terms of the identifications
\begin{multline*}
 \Cohom_X(\cHom_{\cO_X}(\cB^*,\cO_X),{-}) =
 \Cohom_\cA^\cu(\cC_\cu^X(\cA_X,\g,\widetilde\g),{-})\: \\
 \bCom(\cA\Lcth_\bW)\lrarrow\cB^\cu\bLcth_\bW
\end{multline*}
and
$$
 \Cohom_X(\cA,{-}) =
 \Cohom_{\cB^*}^\cu(\cC^\cu_X(\g,\widetilde\g,\cA_X),{-})
 \:\cB^\cu\bLcth_\bW\lrarrow\bCom(\cA\Lcth_\bW).
$$
 Here $\cC^\cu_X(\g,\widetilde\g,\cA_X)$ and
$\cC_\cu^X(\cA_X,\g,\widetilde\g)$ are the Chevalley--Eilenberg
quasi-coherent CDG\+(bi)modules constructed in
Section~\ref{chevalley-eilenberg-cdg-modules-subsecn}.
\end{lem}

\begin{proof}
 Were refer to the proof of
Lemma~\ref{koszul-duality-dg-functors-right-co-side} for
the discussion of the CDG\+bimodules $\cC^\cu_X(\g,\widetilde\g,\cA_X)$
and $\cC_\cu^X(\cA_X,\g,\widetilde\g)$, and only observe that
the notation $\Cohom_X(\bigwedge_X^*(\g),{-})$ can be used as
a shorthand for $\Cohom_X(\cHom_{\cO_X}(\cB^*,\cO_X),{-})$.
 The construction of the $\bW$\+locally contraherent CDG\+modules
$\Cohom^\cu$ from Section~\ref{Cohom-into-lcth-cdg-module} provides
the desired DG\+functors in both directions.  {\hfuzz=1.9pt\par}

 On the level of the underlying (graded) $\bW$\+locally contraherent
cosheaves, the construction making sense of the functor
$\Cohom_X(\bigwedge_X^*(\g),{-})$ was spelled out already
in~\cite[Sections~2.4 and~3.6]{Pcosh}, while to define the functor
$\Cohom_X(\cA,{-})$ one needs to use the construction from
Section~\ref{cohom-from-quasi-modules-subsecn} above.
 The point is that $\bigwedge_X^*(\g)$ is a graded quasi-coherent sheaf
on~$X$ (on which the left and right actions of $\cO_X$ coincide),
while $\cA$ is a quasi-coherent quasi-module/quasi-algebra over~$X$
(so the left and right $\cO_X$\+module structures on $\cA$ are
different).

 Let us explain why the relevant $\bW$\+locally contraherent
cosheaves/modules $\Cohom$ are well-defined.
 Concerning $\Cohom_X$, the point is that the graded
quasi-coherent sheaf $\cHom_{\cO_X}(\cB^*,\cO_X)$ on $X$ is finite
locally free (hence very flat), while the quasi-coherent quasi-module
$\cA$ is locally projective (hence also very flat as a quasi-coherent
sheaf) in its left (as well as right) $\cO_X$\+module structure.
 See Corollary~\ref{twisted-lie-algebroid-pbw-qcoh-cor}(b).

 Concerning $\Cohom_\cA^*$, for the purposes of part~(a) we need to
observe that the graded quasi-coherent $\cA$\+module
$\cC_*^X(\cA_X,\g,\widetilde\g)$ is very flaprojective.
 For the purposes of parts~(b\+-c), we need to point out that
the graded quasi-coherent $\cA$\+$\cO_X$\+bimodule
$\cC_*^X(\cA_X,\g,\widetilde\g)=\cA\ot_{\cO_X}\bigwedge_X^*(\g)$
is robustly flaprojective.

 Concerning $\Cohom_{\cB^*}^*$, for the purposes of part~(a) we need to
observe that the quasi-coherent graded $\cB^*$\+module
$\cC^*_X(\g,\widetilde\g,\cA_X)$ is very flaprojective.
 For the purposes of part~(b), we need to point out that
the quasi-coherent graded $\cB^*$\+$\cO_X$\+bimodule
$\cC^*_X(\g,\widetilde\g,\cA_X)=\cB^*\ot_{\cO_X}\cA$ is robustly
flaprojective.

 Finally, to obtain the DG\+functor in part~(c), we need to explain why
the DG\+functor pointing rightwards in part~(b) actually produces
complexes of $\cA$\+lo\-cally cotorsion (rather than merely $X$\+locally
cotorsion) $\bW$\+locally contraherent $\cA$\+modules.
 For this purpose, it suffices to interpret this functor as
$\Cohom_X(\cA,{-})$ and notice that $\cA$ is a flat quasi-coherent
left $\cO_X$\+module (so
Corollary~\ref{Cohom-contraadjusted-and-more-cor}(c) is applicable).
 Alternatively, one can use
Corollary~\ref{lct-lcth-modules-over-flfrqa-cor} and point out that
$\cC^*_X(\g,\widetilde\g,\cA_X)=\cB^*\ot_{\cO_X}\cA$ is a flat
quasi-coherent graded left $\cB^*$\+module.

 Now let $\P^\bu$ be a complex of $\bW$\+locally contraherent
$\cA$\+modules and $\Q^\cu$ be a $\bW$\+locally contraherent
CDG\+module over~$\cB^\cu$.
 In order to establish the adjunction of the two DG\+functors in
part~(a), one needs to compute that the two complexes of abelian groups
$$ \textstyle
 \Hom^{\cB^*,\bu}(\Cohom_X(\bigwedge_X^*(\g),\P^\bu),\Q^\cu)
 \quad \text{and} \quad
 \Hom^{\cA,\bu}(\P^\bu,\Cohom_X(\cA,\Q^\cu))
$$
have one and the same (up to natural isomorphism) underlying graded
abelian group $\Hom^{X,*}(\P^*,\Q^*)$, and also one and the same
natural differential.
 The differential has three summands: one coming from the differential
on $\P^\bu$, the second one coming from the differential on $\Q^\cu$,
and the third one determined by the action of $\cA$ on $\P^*$ together
with the action of $\cB^*$ on~$\Q^*$.
 The latter one is a kind of Koszul differential.
 The adjunctions of DG\+functors in parts~(b\+-c) follow.

 All the DG\+functors in parts~(a\+-c) are exact because the functor
$\Cohom$ is exact (in both arguments) wherever it is well-defined as
per any one of our sufficient conditions in
Sections~\ref{cohom-from-quasi-modules-subsecn}
and~\ref{Cohom-into-lcth-cdg-module}.
\end{proof}

\subsection{Semiderived Koszul duality on the co side}
\label{semiderived-koszul-duality-co-side-subsecn}
 We refer to Section~\ref{semiderived-quasi-coherent-subsecn}
for the definitions of the Positselski semi(co)derived category
$\sD^\si(\Qcohr\cA)$ and the Becker semi(co)derived category
$\sD^\bsi(\Qcohr\cA)$ of quasi-coherent right $\cA$\+modules.
 The following theorem is (essentially) a generalization
of~\cite[Theorem~B.2(a)]{Pkoszul} to singular/non-Noetherian schemes.
 It is also a nonaffine version of~\cite[Theorem~6.14]{Prel}.

\begin{thm} \label{semiderived-koszul-duality-right-co-side}
 Let $X$ be a scheme and $(\g,\widetilde\g)$ be a quasi-coherent
twisted Lie algebroid over~$X$.
 Assume that\/ $\g$~is a finite locally free sheaf of bounded rank
on~$X$.
 Let $\cA=\cA_X(\g,\widetilde\g)$ be the enveloping quasi-coherent
quasi-algebra and $\cB^\cu=\cC^\cu_X(\g,\widetilde\g)$ be
the Chevalley--Eilenberg quasi-coherent CDG\+quasi-algebra
of $(\g,\widetilde\g)$; so $\cB^n=0$ for $n$~large enough.
 In this setting: \par
\textup{(a)} the pair of adjoint DG\+functors from
Lemma~\ref{koszul-duality-dg-functors-right-co-side} induces
a triangulated equivalence between the Positselski semicoderived
and the Positselski coderived category
$$
 \xymatrix{
  {-}\ot_{\cO_X}\bigwedge\nolimits_X^*(\g)\:
  \sD^\si(\Qcohr\cA) \ar@{=}[r]
  & \sD^\co(\bQcohr\cB^\cu) \,:{-}\ot_{\cO_X}\cA;
 }
$$ \par
\textup{(b)} the pair of adjoint DG\+functors from
Lemma~\ref{koszul-duality-dg-functors-right-co-side} induces
a triangulated equivalence between the Becker semicoderived
and the Becker coderived category
$$
 \xymatrix{
  {-}\ot_{\cO_X}\bigwedge\nolimits_X^*(\g)\:
  \sD^\bsi(\Qcohr\cA) \ar@{=}[r]
  & \sD^\bco(\bQcohr\cB^\cu) \,:{-}\ot_{\cO_X}\cA.
 }
$$
\end{thm}

\begin{proof}
 Four assertions need to be proved in each part~(a) and~(b).
 Firstly, one has to show that the functor
${-}\ot_{\cO_X}\bigwedge_X^*(\g)$ takes semicoacyclic complexes to
coacyclic quasi-coherent CDG\+modules, and that the functor
${-}\ot_{\cO_X}\cA$ takes coacyclic quasi-coherent CDG\+modules to
semicoacyclic complexes.
 Secondly, it needs to be checked that the cone of the adjunction morphism $\M^\bu\ot_{\cO_X}\bigwedge_X^*(\g)\ot_{\cO_X}\cA\rarrow\M^\bu$
is semicoacyclic for every complex $\M^\bu\in\bCom(\Qcohr\cA)$,
and that the cone of the adjunction morphism $\N^\cu\rarrow
\N^\cu\ot_{\cO_X}\cA\ot_{\cO_X}\bigwedge_X^*(\g)$ is coacyclic
for every quasi-coherent CDG\+module $\N^\cu\in\bQcohr\cB^\cu$.

 Let $\M^\bu$ be a complex in the abelian category $\Qcohr\cA$.
 Then the quasi-coherent CDG\+module $\M^\bu\ot_{\cO_X}
\bigwedge_X^*(\g)$ has a natural finite decreasing filtration $F$ by
quasi-coherent CDG\+submodules.
 The filtration $F$ on $\M^\bu\ot_{\cO_X}\bigwedge_X^*(\g)$ is
induced by the filtration $F$ on $\bigwedge_X^*(\g)$ given by the rule
$F^i\bigl(\bigwedge_X^*(\g)\bigr)=\bigoplus_{j\le -i}\bigwedge_X^j(\g)$
for all $i\in\boZ$.
 In fact, the quasi-coherent CDG\+module $\M^\bu\ot_{\cO_X}
\bigwedge_X^*(\g)$ over $\cB^\cu$ is thick in the sense of
Sections~\ref{thick-quasi-coherent-modules-subsecn}
and~\ref{exact-dg-categories-of-thick-qcoh-subsecn}, and
the filtration $F$ defined above is just the filtration $F$
from Corollary~\ref{thick-qcoh-sheaves-are-filtered}(b); but we do not
need to use these observations at this point.

 It is only important that the successive quotients
$$
 F^i\Bigl(\M^\bu\ot_{\cO_X}\bigwedge\nolimits_X^*(\g)\Bigr)\Big/
 F^{i+1}\Bigl(\M^\bu\ot_{\cO_X}\bigwedge\nolimits_X^*(\g)\Bigr)
$$
are trivial quasi-coherent CDG\+modules over~$\cB^\cu$ (in the sense
of Section~\ref{koszul-coresolutions-of-trivial-qcoh-cdg-mods})
corresponding to the complexes of quasi-coherent sheaves
$\M^\bu\ot_{\cO_X}\bigwedge^{-i}_X(\g)$.
 Now $\M^\bu$ being a semicoacyclic complex of quasi-coherent
$\cA$\+modules means precisely that $\M^\bu$ is coacyclic as a complex
of quasi-coherent sheaves on~$X$.
 If this is the case, then the complexes of quasi-coherent sheaves
$\M^\bu\ot_{\cO_X}\bigwedge^{-i}_X(\g)$ on $X$ are coacyclic as well.
 In the context of the Positselski coacyclicity, this holds because
the functor ${-}\ot_{\cO_X}\bigwedge^{-i}_X(\g)\:X\Qcoh\rarrow X\Qcoh$
is exact and preserves infinite direct sums.
 In the context of the Becker coacyclicity, one can refer
to~\cite[Lemma~A.5]{Psemten}, or~\cite[Lemma~B.7.5(a)]{Pcosh}, or any
one of Lemmas~\ref{DG-functors-preserve-Becker-co-contra-acyclicity}(a)
or~\ref{Grothendick-abelian-DG-functor-preserves-Becker-co} above.

 Furthermore, one needs to observe that any coacyclic complex of
quasi-coherent sheaves on $X$ is also coacyclic as a (trivial)
quasi-coherent CDG\+module over~$\cB^\cu$.
 This is clear in the context of the Positselski coacyclicity, as
the related DG\+functor $\bCom(X\Qcoh)\rarrow\bQcohr\cB^\cu$ is exact
and preserves infinite direct sums.
 In the context of the Becker coacyclicity, one can also observe that
the DG\+functor of trivial quasi-coherent CDG\+modules $\bCom(X\Qcoh)
\rarrow\bQcohr\cB^\cu$ has adjoints on both sides; so either one
of Lemmas~\ref{DG-functors-preserve-Becker-co-contra-acyclicity}(a)
or~\ref{Grothendick-abelian-DG-functor-preserves-Becker-co}
is applicable.
 Finally, by Lemma~\ref{filtered-by-coacyclic-is-coacyclic-lemma}
(for a finite filtration~$F$), it follows that $\M^\bu\ot_{\cO_X}
\bigwedge_X^*(\g)$ is a coacyclic quasi-coherent CDG\+module
over $\cB^\cu$ whenever the complex $\M^\bu$ is semicoacyclic.

 Let $\N^\cu$ be a quasi-coherent right CDG\+module over~$\cB^\cu$.
 We claim that the complex of quasi-coherent right $\cA$\+modules
$\N^\cu\ot_{\cO_X}\cA$ is actually coacyclic (and not merely
semicoacyclic) whenever the quasi-coherent CDG\+module $\N^\cu$ is
coacyclic.
 In the case of the Positselski coacyclicity, this holds because
the DG\+functor ${-}\ot_{\cO_X}\cA$ is exact and preserves infinite
direct sums; so it takes coacyclic objects to coacyclic objects.
 In the case of the Becker coacyclicity, we need to refer to
Lemma~\ref{DG-functors-preserve-Becker-co-contra-acyclicity}(a)
or~\ref{Grothendick-abelian-DG-functor-preserves-Becker-co}.

 Next we have to prove that the cone of the adjunction morphism
$\M^\bu\ot_{\cO_X}\bigwedge_X^*(\g)\ot_{\cO_X}\cA\rarrow\M^\bu$
is a semicoacyclic complex for any complex of quasi-coherent
right $\cA$\+modules~$\M^\bu$.
 This simply means that the cone is coacyclic as a complex of
of quasi-coherent sheaves on~$X$.
 All Positselski-coacyclic objects are Becker-coacyclic; so it
suffices to check the Positselski coacyclicity.

 Consider the Poincar\'e--Birkhoff--Witt increasing filtration $F$
on the quasi-coherent enveloping quasi-algebra
$\cA=\cA_X(\g,\widetilde\g)$, as in
Section~\ref{enveloping-algebra-subsecn} and
Corollary~\ref{twisted-lie-algebroid-pbw-qcoh-cor}.
 Define the increasing filtration $F$ on the graded quasi-coherent
sheaf $\bigwedge_X^*(\g)$ by the rule $F_n\bigl(\bigwedge_X^*(\g)\bigr)
=\bigoplus_{j\le n}\bigwedge_X^j(\g)$.
 Let $F$ be the increasing filtration on the tensor product
$\M^\bu\ot_{\cO_X}\bigwedge_X^*(\g)\ot_{\cO_X}\cA$ induced by
the increasing filtrations $F$ on $\bigwedge_X^*(\g)$ and~$\cA$.
 Finally, denote also by $F$ the increasing filtration on
the cone of the adjunction morphism $\M^\bu\ot_{\cO_X}\bigwedge_X^*(\g)
\ot_{\cO_X}\cA\rarrow\M^\bu$ induced by the increasing filtration $F$
on $\M^\bu\ot_{\cO_X}\bigwedge_X^*(\g)\ot_{\cO_X}\cA$.
 Here the filtration $F$ on the complex $\M^\bu$ is concentrated in
the filtration degree $n=0$.

 Then the claim is that the successive quotient pieces of the filtration
$F$ on the cone of the adjunction morphism are absolutely acyclic
complexes of quasi-coherent sheaves on~$X$.
 The point is that the associated graded complex of quasi-coherent
sheaves to the complex
$\M^\bu\ot_{\cO_X}\bigwedge_X^*(\g)\ot_{\cO_X}\cA$ is isomorphic to
the complex obtained by applying the functor of tensor product of
complexes of quasi-coherent sheaves $\M^\bu\ot_{\cO_X}{-}$ to
the Koszul complex
$$
 \bigwedge\nolimits_X^*(\g)\ot_{\cO_X}\Sym_X^*(\g).
$$
 As the graded pieces of the Koszul complex are finite acyclic complexes
of finite locally free sheaves on $X$ (with the exception of
the component of degree~$0$, which is isomorphic to~$\cO_X$),
it remains to point out that the tensor product of any complex of
quasi-coherent sheaves with a finite acyclic complex of flat
quasi-coherent sheaves is absolutely acyclic in $X\Qcoh$ in order
to establish the claim.
 By Lemma~\ref{filtered-by-coacyclic-is-coacyclic-lemma}
(for an infinite increasing filtration~$F$), it follows that
the cone of the adjunction morphism is coacyclic as a complex of
quasi-coherent sheaves on~$X$.

 Finally, we have to show that the cone of the adjunction morphism
$\N^\cu\rarrow\N^\cu\ot_{\cO_X}\cA\ot_{\cO_X}\bigwedge_X^*(\g)$ is 
coacyclic in $\bQcohr\cB^\cu$ for every quasi-coherent right
CDG\+module $\N^\cu$ over~$\cB^\cu$.
 Once again, the Positselski coacyclicity implies the Becker
coacyclicity; so it suffices to check the Positselski coacyclicity.
 In order to prove the desired assertion, we reduce the question to
the case of a trivial quasi-coherent CDG\+module $\N^\cu$ over~$\cB^\cu$
(corresponding to some complex of quasi-coherent sheaves $\N^\bu$
on~$X$).

 The point is that any quasi-coherent CDG\+module over $\cB^\cu$
admits a finite filtration by quasi-coherent CDG\+submodules with
trivial successive quotient quasi-coherent CDG\+modules.
 It suffices to put $F^i\N^*=F^i\cB^*\cdot\N^*=
\im(F^i\cB^*\ot_{\cO_X}\N^*\to\N^*)$, where $F^i\cB^*=
\bigoplus_{j\ge i}\cB^j$ is the decreasing filtration associated with
the grading of the quasi-coherent graded algebra~$\cB^*$.
 Then the filtration components $F^i\N^*\subset\N^*$ are quasi-coherent
graded $\cB^*$\+submodules of~$\N^*$ preserved by the differentials,
so they have natural structures of quasi-coherent CDG\+submodules
$F^i\N^\cu\subset\N^\cu$; and the successive quotient quasi-coherent
CDG\+modules $F^i\N^\cu/F^{i+1}\N^\cu$ are trivial.
 The finite decreasing filtration $F$ on $\N^\cu$ induces a finite
decreasing filtration $F$ on $\N^\cu\ot_{\cO_X}\cA\ot_{\cO_X}
\bigwedge_X^*(\g)$.
 The DG\+functor ${-}\ot_{\cO_X}\cA\ot_{\cO_X}\bigwedge_X^*(\g)$ is
exact, so it commutes with the passage to the successive quotients.

 By Lemma~\ref{filtered-by-coacyclic-is-coacyclic-lemma}
(for a finite decreasing filtration~$F$), we can conclude that it
suffices to check that the cone of the adjunction morphism
$\N^\bu\rarrow\N^\bu\ot_{\cO_X}\cA\ot_{\cO_X}\bigwedge_X^*(\g)$
is coacyclic in $\bQcohr\cB^\cu$ for a trivial quasi-coherent
right CDG\+module $\N^\bu$ over~$\cB^\cu$.
 The latter assertion is the result of
Lemma~\ref{trivial-qcoh-koszul-coresolution-lemma}.
\end{proof}

\begin{rem} \label{loc-Noetherian-right-co-side-b-and-p-agree-remark}
 In the case of a locally Noetherian scheme $X$, there is no
difference between the Positselski and Becker versions of
the exotic derived categories in
Theorem~\ref{semiderived-koszul-duality-right-co-side}.
 So the assertions of parts~(a) and~(b) of the theorem become one
and the same triangulated equivalence
$$
 \xymatrix{
  {-}\ot_{\cO_X}\bigwedge\nolimits_X^*(\g)\:
  \sD^{\si=\bsi}(\Qcohr\cA) \ar@{=}[r]
  & \sD^{\co=\bco}(\bQcohr\cB^\cu) \,:{-}\ot_{\cO_X}\cA.
 }
$$
 Indeed, for a locally Noetherian scheme $X$ we have
$\sD^\si(\Qcohr\cA)=\sD^\bsi(\Qcohr\cA)$ by
Corollary~\ref{noetherian-Positselski=Becker-semicoderived}.
 Then it follows by comparing parts~(a) and~(b) of the theorem
that $\sD^\co(\bQcohr\cB^\cu)=\sD^\bco(\bQcohr\cB^\cu)$.
 In the case of a Noetherian scheme $X$,
Corollary~\ref{noetherian-quasi-algebra-Positselski=Becker} is also
applicable, claiming that
$\sD^\co(\bQcohr\cB^\cu)=\sD^\bco(\bQcohr\cB^\cu)$.
\end{rem}

\subsection{Semiderived Koszul duality on the contra side}
\label{semiderived-koszul-duality-contra-side-subsecn}
 We refer to Section~\ref{semiderived-contraherent-subsecn}
for the definitions of the Positselski semi(contra)derived categories
$\sD^\si(\cA\Lcth_\bW)$, \ $\sD^\si(\cA\Lcth_\bW^{X\dlct})$, and
$\sD^\si(\cA\Lcth_\bW^{\cA\dlct})$.
 The following theorem is our first formulation of
the ``$\cD$\+$\Omega$\+duality on the contra side'' promised
in the title of this paper.
 It is (essentially) a nonaffine, singular/non-Noetherian
generalization of~\cite[Theorem~B.2(b)]{Pkoszul}.
 It is also a nonaffine version of~\cite[Theorem~7.11]{Prel}.

\begin{thm} \label{semiderived-koszul-duality-contra-side}
 Let $X$ be a quasi-compact semi-separated scheme with an open
covering\/ $\bW$ and $(\g,\widetilde\g)$ be a quasi-coherent twisted
Lie algebroid over~$X$.
 Assume that\/ $\g$~is a finite locally free sheaf on~$X$.
 Let $\cA=\cA_X(\g,\widetilde\g)$ be the enveloping quasi-coherent
quasi-algebra and $\cB^\cu=\cC^\cu_X(\g,\widetilde\g)$ be
the Chevalley--Eilenberg quasi-coherent CDG\+quasi-algebra
of~$(\g,\widetilde\g)$.
 In this setting: \par
\textup{(a)} the pair of adjoint DG\+functors from
Lemma~\textup{\ref{koszul-duality-dg-functors-contra-side}(a)} induces
a triangulated equivalence between the Positselski semicontraderived
and the Positselski contraderived category
$$
 \xymatrix{
  \Cohom_X\bigl(\bigwedge\nolimits_X^*(\g),{-}\bigr)\:
  \sD^\si(\cA\Lcth_\bW) \ar@{=}[r]
  & \sD^\ctr(\cB^\cu\bLcth_\bW) \,:\!\Cohom_X(\cA,{-});
 }
$$ \par
\textup{(b)} the pair of adjoint DG\+functors from
Lemma~\textup{\ref{koszul-duality-dg-functors-contra-side}(b)} induces
a triangulated equivalence between the Positselski semicontraderived
and the Positselski contraderived category
$$
 \xymatrix{
  \Cohom_X\bigl(\bigwedge\nolimits_X^*(\g),{-}\bigr)\:
  \sD^\si(\cA\Lcth_\bW^{X\dlct}) \ar@{=}[r]
  & \sD^\ctr(\cB^\cu\bLcth_\bW^{X\dlct}) \,:\!\Cohom_X(\cA,{-});
 }
$$ \par
\textup{(c)} the pair of adjoint DG\+functors from
Lemma~\textup{\ref{koszul-duality-dg-functors-contra-side}(c)} induces
a triangulated equivalence between the Positselski semicontraderived
and the Positselski contraderived category
$$
 \xymatrix{
  \Cohom_X\bigl(\bigwedge\nolimits_X^*(\g),{-}\bigr)\:
  \sD^\si(\cA\Lcth_\bW^{\cA\dlct}) \ar@{=}[r]
  & \sD^\ctr(\cB^\cu\bLcth_\bW^{X\dlct}) \,:\!\Cohom_X(\cA,{-}).
 }
$$
\end{thm}

\begin{proof}
 Four assertions need to be proved in each part~(a), (b), and~(c).
 Firstly, one has to show that the functor
$\Cohom_X\bigl(\bigwedge_X^*(\g),{-}\bigr)$ takes semicontraacyclic
complexes to contraacyclic $\bW$\+locally contraherent CDG\+modules,
and that the functor $\Cohom_X(\cA,{-})$ takes contraacyclic
$\bW$\+locally contraherent CDG\+modules to semicontraacyclic complexes.
 Secondly, it needs to be checked that the cone of the adjunction
morphism $\P^\bu\rarrow\Cohom_X(\cA,\Cohom_X(\bigwedge_X^*(\g),\P^\bu))$
is semicontraacyclic for every complex $\P^\bu\in\Com(\cA\Lcth_\bW)$,
and that the cone of the adjunction morphism
$\Cohom_X(\bigwedge_X^*(\g),\Cohom_X(\cA,\Q^\cu))\rarrow\Q^\cu$ is
contraacyclic for every $\bW$\+locally contraherent CDG\+module
$\Q^\cu\in\cB^\cu\bLcth_\bW$.
 The argument is dual-analogous to, but at one point more complicated
than the proof of
Theorem~\ref{semiderived-koszul-duality-right-co-side}.

 Let us spell out the proof of part~(a).
 Let $\P^\bu$ be a complex in the exact category $\cA\Lcth_\bW$.
 The the $\bW$\+locally contraherent CDG\+module
$\Cohom_X\bigl(\bigwedge_X^*(\g),\P^\bu\bigr)$ has a natural finite
decreasing filtration $F$ by admissible subobjects in the exact
category of $\bW$\+locally contraherent CDG\+modules
$\sZ^0(\cB^\cu\bLcth_\bW)$.
 The filtration $F$ on $\Cohom_X\bigl(\bigwedge_X^*(\g),\P^\bu\bigr)$
is induced by the increasing filtration $F$ on $\bigwedge_X^*(\g)$
given by the rule $F_i\bigl(\bigwedge_X^*(\g)\bigr)=\bigoplus_{j\le i}
\bigwedge_X^j(\g)$ for all $i\in\boZ$; so we have
$$
 F^i\Cohom_X\Bigl(\bigwedge\nolimits_X^*(\g),\P^\bu\Bigr)=
 \Cohom_X\Bigl(\bigwedge\nolimits_X^*(\g)\Big/
 F_{i-1}\Bigl(\bigwedge\nolimits_X^*(\g)\Bigr),\P^\bu\Bigr).
$$
 In fact, the $\bW$\+locally contraherent CDG\+module
$\Cohom_X\bigl(\bigwedge_X^*(\g),\P^\bu\bigr)$ over $\cB^\cu$ is thick
in the sense of Sections~\ref{thick-loc-contraherent-modules-subsecn}
and~\ref{exact-dg-categories-of-thick-lcth-subsecn}, and the filtration
$F$ defined above is just the filtration $F$ from
Corollary~\ref{thick-lcta-lct-cosheaves-are-filtered}(a); but we do not
need to use these observations at this point.

 It is only important that the successive quotients
$$
 F^i\Cohom_X\Bigl(\bigwedge\nolimits_X^*(\g),\P^\bu\Bigr)\Big/
 F^{i+1}\Cohom_X\Bigl(\bigwedge\nolimits_X^*(\g),\P^\bu\Bigr)
$$
are trivial $\bW$\+locally contraherent CDG\+modules over $\cB^\cu$
(in the sense of
Section~\ref{koszul-resolutions-of-trivial-ctrh-cdg-mods})
corresponding to the complexes of $\bW$\+locally contraherent
cosheaves $\Cohom_X\bigl(\bigwedge_X^{-i}(\g),\P^\bu\bigr)$.
 Now $\P^\bu$ being a semicontraacyclic complex of $\bW$\+locally
contraherent $\cA$\+modules means precisely that $\P^\bu$ is
(Positselski) contraacyclic as a complex of $\bW$\+locally contraherent
cosheaves on~$X$.
 If this is the case, then the complexes of $\bW$\+locally contraherent
cosheaves $\Cohom_X\bigl(\bigwedge_X^{-i}(\g),\P^\bu\bigr)$ are
contraacyclic as well (because the functor
$\Cohom_X\bigl(\bigwedge_X^{-i}(\g),{-})\:X\Lcth_\bW\rarrow
X\Lcth_\bW$ is exact and preserves infinite products).

 Furthermore, one needs to observe that any Positselski-contraacyclic
complex of $\bW$\+locally contraherent cosheaves on $X$ is also
Positselski-contraacyclic as a (trivial) $\bW$\+locally contraherent
CDG\+module over~$\cB^\cu$.
 Once again, this holds because the DG\+functor of trivial
$\bW$\+locally contraherent CDG\+modules $\bCom(X\Lcth_\bW)\rarrow
\cB^\cu\bLcth_\bW$ is exact and preserves infinite products.
 Finally, by the dual version of
Lemma~\ref{filtered-by-coacyclic-is-coacyclic-lemma}
(for a finite filtration~$F$), it follows that
$\Cohom_X\bigl(\bigwedge_X^*(\g),\P^\bu\bigr)$ is a contraacyclic
$\bW$\+locally contraherent CDG\+module over $\cB^\cu$ whenever
the complex $\P^\bu$ is semicontraacyclic.

 Let $\Q^\cu$ be a $\bW$\+locally contraherent CDG\+module
over~$\cB^\cu$.
 We claim that the complex of $\bW$\+locally contraherent
$\cA$\+modules $\Cohom_X(\cA,\Q^\cu)$ is actually contraacyclic
(and not merely semicontraacyclic) whenever the $\bW$\+locally
contraherent CDG\+module $\Q^\cu$ is contraacyclic.
 Indeed, the DG\+functor $\Cohom_X(\cA,{-})\:\cB^\cu\bLcth_\bW
\allowbreak\rarrow\bCom(\cA\Lcth_\bW)$ is exact and preserves
infinite products; so it takes Positselski-contraacyclic objects to
Positselski-contraacyclic objects.
{\hfuzz=4pt\par}

 Next we have to prove that the cone of the adjunction morphism
$\P^\bu\rarrow\Cohom_X(\cA,\Cohom_X(\bigwedge_X^*(\g),\P^\bu))$
is a semicontraacyclic complex for any complex of $\bW$\+locally
contraherent $\cA$\+modules~$\P^\bu$.
 This simply means that the cone is Positselski-contraacyclic as
a complex in the exact category $X\Lcth_\bW$ of $\bW$\+locally
contraherent cosheaves on~$X$.

 Consider the Poincar\'e--Birkhoff--Witt increasing filtration $F$
on the quasi-coherent enveloping quasi-algebra
$\cA=\cA_X(\g,\widetilde\g)$, as in
Section~\ref{enveloping-algebra-subsecn} and
Corollary~\ref{twisted-lie-algebroid-pbw-qcoh-cor}.
 Define the increasing filtration $F$ on the graded quasi-coherent
sheaf $\bigwedge_X^*(\g)$ by the rule $F_n\bigl(\bigwedge_X^*(\g)\bigr)
=\bigoplus_{j\le n}\bigwedge_X^j(\g)$.
 Let $F$ be the increasing filtration on the tensor product
$\bigwedge_X^*(\g)\ot_{\cO_X}\cA$ induced by the increasing filtrations
$F$ on $\bigwedge_X^*(\g)$ and~$\cA$.
 Denote also by $F$ the decreasing filtration on the complex of
$\bW$\+locally contraherent cosheaves
$$
 \Cohom_X\Bigl(\cA,
 \Cohom_X\Bigl(\bigwedge\nolimits_X^*(\g),\P^\bu\Bigr)\Bigr) =
 \Cohom_X\Bigl(\bigwedge\nolimits_X^*(\g)\ot_{\cO_X}\cA,\>\P^\bu\Bigr)
$$
induced by the increasing filtration $F$ on
$\bigwedge_X^*(\g)\ot_{\cO_X}\cA$.
 So we have
\begin{multline*}
 F^n\Cohom_X\Bigl(\bigwedge\nolimits_X^*(\g)\ot_{\cO_X}\cA,
 \>\P^\bu\Bigr) \\ =
 \Cohom_X\Bigl(\Bigl(\bigwedge\nolimits_X^*(\g)\ot_{\cO_X}\cA\Bigr)
 \Big/F_{n-1}\Bigl(\bigwedge\nolimits_X^*(\g)\ot_{\cO_X}\cA\Bigr),
 \>\P^\bu\Bigr),
\end{multline*}
and one can easily see that
\begin{multline*}
 \Cohom_X\Bigl(\bigwedge\nolimits_X^*(\g)\ot_{\cO_X}\cA,\>\P^\bu\Bigr)
 \\ = \varprojlim\nolimits_n\Bigl(
 \Cohom_X\Bigl(\bigwedge\nolimits_X^*(\g)\ot_{\cO_X}\cA,\>\P^\bu\Bigr)
 \Big/F^n
 \Cohom_X\Bigl(\bigwedge\nolimits_X^*(\g)\ot_{\cO_X}\cA,\>\P^\bu\Bigr)
 \Bigr),
\end{multline*}
where
\begin{multline*}
 \Cohom_X\Bigl(\bigwedge\nolimits_X^*(\g)\ot_{\cO_X}\cA,\>\P^\bu\Bigr)
 \Big/F^n
 \Cohom_X\Bigl(\bigwedge\nolimits_X^*(\g)\ot_{\cO_X}\cA,\>\P^\bu\Bigr)
 \\ = \Cohom_X\Bigl(
 F_{n-1}\Bigl(\bigwedge\nolimits_X^*(\g)\ot_{\cO_X}\cA\Bigr),
 \>\P^\bu\Bigr).
\end{multline*}
 Finally, denote also by $F$ the decreasing filtration on the cone of
the adjunction morphism $\P^\bu\rarrow
\Cohom_X(\cA,\Cohom_X(\bigwedge_X^*(\g),\P^\bu))$ induced by
the decreasing filtration $F$ on
$\Cohom_X(\cA,\Cohom_X(\bigwedge_X^*(\g),\P^\bu))$.
 Here the filtration $F$ on the complex $\P^\bu$ is concentrated in
the filtration degree $n=0$.

 Then the claim is that the successive quotient pieces of the filtration
$F$ on the cone of the adjunction morphism are absolutely acyclic
complexes of $\bW$\+locally contraherent cosheaves on~$X$.
 The point is that the associated graded complex of $\bW$\+locally
contraherent cosheaves to the complex
$\Cohom_X(\cA,\Cohom_X(\bigwedge_X^*(\g),\P^\bu))$ is isomorphic to
the complex obtained by applying the functor $\Cohom_X({-},\P^\bu)$
to the Koszul complex
$$
 \bigwedge\nolimits_X^*(\g)\ot_{\cO_X}\Sym_X^*(\g).
$$
 The graded pieces of the Koszul complex are finite acyclic complexes
of finite locally free sheaves on $X$ (with the exception of
the component of degree~$0$, which is isomorphic to~$\cO_X$).
 So it remains to point out that the complex $\Cohom_X$ from
a finite acyclic complex of very flat quasi-coherent sheaves into
any complex of $\bW$\+locally contraherent cosheaves on $X$
is absolutely acyclic in $X\Lcth_\bW$ in order to establish the claim.
 By the dual version of
Lemma~\ref{filtered-by-coacyclic-is-coacyclic-lemma}
(for an infinite decreasing filtration~$F$), it follows that
the cone of the adjunction morphism is contraacyclic as a complex of
$\bW$\+locally contraherent cosheaves on~$X$.

 Finally, we have to show that the cone of the adjunction morphism
$\Cohom_X(\bigwedge_X^*(\g),\allowbreak\Cohom_X(\cA,\Q^\cu))
\rarrow\Q^\cu$ is Positselski-contraacyclic in $\cB^\cu\bLcth_\bW$ for
every $\bW$\+lo\-cally contraherent CDG\+module $\Q^\cu$ over~$\cB^\cu$.
 For this purpose, we reduce the question to the case of a trivial
$\bW$\+locally contraherent CDG\+module $\Q^\cu$ over~$\cB^\cu$
(corresponding to some complex of $\bW$\+locally contraherent cosheaves
$\Q^\bu$ on~$X$).

 This is where our argument gets more complicated than the one in
the proof of Theorem~\ref{semiderived-koszul-duality-right-co-side}.
 The difference with the quasi-coherent case is that, given an arbitrary
$\bW$\+locally contraherent graded $\cB^*$\+module $\Q^*$,
a straightforward attempt to define a filtration $F$ on $\Q^*$ induced
(in one way or another) by the decreasing filtration $F$ on $\cB^*$
associated with the grading of $\cB^*$ runs into a problem.

 The problem is that, dealing as we are with $\bW$\+locally contraherent
cosheaves, the filtration $F$ on $\Q^*$ must be, first of all,
a filtration by $\bW$\+locally contraherent subcosheaves.
 So the graded $\cB^*(U)$\+modules $F^i\Q^*[U]$ must satisfy two
conditions~(i) and~(ii) from
Section~\ref{locally-contraherent-cosheaves-subsecn}.
 For an arbitrary $\cO_X(U)$\+con\-tra\-ad\-justed graded
$\cB^*(U)$\+module $\Q^*[U]$, defining a filtration $F$ on $\Q^*[U]$ by
either one of the rules
$F^i\Q^*[U]=\im(F^i\cB^*(U)\ot_{\cO_X(U)}\Q^*[U]\to\Q^*[U])$
or $F^{-i}\Q^*[U]\allowbreak=\ker(\Q^*[U]\to\Hom_{\cO_X(U)}(\cB^*(U)/
F^{i+1}\cB^*(U),\Q^*[U]))$ \emph{cannot} guarantee that the graded
$\cO_X(U)$\+module $F^i\Q^*[U]$ would satisfy these conditions.

 To be more specific, the two constructions run into two different
problems:
\begin{itemize}
\item Setting
$F^i\Q^*[U]=\im(F^i\cB^*(U)\ot_{\cO_X(U)}\Q^*[U]\to\Q^*[U])$ for
affine open subschemes $U\subset X$ subordinate to $\bW$, one
makes the grading components of the graded $\cO(U)$\+module
$F^i\Q^*[U]$ contraadjusted.
 Indeed, the class of contraadjusted modules is closed under quotients
and finite direct sums; so contraadjustedness of
the $\cO_X(U)$\+modules $\Q^n[U]$ for all $n\in\boZ$ implies
contraadjustedness of $F^i\Q^n[U]$.
 Thus condition~(ii) is satisfied.
 But there seems to be just no reason for the maps $F^i\Q^*[V]\rarrow
\Hom_{\cO_X(U)}(\cO_X(V),F^i\Q^*[U])$ to be isomorphisms for
affine open subschemes $V\subset U\subset X$ subordinate to~$\bW$.
\item Setting $F^{-i}\Q^*[U]=\ker(\Q^*[U]\to\Hom_{\cO_X(U)}(\cB^*(U)/
F^{i+1}\cB^*(U),\Q^*[U]))$, one forces the maps $F^i\Q^*[V]\rarrow
\Hom_{\cO_X(U)}(\cO_X(V),F^i\Q^*[U])$ to be isomorphisms for
affine open subschemes $V\subset U\subset X$ subordinate to $\bW$,
as the functors $\Hom_{\cO_X(U)}(\cO_X(V),{-})$ preserve kernels.
 So condition~(i) is satisfied.
 But there seems to be just no reason for the $\cO_X(U)$\+modules
$F^i\Q^n[U]$ to be contraadjusted in this construction, and
condition~(ii) is problematic.
\end{itemize}

 The point is that the construction of the filtration $F$ on
an arbitrary graded $\cB^*(U)$\+module $\Q^*[U]$ by any one of
the two rules above represents a ``brutal'' or ``underived'' approach.
 One can observe that the functor assigning to a $\cB^*(U)$\+module
$\Q^*[U]$ any one of the two filtrations above on $\Q^*[U]$ is
not exact in general (cf.\
Lemmas~\ref{vwr-projective-canonical-filtration-lemma}\+-%
\ref{vwr-injective-canonical-filtration-lemma}, where the very weak
relative projectivity/injectivity condition is imposed).
 The underived, nonexact construction of the filtration was good enough
for the proof of
Theorem~\ref{semiderived-koszul-duality-right-co-side}, but it is
not sufficient for our  present purposes.

 The explanation is that category $\sZ^0(\bQcohr\cB^\cu)$ is abelian,
but the category $\sZ^0(\cB^\cu\bLcth_\bW)$ is only exact.
 As a general rule, whenever a functor into an abelian category is not
exact, an analogous functor into an exact category tends to be only
partially defined (see~\cite[Section~5.6]{Pphil} and~\cite[Section~0.10 
of the Introduction]{Pcosh}).
 In our present context of $\bW$\+locally contraherent CDG\+modules
over $\cB^\cu$, the proper approach now is to ``derive''
the construction of the filtration.

 So we use the upper line of
diagram~\eqref{thick-lcta-cdg-contrader-equivs-diagram} in
Corollary~\ref{thick-cdg-contraderived-equiv-cor}, which tells us that
the inclusion of exact DG\+categories $\cB^\cu\bLcth_\bW^\bth\rarrow
\cB^\cu\bLcth_\bW$ induces a triangulated equivalence
$\sD^\ctr(\cB^\cu\bLcth_\bW^\bth)\simeq\sD^\ctr(\cB^\cu\bLcth_\bW)$.
 More specifically, for any given $\bW$\+locally contraherent
CDG\+module $\Q^\cu$ over $\cB^\cu$ there exists a thick
$\bW$\+locally contraherent CDG\+module $\R^\cu$ over $\cB^\cu$
together with a closed morphism $\R^\cu\rarrow\Q^\cu$ with
a cone contraacyclic in $\cB^\cu\bLcth_\bW$.
 This is where we use the assumption that the scheme $X$ is
quasi-compact and semi-separated.

 Now it follows from the arguments in the first half of this proof that
the cone of the induced morphism
$\Cohom_X(\bigwedge_X^*(\g),\Cohom_X(\cA,\R^\cu))\rarrow
\Cohom_X(\bigwedge_X^*(\g),\allowbreak\Cohom_X(\cA,\Q^\cu))$ is
contraacyclic in $\cB^\cu\bLcth_\bW$.
 Hence, in order to prove that the cone of the morphism
$\Cohom_X(\bigwedge_X^*(\g),\Cohom_X(\cA,\Q^\cu))\rarrow\Q^\cu$
is contraacyclic, it suffices to show that the cone of the morphism
$\Cohom_X(\bigwedge_X^*(\g),\Cohom_X(\cA,\R^\cu))\rarrow\R^\cu$ is
contraacyclic.

 Corollary~\ref{thick-lcta-lct-cosheaves-are-filtered}(a) provides
a finite decreasing filtration $F$ on the $\bW$\+locally contraherent
graded $\cB^*$\+module $\R^*$ by its admissible subobjects in
the exact category of $\bW$\+locally contraherent graded
$\cB^*$\+modules.
 One can easily see that the graded submodules $F^i\R^*[U]\subset
\R^*[U]$ are preserved by the differentials (basically because
the graded ideals $F^i\cB^*(U)\subset\cB^*(U)$ are preserved by
the differentials).
 So we have a filtration of the $\bW$\+locally contraherent
CDG\+module $\R^\cu$ over $\cB^\cu$ by admissible subobjects
$F^i\R^\cu\subset\R^\cu$ in the exact category
$\sZ^0(\cB^\cu\bLcth_\bW)$.

 The successive quotient $\bW$\+locally contraherent CDG\+modules
$F^i\R^\cu/F^{i+1}\R^\cu$ are trivial $\bW$\+locally contraherent
CDG\+modules over~$\cB^\cu$.
 The finite decreasing filtration $F$ on $\R^\cu$ induces a finite
decreasing filtration $F$ on
$\Cohom_X(\bigwedge_X^*(\g),\Cohom_X(\cA,\R^\cu))$.
 The DG\+functor $\Cohom_X(\bigwedge_X^*(\g),\Cohom_X(\cA,{-}))$
is exact, so it commutes with the passage to the successive quotients.
 By the dual version of
Lemma~\ref{filtered-by-coacyclic-is-coacyclic-lemma}
(for an infinite decreasing filtration~$F$), we can conclude that it
suffices to check that the cone of the adjunction morphism
$\Cohom_X(\bigwedge_X^*(\g),\Cohom_X(\cA,\Q^\bu))\rarrow\Q^\bu$
is contraacyclic in $\cB^\cu\bLcth_\bW$ for a trivial
$\bW$\+locally contraherent CDG\+module $\Q^\bu$ over~$\cB^\cu$.
 The latter assertion is the result of
Lemma~\ref{trivial-lcth-koszul-resolution-lemma}(a).

 The proofs of parts~(b\+-c) are similar.
 One needs to use the upper line of
diagram~\eqref{thick-X-lct-cdg-contrader-equivs-diagram} in
Corollary~\ref{thick-cdg-contraderived-equiv-cor}, together with
Corollary~\ref{thick-lcta-lct-cosheaves-are-filtered}(b)
and Lemma~\ref{trivial-lcth-koszul-resolution-lemma}(b).
\end{proof}

\begin{cor} \label{semicontraderived-A-lct=X-lct-corollary}
 Let $X$ be a quasi-compact semi-separated scheme with an open
covering\/ $\bW$ and $(\g,\widetilde\g)$ be a quasi-coherent twisted
Lie algebroid over~$X$.
 Assume that\/ $\g$~is a finite locally free sheaf on~$X$.
 Let $\cA=\cA_X(\g,\widetilde\g)$ be the twisted universal enveloping
quasi-coherent quasi-algebra of~$(\g,\widetilde\g)$.
 Then the inclusion of exact categories $\cA\Lcth_\bW^{\cA\dlct}\rarrow
\cA\Lcth_\bW^{X\dlct}$ induces an equivalence of the Positselski
semicontraderived categories
$$
 \sD^\si(\cA\Lcth_\bW^{\cA\dlct})\simeq\sD^\si(\cA\Lcth_\bW^{X\dlct}).
$$
\end{cor}

\begin{proof}
 We have a commutative diagram of triangulated functors
\begin{equation} \label{semicontraderived-A-lct=X-lct-diagram}
\begin{gathered}
 \qquad\quad\xymatrix{
  \text{\llap{$\Cohom_X\bigl(\bigwedge\nolimits_X^*(\g),{-}\bigr)\:$}}
  \sD^\si(\cA\Lcth_\bW^{X\dlct}) \ar@{=}[r]
  & \sD^\ctr(\cB^\cu\bLcth_\bW^{X\dlct})
  \text{\rlap{$\,\,:\!\Cohom_X(\cA,{-})$}} \\
  \text{\llap{$\Cohom_X\bigl(\bigwedge\nolimits_X^*(\g),{-}\bigr)\:$}}
  \sD^\si(\cA\Lcth_\bW^{\cA\dlct}) \ar@{=}[r] \ar[u]
  & \sD^\ctr(\cB^\cu\bLcth_\bW^{X\dlct})
  \text{\rlap{$\,\,:\!\Cohom_X(\cA,{-})$}} \ar@{=}[u]
 }
\end{gathered}
\end{equation}
 Here the horizontal triangulated equivalences are provided by
Theorem~\ref{semiderived-koszul-duality-contra-side}(b\+-c).
 The rightmost vertical functor is the identity triangulated
equivalence.
 The leftmost vertical functor is induced by the inclusion of exact
categories $\cA\Lcth_\bW^{\cA\dlct}\rarrow\cA\Lcth_\bW^{X\dlct}$.
 Since three sides of the commutative
square~\eqref{semicontraderived-A-lct=X-lct-diagram} are triangulated
equivalences, so is the fourth one.
 Thus the leftmost vertical functor is a triangulated equivalence,
as desired.
\end{proof}

\begin{rem} \label{contraderived-of-lct=of-lcta-remark}
 One expects the contraderived categories of $\cA$\+locally cotorsion
$\bW$\+locally contraherent $\cA$\+modules, $X$\+locally cotorsion
$\bW$\+locally contraherent $\cA$\+modules, and arbitrary ($X$\+locally
contraadjusted) $\bW$\+locally contraherent $\cA$\+mod\-ules to be
equivalent for a quasi-coherent quasi-algebra $\cA$ over a quasi-compact
semi-separated scheme $X$ quite generally.
 The same is expected to hold for the related semicontraderived 
categories.
 Moreover, the same can be conjectured about the contraderived
categories of $\bW$\+locally contraherent CDG\+modules over
any quasi-coherent CDG\+quasi-algebra $\cA^\cu$ over~$X$.
 However, such general expectations apply to the Becker
(semi)contraderived categories and \emph{not} to the Positselski ones.
 For somewhat similar or related (but simpler) results,
see~\cite[Corollary~7.21]{Pphil}, \cite[last line of
Corollary~11.2]{Pbc}, \cite[Theorem~5.5.10]{Pcosh}, and
the paper~\cite{PS7}.
 Under more restrictive assumptions on a scheme $X$, the equivalence
between the exotic derived categories of $X$\+locally cotorsion and
$X$\+locally contraadjusted $\bW$\+locally contraherent CDG\+modules
over $\cA^\cu$ is provable for the Positselski as well as Becker
contraderived categories, and even for the absolute derived categories;
see Theorem~\ref{lcth-lcta-lct-cdg-abs-contraderived-equiv-thm}
below.  \emergencystretch=1em
\end{rem}

\subsection{Semiderived Koszul duality on the contra side,
Becker version}
\label{becker-semiderived-koszul-duality-contra-side-subsecn}
 The definitions of the Becker semi(contra)derived categories
$\sD^\bsi(\cA\Lcth_\bW)$, \ $\sD^\bsi(\cA\Lcth_\bW^{X\dlct})$, and
$\sD^\bsi(\cA\Lcth_\bW^{\cA\dlct})$ were given in
Section~\ref{becker-semicontraderived-defined-subsecn}.
 The aim of this section is to prove a version of
Theorem~\ref{semiderived-koszul-duality-contra-side} for Becker
semicontraderived and contraderived categories.
 We start with three lemmas. {\hbadness=2525\par}

\begin{lem} \label{Cohom-from-finite-loc-free-becker-contraacyclic}
 Let $X$ be a scheme with an open covering\/ $\bW$ and $\E$ be a finite
locally free quasi-coherent sheaf on~$X$.
 Then \par
\textup{(a)} the functor\/ $\Cohom_X(\E,{-})$ takes
Becker-contraacyclic complexes in $X\Lcth_\bW$ to Becker-contraacyclic
complexes in $X\Lcth_\bW$; \par
\textup{(b)} the functor\/ $\Cohom_X(\E,{-})$ takes
Becker-contraacyclic complexes in $X\Lcth_\bW^\lct$ to
Becker-contraacyclic complexes in $X\Lcth_\bW^\lct$.
\end{lem}

\begin{proof}
 Put $\F=\cHom_{\cO_X}(\F,\cO_X)$; so $\F$ is also a finite locally
free sheaf on~$X$.
 The point is that, in part~(a), the functor $\Cohom_X(\F,{-})\:
X\Lcth_\bW\rarrow X\Lcth_\bW$ is adjoint to the functor
$\Cohom_X(\E,{-})\:X\Lcth_\bW\rarrow X\Lcth_\bW$ on both sides.
 The same adjunction holds, in part~(b), for the same functors $\Cohom$
acting on the category $X\Lcth_\bW^\lct$.
 As the functors $\Cohom_X(\E,{-})$ are also exact,
the result of~\cite[Lemma~B.7.5(b)]{Pcosh}
or Lemma~\ref{DG-functors-preserve-Becker-co-contra-acyclicity}(b)
is applicable.
\end{proof}

\begin{lem} \label{trivial-contrah-cdg-module-becker-contraacyclic}
 Let $X$ be a quasi-compact semi-separated scheme with an open
covering\/ $\bW$ and $(\g,\widetilde\g)$ be a quasi-coherent twisted
Lie algebroid over~$X$.
 Assume that\/ $\g$~is a finite locally free sheaf on~$X$.
 Let $\cB^\cu=\cC^\cu_X(\g,\widetilde\g)$ be the Chevalley--Eilenberg
quasi-coherent CDG\+quasi-algebra of~$(\g,\widetilde\g)$.
 Then \par
\textup{(a)} the DG\+functor of trivial\/ $\bW$\+locally contraherent
CDG\+modules\/ $\bCom(X\Lcth_\bW)\allowbreak\rarrow\cB^\cu\bLcth_\bW$
takes Becker-contraacyclic complexes in $X\Lcth_\bW$ to
Becker-contraacyclic objects of $\cB^\cu\bLcth_\bW$;
{\hfuzz=6.5pt \hbadness=1100 \par}
\textup{(b)} the DG\+functor of trivial $X$\+locally cotorsion\/
$\bW$\+locally contraherent CDG\+mod\-ules\/ $\bCom(X\Lcth_\bW^\lct)
\rarrow\cB^\cu\bLcth_\bW^{X\dlct}$ takes Becker-contraacyclic complexes
in $X\Lcth_\bW^\lct$ to Becker-contraacyclic objects of
$\cB^\cu\bLcth_\bW^{X\dlct}$.
\end{lem}

\begin{proof}
 Let us prove part~(a).
 The argument is based on the observation that the existence of
a partially defined left adjoint DG\+functor that is defined on
the graded-projective objects is sufficient for the validity of
Lemma~\ref{DG-functors-preserve-Becker-co-contra-acyclicity}(b).
 Specifically, let $\gM^\bu$ be a Becker-contraacyclic complex in
$X\Lcth_\bW$, viewed as a trivial $\bW$\+locally contraherent
CDG\+module over~$\cB^\cu$.
 Let $\P^\cu$ be a graded-projective object of the exact category
$\cB^\cu\bLcth_\bW$.
 We need to show that the complex of abelian groups
$\Hom^{\cB^*,\bu}(\P^\cu,\gM^\bu)$ is acyclic.

 Consider the graded $\bW$\+locally contraherent $\cB^*$\+module~$\P^*$.
 By assumption, $\P^*$ is a projective object of the exact category
$\cB^*\Lcth_\bW$.
 By Lemma~\ref{enough-thick-antilocal-in-lcth}(a), there exists a thick
(and also antilocal) $\bW$\+locally contraherent graded $\cB^*$\+module
$\Q^*$ together with an admissible epimorphism $\Q^*\rarrow\P^*$ in
$\cB^*\Lcth_\bW$.
 It follows that $\P^*$ is a direct summand of~$\Q^*$.
 So the graded $\bW$\+locally contraherent $\cB^*$\+module $\P^*$ is
thick (and antilocal, but this is not important for us).

 Thus $\P^\cu$ is a thick $\bW$\+locally contraherent CDG\+module
over~$\cB^\cu$.
 As mentioned near the end of the proof of
Theorem~\ref{semiderived-koszul-duality-contra-side},
Corollary~\ref{thick-lcta-lct-cosheaves-are-filtered}(a) provides
a finite decreasing filtration of $\P^\cu$ by its admissible subobjects
in the exact category of $\bW$\+locally contraherent CDG\+modules
$\sZ^0(\cB^\cu\bLcth_\bW)$.
 We will need a version of the filtration from
Corollary~\ref{thick-lcta-lct-cosheaves-are-filtered}(a) based on
Lemma~\ref{vwr-projective-canonical-filtration-lemma} instead of
Lemma~\ref{vwr-injective-canonical-filtration-lemma} (decomposing $X$
into a disjoint union of open subschemes $X_d$ on which the finite
locally free sheaf~$\g$ has constant rank~$d$, the two filtrations only
differ by a shift of the filtration degrees by~$d$ on every~$X_d$).
 So we have a decreasing filtration $F$ of $\P^\cu$ by its admissible
subobjects in $\sZ^0(\cB^\cu\bLcth_\bW)$ such that $\P^\cu=F^0\P^\cu$
and, for every affine open subscheme $U\subset X$,
the graded $\cB^*(U)$\+module $\gr_F^0\P^*[U]=F^0\P^*[U]/F^1\P^*[U]$ is
the maximal quotient $\cB^*(U)$\+module of $\P^*[U]$ on which the action
of $\cB^*(U)$ is trivial.

 Now $\gr_F^0\P^\cu$ is a trivial $\bW$\+locally contraherent
CDG\+module over~$\cB^\cu$ (corresponding to some complex of
$\bW$\+locally contraherent cosheaves on~$X$).
 For any graded $\bW$\+locally contraherent cosheaf $\gN^*$ on $X$,
we have a natural isomorphism of graded abelian groups
$$
 \Hom^{\cB^*,*}(\P^*,\gN^*)\simeq\Hom^{X,*}(\gr_F^0\P^*,\gN^*),
$$
where $\Hom^{X,*}({-},{-})$ denotes the graded abelian group of
homogeneous morphisms of graded objects in $X\Lcth_\bW$.
 For any complex of $\bW$\+locally contraherent cosheaves $\gN^\bu$
on $X$, we have a natural isomorphism of complexes of abelian groups
$$
 \Hom^{\cB^*,\bu}(\P^\cu,\gN^\bu)\simeq
 \Hom^{X,\bu}(\gr_F^0\P^\cu,\gN^\bu),
$$
where $\Hom^{X,\bu}({-},{-})$ denotes the complex of morphisms 
between two complexes in $X\Lcth_\bW$.

 The functor of trivial $\bW$\+locally contraherent graded modules
over~$\cB^*$ (assigning to every graded object of $X\Lcth_\bW$
the corresponding trivial $\bW$\+locally contraherent graded module
over~$\cB^*$) is exact.
 Since the object $\P^*$ is projective in $\cB^*\Lcth_\bW$, it follows
that the functor $\gN^*\longmapsto\Hom^{X,*}(\gr_F^0\P^*,\gN^*)$
is exact on the category of graded $\bW$\+locally contraherent
cosheaves on~$X$.
 So the grading components of the graded $\bW$\+locally contraherent
cosheaf $\gr_F^0\P^*$ are projective, and $\gr_F\P^\cu$ is
a graded-projective complex of $\bW$\+locally contraherent cosheaves
on~$X$.
 Returning to the situation at hand, since $\gM^\bu$ is
a Becker-contraacyclic complex in $X\Lcth_\bW$, the complex of abelian
groups $\Hom^{X,\bu}(\gr_F^0\P^\cu,\gM^\bu)$ is acyclic.
 Thus the complex of abelian groups $\Hom^{\cB^*,\bu}(\P^\cu,\gM^\bu)$
is acyclic as well, and we are done.

 The proof of part~(b) is similar and based
on Lemma~\ref{enough-thick-antilocal-in-lcth}(b)
and Corollary~\ref{thick-lcta-lct-cosheaves-are-filtered}(b).
\end{proof}

\begin{lem} \label{Cohom-X-from-A-becker-contraacyclic}
 Let $X$ be a quasi-compact semi-separated scheme with an open
covering\/ $\bW$ and $(\g,\widetilde\g)$ be a quasi-coherent twisted
Lie algebroid over~$X$.
 Assume that\/ $\g$~is a finite locally free sheaf on~$X$.
 Let $\cA=\cA_X(\g,\widetilde\g)$ be the enveloping quasi-coherent
quasi-algebra and $\cB^\cu=\cC^\cu_X(\g,\widetilde\g)$ be
the Chevalley--Eilenberg quasi-coherent CDG\+quasi-algebra
of~$(\g,\widetilde\g)$.
 Then \par
\textup{(a)} the DG\+functor\/ $\Cohom_X(\cA,{-})\:
\cB^\cu\bLcth_\bW\rarrow\bCom(\cA\Lcth_\bW)$ from
Lemma~\textup{\ref{koszul-duality-dg-functors-contra-side}(a)} takes
Becker-contraacyclic objects of $\cB^\cu\bLcth_\bW$ to
Becker-con\-tra\-acyclic complexes in $\cA\Lcth_\bW$;
{\hbadness=2500\par}
\textup{(b)} the DG\+functor\/ $\Cohom_X(\cA,{-})\:
\cB^\cu\bLcth_\bW^{X\dlct}\rarrow\bCom(\cA\Lcth_\bW^{X\dlct})$ from
Lemma~\textup{\ref{koszul-duality-dg-functors-contra-side}(b)} takes
Becker-contraacyclic objects of $\cB^\cu\bLcth_\bW^{X\dlct}$ to
Becker-con\-tra\-acyclic complexes in $\cA\Lcth_\bW^{X\dlct}$; \par
\textup{(c)} the DG\+functor\/ $\Cohom_X(\cA,{-})\:
\cB^\cu\bLcth_\bW^{X\dlct}\rarrow\bCom(\cA\Lcth_\bW^{\cA\dlct})$ from
Lemma~\textup{\ref{koszul-duality-dg-functors-contra-side}(c)} takes
Becker-contraacyclic objects of $\cB^\cu\bLcth_\bW^{X\dlct}$ to
Becker-con\-tra\-acyclic complexes in $\cA\Lcth_\bW^{\cA\dlct}$.
\end{lem}

\begin{proof}
 All the three DG\+functors in question are exact (as DG\+functors
between exact DG\+categories).
 By Lemma~\ref{koszul-duality-dg-functors-contra-side},
all the three DG\+functors also have left adjoint DG\+functors.
 So Lemma~\ref{DG-functors-preserve-Becker-co-contra-acyclicity}(b)
is applicable.
\end{proof}

\begin{thm} \label{semiderived-koszul-duality-becker-contra-side}
 Let $X$ be a quasi-compact semi-separated scheme with an open
covering\/ $\bW$ and $(\g,\widetilde\g)$ be a quasi-coherent twisted
Lie algebroid over~$X$.
 Assume that\/ $\g$~is a finite locally free sheaf on~$X$.
 Let $\cA=\cA_X(\g,\widetilde\g)$ be the enveloping quasi-coherent
quasi-algebra and $\cB^\cu=\cC^\cu_X(\g,\widetilde\g)$ be
the Chevalley--Eilenberg quasi-coherent CDG\+quasi-algebra
of~$(\g,\widetilde\g)$.
 In this setting: \par
\textup{(a)} the pair of adjoint DG\+functors from
Lemma~\textup{\ref{koszul-duality-dg-functors-contra-side}(a)} induces
a triangulated equivalence between the Becker semicontraderived
and the Becker contraderived category
$$
 \xymatrix{
  \Cohom_X\bigl(\bigwedge\nolimits_X^*(\g),{-}\bigr)\:
  \sD^\bsi(\cA\Lcth_\bW) \ar@{=}[r]
  & \sD^\bctr(\cB^\cu\bLcth_\bW) \,:\!\Cohom_X(\cA,{-});
 }
$$ \par
\textup{(b)} the pair of adjoint DG\+functors from
Lemma~\textup{\ref{koszul-duality-dg-functors-contra-side}(b)} induces
a triangulated equivalence between the Becker semicontraderived
and the Becker contraderived category
$$
 \xymatrix{
  \Cohom_X\bigl(\bigwedge\nolimits_X^*(\g),{-}\bigr)\:
  \sD^\bsi(\cA\Lcth_\bW^{X\dlct}) \ar@{=}[r]
  & \sD^\bctr(\cB^\cu\bLcth_\bW^{X\dlct}) \,:\!\Cohom_X(\cA,{-});
 }
$$ \par
\textup{(c)} the pair of adjoint DG\+functors from
Lemma~\textup{\ref{koszul-duality-dg-functors-contra-side}(c)} induces
a triangulated equivalence between the Becker semicontraderived
and the Becker contraderived category
$$
 \xymatrix{
  \Cohom_X\bigl(\bigwedge\nolimits_X^*(\g),{-}\bigr)\:
  \sD^\bsi(\cA\Lcth_\bW^{\cA\dlct}) \ar@{=}[r]
  & \sD^\bctr(\cB^\cu\bLcth_\bW^{X\dlct}) \,:\!\Cohom_X(\cA,{-}).
 }
$$
\end{thm}

\begin{proof}
 We follow the proof of
Theorem~\ref{semiderived-koszul-duality-contra-side}.
 Let us sketch the proof of part~(a).
 Let $\P^\bu$ be a complex in the exact category $\cA\Lcth_\bW$.
 Then, arguing similarly to the proof of
Theorem~\ref{semiderived-koszul-duality-contra-side} and using
Lemmas~\ref{Cohom-from-finite-loc-free-becker-contraacyclic}(a)
and~\ref{trivial-contrah-cdg-module-becker-contraacyclic}(a),
one proves that the $\bW$\+locally contraherent CDG\+module
$\Cohom_X\bigl(\bigwedge_X^*(\g),\P^\bu\bigr)$ over $\cB^\cu$ is
Becker-contraacyclic whenever the complex $\P^\bu$
is Becker-semicontraacyclic.

 Let $\Q^\cu$ be a Becker-contraacyclic $\bW$\+locally contraherent
CDG\+module over~$\cB^\cu$.
 Then, by Lemma~\ref{Cohom-X-from-A-becker-contraacyclic}(a),
the complex of $\bW$\+locally contraherent $\cA$\+modules
$\Cohom_X(\cA,\Q^\cu)$ is Becker-contraacyclic.
 By Corollary~\ref{forg-functor-preserves-becker-contraacyclicity}(a),
it follows that the complex $\Cohom_X(\cA,\Q^\cu)$ is
Becker-semicontraacyclic.

 Let $\P^\bu$ be a complex of $\bW$\+locally contraherent
$\cA$\+modules.
 According to the proof of
Theorem~\ref{semiderived-koszul-duality-contra-side},
the cone of the adjunction morphism
$\P^\bu\rarrow\Cohom_X(\cA,\Cohom_X(\bigwedge_X^*(\g),\P^\bu))$
is a Positselski-semicontraacyclic complex.
 Since all Positselski-contraacyclic complexes in $X\Lcth_\bW$
(as in any exact category with exact products) are Becker-contraacyclic,
it follows that the same cone is a Becker-contraacyclic complex.
{\hbadness=1350\par}

 Finally, let $\Q^\cu$ be a $\bW$\+locally contraherent CDG\+module
over~$\cB^\cu$.
 According to the proof of
Theorem~\ref{semiderived-koszul-duality-contra-side},
the cone of the adjunction morphism
$\Cohom_X(\bigwedge_X^*(\g),\Cohom_X(\cA,\Q^\cu))\rarrow\Q^\cu$ is
Positselski-contraacyclic in $\cB^\cu\bLcth_\bW$.
 Hence the same cone is Becker-contraacyclic in $\cB^\cu\bLcth_\bW$.
{\hbadness=1375\par}

 The proofs of parts~(b\+-c) are similar and use
Lemmas~\ref{Cohom-from-finite-loc-free-becker-contraacyclic}(b),
\ref{trivial-contrah-cdg-module-becker-contraacyclic}(b),
\ref{Cohom-X-from-A-becker-contraacyclic}(b\+-c), and
Corollary~\ref{forg-functor-preserves-becker-contraacyclicity}(b\+-c).
\end{proof}

\begin{cor} \label{becker-semicontraderived-A-lct=X-lct-corollary}
 Let $X$ be a quasi-compact semi-separated scheme with an open
covering\/ $\bW$ and $(\g,\widetilde\g)$ be a quasi-coherent twisted
Lie algebroid over~$X$.
 Assume that\/ $\g$~is a finite locally free sheaf on~$X$.
 Let $\cA=\cA_X(\g,\widetilde\g)$ be the twisted universal enveloping
quasi-coherent quasi-algebra of~$(\g,\widetilde\g)$.
 Then the inclusion of exact categories $\cA\Lcth_\bW^{\cA\dlct}\rarrow
\cA\Lcth_\bW^{X\dlct}$ induces an equivalence of the Becker
semicontraderived categories
$$
 \sD^\bsi(\cA\Lcth_\bW^{\cA\dlct})\simeq\sD^\bsi(\cA\Lcth_\bW^{X\dlct}).
$$
\end{cor}

\begin{proof}
 The argument is similar to the proof of
Corollary~\ref{semicontraderived-A-lct=X-lct-corollary} and based on
Theorem~\ref{semiderived-koszul-duality-becker-contra-side}(b\+-c).
 The point is that the square diagram of triangulated functors
\begin{equation} \label{becker-semicontraderived-A-lct=X-lct-diagram}
\begin{gathered}
 \qquad\quad\xymatrix{
  \text{\llap{$\Cohom_X\bigl(\bigwedge\nolimits_X^*(\g),{-}\bigr)\:$}}
  \sD^\bsi(\cA\Lcth_\bW^{X\dlct}) \ar@{=}[r]
  & \sD^\bctr(\cB^\cu\bLcth_\bW^{X\dlct})
  \text{\rlap{$\,\,:\!\Cohom_X(\cA,{-})$}} \\
  \text{\llap{$\Cohom_X\bigl(\bigwedge\nolimits_X^*(\g),{-}\bigr)\:$}}
  \sD^\bsi(\cA\Lcth_\bW^{\cA\dlct}) \ar@{=}[r] \ar[u]
  & \sD^\bctr(\cB^\cu\bLcth_\bW^{X\dlct})
  \text{\rlap{$\,\,:\!\Cohom_X(\cA,{-})$}} \ar@{=}[u]
 }
\end{gathered}
\end{equation}
is commutative, so it follows that the leftmost vertical functor is
a triangulated equivalence.
\end{proof}

\begin{rem} \label{Noetherian-contra-side-b-and-p-agree-remark}
 Similarly to
Remark~\ref{loc-Noetherian-right-co-side-b-and-p-agree-remark},
in the case of a semi-separated Noetherian scheme $X$ of finite
Krull dimension there is no difference between the Positselski
versions of the exotic derived categories in
Theorem~\ref{semiderived-koszul-duality-contra-side} and
the respective Becker versions of such categories in
Theorem~\ref{semiderived-koszul-duality-becker-contra-side}.
 So the assertions of
Theorem~\ref{semiderived-koszul-duality-contra-side}(a) and
Theorem~\ref{semiderived-koszul-duality-becker-contra-side}(a) become
one and the same triangulated equivalence
$$
 \xymatrix{
  \Cohom_X\bigl(\bigwedge\nolimits_X^*(\g),{-}\bigr)\:
  \sD^{\si=\bsi}(\cA\Lcth_\bW) \ar@{=}[r]
  & \sD^{\ctr=\bctr}(\cB^\cu\bLcth_\bW) \,:\!\Cohom_X(\cA,{-}).
 }
$$
 The assertions of
Theorem~\ref{semiderived-koszul-duality-contra-side}(b) and
Theorem~\ref{semiderived-koszul-duality-becker-contra-side}(b) become
one and the same triangulated equivalence
$$
 \xymatrixcolsep{0.6em} \xymatrix{
  \Cohom_X\bigl(\bigwedge\nolimits_X^*(\g),{-}\bigr)\:
  \sD^{\si=\bsi}(\cA\Lcth_\bW^{X\dlct}) \ar@{=}[r]
  & \sD^{\ctr=\bctr}(\cB^\cu\bLcth_\bW^{X\dlct}) \,:\!\Cohom_X(\cA,{-}).
 }
$$
 The assertions of
Theorem~\ref{semiderived-koszul-duality-contra-side}(c) and
Theorem~\ref{semiderived-koszul-duality-becker-contra-side}(c) become
one and the same triangulated equivalence
$$
 \xymatrixcolsep{0.6em} \xymatrix{
  \Cohom_X\bigl(\bigwedge\nolimits_X^*(\g),{-}\bigr)\:
  \sD^{\si=\bsi}(\cA\Lcth_\bW^{\cA\dlct}) \ar@{=}[r]
  & \sD^{\ctr=\bctr}(\cB^\cu\bLcth_\bW^{X\dlct}) \,:\!\Cohom_X(\cA,{-}).
 }
$$

 Indeed, for a semi-separated Noetherian scheme $X$ of finite Krull
dimension, we have $\sD^\si(\cA\Lcth_\bW)=\sD^\bsi(\cA\Lcth_\bW)$
by Corollary~\ref{noetherian-Positselski=Becker-semicontraderived}(a).
 More generally, for any semi-separated Noetherian scheme $X$, we have
$\sD^\si(\cA\Lcth_\bW^{X\dlct})=\sD^\bsi(\cA\Lcth_\bW^{X\dlct})$ and
$\sD^\si(\cA\Lcth_\bW^{\cA\dlct})=\sD^\bsi(\cA\Lcth_\bW^{\cA\dlct})$ by
Corollary~\ref{noetherian-Positselski=Becker-semicontraderived}(b\+-c).
 Then it follows by comparing
Theorem~\ref{semiderived-koszul-duality-contra-side}(a)
with Theorem~\ref{semiderived-koszul-duality-becker-contra-side}(a)
that $\Ac^\ctr(\cB^\cu\bLcth_\bW)=\Ac^\bctr(\cB^\cu\bLcth_\bW)$
and $\sD^\ctr(\cB^\cu\bLcth_\bW)=\sD^\bctr(\cB^\cu\bLcth_\bW)$
whenever the scheme $X$ is semi-separated and Noetherian of finite
Krull dimension.
 Comparing 
Theorem~\ref{semiderived-koszul-duality-contra-side}(b) or~(c) with
Theorem~\ref{semiderived-koszul-duality-becker-contra-side}(b) or~(c),
we see that $\Ac^\ctr(\cB^\cu\bLcth_\bW^{X\dlct})=
\Ac^\bctr(\cB^\cu\bLcth_\bW^{X\dlct})$ and
$\sD^\ctr(\cB^\cu\bLcth_\bW^{X\dlct})=
\sD^\bctr(\cB^\cu\bLcth_\bW^{X\dlct})$ whenever the scheme $X$ is
semi-separated and Noetherian.
\end{rem}

\subsection{Reduced Koszul duality on the co side}
\label{reduced-koszul-duality-co-side-subsecn}
 We refer to Section~\ref{reduced-coderived-of-qcoh-subsecn} for
the definition of the reduced (Positselski) coderived category of
quasi-coherent left CDG\+modules $\sD^\co_{X\red}(\cB^\cu\bQcoh)$ over
the Chevalley--Eilenberg quasi-coherent CDG\+quasi-algebra~$\cB^\cu$.
 The reduced coderived category of quasi-coherent right CDG\+modules
$\sD^\co_{X\red}(\bQcohr\cB^\cu)$ is defined similarly.
 The following theorem is (essentially) another generalization
of~\cite[Theorem~B.2(a)]{Pkoszul} to singular/non-Noetherian schemes.
 It is also a nonaffine version of~\cite[Theorem~6.20]{Prel}.

\begin{thm} \label{reduced-koszul-duality-right-co-side}
 Let $X$ be a scheme and $(\g,\widetilde\g)$ be a quasi-coherent
twisted Lie algebroid over~$X$.
 Assume that\/ $\g$~is a finite locally free sheaf of bounded rank
on~$X$.
 Let $\cA=\cA_X(\g,\widetilde\g)$ be the enveloping quasi-coherent
quasi-algebra and $\cB^\cu=\cC^\cu_X(\g,\widetilde\g)$ be
the Chevalley--Eilenberg quasi-coherent CDG\+quasi-algebra
of $(\g,\widetilde\g)$; so $\cB^n=0$ for $n$~large enough.
 Then the pair of adjoint DG\+functors from
Lemma~\ref{koszul-duality-dg-functors-right-co-side} induces
a triangulated equivalence between the conventional derived
category and the reduced Positselski coderived category
$$
 \xymatrix{
  {-}\ot_{\cO_X}\bigwedge\nolimits_X^*(\g)\:
  \sD(\Qcohr\cA) \ar@{=}[r]
  & \sD^\co_{X\red}(\bQcohr\cB^\cu) \,:{-}\ot_{\cO_X}\cA.
 }
$$
 Moreover, there is a natural commutative square diagram of
triangulated equivalences and triangulated Verdier quotient functors
\begin{equation} \label{reduced-koszul-duality-right-co-side-diagram}
\begin{gathered}
 \qquad\quad\xymatrix{
  \text{\llap{${-}\ot_{\cO_X}\bigwedge\nolimits_X^*(\g)\:$}}
  \sD^\si(\Qcohr\cA) \ar@{=}[r] \ar@{->>}[d]
  & \sD^\co(\bQcohr\cB^\cu)
  \text{\rlap{$\,{}:{-}\ot_{\cO_X}\cA$}} \ar@{->>}[d] \\
  \text{\llap{${-}\ot_{\cO_X}\bigwedge\nolimits_X^*(\g)\:$}}
  \sD(\Qcohr\cA) \ar@{=}[r]
  & \sD^\co_{X\red}(\bQcohr\cB^\cu)
  \text{\rlap{$\,{}:{-}\ot_{\cO_X}\cA$}}
 }
\end{gathered}
\end{equation}
\end{thm}

\begin{proof}
 The triangulated equivalence in the upper line
of~\eqref{reduced-koszul-duality-right-co-side-diagram} is provided
by Theorem~\ref{semiderived-koszul-duality-right-co-side}(a).
 In order to prove the theorem, two assertions remain to be checked.
 Firstly, one needs to show that the functor
${-}\ot_{\cO_X}\bigwedge\nolimits_X^*(\g)$ takes acyclic complexes
of quasi-coherent right $\cA$\+modules to reduced-coacyclic
quasi-coherent right CDG\+modules over~$\cB^\cu$.
 Secondly, it needs to be checked that the functor
${-}\ot_{\cO_X}\cA$ takes trivial quasi-coherent right CDG\+modules
over $\cB^\cu$ corresponding to acyclic complexes of quasi-coherent
sheaves on $X$ to acyclic complexes of quasi-coherent right
$\cA$\+modules.

 Both the assertions are easily established.
 For any complex of quasi-coherent right $\cA$\+modules $\M^\bu$,
the quasi-coherent right CDG\+module $\M^\bu\ot_{\cO_X}
\bigwedge_X^*(\g)$ over $\cB^\cu$ has a natural finite decreasing
filtration $F$ by quasi-coherent CDG\+submodules constructed in
the beginning of the proof of
Theorem~\ref{semiderived-koszul-duality-right-co-side}.
 As explained in that proof, the successive quotients to the filtration
$F$ on $\M^\bu\ot_{\cO_X}\bigwedge_X^*(\g)$ are trivial quasi-coherent
CDG\+modules over $\cB^\cu$ corresponding to the complexes of
quasi-coherent sheaves $\M^\bu\ot_{\cO_X}\bigwedge_X^{-i}(\g)$.

 Now if the complex of quasi-coherent $\cA$\+modules $\M^\bu$ is
acyclic, then so are the complexes of quasi-coherent sheaves
$\M^\bu\ot_{\cO_X}\bigwedge_X^{-i}(\g)$.
 It follows that the quasi-coherent CDG\+module
$\M^\bu\ot_{\cO_X}\bigwedge_X^*(\g)$ belongs to the minimal
triangulated subcategory of $\sH^0(\bQcohr\cB^\cu)$ containing
the absolutely acyclic quasi-coherent CDG\+modules over $\cB^\cu$
and the trivial quasi-coherent CDG\+modules corresponding to
acyclic complexes of quasi-coherent sheaves on~$X$.

 Given a complex of quasi-coherent sheaves $\N^\bu$ on $X$ viewed
as a trivial quasi-coherent right CDG\+module over $\cB^\cu$,
the related complex of quasi-coherent right $\cA$\+modules
$\N^\bu\ot_{\cO_X}\cA$ is simply the tensor product of the complex
of quasi-coherent sheaves $\N^\bu$ with the quasi-coherent
quasi-algebra~$\cA$.
 So the differential on $\N^\bu\ot_{\cO_X}\cA$ is induced by
the differential on~$\N^\bu$.
 As $\cA$ is a flat quasi-coherent sheaf in its left (as well as
right) $\cO_X$\+module structure by
Corollary~\ref{twisted-lie-algebroid-pbw-qcoh-cor}(b), it follows
that the complex $\N^\bu\ot_{\cO_X}\cA$ is acyclic whenever
the complex $\N^\bu$~is.
\end{proof}

\subsection{Reduced Koszul duality on the contra side}
\label{reduced-koszul-duality-contra-side-subsecn}
 We refer to Section~\ref{reduced-contraderived-of-lcth-subsecn} for
the definitions of the reduced (Positselski) contraderived categories
$\sD^\ctr_{X\red}(\cB^\cu\bLcth_\bW)$ and
$\sD^\ctr_{X\red}(\cB^\cu\bLcth_\bW^{X\dlct})$.
 The following theorem is our second formulation of
the ``$\cD$\+$\Omega$ duality on the contra side''.
 It is essentially another nonaffine, singular/non-Noetherian
generalization of~\cite[Theorem~B.2(b)]{Pkoszul}.
 It is also a nonaffine version of~\cite[Theorem~7.16]{Prel}.
{\hfuzz=2.5pt\par}

\begin{thm} \label{reduced-koszul-duality-contra-side}
 Let $X$ be a quasi-compact semi-separated scheme with an open
covering\/ $\bW$ and $(\g,\widetilde\g)$ be a quasi-coherent twisted
Lie algebroid over~$X$.
 Assume that\/ $\g$~is a finite locally free sheaf on~$X$.
 Let $\cA=\cA_X(\g,\widetilde\g)$ be the enveloping quasi-coherent
quasi-algebra and $\cB^\cu=\cC^\cu_X(\g,\widetilde\g)$ be
the Chevalley--Eilenberg quasi-coherent CDG\+quasi-algebra
of~$(\g,\widetilde\g)$.
 In this setting: \par
\textup{(a)} The pair of adjoint DG\+functors from
Lemma~\textup{\ref{koszul-duality-dg-functors-contra-side}(a)} induces
a triangulated equivalence between the conventional derived category
and the reduced Positselski contraderived category
$$
 \xymatrix{
  \Cohom_X\bigl(\bigwedge\nolimits_X^*(\g),{-}\bigr)\:
  \sD(\cA\Lcth_\bW) \ar@{=}[r]
  & \sD^\ctr_{X\red}(\cB^\cu\bLcth_\bW) \,:\!\Cohom_X(\cA,{-}).
 }
$$
 Moreover, there is a natural commutative square diagram of
triangulated equivalences and triangulated Verdier quotient functors
\begin{equation} \label{reduced-koszul-duality-lcta-contra-side-diagram}
\begin{gathered}
 \qquad\quad\xymatrix{
  \text{\llap{$\Cohom_X\bigl(\bigwedge\nolimits_X^*(\g),{-}\bigr)\:$}}
  \sD^\si(\cA\Lcth_\bW) \ar@{=}[r] \ar@{->>}[d]
  & \sD^\ctr(\cB^\cu\bLcth_\bW)
  \text{\rlap{$\,\,:\!\Cohom_X(\cA,{-})$}} \ar@{->>}[d] \\
  \text{\llap{$\Cohom_X\bigl(\bigwedge\nolimits_X^*(\g),{-}\bigr)\:$}}
  \sD(\cA\Lcth_\bW) \ar@{=}[r]
  & \sD^\ctr_{X\red}(\cB^\cu\bLcth_\bW)
  \text{\rlap{$\,\,:\!\Cohom_X(\cA,{-})$}}
 }
\end{gathered}
\end{equation} \par
\textup{(b)} The pair of adjoint DG\+functors from
Lemma~\textup{\ref{koszul-duality-dg-functors-contra-side}(b)} induces
a triangulated equivalence between the conventional derived category
and the reduced Positselski contraderived category
$$
 \xymatrix{
  \Cohom_X\bigl(\bigwedge\nolimits_X^*(\g),{-}\bigr)\:
  \sD(\cA\Lcth_\bW^{X\dlct}) \ar@{=}[r]
  & \sD^\ctr_{X\red}(\cB^\cu\bLcth_\bW^{X\dlct}) \,:\!\Cohom_X(\cA,{-}).
 }
$$
 Moreover, there is a natural commutative square diagram of
triangulated equivalences and triangulated Verdier quotient functors
\begin{equation} \label{reduced-koszul-duality-lct-contra-side-diagram}
\begin{gathered}
 \qquad\quad\xymatrix{
  \text{\llap{$\Cohom_X\bigl(\bigwedge\nolimits_X^*(\g),{-}\bigr)\:$}}
  \sD^\si(\cA\Lcth_\bW^{X\dlct}) \ar@{=}[r] \ar@{->>}[d]
  & \sD^\ctr(\cB^\cu\bLcth_\bW^{X\dlct})
  \text{\rlap{$\,\,:\!\Cohom_X(\cA,{-})$}} \ar@{->>}[d] \\
  \text{\llap{$\Cohom_X\bigl(\bigwedge\nolimits_X^*(\g),{-}\bigr)\:$}}
  \sD(\cA\Lcth_\bW^{X\dlct}) \ar@{=}[r]
  & \sD^\ctr_{X\red}(\cB^\cu\bLcth_\bW^{X\dlct})
  \text{\rlap{$\,\,:\!\Cohom_X(\cA,{-})$}}
 }
\end{gathered}
\end{equation} \par
\textup{(c)} The pair of adjoint DG\+functors from
Lemma~\textup{\ref{koszul-duality-dg-functors-contra-side}(c)} induces
a triangulated equivalence between the conventional derived category
and the reduced Positselski contraderived category
$$
 \xymatrix{
  \Cohom_X\bigl(\bigwedge\nolimits_X^*(\g),{-}\bigr)\:
  \sD(\cA\Lcth_\bW^{\cA\dlct}) \ar@{=}[r]
  & \sD^\ctr_{X\red}(\cB^\cu\bLcth_\bW^{X\dlct}) \,:\!\Cohom_X(\cA,{-}).
 }
$$
 Moreover, there is a natural commutative square diagram of
triangulated equivalences and triangulated Verdier quotient functors
\begin{equation} \label{reduced-koszul-duality-lct-A-contra-side-diagr}
\begin{gathered}
 \qquad\quad\xymatrix{
  \text{\llap{$\Cohom_X\bigl(\bigwedge\nolimits_X^*(\g),{-}\bigr)\:$}}
  \sD^\si(\cA\Lcth_\bW^{\cA\dlct}) \ar@{=}[r] \ar@{->>}[d]
  & \sD^\ctr(\cB^\cu\bLcth_\bW^{X\dlct})
  \text{\rlap{$\,\,:\!\Cohom_X(\cA,{-})$}} \ar@{->>}[d] \\
  \text{\llap{$\Cohom_X\bigl(\bigwedge\nolimits_X^*(\g),{-}\bigr)\:$}}
  \sD(\cA\Lcth_\bW^{\cA\dlct}) \ar@{=}[r]
  & \sD^\ctr_{X\red}(\cB^\cu\bLcth_\bW^{X\dlct})
  \text{\rlap{$\,\,:\!\Cohom_X(\cA,{-})$}}
 }
\end{gathered}
\end{equation}
\end{thm}

\begin{proof}
 The triangulated equivalences in the upper lines of all the three
square diagrams are provided by
Theorem~\ref{semiderived-koszul-duality-contra-side}.
 In order to prove the theorem, two assertions remain to be checked in
each part~(a), (b), and~(c).
 Let us spell out part~(a).

 Firsly, one needs to show that the functor
$\Cohom_X\bigl(\bigwedge_X^*(\g),{-}\bigr)$ takes acyclic complexes in
the exact category $\cA\Lcth_\bW$ to reduced-contraacyclic
$\bW$\+locally contraherent CDG\+modules over~$\cB^\cu$.
 Secondly, it needs to be checked that the functor $\Cohom_X(\cA,{-})$
takes trivial $\bW$\+locally contraherent CDG\+modules over $\cB^\cu$
corresponding to acyclic complexes in the exact category $X\Lcth_\bW$
to acyclic complexes in $\cA\Lcth_\bW$.

 Both the assertions are easily established.
 For any complex of $\bW$\+locally contraherent $\cA$\+modules $\P^\bu$,
the $\bW$\+locally contraherent CDG\+module
$\Cohom_X\bigl(\bigwedge_X^*(\g),\P^\bu\bigr)$ over $\cB^\cu$ has
a natural finite decreasing filtration by admissible subobjects in
$\sZ^0(\cB^\cu\bLcth_\bW)$ constructed in the beginning of the proof of
Theorem~\ref{semiderived-koszul-duality-contra-side}.
 As explained in that proof, the successive quotients to the filtration
$F$ on $\Cohom_X\bigl(\bigwedge_X^*(\g),\allowbreak\P^\bu\bigr)$ are
trivial $\bW$\+locally contraherent CDG\+modules over $\cB^\cu$
corresponding to the complexes of $\bW$\+locally contraherent cosheaves
$\Cohom_X\bigl(\bigwedge_X^{-i}(\g),\P^\bu\bigr)$.

 Now if the complex of $\bW$\+locally contraherent $\cA$\+modules is
acyclic in $\cA\Lcth_\bW$, then the complexes of $\bW$\+locally
contraherent  cosheaves
$\Cohom_X\bigl(\bigwedge_X^{-i}(\g),\P^\bu\bigr)$ are acyclic
in $X\Lcth_\bW$.
 It follows that the $\bW$\+locally contraherent CDG\+module
$\Cohom_X\bigl(\bigwedge_X^*(\g),\P^\bu\bigr)$ belongs to the minimal
triangulated subcategory of $\sH^0(\cB^\cu\allowbreak\bLcth_\bW)$
containing the absolutely acyclic objects of $\cB^\cu\bLcth_\bW$ and
the trivial $\bW$\+locally contraherent CDG\+modules corresponding to
acyclic complexes in $X\Lcth_\bW$.

 Given a complex of $\bW$\+locally contraherent cosheaves $\Q^\bu$ on
$X$ viewed as a trivial $\bW$\+locally contraherent CDG\+module over
$\cB^\cu$, the related complex of $\bW$\+locally contraherent
$\cA$\+modules $\Cohom_X(\cA,\Q^\bu)$ is simply the $\Cohom_X$
complex from the quasi-coherent quasi-algebra $\cA$ to the complex
of $\bW$\+locally contraherent cosheaves~$\Q^\bu$.
 So the differential on $\Cohom_X(\cA,\Q^\bu)$ is induced by
the differential on~$\Q^\bu$.
 As the functor $\Cohom_X$ is exact in both its arguments wherever it
is well-defined (we recall once again that $\cA$ is a very flat
quasi-coherent sheaf on $X$ in its left, as well as right
$\cO_X$\+module structure by
Corollary~\ref{twisted-lie-algebroid-pbw-qcoh-cor}(b)), it follows that
the complex $\Cohom_X(\cA,\Q^\bu)$ is acyclic in $\cA\Lcth_\bW$ whenever
the complex $\Q^\bu$ is acyclic in $X\Lcth_\bW$.

 The proofs of parts~(b\+-c) are similar.
\end{proof}

 The result of the following corollary was used in the proofs of
Corollaries~\ref{lcth-reduced-X-lcta-X-lct-equivs-cor}
and~\ref{reduced-thick-lcth-al-ctrder-lcta-lct-equiv-cor}, and played
an important role in Section~\ref{co-contra-cot-lct-subsecn} above.

\begin{cor} \label{reduced-koszul-duality-contra-side-all-equiv}
 Let $X$ be a quasi-compact semi-separated scheme with an open
covering\/ $\bW$ and $(\g,\widetilde\g)$ be a quasi-coherent twisted
Lie algebroid over~$X$.
 Assume that\/ $\g$~is a finite locally free sheaf on~$X$.
 Let $\cA=\cA_X(\g,\widetilde\g)$ be the enveloping quasi-coherent
quasi-algebra and $\cB^\cu=\cC^\cu_X(\g,\widetilde\g)$ be
the Chevalley--Eilenberg quasi-coherent CDG\+quasi-algebra
of~$(\g,\widetilde\g)$.
 Then there is a natural commutative diagram of triangulated
equivalences
\begin{equation} \label{reduced-koszul-duality-contra-side-big-diagram}
\
\end{equation}
$$
 \qquad\ \ \xymatrix{
  \text{\llap{$\Cohom_X\bigl(\bigwedge\nolimits_X^*(\g),{-}\bigr)\:$}}
  \sD(\cA\Lcth_\bW) \ar@{=}[r] \ar@{-}@<-2pt>[d]
  & \sD^\ctr_{X\red}(\cB^\cu\bLcth_\bW)
  \text{\rlap{$\,\,:\!\Cohom_X(\cA,{-})$}} \ar@{-}@<-2pt>[d] \\
  \text{\llap{$\Cohom_X\bigl(\bigwedge\nolimits_X^*(\g),{-}\bigr)\:$}}
  \sD(\cA\Lcth_\bW^{X\dlct}) \ar@{=}[r] \ar@<-2pt>[u] \ar@{-}@<-2pt>[d]
  & \sD^\ctr_{X\red}(\cB^\cu\bLcth_\bW^{X\dlct})
  \text{\rlap{$\,\,:\!\Cohom_X(\cA,{-})$}}
  \ar@<-2pt>[u] \ar@{=}[d] \\
  \text{\llap{$\Cohom_X\bigl(\bigwedge\nolimits_X^*(\g),{-}\bigr)\:$}}
  \sD(\cA\Lcth_\bW^{\cA\dlct}) \ar@{=}[r] \ar@<-2pt>[u]
  & \sD^\ctr_{X\red}(\cB^\cu\bLcth_\bW^{X\dlct})
  \text{\rlap{$\,\,:\!\Cohom_X(\cA,{-})$}}
 }
$$
with the horizontal equivalences provided by
Theorem~\ref{reduced-koszul-duality-contra-side}, the vertical
equivalences in the leftmost column induced by the inclusions of
exact categories $\cA\Lcth_\bW^{\cA\dlct}\rarrow\cA\Lcth_\bW^{X\dlct}
\rarrow\cA\Lcth_\bW$, and the upper vertical equivalence in
the rightmost column induced by the inclusion of exact DG\+categories
$\cB^\cu\bLcth_\bW^{X\dlct}\rarrow\cB^\cu\bLcth_\bW$.
 The lower vertical equivalence in the rightmost column is
the identity triangulated functor.
\end{cor}

\begin{proof}
 Let us point out once again that $\cB^\cu\bLcth_\bW^{\cB^*\dlct}=
\cB^\cu\bLcth_\bW^{X\dlct}$ by
Corollary~\ref{lct-lcth-modules-over-flfrqa-cor}; this is the reason
why we have an identity lower vertical functor in the rightmost
column.
 It is clear from the constructions that the diagram is commutative.
 The vertical triangulated functors in the leftmost column are
triangulated equivalences by
Theorem~\ref{A-lcth-X-lct-A-lct-derived-equivalence}.
 Taking Theorem~\ref{reduced-koszul-duality-contra-side}(a\+b) into
account, it follows that the upper vertical triangulated functor in
the rightmost column is a triangulated equivalence, too.
\end{proof}

\begin{rem} \label{weakly-smooth-de-Rham-and-diffoperators-remark}
 Continuing the discussion in Remark~\ref{weakly-smooth-de-Rham-remark},
let $X\rarrow T$ be a weakly smooth morphism of schemes (of bounded
relative dimension).
 Then the sheaf of rings of crystalline fiberwise differential operators
$\cD^\cry_{X/T}$ it the quasi-coherent universal enveloping
quasi-algebra of the quasi-coherent Lie algebroid of polyvector fields
$\vect_{X/T}$.
 So we have $\cD^\cry_{X/T}=\cA_X(\vect_{X/T})$ (see
Section~\ref{crystalline-diffoperators-subsecn}).
 Therefore, the results of
Theorems~\ref{semiderived-koszul-duality-right-co-side}
and~\ref{reduced-koszul-duality-right-co-side}
are applicable to the quasi-coherent quasi-algebra of crystalline
differential operators $\cA=\cD^\cry_{X/T}$ and the de~Rham
quasi-coherent DG\+quasi-algebra $\cB^\cu=\Omega^\bu_{X/T}$ on~$X$.
 Assuming that the scheme $X$ is quasi-compact and semi-separated,
the results of
Theorems~\ref{semiderived-koszul-duality-contra-side},
\ref{semiderived-koszul-duality-becker-contra-side},
and~\ref{reduced-koszul-duality-contra-side} are applicable to
$\cA=\cD^\cry_{X/T}$ and $\cB^\cu=\Omega^\bu_{X/T}$ as well.
 In particular, when $X$ is a scheme over a field of
characteristic~$0$, the crystalline differential operators coincide
with the conventional (Grothendieck) algebraic differential operators,
so $\cD^\cry_{X/T}=\cD^\str_{X/T}$ (see
Section~\ref{diffoperators-characteristic-zero-subsecn}).
\end{rem}

\Section{Quadrality Diagram} \label{quadrality-diagram-secn}

\subsection{Dual Koszul coresolutions} \label{dual-koszul-subsecn}
 In the next three Sections~\ref{dual-koszul-subsecn}\+-%
\ref{conversion-functor-subsecn}, we roughly follow the exposition
in~\cite[Sections~9.5\+-9.6]{Prel}.

 Let $A$ and $B$ be two rings, and let $E$ be an $A$\+$B$\+bimodule.
 Then the opposite $B^\rop$\+$A^\rop$\+bimodule $E^\rop$ is defined by
the obvious rules $E^\rop=E$ and $b^\rop e^\rop a^\rop=(aeb)^\rop$
for all $a\in A$, \,$b\in B$, and $e\in E$.

 Let $A^\cu=(A^*,d_A,h_A)$ and $B^\cu=(B^*,d_B,h_B)$ be two
CDG\+rings, and let $E^\cu=(E^*,d_E)$ be a CDG\+bimodule over
$A^\cu$ and $B^\cu$ (as defined in
Section~\ref{cdg-rings-cdg-modules-subsecn}).
 The definition of the opposite CDG\+rings $A^\rop{}^\cu$ and
$B^\rop{}^\cu$ can be found in
Section~\ref{opposite-cdg-rings-and-twisted-lie-subsecn}.
 The CDG\+bimodule $E^\rop{}^\cu$ over $B^\cu$ and $A^\cu$ is
defined by the rules
\begin{itemize}
\item the elements of $E^\rop{}^*$ are the symbols $e^\rop$, where
$e\in E^*$; the grading on $E^\rop{}^*$ agrees with the grading
on~$E^*$; the right action of $A^\rop{}^*$ and left action of
$B^\rop{}^*$ on $E^\rop{}^*$ are defined using the sign rules from
Section~\ref{opposite-cdg-rings-and-twisted-lie-subsecn}, so
$b^\rop e^\rop a^\rop=(-1)^{|a||b|+|a||e|+|b||e|}(aeb)^\rop$ for
all $a\in A^{|a|}$, \ $b\in B^{|b|}$ and $e\in E^{|e|}$;
\item the differential $d_E^\rop$ on $E^\rop{}^*$ is equal to
the differential~$d_E$ on~$E$; more precisely, $d^\rop_E(e^\rop)=
(d_E(e))^\rop$ for all $e\in E$.
\end{itemize}

 For a scheme $X$ and two quasi-coherent CDG\+quasi-algebras $\cA^\cu$
and $\cB^\cu$ over $X$, the opposite quasi-coherent CDG\+quasi-algebras
$\cA^\rop{}^\cu$ and $\cB^\rop{}^\cu$ were defined in
Section~\ref{opposite-cdg-rings-and-twisted-lie-subsecn}.
 Given a quasi-coherent CDG\+bimodule $\cE^\cu$
over $\cA^\cu$ and $\cB^\cu$ (as defined in
Section~\ref{fHom-and-contratensor-of-cdg-modules}), one can apply
the rules above in order to produce the opposite quasi-coherent
CDG\+bimodule $\cE^\rop{}^\cu$ over $\cB^\rop{}^\cu$
and~$\cA^\rop{}^\cu$.

 Let $X$ be a scheme and $(\g,\widetilde\g)$ be a quasi-coherent
twisted Lie algebroid over~$X$.
 Assume that the quasi-coherent sheaf~$\g$ over $X$ is finite locally
free.
 Denote by $\cA=\cA_X(\g,\widetilde\g)$ the quasi-coherent twisted
enveloping quasi-algebra of $(\g,\widetilde\g)$, and by
$\cB^\cu=\cC^\cu_X(\g,\widetilde\g)$ the Chevalley--Eilenberg
quasi-coherent CDG\+quasi-algebra of~$(\g,\widetilde\g)$.
 See Sections~\ref{enveloping-algebra-subsecn}
and~\ref{chevalley-eilenberg-qcoh-cdg-quasi-algebra-subsecn}.

 Recall the notation $(\g,\widetilde\g^\circ)$ from
Section~\ref{opposite-cdg-rings-and-twisted-lie-subsecn} for
the opposite quasi-coherent twisted Lie algebroid
to~$(\g,\widetilde\g)$. 
 Put $\cA^\circ=\cA_X(\g,\widetilde\g^\circ)$ and $\cB^\circ{}^\cu
=\cC^\cu_X(\g,\widetilde\g^\circ)$.

 According to Section~\ref{opposite-cdg-rings-and-twisted-lie-subsecn},
we have a natural (identity) isomorphism of quasi-coherent
CDG\+quasi-algebras $\cB^\circ{}^\cu=\cB^\rop{}^\cu$.
 Let us \emph{warn} the reader once again, however, that
the quasi-coherent quasi-algebra $\cA^\circ$ is usually \emph{not}
isomorphic to the quasi-coherent quasi-algebra~$\cA^\rop$.
 For example, when $X\rarrow T$ is a (weakly) smooth morphism of
schemes over a field of characteristic~$0$ and $\g=\vect_{X/T}$ is
the split quasi-coherent twisted Lie algebroid of vector fields on $X$,
one has $\cA^\circ=\cA$ (as for any split quasi-coherent twisted Lie 
algebroid $\widetilde\g=\cO_X\oplus\g$).
 Nevertheless, the quasi-coherent quasi-algebra of differential
operators $\cD^\str_{X/T}=\cA_X(\vect_{X/T})$ (see
Sections~\ref{crystalline-diffoperators-subsecn}
and~\ref{diffoperators-characteristic-zero-subsecn}) is \emph{not}
isomorphic to its opposite quasi-coherent quasi-algebra
$(\cD^\str_{X/T})^\rop$ in any natural way.

 According to Section~\ref{chevalley-eilenberg-cdg-modules-subsecn},
we have a quasi-coherent CDG\+bimodule
$\cC_X^\cu(\g,\widetilde\g,\cA_X)$ over $\cB^\cu$ and $\cA$,
as well as a quasi-coherent CDG\+bimodule
$\cC^X_\cu(\cA_X,\g,\widetilde\g)$ over $\cA$ and~$\cB^\cu$.
 Applying the same construction to the quasi-coherent twisted Lie
algebroid $(\g,\widetilde\g^\circ)$ over $X$ instead of
$(\g,\widetilde\g)$, we obtain a quasi-coherent CDG\+bimodule
$\cC_X^\cu(\g,\widetilde\g^\circ,\cA_X)$ over $\cB^\circ{}^\cu$ and
$\cA^\circ$, as well as a quasi-coherent CDG\+bimodule
$\cC^X_\cu(\cA_X,\g,\widetilde\g^\circ)$ over $\cA^\circ$
and~$\cB^\circ{}^\cu$.
 Hence we also have the opposite quasi-coherent CDG\+bimodule
$\cC_X^\cu(\g,\widetilde\g^\circ,\cA_X)^\rop$ over
$\cA^{\circ,\rop}$ and $\cB^{\circ,\rop}{}^\cu=\cB^\cu$, as well
as the opposite quasi-coherent CDG\+bimodule
$\cC^X_\cu(\cA_X,\g,\widetilde\g^\circ)^\rop$ over
$\cB^{\circ,\rop}{}^\cu=\cB^\cu$ and~$\cA^{\circ,\rop}$.

 We put $\cC_X^\cu(\cA_X^{\circ,\rop},\g,\widetilde\g^\circ)=
\cC_X^\cu(\g,\widetilde\g^\circ,\cA_X)^\rop$ and
$\cC^X_\cu(\g,\widetilde\g^\circ,\cA_X^{\circ,\rop})=
\cC^X_\cu(\cA_X,\g,\widetilde\g^\circ)^\rop$.
 So $\cC_X^\cu(\cA_X^{\circ,\rop},\g,\widetilde\g^\circ)$ is
a quasi-coherent CDG\+bimodule over $\cA^{\circ,\rop}$ and $\cB^\cu$,
while $\cC^X_\cu(\g,\widetilde\g^\circ,\cA_X^{\circ,\rop})$ is
a quasi-coherent CDG\+bimodule over $\cB^\cu$ and~$\cA^{\circ,\rop}$.
 The underlying quasi-coherent graded bimodule of
$\cC_X^\cu(\cA_X^{\circ,\rop},\g,\widetilde\g^\circ)$ is
$\cC_X^*(\cA_X^{\circ,\rop},\g,\widetilde\g^\circ)=
\cA^{\circ,\rop}\ot_{\cO_X}\cB^*$, while the underlying
quasi-coherent graded bimodule of
$\cC^X_\cu(\g,\widetilde\g^\circ,\cA_X^{\circ,\rop})$
is $\cC^X_*(\g,\widetilde\g^\circ,\cA_X^{\circ,\rop})=
\bigwedge_X^*(\g)\ot_{\cO_X}\cA^{\circ,\rop}$.

 The following lemma should be compared with
Lemma~\ref{trivial-qcoh-koszul-coresolution-lemma}
and~\cite[Lemma~9.7]{Prel}.

\begin{lem} \label{trivial-qcoh-dual-koszul-coresolution}
 Let $X$ be a scheme and $(\g,\widetilde\g)$ be a quasi-coherent
twisted Lie algebroid over~$X$.
 Assume that the quasi-coherent sheaf\/~$\g$ on $X$ is finite locally
free of constant rank~$m$.
 In this setting: \par
\textup{(a)} Let $\M^\bu$ be a complex of quasi-coherent sheaves on~$X$.
 Consider the complex of quasi-coherent sheaves
$\cB^m[-m]\ot_{\cO_X}\M^\bu$, and interpret it as a trivial
quasi-coherent left CDG\+module over~$\cB^\cu$.
 Then the tensor product $\cC_X^\cu(\g,\widetilde\g,\cA_X)
\ot_{\cO_X}\M^\bu$ is a quasi-coherent left CDG\+module over $\cB^\cu$
endowed with a natural closed morphism of quasi-coherent left
CDG\+modules $\cB^m[-m]\ot_{\cO_X}\M^\bu\rarrow
\cC_X^\cu(\g,\widetilde\g,\cA_X)\ot_{\cO_X}\M^\bu$ whose cone is
Positselski-coacyclic in the abelian DG\+category $\cB^\cu\bQcoh$. \par
\textup{(b)} Let $\N^\bu$ be a complex of quasi-coherent sheaves on~$X$.
 Consider the complex of quasi-coherent sheaves
$\N^\bu\ot_{\cO_X}\cB^m[-m]$, and interpret it as a trivial
quasi-coherent right CDG\+module over~$\cB^\cu$.
 Then the tensor product $\N^\bu\ot_{\cO_X}
\cC_X^\cu(\cA_X^{\circ,\rop},\g,\widetilde\g^\circ)$ is
a quasi-coherent right CDG\+module over $\cB^\cu$ endowed with
a natural closed morphism of quasi-coherent right CDG\+modules
$\N^\bu\ot_{\cO_X}\cB^m[-m]\rarrow\N^\bu\ot_{\cO_X}
\cC_X^\cu(\cA_X^{\circ,\rop},\g,\widetilde\g^\circ)$ whose cone is
Positselski-coacyclic in the abelian DG\+category\/ $\bQcohr\cB^\cu$.
\end{lem}

\begin{proof}
 Let us prove part~(a); part~(b) is the opposite version.
 The quasi-coherent left CDG\+module structure on the tensor product
$\cC_X^\cu(\g,\widetilde\g,\cA_X)\ot_{\cO_X}\M^\bu$ is induced by
the quasi-coherent left CDG\+module structure on
the Chevalley--Eilenberg CDG\+module
$\cC_X^\cu(\g,\widetilde\g,\cA_X)$ (which is a CDG\+bimodule over
$\cB^\cu$ and $\cA$, hence also over $\cB^\cu$ and~$\cO_X$;
so the construction of tensor product from
Section~\ref{tensor-product-of-qcoh-cdg-bimods} is applicable).
 The map of quasi-coherent graded left modules
$\cB^m[-m]\ot_{\cO_X}\M^*\rarrow\cC_X^*(\g,\widetilde\g,\cA_X)
\ot_{\cO_X}\M^*$ over the quasi-coherent graded algebra $\cB^*$ is
obtained by taking the tensor product of the graded quasi-coherent
sheaf $\M^*$ with the natural map of quasi-coherent graded
quasi-modules $\cB^m[-m]\rarrow\cC_X^*(\g,\widetilde\g,\cA_X)=
\cB^*\ot_{\cO_X}\cA$ over~$X$.
 The latter map is the tensor product of the direct summand/degree~$m$
component inclusion $\cB^m[-m]\rarrow\cB^*$ and the natural injective
morphism of quasi-coherent graded quasi-modules $\cO_X\rarrow\cA$.
 One can easily check that the resulting map is a closed morphism
of quasi-coherent left CDG\+modules $\cB^m[-m]\ot_{\cO_X}\M^\bu\rarrow
\cC_X^\cu(\g,\widetilde\g,\cA_X)\ot_{\cO_X}\M^\bu$.

 The closed morphism of quasi-coherent CDG\+modules
$\cB^m[-m]\ot_{\cO_X}\M^\bu\rarrow\cC_X^\cu(\g,\widetilde\g,\cA_X)
\ot_{\cO_X}\M^\bu$ is injective, so showing that its cone is coacyclic
in $\cB^\cu\bQcoh$ is equivalent to checking that its cokernel is
coacyclic in $\cB^\cu\bQcoh$.
 Similarly to the proof of
Lemma~\ref{trivial-qcoh-koszul-coresolution-lemma}, we use
Lemma~\ref{filtered-by-coacyclic-is-coacyclic-lemma}.
 Let $F$ be the natural increasing filtration on the quasi-coherent
quasi-algebra $\cA=\cA_X(\g,\widetilde\g)$ defined in
Section~\ref{enveloping-algebra-subsecn}.
 Define a natural increasing filtration of the quasi-coherent graded
algebra $\cB^*$ by the rule $F_n\cB^*=\bigoplus_{j\ge -n}\cB^j$.
 Let the increasing filtration $F$ on the tensor product
$\cC_X^*(\g,\widetilde\g,\cA_X)=\cB^*\ot_{\cO_X}\cA$ be obtained as
the tensor product of the filtrations $F$ on $\cB^*$ and~$\cA$.
 Finally, denote also by $F$ the increasing filtration on
the tensor product $\cC_X^\cu(\g,\widetilde\g,\cA_X)\ot_{\cO_X}\M^\bu$
induced by the increasing filtration $F$ on the Chevalley--Eilenberg
quasi-coherent CDG\+module $\cC_X^\cu(\g,\widetilde\g,\cA_X)$.

 One can easily check that $F$ is a filtration of
$\cC_X^\cu(\g,\widetilde\g,\cA_X)\ot_{\cO_X}\M^\bu$ by quasi-coherent
CDG\+submodules over~$\cB^\cu$.
 The associated graded quasi-coherent CDG\+module to the filtration
$F$ on $\cC_X^\cu(\g,\widetilde\g,\cA_X)\ot_{\cO_X}\M^\bu$ is a trivial
left CDG\+module over $\cB^\cu$ corresponding to the complex of
quasi-coherent sheaves
$$
 \gr^F_*\cC_X^\cu(\g,\widetilde\g,\cA_X)\ot_{\cO_X}\M^\bu \simeq
 \bigwedge\nolimits^*_X(\cHom_{\cO_X}(\g,\cO_X))\ot_{\cO_X}
 \Sym_X^*(\g)\ot_{\cO_X}\M^\bu
$$
 Here $\bigwedge^*_X(\cHom_{\cO_X}(\g,\cO_X))\ot_{\cO_X}\Sym_X^*(\g)
\ot_{\cO_X}\M^\bu$ is the tensor product of two complexes of
quasi-coherent sheaves on $X$, viz., the complex $\M^\bu$ and
the dual Koszul complex
$$
 \cB^*\ot_{\cO_X}\gr^F_*\cA=
 \bigwedge\nolimits^*_X(\cHom_{\cO_X}(\g,\cO_X))\ot_{\cO_X}\Sym_X^*(\g)
$$
with its usual dual Koszul differential.

 Now the dual Koszul complex is bigraded; it has the cohomological
grading, corresponding to the cohomological grading on
the Chevalley--Eilenberg CDG\+module $\cC_X^\cu(\g,\widetilde\g,\cA_X)$,
and the internal grading, corresponding to the degrees of
the filtration~$F$.
 The dual Koszul complex is the direct sum of its internal grading
components; the internal grading component of degree~$-m$ is isomorphic
to $\cB^m[-m]$, while the internal grading components of all other
degrees are finite acyclic complexes of finite locally free sheaves
on~$X$.

 Continuing to argue as in the proof of
Lemma~\ref{trivial-qcoh-koszul-coresolution-lemma}, we conclude that
the map $\cB^m[-m]\ot_{\cO_X}\M^\bu\rarrow
\cC_X^\cu(\g,\widetilde\g,\cA_X)\ot_{\cO_X}\M^\bu$ is a closed
isomorphism of quasi-coherent CDG\+modules
$\cB^m[-m]\ot_{\cO_X}\M^\bu\simeq F_{-m}
\bigl(\cC_X^\cu(\g,\widetilde\g,\cA_X)\ot_{\cO_X}\M^\bu\bigr)$,
while the successive quotient CDG\+module
$$
 F_n\bigl(\cC_X^\cu(\g,\widetilde\g,\cA_X)\ot_{\cO_X}\M^\bu\bigr)
 \big/F_{n-1}
 \bigl(\cC_X^\cu(\g,\widetilde\g,\cA_X)\ot_{\cO_X}\M^\bu\bigr)
$$
is Positselski-coacyclic (in fact, absolutely acyclic) in
$\cB^\cu\bQcoh$ for every $n\ge-m+1$.
 By Lemma~\ref{filtered-by-coacyclic-is-coacyclic-lemma}, it follows
that the cone of the closed morphism $\cB^m[-m]\ot_{\cO_X}\M^\bu\rarrow
\cC_X^\cu(\g,\widetilde\g,\cA_X)\ot_{\cO_X}\M^\bu$ is
Positselski-coacyclic in $\cB^\cu\bQcoh$.
\end{proof}

\subsection{Conversion bimodule} \label{conversion-bimodule-subsecn}
 The construction of the tensor product of quasi-coherent
CDG\+bimodules was spelled out in
Section~\ref{tensor-product-of-qcoh-cdg-bimods}.
 In the notation of Section~\ref{dual-koszul-subsecn}, we put
$$
 \cC_X^\bu(\cA_X^{\circ,\rop},\g,\widetilde\g,\cA_X)=
 \cC_X^\cu(\cA_X^{\circ,\rop},\g,\widetilde\g^\circ)
 \ot_{\cB^*}\cC_X^\cu(\g,\widetilde\g,\cA_X).
$$
 So, by construction,
$\cC_X^\bu(\cA_X^{\circ,\rop},\g,\widetilde\g,\cA_X)$ is a complex
of quasi-coherent $\cA^{\circ,\rop}$\+$\cA$\+bi\-mod\-ules on~$X$
(in the sense of Section~\ref{fHom-over-qcoh-quasi-algebra-subsecn}).
 The underlying graded $\cA^{\circ,\rop}$\+$\cA$\+bimodule
$\cC_X^*(\cA_X^{\circ,\rop},\g,\widetilde\g,\cA_X)$ of the complex
$\cC_X^\bu(\cA_X^{\circ,\rop},\g,\widetilde\g,\cA_X)$ is
$$
 \cC_X^*(\cA_X^{\circ,\rop},\g,\widetilde\g,\cA_X)=
 \cA_X^{\circ,\rop}\ot_{\cO_X}\cB^*\ot_{\cO_X}\cA.
$$

\begin{lem} \label{conversion-functor-for-complexes-lemma}
 Let $X$ be a scheme and $(\g,\widetilde\g)$ be a quasi-coherent
twisted Lie algebroid over~$X$.
 Assume that the quasi-coherent sheaf\/~$\g$ on $X$ is finite locally
free of constant rank~$m$.
 Then \par
\textup{(a)} for any complex of quasi-coherent right
$\cA^{\circ,\rop}$\+modules $\N^\bu$, there is a natural
quasi-isomorphism of complexes of quasi-coherent quasi-modules on~$X$
$$
 \N^\bu\ot_{\cO_X}\cB^m[-m] \lrarrow \N^\bu\ot_{\cA^{\circ,\rop}}
 \cC_X^\bu(\cA_X^{\circ,\rop},\g,\widetilde\g,\cA_X);
$$ \par
\textup{(b)} for any complex of quasi-coherent left
$\cA$\+modules $\M^\bu$, there is a natural quasi-isomorphism 
of complexes of quasi-coherent quasi-modules on~$X$
$$
 \cB^m[-m]\ot_{\cO_X}\M^\bu \lrarrow
 \cC_X^\bu(\cA_X^{\circ,\rop},\g,\widetilde\g,\cA_X)
 \ot_\cA\M^\bu.
$$
\end{lem}

\begin{proof}
 This is our version of~\cite[Lemma~9.8]{Prel}.
 Let us explain part~(a); part~(b) is opposite.
 Applying Lemma~\ref{trivial-qcoh-dual-koszul-coresolution}(a) to
the complex of quasi-coherent sheaves $\M^\bu=\cO_X$, we obtain
a natural closed morphism of quasi-coherent left CDG\+modules
$\cB^m[-m]\rarrow\cC_X^\cu(\g,\widetilde\g,\cA_X)$ over $\cB^\cu$
whose cone is coacyclic in $\cB^\cu\bQcoh$.
 In fact, $\cB^m[-m]\rarrow\cC_X^\cu(\g,\widetilde\g,\cA_X)$ is
a morphism of quasi-coherent CDG\+bimodules over $\cB^\cu$
and~$\cO_X$.
 Following the proof of
Lemma~\ref{trivial-qcoh-dual-koszul-coresolution}, one can see that
the cone of the morphism
$\cB^m[-m]\rarrow\cC_X^\cu(\g,\widetilde\g,\cA_X)$ is actually
Positselski-coacyclic \emph{in the DG\+category of quasi-coherent
CDG\+bimodules over $\cB^\cu$ and $\cO_X$ that are flat as
graded left and right $\cO_X$\+modules}.

 On the other hand,
$\cC_X^\cu(\cA_X^{\circ,\rop},\g,\widetilde\g^\circ)$ is
a quasi-coherent CDG\+bimodule over $\cA^{\circ,\rop}$ and $\cB^\cu$
whose underlying quasi-coherent graded bimodule
$\cC_X^*(\cA_X^{\circ,\rop},\g,\widetilde\g^\circ)$ is flat as
a quasi-coherent graded right $\cB^*$\+module.
 Furthermore, the tensor product functor
$\cC_X^*(\cA_X^{\circ,\rop},\g,\widetilde\g^\circ)\ot_{\cB^*}{-}$
takes quasi-coherent graded left $\cB^*$\+modules that are flat
as graded left $\cO_X$\+modules to flat graded quasi-coherent
left $\cA^{\circ,\rop}$\+modules.
 Applying the functor
$\cC_X^\cu(\cA_X^{\circ,\rop},\g,\widetilde\g^\circ)\ot_{\cB^*}{-}$
to the morphism $\cB^m[-m]\rarrow\cC_X^\cu(\g,\widetilde\g,\cA_X)$,
we obtain a morphism $\cA^{\circ,\rop}\ot_{\cO_X}\cB^m[-m]=
\cC_X^\cu(\cA_X^{\circ,\rop},\g,\widetilde\g^\circ)\ot_{\cB*}\cB^m[-m]
\rarrow\cC_X^\bu(\cA_X^{\circ,\rop},\g,\widetilde\g,\cA_X)$ of complexes
of quasi-coherent $\cA^{\circ,\rop}$\+$\cO_X$\+bimodules.
 It follows that the cone of the latter morphism is
Positselski-coacyclic in the exact category of quasi-coherent
$\cA^{\circ,\rop}$\+$\cO_X$\+bimodules that are flat as quasi-coherent
left $\cA^{\circ,\rop}$\+modules and as quasi-coherent right
$\cO_X$\+modules.

 Thus the tensor product functor $\N^\bu\ot_{\cA^{\circ,\rop}}{-}$
preserves coacyclicity of the latter cone.
 So the cone of the resulting morphism
$\N^\bu\ot_{\cO_X}\cB^m[-m]=\N^\bu\ot_{\cA^{\circ,\rop}}
\cA^{\circ,\rop}\ot_{\cO_X}\cB^m[-m]\rarrow\N^\bu\ot_{\cA^{\circ,\rop}}
\cC_X^\bu(\cA_X^{\circ,\rop},\g,\widetilde\g,\cA_X)$ is
not only acyclic but even coacyclic as a complex of quasi-coherent
quasi-modules over~$X$.
\end{proof}

 Keeping the assumption that $\g$~is a finite locally free sheaf of
constant rank~$m$ on $X$, we are interested in the top cohomology
bimodule
$\cE=\cH^m(\cC_X^\bu(\cA_X^{\circ,\rop},\g,\widetilde\g,\cA_X))$
of the finite complex of quasi-coherent
$\cA^{\circ,\rop}$\+$\cA$\+bimodules
$\cC_X^\bu(\cA_X^{\circ,\rop},\g,\widetilde\g,\cA_X)$.
 According to Lemma~\ref{conversion-functor-for-complexes-lemma}
applied to the one-term complexes $\N^\bu=\cA^{\circ,\rop}$ and
$\M^\bu=\cA$, have
$\cH^n(\cC_X^\bu(\cA_X^{\circ,\rop},\g,\widetilde\g,\cA_X))=0$ for all
$n\ne m$, while the quasi-coherent $\cA^{\circ,\rop}$\+$\cA$\+bimodule
$\cE$ is naturally isomorphic to $\cA^{\circ,\rop}\ot_{\cO_X}\cB^m$ as
a quasi-coherent $\cA^{\circ,\rop}$\+$\cO_X$\+bimodule and naturally
isomorphic to $\cB^m\ot_{\cO_X}\cA$ as a quasi-coherent
$\cO_X$\+$\cA$\+bimodule.
 Moreover, it is clear from the constructions that both the isomorphisms
$\cA^{\circ,\rop}\ot_{\cO_X}\cB^m\simeq\cE$ and $\cB^m\ot_{\cO_X}\cA
\simeq\cE$ are induced by one and the same natural morphism of
quasi-coherent quasi-modules $\cB^m\rarrow\cE$ over~$X$.

\subsection{Conversion functor} \label{conversion-functor-subsecn}
 Let $X$ be a scheme, and let $\cA$ and $\cA^\sharp$ be quasi-coherent
quasi-algebras over~$X$.
 Let $\cE$ be a quasi-coherent $\cA^\sharp$\+$\cA$\+bimodule (as
defined in Section~\ref{fHom-over-qcoh-quasi-algebra-subsecn}) and $\M$
be a quasi-coherent left $\cA^\sharp$\+module (in the sense of
Section~\ref{cosheaves-of-A-modules-subsecn}).

 Denote by $\bB$ the base of neighborhoods of zero in $X$
consisting of all affine open subschemes.
 Define a presheaf of left $\cA$\+modules $\cHom_{\cA^\sharp}(\cE,\M)$
on $\bB$ by the rule $\cHom_{\cA^\sharp}(\cE,\M)[U]=
\Hom_{\cA^\sharp(U)}(\cE(U),\M(U))$ for all affine open subschemes
$U\subset X$.
 For any pair of affine open subschemes $V\subset U\subset X$,
the natural isomorphism of $\cA^\sharp(V)$\+$\cA(U)$\+bimodules
$$
 \cE(V)\simeq\cA^\sharp(V)\ot_{\cA^\sharp(U)}\cE(U)
$$
from the alternative proof of
Lemma~\ref{quasi-algebra-co-extension-of-scalars}(b) allows to define
the restriction map $\Hom_{\cA^\sharp(U)}(\cE(U),\M(U))\rarrow
\Hom_{\cA^\sharp(V)}(\cE(V),\M(V))$ of $\cA(U)$\+modules as the map
induced by the restriction map of $\cA^\sharp(U)$\+modules
$\M(U)\rarrow\M(V)$.

 One can check that the presheaf $\cHom_{\cA^\sharp}(\cE,\M)$ is
actually a sheaf on $\bB$, so it extends uniquely to a sheaf of left
$\cA$\+modules $\cHom_{\cA^\sharp}(\cE,\M)$ on the whole scheme~$X$
(see Theorem~\ref{extension-of-co-sheaves-from-topology-base}(a)).
 Basically, the sheaf axiom for affine open coverings of affine open
subschemes of $X$ holds for $\cHom_{\cA^\sharp}(\cE,\M)$ because
it holds for~$\M$.

 For any quasi-coherent $\cA$\+module $\cL$ over $X$, there is
an obvious adjunction isomorphism
\begin{equation} \label{cHom-A-tensor-adjunction}
 \Hom_\cA(\cL,\cHom_{\cA^\sharp}(\cE,\M))
 \simeq\Hom_{\cA^\sharp}(\cE\ot_\cA\cL,\>\M).
\end{equation}
 Here $\Hom_{\cA^\sharp}$ denotes the abelian group of morphisms in
the category of quasi-coherent left $\cA^\sharp$\+modules, while
$\Hom_\cA$ stands for the abelian group of morphisms in the category
of sheaves of left $\cA$\+modules on~$X$.

\begin{lem} \label{cHom-over-quasi-algebra-quasi-coherent}
 Assume that there exists a finite locally free quasi-coherent
sheaf $\cT$ on $X$ such that the underlying quasi-coherent
$\cA^\sharp$\+$\cO_X$\+bimodule of $\cE$ is isomorphic to a direct
summand of the quasi-coherent $\cA^\sharp$\+$\cO_X$\+bimodule
$\cA^\sharp\ot_{\cO_X}\cT$.
 Then the sheaf of left $\cA$\+modules $\cHom_{\cA^\sharp}(\cE,\M)$
on $X$ is quasi-coherent.
\end{lem}

\begin{proof}
 Following the definition in
Section~\ref{cosheaves-of-A-modules-subsecn}, saying that a sheaf
of left $\cA$\+modules is quasi-coherent means simply that its
underlying sheaf of $\cO_X$\+modules is quasi-coherent.
 By construction, the underlying sheaf of $\cO_X$\+modules of
$\cHom_{\cA^\sharp}(\cE,\M)$ only depends on the underlying
quasi-coherent $\cA^\sharp$\+$\cO_X$\+bimodule of~$\cE$.
 Now we have $\cHom_{\cA^\sharp}(\cA^\sharp\ot_{\cO_X}\cT,\>\M)
\simeq\cHom_{\cO_X}(\cT,\M)$, which is a quasi-coherent sheaf on $X$
whenever $\cT$ is a finite locally free sheaf (see
Section~\ref{loc-free-sheaves-subsecn}).
\end{proof}

 An \emph{invertible sheaf} $\cT$ on a scheme $X$ can be defined as
a finite locally free sheaf of constant rank~$1$, or equivalently,
as a finite locally free sheaf such that the natural map
$\cHom_{\cO_X}(\cT,\cO_X)\ot_{\cO_X}\cT\rarrow\cO_X$ is an isomorphism
of quasi-coherent sheaves.

\begin{lem} \label{abstract-conversion-lemma}
 Let $X$ be a scheme, and let $\cA$ and $\cA^\sharp$ be quasi-coherent
quasi-algebras over~$X$.
 Let $\cT$ be an invertible sheaf on $X$ and $\cE$ be a quasi-coherent
$\cA^\sharp$\+$\cA$\+bimodule endowed with a morphism of quasi-coherent
quasi-modules $\cT\rarrow\cE$.
 Assume that the induced morphism of quasi-coherent
$\cA^\sharp$\+$\cO_X$\+bimodules $\cA^\sharp\ot_{\cO_X}\cT\rarrow\cE$
and the induced morphism of quasi-coherent $\cO_X$\+$\cA$\+bimodules
$\cT\ot_{\cO_X}\cA\rarrow\cE$ are isomorphisms.
 Then there is an equivalence of abelian categories
$$
 \cE\ot_\cA{-}\,\:\cA\Qcoh\simeq\cA^\sharp\Qcoh\,:\!
 \cHom_{\cA^\sharp}(\cE,{-}).
$$
\end{lem}

\begin{proof}
 The functor $\cE\ot_\cA{-}\,\:\cA\Qcoh\rarrow\cA^\sharp\Qcoh$ was
constructed in Section~\ref{tensor-product-of-qcoh-cdg-bimods}.
 The functor $\cHom_{\cA^\sharp}(\cE,{-})\:\cA^\sharp\Qcoh\rarrow
\cA\Qcoh$ is well-defined by
Lemma~\ref{cHom-over-quasi-algebra-quasi-coherent}.
 By formula~\eqref{cHom-A-tensor-adjunction}, the functor
$\cE\ot_{\cA}{-}$ is left adjoint to the functor
$\cHom_{\cA^\sharp}(\cE,{-})$.
 Furthermore, it follows from the assumptions of the lemma that
there are commutative square diagrams of adjoint functors and
forgetful functors
\begin{equation} \label{abstract-conversion-adjunctions}
\begin{gathered}
 \xymatrix{
  \cA\Qcoh\ar[rr]^-{\cE\ot_{\cA}{-}} \ar[d]
  && \cA^\sharp\Qcoh \ar[d] \\
  X\Qcoh \ar[rr]^{\cT\ot_{\cO_X}{-}} && X\Qcoh
 }
 \qquad\qquad
 \xymatrix{
  \cA\Qcoh \ar[d]
  && \cA^\sharp\Qcoh \ar[ll]_-{\cHom_{\cA^\sharp}(\cE,{-})} \ar[d] \\
  X\Qcoh && X\Qcoh \ar[ll]_{\cHom_{\cO_X}(\cT,{-})}
 }
\end{gathered}
\end{equation}
 The forgetful functors take the adjunction morphisms of the adjoint
pair in the upper lines to the adjunction morphisms of the adjoint
pair in the lower lines of the diagrams.
 The adjoint functors $\cT\ot_{\cO_X}{-}\,\:X\Qcoh\rarrow X\Qcoh$
and $\cHom_{\cO_X}(\cT,{-})\:X\Qcoh\rarrow X\Qcoh$ are mutually inverse
equivalences, so their adjunction morphisms are isomorphisms.
 As the forgetful functors between abelian categories
$\cA\Qcoh\rarrow X\Qcoh$ and $\cA^\sharp\Qcoh\rarrow X\Qcoh$ are exact
and faithful, hence conservative, it follows that the adjunction
morphisms of the adjoint functors $\cE\ot_{\cA}{-}$ and
$\cHom_{\cA^\sharp}(\cE,{-})$ are isomorphisms, too.
\end{proof}

\begin{cor} \label{conversion-functor-corollary}
 Let $X$ be a scheme and $(\g,\widetilde\g)$ be a quasi-coherent
twisted Lie algebroid over~$X$.
 Assume that the quasi-coherent sheaf\/~$\g$ on $X$ is finite locally
free of constant rank~$m$.
 Consider the quasi-coherent universal enveloping quasi-algebras
$\cA=\cA_X(\g,\widetilde\g)$ and
$\cA^\circ=\cA_X(\g,\widetilde\g^\circ)$, the invertible sheaf
$\cB^m=\bigwedge_X^m(\cHom_{\cO_X}(\g,\cO_X))$, and the quasi-coherent
$\cA^{\circ,\rop}$\+$\cA$\+bimodule $\cE$ constructed in
Section~\ref{conversion-bimodule-subsecn}.
 Put $\cA^\sharp=\cA^{\circ,\rop}$.
 Then there is a natural equivalence of abelian categories
$\cA\Qcoh\simeq\Qcohr\cA^\circ$ forming a commutative square diagram
of exact, faithful forgetful functors and abelian category equivalences
\begin{equation} \label{conversion-functor-diagram}
\begin{gathered}
 \xymatrix{
  \text{\llap{$\cE\ot_\cA{-}\,\:$}} \cA\Qcoh
  \ar@{=}[rr] \ar[d]
  && \cA^\sharp\Qcoh \text{\rlap{$\,\,:\!\cHom_{\cA^\sharp}(\cE,{-})$}}
  \ar[d] \\
  \text{\llap{$\cB^m\ot_{\cO_X}{-}\,\:$}} X\Qcoh \ar@{=}[rr] 
  && X\Qcoh \text{\rlap{$\,\,:\!\cHom_{\cO_X}(\cB^m,{-})$}}
 }
\end{gathered}
\end{equation}
\end{cor}

\begin{proof}
 It suffices to compare Lemma~\ref{abstract-conversion-lemma} with
the discussion at the end of
Section~\ref{conversion-bimodule-subsecn}.
 Notice that we have $\cA^\sharp\Qcoh=\cA^{\circ,\rop}\Qcoh\simeq
\Qcohr\cA^\circ$ according to the discussion in
Section~\ref{opposite-cdg-rings-and-twisted-lie-subsecn}
(by construction, this equivalence commutes with the forgetful
functors $\cA^\sharp\Qcoh\rarrow X\Qcoh$ and $\Qcohr\cA^\circ
\rarrow X\Qcoh$).
\end{proof}

 The result of Corollary~\ref{conversion-functor-corollary} can be
rephrased by saying that, for any quasi-coherent left $\cA$\+module
$\M$, the quasi-coherent sheaf $\cB^m\ot_{\cO_X}\M$ on $X$ has
a natural structure of quasi-coherent right $\cA^\circ$\+module.
 Conversely, for any quasi-coherent right $\cA^\circ$\+module $\N$,
the quasi-coherent sheaf $\cHom_{\cO_X}(\cB^m,\N)$ on $X$ has
a natural structure of quasi-coherent left $\cA$\+module.

\subsection{Comparison of left and right quasi-coherent module Koszul
duality functors} \label{comparison-of-left-and-right-subsecn}
 We start with a partial left version of
Lemma~\ref{koszul-duality-dg-functors-right-co-side}, providing one
Koszul duality DG\+functor for quasi-coherent left $\cA$\+modules and
quasi-coherent left CDG\+modules over~$\cB^\cu$.

\begin{lem} \label{one-koszul-duality-dg-functor-left-co-side}
 Let $X$ be a scheme and $(\g,\widetilde\g)$ be a quasi-coherent
twisted Lie algebroid over $X$ such that the quasi-coherent sheaf\/~$\g$
on $X$ is finite locally free.
 Let $\cA=\cA_X(\g,\widetilde\g)$ be the twisted universal enveloping
quasi-coherent quasi-algebra and $\cB^\cu=\cC_X^\cu(\g,\widetilde\g)$
be the Chevalley--Eilenberg quasi-coherent CDG\+quasi-algebra
of~$(\g,\widetilde\g)$.
 Then there is a natural exact DG\+functor
$$
 \cB^*\ot_{\cO_X}\,\:\bCom(\cA\Qcoh)\lrarrow
 \cB^\cu\bQcoh.
$$
 To be more precise, the differential on the tensor product over
$\cO_X$ is defined in terms of the identification
$$
 \cB^*\ot_{\cO_X}{-}\,=\cC_X^\cu(\g,\widetilde\g,\cA_X)\ot_{\cA}{-}\,\:
 \bCom(\cA\Qcoh)\lrarrow\cB^\cu\bQcoh.
$$
\end{lem}

\begin{proof}
 Recall once again that the Chevalley--Eilenberg quasi-coherent
CDG\+module $\cC_X^\cu(\g,\widetilde\g,\cA_X)$ constructed in
Section~\ref{chevalley-eilenberg-cdg-modules-subsecn} is
a quasi-coherent CDG\+bimodule over $\cB^\cu$ and $\cA$ with
the underlying quasi-coherent graded bimodule
$\cC_X^*(\g,\widetilde\g,\cA_X)=\cB^*\ot_{\cO_X}\nobreak\cA$.
 The construction of the tensor product of quasi-coherent
CDG\+(bi)modules from Section~\ref{tensor-product-of-qcoh-cdg-bimods}
defines the quasi-coherent left CDG\+module structure on the tensor
product of any complex of quasi-coherent left $\cA$\+modules with this
quasi-coherent CDG\+bimodule, providing the desired DG\+functor.
\end{proof}

 The following proposition is our version of~\cite[Lemma~9.11]{Prel}.

\begin{prop} \label{left-right-co-side-koszul-functors-prop}
 Let $X$ be a scheme and $(\g,\widetilde\g)$ be a quasi-coherent
twisted Lie algebroid over $X$ such that the quasi-coherent sheaf\/~$\g$
on $X$ is finite locally free of constant rank~$m$.
 Let $\cA=\cA_X(\g,\widetilde\g)$ and
$\cA^\circ=\cA_X(\g,\widetilde\g^\circ)$ be the two related twisted
universal enveloping quasi-coherent quasi-algebras, and let
$\cB^\cu=\cC_X^\cu(\g,\widetilde\g)$ and
$\cB^\circ{}^\cu=\cC_X^\cu(\g,\widetilde\g^\circ)$ be the two related
Chevalley--Eilenberg quasi-coherent CDG\+quasi-algebras over~$X$.
 Then there is a natural commutative square diagram of exact
DG\+functors and abelian DG\+category equivalences
\begin{equation} \label{left-right-co-side-koszul-functors-diagram}
\begin{gathered}
 \xymatrix{
  \bCom(\Qcohr\cA^\circ)
  \ar[rr]^-{{-}\ot_{\cO_X}\bigwedge\nolimits_X^*(\g)} \ar@{=}[d]
  && \bQcohr\cB^\circ{}^\cu \ar@{=}[d] \\
  \bCom(\cA\Qcoh)
  \ar[rr]^-{\cB^*\ot_{\cO_X}{-}} && \cB^\cu\bQcoh
 }
\end{gathered}
\end{equation}
 Here the upper horizontal DG\+functor is provided by
Lemma~\ref{koszul-duality-dg-functors-right-co-side} applied to
the quasi-coherent twisted Lie algebroid~$(\g,\widetilde\g^\circ)$
over~$X$.
 The lower horizontal DG\+functor is provided by
Lemma~\ref{one-koszul-duality-dg-functor-left-co-side} applied to
the quasi-coherent twisted Lie algebroid~$(\g,\widetilde\g)$.
 The equivalence of the abelian DG\+categories of complexes in
the leftmost column is the composition of the DG\+equivalence\/
$\Com(\cA\Qcoh)\rarrow\Com(\Qcohr\cA^\circ)$ induced by the equivalence
of abelian categories $\cA\Qcoh\simeq\Qcohr\cA^\circ$ from
Corollary~\ref{conversion-functor-corollary} with the cohomological
degree shift\/
$[-m]\:\Com(\Qcohr\cA^\circ)\rarrow\Com(\Qcohr\cA^\circ)$.
 The isomorphism of abelian DG\+categories in the rightmost column
is induced by the natural (identity) isomorphism of quasi-coherent
CDG\+quasi-algebras $\cB^{\circ,\rop}{}^\cu\simeq\cB^\cu$,
as per the discussions in
Sections~\ref{opposite-cdg-rings-and-twisted-lie-subsecn}
and~\ref{dual-koszul-subsecn}.
\end{prop}

\begin{proof}
 For any given complex of quasi-coherent left $\cA$\+modules $\M^\bu$,
we need to construct a closed isomorphism
$$
 (\cB^m[-m]\ot_{\cO_X}\M^\bu)\ot_{\cO_X}
 \bigwedge\nolimits_X^*(\g) \simeq
 \cB^*\ot_{\cO_X}\M^\bu
$$
of quasi-coherent left CDG\+modules over~$\cB^\cu$.
 Of course, the desired map is induced by the natural isomorphism of
graded quasi-coherent sheaves $\cB^m[-m]\ot_{\cO_X}\bigwedge_X^*(\g)
\simeq\cB^*$ on~$X$.
 The following construction of the same map makes it clear that it is
a closed isomorphism of quasi-coherent CDG\+modules.

 Let us rewrite the left-hand side and the right-hand side of
the desired closed isomorphism as follows using the notation of
Section~\ref{dual-koszul-subsecn}.
 We need to construct a closed isomorphism of quasi-coherent
left CDG\+modules
$$
 \cC^X_\cu(\g,\widetilde\g^\circ,\cA^{\circ,\rop}_X)
 \ot_{\cA^{\circ,\rop}}\cE\ot_\cA\M^\bu\simeq
 \cC_X^\cu(\g,\widetilde\g,\cA_X)\ot_\cA\M^\bu
$$
over~$\cB^\cu$.
 For this purpose, it suffices to construct a closed isomorphism of
quasi-coherent CDG\+bimodules
\begin{equation} \label{left-right-comparison-cdg-bimodule-isom}
 \cC^X_\cu(\g,\widetilde\g^\circ,\cA^{\circ,\rop}_X)
 \ot_{\cA^{\circ,\rop}}\cE\simeq
 \cC_X^\cu(\g,\widetilde\g,\cA_X)
\end{equation}
over $\cB^\cu$ and~$\cA$ (then it would remain to apply the DG\+functor
of tensor product ${-}\ot_\cA\M^\bu$).

 Recall that, by the definition, the quasi-coherent
$\cA^{\circ,\rop}$\+$\cA$\+bimodule $\cE$ is constructed as
$\cE=\cH^m(\cC_X^\bu(\cA_X^{\circ,\rop},\g,\widetilde\g,\cA_X))$,
where $\cC_X^\bu(\cA_X^{\circ,\rop},\g,\widetilde\g,\cA_X)=
\cC_X^\cu(\cA_X^{\circ,\rop},\g,\widetilde\g^\circ)
\ot_{\cB^*}\cC_X^\cu(\g,\widetilde\g,\cA_X)$.
 Moreover, there is a natural quasi-isomorphism of finite complexes of
$\cA^{\circ,\rop}$\+$\cA$\+bimodules
$\cC_X^\bu(\cA_X^{\circ,\rop},\g,\widetilde\g,\cA_X)\rarrow\cE[-m]$.
 Applying the DG\+functor
$\cC_\cu^X(\g,\widetilde\g^\circ,\cA_X^{\circ,\rop})
\ot_{\cA^{\circ,\rop}}{-}$, we obtain a closed morphism
\begin{multline} \label{first-composed-CDG-bimodule-morphism}
 \cC_\cu^X(\g,\widetilde\g^\circ,\cA_X^{\circ,\rop})
 \ot_{\cA^{\circ,\rop}}
 \cC_X^\cu(\cA_X^{\circ,\rop},\g,\widetilde\g^\circ)
 \ot_{\cB^*}\cC_X^\cu(\g,\widetilde\g,\cA_X) \\
 \lrarrow\cC_\cu^X(\g,\widetilde\g^\circ,\cA_X^{\circ,\rop})
 \ot_{\cA^{\circ,\rop}}\cE[-m]
\end{multline}
of quasi-coherent CDG\+bimodules over $\cB^\cu$ and~$\cA$.

 For any quasi-coherent left CDG\+module $\N^\cu$ over $\cB^\cu$,
the left-right opposite version of
Lemma~\ref{koszul-duality-dg-functors-right-co-side} for
$(\g,\widetilde\g^\circ)$ provides a natural adjunction morphism
$$
 \N^\cu\lrarrow
 \cC_\cu^X(\g,\widetilde\g^\circ,\cA_X^{\circ,\rop})
 \ot_{\cA^{\circ,\rop}}
 \cC_X^\cu(\cA_X^{\circ,\rop},\g,\widetilde\g^\circ)
 \ot_{\cB^*}\N^\cu,
$$
which is a closed morphism of quasi-coherent left CDG\+modules
over~$\cB^\cu$.
 In particular, we are interested in the adjunction morphism
\begin{equation} \label{second-composed-CDG-bimodule-morphism}
 \cC_X^\cu(\g,\widetilde\g,\cA_X)\lrarrow
 \cC_\cu^X(\g,\widetilde\g^\circ,\cA_X^{\circ,\rop})
 \ot_{\cA^{\circ,\rop}}
 \cC_X^\cu(\cA_X^{\circ,\rop},\g,\widetilde\g^\circ)
 \ot_{\cB^*}\cC_X^\cu(\g,\widetilde\g,\cA_X),
\end{equation}
which is a closed morphism of quasi-coherent CDG\+bimodules over
$\cB^\cu$ and~$\cA$.
 The composition of~\eqref{second-composed-CDG-bimodule-morphism}
and~\eqref{first-composed-CDG-bimodule-morphism} provides
the desired closed isomorphism of quasi-coherent
CDG\+bimodules~\eqref{left-right-comparison-cdg-bimodule-isom}.
\end{proof}

 The Positselski semicoderived category of quasi-coherent left modules
$\sD^\si(\cA\Qcoh)$ for a quasi-coherent quasi-algebra $\cA$ over
a scheme $X$ was defined in
Section~\ref{semiderived-quasi-coherent-subsecn}.
 The Positselski semicoderived category of quasi-coherent right modules
$\sD^\si(\Qcohr\cA)$ is defined similarly.
 The same applies to the Becker semicoderived categories
$\sD^\bsi(\cA\Qcoh)$ and $\sD^\bsi(\Qcohr\cA)$.

\begin{cor} \label{left-right-co-side-semiderived-koszul-duality}
 In the assumptions and notation of
Proposition~\ref{left-right-co-side-koszul-functors-prop},
there is a natural commutative square diagram of triangulated
equivalences
\begin{equation} \label{left-right-co-side-semiderived-koszul-diagram}
\begin{gathered}
 \qquad\quad\xymatrix{
  \text{\llap{${-}\ot_{\cO_X}\bigwedge\nolimits_X^*(\g)\:$}}
  \sD^\si(\Qcohr\cA^\circ) \ar@{=}[r] \ar@{=}[d]
  & \sD^\co(\bQcohr\cB^\circ{}^\cu)
  \text{\rlap{$\,{}:{-}\ot_{\cO_X}\cA^\circ$}} \ar@{=}[d] \\
  \text{\llap{$\cB^*\ot_{\cO_X}{-}\,\:$}}
  \sD^\si(\cA\Qcoh) \ar@<-2pt>[r]
  & \sD^\co(\cB^\cu\bQcoh) \ar@<-2pt>@{-}[l]
 }
\end{gathered}
\end{equation}
 There is also a similar commutative square diagram of triangulated
equivalences of the Becker semicoderived and coderived categories\/
$\sD^\bsi$ and\/~$\sD^\co$.
\end{cor}

\begin{proof}
 All the triangulated functors are induced by the respective
DG\+functors from Lemma~\ref{koszul-duality-dg-functors-right-co-side}
and Proposition~\ref{left-right-co-side-koszul-functors-prop}.
 Then it is clear that the vertical triangulated equivalences hold
(for the triangulated equivalence in the leftmost column,
one needs to take into account the commutative square diagram from
Corollary~\ref{conversion-functor-corollary}).

 The functors shown by the upper horizontal line
in~\eqref{left-right-co-side-semiderived-koszul-diagram} are
mutually inverse triangulated equivalences by
Theorem~\ref{semiderived-koszul-duality-right-co-side}.
 It is clear from the commutativity of the square
diagram~\eqref{left-right-co-side-koszul-functors-diagram}
that the lower horizontal triangulated functor
in~\eqref{left-right-co-side-semiderived-koszul-diagram}
is well-defined (because the upper horizontal functor pointing
rightwards is well-defined by
Theorem~\ref{semiderived-koszul-duality-right-co-side}).
 Alternatively, one can check that the lower horizontal triangulated
functor is well-defined (i.~e., takes semicoacyclic complexes to
coacyclic CDG\+modules) by a direct argument similar to the one in
the proof of Theorem~\ref{semiderived-koszul-duality-right-co-side}.

 The square
diagram~\eqref{left-right-co-side-semiderived-koszul-diagram} (with
the upper horizontal functor pointing rightwards) is commutative,
because so is~\eqref{left-right-co-side-koszul-functors-diagram}.
 In view of the commutativity
of~\eqref{left-right-co-side-semiderived-koszul-diagram}, it follows
that the lower horizontal triangulated functor is a triangulated
equivalence, too.
\end{proof}

 The definition of the reduced (Positselski) coderived categories
appearing in the next corollary can be found in
Section~\ref{reduced-coderived-of-qcoh-subsecn}; see also
Section~\ref{reduced-koszul-duality-co-side-subsecn}.

\begin{cor} \label{left-right-co-side-reduced-koszul-duality}
 In the assumptions and notation of
Proposition~\ref{left-right-co-side-koszul-functors-prop},
there is a natural commutative square diagram of triangulated
equivalences
\begin{equation} \label{left-right-co-side-reduced-koszul-diagram}
\begin{gathered}
 \qquad\quad\xymatrix{
  \text{\llap{${-}\ot_{\cO_X}\bigwedge\nolimits_X^*(\g)\:$}}
  \sD(\Qcohr\cA^\circ) \ar@{=}[r] \ar@{=}[d]
  & \sD^\co_{X\red}(\bQcohr\cB^\circ{}^\cu)
  \text{\rlap{$\,{}:{-}\ot_{\cO_X}\cA^\circ$}} \ar@{=}[d] \\
  \text{\llap{$\cB^*\ot_{\cO_X}{-}\,\:$}}
  \sD(\cA\Qcoh) \ar@<-2pt>[r]
  & \sD^\co_{X\red}(\cB^\cu\bQcoh) \ar@<-2pt>@{-}[l]
 }
\end{gathered}
\end{equation}
 Moreover, there is a natural commutative square diagram of
triangulated equivalences and triangulated Verdier quotient functors
\begin{equation} \label{reduced-koszul-duality-left-co-side-diagram}
\begin{gathered}
 \xymatrix{
  \sD^\si(\cA\Qcoh)
  \ar@<2pt>[rr]^-{\cB^*\ot_{\cO_X}{-}} \ar@{->>}[d]
  && \sD^\co(\cB^\cu\bQcoh)
  \ar@{->>}[d] \ar@<2pt>@{-}[ll] \\
  \sD(\cA\Qcoh)
  \ar@<2pt>[rr]^-{\cB^*\ot_{\cO_X}{-}}
  && \sD^\co_{X\red}(\cB^\cu\bQcoh) \ar@<2pt>@{-}[ll]
 }
\end{gathered}
\end{equation}
\end{cor}

\begin{proof}
 All the triangulated functors
in~\eqref{left-right-co-side-reduced-koszul-diagram} are induced by
the respective DG\+func\-tors from
Lemma~\ref{koszul-duality-dg-functors-right-co-side} and
Proposition~\ref{left-right-co-side-koszul-functors-prop}.
 Then it is clear that the triangulated equivalences shown by
the vertical double lines hold.

 The functors shown by the upper horizontal line
in~\eqref{left-right-co-side-semiderived-koszul-diagram} are
mutually inverse triangulated equivalences by
Theorem~\ref{reduced-koszul-duality-right-co-side}.
 It is clear from the commutativity of the square
diagram~\eqref{left-right-co-side-koszul-functors-diagram}
that the lower horizontal triangulated functor
in~\eqref{left-right-co-side-reduced-koszul-diagram} is well-defined
(because the upper horizontal functor pointing rightwards is
well-defined by Theorem~\ref{reduced-koszul-duality-right-co-side}).
 Alternatively, one can check that the lower horizontal triangulated
functor is well-defined (i.~e., takes acyclic complexes to
reduced-coacyclic CDG\+modules) by a direct argument similar to the one
in the proof of Theorem~\ref{reduced-koszul-duality-right-co-side}.

 The square
diagram~\eqref{left-right-co-side-reduced-koszul-diagram} (with
the upper horizontal functor pointing rightwards) is commutative,
because so is~\eqref{left-right-co-side-koszul-functors-diagram}.
 In view of the commutativity
of~\eqref{left-right-co-side-reduced-koszul-diagram}, it follows that
the lower horizontal triangulated functor is a triangulated
equivalence, too.

 In diagram~\eqref{reduced-koszul-duality-left-co-side-diagram},
the triangulated equivalence in the upper line is provided by
the lower line of
diagram~\eqref{left-right-co-side-semiderived-koszul-diagram}
from Corollary~\ref{left-right-co-side-semiderived-koszul-duality}.
 The triangulated equivalence in the lower line
of~\eqref{reduced-koszul-duality-left-co-side-diagram} is provided
by the lower line of
diagram~\eqref{left-right-co-side-reduced-koszul-diagram}.
 The functors shown by vertical arrows with double heads are
the obvious triangulated Verdier quotient functors.
 The commutativity
of~\eqref{reduced-koszul-duality-left-co-side-diagram} is also
obvious from the constructions.
\end{proof}

\subsection{Internal $\cHom$ of quasi-coherent CDG-modules}
\label{internal-cHom-of-quasi-coherent-cdg-modules-subsecn}
 Let $X$ be a scheme.
 A graded quasi-coherent sheaf $\cT^*=\bigoplus_{i\in\boZ}\cT^i$ on $X$
is said to be \emph{finite locally free} if, for every affine open
subscheme $U\subset X$, the graded $\cO_X(U)$\+module $\cT^*(U)$ is
finitely generated and projective (as a graded module).
 The latter condition means that $\cT^i(U)=0$ for all but a finite
subset of degrees $i\in\boZ$, and $\cT^i(U)$ is a finitely generated
projective $\cO_X(U)$\+module for all $i\in\boZ$.

 Let $\cA^*$, $\cB^*$, and $\cC^*$ be quasi-coherent graded
quasi-algebras over $X$, and let $\cE^*$ be a quasi-coherent graded
$\cA^*$\+$\cB^*$\+bimodule.
 Let $\M^*$ be a quasi-coherent graded left $\cA^*$\+module and $\N^*$
be a quasi-coherent graded $\cA^*$\+$\cC^*$\+bimodule.
 Then the graded version of the construction of
Section~\ref{conversion-functor-subsecn} defines a sheaf of graded
$\cB^*$\+modules $\cHom^*_{\cA^*}(\cE^*,\M^*)$ on~$X$.
 Furthermore, the same construction defines a sheaf of graded
$\cB^*$\+$\cC^*$\+bimodules $\cHom^*_{\cA^*}(\cE^*,\N^*)$ on~$X$.

\begin{lem} \label{cHom-over-graded-quasi-algebra-quasi-coherent-graded}
 Assume that there exists a finite locally free graded quasi-coherent
sheaf $\cT^*$ on $X$ such that the underlying quasi-coherent graded
$\cA^*$\+$\cO_X$\+bimodule of $\cE^*$ is isomorphic to a direct summand
of the quasi-coherent graded $\cA^*$\+$\cO_X$\+bimodule
$\cA^*\ot_{\cO_X}\cT^*$.
 Then \par
\textup{(a)} the sheaf of graded left $\cB^*$\+modules
$\cHom^*_{\cA^*}(\cE^*,\M^*)$ on $X$ is quasi-coherent; \par
\textup{(b)} the sheaf of graded $\cB^*$\+$\cC^*$\+bimodules
$\cHom^*_{\cA^*}(\cE^*,\N^*)$ on $X$ is a quasi-coherent graded
$\cB^*$\+$\cC^*$\+bimodule.
\end{lem}

\begin{proof}
 Part~(a) is just the graded version of
Lemma~\ref{cHom-over-quasi-algebra-quasi-coherent}, and the proof is
similar.
 In part~(b), we need to show that the underlying sheaf of graded
$\cO_X$\+$\cO_X$\+bimodules of $\cHom^*_{\cA^*}(\cE^*,\N^*)$ is
a quasi-coherent graded quasi-module over~$X$.
 Indeed, we have $\cHom^*_{\cA^*}(\cE^*,\N^*)\simeq
\cHom^*_{\cO_X}(\cT^*,\N^*)\simeq\cHom_{\cO_X}(\cT^*,\cO_X)
\ot_{\cO_X}\N^*$, and it remains to point out that the tensor products
(and finite direct sums) of quasi-coherent quasi-modules are
quasi-coherent quasi-modules by
Lemma~\ref{quasi-coherent-quasi-modules-tensor-product}.
\end{proof}

 Let $\cA^\cu$ and $\cB^\cu$ be quasi-coherent CDG\+quasi-algebras
over $X$, and let $\cE^\cu$ be a quasi-coherent CDG\+bimodule over
$\cA^\cu$ and~$\cB^\cu$.
 Let $\M^\cu$ be a quasi-coherent left CDG\+module over~$\cA^\cu$.
 Assume that the quasi-coherent graded $\cA^*$\+$\cB^*$\+bimodule
$\cE^*$ satisfies the assumption of
Lemma~\ref{cHom-over-graded-quasi-algebra-quasi-coherent-graded}.
 Then the rule
$$
 \cHom_{\cA^*}^\cu(\cE^\cu,\M^\cu)(U)=
 \Hom_{\cA^*(U)}^\cu(\cE^\cu(U),\M^\cu(U)),
$$
with the reference to the definition of the Hom CDG\+module in
Section~\ref{cdg-rings-cdg-modules-subsecn}, produces a differential on
the graded left $\cB^*(U)$\+module $\cHom_{\cA^*}^\cu(\cE^*,\M^*)(U)$
for every affine open subscheme $U\subset X$.
 The collection of all such differentials defines a structure of
quasi-coherent CDG\+module over $\cB^\cu$ on the quasi-coherent
graded left $\cB^*$\+module $\cHom^*_{\cA^*}(\cE^*,\M^*)$.
 We denote the resulting quasi-coherent left CDG\+module over $\cB^\cu$
by $\cHom_{\cA^*}^\cu(\cE^\cu,\M^\cu)$.

 Similarly, let $\cC^\cu$ be another quasi-coherent CDG\+quasi-algebra
over $X$, and let $\N^\cu$ be a quasi-coherent CDG\+bimodule over
$\cA^\cu$ and~$\cC^\cu$.
 Assume that the quasi-coherent graded $\cA^*$\+$\cB^*$\+bimodule
$\cE^*$ satisfies the assumption of
Lemma~\ref{cHom-over-graded-quasi-algebra-quasi-coherent-graded}.
 Then the rule
$$
 \cHom_{\cA^*}^\cu(\cE^\cu,\N^\cu)(U)=
 \Hom_{\cA^*(U)}^\cu(\cE^\cu(U),\N^\cu(U)),
$$
with the reference to the definition of the Hom CDG\+bimodule in
Section~\ref{cdg-rings-cdg-modules-subsecn}, produces a differential
on the graded $\cB^*(U)$\+$\cC^*(U$)\+bimodule
$\cHom_{\cA^*}^\cu(\cE^*,\N^*)(U)$ for every affine open subscheme
$U\subset X$.
 The collection of all such differentials defines a structure of
quasi-coherent CDG\+bimodule over $\cB^\cu$ and $\cC^\cu$ on
the quasi-coherent graded $\cB^*$\+$\cC^*$\+bimodule
$\cHom^*_{\cA^*}(\cE^*,\N^*)$.
 We denote the resulting quasi-coherent CDG\+bimodule over $\cB^\cu$
and $\cC^\cu$ by $\cHom_{\cA^*}^\cu(\cE^\cu,\N^\cu)$.

\begin{lem} \label{cdg-bimodule-fin-gen-proj-Hom-tensor-lemma}
 Let $\cA^\cu$ and $\cB^\cu$ be quasi-coherent CDG\+quasi-algebras
over a scheme~$X$.
 Let $\cE^\cu$ be a quasi-coherent CDG\+bimodule over $\cA^\cu$ and
$\cB^\cu$ such that the graded $\cA^*$\+$\cB^*$\+bimodule $\cE^*$
satisfies the assumption of
Lemma~\ref{cHom-over-graded-quasi-algebra-quasi-coherent-graded},
and let $\M^\cu$ be a quasi-coherent left CDG\+module over~$\cA^\cu$.
 Then there is a natural closed isomorphism
$$
 \cHom_{\cA^*}^\cu(\cE^\cu,\cA^\cu)\ot_{\cA^*}\M^\cu
 \simeq\cHom_{\cA^*}^\cu(\cE^\cu,\M^\cu)
$$
of quasi-coherent left CDG\+modules over~$\cB^\cu$.
\end{lem}

\begin{proof}
 It follows from our assumptions that the graded
$\cA^*(U)$\+$\cB^*(U)$\+bimodule $\cE^*(U)$ is finitely generated
and projective as a graded left $\cA^*(U)$\+module for every affine
open subscheme $U\subset X$.
 Now the collection of natural closed isomorphisms
$$
 \Hom_{\cA^*(U)}^\cu(\cE^\cu(U),\cA^\cu(U))\ot_{\cA^*(U)}\M^\cu(U)
 \simeq\Hom_{\cA^*(U)}^\cu(\cE^\cu(U),\M^\cu(U))
$$
of left CDG\+modules over the CDG\+rings $\cB^\cu(U)$ provides
the desired closed isomorphism of quasi-coherent CDG\+modules.
\end{proof}

\begin{cor} \label{tensor-Hom-koszul-duality-functor-isom-cor}
 Let $X$ be a scheme and $(\g,\widetilde\g)$ be a quasi-coherent
twisted Lie algebroid over $X$ such that the quasi-coherent sheaf\/~$\g$
on $X$ is finite locally free.
 Let $\cA=\cA_X(\g,\widetilde\g)$ be the twisted universal enveloping
quasi-coherent quasi-algebra and $\cB^\cu=\cC_X^\cu(\g,\widetilde\g)$
be the Chevalley--Eilenberg quasi-coherent CDG\+quasi-algebra
of~$(\g,\widetilde\g)$.
 Then there is a natural closed isomorphism of DG\+functors
$$
 \cC^\cu_X(\g,\widetilde\g,\cA_X)\ot_\cA{-}\,\simeq
 \cHom_\cA^\cu(\cC_\cu^X(\cA_X,\g,\widetilde\g),{-})\:
 \bCom(\cA\Qcoh)\lrarrow\cB^\cu\bQcoh.
$$
 Here $\cC^\cu_X(\g,\widetilde\g,\cA_X)$ and
$\cC_\cu^X(\cA_X,\g,\widetilde\g)$ are the Chevalley--Eilenberg
quasi-coherent CDG\+bimodules constructed in
Section~\ref{chevalley-eilenberg-cdg-modules-subsecn}.
\end{cor}

\begin{proof}
 On the level of the underlying quasi-coherent graded modules, we
have $\cC^\cu_X(\g,\widetilde\g,\cA_X)\ot_\cA{-}\,=
\cB^*\ot_{\cO_X}{-}$ and
$\cHom_\cA^\cu(\cC_\cu^X(\cA_X,\g,\widetilde\g),{-})=
\cHom_{\cO_X}(\bigwedge_X^*(\g),{-})$; so there is an obvious
isomorphism of quasi-coherent graded $\cB^*$\+modules
\begin{multline} \label{tensor-Hom-koszul-duality-functor-isom}
 \cC^*_X(\g,\widetilde\g,\cA_X)\ot_\cA\M^*=\cB^*\ot_{\cO_X}\M^* \\
 \,\simeq\,\cHom_{\cO_X}^*\Bigl(\bigwedge\nolimits_X^*(\g),\M^*\Bigr)
 =\cHom_\cA^*(\cC_*^X(\cA_X,\g,\widetilde\g),\M^*)
\end{multline}
for any quasi-coherent graded left $\cA$\+module~$\M^*$.
 It remains to check that~\eqref{tensor-Hom-koszul-duality-functor-isom}
underlies a closed isomorphism of quasi-coherent CDG\+modules over
$\cB^\cu$ for any complex of quasi-coherent left
$\cA$\+modules~$\M^\cu$.

 Recall that, according to the construction in
Section~\ref{chevalley-eilenberg-cdg-modules-subsecn}, we have
a closed isomorphism
\begin{equation} \label{homol-ce-cdg-bimod-dualization-of-cohomol-one}
 \cC_\cu^X(\cA_X,\g,\widetilde\g)(U)=
 \Hom_{\cA(U)^\rop}^\cu(\cC^\cu_X(\g,\widetilde\g,\cA_X)(U),\cA(U))
\end{equation}
of CDG\+bimodules over $\cA(U)$ and $\cB^\cu(U)$ for all affine open
subschemes $U\subset X$.
 As the graded right $\cA_X(U)$\+module
$\cC^\cu_X(\g,\widetilde\g,\cA_X)(U)$ is finitely generated and
projective,
formula~\eqref{homol-ce-cdg-bimod-dualization-of-cohomol-one} can be
equivalently rewritten as a closed isomorphism
\begin{equation} \label{cohomol-ce-cdg-bimod-dualization-of-homol-one}
 \cC^\cu_X(\g,\widetilde\g,\cA_X)(U)\simeq
 \Hom_{\cA(U)}^\cu(\cC_\cu^X(\cA_X,\g,\widetilde\g)(U),\cA(U))
\end{equation}
of CDG\+bimodules over $\cB^\cu(U)$ and~$\cA(U)$.
 In the notation of the present
Section~\ref{internal-cHom-of-quasi-coherent-cdg-modules-subsecn},
the collection of closed isomorphisms of
CDG\+bimodules~\eqref{cohomol-ce-cdg-bimod-dualization-of-homol-one}\
for all affine open subschemes $U\subset X$ defines a closed
isomorphism
\begin{equation} \label{ce-qcoh-cdg-bimodules-dualization}
 \cC^\cu_X(\g,\widetilde\g,\cA_X)\simeq
 \cHom_\cA^\cu(\cC_\cu^X(\cA_X,\g,\widetilde\g),\cA).
\end{equation}
 It remains to refer to
Lemma~\ref{cdg-bimodule-fin-gen-proj-Hom-tensor-lemma} in order
to obtain the desired closed isomorphism
$$
 \cC_X^\cu(\g,\widetilde\g,\cA_X)\ot_\cA\M^\bu\simeq
 \cHom_\cA^\cu(\cC^X_\cu(\cA_X,\g,\widetilde\g),\M^\bu)
$$
for any complex of quasi-coherent left $\cA$\+modules~$\M^\bu$.
\end{proof}

\subsection{$\Cohom$ into $\fHom$} \label{Cohom-into-fHom-subsecn}
 The aim of this section is to construct a quasi-coherent CDG\+bimodule
version of the natural isomorphism~\cite[formula~(2.17) in
Section~2.5]{Pcosh}.

\begin{lem} \label{Cohom-into-fHom-lemma}
 Let $\cA$, $\cB$, and $\cC$ be three quasi-coherent quasi-algebras
over a quasi-compact semi-separated scheme~$X$.
 Let $\cE$ be a quasi-coherent $\cA$\+$\cB$\+bimodule, $\F$ be
a quasi-coherent $\cB$\+$\cC$\+bimodule, and $\G$ be a quasi-coherent
left $\cA$\+module on~$X$.
 In this setting: \par
\textup{(a)} if the quasi-coherent $\cA$\+$\cB$\+bimodule $\cE$ is
$\cA$\+very flaprojective, the quasi-coherent $\cB$\+$\cC$\+bimodule
$\F$ is $\cB$\+very flaprojective, and the quasi-coherent left
$\cA$\+module $\G$ is $X$\+contraadjusted, then there is a natural
isomorphism of contraherent $\cC$\+modules
\begin{equation} \label{Cohom-into-fHom-formula}
 \Cohom_\cB(\F,\fHom_\cA(\cE,\G))\simeq\fHom_\cA(\cE\ot_\cB\F,\>\G)
\end{equation}
on~$X$; \par
\textup{(b)} if the quasi-coherent $\cA$\+$\cB$\+bimodule $\cE$ is
$\cA$\+robustly flaprojective, the quasi-coherent
$\cB$\+$\cC$\+bimodule $\F$ is $\cB$\+robustly flaprojective, and
the quasi-coherent left $\cA$\+module $\G$ is $X$\+cotorsion, then
there is a natural isomorphism of locally cotorsion contraherent
$\cC$\+modules~\eqref{Cohom-into-fHom-formula} on~$X$.
\end{lem}

\begin{proof}
 Notice first of all that, in the context of part~(a),
the left-hand side of~\eqref{Cohom-into-fHom-formula} is
well-defined by Lemmas~\ref{fHom-contraadjusted-and-more-lemma}(a)
and~\ref{module-Hom-contraadjusted-and-more-lemma}(a).
 The right-hand side of~\eqref{Cohom-into-fHom-formula} is
well-defined because the quasi-coherent $\cA$\+$\cC$\+bimodule
$\cE\ot_\cB\F$ (constructed in
Section~\ref{tensor-product-of-qcoh-cdg-bimods}) is $\cA$\+very
flaprojective by Lemma~\ref{very-flaprojective-Hom-tensor-lemma}(a).
 Similarly, in the context of part~(b), the left-hand side
of~\eqref{Cohom-into-fHom-formula} is well-defined by
Lemmas~\ref{fHom-contraadjusted-and-more-lemma}(b)
and~\ref{module-Hom-contraadjusted-and-more-lemma}(b).
 The right-hand side of~\eqref{Cohom-into-fHom-formula} is
well-defined because the quasi-coherent $\cA$\+$\cC$\+bimodule
$\cE\ot_\cB\F$ is $\cA$\+robustly flaprojective by
Corollary~\ref{robustly-flaprojective-tensor-product-cor}.

 In both the cases (a) and~(b), the natural
isomorphism~\eqref{Cohom-into-fHom-formula} holds because for
every affine open subscheme $U\subset X$ with the open immersion
morphism $j\:U\rarrow X$ one has
\begin{multline*}
 \Cohom_\cB(\F,\fHom_\cA(\cE,\G))[U]=
 \Hom_{\cB(U)}(\F(U),\fHom_\cA(\cE,\G)[U]) \\ =
 \Hom_{\cB(U)}(\F(U),\Hom_\cA(j_*j^*\cE,\G))\simeq
 \Hom_\cA((j_*j^*\cE)\ot_{\cB(U)}\F(U),\>\G) \\
 \Hom_\cA((j_*j^*(\cE\ot_\cB\F),\>\G) =
 \fHom_\cA(\cE\ot_\cB\F,\>\G)[U].
\end{multline*}
 Indeed, the natural morphism
$$
 (j_*j^*\cE)\ot_{\cB(U)}\F(U)\lrarrow j_*j^*(\cE\ot_\cB\F)
$$
is an isomorphism of quasi-coherent quasi-modules (in fact, of
quasi-coherent $\cA$\+$\cC$\+bi\-mod\-ules) on $X$ because, for
every affine open subscheme $V\subset X$, one has
\begin{multline*}
 ((j_*j^*\cE)\ot_{\cB(U)}\F(U))(V)=
 (j_*j^*\cE)(V)\ot_{\cB(U)}\F(U)=
 \cE(U\cap V)\ot_{\cB(U)}\F(U) \\ \simeq
 \cE(U\cap V)\ot_{\cB(U\cap V)}(\cB(U\cap V)\ot_{\cB(U)}\F(U))
 \\ \simeq \cE(U\cap V)\ot_{\cB(U\cap V)}\F(U\cap V) =
 (j_*j^*(\cE\ot_\cB\F))(V)
$$
\end{multline*}
(cf.~\cite[formula~(2.14) in Section~2.5]{Pcosh}).
\end{proof}

\begin{lem} \label{cdg-bimodule-Cohom-into-fHom-lemma}
 Let $\cA^\cu$, $\cB^\cu$, and $\cC^\cu$ be three quasi-coherent
CDG\+quasi-algebras over a quasi-compact semi-separated scheme~$X$.
 Let $\cE^\cu$ be a quasi-coherent CDG\+bimodule over $\cA^\cu$
and $\cB^\cu$, let $\F^\cu$ be a quasi-coherent CDG\+bimodule over
$\cB^\cu$ and $\cC^\cu$, and left $\G^\cu$ be a quasi-coherent
left CDG\+module over~$\cA^\cu$.
 In this setting: \par
\textup{(a)} if the quasi-coherent graded $\cA^*$\+$\cB^*$\+bimodule
$\cE^*$ is $\cA^*$\+very flaprojective, the quasi-coherent graded
$\cB^*$\+$\cC^*$\+bimodule $\F^*$ is $\cB^*$\+very flaprojective, and
the quasi-coherent graded left $\cA^*$\+module $\G^*$ is
$X$\+contraadjusted, then there is a natural closed isomorphism
\begin{equation} \label{cdg-bimodule-Cohom-into-fHom-formula}
 \Cohom_{\cB^*}^\cu(\F^\cu,\fHom_{\cA^*}^\cu(\cE^\cu,\G^\cu))
 \simeq\fHom_{\cA^*}^\cu(\cE^\cu\ot_{\cB^*}\F^\cu,\>\G^\cu)
\end{equation}
of contraherent CDG\+modules over~$\cC^\cu$; \par
\textup{(b)} if the quasi-coherent graded $\cA^*$\+$\cB^*$\+bimodule
$\cE^*$ is $\cA^*$\+robustly flaprojective, the quasi-coherent graded
$\cB^*$\+$\cC^*$\+bimodule $\F^*$ is $\cB^*$\+robustly flaprojective,
and the quasi-coherent graded left $\cA^*$\+module $\G^*$ is
$X$\+cotorsion, then there is a natural
isomorphism~\eqref{cdg-bimodule-Cohom-into-fHom-formula} of locally
cotorsion contraherent CDG\+modules over~$\cC^\cu$.
\end{lem}

\begin{proof}
 In both the contexts~(a) and~(b), the contraherent CDG\+module
$\fHom_{\cA^*}^\cu(\cE^\cu,\G^\cu)$ over $\cB^\cu$ was constructed in
Section~\ref{fHom-and-contratensor-of-cdg-modules}, the quasi-coherent
CDG\+bimodule $\cE^\cu\ot_{\cB^*}\F^\cu$ over $\cA^\cu$ and $\cC^\cu$
was constructed in Section~\ref{tensor-product-of-qcoh-cdg-bimods},
and the contraherent CDG\+module $\Cohom_{\cB^*}^\cu(\F^\cu,{-})$
over $\cC^\cu$ was constructed in
Section~\ref{Cohom-into-lcth-cdg-module}.
 The isomorphism of contraherent graded $\cC^*$\+modules
$$
 \Cohom_{\cB^*}^*(\F^*,\fHom_{\cA^*}^*(\cE^*,\G^*))
 \simeq\fHom_{\cA^*}^*(\cE^*\ot_{\cB^*}\F^*,\>\G^*)
$$
underlying~\eqref{cdg-bimodule-Cohom-into-fHom-formula} is just
the graded version of the isomorphism~\eqref{Cohom-into-fHom-formula}
from Lemma~\ref{Cohom-into-fHom-lemma}.
 Checking that the resulting isomorphism of contraherent
CDG\+mod\-ules~\eqref{cdg-bimodule-Cohom-into-fHom-formula} is closed
(i.~e., respects the differentials) is straightforward.
\hbadness=1525
\end{proof}

\subsection{$\fHom$ into $\cHom$} \label{fHom-into-cHom-subsecn}
 The aim of this section is to construct a quasi-coherent CDG\+bimodule
version of the natural isomorphism~\cite[formula~(2.15) in
Section~2.5]{Pcosh}.

\begin{lem} \label{fHom-into-cHom-lemma}
 Let $\cA$, $\cB$, and $\cC$ be three quasi-coherent quasi-algebras
over a quasi-compact semi-separated scheme~$X$.
 Let $\cE$ be a quasi-coherent $\cA$\+$\cB$\+bimodule, $\F$ be
a quasi-coherent $\cB$\+$\cC$\+bimodule, and $\G$ be a quasi-coherent
left $\cA$\+module on~$X$.
 In this setting: \par
\textup{(a)} if the quasi-coherent $\cA$\+$\cB$\+bimodule $\cE$
satisfies the assumption of
Lemma~\ref{cHom-over-quasi-algebra-quasi-coherent}, the quasi-coherent
$\cB$\+$\cC$\+bimodule $\F$ is $\cB$\+very flaprojective, and
the quasi-coherent left $\cA$\+module $\G$ is $X$\+contraadjusted,
then there is a natural isomorphism of contraherent $\cC$\+modules
\begin{equation} \label{fHom-into-cHom-formula}
 \fHom_\cB(\F,\cHom_\cA(\cE,\G))\simeq\fHom_\cA(\cE\ot_\cB\F,\>\G)
\end{equation}
on~$X$; \par
\textup{(b)} if the quasi-coherent $\cA$\+$\cB$\+bimodule $\cE$
satisfies the assumption of
Lemma~\ref{cHom-over-quasi-algebra-quasi-coherent}, the quasi-coherent
$\cB$\+$\cC$\+bimodule $\F$ is $\cB$\+robustly flaprojective, and
the quasi-coherent left $\cA$\+module $\G$ is $X$\+cotorsion, then
there is a natural isomorphism of locally cotorsion contraherent
$\cC$\+modules~\eqref{fHom-into-cHom-formula} on~$X$.
\end{lem}

\begin{proof}
 The notation in Lemma~\ref{cHom-over-quasi-algebra-quasi-coherent}
was slightly different, so let us restate the assumption of that
lemma: it says that there exists a finite locally free quasi-coherent
sheaf $\cT$ on $X$ such that the underlying quasi-coherent
$\cA$\+$\cO_X$\+bimodule of $\cE$ is isomorphic to a direct summand
of the quasi-coherent $\cA$\+$\cO_X$\+bimodule $\cA\ot_{\cO_X}\cT$.
 If this is the case, then the quasi-coherent $\cA$\+$\cB$\+bimodule
$\cE$ is obviously $\cA$\+very flaprojective and $\cA$\+robustly
flaprojective.
 Hence the quasi-coherent $\cA$\+$\cC$\+bimodule $\cE\ot_\cB\F$ is
$\cA$\+very flaprojective in the context of part~(a) and
$\cA$\+robustly flaprojective in the context of part~(b) as explained
in the proof of Lemma~\ref{Cohom-into-fHom-lemma}.
 Hence the right-hand side of~\eqref{fHom-into-cHom-formula} is
well-defined in both cases.

 To show that the left-hand side of~\eqref{fHom-into-cHom-formula} is
well-defined, one needs to check that the underlying quasi-coherent 
sheaf of the quasi-coherent left $\cB$\+module $\cHom_\cA(\cE,\G)$ is
contraadjusted in part~(a) and cotorsion in part~(b).
 This holds due to the isomorphism $\cHom_\cA(\cE,\G)\simeq
\cHom_{\cO_X}(\cT,\G)$ from the proof of
Lemma~\ref{cHom-over-quasi-algebra-quasi-coherent} and by virtue of
Lemma~\ref{qcoh-internal-Hom-sheaf-cta-cot}.
 So Lemma~\ref{fHom-contraadjusted-and-more-lemma}(a) or~(b) is
applicable.

 In both the cases~(a) and~(b), the natural
isomorphism~\eqref{fHom-into-cHom-formula} holds because for
every affine open subscheme $U\subset X$ with the open immersion
morphism $j\:U\rarrow X$ one has
\begin{multline*}
 \fHom_\cB(\F,\cHom_\cA(\cE,\G))[U]=
 \Hom_\cB(j_*j^*\F,\cHom_\cA(\cE,\G)) \\ \simeq
 \Hom_\cA(\cE\ot_\cB(j_*j^*\F),\>\G)\simeq
 \Hom_\cA(j_*j^*(\cE\ot_\cB\F),\>\G)=
 \fHom_\cA(\cE\ot_\cB\F,\>\G)[U].
\end{multline*}
 Indeed, the natural morphism
$$
 \cE\ot_\cB(j_*j^*\F)\lrarrow j_*j^*(\cE\ot_\cB\F)
$$
is an isomorphism of quasi-coherent quasi-modules (in fact, of
quasi-coherent $\cA$\+$\cC$\+bi\-mod\-ules) on $X$ because, for
every affine open subscheme $V\subset X$, one has
\begin{multline*}
 (\cE\ot_\cB(j_*j^*\F))(V)=\cE(V)\ot_{\cB(V)}((j_*j^*\F)(V))=
 \cE(V)\ot_{\cB(V)}\F(U\cap V) \\
 \simeq \cE(V)\ot_{\cB(V)}\cB(U\cap V)\ot_{\cB(U\cap V)}\F(U\cap V) \\
 \simeq \cE(U\cap V)\ot_{\cB(U\cap V)}\F(U\cap V)=
 (j_*j^*(\cE\ot_\cB\F))(V)
\end{multline*}
(cf.~\cite[formula~(2.13) in Section~2.5]{Pcosh}).
\end{proof}

\begin{lem} \label{cdg-bimodule-fHom-into-cHom-lemma}
 Let $\cA^\cu$, $\cB^\cu$, and $\cC^\cu$ be three quasi-coherent
CDG\+quasi-algebras over a quasi-compact semi-separated scheme~$X$.
 Let $\cE^\cu$ be a quasi-coherent CDG\+bimodule over $\cA^\cu$
and $\cB^\cu$, let $\F^\cu$ be a quasi-coherent CDG\+bimodule over
$\cB^\cu$ and $\cC^\cu$, and left $\G^\cu$ be a quasi-coherent
left CDG\+module over~$\cA^\cu$.
 In this setting: \par
\textup{(a)} if the quasi-coherent graded $\cA^*$\+$\cB^*$\+bimodule
$\cE^*$ satisfies the assumption of
Lemma~\ref{cHom-over-graded-quasi-algebra-quasi-coherent-graded},
the quasi-coherent graded $\cB^*$\+$\cC^*$\+bimodule $\F^*$ is
$\cB^*$\+very flaprojective, and the quasi-coherent graded left
$\cA^*$\+module $\G^*$ is $X$\+contraadjusted, then there is a natural
closed isomorphism
\begin{equation} \label{cdg-bimodule-fHom-into-cHom-formula}
 \fHom_{\cB^*}^\cu(\F^\cu,\cHom_{\cA^*}^\cu(\cE^\cu,\G^\cu))
 \simeq\fHom_{\cA^*}^\cu(\cE^\cu\ot_{\cB^*}\F^\cu,\>\G^\cu)
\end{equation}
of contraherent CDG\+modules over~$\cC^\cu$; \par
\textup{(b)} if the quasi-coherent graded $\cA^*$\+$\cB^*$\+bimodule
$\cE^*$ satisfies the assumption of
Lemma~\ref{cHom-over-graded-quasi-algebra-quasi-coherent-graded},
the quasi-coherent graded $\cB^*$\+$\cC^*$\+bimodule $\F^*$ is
$\cB^*$\+robustly flaprojective, and the quasi-coherent graded left
$\cA^*$\+module $\G^*$ is $X$\+cotorsion, then there is a natural
isomorphism~\eqref{cdg-bimodule-fHom-into-cHom-formula} of locally
cotorsion contraherent CDG\+modules over~$\cC^\cu$.
\end{lem}

\begin{proof}
 In both the contexts (a) and~(b), the quasi-coherent left CDG\+module
$\cHom_{\cA^*}^\cu(\cE^\cu,\G^\cu)$ over $\cB^\cu$ was constructed in
Section~\ref{internal-cHom-of-quasi-coherent-cdg-modules-subsecn},
the quasi-coherent CDG\+bi\-mod\-ule $\cE^\cu\ot_{\cB^*}\F^\cu$ over
$\cA^\cu$ and $\cC^\cu$ was constructed in
Section~\ref{tensor-product-of-qcoh-cdg-bimods}, and
the contraherent CDG\+module $\fHom^\cu$ was constructed in
Section~\ref{fHom-and-contratensor-of-cdg-modules}.
 The isomorphism of contraherent graded $\cC^*$\+modules
$$
 \fHom_{\cB^*}^*(\F^*,\cHom_{\cA^*}^*(\cE^*,\G^*))
 \simeq\fHom_{\cA^*}^*(\cE^*\ot_{\cB^*}\F^*,\>\G^*)
$$
underlying~\eqref{cdg-bimodule-fHom-into-cHom-formula} is just
the graded version of the isomorphism~\eqref{fHom-into-cHom-formula}
from Lemma~\ref{fHom-into-cHom-lemma}.
 Checking that the resulting isomorphism of contraherent
CDG\+mod\-ules~\eqref{cdg-bimodule-fHom-into-cHom-formula} is closed
(i.~e., respects the differentials) is straightforward.
\hbadness=1525
\end{proof}

\subsection{Quadrality diagram}
 We start with stating a lemma about commutativity of the main square
diagram on the level of DG\+functors.

\begin{lem} \label{main-commutative-square-of-DG-functors-lemma}
 Let $X$ be a quasi-compact semi-separated scheme and
$(\g,\widetilde\g)$ be a quasi-coherent twisted Lie algebroid over~$X$.
 Assume that\/ $\g$~is a finite locally free sheaf on~$X$.
 Let $\cA=\cA_X(\g,\widetilde\g)$ be the twisted univeral enveloping
quasi-coherent quasi-algebra and $\cB^\cu=\cC^\cu_X(\g,\widetilde\g)$
be the Chevalley--Eilenberg quasi-coherent CDG\+quasi-algebra
of~$(\g,\widetilde\g)$.
 Then there is the following commutative diagram of exact DG\+functors
between exact DG\+categories:
\begin{equation} \label{main-commutative-square-of-DG-functors-diagram}
\begin{gathered}
 \xymatrixcolsep{5em}\xymatrix{
  & \cB^\cu\bQcoh^{X\dcta}
  \ar[rd]^-{\quad\fHom_{\cB^*}(\cB^\subcu,{-})} \\
  \bCom(\cA\Qcoh^{X\dcta})
  \ar[ru]^-{\cB^*\ot_{\cO_X}{-}\quad}
  \ar[rd]_-{\fHom_\cA(\cA,{-})\qquad}
  \ar[rr]^-{\fHom_\cA(\cC^X_\subcu(\cA,\g,\widetilde\g),{-})}
  && \cB^\cu\bCtrh_\al \\
  & \bCom(\cA\bCtrh_\al)
  \ar[ru]_-{\qquad\Cohom_X\left(\bigwedge\nolimits_X^*(\g),{-}\right)}
 }
\end{gathered}
\end{equation}
 Here the commutativity of the diagram of DG\+functors means natural
closed isomorphisms between the compositions of DG\+functors along
the two paths around the diamond and the diagonal DG\+functor.
\end{lem}

\begin{proof}
 Let us start with a discussion of the constructions of
the DG\+functors involved.
 The exact DG\+functor
\begin{equation} \label{B-tensor-DG-functor}
 \cB^*\ot_{\cO_X}{-}\,=\cC_X^\cu(\g,\widetilde\g,\cA_X)\ot_{\cA}{-}\,\:
 \bCom(\cA\Qcoh)\lrarrow\cB^\cu\bQcoh.
\end{equation}
was constructed in
Lemma~\ref{one-koszul-duality-dg-functor-left-co-side}.
 According to
Corollary~\ref{tensor-Hom-koszul-duality-functor-isom-cor},
the DG\+functor~\eqref{B-tensor-DG-functor} is naturally isomorphic
to the DG\+functor
\begin{equation} \label{cHom-from-wedge-g-DG-functor}
 \cHom_{\cO_X}\Bigl(\bigwedge\nolimits_X^*(\g),{-}\Bigr)
 =\cHom_\cA^\cu(\cC_\cu^X(\cA_X,\g,\widetilde\g),{-})\:
 \bCom(\cA\Qcoh)\lrarrow\cB^\cu\bQcoh.
\end{equation}
 Now it is claimed that the DG\+functor~\eqref{B-tensor-DG-functor}
(or equivalently, \eqref{cHom-from-wedge-g-DG-functor}) restricts to
a DG\+functor $\bCom(\cA\Qcoh^{X\dcta})\rarrow\cB^\cu\bQcoh^{X\dcta}$.
 It other words, this means existence of a commutative diagram of
DG\+functors
\begin{equation} \label{koszul-duality-restricts-to-X-cta}
\begin{gathered}
 \xymatrix{
  \bCom(\cA\Qcoh)
  \ar[rrr]^-{\cHom_{\cO_X}\left(\bigwedge\nolimits_X^*(\g),{-}\right)}
  &&& \cB^\cu\bQcoh \\
  \bCom(\cA\Qcoh^{X\dcta}) \ar@{>->}[u]
  \ar[rrr]^-{\cHom_{\cO_X}\left(\bigwedge\nolimits_X^*(\g),{-}\right)}
  &&& \cB^\cu\bQcoh^{X\dcta} \ar@{>->}[u]
 }
\end{gathered}
\end{equation}
where the vertical arrows with tails denote the natural fully faithful
inclusions.
 Indeed, the DG\+functor
$\cHom_{\cO_X}\bigl(\bigwedge_X^*(\g),{-}\bigr)$ takes
$X$\+contraadjusted quasi-coherent $\cA$\+modules to
$X$\+contraadjusted quasi-coherent CDG\+modules over $\cB^\cu$
by Lemma~\ref{qcoh-internal-Hom-sheaf-cta-cot}(a).
{\hbadness=1350\par}

 Let $\bW$ be an open covering of~$X$.
 The exact DG\+functor
\begin{multline} \label{Cohom-from-wedge-g-DG-functor}
 \Cohom_X\Bigl(\bigwedge\nolimits_X^*(\g),{-}\Bigr) =
 \Cohom_\cA^\cu(\cC_\cu^X(\cA_X,\g,\widetilde\g),{-})\: \\
 \bCom(\cA\Lcth_\bW)\lrarrow\cB^\cu\bLcth_\bW
\end{multline}
was constructed in
Lemma~\ref{koszul-duality-dg-functors-contra-side}(a).
 Now it is claimed that
the DG\+functor~\eqref{Cohom-from-wedge-g-DG-functor} restricts to
a DG\+functor $\bCom(\cA\Ctrh_\al)\rarrow\cB^\cu\bCtrh_\al$.
 It other words, this means existence of a commutative diagram of
DG\+functors
\begin{equation} \label{koszul-duality-restricts-to-antilocal}
\begin{gathered}
 \xymatrix{
  \bCom(\cA\Lcth_\bW)
  \ar[rrr]^-{\Cohom_X\left(\bigwedge\nolimits_X^*(\g),{-}\right)}
  &&& \cB^\cu\bLcth_\bW \\
  \bCom(\cA\Ctrh_\al) \ar@{>->}[u]
  \ar[rrr]^-{\Cohom_X\left(\bigwedge\nolimits_X^*(\g),{-}\right)}
  &&& \cB^\cu\bCtrh_\al \ar@{>->}[u]
 }
\end{gathered}
\end{equation}
where the vertical arrows with tails denote the natural fully faithful
inclusions.
 Indeed, the DG\+functor
$\Cohom_X\bigl(\bigwedge_X^*(\g),{-}\bigr)$ takes antilocal contraherent
$\cA$\+modules to antilocal contraherent CDG\+modules over $\cB^\cu$
by Lemma~\ref{Cohom-into-antilocal-cosheaf-is-antilocal}(a).

 The exact DG\+functor
$$
 \fHom_\cA(\cA,{-})\:\bCom(\cA\Qcoh^{X\dcta})\lrarrow
 \bCom(\cA\bCtrh_\al)
$$
is constructed by applying the exact functor
$$
 \fHom_\cA(\cA,{-})\:\cA\Qcoh^{X\dcta}\lrarrow\cA\bCtrh_\al
$$
to complexes termwise.
 The latter exact functor was constructed
in Section~\ref{fHom-over-qcoh-quasi-algebra-subsecn} and appeared in
Lemma~\ref{quasi-algebra-underived-naive-co-contra}(a).
 The exact DG\+functor
\begin{equation} \label{fHom-from-B-DG-functor}
 \fHom_{\cB^*}^\cu(\cB^\cu,{-})\:\cB^\cu\bQcoh^{X\dcta}\lrarrow
 \cB^\cu\bCtrh_\al
\end{equation}
was constructed in Section~\ref{fHom-and-contratensor-of-cdg-modules}
and appeared in Lemma~\ref{cta-cot-al-co-contra-dg-equivalence}(a).

 The exact DG\+functor
\begin{equation} \label{fHom-from-ce-DG-functor}
 \fHom_\cA^\cu(\cC^X_\cu(\cA,\g,\widetilde\g),{-})\:
 \bCom(\cA\Qcoh^{X\dcta})\lrarrow\cB^\cu\bCtrh
\end{equation}
is also provided by the construction of
Section~\ref{fHom-and-contratensor-of-cdg-modules}.
 In view of the graded version of
Corollary~\ref{fHom-antilocal-corollary},
the DG\+functor~\eqref{fHom-from-ce-DG-functor} lands in
the full DG\+subcategory $\cB^\cu\bCtrh_\al\subset\cB^\cu\bCtrh$,
so we have the exact DG\+functor
\begin{equation} \label{fHom-from-ce-DG-functor-into-al}
 \fHom_\cA^\cu(\cC^X_\cu(\cA,\g,\widetilde\g),{-})\:
 \bCom(\cA\Qcoh^{X\dcta})\lrarrow\cB^\cu\bCtrh_\al.
\end{equation}

 It remains to establish the commutativity of the two triangles forming
the diagram of
DG\+functors~\eqref{main-commutative-square-of-DG-functors-diagram}.
 Indeed, let $\M^\bu$ be a complex of $X$\+contraadjusted quasi-coherent
left $\cA$\+modules on~$X$.
 Then there are natural closed isomorphisms of (antilocal) contraherent
CDG\+modules over~$\cB^\cu$
\begin{multline*}
 \fHom_{\cB^*}^\cu(\cB^\cu,
 \cHom_\cA^\cu(\cC^X_\cu(\cA_X,\g,\widetilde\g),\M^\bu)) \\ \simeq
 \fHom_\cA^\cu(\cC^X_\cu(\cA_X,\g,\widetilde\g)\ot_{\cB^*}\cB^\cu,
 \>\M^\bu)\simeq\fHom_\cA^\cu(\cC^X_\cu(\cA_X,\g,\widetilde\g),\M^\bu)
\end{multline*}
by Lemma~\ref{cdg-bimodule-fHom-into-cHom-lemma}(a).
 This proves the commutativity of the upper triangle
in~\eqref{main-commutative-square-of-DG-functors-diagram}.

 Similarly, there are natural closed isomorphisms of (antilocal) 
contraherent CDG\+modules over~$\cB^\cu$
\begin{multline*}
 \Cohom_\cA^\cu(\cC^X_\cu(\cA_X,\g,\widetilde\g),
 \fHom_\cA(\cA,\M^\bu)) \\ \simeq
 \fHom_\cA^\cu(\cA\ot_\cA\cC^X_\cu(\cA_X,\g,\widetilde\g),\>\M^\bu)
 \simeq\fHom_\cA^\cu(\cC^X_\cu(\cA_X,\g,\widetilde\g),\M^\bu)
\end{multline*}
by Lemma~\ref{cdg-bimodule-Cohom-into-fHom-lemma}(a).
 This proves the commutativity of the lower triangle
in~\eqref{main-commutative-square-of-DG-functors-diagram}.
\end{proof}

 The following theorem is the main result of
Section~\ref{quadrality-diagram-secn}, and also of the whole paper.

\begin{thm} \label{main-quadrality-theorem}
 Let $X$ be a quasi-compact semi-separated scheme with an open
covering\/ $\bW$ and $(\g,\widetilde\g)$ be a quasi-coherent twisted
Lie algebroid over~$X$.
 Assume that\/ $\g$~is a finite locally free sheaf on~$X$.
 Let $\cA=\cA_X(\g,\widetilde\g)$ be the twisted universal enveloping
quasi-coherent quasi-algebra and $\cB^\cu=\cC^\cu_X(\g,\widetilde\g)$
be the Chevalley--Eilenberg quasi-coherent CDG\+quasi-algebra
of~$(\g,\widetilde\g)$.
 Then there is a natural commutative square diagram of triangulated
category equivalences
\begin{equation} \label{main-quadrality-diagram}
\begin{gathered}
 \xymatrix{
  \sD(\cA\Qcoh) \ar@<3pt>[rrrr]^{\cB^*\ot_{\cO_X}{-}}
  \ar@<3pt>[dd]^{\boR\fHom_\cA(\cA,{-})}
  &&&& \sD^\co_{X\red}(\cB^\cu\bQcoh)
  \ar@{-}@<3pt>[llll]
  \ar@<3pt>[dd]^{\boR\fHom_{\cB^*}(\cB^\subcu,{-})} \\ \\
  \sD(\cA\Lcth_\bW)
  \ar@<3pt>[rrrr]^{\Cohom_X\left(\bigwedge\nolimits_X^*(\g),{-}\right)}
  \ar@<3pt>[uu]^{\cA\ocn_\cA^\boL{-}}
  &&&& \sD^\ctr_{X\red}(\cB^\cu\bLcth_\bW)
  \ar@<3pt>[llll]^{\Cohom_X(\cA,{-})}
  \ar@<3pt>[uu]^{\cB^\subcu\ocn_{\cB^*}^\boL{-}}
 }
\end{gathered}
\end{equation}
\end{thm}

\begin{proof}
 The leftmost vertical triangulated equivalence is provided by
Theorem~\ref{quasi-algebra-bounded-derived-naive-co-contra}
or~\ref{quasi-algebra-unbounded-derived-naive-co-contra};
see formula~\eqref{bounded-derived-naive-co-contra-diagram-II}.
 The rightmost vertical triangulated equivalence is provided by
Corollary~\ref{cdg-module-co-contra-correspondence-cor}.
 The lower horizontal triangulated equivalence is provided by
Theorem~\ref{reduced-koszul-duality-contra-side}(a).
 The upper horizontal triangulated equivalence is provided by
Corollary~\ref{left-right-co-side-reduced-koszul-duality}.

 Notice that, strictly speaking, the upper horizontal triangulated
equivalence was proved in
Corollary~\ref{left-right-co-side-reduced-koszul-duality} under
the additional assumption that $\g$~is a finite locally free sheaf
\emph{of constant rank} on~$X$.
 Decomposing the scheme $X$ into the disjoint union of open
subschemes on which the rank of~$\g$ is constant, one obtains
the upper horizontal triangulated equivalence for any quasi-coherent
twisted Lie algebroid $(\g,\widetilde\g)$ with a finite locally
free sheaf~$\g$ on~$X$.

 It remains to prove commutativity of the square diagram of triangulated
equivalences~\eqref{main-quadrality-diagram}.
 For this purpose, we check commutativity of the related square
diagram of DG\+functors or triangulated functors on the underived level,
and also check that the functors involved in such diagram take objects
adjusted to the respective derived functors to adjusted objects.

 In addition to what has been said in
Lemma~\ref{main-commutative-square-of-DG-functors-lemma}, we only need
make some observations about the thickness condition.
 It is clear from the construction that
the DG\+functor~\eqref{B-tensor-DG-functor}
or~\eqref{cHom-from-wedge-g-DG-functor} lands in the full
DG\+subcategory of thick quasi-coherent CDG\+modules
$\cB^\cu\bQcoh_\bth\subset\cB^\cu\bQcoh$.
 So the commutative diagram~\eqref{koszul-duality-restricts-to-X-cta}
can be redrawn as
\begin{equation} \label{koszul-duality-restricts-to-X-cta-thick}
\begin{gathered}
 \xymatrix{
  \bCom(\cA\Qcoh)
  \ar[rrrrr]^-{\cB^*\ot_{\cO_X}{-}\,\,=\,
  \cHom_{\cO_X}\left(\bigwedge\nolimits_X^*(\g),{-}\right)}
  &&&&& \cB^\cu\bQcoh_\bth \\
  \bCom(\cA\Qcoh^{X\dcta}) \ar@{>->}[u]
  \ar[rrrrr]^-{\cB^*\ot_{\cO_X}{-}\,\,=\,
  \cHom_{\cO_X}\left(\bigwedge\nolimits_X^*(\g),{-}\right)}
  &&&&& \cB^\cu\bQcoh^{X\dcta}_\bth \ar@{>->}[u]
 }
\end{gathered}
\end{equation}

 Similarly, the DG\+functor~\eqref{Cohom-from-wedge-g-DG-functor}
lands in the full DG\+subcategory of thick $\bW$\+locally contraherent
CDG\+modules $\cB^\cu\bLcth_\bW^\bth\subset\cB^\cu\bLcth_\bW$.
 So the commutative
diagram~\eqref{koszul-duality-restricts-to-antilocal} can be redrawn as
\begin{equation} \label{koszul-duality-restricts-to-antilocal-thick}
\begin{gathered}
 \xymatrix{
  \bCom(\cA\Lcth_\bW)
  \ar[rrr]^-{\Cohom_X\left(\bigwedge\nolimits_X^*(\g),{-}\right)}
  &&& \cB^\cu\bLcth_\bW^\bth \\
  \bCom(\cA\Ctrh_\al) \ar@{>->}[u]
  \ar[rrr]^-{\Cohom_X\left(\bigwedge\nolimits_X^*(\g),{-}\right)}
  &&& \cB^\cu\bCtrh_\al^\bth \ar@{>->}[u]
 }
\end{gathered}
\end{equation}

 Finally, the DG\+functor~\eqref{fHom-from-B-DG-functor} takes
thick quasi-coherent CDG\+modules to thick (antilocal) contraherent
CDG\+modules by
Theorem~\ref{thick-corresponds-to-thick-under-co-contra}.
 So we have a commutative diagram of DG\+functors
\begin{equation} \label{co-contra-over-B-restricts-to-thick}
\begin{gathered}
 \xymatrix{
  \cB^\cu\bQcoh^{X\dcta}
  \ar[rrr]^-{\fHom_{\cB^*}(\cB^\subcu,{-})}
  &&& \cB^\cu\bCtrh_\al \\
  \cB^\cu\bQcoh^{X\dcta}_\bth \ar@{>->}[u]
  \ar[rrr]^-{\fHom_{\cB^*}(\cB^\subcu,{-})}
  &&& \cB^\cu\bCtrh_\al^\bth \ar@{>->}[u]
 }
\end{gathered}
\end{equation}
where the vertical arrows with tails denote the natural fully faithful
inclusions.

 We have arrived to a version of the commutative diagram of
DG\+functors~\eqref{main-commutative-square-of-DG-functors-diagram}
taking proper care of the thickness conditions:
\begin{equation} \label{main-commsq-of-DG-functors-thickness-diagram}
\begin{gathered}
 \xymatrixcolsep{5em}\xymatrix{
  & \cB^\cu\bQcoh^{X\dcta}_\bth
  \ar[rd]^-{\quad\fHom_{\cB^*}(\cB^\subcu,{-})} \\
  \bCom(\cA\Qcoh^{X\dcta})
  \ar[ru]^-{\cB^*\ot_{\cO_X}{-}\quad}
  \ar[rd]_-{\fHom_\cA(\cA,{-})\qquad}
  \ar[rr]^-{\fHom_\cA(\cC^X_\subcu(\cA,\g,\widetilde\g),{-})}
  && \cB^\cu\bCtrh_\al^\bth \\
  & \bCom(\cA\bCtrh_\al)
  \ar[ru]_-{\qquad\Cohom_X\left(\bigwedge\nolimits_X^*(\g),{-}\right)}
 }
\end{gathered}
\end{equation}
 The commutativity of the square diagram of
triangulated functors~\eqref{main-quadrality-diagram} now follows
immediately from the commutativity of the diagram of
DG\+functors~\eqref{main-commsq-of-DG-functors-thickness-diagram}
in view of the constructions of the derived functors involved
in~\eqref{main-quadrality-diagram}.
\end{proof}

\begin{cor} \label{main-hexagonality-corollary}
 Let $X$ be a quasi-compact semi-separated scheme with an open
covering\/ $\bW$ and $(\g,\widetilde\g)$ be a quasi-coherent twisted
Lie algebroid over~$X$.
 Assume that the quasi-coherent sheaf\/~$\g$ on $X$ is finite locally
free of constant rank~$m$.
 Let $\cA=\cA_X(\g,\widetilde\g)$ and
$\cA^\circ=\cA_X(\g,\widetilde\g^\circ)$ be the two related twisted
universal enveloping quasi-coherent quasi-algebras, and let
$\cB^\cu=\cC_X^\cu(\g,\widetilde\g)$ and
$\cB^\circ{}^\cu=\cC_X^\cu(\g,\widetilde\g^\circ)$ be the two related
Chevalley--Eilenberg quasi-coherent CDG\+quasi-algebras over~$X$.
 Then there is a natural commutative diagram of triangulated category
equivalences
\begin{equation} \label{main-hexagonality-diagram}
\begin{gathered}
 \xymatrix{
  \sD(\Qcohr\cA^\circ)
  \ar@<3pt>[rrrr]^{{-}\ot_{\cO_X}\bigwedge\nolimits_X^*(\g)}
  \ar@<3pt>[dd]^{\cHom_{\cO_X}(\cB^m[-m],{-})}
  &&&& \sD^\co_{X\red}(\bQcohr\cB^\circ{}^\cu)
  \ar@<3pt>[llll]^{{-}\ot_{\cO_X}\cA^\circ}
  \ar@{=}[dd]^{\cB^{\circ,\rop}{}^\subcu\simeq\cB^\subcu}
  \\ \\
  \sD(\cA\Qcoh) \ar@<3pt>[rrrr]^{\cB^*\ot_{\cO_X}{-}}
  \ar@<3pt>[uu]^{\cB^m[-m]\ot_{\cO_X}{-}}
  \ar@<3pt>[dd]^{\boR\fHom_\cA(\cA,{-})}
  &&&& \sD^\co_{X\red}(\cB^\cu\bQcoh)
  \ar@{-}@<3pt>[llll]
  \ar@<3pt>[dd]^{\boR\fHom_{\cB^*}(\cB^\subcu,{-})} \\ \\
  \sD(\cA\Lcth_\bW)
  \ar@<3pt>[rrrr]^{\Cohom_X\left(\bigwedge\nolimits_X^*(\g),{-}\right)}
  \ar@<3pt>[uu]^{\cA\ocn_\cA^\boL{-}}
  &&&& \sD^\ctr_{X\red}(\cB^\cu\bLcth_\bW)
  \ar@<3pt>[llll]^{\Cohom_X(\cA,{-})}
  \ar@<3pt>[uu]^{\cB^\subcu\ocn_{\cB^*}^\boL{-}}
 }
\end{gathered}
\end{equation}
\end{cor}

\begin{proof}
 The upper commutative square of triangulated equivalences is provided
by Corollary~\ref{left-right-co-side-reduced-koszul-duality}.
 The lower commutative square of triangulated equivalences is
the result of Theorem~\ref{main-quadrality-theorem}.
\end{proof}

\subsection{Quadrality diagrams with locally cotorsion contraherent
modules} \label{quadrality-loc-cotors-subsecn}
 The following lemma is an $X$\+cotorsion and $\cA$\+cotorsion version
of Lemma~\ref{main-commutative-square-of-DG-functors-lemma}.

\begin{lem} \label{main-lct-commutative-square-of-DG-functors-lemma}
 Let $X$ be a quasi-compact semi-separated scheme and
$(\g,\widetilde\g)$ be a quasi-coherent twisted Lie algebroid over~$X$.
 Assume that\/ $\g$~is a finite locally free sheaf on~$X$.
 Let $\cA=\cA_X(\g,\widetilde\g)$ be the twisted univeral enveloping
quasi-coherent quasi-algebra and $\cB^\cu=\cC^\cu_X(\g,\widetilde\g)$
be the Chevalley--Eilenberg quasi-coherent CDG\+quasi-algebra
of~$(\g,\widetilde\g)$.
 Then there are the following commutative diagrams of exact DG\+functors
between exact DG\+categories:
\begin{equation} \label{main-X-cot-square-of-DG-functors-diagram}
\begin{gathered}
 \xymatrixcolsep{4em}\xymatrix{
  & \cB^\cu\bQcoh^{X\dcot}
  \ar[rd]^-{\quad\fHom_{\cB^*}(\cB^\subcu,{-})} \\
  \bCom(\cA\Qcoh^{X\dcot})
  \ar[ru]^-{\cB^*\ot_{\cO_X}{-}\quad}
  \ar[rd]_-{\fHom_\cA(\cA,{-})\qquad}
  \ar[rr]^-{\fHom_\cA(\cC^X_\subcu(\cA,\g,\widetilde\g),{-})}
  && \cB^\cu\bCtrh_\al^{X\dlct} \\
  & \bCom(\cA\bCtrh_\al^{X\dlct})
  \ar[ru]_-{\qquad\Cohom_X\left(\bigwedge\nolimits_X^*(\g),{-}\right)}
 }
\end{gathered}
\end{equation}
\begin{equation} \label{main-A-cot-square-of-DG-functors-diagram}
\begin{gathered}
 \xymatrixcolsep{4em}\xymatrix{
  & \cB^\cu\bQcoh^{X\dcot}
  \ar[rd]^-{\quad\fHom_{\cB^*}(\cB^\subcu,{-})} \\
  \bCom(\cA\Qcoh^{\cA\dcot})
  \ar[ru]^-{\cB^*\ot_{\cO_X}{-}\quad}
  \ar[rd]_-{\fHom_\cA(\cA,{-})\qquad}
  \ar[rr]^-{\fHom_\cA(\cC^X_\subcu(\cA,\g,\widetilde\g),{-})}
  && \cB^\cu\bCtrh_\al^{X\dlct} \\
  & \bCom(\cA\bCtrh_\al^{\cA\dlct})
  \ar[ru]_-{\qquad\Cohom_X\left(\bigwedge\nolimits_X^*(\g),{-}\right)}
 }
\end{gathered}
\end{equation}
\end{lem}

\begin{proof}
 The DG\+functor~\eqref{B-tensor-DG-functor} (or equivalently,
\eqref{cHom-from-wedge-g-DG-functor}) restricts to a DG\+functor
$\bCom(\cA\Qcoh^{X\dcot})\rarrow\cB^\cu(\bQcoh^{X\dcot})$ by
Lemma~\ref{qcoh-internal-Hom-sheaf-cta-cot}(b).
 The exact DG\+functor
\begin{multline} \label{Cohom-from-wedge-g-X-lct-DG-functor}
 \Cohom_X\Bigl(\bigwedge\nolimits_X^*(\g),{-}\Bigr) =
 \Cohom_\cA^\cu(\cC_\cu^X(\cA_X,\g,\widetilde\g),{-})\: \\
 \bCom(\cA\Lcth_\bW^{X\dlct})\lrarrow\cB^\cu\bLcth_\bW^{X\dlct}
\end{multline}
from Lemma~\ref{koszul-duality-dg-functors-contra-side}(b)
restricts to a DG\+functor $\bCom(\cA\Ctrh_\al^{X\dlct})\rarrow
\cB^\cu\bCtrh_\al^{X\dlct}$ by
Lemma~\ref{Cohom-into-antilocal-cosheaf-is-antilocal}(b).

 The exact DG\+functor
$$
 \fHom_\cA(\cA,{-})\:\bCom(\cA\Qcoh^{X\dcot})\lrarrow
 \bCom(\cA\bCtrh_\al^{X\dlct})
$$
is constructed by applying the exact functor
$$
 \fHom_\cA(\cA,{-})\:\cA\Qcoh^{X\dcot}\lrarrow\cA\bCtrh_\al^{X\dlct}
$$
to complexes termwise.
 The latter exact functor was constructed
in Section~\ref{fHom-over-qcoh-quasi-algebra-subsecn} and appeared in
Lemma~\ref{quasi-algebra-underived-naive-co-contra}(b).
 The exact DG\+functor
$$
 \fHom_{\cB^*}^\cu(\cB^\cu,{-})\:\cB^\cu\bQcoh^{X\dcot}\lrarrow
 \cB^\cu\bCtrh_\al^{X\dlct}
$$
was constructed in Section~\ref{fHom-and-contratensor-of-cdg-modules}
and appeared in Lemma~\ref{cta-cot-al-co-contra-dg-equivalence}(b).

 The exact DG\+functor
\begin{equation} \label{fHom-from-ce-X-cot-X-lct-DG-functor}
 \fHom_\cA^\cu(\cC^X_\cu(\cA,\g,\widetilde\g),{-})\:
 \bCom(\cA\Qcoh^{X\dcot})\lrarrow\cB^\cu\bCtrh^{X\dlct}
\end{equation}
is also provided by the construction of
Section~\ref{fHom-and-contratensor-of-cdg-modules}.
 By the graded version of
Corollary~\ref{fHom-antilocal-corollary},
the DG\+functor~\eqref{fHom-from-ce-X-cot-X-lct-DG-functor} lands in
the full DG\+subcategory $\cB^\cu\bCtrh_\al^{X\dlct}\subset
\cB^\cu\bCtrh^{X\dlct}$, so we have the exact DG\+functor
$$
 \fHom_\cA^\cu(\cC^X_\cu(\cA,\g,\widetilde\g),{-})\:
 \bCom(\cA\Qcoh^{X\dcot})\lrarrow\cB^\cu\bCtrh_\al^{X\dlct}.
$$

 The commutativity of
the diagram~\eqref{main-X-cot-square-of-DG-functors-diagram}
follows from the commutativity of
the diagram~\eqref{main-commutative-square-of-DG-functors-diagram},
as the DG\+functors
in~\eqref{main-X-cot-square-of-DG-functors-diagram} are the restictions
of the DG\+functors
in~\eqref{main-commutative-square-of-DG-functors-diagram} to
the respective full exact DG\+subcategories.

 The commutative
diagram~\eqref{main-A-cot-square-of-DG-functors-diagram} is obtained
by a further restriction of the commutative
diagram~\eqref{main-X-cot-square-of-DG-functors-diagram} to
the respective full exact DG\+subcategories.
 It is only the lower leftmost diagonal functor that needs to be
explained.
 The exact DG\+functor
$$
 \fHom_\cA(\cA,{-})\:\bCom(\cA\Qcoh^{\cA\dcot})\lrarrow
 \bCom(\cA\bCtrh_\al^{\cA\dlct})
$$
is constructed by applying the exact functor
$$
 \fHom_\cA(\cA,{-})\:\cA\Qcoh^{\cA\dcot}\lrarrow\cA\bCtrh_\al^{\cA\dlct}
$$
to complexes termwise.
 The latter exact functor was constructed
in Section~\ref{fHom-over-qcoh-quasi-algebra-subsecn} and appeared in
Lemma~\ref{quasi-algebra-underived-naive-co-contra}(c).
\end{proof}

 Now we can formulate the $X$\+cotorsion and $\cA$\+cotorsion versions
of Theorem~\ref{main-quadrality-theorem}.

\begin{thm} \label{main-lct-quadrality-theorem}
 Let $X$ be a quasi-compact semi-separated scheme with an open
covering\/ $\bW$ and $(\g,\widetilde\g)$ be a quasi-coherent twisted
Lie algebroid over~$X$.
 Assume that\/ $\g$~is a finite locally free sheaf on~$X$.
 Let $\cA=\cA_X(\g,\widetilde\g)$ be the twisted universal enveloping
quasi-coherent quasi-algebra and $\cB^\cu=\cC^\cu_X(\g,\widetilde\g)$
be the Chevalley--Eilenberg quasi-coherent CDG\+quasi-algebra
of~$(\g,\widetilde\g)$.
 Then there are natural commutative square diagrams of triangulated
category equivalences
\begin{equation} \label{main-X-lct-quadrality-diagram}
\begin{gathered}
 \xymatrix{
  \sD(\cA\Qcoh) \ar@<3pt>[rrrr]^{\cB^*\ot_{\cO_X}{-}}
  \ar@<3pt>[dd]^{\boR\fHom_\cA(\cA,{-})}
  &&&& \sD^\co_{X\red}(\cB^\cu\bQcoh)
  \ar@{-}@<3pt>[llll]
  \ar@<3pt>[dd]^{\boR\fHom_{\cB^*}(\cB^\subcu,{-})} \\ \\
  \sD(\cA\Lcth_\bW^{X\dlct})
  \ar@<3pt>[rrrr]^{\Cohom_X\left(\bigwedge\nolimits_X^*(\g),{-}\right)}
  \ar@<3pt>[uu]^{\cA\ocn_\cA^\boL{-}}
  &&&& \sD^\ctr_{X\red}(\cB^\cu\bLcth_\bW^{X\dlct})
  \ar@<3pt>[llll]^{\Cohom_X(\cA,{-})}
  \ar@<3pt>[uu]^{\cB^\subcu\ocn_{\cB^*}^\boL{-}}
 }
\end{gathered}
\end{equation}
\begin{equation} \label{main-A-lct-quadrality-diagram}
\begin{gathered}
 \xymatrix{
  \sD(\cA\Qcoh) \ar@<3pt>[rrrr]^{\cB^*\ot_{\cO_X}{-}}
  \ar@<3pt>[dd]^{\boR\fHom_\cA(\cA,{-})}
  &&&& \sD^\co_{X\red}(\cB^\cu\bQcoh)
  \ar@{-}@<3pt>[llll]
  \ar@<3pt>[dd]^{\boR\fHom_{\cB^*}(\cB^\subcu,{-})} \\ \\
  \sD(\cA\Lcth_\bW^{\cA\dlct})
  \ar@<3pt>[rrrr]^{\Cohom_X\left(\bigwedge\nolimits_X^*(\g),{-}\right)}
  \ar@<3pt>[uu]^{\cA\ocn_\cA^\boL{-}}
  &&&& \sD^\ctr_{X\red}(\cB^\cu\bLcth_\bW^{X\dlct})
  \ar@<3pt>[llll]^{\Cohom_X(\cA,{-})}
  \ar@<3pt>[uu]^{\cB^\subcu\ocn_{\cB^*}^\boL{-}}
 }
\end{gathered}
\end{equation}
\end{thm}

\begin{proof}
 In diagram~\eqref{main-X-lct-quadrality-diagram}, the leftmost
vertical triangulated equivalence is provided by
Theorem~\ref{quasi-algebra-unbounded-derived-naive-co-contra}.
 The rightmost vertical triangulated equivalence is provided by
Corollary~\ref{X-cot-X-lct-cdg-module-co-contra-corresp-cor}.
 The lower horizontal triangulated equivalence is provided by
Theorem~\ref{reduced-koszul-duality-contra-side}(b).
 The upper horizontal triangulated equivalence was explained in
the proof of Theorem~\ref{main-quadrality-theorem}.

 In diagram~\ref{main-A-lct-quadrality-diagram}, the leftmost
vertical triangulated equivalence is again provided by
Theorem~\ref{quasi-algebra-unbounded-derived-naive-co-contra}.
 The lower horizontal triangulated equivalence is provided by
Theorem~\ref{reduced-koszul-duality-contra-side}(c).
 The remaining equivalences are the same as in
diagram~\ref{main-X-lct-quadrality-diagram}.

 The proof of the commutativity of both the diagrams is based on
the commutativity of the respective diagrams in
Lemma~\ref{main-lct-commutative-square-of-DG-functors-lemma}.
 The discussion of the thickness condition, which is also relevant,
can be found in the proof of Theorem~\ref{main-quadrality-theorem}.
\end{proof}

\begin{cor} \label{main-lct-hexagonality-corollary}
 Let $X$ be a quasi-compact semi-separated scheme with an open
covering\/ $\bW$ and $(\g,\widetilde\g)$ be a quasi-coherent twisted
Lie algebroid over~$X$.
 Assume that the quasi-coherent sheaf\/~$\g$ on $X$ is finite locally
free of constant rank~$m$.
 Let $\cA=\cA_X(\g,\widetilde\g)$ and
$\cA^\circ=\cA_X(\g,\widetilde\g^\circ)$ be the two related twisted
universal enveloping quasi-coherent quasi-algebras, and let
$\cB^\cu=\cC_X^\cu(\g,\widetilde\g)$ and
$\cB^\circ{}^\cu=\cC_X^\cu(\g,\widetilde\g^\circ)$ be the two related
Chevalley--Eilenberg quasi-coherent CDG\+quasi-algebras over~$X$.
 Then there are natural commutative diagrams of triangulated category
equivalences
\begin{equation} \label{main-X-lct-hexagonality-diagram}
\begin{gathered}
 \xymatrix{
  \sD(\Qcohr\cA^\circ)
  \ar@<3pt>[rrrr]^{{-}\ot_{\cO_X}\bigwedge\nolimits_X^*(\g)}
  \ar@<3pt>[dd]^{\cHom_{\cO_X}(\cB^m[-m],{-})}
  &&&& \sD^\co_{X\red}(\bQcohr\cB^\circ{}^\cu)
  \ar@<3pt>[llll]^{{-}\ot_{\cO_X}\cA^\circ}
  \ar@{=}[dd]^{\cB^{\circ,\rop}{}^\subcu\simeq\cB^\subcu}
  \\ \\
  \sD(\cA\Qcoh) \ar@<3pt>[rrrr]^{\cB^*\ot_{\cO_X}{-}}
  \ar@<3pt>[uu]^{\cB^m[-m]\ot_{\cO_X}{-}}
  \ar@<3pt>[dd]^{\boR\fHom_\cA(\cA,{-})}
  &&&& \sD^\co_{X\red}(\cB^\cu\bQcoh)
  \ar@{-}@<3pt>[llll]
  \ar@<3pt>[dd]^{\boR\fHom_{\cB^*}(\cB^\subcu,{-})} \\ \\
  \sD(\cA\Lcth_\bW^{X\dlct})
  \ar@<3pt>[rrrr]^{\Cohom_X\left(\bigwedge\nolimits_X^*(\g),{-}\right)}
  \ar@<3pt>[uu]^{\cA\ocn_\cA^\boL{-}}
  &&&& \sD^\ctr_{X\red}(\cB^\cu\bLcth_\bW^{X\dlct})
  \ar@<3pt>[llll]^{\Cohom_X(\cA,{-})}
  \ar@<3pt>[uu]^{\cB^\subcu\ocn_{\cB^*}^\boL{-}}
 }
\end{gathered}
\end{equation}
\begin{equation} \label{main-A-lct-hexagonality-diagram}
\begin{gathered}
 \xymatrix{
  \sD(\Qcohr\cA^\circ)
  \ar@<3pt>[rrrr]^{{-}\ot_{\cO_X}\bigwedge\nolimits_X^*(\g)}
  \ar@<3pt>[dd]^{\cHom_{\cO_X}(\cB^m[-m],{-})}
  &&&& \sD^\co_{X\red}(\bQcohr\cB^\circ{}^\cu)
  \ar@<3pt>[llll]^{{-}\ot_{\cO_X}\cA^\circ}
  \ar@{=}[dd]^{\cB^{\circ,\rop}{}^\subcu\simeq\cB^\subcu}
  \\ \\
  \sD(\cA\Qcoh) \ar@<3pt>[rrrr]^{\cB^*\ot_{\cO_X}{-}}
  \ar@<3pt>[uu]^{\cB^m[-m]\ot_{\cO_X}{-}}
  \ar@<3pt>[dd]^{\boR\fHom_\cA(\cA,{-})}
  &&&& \sD^\co_{X\red}(\cB^\cu\bQcoh)
  \ar@{-}@<3pt>[llll]
  \ar@<3pt>[dd]^{\boR\fHom_{\cB^*}(\cB^\subcu,{-})} \\ \\
  \sD(\cA\Lcth_\bW^{\cA\dlct})
  \ar@<3pt>[rrrr]^{\Cohom_X\left(\bigwedge\nolimits_X^*(\g),{-}\right)}
  \ar@<3pt>[uu]^{\cA\ocn_\cA^\boL{-}}
  &&&& \sD^\ctr_{X\red}(\cB^\cu\bLcth_\bW^{X\dlct})
  \ar@<3pt>[llll]^{\Cohom_X(\cA,{-})}
  \ar@<3pt>[uu]^{\cB^\subcu\ocn_{\cB^*}^\boL{-}}
 }
\end{gathered}
\end{equation}
\end{cor}

\begin{proof}
 The upper commutative square of triangulated equivalences in both
the diagrams is provided by
Corollary~\ref{left-right-co-side-reduced-koszul-duality}.
 The lower commutative squares of triangulated equivalences are
the results of Theorem~\ref{main-lct-quadrality-theorem}.
\end{proof}

 Once again we recall that, in the situation at hand, the classes
of $X$\+locally cotorsion and $\cB^*$\+locally cotorsion
$\bW$\+locally contraherent graded modules over $\cB^*$ coincide
by Corollary~\ref{lct-lcth-modules-over-flfrqa-cor} (see also
Corollary~\ref{cot-qcoh-modules-over-flfrqa-cor} for a quasi-coherent
version).
 This is the reason why $\cB^*$\+locally cotorsion $\bW$\+locally
contraherent CDG\+modules over $\cB^\cu$ are not mentioned in
the statement of the main results of this paper in the theorem
and corollary above.

\begin{rem} \label{becker=positselski-at-the-end-remark}
 In the commutative diagrams of
Corollaries~\ref{main-hexagonality-corollary}
and~\ref{main-lct-hexagonality-corollary}, the leftmost columns
consist of the conventional derived categories, while the rightmost
columns are formed of the reduced versions of coderived and
contraderived categories.
 In the context of the latter, the distinction between the Positselski
and Becker definitions might become relevant.
 As the corollaries (and the preceding
Theorems~\ref{main-quadrality-theorem}
and~\ref{main-lct-quadrality-theorem}) are stated, they involve
the reduced versions of the Positselski co/contraderived categories.
 Let us explain why there is no actual difference.

 The definitions of the reduced coderived and reduced contraderived
categories appearing in the diagrams~\eqref{main-hexagonality-diagram},
\eqref{main-X-lct-hexagonality-diagram},
\eqref{main-A-lct-hexagonality-diagram}
were given in Sections~\ref{reduced-coderived-of-qcoh-subsecn}
and~\ref{reduced-contraderived-of-lcth-subsecn}
in the context of the Positselski co/contraderived categories only.
 We have \emph{never} defined reduced Becker co/contraderived
categories.
 The reason was that the main results of
Sections~\ref{reduced-coderived-of-qcoh-subsecn}
and~\ref{reduced-contraderived-of-lcth-subsecn}, viz.,
Theorems~\ref{qcoh-reduced-coderived-category-equivalence-thm}
and~\ref{lcth-reduced-contrader-category-equivalences-thm},
we were unable to prove in the context of the Becker definitions,
for the reasons explained in
Remarks~\ref{proving-becker-coacyclicity-preserved-remark}
and~\ref{proving-becker-contraacyclicity-preserved-remark}.
 Basically, the two theorems established the comparison of
the reduced co/contraderived categories of arbitrary CDG\+modules
with those of \emph{thick} CDG\+modules, and the respective
definitions for the thick CDG\+modules were difficult to work out
in the Becker context.

 Nevertheless, the reduced Becker co/contraderived categories of
arbitrary (\emph{not} necessarily thick) CDG\+modules, which would be
denoted by $\sD^\bco_{X\red}(\bQcohr\cB^\circ{}^\cu)$, \
$\sD^\bco_{X\red}(\cB^\cu\bQcoh)$, \
$\sD^\bctr_{X\red}(\cB^\cu\bLcth_\bW)$, and
$\sD^\bctr_{X\red}(\cB^\cu\bLcth_\bW^{X\dlct})$, can be easily
defined in the way very similar to the definitions in
Sections~\ref{reduced-coderived-of-qcoh-subsecn}
and~\ref{reduced-contraderived-of-lcth-subsecn}.
 Then the claim is that such reduced Becker co/contraderived
categories \emph{coincide} with the respective Positselski versions.
 This is provable using the commutative
diagrams~\eqref{main-hexagonality-diagram},
\eqref{main-X-lct-hexagonality-diagram},
\eqref{main-A-lct-hexagonality-diagram}, or more precisely,
the results of the reduced Koszul duality theorems,
Theorems~\ref{reduced-koszul-duality-right-co-side}
and~\ref{reduced-koszul-duality-contra-side}, by comparing them with
the Becker versions of the semiderived Koszul duality theorems,
Theorems~\ref{semiderived-koszul-duality-right-co-side}(b)
and~\ref{semiderived-koszul-duality-becker-contra-side}.

 The point is that all Positselski-coacyclic objects are always
Becker-coacyclic, and similarly for the contraacyclicity.
 So, in order to prove our claims that, in the respective versions
of the notation, $\Ac^\bco_{X\red}(\bQcohr\cB^\cu)=
\Ac^\co_{X\red}(\bQcohr\cB^\cu)$, \
$\Ac^\bctr_{X\red}(\cB^\cu\bLcth_\bW)=
\Ac^\ctr_{X\red}(\cB^\cu\bLcth_\bW)$, and
$\Ac^\bctr_{X\red}(\cB^\cu\bLcth_\bW^{X\dlct})=
\Ac^\ctr_{X\red}(\cB^\cu\bLcth_\bW^{X\dlct})$, it suffices to
check that the horizontal functors going leftwards in
the diagrams~\eqref{main-hexagonality-diagram},
\eqref{main-X-lct-hexagonality-diagram},
\eqref{main-A-lct-hexagonality-diagram} annihilate
the Becker-coacyclic and Becker-contraacyclic CDG\+modules
over~$\cB^\cu$.
 But these assertions are a part of the results of
Theorems~\ref{semiderived-koszul-duality-right-co-side}(b)
and~\ref{semiderived-koszul-duality-becker-contra-side}.
 We recall that all the Becker-semicoacyclic and
Becker-semicontraacyclic complexes of $\cA$\+modules in
the respective contexts are acyclic, as per the discussions
in Sections~\ref{semiderived-quasi-coherent-subsecn}
and~\ref{becker-semicontraderived-defined-subsecn}.
\emergencystretch=3em \hbadness=8500
\end{rem}

\begin{rem}
 Following the discussion in
Remarks~\ref{naive-co-contra-applies-to-diffoperators-remark},
\ref{weakly-smooth-de-Rham-remark},
and~\ref{weakly-smooth-de-Rham-and-diffoperators-remark},
the results of Theorems~\ref{main-quadrality-theorem},
\ref{main-lct-quadrality-theorem} and
Corollaries~\ref{main-hexagonality-corollary},
\ref{main-lct-hexagonality-corollary} are applicable, in
particular, to the sheaf of rings of crystalline fiberwise differential
operators $\cA=\cD^\cry_{X/T}$ and the sheaf of de~Rham DG\+algebras of
fiberwise differential forms $\cB^\cu=\Omega^\bu_{X/T}$.
 Here $X\rarrow T$ is a weakly smooth morphism of schemes of constant
relative dimension~$m$, and the scheme $X$ is assumed to be
quasi-compact and semi-separated.
 In this context, the results of Theorems~\ref{main-quadrality-theorem}
and/or~\ref{main-lct-quadrality-theorem} can be called
the ``$\cD$\+$\Omega$ quadrality theorem''.
\end{rem}

\Section{Special Classes of Particularly Nice Schemes}
\label{well-behaved-schemes-secn}

 In this section we discuss some improvements and simplifications
of the main results of the preceding sections that can be obtained
under additional assumtions on the scheme~$X$.
 There are two special classes of schemes which we consider:
the semi-separated Noetherian schemes of finite Krull dimension
(this includes all algebraic varieties over fields), and
the semi-separated regular Noetherian schemes of finite Krull
dimension (this includes all smooth algebraic varieties over fields).

\subsection{Finite cotorsion and locally cotorsion coresolution
dimensions}
 We start with an equivalent reformulation of a classical result of
Raynaud and Gruson~\cite[Corollaire~II.3.2.7]{RG}.

\begin{lem} \label{wrt-R-cotorsion-A-modules-coresol-dim}
 Let $R$ be a commutative Noetherian ring of finite Krull dimension~$D$.
 Let $A$ be an associative ring and $R\rarrow A$ be a ring homomorphism.
 Then all $A$\+modules have finite coresolution dimension not
exceeding $D$ with respect to the coresolving subcategory of
$R$\+cotorsion $A$\+modules $A\Modl^{R\dcot}\subset A\Modl$.
\end{lem}

\begin{proof}
 The full subcategory $A\Modl^{R\dcot}$ is coresolving in $A\Modl$ by
Lemmas~\ref{flat-cotorsion-pair-hereditary}(b)
and~\ref{restriction-coextension-injective-cotorsion}(a).
 The bound on the coresolution dimension follows from the facts that
the coresolution dimension does not depend on the choice of
a coresolution (see Section~\ref{co-resolution-dimensions-subsecn})
and that all $R$\+modules have coresolution dimensions not exceeding $D$
with respect to the coresolving subcategory of cotorsion $R$\+modules
$R\Modl^\cot\subset R\Modl$ \,\cite[Corollary~1.5.7]{Pcosh}.
\end{proof}

 The following corollary should be compared
with~\cite[Lemma~6.2.5(b)]{Pcosh}, which does not assume the scheme
to be semi-separated.

\begin{cor} \label{A-lcth-al-X-lct-coresolution-dimension}
 Let $X$ be a semi-separated Noetherian scheme of finite Krull
dimension $D$ with an open covering\/~$\bW$.
 Let $\cA$ be a quasi-coherent quasi-algebra over~$X$.
 Then \par
\textup{(a)} all\/ $\bW$\+locally contraherent $\cA$\+modules have
finite coresolution dimensions not exceeding $D$ with respect to
the coresolving subcategory of $X$\+locally cotorsion $\bW$\+locally
contraherent $\cA$\+modules $\cA\Lcth_\bW^{X\dlct}\subset\cA\Lcth_\bW$;
\par
\textup{(b)} all antilocal contraherent $\cA$\+modules have finite
coresolution dimensions not exceeding $D$ with respect to
the coresolving subcategory of antilocal $X$\+locally cotorsion
contraherent $\cA$\+modules $\cA\Ctrh^{X\dlct}_\al\subset\cA\Ctrh_\al$.
\end{cor}

\begin{proof}
 Any $\bW$\+locally contraherent $\cA$\+module has an admissible
monomorphism into an $X$\+locally cotorsion, and in fact, even
into an $\cA$\+locally injective $\bW$\+locally contraherent
$\cA$\+module by Lemma~\ref{qcomp-qsep-A-loc-inj-preenvelope}.
 Moreover, the cokernel of such an admissible monomorphism can be
made antilocal by
Theorem~\ref{qcomp-qsep-antilocal-complete-cotorsion-pair-thm}(b).
 In view of these observations, it follows from
Lemmas~\ref{flat-cotorsion-pair-hereditary}(b) that the full
subcategory $\cA\Lcth_\bW^{X\dlct}$ is coresolving in
$\cA\Lcth_\bW$ and the full subcategory $\cA\Ctrh^{X\dlct}_\al$
is coresolving in $\cA\Ctrh_\al$.

 Now the bounds on the coresolution dimension follow from the fact
that the coresolution dimension does not depend on the choice of
the coresolution together with
Lemma~\ref{wrt-R-cotorsion-A-modules-coresol-dim}.
\end{proof}

 The following lemma builds on top of
Lemma~\ref{X-A-cta-cot-coresolving-coresolution-dimension}.

\begin{lem} \label{A-qcoh-X-cot-coresolution-dimension}
 Let $X$ be a semi-separated Noetherian scheme of finite Krull
dimension $D$ with a finite affine open covering
$X=\bigcup_{\alpha=1}^N U_\alpha$.
 Let $\cA$ be a quasi-coherent quasi-algebra over~$X$.
 Then all quasi-coherent $\cA$\+modules have finite coresolution
dimension not exceeding $N-1+D$ with respect to the coresolving
subcategory of $X$\+cotorsion $\cA$\+modules
$\cA\Qcoh^{X\dcot}\subset\cA\Qcoh$.
\end{lem}

\begin{proof}
 As in the proofs of
Lemmas~\ref{X-A-cta-cot-coresolving-coresolution-dimension}(b)
and~\ref{wrt-R-cotorsion-A-modules-coresol-dim}, one needs to use
the fact that the coresolution dimension does not depend on the choice
of a coresolution.
 Then the assertion follows from the fact that all quasi-coherent
sheaves on $X$ have finite coresolution dimension not exceeding $N-1+D$ with respect to the full subcategory of cotorsion quasi-coherent sheaves
$X\Qcoh^\cot\subset X\Qcoh$ \,\cite[Lemma~6.4.1(c)]{Pcosh}.
\end{proof}

\subsection{Na\"\i ve co-contra correspondence via $X$-cotorsion
sheaves and $X$-lo\-cally cotorsion cosheaves}
 This section complements the exposition in
Section~\ref{naive-co-contra-secn}.

\begin{cor} \label{qcoh-X-cot-derived-equivalence}
 Let $X$ be a semi-separated Noetherian scheme of finite Krull dimension
and $\cA$ be a quasi-coherent quasi-algebra over~$X$.
 Then, for any conventional derived category symbol\/ $\st=\bb$, $+$,
$-$, or\/~$\varnothing$, the inclusions of exact/abelian categories
$\cA\Qcoh^{X\dcot}\rarrow\cA\Qcoh^{X\dcta}\rarrow\cA\Qcoh$ induce
equivalences of the derived categories
$$
 \sD^\st(\cA\Qcoh^{X\dcot})\simeq
 \sD^\st(\cA\Qcoh^{X\dcta})\simeq\sD^\st(\cA\Qcoh).
$$
\end{cor}

\begin{proof}
 The equivalence $\sD^\st(\cA\Qcoh^{X\dcta})\simeq\sD^\st(\cA\Qcoh)$
holds in greater generality by
Corollary~\ref{qcoh-X-cta-derived-equivalence}.
 The equivalence $\sD^\st(\cA\Qcoh^{X\dcot})\simeq\sD^\st(\cA\Qcoh)$
follows from Lemma~\ref{A-qcoh-X-cot-coresolution-dimension} in view of
the dual version of Proposition~\ref{finite-resolutions}.
\end{proof}

\begin{cor} \label{A-lcth-al-X-lct-X-lcta-derived-equivalences}
 Let $X$ be a semi-separated Noetherian scheme of finite Krull
dimension with an open covering\/ $\bW$, and let $\cA$ be
a quasi-coherent quasi-algebra over~$X$.
 Let\/ $\st=\bb$, $+$, $-$, or\/~$\varnothing$ be a conventional
derived category symbol.
 Then \par
\textup{(a)} the inclusion of exact categories\/ $\cA\Lcth_\bW^{X\dlct}
\rarrow\cA\Lcth_\bW$ induces an equivalence of the derived categories
$$
 \sD^\st(\cA\Lcth_\bW^{X\dlct})\simeq\sD^\st(\cA\Lcth_\bW);
$$ \par
\textup{(b)} the inclusion of exact categories\/ $\cA\Ctrh^{X\dlct}_\al
\rarrow\cA\Ctrh_\al$ induces an equivalence of the derived categories
$$
 \sD^\st(\cA\Ctrh^{X\dlct}_\al)\simeq\sD^\st(\cA\Ctrh_\al).
$$
\end{cor}

\begin{proof}
 The assertions follow from the respective parts of
Corollary~\ref{A-lcth-al-X-lct-coresolution-dimension} in view of
the dual version of Proposition~\ref{finite-resolutions}.
\end{proof}

 The following theorem partly extends the result of
Theorem~\ref{quasi-algebra-unbounded-derived-naive-co-contra} to
the (half\+)boun\-ded derived category symbols $\st=\bb$ and~$-$.

\begin{thm} \label{quasi-algebra-cot-bounded-derived-naive-co-contra}
 Let $X$ be a semi-separated scheme Noetherian scheme of finite Krull
dimension with an open covering\/ $\bW$, and let $\cA$ be
a quasi-coherent quasi-algebra over~$X$.
 Then, for any conventional derived category symbol\/ $\st=\bb$, $+$,
$-$, or\/~$\varnothing$ there is a commutative diagram of natural
triangulated equivalences between the derived categories of various
abelian/exact categories of quasi-coherent and\/ $\bW$\+locally
contraherent $\cA$\+modules on~$X$,
\begin{equation} \label{cot-bounded-derived-naive-co-contra-diagram}
\begin{gathered}
 \xymatrix{
  \sD^\st(\cA\Qcoh)\ar@{-}@<2pt>[r]
  & \sD^\st(\cA\Qcoh^{X\dcta}) \ar@<2pt>[l] \ar@{=}[r] \ar@{-}@<2pt>[d]
  & \sD^\st(\cA\Ctrh_\al) \ar@<-2pt>[r] \ar@{-}@<-2pt>[d]
  & \sD^\st(\cA\Lcth_\bW) \ar@{-}@<-2pt>[l] \ar@{-}@<-2pt>[d] \\
  & \sD^\st(\cA\Qcoh^{X\dcot}) \ar@{=}[r] \ar@<2pt>[u]
  & \sD^\st(\cA\Ctrh_\al^{X\dlct}) \ar@<-2pt>[r] \ar@<-2pt>[u]
  & \sD^\st(\cA\Lcth_\bW^{X\dlct}) \ar@{-}@<-2pt>[l] \ar@<-2pt>[u]
 }
\end{gathered}
\end{equation}
\end{thm}

\begin{proof}
 The two triangulated equivalences shown by double horizontal lines
without arrows (in the middle column of equivalences) are induced by
the equivalences of exact categories provided by
Lemma~\ref{quasi-algebra-underived-naive-co-contra}(a\+-b).
 All the other triangulated functors (shown by arrows) are induced by
the respective exact inclusions of exact/abelian categories.

 The upper horizontal line of triangulated equivalences is provided by
Theorem~\ref{quasi-algebra-bounded-derived-naive-co-contra}.
 The triangulated equivalences in the left-hand part of the diagram
are the result of Corollary~\ref{qcoh-X-cot-derived-equivalence}.
 The horizontal triangulated equivalences in the right-hand part of
the diagram are the result of
Corollary~\ref{A-lcth-ctrh-al-derived-equivalences}(a\+-b).
 The vertical triangulated equivalences in the right-hand part of
the diagram are the result of
Corollary~\ref{A-lcth-al-X-lct-X-lcta-derived-equivalences}.
\end{proof}

\subsection{Antilocally flaprojective contraherent $\cA$-modules}
\label{antilocally-flaprojective-subsecn}
 This section, in which the scheme $X$ is only assumed to be
quasi-compact and semi-separated, is an $\cA$\+module
version/generalization of~\cite[Section~4.4]{Pcosh}, complementing
the exposition in Section~\ref{antilocal-classes-secn}.
 As in Section~\ref{becker-semicontraderived-secn}, in this section
and in Section~\ref{becker-contraderived-finite-coresolutions-subsecn}
below we use more knowledge of the machinery of cotorsion pairs than
in the rest of this paper.
 The results of this section will be used in
Section~\ref{co-contra-derived-loc-cotorsion-subsecn}.
{\hbadness=1225\par}

 Let $X$ be a scheme with an open covering $\bW$ and $\cA$ be
a quasi-coherent quasi-algebra over~$X$.
 We will say that a $\bW$\+locally contraherent $\cA$\+module $\gF$
on $X$ is \emph{antilocally flaprojective} if the functor
$\Hom^\cA(\gF,{-})$ takes short exact sequences of $X$\+locally
cotorsion $\bW$\+locally contraherent $\cA$\+modules on $X$ to
short exact sequences of abelian groups.

\begin{lem} \label{antilocally-flaprojective-on-affines}
 Let $U$ be an affine scheme with the open covering\/ $\bW=\{U\}$ and
$\cA$ be a quasi-coherent quasi-algebra over~$U$.
 Then a contraherent $\cA$\+module\/ $\gF$ on $U$ is antilocally
flaprojective if and only if the $\cA(U)$\+module\/ $\gF[U]$ is
$\cA(U)/\cO(U)$\+flaprojective.
\end{lem}

\begin{proof}
 Put $R=\cO(U)$ and $A=\cA(U)$.
 The rule $\P\longmapsto\P[U]$ establishes an equivalence between
the exact category of contraherent $\cA$\+modules on $U$ and
the exact category of $R$\+contraadjusted left $A$\+modules,
$\cA\Ctrh\simeq A\Modl^{R\dcta}$ (see
Section~\ref{cosheaves-of-A-modules-subsecn}).
 The $U$\+locally cotorsion contraherent $\cA$\+modules on $X$
correspond to $R$\+cotorsion $A$\+modules.
 Put $F=\gF[U]$.
 Using the facts that all injective $A$\+modules are $R$\+cotorsion
and the cokernels of injective morphisms of $R$\+cotorsion $A$\+modules
are $R$\+cotorsion (Lemmas~\ref{flat-cotorsion-pair-hereditary}(b)
and~\ref{restriction-coextension-injective-cotorsion}), one can show
easily that the contraherent $\cA$\+module $\gF$ is antilocally
flaprojective if and only if $\Ext_A^1(F,C)=0$ for all
$R$\+cotorsion $A$\+modules~$C$.
 By the definition, this means precisely that the $A$\+module $F$ is
$A/R$\+flaprojective (see Section~\ref{prelim-flaprojective-subsecn}).
\end{proof}

\begin{lem} \label{antilocally-flaprojective-affine-direct-image}
 Let $X$ be a scheme with an open covering\/ $\bW$ and $\cA$ be
a quasi-coherent quasi-algebra over~$X$.
 Let $Y\subset X$ be an open subscheme such that the open immersion
morphism $j\:Y\rarrow X$ is affine.
 Then the direct image functor $j_!\:\cA|_Y\Lcth_{\bW|_Y}\rarrow
\cA\Lcth_\bW$ takes antilocally flaprojective\/ $\bW|_Y$\+locally
contraherent $\cA|_Y$\+modules on $Y$ to antilocally flaprojective\/
$\bW$\+locally contraherent $\cA$\+modules on~$X$.
\end{lem}

\begin{proof}
 The assertion follows immediately from
the adjunction~\eqref{affine-open-immers-contrah-module-adjunction}
together with the fact that the restrictions of locally cotorsion
contraherent cosheaves to open subschemes are locally cotorsion
contraherent cosheaves again.
\end{proof}

\begin{lem} \label{qcomp-qsep-X-loc-cot-preenvelope}
 Let $X$ be a quasi-compact semi-separated scheme with an open
covering\/ $\bW$ and $X=\bigcup_\alpha U_\alpha$ be a finite affine
open covering of $X$ subordinate to\/~$\bW$.
 Let $\cA$ be a quasi-coherent quasi-algebra over $X$ and\/ $\gM$ be
a\/ $\bW$\+locally contraherent $\cA$\+module on~$X$.
 Then there exists a short exact sequence of\/ $\bW$\+locally
contraherent $\cA$\+modules\/ $0\rarrow\gM\rarrow\Q\rarrow\gF\rarrow0$
on $X$ such that\/ $\Q$ is an $X$\+locally cotorsion\/
$\bW$\+locally contraherent $\cA$\+module and\/ $\gF$ is a finitely
iterated extension of the direct images of antilocally flaprojective
contraherent $\cA|_{U_\alpha}$\+modules from $U_\alpha$ in the exact
category $\cA\Lcth_\bW$ (or $\cA\Ctrh$).
\end{lem}

\begin{proof}
 This is another version of~\cite[Proposition~5.2]{Pal}
and~\cite[Lemma~4.4.1 or Proposition~B.6.4]{Pcosh}.
 The same construction as in
Lemmas~\ref{qcomp-qsep-A-loc-inj-preenvelope}
and~\ref{qcomp-qsep-X-lct-A-loc-inj-preenvelope} is applicable in
the situation at hand.
 The local nature of the local cotorsion property of locally
contraherent cosheaves (see~\cite[Lemma~1.3.6(a)]{Pcosh} and
Section~\ref{locally-contraherent-cosheaves-subsecn} above) plays
a key role.
 It is also important that the class of locally cotorsion
$\bW$\+locally contraherent cosheaves is preserved by the direct images
with respect to $(\bW,\bT)$\+affine open immersions of schemes (see
Section~\ref{direct-images-of-O-co-sheaves-subsecn}), as well as
by extensions and cokernels of admissible monomorphisms in $X\Lcth_\bW$.
 Of course, one also has to use the result of
Theorem~\ref{flaprojective-cotorsion-pair-complete}(b) together with
Lemma~\ref{antilocally-flaprojective-on-affines} in order to resolve
a contraherent cosheaf on an affine open subscheme $U\subset X$.
\end{proof}

 The notation $\Ext^{\cA,*}({-},{-})$ appearing in the following
corollary was introduced in
Section~\ref{A-antilocal=X-antilocal-subsecn} and explained in
Remark~\ref{Ext-unambiguous-remark}.

\begin{cor} \label{antilocally-flaprojective-hereditary}
 Let $X$ be a quasi-compact semi-separated scheme with an open
covering\/ $\bW$, and $\cA$ be a quasi-coherent quasi-algebra over~$X$.
 Then \par
\textup{(a)} A $\bW$\+locally contraherent $\cA$\+module\/ $\gF$ on $X$
is antilocally flaprojective if and only if\/ $\Ext^{\cA,1}(\gF,\Q)=0$
for all $X$\+locally cotorsion\/ $\bW$\+locally contraherent
$\cA$\+modules\/ $\Q$ on~$X$. \par
\textup{(b)} A $\bW$\+locally contraherent $\cA$\+module\/ $\gF$ on $X$
is antilocally flaprojective if and only if\/ $\Ext^{\cA,n}(\gF,\Q)=0$
for all $X$\+locally cotorsion\/ $\bW$\+locally contraherent
$\cA$\+modules\/ $\Q$ on $X$ and all $n\ge1$. \par
\textup{(c)} The class of antilocally flaprojective\/ $\bW$\+locally
contraherent $\cA$\+modules on $X$ is closed under extensions and
kernels of admissible epimorphisms in $\cA\Lcth_\bW$.
\end{cor}

\begin{proof}
 This is a generalization of~\cite[Corollary~4.4.2]{Pcosh}.
 Parts~(a\+-b) follow from the fact that every $\bW$\+contraherent
$\cA$\+module has an admissible monomorphism into an $X$\+locally
cotorsion $\bW$\+locally contraherent $\cA$\+module (which is a weak
form of Lemma~\ref{qcomp-qsep-X-loc-cot-preenvelope} and follows also
from Lemma~\ref{qcomp-qsep-A-loc-inj-preenvelope}).
 One also has to use the fact that the class of locally cotorsion
$\bW$\+locally cosheaves on $X$ is closed under the cokernels of
admissible monomorphisms in $X\Lcth_\bW$.
 The relevant argument can be found in~\cite[Lemma~7.1]{PS6} (based
on~\cite[proof of Lemma~6.17]{Sto-ICRA}), \cite[Lemma~1.3]{BHP},
or~\cite[Lemma~1.1(b)]{Pal}.
 Part~(c) follows from parts~(a\+-b).
\end{proof}

\begin{lem} \label{antilocally-flaprojective-modules-generating-class}
 Let $X$ be a quasi-compact semi-separated scheme with an open
covering\/ $\bW$ and $X=\bigcup_\alpha U_\alpha$ be a finite affine
open covering of $X$ subordinate to\/~$\bW$.
 Denote by $j_\alpha\:U_\alpha\rarrow X$ the open immersion morphisms.
 Let $\cA$ be a quasi-coherent quasi-algebra over $X$ and\/ $\gM$ be
a\/ $\bW$\+locally contraherent $\cA$\+module on~$X$.
 Then there exists an admissible epimorphism onto\/ $\gM$ in
$\cA\Lcth_\bW$ from a finite direct sum\/ $\bigoplus_\alpha
j_\alpha{}_!\gF_\alpha$ for some antilocally flaprojective contraherent
$\cA|_{U_\alpha}$\+modules\/ $\gF_\alpha$ on~$U_\alpha$.
\end{lem}

\begin{proof}
 This is a stronger version of
Lemma~\ref{antilocal-modules-generating-class}, provable
in a similar way.
 Put $R_\alpha=\cO_X(U_\alpha)$ and $A_\alpha=\cA(U_\alpha)$.
 Using Theorem~\ref{flaprojective-cotorsion-pair-complete}(a)
or~\ref{very-flaprojective-cotorsion-pair-complete}(a), one can
construct for any $R_\alpha$\+contraadjusted $A_\alpha$\+module
$M_\alpha$ a short exact sequence of $A_\alpha$\+modules
$0\rarrow C_\alpha\rarrow F_\alpha\rarrow M_\alpha\rarrow0$ with
an $A_\alpha/R_\alpha$\+flaprojective $A_\alpha$\+module $F_\alpha$
and an $R_\alpha$\+contraadjusted $A_\alpha$\+module~$C_\alpha$.
 Then it follows that the $A_\alpha$\+module $F_\alpha$ is
$R_\alpha$\+contraadjusted as well.
 This allows one to produce, for every index~$\alpha$, an admissible
epimorphism $\gF_\alpha\rarrow j_\alpha^!\gM$ in the exact category
$\cA|_{U\alpha}\Ctrh$ with an antilocally flaprojective
contraherent $\cA|_{U_\alpha}$\+module~$\gF_\alpha$.
 Now the composition of admissible epimorphisms
$\bigoplus_\alpha j_\alpha{}_!\gF_\alpha\rarrow
\bigoplus_\alpha j_\alpha{}_!j_\alpha^!\gM\rarrow\gM$ provides
the desired admissible epimorphism
$\bigoplus_\alpha j_\alpha{}_!\gF_\alpha\rarrow\gM$.
\end{proof}

\begin{thm} \label{antilocally-flaproj-complete-cotorsion-pair-thm}
 Let $X$ be a quasi-compact semi-separated scheme with an open
covering\/ $\bW$ and $X=\bigcup_\alpha U_\alpha$ be a finite affine
open covering of $X$ subordinate to\/~$\bW$.
 Let $\cA$ be a quasi-coherent quasi-algebra over $X$ and\/ $\gM$ be
a\/ $\bW$\+locally contraherent $\cA$\+module on~$X$.
 In this context: \par
\textup{(a)} there exists a short exact sequence of\/ $\bW$\+locally
contraherent $\cA$\+modules\/ $0\rarrow\Q'\rarrow\gF\rarrow\gM\rarrow0$
on $X$ such that\/ $\gF$ is an antilocally flaprojective\/
$\bW$\+locally contraherent $\cA$\+module and\/ $\Q'$ is an $X$\+locally
cotorsion\/ $\bW$\+locally contraherent $\cA$\+module; \par
\textup{(b)} there exists a short exact sequence of\/ $\bW$\+locally
contraherent $\cA$\+modules\/ $0\rarrow\gM\rarrow\Q\rarrow\gF'\rarrow0$
on $X$ such that\/ $\Q$ is an $X$\+locally cotorsion\/ $\bW$\+locally
contraherent $\cA$\+module and\/ $\gF'$ is an antilocally
flaprojective\/ $\bW$\+locally contraherent $\cA$\+module; \par
\textup{(c)} a\/ $\bW$\+locally contraherent $\cA$\+module on $X$ is
antilocally flaprojective if and only if it is a direct summand of
a finitely iterated extension of the direct images of antilocally
flaprojective contraherent $\cA|_{U_\alpha}$\+modules from $U_\alpha$
in the exact category $\cA\Lcth_\bW$ or $\cA\Ctrh$.
\end{thm}

\begin{proof}
 This is a generalization of~\cite[Corollary~4.4.4]{Pcosh}.
 The proof is similar to those of
Theorems~\ref{qcomp-qsep-antilocal-complete-cotorsion-pair-thm}
and~\ref{qcomp-qsep-X-lct-aloc-complete-cotorsion-pair-thm}.
 Part~(b) follows from Lemma~\ref{qcomp-qsep-X-loc-cot-preenvelope}
in view of Lemma~\ref{antilocally-flaprojective-affine-direct-image}
and Corollary~\ref{antilocally-flaprojective-hereditary}(c).
 Part~(a) follows from part~(b) and
Lemma~\ref{antilocally-flaprojective-modules-generating-class} by
virtue of the Salce lemma.
 Part~(c) is deduced from the more precise version of part~(a) provided
by its proof based on Lemma~\ref{qcomp-qsep-X-loc-cot-preenvelope}
(one also needs to use
Corollary~\ref{antilocally-flaprojective-hereditary}(a)).
\end{proof}

 It follows from
Theorem~\ref{antilocally-flaproj-complete-cotorsion-pair-thm} that
all antilocally flaprojective $\bW$\+locally contraherent $\cA$\+modules
on a quasi-compact semi-separated scheme $X$ are globally contraherent
(and antilocal), and moreover, the class of all such contraherent
$\cA$\+modules does not depend on the open covering~$\bW$.
 So we will speak simply of \emph{antilocally flaprojective contraherent
$\cA$\+modules} on $X$, and denote the class of such contraherent
$\cA$\+modules by $\cA\Ctrh_\alfp\subset\cA\Lcth_\bW$.

\begin{cor} \label{antilocally-flaprojective-cotorsion-pair-cor}
 Let $X$ be a quasi-compact semi-separated scheme with an open
covering\/ $\bW$ and $\cA$ be a quasi-coherent quasi-algebra over~$X$.
 Then the pair of classes of antilocally flaprojective contraherent
$\cA$\+modules\/ $\sF=\cA\Ctrh_\alfp$ and $X$\+locally cotorsion\/
$\bW$\+locally contraherent $\cA$\+modules\/
$\sC=\cA\Lcth^{X\dlct}_\bW$ is a hereditary complete cotorsion pair
$(\sF,\sC)$ in the exact category $\cA\Lcth_\bW$.
\end{cor}

\begin{proof}
 This is similar to Remarks~\ref{antilocal-A-modules-concluding-remark}
and~\ref{antilocal-X-lct-A-modules-concluding-remark}.
 The assertion follows from
Theorem~\ref{antilocally-flaproj-complete-cotorsion-pair-thm}(a\+-b)
and Corollary~\ref{antilocally-flaprojective-hereditary} by
virtue of the direct summand lemma~\cite[Lemma~B.1.2]{Pcosh}.
 One needs to use the fact that the class $\cA\Lcth^{X\dlct}_\bW$ is
closed under direct summands in $\cA\Lcth_\bW$ (since the class
$X\Lcth^\lct$ is closed under direct summands in $X\Lcth$).
\end{proof}

 The assertion of
Theorem~\ref{antilocally-flaproj-complete-cotorsion-pair-thm}(c)
provides an additional piece of information on top of the result of
Corollary~\ref{antilocally-flaprojective-cotorsion-pair-cor}.
 This additional aspect can be expressed by saying that the class of
all antilocally flaprojective contraherent $\cA$\+modules on $X$ is
\emph{antilocal} in the sense of~\cite[Section~4]{Pal}.

 For our purposes below in
Section~\ref{co-contra-derived-loc-cotorsion-subsecn}, the important
conclusion from the results of this
Section~\ref{antilocally-flaprojective-subsecn} is that the class
$\cA\Lcth_\bW^{X\dlct}$ is the right-hand side of a (hereditary)
complete cotorsion pair in $\cA\Lcth_\bW$.

\subsection{Exact DG-pairs} \label{exact-dg-pairs-subsecn}
 The aim of this section is to formulate a generalization
of~\cite[Theorem~5.6]{Pedg} promised in~\cite[Remark~5.12]{Pedg}.
 We will need it in the next
Section~\ref{becker-contraderived-finite-coresolutions-subsecn}.

 Let $\bE$ be an exact DG\+category.
 An \emph{exact DG\+pair}~\cite[Section~6.1]{Pedg} is a pair
$(\bE,\sK)$, where $\sK$ is a full additive subcategory in the exact
category $\sZ^0(\bE^\bec)$ such that $\sK$ is preserved by the shift
functors $[1]$ and~$[-1]$, the image of the functor
$\Phi_\bE\:\sZ^0(\bE)\rarrow\sZ^0(\bE^\bec)$ is contained in $\sK$,
and the full subcategory $\sK$ inherits an exact category structure
from $\sZ^0(\bE^\bec)$ (in the sense of the definition in
Section~\ref{derived-second-kind-subsecn}).

 Following~\cite[Section~6.1]{Pedg}, we denote the functor $\Phi_\bE$,
viewed as taking values in $\sK$, by $\Phi_\bE^\sK\:\sZ^0(\bE)\rarrow
\sK$, while the restrictions of the two functors $\Phi_\bE^{\pm}\:
\sZ^0(\bE^\bec)\rarrow\sZ^0(\bE)$ to the full subcategory
$\sK\subset\sZ^0(\bE^\bec)$ are denoted by $\Phi_\bE^{\sK,+}$
and $\Phi_\bE^{\sK,-}\:\sK\rarrow\sZ^0(\bE)$.

\begin{lem} \label{exact-dg-pair-co-generating-property}
 Let $(\bE,\sK)$ be an exact DG\+pair, and let $E\in\sZ^0(\sE^\bec)$
and $K\in\sK$ be two objects.
 In this context: \par
\textup{(a)} for any admissible monomorphism $K\rarrow E$ in
the exact category\/ $\sZ^0(\bE^\bec)$ there exists an object
$L\in\sK$ and an admissible monomorphism $E\rarrow L$ in\/
$\sZ^0(\bE^\bec)$ such that the composition $K\rarrow E\rarrow L$
is an admissible monomorphism in\/~$\sK$; \par
\textup{(b)} for any admissible epimorphism $E\rarrow K$ in
the exact category\/ $\sZ^0(\bE^\bec)$ there exists an object
$L\in\sK$ and an admissible epimorphism $L\rarrow E$ in\/
$\sZ^0(\bE^\bec)$ such that the composition $L\rarrow E\rarrow K$
is an admissible epimorphism in\/~$\sK$.
\end{lem}

\begin{proof}
 Let us prove part~(b).
 The object $L\in\sK$ and the admissible epimorphism $L\rarrow E$ will
only depend on the object $E\in\sZ^0(\bE^\bec)$ and not on~$K$.
 We simply put $L=\Phi_\bE^\sK(\Psi_\bE^-(E))$, and let $L\rarrow E$ be
the adjunction morphism for the pair of adjoint functors $\Phi_\bE$
and~$\Psi_\bE^-$.
 By~\cite[Lemmas~3.8 and~4.11]{Pedg}, we have an admissible short
exact sequence $0\rarrow E[-1]\rarrow\Phi_\bE^\sK(\Psi_\bE^-(E))
\rarrow E\rarrow0$ in $\sZ^0(\bE^\bec)$; so $L\rarrow E$ is
an admissible epimorphism in $\sZ^0(\bE^\bec)$.
 Furthermore, let $F\in\sZ^0(\bE^\bec)$ be the kernel of
the admissible epimorphism $E\rarrow K$.
 Then we have a commutative diagram
$$
 \xymatrix{
  F[-1] \ar@{>->}[r] \ar@{>->}[d]
  & E[-1] \ar@{->>}[r] \ar@{>->}[d] &  K[-1] \ar@{>->}[d] \\
  \Phi_\bE^\sK(\Psi_\bE^-(F)) \ar@{>->}[r] \ar@{->>}[d]
  & \Phi_\bE^\sK(\Psi_\bE^-(E)) \ar@{->>}[r] \ar@{->>}[d]
  & \Phi_\bE^\sK(\Psi_\bE^{\sK,-}(K)) \ar@{->>}[d] \\
  F \ar@{>->}[r] & E \ar@{->>}[r] & K
 }
$$
where all the rows and colums are admissible short exact sequences
in $\sZ^0(\bE^\bec)$.
 Moreover, all the terms of the middle row and the rightmost column
belong to~$\sK$; so the middle row and the rightmost column are
admissible short exact sequences in~$\sK$.
 Thus the morphisms $\Phi_\bE^\sK(\Psi_\bE^-(E))\rarrow
\Phi_\bE^\sK(\Psi_\bE^{\sK,-}(K))$ and
$\Phi_\bE^\sK(\Psi_\bE^{\sK,-}(K))\rarrow K$ are admissible
epimorphisms in $\sK$, and it follows that so is their composition
$L=\Phi_\bE^\sK(\Psi_\bE^-(E))\rarrow K$.
 By commutativity of the lower rightmost square of the diagram,
the latter composition is equal to the composition
$L\rarrow E\rarrow K$.
\end{proof}

\begin{cor} \label{exact-dg-pair-inj-projectivity-in-K-and-E-agree}
 Let $(\bE,\sK)$ be an exact DG\+pair.
 In this context: \par
\textup{(a)} An object of\/ $\sK$ is injective in\/ $\sK$ if and only if
it is injective in\/ $\sZ^0(\bE^\bec)$.
 Every injective object of\/ $\sZ^0(\bE^\bec)$ is a direct summand of
an injective object of\/~$\sK$. \par
\textup{(b)} An object of\/ $\sK$ is projective in\/ $\sK$ if and only
if it is projective in\/ $\sZ^0(\bE^\bec)$.
 Every projective object of\/ $\sZ^0(\bE^\bec)$ is a direct summand of
an projective object of\/~$\sK$.
\end{cor}

\begin{proof}
 Let us prove part~(b).
 Any object of $\sK$ that is projective in $\sZ^0(\bE^\bec)$ is
projective in $\sK$, because the inclusion $\sK\rarrow\sZ^0(\bE^\bec)$
is a fully faithful exact functor.
 Conversely, let $P\in\sK$ be a projective object.
 In order to prove that $P$ is projective in $\sZ^0(\bE^\bec)$, it
suffuces to check that every admissible epimorphism $E\rarrow P$ in
$\sZ^0(\bE^\bec)$ is split.
 Indeed, by Lemma~\ref{exact-dg-pair-co-generating-property}(b)
we have an object $L\in\sK$ and a morphism $L\rarrow E$ in
$\sZ^0(\bE^\bec)$ such that the composition $L\rarrow P$ is
an admissible epimorphism in~$\sK$.
 Now the epimorphism $L\rarrow P$ is split, and it follows that so is
the epimorphism $E\rarrow P$.

 Let $Q\in\sZ^0(\bE^\bec)$ be a projective object.
 Notice that the functors $\Phi_\bE^\sK$, \,$\Psi_\bE^{\sK,-}=
\Psi_\bE^{\sK,+}[1]$, and their shifts form an infinite ladder of
adjoint exact functors between the exact categories $\sZ^0(\bE)$
and~$\sK$.
 Hence all these functors take projective objects to projective
objects.
 The same applies to the functors $\Phi_\bE$ and $\Psi_\bE^\pm$,
and their shifts, acting between the exact categories $\sZ^0(\bE)$
and $\sZ^0(\bE^\bec)$.
 Thus the adjunction morphism $\Phi_\bE^\sK(\Psi_\bE^-(Q))\rarrow Q$
from the proof of Lemma~\ref{exact-dg-pair-co-generating-property}(b)
is an admissible epimorphism onto $Q$ from a projective object
of~$\sK$.
\end{proof}

 The next corollary is a generalization of~\cite[Lemmas~6.2
and~7.1]{PS5}.

\begin{cor} \label{exact-dg-pair-enough-injectives-projectives}
 Let $(\bE,\sK)$ be an exact DG\+pair.
 In this context: \par
\textup{(a)} The following conditions are equivalent:
\begin{enumerate}
\item there are enough injective objects in the exact category\/
$\sZ^0(\bE)$; \par
\item there are enough injective objects in the exact category\/
$\sZ^0(\bE^\bec)$; \par
\item there are enough injective objects in the exact category\/~$\sK$.
\end{enumerate} \par
\textup{(b)} The following conditions are equivalent:
\begin{enumerate}
\item there are enough projective objects in the exact category\/
$\sZ^0(\bE)$; \par
\item there are enough projective objects in the exact category\/
$\sZ^0(\bE^\bec)$; \par
\item there are enough projective objects in the exact category\/~$\sK$.
\end{enumerate}
\end{cor}

\begin{proof}
 Let us prove part~(b).

 (3)~$\Longrightarrow$~(1)
 Let $B\in\bE$ be an object.
 According to~(3), there exists a projective object $P\in\sK$ together
with an admissible epimorphism $P\rarrow\Phi_\bE^\sK(B)$ in~$\sK$.
 By~\cite[Lemma~5.7(b)]{Pedg}, it follows that the adjoint morphism
$\Phi_\bE^{\sK,+}(P)\rarrow B$ is an admissible epimorphism in
$\sZ^0(\bE)$.
 According to the second paragraph of the proof of
Corollary~\ref{exact-dg-pair-inj-projectivity-in-K-and-E-agree},
the object $\Phi_\bE^{\sK,+}(P)$ is projective in $\sZ^0(\bE)$.

 (1)~$\Longrightarrow$~(3)
 Let $K\in\sK$ be an object.
 According to~(1), there exists a projective object $P\in\sZ^0(\bE)$
together with an admissible epimorphism $P\rarrow\Psi_\bE^{\sK,-}(K)$
in $\sZ^0(\bE)$.
 The adjoint morphism $\Phi_\bE^\sK(P)\rarrow K$ is the composition
$\Phi_\bE^\sK(P)\rarrow\Phi_\bE^\sK\Psi_\bE^{\sK,-}(K)\rarrow K$.
 The functor $\Phi_\bE^\sK$ is exact, so it takes admissible
epimorphisms to admissible epimorphisms; hence
$\Phi_\bE^\sK(P)\rarrow\Phi_\bE^\sK\Psi_\bE^{\sK,-}(K)$ is
an admissible epimorphism in~$\sK$.
 Following the proof of
Lemma~\ref{exact-dg-pair-co-generating-property}(b), the morphism
$\Phi_\bE^\sK\Psi_\bE^{\sK,-}(K)\rarrow K$ is also an admissible
epimorphism in~$\sK$.
 Thus the composition $\Phi_\bE^\sK(P)\rarrow K$ is an admissible
epimorphism in~$\sK$.
 It remains to point out that the object $\Phi_\bE^\sK(P)$ is
projective in $\sK$ according to the second paragraph of the proof of
Corollary~\ref{exact-dg-pair-inj-projectivity-in-K-and-E-agree}.

 (2)~$\Longrightarrow$~(1) is a special case of (3)~$\Rightarrow$~(1).

 (1)~$\Longrightarrow$~(2) is a special case of (1)~$\Rightarrow$~(3).
 
 (3)~$\Longrightarrow$~(2)
 Let $E\in\sZ^0(\bE^\bec)$ be an object.
 By Lemma~\ref{exact-dg-pair-co-generating-property}(b), there exists
an object $L\in\sK$ together with an admissible epimorphism
$L\rarrow E$ in $\sZ^0(\bE^\bec)$.
 According to~(3), there exists a projective object $P\in\sK$
together with an admissible epimorphism $P\rarrow L$ in~$\sK$.
 Then the composition $P\rarrow L\rarrow E$ is an admissible epimorphism
in $\sZ^0(\bE^\bec)$.
 It remains to point out that the object $P$ is also projective in
$\sZ^0(\bE^\bec)$ by
Corollary~\ref{exact-dg-pair-inj-projectivity-in-K-and-E-agree}(b).

 (2)~$\Longrightarrow$~(3)
 Let $K\in\sK$ be an object.
 According to~(2), there exists a projective object
$Q\in\sZ^0(\bE^\bec)$ together with an admissible epimorphism
$Q\rarrow K$ in $\sZ^0(\bE^\bec)$.
 According to (the proof of)
Lemma~\ref{exact-dg-pair-co-generating-property}(b), there is
an admissible epimorphism $\Phi_\bE^\sK\Psi_\bE^-(Q)\rarrow K$ in~$\sK$.
 It remains to point out that the object $\Phi_\bE^\sK\Psi_\bE^-(Q)$
is projective in $\sK$ according to the second paragraph of the proof of
Corollary~\ref{exact-dg-pair-inj-projectivity-in-K-and-E-agree}.
\end{proof}

 It follows from
Corollary~\ref{exact-dg-pair-inj-projectivity-in-K-and-E-agree}
that an object $J\in\bE$ is graded-injective (in the sense of
the definition in Section~\ref{derived-second-kind-subsecn}) if and
only if the object $\Phi_\bE^\sK(J)$ is injective in~$\sK$.
 Similarly, an object $P\in\bE$ is graded-projective if and only if
the object $\Phi_\bE^\sK(P)$ is projective in~$\sK$.
 So the graded-injectivity and graded-projectivity with respect to
an exact DG\+pair $(\bE,\sK)$ agree with the same notions for
the underlying exact DG\+category $\bE$ of the exact DG\+pair
$(\bE,\sK)$, and do not depend on~$\sK$.

 For any exact category $\sE$ we denote by $\sE_\proj$ the class
of all projective objects in $\sE$ and by $\sE^\inj$ the class of
all injective objects in~$\sE$.
 An exact category $\sK$ is said to have \emph{homological
dimension\/~$\le d$} (where $d\ge-1$ is an integer) if
$\Ext_\sK^{d+1}(X,Y)=0$ for all $X$, $Y\in\sK$.
 The following proposition is our version
of~\cite[Theorem~5.6]{Pedg} (and a partial generalization
of~\cite[Proposition~7.5]{Pphil}).

\begin{prop} \label{exact-dg-pair-finite-homol-dim-prop}
 Let\/ $(\bE,\sK)$ be an exact DG\+pair such that the exact category\/
$\sK$ has finite homological dimension.  Then \par
\textup{(a)} Assume that the exact category\/ $\sK$ has enough
injective objects.
 Then the classes of absolutely acyclic and Becker-coacyclic objects
coincide in\/ $\bE$, that is\/ $\Ac^\bco(\bE)=\Ac^\abs(\bE)$.
 Furthermore, the composition of the triangulated inclusion functor\/
$\sH^0(\bE^\binj)\rarrow\sH^0(\bE)$ and the Verdier quotient functor\/
$\sH^0(\bE)\rarrow\sD^\bco(\bE)$ is a triangulated equivalence.
 So one has
$$
 \sH^0(\bE^\binj)\simeq\sD^\abs(\bE)=\sD^\bco(\bE).
$$ \par
\textup{(b)} Assume that the exact category\/ $\sK$ has enough
projective objects.
 Then the classes of absolutely acyclic and Becker-contraacyclic objects
coincide in\/ $\bE$, that is\/ $\Ac^\bctr(\bE)=\Ac^\abs(\bE)$.
 Furthermore, the composition of the triangulated inclusion functor\/
$\sH^0(\bE_\bproj)\rarrow\sH^0(\bE)$ and the Verdier quotient functor\/
$\sH^0(\bE)\rarrow\sD^\bctr(\bE)$ is a triangulated equivalence.
 So one has
$$
 \sH^0(\bE_\bproj)\simeq\sD^\abs(\bE)=\sD^\bctr(\bE).
$$
\end{prop}

\begin{proof}
 Part~(b): the proof of~\cite[Theorem~5.6(b)]{Pedg} is applicable.
 The only changes one needs to make are that $\Phi_\bE^\sK(P_i)\in
\sK\cap\sZ^0(\bE^\bec)_\proj=\sK_\proj$ by
Corollary~\ref{exact-dg-pair-inj-projectivity-in-K-and-E-agree}(b),
the complex $0\rarrow\Phi_\bE^\sK(Q)\rarrow\Phi_\bE^\sK(P_{n-1})
\rarrow\dotsb\rarrow\Phi_\bE^\sK(P_0)\rarrow\Phi_\bE^\sK(B)\rarrow0$
is exact in $\sK$, and it follows that $\Phi_\bE(Q)\in\sK_\proj$.
\end{proof}

 Given an exact DG\+pair $(\bE,\sK)$, an \emph{exact DG\+subpair}
$(\bF,\sL)\subset(\bE,\sK)$ consists of a full DG\+subcategory
$\bF\subset\bE$ inheriting an exact DG\+category structure from $\bE$
and a full additive subcategory $\sL\subset\sK\cap\sZ^0(\bF^\bec)
\subset\sZ^0(\bE^\bec)$ such that $(\bF,\sL)$ is an exact
DG\+pair~\cite[Section~6.1]{Pedg}.
 An exact DG\+subpair $(\bF,\sL)\subset(\bE,\sK)$ is said to be
\emph{strict} if $\sL$ is a strictly full subcategory in $\sK$ and,
for any given object $F\in\bE$, one has $F\in\bF$ whenever
$\Phi_\bE^\sK(F)\in\sL$.
 The definition of a \emph{strict exact DG\+subcategory} was given
in Section~\ref{derived-second-kind-subsecn}.

\begin{lem} \label{strict-exact-dg-subpair-constructed}
 Let $(\bE,\sK)$ be an exact DG\+pair and\/ $\sL\subset\sK$ be a full
additive subcategory preserved by the shift functors and closed under
extensions.
 Denote by\/ $\bF\subset\bE$ the full DG\+subcategory consisting of
all objects $F\in\bE$ such that\/ $\Phi_\bE^\sK(F)\in\sL$.
 Then $(\bF,\sL)$ is a strict exact DG\+subpair in $(\bE,\sK)$.
\end{lem}

\begin{proof}
 This is~\cite[Examples~6.1(1\+-2)]{Pedg}.
\end{proof}

\subsection{Becker contraderived category for exact DG-subcategory
with finite coresolution dimension}
\label{becker-contraderived-finite-coresolutions-subsecn}
 The aim of this section, which complements the exposition in
Section~\ref{derived-second-kind-subsecn}, is to prove the following
extension of the result~\cite[Proposition~B.7.10]{Pcosh} from
the context of exact categories to that of exact DG\+categories.

\begin{prop} \label{becker-contraderived-finite-coresolutions}
 Let\/ $\bE$ be an exact DG\+category and\/ $\bC\subset\bE$ be
a strict exact DG\+subcategory.
 Assume that the exact category\/ $\sE=\sZ^0(\bE^\bec)$ is weakly
idempotent-complete and has enough projective objects, the full
subcategory\/ $\sC=\sZ^0(\bC^\bec)\subset\sZ^0(\bE^\bec)$ is
the right-hand class of a hereditary complete cotorsion pair $(\sF,\sC)$
in the exact category\/ $\sE$, and the\/ $\sC$\+coresolution dimensions
of all the objects of\/ $\sE$ do not exceed a fixed constant~$n$.
 Then an object of\/ $\bC$ is Becker-contraacyclic in\/ $\bC$ if and
only if it is Becker-contraacyclic as an object of\/ $\bE$,
$$
 \Ac^\bctr(\bC)=\bC\cap\Ac^\bctr(\bE).
$$
 The inclusion of exact DG\+categories\/ $\bC\rarrow\bE$ induces
an equivalence of the Becker contraderived categories
$$
 \sD^\bctr(\bC)\simeq\sD^\bctr(\bE).
$$
\end{prop}

\begin{proof}
 By~\cite[Lemma~B.1.9]{Pcosh} or~\cite[Lemma~6.4]{Pfltp},
the homological dimension of the exact category~$\sF$ (with the exact
structure inherited from~$\sE$) is finite and does not exceed~$n$.
 Furthermore, by~\cite[Lemmas~1.9 and~7.1(a)]{Pfltp}, there are
enough projective and injective objects in the exact category $\sF$,
and also enough projective objects in the exact category~$\sC$;
the classes of such objects are $\sF_\proj=\sE_\proj$ and
$\sF^\inj=\sF\cap\sC=\sC_\proj$.

 The full subcategory $\sF\subset\sZ^0(\bE^\bec)$ is preserved by
the shift functors (because the full subcategory $\sC=\sZ^0(\bC^\bec)
\subset\sZ^0(\bE^\bec)$ is), and closed under extensions in
$\sZ^0(\bE^\bec)$.
 So Lemma~\ref{strict-exact-dg-subpair-constructed} is applicable,
and we have a strict exact DG\+subpair $(\bF,\sF)$ in
$(\bE,\sZ^0(\bE^\bec))$.
 By Proposition~\ref{exact-dg-pair-finite-homol-dim-prop}, we have
$\sH^0(\bF^\binj)\simeq\sD^\abs(\bF)\simeq\sH^0(\bF_\bproj)$.
 Specifically, following the proof of the proposition, this means that,
for every object $F\in\bF$, there exists a graded-projective object
$P\in\bF_\bproj$ together with a closed morphism $P\rarrow F$ with
an absolutely acyclic cone, and there also exists a graded-injective
object $Q\in\bF^\binj$ together with a closed morphism $F\rarrow Q$
with an absolutely acyclic cone.
 Following the previous paragraph, we have $\bF_\bproj=\bE_\bproj$
and $\bF^\binj=\bC_\bproj$.

 Let $A\in\Ac^\abs(\bF)$ be an absolutely acyclic object in~$\bF$.
 Then, for every object $C\in\bC$, the complex of abelian groups
$\Hom_\bE^\bu(A,C)$ is acyclic.
 Indeed, for any short exact sequence $0\rarrow F\rarrow G\rarrow H
\rarrow0$ in $\sZ^0(\bF)$, the short sequence of complexes of abelian
groups $0\rarrow\Hom_\bE^\bu(H,C)\rarrow\Hom_\bE^\bu(G,C)\rarrow
\Hom_\bE^\bu(F,C)\rarrow0$ is exact in view of~\cite[Lemma~3.9]{Pedg}.
 The point is that $0\rarrow\Phi_\bE(F)\rarrow\Phi_\bE(G)\rarrow
\Phi_\bE(H)\rarrow0$ is an admissible short exact sequence in $\sF$,
and the functor $\Hom_{\sZ^0(\bE^\bec)}({-},\Phi_\bE(C))$ preserves
exactness of any such short exact sequence
(cf.~\cite[proof of Theorem~5.5]{Pedg}).

 Now for every object $Q\in\bC_\bproj$ we have an object
$P\in\bE_\bproj$ together with a closed morphism $P\rarrow Q$ in $\bE$
with the cone absolutely acyclic in $\bF$, and conversely, for every
object $P\in\bE_\bproj$ we have an object $Q\in\bC_\bproj$ together
with a closed morphism $P\rarrow Q$ in $\bE$ with the cone absolutely
acyclic in~$\bF$.
 For any object $C\in\bC$, the cone of the morphism of complexes of
abelian groups $\Hom_\bE^\bu(Q,C)\rarrow\Hom_\bE^\bu(P,C)$ is acyclic;
so this morphism of complexes of abelian groups is a quasi-isomorphism.
 Thus the complex $\Hom_\bE^\bu(Q,C)$ is acyclic if and only if
the complex $\Hom_\bE^\bu(P,C)$ is acyclic.
 We have proved that the object $C$ is Becker-contraacyclic in $\bC$ if
and only if it is Becker-contraacyclic in~$\bE$.

 Finally, for any object $E\in\bE$ there exists an object $C\in\bC$
together with a closed morphism $E\rarrow C$ whose cone in absolutely
acyclic in~$\bE$.
 This holds because all objects of $\sZ^0(\bE^\bec)$ have finite
coresolution dimensions with respect to the coresoloving subcategory
$\sZ^0(\bC^\bec)$ (use the dual version of the construction
in~\cite[beginning of the proof of Theorem~6.6]{Pedg}).
 As any absolutely acyclic object is Becker-coacyclic in $\bE$,
the triangulated equivalence $\sD^\bctr(\bC)\simeq\sD^\bctr(\bE)$
follows.
\end{proof}

\subsection{Coderived and contraderived categories of $X$-cotorsion
and $X$-locally cotorsion CDG-modules}
\label{co-contra-derived-loc-cotorsion-subsecn}
 This section complements
Sections~\ref{coderived-qcoh-cdg-modules-subsecn}\+-%
\ref{contraderived-lcth-cdg-modules-subsecn}.

\begin{cor} \label{qcoh-cot-cdg-absolute-derived-equivalence}
 Let $X$ be a semi-separated Noetherian scheme of finite Krull dimension
and $\cB^\cu$ be a quasi-coherent CDG\+quasi-algebra over~$X$.
 Then the inclusions of exact/abelian DG\+categories
$\cB^\cu\bQcoh^{X\dcot}\rarrow\cB^\cu\bQcoh^{X\dcta}\rarrow
\cB^\cu\bQcoh$ induce equivalences of the absolute derived categories
$$
 \sD^\abs(\cB^\cu\bQcoh^{X\dcot})\simeq
 \sD^\abs(\cB^\cu\bQcoh^{X\dcta})\simeq\sD^\abs(\cB^\cu\bQcoh).
$$
\end{cor}

\begin{proof}
 The equivalence $\sD^\abs(\cB^\cu\bQcoh^{X\dcta})\simeq
\sD^\abs(\cB^\cu\bQcoh)$ holds in greater generality by
Corollary~\ref{qcoh-cta-cdg-absolute-derived-equivalence}.
 The discussion in the proof of that corollary is relevant.
 The equivalence $\sD^\abs(\cB^\cu\bQcoh^{X\dcot})\simeq
\sD^\abs(\cB^\cu\bQcoh)$ follows from the graded version of
Lemma~\ref{A-qcoh-X-cot-coresolution-dimension} in view of the dual
version of Proposition~\ref{second-kind-finite-resolutions}(a).
\end{proof}

 The following lemma is our version of the main results
of~\cite[Section~4.5]{Pcosh}.

\begin{lem} \label{enough-projectives-in-lcth-A-modules}
 Let $X$ be a quasi-compact semi-separated scheme with an open
covering\/ $\bW$ and $\cA$ be a quasi-coherent quasi-algebra
over~$X$.
 Then \par
\textup{(a)} the exact category $\cA\Lcth_\bW$ has enough projective
objects, all such objects are (globally) contraherent on $X$, and
the class of all projective objects in $\cA\Lcth_\bW$ does not depend
on an open covering\/~$\bW$; \par
\textup{(b)} the exact category $\cA\Lcth_\bW^{X\dlct}$ has enough
projective objects, all such objects are (globally) contraherent on $X$,
and the class of all projective objects in $\cA\Lcth_\bW^{X\dlct}$
does not depend on an open covering\/~$\bW$; \par
\textup{(c)} the exact category $\cA\Lcth_\bW^{\cA\dlct}$ has enough
projective objects, all such objects are (globally) contraherent on $X$,
and the class of all projective objects in $\cA\Lcth_\bW^{\cA\dlct}$
does not depend on an open covering\/~$\bW$.
\end{lem}

\begin{proof}
 We will spell out a detailed proof of the first assertion of part~(a).
 In the case of an affine scheme $U$ with the open covering $\{U\}$,
the claim is that there are enough projective objects in
the exact category of $R$\+contraadjusted $A$\+modules
$A\Modl^{R\dcta}$ for $R=\cO(U)$ and $A=\cA(U)$.
 Indeed, by Remark~\ref{very-flaprojective-cotorsion-pair-generated-by}
with Lemma~\ref{very-flaprojective-cotorsion-pair-hereditary} and
Theorem~\ref{very-flaprojective-cotorsion-pair-complete},
the pair of classes of $A/R$\+very flaprojective $A$\+modules
$A\Modl_{A/R\dvflp}$ and $R$\+contraadjusted $A$\+modules
$A\Modl^{R\dcta}$ is a hereditary complete cotorsion pair
($A\Modl_{A/R\dvflp}$, $A\Modl^{R\dcta}$) in $A\Modl$.
 By~\cite[Lemma~1.9(b)]{Pfltp}, it follows that there are enough
projective objects in $A\Modl^{R\dcta}$, and the class of such
projective objects is precisely the intersection of two classes
$A\Modl_{A/R\dvflp}^{R\dcta}=A\Modl_{A/R\dvflp}\cap A\Modl^{R\dcta}$.

 Now let $\gM$ be a locally contraherent $\cA$\+module on $X$, and
let $X=\bigcup_\alpha U_\alpha$ be a finite affine open covering
of $X$ subordinate to~$\bW$.
 Denote by $j_\alpha\:U_\alpha\rarrow X$ the open immersion morphisms.
 Arguing as in the proof of
Lemma~\ref{antilocally-flaprojective-modules-generating-class} and
using Theorem~\ref{very-flaprojective-cotorsion-pair-complete}(a),
we produce an admissible epimorphism
$\bigoplus_\alpha j_\alpha{}_!\P_\alpha\rarrow\gM$ in $\cA\Lcth_\bW$,
where $\P_\alpha$ is a projective object of the exact category
$\cA|_{U_\alpha}\Ctrh$ of contraherent $\cA|_{U_\alpha}$\+modules
on~$U_\alpha$.
 It remains to observe that, for any scheme $X$ with an open covering
$\bW$, any affine open immersion morphism $j\:Y\rarrow X$ and any
quasi-coherent quasi-algebra $\cA$ on $X$, the direct image
functor~$j_!$ takes projective objects of $\cA|_Y\Lcth_{\bW|_Y}$ to
projective objects of $\cA_Y\Lcth_\bW$ in view of
the adjunction~\eqref{affine-open-immers-contrah-module-adjunction}.

 It follows that, for any choice of a finite affine open covering
$X=\bigcup_\alpha U_\alpha$ subordinate to $\bW$, the projective
objects of $\cA\Lcth_\bW$ are precisely the direct summands of
the finite direct sums $\bigoplus_\alpha j_\alpha{}_!\P_\alpha$,
where $\P_\alpha$ are projective objects of $\cA|_{U_\alpha}\Ctrh$.
 This implies all the assertions of part~(a).
 To prove part~(b), one needs to use
Theorem~\ref{flaprojective-cotorsion-pair-complete}, and for part~(c),
Theorem~\ref{flat-cotorsion-pair-complete}.
\end{proof}

 The following theorem should be compared with~\cite[Theorem~5.5.10
and Corollary~6.4.4(b)]{Pcosh} (see also
Remark~\ref{contraderived-of-lct=of-lcta-remark} above).

\begin{thm} \label{lcth-lcta-lct-cdg-abs-contraderived-equiv-thm}
 Let $X$ be a semi-separated Noetherian scheme of finite Krull dimension
with an open covering\/ $\bW$ and $\cB^\cu$ be a quasi-coherent
CDG\+quasi-algebra over~$X$.
 Then, for any absolute or contraderived category symbol\/ $\st=\abs$,
$\ctr$, or\/~$\bctr$, the inclusion of exact DG\+categories
$\cB^\cu\bLcth_\bW^{X\dlct}\rarrow\cB^\cu\bLcth_\bW$ induces
an equivalence of triangulated categories
$$
 \sD^\st(\cB^\cu\bLcth_\bW^{X\dlct})\simeq\sD^\st(\cB^\cu\bLcth_\bW).
$$
 There is a commutative diagram of triangulated equivalences induced
by the inclusions of exact DG\+categories
\begin{equation} \label{lcth-al-lcta-lct-cdg-abs-contrader-diagram}
\begin{gathered}
 \xymatrix{
  \sD^\st(\cB^\cu\bCtrh_\al) \ar@<-2pt>[r] \ar@{-}@<-2pt>[d]
  & \sD^\st(\cB^\cu\bCtrh)
  \ar@<-2pt>[r] \ar@{-}@<-2pt>[d] \ar@{-}@<-2pt>[l]
  & \sD^\st(\cB^\cu\bLcth_\bW) \ar@{-}@<-2pt>[l] \ar@{-}@<-2pt>[d] \\
  \sD^\st(\cB^\cu\bCtrh^{X\dlct}_\al) \ar@<-2pt>[r] \ar@<-2pt>[u]
  & \sD^\st(\cB^\cu\bCtrh^{X\dlct})
  \ar@<-2pt>[r] \ar@<-2pt>[u] \ar@{-}@<-2pt>[l]
  & \sD^\st(\cB^\cu\bLcth_\bW^{X\dlct}) \ar@{-}@<-2pt>[l] \ar@<-2pt>[u]
 }
\end{gathered}
\end{equation}
\end{thm}

\begin{proof}
 By Corollary~\ref{exact-dg-categories-of-lcth-cdg-modules},
we have $\sZ^0((\cB^\cu\bLcth_\bW)^\bec)\simeq\cB^*\Lcth_\bW$ and
similarly $\sZ^0((\cB^\cu\bLcth_\bW^{X\dlct})^\bec)\simeq
\cB^*\Lcth_\bW^{X\dlct}$, so $\cB^\cu\bLcth_\bW^{X\dlct}$ is a strict
exact DG\+subcategory in $\cB^\cu\bLcth_\bW$.
 By the graded version of
Corollary~\ref{A-lcth-al-X-lct-coresolution-dimension}(a),
the coresolution dimensions of all the objects of $\cB^*\Lcth_\bW$
with respect to the coresolving subcategory $\cB^*\Lcth_\bW^{X\dlct}$
do not exceed the Krull dimension $D$ of the scheme~$X$.
 By the dual version of
Proposition~\ref{second-kind-finite-resolutions}(a\+-b),
the triangulated equivalence $\sD^\st(\cB^\cu\bLcth_\bW^{X\dlct})\simeq
\sD^\st(\cB^\cu\bLcth_\bW)$ for $\st=\abs$ or~$\ctr$ follows.
{\emergencystretch=1em\par}

 The case of the Becker contraderived categories, $\st=\bctr$,
is a bit harder to prove.
 By the graded version of
Corollary~\ref{antilocally-flaprojective-cotorsion-pair-cor},
the full subcategory $\cB^*\Lcth_\bW^{X\dlct}$ is the right-hand
part of a hereditary cotorsion par in $\cB^*\Lcth_\bW$.
 By the graded version of
Lemma~\ref{enough-projectives-in-lcth-A-modules}(a), the exact
category $\cB^*\Lcth_\bW$ has enough projective objects.
 So the assumptions of
Proposition~\ref{becker-contraderived-finite-coresolutions} are
satisfied, and we can conclude that
$\sD^\bctr(\cB^\cu\bLcth_\bW^{X\dlct})\simeq
\sD^\bctr(\cB^\cu\bLcth_\bW)$.

 This proves the rightmost and the middle vertical triangulated
equivalences in
the diagram~\eqref{lcth-al-lcta-lct-cdg-abs-contrader-diagram},
while the horizontal equivalences are provided by
Corollary~\ref{contraderived-indep-of-covering-or-antilocal}(a\+-b).
 It remains to explain why the leftmost vertical triangulated functor
in~\eqref{lcth-al-lcta-lct-cdg-abs-contrader-diagram} is well-defined
(then it would follow from the commutativity of the diagram that it
is a triangulated equivalences).
 In the cases of the absolute derived or Positselski contraderived
categories, this holds because the DG\+functor
$\cB^\cu\bCtrh^{X\dlct}_\al\rarrow\cB^\cu\bCtrh_\al$
is exact and preserves infinite products.

 In the case of the Becker contraderived categories, the point is
that the classes of projective objects in $\cB^*\Ctrh_\al$ and in
$\cB^*\Lcth_\bW$ coincide, and similarly, the classes of projective
objects in $\cB^*\Ctrh^{X\dlct}_\al$ and in $\cB^*\Lcth_\bW^{X\dlct}$
coincide (e.~g., by the proofs of
Corollary~\ref{contraderived-indep-of-covering-or-antilocal}(a\+-b)
and Proposition~\ref{second-kind-finite-resolutions}(d)).
 So an object of $\cB^\cu\bCtrh_\al$ is Becker-contraacyclic if and
only if it is Becker-contraacyclic in $\cB^\cu\bLcth_\bW$, and
similarly, an object of $\cB^\cu\bCtrh^{X\dlct}_\al$ is
Becker-contraacyclic if and only if it is Becker-contraacyclic
in $\cB^\cu\bLcth_\bW^{X\dlct}$.
 This can be also seen from the assertions of
Corollary~\ref{contraderived-indep-of-covering-or-antilocal}(a\+-b).
 It remains to recall that an object of $\cB^\cu\bLcth_\bW^{X\dlct}$
is Becker-contraacyclic if and only if it is Becker-contraacyclic
in $\cB^\cu\bLcth_\bW$, as per
Proposition~\ref{becker-contraderived-finite-coresolutions}.
 Thus an object of $\cB^\cu\bCtrh^{X\dlct}_\al$ is Becker-contraacyclic
if and only if it is Becker-contraacyclic in $\cB^\cu\bCtrh_\al$.
\end{proof}

\subsection{Semicontraderived categories of $X$-locally cotorsion
contraherent modules} \label{semicontraderived-loc-cotorsion-subsecn}
 The following corollary complements
Sections~\ref{semiderived-contraherent-subsecn}
and~\ref{becker-semicontraderived-defined-subsecn}.

\begin{cor} \label{semicontraderived-X-lcta-X-lct-A-mod-equivalence}
 Let $X$ be a semi-separated Noetherian scheme of finite Krull dimension
with an open covering\/ $\bW$ and $\cA$ be a quasi-coherent
quasi-algebra over~$X$.
 Then, for any semi(contra)derived category symbol\/ $\st=\si$
or\/~$\bsi$, the inclusion of exact categories $\cA\Lcth_\bW^{X\dlct}
\rarrow\cA\Lcth_\bW$ induces an equivalence of triangulated categories
$$
 \sD^\st(\cA\Lcth_\bW^{X\dlct})\simeq\sD^\st(\cA\Lcth_\bW).
$$
\end{cor}

\begin{proof}
 By~\cite[Theorem~5.5.10 or Corollary~6.4.4(b)]{Pcosh}, or by
Theorem~\ref{lcth-lcta-lct-cdg-abs-contraderived-equiv-thm} above,
a complex in $X\Lcth_\bW^\lct$ is Becker (respectively, Positselski)
contraacyclic if and only if it is Becker (resp., Positselski)
contraacyclic as a complex in $X\Lcth_\bW$.
 Consequently, a complex in $\cA\Lcth_\bW^{X\dlct}$ is Becker
(resp., Positselski) semicontraacyclic if and only if its is Becker
(resp., Positselski) semicontraacyclic as a complex in $\cA\Lcth_\bW$.
 Hence the triangulated functor
$\sD^\st(\cA\Lcth_\bW^{X\dlct})\simeq\sD^\st(\cA\Lcth_\bW)$ induced
by the inclusion of exact categories $\cA\Lcth_\bW^{X\dlct}
\rarrow\cA\Lcth_\bW$ is well-defined.

 It remains to mention that, for any complex of $\bW$\+locally
contraherent $\cA$\+modules $\P^\bu$ on $X$ there exists a complex of
$X$\+locally cotorsion $\bW$\+locally contraherent $\cA$\+modules
$Q^\bu$ together with a morphism of complexes $\P^\bu\rarrow\Q^\bu$
whose cone is absolutely acyclic in $\cA\Lcth_\bW$.
 This follows from
Corollary~\ref{A-lcth-al-X-lct-coresolution-dimension}(a)
(use the dual version of the construction in the beginning of the proof
of~\cite[Proposition~A.5.8]{Pcosh} or the last paragraph of the proof
of Proposition~\ref{becker-contraderived-finite-coresolutions} above).
 As all absolutely acyclic complexes in $\cA\Lcth_\bW$ are obviously
semicontraacyclic, the desired triangulated equivalence follows.
\end{proof}

 In fact, there is no difference between the Positselski and
Becker semicontraderived categories under the assumptions of
Corollary~\ref{semicontraderived-X-lcta-X-lct-A-mod-equivalence};
see Corollary~\ref{noetherian-Positselski=Becker-semicontraderived}.

\subsection{Coderived and contraderived categories of thick
$X$-cotorsion and $X$-locally cotorsion CDG-modules}
\label{co-contra-derived-thick-loc-cotorsion-subsecn}
 This section complements
Sections~\ref{coderived-of-thick-cdg-modules-subsecn}\+-%
\ref{contraderived-of-thick-cdg-modules-subsecn}.

\begin{cor} \label{thick-cdg-qcoh-cot-abs-derived-equiv-cor}
 Let $X$ be a semi-separated Noetherian scheme of finite Krull dimension
and $(\g,\widetilde\g)$ be a quasi-coherent twisted Lie algebroid
over~$X$.
 Assume that\/ $\g$~is a finite locally free sheaf on $X$, and let
$\cB^\cu=\cC^\cu_X(\g,\widetilde\g)$ be the related Chevalley--Eilenberg
quasi-coherent CDG\+quasi-algebra over~$X$.
 Then the inclusions of exact DG\+categories
$\cB^\cu\bQcoh^{X\dcot}_\bth\rarrow\cB^\cu\bQcoh^{X\dcta}_\bth\rarrow
\cB^\cu\bQcoh_\bth$ induce equivalences of the absolute derived
categories
$$
 \sD^\abs(\cB^\cu\bQcoh^{X\dcot}_\bth)\simeq
 \sD^\abs(\cB^\cu\bQcoh^{X\dcta}_\bth)\simeq
 \sD^\abs(\cB^\cu\bQcoh_\bth).
$$
\end{cor}

\begin{proof}
 The equivalence $\sD^\abs(\cB^\cu\bQcoh^{X\dcta}_\bth)\simeq
\sD^\abs(\cB^\cu\bQcoh_\bth)$ holds in greater generality by
Corollary~\ref{thick-cdg-qcoh-cta-abs-derived-equiv-cor}.
 The proof of the equivalence $\sD^\abs(\cB^\cu\bQcoh^{X\dcot}_\bth)
\simeq\sD^\abs(\cB^\cu\bQcoh_\bth)$ is similar and based on the graded
version of Lemma~\ref{A-qcoh-X-cot-coresolution-dimension} together
with the dual version of
Proposition~\ref{second-kind-finite-resolutions}(a).
 The upper line of
diagram~\eqref{Upsilon-qcoh-X-cot-thick-diagram}
in Corollary~\ref{exact-dg-categories-of-qcoh-thick-cdg-modules}
is also relevant.
\end{proof}

\begin{cor} \label{thick-lcth-al-abs-ctr-derived-lcta-lct-equiv-cor}
 Let $X$ be a semi-separated Noetherian scheme of finite Krull
dimension with an open covering\/ $\bW$ and $(\g,\widetilde\g)$ be
a quasi-coherent twisted Lie algebroid over~$X$.
 Assume that\/ $\g$~is a finite locally free sheaf on $X$, and let
$\cB^\cu=\cC^\cu_X(\g,\widetilde\g)$ be the related Chevalley--Eilenberg
quasi-coherent CDG\+quasi-algebra over~$X$.
 Let\/ $\st=\abs$, $\ctr$, or\/~$\bctr$ be an absolute derived or
contraderived category symbol.
 Then the inclusion of exact DG\+categories
$\cB^\cu\bLcth_\bW^{X\dlct,\bth}\rarrow\cB^\cu\bLcth_\bW^\bth$
induces a triangulated equivalence
$$
 \sD^\st(\cB^\cu\bLcth_\bW^{X\dlct,\bth})\simeq
 \sD^\st(\cB^\cu\bLcth_\bW^\bth).
$$
 There is a commutative diagram of triangulated equivalences
induced by the inclusions of exact DG\+categories
\begin{equation} \label{thick-lcth-al-lcta-lct-cdg-abs-contrader-diag}
\begin{gathered}
 \xymatrix{
  \sD^\st(\cB^\cu\bCtrh_\al^\bth) \ar@<-2pt>[r] \ar@{-}@<-2pt>[d]
  & \sD^\st(\cB^\cu\bCtrh^\bth)
  \ar@<-2pt>[r] \ar@{-}@<-2pt>[d] \ar@{-}@<-2pt>[l]
  & \sD^\st(\cB^\cu\bLcth_\bW^\bth)
  \ar@{-}@<-2pt>[l] \ar@{-}@<-2pt>[d] \\
  \sD^\st(\cB^\cu\bCtrh^{X\dlct,\bth}_\al) \ar@<-2pt>[r] \ar@<-2pt>[u]
  & \sD^\st(\cB^\cu\bCtrh^{X\dlct,\bth})
  \ar@<-2pt>[r] \ar@<-2pt>[u] \ar@{-}@<-2pt>[l]
  & \sD^\st(\cB^\cu\bLcth_\bW^{X\dlct,\bth})
  \ar@{-}@<-2pt>[l] \ar@<-2pt>[u]
 }
\end{gathered}
\end{equation}
\end{cor}

\begin{proof}
 For the contraderived category symbols $\st=\ctr$ or~$\bctr$, we have
a commutative square diagram of triangulated functors induced by
the inclusions of exact DG\+categories
\begin{equation} \label{thick-all-lcth-lcta-lct-cdg-abs-contrader-diag}
\begin{gathered}
 \xymatrix{
  \sD^\st(\cB^\cu\bLcth_\bW^\bth) \ar@<-2pt>[r] \ar@{-}@<-2pt>[d]
  & \sD^\st(\cB^\cu\bLcth_\bW)
  \ar@{-}@<-2pt>[l] \ar@{-}@<-2pt>[d] \\
  \sD^\st(\cB^\cu\bLcth_\bW^{X\dlct,\bth}) \ar@<-2pt>[r] \ar@<-2pt>[u]
  & \sD^\st(\cB^\cu\bLcth_\bW^{X\dlct})
  \ar@{-}@<-2pt>[l] \ar@<-2pt>[u]
 }
\end{gathered}
\end{equation}
 The horizontal functors are triangulated equivalences according to
the upper lines of the two commutative
diagrams~(\ref{thick-lcta-cdg-contrader-equivs-diagram}\+-%
\ref{thick-X-lct-cdg-contrader-equivs-diagram}) in
Corollary~\ref{thick-cdg-contraderived-equiv-cor}.
 The rightmost vertical line is a well-defined triangulated equivalence
by Theorem~\ref{lcth-lcta-lct-cdg-abs-contraderived-equiv-thm}.
 The argument proving that the leftmost vertical functor is well-defined
is similar to the respective argument in the proof of
Theorem~\ref{lcth-lcta-lct-cdg-abs-contraderived-equiv-thm}.
 It follows from the commutativity of the diagram that the leftmost
vertical functor is a triangulated equivalence, too.
 This proves the first assertion of the corollary.

 In the case of the absolute derived or Positselski contraderived
category symbols $\st=\abs$ or~$\ctr$, the triangulated equivalence
$\sD^\st(\cB^\cu\bLcth_\bW^{X\dlct,\bth})\simeq
\sD^\st(\cB^\cu\bLcth_\bW^\bth)$ is provable by applying the dual
version of Proposition~\ref{second-kind-finite-resolutions}(a\+-b).
 According to the upper lines of the two
diagrams~(\ref{Upsilon-lcth-X-lct-thick-diagram}\+-%
\ref{Upsilon-lcth-X-thick-diagram}) in
Corollary~\ref{exact-dg-categories-of-lcth-thick-cdg-modules},
we have $\sZ^0((\cB^\cu\bLcth_\bW^{X\dlct,\bth})^\bec)\simeq
\cB^*\Lcth_\bW^{X\dlct,\thk}$ and
$\sZ^0((\cB^\cu\bLcth_\bW^\bth)^\bec)\simeq\cB^*\Lcth_\bW^\thk$.
 Hence $\cB^\cu\bLcth_\bW^{X\dlct,\bth}$ is a strict exact
DG\+subcategory in $\cB^\cu\bLcth_\bW^\bth$.

 To show that the full subcategory $\cB^*\Lcth_\bW^{X\dlct,\thk}$
is coresolving in $\cB^*\Lcth_\bW^\thk$, one needs to use the graded
version of
Theorem~\ref{antilocally-flaproj-complete-cotorsion-pair-thm}(b)
together with the fact that all antilocally flaprojective contraherent
graded $\cB^*$\+modules are thick
(which is provable using the graded versions of
Theorem~\ref{antilocally-flaproj-complete-cotorsion-pair-thm}(c)
and Lemma~\ref{antilocally-flaprojective-on-affines}).
 Notice that all flat graded $\cB^*(U)$\+modules are thick for any
affine open subscheme $U\subset X$; hence all
$\cB^*(U)/\cO_X(U)$\+flaprojective graded $\cB^*(U)$\+modules are
thick, too.
 The result of Lemma~\ref{lcth-thick-direct-image} also needs to be
used.
 The coresolution dimensions are finite and bounded by a constant by
Corollary~\ref{A-lcth-al-X-lct-coresolution-dimension}(a)
and~\cite[Corollary~A.5.5]{Pcosh}.

 Now the horizontal lines in
the diagram~\eqref{thick-lcth-al-lcta-lct-cdg-abs-contrader-diag} are
triangulated equivalences by
Corollary~\ref{thick-abs-ctr-derived-indep-of-covering-or-al}(a\+-b),
while the rightmost and middle vertical functors are triangulated
equivalences by the previous paragraph.
 Once again, similarly to the respective argument in the proof of
Theorem~\ref{lcth-lcta-lct-cdg-abs-contraderived-equiv-thm}, one shows
that the leftmost vertical functor
in~\eqref{thick-lcth-al-lcta-lct-cdg-abs-contrader-diag} is
well-defined, and it follows that this functor is a triangulated
equivalence.
\end{proof}

\subsection{Reduced co/contraderived categories of $X$-cotorsion
and $X$-locally cotorsion CDG-modules}
\label{reduced-co-contrader-X-cot-X-lct-cdg-modules}
 This section complements
Sections~\ref{reduced-coderived-of-qcoh-X-cta-subsecn}\+-%
\ref{reduced-contraderived-of-antilocal-subsecn},
\ref{co-contra-cot-lct-subsecn},
and~\ref{reduced-koszul-duality-contra-side-subsecn}.
 We start with an $X$\+cotorsion version of
Corollary~\ref{thick-cdg-qcoh-cta-reduced-abs-derived-equiv-cor}.
 This also provides a more direct proof of the result of
Corollary~\ref{thick-reduced-co-contra-cta-cot-equivs-cor}
applicable under more restrictive assumptions on the scheme~$X$.

\begin{cor} \label{thick-cdg-qcoh-cot-reduced-abs-derived-equiv-cor}
 Let $X$ be a semi-separated Noetherian scheme of finite Krull
dimension and $(\g,\widetilde\g)$ be a quasi-coherent twisted
Lie algebroid over~$X$.
 Assume that\/ $\g$~is a finite locally free sheaf on $X$, and let
$\cB^\cu=\cC^\cu_X(\g,\widetilde\g)$ be the related Chevalley--Eilenberg
quasi-coherent CDG\+quasi-algebra over~$X$.
 Then there is a commutative diagram of triangulated equivalences and
triangulated Verdier quotient functors
\begin{equation} \label{thick-cdg-qcoh-cot-cta-reduced-abs-derived-diag}
\begin{gathered}
 \xymatrix{
  \sD^\abs(\cB^\cu\bQcoh^{X\dcot}_\bth) \ar@{=}[r] \ar@{->>}[d]
  & \sD^\abs(\cB^\cu\bQcoh^{X\dcta}_\bth) \ar@{=}[r] \ar@{->>}[d]
  & \sD^\abs(\cB^\cu\bQcoh_\bth) \ar@{->>}[d] \\
  \sD^\abs_{X\red}(\cB^\cu\bQcoh^{X\dcot}_\bth) \ar@{=}[r]
  & \sD^\abs_{X\red}(\cB^\cu\bQcoh^{X\dcta}_\bth) \ar@{=}[r]
  & \sD^{\co=\abs}_{X\red}(\cB^\cu\bQcoh_\bth)
 }
\end{gathered}
\end{equation}
where the vertical arrows with double heads show the natural
triangulated Verdier quotient functors, while the horizontal
triangulated equivalences are induced by the inclusions of exact
DG\+categories $\cB^\cu\bQcoh^{X\dcot}_\bth\rarrow
\cB^\cu\bQcoh^{X\dcta}_\bth\rarrow\cB^\cu\bQcoh_\bth$.
\end{cor}

\begin{proof}
 We follow the proof of
Corollary~\ref{thick-cdg-qcoh-cta-reduced-abs-derived-equiv-cor},
which provides the right-hand side of the diagram.
 The upper horizontal equivalences are provided by
Corollary~\ref{thick-cdg-qcoh-cot-abs-derived-equiv-cor}.

 To construct the leftmost lower horizontal equivalence and the whole
diagram, it remains to check that the upper horizontal equivalences
identify the thick subcategory
$\Ac_{X\red}\sD^\abs(\cB^\cu\bQcoh^{X\dcot}_\bth)
\subset\sD^\abs(\cB^\cu\bQcoh^{X\dcot}_\bth)$ with the thick
subcategory $\Ac_{X\red}\sD^\abs(\cB^\cu\bQcoh^{X\dcta}_\bth)
\subset\sD^\abs(\cB^\cu\bQcoh^{X\dcta}_\bth)$ and/or with the thick
subcategory $\Ac_{X\red}\sD^\abs(\cB^\cu\bQcoh_\bth)\subset
\sD^\abs(\cB^\cu\bQcoh_\bth)$.
 This follows from commutativity of the diagram formed by
the functors $\gr_F^0$ from the proofs of
Lemmas~\ref{reduced-acyclic-thick-qcoh-localizing-subcategory}
and~\ref{reduced-acyclic-thick-X-cta-cot-qcoh-thick-subcat},
$$
 \xymatrix{
  \sD^\abs(\cB^\cu\bQcoh^{X\dcot}_\bth) \ar@{=}[r] \ar[d]_{\gr_F^0}
  & \sD^\abs(\cB^\cu\bQcoh^{X\dcta}_\bth) \ar@{=}[r] \ar[d]_{\gr_F^0}
  & \sD^\abs(\cB^\cu\bQcoh_\bth) \ar[d]^{\gr_F^0} \\
  \sD^\abs(X\Qcoh^\cot) \ar@{=}[r]
  & \sD^\abs(X\Qcoh^\cta) \ar@{=}[r]
  & \sD^\abs(X\Qcoh)
 }
$$
together with the fact that the functors $\sD^\abs(X\Qcoh^\cot)\rarrow
\sD^\abs(X\Qcoh^\cta)\rarrow\sD^\abs(X\Qcoh)$ are triangulated
equivalences (as shown in the lower line of the diagram).
 The latter fact follows from~\cite[Lemma~6.4.1(c)]{Pcosh}
or Lemma~\ref{A-qcoh-X-cot-coresolution-dimension} above
by virtue of the dual version of~\cite[Proposition~A.5.8]{Pcosh} or
Proposition~\ref{second-kind-finite-resolutions}(a) above.

 Furthermore, one needs to use the fact that the triangulated
equivalences $\sD^\abs(X\Qcoh^\cot)\rarrow\sD^\abs(X\Qcoh^\cta)\rarrow
\sD^\abs(X\Qcoh)$ identify the full subcategory of acyclic complexes
in $\sD^\abs(X\Qcoh^\cot)$ with the full subcategory of acyclic
complexes in $\sD^\abs(X\Qcoh^\cta)$ and/or with the full subcategory
of acyclic complexes in $\sD^\abs(X\Qcoh)$.
 This follows from commutativity of the diagram of triangulated
equivalences and triangulated Verdier quotient functors
$$
 \xymatrix{
  \sD^\abs(X\Qcoh^\cot) \ar@{=}[r] \ar@{->>}[d]
  & \sD^\abs(X\Qcoh^\cta) \ar@{=}[r] \ar@{->>}[d]
  & \sD^\abs(X\Qcoh) \ar@{->>}[d] \\
  \sD(X\Qcoh^\cot) \ar@{=}[r] & \sD(X\Qcoh^\cta) \ar@{=}[r]
  & \sD(X\Qcoh)
 }
$$
 The fact that the functors $\sD(X\Qcoh^\cot)\rarrow\sD(X\Qcoh^\cta)
\rarrow\sD(X\Qcoh)$ are triangulated equivalences, which is
the assertion of~\cite[Corollary~10.7]{PS6}, a part
of~\cite[Corollary~4.8.7]{Pcosh}, and a particular case of
Theorem~\ref{qcoh-X-cot-A-cot-derived-equivalence},
also needs to be used here.
\end{proof}

 The following corollary provides a more direct proof of
Corollary~\ref{lcth-reduced-X-lcta-X-lct-equivs-cor} and
the assertion about reduced contraderived categories of CDG\+modules
in Corollary~\ref{reduced-koszul-duality-contra-side-all-equiv}
working under more restrictive assumptions.

\begin{cor} \label{lcth-reduced-X-lcta-X-lct-equivalences-cor}
 Let $X$ be a semi-separated Noetherian scheme of finite Krull
dimension with an open covering\/ $\bW$ and $(\g,\widetilde\g)$ be
a quasi-coherent twisted Lie algebroid over~$X$.
 Assume that\/ $\g$~is a finite locally free sheaf on $X$, and let
$\cB^\cu=\cC^\cu_X(\g,\widetilde\g)$ be the related Chevalley--Eilenberg
quasi-coherent CDG\+quasi-algebra over~$X$.
 In this setting: \par
\textup{(a)} There is a commutative diagram of triangulated equivalences
and triangulated Verdier quotient functors 
\begin{equation} \label{lcth-X-lct-X-lcta-reduced-contrader-diag}
\begin{gathered}
 \xymatrix{
  \sD^\ctr(\cB^\cu\bLcth_\bW^{X\dlct}) \ar@{=}[rr] \ar@{->>}[d]
  && \sD^\ctr(\cB^\cu\bLcth_\bW) \ar@{->>}[d] \\
  \sD^\ctr_{X\red}(\cB^\cu\bLcth_\bW^{X\dlct}) \ar@{=}[rr]
  && \sD^\ctr_{X\red}(\cB^\cu\bLcth_\bW)
 }
\end{gathered}
\end{equation}
where the vertical arrows with double heads show the natural
triangulated Verdier quotient functors, while the horizontal
triangulated equivalences are induced by the inclusion of exact
DG\+categories $\cB^\cu\bLcth^{X\dlct}_\bW\rarrow\cB^\cu\bLcth_\bW$.
\par
\textup{(b)} There is a commutative diagram of triangulated equivalences
and triangulated Verdier quotient functors 
\begin{equation} \label{thick-lcth-X-lct-X-lcta-reduced-contrader-diag}
\begin{gathered}
 \xymatrix{
  \sD^\ctr(\cB^\cu\bLcth_\bW^{X\dlct,\bth}) \ar@{=}[rr] \ar@{->>}[d]
  && \sD^\ctr(\cB^\cu\bLcth_\bW^\bth) \ar@{->>}[d] \\
  \sD^{\ctr=\abs}_{X\red}(\cB^\cu\bLcth_\bW^{X\dlct,\bth}) \ar@{=}[rr]
  && \sD^{\ctr=\abs}_{X\red}(\cB^\cu\bLcth_\bW^\bth)
 }
\end{gathered}
\end{equation}
where the vertical arrows with double heads show the natural
triangulated Verdier quotient functors, while the horizontal
triangulated equivalences are induced by the inclusion of exact
DG\+categories $\cB^\cu\bLcth^{X\dlct,\bth}_\bW\rarrow
\cB^\cu\bLcth_\bW^\bth$.
\end{cor}

\begin{proof}
 To prove part~(a), we argue as in the proof of the rightmost vertical
triangulated equivalences in
Corollary~\ref{reduced-cdg-contraderived-equiv-cor}.
 The assertion follows from commutativity of the diagram of
triangulated functors and triangulated equivalences
$$
 \xymatrix{
  \sD^\ctr(X\Lcth_\bW^\lct) \ar@<-2pt>[rr] \ar[d]
  && \sD^\ctr(X\Lcth_\bW) \ar@{-}@<-2pt>[ll] \ar[d] \\
  \sD^\ctr(\cB^\cu\bLcth_\bW^{X\dlct}) \ar@<-2pt>[rr]
  && \sD^\ctr(\cB^\cu\bLcth_\bW) \ar@{-}@<-2pt>[ll]
 }
$$
where the lower horizontal equivalence is the result of
Theorem~\ref{lcth-lcta-lct-cdg-abs-contraderived-equiv-thm},
while the upper horizontal equivalence is provided
by~\cite[Corollary~6.4.4(b)]{Pcosh} (and can be also viewed as
a particular case of
Theorem~\ref{lcth-lcta-lct-cdg-abs-contraderived-equiv-thm}).

 In addition, one needs to use the commutativity of the diagram of
triangulated Verdier quotient functors and triangulated equivalences
$$
 \xymatrix{
  \sD^\ctr(X\Lcth_\bW^\lct) \ar@<-2pt>[rr] \ar@{->>}[d]
  && \sD^\ctr(X\Lcth_\bW) \ar@{-}@<-2pt>[ll] \ar@{->>}[d] \\
  \sD(X\Lcth_\bW^\lct) \ar@<-2pt>[rr]
  && \sD(X\Lcth_\bW) \ar@{-}@<-2pt>[ll]
 }
$$
where the lower horizontal equivalence is provided
by~\cite[Theorem~4.7.9 or Corollary~6.4.4(b)]{Pcosh}
(and can be also viewed as a particular case of
Theorem~\ref{A-lcth-X-lct-A-lct-derived-equivalence} above).

 Furthermore, we have $\sD^\ctr_{X\red}(\cB^\cu\bLcth_\bW^\bth)
=\sD^\abs_{X\red}(\cB^\cu\bLcth_\bW^\bth)$ by
formula~\eqref{lcth-reduced-contraderived=absolute-derived}
$\sD^\ctr_{X\red}(\cB^\cu\bLcth_\bW^{X\dlct,\bth})=
\sD^\abs_{X\red}(\cB^\cu\bLcth_\bW^{X\dlct,\bth})$ by
formula~\eqref{lcth-lct-reduced-contraderived=absolute-derived}
from Section~\ref{reduced-contraderived-of-lcth-subsecn}.
 This explains the notation $\sD^{\ctr=\abs}_{X\red}$ in
parts~(b) and~(c).

 To prove part~(b), we argue as in the proof of the leftmost vertical
triangulated equivalences in
Corollary~\ref{reduced-cdg-contraderived-equiv-cor}.
 The assertion follows from commutativity of the second diagram in
the argument above, together with the commutativity of the diagram
formed by the functors $\gr_F^0$ from the proof of
Lemma~\ref{reduced-acyclic-thick-lcth-colocalizing-subcategory}(a\+-b),
$$
 \xymatrix{
  \sD^\ctr(\cB^\cu\bLcth_\bW^{X\dlct,\bth})
  \ar@<-2pt>[rr] \ar[d]_{\gr_F^0}
  && \sD^\ctr(\cB^\cu\bLcth_\bW^\bth)
  \ar@{-}@<-2pt>[ll] \ar[d]^{\gr_F^0} \\
  \sD^\ctr(X\Lcth_\bW^\lct) \ar@<-2pt>[rr]
  && \sD^\ctr(X\Lcth_\bW) \ar@{-}@<-2pt>[ll]
 }
$$
 Here the upper horizontal equivalence is provided by
Corollary~\ref{thick-lcth-al-abs-ctr-derived-lcta-lct-equiv-cor},
while the lower horizontal equivalence was mentioned above.
\end{proof}

 The next corollary provides a more direct proof of
Corollary~\ref{reduced-thick-lcth-al-ctrder-lcta-lct-equiv-cor}
applicable under more restrictive assumptions.

\begin{cor} \label{antilocal-reduced-X-lcta-X-lct-equivalences-cor}
 Let $X$ be a semi-separated Noetherian scheme of finite Krull
dimension with an open covering\/ $\bW$ and $(\g,\widetilde\g)$ be
a quasi-coherent twisted Lie algebroid over~$X$.
 Assume that\/ $\g$~is a finite locally free sheaf on $X$, and let
$\cB^\cu=\cC^\cu_X(\g,\widetilde\g)$ be the related Chevalley--Eilenberg
quasi-coherent CDG\+quasi-algebra over~$X$.
 In this setting: \par
\textup{(a)} There is a commutative diagram of triangulated equivalences
and triangulated Verdier quotient functors 
\begin{equation} \label{antiloc-X-lct-X-lcta-reduced-contrader-diag}
\begin{gathered}
 \xymatrix{
  \sD^\ctr(\cB^\cu\bCtrh_\al^{X\dlct}) \ar@{=}[rr] \ar@{->>}[d]
  && \sD^\ctr(\cB^\cu\bCtrh_\al) \ar@{->>}[d] \\
  \sD^\ctr_{X\red}(\cB^\cu\bCtrh_\al^{X\dlct}) \ar@{=}[rr]
  && \sD^\ctr_{X\red}(\cB^\cu\bCtrh_\al)
 }
\end{gathered}
\end{equation}
where the vertical arrows with double heads show the natural
triangulated Verdier quotient functors, while the horizontal
triangulated equivalences are induced by the inclusion of exact
DG\+categories $\cB^\cu\bCtrh_\al^{X\dlct}\rarrow\cB^\cu\bCtrh_\al$.
\par
\textup{(b)} There is a commutative diagram of triangulated equivalences
and triangulated Verdier quotient functors 
\begin{equation} \label{thick-antil-lct-lcta-reduced-contrader-diag}
\begin{gathered}
 \xymatrix{
  \sD^\ctr(\cB^\cu\bCtrh_\al^{X\dlct,\bth}) \ar@{=}[rr] \ar@{->>}[d]
  && \sD^\ctr(\cB^\cu\bCtrh_\al^\bth) \ar@{->>}[d] \\
  \sD^{\ctr=\abs}_{X\red}(\cB^\cu\bCtrh_\al^{X\dlct,\bth}) \ar@{=}[rr]
  && \sD^{\ctr=\abs}_{X\red}(\cB^\cu\bCtrh_\al^\bth)
 }
\end{gathered}
\end{equation}
where the vertical arrows with double heads show the natural
triangulated Verdier quotient functors, while the horizontal
triangulated equivalences are induced by the inclusion of exact
DG\+categories $\cB^\cu\bCtrh_\al^{X\dlct,\bth}\rarrow
\cB^\cu\bCtrh_\al^\bth$.
\end{cor}

\begin{proof}
 Similar to the proof of
Corollary~\ref{lcth-reduced-X-lcta-X-lct-equivalences-cor}.
 One has to use the following triangulated equivalences:
\begin{itemize}
\item $\sD^\ctr(\cB^\cu\bCtrh_\al^{X\dlct})\simeq
\sD^\ctr(\cB^\cu\bCtrh_\al)$ (the result of
Theorem~\ref{lcth-lcta-lct-cdg-abs-contraderived-equiv-thm});
\item $\sD^\ctr(X\Ctrh_\al^\lct)\simeq\sD^\ctr(X\Ctrh_\al)$
(a combination of~\cite[Corollaries~4.7.4(a), 4.7.5(a),
and~6.4.4(b)]{Pcosh}, and also a particular case of
Theorem~\ref{lcth-lcta-lct-cdg-abs-contraderived-equiv-thm});
\item $\sD(X\Ctrh_\al^\lct)\simeq\sD(X\Ctrh_\al)$ (a combination
of~\cite[Corollaries~4.7.4(a) and~4.7.5(a)]{Pcosh}
with~\cite[Theorem~4.7.9 or Corollary~6.4.4(b)]{Pcosh},
or a particular case of the combination
of Corollary~\ref{A-lcth-ctrh-al-derived-equivalences}(a\+-b)
with Theorem~\ref{A-lcth-X-lct-A-lct-derived-equivalence});
\item $\sD^\ctr(\cB^\cu\bCtrh_\al^{X\dlct,\bth})\simeq
\sD^\ctr(\cB^\cu\bCtrh_\al^\bth)$ (the result of
Corollary~\ref{thick-lcth-al-abs-ctr-derived-lcta-lct-equiv-cor}).
\end{itemize}
 One also has to use
the formulas~(\ref{al-reduced-contraderived=absolute-derived}\+-%
\ref{al-lct-reduced-contraderived=absolute-derived}), the functors
$\gr_F^0$ from the (omitted) proof of
Lemma~\ref{reduced-acycl-thick-antiloc-ctrh-colocal-subcat}, and
the result of Theorem~\ref{antiloc-reduced-contrader-categ-equivs-thm}.
\end{proof}

\subsection{Comparison of semiderived Koszul duality equivalences on
the contra side} \label{koszul-duality-contra-side-comparison-subsecn}
 This section complements
Sections~\ref{semiderived-koszul-duality-contra-side-subsecn}\+-%
\ref{becker-semiderived-koszul-duality-contra-side-subsecn}.

\begin{cor} \label{semicontraderived-A-lct=X-lct=X-lcta-corollary}
 Let $X$ be a semi-separated Noetherian scheme of finite Krull
dimension with an open covering\/ $\bW$ and $(\g,\widetilde\g)$
be a quasi-coherent twisted Lie algebroid over~$X$.
 Assume that\/ $\g$~is a finite locally free sheaf on~$X$.
 Let $\cA=\cA_X(\g,\widetilde\g)$ be the enveloping quasi-coherent
quasi-algebra and $\cB^\cu=\cC^\cu_X(\g,\widetilde\g)$ be
the Chevalley--Eilenberg quasi-coherent CDG\+quasi-algebra
of~$(\g,\widetilde\g)$.
 Then the inclusions of exact categories $\cA\Lcth_\bW^{\cA\dlct}
\rarrow\cA\Lcth_\bW^{X\dlct}\rarrow\cA\Lcth_\bW$ induce equivalences
of the semicontraderived categories
$$
 \sD^{\si=\bsi}(\cA\Lcth_\bW^{\cA\dlct})\simeq
 \sD^{\si=\bsi}(\cA\Lcth_\bW^{X\dlct})\simeq
 \sD^{\si=\bsi}(\cA\Lcth_\bW),
$$
while the inclusion of exact DG\+categories $\cB^\cu\bLcth_\bW^{X\dlct}
\rarrow\cB^\cu\bLcth_\bW$ induces an equivalence of the contraderived
categories
$$
 \sD^{\ctr=\bctr}(\cB^\cu\bLcth_\bW^{X\dlct})\simeq
 \sD^{\ctr=\bctr}(\cB^\cu\bLcth_\bW).
$$
 There is a commutative diagram of triangulated equivalences
\begin{equation} \label{semicontraderived-A-lct=X-lct=X-lcta-diagram}
\begin{gathered}
 \xymatrix{
  \sD^{\si=\bsi}(\cA\Lcth_\bW) \ar@{=}[r] \ar@{-}@<-2pt>[d]
  & \sD^{\ctr=\bctr}(\cB^\cu\bLcth_\bW) \ar@{-}@<-2pt>[d] \\
  \sD^{\si=\bsi}(\cA\Lcth_\bW^{X\dlct})
  \ar@{=}[r] \ar@{-}@<-2pt>[d] \ar@<-2pt>[u]
  & \sD^{\ctr=\bctr}(\cB^\cu\bLcth_\bW^{X\dlct}) \ar@<-2pt>[u] \\
  \sD^{\si=\bsi}(\cA\Lcth_\bW^{\cA\dlct}) \ar@{=}[r] \ar@<-2pt>[u]
  & \sD^{\ctr=\bctr}(\cB^\cu\bLcth_\bW^{X\dlct}) \ar@{=}[u]
 }
\end{gathered}
\end{equation}
where the horizontal double lines denote the triangulated equivalences
of Koszul duality provided by
Theorem~\ref{semiderived-koszul-duality-contra-side}(a\+-c) and/or
Theorem~\ref{semiderived-koszul-duality-becker-contra-side}(a\+-c).
 The vertical triangulated equivalences in the leftmost column are
induced by the inclusions of exact categories, the upper vertical
triangulated equivalence in the rightmost column is induced by
the inclusion of exact DG\+categories, and the lower vertical
equivalence in the rightmost column is the identity functor.
\end{cor}

\begin{proof}
 The Positselski and Becker versions of the semiderived categories
in the leftmost column coincide by
Corollary~\ref{noetherian-Positselski=Becker-semicontraderived}.
 The Positselski and Becker versions of the contraderived categories
in the rightmost column coincide by
Remark~\ref{Noetherian-contra-side-b-and-p-agree-remark}.
 The lower vertical functor in the leftmost column is a triangulated
equivalence by Corollary~\ref{semicontraderived-A-lct=X-lct-corollary}
and/or~\ref{becker-semicontraderived-A-lct=X-lct-corollary};
the lower commutative square
in~\eqref{semicontraderived-A-lct=X-lct=X-lcta-diagram}
can be found on diagrams~\eqref{semicontraderived-A-lct=X-lct-diagram}
and/or~\eqref{becker-semicontraderived-A-lct=X-lct-diagram}.
 The upper vertical functor in the leftmost column is a triangulated
equivalence by
Corollary~\ref{semicontraderived-X-lcta-X-lct-A-mod-equivalence}.
 The upper vertical functor in the rightmost column is a triangulated
equivalence by
Theorem~\ref{lcth-lcta-lct-cdg-abs-contraderived-equiv-thm}.
 The commutativity of the diagram is obvious from the constructions.
\end{proof}

\subsection{Regular Noetherian schemes of finite Krull dimension,
technical background}
 The aim of this and the next section is to point out a major
simplification of the main (and many other) results and constructions
of this paper which occurs when the scheme $X$ is Noetherian
\emph{regular} of finite Krull dimension.
 The bottom line is that, under such assumptions on the scheme $X$,
the adjective ``$X$\+reduced'' and the corresponding subindex
$X\red$ (in the terminology and notation for the exotic derived
categories of CDG\+modules over~$\cB^\cu$) can be dropped everywhere.

 On the Koszul dual side, this means that the semiderived categories
of $\cA$\+modules coincide with the derived categories.
 So there is no difference between the semiderived and the reduced
Koszul duality theorems of Section~\ref{D-Omega-duality-secn} over
regular Noetherian schemes $X$ of finite Krull dimension.

\begin{lem} \label{regular-scheme-finite-homological-dimension-lemma}
 Let $X$ be a semi-separated \emph{regular} Noetherian scheme of finite
Krull dimension $D$ with an open covering\/ $\bW$ and a finite affine
open covering $X=\bigcup_{\alpha=1}^N U_\alpha$ subordinate to\/~$\bW$.
 Then \par
\textup{(a)} the homological dimension of the abelian category
of quasi-coherent sheaves $X\Qcoh$ does not exceed $D+N-1$; \par
\textup{(b)} the homological dimension of the exact category of
contraadjusted quasi-coherent sheaves $X\Qcoh^\cta$ does not exceed~$D$;
\par
\textup{(c)} the homological dimension of the exact category of
cotorsion quasi-coherent sheaves $X\Qcoh^\cot$ does not exceed~$D$; \par
\textup{(d)} the homological dimension of the exact category of\/
$\bW$\+locally contraherent cosheaves $X\Lcth_\bW$ does not exceed
$D+N-1$; \par
\textup{(e)} the homological dimension of the exact category of
locally cotorsion\/ $\bW$\+locally contraherent cosheaves
$X\Lcth_\bW^\lct$ does not exceed $D+N-1$; \par
\textup{(f)} the homological dimension of the exact category of
antilocal contraherent cosheaves $X\Ctrh_\al$ does not exceed~$D$; \par
\textup{(g)} the homological dimension of the exact category of
antilocal locally cotorsion contraherent cosheaves $X\Ctrh^\lct_\al$
does not exceed~$D$.
\end{lem}

\begin{proof}
 Notice first of all that, for any quasi-compact semi-separated
scheme $X$, the full subcategories $X\Qcoh^\cot\subset X\Qcoh^\cta$
are coresolving in $X\Qcoh$ \,\cite[Corollaries~4.1.2(c)
and~4.1.9(c)]{Pcosh}.
 Similarly, the full subcategory $X\Lcth_\bW^\lct$ is coresolving
in $X\Lcth_\bW$ (\cite[Lemma~4.3.2(a) or~4.4.1]{Pcosh} together with
Lemma~{flat-cotorsion-pair-hereditary}(b) above), while the full
subcategory $X\Ctrh_\al^\lct$ is coresolving in $X\Ctrh_\al$
(as one can see from~\cite[Corollary~4.3.5(a) or~4.4.4(a)]{Pcosh}).
 On the other hand, the full subcategory $X\Ctrh_\al$ is resolving
in $X\Lcth_\bW$ (by~\cite[Corollaries~4.3.3(b), 4.3.5(b\+-c),
and~4.3.7]{Pcosh}), while the full subcategory $X\Ctrh_\al^\lct$ is
resolving in $X\Lcth_\bW^\lct$ (by~\cite[Corollaries~4.3.3(b),
4.3.6(b\+-c), and~4.3.7]{Pcosh}).
 So all the categories in parts~(a\+-c) are coresolving in each other,
while all the categories in parts~(d\+-g) are either resolving or
coresolving in each other.

 By~\cite[Proposition~13.2.2(i)]{KS}, \cite[Theorem~12.1(b)]{Kel0},
or~\cite[Proposition~A.2.1 or~A.3.1(a)]{Pcosh}, the $\Ext$ groups
computed in an exact category and its (co)resolving subcategory agree.
 So the $\Ext$ groups in all the categories~(a\+-c) agree with each
other, and dual-analogously, the $\Ext$ groups in all
the categories~(d\+-g) agree with each other.
 In the case of an affine scheme $U=\Spec R$, the full subcategories
$U\Ctrh^\lct=R\Modl^\cot$ and $U\Ctrh=R\Modl^\cta$ are also coresolving
in $R\Modl$, so the $\Ext$ groups computed in $U\Ctrh^\lct\subset
U\Ctrh$ and in $R\Modl$ also agree.
 For a regular Noetherian commutative ring of Krull dimension~$D$, it
is well-known that the homological dimension of the abelian category
$R\Modl$ does not exceed~$D$.

 Furthermore, for any affine open subscheme $U\subset X$ subordinate
to~$\bW$, we have an adjoint pair of exact functors of direct and
inverse image of quasi-coherent sheaves $U\Qcoh\leftrightarrows X\Qcoh$
and an adjoint pair of exact functors of direct and inverse image of
($\bW$\+locally) contraherent cosheaves $U\Ctrh\rightleftarrows
X\Lcth_\bW$.
 An exact functor with an exact right adjoint does not increase
the projective dimension of objects, while an exact functor with
an exact left adjoint does not increase the injective
dimension~\cite[Lemma~1.7(a)]{Pal}.
 In view of the antilocality of the classes of contraadusted
quasi-coherent sheaves~\cite[Corollary~4.1.4(c)]{Pcosh},
cotorsion quasi-coherent sheaves~\cite[Corollary~4.1.11(c)]{Pcosh},
antilocal contraherent cosheaves~\cite[Corollary~4.3.5(c)]{Pcosh},
and antilocal locally cotorsion contraherent
cosheaves~\cite[Corollary~4.3.6(c)]{Pcosh}, the assertions of
parts~(b\+-c) and~(f\+-g) follow.

 Parts~(d\+-e) follow from parts~(f\+-g), because the respective
resolution dimensions do not exceed $N-1$
(\cite[Lemmas~4.7.1(c) and~4.7.2(a)]{Pcosh} or
Lemma~\ref{antilocal-resolving-resolution-dimension} above).
 Part~(a) can be obtained as particular case
of~\cite[Theorem~6.3(b)]{PS6}.
\end{proof}

\begin{cor} \label{all-acyclic-classes-coincide-over-regular-X}
 Let $X$ be a semi-separated \emph{regular} Noetherian scheme of finite
Krull dimension.
 Then one has: \par
\textup{(a)} $\Ac^\co(X\Qcoh)=\Ac^\bco(X\Qcoh)=\Ac^\abs(X\Qcoh)=
\Ac(X\Qcoh)$; \par
\textup{(b)} $\Ac^\bco(X\Qcoh^\cta)=\Ac^\abs(X\Qcoh^\cta)=
\Ac(X\Qcoh^\cta)$; \par
\textup{(c)} $\Ac^\bco(X\Qcoh^\cot)=\Ac^\abs(X\Qcoh^\cot)=
\Ac(X\Qcoh^\cot)$; \par
\textup{(d)} $\Ac^\ctr(X\Lcth_\bW)=\Ac^\bctr(X\Lcth_\bW)=
\Ac^\abs(X\Lcth_\bW)=\Ac(X\Lcth_\bW)$; \par
\textup{(e)} $\Ac^\ctr(X\Lcth_\bW^\lct)=\Ac^\bctr(X\Lcth_\bW^\lct)=
\Ac^\abs(X\Lcth_\bW^\lct)=\Ac(X\Lcth_\bW^\lct)$; \par
\textup{(f)} $\Ac^\ctr(X\Ctrh_\al)=\Ac^\bctr(X\Ctrh_\al)=
\Ac^\abs(X\Ctrh_\al)=\Ac(X\Ctrh_\al)$; \par
\textup{(g)} $\Ac^\ctr(X\Ctrh_\al^\lct)=\Ac^\bctr(X\Ctrh_\al^\lct)=
\Ac^\abs(X\Ctrh_\al^\lct)=\Ac(X\Ctrh_\al^\lct)$.
\end{cor}

\begin{proof}
 In an exact category of finite homological dimension, the classes of
acyclic and absolutely acyclic complexes
coincide~\cite[Remark~2.1]{Psemi}, \cite[Lemma~A.1.1]{Pcosh}.
 It follows that, in an exact category of finite homological dimension
with exact infinite direct sums, the classes of acyclic and
Positselski-coacyclic complexes coincide, while in an exact category
of finite homological dimension with exact products, the classes of
acyclic and Positselski-contraacyclic complexes coincide.

 Furthermore, in any exact category, all absolutely acyclic complexes
are both Becker-coacyclic and Becker-contraacyclic, and also acyclic.
 In any exact category with exact infinite direct sums, all
Positselski-coacyclic complexes are acyclic and Becker-coacyclic;
and dually, in any exact category with exact infinite products,
all Positselski-contraacyclic complexes are acyclic and
Becker-contraacyclic.

 The question of acyclicity of Becker-coacyclic and
Becker-contraacyclic complexes in exact categories may be more
delicate~\cite[Remark~B.7.4]{Pcosh}.
 We know that, for any quasi-compact semi-separated scheme $X$,
all Becker-coacyclic complexes in $X\Qcoh$ are acyclic in $X\Qcoh$,
all Becker-coacyclic complexes in $X\Qcoh^\cta$ are acyclic in
$X\Qcoh^\cta$, and all Becker-coacyclic complexes in $X\Qcoh^\cot$
are acyclic in $X\Qcoh^\cot$,
by Lemma~\ref{becker-coacyclic-complexes-are-acyclic}.
 Under the same assumptions, all Becker-contraacyclic complexes
in $X\Lcth_\bW$ are acyclic in $X\Lcth_\bW$, and all
Becker-contraacyclic complexes in $X\Lcth_\bW^\lct$ are acyclic
in $X\Lcth_\bW^\lct$ by
Corollary~\ref{becker-contraacyclic-complexes-are-acyclic}.

 In any exact category of finite homological dimension with enough
injective objects, the classes of acyclic, absolutely acyclic,
and Becker-coacyclic complexes coincide; and dually, in any exact
category of finite homological dimension with enough projective
objects, the classes of acyclic, absolutely acyclic, and
Becker-contraacyclic complexes coincide~\cite[Proposition~7.5]{Pphil},
\cite[Theorem~B.7.6]{Pcosh}.

 In view of the mentioned observations, the assertions parts~(a\+-g)
of the corollary follow from the respective parts of
Lemma~\ref{regular-scheme-finite-homological-dimension-lemma}.
\end{proof}

 The definitions of the semiderived categories appearing in
the next corollary were given in
Sections~\ref{semiderived-quasi-coherent-subsecn}\+-%
\ref{semiderived-contraherent-subsecn}
and~\ref{becker-semicontraderived-defined-subsecn}.

\begin{cor} \label{regular-reduced-semiderived=derived-cor}
 Let $X$ be a semi-separated \emph{regular} Noetherian scheme of finite
Krull dimension with an open covering\/ $\bW$ and $\cA$ be
a quasi-coherent quasi-algebra over~$X$.
 Then one has \par \hfuzz=15pt
\textup{(a)} $\Ac^\si(\cA\Qcoh)=\Ac(\cA\Qcoh)$ and\/
$\sD^\si(\cA\Qcoh)=\sD(\cA\Qcoh)$; \par
\textup{(b)} $\Ac^\bsi(\cA\Qcoh)=\Ac(\cA\Qcoh)$ and\/
$\sD^\bsi(\cA\Qcoh)=\sD(\cA\Qcoh)$; \par
\textup{(c)} $\Ac^\bsi(\cA\Qcoh^{X\dcta})=\Ac(\cA\Qcoh^{X\dcta})$ and\/
$\sD^\bsi(\cA\Qcoh^{X\dcta})=\sD(\cA\Qcoh^{X\dcta})$; \par
\textup{(d)} $\Ac^\bsi(\cA\Qcoh^{X\dcot})=\Ac(\cA\Qcoh^{X\dcot})$ and\/
$\sD^\bsi(\cA\Qcoh^{X\dcot})=\sD(\cA\Qcoh^{X\dcot})$; \par
\textup{(e)} $\Ac^\bsi(\cA\Qcoh^{\cA\dcot})=\Ac(\cA\Qcoh^{\cA\dcot})$
and\/ $\sD^\bsi(\cA\Qcoh^{\cA\dcot})=\sD(\cA\Qcoh^{\cA\dcot})$; \par
\textup{(f)} $\Ac^\si(\cA\Lcth_\bW)=\Ac(\cA\Lcth_\bW)$ and\/
$\sD^\si(\cA\Lcth_\bW)=\sD(\cA\Lcth_\bW)$; \par
\textup{(g)} $\Ac^\bsi(\cA\Lcth_\bW)=\Ac(\cA\Lcth_\bW)$ and\/
$\sD^\bsi(\cA\Lcth_\bW)=\sD(\cA\Lcth_\bW)$; \par
\textup{(h)} $\Ac^\si(\cA\Lcth_\bW^{X\dlct})=\Ac(\cA\Lcth_\bW^{X\dlct})$
and\/ $\sD^\si(\cA\Lcth_\bW^{X\dlct})=\sD(\cA\Lcth_\bW^{X\dlct})$; \par
\textup{(i)} $\Ac^\bsi(\cA\Lcth_\bW^{X\dlct})=
\Ac(\cA\Lcth_\bW^{X\dlct})$ and\/ $\sD^\bsi(\cA\Lcth_\bW^{X\dlct})=
\sD(\cA\Lcth_\bW^{X\dlct})$; \par
\textup{(j)} $\Ac^\si(\cA\Lcth_\bW^{\cA\dlct})=
\Ac(\cA\Lcth_\bW^{\cA\dlct})$ and\/ $\sD^\si(\cA\Lcth_\bW^{\cA\dlct})=
\sD(\cA\Lcth_\bW^{\cA\dlct})$; \par
\textup{(k)} $\Ac^\bsi(\cA\Lcth_\bW^{\cA\dlct})=
\Ac(\cA\Lcth_\bW^{\cA\dlct})$ and\/ $\sD^\bsi(\cA\Lcth_\bW^{\cA\dlct})=
\sD(\cA\Lcth_\bW^{\cA\dlct})$.
\end{cor}

\begin{proof}
 All the claims follow from the relevant assertions of
Corollary~\ref{all-acyclic-classes-coincide-over-regular-X} in view
of the results of
Lemmas~\ref{qcoh-forgetful-functor-reflects-acyclicity}
and~\ref{lcth-forgetful-functor-reflects-acyclicity}.
 Specifically:
\begin{itemize}
\item parts~(a\+-b) follow from
Corollary~\ref{all-acyclic-classes-coincide-over-regular-X}(a)
and Lemma~\ref{qcoh-forgetful-functor-reflects-acyclicity}(a);
\item part~(c) follows from
Corollary~\ref{all-acyclic-classes-coincide-over-regular-X}(b)
and Lemma~\ref{qcoh-forgetful-functor-reflects-acyclicity}(b);
\item part~(d) follows from
Corollary~\ref{all-acyclic-classes-coincide-over-regular-X}(c)
and Lemma~\ref{qcoh-forgetful-functor-reflects-acyclicity}(c);
\item part~(e) follows from
Corollary~\ref{all-acyclic-classes-coincide-over-regular-X}(c)
and Lemma~\ref{qcoh-forgetful-functor-reflects-acyclicity}(d);
\item parts~(f\+-g) follow from
Corollary~\ref{all-acyclic-classes-coincide-over-regular-X}(d)
and Lemma~\ref{lcth-forgetful-functor-reflects-acyclicity}(a).
\item parts~(h\+-i) follow from
Corollary~\ref{all-acyclic-classes-coincide-over-regular-X}(e)
and Lemma~\ref{lcth-forgetful-functor-reflects-acyclicity}(b).
\item parts~(j\+-k) follow from
Corollary~\ref{all-acyclic-classes-coincide-over-regular-X}(e)
and Lemma~\ref{lcth-forgetful-functor-reflects-acyclicity}(c).
\end{itemize}
\end{proof}

 The definitions of the classes of reduced-coacyclic CDG\+modules,
reduced-con\-tra\-acyclic CDG\+modules, reduced-acyclic thick
CDG\+modules, and the reduced coderived and reduced contraderived
categories appearing in the following corollary were given in
Sections~\ref{reduced-coderived-of-qcoh-subsecn}\+-%
\ref{reduced-contraderived-of-antilocal-subsecn}.

\begin{cor} \label{regular-reduced-dropped-corollary}
 Let $X$ be a semi-separated \emph{regular} Noetherian scheme of finite
Krull dimension with an open covering\/ $\bW$ and $(\g,\widetilde\g)$
be a quasi-coherent twisted Lie algebroid over~$X$.
 Assume that\/ $\g$~is a finite locally free sheaf on $X$, and let
$\cB^\cu=\cC^\cu_X(\g,\widetilde\g)$ be the related Chevalley--Eilenberg
quasi-coherent CDG\+quasi-algebra over~$X$.
 Then one has \par
\hfuzz=18pt \emergencystretch=3em \hbadness=5500
\textup{(a)} $\Ac^\co_{X\red}(\cB^\cu\bQcoh)=\Ac^\co(\cB^\cu\bQcoh)$
and\/ $\sD^\co_{X\red}(\cB^\cu\bQcoh)=\sD^\co(\cB^\cu\bQcoh)$; \par
\textup{(b)} $\Ac_{X\red}\sD^\co(\cB^\cu\bQcoh_\bth)=0
=\Ac_{X\red}\sD^\abs(\cB^\cu\bQcoh_\bth)$ and\/
$\sD^\co(\cB^\cu\bQcoh_\bth)=\sD^\co_{X\red}(\cB^\cu\bQcoh_\bth)=
\sD^\abs_{X\red}(\cB^\cu\bQcoh_\bth)=\sD^\abs(\cB^\cu\bQcoh_\bth)$; \par
\textup{(c)} $\Ac_{X\red}\sD^\abs(\cB^\cu\bQcoh^{X\dcta}_\bth)=0$ and\/
$\sD^\abs_{X\red}(\cB^\cu\bQcoh^{X\dcta}_\bth)=
\sD^\abs(\cB^\cu\bQcoh^{X\dcta}_\bth)$; \par
\textup{(d)} $\Ac_{X\red}\sD^\abs(\cB^\cu\bQcoh^{X\dcot}_\bth)=0$ and\/
$\sD^\abs_{X\red}(\cB^\cu\bQcoh^{X\dcot}_\bth)=
\sD^\abs(\cB^\cu\bQcoh^{X\dcot}_\bth)$; \par
\textup{(e)} $\Ac^\ctr_{X\red}(\cB^\cu\bLcth_\bW)=
\Ac^\ctr(\cB^\cu\bLcth_\bW)$ and\/ $\sD^\ctr_{X\red}(\cB^\cu\bLcth_\bW)
=\sD^\ctr(\cB^\cu\bLcth_\bW)$; \par
\textup{(f)} $\Ac^\ctr_{X\red}(\cB^\cu\bLcth_\bW^{X\dlct})=
\Ac^\ctr(\cB^\cu\bLcth_\bW^{X\dlct})$ and\/
$\sD^\ctr_{X\red}(\cB^\cu\bLcth_\bW^{X\dlct})=
\sD^\ctr(\cB^\cu\bLcth_\bW^{X\dlct})$; \par
\textup{(g)} $\Ac_{X\red}\sD^\ctr(\cB^\cu\bLcth_\bW^\bth)=0
=\Ac_{X\red}\sD^\abs(\cB^\cu\bLcth_\bW^\bth)$ and\/
$\sD^\ctr(\cB^\cu\bLcth_\bW^\bth)=
\sD^\ctr_{X\red}(\cB^\cu\bLcth_\bW^\bth)=
\sD^\abs_{X\red}(\cB^\cu\bLcth_\bW^\bth)=
\sD^\abs(\cB^\cu\bLcth_\bW^\bth)$; \par
\textup{(h)} $\Ac_{X\red}\sD^\ctr(\cB^\cu\bLcth_\bW^{X\dlct,\bth})=0
=\Ac_{X\red}\sD^\abs(\cB^\cu\bLcth_\bW^{X\dlct,\bth})$ and\/
$\sD^\ctr(\cB^\cu\bLcth_\bW^{X\dlct,\bth})=
\sD^\ctr_{X\red}(\cB^\cu\bLcth_\bW^{X\dlct,\bth})=
\sD^\abs_{X\red}(\cB^\cu\bLcth_\bW^{X\dlct,\bth})=
\sD^\abs(\cB^\cu\bLcth_\bW^{X\dlct,\bth})$; \par
\textup{(i)} $\Ac^\ctr_{X\red}(\cB^\cu\bCtrh_\al)=
\Ac^\ctr(\cB^\cu\bCtrh_\al)$ and\/ $\sD^\ctr_{X\red}(\cB^\cu\bCtrh_\al)
=\sD^\ctr(\cB^\cu\bCtrh_\al)$; \par
\textup{(j)} $\Ac^\ctr_{X\red}(\cB^\cu\bCtrh_\al^{X\dlct})=
\Ac^\ctr(\cB^\cu\bCtrh_\al^{X\dlct})$ and\/
$\sD^\ctr_{X\red}(\cB^\cu\bCtrh_\al^{X\dlct})=
\sD^\ctr(\cB^\cu\bCtrh_\al^{X\dlct})$; \par
\textup{(k)} $\Ac_{X\red}\sD^\ctr(\cB^\cu\bCtrh_\al^\bth)=0
=\Ac_{X\red}\sD^\abs(\cB^\cu\bCtrh_\al^\bth)$ and\/
$\sD^\ctr(\cB^\cu\bCtrh_\al^\bth)=
\sD^\ctr_{X\red}(\cB^\cu\bCtrh_\al^\bth)=
\sD^\abs_{X\red}(\cB^\cu\bCtrh_\al^\bth)=
\sD^\abs(\cB^\cu\bCtrh_\al^\bth)$; \par
\textup{(l)} $\Ac_{X\red}\sD^\ctr(\cB^\cu\bCtrh_\al^{X\dlct,\bth})=0
=\Ac_{X\red}\sD^\abs(\cB^\cu\bCtrh_\al^{X\dlct,\bth})$ and\/
$\sD^\ctr(\cB^\cu\bCtrh_\al^{X\dlct,\bth})=
\sD^\ctr_{X\red}(\cB^\cu\bCtrh_\al^{X\dlct,\bth})=
\sD^\abs_{X\red}(\cB^\cu\bCtrh_\al^{X\dlct,\bth})=
\sD^\abs(\cB^\cu\bCtrh_\al^{X\dlct,\bth})$.
\end{cor}

\begin{proof}
 All the claims follow immediately from the definitions in view
of the relevant assertions of
Corollary~\ref{all-acyclic-classes-coincide-over-regular-X}.
 Specifically:
\begin{itemize}
\item part~(a) follows from
Corollary~\ref{all-acyclic-classes-coincide-over-regular-X}(a);
\item part~(b) follows from
Corollary~\ref{all-acyclic-classes-coincide-over-regular-X}(a)
together with formula~\eqref{reduced-coderived=absolute-derived};
\item part~(c) follows from
Corollary~\ref{all-acyclic-classes-coincide-over-regular-X}(b);
\item part~(d) follows from
Corollary~\ref{all-acyclic-classes-coincide-over-regular-X}(c);
\item part~(e) follows from
Corollary~\ref{all-acyclic-classes-coincide-over-regular-X}(d);
\item part~(f) follows from
Corollary~\ref{all-acyclic-classes-coincide-over-regular-X}(e);
\item part~(g) follows from
Corollary~\ref{all-acyclic-classes-coincide-over-regular-X}(d)
together with
formula~\eqref{lcth-reduced-contraderived=absolute-derived};
\item part~(h) follows from
Corollary~\ref{all-acyclic-classes-coincide-over-regular-X}(e)
together with
formula~\eqref{lcth-lct-reduced-contraderived=absolute-derived};
\item part~(i) follows from
Corollary~\ref{all-acyclic-classes-coincide-over-regular-X}(f);
\item part~(j) follows from
Corollary~\ref{all-acyclic-classes-coincide-over-regular-X}(g);
\item part~(k) follows from
Corollary~\ref{all-acyclic-classes-coincide-over-regular-X}(f)
together with
formula~\eqref{al-reduced-contraderived=absolute-derived};
\item part~(l) follows from
Corollary~\ref{all-acyclic-classes-coincide-over-regular-X}(g)
together with
formula~\eqref{al-lct-reduced-contraderived=absolute-derived}.
\end{itemize}
\end{proof}

\subsection{Regular Noetherian schemes of finite Krull dimension,
main results} 
 In this section collects some of the main results of this paper in
the form of commutative diagrams of triangulated equivalences,
redrawn with the simplifications from
Corollaries~\ref{regular-reduced-semiderived=derived-cor}
and~\ref{regular-reduced-dropped-corollary} taken into account.
 We start with the co-contra correpondence for CDG\+modules over
the Chevalley--Eilenberg quasi-coherent CDG\+quasi-algebra before
passing to the Koszul duality and quadrality.

\begin{cor} \label{regular-cdg-co-contra-via-lcta-and-lct}
 Let $X$ be a semi-separated \emph{regular} Noetherian scheme of
finite Krull dimension with an open covering\/ $\bW$ and
$(\g,\widetilde\g)$ be a quasi-coherent twisted Lie algebroid over~$X$.
 Assume that\/ $\g$~is a finite locally free sheaf on $X$, and let
$\cB^\cu=\cC^\cu_X(\g,\widetilde\g)$ be the related Chevalley--Eilenberg
quasi-coherent CDG\+quasi-algebra over~$X$.
 Then there is a commutative diagram of triangulated equivalences
\begin{equation} \label{regular-cdg-co-contra-via-lcta-and-lct-diag}
\begin{gathered}
 \xymatrix{
  & \sD^{\ctr=\bctr}(\cB^\cu\bLcth_\bW)
  \ar@{=}[ld] \ar@{-}@<-2pt>[dd] \\
  \sD^{\co=\bco}(\cB^\cu\bQcoh) \ar@{=}[rd] \\
  & \sD^{\ctr=\bctr}(\cB^\cu\bLcth_\bW^{X\dlct}) \ar@<-2pt>[uu]
 }
\end{gathered}
\end{equation}
 The vertical triangulated equivalence is induced by the inclusion
of exact DG\+cat\-e\-gories $\cB^\cu\bLcth_\bW^{X\dlct}\rarrow
\cB^\cu\bLcth_\bW$, as per
Theorem~\textup{\ref{lcth-lcta-lct-cdg-abs-contraderived-equiv-thm}} or
Corollary~\textup{\ref{lcth-reduced-X-lcta-X-lct-equivalences-cor}(a)}.
 The diagonal triangulated equivalences are provided by
Corollaries~\textup{\ref{cdg-module-co-contra-correspondence-cor}}
and~\textup{\ref{X-cot-X-lct-cdg-module-co-contra-corresp-cor}}.
\end{cor}

\begin{proof}
 This is the commutative diagram of triangulated equivalences from
Corollary~\ref{cdg-co-contra-via-lcta-and-lct-agree-cor}.
 For a semi-separated regular Noetherian scheme $X$ of finite Krull
dimension, the subindices $X\red$ in the notation of
diagram~\eqref{cdg-co-contra-via-lcta-and-lct-diagram} can be dropped
by Corollary~\ref{regular-reduced-dropped-corollary}(a,e,f).
 The Positselski co/contraderived categories coincide with
the Becker ones for the exact DG\+categories at hand by
Corollary~\ref{noetherian-quasi-algebra-Positselski=Becker}
and Remark~\ref{becker=positselski-at-the-end-remark}, or by
Remarks~\ref{loc-Noetherian-right-co-side-b-and-p-agree-remark}
and~\ref{Noetherian-contra-side-b-and-p-agree-remark};
hence the notation $\co=\bco$ and $\ctr=\bctr$.
\end{proof}

 The following corollary states the Koszul ($\cD$\+$\Omega$) duality
on the co side for regular schemes.

\begin{cor} \label{regular-koszul-duality-right-co-side}
 Let $X$ be a semi-separated \emph{regular} Noetherian scheme of
finite Krull dimension and $(\g,\widetilde\g)$ be a quasi-coherent
twisted Lie algebroid over~$X$.
 Assume that\/ $\g$~is a finite locally free sheaf  on~$X$.
 Let $\cA=\cA_X(\g,\widetilde\g)$ be the enveloping quasi-coherent
quasi-algebra and $\cB^\cu=\cC^\cu_X(\g,\widetilde\g)$ be
the Chevalley--Eilenberg quasi-coherent CDG\+quasi-algebra
of~$(\g,\widetilde\g)$.
 Then the pair of adjoint DG\+functors from
Lemma~\ref{koszul-duality-dg-functors-right-co-side} induces
a triangulated equivalence
$$
 \xymatrix{
  {-}\ot_{\cO_X}\bigwedge\nolimits_X^*(\g)\:
  \sD(\Qcohr\cA) \ar@{=}[r]
  & \sD^{\co=\bco}(\bQcohr\cB^\cu) \,:{-}\ot_{\cO_X}\cA.
 }
$$
\end{cor}

\begin{proof}
 This is a common special case of
Theorem~\ref{semiderived-koszul-duality-right-co-side}(a) or~(b)
and Theorem~\ref{reduced-koszul-duality-right-co-side}.
 For a semi-separated regular Noetherian scheme $X$ of finite Krull
dimension, the superindices $\si$ or~$\bsi$ in the notation of
Theorem~\ref{semiderived-koszul-duality-right-co-side}(a\+-b)
can be dropped by
Corollary~\ref{regular-reduced-semiderived=derived-cor}(a\+-b),
while the subindex $X\red$ in the notation of
Theorem~\ref{reduced-koszul-duality-right-co-side} can be dropped by
Corollary~\ref{regular-reduced-dropped-corollary}(a).
 Under the assumptions of the present corollary, the Positselski
and Becker coderived categories for the exact DG\+category
$\bQcohr\cB^\cu$ coincide by
Corollary~\ref{noetherian-quasi-algebra-Positselski=Becker},
or Remark~\ref{loc-Noetherian-right-co-side-b-and-p-agree-remark},
or Remark~\ref{becker=positselski-at-the-end-remark};
hence the notation $\co=\bco$.
\end{proof}

 The next corollary is a summary of the Koszul ($\cD$\+$\Omega$)
duality on the contra side for regular schemes.

\begin{cor} \label{regular-koszul-duality-contra-side-all-equiv}
 Let $X$ be a semi-separated \emph{regular} Noetherian scheme of
finite Krull dimension with an open covering\/ $\bW$ and
$(\g,\widetilde\g)$ be a quasi-coherent twisted Lie algebroid over~$X$.
 Assume that\/ $\g$~is a finite locally free sheaf on~$X$.
 Let $\cA=\cA_X(\g,\widetilde\g)$ be the enveloping quasi-coherent
quasi-algebra and $\cB^\cu=\cC^\cu_X(\g,\widetilde\g)$ be
the Chevalley--Eilenberg quasi-coherent CDG\+quasi-algebra
of~$(\g,\widetilde\g)$.
 Then there is a natural commutative diagram of triangulated
equivalences
\begin{equation} \label{regular-koszul-duality-contra-side-big-diagram}
\
\end{equation}
$$
 \qquad\ \xymatrix{
  \text{\llap{$\Cohom_X\bigl(\bigwedge\nolimits_X^*(\g),{-}\bigr)\:$}}
  \sD(\cA\Lcth_\bW) \ar@{=}[r] \ar@{-}@<-2pt>[d]
  & \sD^{\ctr=\bctr}(\cB^\cu\bLcth_\bW)
  \text{\rlap{$\,\,:\!\Cohom_X(\cA,{-})$}} \ar@{-}@<-2pt>[d] \\
  \text{\llap{$\Cohom_X\bigl(\bigwedge\nolimits_X^*(\g),{-}\bigr)\:$}}
  \sD(\cA\Lcth_\bW^{X\dlct}) \ar@{=}[r] \ar@<-2pt>[u] \ar@{-}@<-2pt>[d]
  & \sD^{\ctr=\bctr}(\cB^\cu\bLcth_\bW^{X\dlct})
  \text{\rlap{$\,\,:\!\Cohom_X(\cA,{-})$}}
  \ar@<-2pt>[u] \ar@{=}[d] \\
  \text{\llap{$\Cohom_X\bigl(\bigwedge\nolimits_X^*(\g),{-}\bigr)\:$}}
  \sD(\cA\Lcth_\bW^{\cA\dlct}) \ar@{=}[r] \ar@<-2pt>[u]
  & \sD^{\ctr=\bctr}(\cB^\cu\bLcth_\bW^{X\dlct})
  \text{\rlap{$\,\,:\!\Cohom_X(\cA,{-})$}}
 }
$$
with the horizontal equivalences provided by
Theorem~\ref{semiderived-koszul-duality-contra-side},
Theorem~\ref{semiderived-koszul-duality-becker-contra-side},
or~\ref{reduced-koszul-duality-contra-side}, the vertical
equivalences in the leftmost column induced by the inclusions of
exact categories $\cA\Lcth_\bW^{\cA\dlct}\rarrow\cA\Lcth_\bW^{X\dlct}
\rarrow\cA\Lcth_\bW$, and the upper vertical equivalence in
the rightmost column induced by the inclusion of exact DG\+categories
$\cB^\cu\bLcth_\bW^{X\dlct}\rarrow\cB^\cu\bLcth_\bW$.
 The lower vertical equivalence in the rightmost column is
the identity triangulated functor.
\end{cor}

\begin{proof}
 This is a common special case of
Corollaries~\ref{reduced-koszul-duality-contra-side-all-equiv}
and~\ref{semicontraderived-A-lct=X-lct=X-lcta-corollary}.
 For a semi-separated regular Noetherian scheme $X$ of finite Krull
dimension, the superindices $\si$ or~$\bsi$ in the notation of
diagram~\eqref{semicontraderived-A-lct=X-lct=X-lcta-diagram}
can be dropped by
Corollary~\ref{regular-reduced-semiderived=derived-cor}(f\+-k),
while the subindices $X\red$ in the notation of
diagram~\eqref{reduced-koszul-duality-contra-side-big-diagram}
can be dropped by
Corollary~\ref{regular-reduced-dropped-corollary}(e\+-f).
 Under the assumptions of the present corollary, the Positselski
and Becker contraderived categories for the exact DG\+categories
of CDG\+modules in the rightmost column coincide by
Remark~\ref{Noetherian-contra-side-b-and-p-agree-remark}
or~\ref{becker=positselski-at-the-end-remark};
hence the notation $\ctr=\bctr$.
\end{proof}

 Finally, we draw three versions of the big quadrality diagram
for regular schemes.

\begin{cor} \label{regular-main-hexagonality-corollary}
 Let $X$ be a semi-separated \emph{regular} Noetherian scheme of
finite Krull dimension with an open covering\/ $\bW$ and
$(\g,\widetilde\g)$ be a quasi-coherent twisted Lie algebroid over~$X$.
 Assume that the quasi-coherent sheaf\/~$\g$ on $X$ is finite locally
free of constant rank~$m$.
 Let $\cA=\cA_X(\g,\widetilde\g)$ and
$\cA^\circ=\cA_X(\g,\widetilde\g^\circ)$ be the two related twisted
universal enveloping quasi-coherent quasi-algebras, and let
$\cB^\cu=\cC_X^\cu(\g,\widetilde\g)$ and
$\cB^\circ{}^\cu=\cC_X^\cu(\g,\widetilde\g^\circ)$ be the two related
Chevalley--Eilenberg quasi-coherent CDG\+quasi-algebras over~$X$.
 Then there are natural commutative diagrams of triangulated category
equivalences
\begin{equation} \label{regular-main-hexagonality-diagram}
\begin{gathered}
 \xymatrix{
  \sD(\Qcohr\cA^\circ)
  \ar@<3pt>[rrrr]^{{-}\ot_{\cO_X}\bigwedge\nolimits_X^*(\g)}
  \ar@<3pt>[dd]^{\cHom_{\cO_X}(\cB^m[-m],{-})}
  &&&& \sD^{\co=\bco}(\bQcohr\cB^\circ{}^\cu)
  \ar@<3pt>[llll]^{{-}\ot_{\cO_X}\cA^\circ}
  \ar@{=}[dd]^{\cB^{\circ,\rop}{}^\subcu\simeq\cB^\subcu}
  \\ \\
  \sD(\cA\Qcoh) \ar@<3pt>[rrrr]^{\cB^*\ot_{\cO_X}{-}}
  \ar@<3pt>[uu]^{\cB^m[-m]\ot_{\cO_X}{-}}
  \ar@<3pt>[dd]^{\boR\fHom_\cA(\cA,{-})}
  &&&& \sD^{\co=\bco}(\cB^\cu\bQcoh)
  \ar@{-}@<3pt>[llll]
  \ar@<3pt>[dd]^{\boR\fHom_{\cB^*}(\cB^\subcu,{-})} \\ \\
  \sD(\cA\Lcth_\bW)
  \ar@<3pt>[rrrr]^{\Cohom_X\left(\bigwedge\nolimits_X^*(\g),{-}\right)}
  \ar@<3pt>[uu]^{\cA\ocn_\cA^\boL{-}}
  &&&& \sD^{\ctr=\bctr}(\cB^\cu\bLcth_\bW)
  \ar@<3pt>[llll]^{\Cohom_X(\cA,{-})}
  \ar@<3pt>[uu]^{\cB^\subcu\ocn_{\cB^*}^\boL{-}}
 }
\end{gathered}
\end{equation}
\begin{equation} \label{regular-main-X-lct-hexagonality-diagram}
\begin{gathered}
 \xymatrix{
  \sD(\Qcohr\cA^\circ)
  \ar@<3pt>[rrrr]^{{-}\ot_{\cO_X}\bigwedge\nolimits_X^*(\g)}
  \ar@<3pt>[dd]^{\cHom_{\cO_X}(\cB^m[-m],{-})}
  &&&& \sD^{\co=\bco}(\bQcohr\cB^\circ{}^\cu)
  \ar@<3pt>[llll]^{{-}\ot_{\cO_X}\cA^\circ}
  \ar@{=}[dd]^{\cB^{\circ,\rop}{}^\subcu\simeq\cB^\subcu}
  \\ \\
  \sD(\cA\Qcoh) \ar@<3pt>[rrrr]^{\cB^*\ot_{\cO_X}{-}}
  \ar@<3pt>[uu]^{\cB^m[-m]\ot_{\cO_X}{-}}
  \ar@<3pt>[dd]^{\boR\fHom_\cA(\cA,{-})}
  &&&& \sD^{\co=\bco}(\cB^\cu\bQcoh)
  \ar@{-}@<3pt>[llll]
  \ar@<3pt>[dd]^{\boR\fHom_{\cB^*}(\cB^\subcu,{-})} \\ \\
  \sD(\cA\Lcth_\bW^{X\dlct})
  \ar@<3pt>[rrrr]^{\Cohom_X\left(\bigwedge\nolimits_X^*(\g),{-}\right)}
  \ar@<3pt>[uu]^{\cA\ocn_\cA^\boL{-}}
  &&&& \sD^{\ctr=\bctr}(\cB^\cu\bLcth_\bW^{X\dlct})
  \ar@<3pt>[llll]^{\Cohom_X(\cA,{-})}
  \ar@<3pt>[uu]^{\cB^\subcu\ocn_{\cB^*}^\boL{-}}
 }
\end{gathered}
\end{equation}
\begin{equation} \label{regular-main-A-lct-hexagonality-diagram}
\begin{gathered}
 \xymatrix{
  \sD(\Qcohr\cA^\circ)
  \ar@<3pt>[rrrr]^{{-}\ot_{\cO_X}\bigwedge\nolimits_X^*(\g)}
  \ar@<3pt>[dd]^{\cHom_{\cO_X}(\cB^m[-m],{-})}
  &&&& \sD^{\co=\bco}(\bQcohr\cB^\circ{}^\cu)
  \ar@<3pt>[llll]^{{-}\ot_{\cO_X}\cA^\circ}
  \ar@{=}[dd]^{\cB^{\circ,\rop}{}^\subcu\simeq\cB^\subcu}
  \\ \\
  \sD(\cA\Qcoh) \ar@<3pt>[rrrr]^{\cB^*\ot_{\cO_X}{-}}
  \ar@<3pt>[uu]^{\cB^m[-m]\ot_{\cO_X}{-}}
  \ar@<3pt>[dd]^{\boR\fHom_\cA(\cA,{-})}
  &&&& \sD^{\co=\bco}(\cB^\cu\bQcoh)
  \ar@{-}@<3pt>[llll]
  \ar@<3pt>[dd]^{\boR\fHom_{\cB^*}(\cB^\subcu,{-})} \\ \\
  \sD(\cA\Lcth_\bW^{\cA\dlct})
  \ar@<3pt>[rrrr]^{\Cohom_X\left(\bigwedge\nolimits_X^*(\g),{-}\right)}
  \ar@<3pt>[uu]^{\cA\ocn_\cA^\boL{-}}
  &&&& \sD^{\ctr=\bctr}(\cB^\cu\bLcth_\bW^{X\dlct})
  \ar@<3pt>[llll]^{\Cohom_X(\cA,{-})}
  \ar@<3pt>[uu]^{\cB^\subcu\ocn_{\cB^*}^\boL{-}}
 }
\end{gathered}
\end{equation}
\end{cor}

\begin{proof}
 These are the results of
Corollaries~\ref{main-hexagonality-corollary}
and~\ref{main-lct-hexagonality-corollary}, restated with
the simplifications for the case of a regular scheme.
 For a semi-separated regular Noetherian scheme $X$ of finite Krull
dimension, the subindices $X\red$ in the notation of
diagrams~\eqref{main-hexagonality-diagram}
and~(\ref{main-X-lct-hexagonality-diagram}\+-%
\ref{main-A-lct-hexagonality-diagram}) can be dropped by
Corollary~\ref{regular-reduced-dropped-corollary}(a,e,f).
 The Positselski co/contraderived categories coincide with
the Becker ones for the exact DG\+categories at hand, as per
the proof of Corollary~\ref{regular-cdg-co-contra-via-lcta-and-lct}
and the other proofs above in this section; hence the notation
$\co=\bco$ and $\ctr=\bctr$.
\end{proof}

\end{document}